\documentclass[12pt]{article}

\usepackage{latexsym}
\usepackage{amssymb}
\usepackage{amsmath}
\usepackage{euscript}
\usepackage{authblk}

\newtheorem{lemma}{\bf Lemma}
\newtheorem{theorem}[lemma]{\bf Theorem}

\title{Hasse diagrams of  posets with up to 7 elements, and the number of posets with $10$ elements, without the use of computer programs}
\author[1]{Luiz F. Monteiro}
\author[2]{ Sonia Savini} 
\author[2,3]{Ignacio Viglizzo} 
\affil[1]{Universidad Nacional del Sur}
\affil[2] {Departamento de Matem\'atica, Universidad Nacional del Sur}
\affil[3] {INMABB-CONICET-Universidad Nacional del Sur}
\affil[ ]{\tt lfmonteiro0510@gmail.com, ssavini@criba.edu.ar, viglizzo@gmail.com}

\begin{document}
\font\fivrm=cmr5 \relax
\input{prepictex}
\input{pictex}
\input{postpictex}

\date {}
\maketitle

\begin{abstract}
	Let $P(n)$ be the set of all posets with $n$ elements. Let  $P^{(j)}(n)$, $1\leq j\leq 2^n,$ be the number of all posets with $n$ elements possessing exactly $j$ antichains. We have determined the numbers $P^{(j)}(7),$ $1\leq j\leq 128$, and 	using a result of M.~Ern\'e \cite{EM4}, we compute $|P(10)|$ without the aid of any computer program. We include the Hasse diagrams of all the non-isomorphic posets of $P(7)$. We also present formulas for the number of connected posets of certain forms, and use them to compute $|P(n)|$ with $1\le n\le 8$ by a different method.
	\end{abstract}
	
	\section{Introduction}
	
 Let $P(n)$ be the set of all posets with $n$ elements and   $NIP(n)$ the set of non-isomorphic posets with $n$ elements. It is a long-standing combinatorial problem to determine the numbers $|P(n)|$ and $|NIP(n)|$.
 
No simple explicit  or recursive formula for $|P(n)|,(|NIP(n)|)$ has been found yet.
The numbers determined so far are $|P(n)|$,  $1\leq n\leq 18$ and $|NIP(n)|$,  $1\leq n\leq 16$, see \cite{BG}, \cite{CL3}, \cite{CSD}, \cite{KV}, \cite{SA}, \cite{SH1}, \cite{EHL}, \cite{EM1}, \cite{EM2, EM3}, \cite{DS}, \cite{BVDR}, \cite{BBR}, \cite{RVI1}, \cite{RVI2}, \cite{RRA}, \cite{CR}, \cite{EM5}, \cite{EM6}, \cite{HR}, \cite{BM}.

 Let  $P^{(j)}(n)$, $1\leq j\leq 2^n,$ be the number of all posets with $n$ elements possessing exactly $j$ antichains. Note that $P^{(j)}(n)=0$, $1\leq j\leq n$, then  
 \begin{equation}\label{PnSum}
 |P(n)|=\sum\limits_{j=n+1}^{2^n} P^{(j)}(n).
 \end{equation}

  For the determination of the number of antichians of a poset $X$ we use some results of L.~Monteiro \cite{ML} (section 3). S.~K.~Das \cite{DS} has found, with the help of a computer program (see \cite{EM4}, page 128) $|P(n)|$, for $n\leq 11$, in particular:
$$ |P(10)|=6{,}611{,}065{,}248{,}783.$$
Note that M.~Ern\'e  determines in \cite{EM4} $|P(n)|$ for $1\leq n\leq 9$ without the aid of any computer program. In order to use this method to calculate $|P(10)|$, we have determined the numbers $P^{(j)}(7),$ $1\leq j\leq 128$  using for this the  Hasse diagrams of all the non-isomorphic posets of $P(7)$.

An extended version of this work, in Spanish, can be found in \cite{MSV}.

\section{Hasse diagrams of posets with up to 7 elements}\label{diagrams}

 In the Hasse diagrams below, for each poset $X$ we indicate in boldface the number of antichains of $X$ and under each one of them we indicate the number of posets which are isomorphic to $X$.
 
 For $n=2,3,4$ the Hasse diagrams are well known.

  ${\bf NIP(2)}$

$$
\begin{minipage}{4cm}
\beginpicture
\setcoordinatesystem units <1.5mm,2mm>
\setplotarea x from 0 to 16, y from 0 to 10
\put{${\bf 3}$} [l] at 0 10
\put{1)} [l] at 0 7
\put {$\scriptstyle \bullet$} [c] at  8  8
\put {$ \scriptstyle \bullet$} [c] at  8  2
\setlinear \plot 8 2  8 8 /
\put{$2$} [c] at 8 0
\endpicture
\end{minipage}
\begin{minipage}{4cm}
\beginpicture
\setcoordinatesystem units <1.5mm,2mm>
\setplotarea x from 0 to 16, y from 0 to 10
\put{${\bf 4}$} [l] at 0 10
\put{2)} [l] at 0 7
\put {$ \scriptstyle \bullet$} [c] at  6  2
\put {$ \scriptstyle \bullet$} [c] at  10 2
\put{$1$} [c] at 8 0
\endpicture
\end{minipage}
$$
\medskip

  ${\bf NIP(3)}$

$$
\begin{minipage}{4cm}
\beginpicture
\setcoordinatesystem units <1.5mm,2mm>
\setplotarea x from 0 to 16, y from 0 to 11
\put{${\bf 4}$} [l] at 0 11
\put{1)} [l] at 0 8
\put {$ \scriptstyle \bullet$} [c] at  8  2
\put {$ \scriptstyle \bullet$} [c] at  8  5
\put {$ \scriptstyle \bullet$} [c] at  8  8
\setlinear \plot 8 2 8 8 /
\put{$6$} [c] at 8 0
\endpicture
\end{minipage}
\begin{minipage}{4cm}
\beginpicture
\setcoordinatesystem units <1.5mm,2mm>
\setplotarea x from 0 to 16, y from 0 to 11
\put{${\bf 5}$} [l] at 0 11
\put{2)} [l] at 0 8
\put {$ \scriptstyle \bullet$} [c] at  6 8
\put {$ \scriptstyle \bullet$} [c] at  10 8
\put {$ \scriptstyle \bullet$} [c] at  8 2
\setlinear \plot 6 8 8 2 10 8  /
\put{$3$} [c] at 8 0
\endpicture
\end{minipage}
\begin{minipage}{4cm}
\beginpicture
\setcoordinatesystem units <1.5mm,2mm>
\setplotarea x from 0 to 16, y from 0 to 11
\put{3)} [l] at 0 8
\put {$\scriptstyle \bullet$} [c] at  6 2
\put {$\scriptstyle \bullet$} [c] at  10 2
\put {$\scriptstyle \bullet$} [c] at  8 8
\setlinear \plot 6 2 8 8 10 2  /
\put{$3$} [c] at 8 0
\endpicture
\end{minipage}
\begin{minipage}{4cm}
\beginpicture
\setcoordinatesystem units <1.5mm,2mm>
\setplotarea x from 0 to 16, y from 0 to 11
\put{${\bf 6}$} [l] at 0 11
\put{4)} [l] at 0 8
\put {$ \scriptstyle \bullet$} [c] at 6 2
\put {$ \scriptstyle \bullet$} [c] at 6 8
\put {$ \scriptstyle \bullet$} [c] at 10 2
\setlinear \plot 6 2 6 8 /
\put{$6$} [c] at 8 0
\endpicture
\end{minipage}
\begin{minipage}{4cm}
\beginpicture
\setcoordinatesystem units <1.5mm,2mm>
\setplotarea x from 0 to 16, y from 0 to 11
\put{${\bf 8}$} [l] at 0 11
\put{5)} [l] at 0 8
\put {$ \scriptstyle \bullet$} [c] at 6 2
\put {$ \scriptstyle \bullet$} [c] at 8 2
\put {$ \scriptstyle \bullet$} [c] at 10 2
\put{$1$} [c] at 8 0
\endpicture
\end{minipage}
$$

  ${\bf NIP(4)}$

$$
\begin{minipage}{4cm}
\beginpicture
\setcoordinatesystem units <1.5mm,2mm>
\setplotarea x from 0 to 16, y from  0 to 11
\put{${\bf 5}$} [l] at 0 11
\put{1)} [l] at 0 8
\put {$ \scriptstyle \bullet$} [c] at 8 8
\put {$ \scriptstyle \bullet$} [c] at 8 6
\put {$ \scriptstyle \bullet$} [c] at 8 4
\put {$ \scriptstyle \bullet$} [c] at 8 2
\setlinear \plot 8 2 8 8 /
\put{$24$} [c] at 8 0
\endpicture
\end{minipage}
\begin{minipage}{4cm}
\beginpicture
\setcoordinatesystem units <1.5mm,2mm>
\setplotarea x from 0 to 16, y from 0 to  11
\put{${\bf 6}$} [l] at 0 11
\put{2)} [l] at 0 8
\put {$ \scriptstyle \bullet$} [c] at 8 2
\put {$ \scriptstyle \bullet$} [c] at 8 5
\put {$ \scriptstyle \bullet$} [c] at 6 8
\put {$ \scriptstyle \bullet$} [c] at 10 8
\setlinear \plot 8 2 8 5 6 8  /
\setlinear \plot 8 5 10 8  /
\put{$12$} [c] at 8 0
\endpicture
\end{minipage}
\begin{minipage}{4cm}
\beginpicture
\setcoordinatesystem units <1.5mm,2mm>
\setplotarea x from 0 to 16, y from 0 to  11
\put{3)} [l] at 0 8
\put {$ \scriptstyle \bullet$} [c] at 8 8
\put {$ \scriptstyle \bullet$} [c] at 8 5
\put {$ \scriptstyle \bullet$} [c] at 6 2
\put {$ \scriptstyle \bullet$} [c] at 10 2
\setlinear \plot 6 2 8 5 8 8  /
\setlinear \plot 8 5 10 2  /
\put{$12$} [c] at 8 0
\endpicture
\end{minipage}
\begin{minipage}{4cm}
\beginpicture
\setcoordinatesystem units <1.5mm,2mm>
\setplotarea x from 0 to 16, y from 0 to  11
\put{4)} [l] at 0 8
\put {$ \scriptstyle \bullet$} [c] at 8 8
\put {$ \scriptstyle \bullet$} [c] at 8 2
\put {$ \scriptstyle \bullet$} [c] at 6 5
\put {$ \scriptstyle \bullet$} [c] at 10 5
\setlinear \plot 8 2 6 5 8 8 10 5 8 2 /
\put{$12$} [c] at 8 0
\endpicture
\end{minipage}
\begin{minipage}{4cm}
\beginpicture
\setcoordinatesystem units <1.5mm,2mm>
\setplotarea x from 0 to 16, y from 0 to  11
\put{${\bf 7}$} [l] at 0 11
\put{5)} [l] at 0 8
\put {$ \scriptstyle \bullet$} [c] at 8 2
\put {$ \scriptstyle \bullet$} [c] at 6 8
\put {$ \scriptstyle \bullet$} [c] at 10 8
\put {$ \scriptstyle \bullet$} [c] at 7 5
\setlinear \plot 6 8 8 2 10 8  /
\put{$24$} [c] at 8 0
\endpicture
\end{minipage}
\begin{minipage}{4cm}
\beginpicture
\setcoordinatesystem units <1.5mm,2mm>
\setplotarea x from 0 to 16, y from 0 to  11
\put{6)} [l] at 0 8
\put {$ \scriptstyle \bullet$} [c] at 8 8
\put {$ \scriptstyle \bullet$} [c] at 6 2
\put {$ \scriptstyle \bullet$} [c] at 10 2
\put {$ \scriptstyle \bullet$} [c] at 7 5
\setlinear \plot 6 2 8 8 10 2  /
\put{$24$} [c] at 8 0
\endpicture
\end{minipage}
$$
$$
\begin{minipage}{4cm}
\beginpicture
\setcoordinatesystem units <1.5mm,2mm>
\setplotarea x from 0 to 16, y from 0 to  11
\put{7)} [l] at 0 8
\put {$ \scriptstyle \bullet$} [c] at 6 2
\put {$ \scriptstyle \bullet$} [c] at 10 2
\put {$ \scriptstyle \bullet$} [c] at 6 8
\put {$ \scriptstyle \bullet$} [c] at 10 8
\setlinear \plot 6 2 6 8 10 2 10 8 6 2 /
\put{$6$} [c] at 8 0
\endpicture
\end{minipage}
\begin{minipage}{4cm}
\beginpicture
\setcoordinatesystem units <1.5mm,2mm>
\setplotarea x from 0 to 16, y from 0 to  11
\put{${\bf 8}$} [l] at 0 11
\put{8)} [l] at 0 8
\put {$ \scriptstyle \bullet$} [c] at 6 2
\put {$ \scriptstyle \bullet$} [c] at 10 2
\put {$ \scriptstyle \bullet$} [c] at 6 8
\put {$ \scriptstyle \bullet$} [c] at 10 8
\setlinear \plot 6 8 6 2 10 8 10 2  /
\put{$24$} [c] at 8 0
\endpicture
\end{minipage}
\begin{minipage}{4cm}
\beginpicture
\setcoordinatesystem units <1.5mm,2mm>
\setplotarea x from 0 to 16, y from 0 to  11
\put{9)} [l] at 0 8
\put {$ \scriptstyle \bullet$} [c] at 6 2
\put {$ \scriptstyle \bullet$} [c] at 6 5
\put {$ \scriptstyle \bullet$} [c] at 6 8
\put {$ \scriptstyle \bullet$} [c] at 10 2
\setlinear \plot 6 2 6  8  /
\put{$24$} [c] at 8 0
\endpicture
\end{minipage}
\begin{minipage}{4cm}
\beginpicture
\setcoordinatesystem units <1.5mm,2mm>
\setplotarea x from 0 to 16, y from 0 to  11
\put{${\bf 9}$} [l] at 0 11
\put{10)} [l] at 0 8
\put {$ \scriptstyle \bullet$} [c] at 8 2
\put {$ \scriptstyle \bullet$} [c] at 8 8
\put {$ \scriptstyle \bullet$} [c] at 6  8
\put {$ \scriptstyle \bullet$} [c] at 10 8
\setlinear \plot 6 8 8 2 10 8  /
\setlinear \plot 8 2 8 8  /
\put{$4$} [c] at 8 0
\endpicture
\end{minipage}
\begin{minipage}{4cm}
\beginpicture
\setcoordinatesystem units <1.5mm,2mm>
\setplotarea x from 0 to 16, y from 0 to  11
\put{11)} [l] at 0 8
\put {$ \scriptstyle \bullet$} [c] at 8 2
\put {$ \scriptstyle \bullet$} [c] at 8 8
\put {$ \scriptstyle \bullet$} [c] at 6  2
\put {$ \scriptstyle \bullet$} [c] at 10 2
\setlinear \plot 6 2 8 8 10 2  /
\setlinear \plot 8 2 8 8  /
\put{$4$} [c] at 8 0
\endpicture
\end{minipage}
\begin{minipage}{4cm}
\beginpicture
\setcoordinatesystem units <1.5mm,2mm>
\setplotarea x from 0 to 16, y from 0 to  11
\put{12)} [l] at 0 8
\put {$ \scriptstyle \bullet$} [c] at 6 8
\put {$ \scriptstyle \bullet$} [c] at 6 2
\put {$ \scriptstyle \bullet$} [c] at 10 2
\put {$ \scriptstyle \bullet$} [c] at 10 8
\setlinear \plot 6 8 6 2   /
\setlinear \plot 10  2 10 8  /
\put{$12$} [c] at 8 0
\endpicture
\end{minipage}
$$
$$
\begin{minipage}{4cm}
\beginpicture
\setcoordinatesystem units <1.5mm,2mm>
\setplotarea x from 0 to 16, y from 0 to  11
\put{${\bf 10}$} [l] at 0 11
\put{13)} [l] at 0 8
\put {$ \scriptstyle \bullet$} [c] at 8 2
\put {$ \scriptstyle \bullet$} [c] at 6 8
\put {$ \scriptstyle \bullet$} [c] at 10 8
\put {$ \scriptstyle \bullet$} [c] at 12 2
\setlinear \plot 6 8 8 2 10 8   /
\put{$12$} [c] at 8 0
\endpicture
\end{minipage}
\begin{minipage}{4cm}
\beginpicture
\setcoordinatesystem units <1.5mm,2mm>
\setplotarea x from 0 to 16, y from 0 to  11
\put{14)} [l] at 0 8
\put {$ \scriptstyle \bullet$} [c] at 8 8
\put {$ \scriptstyle \bullet$} [c] at 6 2
\put {$ \scriptstyle \bullet$} [c] at 10 2
\put {$ \scriptstyle \bullet$} [c] at 12 2
\setlinear \plot 6 2 8 8 10 2   /
\put{$12$} [c] at 8 0
\endpicture
\end{minipage}
\begin{minipage}{4cm}
\beginpicture
\setcoordinatesystem units <1.5mm,2mm>
\setplotarea x from 0 to 16, y from 0 to  11
\put{${\bf 12}$} [l] at 0 11
\put{15)} [l] at 0 8
\put {$ \scriptstyle \bullet$} [c] at 6 2
\put {$ \scriptstyle \bullet$} [c] at 6 8
\put {$ \scriptstyle \bullet$} [c] at 8  2
\put {$ \scriptstyle \bullet$} [c] at 10 2
\setlinear \plot 6 2 6 8   /
\put{$12$} [c] at 8 0

\endpicture
\end{minipage}
\begin{minipage}{4cm}
\beginpicture
\setcoordinatesystem units <1.5mm,2mm>
\setplotarea x from 0 to 16, y from 0 to  11
\put{${\bf 16}$} [l] at 0 11

\put{16)} [l] at 0 8
\put {$ \scriptstyle \bullet$} [c] at 6 2
\put {$ \scriptstyle \bullet$} [c] at 8 2
\put {$ \scriptstyle \bullet$} [c] at 10 2
\put {$ \scriptstyle \bullet$} [c] at 12 2
\put{$1$} [c] at 9 0
\endpicture
\end{minipage}
$$

\medskip

  ${\bf NIP(5)}$

$$
\begin{minipage}{4cm}
\beginpicture
\setcoordinatesystem units <1.5mm,2mm>
\setplotarea x from 0 to 16, y from 0 to 13
\put{\bf 6} [l] at 0 13
\put{1)} [l] at 0 10
\put {$ \scriptstyle \bullet$} [c] at 8 2
\put {$ \scriptstyle \bullet$} [c] at 8 4
\put {$ \scriptstyle \bullet$} [c] at 8 6
\put {$ \scriptstyle \bullet$} [c] at 8 8
\put {$ \scriptstyle \bullet$} [c] at 8 10
\setlinear \plot 8 2 8 10   /
\put{$120$} [c] at 8 0
\endpicture
\end{minipage}
\begin{minipage}{4cm}
\beginpicture
\setcoordinatesystem units <1.5mm,2mm>
\setplotarea x from 0 to 16, y from 0 to 13
\put{\bf 7} [l] at 0 13
\put{2)} [l] at 0 10
\put {$ \scriptstyle \bullet$} [c] at 8 2
\put {$ \scriptstyle \bullet$} [c] at 8 6
\put {$ \scriptstyle \bullet$} [c] at 8 10
\put {$ \scriptstyle \bullet$} [c] at 6 8
\put {$ \scriptstyle \bullet$} [c] at 10 8
\setlinear \plot 8 2 8 6  10 8 8  10  6 8 8 6 /
\put{$60$} [c] at 8 0
\endpicture
\end{minipage}
\begin{minipage}{4cm}
\beginpicture
\setcoordinatesystem units <1.5mm,2mm>
\setplotarea x from 0 to 16, y from 0 to 13
\put{3)} [l] at 0 10
\put {$ \scriptstyle \bullet$} [c] at 8 2
\put {$ \scriptstyle \bullet$} [c] at 8 6
\put {$ \scriptstyle \bullet$} [c] at 8 10
\put {$ \scriptstyle \bullet$} [c] at 6 4
\put {$ \scriptstyle \bullet$} [c] at 10 4
\setlinear \plot 8 10 8 6  10 4 8  2  6 4 8 6 /
\put{$60$} [c] at 8 0
\endpicture
\end{minipage}
\begin{minipage}{4cm}
\beginpicture
\setcoordinatesystem units <1.5mm,2mm>
\setplotarea x from 0 to 16, y from 0 to 13
\put{4)} [l] at 0 10
\put {$ \scriptstyle \bullet$} [c] at 8 2
\put {$ \scriptstyle \bullet$} [c] at 8 6
\put {$ \scriptstyle \bullet$} [c] at 8 4
\put {$ \scriptstyle \bullet$} [c] at 6 10
\put {$ \scriptstyle \bullet$} [c] at 10 10
\setlinear \plot 8 2  8 6  6 10   /
\setlinear \plot 8 6  10 10   /
\put{$60$} [c] at 8 0
\endpicture
\end{minipage}
\begin{minipage}{4cm}
\beginpicture
\setcoordinatesystem units <1.5mm,2mm>
\setplotarea x from 0 to 16, y from 0 to 13
\put{5)} [l] at 0 10
\put {$ \scriptstyle \bullet$} [c] at 8 10
\put {$ \scriptstyle \bullet$} [c] at 8 8
\put {$ \scriptstyle \bullet$} [c] at 8 6
\put {$ \scriptstyle \bullet$} [c] at 6 2
\put {$ \scriptstyle \bullet$} [c] at 10 2
\setlinear \plot 8 10  8 6  6 2   /
\setlinear \plot 8 6  10 2   /
\put{$60$} [c] at 8 0
\endpicture
\end{minipage}
\begin{minipage}{4cm}
\beginpicture
\setcoordinatesystem units <1.5mm,2mm>
\setplotarea x from 0 to 16, y from 0 to 13
\put{\bf 8} [l] at 0 13
\put{6)} [l] at 0 10
\put {$ \scriptstyle \bullet$} [c] at 6 10
\put {$ \scriptstyle \bullet$} [c] at 10 10
\put {$ \scriptstyle \bullet$} [c] at 7 8
\put {$ \scriptstyle \bullet$} [c] at 8 6
\put {$ \scriptstyle \bullet$} [c] at 8 2
\setlinear \plot 6 10 8 6  8 2   /
\setlinear \plot 8 6  10 10   /
\put{$120$} [c] at 8 0
\endpicture
\end{minipage}
$$

$$
\begin{minipage}{4cm}
\beginpicture
\setcoordinatesystem units <1.5mm,2mm>
\setplotarea x from 0 to 16, y from 0 to 13
\put{7)} [l] at 0 10
\put {$ \scriptstyle \bullet$} [c] at 6 2
\put {$ \scriptstyle \bullet$} [c] at 10 2
\put {$ \scriptstyle \bullet$} [c] at 7 4
\put {$ \scriptstyle \bullet$} [c] at 8 6
\put {$ \scriptstyle \bullet$} [c] at 8 10
\setlinear \plot 6 2 8 6  8 10   /
\setlinear \plot 8 6  10 2   /
\put{$120$} [c] at 8 0
\endpicture
\end{minipage}
\begin{minipage}{4cm}
\beginpicture
\setcoordinatesystem units <1.5mm,2mm>
\setplotarea x from 0 to 16, y from 0 to 13
\put{8)} [l] at 0 10
\put {$ \scriptstyle \bullet$} [c] at 6 4
\put {$ \scriptstyle \bullet$} [c] at 6 8
\put {$ \scriptstyle \bullet$} [c] at 8 2
\put {$ \scriptstyle \bullet$} [c] at 8 10
\put {$ \scriptstyle \bullet$} [c] at 10 6
\setlinear \plot 8 2 6 4 6 8 8 10 10 6 8 2 /
\put{$120$} [c] at 8 0
\endpicture
\end{minipage}
\begin{minipage}{4cm}
\beginpicture
\setcoordinatesystem units <1.5mm,2mm>
\setplotarea x from 0 to 16, y from 0 to 13
\put{9)} [l] at 0 10
\put {$ \scriptstyle \bullet$} [c] at 8 2
\put {$ \scriptstyle \bullet$} [c] at 6 6
\put {$ \scriptstyle \bullet$} [c] at 6 10
\put {$ \scriptstyle \bullet$} [c] at 10 6
\put {$ \scriptstyle \bullet$} [c] at 10 10
\setlinear \plot  6 6  8 2  10 6 10 10  6 6 6 10 10 6   /
\put{$30$} [c] at 8 0
\endpicture
\end{minipage}
\begin{minipage}{4cm}
\beginpicture
\setcoordinatesystem units <1.5mm,2mm>
\setplotarea x from 0 to 16, y from 0 to 13
\put{10)} [l] at 0 10
\put {$ \scriptstyle \bullet$} [c] at 8 10
\put {$ \scriptstyle \bullet$} [c] at 6 6
\put {$ \scriptstyle \bullet$} [c] at 6 2
\put {$ \scriptstyle \bullet$} [c] at 10 6
\put {$ \scriptstyle \bullet$} [c] at 10 2
\setlinear \plot  6 6  8 10  10 6 10 2  6 6 6 2 10 6   /
\put{$30$} [c] at 8 0
\endpicture
\end{minipage}
\begin{minipage}{4cm}
\beginpicture
\setcoordinatesystem units <1.5mm,2mm>
\setplotarea x from 0 to 16, y from 0 to 13
\put{11)} [l] at 0 10
\put {$ \scriptstyle \bullet$} [c] at 6 2
\put {$ \scriptstyle \bullet$} [c] at 6 10
\put {$ \scriptstyle \bullet$} [c] at 8 6
\put {$ \scriptstyle \bullet$} [c] at 10 2
\put {$ \scriptstyle \bullet$} [c] at 10 10
\setlinear \plot 6 2 10 10   /
\setlinear \plot 6 10  10 2  /
\put{$30$} [c] at 8 0
\endpicture
\end{minipage}
\begin{minipage}{4cm}
\beginpicture
\setcoordinatesystem units <1.5mm,2mm>
\setplotarea x from 0 to 16, y from 0 to 13
\put{\bf 9} [l] at 0 13
\put{12)} [l] at 0 10
\put {$ \scriptstyle \bullet$} [c] at 6 10
\put {$ \scriptstyle \bullet$} [c] at 6.5 8
\put {$ \scriptstyle \bullet$} [c] at 7.3 4.7
\put {$ \scriptstyle \bullet$} [c] at 10 10
\put {$ \scriptstyle \bullet$} [c] at 8 2
\setlinear \plot  6 10  8 2   10 10    /
\put{$120$} [c] at 8 0
\endpicture
\end{minipage}
$$
$$
\begin{minipage}{4cm}
\beginpicture
\setcoordinatesystem units <1.5mm,2mm>
\setplotarea x from 0 to 16, y from 0 to 13
\put{13)} [l] at 0 10
\put {$ \scriptstyle \bullet$} [c] at 6 2
\put {$ \scriptstyle \bullet$} [c] at 6.6 4.3
\put {$ \scriptstyle \bullet$} [c] at 7.4 7.6
\put {$ \scriptstyle \bullet$} [c] at 10 2
\put {$ \scriptstyle \bullet$} [c] at 8 10
\setlinear \plot  6 2  8 10   10 2    /
\put{$120$} [c] at 8 0
\endpicture
\end{minipage}
\begin{minipage}{4cm}
\beginpicture
\setcoordinatesystem units <1.5mm,2mm>
\setplotarea x from 0 to 16, y from 0 to 13
\put{14)} [l] at 0 10
\put {$ \scriptstyle \bullet$} [c] at 6 10
\put {$ \scriptstyle \bullet$} [c] at 6 6
\put {$ \scriptstyle \bullet$} [c] at 8 2
\put {$ \scriptstyle \bullet$} [c] at 8 10
\put {$ \scriptstyle \bullet$} [c] at 10 6
\setlinear \plot 6 10 6 6  8 2 10 6 8 10 6 6  /
\put{$120$} [c] at 8 0
\endpicture
\end{minipage}
\begin{minipage}{4cm}
\beginpicture
\setcoordinatesystem units <1.5mm,2mm>
\setplotarea x from 0 to 16, y from 0 to 13
\put{15)} [l] at 0 10
\put {$ \scriptstyle \bullet$} [c] at 6 2
\put {$ \scriptstyle \bullet$} [c] at 6 6
\put {$ \scriptstyle \bullet$} [c] at 8 2
\put {$ \scriptstyle \bullet$} [c] at 8 10
\put {$ \scriptstyle \bullet$} [c] at 10 6
\setlinear \plot 6 2 6 6  8 10 10 6 8 2 6 6  /
\put{$120$} [c] at 8 0
\endpicture
\end{minipage}
\begin{minipage}{4cm}
\beginpicture
\setcoordinatesystem units <1.5mm,2mm>
\setplotarea x from 0 to 16, y from 0 to 13
\put{16)} [l] at 0 10
\put {$ \scriptstyle \bullet$} [c] at 6 2
\put {$ \scriptstyle \bullet$} [c] at 6 6
\put {$ \scriptstyle \bullet$} [c] at 6 10
\put {$ \scriptstyle \bullet$} [c] at 10 2
\put {$ \scriptstyle \bullet$} [c] at 10 10
\setlinear \plot  6 10  6 2   10 10 10 2 6 6   /
\put{$60$} [c] at 8 0
\endpicture
\end{minipage}
\begin{minipage}{4cm}
\beginpicture
\setcoordinatesystem units <1.5mm,2mm>
\setplotarea x from 0 to 16, y from 0 to 13
\put{17)} [l] at 0 10
\put {$ \scriptstyle \bullet$} [c] at 6 2
\put {$ \scriptstyle \bullet$} [c] at 6 6
\put {$ \scriptstyle \bullet$} [c] at 6 10
\put {$ \scriptstyle \bullet$} [c] at 10 2
\put {$ \scriptstyle \bullet$} [c] at 10 10
\setlinear \plot  6 2  6 10   10 2 10 10 6 6   /
\put{$60$} [c] at 8 0
\endpicture
\end{minipage}
\begin{minipage}{4cm}
\beginpicture
\setcoordinatesystem units <1.5mm,2mm>
\setplotarea x from 0 to 16, y from 0 to 13
\put{\bf 10} [l] at 0 13
\put{18)} [l] at 0 10
\put {$ \scriptstyle \bullet$} [c] at 6 2
\put {$ \scriptstyle \bullet$} [c] at 6  10
\put {$ \scriptstyle \bullet$} [c] at 8 2
\put {$ \scriptstyle \bullet$} [c] at 8 6
\put {$ \scriptstyle \bullet$} [c] at 10 10
\setlinear \plot 6 2 6  10    8 6 8 2    /
\setlinear \plot  8 6 10 10    /
\put{$120$} [c] at 8 0
\endpicture
\end{minipage}
$$
$$
\begin{minipage}{4cm}
\beginpicture
\setcoordinatesystem units <1.5mm,2mm>
\setplotarea x from 0 to 16, y from 0 to 13
\put{19)} [l] at 0 10
\put {$ \scriptstyle \bullet$} [c] at 6 2
\put {$ \scriptstyle \bullet$} [c] at 6  10
\put {$ \scriptstyle \bullet$} [c] at 8 10
\put {$ \scriptstyle \bullet$} [c] at 8 6
\put {$ \scriptstyle \bullet$} [c] at 10 2
\setlinear \plot 6 10 6  2    8 6 8 10    /
\setlinear \plot  8 6 10 2    /
\put{$120$} [c] at 8 0
\endpicture
\end{minipage}
\begin{minipage}{4cm}
\beginpicture
\setcoordinatesystem units <1.5mm,2mm>
\setplotarea x from 0 to 16, y from 0 to 13
\put{20)} [l] at 0 10
\put {$ \scriptstyle \bullet$} [c] at 6 2
\put {$ \scriptstyle \bullet$} [c] at 6 6
\put {$ \scriptstyle \bullet$} [c] at 6 10
\put {$ \scriptstyle \bullet$} [c] at 10 2
\put {$ \scriptstyle \bullet$} [c] at 10 10
\setlinear \plot  6 2  6 10   10 2 10 10 6 2   /
\put{$120$} [c] at 8 0
\endpicture
\end{minipage}
\begin{minipage}{4cm}
\beginpicture
\setcoordinatesystem units <1.5mm,2mm>
\setplotarea x from 0 to 16, y from 0 to 13
\put{21)} [l] at 0 10
\put {$ \scriptstyle \bullet$} [c] at 6 10
\put {$ \scriptstyle \bullet$} [c] at 6 6
\put {$ \scriptstyle \bullet$} [c] at 8 2
\put {$ \scriptstyle \bullet$} [c] at 10 10
\put {$ \scriptstyle \bullet$} [c] at 10 6
\setlinear \plot 6 10 6 6 8 2 10 6 10 10   /
\put{$60$} [c] at 8 0
\endpicture
\end{minipage}
\begin{minipage}{4cm}
\beginpicture
\setcoordinatesystem units <1.5mm,2mm>
\setplotarea x from 0 to 16, y from 0 to 13
\put{22)} [l] at 0 10
\put {$ \scriptstyle \bullet$} [c] at 6 2
\put {$ \scriptstyle \bullet$} [c] at 6 6
\put {$ \scriptstyle \bullet$} [c] at 8 10
\put {$ \scriptstyle \bullet$} [c] at 10 2
\put {$ \scriptstyle \bullet$} [c] at 10 6
\setlinear \plot 6 2 6 6 8 10 10 6 10 2   /
\put{$60$} [c] at 8 0
\endpicture
\end{minipage}
\begin{minipage}{4cm}
\beginpicture
\setcoordinatesystem units <1.5mm,2mm>
\setplotarea x from 0 to 16, y from 0 to 13
\put{23)} [l] at 0 10
\put {$ \scriptstyle \bullet$} [c] at 6 6
\put {$ \scriptstyle \bullet$} [c] at 8 2
\put {$ \scriptstyle \bullet$} [c] at 8 6
\put {$ \scriptstyle \bullet$} [c] at 8 10
\put {$ \scriptstyle \bullet$} [c] at 10 6
\setlinear \plot  6 6   8 10   10 6 8 2 6 6    /
\setlinear \plot  8 2 8 10    /
\put{$20$} [c] at 8 0
\endpicture
\end{minipage}
\begin{minipage}{4cm}
\beginpicture
\setcoordinatesystem units <1.5mm,2mm>
\setplotarea x from 0 to 16, y from 0 to 13
\put{24)} [l] at 0 10
\put {$ \scriptstyle \bullet$} [c] at 6 10
\put {$ \scriptstyle \bullet$} [c] at 8 10
\put {$ \scriptstyle \bullet$} [c] at 10 10
\put {$ \scriptstyle \bullet$} [c] at 8 6
\put {$ \scriptstyle \bullet$} [c] at 8 2
\setlinear \plot  8 2  8 10     /
\setlinear \plot  6 10 8 6  10 10     /
\put{$20$} [c] at 8 0
\endpicture
\end{minipage}
$$
$$
\begin{minipage}{4cm}
\beginpicture
\setcoordinatesystem units <1.5mm,2mm>
\setplotarea x from 0 to 16, y from 0 to 13
\put{25)} [l] at 0 10
\put {$ \scriptstyle \bullet$} [c] at 6 2
\put {$ \scriptstyle \bullet$} [c] at 8 2
\put {$ \scriptstyle \bullet$} [c] at 10 2
\put {$ \scriptstyle \bullet$} [c] at 8 6
\put {$ \scriptstyle \bullet$} [c] at 8 10
\setlinear \plot  8 2  8 10     /
\setlinear \plot  6 2 8 6  10 2     /
\put{$20$} [c] at 8 0
\endpicture
\end{minipage}
\begin{minipage}{4cm}
\beginpicture
\setcoordinatesystem units <1.5mm,2mm>
\setplotarea x from 0 to 16, y from 0 to 13
\put{26)} [l] at 0 10
\put {$ \scriptstyle \bullet$} [c] at 6 2
\put {$ \scriptstyle \bullet$} [c] at 6 4.5
\put {$ \scriptstyle \bullet$} [c] at 6 7.5
\put {$ \scriptstyle \bullet$} [c] at 6 10
\put {$ \scriptstyle \bullet$} [c] at 10 2
\setlinear \plot 6 2  6 10    /
\put{$120$} [c] at 8 0
\endpicture
\end{minipage}
\begin{minipage}{4cm}
\beginpicture
\setcoordinatesystem units <1.5mm,2mm>
\setplotarea x from 0 to 16, y from 0 to 13
\put{\bf 11} [l] at 0 13
\put{27)} [l] at 0 10
\put {$ \scriptstyle \bullet$} [c] at 6 2
\put {$ \scriptstyle \bullet$} [c] at 6 6
\put {$ \scriptstyle \bullet$} [c] at 6 10
\put {$ \scriptstyle \bullet$} [c] at 10 2
\put {$ \scriptstyle \bullet$} [c] at 10 10
\setlinear \plot  6  10  6 2   10 10 10 2   /
\put{$120$} [c] at 8 0
\endpicture
\end{minipage}
\begin{minipage}{4cm}
\beginpicture
\setcoordinatesystem units <1.5mm,2mm>
\setplotarea x from 0 to 16, y from 0 to 13
\put{28)} [l] at 0 10
\put {$ \scriptstyle \bullet$} [c] at 6 2
\put {$ \scriptstyle \bullet$} [c] at 6 6
\put {$ \scriptstyle \bullet$} [c] at 6 10
\put {$ \scriptstyle \bullet$} [c] at 10 2
\put {$ \scriptstyle \bullet$} [c] at 10 10
\setlinear \plot  6  2  6 10   10 2 10 10   /
\put{$120$} [c] at 8 0
\endpicture
\end{minipage}
\begin{minipage}{4cm}
\beginpicture
\setcoordinatesystem units <1.5mm,2mm>
\setplotarea x from 0 to 16, y from 0 to 13
\put{29)} [l] at 0 10
\put {$ \scriptstyle \bullet$} [c] at 6 10
\put {$ \scriptstyle \bullet$} [c] at 8.5 6
\put {$ \scriptstyle \bullet$} [c] at 9.5 2
\put {$ \scriptstyle \bullet$} [c] at 9.5 10
\put {$ \scriptstyle \bullet$} [c] at 10.5 6
\setlinear \plot  6 10   9.5 2 8.5 6 9.5 10 10.5 6 9.5 2   /
\put{$60$} [c] at 8 0
\endpicture
\end{minipage}
\begin{minipage}{4cm}
\beginpicture
\setcoordinatesystem units <1.5mm,2mm>
\setplotarea x from 0 to 16, y from 0 to 13
\put{30)} [l] at 0 10
\put {$ \scriptstyle \bullet$} [c] at 6 2
\put {$ \scriptstyle \bullet$} [c] at 8.5 6
\put {$ \scriptstyle \bullet$} [c] at 9.5 2
\put {$ \scriptstyle \bullet$} [c] at 9.5 10
\put {$ \scriptstyle \bullet$} [c] at 10.5 6
\setlinear \plot  6 2   9.5 10 8.5 6 9.5 2 10.5 6 9.5 10   /
\put{$60$} [c] at 8 0
\endpicture
\end{minipage}
$$
$$
\begin{minipage}{4cm}
\beginpicture
\setcoordinatesystem units <1.5mm,2mm>
\setplotarea x from 0 to 16, y from 0 to 13
\put{31)} [l] at 0 10
\put {$ \scriptstyle \bullet$} [c] at 6 10
\put {$ \scriptstyle \bullet$} [c] at 9 2
\put {$ \scriptstyle \bullet$} [c] at 9 6
\put {$ \scriptstyle \bullet$} [c] at 8 10
\put {$ \scriptstyle \bullet$} [c] at 10 10
\setlinear \plot 6 10 9 2 9  6   8 10    /
\setlinear \plot 9  6   10 10    /
\put{$60$} [c] at 8 0
\endpicture
\end{minipage}
\begin{minipage}{4cm}
\beginpicture
\setcoordinatesystem units <1.5mm,2mm>
\setplotarea x from 0 to 16, y from 0 to 13
\put{32)} [l] at 0 10
\put {$ \scriptstyle \bullet$} [c] at 6 2
\put {$ \scriptstyle \bullet$} [c] at 9 10
\put {$ \scriptstyle \bullet$} [c] at 9 6
\put {$ \scriptstyle \bullet$} [c] at 8 2
\put {$ \scriptstyle \bullet$} [c] at 10 2
\setlinear \plot 6 2 9 10 9  6   8 2    /
\setlinear \plot 9  6   10 2    /
\put{$60$} [c] at 8 0
\endpicture
\end{minipage}
\begin{minipage}{4cm}
\beginpicture
\setcoordinatesystem units <1.5mm,2mm>
\setplotarea x from 0 to 16, y from 0 to 13
\put{33)} [l] at 0 10
\put {$ \scriptstyle \bullet$} [c] at 6 2
\put {$ \scriptstyle \bullet$} [c] at 6 10
\put {$ \scriptstyle \bullet$} [c] at 8 10
\put {$ \scriptstyle \bullet$} [c] at 10 10
\put {$ \scriptstyle \bullet$} [c] at 10 2
\setlinear \plot 10 2 6 10 6 2 8 10 10 2 10 10 6 2  /
\put{$10$} [c] at 8 0
\endpicture
\end{minipage}
\begin{minipage}{4cm}
\beginpicture
\setcoordinatesystem units <1.5mm,2mm>
\setplotarea x from 0 to 16, y from 0 to 13
\put{34)} [l] at 0 10
\put {$ \scriptstyle \bullet$} [c] at 6 2
\put {$ \scriptstyle \bullet$} [c] at 6 10
\put {$ \scriptstyle \bullet$} [c] at 8 2
\put {$ \scriptstyle \bullet$} [c] at 10 10
\put {$ \scriptstyle \bullet$} [c] at 10 2
\setlinear \plot 10 10 6 2 6 10 8 2  10 10 10 2  6 10  /
\put{$10$} [c] at 8 0
\endpicture
\end{minipage}
\begin{minipage}{4cm}
\beginpicture
\setcoordinatesystem units <1.5mm,2mm>
\setplotarea x from 0 to 16, y from 0 to 13
\put{\bf 12} [l] at 0 13
\put{35)} [l] at 0 10
\put {$ \scriptstyle \bullet$} [c] at 6 2
\put {$ \scriptstyle \bullet$} [c] at 6 10
\put {$ \scriptstyle \bullet$} [c] at 8 6
\put {$ \scriptstyle \bullet$} [c] at 10 2
\put {$ \scriptstyle \bullet$} [c] at 10 10
\setlinear \plot  6 10   6 2 10  10 10 2   /
\put{$120$} [c] at 8 0
\endpicture
\end{minipage}
\begin{minipage}{4cm}
\beginpicture
\setcoordinatesystem units <1.5mm,2mm>
\setplotarea x from 0 to 16, y from 0 to 13
\put{36)} [l] at 0 10
\put {$ \scriptstyle \bullet$} [c] at 6 10
\put {$ \scriptstyle \bullet$} [c] at 8 2
\put {$ \scriptstyle \bullet$} [c] at 8 10
\put {$ \scriptstyle \bullet$} [c] at 10 2
\put {$ \scriptstyle \bullet$} [c] at 10 10
\setlinear \plot  6 10  8 2 8 10 10 2 10 10 8 2 /
\put{$60$} [c] at 8 0
\endpicture
\end{minipage}
$$
$$
\begin{minipage}{4cm}
\beginpicture
\setcoordinatesystem units <1.5mm,2mm>
\setplotarea x from 0 to 16, y from 0 to 13
\put{37)} [l] at 0 10
\put {$ \scriptstyle \bullet$} [c] at 6 2
\put {$ \scriptstyle \bullet$} [c] at 8 2
\put {$ \scriptstyle \bullet$} [c] at 8 10
\put {$ \scriptstyle \bullet$} [c] at 10 2
\put {$ \scriptstyle \bullet$} [c] at 10 10
\setlinear \plot  6 2  8 10 8 2 10 10 10 2 8 10 /
\put{$60$} [c] at 8 0
\endpicture
\end{minipage}
\begin{minipage}{4cm}
\beginpicture
\setcoordinatesystem units <1.5mm,2mm>
\setplotarea x from 0 to 16, y from 0 to 13
\put{38)} [l] at 0 10
\put {$ \scriptstyle \bullet$} [c] at 6 2
\put {$ \scriptstyle \bullet$} [c] at 6  6
\put {$ \scriptstyle \bullet$} [c] at 6 10
\put {$ \scriptstyle \bullet$} [c] at 10 2
\put {$ \scriptstyle \bullet$} [c] at 10 10
\setlinear \plot 6 2 6 10     /
\setlinear \plot 10 2  10 10    /
\put{$120$} [c] at 8 0
\endpicture
\end{minipage}
\begin{minipage}{4cm}
\beginpicture
\setcoordinatesystem units <1.5mm,2mm>
\setplotarea x from 0 to 16, y from 0 to 10
\put{39)} [l] at 0 10
\put {$ \scriptstyle \bullet$} [c] at 6 10
\put {$ \scriptstyle \bullet$} [c] at 7 2
\put {$ \scriptstyle \bullet$} [c] at 7 6
\put {$ \scriptstyle \bullet$} [c] at 8 10
\put {$ \scriptstyle \bullet$} [c] at 10 2
\setlinear \plot  6 10  7 6  7 2    /
\setlinear \plot   7 6  8 10     /
\put{$60$} [c] at 8 0
\endpicture
\end{minipage}
\begin{minipage}{4cm}
\beginpicture
\setcoordinatesystem units <1.5mm,2mm>
\setplotarea x from 0 to 16, y from 0 to 13
\put{40)} [l] at 0 10
\put {$ \scriptstyle \bullet$} [c] at 6 2
\put {$ \scriptstyle \bullet$} [c] at 7 10
\put {$ \scriptstyle \bullet$} [c] at 7 6
\put {$ \scriptstyle \bullet$} [c] at 8 2
\put {$ \scriptstyle \bullet$} [c] at 10 2
\setlinear \plot  6 2  7 6  7 10    /
\setlinear \plot   7 6  8 2     /
\put{$60$} [c] at 8 0
\endpicture
\end{minipage}
\begin{minipage}{4cm}
\beginpicture
\setcoordinatesystem units <1.5mm,2mm>
\setplotarea x from 0 to 16, y from 0 to 13
\put{41)} [l] at 0 10
\put {$ \scriptstyle \bullet$} [c] at 6 6
\put {$ \scriptstyle \bullet$} [c] at 7 2
\put {$ \scriptstyle \bullet$} [c] at 7 10
\put {$ \scriptstyle \bullet$} [c] at 10 2
\put {$ \scriptstyle \bullet$} [c] at 8 6
\setlinear \plot  6 6  7 2   8 6 7 10 6 6   /
\put{$60$} [c] at 8 0
\endpicture
\end{minipage}
\begin{minipage}{4cm}
\beginpicture
\setcoordinatesystem units <1.5mm,2mm>
\setplotarea x from 0 to 16, y from 0 to 13
\put{\bf 13} [l] at 0 13
\put{42)} [l] at 0 10
\put {$ \scriptstyle \bullet$} [c] at 6 10
\put {$ \scriptstyle \bullet$} [c] at 8 10
\put {$ \scriptstyle \bullet$} [c] at 10 10
\put {$ \scriptstyle \bullet$} [c] at 8 6
\put {$ \scriptstyle \bullet$} [c] at 8 2
\setlinear \plot  8 10   8 2     /
\setlinear \plot  6 10   8 2  10 10    /
\put{$60$} [c] at 8 0
\endpicture
\end{minipage}
$$
$$
\begin{minipage}{4cm}
\beginpicture
\setcoordinatesystem units <1.5mm,2mm>
\setplotarea x from 0 to 16, y from 0 to 13
\put{43)} [l] at 0 10
\put {$ \scriptstyle \bullet$} [c] at 6 2
\put {$ \scriptstyle \bullet$} [c] at 8 2
\put {$ \scriptstyle \bullet$} [c] at 10 2
\put {$ \scriptstyle \bullet$} [c] at 8 6
\put {$ \scriptstyle \bullet$} [c] at 8 10
\setlinear \plot  8 10   8 2     /
\setlinear \plot  6 2   8 10 10 2    /
\put{$60$} [c] at 8 0
\endpicture
\end{minipage}
\begin{minipage}{4cm}
\beginpicture
\setcoordinatesystem units <1.5mm,2mm>
\setplotarea x from 0 to 16, y from 0 to 13
\put{44)} [l] at 0 10
\put {$ \scriptstyle \bullet$} [c] at 6 2
\put {$ \scriptstyle \bullet$} [c] at 6 10
\put {$ \scriptstyle \bullet$} [c] at 8 10
\put {$ \scriptstyle \bullet$} [c] at 10 2
\put {$ \scriptstyle \bullet$} [c] at 10 10
\setlinear \plot  6 10 6 2 8 10 10 2   10 10  /
\put{$60$} [c] at 8 0

\endpicture
\end{minipage}
\begin{minipage}{4cm}
\beginpicture
\setcoordinatesystem units <1.5mm,2mm>
\setplotarea x from 0 to 16, y from 0 to 13
\put{45)} [l] at 0 10
\put {$ \scriptstyle \bullet$} [c] at 6 2
\put {$ \scriptstyle \bullet$} [c] at 6 10
\put {$ \scriptstyle \bullet$} [c] at 8 2
\put {$ \scriptstyle \bullet$} [c] at 10 2
\put {$ \scriptstyle \bullet$} [c] at 10 10
\setlinear \plot  6 2 6 10 8 2 10 10  10 2  /
\put{$60$} [c] at 8 0
\endpicture
\end{minipage}
\begin{minipage}{4cm}
\beginpicture
\setcoordinatesystem units <1.5mm,2mm>
\setplotarea x from 0 to 16, y from 0 to 13
\put{\bf 14} [l] at 0 13
\put{46)} [l] at 0 10
\put {$ \scriptstyle \bullet$} [c] at 6 10
\put {$ \scriptstyle \bullet$} [c] at 8 10
\put {$ \scriptstyle \bullet$} [c] at 10 10
\put {$ \scriptstyle \bullet$} [c] at 8 2
\put {$ \scriptstyle \bullet$} [c] at 10 2
\setlinear \plot  6  10   8 2 10 10 10 2    /
\setlinear \plot  8 10   8 2      /
\put{$60$} [c] at 8 0
\endpicture
\end{minipage}
\begin{minipage}{4cm}
\beginpicture
\setcoordinatesystem units <1.5mm,2mm>
\setplotarea x from 0 to 16, y from 0 to 13
\put{47)} [l] at 0 10
\put {$ \scriptstyle \bullet$} [c] at 6 2
\put {$ \scriptstyle \bullet$} [c] at 8 2
\put {$ \scriptstyle \bullet$} [c] at 10 2
\put {$ \scriptstyle \bullet$} [c] at 8 10
\put {$ \scriptstyle \bullet$} [c] at 10 10
\setlinear \plot  6  2   8 10 10 2 10 10    /
\setlinear \plot  8 10   8 2      /
\put{$60$} [c] at 8 0
\endpicture
\end{minipage}
\begin{minipage}{4cm}
\beginpicture
\setcoordinatesystem units <1.5mm,2mm>
\setplotarea x from 0 to 16, y from 0 to 13
\put{48)} [l] at 0 10
\put {$ \scriptstyle \bullet$} [c] at 6 10
\put {$ \scriptstyle \bullet$} [c] at 7 2
\put {$ \scriptstyle \bullet$} [c] at 8 10
\put {$ \scriptstyle \bullet$} [c] at 7.5 6
\put {$ \scriptstyle \bullet$} [c] at 10 2
\setlinear \plot  6 10 7 2 8 10    /
\put{$120$} [c] at 8 0

\endpicture
\end{minipage}
$$
$$
\begin{minipage}{4cm}
\beginpicture
\setcoordinatesystem units <1.5mm,2mm>
\setplotarea x from 0 to 16, y from 0 to 13
\put{49)} [l] at 0 10
\put {$ \scriptstyle \bullet$} [c] at 6 2
\put {$ \scriptstyle \bullet$} [c] at 7 10
\put {$ \scriptstyle \bullet$} [c] at 8 2
\put {$ \scriptstyle \bullet$} [c] at 7.5 6
\put {$ \scriptstyle \bullet$} [c] at 10 2
\setlinear \plot  6 2 7 10 8 2    /
\put{$120$} [c] at 8 0
\endpicture
\end{minipage}
\begin{minipage}{4cm}
\beginpicture
\setcoordinatesystem units <1.5mm,2mm>
\setplotarea x from 0 to 16, y from 0 to 13
\put{50)} [l] at 0 10
\put {$ \scriptstyle \bullet$} [c] at 6 2
\put {$ \scriptstyle \bullet$} [c] at 6 10
\put {$ \scriptstyle \bullet$} [c] at 8 2
\put {$ \scriptstyle \bullet$} [c] at 8 10
\put {$ \scriptstyle \bullet$} [c] at 10 2
\setlinear \plot  6 2 6 10 8 2  8 10  6 2  /
\put{$30$} [c] at 8 0
\endpicture
\end{minipage}
\begin{minipage}{4cm}
\beginpicture
\setcoordinatesystem units <1.5mm,2mm>
\setplotarea x from 0 to 16, y from 0 to 13
\put{\bf 15} [l] at 0 13
\put{51)} [l] at 0 10
\put {$ \scriptstyle \bullet$} [c] at 6 10
\put {$ \scriptstyle \bullet$} [c] at 7 2
\put {$ \scriptstyle \bullet$} [c] at 8 10
\put {$ \scriptstyle \bullet$} [c] at 10 2
\put {$ \scriptstyle \bullet$} [c] at 10 10
\setlinear \plot  6 10   7 2 8 10      /
\setlinear \plot   10 10 10 2       /
\put{$60$} [c] at 8 0
\endpicture
\end{minipage}
\begin{minipage}{4cm}
\beginpicture
\setcoordinatesystem units <1.5mm,2mm>
\setplotarea x from 0 to 16, y from 0 to 13
\put{52)} [l] at 0 10
\put {$ \scriptstyle \bullet$} [c] at 6 2
\put {$ \scriptstyle \bullet$} [c] at 7 10
\put {$ \scriptstyle \bullet$} [c] at 8 2
\put {$ \scriptstyle \bullet$} [c] at 10 2
\put {$ \scriptstyle \bullet$} [c] at 10 10
\setlinear \plot  6 2   7 10 8 2      /
\setlinear \plot   10 10 10 2       /
\put{$60$} [c] at 8 0
\endpicture
\end{minipage}
\begin{minipage}{4cm}
\beginpicture
\setcoordinatesystem units <1.5mm,2mm>
\setplotarea x from 0 to 16, y from 0 to 13
\put{\bf 16} [l] at 0 13
\put{53)} [l] at 0 10
\put {$ \scriptstyle \bullet$} [c] at 6 10
\put {$ \scriptstyle \bullet$} [c] at 6 2
\put {$ \scriptstyle \bullet$} [c] at 8 10
\put {$ \scriptstyle \bullet$} [c] at 8 2
\put {$ \scriptstyle \bullet$} [c] at 10 2
\setlinear \plot  6 10   6 2 8 10 8 2     /
\put{$120$} [c] at 8 0
\endpicture
\end{minipage}
\begin{minipage}{4cm}
\beginpicture
\setcoordinatesystem units <1.5mm,2mm>
\setplotarea x from 0 to 16, y from 0 to 13
\put{54)} [l] at 0 10
\put {$ \scriptstyle \bullet$} [c] at 6 2
\put {$ \scriptstyle \bullet$} [c] at 6 6
\put {$ \scriptstyle \bullet$} [c] at 6 10
\put {$ \scriptstyle \bullet$} [c] at 8 2
\put {$ \scriptstyle \bullet$} [c] at 10 2
\setlinear \plot  6  2  6  10      /
\put{$60$} [c] at 8 0
\endpicture
\end{minipage}
$$
$$
\begin{minipage}{4cm}
\beginpicture
\setcoordinatesystem units <1.5mm,2mm>
\setplotarea x from 0 to 16, y from 0 to 13
\put{\bf 17} [l] at 0 13
\put{55)} [l] at 0 10
\put {$ \scriptstyle \bullet$} [c] at 5.5 10
\put {$ \scriptstyle \bullet$} [c] at 7 10
\put {$ \scriptstyle \bullet$} [c] at 9 10
\put {$ \scriptstyle \bullet$} [c] at 10.5 10
\put {$ \scriptstyle \bullet$} [c] at 8 2
\setlinear \plot  5.5 10   8 2 10.5 10     /
\setlinear \plot  7 10   8 2 9 10      /
\put{$5$} [c] at 8 0
\endpicture
\end{minipage}
\begin{minipage}{4cm}
\beginpicture
\setcoordinatesystem units <1.5mm,2mm>
\setplotarea x from 0 to 16, y from 0 to 13
\put{56)} [l] at 0 10
\put {$ \scriptstyle \bullet$} [c] at 5.5 2
\put {$ \scriptstyle \bullet$} [c] at 7 2
\put {$ \scriptstyle \bullet$} [c] at 9 2
\put {$ \scriptstyle \bullet$} [c] at 10.5 2
\put {$ \scriptstyle \bullet$} [c] at 8 10
\setlinear \plot  5.5 2   8 10 10.5 2     /
\setlinear \plot  7 2   8 10 9 2      /
\put{$5$} [c] at 8 0
\endpicture
\end{minipage}
\begin{minipage}{4cm}
\beginpicture
\setcoordinatesystem units <1.5mm,2mm>
\setplotarea x from 0 to 16, y from 0 to 13
\put{\bf 18} [l] at 0 13
\put{57)} [l] at 0 10
\put {$ \scriptstyle \bullet$} [c] at 6 2
\put {$ \scriptstyle \bullet$} [c] at 6 10
\put {$ \scriptstyle \bullet$} [c] at 8 2
\put {$ \scriptstyle \bullet$} [c] at 8 10
\put {$ \scriptstyle \bullet$} [c] at 10 2
\setlinear \plot  6 2   6 10      /
\setlinear \plot  8 2 8 10     /
\put{$60$} [c] at 8 0
\endpicture
\end{minipage}
\begin{minipage}{4cm}
\beginpicture
\setcoordinatesystem units <1.5mm,2mm>
\setplotarea x from 0 to 16, y from 0 to 13
\put{58)} [l] at 0 10
\put {$ \scriptstyle \bullet$} [c] at 5.5 10
\put {$ \scriptstyle \bullet$} [c] at 7 10
\put {$ \scriptstyle \bullet$} [c] at 8.5 10
\put {$ \scriptstyle \bullet$} [c] at 7 2
\put {$ \scriptstyle \bullet$} [c] at 10  2
\setlinear \plot  5.5 10   7 2  8.5 10     /
\setlinear \plot   7 2 7 10     /
\put{$20$} [c] at 8 0
\endpicture
\end{minipage}
\begin{minipage}{4cm}
\beginpicture
\setcoordinatesystem units <1.5mm,2mm>
\setplotarea x from 0 to 16, y from 0 to 13
\put{59)} [l] at 0 10
\put {$ \scriptstyle \bullet$} [c] at 5.5 2
\put {$ \scriptstyle \bullet$} [c] at 7 2
\put {$ \scriptstyle \bullet$} [c] at 8.5 2
\put {$ \scriptstyle \bullet$} [c] at 7 10
\put {$ \scriptstyle \bullet$} [c] at 10  2
\setlinear \plot  5.5 2   7 10  8.5 2     /
\setlinear \plot   7 2 7 10     /
\put{$20$} [c] at 8 0
\endpicture
\end{minipage}
\begin{minipage}{4cm}
\beginpicture
\setcoordinatesystem units <1.5mm,2mm>
\setplotarea x from 0 to 16, y from 0 to 13
\put{\bf 20} [l] at 0 13
\put{60)} [l] at 0 10
\put {$ \scriptstyle \bullet$} [c] at 6 10
\put {$ \scriptstyle \bullet$} [c] at 7 2
\put {$ \scriptstyle \bullet$} [c] at 8 10
\put {$ \scriptstyle \bullet$} [c] at 9 2
\put {$ \scriptstyle \bullet$} [c] at 11  2
\setlinear \plot  6 10   7 2 8 10     /
\put{$30$} [c] at 8 0
\endpicture
\end{minipage}
$$
$$
\begin{minipage}{4cm}
\beginpicture
\setcoordinatesystem units <1.5mm,2mm>
\setplotarea x from 0 to 16, y from 0 to 13
\put{61)} [l] at 0 10
\put {$ \scriptstyle \bullet$} [c] at 6 2
\put {$ \scriptstyle \bullet$} [c] at 10  2
\put {$ \scriptstyle \bullet$} [c] at 7 10
\put {$ \scriptstyle \bullet$} [c] at 8 2
\put {$ \scriptstyle \bullet$} [c] at 10 2
\setlinear \plot  6 2   7 10 8 2     /
\put{$30$} [c] at 8 0
\endpicture
\end{minipage}
\begin{minipage}{4cm}
\beginpicture
\setcoordinatesystem units <1.5mm,2mm>
\setplotarea x from 0 to 16, y from 0 to 13
\put{\bf 24} [l] at 0 13
\put{62)} [l] at 0 10
\put {$ \scriptstyle \bullet$} [c] at 6 10
\put {$ \scriptstyle \bullet$} [c] at 6 2
\put {$ \scriptstyle \bullet$} [c] at 8 2
\put {$ \scriptstyle \bullet$} [c] at 10 2
\put {$ \scriptstyle \bullet$} [c] at 12  2
\setlinear \plot  6 2   6 10     /
\put{$20$} [c] at 8 0
\endpicture
\end{minipage}
\begin{minipage}{4cm}
\beginpicture
\setcoordinatesystem units <1.5mm,2mm>
\setplotarea x from 0 to 16, y from 0 to 13
\put{\bf 32} [l] at 0 13
\put{63)} [l] at 0 10
\put {$ \scriptstyle \bullet$} [c] at 6 2
\put {$ \scriptstyle \bullet$} [c] at 8 2
\put {$ \scriptstyle \bullet$} [c] at 10 2
\put {$ \scriptstyle \bullet$} [c] at 12 2
\put {$ \scriptstyle \bullet$} [c] at 14  2
\put{$1$} [c] at 10 0
\endpicture
\end{minipage}
$$

  ${\bf NIP(6)}$
\medskip

$$
\begin{minipage}{4cm}
\beginpicture
\setcoordinatesystem units <1.5mm,2mm>
\setplotarea x from 0 to 16, y from -2 to 15
\put{\bf 7} [l] at 0 15
\put{1)} [l] at 0 12
\put {$ \scriptstyle \bullet$} [c] at 12 0
\put {$ \scriptstyle \bullet$} [c] at 12 2.4
\put {$ \scriptstyle \bullet$} [c] at 12 4.8
\put {$ \scriptstyle \bullet$} [c] at 12 7.2
\put {$ \scriptstyle \bullet$} [c] at 12 9.6
\put {$ \scriptstyle \bullet$} [c] at 12  12
\setlinear \plot  12 0   12 12     /
\put{$720$} [c] at 12 -2
\endpicture
\end{minipage}
\begin{minipage}{4cm}
\beginpicture
\setcoordinatesystem units <1.5mm,2mm>
\setplotarea x from 0 to 16, y from -2 to 15
\put{\bf 8} [l] at 0 15
\put{2)} [l] at 0 12
\put {$ \scriptstyle \bullet$} [c] at 10 12
\put {$ \scriptstyle \bullet$} [c] at 10 6
\put {$ \scriptstyle \bullet$} [c] at 10 3
\put {$ \scriptstyle \bullet$} [c] at 10 0
\put {$ \scriptstyle \bullet$} [c] at 6 9
\put {$ \scriptstyle \bullet$} [c] at 14 9
\setlinear \plot  10 0 10 6  6 9 10 12 14 9 10  6     /
\put{$360$} [c] at 10 -2
\endpicture
\end{minipage}
\begin{minipage}{4cm}
\beginpicture
\setcoordinatesystem units <1.5mm,2mm>
\setplotarea x from 0 to 16, y from -2 to 15
\put{3)} [l] at 0 12
\put {$ \scriptstyle \bullet$} [c] at 10 12
\put {$ \scriptstyle \bullet$} [c] at 10 9
\put {$ \scriptstyle \bullet$} [c] at 10 6
\put {$ \scriptstyle \bullet$} [c] at 10 0
\put {$ \scriptstyle \bullet$} [c] at 6 3
\put {$ \scriptstyle \bullet$} [c] at 14 3
\setlinear \plot  10 12 10 6  6 3 10 0 14 3 10  6     /
\put{$360$} [c] at 10 -2
\endpicture
\end{minipage}
\begin{minipage}{4cm}
\beginpicture
\setcoordinatesystem units <1.5mm,2mm>
\setplotarea x from 0 to 16, y from -2 to 15
\put{4)} [l] at 0 12
\put {$ \scriptstyle \bullet$} [c] at 10 0
\put {$ \scriptstyle \bullet$} [c] at 10 2
\put {$ \scriptstyle \bullet$} [c] at 10 4
\put {$ \scriptstyle \bullet$} [c] at 10 6
\put {$ \scriptstyle \bullet$} [c] at 6 12
\put {$ \scriptstyle \bullet$} [c] at 14 12
\setlinear \plot  10 0 10 6  6 12      /
\setlinear \plot  10 6  14 12      /
\put{$360$} [c] at 10 -2
\endpicture
\end{minipage}
\begin{minipage}{4cm}
\beginpicture
\setcoordinatesystem units <1.5mm,2mm>
\setplotarea x from 0 to 16, y from -2 to 15
\put{5)} [l] at 0 12
\put {$ \scriptstyle \bullet$} [c] at 10 12
\put {$ \scriptstyle \bullet$} [c] at 10 10
\put {$ \scriptstyle \bullet$} [c] at 10 8
\put {$ \scriptstyle \bullet$} [c] at 10 6
\put {$ \scriptstyle \bullet$} [c] at 6 0
\put {$ \scriptstyle \bullet$} [c] at 14 0
\setlinear \plot  6 0  10 6  10 12      /
\setlinear \plot  10 6  14 0      /
\put{$360$} [c] at 10 -2

\endpicture
\end{minipage}
\begin{minipage}{4cm}
\beginpicture
\setcoordinatesystem units <1.5mm,2mm>
\setplotarea x from 0 to 16, y from -2 to 15
\put{6)} [l] at 0 12
\put {$ \scriptstyle \bullet$} [c] at 10 12
\put {$ \scriptstyle \bullet$} [c] at 10 10
\put {$ \scriptstyle \bullet$} [c] at 10 2
\put {$ \scriptstyle \bullet$} [c] at 10 0
\put {$ \scriptstyle \bullet$} [c] at 6  6
\put {$ \scriptstyle \bullet$} [c] at 14 6
\setlinear \plot  10 0 10 2 6 6 10  10  10 12     /
\setlinear \plot  10 10   14 6 10 2      /
\put{$360$} [c] at 8 -2
\endpicture
\end{minipage}
$$

$$
\begin{minipage}{4cm}
\beginpicture
\setcoordinatesystem units <1.5mm,2mm>
\setplotarea x from 0 to 16, y from -2 to 15
\put{\bf 9} [l] at 0 15
\put{7)} [l] at 0 12
\put {$ \scriptstyle \bullet$} [c] at 10 0
\put {$ \scriptstyle \bullet$} [c] at 10 2
\put {$ \scriptstyle \bullet$} [c] at 10  4
\put {$ \scriptstyle \bullet$} [c] at 6 12
\put {$ \scriptstyle \bullet$} [c] at 8 8
\put {$ \scriptstyle \bullet$} [c] at 14 12
\setlinear \plot  10 0  10 4  6 12      /
\setlinear \plot  10 4  14 12      /
\put{$720$} [c] at 10 -2
\endpicture
\end{minipage}
\begin{minipage}{4cm}
\beginpicture
\setcoordinatesystem units <1.5mm,2mm>
\setplotarea x from 0 to 16, y from -2 to 15
\put{8)} [l] at 0 12
\put {$ \scriptstyle \bullet$} [c] at 10 12
\put {$ \scriptstyle \bullet$} [c] at 10 10
\put {$ \scriptstyle \bullet$} [c] at 10 8
\put {$ \scriptstyle \bullet$} [c] at 6 0
\put {$ \scriptstyle \bullet$} [c] at 8 4
\put {$ \scriptstyle \bullet$} [c] at 14 0
\setlinear \plot  10 12  10 8  6 0      /
\setlinear \plot  10 8  14 0      /
\put{$720$} [c] at 10 -2
\endpicture
\end{minipage}
\begin{minipage}{4cm}
\beginpicture
\setcoordinatesystem units <1.5mm,2mm>
\setplotarea x from 0 to 16, y from -2 to 15
\put{9)} [l] at 0 12
\put {$ \scriptstyle \bullet$} [c] at 10 12
\put {$ \scriptstyle \bullet$} [c] at 10 8
\put {$ \scriptstyle \bullet$} [c] at 6 4
\put {$ \scriptstyle \bullet$} [c] at 10 0
\put {$ \scriptstyle \bullet$} [c] at 14 4
\put {$ \scriptstyle \bullet$} [c] at 8 6
\setlinear \plot   10 12 10 8 6 4 10 0 14 4 10 8     /
\put{$720$} [c] at 10 -2
\endpicture
\end{minipage}
\begin{minipage}{4cm}
\beginpicture
\setcoordinatesystem units <1.5mm,2mm>
\setplotarea x from 0 to 16, y from -2 to 15
\put{10)} [l] at 0 12
\put {$ \scriptstyle \bullet$} [c] at 10 12
\put {$ \scriptstyle \bullet$} [c] at 10 4
\put {$ \scriptstyle \bullet$} [c] at 6 8
\put {$ \scriptstyle \bullet$} [c] at 10 0
\put {$ \scriptstyle \bullet$} [c] at 14 8
\put {$ \scriptstyle \bullet$} [c] at 8 10
\setlinear \plot    10 0 10 4 6 8 10 12 14 8  10 4     /
\put{$720$} [c] at 10 -2
\endpicture
\end{minipage}
\begin{minipage}{4cm}
\beginpicture
\setcoordinatesystem units <1.5mm,2mm>
\setplotarea x from 0 to 16, y from -2 to 15
\put{11)} [l] at 0 12
\put {$ \scriptstyle \bullet$} [c] at 10 0
\put {$ \scriptstyle \bullet$} [c] at 10 4
\put {$ \scriptstyle \bullet$} [c] at 6 12
\put {$ \scriptstyle \bullet$} [c] at 6 8
\put {$ \scriptstyle \bullet$} [c] at 14 12
\put {$ \scriptstyle \bullet$} [c] at 14 8
\setlinear \plot    10 0 10 4 6 8 6 12 14 8 14 12 6 8  /
\setlinear \plot    10 4 14 8 /
\put{$180$} [c] at 10 -2
\endpicture
\end{minipage}
\begin{minipage}{4cm}
\beginpicture
\setcoordinatesystem units <1.5mm,2mm>
\setplotarea x from 0 to 16, y from -2 to 15
\put{12)} [l] at 0 12
\put {$ \scriptstyle \bullet$} [c] at 10 12
\put {$ \scriptstyle \bullet$} [c] at 10  8
\put {$ \scriptstyle \bullet$} [c] at 6 0
\put {$ \scriptstyle \bullet$} [c] at 6 4
\put {$ \scriptstyle \bullet$} [c] at 14 0
\put {$ \scriptstyle \bullet$} [c] at 14 4
\setlinear \plot    10 12 10 8 6 4 6 0 14 4 14 0 6 4  /
\setlinear \plot    10 8 14 4 /
\put{$180$} [c] at 10 -2
\endpicture
\end{minipage}
$$
$$
\begin{minipage}{4cm}
\beginpicture
\setcoordinatesystem units <1.5mm,2mm>
\setplotarea x from 0 to 16, y from -2 to 15
\put{13)} [l] at 0 12
\put {$ \scriptstyle \bullet$} [c] at 10 0
\put {$ \scriptstyle \bullet$} [c] at 10 10
\put {$ \scriptstyle \bullet$} [c] at 6 5
\put {$ \scriptstyle \bullet$} [c] at 6 12
\put {$ \scriptstyle \bullet$} [c] at 14 5
\put {$ \scriptstyle \bullet$} [c] at 14 12
\setlinear \plot  10 0  6 5 10  10  14  12      /
\setlinear \plot  10 0  14 5 10 10 6 12      /
\put{$180$} [c] at 10 -2
\endpicture
\end{minipage}
\begin{minipage}{4cm}
\beginpicture
\setcoordinatesystem units <1.5mm,2mm>
\setplotarea x from 0 to 16, y from -2 to 15
\put{14)} [l] at 0 12
\put {$ \scriptstyle \bullet$} [c] at 10 2
\put {$ \scriptstyle \bullet$} [c] at 10 12
\put {$ \scriptstyle \bullet$} [c] at 6 7
\put {$ \scriptstyle \bullet$} [c] at 6 0
\put {$ \scriptstyle \bullet$} [c] at 14 7
\put {$ \scriptstyle \bullet$} [c] at 14 0
\setlinear \plot  10 12  6 7 10 2  14  0      /
\setlinear \plot  10 12  14 7 10 2 6 0      /
\put{$180$} [c] at 10 -2
\endpicture
\end{minipage}
\begin{minipage}{4cm}
\beginpicture
\setcoordinatesystem units <1.5mm,2mm>
\setplotarea x from 0 to 16, y from -2 to 15
\put{15)} [l] at 0 12
\put {$ \scriptstyle \bullet$} [c] at 6 4
\put {$ \scriptstyle \bullet$} [c] at 6 8
\put {$ \scriptstyle \bullet$} [c] at 14 4
\put {$ \scriptstyle \bullet$} [c] at 14 8
\put {$ \scriptstyle \bullet$} [c] at 10 0
\put {$ \scriptstyle \bullet$} [c] at 10 12
\setlinear \plot 6 8 14 4  10 0  6 4 14  8 10 12  6 8 6 4  14 8   14 4 /
\put{$180$} [c] at 10 -2
\endpicture
\end{minipage}
\begin{minipage}{4cm}
\beginpicture
\setcoordinatesystem units <1.5mm,2mm>
\setplotarea x from 0 to 16, y from -2 to 15
\put{16)} [l] at 0 12
\put {$ \scriptstyle \bullet$} [c] at 6 0
\put {$ \scriptstyle \bullet$} [c] at 6 12
\put {$ \scriptstyle \bullet$} [c] at 10 4
\put {$ \scriptstyle \bullet$} [c] at 10 8
\put {$ \scriptstyle \bullet$} [c] at 14 0
\put {$ \scriptstyle \bullet$} [c] at 14 12
\setlinear \plot   6 0 10 4 10 8 6 12     /
\setlinear \plot   10 8 14 12     /
\setlinear \plot   14 0 10 4     /
\put{$180$} [c] at 10 -2
\endpicture
\end{minipage}
\begin{minipage}{4cm}
\beginpicture
\setcoordinatesystem units <1.5mm,2mm>
\setplotarea x from 0 to 16, y from -2 to 15
\put{\bf 10} [l] at 0 15
\put{17)} [l] at 0 12
\put {$ \scriptstyle \bullet$} [c] at 10 0
\put {$ \scriptstyle \bullet$} [c] at 10 2
\put {$ \scriptstyle \bullet$} [c] at 7.15 9
\put {$ \scriptstyle \bullet$} [c] at 8.35 6
\put {$ \scriptstyle \bullet$} [c] at 6 12
\put {$ \scriptstyle \bullet$} [c] at 14 12
\setlinear \plot    10 0 10 2  6 12  /
\setlinear \plot    10 2 14 12 /
\put{$720$} [c] at 10 -2
\endpicture
\end{minipage}
\begin{minipage}{4cm}
\beginpicture
\setcoordinatesystem units <1.5mm,2mm>
\setplotarea x from 0 to 16, y from -2 to 15
\put{18)} [l] at 0 12
\put {$ \scriptstyle \bullet$} [c] at 10 12
\put {$ \scriptstyle \bullet$} [c] at 10 10
\put {$ \scriptstyle \bullet$} [c] at 7.15 3
\put {$ \scriptstyle \bullet$} [c] at 8.35 6
\put {$ \scriptstyle \bullet$} [c] at 6 0
\put {$ \scriptstyle \bullet$} [c] at 14 0
\setlinear \plot    10 12 10 10  6 0  /
\setlinear \plot    10 10 14 0 /
\put{$720$} [c] at 10 -2
\endpicture
\end{minipage}
$$
$$
\begin{minipage}{4cm}
\beginpicture
\setcoordinatesystem units <1.5mm,2mm>
\setplotarea x from 0 to 16, y from -2 to 15
\put{19)} [l] at 0 12
\put {$ \scriptstyle \bullet$} [c] at 6 7
\put {$ \scriptstyle \bullet$} [c] at 6 12
\put {$ \scriptstyle \bullet$} [c] at 10 0
\put {$ \scriptstyle \bullet$} [c] at 10 2
\put {$ \scriptstyle \bullet$} [c] at 10 12
\put {$ \scriptstyle \bullet$} [c] at 14 7
\setlinear \plot  10 0  10 2 6 7  10 12  14  7 10 2     /
\setlinear \plot  6 7 6 12 /
\put{$720$} [c] at 10 -2
\endpicture
\end{minipage}
\begin{minipage}{4cm}
\beginpicture
\setcoordinatesystem units <1.5mm,2mm>
\setplotarea x from 0 to 16, y from -2 to 15
\put{20)} [l] at 0 12
\put {$ \scriptstyle \bullet$} [c] at 6 5
\put {$ \scriptstyle \bullet$} [c] at 6 0
\put {$ \scriptstyle \bullet$} [c] at 10 0
\put {$ \scriptstyle \bullet$} [c] at 10 10
\put {$ \scriptstyle \bullet$} [c] at 10 12
\put {$ \scriptstyle \bullet$} [c] at 14 5
\setlinear \plot  10 12  10 10 6 5  10 0  14  5 10 10     /
\setlinear \plot  6 5 6 0 /
\put{$720$} [c] at 10 -2
\endpicture
\end{minipage}
\begin{minipage}{4cm}
\beginpicture
\setcoordinatesystem units <1.5mm,2mm>
\setplotarea x from 0 to 16, y from -2 to 15
\put{21)} [l] at 0 12
\put {$ \scriptstyle \bullet$} [c] at 10 0
\put {$ \scriptstyle \bullet$} [c] at 10 12
\put {$ \scriptstyle \bullet$} [c] at 6 6
\put {$ \scriptstyle \bullet$} [c] at 14 6
\put {$ \scriptstyle \bullet$} [c] at 8 3
\put {$ \scriptstyle \bullet$} [c] at 12 9
\setlinear \plot  10 0  6 6 10 12  14  6 10 0     /
\setlinear \plot  8 3 12 9 /
\put{$720$} [c] at 10 -2
\endpicture
\end{minipage}
\begin{minipage}{4cm}
\beginpicture
\setcoordinatesystem units <1.5mm,2mm>
\setplotarea x from 0 to 16, y from -2 to 15
\put{22)} [l] at 0 12
\put {$ \scriptstyle \bullet$} [c] at 10 0
\put {$ \scriptstyle \bullet$} [c] at 10 12
\put {$ \scriptstyle \bullet$} [c] at 6 4
\put {$ \scriptstyle \bullet$} [c] at 6 6
\put {$ \scriptstyle \bullet$} [c] at 6 8
\put {$ \scriptstyle \bullet$} [c] at 14 6
\setlinear \plot    10 0 6 4 6 8 10 12 14 6 10 0  /
\put{$720$} [c] at 10 -2
\endpicture
\end{minipage}
\begin{minipage}{4cm}
\beginpicture
\setcoordinatesystem units <1.5mm,2mm>
\setplotarea x from 0 to 16, y from -2 to 15
\put{23)} [l] at 0 12
\put {$ \scriptstyle \bullet$} [c] at 6 0
\put {$ \scriptstyle \bullet$} [c] at 6 12
\put {$ \scriptstyle \bullet$} [c] at 10 6
\put {$ \scriptstyle \bullet$} [c] at 8 9
\put {$ \scriptstyle \bullet$} [c] at 14 0
\put {$ \scriptstyle \bullet$} [c] at 14 12
\setlinear \plot    6 0 14  12 /
\setlinear \plot    6 12 14 0 /
\put{$360$} [c] at 10 -2
\endpicture
\end{minipage}
\begin{minipage}{4cm}
\beginpicture
\setcoordinatesystem units <1.5mm,2mm>
\setplotarea x from 0 to 16, y from -2 to 15
\put{24)} [l] at 0 12
\put {$ \scriptstyle \bullet$} [c] at 6 0
\put {$ \scriptstyle \bullet$} [c] at 6 12
\put {$ \scriptstyle \bullet$} [c] at 10 6
\put {$ \scriptstyle \bullet$} [c] at 8 3
\put {$ \scriptstyle \bullet$} [c] at 14 0
\put {$ \scriptstyle \bullet$} [c] at 14 12
\setlinear \plot    6 0 14  12 /
\setlinear \plot    6 12 14 0 /
\put{$360$} [c] at 10 -2
\endpicture
\end{minipage}
$$

$$
\begin{minipage}{4cm}
\beginpicture
\setcoordinatesystem units <1.5mm,2mm>
\setplotarea x from 0 to 16, y from -2 to 15
\put{25)} [l] at 0 12
\put {$ \scriptstyle \bullet$} [c] at 10 0
\put {$ \scriptstyle \bullet$} [c] at 6 2
\put {$ \scriptstyle \bullet$} [c] at 6 7
\put {$ \scriptstyle \bullet$} [c] at 6 12
\put {$ \scriptstyle \bullet$} [c] at 14 2
\put {$ \scriptstyle \bullet$} [c] at 14 12
\setlinear \plot    6 12 6 2  10 0 14 2 14 12 6 2 /
\setlinear \plot    6 7 14 2 /
\put{$360$} [c] at 10 -2
\endpicture
\end{minipage}
\begin{minipage}{4cm}
\beginpicture
\setcoordinatesystem units <1.5mm,2mm>
\setplotarea x from 0 to 16, y from -2 to 15
\put{26)} [l] at 0 12
\put {$ \scriptstyle \bullet$} [c] at 10 12
\put {$ \scriptstyle \bullet$} [c] at 6 10
\put {$ \scriptstyle \bullet$} [c] at 6 5
\put {$ \scriptstyle \bullet$} [c] at 6 0
\put {$ \scriptstyle \bullet$} [c] at 14 10
\put {$ \scriptstyle \bullet$} [c] at 14 0
\setlinear \plot    6 0 6 10  10 12 14 10 14 0 6 10 /
\setlinear \plot    6 5 14 10 /
\put{$360$} [c] at 10 -2
\endpicture
\end{minipage}
\begin{minipage}{4cm}
\beginpicture
\setcoordinatesystem units <1.5mm,2mm>
\setplotarea x from 0 to 16, y from -2 to 15
\put{27)} [l] at 0 12
\put {$ \scriptstyle \bullet$} [c] at 10 0
\put {$ \scriptstyle \bullet$} [c] at 8 4
\put {$ \scriptstyle \bullet$} [c] at 6 12
\put {$ \scriptstyle \bullet$} [c] at 6 8
\put {$ \scriptstyle \bullet$} [c] at 14 12
\put {$ \scriptstyle \bullet$} [c] at 14 8
\setlinear \plot    10 0 6 8 6 12 14 8 14 12 6 8  /
\setlinear \plot    10 0 14 8 /
\put{$360$} [c] at 10 -2
\endpicture
\end{minipage}
\begin{minipage}{4cm}
\beginpicture
\setcoordinatesystem units <1.5mm,2mm>
\setplotarea x from 0 to 16, y from -2 to 15
\put{28)} [l] at 0 12
\put {$ \scriptstyle \bullet$} [c] at 10 12
\put {$ \scriptstyle \bullet$} [c] at 8 8
\put {$ \scriptstyle \bullet$} [c] at 6 0
\put {$ \scriptstyle \bullet$} [c] at 6 4
\put {$ \scriptstyle \bullet$} [c] at 14 0
\put {$ \scriptstyle \bullet$} [c] at 14 4
\setlinear \plot    10 12 6 4 6 0 14 4 14 0 6 4  /
\setlinear \plot    10 12 14 4 /
\put{$360$} [c] at 10 -2
\endpicture
\end{minipage}
\begin{minipage}{4cm}
\beginpicture
\setcoordinatesystem units <1.5mm,2mm>
\setplotarea x from 0 to 16, y from -2 to 15
\put{29)} [l] at 0 12
\put {$ \scriptstyle \bullet$} [c] at 6 0
\put {$ \scriptstyle \bullet$} [c] at 6 6
\put {$ \scriptstyle \bullet$} [c] at 6 12
\put {$ \scriptstyle \bullet$} [c] at 14 0
\put {$ \scriptstyle \bullet$} [c] at 14 6
\put {$ \scriptstyle \bullet$} [c] at 14 12
\setlinear \plot    6 12 6 0 14 6  6 12  /
\setlinear \plot    14 12 14 0 6 6 14 12 /
\put{$90$} [c] at 10 -2
\endpicture
\end{minipage}
\begin{minipage}{4cm}
\beginpicture
\setcoordinatesystem units <1.5mm,2mm>
\setplotarea x from 0 to 16, y from -2 to 15
\put{\bf 11} [l] at 0 15
\put{30)} [l] at 0 12
\put {$ \scriptstyle \bullet$} [c] at 6 12
\put {$ \scriptstyle \bullet$} [c] at 14 12
\put {$ \scriptstyle \bullet$} [c] at 10 0
\put {$ \scriptstyle \bullet$} [c] at 7.2 8.5
\put {$ \scriptstyle \bullet$} [c] at 8.2 5.5
\put {$ \scriptstyle \bullet$} [c] at 9 2.8
\setlinear \plot  6 12  10 0 14 12      /
\put{$720$} [c] at 10 -2
\endpicture
\end{minipage}
$$
$$
\begin{minipage}{4cm}
\beginpicture
\setcoordinatesystem units <1.5mm,2mm>
\setplotarea x from 0 to 16, y from -2 to 15
\put{31)} [l] at 0 12
\put {$ \scriptstyle \bullet$} [c] at 6 0
\put {$ \scriptstyle \bullet$} [c] at 14 0
\put {$ \scriptstyle \bullet$} [c] at 10 12
\put {$ \scriptstyle \bullet$} [c] at 7 2.8
\put {$ \scriptstyle \bullet$} [c] at 7.8 5.5
\put {$ \scriptstyle \bullet$} [c] at 8.8 8.5
\setlinear \plot  6 0  10 12 14 0      /
\put{$720$} [c] at 10 -2
\endpicture
\end{minipage}
\begin{minipage}{4cm}
\beginpicture
\setcoordinatesystem units <1.5mm,2mm>
\setplotarea x from 0 to 16, y from -2 to 15
\put{32)} [l] at 0 12
\put {$ \scriptstyle \bullet$} [c] at 6 6
\put {$ \scriptstyle \bullet$} [c] at 6 12
\put {$ \scriptstyle \bullet$} [c] at 10 0
\put {$ \scriptstyle \bullet$} [c] at 10 12
\put {$ \scriptstyle \bullet$} [c] at 14 6
\put {$ \scriptstyle \bullet$} [c] at 8 3
\setlinear \plot  6 12  6 6 10 0 14 6 10 12 6 6      /
\put{$720$} [c] at 10 -2
\endpicture
\end{minipage}
\begin{minipage}{4cm}
\beginpicture
\setcoordinatesystem units <1.5mm,2mm>
\setplotarea x from 0 to 16, y from -2 to 15
\put{33)} [l] at 0 12
\put {$ \scriptstyle \bullet$} [c] at 6 6
\put {$ \scriptstyle \bullet$} [c] at 6 0
\put {$ \scriptstyle \bullet$} [c] at 10 0
\put {$ \scriptstyle \bullet$} [c] at 10 12
\put {$ \scriptstyle \bullet$} [c] at 14 6
\put {$ \scriptstyle \bullet$} [c] at 8 9
\setlinear \plot  6 0  6 6 10 0 14 6 10 12 6 6      /
\put{$720$} [c] at 10 -2
\endpicture
\end{minipage}
\begin{minipage}{4cm}
\beginpicture
\setcoordinatesystem units <1.5mm,2mm>
\setplotarea x from 0 to 16, y from -2 to 15
\put{34)} [l] at 0 12
\put {$ \scriptstyle \bullet$} [c] at 10 0
\put {$ \scriptstyle \bullet$} [c] at 10 10
\put {$ \scriptstyle \bullet$} [c] at 10 12
\put {$ \scriptstyle \bullet$} [c] at 6 5
\put {$ \scriptstyle \bullet$} [c] at 14  5
\put {$ \scriptstyle \bullet$} [c] at 14 12
\setlinear \plot  14 12 14 5  10 0 6 5 10 10 10 12    /
\setlinear \plot   10 10 14 5  /
\put{$720$} [c] at 10 -2
\endpicture
\end{minipage}
\begin{minipage}{4cm}
\beginpicture
\setcoordinatesystem units <1.5mm,2mm>
\setplotarea x from 0 to 16, y from -2 to 15
\put{35)} [l] at 0 12
\put {$ \scriptstyle \bullet$} [c] at 10 0
\put {$ \scriptstyle \bullet$} [c] at 10 2
\put {$ \scriptstyle \bullet$} [c] at 10 12
\put {$ \scriptstyle \bullet$} [c] at 6 7
\put {$ \scriptstyle \bullet$} [c] at 14  7
\put {$ \scriptstyle \bullet$} [c] at 14 0
\setlinear \plot  14 0 14 7  10 12 6 7 10 2 10 0    /
\setlinear \plot   10 2 14 7  /
\put{$720$} [c] at 10 -2
\endpicture
\end{minipage}
\begin{minipage}{4cm}
\beginpicture
\setcoordinatesystem units <1.5mm,2mm>
\setplotarea x from 0 to 16, y from -2 to 12
\put{36)} [l] at 0 12
\put {$ \scriptstyle \bullet$} [c] at 6 2
\put {$ \scriptstyle \bullet$} [c] at 6 7
\put {$ \scriptstyle \bullet$} [c] at 6 12
\put {$ \scriptstyle \bullet$} [c] at 10 0
\put {$ \scriptstyle \bullet$} [c] at 14 2
\put {$ \scriptstyle \bullet$} [c] at 14 12
\setlinear \plot 14 2 10 0 6 2 6 12 14 2 14 12 6 2     /
\put{$720$} [c] at 10 -2
\endpicture
\end{minipage}
$$
$$
\begin{minipage}{4cm}
\beginpicture
\setcoordinatesystem units <1.5mm,2mm>
\setplotarea x from 0 to 16, y from -2 to 15
\put{37)} [l] at 0 12
\put {$ \scriptstyle \bullet$} [c] at 6 10
\put {$ \scriptstyle \bullet$} [c] at 6 5
\put {$ \scriptstyle \bullet$} [c] at 6 0
\put {$ \scriptstyle \bullet$} [c] at 10 12
\put {$ \scriptstyle \bullet$} [c] at 14 10
\put {$ \scriptstyle \bullet$} [c] at 14 0
\setlinear \plot 14 10 10 12 6 10 6 0 14 10 14 0 6 10     /
\put{$720$} [c] at 10 -2
\endpicture
\end{minipage}
\begin{minipage}{4cm}
\beginpicture
\setcoordinatesystem units <1.5mm,2mm>
\setplotarea x from 0 to 16, y from -2 to 15
\put{38)} [l] at 0 12
\put {$ \scriptstyle \bullet$} [c] at 6 0
\put {$ \scriptstyle \bullet$} [c] at 6 4
\put {$ \scriptstyle \bullet$} [c] at 6 8
\put {$ \scriptstyle \bullet$} [c] at 6 12
\put {$ \scriptstyle \bullet$} [c] at 14 0
\put {$ \scriptstyle \bullet$} [c] at 14 12
\setlinear \plot    6 0 6 12   /
\setlinear \plot    6 4 14 12 14 0 6 8   /
\put{$720$} [c] at 10 -2
\endpicture
\end{minipage}
\begin{minipage}{4cm}
\beginpicture
\setcoordinatesystem units <1.5mm,2mm>
\setplotarea x from 0 to 16, y from -2 to 15
\put{39)} [l] at 0 12
\put {$ \scriptstyle \bullet$} [c] at 10 0
\put {$ \scriptstyle \bullet$} [c] at 10 4
\put {$ \scriptstyle \bullet$} [c] at 6 8
\put {$ \scriptstyle \bullet$} [c] at 6 12
\put {$ \scriptstyle \bullet$} [c] at 14  8
\put {$ \scriptstyle \bullet$} [c] at 14 12
\setlinear \plot    6 12 6 8 10 4 14 8 14 12    /
\setlinear \plot   10 0 10 4  /
\put{$360$} [c] at 10 -2
\endpicture
\end{minipage}
\begin{minipage}{4cm}
\beginpicture
\setcoordinatesystem units <1.5mm,2mm>
\setplotarea x from 0 to 16, y from -2 to 15
\put{40)} [l] at 0 12
\put {$ \scriptstyle \bullet$} [c] at 10 12
\put {$ \scriptstyle \bullet$} [c] at 10 8
\put {$ \scriptstyle \bullet$} [c] at 6 4
\put {$ \scriptstyle \bullet$} [c] at 6 0
\put {$ \scriptstyle \bullet$} [c] at 14 4
\put {$ \scriptstyle \bullet$} [c] at 14 0
\setlinear \plot    6 0 6 4 10 8 14 4 14 0    /
\setlinear \plot  10 8 10 12  /
\put{$360$} [c] at 10 -2
\endpicture
\end{minipage}
\begin{minipage}{4cm}
\beginpicture
\setcoordinatesystem units <1.5mm,2mm>
\setplotarea x from 0 to 16, y from -2 to 15
\put{41)} [l] at 0 12
\put {$ \scriptstyle \bullet$} [c] at 6 0
\put {$ \scriptstyle \bullet$} [c] at 6 4
\put {$ \scriptstyle \bullet$} [c] at 6 8
\put {$ \scriptstyle \bullet$} [c] at 6 12
\put {$ \scriptstyle \bullet$} [c] at 14 0
\put {$ \scriptstyle \bullet$} [c] at 14 12
\setlinear \plot 6 12 6 0 14 12 14 0 6 4    /
\put{$360$} [c] at 10 -2
\endpicture
\end{minipage}
\begin{minipage}{4cm}
\beginpicture
\setcoordinatesystem units <1.5mm,2mm>
\setplotarea x from 0 to 16, y from -2 to 15
\put{42)} [l] at 0 12
\put {$ \scriptstyle \bullet$} [c] at 6 0
\put {$ \scriptstyle \bullet$} [c] at 6 4
\put {$ \scriptstyle \bullet$} [c] at 6 8
\put {$ \scriptstyle \bullet$} [c] at 6 12
\put {$ \scriptstyle \bullet$} [c] at 14 0
\put {$ \scriptstyle \bullet$} [c] at 14 12
\setlinear \plot 6 0 6 12 14  0 14 12 6 8    /
\put{$360$} [c] at 10 -2
\endpicture
\end{minipage}
$$
$$
\begin{minipage}{4cm}
\beginpicture
\setcoordinatesystem units <1.5mm,2mm>
\setplotarea x from 0 to 16, y from -2 to 15
\put{43)} [l] at 0 12
\put {$ \scriptstyle \bullet$} [c] at 6 4
\put {$ \scriptstyle \bullet$} [c] at 6 12
\put {$ \scriptstyle \bullet$} [c] at 10 0
\put {$ \scriptstyle \bullet$} [c] at 14 0
\put {$ \scriptstyle \bullet$} [c] at 14 4
\put {$ \scriptstyle \bullet$} [c] at 14 12
\setlinear \plot 14 0 14 12 6 4 6 12 14 4 10 0 6 4    /
\put{$360$} [c] at 10 -2
\endpicture
\end{minipage}
\begin{minipage}{4cm}
\beginpicture
\setcoordinatesystem units <1.5mm,2mm>
\setplotarea x from 0 to 16, y from -2 to 15
\put{44)} [l] at 0 12
\put {$ \scriptstyle \bullet$} [c] at 6 8
\put {$ \scriptstyle \bullet$} [c] at 6 0
\put {$ \scriptstyle \bullet$} [c] at 10 12
\put {$ \scriptstyle \bullet$} [c] at 14 0
\put {$ \scriptstyle \bullet$} [c] at 14 8
\put {$ \scriptstyle \bullet$} [c] at 14 12
\setlinear \plot 14 12 14 0 6 8 6 0 14 8 10 12 6 8    /
\put{$360$} [c] at 10 -2
\endpicture
\end{minipage}
\begin{minipage}{4cm}
\beginpicture
\setcoordinatesystem units <1.5mm,2mm>
\setplotarea x from 0 to 16, y from -2 to 15
\put{45)} [l] at 0 12
\put {$ \scriptstyle \bullet$} [c] at 6 4
\put {$ \scriptstyle \bullet$} [c] at 6 8
\put {$ \scriptstyle \bullet$} [c] at 10 0
\put {$ \scriptstyle \bullet$} [c] at 10 12
\put {$ \scriptstyle \bullet$} [c] at 14 4
\put {$ \scriptstyle \bullet$} [c] at 14 8
\setlinear \plot  10 0  6 4 6 8  10 12  14 8 14 4  10 0     /
\put{$360$} [c] at 10 -2
\endpicture
\end{minipage}
\begin{minipage}{4cm}
\beginpicture
\setcoordinatesystem units <1.5mm,2mm>
\setplotarea x from 0 to 16, y from -2 to 15
\put{46)} [l] at 0 12
\put {$ \scriptstyle \bullet$} [c] at 10 12
\put {$ \scriptstyle \bullet$} [c] at 10 8
\put {$ \scriptstyle \bullet$} [c] at 10 4
\put {$ \scriptstyle \bullet$} [c] at 10 0
\put {$ \scriptstyle \bullet$} [c] at 6 4
\put {$ \scriptstyle \bullet$} [c] at 14 4
\setlinear \plot    10 0 6 4 10 8 14 4 10 0   /
\setlinear \plot   10 0 10 12  /
\put{$120$} [c] at 10 -2
\endpicture
\end{minipage}
\begin{minipage}{4cm}
\beginpicture
\setcoordinatesystem units <1.5mm,2mm>
\setplotarea x from 0 to 16, y from -2 to 15
\put{47)} [l] at 0 12
\put {$ \scriptstyle \bullet$} [c] at 10 0
\put {$ \scriptstyle \bullet$} [c] at 10 4
\put {$ \scriptstyle \bullet$} [c] at 10 8
\put {$ \scriptstyle \bullet$} [c] at 10 12
\put {$ \scriptstyle \bullet$} [c] at 6 8
\put {$ \scriptstyle \bullet$} [c] at 14 8
\setlinear \plot    10 12 6 8 10 4 14 8 10 12   /
\setlinear \plot   10 0 10 12  /
\put{$120$} [c] at 10 -2
\endpicture
\end{minipage}
\begin{minipage}{4cm}
\beginpicture
\setcoordinatesystem units <1.5mm,2mm>
\setplotarea x from 0 to 16, y from -2 to 15
\put{48)} [l] at 0 12
\put {$ \scriptstyle \bullet$} [c] at 10 0
\put {$ \scriptstyle \bullet$} [c] at 10 4
\put {$ \scriptstyle \bullet$} [c] at 10 8
\put {$ \scriptstyle \bullet$} [c] at 10 12
\put {$ \scriptstyle \bullet$} [c] at 6 12
\put {$ \scriptstyle \bullet$} [c] at 14 12
\setlinear \plot  10 0 10 12     /
\setlinear \plot  6 12 10 8 14 12  /
\put{$120$} [c] at 10 -2
\endpicture
\end{minipage}
$$
$$
\begin{minipage}{4cm}
\beginpicture
\setcoordinatesystem units <1.5mm,2mm>
\setplotarea x from 0 to 16, y from -2 to 15
\put{49)} [l] at 0 12
\put {$ \scriptstyle \bullet$} [c] at 10 0
\put {$ \scriptstyle \bullet$} [c] at 10 4
\put {$ \scriptstyle \bullet$} [c] at 10 8
\put {$ \scriptstyle \bullet$} [c] at 10 12
\put {$ \scriptstyle \bullet$} [c] at 6 0
\put {$ \scriptstyle \bullet$} [c] at 14 0
\setlinear \plot  10 0 10 12     /
\setlinear \plot  6 0  10 4 14 0  /
\put{$120$} [c] at 10 -2
\endpicture
\end{minipage}
\begin{minipage}{4cm}
\beginpicture
\setcoordinatesystem units <1.5mm,2mm>
\setplotarea x from 0 to 16, y from -2 to 15
\put{\bf 12} [l] at 0 15
\put{50)} [l] at 0 12
\put {$ \scriptstyle \bullet$} [c] at 6 6
\put {$ \scriptstyle \bullet$} [c] at 8 9
\put {$ \scriptstyle \bullet$} [c] at 10 0
\put {$ \scriptstyle \bullet$} [c] at 10 12
\put {$ \scriptstyle \bullet$} [c] at 14 6
\put {$ \scriptstyle \bullet$} [c] at 14 12
\setlinear \plot 14 6 10 12  6 6 10  0  14 6 14 12 /
\put{$720$} [c] at 10 -2
\endpicture
\end{minipage}
\begin{minipage}{4cm}
\beginpicture
\setcoordinatesystem units <1.5mm,2mm>
\setplotarea x from 0 to 16, y from -2 to 15
\put{51)} [l] at 0 12
\put {$ \scriptstyle \bullet$} [c] at 6 6
\put {$ \scriptstyle \bullet$} [c] at 8 3
\put {$ \scriptstyle \bullet$} [c] at 10 0
\put {$ \scriptstyle \bullet$} [c] at 10 12
\put {$ \scriptstyle \bullet$} [c] at 14 6
\put {$ \scriptstyle \bullet$} [c] at 14 0
\setlinear \plot 14 6 10 0  6 6 10  12  14 6 14 0 /
\put{$720$} [c] at 10 -2
\endpicture
\end{minipage}
\begin{minipage}{4cm}
\beginpicture
\setcoordinatesystem units <1.5mm,2mm>
\setplotarea x from 0 to 16, y from -2 to 15
\put{52)} [l] at 0 12
\put {$ \scriptstyle \bullet$} [c] at 6 12
\put {$ \scriptstyle \bullet$} [c] at 8 9
\put {$ \scriptstyle \bullet$} [c] at 10 6
\put {$ \scriptstyle \bullet$} [c] at 12 0
\put {$ \scriptstyle \bullet$} [c] at 12 12
\put {$ \scriptstyle \bullet$} [c] at 14 6
\setlinear \plot 6 12 10 6 12 0 14  6 12  12 10 6  /
\put{$720$} [c] at 10 -2
\endpicture
\end{minipage}
\begin{minipage}{4cm}
\beginpicture
\setcoordinatesystem units <1.5mm,2mm>
\setplotarea x from 0 to 16, y from -2 to 15
\put{53)} [l] at 0 12
\put {$ \scriptstyle \bullet$} [c] at 6 0
\put {$ \scriptstyle \bullet$} [c] at 8 3
\put {$ \scriptstyle \bullet$} [c] at 10 6
\put {$ \scriptstyle \bullet$} [c] at 12 0
\put {$ \scriptstyle \bullet$} [c] at 12 12
\put {$ \scriptstyle \bullet$} [c] at 14 6
\setlinear \plot 6 0 10 6 12 0 14  6 12  12 10 6  /
\put{$720$} [c] at 10 -2
\endpicture
\end{minipage}
\begin{minipage}{4cm}
\beginpicture
\setcoordinatesystem units <1.5mm,2mm>
\setplotarea x from 0 to 16, y from -2 to 15
\put{54)} [l] at 0 12
\put {$ \scriptstyle \bullet$} [c] at 6 0
\put {$ \scriptstyle \bullet$} [c] at 6 4
\put {$ \scriptstyle \bullet$} [c] at 6 8
\put {$ \scriptstyle \bullet$} [c] at 6 12
\put {$ \scriptstyle \bullet$} [c] at 14 0
\put {$ \scriptstyle \bullet$} [c] at 14 12
\setlinear \plot 6 12 6 0  /
\setlinear \plot 6 8 14 12 14 0  /
\put{$720$} [c] at 10 -2
\endpicture
\end{minipage}
$$
$$
\begin{minipage}{4cm}
\beginpicture
\setcoordinatesystem units <1.5mm,2mm>
\setplotarea x from 0 to 16, y from -2 to 15
\put{55)} [l] at 0 12
\put {$ \scriptstyle \bullet$} [c] at 6 0
\put {$ \scriptstyle \bullet$} [c] at 6 4
\put {$ \scriptstyle \bullet$} [c] at 6 8
\put {$ \scriptstyle \bullet$} [c] at 6 12
\put {$ \scriptstyle \bullet$} [c] at 14 0
\put {$ \scriptstyle \bullet$} [c] at 14 12
\setlinear \plot 6 12 6 0  /
\setlinear \plot 6 4 14 0 14 12  /
\put{$720$} [c] at 10 -2
\endpicture
\end{minipage}
\begin{minipage}{4cm}
\beginpicture
\setcoordinatesystem units <1.5mm,2mm>
\setplotarea x from 0 to 16, y from -2 to 15
\put{56)} [l] at 0 12
\put {$ \scriptstyle \bullet$} [c] at 6 0
\put {$ \scriptstyle \bullet$} [c] at 6 4
\put {$ \scriptstyle \bullet$} [c] at 6 8
\put {$ \scriptstyle \bullet$} [c] at 6 12
\put {$ \scriptstyle \bullet$} [c] at 14 0
\put {$ \scriptstyle \bullet$} [c] at 14 12
\setlinear \plot 6 12 6 0  14 12 14 0 6 8 /
\put{$720$} [c] at 10 -2
\endpicture
\end{minipage}
\begin{minipage}{4cm}
\beginpicture
\setcoordinatesystem units <1.5mm,2mm>
\setplotarea x from 0 to 16, y from -2 to 15
\put{57)} [l] at 0 12
\put {$ \scriptstyle \bullet$} [c] at 6 0
\put {$ \scriptstyle \bullet$} [c] at 6 4
\put {$ \scriptstyle \bullet$} [c] at 6 8
\put {$ \scriptstyle \bullet$} [c] at 6 12
\put {$ \scriptstyle \bullet$} [c] at 14 0
\put {$ \scriptstyle \bullet$} [c] at 14 12
\setlinear \plot 6 0 6 12  14 0 14 12 6 4 /
\put{$720$} [c] at 10 -2
\endpicture
\end{minipage}
\begin{minipage}{4cm}
\beginpicture
\setcoordinatesystem units <1.5mm,2mm>
\setplotarea x from 0 to 16, y from -2 to 15
\put{58)} [l] at 0 12
\put {$ \scriptstyle \bullet$} [c] at 6 6
\put {$ \scriptstyle \bullet$} [c] at 8 0
\put {$ \scriptstyle \bullet$} [c] at 8 12
\put {$ \scriptstyle \bullet$} [c] at 10 6
\put {$ \scriptstyle \bullet$} [c] at 14 0
\put {$ \scriptstyle \bullet$} [c] at 14 12
\setlinear \plot 8 0 6 6 8 12 10 6 8 0 /
\setlinear \plot 14 0 10 6  14 12 /
\put{$720$} [c] at 10 -2
\endpicture
\end{minipage}
\begin{minipage}{4cm}
\beginpicture
\setcoordinatesystem units <1.5mm,2mm>
\setplotarea x from 0 to 16, y from -2 to 15
\put{59)} [l] at 0 12
\put {$ \scriptstyle \bullet$} [c] at 6 0
\put {$ \scriptstyle \bullet$} [c] at 6 8
\put {$ \scriptstyle \bullet$} [c] at 6 12
\put {$ \scriptstyle \bullet$} [c] at 14 0
\put {$ \scriptstyle \bullet$} [c] at 14 4
\put {$ \scriptstyle \bullet$} [c] at 14 12
\setlinear \plot 6 12 6 0   /
\setlinear \plot 14 12 14 0 /
\setlinear \plot 6 8 14 4 /
\put{$720$} [c] at 10 -2
\endpicture
\end{minipage}
\begin{minipage}{4cm}
\beginpicture
\setcoordinatesystem units <1.5mm,2mm>
\setplotarea x from 0 to 16, y from -2 to 15
\put{60)} [l] at 0 12
\put {$ \scriptstyle \bullet$} [c] at 6 0
\put {$ \scriptstyle \bullet$} [c] at 6 6
\put {$ \scriptstyle \bullet$} [c] at 6 12
\put {$ \scriptstyle \bullet$} [c] at 14 0
\put {$ \scriptstyle \bullet$} [c] at 14 6
\put {$ \scriptstyle \bullet$} [c] at 14 12
\setlinear \plot 6 0 6 12 14 6  14 12 6 0 /
\setlinear \plot 6 6 14 0 14 6 /
\put{$720$} [c] at 10 -2
\endpicture
\end{minipage}
$$
$$
\begin{minipage}{4cm}
\beginpicture
\setcoordinatesystem units <1.5mm,2mm>
\setplotarea x from 0 to 16, y from -2 to 15
\put{61)} [l] at 0 12
\put {$ \scriptstyle \bullet$} [c] at 6 8
\put {$ \scriptstyle \bullet$} [c] at 8 0
\put {$ \scriptstyle \bullet$} [c] at 8 4
\put {$ \scriptstyle \bullet$} [c] at 8 12
\put {$ \scriptstyle \bullet$} [c] at 10 8
\put {$ \scriptstyle \bullet$} [c] at 14 6
\setlinear \plot 8 4 6 8 8 12 10 8 8 4 8 0 /
\setlinear \plot  8 0 14 6 8 12  /
\put{$360$} [c] at 10 -2
\endpicture
\end{minipage}
\begin{minipage}{4cm}
\beginpicture
\setcoordinatesystem units <1.5mm,2mm>
\setplotarea x from 0 to 16, y from -2 to 15
\put{62)} [l] at 0 12
\put {$ \scriptstyle \bullet$} [c] at 6 4
\put {$ \scriptstyle \bullet$} [c] at 8 0
\put {$ \scriptstyle \bullet$} [c] at 8 8
\put {$ \scriptstyle \bullet$} [c] at 8 12
\put {$ \scriptstyle \bullet$} [c] at 10 4
\put {$ \scriptstyle \bullet$} [c] at 14 6
\setlinear \plot 8 8  6 4 8 0 10 4 8 8 8 12 /
\setlinear \plot  8 0 14 6 8 12  /
\put{$360$} [c] at 10 -2
\endpicture
\end{minipage}
\begin{minipage}{4cm}
\beginpicture
\setcoordinatesystem units <1.5mm,2mm>
\setplotarea x from 0 to 16, y from -2 to 15
\put{63)} [l] at 0 12
\put {$ \scriptstyle \bullet$} [c] at 6 8
\put {$ \scriptstyle \bullet$} [c] at 8 0
\put {$ \scriptstyle \bullet$} [c] at 8 4
\put {$ \scriptstyle \bullet$} [c] at 8 12
\put {$ \scriptstyle \bullet$} [c] at 10 8
\put {$ \scriptstyle \bullet$} [c] at 14 12
\setlinear \plot 8 4 6 8 8 12 10 8 8 4 8 0 /
\setlinear \plot  8 4 14 12  /
\put{$360$} [c] at 10 -2
\endpicture
\end{minipage}
\begin{minipage}{4cm}
\beginpicture
\setcoordinatesystem units <1.5mm,2mm>
\setplotarea x from 0 to 16, y from -2 to 15
\put{64)} [l] at 0 12
\put {$ \scriptstyle \bullet$} [c] at 6 4
\put {$ \scriptstyle \bullet$} [c] at 8 0
\put {$ \scriptstyle \bullet$} [c] at 8 8
\put {$ \scriptstyle \bullet$} [c] at 8 12
\put {$ \scriptstyle \bullet$} [c] at 10 4
\put {$ \scriptstyle \bullet$} [c] at 14 0
\setlinear \plot 8 8 6 4 8 0 10 4 8 8 8 12 /
\setlinear \plot  8 8 14 0  /
\put{$360$} [c] at 10 -2
\endpicture
\end{minipage}
\begin{minipage}{4cm}
\beginpicture
\setcoordinatesystem units <1.5mm,2mm>
\setplotarea x from 0 to 16, y from -2 to 15
\put{65)} [l] at 0 12
\put {$ \scriptstyle \bullet$} [c] at 6 12
\put {$ \scriptstyle \bullet$} [c] at 10 0
\put {$ \scriptstyle \bullet$} [c] at 10 4
\put {$ \scriptstyle \bullet$} [c] at 10 8
\put {$ \scriptstyle \bullet$} [c] at 10 12
\put {$ \scriptstyle \bullet$} [c] at 14 12
\setlinear \plot 10 0 10 12  /
\setlinear \plot 6 12 10 4 /
\setlinear \plot 14 12 10 8 /
\put{$360$} [c] at 10 -2
\endpicture
\end{minipage}
\begin{minipage}{4cm}
\beginpicture
\setcoordinatesystem units <1.5mm,2mm>
\setplotarea x from 0 to 16, y from -2 to 15
\put{66)} [l] at 0 12
\put {$ \scriptstyle \bullet$} [c] at 6 0
\put {$ \scriptstyle \bullet$} [c] at 10 0
\put {$ \scriptstyle \bullet$} [c] at 10 4
\put {$ \scriptstyle \bullet$} [c] at 10 8
\put {$ \scriptstyle \bullet$} [c] at 10 12
\put {$ \scriptstyle \bullet$} [c] at 14 0
\setlinear \plot 10 0 10 12  /
\setlinear \plot 6 0 10 8 /
\setlinear \plot 14 0 10 4 /
\put{$360$} [c] at 10 -2
\endpicture
\end{minipage}
$$
$$
\begin{minipage}{4cm}
\beginpicture
\setcoordinatesystem units <1.5mm,2mm>
\setplotarea x from 0 to 16, y from -2 to 15
\put{67)} [l] at 0 12
\put {$ \scriptstyle \bullet$} [c] at 6 0
\put {$ \scriptstyle \bullet$} [c] at 6 6
\put {$ \scriptstyle \bullet$} [c] at 6 12
\put {$ \scriptstyle \bullet$} [c] at 14 0
\put {$ \scriptstyle \bullet$} [c] at 14 6
\put {$ \scriptstyle \bullet$} [c] at 14 12
\setlinear \plot 6 0 6 12 14 6 14 12 6  6  /
\setlinear \plot 14 0 14 6 /
\put{$180$} [c] at 10 -2
\endpicture
\end{minipage}
\begin{minipage}{4cm}
\beginpicture
\setcoordinatesystem units <1.5mm,2mm>
\setplotarea x from 0 to 16, y from -2 to 15
\put{68)} [l] at 0 12
\put {$ \scriptstyle \bullet$} [c] at 6 0
\put {$ \scriptstyle \bullet$} [c] at 6 6
\put {$ \scriptstyle \bullet$} [c] at 6 12
\put {$ \scriptstyle \bullet$} [c] at 14 0
\put {$ \scriptstyle \bullet$} [c] at 14 6
\put {$ \scriptstyle \bullet$} [c] at 14 12
\setlinear \plot 6 12 6 0 14 6 14 0 6  6  /
\setlinear \plot 14 12 14 6 /
\put{$180$} [c] at 10 -2
\endpicture
\end{minipage}
\begin{minipage}{4cm}
\beginpicture
\setcoordinatesystem units <1.5mm,2mm>
\setplotarea x from 0 to 16, y from -2 to 15
\put{69)} [l] at 0 12
\put {$ \scriptstyle \bullet$} [c] at 6 6
\put {$ \scriptstyle \bullet$} [c] at  6 12
\put {$ \scriptstyle \bullet$} [c] at 10 0
\put {$ \scriptstyle \bullet$} [c] at 10 12
\put {$ \scriptstyle \bullet$} [c] at 14 6
\put {$ \scriptstyle \bullet$} [c] at 14 12
\setlinear \plot 6 12 6 6  10 0 14  6 14  12 6 6  10 12 14 6 6 12 /
\put{$60$} [c] at 10 -2
\endpicture
\end{minipage}
\begin{minipage}{4cm}
\beginpicture
\setcoordinatesystem units <1.5mm,2mm>
\setplotarea x from 0 to 16, y from -2 to 15
\put{70)} [l] at 0 12
\put {$ \scriptstyle \bullet$} [c] at 6 6
\put {$ \scriptstyle \bullet$} [c] at  6 0
\put {$ \scriptstyle \bullet$} [c] at 10 0
\put {$ \scriptstyle \bullet$} [c] at 10 12
\put {$ \scriptstyle \bullet$} [c] at 14 6
\put {$ \scriptstyle \bullet$} [c] at 14 0
\setlinear \plot 6 0 6 6  10 12 14  6 14  0 6 6  10 0 14 6 6 0 /
\put{$60$} [c] at 10 -2
\endpicture
\end{minipage}
\begin{minipage}{4cm}
\beginpicture
\setcoordinatesystem units <1.5mm,2mm>
\setplotarea x from 0 to 16, y from -2 to 15
\put{71)} [l] at 0 12
\put {$ \scriptstyle \bullet$} [c] at 6 12
\put {$ \scriptstyle \bullet$} [c] at 10 12
\put {$ \scriptstyle \bullet$} [c] at 14 12
\put {$ \scriptstyle \bullet$} [c] at 10 6
\put {$ \scriptstyle \bullet$} [c] at 6 0
\put {$ \scriptstyle \bullet$} [c] at 14 0
\setlinear \plot 6 12 14 0  /
\setlinear \plot 6 0 14 12   /
\setlinear \plot 10 6 10 12   /
\put{$60$} [c] at 10 -2
\endpicture
\end{minipage}
\begin{minipage}{4cm}
\beginpicture
\setcoordinatesystem units <1.5mm,2mm>
\setplotarea x from 0 to 16, y from -2 to 15
\put{72)} [l] at 0 12
\put {$ \scriptstyle \bullet$} [c] at 6 0
\put {$ \scriptstyle \bullet$} [c] at 10 0
\put {$ \scriptstyle \bullet$} [c] at 14 0
\put {$ \scriptstyle \bullet$} [c] at 10 6
\put {$ \scriptstyle \bullet$} [c] at 6 12
\put {$ \scriptstyle \bullet$} [c] at 14 12
\setlinear \plot 6 12 14 0  /
\setlinear \plot 6 0 14 12   /
\setlinear \plot 10 6 10 0   /
\put{$60$} [c] at 10 -2
\endpicture
\end{minipage}
$$
$$
\begin{minipage}{4cm}
\beginpicture
\setcoordinatesystem units <1.5mm,2mm>
\setplotarea x from 0 to 16, y from -2 to 15
\put{73)} [l] at 0 12
\put {$ \scriptstyle \bullet$} [c] at 6 6
\put {$ \scriptstyle \bullet$} [c] at 6 12
\put {$ \scriptstyle \bullet$} [c] at 8 0
\put {$ \scriptstyle \bullet$} [c] at 10 12
\put {$ \scriptstyle \bullet$} [c] at 10 6
\put {$ \scriptstyle \bullet$} [c] at 14 6
\setlinear \plot 6 12 6 6 8 0 14 6 10  12 10 6 8 0 /
\setlinear \plot  14 6 6 12 10 6 /
\setlinear \plot  6 6  10 12  /
\put{$60$} [c] at 10 -2
\endpicture
\end{minipage}
\begin{minipage}{4cm}
\beginpicture
\setcoordinatesystem units <1.5mm,2mm>
\setplotarea x from 0 to 16, y from -2 to 15
\put{74)} [l] at 0 12
\put {$ \scriptstyle \bullet$} [c] at 6 6
\put {$ \scriptstyle \bullet$} [c] at 6 0
\put {$ \scriptstyle \bullet$} [c] at 8 12
\put {$ \scriptstyle \bullet$} [c] at 10 0
\put {$ \scriptstyle \bullet$} [c] at 10 6
\put {$ \scriptstyle \bullet$} [c] at 14 6
\setlinear \plot 6 0 6 6 8 12 14 6 10  0 10 6 8 12 /
\setlinear \plot  14 6 6 0 10 6 /
\setlinear \plot  6 6  10 0  /
\put{$60$} [c] at 10 -2
\endpicture
\end{minipage}
\begin{minipage}{4cm}
\beginpicture
\setcoordinatesystem units <1.5mm,2mm>
\setplotarea x from 0 to 16, y from -2 to 15
\put{75)} [l] at 0 12
\put {$ \scriptstyle \bullet$} [c] at 6 0
\put {$ \scriptstyle \bullet$} [c] at 6 3
\put {$ \scriptstyle \bullet$} [c] at 6 6
\put {$ \scriptstyle \bullet$} [c] at 6 9
\put {$ \scriptstyle \bullet$} [c] at 6 12
\put {$ \scriptstyle \bullet$} [c] at 14 0
\setlinear \plot 6 12 6 0  /
\put{$720$} [c] at 10 -2
\endpicture
\end{minipage}
\begin{minipage}{4cm}
\beginpicture
\setcoordinatesystem units <1.5mm,2mm>
\setplotarea x from 0 to 16, y from -2 to 15
\put{${\bf 13}$} [l] at 0 15
\put{76)} [l] at 0 12
\put {$ \scriptstyle \bullet$} [c] at 6 4
\put {$ \scriptstyle \bullet$} [c] at 6 8
\put {$ \scriptstyle \bullet$} [c] at 6 12
\put {$ \scriptstyle \bullet$} [c] at 10 0
\put {$ \scriptstyle \bullet$} [c] at 14 4
\put {$ \scriptstyle \bullet$} [c] at 14 12
\setlinear \plot 6 12 6 4 10 0 14 4 14 12  /
\put{$720$} [c] at 10 -2
\endpicture
\end{minipage}
\begin{minipage}{4cm}
\beginpicture
\setcoordinatesystem units <1.5mm,2mm>
\setplotarea x from 0 to 16, y from -2 to 15
\put{77)} [l] at 0 12
\put {$ \scriptstyle \bullet$} [c] at 6 0
\put {$ \scriptstyle \bullet$} [c] at 6 4
\put {$ \scriptstyle \bullet$} [c] at 6 8
\put {$ \scriptstyle \bullet$} [c] at 10 12
\put {$ \scriptstyle \bullet$} [c] at 14 0
\put {$ \scriptstyle \bullet$} [c] at 14 8
\setlinear \plot 6 0 6 8 10 12 14 8 14 0  /
\put{$720$} [c] at 10 -2
\endpicture
\end{minipage}
\begin{minipage}{4cm}
\beginpicture
\setcoordinatesystem units <1.5mm,2mm>
\setplotarea x from 0 to 16, y from -2 to 15
\put{78)} [l] at 0 12
\put {$ \scriptstyle \bullet$} [c] at 6 6
\put {$ \scriptstyle \bullet$} [c] at 6 12
\put {$ \scriptstyle \bullet$} [c] at 8 9
\put {$ \scriptstyle \bullet$} [c] at 10 0
\put {$ \scriptstyle \bullet$} [c] at 10 12
\put {$ \scriptstyle \bullet$} [c] at 14 6
\setlinear \plot 6 12 6 6 10 0 14 6 10 12 6 6   /
\put{$720$} [c] at 10 -2
\endpicture
\end{minipage}
$$
$$
\begin{minipage}{4cm}
\beginpicture
\setcoordinatesystem units <1.5mm,2mm>
\setplotarea x from 0 to 16, y from -2 to 15
\put{79)} [l] at 0 12
\put {$ \scriptstyle \bullet$} [c] at 6 6
\put {$ \scriptstyle \bullet$} [c] at 6 0
\put {$ \scriptstyle \bullet$} [c] at 8 3
\put {$ \scriptstyle \bullet$} [c] at 10 0
\put {$ \scriptstyle \bullet$} [c] at 10 12
\put {$ \scriptstyle \bullet$} [c] at 14 6
\setlinear \plot 6 0 6 6 10 0 14 6 10 12 6 6   /
\put{$720$} [c] at 10 -2
\endpicture
\end{minipage}
\begin{minipage}{4cm}
\beginpicture
\setcoordinatesystem units <1.5mm,2mm>
\setplotarea x from 0 to 16, y from -2 to 15
\put{80)} [l] at 0 12
\put {$ \scriptstyle \bullet$} [c] at 6 0
\put {$ \scriptstyle \bullet$} [c] at 6 6
\put {$ \scriptstyle \bullet$} [c] at 6 12
\put {$ \scriptstyle \bullet$} [c] at 14 0
\put {$ \scriptstyle \bullet$} [c] at 14 6
\put {$ \scriptstyle \bullet$} [c] at 14 12
\setlinear \plot 6 12 6 0 14 12 14 0 6 6  /
\put{$720$} [c] at 10 -2
\endpicture
\end{minipage}
\begin{minipage}{4cm}
\beginpicture
\setcoordinatesystem units <1.5mm,2mm>
\setplotarea x from 0 to 16, y from -2 to 15
\put{81)} [l] at 0 12
\put {$ \scriptstyle \bullet$} [c] at 6 0
\put {$ \scriptstyle \bullet$} [c] at 6 6
\put {$ \scriptstyle \bullet$} [c] at 6 12
\put {$ \scriptstyle \bullet$} [c] at 14 0
\put {$ \scriptstyle \bullet$} [c] at 14 6
\put {$ \scriptstyle \bullet$} [c] at 14 12
\setlinear \plot 6 0 6 12 14 0 14 12 6 6  /
\put{$720$} [c] at 10 -2
\endpicture
\end{minipage}
\begin{minipage}{4cm}
\beginpicture
\setcoordinatesystem units <1.5mm,2mm>
\setplotarea x from 0 to 16, y from -2 to 15
\put{82)} [l] at 0 12
\put {$ \scriptstyle \bullet$} [c] at 6 0
\put {$ \scriptstyle \bullet$} [c] at 6 4
\put {$ \scriptstyle \bullet$} [c] at 6 8
\put {$ \scriptstyle \bullet$} [c] at 6 12
\put {$ \scriptstyle \bullet$} [c] at 14 0
\put {$ \scriptstyle \bullet$} [c] at 14 12
\setlinear \plot 6 0 6 12  /
\setlinear \plot 14  0 14  12 6 4  /
\put{$720$} [c] at 10 -2
\endpicture
\end{minipage}
\begin{minipage}{4cm}
\beginpicture
\setcoordinatesystem units <1.5mm,2mm>
\setplotarea x from 0 to 16, y from -2 to 15
\put{83)} [l] at 0 12
\put {$ \scriptstyle \bullet$} [c] at 6 0
\put {$ \scriptstyle \bullet$} [c] at 6 4
\put {$ \scriptstyle \bullet$} [c] at 6 8
\put {$ \scriptstyle \bullet$} [c] at 6 12
\put {$ \scriptstyle \bullet$} [c] at 14 0
\put {$ \scriptstyle \bullet$} [c] at 14 12
\setlinear \plot 6 0 6 12  /
\setlinear \plot 14  12 14  0 6 8  /
\put{$720$} [c] at 10 -2
\endpicture
\end{minipage}
\begin{minipage}{4cm}
\beginpicture
\setcoordinatesystem units <1.5mm,2mm>
\setplotarea x from 0 to 16, y from -2 to 15
\put{84)} [l] at 0 12
\put {$ \scriptstyle \bullet$} [c] at 6 0
\put {$ \scriptstyle \bullet$} [c] at 6 4
\put {$ \scriptstyle \bullet$} [c] at 6 8
\put {$ \scriptstyle \bullet$} [c] at 6 12
\put {$ \scriptstyle \bullet$} [c] at 14 0
\put {$ \scriptstyle \bullet$} [c] at 14 12
\setlinear \plot 6 12 6 0 14 12 14 0 6 12  /
\put{$720$} [c] at 10 -2
\endpicture
\end{minipage}
$$
$$
\begin{minipage}{4cm}
\beginpicture
\setcoordinatesystem units <1.5mm,2mm>
\setplotarea x from 0 to 16, y from -2 to 15
\put{85)} [l] at 0 12
\put {$ \scriptstyle \bullet$} [c] at 6 0
\put {$ \scriptstyle \bullet$} [c] at 6 6
\put {$ \scriptstyle \bullet$} [c] at 10 0
\put {$ \scriptstyle \bullet$} [c] at 10 12
\put {$ \scriptstyle \bullet$} [c] at 14 6
\put {$ \scriptstyle \bullet$} [c] at 14 12
\setlinear \plot  6 0 6 6 10 12 14 6 10 0 6 6  /
\setlinear \plot 14 6  14 12   /
\put{$720$} [c] at 10 -2
\endpicture
\end{minipage}
\begin{minipage}{4cm}
\beginpicture
\setcoordinatesystem units <1.5mm,2mm>
\setplotarea x from 0 to 16, y from -2 to 15
\put{86)} [l] at 0 12
\put {$ \scriptstyle \bullet$} [c] at 6 8
\put {$ \scriptstyle \bullet$} [c] at 8 0
\put {$ \scriptstyle \bullet$} [c] at 8 4
\put {$ \scriptstyle \bullet$} [c] at 8 12
\put {$ \scriptstyle \bullet$} [c] at 10 8
\put {$ \scriptstyle \bullet$} [c] at 14 12
\setlinear \plot 14  12 8 0 8 4 6 8 8 12 10 8 8 4 /
\put{$360$} [c] at 10 -2
\endpicture
\end{minipage}
\begin{minipage}{4cm}
\beginpicture
\setcoordinatesystem units <1.5mm,2mm>
\setplotarea x from 0 to 16, y from -2 to 15
\put{87)} [l] at 0 12
\put {$ \scriptstyle \bullet$} [c] at 6 4
\put {$ \scriptstyle \bullet$} [c] at 8 0
\put {$ \scriptstyle \bullet$} [c] at 8 8
\put {$ \scriptstyle \bullet$} [c] at 8 12
\put {$ \scriptstyle \bullet$} [c] at 10 4
\put {$ \scriptstyle \bullet$} [c] at 14 0
\setlinear \plot 14  0 8 12 8 8 6 4 8 0 10 4 8 8 /
\put{$360$} [c] at 10 -2
\endpicture
\end{minipage}
\begin{minipage}{4cm}
\beginpicture
\setcoordinatesystem units <1.5mm,2mm>
\setplotarea x from 0 to 16, y from -2 to 15
\put{88)} [l] at 0 12
\put {$ \scriptstyle \bullet$} [c] at 6 5
\put {$ \scriptstyle \bullet$} [c] at 7 0
\put {$ \scriptstyle \bullet$} [c] at 7 10
\put {$ \scriptstyle \bullet$} [c] at 7 12
\put {$ \scriptstyle \bullet$} [c] at 8 5
\put {$ \scriptstyle \bullet$} [c] at 14 12
\setlinear \plot 14  12 7 0 6 5 7 10 7 12 /
\setlinear \plot 7 10 8 5  7 0 /
\put{$360$} [c] at 10 -2
\endpicture
\end{minipage}
\begin{minipage}{4cm}
\beginpicture
\setcoordinatesystem units <1.5mm,2mm>
\setplotarea x from 0 to 16, y from -2 to 15
\put{89)} [l] at 0 12
\put {$ \scriptstyle \bullet$} [c] at 6 7
\put {$ \scriptstyle \bullet$} [c] at 7 0
\put {$ \scriptstyle \bullet$} [c] at 7 2
\put {$ \scriptstyle \bullet$} [c] at 7 12
\put {$ \scriptstyle \bullet$} [c] at 8 7
\put {$ \scriptstyle \bullet$} [c] at 14 0
\setlinear \plot 14  0 7 12 6 7 7 2 7 0 /
\setlinear \plot 7 2 8 7  7 12 /
\put{$360$} [c] at 10 -2
\endpicture
\end{minipage}
\begin{minipage}{4cm}
\beginpicture
\setcoordinatesystem units <1.5mm,2mm>
\setplotarea x from 0 to 16, y from -2 to 15
\put{90)} [l] at 0 12
\put {$ \scriptstyle \bullet$} [c] at 6 12
\put {$ \scriptstyle \bullet$} [c] at 10 12
\put {$ \scriptstyle \bullet$} [c] at 10 8
\put {$ \scriptstyle \bullet$} [c] at 10 4
\put {$ \scriptstyle \bullet$} [c] at 10 0
\put {$ \scriptstyle \bullet$} [c] at 14 12
\setlinear \plot 6 12 10 0 10 12 /
\setlinear \plot 10 8  14 12   /
\put{$360$} [c] at 10 -2
\endpicture
\end{minipage}
$$
$$
\begin{minipage}{4cm}
\beginpicture
\setcoordinatesystem units <1.5mm,2mm>
\setplotarea x from 0 to 16, y from -2 to 15
\put{91)} [l] at 0 12
\put {$ \scriptstyle \bullet$} [c] at 6 0
\put {$ \scriptstyle \bullet$} [c] at 10 12
\put {$ \scriptstyle \bullet$} [c] at 10 8
\put {$ \scriptstyle \bullet$} [c] at 10 4
\put {$ \scriptstyle \bullet$} [c] at 10 0
\put {$ \scriptstyle \bullet$} [c] at 14 0
\setlinear \plot 6 0 10 12  10 0 /
\setlinear \plot 10 4  14 0   /
\put{$360$} [c] at 10 -2
\endpicture
\end{minipage}
\begin{minipage}{4cm}
\beginpicture
\setcoordinatesystem units <1.5mm,2mm>
\setplotarea x from 0 to 16, y from -2 to 15
\put{92)} [l] at 0 12
\put {$ \scriptstyle \bullet$} [c] at 6 6
\put {$ \scriptstyle \bullet$} [c] at 6 12
\put {$ \scriptstyle \bullet$} [c] at 9 0
\put {$ \scriptstyle \bullet$} [c] at 12 6
\put {$ \scriptstyle \bullet$} [c] at 12 12
\put {$ \scriptstyle \bullet$} [c] at 14 12
\setlinear \plot 14 12 12 6 9 0 6 6 6 12 12 6 12 12 6 6 /
\put{$360$} [c] at 10 -2
\endpicture
\end{minipage}
\begin{minipage}{4cm}
\beginpicture
\setcoordinatesystem units <1.5mm,2mm>
\setplotarea x from 0 to 16, y from -2 to 15
\put{93)} [l] at 0 12
\put {$ \scriptstyle \bullet$} [c] at 6 6
\put {$ \scriptstyle \bullet$} [c] at 6 0
\put {$ \scriptstyle \bullet$} [c] at 9 12
\put {$ \scriptstyle \bullet$} [c] at 12 6
\put {$ \scriptstyle \bullet$} [c] at 12 0
\put {$ \scriptstyle \bullet$} [c] at 14 0
\setlinear \plot 14 0 12 6 9 12 6 6 6 0 12 6 12 0 6 6 /
\put{$360$} [c] at 10 -2
\endpicture
\end{minipage}
\begin{minipage}{4cm}
\beginpicture
\setcoordinatesystem units <1.5mm,2mm>
\setplotarea x from 0 to 16, y from -2 to 15
\put{94)} [l] at 0 12
\put {$ \scriptstyle \bullet$} [c] at 6 6
\put {$ \scriptstyle \bullet$} [c] at 6 12
\put {$ \scriptstyle \bullet$} [c] at 10 0
\put {$ \scriptstyle \bullet$} [c] at 10 6
\put {$ \scriptstyle \bullet$} [c] at 10 12
\put {$ \scriptstyle \bullet$} [c] at 14 6
\setlinear \plot 10 6 6 12 6 6 10 0 14 6 10 12 10 0 /
\setlinear \plot 6 6 10 12 /
\put{$360$} [c] at 10 -2
\endpicture
\end{minipage}
\begin{minipage}{4cm}
\beginpicture
\setcoordinatesystem units <1.5mm,2mm>
\setplotarea x from 0 to 16, y from -2 to 15
\put{95)} [l] at 0 12
\put {$ \scriptstyle \bullet$} [c] at 6 6
\put {$ \scriptstyle \bullet$} [c] at 6 0
\put {$ \scriptstyle \bullet$} [c] at 10 0
\put {$ \scriptstyle \bullet$} [c] at 10 6
\put {$ \scriptstyle \bullet$} [c] at 10 12
\put {$ \scriptstyle \bullet$} [c] at 14 6
\setlinear \plot 10 6 6 0 6 6 10 0 14 6 10 12 10 0 /
\setlinear \plot 6 6 10 12 /
\put{$360$} [c] at 10 -2
\endpicture
\end{minipage}
\begin{minipage}{4cm}
\beginpicture
\setcoordinatesystem units <1.5mm,2mm>
\setplotarea x from 0 to 16, y from -2 to 15
\put{96)} [l] at 0 12
\put {$ \scriptstyle \bullet$} [c] at 6 6
\put {$ \scriptstyle \bullet$} [c] at 6 12
\put {$ \scriptstyle \bullet$} [c] at 7 0
\put {$ \scriptstyle \bullet$} [c] at 8 6
\put {$ \scriptstyle \bullet$} [c] at 8 12
\put {$ \scriptstyle \bullet$} [c] at 14 0
\setlinear \plot 8 6 7 0 6 6 6 12 8 6 8 12 6  6  /
\setlinear \plot 6 12 14  0 8 12  /
\put{$180$} [c] at 10 -2
\endpicture
\end{minipage}
$$
$$
\begin{minipage}{4cm}
\beginpicture
\setcoordinatesystem units <1.5mm,2mm>
\setplotarea x from 0 to 16, y from -2 to 15
\put{97)} [l] at 0 12
\put {$ \scriptstyle \bullet$} [c] at 6 6
\put {$ \scriptstyle \bullet$} [c] at 6 0
\put {$ \scriptstyle \bullet$} [c] at 7 12
\put {$ \scriptstyle \bullet$} [c] at 8 6
\put {$ \scriptstyle \bullet$} [c] at 8 0
\put {$ \scriptstyle \bullet$} [c] at 14 12
\setlinear \plot 8 6 7 12 6 6 6 0 8 6 8 0 6  6  /
\setlinear \plot 6 0 14  12 8 0  /
\put{$180$} [c] at 10 -2
\endpicture
\end{minipage}
\begin{minipage}{4cm}
\beginpicture
\setcoordinatesystem units <1.5mm,2mm>
\setplotarea x from 0 to 16, y from -2 to 15
\put{98)} [l] at 0 12
\put {$ \scriptstyle \bullet$} [c] at 6 0
\put {$ \scriptstyle \bullet$} [c] at 6 6
\put {$ \scriptstyle \bullet$} [c] at 6 12
\put {$ \scriptstyle \bullet$} [c] at 10 12
\put {$ \scriptstyle \bullet$} [c] at 14 0
\put {$ \scriptstyle \bullet$} [c] at 14 12
\setlinear \plot  6 12 6 0 14 12 14  0 6 6 10 12  /
\put{$180$} [c] at 10 -2
\endpicture
\end{minipage}
\begin{minipage}{4cm}
\beginpicture
\setcoordinatesystem units <1.5mm,2mm>
\setplotarea x from 0 to 16, y from -2 to 15
\put{99)} [l] at 0 12
\put {$ \scriptstyle \bullet$} [c] at 6 0
\put {$ \scriptstyle \bullet$} [c] at 6 6
\put {$ \scriptstyle \bullet$} [c] at 6 12
\put {$ \scriptstyle \bullet$} [c] at 10 0
\put {$ \scriptstyle \bullet$} [c] at 14 0
\put {$ \scriptstyle \bullet$} [c] at 14 12
\setlinear \plot  6 0 6 12 14 0 14 12 6 6 10 0  /
\put{$180$} [c] at 10 -2
\endpicture
\end{minipage}
\begin{minipage}{4cm}
\beginpicture
\setcoordinatesystem units <1.5mm,2mm>
\setplotarea x from 0 to 16, y from -2 to 15
\put{100)} [l] at 0 12
\put {$ \scriptstyle \bullet$} [c] at 6 0
\put {$ \scriptstyle \bullet$} [c] at 6 12
\put {$ \scriptstyle \bullet$} [c] at 10 12
\put {$ \scriptstyle \bullet$} [c] at 14 0
\put {$ \scriptstyle \bullet$} [c] at 14 6
\put {$ \scriptstyle \bullet$} [c] at 14 12
\setlinear \plot 14 0 14 12  6 0  6 12 14 6 10 12 6 0   /
\put{$120$} [c] at 10 -2
\endpicture
\end{minipage}
\begin{minipage}{4cm}
\beginpicture
\setcoordinatesystem units <1.5mm,2mm>
\setplotarea x from 0 to 16, y from -2 to 15
\put{101)} [l] at 0 12
\put {$ \scriptstyle \bullet$} [c] at 6 0
\put {$ \scriptstyle \bullet$} [c] at 6 12
\put {$ \scriptstyle \bullet$} [c] at 10 0
\put {$ \scriptstyle \bullet$} [c] at 14 0
\put {$ \scriptstyle \bullet$} [c] at 14 6
\put {$ \scriptstyle \bullet$} [c] at 14 12
\setlinear \plot 14 12 14 0  6 12  6 0 14 6 10 0 6 12   /
\put{$120$} [c] at 10 -2
\endpicture
\end{minipage}
\begin{minipage}{4cm}
\beginpicture
\setcoordinatesystem units <1.5mm,2mm>
\setplotarea x from 0 to 16, y from -2 to 15
\put{${\bf 14}$} [l] at 0 15
\put{102)} [l] at 0 12
\put {$ \scriptstyle \bullet$} [c] at 6 6
\put {$ \scriptstyle \bullet$} [c] at 6 0
\put {$ \scriptstyle \bullet$} [c] at 6 12
\put {$ \scriptstyle \bullet$} [c] at 14 6
\put {$ \scriptstyle \bullet$} [c] at 14 12
\put {$ \scriptstyle \bullet$} [c] at 14 0
\setlinear \plot 6 12  6 0     /
\setlinear \plot 14 12  14 0     /
\setlinear \plot 6 0  14 6     /
\put{$720$} [c] at 10 -2
\endpicture
\end{minipage}
$$
$$
\begin{minipage}{4cm}
\beginpicture
\setcoordinatesystem units <1.5mm,2mm>
\setplotarea x from 0 to 16, y from -2 to 15
\put{103)} [l] at 0 12
\put {$ \scriptstyle \bullet$} [c] at 6 6
\put {$ \scriptstyle \bullet$} [c] at 6 0
\put {$ \scriptstyle \bullet$} [c] at 6 12
\put {$ \scriptstyle \bullet$} [c] at 14 6
\put {$ \scriptstyle \bullet$} [c] at 14 12
\put {$ \scriptstyle \bullet$} [c] at 14 0
\setlinear \plot 6 12  6 0     /
\setlinear \plot 14 12  14 0     /
\setlinear \plot 6 12  14 6     /
\put{$720$} [c] at 10 -2
\endpicture
\end{minipage}
\begin{minipage}{4cm}
\beginpicture
\setcoordinatesystem units <1.5mm,2mm>
\setplotarea x from 0 to 16, y from -2 to 15
\put{104)} [l] at 0 12
\put {$ \scriptstyle \bullet$} [c] at 6 0
\put {$ \scriptstyle \bullet$} [c] at 6 4
\put {$ \scriptstyle \bullet$} [c] at 6 8
\put {$ \scriptstyle \bullet$} [c] at 6 12
\put {$ \scriptstyle \bullet$} [c] at 14 12
\put {$ \scriptstyle \bullet$} [c] at 14 0
\setlinear \plot 6 12  6 0  14 12 14 0    /
\put{$720$} [c] at 10 -2
\endpicture
\end{minipage}
\begin{minipage}{4cm}
\beginpicture
\setcoordinatesystem units <1.5mm,2mm>
\setplotarea x from 0 to 16, y from -2 to 15
\put{105)} [l] at 0 12
\put {$ \scriptstyle \bullet$} [c] at 6 0
\put {$ \scriptstyle \bullet$} [c] at 6 4
\put {$ \scriptstyle \bullet$} [c] at 6 8
\put {$ \scriptstyle \bullet$} [c] at 6 12
\put {$ \scriptstyle \bullet$} [c] at 14 12
\put {$ \scriptstyle \bullet$} [c] at 14 0
\setlinear \plot 6 0  6 12  14 0 14 12    /
\put{$720$} [c] at 10 -2
\endpicture
\end{minipage}
\begin{minipage}{4cm}
\beginpicture
\setcoordinatesystem units <1.5mm,2mm>
\setplotarea x from 0 to 16, y from -2 to 15
\put{106)} [l] at 0 12
\put {$ \scriptstyle \bullet$} [c] at 6 6
\put {$ \scriptstyle \bullet$} [c] at 6 0
\put {$ \scriptstyle \bullet$} [c] at 6 12
\put {$ \scriptstyle \bullet$} [c] at 10 9
\put {$ \scriptstyle \bullet$} [c] at 14 12
\put {$ \scriptstyle \bullet$} [c] at 14 0
\setlinear \plot 6 12  6 0     /
\setlinear \plot 6 6  14 12 14 0     /
\put{$720$} [c] at 10 -2
\endpicture
\end{minipage}
\begin{minipage}{4cm}
\beginpicture
\setcoordinatesystem units <1.5mm,2mm>
\setplotarea x from 0 to 16, y from -2 to 15
\put{107)} [l] at 0 12
\put {$ \scriptstyle \bullet$} [c] at 6 6
\put {$ \scriptstyle \bullet$} [c] at 6 0
\put {$ \scriptstyle \bullet$} [c] at 6 12
\put {$ \scriptstyle \bullet$} [c] at 10 3
\put {$ \scriptstyle \bullet$} [c] at 14 12
\put {$ \scriptstyle \bullet$} [c] at 14 0
\setlinear \plot 6 12  6 0     /
\setlinear \plot 6 6  14 0 14 12     /
\put{$720$} [c] at 10 -2
\endpicture
\end{minipage}
\begin{minipage}{4cm}
\beginpicture
\setcoordinatesystem units <1.5mm,2mm>
\setplotarea x from 0 to 16, y from -2 to 15
\put{108)} [l] at 0 12
\put {$ \scriptstyle \bullet$} [c] at 6 6
\put {$ \scriptstyle \bullet$} [c] at 8 0
\put {$ \scriptstyle \bullet$} [c] at 8 12
\put {$ \scriptstyle \bullet$} [c] at 10 6
\put {$ \scriptstyle \bullet$} [c] at 14 12
\put {$ \scriptstyle \bullet$} [c] at 14 0
\setlinear \plot 8 12  14 0 14 12 10 6 8 0 6 6 8 12 10 6     /
\put{$720$} [c] at 10 -2
\endpicture
\end{minipage}
$$
$$
\begin{minipage}{4cm}
\beginpicture
\setcoordinatesystem units <1.5mm,2mm>
\setplotarea x from 0 to 16, y from -2 to 15
\put{109)} [l] at 0 12
\put {$ \scriptstyle \bullet$} [c] at 6 6
\put {$ \scriptstyle \bullet$} [c] at 8 0
\put {$ \scriptstyle \bullet$} [c] at 8 12
\put {$ \scriptstyle \bullet$} [c] at 10 6
\put {$ \scriptstyle \bullet$} [c] at 14 12
\put {$ \scriptstyle \bullet$} [c] at 14 0
\setlinear \plot 8 0  14 12 14 0 10 6 8 12 6 6 8 0 10 6     /
\put{$720$} [c] at 10 -2
\endpicture
\end{minipage}
\begin{minipage}{4cm}
\beginpicture
\setcoordinatesystem units <1.5mm,2mm>
\setplotarea x from 0 to 16, y from -2 to 15
\put{110)} [l] at 0 12
\put {$ \scriptstyle \bullet$} [c] at 6 4
\put {$ \scriptstyle \bullet$} [c] at 6 12
\put {$ \scriptstyle \bullet$} [c] at 10 0
\put {$ \scriptstyle \bullet$} [c] at 10 8
\put {$ \scriptstyle \bullet$} [c] at 14 4
\put {$ \scriptstyle \bullet$} [c] at 14 12
\setlinear \plot 6 12 10 8 14 12 14 4 10 0 6 4 6 12  /
\setlinear \plot 10 0 10 8  /
\put{$360$} [c] at 10 -2
\endpicture
\end{minipage}
\begin{minipage}{4cm}
\beginpicture
\setcoordinatesystem units <1.5mm,2mm>
\setplotarea x from 0 to 16, y from -2 to 15
\put{111)} [l] at 0 12
\put {$ \scriptstyle \bullet$} [c] at 6 0
\put {$ \scriptstyle \bullet$} [c] at 6 8
\put {$ \scriptstyle \bullet$} [c] at 10 4
\put {$ \scriptstyle \bullet$} [c] at 10 12
\put {$ \scriptstyle \bullet$} [c] at 14 8
\put {$ \scriptstyle \bullet$} [c] at 14 0
\setlinear \plot 6 0 10 4 14 0 14 8 10 12 6 8 6 0  /
\setlinear \plot 10 4 10 12  /
\put{$360$} [c] at 10 -2
\endpicture
\end{minipage}
\begin{minipage}{4cm}
\beginpicture
\setcoordinatesystem units <1.5mm,2mm>
\setplotarea x from 0 to 16, y from -2 to 15
\put{112)} [l] at 0 12
\put {$ \scriptstyle \bullet$} [c] at 6 12
\put {$ \scriptstyle \bullet$} [c] at 10 0
\put {$ \scriptstyle \bullet$} [c] at 10 4
\put {$ \scriptstyle \bullet$} [c] at 10 8
\put {$ \scriptstyle \bullet$} [c] at 10 12
\put {$ \scriptstyle \bullet$} [c] at 14 12
\setlinear \plot 6 12 10 4 14 12   /
\setlinear \plot 10 0 10 12  /
\put{$360$} [c] at 10 -2
\endpicture
\end{minipage}
\begin{minipage}{4cm}
\beginpicture
\setcoordinatesystem units <1.5mm,2mm>
\setplotarea x from 0 to 16, y from -2 to 15
\put{113)} [l] at 0 12
\put {$ \scriptstyle \bullet$} [c] at 6 0
\put {$ \scriptstyle \bullet$} [c] at 10 0
\put {$ \scriptstyle \bullet$} [c] at 10 4
\put {$ \scriptstyle \bullet$} [c] at 10 8
\put {$ \scriptstyle \bullet$} [c] at 10 12
\put {$ \scriptstyle \bullet$} [c] at 14 0
\setlinear \plot 6 0 10 8 14 0   /
\setlinear \plot 10 0 10 12  /
\put{$360$} [c] at 10 -2
\endpicture
\end{minipage}
\begin{minipage}{4cm}
\beginpicture
\setcoordinatesystem units <1.5mm,2mm>
\setplotarea x from 0 to 16, y from -2 to 15
\put{114)} [l] at 0 12
\put {$ \scriptstyle \bullet$} [c] at 6 6
\put {$ \scriptstyle \bullet$} [c] at 6 12
\put {$ \scriptstyle \bullet$} [c] at 10 0
\put {$ \scriptstyle \bullet$} [c] at 10 12
\put {$ \scriptstyle \bullet$} [c] at 14 6
\put {$ \scriptstyle \bullet$} [c] at 14 12
\setlinear \plot 6 12 6 6 10 0 14 6 14 12    /
\setlinear \plot 6 6 10 12 14 6    /
\put{$360$} [c] at 10 -2
\endpicture
\end{minipage}
$$
$$
\begin{minipage}{4cm}
\beginpicture
\setcoordinatesystem units <1.5mm,2mm>
\setplotarea x from 0 to 16, y from -2 to 15
\put{115)} [l] at 0 12
\put {$ \scriptstyle \bullet$} [c] at 6 6
\put {$ \scriptstyle \bullet$} [c] at 6 0
\put {$ \scriptstyle \bullet$} [c] at 10 0
\put {$ \scriptstyle \bullet$} [c] at 10 12
\put {$ \scriptstyle \bullet$} [c] at 14 6
\put {$ \scriptstyle \bullet$} [c] at 14 0
\setlinear \plot 6 0 6 6 10 0 14 6 14 0    /
\setlinear \plot 6 6 10 12 14 6    /
\put{$360$} [c] at 10 -2
\endpicture
\end{minipage}
\begin{minipage}{4cm}
\beginpicture
\setcoordinatesystem units <1.5mm,2mm>
\setplotarea x from 0 to 16, y from -2 to 15
\put{116)} [l] at 0 12
\put {$ \scriptstyle \bullet$} [c] at 6 6
\put {$ \scriptstyle \bullet$} [c] at 6 12
\put {$ \scriptstyle \bullet$} [c] at 8 0
\put {$ \scriptstyle \bullet$} [c] at 10 12
\put {$ \scriptstyle \bullet$} [c] at 10 6
\put {$ \scriptstyle \bullet$} [c] at 14 0
\setlinear \plot 14 0 10 12 10 6 6 12  6 6 8 0 10 6    /
\setlinear \plot 6 6  10 12    /
\put{$360$} [c] at 10 -2
\endpicture
\end{minipage}
\begin{minipage}{4cm}
\beginpicture
\setcoordinatesystem units <1.5mm,2mm>
\setplotarea x from 0 to 16, y from -2 to 15
\put{117)} [l] at 0 12
\put {$ \scriptstyle \bullet$} [c] at 6 6
\put {$ \scriptstyle \bullet$} [c] at 6 0
\put {$ \scriptstyle \bullet$} [c] at 8 12
\put {$ \scriptstyle \bullet$} [c] at 10 0
\put {$ \scriptstyle \bullet$} [c] at 10 6
\put {$ \scriptstyle \bullet$} [c] at 14 12
\setlinear \plot 14 12 10 0 10 6 6 0  6 6  10 0    /
\setlinear \plot 6 6  8 12 10 6    /
\put{$360$} [c] at 10 -2
\endpicture
\end{minipage}
\begin{minipage}{4cm}
\beginpicture
\setcoordinatesystem units <1.5mm,2mm>
\setplotarea x from 0 to 16, y from -2 to 15
\put{118)} [l] at 0 12
\put {$ \scriptstyle \bullet$} [c] at 6 12
\put {$ \scriptstyle \bullet$} [c] at 8 0
\put {$ \scriptstyle \bullet$} [c] at 8 12
\put {$ \scriptstyle \bullet$} [c] at 11 6
\put {$ \scriptstyle \bullet$} [c] at 14 12
\put {$ \scriptstyle \bullet$} [c] at 14 0
\setlinear \plot 6 12  8 0 14 12    /
\setlinear \plot 8 12  14 0     /
\put{$360$} [c] at 10 -2
\endpicture
\end{minipage}
\begin{minipage}{4cm}
\beginpicture
\setcoordinatesystem units <1.5mm,2mm>
\setplotarea x from 0 to 16, y from -2 to 15
\put{119)} [l] at 0 12
\put {$ \scriptstyle \bullet$} [c] at 6 0
\put {$ \scriptstyle \bullet$} [c] at 8 0
\put {$ \scriptstyle \bullet$} [c] at 8 12
\put {$ \scriptstyle \bullet$} [c] at 11 6
\put {$ \scriptstyle \bullet$} [c] at 14 12
\put {$ \scriptstyle \bullet$} [c] at 14 0
\setlinear \plot 6 0  8 12 14 0    /
\setlinear \plot 8 0  14 12     /
\put{$360$} [c] at 10 -2
\endpicture
\end{minipage}
\begin{minipage}{4cm}
\beginpicture
\setcoordinatesystem units <1.5mm,2mm>
\setplotarea x from 0 to 16, y from -2 to 15
\put{120)} [l] at 0 12
\put {$ \scriptstyle \bullet$} [c] at 6 0
\put {$ \scriptstyle \bullet$} [c] at 6 12
\put {$ \scriptstyle \bullet$} [c] at 10 12
\put {$ \scriptstyle \bullet$} [c] at 14 6
\put {$ \scriptstyle \bullet$} [c] at 14 12
\put {$ \scriptstyle \bullet$} [c] at 14 0
\setlinear \plot 6 0 6 12  14 0 14 12 6 0 14 12  /
\setlinear \plot 6 0 10  12  14  6    /
\put{$360$} [c] at 10 -2
\endpicture
\end{minipage}
$$
$$
\begin{minipage}{4cm}
\beginpicture
\setcoordinatesystem units <1.5mm,2mm>
\setplotarea x from 0 to 16, y from -2 to 15
\put{121)} [l] at 0 12
\put {$ \scriptstyle \bullet$} [c] at 6 0
\put {$ \scriptstyle \bullet$} [c] at 6 12
\put {$ \scriptstyle \bullet$} [c] at 10 0
\put {$ \scriptstyle \bullet$} [c] at 14 6
\put {$ \scriptstyle \bullet$} [c] at 14 12
\put {$ \scriptstyle \bullet$} [c] at 14 0
\setlinear \plot 6 12 6 0  14 12 14 0 6 12 14 0  /
\setlinear \plot 6 12 10  0  14  6    /
\put{$360$} [c] at 10 -2
\endpicture
\end{minipage}
\begin{minipage}{4cm}
\beginpicture
\setcoordinatesystem units <1.5mm,2mm>
\setplotarea x from 0 to 16, y from -2 to 15
\put{122)} [l] at 0 12
\put {$ \scriptstyle \bullet$} [c] at 6 0
\put {$ \scriptstyle \bullet$} [c] at 6 12
\put {$ \scriptstyle \bullet$} [c] at 10 12
\put {$ \scriptstyle \bullet$} [c] at 14 6
\put {$ \scriptstyle \bullet$} [c] at 14 12
\put {$ \scriptstyle \bullet$} [c] at 14 0
\setlinear \plot 6 0 6 12  14 6 10 12 6 0  /
\setlinear \plot 14 0 14  12     /
\put{$360$} [c] at 10 -2
\endpicture
\end{minipage}
\begin{minipage}{4cm}
\beginpicture
\setcoordinatesystem units <1.5mm,2mm>
\setplotarea x from 0 to 16, y from -2 to 15
\put{123)} [l] at 0 12
\put {$ \scriptstyle \bullet$} [c] at 6 0
\put {$ \scriptstyle \bullet$} [c] at 6 12
\put {$ \scriptstyle \bullet$} [c] at 10 0
\put {$ \scriptstyle \bullet$} [c] at 14 6
\put {$ \scriptstyle \bullet$} [c] at 14 12
\put {$ \scriptstyle \bullet$} [c] at 14 0
\setlinear \plot 6 12 6 0  14 6 14 12  /
\setlinear \plot 6 12 10 0 14 6 14 0     /
\put{$360$} [c] at 10 -2
\endpicture
\end{minipage}
\begin{minipage}{4cm}
\beginpicture
\setcoordinatesystem units <1.5mm,2mm>
\setplotarea x from 0 to 16, y from -2 to 15
\put{124)} [l] at 0 12
\put {$ \scriptstyle \bullet$} [c] at 6 4
\put {$ \scriptstyle \bullet$} [c] at 6 8
\put {$ \scriptstyle \bullet$} [c] at 10 0
\put {$ \scriptstyle \bullet$} [c] at 10 6
\put {$ \scriptstyle \bullet$} [c] at 10 12
\put {$ \scriptstyle \bullet$} [c] at 14 6
\setlinear \plot 10 0 6 4 6 8 10 12 14 6 10 0   /
\setlinear \plot 10 0 10 12   /
\put{$360$} [c] at 10 -2
\endpicture
\end{minipage}
\begin{minipage}{4cm}
\beginpicture
\setcoordinatesystem units <1.5mm,2mm>
\setplotarea x from 0 to 16, y from -2 to 15
\put{125)} [l] at 0 12
\put {$ \scriptstyle \bullet$} [c] at 6 6
\put {$ \scriptstyle \bullet$} [c] at 6 0
\put {$ \scriptstyle \bullet$} [c] at 6 12
\put {$ \scriptstyle \bullet$} [c] at 14 12
\put {$ \scriptstyle \bullet$} [c] at 14 6
\put {$ \scriptstyle \bullet$} [c] at 14 0
\setlinear \plot 6 12  6 0 14 12 14 0 6 12    /
\put{$360$} [c] at 10 -2
\endpicture
\end{minipage}
\begin{minipage}{4cm}
\beginpicture
\setcoordinatesystem units <1.5mm,2mm>
\setplotarea x from 0 to 16, y from -2 to 15
\put{126)} [l] at 0 12
\put {$ \scriptstyle \bullet$} [c] at 6 8
\put {$ \scriptstyle \bullet$} [c] at 9 0
\put {$ \scriptstyle \bullet$} [c] at 9 4
\put {$ \scriptstyle \bullet$} [c] at 9 12
\put {$ \scriptstyle \bullet$} [c] at 12 8
\put {$ \scriptstyle \bullet$} [c] at 14 0
\setlinear \plot 9 0 9 4 6 8 9 12 12 8 9 4   /
\put{$360$} [c] at 10 -2
\endpicture
\end{minipage}
$$
$$
\begin{minipage}{4cm}
\beginpicture
\setcoordinatesystem units <1.5mm,2mm>
\setplotarea x from 0 to 16, y from -2 to 15
\put{127)} [l] at 0 12
\put {$ \scriptstyle \bullet$} [c] at 6 4
\put {$ \scriptstyle \bullet$} [c] at 9 0
\put {$ \scriptstyle \bullet$} [c] at 9 8
\put {$ \scriptstyle \bullet$} [c] at 9 12
\put {$ \scriptstyle \bullet$} [c] at 12 4
\put {$ \scriptstyle \bullet$} [c] at 14 0
\setlinear \plot 9 12 9 8 6 4 9 0 12 4 9 8   /
\put{$360$} [c] at 10 -2
\endpicture
\end{minipage}
\begin{minipage}{4cm}
\beginpicture
\setcoordinatesystem units <1.5mm,2mm>
\setplotarea x from 0 to 16, y from -2 to 15
\put{128)} [l] at 0 12
\put {$ \scriptstyle \bullet$} [c] at 6 12
\put {$ \scriptstyle \bullet$} [c] at 9 0
\put {$ \scriptstyle \bullet$} [c] at 9 4
\put {$ \scriptstyle \bullet$} [c] at 9 8
\put {$ \scriptstyle \bullet$} [c] at 12 12
\put {$ \scriptstyle \bullet$} [c] at 14 0
\setlinear \plot 9 0 9 8 6  12    /
\setlinear \plot 9 8 12  12    /
\put{$360$} [c] at 10 -2
\endpicture
\end{minipage}
\begin{minipage}{4cm}
\beginpicture
\setcoordinatesystem units <1.5mm,2mm>
\setplotarea x from 0 to 16, y from -2 to 15
\put{129)} [l] at 0 12
\put {$ \scriptstyle \bullet$} [c] at 6 0
\put {$ \scriptstyle \bullet$} [c] at 9 12
\put {$ \scriptstyle \bullet$} [c] at 9 4
\put {$ \scriptstyle \bullet$} [c] at 9 8
\put {$ \scriptstyle \bullet$} [c] at 12 0
\put {$ \scriptstyle \bullet$} [c] at 14 0
\setlinear \plot 9 12 9 4 6  0    /
\setlinear \plot 9 4 12  0    /
\put{$360$} [c] at 10 -2
\endpicture
\end{minipage}
\begin{minipage}{4cm}
\beginpicture
\setcoordinatesystem units <1.5mm,2mm>
\setplotarea x from 0 to 16, y from -2 to 15
\put{${\bf  15}$} [l] at 0 15
\put{130)} [l] at 0 12
\put {$ \scriptstyle \bullet$} [c] at 6 6
\put {$ \scriptstyle \bullet$} [c] at 8 12
\put {$ \scriptstyle \bullet$} [c] at 8 0
\put {$ \scriptstyle \bullet$} [c] at 10 6
\put {$ \scriptstyle \bullet$} [c] at 7 3
\put {$ \scriptstyle \bullet$} [c] at 14 12
\setlinear \plot 8 0 6 6  8 12 10 6 8 0 14 12 /
\put{$720$} [c] at 10 -2
\endpicture
\end{minipage}
\begin{minipage}{4cm}
\beginpicture
\setcoordinatesystem units <1.5mm,2mm>
\setplotarea x from 0 to 16, y from -2 to 15
\put{131)} [l] at 0 12
\put {$ \scriptstyle \bullet$} [c] at 6 6
\put {$ \scriptstyle \bullet$} [c] at 8 12
\put {$ \scriptstyle \bullet$} [c] at 8 0
\put {$ \scriptstyle \bullet$} [c] at 10 6
\put {$ \scriptstyle \bullet$} [c] at 7 9
\put {$ \scriptstyle \bullet$} [c] at 14 0
\setlinear \plot 8 12 6 6  8 0 10 6 8 12 14 0 /
\put{$720$} [c] at 10 -2
\endpicture
\end{minipage}
\begin{minipage}{4cm}
\beginpicture
\setcoordinatesystem units <1.5mm,2mm>
\setplotarea x from 0 to 16, y from -2 to 15
\put{132)} [l] at 0 12
\put {$ \scriptstyle \bullet$} [c] at 6 12
\put {$ \scriptstyle \bullet$} [c] at 10  0
\put {$ \scriptstyle \bullet$} [c] at 10 4
\put {$ \scriptstyle \bullet$} [c] at 10 8
\put {$ \scriptstyle \bullet$} [c] at 10 12
\put {$ \scriptstyle \bullet$} [c] at 14 12
\setlinear \plot 6 12  10 0 10 12  /
\setlinear \plot 10 4  14  12     /
\put{$720$} [c] at 10 -2
\endpicture
\end{minipage}
$$
$$
\begin{minipage}{4cm}
\beginpicture
\setcoordinatesystem units <1.5mm,2mm>
\setplotarea x from 0 to 16, y from -2 to 15
\put{133)} [l] at 0 12
\put {$ \scriptstyle \bullet$} [c] at 6 0
\put {$ \scriptstyle \bullet$} [c] at 10  0
\put {$ \scriptstyle \bullet$} [c] at 10 4
\put {$ \scriptstyle \bullet$} [c] at 10 8
\put {$ \scriptstyle \bullet$} [c] at 10 12
\put {$ \scriptstyle \bullet$} [c] at 14 0
\setlinear \plot 6 0  10 12 10 0  /
\setlinear \plot 10 8  14  0     /
\put{$720$} [c] at 10 -2
\endpicture
\end{minipage}
\begin{minipage}{4cm}
\beginpicture
\setcoordinatesystem units <1.5mm,2mm>
\setplotarea x from 0 to 16, y from -2 to 15
\put{134)} [l] at 0 12
\put {$ \scriptstyle \bullet$} [c] at 6 6
\put {$ \scriptstyle \bullet$} [c] at 8 0
\put {$ \scriptstyle \bullet$} [c] at 8 12
\put {$ \scriptstyle \bullet$} [c] at 10 6
\put {$ \scriptstyle \bullet$} [c] at 14 12
\put {$ \scriptstyle \bullet$} [c] at 14 0
\setlinear \plot 10 6 8 0  6 6 8 12 10 6 14 12 14 0  /
\put{$720$} [c] at 10 -2
\endpicture
\end{minipage}
\begin{minipage}{4cm}
\beginpicture
\setcoordinatesystem units <1.5mm,2mm>
\setplotarea x from 0 to 16, y from -2 to 15
\put{135)} [l] at 0 12
\put {$ \scriptstyle \bullet$} [c] at 6 6
\put {$ \scriptstyle \bullet$} [c] at 8 0
\put {$ \scriptstyle \bullet$} [c] at 8 12
\put {$ \scriptstyle \bullet$} [c] at 10 6
\put {$ \scriptstyle \bullet$} [c] at 14 12
\put {$ \scriptstyle \bullet$} [c] at 14 0
\setlinear \plot 10 6 8 0  6 6 8 12 10 6 14 0 14 12  /
\put{$720$} [c] at 10 -2
\endpicture
\end{minipage}
\begin{minipage}{4cm}
\beginpicture
\setcoordinatesystem units <1.5mm,2mm>
\setplotarea x from 0 to 16, y from -2 to 15
\put{136)} [l] at 0 12
\put {$ \scriptstyle \bullet$} [c] at 6 0
\put {$ \scriptstyle \bullet$} [c] at 6 12
\put {$ \scriptstyle \bullet$} [c] at 10 0
\put {$ \scriptstyle \bullet$} [c] at 14 0
\put {$ \scriptstyle \bullet$} [c] at 14 6
\put {$ \scriptstyle \bullet$} [c] at 14 12
\setlinear \plot 6 12 6  0 14 12 14 0   /
\setlinear \plot 6 12 10 0 14 6     /
\put{$720$} [c] at 10 -2
\endpicture
\end{minipage}
\begin{minipage}{4cm}
\beginpicture
\setcoordinatesystem units <1.5mm,2mm>
\setplotarea x from 0 to 16, y from -2 to 15
\put{137)} [l] at 0 12
\put {$ \scriptstyle \bullet$} [c] at 6 0
\put {$ \scriptstyle \bullet$} [c] at 6 12
\put {$ \scriptstyle \bullet$} [c] at 10 12
\put {$ \scriptstyle \bullet$} [c] at 14 0
\put {$ \scriptstyle \bullet$} [c] at 14 6
\put {$ \scriptstyle \bullet$} [c] at 14 12
\setlinear \plot 6 0 6  12 14 0 14 12   /
\setlinear \plot 6 0 10 12 14 6     /
\put{$720$} [c] at 10 -2
\endpicture
\end{minipage}
\begin{minipage}{4cm}
\beginpicture
\setcoordinatesystem units <1.5mm,2mm>
\setplotarea x from 0 to 16, y from -2 to 15
\put{138)} [l] at 0 12
\put {$ \scriptstyle \bullet$} [c] at 6 0
\put {$ \scriptstyle \bullet$} [c] at 6 6
\put {$ \scriptstyle \bullet$} [c] at 6 12
\put {$ \scriptstyle \bullet$} [c] at 14 0
\put {$ \scriptstyle \bullet$} [c] at 14 6
\put {$ \scriptstyle \bullet$} [c] at 14 12
\setlinear \plot 6 12 6 0 14 12 14 0   /
\put{$720$} [c] at 10 -2
\endpicture
\end{minipage}
$$
$$
\begin{minipage}{4cm}
\beginpicture
\setcoordinatesystem units <1.5mm,2mm>
\setplotarea x from 0 to 16, y from -2 to 15
\put{139)} [l] at 0 12
\put {$ \scriptstyle \bullet$} [c] at 6 6
\put {$ \scriptstyle \bullet$} [c] at 10 0
\put {$ \scriptstyle \bullet$} [c] at 10 6
\put {$ \scriptstyle \bullet$} [c] at 10 12
\put {$ \scriptstyle \bullet$} [c] at 14 6
\put {$ \scriptstyle \bullet$} [c] at 14 12
\setlinear \plot 14 12 14 6 10 12 6 6 10  0 14 6    /
\setlinear \plot 10 12 10  0     /
\put{$360$} [c] at 10 -2
\endpicture
\end{minipage}
\begin{minipage}{4cm}
\beginpicture
\setcoordinatesystem units <1.5mm,2mm>
\setplotarea x from 0 to 16, y from -2 to 15
\put{140)} [l] at 0 12
\put {$ \scriptstyle \bullet$} [c] at 6 6
\put {$ \scriptstyle \bullet$} [c] at 10 0
\put {$ \scriptstyle \bullet$} [c] at 10 6
\put {$ \scriptstyle \bullet$} [c] at 10 12
\put {$ \scriptstyle \bullet$} [c] at 14 6
\put {$ \scriptstyle \bullet$} [c] at 14 0
\setlinear \plot 14 0 14 6 10 12 6 6 10  0 14 6    /
\setlinear \plot 10 12 10  0     /
\put{$360$} [c] at 10 -2
\endpicture
\end{minipage}
\begin{minipage}{4cm}
\beginpicture
\setcoordinatesystem units <1.5mm,2mm>
\setplotarea x from 0 to 16, y from -2 to 15
\put{141)} [l] at 0 12
\put {$ \scriptstyle \bullet$} [c] at 6 6
\put {$ \scriptstyle \bullet$} [c] at 9 0
\put {$ \scriptstyle \bullet$} [c] at 9 12
\put {$ \scriptstyle \bullet$} [c] at 12 6
\put {$ \scriptstyle \bullet$} [c] at 12 12
\put {$ \scriptstyle \bullet$} [c] at 14 12
\setlinear \plot 12 12 12 6 9 12  6 6 9 0 12 6 14 12    /
\put{$360$} [c] at 10 -2
\endpicture
\end{minipage}
\begin{minipage}{4cm}
\beginpicture
\setcoordinatesystem units <1.5mm,2mm>
\setplotarea x from 0 to 16, y from -2 to 15
\put{142)} [l] at 0 12
\put {$ \scriptstyle \bullet$} [c] at 6 6
\put {$ \scriptstyle \bullet$} [c] at 9 0
\put {$ \scriptstyle \bullet$} [c] at 9 12
\put {$ \scriptstyle \bullet$} [c] at 12 6
\put {$ \scriptstyle \bullet$} [c] at 12 0
\put {$ \scriptstyle \bullet$} [c] at 14 0
\setlinear \plot 12 0 12 6 9 12  6 6 9 0 12 6 14 0    /
\put{$360$} [c] at 10 -2
\endpicture
\end{minipage}
\begin{minipage}{4cm}
\beginpicture
\setcoordinatesystem units <1.5mm,2mm>
\setplotarea x from 0 to 16, y from -2 to 15
\put{143)} [l] at 0 12
\put {$ \scriptstyle \bullet$} [c] at 6 0
\put {$ \scriptstyle \bullet$} [c] at 6 6
\put {$ \scriptstyle \bullet$} [c] at 6 12
\put {$ \scriptstyle \bullet$} [c] at 8 3
\put {$ \scriptstyle \bullet$} [c] at 14 0
\put {$ \scriptstyle \bullet$} [c] at 14 12
\setlinear \plot 6 12 6 0 14 12 14 0 6 12  /
\put{$360$} [c] at 10 -2
\endpicture
\end{minipage}
\begin{minipage}{4cm}
\beginpicture
\setcoordinatesystem units <1.5mm,2mm>
\setplotarea x from 0 to 16, y from -2 to 15
\put{144)} [l] at 0 12
\put {$ \scriptstyle \bullet$} [c] at 6 0
\put {$ \scriptstyle \bullet$} [c] at 6 6
\put {$ \scriptstyle \bullet$} [c] at 6 12
\put {$ \scriptstyle \bullet$} [c] at 8 9
\put {$ \scriptstyle \bullet$} [c] at 14 0
\put {$ \scriptstyle \bullet$} [c] at 14 12
\setlinear \plot 6 12 6 0 14 12 14 0 6 12  /
\put{$360$} [c] at 10 -2
\endpicture
\end{minipage}
$$
$$
\begin{minipage}{4cm}
\beginpicture
\setcoordinatesystem units <1.5mm,2mm>
\setplotarea x from 0 to 16, y from -2 to 15
\put{145)} [l] at 0 12
\put {$ \scriptstyle \bullet$} [c] at 6 0
\put {$ \scriptstyle \bullet$} [c] at 6 12
\put {$ \scriptstyle \bullet$} [c] at 10 12
\put {$ \scriptstyle \bullet$} [c] at 14 0
\put {$ \scriptstyle \bullet$} [c] at 14 6
\put {$ \scriptstyle \bullet$} [c] at 14 12
\setlinear \plot 6 12 6 0  10 12 14 6 14  12 6 0 /
\setlinear \plot 14 0 14  6     /
\put{$360$} [c] at 10 -2
\endpicture
\end{minipage}
\begin{minipage}{4cm}
\beginpicture
\setcoordinatesystem units <1.5mm,2mm>
\setplotarea x from 0 to 16, y from -2 to 15
\put{146)} [l] at 0 12
\put {$ \scriptstyle \bullet$} [c] at 6 0
\put {$ \scriptstyle \bullet$} [c] at 6 12
\put {$ \scriptstyle \bullet$} [c] at 10 0
\put {$ \scriptstyle \bullet$} [c] at 14 0
\put {$ \scriptstyle \bullet$} [c] at 14 6
\put {$ \scriptstyle \bullet$} [c] at 14 12
\setlinear \plot 6 0 6 12  10 0 14 6 14  0 6 12 /
\setlinear \plot 14 12 14  6     /
\put{$360$} [c] at 10 -2
\endpicture
\end{minipage}
\begin{minipage}{4cm}
\beginpicture
\setcoordinatesystem units <1.5mm,2mm>
\setplotarea x from 0 to 16, y from -2 to 15
\put{147)} [l] at 0 12
\put {$ \scriptstyle \bullet$} [c] at 6 6
\put {$ \scriptstyle \bullet$} [c] at 6 12
\put {$ \scriptstyle \bullet$} [c] at 10 12
\put {$ \scriptstyle \bullet$} [c] at 10 6
\put {$ \scriptstyle \bullet$} [c] at 8 0
\put {$ \scriptstyle \bullet$} [c] at 14 12
\setlinear \plot 14 12 8 0  6 6 6 12 10 6 10 12 6 6  /
\setlinear \plot 8 0 10  6     /
\put{$180$} [c] at 10 -2
\endpicture
\end{minipage}
\begin{minipage}{4cm}
\beginpicture
\setcoordinatesystem units <1.5mm,2mm>
\setplotarea x from 0 to 16, y from -2 to 15
\put{148)} [l] at 0 12
\put {$ \scriptstyle \bullet$} [c] at 6 6
\put {$ \scriptstyle \bullet$} [c] at 6 0
\put {$ \scriptstyle \bullet$} [c] at 10 0
\put {$ \scriptstyle \bullet$} [c] at 10 6
\put {$ \scriptstyle \bullet$} [c] at 8 12
\put {$ \scriptstyle \bullet$} [c] at 14 0
\setlinear \plot 14 0 8 12  6 6 6 0 10 6 10 0 6 6  /
\setlinear \plot 8 12 10  6     /
\put{$180$} [c] at 10 -2
\endpicture
\end{minipage}
\begin{minipage}{4cm}
\beginpicture
\setcoordinatesystem units <1.5mm,2mm>
\setplotarea x from 0 to 16, y from -2 to 15
\put{149)} [l] at 0 12
\put {$ \scriptstyle \bullet$} [c] at 6 0
\put {$ \scriptstyle \bullet$} [c] at 6 12
\put {$ \scriptstyle \bullet$} [c] at 10 0
\put {$ \scriptstyle \bullet$} [c] at 10 12
\put {$ \scriptstyle \bullet$} [c] at 14 6
\put {$ \scriptstyle \bullet$} [c] at 14 12
\setlinear \plot 14 12 14 6 10 0 10 12 6 0 6 12 10 0  /
\setlinear \plot 6 0  14 6 /
\put{$180$} [c] at 10 -2
\endpicture
\end{minipage}
\begin{minipage}{4cm}
\beginpicture
\setcoordinatesystem units <1.5mm,2mm>
\setplotarea x from 0 to 16, y from -2 to 15
\put{150)} [l] at 0 12
\put {$ \scriptstyle \bullet$} [c] at 6 0
\put {$ \scriptstyle \bullet$} [c] at 6 12
\put {$ \scriptstyle \bullet$} [c] at 10 0
\put {$ \scriptstyle \bullet$} [c] at 10 12
\put {$ \scriptstyle \bullet$} [c] at 14 6
\put {$ \scriptstyle \bullet$} [c] at 14 0
\setlinear \plot 14 0 14 6 10 12 10 0 6 12 6 0 10 12    /
\setlinear \plot 6 12  14 6 /
\put{$180$} [c] at 10 -2
\endpicture
\end{minipage}
$$
$$
\begin{minipage}{4cm}
\beginpicture
\setcoordinatesystem units <1.5mm,2mm>
\setplotarea x from 0 to 16, y from -2 to 15
\put{151)} [l] at 0 12
\put {$ \scriptstyle \bullet$} [c] at 6 0
\put {$ \scriptstyle \bullet$} [c] at 6 12
\put {$ \scriptstyle \bullet$} [c] at 10 12
\put {$ \scriptstyle \bullet$} [c] at 10 0
\put {$ \scriptstyle \bullet$} [c] at 14 0
\put {$ \scriptstyle \bullet$} [c] at 14 12
\setlinear \plot 6 12 6 0  10 12 14 0 14  12 10 0 6 12 14 0 /
\setlinear \plot 6 0 14  12     /
\setlinear \plot 10 0 10  12     /
\put{$20$} [c] at 10 -2
\endpicture
\end{minipage}
\begin{minipage}{4cm}
\beginpicture
\setcoordinatesystem units <1.5mm,2mm>
\setplotarea x from 0 to 16, y from -2 to 15
\put{152)} [l] at 0 12
\put {$ \scriptstyle \bullet$} [c] at 6 0
\put {$ \scriptstyle \bullet$} [c] at 6 4
\put {$ \scriptstyle \bullet$} [c] at 6 8
\put {$ \scriptstyle \bullet$} [c] at 6 12
\put {$ \scriptstyle \bullet$} [c] at 14 0
\put {$ \scriptstyle \bullet$} [c] at 14 12
\setlinear \plot 6 12 6  0    /
\setlinear \plot 14 0 14 12     /
\put{$720$} [c] at 10 -2
\endpicture
\end{minipage}
\begin{minipage}{4cm}
\beginpicture
\setcoordinatesystem units <1.5mm,2mm>
\setplotarea x from 0 to 16, y from -2 to 15
\put{ ${\bf  16}$} [l] at 0 15
\put{153)} [l] at 0 12
\put {$ \scriptstyle \bullet$} [c] at 6 0
\put {$ \scriptstyle \bullet$} [c] at 10 6
\put {$ \scriptstyle \bullet$} [c] at 12 12
\put {$ \scriptstyle \bullet$} [c] at 12 0
\put {$ \scriptstyle \bullet$} [c] at 14 6
\put {$ \scriptstyle \bullet$} [c] at 14 12
\setlinear \plot 6 0 12 12 10 6 12 0 14 6 14 12     /
\setlinear \plot 14 6  12 12 /
\put{$720$} [c] at 10 -2
\endpicture
\end{minipage}
\begin{minipage}{4cm}
\beginpicture
\setcoordinatesystem units <1.5mm,2mm>
\setplotarea x from 0 to 16, y from -2 to 15
\put{154)} [l] at 0 12
\put {$ \scriptstyle \bullet$} [c] at 6 12
\put {$ \scriptstyle \bullet$} [c] at 10 6
\put {$ \scriptstyle \bullet$} [c] at 12 12
\put {$ \scriptstyle \bullet$} [c] at 12 0
\put {$ \scriptstyle \bullet$} [c] at 14 6
\put {$ \scriptstyle \bullet$} [c] at 14 0
\setlinear \plot 6 12 12 0 10 6 12 12 14 6 14 0     /
\setlinear \plot 14 6  12 0 /
\put{$720$} [c] at 10 -2
\endpicture
\end{minipage}
\begin{minipage}{4cm}
\beginpicture
\setcoordinatesystem units <1.5mm,2mm>
\setplotarea x from 0 to 16, y from -2 to 15
\put{155)} [l] at 0 12
\put {$ \scriptstyle \bullet$} [c] at 6 0
\put {$ \scriptstyle \bullet$} [c] at 6 12
\put {$ \scriptstyle \bullet$} [c] at 10 0
\put {$ \scriptstyle \bullet$} [c] at 10 6
\put {$ \scriptstyle \bullet$} [c] at 10 12
\put {$ \scriptstyle \bullet$} [c] at 14 12
\setlinear \plot 10 6 6 0 6 12 10 0 10 12   /
\setlinear \plot 10 0 14 12   /
\put{$720$} [c] at 10 -2
\endpicture
\end{minipage}
\begin{minipage}{4cm}
\beginpicture
\setcoordinatesystem units <1.5mm,2mm>
\setplotarea x from 0 to 16, y from -2 to 15
\put{156)} [l] at 0 12
\put {$ \scriptstyle \bullet$} [c] at 6 0
\put {$ \scriptstyle \bullet$} [c] at 6 12
\put {$ \scriptstyle \bullet$} [c] at 10 0
\put {$ \scriptstyle \bullet$} [c] at 10 4
\put {$ \scriptstyle \bullet$} [c] at 10 12
\put {$ \scriptstyle \bullet$} [c] at 14 0
\setlinear \plot 10 4 6 12 6 0 10 12 10 0   /
\setlinear \plot 10 12 14 0   /
\put{$720$} [c] at 10 -2
\endpicture
\end{minipage}
$$
$$
\begin{minipage}{4cm}
\beginpicture
\setcoordinatesystem units <1.5mm,2mm>
\setplotarea x from 0 to 16, y from -2 to 15
\put{157)} [l] at 0 12
\put {$ \scriptstyle \bullet$} [c] at 6 12
\put {$ \scriptstyle \bullet$} [c] at 6 0
\put {$ \scriptstyle \bullet$} [c] at 10 12
\put {$ \scriptstyle \bullet$} [c] at 14 0
\put {$ \scriptstyle \bullet$} [c] at 14  6
\put {$ \scriptstyle \bullet$} [c] at 14 12
\setlinear \plot 6 12 6 0 10 12 14 6 14 12   /
\setlinear \plot 14 0 14 6  /
\put{$720$} [c] at 10 -2
\endpicture
\end{minipage}
\begin{minipage}{4cm}
\beginpicture
\setcoordinatesystem units <1.5mm,2mm>
\setplotarea x from 0 to 16, y from -2 to 15
\put{158)} [l] at 0 12
\put {$ \scriptstyle \bullet$} [c] at 6 12
\put {$ \scriptstyle \bullet$} [c] at 6 0
\put {$ \scriptstyle \bullet$} [c] at 10 0
\put {$ \scriptstyle \bullet$} [c] at 14 0
\put {$ \scriptstyle \bullet$} [c] at 14  6
\put {$ \scriptstyle \bullet$} [c] at 14 12
\setlinear \plot 6 0  6 12 10 0 14 6 14 0   /
\setlinear \plot 14 12 14 6  /
\put{$720$} [c] at 10 -2
\endpicture
\end{minipage}
\begin{minipage}{4cm}
\beginpicture
\setcoordinatesystem units <1.5mm,2mm>
\setplotarea x from 0 to 16, y from -2 to 15
\put{159)} [l] at 0 12
\put {$ \scriptstyle \bullet$} [c] at 6 0
\put {$ \scriptstyle \bullet$} [c] at 6 12
\put {$ \scriptstyle \bullet$} [c] at 8 9
\put {$ \scriptstyle \bullet$} [c] at 12 3
\put {$ \scriptstyle \bullet$} [c] at 14  12
\put {$ \scriptstyle \bullet$} [c] at 14 0
\setlinear \plot 6 0 6 12  14 0 14 12   /
\put{$720$} [c] at 10 -2
\endpicture
\end{minipage}
\begin{minipage}{4cm}
\beginpicture
\setcoordinatesystem units <1.5mm,2mm>
\setplotarea x from 0 to 16, y from -2 to 15
\put{160)} [l] at 0 12
\put {$ \scriptstyle \bullet$} [c] at 6 6
\put {$ \scriptstyle \bullet$} [c] at 9 0
\put {$ \scriptstyle \bullet$} [c] at 9 12
\put {$ \scriptstyle \bullet$} [c] at 12 6
\put {$ \scriptstyle \bullet$} [c] at 14 6
\put {$ \scriptstyle \bullet$} [c] at 14 12
\setlinear \plot 9  0 12 6 9 12 6 6 9 0 14 6 14 12    /
\put{$360$} [c] at 10 -2
\endpicture
\end{minipage}
\begin{minipage}{4cm}
\beginpicture
\setcoordinatesystem units <1.5mm,2mm>
\setplotarea x from 0 to 16, y from -2 to 15
\put{161)} [l] at 0 12
\put {$ \scriptstyle \bullet$} [c] at 6 6
\put {$ \scriptstyle \bullet$} [c] at 9 0
\put {$ \scriptstyle \bullet$} [c] at 9 12
\put {$ \scriptstyle \bullet$} [c] at 12 6
\put {$ \scriptstyle \bullet$} [c] at 14 6
\put {$ \scriptstyle \bullet$} [c] at 14 0
\setlinear \plot 9  12 12 6 9 0 6 6 9 12 14 6 14 0    /
\put{$360$} [c] at 10 -2
\endpicture
\end{minipage}
\begin{minipage}{4cm}
\beginpicture
\setcoordinatesystem units <1.5mm,2mm>
\setplotarea x from 0 to 16, y from -2 to 15
\put{162)} [l] at 0 12
\put {$ \scriptstyle \bullet$} [c] at 6 6
\put {$ \scriptstyle \bullet$} [c] at 6 12
\put {$ \scriptstyle \bullet$} [c] at 10 0
\put {$ \scriptstyle \bullet$} [c] at 10 12
\put {$ \scriptstyle \bullet$} [c] at 14 6
\put {$ \scriptstyle \bullet$} [c] at 14 12
\setlinear \plot 6 12 6 6 10 0 14 6 14 12     /
\setlinear \plot 10 12 14 6      /
\put{$360$} [c] at 10 -2
\endpicture
\end{minipage}
$$
$$
\begin{minipage}{4cm}
\beginpicture
\setcoordinatesystem units <1.5mm,2mm>
\setplotarea x from 0 to 16, y from -2 to 15
\put{163)} [l] at 0 12
\put {$ \scriptstyle \bullet$} [c] at 6 6
\put {$ \scriptstyle \bullet$} [c] at 6 0
\put {$ \scriptstyle \bullet$} [c] at 10 0
\put {$ \scriptstyle \bullet$} [c] at 10 12
\put {$ \scriptstyle \bullet$} [c] at 14 6
\put {$ \scriptstyle \bullet$} [c] at 14 0
\setlinear \plot 6 0 6 6 10 12 14 6 14 0     /
\setlinear \plot 10 0 14 6      /
\put{$360$} [c] at 10 -2
\endpicture
\end{minipage}
\begin{minipage}{4cm}
\beginpicture
\setcoordinatesystem units <1.5mm,2mm>
\setplotarea x from 0 to 16, y from -2 to 15
\put{164)} [l] at 0 12
\put {$ \scriptstyle \bullet$} [c] at 6 12
\put {$ \scriptstyle \bullet$} [c] at 6 0
\put {$ \scriptstyle \bullet$} [c] at 10 12
\put {$ \scriptstyle \bullet$} [c] at 14 0
\put {$ \scriptstyle \bullet$} [c] at 14 6
\put {$ \scriptstyle \bullet$} [c] at 14 12
\setlinear \plot 14 0 6 12 6 0 14 12 14 0 10 12 6 0   /
\put{$360$} [c] at 10 -2
\endpicture
\end{minipage}
\begin{minipage}{4cm}
\beginpicture
\setcoordinatesystem units <1.5mm,2mm>
\setplotarea x from 0 to 16, y from -2 to 15
\put{165)} [l] at 0 12
\put {$ \scriptstyle \bullet$} [c] at 6 12
\put {$ \scriptstyle \bullet$} [c] at 6 0
\put {$ \scriptstyle \bullet$} [c] at 10 0
\put {$ \scriptstyle \bullet$} [c] at 14 0
\put {$ \scriptstyle \bullet$} [c] at 14 6
\put {$ \scriptstyle \bullet$} [c] at 14 12
\setlinear \plot 14 12 6 0 6 12 14 0 14 12 10 0 6 12   /
\put{$360$} [c] at 10 -2
\endpicture
\end{minipage}
\begin{minipage}{4cm}
\beginpicture
\setcoordinatesystem units <1.5mm,2mm>
\setplotarea x from 0 to 16, y from -2 to 15
\put{166)} [l] at 0 12
\put {$ \scriptstyle \bullet$} [c] at 6 12
\put {$ \scriptstyle \bullet$} [c] at 6 0
\put {$ \scriptstyle \bullet$} [c] at 10 12
\put {$ \scriptstyle \bullet$} [c] at 14 0
\put {$ \scriptstyle \bullet$} [c] at 14  6
\put {$ \scriptstyle \bullet$} [c] at 14 12
\setlinear \plot 6 0 6 12 14 6 10 12   /
\setlinear \plot 14 0 14 12  /
\put{$360$} [c] at 10 -2
\endpicture
\end{minipage}
\begin{minipage}{4cm}
\beginpicture
\setcoordinatesystem units <1.5mm,2mm>
\setplotarea x from 0 to 16, y from -2 to 15
\put{167)} [l] at 0 12
\put {$ \scriptstyle \bullet$} [c] at 6 12
\put {$ \scriptstyle \bullet$} [c] at 6 0
\put {$ \scriptstyle \bullet$} [c] at 10 0
\put {$ \scriptstyle \bullet$} [c] at 14 0
\put {$ \scriptstyle \bullet$} [c] at 14  6
\put {$ \scriptstyle \bullet$} [c] at 14 12
\setlinear \plot 6 12 6 0 14 6 10 0   /
\setlinear \plot 14 0 14 12  /
\put{$360$} [c] at 10 -2
\endpicture
\end{minipage}
\begin{minipage}{4cm}
\beginpicture
\setcoordinatesystem units <1.5mm,2mm>
\setplotarea x from 0 to 16, y from -2 to 15
\put{168)} [l] at 0 12
\put {$ \scriptstyle \bullet$} [c] at 6 12
\put {$ \scriptstyle \bullet$} [c] at 6 0
\put {$ \scriptstyle \bullet$} [c] at 10 12
\put {$ \scriptstyle \bullet$} [c] at 14 0
\put {$ \scriptstyle \bullet$} [c] at 10  6
\put {$ \scriptstyle \bullet$} [c] at 14 12
\setlinear \plot 14 12 14  0 6 12 6 0 10 12 10 6    /
\put{$360$} [c] at 10 -2
\endpicture
\end{minipage}
$$
$$
\begin{minipage}{4cm}
\beginpicture
\setcoordinatesystem units <1.5mm,2mm>
\setplotarea x from 0 to 16, y from -2 to 15
\put{169)} [l] at 0 12
\put {$ \scriptstyle \bullet$} [c] at 6 12
\put {$ \scriptstyle \bullet$} [c] at 6 0
\put {$ \scriptstyle \bullet$} [c] at 10 0
\put {$ \scriptstyle \bullet$} [c] at 14 0
\put {$ \scriptstyle \bullet$} [c] at 10  6
\put {$ \scriptstyle \bullet$} [c] at 14 12
\setlinear \plot 14 0 14  12 6 0 6 12 10 0 10 6    /
\put{$360$} [c] at 10 -2
\endpicture
\end{minipage}
\begin{minipage}{4cm}
\beginpicture
\setcoordinatesystem units <1.5mm,2mm>
\setplotarea x from 0 to 16, y from -2 to 15
\put{170)} [l] at 0 12
\put {$ \scriptstyle \bullet$} [c] at 6 12
\put {$ \scriptstyle \bullet$} [c] at 6 0
\put {$ \scriptstyle \bullet$} [c] at 10 6
\put {$ \scriptstyle \bullet$} [c] at 12 12
\put {$ \scriptstyle \bullet$} [c] at 12 0
\put {$ \scriptstyle \bullet$} [c] at 14 6
\setlinear \plot 6 0 6 12  12 0 10 6 12 12 14 6 12 0 6 12 /
\setlinear \plot 6 0 12 12 /
\put{$360$} [c] at 10 -2
\endpicture
\end{minipage}
\begin{minipage}{4cm}
\beginpicture
\setcoordinatesystem units <1.5mm,2mm>
\setplotarea x from 0 to 16, y from -2 to 15
\put{171)} [l] at 0 12
\put {$ \scriptstyle \bullet$} [c] at 6 12
\put {$ \scriptstyle \bullet$} [c] at 6 0
\put {$ \scriptstyle \bullet$} [c] at 10 12
\put {$ \scriptstyle \bullet$} [c] at 10 0
\put {$ \scriptstyle \bullet$} [c] at 14  0
\put {$ \scriptstyle \bullet$} [c] at 14 12
\setlinear \plot  6 0 6 12 14 0 10 12 6 0 14 12 10 0 10  12 /
\setlinear \plot  6 12 10 0 /
\put{$180$} [c] at 10 -2
\endpicture
\end{minipage}
\begin{minipage}{4cm}
\beginpicture
\setcoordinatesystem units <1.5mm,2mm>
\setplotarea x from 0 to 16, y from -2 to 15
\put{172)} [l] at 0 12
\put {$ \scriptstyle \bullet$} [c] at 6 12
\put {$ \scriptstyle \bullet$} [c] at 7.2 9
\put {$ \scriptstyle \bullet$} [c] at 9 4
\put {$ \scriptstyle \bullet$} [c] at 9 0
\put {$ \scriptstyle \bullet$} [c] at 12 12
\put {$ \scriptstyle \bullet$} [c] at 14 0
\setlinear \plot  9 0 9 4 6 12  /
\setlinear \plot  9 4 12 12  /
\put{$720$} [c] at 10 -2
\endpicture
\end{minipage}
\begin{minipage}{4cm}
\beginpicture
\setcoordinatesystem units <1.5mm,2mm>
\setplotarea x from 0 to 16, y from -2 to 15
\put{173)} [l] at 0 12
\put {$ \scriptstyle \bullet$} [c] at 6 0
\put {$ \scriptstyle \bullet$} [c] at 7.2 3
\put {$ \scriptstyle \bullet$} [c] at 9 12
\put {$ \scriptstyle \bullet$} [c] at 9 8
\put {$ \scriptstyle \bullet$} [c] at 12 0
\put {$ \scriptstyle \bullet$} [c] at 14 0
\setlinear \plot  9 12 9 8 6 0  /
\setlinear \plot  9 8 12 0  /
\put{$720$} [c] at 10 -2
\endpicture
\end{minipage}
\begin{minipage}{4cm}
\beginpicture
\setcoordinatesystem units <1.5mm,2mm>
\setplotarea x from 0 to 16, y from -2 to 15
\put{174)} [l] at 0 12
\put {$ \scriptstyle \bullet$} [c] at 6 8
\put {$ \scriptstyle \bullet$} [c] at 6 4
\put {$ \scriptstyle \bullet$} [c] at 9 12
\put {$ \scriptstyle \bullet$} [c] at 9 0
\put {$ \scriptstyle \bullet$} [c] at 12  6
\put {$ \scriptstyle \bullet$} [c] at 14 0
\setlinear \plot 9 0 6 4  6 8  9 12 12 6  9 0   /
\put{$720$} [c] at 10 -2
\endpicture
\end{minipage}
$$
$$
\begin{minipage}{4cm}
\beginpicture
\setcoordinatesystem units <1.5mm,2mm>
\setplotarea x from 0 to 16, y from -2 to 15
\put{175)} [l] at 0 12
\put {$ \scriptstyle \bullet$} [c] at 6 0
\put {$ \scriptstyle \bullet$} [c] at 6 6
\put {$ \scriptstyle \bullet$} [c] at 6 12
\put {$ \scriptstyle \bullet$} [c] at 14 6
\put {$ \scriptstyle \bullet$} [c] at 14 0
\put {$ \scriptstyle \bullet$} [c] at 14 12
\setlinear \plot  6 0 6 12  /
\setlinear \plot  14 0 14 12  /
\put{$360$} [c] at 10 -2
\endpicture
\end{minipage}
\begin{minipage}{4cm}
\beginpicture
\setcoordinatesystem units <1.5mm,2mm>
\setplotarea x from 0 to 16, y from -2 to 15
\put{176)} [l] at 0 12
\put {$ \scriptstyle \bullet$} [c] at 6 12
\put {$ \scriptstyle \bullet$} [c] at 6 6
\put {$ \scriptstyle \bullet$} [c] at 9 0
\put {$ \scriptstyle \bullet$} [c] at 12 6
\put {$ \scriptstyle \bullet$} [c] at 12 12
\put {$ \scriptstyle \bullet$} [c] at 14 0
\setlinear \plot  9 0 6 6 6 12 12 6 12 12 6 6 9 0 12 6 /
\put{$180$} [c] at 10 -2
\endpicture
\end{minipage}
\begin{minipage}{4cm}
\beginpicture
\setcoordinatesystem units <1.5mm,2mm>
\setplotarea x from 0 to 16, y from -2 to 15
\put{177)} [l] at 0 12
\put {$ \scriptstyle \bullet$} [c] at 6 0
\put {$ \scriptstyle \bullet$} [c] at 6 6
\put {$ \scriptstyle \bullet$} [c] at 9 12
\put {$ \scriptstyle \bullet$} [c] at 12 6
\put {$ \scriptstyle \bullet$} [c] at 12 0
\put {$ \scriptstyle \bullet$} [c] at 14 0
\setlinear \plot  9 12 6 6 6 0 12 6 12 0 6 6 9 12 12 6 /
\put{$180$} [c] at 10 -2
\endpicture
\end{minipage}
\begin{minipage}{4cm}
\beginpicture
\setcoordinatesystem units <1.5mm,2mm>
\setplotarea x from 0 to 16, y from -2 to 15
\put{178)} [l] at 0 12
\put {$ \scriptstyle \bullet$} [c] at 6 0
\put {$ \scriptstyle \bullet$} [c] at 6 12
\put {$ \scriptstyle \bullet$} [c] at 9 6
\put {$ \scriptstyle \bullet$} [c] at 12 12
\put {$ \scriptstyle \bullet$} [c] at 12 0
\put {$ \scriptstyle \bullet$} [c] at 14 0
\setlinear \plot  6 0 12 12  /
\setlinear \plot  6 12 12 0 /
\put{$180$} [c] at 10 -2
\endpicture
\end{minipage}
\begin{minipage}{4cm}
\beginpicture
\setcoordinatesystem units <1.5mm,2mm>
\setplotarea x from 0 to 16, y from -2 to 15
\put{ ${\bf  17}$} [l] at 0 15
\put{179)} [l] at 0 12
\put {$ \scriptstyle \bullet$} [c] at 6 12
\put {$ \scriptstyle \bullet$} [c] at 9 6
\put {$ \scriptstyle \bullet$} [c] at 10 0
\put {$ \scriptstyle \bullet$} [c] at 10 12
\put {$ \scriptstyle \bullet$} [c] at 11 6
\put {$ \scriptstyle \bullet$} [c] at 14 12
\setlinear \plot 6 12 10 0  9 6 10 12 11 6 10 0  /
\setlinear \plot 11 6 14 12 /
\put{$720$} [c] at 10 -2
\endpicture
\end{minipage}
\begin{minipage}{4cm}
\beginpicture
\setcoordinatesystem units <1.5mm,2mm>
\setplotarea x from 0 to 16, y from -2 to 15
\put{180)} [l] at 0 12
\put {$ \scriptstyle \bullet$} [c] at 6 0
\put {$ \scriptstyle \bullet$} [c] at 9 6
\put {$ \scriptstyle \bullet$} [c] at 10 0
\put {$ \scriptstyle \bullet$} [c] at 10 12
\put {$ \scriptstyle \bullet$} [c] at 11 6
\put {$ \scriptstyle \bullet$} [c] at 14 0
\setlinear \plot 6 0 10 12  9 6 10 0 11 6 10 12  /
\setlinear \plot 11 6 14 0 /
\put{$720$} [c] at 10 -2
\endpicture
\end{minipage}
$$
$$
\begin{minipage}{4cm}
\beginpicture
\setcoordinatesystem units <1.5mm,2mm>
\setplotarea x from 0 to 16, y from -2 to 15
\put{181)} [l] at 0 12
\put {$ \scriptstyle \bullet$} [c] at 6 0
\put {$ \scriptstyle \bullet$} [c] at 6 6
\put {$ \scriptstyle \bullet$} [c] at 6 12
\put {$ \scriptstyle \bullet$} [c] at 10 6
\put {$ \scriptstyle \bullet$} [c] at 14 12
\put {$ \scriptstyle \bullet$} [c] at 14 0
\setlinear \plot  6 12 6 0 14 12 14 0   /
\put{$720$} [c] at 10 -2
\endpicture
\end{minipage}
\begin{minipage}{4cm}
\beginpicture
\setcoordinatesystem units <1.5mm,2mm>
\setplotarea x from 0 to 16, y from -2 to 15
\put{182)} [l] at 0 12
\put {$ \scriptstyle \bullet$} [c] at 6 0
\put {$ \scriptstyle \bullet$} [c] at 6 6
\put {$ \scriptstyle \bullet$} [c] at 6 12
\put {$ \scriptstyle \bullet$} [c] at 10 6
\put {$ \scriptstyle \bullet$} [c] at 14 12
\put {$ \scriptstyle \bullet$} [c] at 14 0
\setlinear \plot  6 0 6 12 14 0 14 12   /
\put{$720$} [c] at 10 -2
\endpicture
\end{minipage}
\begin{minipage}{4cm}
\beginpicture
\setcoordinatesystem units <1.5mm,2mm>
\setplotarea x from 0 to 16, y from -2 to 15
\put{183)} [l] at 0 12
\put {$ \scriptstyle \bullet$} [c] at 6 0
\put {$ \scriptstyle \bullet$} [c] at 6 12
\put {$ \scriptstyle \bullet$} [c] at 10 12
\put {$ \scriptstyle \bullet$} [c] at 14 0
\put {$ \scriptstyle \bullet$} [c] at 14 6
\put {$ \scriptstyle \bullet$} [c] at 14 12
\setlinear \plot 6 12 6 0  10 12 14 0 14 12 6 0 /
\put{$720$} [c] at 10 -2
\endpicture
\end{minipage}
\begin{minipage}{4cm}
\beginpicture
\setcoordinatesystem units <1.5mm,2mm>
\setplotarea x from 0 to 16, y from -2 to 15
\put{184)} [l] at 0 12
\put {$ \scriptstyle \bullet$} [c] at 6 0
\put {$ \scriptstyle \bullet$} [c] at 6 12
\put {$ \scriptstyle \bullet$} [c] at 10 0
\put {$ \scriptstyle \bullet$} [c] at 14 0
\put {$ \scriptstyle \bullet$} [c] at 14 6
\put {$ \scriptstyle \bullet$} [c] at 14 12
\setlinear \plot 6 0 6 12  10 0 14 12 14 0 6 12 /
\put{$720$} [c] at 10 -2
\endpicture
\end{minipage}
\begin{minipage}{4cm}
\beginpicture
\setcoordinatesystem units <1.5mm,2mm>
\setplotarea x from 0 to 16, y from -2 to 15
\put{185)} [l] at 0 12
\put {$ \scriptstyle \bullet$} [c] at 6 12
\put {$ \scriptstyle \bullet$} [c] at 10 0
\put {$ \scriptstyle \bullet$} [c] at 10 4
\put {$ \scriptstyle \bullet$} [c] at 10 8
\put {$ \scriptstyle \bullet$} [c] at 10 12
\put {$ \scriptstyle \bullet$} [c] at 14 12
\setlinear \plot  6  12 10 0 14 12    /
\setlinear \plot  10  0  10 12    /
\put{$360$} [c] at 10 -2
\endpicture
\end{minipage}
\begin{minipage}{4cm}
\beginpicture
\setcoordinatesystem units <1.5mm,2mm>
\setplotarea x from 0 to 16, y from -2 to 15
\put{186)} [l] at 0 12
\put {$ \scriptstyle \bullet$} [c] at 6 0
\put {$ \scriptstyle \bullet$} [c] at 10 0
\put {$ \scriptstyle \bullet$} [c] at 10 4
\put {$ \scriptstyle \bullet$} [c] at 10 8
\put {$ \scriptstyle \bullet$} [c] at 10 12
\put {$ \scriptstyle \bullet$} [c] at 14 0
\setlinear \plot  6  0 10 12 14 0    /
\setlinear \plot  10  0  10 12    /
\put{$360$} [c] at 10 -2
\endpicture
\end{minipage}
$$

$$
\begin{minipage}{4cm}
\beginpicture
\setcoordinatesystem units <1.5mm,2mm>
\setplotarea x from 0 to 16, y from -2 to 15
\put{187)} [l] at 0 12
\put {$ \scriptstyle \bullet$} [c] at 6 6
\put {$ \scriptstyle \bullet$} [c] at 8 0
\put {$ \scriptstyle \bullet$} [c] at 8 12
\put {$ \scriptstyle \bullet$} [c] at 10 6
\put {$ \scriptstyle \bullet$} [c] at 14 12
\put {$ \scriptstyle \bullet$} [c] at 14 0
\setlinear \plot 14 0 14  12 8 0  10 6 8 12 6 6 8 0  /
\put{$360$} [c] at 10 -2
\endpicture
\end{minipage}
\begin{minipage}{4cm}
\beginpicture
\setcoordinatesystem units <1.5mm,2mm>
\setplotarea x from 0 to 16, y from -2 to 15
\put{188)} [l] at 0 12
\put {$ \scriptstyle \bullet$} [c] at 6 6
\put {$ \scriptstyle \bullet$} [c] at  8 0
\put {$ \scriptstyle \bullet$} [c] at 8 12
\put {$ \scriptstyle \bullet$} [c] at 10 6
\put {$ \scriptstyle \bullet$} [c] at 14 12
\put {$ \scriptstyle \bullet$} [c] at 14 0
\setlinear \plot 14 12 14  0 8 12  10 6 8 0 6 6 8 12  /
\put{$360$} [c] at 10 -2
\endpicture
\end{minipage}
\begin{minipage}{4cm}
\beginpicture
\setcoordinatesystem units <1.5mm,2mm>
\setplotarea x from 0 to 16, y from -2 to 15
\put{189)} [l] at 0 12
\put {$ \scriptstyle \bullet$} [c] at 6 0
\put {$ \scriptstyle \bullet$} [c] at 6 12
\put {$ \scriptstyle \bullet$} [c] at 10 6
\put {$ \scriptstyle \bullet$} [c] at 10 12
\put {$ \scriptstyle \bullet$} [c] at 14 12
\put {$ \scriptstyle \bullet$} [c] at 14 0
\setlinear \plot  6 12 6  0  10 6 10  12    /
\setlinear \plot  10 6 14  0  14 12    /
\put{$360$} [c] at 10 -2
\endpicture
\end{minipage}
\begin{minipage}{4cm}
\beginpicture
\setcoordinatesystem units <1.5mm,2mm>
\setplotarea x from 0 to 16, y from -2 to 15
\put{190)} [l] at 0 12
\put {$ \scriptstyle \bullet$} [c] at 6 0
\put {$ \scriptstyle \bullet$} [c] at 6 12
\put {$ \scriptstyle \bullet$} [c] at 10 6
\put {$ \scriptstyle \bullet$} [c] at 10 0
\put {$ \scriptstyle \bullet$} [c] at 14 12
\put {$ \scriptstyle \bullet$} [c] at 14 0
\setlinear \plot  6 0 6  12  10 6 10 0    /
\setlinear \plot  10 6 14  12  14 0    /
\put{$360$} [c] at 10 -2
\endpicture
\end{minipage}
\begin{minipage}{4cm}
\beginpicture
\setcoordinatesystem units <1.5mm,2mm>
\setplotarea x from 0 to 16, y from -2 to 15
\put{191)} [l] at 0 12
\put {$ \scriptstyle \bullet$} [c] at 6 0
\put {$ \scriptstyle \bullet$} [c] at 6 12
\put {$ \scriptstyle \bullet$} [c] at 10 0
\put {$ \scriptstyle \bullet$} [c] at 10 12
\put {$ \scriptstyle \bullet$} [c] at 12 6
\put {$ \scriptstyle \bullet$} [c] at 14 12
\setlinear \plot  14 12 10  0 10 12 6 0  6 12 10 0   /
\put{$360$} [c] at 10 -2
\endpicture
\end{minipage}
\begin{minipage}{4cm}
\beginpicture
\setcoordinatesystem units <1.5mm,2mm>
\setplotarea x from 0 to 16, y from -2 to 15
\put{192)} [l] at 0 12
\put {$ \scriptstyle \bullet$} [c] at 6 0
\put {$ \scriptstyle \bullet$} [c] at 6 12
\put {$ \scriptstyle \bullet$} [c] at 10 0
\put {$ \scriptstyle \bullet$} [c] at 10 12
\put {$ \scriptstyle \bullet$} [c] at 12 6
\put {$ \scriptstyle \bullet$} [c] at 14 0
\setlinear \plot  14 0 10 12 10 0 6 12  6 0 10 12   /
\put{$360$} [c] at 10 -2
\endpicture
\end{minipage}
$$
$$
\begin{minipage}{4cm}
\beginpicture
\setcoordinatesystem units <1.5mm,2mm>
\setplotarea x from 0 to 16, y from -2 to 15
\put{193)} [l] at 0 12
\put {$ \scriptstyle \bullet$} [c] at 6 0
\put {$ \scriptstyle \bullet$} [c] at 6 12
\put {$ \scriptstyle \bullet$} [c] at 12 6
\put {$ \scriptstyle \bullet$} [c] at 12 0
\put {$ \scriptstyle \bullet$} [c] at 10 12
\put {$ \scriptstyle \bullet$} [c] at 14 12
\setlinear \plot 6 0 6 12 12 0  12 6 10 12  /
\setlinear \plot 12 6 14 12  /
\put{$360$} [c] at 10 -2
\endpicture
\end{minipage}
\begin{minipage}{4cm}
\beginpicture
\setcoordinatesystem units <1.5mm,2mm>
\setplotarea x from 0 to 16, y from -2 to 15
\put{194)} [l] at 0 12
\put {$ \scriptstyle \bullet$} [c] at 6 0
\put {$ \scriptstyle \bullet$} [c] at 6 12
\put {$ \scriptstyle \bullet$} [c] at 12 6
\put {$ \scriptstyle \bullet$} [c] at 12 12
\put {$ \scriptstyle \bullet$} [c] at 10 0
\put {$ \scriptstyle \bullet$} [c] at 14 0
\setlinear \plot 6 12 6 0 12 12  12 6 10 0  /
\setlinear \plot 12 6 14 0  /
\put{$360$} [c] at 10 -2
\endpicture
\end{minipage}
\begin{minipage}{4cm}
\beginpicture
\setcoordinatesystem units <1.5mm,2mm>
\setplotarea x from 0 to 16, y from -2 to 15
\put{195)} [l] at 0 12
\put {$ \scriptstyle \bullet$} [c] at 6 0
\put {$ \scriptstyle \bullet$} [c] at 6 12
\put {$ \scriptstyle \bullet$} [c] at 10 12
\put {$ \scriptstyle \bullet$} [c] at 10 0
\put {$ \scriptstyle \bullet$} [c] at 14 12
\put {$ \scriptstyle \bullet$} [c] at 14 0
\setlinear \plot  6  0  6  12 10 0 14 12 14 0  10 12  6 0  /
\setlinear \plot  10  0  10 12    /
\put{$360$} [c] at 10 -2
\endpicture
\end{minipage}
\begin{minipage}{4cm}
\beginpicture
\setcoordinatesystem units <1.5mm,2mm>
\setplotarea x from 0 to 16, y from -2 to 15
\put{ ${\bf  18}$} [l] at 0 15
\put{196)} [l] at 0 12
\put {$ \scriptstyle \bullet$} [c] at 6 0
\put {$ \scriptstyle \bullet$} [c] at 6 12
\put {$ \scriptstyle \bullet$} [c] at 10 12
\put {$ \scriptstyle \bullet$} [c] at 14 12
\put {$ \scriptstyle \bullet$} [c] at 14 0
\put {$ \scriptstyle \bullet$} [c] at 14 6
\setlinear \plot 6 12 6 0 10 12 14 0 14 12   /
\put{$720$} [c] at 10 -2
\endpicture
\end{minipage}
\begin{minipage}{4cm}
\beginpicture
\setcoordinatesystem units <1.5mm,2mm>
\setplotarea x from 0 to 16, y from -2 to 15
\put{197)} [l] at 0 12
\put {$ \scriptstyle \bullet$} [c] at 6 0
\put {$ \scriptstyle \bullet$} [c] at 6 12
\put {$ \scriptstyle \bullet$} [c] at 10 0
\put {$ \scriptstyle \bullet$} [c] at 14 12
\put {$ \scriptstyle \bullet$} [c] at 14 0
\put {$ \scriptstyle \bullet$} [c] at 14 6
\setlinear \plot 6 0 6 12 10 0 14 12 14 0   /
\put{$720$} [c] at 10 -2
\endpicture
\end{minipage}
\begin{minipage}{4cm}
\beginpicture
\setcoordinatesystem units <1.5mm,2mm>
\setplotarea x from 0 to 16, y from -2 to 15
\put{198)} [l] at 0 12
\put {$ \scriptstyle \bullet$} [c] at 6 0
\put {$ \scriptstyle \bullet$} [c] at 6 12
\put {$ \scriptstyle \bullet$} [c] at 10 6
\put {$ \scriptstyle \bullet$} [c] at 10 12
\put {$ \scriptstyle \bullet$} [c] at 14 0
\put {$ \scriptstyle \bullet$} [c] at 14 12
\setlinear \plot 6 0  6 12 14 0 14 12  /
\setlinear \plot 10 6 10  12  /
\put{$720$} [c] at 10 -2
\endpicture
\end{minipage}
$$
$$
\begin{minipage}{4cm}
\beginpicture
\setcoordinatesystem units <1.5mm,2mm>
\setplotarea x from 0 to 16, y from -2 to 15
\put{199)} [l] at 0 12
\put {$ \scriptstyle \bullet$} [c] at 6 0
\put {$ \scriptstyle \bullet$} [c] at 6 12
\put {$ \scriptstyle \bullet$} [c] at 10 6
\put {$ \scriptstyle \bullet$} [c] at 10 0
\put {$ \scriptstyle \bullet$} [c] at 14 0
\put {$ \scriptstyle \bullet$} [c] at 14 12
\setlinear \plot 6 12  6 0 14 12 14 0  /
\setlinear \plot 10 6 10  0  /
\put{$720$} [c] at 10 -2
\endpicture
\end{minipage}
\begin{minipage}{4cm}
\beginpicture
\setcoordinatesystem units <1.5mm,2mm>
\setplotarea x from 0 to 16, y from -2 to 15
\put{200)} [l] at 0 12
\put {$ \scriptstyle \bullet$} [c] at 6 0
\put {$ \scriptstyle \bullet$} [c] at 6 12
\put {$ \scriptstyle \bullet$} [c] at 12 0
\put {$ \scriptstyle \bullet$} [c] at 12 6
\put {$ \scriptstyle \bullet$} [c] at 12 12
\put {$ \scriptstyle \bullet$} [c] at 14 12
\setlinear \plot 14 12 12 0 12  12 6 0 6 12 12 0 /
\put{$720$} [c] at 10 -2
\endpicture
\end{minipage}
\begin{minipage}{4cm}
\beginpicture
\setcoordinatesystem units <1.5mm,2mm>
\setplotarea x from 0 to 16, y from -2 to 15
\put{201)} [l] at 0 12
\put {$ \scriptstyle \bullet$} [c] at 6 0
\put {$ \scriptstyle \bullet$} [c] at 6 12
\put {$ \scriptstyle \bullet$} [c] at 12 0
\put {$ \scriptstyle \bullet$} [c] at 12 6
\put {$ \scriptstyle \bullet$} [c] at 12 12
\put {$ \scriptstyle \bullet$} [c] at 14 0
\setlinear \plot 14 0 12 12 12  0 6 12 6 0 12 12 /
\put{$720$} [c] at 10 -2
\endpicture
\end{minipage}
\begin{minipage}{4cm}
\beginpicture
\setcoordinatesystem units <1.5mm,2mm>
\setplotarea x from 0 to 16, y from -2 to 15
\put{202)} [l] at 0 12
\put {$ \scriptstyle \bullet$} [c] at 6 12
\put {$ \scriptstyle \bullet$} [c] at 10 0
\put {$ \scriptstyle \bullet$} [c] at 10 12
\put {$ \scriptstyle \bullet$} [c] at 12 6
\put {$ \scriptstyle \bullet$} [c] at 12 12
\put {$ \scriptstyle \bullet$} [c] at 14 0
\setlinear \plot 6 12 10 0 10 12  /
\setlinear \plot 10 0 12 6 12 12  /
\setlinear \plot 14 0 12 6 /
\put{$360$} [c] at 10 -2
\endpicture
\end{minipage}
\begin{minipage}{4cm}
\beginpicture
\setcoordinatesystem units <1.5mm,2mm>
\setplotarea x from 0 to 16, y from -2 to 15
\put{203)} [l] at 0 12
\put {$ \scriptstyle \bullet$} [c] at 6 0
\put {$ \scriptstyle \bullet$} [c] at 10 0
\put {$ \scriptstyle \bullet$} [c] at 10 12
\put {$ \scriptstyle \bullet$} [c] at 12 6
\put {$ \scriptstyle \bullet$} [c] at 12 0
\put {$ \scriptstyle \bullet$} [c] at 14 12
\setlinear \plot 6 0 10 12 10 0  /
\setlinear \plot 10 12 12 6 12 0  /
\setlinear \plot 14 12 12 6 /
\put{$360$} [c] at 10 -2
\endpicture
\end{minipage}
\begin{minipage}{4cm}
\beginpicture
\setcoordinatesystem units <1.5mm,2mm>
\setplotarea x from 0 to 16, y from -2 to 15
\put{204)} [l] at 0 12
\put {$ \scriptstyle \bullet$} [c] at 6 0
\put {$ \scriptstyle \bullet$} [c] at 6 12
\put {$ \scriptstyle \bullet$} [c] at 10 0
\put {$ \scriptstyle \bullet$} [c] at 10 12
\put {$ \scriptstyle \bullet$} [c] at 14 0
\put {$ \scriptstyle \bullet$} [c] at 14 12
\setlinear \plot 14 0 10 12 10 0  6 12 6  0 10 12  /
\setlinear \plot 6 12 14  0  14 12   /
\put{$180$} [c] at 10 -2
\endpicture
\end{minipage}
$$
$$
\begin{minipage}{4cm}
\beginpicture
\setcoordinatesystem units <1.5mm,2mm>
\setplotarea x from 0 to 16, y from -2 to 15
\put{205)} [l] at 0 12
\put {$ \scriptstyle \bullet$} [c] at 6 0
\put {$ \scriptstyle \bullet$} [c] at 6 12
\put {$ \scriptstyle \bullet$} [c] at 10 0
\put {$ \scriptstyle \bullet$} [c] at 10 12
\put {$ \scriptstyle \bullet$} [c] at 14 0
\put {$ \scriptstyle \bullet$} [c] at 14 12
\setlinear \plot 14 12 10 0 10 12  6 0 6  12 10 0  /
\setlinear \plot 6 0 14 12  14 0   /
\put{$180$} [c] at 10 -2
\endpicture
\end{minipage}
\begin{minipage}{4cm}
\beginpicture
\setcoordinatesystem units <1.5mm,2mm>
\setplotarea x from 0 to 16, y from -2 to 15
\put{206)} [l] at 0 12
\put {$ \scriptstyle \bullet$} [c] at 6 0
\put {$ \scriptstyle \bullet$} [c] at 6 12
\put {$ \scriptstyle \bullet$} [c] at 10 0
\put {$ \scriptstyle \bullet$} [c] at 10 12
\put {$ \scriptstyle \bullet$} [c] at 14 0
\put {$ \scriptstyle \bullet$} [c] at 14 12
\setlinear \plot 6 0 6 12  10 0 14 12 14 0 10 12 6 0 /
\put{$120$} [c] at 10 -2
\endpicture
\end{minipage}
\begin{minipage}{4cm}
\beginpicture
\setcoordinatesystem units <1.5mm,2mm>
\setplotarea x from 0 to 16, y from -2 to 15
\put{207)} [l] at 0 12
\put {$ \scriptstyle \bullet$} [c] at 6 12
\put {$ \scriptstyle \bullet$} [c] at 8 12
\put {$ \scriptstyle \bullet$} [c] at 12 12
\put {$ \scriptstyle \bullet$} [c] at 14 12
\put {$ \scriptstyle \bullet$} [c] at 10 6
\put {$ \scriptstyle \bullet$} [c] at 10 0
\setlinear \plot 6 12 10 6 14 12    /
\setlinear \plot 8 12 10 6 12 12    /
\setlinear \plot 10 0 10 6    /
\put{$30$} [c] at 10 -2
\endpicture
\end{minipage}
\begin{minipage}{4cm}
\beginpicture
\setcoordinatesystem units <1.5mm,2mm>
\setplotarea x from 0 to 16, y from -2 to 15
\put{208)} [l] at 0 12
\put {$ \scriptstyle \bullet$} [c] at 6 0
\put {$ \scriptstyle \bullet$} [c] at 8 0
\put {$ \scriptstyle \bullet$} [c] at 12 0
\put {$ \scriptstyle \bullet$} [c] at 14 0
\put {$ \scriptstyle \bullet$} [c] at 10 6
\put {$ \scriptstyle \bullet$} [c] at 10 12
\setlinear \plot 6 0 10 6 14 0    /
\setlinear \plot 8 0 10 6 12 0    /
\setlinear \plot 10 12 10 6    /
\put{$30$} [c] at 10 -2
\endpicture
\end{minipage}
\begin{minipage}{4cm}
\beginpicture
\setcoordinatesystem units <1.5mm,2mm>
\setplotarea x from 0 to 16, y from -2 to 15
\put{209)} [l] at 0 12
\put {$ \scriptstyle \bullet$} [c] at 6 6
\put {$ \scriptstyle \bullet$} [c] at 8 6
\put {$ \scriptstyle \bullet$} [c] at 12 6
\put {$ \scriptstyle \bullet$} [c] at 14 6
\put {$ \scriptstyle \bullet$} [c] at 10 0
\put {$ \scriptstyle \bullet$} [c] at 10 12
\setlinear \plot  10 0 6  6  10 12 14 6 10  0   /
\setlinear \plot  10 0 8  6  10 12 12 6 10  0   /
\put{$30$} [c] at 10 -2
\endpicture
\end{minipage}
\begin{minipage}{4cm}
\beginpicture
\setcoordinatesystem units <1.5mm,2mm>
\setplotarea x from 0 to 16, y from -2 to 15
\put{210)} [l] at 0 12
\put {$ \scriptstyle \bullet$} [c] at 6 12
\put {$ \scriptstyle \bullet$} [c] at 7 8
\put {$ \scriptstyle \bullet$} [c] at 8 4
\put {$ \scriptstyle \bullet$} [c] at 9 0
\put {$ \scriptstyle \bullet$} [c] at 12 12
\put {$ \scriptstyle \bullet$} [c] at 14 0
\setlinear \plot 6 12 9 0  12 12  /
\put{$720$} [c] at 10 -2
\endpicture
\end{minipage}
$$
$$
\begin{minipage}{4cm}
\beginpicture
\setcoordinatesystem units <1.5mm,2mm>
\setplotarea x from 0 to 16, y from -2 to 15
\put{211)} [l] at 0 12
\put {$ \scriptstyle \bullet$} [c] at 6 0
\put {$ \scriptstyle \bullet$} [c] at 7 4
\put {$ \scriptstyle \bullet$} [c] at 8 8
\put {$ \scriptstyle \bullet$} [c] at 9 12
\put {$ \scriptstyle \bullet$} [c] at 12 0
\put {$ \scriptstyle \bullet$} [c] at 14 0
\setlinear \plot 6 0 9 12  12 0  /
\put{$720$} [c] at 10 -2
\endpicture
\end{minipage}
\begin{minipage}{4cm}
\beginpicture
\setcoordinatesystem units <1.5mm,2mm>
\setplotarea x from 0 to 16, y from -2 to 15
\put{212)} [l] at 0 12
\put {$ \scriptstyle \bullet$} [c] at 6 0
\put {$ \scriptstyle \bullet$} [c] at 8 6
\put {$ \scriptstyle \bullet$} [c] at 10 0
\put {$ \scriptstyle \bullet$} [c] at 10 12
\put {$ \scriptstyle \bullet$} [c] at 12 6
\put {$ \scriptstyle \bullet$} [c] at 14 0
\setlinear \plot 6 0 8 6 10 12 12 6 10  0 8 6  /
\put{$720$} [c] at 10 -2
\endpicture
\end{minipage}
\begin{minipage}{4cm}
\beginpicture
\setcoordinatesystem units <1.5mm,2mm>
\setplotarea x from 0 to 16, y from -2 to 15
\put{213)} [l] at 0 12
\put {$ \scriptstyle \bullet$} [c] at 6 12
\put {$ \scriptstyle \bullet$} [c] at 8 6
\put {$ \scriptstyle \bullet$} [c] at 10 0
\put {$ \scriptstyle \bullet$} [c] at 10 12
\put {$ \scriptstyle \bullet$} [c] at 12 6
\put {$ \scriptstyle \bullet$} [c] at 14 0
\setlinear \plot 6 12 8 6 10 12 12 6 10  0 8 6  /
\put{$720$} [c] at 10 -2
\endpicture
\end{minipage}
\begin{minipage}{4cm}
\beginpicture
\setcoordinatesystem units <1.5mm,2mm>
\setplotarea x from 0 to 16, y from -2 to 15
\put{214)} [l] at 0 12
\put {$ \scriptstyle \bullet$} [c] at 6 12
\put {$ \scriptstyle \bullet$} [c] at 9 0
\put {$ \scriptstyle \bullet$} [c] at 9 8
\put {$ \scriptstyle \bullet$} [c] at 12 12
\put {$ \scriptstyle \bullet$} [c] at 14 12
\put {$ \scriptstyle \bullet$} [c] at 14 0
\setlinear \plot 6 12 9 8 9 0   /
\setlinear \plot 12 12 9 8    /
\setlinear \plot 14  0 14 12   /
\put{$360$} [c] at 10 -2
\endpicture
\end{minipage}
\begin{minipage}{4cm}
\beginpicture
\setcoordinatesystem units <1.5mm,2mm>
\setplotarea x from 0 to 16, y from -2 to 15
\put{215)} [l] at 0 12
\put {$ \scriptstyle \bullet$} [c] at 6 0
\put {$ \scriptstyle \bullet$} [c] at 9 12
\put {$ \scriptstyle \bullet$} [c] at 9 4
\put {$ \scriptstyle \bullet$} [c] at 12 0
\put {$ \scriptstyle \bullet$} [c] at 14 12
\put {$ \scriptstyle \bullet$} [c] at 14 0
\setlinear \plot  9 12 9 4  6  0  /
\setlinear \plot  12 0 9 4   /
\setlinear \plot  14  0  14 12    /
\put{$360$} [c] at 10 -2
\endpicture
\end{minipage}
\begin{minipage}{4cm}
\beginpicture
\setcoordinatesystem units <1.5mm,2mm>
\setplotarea x from 0 to 16, y from -2 to 15
\put{216)} [l] at 0 12
\put {$ \scriptstyle \bullet$} [c] at 6 0
\put {$ \scriptstyle \bullet$} [c] at 6 6
\put {$ \scriptstyle \bullet$} [c] at 6  12
\put {$ \scriptstyle \bullet$} [c] at 12 0
\put {$ \scriptstyle \bullet$} [c] at 12 12
\put {$ \scriptstyle \bullet$} [c] at 14 0
\setlinear \plot  6 0  6 12 12 0 12 12 6 6   /
\put{$360$} [c] at 10 -2
\endpicture
\end{minipage}
$$
$$
\begin{minipage}{4cm}
\beginpicture
\setcoordinatesystem units <1.5mm,2mm>
\setplotarea x from 0 to 16, y from -2 to 15
\put{217)} [l] at 0 12
\put {$ \scriptstyle \bullet$} [c] at 6 0
\put {$ \scriptstyle \bullet$} [c] at 6 6
\put {$ \scriptstyle \bullet$} [c] at 6  12
\put {$ \scriptstyle \bullet$} [c] at 12 0
\put {$ \scriptstyle \bullet$} [c] at 12 12
\put {$ \scriptstyle \bullet$} [c] at 14 0
\setlinear \plot  6 12  6 0 12 12 12 0 6 6   /
\put{$360$} [c] at 10 -2
\endpicture
\end{minipage}
\begin{minipage}{4cm}
\beginpicture
\setcoordinatesystem units <1.5mm,2mm>
\setplotarea x from 0 to 16, y from -2 to 15
\put{218)} [l] at 0 12
\put {$ \scriptstyle \bullet$} [c] at 6 6
\put {$ \scriptstyle \bullet$} [c] at 9 0
\put {$ \scriptstyle \bullet$} [c] at 9 12
\put {$ \scriptstyle \bullet$} [c] at 12 6
\put {$ \scriptstyle \bullet$} [c] at 14 12
\put {$ \scriptstyle \bullet$} [c] at 14 0
\setlinear \plot  9 0 6  6  9  12 12 6 9 0   /
\setlinear \plot  14  0  14 12    /
\put{$360$} [c] at 10 -2
\endpicture
\end{minipage}
\begin{minipage}{4cm}
\beginpicture
\setcoordinatesystem units <1.5mm,2mm>
\setplotarea x from 0 to 16, y from -2 to 15
\put{ ${\bf  19}$} [l] at 0 15
\put{219)} [l] at 0 12
\put {$ \scriptstyle \bullet$} [c] at 6 0
\put {$ \scriptstyle \bullet$} [c] at 6 12
\put {$ \scriptstyle \bullet$} [c] at 10 12
\put {$ \scriptstyle \bullet$} [c] at 12 6
\put {$ \scriptstyle \bullet$} [c] at 14 0
\put {$ \scriptstyle \bullet$} [c] at 14 12
\setlinear \plot 6 12  6 0 10 12 14 0  14 12 /
\put{$720$} [c] at 10 -2
\endpicture
\end{minipage}
\begin{minipage}{4cm}
\beginpicture
\setcoordinatesystem units <1.5mm,2mm>
\setplotarea x from 0 to 16, y from -2 to 15
\put{220)} [l] at 0 12
\put {$ \scriptstyle \bullet$} [c] at 6 0
\put {$ \scriptstyle \bullet$} [c] at 6 12
\put {$ \scriptstyle \bullet$} [c] at 10 0
\put {$ \scriptstyle \bullet$} [c] at 12 6
\put {$ \scriptstyle \bullet$} [c] at 14 0
\put {$ \scriptstyle \bullet$} [c] at 14 12
\setlinear \plot 6 0  6 12 10 0 14 12  14 0 /
\put{$720$} [c] at 10 -2
\endpicture
\end{minipage}
\begin{minipage}{4cm}
\beginpicture
\setcoordinatesystem units <1.5mm,2mm>
\setplotarea x from 0 to 16, y from -2 to 15
\put{221)} [l] at 0 12
\put {$ \scriptstyle \bullet$} [c] at 6 6
\put {$ \scriptstyle \bullet$} [c] at 6 12
\put {$ \scriptstyle \bullet$} [c] at  10 12
\put {$ \scriptstyle \bullet$} [c] at 10 0
\put {$ \scriptstyle \bullet$} [c] at 14 6
\put {$ \scriptstyle \bullet$} [c] at 14 12
\setlinear \plot 6 12 6 6 10 0 14 6 14 12  /
\setlinear \plot  10 0 10 12  /
\put{$360$} [c] at 10 -2
\endpicture
\end{minipage}
\begin{minipage}{4cm}
\beginpicture
\setcoordinatesystem units <1.5mm,2mm>
\setplotarea x from 0 to 16, y from -2 to 15
\put{222)} [l] at 0 12
\put {$ \scriptstyle \bullet$} [c] at 6 6
\put {$ \scriptstyle \bullet$} [c] at 6 0
\put {$ \scriptstyle \bullet$} [c] at  10 12
\put {$ \scriptstyle \bullet$} [c] at 10 0
\put {$ \scriptstyle \bullet$} [c] at 14 6
\put {$ \scriptstyle \bullet$} [c] at 14 0
\setlinear \plot 6 0 6 6 10 12 14 6 14 0  /
\setlinear \plot  10 0 10 12  /
\put{$360$} [c] at 10 -2
\endpicture
\end{minipage}
$$
$$
\begin{minipage}{4cm}
\beginpicture
\setcoordinatesystem units <1.5mm,2mm>
\setplotarea x from 0 to 16, y from -2 to 15
\put{223)} [l] at 0 12
\put {$ \scriptstyle \bullet$} [c] at 6 0
\put {$ \scriptstyle \bullet$} [c] at 6 12
\put {$ \scriptstyle \bullet$} [c] at  9 12
\put {$ \scriptstyle \bullet$} [c] at 14 0
\put {$ \scriptstyle \bullet$} [c] at 14 6
\put {$ \scriptstyle \bullet$} [c] at 14 12
\setlinear \plot 6 12 6 0 14 12 14 0 /
\setlinear \plot  6 0 9 12  /
\put{$360$} [c] at 10 -2
\endpicture
\end{minipage}
\begin{minipage}{4cm}
\beginpicture
\setcoordinatesystem units <1.5mm,2mm>
\setplotarea x from 0 to 16, y from -2 to 15
\put{224)} [l] at 0 12
\put {$ \scriptstyle \bullet$} [c] at 6 0
\put {$ \scriptstyle \bullet$} [c] at 6 12
\put {$ \scriptstyle \bullet$} [c] at  9 0
\put {$ \scriptstyle \bullet$} [c] at 14 0
\put {$ \scriptstyle \bullet$} [c] at 14 6
\put {$ \scriptstyle \bullet$} [c] at 14 12
\setlinear \plot 6 0 6 12 14 0 14 12 /
\setlinear \plot  6 12 9 0  /
\put{$360$} [c] at 10 -2
\endpicture
\end{minipage}
\begin{minipage}{4cm}
\beginpicture
\setcoordinatesystem units <1.5mm,2mm>
\setplotarea x from 0 to 16, y from -2 to 15
\put{225)} [l] at 0 12
\put {$ \scriptstyle \bullet$} [c] at 6 0
\put {$ \scriptstyle \bullet$} [c] at 6 12
\put {$ \scriptstyle \bullet$} [c] at 12 0
\put {$ \scriptstyle \bullet$} [c] at 12 12
\put {$ \scriptstyle \bullet$} [c] at 14 0
\put {$ \scriptstyle \bullet$} [c] at 14 12
\setlinear \plot 14  0 14  12 12 0  12 12 6 0 6 12 12 0 /
\put{$360$} [c] at 10 -2
\endpicture
\end{minipage}
\begin{minipage}{4cm}
\beginpicture
\setcoordinatesystem units <1.5mm,2mm>
\setplotarea x from 0 to 16, y from -2 to 15
\put{226)} [l] at 0 12
\put {$ \scriptstyle \bullet$} [c] at 6 0
\put {$ \scriptstyle \bullet$} [c] at 6 12
\put {$ \scriptstyle \bullet$} [c] at 12 0
\put {$ \scriptstyle \bullet$} [c] at 12 12
\put {$ \scriptstyle \bullet$} [c] at 14 0
\put {$ \scriptstyle \bullet$} [c] at 14 12
\setlinear \plot 14  12 14  0 12 12  12 0 6 12 6 0 12 12 /
\put{$360$} [c] at 10 -2
\endpicture
\end{minipage}
\begin{minipage}{4cm}
\beginpicture
\setcoordinatesystem units <1.5mm,2mm>
\setplotarea x from 0 to 16, y from -2 to 15
\put{227)} [l] at 0 12
\put {$ \scriptstyle \bullet$} [c] at 6 6
\put {$ \scriptstyle \bullet$} [c] at 8 0
\put {$ \scriptstyle \bullet$} [c] at  8 6
\put {$ \scriptstyle \bullet$} [c] at 8 12
\put {$ \scriptstyle \bullet$} [c] at 10 6
\put {$ \scriptstyle \bullet$} [c] at 14 12
\setlinear \plot 14 12 8 0 6 6 8 12 10 6 8 0 8 12 /
\put{$120$} [c] at 10 -2
\endpicture
\end{minipage}
\begin{minipage}{4cm}
\beginpicture
\setcoordinatesystem units <1.5mm,2mm>
\setplotarea x from 0 to 16, y from -2 to 15
\put{228)} [l] at 0 12
\put {$ \scriptstyle \bullet$} [c] at 6 6
\put {$ \scriptstyle \bullet$} [c] at 8 0
\put {$ \scriptstyle \bullet$} [c] at  8 6
\put {$ \scriptstyle \bullet$} [c] at 8 12
\put {$ \scriptstyle \bullet$} [c] at 10 6
\put {$ \scriptstyle \bullet$} [c] at 14 0
\setlinear \plot 14 0 8 12 6 6 8 0 10 6 8 12 8 0 /
\put{$120$} [c] at 10 -2
\endpicture
\end{minipage}
$$
$$
\begin{minipage}{4cm}
\beginpicture
\setcoordinatesystem units <1.5mm,2mm>
\setplotarea x from 0 to 16, y from -2 to 15
\put{229)} [l] at 0 12
\put {$ \scriptstyle \bullet$} [c] at 6 12
\put {$ \scriptstyle \bullet$} [c] at 10 12
\put {$ \scriptstyle \bullet$} [c] at 12 12
\put {$ \scriptstyle \bullet$} [c] at 14 12
\put {$ \scriptstyle \bullet$} [c] at 12 0
\put {$ \scriptstyle \bullet$} [c] at 12 6
\setlinear \plot 6 12  12 0  12 12 /
\setlinear \plot 10 12  12 6  14 12 /
\put{$120$} [c] at 10 -2
\endpicture
\end{minipage}
\begin{minipage}{4cm}
\beginpicture
\setcoordinatesystem units <1.5mm,2mm>
\setplotarea x from 0 to 16, y from -2 to 15
\put{230)} [l] at 0 12
\put {$ \scriptstyle \bullet$} [c] at 6 0
\put {$ \scriptstyle \bullet$} [c] at 10 0
\put {$ \scriptstyle \bullet$} [c] at 12 0
\put {$ \scriptstyle \bullet$} [c] at 14 0
\put {$ \scriptstyle \bullet$} [c] at 12 12
\put {$ \scriptstyle \bullet$} [c] at 12 6
\setlinear \plot 6 0  12 12  12 0 /
\setlinear \plot 10 0  12 6  14 0 /
\put{$120$} [c] at 10 -2
\endpicture
\end{minipage}
\begin{minipage}{4cm}
\beginpicture
\setcoordinatesystem units <1.5mm,2mm>
\setplotarea x from 0 to 16, y from -2 to 15
\put{231)} [l] at 0 12
\put {$ \scriptstyle \bullet$} [c] at 6 12
\put {$ \scriptstyle \bullet$} [c] at 8 12
\put {$ \scriptstyle \bullet$} [c] at 12 12
\put {$ \scriptstyle \bullet$} [c] at 14 12
\put {$ \scriptstyle \bullet$} [c] at 6 0
\put {$ \scriptstyle \bullet$} [c] at 14 0
\setlinear \plot 6 12 6 0 8 12 14 0 14 12 6 0 12 12 /
\setlinear \plot  6 12 14 0 12 12 /
\put{$15$} [c] at 10 -2
\endpicture
\end{minipage}
\begin{minipage}{4cm}
\beginpicture
\setcoordinatesystem units <1.5mm,2mm>
\setplotarea x from 0 to 16, y from -2 to 15
\put{232)} [l] at 0 12
\put {$ \scriptstyle \bullet$} [c] at 6 0
\put {$ \scriptstyle \bullet$} [c] at 8 0
\put {$ \scriptstyle \bullet$} [c] at 12 0
\put {$ \scriptstyle \bullet$} [c] at 14 0
\put {$ \scriptstyle \bullet$} [c] at 6 12
\put {$ \scriptstyle \bullet$} [c] at 14 12
\setlinear \plot 6 0 6 12 8 0 14 12 14 0 6 12 12 0 /
\setlinear \plot  6 0 14 12 12 0 /
\put{$15$} [c] at 10 -2
\endpicture
\end{minipage}
\begin{minipage}{4cm}
\beginpicture
\setcoordinatesystem units <1.5mm,2mm>
\setplotarea x from 0 to 16, y from -2 to 15
\put{ ${\bf  20}$} [l] at 0 15
\put{233)} [l] at 0 12
\put {$ \scriptstyle \bullet$} [c] at 6 0
\put {$ \scriptstyle \bullet$} [c] at 6 12
\put {$ \scriptstyle \bullet$} [c] at 10 0
\put {$ \scriptstyle \bullet$} [c] at 14 0
\put {$ \scriptstyle \bullet$} [c] at 14 6
\put {$ \scriptstyle \bullet$} [c] at 14 12
\setlinear \plot 6 12 6 0 14 12  14 0 /
\setlinear \plot 10 0 14 12 /
\put{$720$} [c] at 10 -2
\endpicture
\end{minipage}
\begin{minipage}{4cm}
\beginpicture
\setcoordinatesystem units <1.5mm,2mm>
\setplotarea x from 0 to 16, y from -2 to 15
\put{234)} [l] at 0 12
\put {$ \scriptstyle \bullet$} [c] at 6 0
\put {$ \scriptstyle \bullet$} [c] at 6 12
\put {$ \scriptstyle \bullet$} [c] at 10 12
\put {$ \scriptstyle \bullet$} [c] at 14 0
\put {$ \scriptstyle \bullet$} [c] at 14 6
\put {$ \scriptstyle \bullet$} [c] at 14 12
\setlinear \plot 6 0 6 12 14 0  14 12 /
\setlinear \plot 10 12 14 0 /
\put{$720$} [c] at 10 -2
\endpicture
\end{minipage}
$$
$$
\begin{minipage}{4cm}
\beginpicture
\setcoordinatesystem units <1.5mm,2mm>
\setplotarea x from 0 to 16, y from -2 to 15
\put{235)} [l] at 0 12
\put {$ \scriptstyle \bullet$} [c] at 6 0
\put {$ \scriptstyle \bullet$} [c] at 9 0
\put {$ \scriptstyle \bullet$} [c] at 9 12
\put {$ \scriptstyle \bullet$} [c] at 11 12
\put {$ \scriptstyle \bullet$} [c] at 11 0
\put {$ \scriptstyle \bullet$} [c] at 14 12
\setlinear \plot 6 0 9 12 9 0 11 12  11 0 9 12  /
\setlinear \plot  11 0 14 12  /
\put{$720$} [c] at 10 -2
\endpicture
\end{minipage}
\begin{minipage}{4cm}
\beginpicture
\setcoordinatesystem units <1.5mm,2mm>
\setplotarea x from 0 to 16, y from -2 to 15
\put{236)} [l] at 0 12
\put {$ \scriptstyle \bullet$} [c] at 6 12
\put {$ \scriptstyle \bullet$} [c] at 9 6
\put {$ \scriptstyle \bullet$} [c] at 10 0
\put {$ \scriptstyle \bullet$} [c] at 10 12
\put {$ \scriptstyle \bullet$} [c] at 11 6
\put {$ \scriptstyle \bullet$} [c] at 14 0
\setlinear \plot 6 12 10 0 9 6 10  12 11 6 10 0 /
\setlinear \plot  10 12 14 0 /
\put{$360$} [c] at 10 -2
\endpicture
\end{minipage}
\begin{minipage}{4cm}
\beginpicture
\setcoordinatesystem units <1.5mm,2mm>
\setplotarea x from 0 to 16, y from -2 to 15
\put{237)} [l] at 0 12
\put {$ \scriptstyle \bullet$} [c] at 6 0
\put {$ \scriptstyle \bullet$} [c] at 8 12
\put {$ \scriptstyle \bullet$} [c] at 8 0
\put {$ \scriptstyle \bullet$} [c] at 11 0
\put {$ \scriptstyle \bullet$} [c] at 14 12
\put {$ \scriptstyle \bullet$} [c] at 14 0
\setlinear \plot 6 0 8 12 8 0 14 12 14 0 8 12 /
\setlinear \plot 8 12 11  0 14 12    /
\put{$120$} [c] at 10 -2
\endpicture
\end{minipage}
\begin{minipage}{4cm}
\beginpicture
\setcoordinatesystem units <1.5mm,2mm>
\setplotarea x from 0 to 16, y from -2 to 15
\put{238)} [l] at 0 12
\put {$ \scriptstyle \bullet$} [c] at 6 12
\put {$ \scriptstyle \bullet$} [c] at 8 12
\put {$ \scriptstyle \bullet$} [c] at 8 0
\put {$ \scriptstyle \bullet$} [c] at 11 12
\put {$ \scriptstyle \bullet$} [c] at 14 12
\put {$ \scriptstyle \bullet$} [c] at 14 0
\setlinear \plot 6 12 8 0 8 12 14 0 14 12 8 0 /
\setlinear \plot 8 0 11 12 14 0    /
\put{$120$} [c] at 10 -2
\endpicture
\end{minipage}
\begin{minipage}{4cm}
\beginpicture
\setcoordinatesystem units <1.5mm,2mm>
\setplotarea x from 0 to 16, y from -2 to 15
\put{239)} [l] at 0 12
\put {$ \scriptstyle \bullet$} [c] at 6 0
\put {$ \scriptstyle \bullet$} [c] at 6 12
\put {$ \scriptstyle \bullet$} [c] at 9 12
\put {$ \scriptstyle \bullet$} [c] at 9 6
\put {$ \scriptstyle \bullet$} [c] at 12 0
\put {$ \scriptstyle \bullet$} [c] at 14 0
\setlinear \plot 6 12 6 0 9 6 9 12  /
\setlinear \plot 9 6 12 0  /
\put{$720$} [c] at 10 -2
\endpicture
\end{minipage}
\begin{minipage}{4cm}
\beginpicture
\setcoordinatesystem units <1.5mm,2mm>
\setplotarea x from 0 to 16, y from -2 to 12
\put{240)} [l] at 0 12
\put {$ \scriptstyle \bullet$} [c] at 6 0
\put {$ \scriptstyle \bullet$} [c] at 6 12
\put {$ \scriptstyle \bullet$} [c] at 9 0
\put {$ \scriptstyle \bullet$} [c] at 9 6
\put {$ \scriptstyle \bullet$} [c] at 12 12
\put {$ \scriptstyle \bullet$} [c] at 14 0
\setlinear \plot 6 0 6 12 9 6 9 0  /
\setlinear \plot 9 6 12 12  /
\put{$720$} [c] at 10 -2
\endpicture
\end{minipage}
$$

$$
\begin{minipage}{4cm}
\beginpicture
\setcoordinatesystem units <1.5mm,2mm>
\setplotarea x from 0 to 16, y from -2 to 15
\put{241)} [l] at 0 12
\put {$ \scriptstyle \bullet$} [c] at 6 0
\put {$ \scriptstyle \bullet$} [c] at 6 6
\put {$ \scriptstyle \bullet$} [c] at 6 12
\put {$ \scriptstyle \bullet$} [c] at 12 0
\put {$ \scriptstyle \bullet$} [c] at 12 12
\put {$ \scriptstyle \bullet$} [c] at 14 0
\setlinear \plot 6 12 6 0 12  12 12  0 6 12 /
\put{$720$} [c] at 10 -2
\endpicture
\end{minipage}
\begin{minipage}{4cm}
\beginpicture
\setcoordinatesystem units <1.5mm,2mm>
\setplotarea x from 0 to 16, y from -2 to 15
\put{242)} [l] at 0 12
\put {$ \scriptstyle \bullet$} [c] at 6 12
\put {$ \scriptstyle \bullet$} [c] at 9 0
\put {$ \scriptstyle \bullet$} [c] at 12 12
\put {$ \scriptstyle \bullet$} [c] at 14 0
\put {$ \scriptstyle \bullet$} [c] at 14 6
\put {$ \scriptstyle \bullet$} [c] at 14 12
\setlinear \plot 6 12 9 0 12 12 /
\setlinear \plot 14 12 14 0 /
\put{$360$} [c] at 10 -2
\endpicture
\end{minipage}
\begin{minipage}{4cm}
\beginpicture
\setcoordinatesystem units <1.5mm,2mm>
\setplotarea x from 0 to 16, y from -2 to 15
\put{243)} [l] at 0 12
\put {$ \scriptstyle \bullet$} [c] at 6 0
\put {$ \scriptstyle \bullet$} [c] at 9 12
\put {$ \scriptstyle \bullet$} [c] at 12 0
\put {$ \scriptstyle \bullet$} [c] at 14 0
\put {$ \scriptstyle \bullet$} [c] at 14 6
\put {$ \scriptstyle \bullet$} [c] at 14 12
\setlinear \plot 6 0 9 12 12 0 /
\setlinear \plot 14 12 14 0 /
\put{$360$} [c] at 10 -2
\endpicture
\end{minipage}
\begin{minipage}{4cm}
\beginpicture
\setcoordinatesystem units <1.5mm,2mm>
\setplotarea x from 0 to 16, y from -2 to 15
\put{244)} [l] at 0 12
\put {$ \scriptstyle \bullet$} [c] at 6 0
\put {$ \scriptstyle \bullet$} [c] at 6 6
\put {$ \scriptstyle \bullet$} [c] at 9 12
\put {$ \scriptstyle \bullet$} [c] at 12 6
\put {$ \scriptstyle \bullet$} [c] at 12 0
\put {$ \scriptstyle \bullet$} [c] at 14 0
\setlinear \plot 6 0 6 6 9 12 12 6 12 0  /
\put{$360$} [c] at 10 -2
\endpicture
\end{minipage}
\begin{minipage}{4cm}
\beginpicture
\setcoordinatesystem units <1.5mm,2mm>
\setplotarea x from 0 to 16, y from -2 to 15
\put{245)} [l] at 0 12
\put {$ \scriptstyle \bullet$} [c] at 6 6
\put {$ \scriptstyle \bullet$} [c] at 6 12
\put {$ \scriptstyle \bullet$} [c] at 9 0
\put {$ \scriptstyle \bullet$} [c] at 12 12
\put {$ \scriptstyle \bullet$} [c] at 12 6
\put {$ \scriptstyle \bullet$} [c] at 14 0
\setlinear \plot 6 12 6 6 9 0 12 6 12 12  /
\put{$360$} [c] at 10 -2
\endpicture
\end{minipage}
\begin{minipage}{4cm}
\beginpicture
\setcoordinatesystem units <1.5mm,2mm>
\setplotarea x from 0 to 16, y from -2 to 15
\put{246)} [l] at 0 12
\put {$ \scriptstyle \bullet$} [c] at 6 0
\put {$ \scriptstyle \bullet$} [c] at 6 4
\put {$ \scriptstyle \bullet$} [c] at 6 8
\put {$ \scriptstyle \bullet$} [c] at 6 12
\put {$ \scriptstyle \bullet$} [c] at 10 0
\put {$ \scriptstyle \bullet$} [c] at 14 0
\setlinear \plot 6 0 6 12 /
\put{$360$} [c] at 10 -2
\endpicture
\end{minipage}
$$
$$
\begin{minipage}{4cm}
\beginpicture
\setcoordinatesystem units <1.5mm,2mm>
\setplotarea x from 0 to 16, y from -2 to 15
\put{247)} [l] at 0 12
\put {$ \scriptstyle \bullet$} [c] at 6 0
\put {$ \scriptstyle \bullet$} [c] at 9 6
\put {$ \scriptstyle \bullet$} [c] at 12 0
\put {$ \scriptstyle \bullet$} [c] at 9 0
\put {$ \scriptstyle \bullet$} [c] at 9 12
\put {$ \scriptstyle \bullet$} [c] at 14 0
\setlinear \plot 6 0  9 6 12 0  /
\setlinear \plot 9 0  9  12   /
\put{$120$} [c] at 10 -2
\endpicture
\end{minipage}
\begin{minipage}{4cm}
\beginpicture
\setcoordinatesystem units <1.5mm,2mm>
\setplotarea x from 0 to 16, y from -2 to 15
\put{248)} [l] at 0 12
\put {$ \scriptstyle \bullet$} [c] at 6 12
\put {$ \scriptstyle \bullet$} [c] at 9 6
\put {$ \scriptstyle \bullet$} [c] at 12 12
\put {$ \scriptstyle \bullet$} [c] at 9 0
\put {$ \scriptstyle \bullet$} [c] at 9 12
\put {$ \scriptstyle \bullet$} [c] at 14 0
\setlinear \plot 6 12  9 6 12 12  /
\setlinear \plot 9 0  9  12   /
\put{$120$} [c] at 10 -2
\endpicture
\end{minipage}
\begin{minipage}{4cm}
\beginpicture
\setcoordinatesystem units <1.5mm,2mm>
\setplotarea x from 0 to 16, y from -2 to 15
\put{249)} [l] at 0 12
\put {$ \scriptstyle \bullet$} [c] at 6 6
\put {$ \scriptstyle \bullet$} [c] at 9 0
\put {$ \scriptstyle \bullet$} [c] at 9 6
\put {$ \scriptstyle \bullet$} [c] at 9 12
\put {$ \scriptstyle \bullet$} [c] at 12 6
\put {$ \scriptstyle \bullet$} [c] at 14 0
\setlinear \plot 9 0 6 6 9 12  12 6 9 0  /
\setlinear \plot 9 0  9  12   /
\put{$120$} [c] at 10 -2
\endpicture
\end{minipage}
\begin{minipage}{4cm}
\beginpicture
\setcoordinatesystem units <1.5mm,2mm>
\setplotarea x from 0 to 16, y from -2 to 15
\put{ ${\bf  21}$} [l] at 0 15

\put{250)} [l] at 0 12
\put {$ \scriptstyle \bullet$} [c] at 6 12
\put {$ \scriptstyle \bullet$} [c] at 6 0
\put {$ \scriptstyle \bullet$} [c] at 9 0
\put {$ \scriptstyle \bullet$} [c] at 11 12
\put {$ \scriptstyle \bullet$} [c] at 14 0
\put {$ \scriptstyle \bullet$} [c] at 14 12
\setlinear \plot 6 0 6 12 9 0 11  12 14 0 14 12 /
\put{$720$} [c] at 10 -2
\endpicture
\end{minipage}
\begin{minipage}{4cm}
\beginpicture
\setcoordinatesystem units <1.5mm,2mm>
\setplotarea x from 0 to 16, y from -2 to 15
\put{251)} [l] at 0 12
\put {$ \scriptstyle \bullet$} [c] at 6 0
\put {$ \scriptstyle \bullet$} [c] at 9 6
\put {$ \scriptstyle \bullet$} [c] at 10 12
\put {$ \scriptstyle \bullet$} [c] at 10 0
\put {$ \scriptstyle \bullet$} [c] at 11 6
\put {$ \scriptstyle \bullet$} [c] at 14 0
\setlinear \plot 6 0 10 12 9 6 10 0 11 6 10 12 14 0 /
\put{$180$} [c] at 10 -2
\endpicture
\end{minipage}
\begin{minipage}{4cm}
\beginpicture
\setcoordinatesystem units <1.5mm,2mm>
\setplotarea x from 0 to 16, y from -2 to 15
\put{252)} [l] at 0 12
\put {$ \scriptstyle \bullet$} [c] at 6 12
\put {$ \scriptstyle \bullet$} [c] at 9 6
\put {$ \scriptstyle \bullet$} [c] at 10 12
\put {$ \scriptstyle \bullet$} [c] at 10 0
\put {$ \scriptstyle \bullet$} [c] at 11 6
\put {$ \scriptstyle \bullet$} [c] at 14 12
\setlinear \plot 6 12 10 0 9 6 10 12 11 6 10 0 14 12 /
\put{$180$} [c] at 10 -2
\endpicture
\end{minipage}
$$
$$
\begin{minipage}{4cm}
\beginpicture
\setcoordinatesystem units <1.5mm,2mm>
\setplotarea x from 0 to 16, y from -2 to 15
\put{253)} [l] at 0 12
\put {$ \scriptstyle \bullet$} [c] at 6 0
\put {$ \scriptstyle \bullet$} [c] at 8 12
\put {$ \scriptstyle \bullet$} [c] at 8 0
\put {$ \scriptstyle \bullet$} [c] at 12 0
\put {$ \scriptstyle \bullet$} [c] at 12 12
\put {$ \scriptstyle \bullet$} [c] at 14 0
\setlinear \plot 6 0 8 12 12 0 12 12 14 0 /
\setlinear \plot 8 12 8  0 12 12    /
\put{$180$} [c] at 10 -2
\endpicture
\end{minipage}
\begin{minipage}{4cm}
\beginpicture
\setcoordinatesystem units <1.5mm,2mm>
\setplotarea x from 0 to 16, y from -2 to 15
\put{254)} [l] at 0 12
\put {$ \scriptstyle \bullet$} [c] at 6 12
\put {$ \scriptstyle \bullet$} [c] at 8 12
\put {$ \scriptstyle \bullet$} [c] at 8 0
\put {$ \scriptstyle \bullet$} [c] at 12 0
\put {$ \scriptstyle \bullet$} [c] at 12 12
\put {$ \scriptstyle \bullet$} [c] at 14 12
\setlinear \plot 6 12 8 0 12 12 12 0 14 12 /
\setlinear \plot 8 0 8  12 12 0    /
\put{$180$} [c] at 10 -2
\endpicture
\end{minipage}
\begin{minipage}{4cm}
\beginpicture
\setcoordinatesystem units <1.5mm,2mm>
\setplotarea x from 0 to 16, y from -2 to 15
\put{255)} [l] at 0 12
\put {$ \scriptstyle \bullet$} [c] at 6 12
\put {$ \scriptstyle \bullet$} [c] at 8 12
\put {$ \scriptstyle \bullet$} [c] at 12 12
\put {$ \scriptstyle \bullet$} [c] at 14 12
\put {$ \scriptstyle \bullet$} [c] at 8 0
\put {$ \scriptstyle \bullet$} [c] at 13 6
\setlinear \plot 6 12 8 0 13 6 14 12   /
\setlinear \plot 13 6 12  12   /
\setlinear \plot 8 0 8 12   /
\put{$180$} [c] at 10 -2
\endpicture
\end{minipage}
\begin{minipage}{4cm}
\beginpicture
\setcoordinatesystem units <1.5mm,2mm>
\setplotarea x from 0 to 16, y from -2 to 15
\put{256)} [l] at 0 12
\put {$ \scriptstyle \bullet$} [c] at 6 0
\put {$ \scriptstyle \bullet$} [c] at 8 12
\put {$ \scriptstyle \bullet$} [c] at 12 0
\put {$ \scriptstyle \bullet$} [c] at 14 0
\put {$ \scriptstyle \bullet$} [c] at 8 0
\put {$ \scriptstyle \bullet$} [c] at 13 6
\setlinear \plot 6 0 8 12 13 6 14 0   /
\setlinear \plot 13 6 12  0   /
\setlinear \plot 8 0 8 12   /
\put{$180$} [c] at 10 -2
\endpicture
\end{minipage}
\begin{minipage}{4cm}
\beginpicture
\setcoordinatesystem units <1.5mm,2mm>
\setplotarea x from 0 to 16, y from -2 to 15
\put{257)} [l] at 0 12
\put {$ \scriptstyle \bullet$} [c] at 6 0
\put {$ \scriptstyle \bullet$} [c] at 12 0
\put {$ \scriptstyle \bullet$} [c] at 9 12
\put {$ \scriptstyle \bullet$} [c] at 7.5 6
\put {$ \scriptstyle \bullet$} [c] at 14 0
\put {$ \scriptstyle \bullet$} [c] at 14 12
\setlinear \plot 6 0 9 12 12 0   /
\setlinear \plot  14 12 14 0   /
\put{$720$} [c] at 10 -2
\endpicture
\end{minipage}
\begin{minipage}{4cm}
\beginpicture
\setcoordinatesystem units <1.5mm,2mm>
\setplotarea x from 0 to 16, y from -2 to 15
\put{258)} [l] at 0 12
\put {$ \scriptstyle \bullet$} [c] at 6 12
\put {$ \scriptstyle \bullet$} [c] at 12 12
\put {$ \scriptstyle \bullet$} [c] at 9 0
\put {$ \scriptstyle \bullet$} [c] at 7.5 6
\put {$ \scriptstyle \bullet$} [c] at 14 0
\put {$ \scriptstyle \bullet$} [c] at 14 12
\setlinear \plot 6 12 9 0 12 12   /
\setlinear \plot  14 12 14 0   /
\put{$720$} [c] at 10 -2
\endpicture
\end{minipage}
$$
$$
\begin{minipage}{4cm}
\beginpicture
\setcoordinatesystem units <1.5mm,2mm>
\setplotarea x from 0 to 16, y from -2 to 15
\put{259)} [l] at 0 12
\put {$ \scriptstyle \bullet$} [c] at 6 12
\put {$ \scriptstyle \bullet$} [c] at 6 0
\put {$ \scriptstyle \bullet$} [c] at 12 12
\put {$ \scriptstyle \bullet$} [c] at 12 0
\put {$ \scriptstyle \bullet$} [c] at 14  0
\put {$ \scriptstyle \bullet$} [c] at 14 12
\setlinear \plot 6 12 6 0  12 12 12 0 6 12  /
\setlinear \plot 14 0 14 12    /
\put{$180$} [c] at 10 -2
\endpicture
\end{minipage}
\begin{minipage}{4cm}
\beginpicture
\setcoordinatesystem units <1.5mm,2mm>
\setplotarea x from 0 to 16, y from -2 to 15
\put{ ${\bf  22}$} [l] at 0 15
\put{260)} [l] at 0 12
\put {$ \scriptstyle \bullet$} [c] at 6 12
\put {$ \scriptstyle \bullet$} [c] at 10 12
\put {$ \scriptstyle \bullet$} [c] at 14 12
\put {$ \scriptstyle \bullet$} [c] at 10 0
\put {$ \scriptstyle \bullet$} [c] at 14 0
\put {$ \scriptstyle \bullet$} [c] at 12 6
\setlinear \plot 6 12 10 0 10  12  /
\setlinear \plot 10 0 14 12 14 0   /
\put{$360$} [c] at 10 -2
\endpicture
\end{minipage}
\begin{minipage}{4cm}
\beginpicture
\setcoordinatesystem units <1.5mm,2mm>
\setplotarea x from 0 to 16, y from -2 to 15
\put{261)} [l] at 0 12
\put {$ \scriptstyle \bullet$} [c] at 6 0
\put {$ \scriptstyle \bullet$} [c] at 10 12
\put {$ \scriptstyle \bullet$} [c] at 14 12
\put {$ \scriptstyle \bullet$} [c] at 10 0
\put {$ \scriptstyle \bullet$} [c] at 14 0
\put {$ \scriptstyle \bullet$} [c] at 12 6
\setlinear \plot 6 0 10 12 10  0  /
\setlinear \plot 10 12 14 0 14 12   /
\put{$360$} [c] at 10 -2
\endpicture
\end{minipage}
\begin{minipage}{4cm}
\beginpicture
\setcoordinatesystem units <1.5mm,2mm>
\setplotarea x from 0 to 16, y from -2 to 15
\put{262)} [l] at 0 12
\put {$ \scriptstyle \bullet$} [c] at 6 12
\put {$ \scriptstyle \bullet$} [c] at 6 0
\put {$ \scriptstyle \bullet$} [c] at 10 12
\put {$ \scriptstyle \bullet$} [c] at 10 0
\put {$ \scriptstyle \bullet$} [c] at 14  0
\put {$ \scriptstyle \bullet$} [c] at 14 12
\setlinear \plot 6 0 6 12  10 0 14 12 14 0  /
\setlinear \plot 10 0 10 12    /
\put{$360$} [c] at 10 -2
\endpicture
\end{minipage}
\begin{minipage}{4cm}
\beginpicture
\setcoordinatesystem units <1.5mm,2mm>
\setplotarea x from 0 to 16, y from -2 to 15
\put{263)} [l] at 0 12
\put {$ \scriptstyle \bullet$} [c] at 6 12
\put {$ \scriptstyle \bullet$} [c] at 6 0
\put {$ \scriptstyle \bullet$} [c] at 10 12
\put {$ \scriptstyle \bullet$} [c] at 10 0
\put {$ \scriptstyle \bullet$} [c] at 14  0
\put {$ \scriptstyle \bullet$} [c] at 14 12
\setlinear \plot 6 12 6 0  10 12 14 0 14 12  /
\setlinear \plot 10 0 10 12    /
\put{$360$} [c] at 10 -2
\endpicture
\end{minipage}
\begin{minipage}{4cm}
\beginpicture
\setcoordinatesystem units <1.5mm,2mm>
\setplotarea x from 0 to 16, y from -2 to 15
\put{264)} [l] at 0 12
\put {$ \scriptstyle \bullet$} [c] at 6 12
\put {$ \scriptstyle \bullet$} [c] at 8 12
\put {$ \scriptstyle \bullet$} [c] at 10 12
\put {$ \scriptstyle \bullet$} [c] at 14 12
\put {$ \scriptstyle \bullet$} [c] at 10 0
\put {$ \scriptstyle \bullet$} [c] at 14 0
\setlinear \plot  10 12 14 0 14 12 10 0 10 12 /
\setlinear \plot  6 12 10 0 8 12  /

\put{$180$} [c] at 10 -2
\endpicture
\end{minipage}
$$
$$
\begin{minipage}{4cm}
\beginpicture
\setcoordinatesystem units <1.5mm,2mm>
\setplotarea x from 0 to 16, y from -2 to 15
\put{265)} [l] at 0 12
\put {$ \scriptstyle \bullet$} [c] at 6 0
\put {$ \scriptstyle \bullet$} [c] at 8 0
\put {$ \scriptstyle \bullet$} [c] at 10 12
\put {$ \scriptstyle \bullet$} [c] at 14 12
\put {$ \scriptstyle \bullet$} [c] at 10 0
\put {$ \scriptstyle \bullet$} [c] at 14 0
\setlinear \plot  10 12 14 0 14 12 10 0 10 12 /
\setlinear \plot  6 0 10 12 8 0  /
\put{$180$} [c] at 10 -2
\endpicture
\end{minipage}
\begin{minipage}{4cm}
\beginpicture
\setcoordinatesystem units <1.5mm,2mm>
\setplotarea x from 0 to 16, y from -2 to 15
\put{266)} [l] at 0 12
\put {$ \scriptstyle \bullet$} [c] at 6 0
\put {$ \scriptstyle \bullet$} [c] at 6 6
\put {$ \scriptstyle \bullet$} [c] at 6 12
\put {$ \scriptstyle \bullet$} [c] at 12 12
\put {$ \scriptstyle \bullet$} [c] at 12 0
\put {$ \scriptstyle \bullet$} [c] at 14 0
\setlinear \plot 6 12 6 0  12 12 12 0  /
\put{$720$} [c] at 10 -2
\endpicture
\end{minipage}
\begin{minipage}{4cm}
\beginpicture
\setcoordinatesystem units <1.5mm,2mm>
\setplotarea x from 0 to 16, y from -2 to 15
\put{267)} [l] at 0 12
\put {$ \scriptstyle \bullet$} [c] at 6 0
\put {$ \scriptstyle \bullet$} [c] at 6 6
\put {$ \scriptstyle \bullet$} [c] at 6 12
\put {$ \scriptstyle \bullet$} [c] at 12 12
\put {$ \scriptstyle \bullet$} [c] at 12 0
\put {$ \scriptstyle \bullet$} [c] at 14 0
\setlinear \plot 6 0 6 12  12 0 12 12  /
\put{$720$} [c] at 10 -2
\endpicture
\end{minipage}
\begin{minipage}{4cm}
\beginpicture
\setcoordinatesystem units <1.5mm,2mm>
\setplotarea x from 0 to 16, y from -2 to 15
\put{268)} [l] at 0 12
\put {$ \scriptstyle \bullet$} [c] at 6 0
\put {$ \scriptstyle \bullet$} [c] at 8 0
\put {$ \scriptstyle \bullet$} [c] at 12 0
\put {$ \scriptstyle \bullet$} [c] at 14 0
\put {$ \scriptstyle \bullet$} [c] at 10 6
\put {$ \scriptstyle \bullet$} [c] at 10 12
\setlinear \plot 6 0 10  12 10 6 8 0 /
\setlinear \plot 10 6  12 0   /
\put{$360$} [c] at 10 -2
\endpicture
\end{minipage}
\begin{minipage}{4cm}
\beginpicture
\setcoordinatesystem units <1.5mm,2mm>
\setplotarea x from 0 to 16, y from -2 to 15
\put{269)} [l] at 0 12
\put {$ \scriptstyle \bullet$} [c] at 6 12
\put {$ \scriptstyle \bullet$} [c] at 8 12
\put {$ \scriptstyle \bullet$} [c] at 12 12
\put {$ \scriptstyle \bullet$} [c] at 14 0
\put {$ \scriptstyle \bullet$} [c] at 10 6
\put {$ \scriptstyle \bullet$} [c] at 10 0
\setlinear \plot 6 12 10  0 10 6 8 12 /
\setlinear \plot 10 6  12 12   /
\put{$360$} [c] at 10 -2
\endpicture
\end{minipage}
\begin{minipage}{4cm}
\beginpicture
\setcoordinatesystem units <1.5mm,2mm>
\setplotarea x from 0 to 16, y from -2 to 15
\put{270)} [l] at 0 12
\put {$ \scriptstyle \bullet$} [c] at 6 12
\put {$ \scriptstyle \bullet$} [c] at 6 0
\put {$ \scriptstyle \bullet$} [c] at 9 0
\put {$ \scriptstyle \bullet$} [c] at 12 0
\put {$ \scriptstyle \bullet$} [c] at 12 12
\put {$ \scriptstyle \bullet$} [c] at 14 0
\setlinear \plot 6 12 6 0 12 12 12 0 6 12 /
\setlinear \plot 6 12 9 0 12 12    /
\put{$60$} [c] at 10 -2
\endpicture
\end{minipage}
$$
$$
\begin{minipage}{4cm}
\beginpicture
\setcoordinatesystem units <1.5mm,2mm>
\setplotarea x from 0 to 16, y from -2 to 15
\put{271)} [l] at 0 12
\put {$ \scriptstyle \bullet$} [c] at 6 12
\put {$ \scriptstyle \bullet$} [c] at 6 0
\put {$ \scriptstyle \bullet$} [c] at 9 12
\put {$ \scriptstyle \bullet$} [c] at 12 0
\put {$ \scriptstyle \bullet$} [c] at 12 12
\put {$ \scriptstyle \bullet$} [c] at 14 0
\setlinear \plot 6 12 6 0 12 12 12 0 6 12 /
\setlinear \plot 6 0 9 12 12 0    /

\put{$60$} [c] at 10 -2
\endpicture
\end{minipage}
\begin{minipage}{4cm}
\beginpicture
\setcoordinatesystem units <1.5mm,2mm>
\setplotarea x from 0 to 16, y from -2 to 15
\put{272)} [l] at 0 12
\put {$ \scriptstyle \bullet$} [c] at 6 0
\put {$ \scriptstyle \bullet$} [c] at 11 6
\put {$ \scriptstyle \bullet$} [c] at 12 12
\put {$ \scriptstyle \bullet$} [c] at 12 0
\put {$ \scriptstyle \bullet$} [c] at 13 6
\put {$ \scriptstyle \bullet$} [c] at 14 0
\setlinear \plot 6 0 12  12 11 6 12 0 13 6 12 12 /
\put{$360$} [c] at 10 -2
\endpicture
\end{minipage}
\begin{minipage}{4cm}
\beginpicture
\setcoordinatesystem units <1.5mm,2mm>
\setplotarea x from 0 to 16, y from -2 to 15
\put{273)} [l] at 0 12
\put {$ \scriptstyle \bullet$} [c] at 6 12
\put {$ \scriptstyle \bullet$} [c] at 11 6
\put {$ \scriptstyle \bullet$} [c] at 12 12
\put {$ \scriptstyle \bullet$} [c] at 12 0
\put {$ \scriptstyle \bullet$} [c] at 13 6
\put {$ \scriptstyle \bullet$} [c] at 14 0
\setlinear \plot 6 12 12  0 11 6 12 12 13 6 12 0 /
\put{$360$} [c] at 10 -2
\endpicture
\end{minipage}
\begin{minipage}{4cm}
\beginpicture
\setcoordinatesystem units <1.5mm,2mm>
\setplotarea x from 0 to 16, y from -2 to 15
\put{ ${\bf  23}$} [l] at 0 15

\put{274)} [l] at 0 12
\put {$ \scriptstyle \bullet$} [c] at 6 12
\put {$ \scriptstyle \bullet$} [c] at 8 12
\put {$ \scriptstyle \bullet$} [c] at 10 12
\put {$ \scriptstyle \bullet$} [c] at 14 12
\put {$ \scriptstyle \bullet$} [c] at 8 0
\put {$ \scriptstyle \bullet$} [c] at 12 0
\setlinear \plot 6 12 8 0 10  12 12 0 14 12 /
\setlinear \plot 8 0 8 12    /

\put{$360$} [c] at 10 -2
\endpicture
\end{minipage}
\begin{minipage}{4cm}
\beginpicture
\setcoordinatesystem units <1.5mm,2mm>
\setplotarea x from 0 to 16, y from -2 to 15
\put{275)} [l] at 0 12
\put {$ \scriptstyle \bullet$} [c] at 6 0
\put {$ \scriptstyle \bullet$} [c] at 8 12
\put {$ \scriptstyle \bullet$} [c] at 10 0
\put {$ \scriptstyle \bullet$} [c] at 14 0
\put {$ \scriptstyle \bullet$} [c] at 8 0
\put {$ \scriptstyle \bullet$} [c] at 12 12
\setlinear \plot 6 0 8 12 10  0 12 12 14 0 /
\setlinear \plot 8 0 8 12    /
\put{$360$} [c] at 10 -2
\endpicture
\end{minipage}
\begin{minipage}{4cm}
\beginpicture
\setcoordinatesystem units <1.5mm,2mm>
\setplotarea x from 0 to 16, y from -2 to 15
\put{ ${\bf  24}$} [l] at 0 15
\put{276)} [l] at 0 12
\put {$ \scriptstyle \bullet$} [c] at 6 12
\put {$ \scriptstyle \bullet$} [c] at 8 12
\put {$ \scriptstyle \bullet$} [c] at 12 12
\put {$ \scriptstyle \bullet$} [c] at 14 0
\put {$ \scriptstyle \bullet$} [c] at 8 0
\put {$ \scriptstyle \bullet$} [c] at 12 0
\setlinear \plot 6 12 8 0 12  12 14 0  /
\setlinear \plot 8 0 8 12    /
\setlinear \plot 12 0 12 12    /
\put{$180$} [c] at 10 -2
\endpicture
\end{minipage}
$$
$$
\begin{minipage}{4cm}
\beginpicture
\setcoordinatesystem units <1.5mm,2mm>
\setplotarea x from 0 to 16, y from -2 to 15
\put{277)} [l] at 0 12
\put {$ \scriptstyle \bullet$} [c] at 6 0
\put {$ \scriptstyle \bullet$} [c] at 6 12
\put {$ \scriptstyle \bullet$} [c] at 9 6
\put {$ \scriptstyle \bullet$} [c] at 12 0
\put {$ \scriptstyle \bullet$} [c] at 12 12
\put {$ \scriptstyle \bullet$} [c] at 14 0
\setlinear \plot 6 12 6 0 12 12 12 0    /
\put{$720$} [c] at 10 -2
\endpicture
\end{minipage}
\begin{minipage}{4cm}
\beginpicture
\setcoordinatesystem units <1.5mm,2mm>
\setplotarea x from 0 to 16, y from -2 to 15
\put{278)} [l] at 0 12
\put {$ \scriptstyle \bullet$} [c] at 6 0
\put {$ \scriptstyle \bullet$} [c] at 6 6
\put {$ \scriptstyle \bullet$} [c] at 6 12
\put {$ \scriptstyle \bullet$} [c] at 10 0
\put {$ \scriptstyle \bullet$} [c] at 10 12
\put {$ \scriptstyle \bullet$} [c] at 14 0
\setlinear \plot 6 0 6 12     /
\setlinear \plot 10 0 10 12   /
\put{$720$} [c] at 10 -2
\endpicture
\end{minipage}
\begin{minipage}{4cm}
\beginpicture
\setcoordinatesystem units <1.5mm,2mm>
\setplotarea x from 0 to 16, y from -2 to 15
\put{279)} [l] at 0 12
\put {$ \scriptstyle \bullet$} [c] at 6 0
\put {$ \scriptstyle \bullet$} [c] at 6 12
\put {$ \scriptstyle \bullet$} [c] at 12 0
\put {$ \scriptstyle \bullet$} [c] at 12 12
\put {$ \scriptstyle \bullet$} [c] at 14 12
\put {$ \scriptstyle \bullet$} [c] at 14 0
\setlinear \plot 6 12 6 0 12 12 12 0    /
\setlinear \plot 14  0 14 12      /
\put{$720$} [c] at 10 -2
\endpicture
\end{minipage}
\begin{minipage}{4cm}
\beginpicture
\setcoordinatesystem units <1.5mm,2mm>
\setplotarea x from 0 to 16, y from -2 to 15
\put{280)} [l] at 0 12
\put {$ \scriptstyle \bullet$} [c] at 6 0
\put {$ \scriptstyle \bullet$} [c] at 8 12
\put {$ \scriptstyle \bullet$} [c] at 8 0
\put {$ \scriptstyle \bullet$} [c] at 12 12
\put {$ \scriptstyle \bullet$} [c] at 12 0
\put {$ \scriptstyle \bullet$} [c] at 14 0
\setlinear \plot 6 0 8 12 8 0 12 12 12 0 8 12    /
\put{$360$} [c] at 10 -2
\endpicture
\end{minipage}
\begin{minipage}{4cm}
\beginpicture
\setcoordinatesystem units <1.5mm,2mm>
\setplotarea x from 0 to 16, y from -2 to 15
\put{281)} [l] at 0 12
\put {$ \scriptstyle \bullet$} [c] at 6 12
\put {$ \scriptstyle \bullet$} [c] at 8 12
\put {$ \scriptstyle \bullet$} [c] at 8 0
\put {$ \scriptstyle \bullet$} [c] at 12 12
\put {$ \scriptstyle \bullet$} [c] at 12 0
\put {$ \scriptstyle \bullet$} [c] at 14 0
\setlinear \plot 6 12 8 0 8 12 12 0 12 12 8 0    /
\put{$360$} [c] at 10 -2
\endpicture
\end{minipage}
\begin{minipage}{4cm}
\beginpicture
\setcoordinatesystem units <1.5mm,2mm>
\setplotarea x from 0 to 16, y from -2 to 15
\put{282)} [l] at 0 12
\put {$ \scriptstyle \bullet$} [c] at 6 6
\put {$ \scriptstyle \bullet$} [c] at 8 0
\put {$ \scriptstyle \bullet$} [c] at 8 12
\put {$ \scriptstyle \bullet$} [c] at 10 6
\put {$ \scriptstyle \bullet$} [c] at 12 0
\put {$ \scriptstyle \bullet$} [c] at 14 0
\setlinear \plot 8 0 6 6 8 12 10  6 8 0   /
\put{$180$} [c] at 10 -2
\endpicture
\end{minipage}
$$
$$
\begin{minipage}{4cm}
\beginpicture
\setcoordinatesystem units <1.5mm,2mm>
\setplotarea x from 0 to 16, y from -2 to 15
\put{283)} [l] at 0 12
\put {$ \scriptstyle \bullet$} [c] at 6 12
\put {$ \scriptstyle \bullet$} [c] at 10  12
\put {$ \scriptstyle \bullet$} [c] at 8 6
\put {$ \scriptstyle \bullet$} [c] at 8  0
\put {$ \scriptstyle \bullet$} [c] at 12 0
\put {$ \scriptstyle \bullet$} [c] at 14 0
\setlinear \plot 6 12  8 6 8 0    /
\setlinear \plot 8 6 10  12      /
\put{$180$} [c] at 10 -2
\endpicture
\end{minipage}
\begin{minipage}{4cm}
\beginpicture
\setcoordinatesystem units <1.5mm,2mm>
\setplotarea x from 0 to 16, y from -2 to 15
\put{284)} [l] at 0 12
\put {$ \scriptstyle \bullet$} [c] at 6 0
\put {$ \scriptstyle \bullet$} [c] at 10  0
\put {$ \scriptstyle \bullet$} [c] at 8 6
\put {$ \scriptstyle \bullet$} [c] at 8  12
\put {$ \scriptstyle \bullet$} [c] at 12 0
\put {$ \scriptstyle \bullet$} [c] at 14 0
\setlinear \plot 6 0  8 6 8 12    /
\setlinear \plot 8 6 10  0      /
\put{$180$} [c] at 10 -2
\endpicture
\end{minipage}
\begin{minipage}{4cm}
\beginpicture
\setcoordinatesystem units <1.5mm,2mm>
\setplotarea x from 0 to 16, y from -2 to 15
\put{ ${\bf  25}$} [l] at 0 15

\put{285)} [l] at 0 12
\put {$ \scriptstyle \bullet$} [c] at 6 12
\put {$ \scriptstyle \bullet$} [c] at 9 12
\put {$ \scriptstyle \bullet$} [c] at 11 12
\put {$ \scriptstyle \bullet$} [c] at 14 12
\put {$ \scriptstyle \bullet$} [c] at 10 0
\put {$ \scriptstyle \bullet$} [c] at 12 6
\setlinear \plot 6 12 10 0 14  12   /
\setlinear \plot 9 12 10 0 11 12   /
\put{$120$} [c] at 10 -2
\endpicture
\end{minipage}
\begin{minipage}{4cm}
\beginpicture
\setcoordinatesystem units <1.5mm,2mm>
\setplotarea x from 0 to 16, y from -2 to 15
\put{286)} [l] at 0 12
\put {$ \scriptstyle \bullet$} [c] at 6 0
\put {$ \scriptstyle \bullet$} [c] at 9 0
\put {$ \scriptstyle \bullet$} [c] at 11 0
\put {$ \scriptstyle \bullet$} [c] at 14 0
\put {$ \scriptstyle \bullet$} [c] at 12 6
\put {$ \scriptstyle \bullet$} [c] at 10 12
\setlinear \plot 6 0 10 12 14  0   /
\setlinear \plot 9 0 10 12 11 0   /
\put{$120$} [c] at 10 -2
\endpicture
\end{minipage}
\begin{minipage}{4cm}
\beginpicture
\setcoordinatesystem units <1.5mm,2mm>
\setplotarea x from 0 to 16, y from -2 to 15
\put{287)} [l] at 0 12
\put {$ \scriptstyle \bullet$} [c] at 6 12
\put {$ \scriptstyle \bullet$} [c] at 7.5 0
\put {$ \scriptstyle \bullet$} [c] at 9 12
\put {$ \scriptstyle \bullet$} [c] at 11 12
\put {$ \scriptstyle \bullet$} [c] at 12.5 0
\put {$ \scriptstyle \bullet$} [c] at 14 12
\setlinear \plot 6 12 7.5 0 9 12    /
\setlinear \plot 11 12 12.5 0 14 12     /
\put{$90$} [c] at 10 -2
\endpicture
\end{minipage}
\begin{minipage}{4cm}
\beginpicture
\setcoordinatesystem units <1.5mm,2mm>
\setplotarea x from 0 to 16, y from -2 to 15
\put{288)} [l] at 0 12
\put {$ \scriptstyle \bullet$} [c] at 6 0
\put {$ \scriptstyle \bullet$} [c] at 7.5 12
\put {$ \scriptstyle \bullet$} [c] at 9 0
\put {$ \scriptstyle \bullet$} [c] at 11 0
\put {$ \scriptstyle \bullet$} [c] at 12.5 12
\put {$ \scriptstyle \bullet$} [c] at 14 0
\setlinear \plot 6 0 7.5 12 9 0    /
\setlinear \plot 11 0 12.5 12 14 0     /
\put{$90$} [c] at 10 -2
\endpicture
\end{minipage}
$$
$$
\begin{minipage}{4cm}
\beginpicture
\setcoordinatesystem units <1.5mm,2mm>
\setplotarea x from 0 to 16, y from -2 to 15
\put{289)} [l] at 0 12
\put {$ \scriptstyle \bullet$} [c] at 6 0
\put {$ \scriptstyle \bullet$} [c] at 7.5 12
\put {$ \scriptstyle \bullet$} [c] at 9 0
\put {$ \scriptstyle \bullet$} [c] at 11 12
\put {$ \scriptstyle \bullet$} [c] at 12.5 0
\put {$ \scriptstyle \bullet$} [c] at 14 12
\setlinear \plot 6 0 7.5 12 9 0    /
\setlinear \plot 11 12 12.5 0 14 12     /
\put{$180$} [c] at 10 -2
\endpicture
\end{minipage}
\begin{minipage}{4cm}
\beginpicture
\setcoordinatesystem units <1.5mm,2mm>
\setplotarea x from 0 to 16, y from -2 to 15
\put{ ${\bf  26}$} [l] at 0 15
\put{290)} [l] at 0 12
\put {$ \scriptstyle \bullet$} [c] at 6 12
\put {$ \scriptstyle \bullet$} [c] at 8 12
\put {$ \scriptstyle \bullet$} [c] at 12 12
\put {$ \scriptstyle \bullet$} [c] at 14 12
\put {$ \scriptstyle \bullet$} [c] at 14 0
\put {$ \scriptstyle \bullet$} [c] at 10 0
\setlinear \plot 6 12 10 0 14  12 14 0  /
\setlinear \plot 8 12 10 0 12 12   /
\put{$120$} [c] at 10 -2
\endpicture
\end{minipage}
\begin{minipage}{4cm}
\beginpicture
\setcoordinatesystem units <1.5mm,2mm>
\setplotarea x from 0 to 16, y from -2 to 15
\put{291)} [l] at 0 12
\put {$ \scriptstyle \bullet$} [c] at 6 0
\put {$ \scriptstyle \bullet$} [c] at 8 0
\put {$ \scriptstyle \bullet$} [c] at 12 0
\put {$ \scriptstyle \bullet$} [c] at 14 12
\put {$ \scriptstyle \bullet$} [c] at 14 0
\put {$ \scriptstyle \bullet$} [c] at 10 12
\setlinear \plot 6 0 10 12 14  0 14 12  /
\setlinear \plot 8 0 10 12 12 0   /
\put{$120$} [c] at 10 -2
\endpicture
\end{minipage}
\begin{minipage}{4cm}
\beginpicture
\setcoordinatesystem units <1.5mm,2mm>
\setplotarea x from 0 to 16, y from -2 to 15
\put{292)} [l] at 0 12
\put {$ \scriptstyle \bullet$} [c] at 6 12
\put {$ \scriptstyle \bullet$} [c] at 9 12
\put {$ \scriptstyle \bullet$} [c] at 9 6
\put {$ \scriptstyle \bullet$} [c] at 9  0
\put {$ \scriptstyle \bullet$} [c] at 12 12
\put {$ \scriptstyle \bullet$} [c] at 14 0
\setlinear \plot 6 12 9 0 12 12     /
\setlinear \plot 9  0 9 12     /
\put{$360$} [c] at 10 -2
\endpicture
\end{minipage}
\begin{minipage}{4cm}
\beginpicture
\setcoordinatesystem units <1.5mm,2mm>
\setplotarea x from 0 to 16, y from -2 to 15
\put{293)} [l] at 0 12
\put {$ \scriptstyle \bullet$} [c] at 6 0
\put {$ \scriptstyle \bullet$} [c] at 9 12
\put {$ \scriptstyle \bullet$} [c] at 9 6
\put {$ \scriptstyle \bullet$} [c] at 9  0
\put {$ \scriptstyle \bullet$} [c] at 12 0
\put {$ \scriptstyle \bullet$} [c] at 14 0
\setlinear \plot 6 0 9 12 12 0     /
\setlinear \plot 9  0 9 12     /
\put{$360$} [c] at 10 -2
\endpicture
\end{minipage}
\begin{minipage}{4cm}
\beginpicture
\setcoordinatesystem units <1.5mm,2mm>
\setplotarea x from 0 to 16, y from -2 to 15
\put{294)} [l] at 0 12
\put {$ \scriptstyle \bullet$} [c] at 6 12
\put {$ \scriptstyle \bullet$} [c] at 9 12
\put {$ \scriptstyle \bullet$} [c] at 12 12
\put {$ \scriptstyle \bullet$} [c] at 6 0
\put {$ \scriptstyle \bullet$} [c] at 12 0
\put {$ \scriptstyle \bullet$} [c] at 14 0
\setlinear \plot 6 12 6 0 9 12 12 0 12 12    /
\put{$360$} [c] at 10 -2
\endpicture
\end{minipage}
$$
$$
\begin{minipage}{4cm}
\beginpicture
\setcoordinatesystem units <1.5mm,2mm>
\setplotarea x from 0 to 16, y from -2 to 15
\put{295)} [l] at 0 12
\put {$ \scriptstyle \bullet$} [c] at 6 0
\put {$ \scriptstyle \bullet$} [c] at 9 0
\put {$ \scriptstyle \bullet$} [c] at 12 0
\put {$ \scriptstyle \bullet$} [c] at 6 12
\put {$ \scriptstyle \bullet$} [c] at 12 12
\put {$ \scriptstyle \bullet$} [c] at 14 0
\setlinear \plot 6 0 6 12 9 0 12 12 12 0    /
\put{$360$} [c] at 10 -2
\endpicture
\end{minipage}
\begin{minipage}{4cm}
\beginpicture
\setcoordinatesystem units <1.5mm,2mm>
\setplotarea x from 0 to 16, y from -2 to 15
\put{ ${\bf  27}$} [l] at 0 15
\put{296)} [l] at 0 12
\put {$ \scriptstyle \bullet$} [c] at 6 12
\put {$ \scriptstyle \bullet$} [c] at 8 12
\put {$ \scriptstyle \bullet$} [c] at 10 12
\put {$ \scriptstyle \bullet$} [c] at 8 0
\put {$ \scriptstyle \bullet$} [c] at 14 12
\put {$ \scriptstyle \bullet$} [c] at 14 0
\setlinear \plot 6 12 8 0 10 12    /
\setlinear \plot 8  0 8 12     /
\setlinear \plot 14  0 14 12     /
\put{$120$} [c] at 10 -2
\endpicture
\end{minipage}
\begin{minipage}{4cm}
\beginpicture
\setcoordinatesystem units <1.5mm,2mm>
\setplotarea x from 0 to 16, y from -2 to 15
\put{297)} [l] at 0 12
\put {$ \scriptstyle \bullet$} [c] at 6 0
\put {$ \scriptstyle \bullet$} [c] at 8 0
\put {$ \scriptstyle \bullet$} [c] at 10 0
\put {$ \scriptstyle \bullet$} [c] at 8 12
\put {$ \scriptstyle \bullet$} [c] at 14 12
\put {$ \scriptstyle \bullet$} [c] at 14 0
\setlinear \plot 6 0 8 12 10 0    /
\setlinear \plot 8  0 8 12     /
\setlinear \plot 14  0 14 12     /
\put{$120$} [c] at 10 -2
\endpicture
\end{minipage}
\begin{minipage}{4cm}
\beginpicture
\setcoordinatesystem units <1.5mm,2mm>
\setplotarea x from 0 to 16, y from -2 to 15
\put{298)} [l] at 0 12
\put {$ \scriptstyle \bullet$} [c] at 6 0
\put {$ \scriptstyle \bullet$} [c] at 6 12
\put {$ \scriptstyle \bullet$} [c] at 10 0
\put {$ \scriptstyle \bullet$} [c] at 10 12
\put {$ \scriptstyle \bullet$} [c] at 14 12
\put {$ \scriptstyle \bullet$} [c] at 14 0
\setlinear \plot 6 0 6 12   /
\setlinear \plot 10 0 10 12   /
\setlinear \plot 14 0 14 12   /
\put{$120$} [c] at 10 -2
\endpicture
\end{minipage}
\begin{minipage}{4cm}
\beginpicture
\setcoordinatesystem units <1.5mm,2mm>
\setplotarea x from 0 to 16, y from -2 to 15
\put{ ${\bf  28}$} [l] at 0 15

\put{299)} [l] at 0 12
\put {$ \scriptstyle \bullet$} [c] at 6 12
\put {$ \scriptstyle \bullet$} [c] at 8 12
\put {$ \scriptstyle \bullet$} [c] at 10 12
\put {$ \scriptstyle \bullet$} [c] at 8 0
\put {$ \scriptstyle \bullet$} [c] at 10 0
\put {$ \scriptstyle \bullet$} [c] at 14 0
\setlinear \plot 6 12 8 0 10 12 10 0   /
\setlinear \plot 8  0 8 12     /
\put{$360$} [c] at 10 -2
\endpicture
\end{minipage}
\begin{minipage}{4cm}
\beginpicture
\setcoordinatesystem units <1.5mm,2mm>
\setplotarea x from 0 to 16, y from -2 to 15
\put{300)} [l] at 0 12
\put {$ \scriptstyle \bullet$} [c] at 6 0
\put {$ \scriptstyle \bullet$} [c] at 8 0
\put {$ \scriptstyle \bullet$} [c] at 10 0
\put {$ \scriptstyle \bullet$} [c] at 8 12
\put {$ \scriptstyle \bullet$} [c] at 10 12
\put {$ \scriptstyle \bullet$} [c] at 14 0
\setlinear \plot 6 0 8 12 10 0 10 12   /
\setlinear \plot 8  0 8 12     /
\put{$360$} [c] at 10 -2
\endpicture
\end{minipage}
$$
$$
\begin{minipage}{4cm}
\beginpicture
\setcoordinatesystem units <1.5mm,2mm>
\setplotarea x from 0 to 16, y from -2 to 15
\put{301)} [l] at 0 12
\put {$ \scriptstyle \bullet$} [c] at 6 12
\put {$ \scriptstyle \bullet$} [c] at 7 6
\put {$ \scriptstyle \bullet$} [c] at 10 12
\put {$ \scriptstyle \bullet$} [c] at 8 0
\put {$ \scriptstyle \bullet$} [c] at 12 0
\put {$ \scriptstyle \bullet$} [c] at 14 0
\setlinear \plot 6 12 8 0 10 12    /
\put{$360$} [c] at 10 -2
\endpicture
\end{minipage}
\begin{minipage}{4cm}
\beginpicture
\setcoordinatesystem units <1.5mm,2mm>
\setplotarea x from 0 to 16, y from -2 to 15
\put{302)} [l] at 0 12
\put {$ \scriptstyle \bullet$} [c] at 6 0
\put {$ \scriptstyle \bullet$} [c] at 7 6
\put {$ \scriptstyle \bullet$} [c] at 10 0
\put {$ \scriptstyle \bullet$} [c] at 8 12
\put {$ \scriptstyle \bullet$} [c] at 12 0
\put {$ \scriptstyle \bullet$} [c] at 14 0
\setlinear \plot 6 0 8 12 10 0    /
\put{$360$} [c] at 10 -2
\endpicture
\end{minipage}
\begin{minipage}{4cm}
\beginpicture
\setcoordinatesystem units <1.5mm,2mm>
\setplotarea x from 0 to 16, y from -2 to 15
\put{303)} [l] at 0 12
\put {$ \scriptstyle \bullet$} [c] at 6 0
\put {$ \scriptstyle \bullet$} [c] at 6 12
\put {$ \scriptstyle \bullet$} [c] at 10 0
\put {$ \scriptstyle \bullet$} [c] at 10 12
\put {$ \scriptstyle \bullet$} [c] at 12 0
\put {$ \scriptstyle \bullet$} [c] at 14 0
\setlinear \plot 6 0 6 12 10 0 10 12 6 0   /
\put{$90$} [c] at 10 -2
\endpicture
\end{minipage}
\begin{minipage}{4cm}
\beginpicture
\setcoordinatesystem units <1.5mm,2mm>
\setplotarea x from 0 to 16, y from -2 to 15
\put{ ${\bf  30}$} [l] at 0 15

\put{304)} [l] at 0 12
\put {$ \scriptstyle \bullet$} [c] at 6 0
\put {$ \scriptstyle \bullet$} [c] at 8 12
\put {$ \scriptstyle \bullet$} [c] at 10 0
\put {$ \scriptstyle \bullet$} [c] at 12 12
\put {$ \scriptstyle \bullet$} [c] at 12 0
\put {$ \scriptstyle \bullet$} [c] at 14 0
\setlinear \plot 6 0 8 12 10 0    /
\setlinear \plot 12  0 12 12     /
\put{$360$} [c] at 10 -2
\endpicture
\end{minipage}
\begin{minipage}{4cm}
\beginpicture
\setcoordinatesystem units <1.5mm,2mm>
\setplotarea x from 0 to 16, y from -2 to 15
\put{305)} [l] at 0 12
\put {$ \scriptstyle \bullet$} [c] at 6 12
\put {$ \scriptstyle \bullet$} [c] at 8 0
\put {$ \scriptstyle \bullet$} [c] at 10 12
\put {$ \scriptstyle \bullet$} [c] at 12 12
\put {$ \scriptstyle \bullet$} [c] at 12 0
\put {$ \scriptstyle \bullet$} [c] at 14 0
\setlinear \plot 6 12 8 0 10 12    /
\setlinear \plot 12  0 12 12     /
\put{$360$} [c] at 10 -2
\endpicture
\end{minipage}
\begin{minipage}{4cm}
\beginpicture
\setcoordinatesystem units <1.5mm,2mm>
\setplotarea x from 0 to 16, y from -2 to 15
\put{ ${\bf  32}$} [l] at 0 15

\put{306)} [l] at 0 12
\put {$ \scriptstyle \bullet$} [c] at 6 0
\put {$ \scriptstyle \bullet$} [c] at 6 12
\put {$ \scriptstyle \bullet$} [c] at 10 0
\put {$ \scriptstyle \bullet$} [c] at 10 12
\put {$ \scriptstyle \bullet$} [c] at 12 0
\put {$ \scriptstyle \bullet$} [c] at 14 0
\setlinear \plot 6 12 6 0 10 12 10 0    /

\put{$360$} [c] at 10 -2
\endpicture
\end{minipage}
$$
$$
\begin{minipage}{4cm}
\beginpicture
\setcoordinatesystem units <1.5mm,2mm>
\setplotarea x from 0 to 16, y from -2 to 15
\put{307)} [l] at 0 12
\put {$ \scriptstyle \bullet$} [c] at 6 0
\put {$ \scriptstyle \bullet$} [c] at 6 6
\put {$ \scriptstyle \bullet$} [c] at 6 12
\put {$ \scriptstyle \bullet$} [c] at 10 0
\put {$ \scriptstyle \bullet$} [c] at 12 0
\put {$ \scriptstyle \bullet$} [c] at 14 0
\setlinear \plot 6 0 6 12      /
\put{$120$} [c] at 10 -2
\endpicture
\end{minipage}
\begin{minipage}{4cm}
\beginpicture
\setcoordinatesystem units <1.5mm,2mm>
\setplotarea x from 0 to 16, y from -2 to 15
\put{ ${\bf  33}$} [l] at 0 15

\put{308)} [l] at 0 12
\put {$ \scriptstyle \bullet$} [c] at 6 12
\put {$ \scriptstyle \bullet$} [c] at 8 12
\put {$ \scriptstyle \bullet$} [c] at 10 12
\put {$ \scriptstyle \bullet$} [c] at 12 12
\put {$ \scriptstyle \bullet$} [c] at 14 12
\put {$ \scriptstyle \bullet$} [c] at 10 0
\setlinear \plot 6 12 10 0 14 12     /
\setlinear \plot 8 12 10 0 12 12     /
\setlinear \plot 10 0 10 12     /

\put{$6$} [c] at 10 -2
\endpicture
\end{minipage}
\begin{minipage}{4cm}
\beginpicture
\setcoordinatesystem units <1.5mm,2mm>
\setplotarea x from 0 to 16, y from -2 to 15
\put{309)} [l] at 0 12
\put {$ \scriptstyle \bullet$} [c] at 6 0
\put {$ \scriptstyle \bullet$} [c] at 8 0
\put {$ \scriptstyle \bullet$} [c] at 10 0
\put {$ \scriptstyle \bullet$} [c] at 12 0
\put {$ \scriptstyle \bullet$} [c] at 14 0
\put {$ \scriptstyle \bullet$} [c] at 10 12
\setlinear \plot 6 0 10 12 14 0     /
\setlinear \plot 8 0 10 12 12 0     /
\setlinear \plot 10 0 10 12     /
\put{$6$} [c] at 10 -2
\endpicture
\end{minipage}
\begin{minipage}{4cm}
\beginpicture
\setcoordinatesystem units <1.5mm,2mm>
\setplotarea x from 0 to 16, y from -2 to 15
\put{ ${\bf  34}$} [l] at 0 15

\put{310)} [l] at 0 12
\put {$ \scriptstyle \bullet$} [c] at 6 12
\put {$ \scriptstyle \bullet$} [c] at 8 12
\put {$ \scriptstyle \bullet$} [c] at 10 12
\put {$ \scriptstyle \bullet$} [c] at 12 12
\put {$ \scriptstyle \bullet$} [c] at 9 0
\put {$ \scriptstyle \bullet$} [c] at 14  0
\setlinear \plot 6 12 9 0 12 12     /
\setlinear \plot 8 12 9 0 10 12     /
\put{$30$} [c] at 10 -2
\endpicture
\end{minipage}
\begin{minipage}{4cm}
\beginpicture
\setcoordinatesystem units <1.5mm,2mm>
\setplotarea x from 0 to 16, y from -2 to 15
\put{311)} [l] at 0 12
\put {$ \scriptstyle \bullet$} [c] at 6 0
\put {$ \scriptstyle \bullet$} [c] at 8 0
\put {$ \scriptstyle \bullet$} [c] at 10 0
\put {$ \scriptstyle \bullet$} [c] at 12 0
\put {$ \scriptstyle \bullet$} [c] at 9 12
\put {$ \scriptstyle \bullet$} [c] at 14  0
\setlinear \plot 6 0 9 12 12 0     /
\setlinear \plot 8 0 9 12 10 0     /
\put{$30$} [c] at 10 -2
\endpicture
\end{minipage}
\begin{minipage}{4cm}
\beginpicture
\setcoordinatesystem units <1.5mm,2mm>
\setplotarea x from 0 to 16, y from -2 to 15
\put{ ${\bf  36}$} [l] at 0 15

\put{312)} [l] at 0 12
\put {$ \scriptstyle \bullet$} [c] at 6 12
\put {$ \scriptstyle \bullet$} [c] at 6 0
\put {$ \scriptstyle \bullet$} [c] at 9 12
\put {$ \scriptstyle \bullet$} [c] at 9 0
\put {$ \scriptstyle \bullet$} [c] at 12 0
\put {$ \scriptstyle \bullet$} [c] at 14  0
\setlinear \plot 6 0 6 12     /
\setlinear \plot 9 0 9 12     /
\put{$180$} [c] at 10 -2
\endpicture
\end{minipage}
$$
$$
\begin{minipage}{4cm}
\beginpicture
\setcoordinatesystem units <1.5mm,2mm>
\setplotarea x from 0 to 16, y from -2 to 15
\put{313)} [l] at 0 12
\put {$ \scriptstyle \bullet$} [c] at 6 12
\put {$ \scriptstyle \bullet$} [c] at 8 12
\put {$ \scriptstyle \bullet$} [c] at 10 12
\put {$ \scriptstyle \bullet$} [c] at 8 0
\put {$ \scriptstyle \bullet$} [c] at 12 0
\put {$ \scriptstyle \bullet$} [c] at 14  0
\setlinear \plot 6 12 8 0 10  12     /
\setlinear \plot  8 0 8  12     /
\put{$60$} [c] at 10 -2
\endpicture
\end{minipage}
\begin{minipage}{4cm}
\beginpicture
\setcoordinatesystem units <1.5mm,2mm>
\setplotarea x from 0 to 16, y from -2 to 15
\put{314)} [l] at 0 12
\put {$ \scriptstyle \bullet$} [c] at 6 0
\put {$ \scriptstyle \bullet$} [c] at 8 0
\put {$ \scriptstyle \bullet$} [c] at 10 0
\put {$ \scriptstyle \bullet$} [c] at 8 12
\put {$ \scriptstyle \bullet$} [c] at 12 0
\put {$ \scriptstyle \bullet$} [c] at 14  0
\setlinear \plot 6 0 8 12 10  0     /
\setlinear \plot  8 0 8  12     /
\put{$60$} [c] at 10 -2
\endpicture
\end{minipage}
\begin{minipage}{4cm}
\beginpicture
\setcoordinatesystem units <1.5mm,2mm>
\setplotarea x from 0 to 16, y from -2 to 15
\put{ ${\bf  40}$} [l] at 0 15
\put{315)} [l] at 0 12
\put {$ \scriptstyle \bullet$} [c] at 6 12
\put {$ \scriptstyle \bullet$} [c] at 7 0
\put {$ \scriptstyle \bullet$} [c] at 8 12
\put {$ \scriptstyle \bullet$} [c] at 10 0
\put {$ \scriptstyle \bullet$} [c] at 12 0
\put {$ \scriptstyle \bullet$} [c] at 14  0
\setlinear \plot 6 12 7 0 8 12     /
\put{$60$} [c] at 10 -2
\endpicture
\end{minipage}
\begin{minipage}{4cm}
\beginpicture
\setcoordinatesystem units <1.5mm,2mm>
\setplotarea x from 0 to 16, y from -2 to 15
\put{316)} [l] at 0 12
\put {$ \scriptstyle \bullet$} [c] at 6 0
\put {$ \scriptstyle \bullet$} [c] at 7 12
\put {$ \scriptstyle \bullet$} [c] at 8 0
\put {$ \scriptstyle \bullet$} [c] at 10 0
\put {$ \scriptstyle \bullet$} [c] at 12 0
\put {$ \scriptstyle \bullet$} [c] at 14  0
\setlinear \plot 6 0 7 12 8 0     /
\put{$60$} [c] at 10 -2
\endpicture
\end{minipage}
\begin{minipage}{4cm}
\beginpicture
\setcoordinatesystem units <1.5mm,2mm>
\setplotarea x from 0 to 16, y from -2 to 15
\put{ ${\bf  48}$} [l] at 0 15

\put{317)} [l] at 0 12
\put {$ \scriptstyle \bullet$} [c] at 6 0
\put {$ \scriptstyle \bullet$} [c] at 6 12
\put {$ \scriptstyle \bullet$} [c] at 8 0
\put {$ \scriptstyle \bullet$} [c] at 10 0
\put {$ \scriptstyle \bullet$} [c] at 12 0
\put {$ \scriptstyle \bullet$} [c] at 14  0
\setlinear \plot 6 0 6 12     /
\put{$30$} [c] at 10 -2
\endpicture
\end{minipage}
\begin{minipage}{4cm}
\beginpicture
\setcoordinatesystem units <1.5mm,2mm>
\setplotarea x from 0 to 16, y from -2 to 15
\put{ ${\bf  64}$} [l] at 0 15
\put{318)} [l] at 0 12
\put {$ \scriptstyle \bullet$} [c] at 6 0
\put {$ \scriptstyle \bullet$} [c] at 8 0
\put {$ \scriptstyle \bullet$} [c] at 10 0
\put {$ \scriptstyle \bullet$} [c] at 12 0
\put {$ \scriptstyle \bullet$} [c] at 14 0
\put {$ \scriptstyle \bullet$} [c] at 16 0
\put{$1$} [c] at 11 -2
\endpicture
\end{minipage}
$$
\newpage

   ${\bf  NIP(7)}$

$$
\begin{minipage}{4cm}
\beginpicture
\setcoordinatesystem units   <1.5mm,2mm>
\setplotarea x from 0 to 16, y from -2 to 15
\put{${\bf  8}$} [l] at 2 15
\put{1)} [l] at 2 12
\put {$ \scriptstyle \bullet$} [c] at 12 0
\put {$ \scriptstyle \bullet$} [c] at 12 2
\put {$ \scriptstyle \bullet$} [c] at 12 4
\put {$ \scriptstyle \bullet$} [c] at 12 6
\put {$ \scriptstyle \bullet$} [c] at 12 8
\put {$ \scriptstyle \bullet$} [c] at 12 10
\put {$ \scriptstyle \bullet$} [c] at 12  12
\setlinear \plot  12 0   12 12     /
\put{$5{,}040$} [c] at 12 -2
\endpicture
\end{minipage}
\begin{minipage}{4cm}
\beginpicture
\setcoordinatesystem units   <1.5mm,2mm>
\setplotarea x from 0 to 16, y from -2 to 15
\put{${\bf  9}$} [l] at 2 15

\put{2)} [l] at 2 12
\put {$ \scriptstyle \bullet$} [c] at 13 0
\put {$ \scriptstyle \bullet$} [c] at 13 2
\put {$ \scriptstyle \bullet$} [c] at 13 4
\put {$ \scriptstyle \bullet$} [c] at 13 6
\put {$ \scriptstyle \bullet$} [c] at 10 9
\put {$ \scriptstyle \bullet$} [c] at 16 9
\put {$ \scriptstyle \bullet$} [c] at 13 12
\setlinear \plot  13 0 13 6  10 9 13 12 16 9 13  6     /
\put{$2{,}520$} [c] at 13 -2
\endpicture
\end{minipage}
\begin{minipage}{4cm}
\beginpicture
\setcoordinatesystem units   <1.5mm,2mm>
\setplotarea x from 0 to 16, y from -2 to 15
\put{3)} [l] at 2 12
\put {$ \scriptstyle \bullet$} [c] at 13 0
\put {$ \scriptstyle \bullet$} [c] at 10 3
\put {$ \scriptstyle \bullet$} [c] at 16 3
\put {$ \scriptstyle \bullet$} [c] at 13 6
\put {$ \scriptstyle \bullet$} [c] at 13 8
\put {$ \scriptstyle \bullet$} [c] at 13 10
\put {$ \scriptstyle \bullet$} [c] at 13 12
\setlinear \plot  13 12 13 6  10 3 13 0 16 3 13  6     /
\put{$2{,}520$} [c] at 13 -2
\endpicture
\end{minipage}
\begin{minipage}{4cm}
\beginpicture
\setcoordinatesystem units   <1.5mm,2mm>
\setplotarea x from 0 to 16, y from -2 to 15
\put{4)} [l] at 2 12
\put {$ \scriptstyle \bullet$} [c] at 13 0
\put {$ \scriptstyle \bullet$} [c] at 13 2
\put {$ \scriptstyle \bullet$} [c] at 13 4
\put {$ \scriptstyle \bullet$} [c] at 10 7
\put {$ \scriptstyle \bullet$} [c] at 16 7
\put {$ \scriptstyle \bullet$} [c] at 13 10
\put {$ \scriptstyle \bullet$} [c] at 13 12
\setlinear \plot  13 0 13 4 10 7 13 10 13 12    /
\setlinear \plot  13 4 16 7 13 10  /
\put{$2{,}520$} [c] at 13 -2
\endpicture
\end{minipage}
\begin{minipage}{4cm}
\beginpicture
\setcoordinatesystem units   <1.5mm,2mm>
\setplotarea x from 0 to 16, y from -2 to 15
\put{5)} [l] at 2 12
\put {$ \scriptstyle \bullet$} [c] at 13 0
\put {$ \scriptstyle \bullet$} [c] at 13 2
\put {$ \scriptstyle \bullet$} [c] at 10 5
\put {$ \scriptstyle \bullet$} [c] at 16 5
\put {$ \scriptstyle \bullet$} [c] at 13 8
\put {$ \scriptstyle \bullet$} [c] at 13 10
\put {$ \scriptstyle \bullet$} [c] at 13 12
\setlinear \plot  13 0 13 2 10 5 13 8 13 12    /
\setlinear \plot  13 2 16 5 13 8  /
\put{$2{,}520$} [c] at 13 -2
\endpicture
\end{minipage}
\begin{minipage}{4cm}
\beginpicture
\setcoordinatesystem units   <1.5mm,2mm>
\setplotarea x from 0 to 16, y from -2 to 15
\put{6)} [l] at 2 12
\put {$ \scriptstyle \bullet$} [c] at 13 0
\put {$ \scriptstyle \bullet$} [c] at 13 2
\put {$ \scriptstyle \bullet$} [c] at 13 4
\put {$ \scriptstyle \bullet$} [c] at 13 6
\put {$ \scriptstyle \bullet$} [c] at 13 8
\put {$ \scriptstyle \bullet$} [c] at 10 12
\put {$ \scriptstyle \bullet$} [c] at 16 12
\setlinear \plot  13 0 13 8  10 12      /
\setlinear \plot  13 8  16 12      /
\put{$2{,}520$} [c] at 13 -2
\endpicture
\end{minipage}
$$
$$
\begin{minipage}{4cm}
\beginpicture
\setcoordinatesystem units   <1.5mm,2mm>
\setplotarea x from 0 to 16, y from -2 to 15
\put{7)} [l] at 2 12
\put {$ \scriptstyle \bullet$} [c] at 13 12
\put {$ \scriptstyle \bullet$} [c] at 13 10
\put {$ \scriptstyle \bullet$} [c] at 13 8
\put {$ \scriptstyle \bullet$} [c] at 13 6
\put {$ \scriptstyle \bullet$} [c] at 13 4
\put {$ \scriptstyle \bullet$} [c] at 10 0
\put {$ \scriptstyle \bullet$} [c] at 16 0
\setlinear \plot  13 12 13 4  10 0      /
\setlinear \plot  16 0  13 4     /
\put{$2{,}520$} [c] at 13 -2
\endpicture
\end{minipage}
\begin{minipage}{4cm}
\beginpicture
\setcoordinatesystem units   <1.5mm,2mm>
\setplotarea x from 0 to 16, y from -2 to 15
\put{${\bf  10}$} [l] at 1 15

\put{8)} [l] at 2 12
\put {$ \scriptstyle \bullet$} [c] at 13 0
\put {$ \scriptstyle \bullet$} [c] at 13 8
\put {$ \scriptstyle \bullet$} [c] at 13 10
\put {$ \scriptstyle \bullet$} [c] at 13 12
\put {$ \scriptstyle \bullet$} [c] at 10 4
\put {$ \scriptstyle \bullet$} [c] at 16 4
\put {$ \scriptstyle \bullet$} [c] at 11.5 6
\setlinear \plot  13 0 10 4 13 8 13 12      /
\setlinear \plot  13 0  16 4 13 8      /
\put{$5{,}040$} [c] at 13 -2
\endpicture
\end{minipage}
\begin{minipage}{4cm}
\beginpicture
\setcoordinatesystem units   <1.5mm,2mm>
\setplotarea x from 0 to 16, y from -2 to 15
\put{9)} [l] at 2 12
\put {$ \scriptstyle \bullet$} [c] at 13 0
\put {$ \scriptstyle \bullet$} [c] at 13 2
\put {$ \scriptstyle \bullet$} [c] at 13 4
\put {$ \scriptstyle \bullet$} [c] at 10 8
\put {$ \scriptstyle \bullet$} [c] at 16 8
\put {$ \scriptstyle \bullet$} [c] at 13 12
\put {$ \scriptstyle \bullet$} [c] at 11.5 6
\setlinear \plot  13 0 13 4 10 8 13 12      /
\setlinear \plot  13 12  16 8 13 4      /
\put{$5{,}040$} [c] at 13 -2
\endpicture
\end{minipage}
\begin{minipage}{4cm}
\beginpicture
\setcoordinatesystem units   <1.5mm,2mm>
\setplotarea x from 0 to 16, y from -2 to 15
\put{10)} [l] at 2 12
\put {$ \scriptstyle \bullet$} [c] at 13 0
\put {$ \scriptstyle \bullet$} [c] at 13 3
\put {$ \scriptstyle \bullet$} [c] at 10 5
\put {$ \scriptstyle \bullet$} [c] at 16 6
\put {$ \scriptstyle \bullet$} [c] at 10 7
\put {$ \scriptstyle \bullet$} [c] at 13 9
\put {$ \scriptstyle \bullet$} [c] at 13 12
\setlinear \plot  13 0  13 3 10 5 10 7 13 9 13 12      /
\setlinear \plot  13 3 16 6 13 9      /
\put{$5{,}040$} [c] at 13 -2
\endpicture
\end{minipage}
\begin{minipage}{4cm}
\beginpicture
\setcoordinatesystem units   <1.5mm,2mm>
\setplotarea x from 0 to 16, y from -2 to 15
\put{11)} [l] at 2 12
\put {$ \scriptstyle \bullet$} [c] at 13 0
\put {$ \scriptstyle \bullet$} [c] at 13 2
\put {$ \scriptstyle \bullet$} [c] at 13 4
\put {$ \scriptstyle \bullet$} [c] at 13 6
\put {$ \scriptstyle \bullet$} [c] at 11.5 9
\put {$ \scriptstyle \bullet$} [c] at 10 12
\put {$ \scriptstyle \bullet$} [c] at 16 12
\setlinear \plot    13 0 13 6 10 12  /
\setlinear \plot    13 6 16 12 /
\put{$5{,}040$} [c] at 13 -2
\endpicture
\end{minipage}
\begin{minipage}{4cm}
\beginpicture
\setcoordinatesystem units   <1.5mm,2mm>
\setplotarea x from 0 to 16, y from -2 to 15
\put{12)} [l] at 2 12
\put {$ \scriptstyle \bullet$} [c] at 13 12
\put {$ \scriptstyle \bullet$} [c] at 13 10
\put {$ \scriptstyle \bullet$} [c] at 13 8
\put {$ \scriptstyle \bullet$} [c] at 13 6
\put {$ \scriptstyle \bullet$} [c] at 11.5 3
\put {$ \scriptstyle \bullet$} [c] at 10 0
\put {$ \scriptstyle \bullet$} [c] at 16 0
\setlinear \plot    13 12 13 6 10 0  /
\setlinear \plot    13 6 16 0 /
\put{$5{,}040$} [c] at 13 -2
\endpicture
\end{minipage}
$$
$$
\begin{minipage}{4cm}
\beginpicture
\setcoordinatesystem units   <1.5mm,2mm>
\setplotarea x from 0 to 16, y from -2 to 15
\put{13)} [l] at 2 12
\put {$ \scriptstyle \bullet$} [c] at 13 0
\put {$ \scriptstyle \bullet$} [c] at 10 2
\put {$ \scriptstyle \bullet$} [c] at 16 2
\put {$ \scriptstyle \bullet$} [c] at 10 8
\put {$ \scriptstyle \bullet$} [c] at 16 8
\put {$ \scriptstyle \bullet$} [c] at 13 10
\put {$ \scriptstyle \bullet$} [c] at 13 12
\setlinear \plot  16 8 16 2 13 0 10 2  16 8 13 10 13 12      /
\setlinear \plot  10 2  10 8 13 10      /
\setlinear \plot  10 8  16 2      /
\put{$1{,}260$} [c] at 13 -2
\endpicture
\end{minipage}
\begin{minipage}{4cm}
\beginpicture
\setcoordinatesystem units   <1.5mm,2mm>
\setplotarea x from 0 to 16, y from -2 to 15
\put{14)} [l] at 2 12
\put {$ \scriptstyle \bullet$} [c] at 13 0
\put {$ \scriptstyle \bullet$} [c] at 13 2
\put {$ \scriptstyle \bullet$} [c] at 10 4
\put {$ \scriptstyle \bullet$} [c] at 10 10
\put {$ \scriptstyle \bullet$} [c] at 16 4
\put {$ \scriptstyle \bullet$} [c] at 16 10
\put {$ \scriptstyle \bullet$} [c] at 13 12
\setlinear \plot 13 2 10 4 10 10 13 12 16 10 16 4   13 2 13 0      /
\setlinear \plot  10 4 16 10      /
\setlinear \plot  10 10 16 4       /
\put{$1{,}260$} [c] at 13 -2
\endpicture
\end{minipage}
\begin{minipage}{4cm}
\beginpicture
\setcoordinatesystem units   <1.5mm,2mm>
\setplotarea x from 0 to 16, y from -2 to 15
\put{15)} [l] at 2 12
\put {$ \scriptstyle \bullet$} [c] at 10 3
\put {$ \scriptstyle \bullet$} [c] at 10 9
\put {$ \scriptstyle \bullet$} [c] at 13 0
\put {$ \scriptstyle \bullet$} [c] at 13 6
\put {$ \scriptstyle \bullet$} [c] at 13 12
\put {$ \scriptstyle \bullet$} [c] at 16 3
\put {$ \scriptstyle \bullet$} [c] at 16 9
\setlinear \plot   13 0 10 3  16 9 13 12 10 9 16 3 13 0     /
\put{$1{,}260$} [c] at 13 -2
\endpicture
\end{minipage}
\begin{minipage}{4cm}
\beginpicture
\setcoordinatesystem units   <1.5mm,2mm>
\setplotarea x from 0 to 16, y from -2 to 15
\put{16)} [l] at 2 12
\put {$ \scriptstyle \bullet$} [c] at 13 0
\put {$ \scriptstyle \bullet$} [c] at 13 2
\put {$ \scriptstyle \bullet$} [c] at 13 4
\put {$ \scriptstyle \bullet$} [c] at 10 8
\put {$ \scriptstyle \bullet$} [c] at 10 12
\put {$ \scriptstyle \bullet$} [c] at 16 8
\put {$ \scriptstyle \bullet$} [c] at 16 12
\setlinear \plot    13 0 13 4 10 8 10 12 16 8  16 12 10 8     /
\setlinear \plot    13 4 16  8     /
\put{$1{,}260$} [c] at 13 -2
\endpicture
\end{minipage}
\begin{minipage}{4cm}
\beginpicture
\setcoordinatesystem units   <1.5mm,2mm>
\setplotarea x from 0 to 16, y from -2 to 15
\put{17)} [l] at 2 12
\put {$ \scriptstyle \bullet$} [c] at 13 12
\put {$ \scriptstyle \bullet$} [c] at 13 10
\put {$ \scriptstyle \bullet$} [c] at 13 8
\put {$ \scriptstyle \bullet$} [c] at 10 4
\put {$ \scriptstyle \bullet$} [c] at 10 0
\put {$ \scriptstyle \bullet$} [c] at 16 4
\put {$ \scriptstyle \bullet$} [c] at 16 0
\setlinear \plot    13 12 13 8 10 4 10 0 16 4  16 0 10 4     /
\setlinear \plot    13 8 16  4     /
\put{$1{,}260$} [c] at 13 -2
\endpicture
\end{minipage}
\begin{minipage}{4cm}
\beginpicture
\setcoordinatesystem units   <1.5mm,2mm>
\setplotarea x from 0 to 16, y from -2 to 15
\put{18)} [l] at 2 12
\put {$ \scriptstyle \bullet$} [c] at 13 0
\put {$ \scriptstyle \bullet$} [c] at 13 2
\put {$ \scriptstyle \bullet$} [c] at 10 6
\put {$ \scriptstyle \bullet$} [c] at 16 6
\put {$ \scriptstyle \bullet$} [c] at 13 10
\put {$ \scriptstyle \bullet$} [c] at 10 12
\put {$ \scriptstyle \bullet$} [c] at 16 12
\setlinear \plot   13 0 13 2 10 6 13 10 10 12     /
\setlinear \plot   16 12 13 10 16 6 13 2     /
\put{$1{,}260$} [c] at 13 -2
\endpicture
\end{minipage}
$$
$$
\begin{minipage}{4cm}
\beginpicture
\setcoordinatesystem units   <1.5mm,2mm>
\setplotarea x from 0 to 16, y from -2 to 15
\put{19)} [l] at 2 12
\put {$ \scriptstyle \bullet$} [c] at 13 12
\put {$ \scriptstyle \bullet$} [c] at 13 10
\put {$ \scriptstyle \bullet$} [c] at 10 6
\put {$ \scriptstyle \bullet$} [c] at 16 6
\put {$ \scriptstyle \bullet$} [c] at 13 2
\put {$ \scriptstyle \bullet$} [c] at 10 0
\put {$ \scriptstyle \bullet$} [c] at 16 0
\setlinear \plot   13 12 13 10 10 6 13 2 10 0     /
\setlinear \plot   16 0 13 2 16 6 13 10     /
\put{$1{,}260$} [c] at 13 -2
\endpicture
\end{minipage}
\begin{minipage}{4cm}
\beginpicture
\setcoordinatesystem units   <1.5mm,2mm>
\setplotarea x from 0 to 16, y from -2 to 15
\put{20)} [l] at 2 12
\put {$ \scriptstyle \bullet$} [c] at 13 0
\put {$ \scriptstyle \bullet$} [c] at 10 4
\put {$ \scriptstyle \bullet$} [c] at 16 4
\put {$ \scriptstyle \bullet$} [c] at 13 8
\put {$ \scriptstyle \bullet$} [c] at 13 10
\put {$ \scriptstyle \bullet$} [c] at 10 12
\put {$ \scriptstyle \bullet$} [c] at 16 12
\setlinear \plot   13 0 10 4  13 8 13 10 16 12     /
\setlinear \plot   13 0 16 4 13 8     /
\setlinear \plot   10 12 13 10     /
\put{$1{,}260$} [c] at 13 -2
\endpicture
\end{minipage}
\begin{minipage}{4cm}
\beginpicture
\setcoordinatesystem units   <1.5mm,2mm>
\setplotarea x from 0 to 16, y from -2 to 15
\put{21)} [l] at 2 12
\put {$ \scriptstyle \bullet$} [c] at 10 0
\put {$ \scriptstyle \bullet$} [c] at 16 0
\put {$ \scriptstyle \bullet$} [c] at 13 2
\put {$ \scriptstyle \bullet$} [c] at 13 4
\put {$ \scriptstyle \bullet$} [c] at 10 8
\put {$ \scriptstyle \bullet$} [c] at 16 8
\put {$ \scriptstyle \bullet$} [c] at 13 12
\setlinear \plot   10 0 13 2 13 4 10 8 13 12 16 8 13 4     /
\setlinear \plot   16 0  13 2     /
\put{$1{,}260$} [c] at 13 -2
\endpicture
\end{minipage}
\begin{minipage}{4cm}
\beginpicture
\setcoordinatesystem units   <1.5mm,2mm>
\setplotarea x from 0 to 16, y from -2 to 15
\put{22)} [l] at 2 12
\put {$ \scriptstyle \bullet$} [c] at  10 0
\put {$ \scriptstyle \bullet$} [c] at  16 0
\put {$ \scriptstyle \bullet$} [c] at  13 4
\put {$ \scriptstyle \bullet$} [c] at  13 6
\put {$ \scriptstyle \bullet$} [c] at  13 8
\put {$ \scriptstyle \bullet$} [c] at  10  12
\put {$ \scriptstyle \bullet$} [c] at  16  12
\setlinear \plot  10  0 13 4 13 8 10 12    /
\setlinear \plot  16 12 13 8 /
\setlinear \plot  16 0 13 4   /
\put{$1{,}260$} [c] at 13 -2
\endpicture
\end{minipage}
\begin{minipage}{4cm}
\beginpicture
\setcoordinatesystem units   <1.5mm,2mm>
\setplotarea x from 0 to 16, y from -2 to 15
\put{${\bf  11}$} [l] at 2 15
\put{23)} [l] at 2 12
\put {$ \scriptstyle \bullet$} [c] at 13 0
\put {$ \scriptstyle \bullet$} [c] at 13 2
\put {$ \scriptstyle \bullet$} [c] at 10 8
\put {$ \scriptstyle \bullet$} [c] at 16 6
\put {$ \scriptstyle \bullet$} [c] at 13 12
\put {$ \scriptstyle \bullet$} [c] at 11.5 5
\put {$ \scriptstyle \bullet$} [c] at 14.5 9
\setlinear \plot  13 0 13 2 10 8 13 12 16 6 13 2   /
\setlinear \plot  11.5 5 14.5 9 /
\put{$5{,}040$} [c] at 13 -2
\endpicture
\end{minipage}
\begin{minipage}{4cm}
\beginpicture
\setcoordinatesystem units   <1.5mm,2mm>
\setplotarea x from 0 to 16, y from -2 to 15
\put{24)} [l] at 2 12
\put {$ \scriptstyle \bullet$} [c] at 13 10
\put {$ \scriptstyle \bullet$} [c] at 13 12
\put {$ \scriptstyle \bullet$} [c] at 16 6
\put {$ \scriptstyle \bullet$} [c] at 13 0
\put {$ \scriptstyle \bullet$} [c] at 10 4
\put {$ \scriptstyle \bullet$} [c] at 11.5 7
\put {$ \scriptstyle \bullet$} [c] at 14.5 3
\setlinear \plot  13 12 13 10 16 6 13 0 10 4 13 10   /
\setlinear \plot  11.5 7 14.5 3 /
\put{$5{,}040$} [c] at 13 -2
\endpicture
\end{minipage}
$$
$$
\begin{minipage}{4cm}
\beginpicture
\setcoordinatesystem units   <1.5mm,2mm>
\setplotarea x from 0 to 16, y from -2 to 15
\put{25)} [l] at 2 12
\put {$ \scriptstyle \bullet$} [c] at 13 0
\put {$ \scriptstyle \bullet$} [c] at 13 2
\put {$ \scriptstyle \bullet$} [c] at 10 4
\put {$ \scriptstyle \bullet$} [c] at 10 6
\put {$ \scriptstyle \bullet$} [c] at 10 8
\put {$ \scriptstyle \bullet$} [c] at 13 12
\put {$ \scriptstyle \bullet$} [c] at 16 6
\setlinear \plot    13 0 13 2 10 4 10 8 13 12 16 6 13 2  /
\put{$5{,}040$} [c] at 13 -2
\endpicture
\end{minipage}
\begin{minipage}{4cm}
\beginpicture
\setcoordinatesystem units   <1.5mm,2mm>
\setplotarea x from 0 to 16, y from -2 to 15
\put{26)} [l] at 2 12
\put {$ \scriptstyle \bullet$} [c] at 13 12
\put {$ \scriptstyle \bullet$} [c] at 13 10
\put {$ \scriptstyle \bullet$} [c] at 10 8
\put {$ \scriptstyle \bullet$} [c] at 10 6
\put {$ \scriptstyle \bullet$} [c] at 10 4
\put {$ \scriptstyle \bullet$} [c] at 13 0
\put {$ \scriptstyle \bullet$} [c] at 16 6
\setlinear \plot    13 12 13 10 10 8 10 4 13 0 16 6 13 10  /
\put{$5{,}040$} [c] at 13 -2
\endpicture
\end{minipage}
\begin{minipage}{4cm}
\beginpicture
\setcoordinatesystem units   <1.5mm,2mm>
\setplotarea x from 0 to 16, y from -2 to 15
\put{27)} [l] at 2 12
\put {$ \scriptstyle \bullet$} [c] at 13 0
\put {$ \scriptstyle \bullet$} [c] at 13 2
\put {$ \scriptstyle \bullet$} [c] at 13 4
\put {$ \scriptstyle \bullet$} [c] at 11.1 9
\put {$ \scriptstyle \bullet$} [c] at 11.85 7
\put {$ \scriptstyle \bullet$} [c] at 10 12
\put {$ \scriptstyle \bullet$} [c] at 16 12
\setlinear \plot    13 0 13 4  10 12  /
\setlinear \plot    13 4 16 12 /
\put{$5{,}040$} [c] at 13 -2
\endpicture
\end{minipage}
\begin{minipage}{4cm}
\beginpicture
\setcoordinatesystem units   <1.5mm,2mm>
\setplotarea x from 0 to 16, y from -2 to 15
\put{28)} [l] at 2 12
\put {$ \scriptstyle \bullet$} [c] at 13 12
\put {$ \scriptstyle \bullet$} [c] at 13 10
\put {$ \scriptstyle \bullet$} [c] at 13 8
\put {$ \scriptstyle \bullet$} [c] at 11.1 3
\put {$ \scriptstyle \bullet$} [c] at 11.85 5
\put {$ \scriptstyle \bullet$} [c] at 10 0
\put {$ \scriptstyle \bullet$} [c] at 16 0
\setlinear \plot    10 0 13 8  13 12  /
\setlinear \plot    13 8 16 0 /
\put{$5{,}040$} [c] at 13 -2
\endpicture
\end{minipage}
\begin{minipage}{4cm}
\beginpicture
\setcoordinatesystem units   <1.5mm,2mm>
\setplotarea x from 0 to 16, y from -2 to 15
\put{29)} [l] at 2 12
\put {$ \scriptstyle \bullet$} [c] at 13 0
\put {$ \scriptstyle \bullet$} [c] at 13 2
\put {$ \scriptstyle \bullet$} [c] at 13 4
\put {$ \scriptstyle \bullet$} [c] at 10 8
\put {$ \scriptstyle \bullet$} [c] at 16 8
\put {$ \scriptstyle \bullet$} [c] at 13 12
\put {$ \scriptstyle \bullet$} [c] at 10 12
\setlinear \plot    10 12 10 8 13 4 13 0  /
\setlinear \plot    13 4 16 8 13 12 10 8 /
\put{$5{,}040$} [c] at 13 -2
\endpicture
\end{minipage}
\begin{minipage}{4cm}
\beginpicture
\setcoordinatesystem units   <1.5mm,2mm>
\setplotarea x from 0 to 16, y from -2 to 15
\put{30)} [l] at 2 12
\put {$ \scriptstyle \bullet$} [c] at 10 0
\put {$ \scriptstyle \bullet$} [c] at 13 0
\put {$ \scriptstyle \bullet$} [c] at 10 4
\put {$ \scriptstyle \bullet$} [c] at 16 4
\put {$ \scriptstyle \bullet$} [c] at 13 8
\put {$ \scriptstyle \bullet$} [c] at 13 10
\put {$ \scriptstyle \bullet$} [c] at 13 12
\setlinear \plot    10 0 10 4 13 8 13 12  /
\setlinear \plot   10 4 13 0 16 4 13 8  /
\put{$5{,}040$} [c] at 13 -2
\endpicture
\end{minipage}
$$

$$
\begin{minipage}{4cm}
\beginpicture
\setcoordinatesystem units   <1.5mm,2mm>
\setplotarea x from 0 to 16, y from -2 to 15
\put{31)} [l] at 2 12
\put {$ \scriptstyle \bullet$} [c] at 13 0
\put {$ \scriptstyle \bullet$} [c] at 10 2
\put {$ \scriptstyle \bullet$} [c] at 10 6
\put {$ \scriptstyle \bullet$} [c] at 10 10
\put {$ \scriptstyle \bullet$} [c] at 13 12
\put {$ \scriptstyle \bullet$} [c] at 16 10
\put {$ \scriptstyle \bullet$} [c] at 16 2
\setlinear \plot  16 2  13 0 10 2  10 10 13 12 16 10 16 2 10 10 /
\setlinear \plot    10 6 16 10  /
\put{$2{,}520$} [c] at 13 -2
\endpicture
\end{minipage}
\begin{minipage}{4cm}
\beginpicture
\setcoordinatesystem units   <1.5mm,2mm>
\setplotarea x from 0 to 16, y from -2 to 15
\put{32)} [l] at 2 12
\put {$ \scriptstyle \bullet$} [c] at 13 0
\put {$ \scriptstyle \bullet$} [c] at 10 2
\put {$ \scriptstyle \bullet$} [c] at 10 6
\put {$ \scriptstyle \bullet$} [c] at 10 10
\put {$ \scriptstyle \bullet$} [c] at 13 12
\put {$ \scriptstyle \bullet$} [c] at 16 10
\put {$ \scriptstyle \bullet$} [c] at 16 2
\setlinear \plot  16 10  13 12 10 10  10 2 13 0 16 2 16 10 10 2 /
\setlinear \plot    10 6 16 2  /
\put{$2{,}520$} [c] at 13 -2
\endpicture
\end{minipage}
\begin{minipage}{4cm}
\beginpicture
\setcoordinatesystem units   <1.5mm,2mm>
\setplotarea x from 0 to 16, y from -2 to 15
\put{33)} [l] at 2 12
\put {$ \scriptstyle \bullet$} [c] at 13 0
\put {$ \scriptstyle \bullet$} [c] at 13 2
\put {$ \scriptstyle \bullet$} [c] at 10 4
\put {$ \scriptstyle \bullet$} [c] at 10 8
\put {$ \scriptstyle \bullet$} [c] at 10 12
\put {$ \scriptstyle \bullet$} [c] at 16 4
\put {$ \scriptstyle \bullet$} [c] at 16 12
\setlinear \plot    13 0 13 2 10 4 10 8 16 4 16 12  10 4  /
\setlinear \plot    13 2 16 4  /
\setlinear \plot    10 8 10 12  /
\put{$2{,}520$} [c] at 13 -2
\endpicture
\end{minipage}
\begin{minipage}{4cm}
\beginpicture
\setcoordinatesystem units   <1.5mm,2mm>
\setplotarea x from 0 to 16, y from -2 to 15
\put{34)} [l] at 2 12
\put {$ \scriptstyle \bullet$} [c] at 13 12
\put {$ \scriptstyle \bullet$} [c] at 13 10
\put {$ \scriptstyle \bullet$} [c] at 10 8
\put {$ \scriptstyle \bullet$} [c] at 10 4
\put {$ \scriptstyle \bullet$} [c] at 10 0
\put {$ \scriptstyle \bullet$} [c] at 16 8
\put {$ \scriptstyle \bullet$} [c] at 16 0
\setlinear \plot    13 12 13 10 10 8 16 0 16 8 10 4 10 0 /
\setlinear \plot    13 10 16 8  /
\setlinear \plot    10 4 10 8  /
\put{$2{,}520$} [c] at 13 -2
\endpicture
\end{minipage}
\begin{minipage}{4cm}
\beginpicture
\setcoordinatesystem units   <1.5mm,2mm>
\setplotarea x from 0 to 16, y from -2 to 15
\put{35)} [l] at 2 12
\put {$ \scriptstyle \bullet$} [c] at 13 0
\put {$ \scriptstyle \bullet$} [c] at 10 3
\put {$ \scriptstyle \bullet$} [c] at 16 3
\put {$ \scriptstyle \bullet$} [c] at 13 7.5
\put {$ \scriptstyle \bullet$} [c] at 11.5 9.7
\put {$ \scriptstyle \bullet$} [c] at 10 12
\put {$ \scriptstyle \bullet$} [c] at 16 12
\setlinear \plot    10 12 16 3 13 0 10 3 16  12 /
\put{$2{,}520$} [c] at 13 -2
\endpicture
\end{minipage}
\begin{minipage}{4cm}
\beginpicture
\setcoordinatesystem units   <1.5mm,2mm>
\setplotarea x from 0 to 16, y from -2 to 15
\put{36)} [l] at 2 12
\put {$ \scriptstyle \bullet$} [c] at 13 12
\put {$ \scriptstyle \bullet$} [c] at 10 9
\put {$ \scriptstyle \bullet$} [c] at 16 9
\put {$ \scriptstyle \bullet$} [c] at 13 4.5
\put {$ \scriptstyle \bullet$} [c] at 11.5 2.3
\put {$ \scriptstyle \bullet$} [c] at 10 0
\put {$ \scriptstyle \bullet$} [c] at 16 0
\setlinear \plot    10 0 16 9 13 12 10 9 16  0 /
\put{$2{,}520$} [c] at 13 -2
\endpicture
\end{minipage}
$$
$$
\begin{minipage}{4cm}
\beginpicture
\setcoordinatesystem units   <1.5mm,2mm>
\setplotarea x from 0 to 16, y from -2 to 15
\put{37)} [l] at 2 12
\put {$ \scriptstyle \bullet$} [c] at 13 0
\put {$ \scriptstyle \bullet$} [c] at 10 3
\put {$ \scriptstyle \bullet$} [c] at 16 3
\put {$ \scriptstyle \bullet$} [c] at 11.5 4.5
\put {$ \scriptstyle \bullet$} [c] at 13 6
\put {$ \scriptstyle \bullet$} [c] at 10 12
\put {$ \scriptstyle \bullet$} [c] at 16 12
\setlinear \plot    13 0 10 3 13 6 10 12 /
\setlinear \plot    13 0 16 3 13 6 16 12 /
\put{$2{,}520$} [c] at 13 -2
\endpicture
\end{minipage}
\begin{minipage}{4cm}
\beginpicture
\setcoordinatesystem units   <1.5mm,2mm>
\setplotarea x from 0 to 16, y from -2 to 15
\put{38)} [l] at 2 12
\put {$ \scriptstyle \bullet$} [c] at 13 12
\put {$ \scriptstyle \bullet$} [c] at 10 9
\put {$ \scriptstyle \bullet$} [c] at 16 9
\put {$ \scriptstyle \bullet$} [c] at 11.5 7.5
\put {$ \scriptstyle \bullet$} [c] at 13 6
\put {$ \scriptstyle \bullet$} [c] at 10 0
\put {$ \scriptstyle \bullet$} [c] at 16 0
\setlinear \plot    13 12 10 9 13 6 10 0 /
\setlinear \plot    13 12 16 9 13 6 16 0 /
\put{$2{,}520$} [c] at 13 -2
\endpicture
\end{minipage}
\begin{minipage}{4cm}
\beginpicture
\setcoordinatesystem units   <1.5mm,2mm>
\setplotarea x from 0 to 16, y from -2 to 15
\put{39)} [l] at 2 12
\put {$ \scriptstyle \bullet$} [c] at 13 0
\put {$ \scriptstyle \bullet$} [c] at 13 2
\put {$ \scriptstyle \bullet$} [c] at 10 8
\put {$ \scriptstyle \bullet$} [c] at 10 12
\put {$ \scriptstyle \bullet$} [c] at 16 12
\put {$ \scriptstyle \bullet$} [c] at 16 8
\put {$ \scriptstyle \bullet$} [c] at 11.5 5
\setlinear \plot    13 0 13 2 10 8 10 12 16 8 16 12 10 8  /
\setlinear \plot    13 2 16 8 /
\put{$2{,}520$} [c] at 13 -2
\endpicture
\end{minipage}
\begin{minipage}{4cm}
\beginpicture
\setcoordinatesystem units   <1.5mm,2mm>
\setplotarea x from 0 to 16, y from -2 to 15
\put{40)} [l] at 2 12
\put {$ \scriptstyle \bullet$} [c] at 13 10
\put {$ \scriptstyle \bullet$} [c] at 13 12
\put {$ \scriptstyle \bullet$} [c] at 10 4
\put {$ \scriptstyle \bullet$} [c] at 10 0
\put {$ \scriptstyle \bullet$} [c] at 16 0
\put {$ \scriptstyle \bullet$} [c] at 16 4
\put {$ \scriptstyle \bullet$} [c] at 11.5 7
\setlinear \plot    13 12 13 10 10 4 10 0 16 4 16 0 10 4  /
\setlinear \plot    13 10 16 4 /
\put{$2{,}520$} [c] at 13 -2
\endpicture
\end{minipage}
\begin{minipage}{4cm}
\beginpicture
\setcoordinatesystem units   <1.5mm,2mm>
\setplotarea x from 0 to 16, y from -2 to 15
\put{41)} [l] at 2 12
\put {$ \scriptstyle \bullet$} [c] at  10 0
\put {$ \scriptstyle \bullet$} [c] at  16 0
\put {$ \scriptstyle \bullet$} [c] at  11.5  2
\put {$ \scriptstyle \bullet$} [c] at  13 4
\put {$ \scriptstyle \bullet$} [c] at  13 8
\put {$ \scriptstyle \bullet$} [c] at  10  12
\put {$ \scriptstyle \bullet$} [c] at  16  12
\setlinear \plot  10 0 13 4 13 8 10 12    /
\setlinear \plot  16 0 13 4 /
\setlinear \plot  13 8 16 12  /
\put{$2{,}520  $} [c] at 13 -2
\endpicture
\end{minipage}
\begin{minipage}{4cm}
\beginpicture
\setcoordinatesystem units   <1.5mm,2mm>
\setplotarea x from 0 to 16, y from -2 to 15
\put{42)} [l] at 2 12
\put {$ \scriptstyle \bullet$} [c] at  10 0
\put {$ \scriptstyle \bullet$} [c] at  16 0
\put {$ \scriptstyle \bullet$} [c] at  11.5  10
\put {$ \scriptstyle \bullet$} [c] at  13 4
\put {$ \scriptstyle \bullet$} [c] at  13 8
\put {$ \scriptstyle \bullet$} [c] at  10  12
\put {$ \scriptstyle \bullet$} [c] at  16  12
\setlinear \plot  10 0 13 4 13 8 10 12    /
\setlinear \plot  16 0 13 4 /
\setlinear \plot  13 8 16 12  /
\put{$2{,}520 $} [c] at 13 -2
\endpicture
\end{minipage}
$$

$$
\begin{minipage}{4cm}
\beginpicture
\setcoordinatesystem units   <1.5mm,2mm>
\setplotarea x from 0 to 16, y from -2 to 15
\put{43)} [l] at 2 12
\put {$ \scriptstyle \bullet$} [c] at 13 0
\put {$ \scriptstyle \bullet$} [c] at 10 2
\put {$ \scriptstyle \bullet$} [c] at 10 7
\put {$ \scriptstyle \bullet$} [c] at 10 12
\put {$ \scriptstyle \bullet$} [c] at 16 2
\put {$ \scriptstyle \bullet$} [c] at 16 7
\put {$ \scriptstyle \bullet$} [c] at 16 12
\setlinear \plot  16 2  13 0 10 2 10 12 16 7 16 12 10 7 16 2 16 7 10 2 /
\put{$630$} [c] at 13 -2
\endpicture
\end{minipage}
\begin{minipage}{4cm}
\beginpicture
\setcoordinatesystem units   <1.5mm,2mm>
\setplotarea x from 0 to 16, y from -2 to 15
\put{44)} [l] at 2 12
\put {$ \scriptstyle \bullet$} [c] at 13 12
\put {$ \scriptstyle \bullet$} [c] at 10 10
\put {$ \scriptstyle \bullet$} [c] at 10 5
\put {$ \scriptstyle \bullet$} [c] at 10 0
\put {$ \scriptstyle \bullet$} [c] at 16 0
\put {$ \scriptstyle \bullet$} [c] at 16 5
\put {$ \scriptstyle \bullet$} [c] at 16 10
\setlinear \plot  16 10  13 12 10 10 10 0 16 5 16 0 10 5 16 10 16 5 10 10 /
\put{$630$} [c] at 13 -2
\endpicture
\end{minipage}
\begin{minipage}{4cm}
\beginpicture
\setcoordinatesystem units   <1.5mm,2mm>
\setplotarea x from 0 to 16, y from -2 to 15
\put{45)} [l] at 2 12
\put {$ \scriptstyle \bullet$} [c] at  10 0
\put {$ \scriptstyle \bullet$} [c] at  10 6
\put {$ \scriptstyle \bullet$} [c] at  16 0
\put {$ \scriptstyle \bullet$} [c] at  16 6
\put {$ \scriptstyle \bullet$} [c] at  13 10
\put {$ \scriptstyle \bullet$} [c] at  10 12
\put {$ \scriptstyle \bullet$} [c] at  16  12
\setlinear \plot 16 6  10 0 10 6 13 10 10 12   /
\setlinear \plot  16 12 13 10 16 6 16 0 10 6  /
\put{$630  $} [c] at 13 -2
\endpicture
\end{minipage}
\begin{minipage}{4cm}
\beginpicture
\setcoordinatesystem units   <1.5mm,2mm>
\setplotarea x from 0 to 16, y from 0 to 15
\put{46)} [l] at 2 12
\put {$ \scriptstyle \bullet$} [c] at  10 0
\put {$ \scriptstyle \bullet$} [c] at  16 0
\put {$ \scriptstyle \bullet$} [c] at  13 2
\put {$ \scriptstyle \bullet$} [c] at  10 6
\put {$ \scriptstyle \bullet$} [c] at  16 6
\put {$ \scriptstyle \bullet$} [c] at  16 12
\put {$ \scriptstyle \bullet$} [c] at  10  12
\setlinear \plot  10 0 13 2 16 6 16 12 10 6 10 12 16 6    /
\setlinear \plot  16 0 13 2 10 6  /
\put{$630  $} [c] at 13 -2
\endpicture
\end{minipage}
\begin{minipage}{4cm}
\beginpicture
\setcoordinatesystem units   <1.5mm,2mm>
\setplotarea x from 0 to 16, y from -2 to 15
\put{${\bf  12}$} [l] at 2 15
\put{47)} [l] at 2 12
\put {$ \scriptstyle \bullet$} [c] at 13 0
\put {$ \scriptstyle \bullet$} [c] at 10 3
\put {$ \scriptstyle \bullet$} [c] at 16 9
\put {$ \scriptstyle \bullet$} [c] at 13  12
\put {$ \scriptstyle \bullet$} [c] at 14 3
\put {$ \scriptstyle \bullet$} [c] at 11 6
\put {$ \scriptstyle \bullet$} [c] at 12 9
\setlinear \plot  13 0 10 3 13 12 16 9 13 0   /
\setlinear \plot  14 3 11 6   /
\put{$5{,}040$} [c] at 13 -2
\endpicture
\end{minipage}
\begin{minipage}{4cm}
\beginpicture
\setcoordinatesystem units   <1.5mm,2mm>
\setplotarea x from 0 to 16, y from -2 to 15
\put{48)} [l] at 2 12
\put {$ \scriptstyle \bullet$} [c] at 13 0
\put {$ \scriptstyle \bullet$} [c] at 10 3
\put {$ \scriptstyle \bullet$} [c] at 16 9
\put {$ \scriptstyle \bullet$} [c] at 13  12
\put {$ \scriptstyle \bullet$} [c] at 15 6
\put {$ \scriptstyle \bullet$} [c] at 14 3
\put {$ \scriptstyle \bullet$} [c] at 12 9
\setlinear \plot  13 0 10 3 13 12 16 9 13 0   /
\setlinear \plot  12 9 15 6   /
\put{$5{,}040$} [c] at 13 -2
\endpicture
\end{minipage}
$$
$$
\begin{minipage}{4cm}
\beginpicture
\setcoordinatesystem units   <1.5mm,2mm>
\setplotarea x from 0 to 16, y from -2 to 15
\put{49)} [l] at 2 12
\put {$ \scriptstyle \bullet$} [c] at 13 0
\put {$ \scriptstyle \bullet$} [c] at 10 2
\put {$ \scriptstyle \bullet$} [c] at 10 5
\put {$ \scriptstyle \bullet$} [c] at 10 7
\put {$ \scriptstyle \bullet$} [c] at 10 10
\put {$ \scriptstyle \bullet$} [c] at 16 6
\put {$ \scriptstyle \bullet$} [c] at 13 12
\setlinear \plot  13 0 10 2  10 10 13  12 16 6 13 0     /
\put{$5{,}040$} [c] at 13 -2
\endpicture
\end{minipage}
\begin{minipage}{4cm}
\beginpicture
\setcoordinatesystem units   <1.5mm,2mm>
\setplotarea x from 0 to 16, y from -2 to 15
\put{50)} [l] at 2 12
\put {$ \scriptstyle \bullet$} [c] at 13 0
\put {$ \scriptstyle \bullet$} [c] at 13 12
\put {$ \scriptstyle \bullet$} [c] at 10 3
\put {$ \scriptstyle \bullet$} [c] at 10 6
\put {$ \scriptstyle \bullet$} [c] at 10 9
\put {$ \scriptstyle \bullet$} [c] at 16 3
\put {$ \scriptstyle \bullet$} [c] at 16 9
\setlinear \plot 13 0 10 3 10 9  13 12 16 9 16 3 13 0     /
\setlinear \plot  10 9  16 3      /
\setlinear \plot  16 9  10 3      /
\put{$5{,}040$} [c] at 13 -2
\endpicture
\end{minipage}
\begin{minipage}{4cm}
\beginpicture
\setcoordinatesystem units   <1.5mm,2mm>
\setplotarea x from 0 to 16, y from -2 to 15
\put{51)} [l] at 2 12
\put {$ \scriptstyle \bullet$} [c] at 13 0
\put {$ \scriptstyle \bullet$} [c] at 13 2
\put {$ \scriptstyle \bullet$} [c] at 13 12
\put {$ \scriptstyle \bullet$} [c] at 13 10
\put {$ \scriptstyle \bullet$} [c] at 10 6
\put {$ \scriptstyle \bullet$} [c] at 16 12
\put {$ \scriptstyle \bullet$} [c] at 16 6
\setlinear \plot  13 0 13 2 10 6 13 10 13 12   /
\setlinear \plot   13 10 16 6 13 2  /
\setlinear \plot   16 12 16 6  /
\put{$5{,}040$} [c] at 13 -2
\endpicture
\end{minipage}
\begin{minipage}{4cm}
\beginpicture
\setcoordinatesystem units   <1.5mm,2mm>
\setplotarea x from 0 to 16, y from -2 to 15
\put{52)} [l] at 2 12
\put {$ \scriptstyle \bullet$} [c] at 13 0
\put {$ \scriptstyle \bullet$} [c] at 13 2
\put {$ \scriptstyle \bullet$} [c] at 13 12
\put {$ \scriptstyle \bullet$} [c] at 13 10
\put {$ \scriptstyle \bullet$} [c] at 10 6
\put {$ \scriptstyle \bullet$} [c] at 16 0
\put {$ \scriptstyle \bullet$} [c] at 16 6
\setlinear \plot  13 0 13 2 10 6 13 10 13 12   /
\setlinear \plot   13 10 16 6 13 2  /
\setlinear \plot   16 0 16 6  /
\put{$5{,}040$} [c] at 13 -2
\endpicture
\end{minipage}
\begin{minipage}{4cm}
\beginpicture
\setcoordinatesystem units   <1.5mm,2mm>
\setplotarea x from 0 to 16, y from -2 to 15
\put{53)} [l] at 2 12
\put {$ \scriptstyle \bullet$} [c] at 13 0
\put {$ \scriptstyle \bullet$} [c] at 13 2
\put {$ \scriptstyle \bullet$} [c] at 10 12
\put {$ \scriptstyle \bullet$} [c] at 16 12
\put {$ \scriptstyle \bullet$} [c] at 10.55 10
\put {$ \scriptstyle \bullet$} [c] at 11.5 7
\put {$ \scriptstyle \bullet$} [c] at 12.3 4
\setlinear \plot  13 0 13 2  10 12      /
\setlinear \plot  13 2  16 12      /
\put{$5{,}040$} [c] at 13 -2
\endpicture
\end{minipage}
\begin{minipage}{4cm}
\beginpicture
\setcoordinatesystem units   <1.5mm,2mm>
\setplotarea x from 0 to 16, y from -2 to 15
\put{54)} [l] at 2 12
\put {$ \scriptstyle \bullet$} [c] at 13 12
\put {$ \scriptstyle \bullet$} [c] at 13 10
\put {$ \scriptstyle \bullet$} [c] at 10 0
\put {$ \scriptstyle \bullet$} [c] at 16 0
\put {$ \scriptstyle \bullet$} [c] at 10.55 2
\put {$ \scriptstyle \bullet$} [c] at 11.5 5
\put {$ \scriptstyle \bullet$} [c] at 12.3 8
\setlinear \plot  13 12 13 10  10 0      /
\setlinear \plot  13 10  16 0      /
\put{$5{,}040$} [c] at 13 -2
\endpicture
\end{minipage}
$$

$$
\begin{minipage}{4cm}
\beginpicture
\setcoordinatesystem units   <1.5mm,2mm>
\setplotarea x from 0 to 16, y from -2 to 15
\put{55)} [l] at 2 12
\put {$ \scriptstyle \bullet$} [c] at 13 0
\put {$ \scriptstyle \bullet$} [c] at 13 4
\put {$ \scriptstyle \bullet$} [c] at 13 12
\put {$ \scriptstyle \bullet$} [c] at 11.5 6
\put {$ \scriptstyle \bullet$} [c] at 10 8
\put {$ \scriptstyle \bullet$} [c] at 10 12
\put {$ \scriptstyle \bullet$} [c] at 16 8
\setlinear \plot  13 0 13 4 10 8 13 12 16 8  13 4      /
\setlinear \plot  10 8 10 12   /
\put{$5{,}040$} [c] at 13 -2
\endpicture
\end{minipage}
\begin{minipage}{4cm}
\beginpicture
\setcoordinatesystem units   <1.5mm,2mm>
\setplotarea x from 0 to 16, y from -2 to 15
\put{56)} [l] at 2 12
\put {$ \scriptstyle \bullet$} [c] at 13 12
\put {$ \scriptstyle \bullet$} [c] at 13 8
\put {$ \scriptstyle \bullet$} [c] at 13 0
\put {$ \scriptstyle \bullet$} [c] at 11.5 6
\put {$ \scriptstyle \bullet$} [c] at 10 4
\put {$ \scriptstyle \bullet$} [c] at 10 0
\put {$ \scriptstyle \bullet$} [c] at 16 4
\setlinear \plot  13 12 13 8 10 4 13 0 16 4  13 8      /
\setlinear \plot  10 4 10 0   /
\put{$5{,}040$} [c] at 13 -2
\endpicture
\end{minipage}
\begin{minipage}{4cm}
\beginpicture
\setcoordinatesystem units   <1.5mm,2mm>
\setplotarea x from 0 to 16, y from -2 to 15
\put{57)} [l] at 2 12
\put {$ \scriptstyle \bullet$} [c] at 13 0
\put {$ \scriptstyle \bullet$} [c] at 13 3
\put {$ \scriptstyle \bullet$} [c] at 10 6
\put {$ \scriptstyle \bullet$} [c] at 10 9
\put {$ \scriptstyle \bullet$} [c] at 10 12
\put {$ \scriptstyle \bullet$} [c] at 16 6
\put {$ \scriptstyle \bullet$} [c] at 16 12
\setlinear \plot 13 0 13 3 10 6 10 12 16 6 16 12 10 6    /
\setlinear \plot 13 3 16 6     /
\put{$5{,}040$} [c] at 13 -2
\endpicture
\end{minipage}
\begin{minipage}{4cm}
\beginpicture
\setcoordinatesystem units   <1.5mm,2mm>
\setplotarea x from 0 to 16, y from -2 to 15
\put{58)} [l] at 2 12
\put {$ \scriptstyle \bullet$} [c] at 13 12
\put {$ \scriptstyle \bullet$} [c] at 13 9
\put {$ \scriptstyle \bullet$} [c] at 10 6
\put {$ \scriptstyle \bullet$} [c] at 10 3
\put {$ \scriptstyle \bullet$} [c] at 10 0
\put {$ \scriptstyle \bullet$} [c] at 16 6
\put {$ \scriptstyle \bullet$} [c] at 16 0
\setlinear \plot 13 12 13 9 10 6 10 0 16 6 16 0 10 6    /
\setlinear \plot 13 9 16 6     /
\put{$5{,}040$} [c] at 13 -2
\endpicture
\end{minipage}
\begin{minipage}{4cm}
\beginpicture
\setcoordinatesystem units   <1.5mm,2mm>
\setplotarea x from 0 to 16, y from -2 to 15
\put{59)} [l] at 2 12
\put {$ \scriptstyle \bullet$} [c] at 13  0
\put {$ \scriptstyle \bullet$} [c] at 10 3
\put {$ \scriptstyle \bullet$} [c] at 10 6
\put {$ \scriptstyle \bullet$} [c] at 10 9
\put {$ \scriptstyle \bullet$} [c] at 10 12
\put {$ \scriptstyle \bullet$} [c] at 16 12
\put {$ \scriptstyle \bullet$} [c] at 16 3
\setlinear \plot 10 6 16 12 16 3 13 0 10 3 10 12   /
\setlinear \plot   10 9 16 3   /
\put{$5{,}040$} [c] at 13 -2
\endpicture
\end{minipage}
\begin{minipage}{4cm}
\beginpicture
\setcoordinatesystem units   <1.5mm,2mm>
\setplotarea x from 0 to 16, y from -2 to 15
\put{60)} [l] at 2 12
\put {$ \scriptstyle \bullet$} [c] at 13  12
\put {$ \scriptstyle \bullet$} [c] at 10 3
\put {$ \scriptstyle \bullet$} [c] at 10 6
\put {$ \scriptstyle \bullet$} [c] at 10 9
\put {$ \scriptstyle \bullet$} [c] at 10 0
\put {$ \scriptstyle \bullet$} [c] at 16 0
\put {$ \scriptstyle \bullet$} [c] at 16 9
\setlinear \plot 10 6 16 0 16 9 13 12 10 9 10 0   /
\setlinear \plot   10 3 16 9   /
\put{$5{,}040$} [c] at 13 -2
\endpicture
\end{minipage}
$$
$$
\begin{minipage}{4cm}
\beginpicture
\setcoordinatesystem units   <1.5mm,2mm>
\setplotarea x from 0 to 16, y from -2 to 15
\put {61)} [l] at 2 12
\put {$ \scriptstyle \bullet$} [c] at  10 0
\put {$ \scriptstyle \bullet$} [c] at  16 0
\put {$ \scriptstyle \bullet$} [c] at  10 12
\put {$ \scriptstyle \bullet$} [c] at  11.5 3
\put {$ \scriptstyle \bullet$} [c] at  11.5 9
\put {$ \scriptstyle \bullet$} [c] at  16 12
\put {$ \scriptstyle \bullet$} [c] at  13 6
\setlinear \plot  10 0 16 12    /
\setlinear \plot  16 0  10 12  /
\put{$5{,}040  $} [c] at 13 -2
\endpicture
\end{minipage}
\begin{minipage}{4cm}
\beginpicture
\setcoordinatesystem units   <1.5mm,2mm>
\setplotarea x from 0 to 16, y from -2 to 15
\put{62)} [l] at 2 12
\put {$ \scriptstyle \bullet$} [c] at 13 12
\put {$ \scriptstyle \bullet$} [c] at 13 9
\put {$ \scriptstyle \bullet$} [c] at 10 6
\put {$ \scriptstyle \bullet$} [c] at 10 3
\put {$ \scriptstyle \bullet$} [c] at 13 0
\put {$ \scriptstyle \bullet$} [c] at 16 3
\put {$ \scriptstyle \bullet$} [c] at 16 6
\setlinear \plot   13 12 13 9 10 6 10 3 13 0 16 3 16 6 13 9  /
\put{$2{,}520$} [c] at 13 -2
\endpicture
\end{minipage}
\begin{minipage}{4cm}
\beginpicture
\setcoordinatesystem units   <1.5mm,2mm>
\setplotarea x from 0 to 16, y from -2 to 15
\put{63)} [l] at 2 12
\put {$ \scriptstyle \bullet$} [c] at 13 0
\put {$ \scriptstyle \bullet$} [c] at 13 3
\put {$ \scriptstyle \bullet$} [c] at 10 6
\put {$ \scriptstyle \bullet$} [c] at 10 9
\put {$ \scriptstyle \bullet$} [c] at 13 12
\put {$ \scriptstyle \bullet$} [c] at 16 9
\put {$ \scriptstyle \bullet$} [c] at 16 6
\setlinear \plot   13 0 13 3 10 6 10 9 13 12 16 9 16 6 13 3  /
\put{$2{,}520$} [c] at 13 -2
\endpicture
\end{minipage}
\begin{minipage}{4cm}
\beginpicture
\setcoordinatesystem units   <1.5mm,2mm>
\setplotarea x from 0 to 16, y from -2 to 15
\put{64)} [l] at 2 12
\put {$ \scriptstyle \bullet$} [c] at 13 0
\put {$ \scriptstyle \bullet$} [c] at 13 3
\put {$ \scriptstyle \bullet$} [c] at 13 6
\put {$ \scriptstyle \bullet$} [c] at 10 9
\put {$ \scriptstyle \bullet$} [c] at 10 12
\put {$ \scriptstyle \bullet$} [c] at 16 12
\put {$ \scriptstyle \bullet$} [c] at 16 9
\setlinear \plot    13 0 13 6 10 9 10 12    /
\setlinear \plot   13 6 16 9 16 12  /
\put{$2{,}520$} [c] at 13 -2
\endpicture
\end{minipage}
\begin{minipage}{4cm}
\beginpicture
\setcoordinatesystem units   <1.5mm,2mm>
\setplotarea x from 0 to 16, y from -2 to 15
\put{65)} [l] at 2 12
\put {$ \scriptstyle \bullet$} [c] at 13 12
\put {$ \scriptstyle \bullet$} [c] at 13 9
\put {$ \scriptstyle \bullet$} [c] at 13 6
\put {$ \scriptstyle \bullet$} [c] at 10 3
\put {$ \scriptstyle \bullet$} [c] at 10 0
\put {$ \scriptstyle \bullet$} [c] at 16 0
\put {$ \scriptstyle \bullet$} [c] at 16 3
\setlinear \plot    13 12 13 6 10 3 10 0    /
\setlinear \plot   13 6 16 3 16 0  /
\put{$2{,}520$} [c] at 13 -2
\endpicture
\end{minipage}
\begin{minipage}{4cm}
\beginpicture
\setcoordinatesystem units   <1.5mm,2mm>
\setplotarea x from 0 to 16, y from -2 to 15
\put{66)} [l] at 2 12
\put {$ \scriptstyle \bullet$} [c] at 13 0
\put {$ \scriptstyle \bullet$} [c] at 10 3
\put {$ \scriptstyle \bullet$} [c] at 10 6
\put {$ \scriptstyle \bullet$} [c] at 10 9
\put {$ \scriptstyle \bullet$} [c] at 10 12
\put {$ \scriptstyle \bullet$} [c] at 16 3
\put {$ \scriptstyle \bullet$} [c] at 16 12
\setlinear \plot 10 3 16 12  16 3 13 0 10 3 10 12    /
\setlinear \plot 10 6  16 3    /
\put{$2{,}520$} [c] at 13 -2
\endpicture
\end{minipage}
$$

$$
\begin{minipage}{4cm}
\beginpicture
\setcoordinatesystem units   <1.5mm,2mm>
\setplotarea x from 0 to 16, y from -2 to 15
\put{67)} [l] at 2 12
\put {$ \scriptstyle \bullet$} [c] at 13 12
\put {$ \scriptstyle \bullet$} [c] at 10 0
\put {$ \scriptstyle \bullet$} [c] at 10 3
\put {$ \scriptstyle \bullet$} [c] at 10 6
\put {$ \scriptstyle \bullet$} [c] at 10 9
\put {$ \scriptstyle \bullet$} [c] at 16 9
\put {$ \scriptstyle \bullet$} [c] at 16 0
\setlinear \plot 10 9 16 0  16 9 13 12 10 9 10 0    /
\setlinear \plot 10 6  16 9    /
\put{$2{,}520$} [c] at 13 -2
\endpicture
\end{minipage}
\begin{minipage}{4cm}
\beginpicture
\setcoordinatesystem units   <1.5mm,2mm>
\setplotarea x from 0 to 16, y from -2 to 15
\put{68)} [l] at 2 12
\put {$ \scriptstyle \bullet$} [c] at 13 0
\put {$ \scriptstyle \bullet$} [c] at 10 3
\put {$ \scriptstyle \bullet$} [c] at 10 6
\put {$ \scriptstyle \bullet$} [c] at 10 9
\put {$ \scriptstyle \bullet$} [c] at 10 12
\put {$ \scriptstyle \bullet$} [c] at 16 3
\put {$ \scriptstyle \bullet$} [c] at 16 12
\setlinear \plot 16 3 13  0 10 3 10 12 16 3 16 12 10 9    /
\put{$2{,}520$} [c] at 13 -2
\endpicture
\end{minipage}
\begin{minipage}{4cm}
\beginpicture
\setcoordinatesystem units   <1.5mm,2mm>
\setplotarea x from 0 to 16, y from -2 to 15
\put{69)} [l] at 2 12
\put {$ \scriptstyle \bullet$} [c] at 13 12
\put {$ \scriptstyle \bullet$} [c] at 10 9
\put {$ \scriptstyle \bullet$} [c] at 10 6
\put {$ \scriptstyle \bullet$} [c] at 10 3
\put {$ \scriptstyle \bullet$} [c] at 10 0
\put {$ \scriptstyle \bullet$} [c] at 16 9
\put {$ \scriptstyle \bullet$} [c] at 16 0
\setlinear \plot 16 9 13  12 10 9 10 0 16 9 16 0 10 3    /
\put{$2{,}520$} [c] at 13 -2
\endpicture
\end{minipage}
\begin{minipage}{4cm}
\beginpicture
\setcoordinatesystem units   <1.5mm,2mm>
\setplotarea x from 0 to 16, y from -2 to 15
\put{70)} [l] at 2 12
\put {$ \scriptstyle \bullet$} [c] at 10 6
\put {$ \scriptstyle \bullet$} [c] at 10 12
\put {$ \scriptstyle \bullet$} [c] at 14 6
\put {$ \scriptstyle \bullet$} [c] at 14 12
\put {$ \scriptstyle \bullet$} [c] at 14 0
\put {$ \scriptstyle \bullet$} [c] at 16 3
\put {$ \scriptstyle \bullet$} [c] at 12 3
\setlinear \plot 14 0 10 6 10 12 14 6 14 12 10 6     /
\setlinear \plot 12 3 14 6 16 3 14 0   /
\put{$2{,}520$} [c] at 13 -2
\endpicture
\end{minipage}
\begin{minipage}{4cm}
\beginpicture
\setcoordinatesystem units   <1.5mm,2mm>
\setplotarea x from 0 to 16, y from -2 to 15
\put{71)} [l] at 2 12
\put {$ \scriptstyle \bullet$} [c] at 10 0
\put {$ \scriptstyle \bullet$} [c] at 10 6
\put {$ \scriptstyle \bullet$} [c] at 12 9
\put {$ \scriptstyle \bullet$} [c] at 14 0
\put {$ \scriptstyle \bullet$} [c] at 14 6
\put {$ \scriptstyle \bullet$} [c] at 14 12
\put {$ \scriptstyle \bullet$} [c] at 16 9
\setlinear \plot 10 0 10 6 14 12 16 9 14 6 14 0 10 6     /
\setlinear \plot 10 0 14 6 12 9   /
\put{$2{,}520$} [c] at 13 -2
\endpicture
\end{minipage}
\begin{minipage}{4cm}
\beginpicture
\setcoordinatesystem units   <1.5mm,2mm>
\setplotarea x from 0 to 16, y from -2 to 15
\put{72)} [l] at 2 12
\put {$ \scriptstyle \bullet$} [c] at 13 0
\put {$ \scriptstyle \bullet$} [c] at 13 12
\put {$ \scriptstyle \bullet$} [c] at 10 3
\put {$ \scriptstyle \bullet$} [c] at 10 9
\put {$ \scriptstyle \bullet$} [c] at 16 3
\put {$ \scriptstyle \bullet$} [c] at 16 9
\put {$ \scriptstyle \bullet$} [c] at 16 12
\setlinear \plot 16 3 13 0 10 3 10 9 13 12 16 9 16 3 10 9    /
\setlinear \plot 16 12 16 9 10 3    /
\put{$2{,}520$} [c] at 13 -2
\endpicture
\end{minipage}
$$
$$
\begin{minipage}{4cm}
\beginpicture
\setcoordinatesystem units   <1.5mm,2mm>
\setplotarea x from 0 to 16, y from -2 to 15
\put{73)} [l] at 2 12
\put {$ \scriptstyle \bullet$} [c] at 13 0
\put {$ \scriptstyle \bullet$} [c] at 13 12
\put {$ \scriptstyle \bullet$} [c] at 10 3
\put {$ \scriptstyle \bullet$} [c] at 10 9
\put {$ \scriptstyle \bullet$} [c] at 16 3
\put {$ \scriptstyle \bullet$} [c] at 16 9
\put {$ \scriptstyle \bullet$} [c] at 16 0
\setlinear \plot 16 0 16 3 13 0 10 3 10 9 13 12 16 9 16 3 10 9    /
\setlinear \plot 10 3 16 9    /
\put{$2{,}520$} [c] at 13 -2
\endpicture
\end{minipage}
\begin{minipage}{4cm}
\beginpicture
\setcoordinatesystem units   <1.5mm,2mm>
\setplotarea x from 0 to 16, y from -2 to 15
\put{74)} [l] at 2 12
\put {$ \scriptstyle \bullet$} [c] at  10 0
\put {$ \scriptstyle \bullet$} [c] at  16 0
\put {$ \scriptstyle \bullet$} [c] at  13 4
\put {$ \scriptstyle \bullet$} [c] at  10 12
\put {$ \scriptstyle \bullet$} [c] at  16 12
\put {$ \scriptstyle \bullet$} [c] at  11.2 9
\put {$ \scriptstyle \bullet$} [c] at  12.2  6
\setlinear \plot  10 0 13 4 16 12    /
\setlinear \plot  16 0 13 4 10 12 /
\put{$2{,}520  $} [c] at 13 -2
\endpicture
\end{minipage}
\begin{minipage}{4cm}
\beginpicture
\setcoordinatesystem units   <1.5mm,2mm>
\setplotarea x from 0 to 16, y from -2 to 15
\put{75)} [l] at 2 12
\put {$ \scriptstyle \bullet$} [c] at  10 0
\put {$ \scriptstyle \bullet$} [c] at  16 0
\put {$ \scriptstyle \bullet$} [c] at  13 7
\put {$ \scriptstyle \bullet$} [c] at  10 12
\put {$ \scriptstyle \bullet$} [c] at  16 12
\put {$ \scriptstyle \bullet$} [c] at  11.2 2.8
\put {$ \scriptstyle \bullet$} [c] at  12.2 5.1
\setlinear \plot  10 0 13 7 16 12    /
\setlinear \plot  16 0 13 7 10 12 /
\put{$2{,}520   $} [c] at 13 -2
\endpicture
\end{minipage}
\begin{minipage}{4cm}
\beginpicture
\setcoordinatesystem units   <1.5mm,2mm>
\setplotarea x from 0 to 16, y from -2 to 15
\put{76)} [l] at 2 12
\put {$ \scriptstyle \bullet$} [c] at  10 0
\put {$ \scriptstyle \bullet$} [c] at  16 0
\put {$ \scriptstyle \bullet$} [c] at  13 6
 \put {$ \scriptstyle \bullet$} [c] at  10 9
 \put {$ \scriptstyle \bullet$} [c] at  16 9
 \put {$ \scriptstyle \bullet$} [c] at  10 12
 \put {$ \scriptstyle \bullet$} [c] at  13 12
\setlinear \plot 10  0 13 6 10 9 10 12  /
\setlinear \plot 16  0 13 6 16 9 13 12 10 9   /
\put{$2{,}520 $} [c] at 13 -2
\endpicture
\end{minipage}
\begin{minipage}{4cm}
\beginpicture
\setcoordinatesystem units   <1.5mm,2mm>
\setplotarea x from 0 to 16, y from -2 to 15
\put{77)} [l] at 2 12
\put {$ \scriptstyle \bullet$} [c] at  10 0
\put {$ \scriptstyle \bullet$} [c] at  16 12
\put {$ \scriptstyle \bullet$} [c] at  13 6
 \put {$ \scriptstyle \bullet$} [c] at  10 3
 \put {$ \scriptstyle \bullet$} [c] at  16 3
 \put {$ \scriptstyle \bullet$} [c] at  10 12
 \put {$ \scriptstyle \bullet$} [c] at  13 0
\setlinear \plot 10  12 13 6 10 3 10 0  /
\setlinear \plot 16  12 13 6 16 3 13 0 10 3   /
\put{$2{,}520 $} [c] at 13 -2
\endpicture
\end{minipage}
\begin{minipage}{4cm}
\beginpicture
\setcoordinatesystem units   <1.5mm,2mm>
\setplotarea x from 0 to 16, y from -2 to 15
\put{78)} [l] at 2 12
\put {$ \scriptstyle \bullet$} [c] at  10 0
\put {$ \scriptstyle \bullet$} [c] at  10 6
\put {$ \scriptstyle \bullet$} [c] at  10 9
\put {$ \scriptstyle \bullet$} [c] at  10 12
\put {$ \scriptstyle \bullet$} [c] at  16 0
\put {$ \scriptstyle \bullet$} [c] at  16 6
\put {$ \scriptstyle \bullet$} [c] at  16 12
\setlinear \plot  10 12 10 0 16 6 16 12 10 6 16 0 16 6 10 9 /
\put{$1{,}260  $} [c] at 13 -2
\endpicture
\end{minipage}
$$

$$
\begin{minipage}{4cm}
\beginpicture
\setcoordinatesystem units   <1.5mm,2mm>
\setplotarea x from 0 to 16, y from -2 to 15
\put{79)} [l] at 2 12
\put {$ \scriptstyle \bullet$} [c] at  10 0
\put {$ \scriptstyle \bullet$} [c] at  10 6
\put {$ \scriptstyle \bullet$} [c] at  10 3
\put {$ \scriptstyle \bullet$} [c] at  10 12
\put {$ \scriptstyle \bullet$} [c] at  16 0
\put {$ \scriptstyle \bullet$} [c] at  16 6
\put {$ \scriptstyle \bullet$} [c] at  16 12
\setlinear \plot  10 0 10 12 16 6 16 0 10 6 16 12 16 6 10 3 /
\put{$1{,}260  $} [c] at 13 -2
\endpicture
\end{minipage}
\begin{minipage}{4cm}
\beginpicture
\setcoordinatesystem units   <1.5mm,2mm>
\setplotarea x from 0 to 16, y from -2 to 15
\put{80)} [l] at 2 12
\put {$ \scriptstyle \bullet$} [c] at  10  0
\put {$ \scriptstyle \bullet$} [c] at  10 3
\put {$ \scriptstyle \bullet$} [c] at  10 9
\put {$ \scriptstyle \bullet$} [c] at  10 12
\put {$ \scriptstyle \bullet$} [c] at  16 0
\put {$ \scriptstyle \bullet$} [c] at  16 6
\put {$ \scriptstyle \bullet$} [c] at  16 12
\setlinear \plot  10 0 10 12 16 6 16 12 10 9   /
\setlinear \plot  10 0 16 6 16 0 10 3    /
\put{$1{,}260  $} [c] at 13 -2
\endpicture
\end{minipage}
\begin{minipage}{4cm}
\beginpicture
\setcoordinatesystem units   <1.5mm,2mm>
\setplotarea x from 0 to 16, y from -2 to 15
\put{81)} [l] at 2 12
\put {$ \scriptstyle \bullet$} [c] at 13 0
\put {$ \scriptstyle \bullet$} [c] at 13 3
\put {$ \scriptstyle \bullet$} [c] at 13 6
\put {$ \scriptstyle \bullet$} [c] at 13 9
\put {$ \scriptstyle \bullet$} [c] at 13 12
\put {$ \scriptstyle \bullet$} [c] at 10 9
\put {$ \scriptstyle \bullet$} [c] at 16 9
\setlinear \plot    13 0 13 12   /
\setlinear \plot   13 6 10 9 13 12 16 9 13 6  /
\put{$840$} [c] at 13 -2
\endpicture
\end{minipage}
\begin{minipage}{4cm}
\beginpicture
\setcoordinatesystem units   <1.5mm,2mm>
\setplotarea x from 0 to 16, y from -2 to 15
\put{82)} [l] at 2 12
\put {$ \scriptstyle \bullet$} [c] at 13 0
\put {$ \scriptstyle \bullet$} [c] at 13 3
\put {$ \scriptstyle \bullet$} [c] at 13 6
\put {$ \scriptstyle \bullet$} [c] at 13 9
\put {$ \scriptstyle \bullet$} [c] at 13 12
\put {$ \scriptstyle \bullet$} [c] at 10 3
\put {$ \scriptstyle \bullet$} [c] at 16 3
\setlinear \plot    13 0 13 12   /
\setlinear \plot   13 6 10 3 13 0 16 3 13 6  /
\put{$840$} [c] at 13 -2
\endpicture
\end{minipage}
\begin{minipage}{4cm}
\beginpicture
\setcoordinatesystem units   <1.5mm,2mm>
\setplotarea x from 0 to 16, y from -2 to 15
\put{83)} [l] at 2 12
\put {$ \scriptstyle \bullet$} [c] at 13 0
\put {$ \scriptstyle \bullet$} [c] at 13 3
\put {$ \scriptstyle \bullet$} [c] at 13 6
\put {$ \scriptstyle \bullet$} [c] at 13 9
\put {$ \scriptstyle \bullet$} [c] at 13 12
\put {$ \scriptstyle \bullet$} [c] at 10 6
\put {$ \scriptstyle \bullet$} [c] at 16 6
\setlinear \plot    13 3 10 6 13 9 16 6 13 3  /
\setlinear \plot   13 0 13 12  /
\put{$840$} [c] at 13 -2
\endpicture
\end{minipage}
\begin{minipage}{4cm}
\beginpicture
\setcoordinatesystem units   <1.5mm,2mm>
\setplotarea x from 0 to 16, y from -2 to 15
\put{84)} [l] at 2 12
\put {$ \scriptstyle \bullet$} [c] at 13 0
\put {$ \scriptstyle \bullet$} [c] at 13 3
\put {$ \scriptstyle \bullet$} [c] at 13 6
\put {$ \scriptstyle \bullet$} [c] at 13 9
\put {$ \scriptstyle \bullet$} [c] at 13 12
\put {$ \scriptstyle \bullet$} [c] at 10 12
\put {$ \scriptstyle \bullet$} [c] at 16 12
\setlinear \plot  13 0 13 12     /
\setlinear \plot  10 12 13 9 16  12   /
\put{$840$} [c] at 13 -2
\endpicture
\end{minipage}
$$
$$
\begin{minipage}{4cm}
\beginpicture
\setcoordinatesystem units   <1.5mm,2mm>
\setplotarea x from 0 to 16, y from -2 to 15
\put{85)} [l] at 2 12
\put {$ \scriptstyle \bullet$} [c] at 13 0
\put {$ \scriptstyle \bullet$} [c] at 13 3
\put {$ \scriptstyle \bullet$} [c] at 13 6
\put {$ \scriptstyle \bullet$} [c] at 13 9
\put {$ \scriptstyle \bullet$} [c] at 13 12
\put {$ \scriptstyle \bullet$} [c] at 10 0
\put {$ \scriptstyle \bullet$} [c] at 16 0
\setlinear \plot  13 0 13 12     /
\setlinear \plot  10 0 13 3 16  0   /
\put{$840$} [c] at 13 -2
\endpicture
\end{minipage}
\begin{minipage}{4cm}
\beginpicture
\setcoordinatesystem units   <1.5mm,2mm>
\setplotarea x from 0 to 16, y from -2 to 15
\put{${\bf  13}$} [l] at 2 15
\put{86)} [l] at 2 12
\put {$ \scriptstyle \bullet$} [c] at 13 0
\put {$ \scriptstyle \bullet$} [c] at 10 3
\put {$ \scriptstyle \bullet$} [c] at 16 9
\put {$ \scriptstyle \bullet$} [c] at 13  12
\put {$ \scriptstyle \bullet$} [c] at 14 3
\put {$ \scriptstyle \bullet$} [c] at 11 6
\put {$ \scriptstyle \bullet$} [c] at 15 6
\setlinear \plot  13 0 10 3 13 12 16 9 13 0   /
\setlinear \plot  14 3 11 6   /
\put{$5{,}040$} [c] at 13 -2
\endpicture
\end{minipage}
\begin{minipage}{4cm}
\beginpicture
\setcoordinatesystem units   <1.5mm,2mm>
\setplotarea x from 0 to 16, y from -2 to 15
\put{87)} [l] at 2 12
\put {$ \scriptstyle \bullet$} [c] at 13 0
\put {$ \scriptstyle \bullet$} [c] at 10 3
\put {$ \scriptstyle \bullet$} [c] at 16 9
\put {$ \scriptstyle \bullet$} [c] at 13  12
\put {$ \scriptstyle \bullet$} [c] at 15 6
\put {$ \scriptstyle \bullet$} [c] at 11 6
\put {$ \scriptstyle \bullet$} [c] at 12 9
\setlinear \plot  13 0 10 3 13 12 16 9 13 0   /
\setlinear \plot  12 9 15 6   /
\put{$5{,}040$} [c] at 13 -2
\endpicture
\end{minipage}
\begin{minipage}{4cm}
\beginpicture
\setcoordinatesystem units   <1.5mm,2mm>
\setplotarea x from 0 to 16, y from -2 to 15
\put{88)} [l] at 2 12
\put {$ \scriptstyle \bullet$} [c] at 13  0
\put {$ \scriptstyle \bullet$} [c] at 13 3
\put {$ \scriptstyle \bullet$} [c] at 10 6
\put {$ \scriptstyle \bullet$} [c] at 16 6
\put {$ \scriptstyle \bullet$} [c] at 11.5 9
\put {$ \scriptstyle \bullet$} [c] at 13 12
\put {$ \scriptstyle \bullet$} [c] at 16 12
\setlinear \plot 13 0 13 3 10 6 13 12 16 6 13 3 /
\setlinear \plot 16 6 16 12 /
\put{$5{,}040$} [c] at 13 -2
\endpicture
\end{minipage}
\begin{minipage}{4cm}
\beginpicture
\setcoordinatesystem units   <1.5mm,2mm>
\setplotarea x from 0 to 16, y from -2 to 15
\put{89)} [l] at 2 12
\put {$ \scriptstyle \bullet$} [c] at 13  12
\put {$ \scriptstyle \bullet$} [c] at 13 9
\put {$ \scriptstyle \bullet$} [c] at 10 6
\put {$ \scriptstyle \bullet$} [c] at 16 6
\put {$ \scriptstyle \bullet$} [c] at 11.5 3
\put {$ \scriptstyle \bullet$} [c] at 13 0
\put {$ \scriptstyle \bullet$} [c] at 16 0
\setlinear \plot 13 12 13 9 10 6 13 0 16 6 13 9 /
\setlinear \plot 16 6 16 0 /
\put{$5{,}040$} [c] at 13 -2
\endpicture
\end{minipage}
\begin{minipage}{4cm}
\beginpicture
\setcoordinatesystem units   <1.5mm,2mm>
\setplotarea x from 0 to 16, y from -2 to 15
\put{90)} [l] at 2 12
\put {$ \scriptstyle \bullet$} [c] at 13 0
\put {$ \scriptstyle \bullet$} [c] at 13 3
\put {$ \scriptstyle \bullet$} [c] at 10 6
\put {$ \scriptstyle \bullet$} [c] at 16 6
\put {$ \scriptstyle \bullet$} [c] at 10 9
\put {$ \scriptstyle \bullet$} [c] at 10 12
\put {$ \scriptstyle \bullet$} [c] at 13 12
\setlinear \plot 13 0 13 3 10 6 10 12  /
\setlinear \plot 13 3 16 6 13 12 10 6 /
\put{$5{,}040$} [c] at 13 -2
\endpicture
\end{minipage}
$$
$$
\begin{minipage}{4cm}
\beginpicture
\setcoordinatesystem units   <1.5mm,2mm>
\setplotarea x from 0 to 16, y from -2 to 15
\put{91)} [l] at 2 12
\put {$ \scriptstyle \bullet$} [c] at 13 0
\put {$ \scriptstyle \bullet$} [c] at 13 10
\put {$ \scriptstyle \bullet$} [c] at 10 6
\put {$ \scriptstyle \bullet$} [c] at 16 6
\put {$ \scriptstyle \bullet$} [c] at 10 3
\put {$ \scriptstyle \bullet$} [c] at 10 0
\put {$ \scriptstyle \bullet$} [c] at 13 12
\setlinear \plot 13 12 13 10 10 6 10 0  /
\setlinear \plot 13 10 16 6 13 0 10 6 /
\put{$5{,}040$} [c] at 13 -2
\endpicture
\end{minipage}
\begin{minipage}{4cm}
\beginpicture
\setcoordinatesystem units   <1.5mm,2mm>
\setplotarea x from 0 to 16, y from -2 to 15
\put{92)} [l] at 2 12
\put {$ \scriptstyle \bullet$} [c] at 10 3
\put {$ \scriptstyle \bullet$} [c] at 10 6
\put {$ \scriptstyle \bullet$} [c] at 10 9
\put {$ \scriptstyle \bullet$} [c] at 10  12
\put {$ \scriptstyle \bullet$} [c] at 13 0
\put {$ \scriptstyle \bullet$} [c] at 16 12
\put {$ \scriptstyle \bullet$} [c] at 16 3
\setlinear \plot 10 12 10 3 13 0 16 3 16 12 10 9 /
\put{$5{,}040$} [c] at 13 -2
\endpicture
\end{minipage}
\begin{minipage}{4cm}
\beginpicture
\setcoordinatesystem units   <1.5mm,2mm>
\setplotarea x from 0 to 16, y from -2 to 15
\put{93)} [l] at 2 12
\put {$ \scriptstyle \bullet$} [c] at 10 3
\put {$ \scriptstyle \bullet$} [c] at 10 6
\put {$ \scriptstyle \bullet$} [c] at 10 9
\put {$ \scriptstyle \bullet$} [c] at 10  0
\put {$ \scriptstyle \bullet$} [c] at 13 12
\put {$ \scriptstyle \bullet$} [c] at 16 9
\put {$ \scriptstyle \bullet$} [c] at 16 0
\setlinear \plot 10 0 10 9 13 12 16 9 16 0 10 3 /
\put{$5{,}040$} [c] at 13 -2
\endpicture
\end{minipage}
\begin{minipage}{4cm}
\beginpicture
\setcoordinatesystem units   <1.5mm,2mm>
\setplotarea x from 0 to 16, y from -2 to 15
\put{94)} [l] at 2 12
\put {$ \scriptstyle \bullet$} [c] at 10 3
\put {$ \scriptstyle \bullet$} [c] at 16 3
\put {$ \scriptstyle \bullet$} [c] at 13 0
\put {$ \scriptstyle \bullet$} [c] at 13 6
\put {$ \scriptstyle \bullet$} [c] at 13 9
\put {$ \scriptstyle \bullet$} [c] at 13 12
\put {$ \scriptstyle \bullet$} [c] at 16 12
\setlinear \plot 16 12 16 3 13 0 10 3 13 6 13 12 /
\setlinear \plot 13 6 16 3 /
\put{$5{,}040$} [c] at 13 -2
\endpicture
\end{minipage}
\begin{minipage}{4cm}
\beginpicture
\setcoordinatesystem units   <1.5mm,2mm>
\setplotarea x from 0 to 16, y from -2 to 15
\put{95)} [l] at 2 12
\put {$ \scriptstyle \bullet$} [c] at 10 9
\put {$ \scriptstyle \bullet$} [c] at 16 9
\put {$ \scriptstyle \bullet$} [c] at 13 12
\put {$ \scriptstyle \bullet$} [c] at 13 6
\put {$ \scriptstyle \bullet$} [c] at 13 3
\put {$ \scriptstyle \bullet$} [c] at 13 0
\put {$ \scriptstyle \bullet$} [c] at 16 0
\setlinear \plot 16 0 16 9 13 12 10 9 13 6 13 0 /
\setlinear \plot 13 6 16 9 /
\put{$5{,}040$} [c] at 13 -2
\endpicture
\end{minipage}
\begin{minipage}{4cm}
\beginpicture
\setcoordinatesystem units   <1.5mm,2mm>
\setplotarea x from 0 to 16, y from -2 to 15
\put{96)} [l] at 2 12
\put {$ \scriptstyle \bullet$} [c] at 13 0
\put {$ \scriptstyle \bullet$} [c] at 10 3
\put {$ \scriptstyle \bullet$} [c] at 10 6
\put {$ \scriptstyle \bullet$} [c] at 10 9
\put {$ \scriptstyle \bullet$} [c] at 10 12
\put {$ \scriptstyle \bullet$} [c] at 16 12
\put {$ \scriptstyle \bullet$} [c] at 16 3
\setlinear \plot 10 12 10 3 13 0  16 3 16 12 10 3 /
\setlinear \plot 10 9 16 3 /
\put{$5{,}040$} [c] at 13 -2
\endpicture
\end{minipage}
$$
$$
\begin{minipage}{4cm}
\beginpicture
\setcoordinatesystem units   <1.5mm,2mm>
\setplotarea x from 0 to 16, y from -2 to 15
\put{97)} [l] at 2 12
\put {$ \scriptstyle \bullet$} [c] at 13 12
\put {$ \scriptstyle \bullet$} [c] at 10 3
\put {$ \scriptstyle \bullet$} [c] at 10 6
\put {$ \scriptstyle \bullet$} [c] at 10 9
\put {$ \scriptstyle \bullet$} [c] at 10 0
\put {$ \scriptstyle \bullet$} [c] at 16 9
\put {$ \scriptstyle \bullet$} [c] at 16 0
\setlinear \plot 10 0 10 9 13 12  16 9 16 0 10 9 /
\setlinear \plot 10 3 16 9 /
\put{$5{,}040$} [c] at 13 -2
\endpicture
\end{minipage}
\begin{minipage}{4cm}
\beginpicture
\setcoordinatesystem units   <1.5mm,2mm>
\setplotarea x from 0 to 16, y from -2 to 15
\put{98)} [l] at 2 12
\put {$ \scriptstyle \bullet$} [c] at 10 3
\put {$ \scriptstyle \bullet$} [c] at 10 6
\put {$ \scriptstyle \bullet$} [c] at 10 9
\put {$ \scriptstyle \bullet$} [c] at 10 12
\put {$ \scriptstyle \bullet$} [c] at 13 0
\put {$ \scriptstyle \bullet$} [c] at 16 12
\put {$ \scriptstyle \bullet$} [c] at 16 3
\setlinear \plot 13 0 10 3 10 12 16 3 16 12  10 6 /
\setlinear \plot 10 12 16 3 13 0 /
\put{$5{,}040$} [c] at 13 -2
\endpicture
\end{minipage}
\begin{minipage}{4cm}
\beginpicture
\setcoordinatesystem units   <1.5mm,2mm>
\setplotarea x from 0 to 16, y from -2 to 15
\put{99)} [l] at 2 12
\put {$ \scriptstyle \bullet$} [c] at 10 3
\put {$ \scriptstyle \bullet$} [c] at 10 6
\put {$ \scriptstyle \bullet$} [c] at 10 9
\put {$ \scriptstyle \bullet$} [c] at 10 0
\put {$ \scriptstyle \bullet$} [c] at 13 12
\put {$ \scriptstyle \bullet$} [c] at 16 0
\put {$ \scriptstyle \bullet$} [c] at 16 9
\setlinear \plot 13 12 10 9 10 0 16 9 16 0  10 6 /
\setlinear \plot 10 0 16 9 13 12 /
\put{$5{,}040$} [c] at 13 -2
\endpicture
\end{minipage}
\begin{minipage}{4cm}
\beginpicture
\setcoordinatesystem units   <1.5mm,2mm>
\setplotarea x from 0 to 16, y from -2 to 15
\put{100)} [l] at 2 12
\put {$ \scriptstyle \bullet$} [c] at 13 0
\put {$ \scriptstyle \bullet$} [c] at 11.5 3
\put {$ \scriptstyle \bullet$} [c] at 10 6
\put {$ \scriptstyle \bullet$} [c] at 13 12
\put {$ \scriptstyle \bullet$} [c] at 14.5 9
\put {$ \scriptstyle \bullet$} [c] at 16 6
\put {$ \scriptstyle \bullet$} [c] at 16 12
\setlinear \plot 13 0 10 6 13 12 16 6 13 0 /
\setlinear \plot 11.5 3 14.5 9 16 12 /
\put{$5{,}040$} [c] at 13 -2
\endpicture
\end{minipage}
\begin{minipage}{4cm}
\beginpicture
\setcoordinatesystem units   <1.5mm,2mm>
\setplotarea x from 0 to 16, y from -2 to 15
\put{101)} [l] at 2 12
\put {$ \scriptstyle \bullet$} [c] at 13 0
\put {$ \scriptstyle \bullet$} [c] at 11.5 9
\put {$ \scriptstyle \bullet$} [c] at 10 6
\put {$ \scriptstyle \bullet$} [c] at 13 12
\put {$ \scriptstyle \bullet$} [c] at 14.5 3
\put {$ \scriptstyle \bullet$} [c] at 16 6
\put {$ \scriptstyle \bullet$} [c] at 16 0
\setlinear \plot 13 12 10 6 13 0 16 6 13 12 /
\setlinear \plot 11.5 9 14.5 3 16 0 /
\put{$5{,}040$} [c] at 13 -2
\endpicture
\end{minipage}
\begin{minipage}{4cm}
\beginpicture
\setcoordinatesystem units   <1.5mm,2mm>
\setplotarea x from 0 to 16, y from -2 to 15
\put{102)} [l] at 2 12
\put {$ \scriptstyle \bullet$} [c] at 13 0
\put {$ \scriptstyle \bullet$} [c] at 10 3
\put {$ \scriptstyle \bullet$} [c] at 10 9
\put {$ \scriptstyle \bullet$} [c] at 10 12
\put {$ \scriptstyle \bullet$} [c] at 16 3
\put {$ \scriptstyle \bullet$} [c] at 16 6
\put {$ \scriptstyle \bullet$} [c] at 16 12
\setlinear \plot 13 0 10 3 10 12    /
\setlinear \plot 13 0 16 3 16 12  /
\setlinear \plot 10 9 16 6  /
\put{$5{,}040$} [c] at 13 -2
\endpicture
\end{minipage}
$$
$$
\begin{minipage}{4cm}
\beginpicture
\setcoordinatesystem units   <1.5mm,2mm>
\setplotarea x from 0 to 16, y from -2 to 15
\put{103)} [l] at 2 12
\put {$ \scriptstyle \bullet$} [c] at 13 12
\put {$ \scriptstyle \bullet$} [c] at 10 0
\put {$ \scriptstyle \bullet$} [c] at 10 3
\put {$ \scriptstyle \bullet$} [c] at 10 9
\put {$ \scriptstyle \bullet$} [c] at 16 9
\put {$ \scriptstyle \bullet$} [c] at 16 6
\put {$ \scriptstyle \bullet$} [c] at 16 0
\setlinear \plot 13 12 10 9 10 0    /
\setlinear \plot 13 12 16 9 16 0  /
\setlinear \plot 10 3 16 6  /
\put{$5{,}040$} [c] at 13 -2
\endpicture
\end{minipage}
\begin{minipage}{4cm}
\beginpicture
\setcoordinatesystem units   <1.5mm,2mm>
\setplotarea x from 0 to 16, y from -2 to 15
\put{104)} [l] at 2 12
\put {$ \scriptstyle \bullet$} [c] at 13 0
\put {$ \scriptstyle \bullet$} [c] at 10 3
\put {$ \scriptstyle \bullet$} [c] at 10 9
\put {$ \scriptstyle \bullet$} [c] at 10 12
\put {$ \scriptstyle \bullet$} [c] at 16  3
\put {$ \scriptstyle \bullet$} [c] at 16 6
\put {$ \scriptstyle \bullet$} [c] at 16 12
\setlinear \plot 16 6 10 12 10 3 13 0 16 3 16 12 10 3 /
\setlinear \plot 10 9 16 3 /
\put{$5{,}040$} [c] at 13 -2
\endpicture
\end{minipage}
\begin{minipage}{4cm}
\beginpicture
\setcoordinatesystem units   <1.5mm,2mm>
\setplotarea x from 0 to 16, y from -2 to 15
\put{105)} [l] at 2 12
\put {$ \scriptstyle \bullet$} [c] at 13 12
\put {$ \scriptstyle \bullet$} [c] at 10 9
\put {$ \scriptstyle \bullet$} [c] at 10 3
\put {$ \scriptstyle \bullet$} [c] at 10 0
\put {$ \scriptstyle \bullet$} [c] at 16  9
\put {$ \scriptstyle \bullet$} [c] at 16 6
\put {$ \scriptstyle \bullet$} [c] at 16 0
\setlinear \plot 16 6 10 0 10 9 13 12 16 9 16 0 10 9 /
\setlinear \plot 10 3 16 9 /
\put{$5{,}040$} [c] at 13 -2
\endpicture
\end{minipage}
\begin{minipage}{4cm}
\beginpicture
\setcoordinatesystem units   <1.5mm,2mm>
\setplotarea x from 0 to 16, y from -2 to 15
\put{106)} [l] at 2 12
\put {$ \scriptstyle \bullet$} [c] at 13 0
\put {$ \scriptstyle \bullet$} [c] at 10 12
\put {$ \scriptstyle \bullet$} [c] at 16 12
\put {$ \scriptstyle \bullet$} [c] at 12.3 2.3
\put {$ \scriptstyle \bullet$} [c] at 11.8 4.8
\put {$ \scriptstyle \bullet$} [c] at 11.2 7.2
\put {$ \scriptstyle \bullet$} [c] at 10.6 9.6
\setlinear \plot 10 12 13 0 16 12  /
\put{$5{,}040$} [c] at 13 -2
\endpicture
\end{minipage}
\begin{minipage}{4cm}
\beginpicture
\setcoordinatesystem units   <1.5mm,2mm>
\setplotarea x from 0 to 16, y from -2 to 15
\put{107)} [l] at 2 12
\put {$ \scriptstyle \bullet$} [c] at 13 12
\put {$ \scriptstyle \bullet$} [c] at 10 0
\put {$ \scriptstyle \bullet$} [c] at 16 0
\put {$ \scriptstyle \bullet$} [c] at 12.3 9.6
\put {$ \scriptstyle \bullet$} [c] at 11.8 7.2
\put {$ \scriptstyle \bullet$} [c] at 11.2 4.8
\put {$ \scriptstyle \bullet$} [c] at 10.6 2.3
\setlinear \plot 10 0 13 12 16 0  /
\put{$5{,}040$} [c] at 13 -2
\endpicture
\end{minipage}
\begin{minipage}{4cm}
\beginpicture
\setcoordinatesystem units   <1.5mm,2mm>
\setplotarea x from 0 to 16, y from -2 to 15
\put{108)} [l] at 2 12
\put {$ \scriptstyle \bullet$} [c] at  10  0
\put {$ \scriptstyle \bullet$} [c] at  10 3
\put {$ \scriptstyle \bullet$} [c] at  10 6
\put {$ \scriptstyle \bullet$} [c] at  10 9
\put {$ \scriptstyle \bullet$} [c] at  10 12
\put {$ \scriptstyle \bullet$} [c] at  16 0
\put {$ \scriptstyle \bullet$} [c] at  16 12
\setlinear \plot  10 0 10  12    /
\setlinear \plot  16  0 16 12  10 6 /
\setlinear \plot  10 9 16 0    /
\put{$5{,}040  $} [c] at 13 -2
\endpicture
\end{minipage}
$$
$$
\begin{minipage}{4cm}
\beginpicture
\setcoordinatesystem units   <1.5mm,2mm>
\setplotarea x from 0 to 16, y from -2 to 15
\put{109)} [l] at 2 12
\put {$ \scriptstyle \bullet$} [c] at  10  0
\put {$ \scriptstyle \bullet$} [c] at  10 3
\put {$ \scriptstyle \bullet$} [c] at  10 6
\put {$ \scriptstyle \bullet$} [c] at  10 9
\put {$ \scriptstyle \bullet$} [c] at  10 12
\put {$ \scriptstyle \bullet$} [c] at  16 0
\put {$ \scriptstyle \bullet$} [c] at  16 12
\setlinear \plot  10 0 10  12    /
\setlinear \plot  16  12 16 0  10 6 /
\setlinear \plot  10 3 16 12    /
\put{$5{,}040 $} [c] at 13 -2
\endpicture
\end{minipage}
\begin{minipage}{4cm}
\beginpicture
\setcoordinatesystem units   <1.5mm,2mm>
\setplotarea x from 0 to 16, y from -2 to 15
\put{110)} [l] at 2 12
\put {$ \scriptstyle \bullet$} [c] at  10 0
\put {$ \scriptstyle \bullet$} [c] at  10 3
\put {$ \scriptstyle \bullet$} [c] at  10 9
\put {$ \scriptstyle \bullet$} [c] at  10 12
\put {$ \scriptstyle \bullet$} [c] at  16 0
\put {$ \scriptstyle \bullet$} [c] at  16 6
\put {$ \scriptstyle \bullet$} [c] at  16 12
\setlinear \plot 10 12 10 9 16 0 16 12 10 9  10 3 16 6   /
\setlinear \plot 10 0 10 3  /
\put{$5{,}040  $} [c] at 13 -2
\endpicture
\end{minipage}
\begin{minipage}{4cm}
\beginpicture
\setcoordinatesystem units   <1.5mm,2mm>
\setplotarea x from 0 to 16, y from -2 to 15
\put{111)} [l] at 2 12
\put {$ \scriptstyle \bullet$} [c] at  10 0
\put {$ \scriptstyle \bullet$} [c] at  10 3
\put {$ \scriptstyle \bullet$} [c] at  10 9
\put {$ \scriptstyle \bullet$} [c] at  10 12
\put {$ \scriptstyle \bullet$} [c] at  16 0
\put {$ \scriptstyle \bullet$} [c] at  16 6
\put {$ \scriptstyle \bullet$} [c] at  16 12
\setlinear \plot 10 0 10 9 16 6 16 0 10 3 16 12 16 6   /
\setlinear \plot 10 9 10 12  /
\put{$5{,}040  $} [c] at 13 -2
\endpicture
\end{minipage}
\begin{minipage}{4cm}
\beginpicture
\setcoordinatesystem units   <1.5mm,2mm>
\setplotarea x from 0 to 16, y from -2 to 15
\put{112)} [l] at 2 12
\put {$ \scriptstyle \bullet$} [c] at 12 0
\put {$ \scriptstyle \bullet$} [c] at 12 3
\put {$ \scriptstyle \bullet$} [c] at 12 6
\put {$ \scriptstyle \bullet$} [c] at 10 9
\put {$ \scriptstyle \bullet$} [c] at 14 9
\put {$ \scriptstyle \bullet$} [c] at 16 9
\put {$ \scriptstyle \bullet$} [c] at 12 12
\setlinear \plot 12 0 12 6 10 9 12 12 14 9 12 6 /
\setlinear \plot  12 12 16 9 12 3  /
\put{$2{,}520$} [c] at 12 -2
\endpicture
\end{minipage}
\begin{minipage}{4cm}
\beginpicture
\setcoordinatesystem units   <1.5mm,2mm>
\setplotarea x from 0 to 16, y from -2 to 15
\put{113)} [l] at 2 12
\put {$ \scriptstyle \bullet$} [c] at 12 12
\put {$ \scriptstyle \bullet$} [c] at 12 9
\put {$ \scriptstyle \bullet$} [c] at 12 6
\put {$ \scriptstyle \bullet$} [c] at 10 3
\put {$ \scriptstyle \bullet$} [c] at 14 3
\put {$ \scriptstyle \bullet$} [c] at 16 3
\put {$ \scriptstyle \bullet$} [c] at 12 0
\setlinear \plot 12 12 12 6 10 3 12 0 14 3 12 6 /
\setlinear \plot  12 0 16 3 12 9  /
\put{$2{,}520$} [c] at 12 -2
\endpicture
\end{minipage}
\begin{minipage}{4cm}
\beginpicture
\setcoordinatesystem units   <1.5mm,2mm>
\setplotarea x from 0 to 16, y from -2 to 15
\put{114)} [l] at 2 12
\put {$ \scriptstyle \bullet$} [c] at 12 0
\put {$ \scriptstyle \bullet$} [c] at 12 3
\put {$ \scriptstyle \bullet$} [c] at 12 9
\put {$ \scriptstyle \bullet$} [c] at 12 12
\put {$ \scriptstyle \bullet$} [c] at 10 6
\put {$ \scriptstyle \bullet$} [c] at 14 6
\put {$ \scriptstyle \bullet$} [c] at 16 6
\setlinear \plot 12 0 12 3 10 6 12 9 14 6 12 3 /
\setlinear \plot  12 9 12 12 16 6 12 3  /
\put{$2{,}520$} [c] at 12 -2
\endpicture
\end{minipage}
$$
$$
\begin{minipage}{4cm}
\beginpicture
\setcoordinatesystem units   <1.5mm,2mm>
\setplotarea x from 0 to 16, y from -2 to 15
\put{115)} [l] at 2 12
\put {$ \scriptstyle \bullet$} [c] at 12 0
\put {$ \scriptstyle \bullet$} [c] at 12 3
\put {$ \scriptstyle \bullet$} [c] at 12 9
\put {$ \scriptstyle \bullet$} [c] at 12 12
\put {$ \scriptstyle \bullet$} [c] at 10 6
\put {$ \scriptstyle \bullet$} [c] at 14 6
\put {$ \scriptstyle \bullet$} [c] at 16 6
\setlinear \plot 12 0 12 3 10 6 12 9 14 6 12 3 /
\setlinear \plot  12 12 12 9 16 6 12 0  /
\put{$2{,}520$} [c] at 12 -2
\endpicture
\end{minipage}
\begin{minipage}{4cm}
\beginpicture
\setcoordinatesystem units   <1.5mm,2mm>
\setplotarea x from 0 to 16, y from -2 to 15
\put{116)} [l] at 2 12
\put {$ \scriptstyle \bullet$} [c] at 11 0
\put {$ \scriptstyle \bullet$} [c] at 11 3
\put {$ \scriptstyle \bullet$} [c] at 11 6
\put {$ \scriptstyle \bullet$} [c] at 11  12
\put {$ \scriptstyle \bullet$} [c] at 10 9
\put {$ \scriptstyle \bullet$} [c] at 12 9
\put {$ \scriptstyle \bullet$} [c] at 16 12
\setlinear \plot 11 0 11 6 10 9 11 12 12 9 11 6 /
\setlinear \plot  11 6 16 12  /
\put{$2{,}520$} [c] at 12 -2
\endpicture
\end{minipage}
\begin{minipage}{4cm}
\beginpicture
\setcoordinatesystem units   <1.5mm,2mm>
\setplotarea x from 0 to 16, y from -2 to 15
\put{117)} [l] at 2 12
\put {$ \scriptstyle \bullet$} [c] at 11 12
\put {$ \scriptstyle \bullet$} [c] at 11 9
\put {$ \scriptstyle \bullet$} [c] at 11 6
\put {$ \scriptstyle \bullet$} [c] at 11  0
\put {$ \scriptstyle \bullet$} [c] at 10 3
\put {$ \scriptstyle \bullet$} [c] at 12 3
\put {$ \scriptstyle \bullet$} [c] at 16 0
\setlinear \plot 11 12 11 6 10 3 11 0 12 3 11 6 /
\setlinear \plot  11 6 16 0  /
\put{$2{,}520$} [c] at 12 -2
\endpicture
\end{minipage}
\begin{minipage}{4cm}
\beginpicture
\setcoordinatesystem units   <1.5mm,2mm>
\setplotarea x from 0 to 16, y from -2 to 15
\put{118)} [l] at 2 12
\put {$ \scriptstyle \bullet$} [c] at 13 0
\put {$ \scriptstyle \bullet$} [c] at 13 3
\put {$ \scriptstyle \bullet$} [c] at 13 6
\put {$ \scriptstyle \bullet$} [c] at 13 9
\put {$ \scriptstyle \bullet$} [c] at 13 12
\put {$ \scriptstyle \bullet$} [c] at 10 12
\put {$ \scriptstyle \bullet$} [c] at 16 12
\setlinear \plot 13 0 13 12  /
\setlinear \plot 13 6 10 12 /
\setlinear \plot 13 9 16 12 /
\put{$2{,}520$} [c] at 13 -2
\endpicture
\end{minipage}
\begin{minipage}{4cm}
\beginpicture
\setcoordinatesystem units   <1.5mm,2mm>
\setplotarea x from 0 to 16, y from -2 to 15
\put{119)} [l] at 2 12
\put {$ \scriptstyle \bullet$} [c] at 13 0
\put {$ \scriptstyle \bullet$} [c] at 13 3
\put {$ \scriptstyle \bullet$} [c] at 13 6
\put {$ \scriptstyle \bullet$} [c] at 13 9
\put {$ \scriptstyle \bullet$} [c] at 13 12
\put {$ \scriptstyle \bullet$} [c] at 10 0
\put {$ \scriptstyle \bullet$} [c] at 16 0
\setlinear \plot 13 0 13 12  /
\setlinear \plot 13 6 10 0 /
\setlinear \plot 13 3 16 0 /
\put{$2{,}520$} [c] at 13 -2
\endpicture
\end{minipage}
\begin{minipage}{4cm}
\beginpicture
\setcoordinatesystem units   <1.5mm,2mm>
\setplotarea x from 0 to 16, y from -2 to 15
\put{120)} [l] at 2 12
\put {$ \scriptstyle \bullet$} [c] at  10 0
\put {$ \scriptstyle \bullet$} [c] at  10 3
\put {$ \scriptstyle \bullet$} [c] at  10 6
\put {$ \scriptstyle \bullet$} [c] at  10 9
\put {$ \scriptstyle \bullet$} [c] at  10 12
\put {$ \scriptstyle \bullet$} [c] at  16 0
\put {$ \scriptstyle \bullet$} [c] at  16 12
\setlinear \plot  10 12 10 3 16 0 16 12 10 0 10 3   /
\put{$2{,}520  $} [c] at 13 -2
\endpicture
\end{minipage}
$$
$$
\begin{minipage}{4cm}
\beginpicture
\setcoordinatesystem units   <1.5mm,2mm>
\setplotarea x from 0 to 16, y from -2 to 15
\put{121)} [l] at 2 12
\put {$ \scriptstyle \bullet$} [c] at  10 0
\put {$ \scriptstyle \bullet$} [c] at  10 3
\put {$ \scriptstyle \bullet$} [c] at  10 6
\put {$ \scriptstyle \bullet$} [c] at  10 9
\put {$ \scriptstyle \bullet$} [c] at  10 12
\put {$ \scriptstyle \bullet$} [c] at  16 0
\put {$ \scriptstyle \bullet$} [c] at  16 12
\setlinear \plot  10 0 10 9 16 12 16 0 10 12 10 9   /
\put{$2{,}520  $} [c] at 13 -2
\endpicture
\end{minipage}
\begin{minipage}{4cm}
\beginpicture
\setcoordinatesystem units   <1.5mm,2mm>
\setplotarea x from 0 to 16, y from -2 to 15
\put{122)} [l] at 2 12
\put {$ \scriptstyle \bullet$} [c] at  10 0
\put {$ \scriptstyle \bullet$} [c] at  16 6
\put {$ \scriptstyle \bullet$} [c] at  10 6
\put {$ \scriptstyle \bullet$} [c] at  10 9
\put {$ \scriptstyle \bullet$} [c] at  10 12
\put {$ \scriptstyle \bullet$} [c] at  16 0
\put {$ \scriptstyle \bullet$} [c] at  16 12
\setlinear \plot  10 12 10 0 16 6 16 0 10 6 16 12 16 6 10 12   /
\put{$2{,}520   $} [c] at 13 -2
\endpicture
\end{minipage}
\begin{minipage}{4cm}
\beginpicture
\setcoordinatesystem units   <1.5mm,2mm>
\setplotarea x from 0 to 16, y from -2 to 15
\put{123)} [l] at 2 12
\put {$ \scriptstyle \bullet$} [c] at  10 0
\put {$ \scriptstyle \bullet$} [c] at  16 6
\put {$ \scriptstyle \bullet$} [c] at  10 6
\put {$ \scriptstyle \bullet$} [c] at  10 3
\put {$ \scriptstyle \bullet$} [c] at  10 12
\put {$ \scriptstyle \bullet$} [c] at  16 0
\put {$ \scriptstyle \bullet$} [c] at  16 12
\setlinear \plot  10 12 10 0 16 6 16 0 10 6 16 12 16 6 10 12   /
\put{$2{,}520 $} [c] at 13 -2
\endpicture
\end{minipage}
\begin{minipage}{4cm}
\beginpicture
\setcoordinatesystem units   <1.5mm,2mm>
\setplotarea x from 0 to 16, y from -2 to 15
\put{124)} [l] at 2 12
\put {$ \scriptstyle \bullet$} [c] at  10 0
\put {$ \scriptstyle \bullet$} [c] at  16 6
\put {$ \scriptstyle \bullet$} [c] at  10 6
\put {$ \scriptstyle \bullet$} [c] at  10 9
\put {$ \scriptstyle \bullet$} [c] at  10 12
\put {$ \scriptstyle \bullet$} [c] at  16 0
\put {$ \scriptstyle \bullet$} [c] at  16 12
\setlinear \plot 10 12 10 0 16 6 16 12 10 9 /
\setlinear \plot 10 6 16 0 16 6 /
\put{$2{,}520   $} [c] at 13 -2
\endpicture
\end{minipage}
\begin{minipage}{4cm}
\beginpicture
\setcoordinatesystem units   <1.5mm,2mm>
\setplotarea x from 0 to 16, y from -2 to 15
\put{125)} [l] at 2 12
\put {$ \scriptstyle \bullet$} [c] at  10 0
\put {$ \scriptstyle \bullet$} [c] at  16 6
\put {$ \scriptstyle \bullet$} [c] at  10 6
\put {$ \scriptstyle \bullet$} [c] at  10 3
\put {$ \scriptstyle \bullet$} [c] at  10 12
\put {$ \scriptstyle \bullet$} [c] at  16 0
\put {$ \scriptstyle \bullet$} [c] at  16 12
\setlinear \plot 10 0 10 12 16 6 16 0 10 3 /
\setlinear \plot 10 6 16 12 16 6 /
\put{$2{,}520   $} [c] at 13 -2
\endpicture
\end{minipage}
\begin{minipage}{4cm}
\beginpicture
\setcoordinatesystem units   <1.5mm,2mm>
\setplotarea x from 0 to 16, y from -2 to 15
\put{126)} [l] at 2 12
\put {$ \scriptstyle \bullet$} [c] at  10 6
\put {$ \scriptstyle \bullet$} [c] at  10 12
\put {$ \scriptstyle \bullet$} [c] at  13 0
\put {$ \scriptstyle \bullet$} [c] at  13 3
\put {$ \scriptstyle \bullet$} [c] at  16 0
\put {$ \scriptstyle \bullet$} [c] at  16 6
\put {$ \scriptstyle \bullet$} [c] at  16 12
\setlinear \plot 13 0 13 3 10 6 10 12 16 6 16 12 10 6  /
\setlinear \plot 13 3 16 6 16 0  /
\put{$2{,}520$} [c] at 13 -2
\endpicture
\end{minipage}
$$
$$
\begin{minipage}{4cm}
\beginpicture
\setcoordinatesystem units   <1.5mm,2mm>
\setplotarea x from 0 to 16, y from -2 to 15
\put{127)} [l] at 2 12
\put {$ \scriptstyle \bullet$} [c] at  10 6
\put {$ \scriptstyle \bullet$} [c] at  10 0
\put {$ \scriptstyle \bullet$} [c] at  13 12
\put {$ \scriptstyle \bullet$} [c] at  13 9
\put {$ \scriptstyle \bullet$} [c] at  16 0
\put {$ \scriptstyle \bullet$} [c] at  16 6
\put {$ \scriptstyle \bullet$} [c] at  16 12
\setlinear \plot 13 12 13 9 10 6 10 0 16 6 16 0 10 6  /
\setlinear \plot 13 9 16 6 16 12  /
\put{$2{,}520$} [c] at 13 -2
\endpicture
\end{minipage}
\begin{minipage}{4cm}
\beginpicture
\setcoordinatesystem units   <1.5mm,2mm>
\setplotarea x from 0 to 16, y from -2 to 15
\put{128)} [l] at 2 12
\put {$ \scriptstyle \bullet$} [c] at 13 0
\put {$ \scriptstyle \bullet$} [c] at 10 3
\put {$ \scriptstyle \bullet$} [c] at 16 3
\put {$ \scriptstyle \bullet$} [c] at 10 6
\put {$ \scriptstyle \bullet$} [c] at 10 12
\put {$ \scriptstyle \bullet$} [c] at 16 6
\put {$ \scriptstyle \bullet$} [c] at 16 12
\setlinear \plot  16 6 16 3 13 0 10 3 10 12 16 6 16 12 10 6  /
\put{$1{,}260$} [c] at 13 -2
\endpicture
\end{minipage}
\begin{minipage}{4cm}
\beginpicture
\setcoordinatesystem units   <1.5mm,2mm>
\setplotarea x from 0 to 16, y from -2 to 15
\put{129)} [l] at 2 12
\put {$ \scriptstyle \bullet$} [c] at 13 12
\put {$ \scriptstyle \bullet$} [c] at 10 9
\put {$ \scriptstyle \bullet$} [c] at 16 9
\put {$ \scriptstyle \bullet$} [c] at 10 6
\put {$ \scriptstyle \bullet$} [c] at 10 0
\put {$ \scriptstyle \bullet$} [c] at 16 6
\put {$ \scriptstyle \bullet$} [c] at 16 0
\setlinear \plot  16 6 16 9 13 12 10 9 10 0 16 6 16 0 10 6  /
\put{$1{,}260$} [c] at 13 -2
\endpicture
\end{minipage}
\begin{minipage}{4cm}
\beginpicture
\setcoordinatesystem units   <1.5mm,2mm>
\setplotarea x from 0 to 16, y from -2 to 15
\put{130)} [l] at 2 12
\put {$ \scriptstyle \bullet$} [c] at 13 0
\put {$ \scriptstyle \bullet$} [c] at 10 3
\put {$ \scriptstyle \bullet$} [c] at 16 3
\put {$ \scriptstyle \bullet$} [c] at 10 6
\put {$ \scriptstyle \bullet$} [c] at 10 12
\put {$ \scriptstyle \bullet$} [c] at 16 6
\put {$ \scriptstyle \bullet$} [c] at 16 12
\setlinear \plot  10 12 10 3 13 0 16 3 16 12  /
\setlinear \plot  10 3 16 6  /
\setlinear \plot  10 6 16 3  /
\put{$1{,}260$} [c] at 13 -2
\endpicture
\end{minipage}
\begin{minipage}{4cm}
\beginpicture
\setcoordinatesystem units   <1.5mm,2mm>
\setplotarea x from 0 to 16, y from -2 to 15
\put{131)} [l] at 2 12
\put {$ \scriptstyle \bullet$} [c] at 13 12
\put {$ \scriptstyle \bullet$} [c] at 10 9
\put {$ \scriptstyle \bullet$} [c] at 16 9
\put {$ \scriptstyle \bullet$} [c] at 10 6
\put {$ \scriptstyle \bullet$} [c] at 10 0
\put {$ \scriptstyle \bullet$} [c] at 16 6
\put {$ \scriptstyle \bullet$} [c] at 16 0
\setlinear \plot  10 0 10 9 13 12 16 9 16 0  /
\setlinear \plot  10 9 16 6  /
\setlinear \plot  10 6 16 9  /
\put{$1{,}260$} [c] at 13 -2
\endpicture
\end{minipage}
\begin{minipage}{4cm}
\beginpicture
\setcoordinatesystem units   <1.5mm,2mm>
\setplotarea x from 0 to 16, y from -2 to 15
\put{132)} [l] at 2 12
\put {$ \scriptstyle \bullet$} [c] at  10 0
\put {$ \scriptstyle \bullet$} [c] at  10  9
\put {$ \scriptstyle \bullet$} [c] at  10 12
\put {$ \scriptstyle \bullet$} [c] at  13 6
\put {$ \scriptstyle \bullet$} [c] at  16 0
\put {$ \scriptstyle \bullet$} [c] at  16 9
\put {$ \scriptstyle \bullet$} [c] at  16 12
\setlinear \plot 10 0 13 6 16 9 16 12  /
\setlinear \plot 16 0 13 6 10 9 10 12  /
\put{$1{,}260 $} [c] at 13 -2
\endpicture
\end{minipage}
$$
$$
\begin{minipage}{4cm}
\beginpicture
\setcoordinatesystem units   <1.5mm,2mm>
\setplotarea x from 0 to 16, y from -2 to 15
\put{133)} [l] at 2 12
\put {$ \scriptstyle \bullet$} [c] at  10 0
\put {$ \scriptstyle \bullet$} [c] at  10  3
\put {$ \scriptstyle \bullet$} [c] at  10 12
\put {$ \scriptstyle \bullet$} [c] at  13 6
\put {$ \scriptstyle \bullet$} [c] at  16 0
\put {$ \scriptstyle \bullet$} [c] at  16 3
\put {$ \scriptstyle \bullet$} [c] at  16 12
\setlinear \plot 10 12 13 6 16 3 16 0  /
\setlinear \plot 16 12 13 6 10 3 10 0  /
\put{$1{,}260 $} [c] at 13 -2
\endpicture
\end{minipage}
\begin{minipage}{4cm}
\beginpicture
\setcoordinatesystem units   <1.5mm,2mm>
\setplotarea x from 0 to 16, y from -2 to 15
\put{134)} [l] at 2 12
\put {$ \scriptstyle \bullet$} [c] at 13 0
\put {$ \scriptstyle \bullet$} [c] at 13 3
\put {$ \scriptstyle \bullet$} [c] at 13 12
\put {$ \scriptstyle \bullet$} [c] at 10 3
\put {$ \scriptstyle \bullet$} [c] at 10 9
\put {$ \scriptstyle \bullet$} [c] at 16 3
\put {$ \scriptstyle \bullet$} [c] at 16 9
\setlinear \plot 13 3 13 0 10 3 10 9 13 12  16 9 16 3 13 0 /
\setlinear \plot 16 3 10 9 13 3 16 9 10 3  /
\put{$420$} [c] at 13 -2
\endpicture
\end{minipage}
\begin{minipage}{4cm}
\beginpicture
\setcoordinatesystem units   <1.5mm,2mm>
\setplotarea x from 0 to 16, y from -2 to 15
\put{135)} [l] at 2 12
\put {$ \scriptstyle \bullet$} [c] at 13 0
\put {$ \scriptstyle \bullet$} [c] at 13 9
\put {$ \scriptstyle \bullet$} [c] at 13 12
\put {$ \scriptstyle \bullet$} [c] at 10 3
\put {$ \scriptstyle \bullet$} [c] at 10 9
\put {$ \scriptstyle \bullet$} [c] at 16 3
\put {$ \scriptstyle \bullet$} [c] at 16 9
\setlinear \plot 13 9 13 12 10 9 10 3 13 0  16 3 16 9 13 12 /
\setlinear \plot 16 9 10 3 13 9 16 3 10 9  /
\put{$420$} [c] at 13 -2
\endpicture
\end{minipage}
\begin{minipage}{4cm}
\beginpicture
\setcoordinatesystem units   <1.5mm,2mm>
\setplotarea x from 0 to 16, y from -2 to 15
\put{136)} [l] at 2 12
\put {$ \scriptstyle \bullet$} [c] at 13 0
\put {$ \scriptstyle \bullet$} [c] at 10 3
\put {$ \scriptstyle \bullet$} [c] at 13 3
\put {$ \scriptstyle \bullet$} [c] at 16 3
\put {$ \scriptstyle \bullet$} [c] at 13 6
\put {$ \scriptstyle \bullet$} [c] at 10 12
\put {$ \scriptstyle \bullet$} [c] at 16 12
\setlinear \plot 10 12 13 6 10 3 13 0 13 6 16 3 13 0  /
\setlinear \plot 13 6 16 12 /
\put{$420$} [c] at 13 -2
\endpicture
\end{minipage}
\begin{minipage}{4cm}
\beginpicture
\setcoordinatesystem units   <1.5mm,2mm>
\setplotarea x from 0 to 16, y from -2 to 15
\put{137)} [l] at 2 12
\put {$ \scriptstyle \bullet$} [c] at 13 12
\put {$ \scriptstyle \bullet$} [c] at 10 9
\put {$ \scriptstyle \bullet$} [c] at 13 9
\put {$ \scriptstyle \bullet$} [c] at 16 9
\put {$ \scriptstyle \bullet$} [c] at 13 6
\put {$ \scriptstyle \bullet$} [c] at 10 0
\put {$ \scriptstyle \bullet$} [c] at 16 0
\setlinear \plot 10 0 13 6 10 9 13 12 13 6 16 9  13 12  /
\setlinear \plot 13 6 16 0 /
\put{$420$} [c] at 13 -2
\endpicture
\end{minipage}
\begin{minipage}{4cm}
\beginpicture
\setcoordinatesystem units   <1.5mm,2mm>
\setplotarea x from 0 to 16, y from -2 to 15
\put{138)} [l] at 2 12
\put {$ \scriptstyle \bullet$} [c] at 13 0
\put {$ \scriptstyle \bullet$} [c] at 13 3
\put {$ \scriptstyle \bullet$} [c] at 10 6
\put {$ \scriptstyle \bullet$} [c] at 13 6
\put {$ \scriptstyle \bullet$} [c] at 16 6
\put {$ \scriptstyle \bullet$} [c] at 10 12
\put {$ \scriptstyle \bullet$} [c] at 16 12
\setlinear \plot 13 0 13 3 10 6 10 12 16 6 16 12 13 6 13 3 16 6 /
\setlinear \plot  10 12 13 6 /
\setlinear \plot  16 12 10 6  /
\put{$420$} [c] at 13 -2
\endpicture
\end{minipage}
$$
$$
\begin{minipage}{4cm}
\beginpicture
\setcoordinatesystem units   <1.5mm,2mm>
\setplotarea x from 0 to 16, y from -2 to 15
\put{139)} [l] at 2 12
\put {$ \scriptstyle \bullet$} [c] at 13 12
\put {$ \scriptstyle \bullet$} [c] at 13 9
\put {$ \scriptstyle \bullet$} [c] at 10 6
\put {$ \scriptstyle \bullet$} [c] at 13 6
\put {$ \scriptstyle \bullet$} [c] at 16 6
\put {$ \scriptstyle \bullet$} [c] at 10 0
\put {$ \scriptstyle \bullet$} [c] at 16 0
\setlinear \plot 13 12  13 9 10 6 10 0 16 6 16 0 13 6 13 9 16 6 /
\setlinear \plot  10 0 13 6 /
\setlinear \plot  16 0 10 6  /
\put{$420$} [c] at 13 -2
\endpicture
\end{minipage}
\begin{minipage}{4cm}
\beginpicture
\setcoordinatesystem units   <1.5mm,2mm>
\setplotarea x from 0 to 16, y from -2 to 15
\put{140)} [l] at 2 12
\put {$ \scriptstyle \bullet$} [c] at 13 0
\put {$ \scriptstyle \bullet$} [c] at 13 3
\put {$ \scriptstyle \bullet$} [c] at 13 12
\put {$ \scriptstyle \bullet$} [c] at 10 6
\put {$ \scriptstyle \bullet$} [c] at 10 12
\put {$ \scriptstyle \bullet$} [c] at 16 6
\put {$ \scriptstyle \bullet$} [c] at 16 12
\setlinear \plot 13 0 13 3 10 6 10 12 16 6 16 12  10 6 13 12 16 6 /
\setlinear \plot 13 3 16 6 /
\put{$420$} [c] at 13 -2
\endpicture
\end{minipage}
\begin{minipage}{4cm}
\beginpicture
\setcoordinatesystem units   <1.5mm,2mm>
\setplotarea x from 0 to 16, y from -2 to 15
\put{141)} [l] at 2 12
\put {$ \scriptstyle \bullet$} [c] at 13 0
\put {$ \scriptstyle \bullet$} [c] at 13 9
\put {$ \scriptstyle \bullet$} [c] at 13 12
\put {$ \scriptstyle \bullet$} [c] at 10 6
\put {$ \scriptstyle \bullet$} [c] at 10 0
\put {$ \scriptstyle \bullet$} [c] at 16 6
\put {$ \scriptstyle \bullet$} [c] at 16 0
\setlinear \plot 13 12 13 9 10 6 10 0 16 6 16 0  10 6 13 0 16 6 /
\setlinear \plot 13 9 16 6 /
\put{$420$} [c] at 13 -2
\endpicture
\end{minipage}
\begin{minipage}{4cm}
\beginpicture
\setcoordinatesystem units   <1.5mm,2mm>
\setplotarea x from 0 to 16, y from -2 to 15
\put{142)} [l] at 2 12
\put {$ \scriptstyle \bullet$} [c] at 13 0
\put {$ \scriptstyle \bullet$} [c] at 10 3
\put {$ \scriptstyle \bullet$} [c] at 16 3
\put {$ \scriptstyle \bullet$} [c] at 13 6
\put {$ \scriptstyle \bullet$} [c] at 13 12
\put {$ \scriptstyle \bullet$} [c] at 10 12
\put {$ \scriptstyle \bullet$} [c] at 16 12
\setlinear \plot 10 12 13 6 10 3 13  0 16 3 13 6 /
\setlinear \plot 13 12 13 6 16 12   /
\put{$420$} [c] at 13 -2
\endpicture
\end{minipage}
\begin{minipage}{4cm}
\beginpicture
\setcoordinatesystem units   <1.5mm,2mm>
\setplotarea x from 0 to 16, y from -2 to 15
\put{143)} [l] at 2 12
\put {$ \scriptstyle \bullet$} [c] at 13 0
\put {$ \scriptstyle \bullet$} [c] at 10 9
\put {$ \scriptstyle \bullet$} [c] at 16 9
\put {$ \scriptstyle \bullet$} [c] at 13 6
\put {$ \scriptstyle \bullet$} [c] at 13 12
\put {$ \scriptstyle \bullet$} [c] at 10 0
\put {$ \scriptstyle \bullet$} [c] at 16 0
\setlinear \plot 10 0 13 6 10 9 13  12 16 9 13 6 /
\setlinear \plot 13 0 13 6 16 0   /
\put{$420$} [c] at 13 -2
\endpicture
\end{minipage}
\begin{minipage}{4cm}
\beginpicture
\setcoordinatesystem units   <1.5mm,2mm>
\setplotarea x from 0 to 16, y from -2 to 15
\put{144)} [l] at 2 12
\put {$ \scriptstyle \bullet$} [c] at 10 0
\put {$ \scriptstyle \bullet$} [c] at 16 0
\put {$ \scriptstyle \bullet$} [c] at 13 3
\put {$ \scriptstyle \bullet$} [c] at 13 9
\put {$ \scriptstyle \bullet$} [c] at 10 12
\put {$ \scriptstyle \bullet$} [c] at 13 12
\put {$ \scriptstyle \bullet$} [c] at 16 12
\setlinear \plot 10 0 13 3 13 9 13 12 /
\setlinear \plot 10 12 13 9 16 12   /
\setlinear \plot 16 0 13 3  /
\put{$420$} [c] at 13 -2
\endpicture
\end{minipage}
$$
$$
\begin{minipage}{4cm}
\beginpicture
\setcoordinatesystem units   <1.5mm,2mm>
\setplotarea x from 0 to 16, y from -2 to 15
\put{145)} [l] at 2 12
\put {$ \scriptstyle \bullet$} [c] at 10 12
\put {$ \scriptstyle \bullet$} [c] at 16 12
\put {$ \scriptstyle \bullet$} [c] at 13 3
\put {$ \scriptstyle \bullet$} [c] at 13 9
\put {$ \scriptstyle \bullet$} [c] at 10 0
\put {$ \scriptstyle \bullet$} [c] at 13 0
\put {$ \scriptstyle \bullet$} [c] at 16 0
\setlinear \plot 10 12 13 9 13 3 13 0 /
\setlinear \plot 10 0 13 3 16 0   /
\setlinear \plot 16 12 13 9  /
\put{$420$} [c] at 13 -2
\endpicture
\end{minipage}
\begin{minipage}{4cm}
\beginpicture
\setcoordinatesystem units   <1.5mm,2mm>
\setplotarea x from 0 to 16, y from -2 to 15
\put{${\bf  14}$} [l] at 2 15
\put{146)} [l] at 2 12
\put {$ \scriptstyle \bullet$} [c] at 10 3
\put {$ \scriptstyle \bullet$} [c] at 10 6
\put {$ \scriptstyle \bullet$} [c] at 10 9
\put {$ \scriptstyle \bullet$} [c] at 13 0
\put {$ \scriptstyle \bullet$} [c] at 13 12
\put {$ \scriptstyle \bullet$} [c] at 16 3
\put {$ \scriptstyle \bullet$} [c] at 16 9
\setlinear \plot 13 0 10 3 10 9 13 12 16 9 16 3 13 0  /
\put{$5{,}040$} [c] at 13 -2
\endpicture
\end{minipage}
\begin{minipage}{4cm}
\beginpicture
\setcoordinatesystem units   <1.5mm,2mm>
\setplotarea x from 0 to 16, y from -2 to 15
\put{147)} [l] at 2 12
\put {$ \scriptstyle \bullet$} [c] at 10 6
\put {$ \scriptstyle \bullet$} [c] at 11.5 9
\put {$ \scriptstyle \bullet$} [c] at 13 0
\put {$ \scriptstyle \bullet$} [c] at 13 6
\put {$ \scriptstyle \bullet$} [c] at 13 12
\put {$ \scriptstyle \bullet$} [c] at 14.5 3
\put {$ \scriptstyle \bullet$} [c] at 16 6
\setlinear \plot 13 0 10 6 13 12 16  6 13 0   /
\setlinear \plot 11.5 9 14.5 3 /
\put{$5{,}040$} [c] at 13 -2
\endpicture
\end{minipage}
\begin{minipage}{4cm}
\beginpicture
\setcoordinatesystem units   <1.5mm,2mm>
\setplotarea x from 0 to 16, y from -2 to 15
\put{148)} [l] at 2 12
\put {$ \scriptstyle \bullet$} [c] at 10 4
\put {$ \scriptstyle \bullet$} [c] at 10 8
\put {$ \scriptstyle \bullet$} [c] at 10 12
\put {$ \scriptstyle \bullet$} [c] at 13 0
\put {$ \scriptstyle \bullet$} [c] at 16 4
\put {$ \scriptstyle \bullet$} [c] at 16 8
\put {$ \scriptstyle \bullet$} [c] at 16 12
\setlinear \plot 10 12 10 4 13 0 16 4 16 12 10 4  /
\setlinear \plot 10 8 16 4 /
\put{$5{,}040$} [c] at 13 -2
\endpicture
\end{minipage}
\begin{minipage}{4cm}
\beginpicture
\setcoordinatesystem units   <1.5mm,2mm>
\setplotarea x from 0 to 16, y from -2 to 15
\put{149)} [l] at 2 12
\put {$ \scriptstyle \bullet$} [c] at 10 4
\put {$ \scriptstyle \bullet$} [c] at 10 8
\put {$ \scriptstyle \bullet$} [c] at 10 0
\put {$ \scriptstyle \bullet$} [c] at 13 12
\put {$ \scriptstyle \bullet$} [c] at 16 4
\put {$ \scriptstyle \bullet$} [c] at 16 8
\put {$ \scriptstyle \bullet$} [c] at 16 0
\setlinear \plot 10 0 10 8 13 12 16 8 16 0 10 8  /
\setlinear \plot 10 4 16 8 /
\put{$5{,}040$} [c] at 13 -2
\endpicture
\end{minipage}
\begin{minipage}{4cm}
\beginpicture
\setcoordinatesystem units   <1.5mm,2mm>
\setplotarea x from 0 to 16, y from -2 to 15
\put{150)} [l] at 2 12
\put {$ \scriptstyle \bullet$} [c] at 10 4
\put {$ \scriptstyle \bullet$} [c] at 10 8
\put {$ \scriptstyle \bullet$} [c] at 10 12
\put {$ \scriptstyle \bullet$} [c] at 13 0
\put {$ \scriptstyle \bullet$} [c] at 16 4
\put {$ \scriptstyle \bullet$} [c] at 16 8
\put {$ \scriptstyle \bullet$} [c] at 16 12
\setlinear \plot 10 12 10 4 13 0 16 4 16 12 10 8  /
\setlinear \plot 10 12 16 4 /
\put{$5{,}040$} [c] at 13 -2
\endpicture
\end{minipage}
$$
$$
\begin{minipage}{4cm}
\beginpicture
\setcoordinatesystem units   <1.5mm,2mm>
\setplotarea x from 0 to 16, y from -2 to 15
\put{151)} [l] at 2 12
\put {$ \scriptstyle \bullet$} [c] at 10 4
\put {$ \scriptstyle \bullet$} [c] at 10 8
\put {$ \scriptstyle \bullet$} [c] at 10 0
\put {$ \scriptstyle \bullet$} [c] at 13 12
\put {$ \scriptstyle \bullet$} [c] at 16 4
\put {$ \scriptstyle \bullet$} [c] at 16 8
\put {$ \scriptstyle \bullet$} [c] at 16 0
\setlinear \plot 10 0 10 8 13 12 16 8 16 0 10 4  /
\setlinear \plot 10 0 16 8 /
\put{$5{,}040$} [c] at 13 -2
\endpicture
\end{minipage}
\begin{minipage}{4cm}
\beginpicture
\setcoordinatesystem units   <1.5mm,2mm>
\setplotarea x from 0 to 16, y from -2 to 15
\put{152)} [l] at 2 12
\put {$ \scriptstyle \bullet$} [c] at 10 3
\put {$ \scriptstyle \bullet$} [c] at 10 6
\put {$ \scriptstyle \bullet$} [c] at 10 9
\put {$ \scriptstyle \bullet$} [c] at 10 12
\put {$ \scriptstyle \bullet$} [c] at 13 0
\put {$ \scriptstyle \bullet$} [c] at 13 12
\put {$ \scriptstyle \bullet$} [c] at 16 6
\setlinear \plot 10 12 10 3 13 0 16 6 13 12 10 6  /
\put{$5{,}040$} [c] at 13 -2
\endpicture
\end{minipage}
\begin{minipage}{4cm}
\beginpicture
\setcoordinatesystem units   <1.5mm,2mm>
\setplotarea x from 0 to 16, y from -2 to 15
\put{153)} [l] at 2 12
\put {$ \scriptstyle \bullet$} [c] at 10 3
\put {$ \scriptstyle \bullet$} [c] at 10 6
\put {$ \scriptstyle \bullet$} [c] at 10 9
\put {$ \scriptstyle \bullet$} [c] at 10 0
\put {$ \scriptstyle \bullet$} [c] at 13 0
\put {$ \scriptstyle \bullet$} [c] at 13 12
\put {$ \scriptstyle \bullet$} [c] at 16 6
\setlinear \plot 10 0 10 9 13 12 16 6 13 0 10 6  /
\put{$5{,}040$} [c] at 13 -2
\endpicture
\end{minipage}
\begin{minipage}{4cm}
\beginpicture
\setcoordinatesystem units   <1.5mm,2mm>
\setplotarea x from 0 to 16, y from -2 to 15
\put{154)} [l] at 2 12
\put {$ \scriptstyle \bullet$} [c] at 10 3
\put {$ \scriptstyle \bullet$} [c] at 11.5 6
\put {$ \scriptstyle \bullet$} [c] at 13 0
\put {$ \scriptstyle \bullet$} [c] at 13 9
\put {$ \scriptstyle \bullet$} [c] at 13 12
\put {$ \scriptstyle \bullet$} [c] at 16 3
\put {$ \scriptstyle \bullet$} [c] at 16 12
\setlinear \plot 16 12 16 3 13 9 10 3 13 0 16 3  /
\setlinear \plot 13  12 13  9  /
\put{$5{,}040$} [c] at 13 -2
\endpicture
\end{minipage}
\begin{minipage}{4cm}
\beginpicture
\setcoordinatesystem units   <1.5mm,2mm>
\setplotarea x from 0 to 16, y from -2 to 15
\put{155)} [l] at 2 12
\put {$ \scriptstyle \bullet$} [c] at 10 9
\put {$ \scriptstyle \bullet$} [c] at 11.5 6
\put {$ \scriptstyle \bullet$} [c] at 13 0
\put {$ \scriptstyle \bullet$} [c] at 13 3
\put {$ \scriptstyle \bullet$} [c] at 13 12
\put {$ \scriptstyle \bullet$} [c] at 16 9
\put {$ \scriptstyle \bullet$} [c] at 16 0
\setlinear \plot 16 0 16 9 13 3 10 9 13 12 16 9  /
\setlinear \plot 13  0 13  3  /
\put{$5{,}040$} [c] at 13 -2
\endpicture
\end{minipage}
\begin{minipage}{4cm}
\beginpicture
\setcoordinatesystem units   <1.5mm,2mm>
\setplotarea x from 0 to 16, y from -2 to 15
\put{156)} [l] at 2 12
\put {$ \scriptstyle \bullet$} [c] at 10 3
\put {$ \scriptstyle \bullet$} [c] at 10 6
\put {$ \scriptstyle \bullet$} [c] at 10 9
\put {$ \scriptstyle \bullet$} [c] at 10 12
\put {$ \scriptstyle \bullet$} [c] at 13 0
\put {$ \scriptstyle \bullet$} [c] at 16 3
\put {$ \scriptstyle \bullet$} [c] at 16 12
\setlinear \plot 16  3 10 12 10 3 13 0 16 3 16 12 10 3  /
\put{$5{,}040$} [c] at 13 -2
\endpicture
\end{minipage}
$$
$$
\begin{minipage}{4cm}
\beginpicture
\setcoordinatesystem units   <1.5mm,2mm>
\setplotarea x from 0 to 16, y from -2 to 15
\put{157)} [l] at 2 12
\put {$ \scriptstyle \bullet$} [c] at 10 3
\put {$ \scriptstyle \bullet$} [c] at 10 6
\put {$ \scriptstyle \bullet$} [c] at 10 9
\put {$ \scriptstyle \bullet$} [c] at 10 0
\put {$ \scriptstyle \bullet$} [c] at 13 12
\put {$ \scriptstyle \bullet$} [c] at 16 9
\put {$ \scriptstyle \bullet$} [c] at 16 0
\setlinear \plot 16  9 10 0 10 9 13 12 16 9 16 0 10 9  /
\put{$5{,}040$} [c] at 13 -2
\endpicture
\end{minipage}
\begin{minipage}{4cm}
\beginpicture
\setcoordinatesystem units   <1.5mm,2mm>
\setplotarea x from 0 to 16, y from -2 to 15
\put{158)} [l] at 2 12
\put {$ \scriptstyle \bullet$} [c] at 10 3
\put {$ \scriptstyle \bullet$} [c] at 12 0
\put {$ \scriptstyle \bullet$} [c] at 12 7.5
\put {$ \scriptstyle \bullet$} [c] at 14 4.5
\put {$ \scriptstyle \bullet$} [c] at 14 12
\put {$ \scriptstyle \bullet$} [c] at 16 12
\put {$ \scriptstyle \bullet$} [c] at 16 9
\setlinear \plot  12 0  10 3 14 12 16 9 12 0   /
\setlinear \plot 16 12  16 9   /
\setlinear \plot 14 4.5 12  7.5    /
\put{$5{,}040$} [c] at 13 -2
\endpicture
\end{minipage}
\begin{minipage}{4cm}
\beginpicture
\setcoordinatesystem units   <1.5mm,2mm>
\setplotarea x from 0 to 16, y from -2 to 15
\put{159)} [l] at 2 12
\put {$ \scriptstyle \bullet$} [c] at 10 9
\put {$ \scriptstyle \bullet$} [c] at 12 12
\put {$ \scriptstyle \bullet$} [c] at 12 4.5
\put {$ \scriptstyle \bullet$} [c] at 14 7.5
\put {$ \scriptstyle \bullet$} [c] at 14 0
\put {$ \scriptstyle \bullet$} [c] at 16 0
\put {$ \scriptstyle \bullet$} [c] at 16 3
\setlinear \plot  12 12  10 9 14 0 16 3 12 12   /
\setlinear \plot 16 0  16 3   /
\setlinear \plot 14 7.5 12  4.5    /
\put{$5{,}040$} [c] at 13 -2
\endpicture
\end{minipage}
\begin{minipage}{4cm}
\beginpicture
\setcoordinatesystem units   <1.5mm,2mm>
\setplotarea x from 0 to 16, y from -2 to 15
\put{160)} [l] at 2 12
\put {$ \scriptstyle \bullet$} [c] at 10 6
\put {$ \scriptstyle \bullet$} [c] at 10 9
\put {$ \scriptstyle \bullet$} [c] at 10 12
\put {$ \scriptstyle \bullet$} [c] at 13 0
\put {$ \scriptstyle \bullet$} [c] at 13 3
\put {$ \scriptstyle \bullet$} [c] at 16 12
\put {$ \scriptstyle \bullet$} [c] at 16 6
\setlinear \plot 10 12 10 6 13 3 16 6 16 12  /
\setlinear \plot 13 0 13 3 /
\put{$5{,}040$} [c] at 13 -2
\endpicture
\end{minipage}
\begin{minipage}{4cm}
\beginpicture
\setcoordinatesystem units   <1.5mm,2mm>
\setplotarea x from 0 to 16, y from -2 to 15
\put{161)} [l] at 2 12
\put {$ \scriptstyle \bullet$} [c] at 10 0
\put {$ \scriptstyle \bullet$} [c] at 10 3
\put {$ \scriptstyle \bullet$} [c] at 10 6
\put {$ \scriptstyle \bullet$} [c] at 13 12
\put {$ \scriptstyle \bullet$} [c] at 13 9
\put {$ \scriptstyle \bullet$} [c] at 16 0
\put {$ \scriptstyle \bullet$} [c] at 16 6
\setlinear \plot 10 0 10 6 13 9 16 6 16 0  /
\setlinear \plot 13 12 13 9 /
\put{$5{,}040$} [c] at 13 -2
\endpicture
\end{minipage}
\begin{minipage}{4cm}
\beginpicture
\setcoordinatesystem units   <1.5mm,2mm>
\setplotarea x from 0 to 16, y from -2 to 15
\put{162)} [l] at 2 12
\put {$ \scriptstyle \bullet$} [c] at 10 6
\put {$ \scriptstyle \bullet$} [c] at 10 12
\put {$ \scriptstyle \bullet$} [c] at 11.5  9
\put {$ \scriptstyle \bullet$} [c] at 13 0
\put {$ \scriptstyle \bullet$} [c] at 13 3
\put {$ \scriptstyle \bullet$} [c] at 13 12
\put {$ \scriptstyle \bullet$} [c] at 16 6
\setlinear \plot 10 12 10 6 13 12 16 6 13 3 10 6   /
\setlinear \plot 13 0 13 3 /
\put{$5{,}040$} [c] at 13 -2
\endpicture
\end{minipage}
$$
$$
\begin{minipage}{4cm}
\beginpicture
\setcoordinatesystem units   <1.5mm,2mm>
\setplotarea x from 0 to 16, y from -2 to 15
\put{163)} [l] at 2 12
\put {$ \scriptstyle \bullet$} [c] at 10 6
\put {$ \scriptstyle \bullet$} [c] at 10 0
\put {$ \scriptstyle \bullet$} [c] at 11.5  3
\put {$ \scriptstyle \bullet$} [c] at 13 0
\put {$ \scriptstyle \bullet$} [c] at 13 9
\put {$ \scriptstyle \bullet$} [c] at 13 12
\put {$ \scriptstyle \bullet$} [c] at 16 6
\setlinear \plot 10 0 10 6 13 0 16 6 13 9 10 6   /
\setlinear \plot 13 9 13 12 /
\put{$5{,}040$} [c] at 13 -2
\endpicture
\end{minipage}
\begin{minipage}{4cm}
\beginpicture
\setcoordinatesystem units   <1.5mm,2mm>
\setplotarea x from 0 to 16, y from -2 to 15
\put{164)} [l] at 2 12
\put {$ \scriptstyle \bullet$} [c] at  10 0
\put {$ \scriptstyle \bullet$} [c] at  10 3
\put {$ \scriptstyle \bullet$} [c] at  10 6
\put {$ \scriptstyle \bullet$} [c] at  10 9
\put {$ \scriptstyle \bullet$} [c] at  10 12
\put {$ \scriptstyle \bullet$} [c] at  16 0
\put {$ \scriptstyle \bullet$} [c] at  16 12
\setlinear \plot  10 12 10 0 16 12 16 0 10 6   /
\put{$5{,}040  $} [c] at 13 -2
\endpicture
\end{minipage}
\begin{minipage}{4cm}
\beginpicture
\setcoordinatesystem units   <1.5mm,2mm>
\setplotarea x from 0 to 16, y from -2 to 15
\put{165)} [l] at 2 12
\put {$ \scriptstyle \bullet$} [c] at  10 0
\put {$ \scriptstyle \bullet$} [c] at  10 3
\put {$ \scriptstyle \bullet$} [c] at  10 6
\put {$ \scriptstyle \bullet$} [c] at  10 9
\put {$ \scriptstyle \bullet$} [c] at  10 12
\put {$ \scriptstyle \bullet$} [c] at  16 0
\put {$ \scriptstyle \bullet$} [c] at  16 12
\setlinear \plot  10 0 10 12 16 0 16 12 10 6   /
\put{$5{,}040  $} [c] at 13 -2
\endpicture
\end{minipage}
\begin{minipage}{4cm}
\beginpicture
\setcoordinatesystem units   <1.5mm,2mm>
\setplotarea x from 0 to 16, y from -2 to 15
\put{166)} [l] at 2 12
\put {$ \scriptstyle \bullet$} [c] at  10 12
\put {$ \scriptstyle \bullet$} [c] at  11.5 9
\put {$ \scriptstyle \bullet$} [c] at  13 6
\put {$ \scriptstyle \bullet$} [c] at  13 3
\put {$ \scriptstyle \bullet$} [c] at  13 0
\put {$ \scriptstyle \bullet$} [c] at  16 12
\put {$ \scriptstyle \bullet$} [c] at  11.5 0
\setlinear \plot  13 0 13 6 16 12   /
\setlinear \plot  10 12  13 6 /
\setlinear \plot  11.5 0 11.5 9 /
\put{$5{,}040  $} [c] at 13 -2
\endpicture
\end{minipage}
\begin{minipage}{4cm}
\beginpicture
\setcoordinatesystem units   <1.5mm,2mm>
\setplotarea x from 0 to 16, y from -2 to 15
\put{167)} [l] at 2 12
\put {$ \scriptstyle \bullet$} [c] at  10 0
\put {$ \scriptstyle \bullet$} [c] at  11.5 3
\put {$ \scriptstyle \bullet$} [c] at  13 12
\put {$ \scriptstyle \bullet$} [c] at  13 9
\put {$ \scriptstyle \bullet$} [c] at  13 6
\put {$ \scriptstyle \bullet$} [c] at  16 0
\put {$ \scriptstyle \bullet$} [c] at  11.5 12
\setlinear \plot  13 12 13 6 16 0   /
\setlinear \plot  10 0  13 6 /
\setlinear \plot  11.5 3 11.5 12 /
\put{$5{,}040  $} [c] at 13 -2
\endpicture
\end{minipage}
\begin{minipage}{4cm}
\beginpicture
\setcoordinatesystem units   <1.5mm,2mm>
\setplotarea x from 0 to 16, y from -2 to 15
\put{168)} [l] at 2 12
\put {$ \scriptstyle \bullet$} [c] at  10 9
\put {$ \scriptstyle \bullet$} [c] at  11 0
\put {$ \scriptstyle \bullet$} [c] at  11 6
\put {$ \scriptstyle \bullet$} [c] at  11 12
\put {$ \scriptstyle \bullet$} [c] at  12 9
\put {$ \scriptstyle \bullet$} [c] at  16 0
\put {$ \scriptstyle \bullet$} [c] at  16 12
\setlinear \plot 10 9 16 0 16 12 12 9 11  12 10 9 11 6 12 9    /
\setlinear \plot  11 0 11 6  /
\put{$5{,}040  $} [c] at 13 -2
\endpicture
\end{minipage}
$$
$$
\begin{minipage}{4cm}
\beginpicture
\setcoordinatesystem units   <1.5mm,2mm>
\setplotarea x from 0 to 16, y from -2 to 15
\put{169)} [l] at 2 12
\put {$ \scriptstyle \bullet$} [c] at  10 3
\put {$ \scriptstyle \bullet$} [c] at  11 0
\put {$ \scriptstyle \bullet$} [c] at  11 6
\put {$ \scriptstyle \bullet$} [c] at  11 12
\put {$ \scriptstyle \bullet$} [c] at  12 3
\put {$ \scriptstyle \bullet$} [c] at  16 0
\put {$ \scriptstyle \bullet$} [c] at  16 12
\setlinear \plot 10 3 16 12 16 0 12 3 11  0 10 3 11 6 12 3    /
\setlinear \plot  11 12 11 6  /
\put{$5{,}040  $} [c] at 13 -2
\endpicture
\end{minipage}
\begin{minipage}{4cm}
\beginpicture
\setcoordinatesystem units   <1.5mm,2mm>
\setplotarea x from 0 to 16, y from -2 to 15
\put{170)} [l] at 2 12
\put {$ \scriptstyle \bullet$} [c] at  10 0
\put {$ \scriptstyle \bullet$} [c] at  10 6
\put {$ \scriptstyle \bullet$} [c] at  10 9
\put {$ \scriptstyle \bullet$} [c] at  10 12
\put {$ \scriptstyle \bullet$} [c] at  13 0
\put {$ \scriptstyle \bullet$} [c] at  16 6
\put {$ \scriptstyle \bullet$} [c] at  16 12
\setlinear \plot  10 0 10 12 16 6  16 12 10 6 13 0 16 6  /
\put{$5{,}040  $} [c] at 13 -2
\endpicture
\end{minipage}
\begin{minipage}{4cm}
\beginpicture
\setcoordinatesystem units   <1.5mm,2mm>
\setplotarea x from 0 to 16, y from -2 to 15
\put{171)} [l] at 2 12
\put {$ \scriptstyle \bullet$} [c] at  10 0
\put {$ \scriptstyle \bullet$} [c] at  10 3
\put {$ \scriptstyle \bullet$} [c] at  10 6
\put {$ \scriptstyle \bullet$} [c] at  10 12
\put {$ \scriptstyle \bullet$} [c] at  13 12
\put {$ \scriptstyle \bullet$} [c] at  16 6
\put {$ \scriptstyle \bullet$} [c] at  16 0
\setlinear \plot  10 12 10 0 16 6  16 0 10 6 13 12 16 6  /
\put{$5{,}040  $} [c] at 13 -2
\endpicture
\end{minipage}
\begin{minipage}{4cm}
\beginpicture
\setcoordinatesystem units   <1.5mm,2mm>
\setplotarea x from 0 to 16, y from -2 to 15
\put{172)} [l] at 2 12
\put {$ \scriptstyle \bullet$} [c] at  10 3
\put {$ \scriptstyle \bullet$} [c] at  13 0
\put {$ \scriptstyle \bullet$} [c] at  13 6
\put {$ \scriptstyle \bullet$} [c] at  13 12
\put {$ \scriptstyle \bullet$} [c] at  16 0
\put {$ \scriptstyle \bullet$} [c] at  16 3
\put {$ \scriptstyle \bullet$} [c] at  16 12
\setlinear \plot  13 12 13 6 16 3 13 0 10 3 13 6   /
\setlinear \plot 16 0  16 12 /
\put{$5{,}040  $} [c] at 13 -2
\endpicture
\end{minipage}
\begin{minipage}{4cm}
\beginpicture
\setcoordinatesystem units   <1.5mm,2mm>
\setplotarea x from 0 to 16, y from -2 to 15
\put{173)} [l] at 2 12
\put {$ \scriptstyle \bullet$} [c] at  10 9
\put {$ \scriptstyle \bullet$} [c] at  13 0
\put {$ \scriptstyle \bullet$} [c] at  13 6
\put {$ \scriptstyle \bullet$} [c] at  13 12
\put {$ \scriptstyle \bullet$} [c] at  16 0
\put {$ \scriptstyle \bullet$} [c] at  16 9
\put {$ \scriptstyle \bullet$} [c] at  16 12
\setlinear \plot  13 0 13 6 16 9 13 12 10 9 13 6   /
\setlinear \plot 16 0  16 12 /
\put{$5{,}040  $} [c] at 13 -2
\endpicture
\end{minipage}
\begin{minipage}{4cm}
\beginpicture
\setcoordinatesystem units   <1.5mm,2mm>
\setplotarea x from 0 to 16, y from -2 to 15
\put{174)} [l] at 2 12
\put {$ \scriptstyle \bullet$} [c] at  10 12
\put {$ \scriptstyle \bullet$} [c] at  13 0
\put {$ \scriptstyle \bullet$} [c] at  13 3
\put {$ \scriptstyle \bullet$} [c] at  13 6
\put {$ \scriptstyle \bullet$} [c] at  13 9
\put {$ \scriptstyle \bullet$} [c] at  16 0
\put {$ \scriptstyle \bullet$} [c] at  16 12
\setlinear \plot  13  0 13 9 10 12    /
\setlinear \plot  16  0 16 12  13 9 /
\put{$5{,}040  $} [c] at 13 -2
\endpicture
\end{minipage}
$$
$$
\begin{minipage}{4cm}
\beginpicture
\setcoordinatesystem units   <1.5mm,2mm>
\setplotarea x from 0 to 16, y from -2 to 15
\put{175)} [l] at 2 12
\put {$ \scriptstyle \bullet$} [c] at  10 0
\put {$ \scriptstyle \bullet$} [c] at  13 12
\put {$ \scriptstyle \bullet$} [c] at  13 9
\put {$ \scriptstyle \bullet$} [c] at  13 6
\put {$ \scriptstyle \bullet$} [c] at  13 3
\put {$ \scriptstyle \bullet$} [c] at  16 0
\put {$ \scriptstyle \bullet$} [c] at  16 12
\setlinear \plot  13  12 13 3 10 0    /
\setlinear \plot  16  12 16 0  13 3 /
\put{$5{,}040  $} [c] at 13 -2
\endpicture
\end{minipage}
\begin{minipage}{4cm}
\beginpicture
\setcoordinatesystem units   <1.5mm,2mm>
\setplotarea x from 0 to 16, y from -2 to 15
\put{176)} [l] at 2 12
\put {$ \scriptstyle \bullet$} [c] at  10 0
\put {$ \scriptstyle \bullet$} [c] at  10 3
\put {$ \scriptstyle \bullet$} [c] at  10 6
\put {$ \scriptstyle \bullet$} [c] at  10 9
\put {$ \scriptstyle \bullet$} [c] at  10 12
\put {$ \scriptstyle \bullet$} [c] at  16 0
\put {$ \scriptstyle \bullet$} [c] at  16 12
\setlinear \plot  10 0 10 12   /
\setlinear \plot 10 9 16 0 16 12 10 3 /
\put{$5{,}040  $} [c] at 13 -2
\endpicture
\end{minipage}
\begin{minipage}{4cm}
\beginpicture
\setcoordinatesystem units   <1.5mm,2mm>
\setplotarea x from 0 to 16, y from -2 to 15
\put{177)} [l] at 2 12
\put {$ \scriptstyle \bullet$} [c] at  10 0
\put {$ \scriptstyle \bullet$} [c] at  10 3
\put {$ \scriptstyle \bullet$} [c] at  10 9
\put {$ \scriptstyle \bullet$} [c] at  10 12
\put {$ \scriptstyle \bullet$} [c] at  13 0
\put {$ \scriptstyle \bullet$} [c] at  13 12
\put {$ \scriptstyle \bullet$} [c] at  16 6
\setlinear \plot  10 0  10 12   /
\setlinear \plot  10  9 13 12  16 6  13 0 10 3    /
\put{$5{,}040  $} [c] at 13 -2
\endpicture
\end{minipage}
\begin{minipage}{4cm}
\beginpicture
\setcoordinatesystem units   <1.5mm,2mm>
\setplotarea x from 0 to 16, y from -2 to 15
\put{178)} [l] at 2 12
\put {$ \scriptstyle \bullet$} [c] at  10 0
\put {$ \scriptstyle \bullet$} [c] at  10 12
\put {$ \scriptstyle \bullet$} [c] at  12 4
\put {$ \scriptstyle \bullet$} [c] at  12 8
\put {$ \scriptstyle \bullet$} [c] at  14 0
\put {$ \scriptstyle \bullet$} [c] at  14 12
\put {$ \scriptstyle \bullet$} [c] at  16 6
\setlinear \plot 12 4 10 12 16 6 10 0 12 8 12 4 14 0 16 6 14 12 12 8   /
\put{$5{,}040  $} [c] at 13 -2
\endpicture
\end{minipage}
\begin{minipage}{4cm}
\beginpicture
\setcoordinatesystem units   <1.5mm,2mm>
\setplotarea x from 0 to 16, y from -2 to 15
\put{179)} [l] at 2 12
\put {$ \scriptstyle \bullet$} [c] at 10 3
\put {$ \scriptstyle \bullet$} [c] at 10 9
\put {$ \scriptstyle \bullet$} [c] at 12 12
\put {$ \scriptstyle \bullet$} [c] at 12 0
\put {$ \scriptstyle \bullet$} [c] at 14 3
\put {$ \scriptstyle \bullet$} [c] at 14 9
\put {$ \scriptstyle \bullet$} [c] at 16 3
\setlinear \plot 14 9 16 3 12 0 14 3 14 9 12 12 10 9 10 3 12 0 /
\setlinear \plot 10 9 14 3  /
\setlinear \plot 14 9 10 3  /
\put{$2{,}520$} [c] at 13 -2
\endpicture
\end{minipage}
\begin{minipage}{4cm}
\beginpicture
\setcoordinatesystem units   <1.5mm,2mm>
\setplotarea x from 0 to 16, y from -2 to 15
\put{180)} [l] at 2 12
\put {$ \scriptstyle \bullet$} [c] at 10 3
\put {$ \scriptstyle \bullet$} [c] at 10 9
\put {$ \scriptstyle \bullet$} [c] at 12 12
\put {$ \scriptstyle \bullet$} [c] at 12 0
\put {$ \scriptstyle \bullet$} [c] at 14 3
\put {$ \scriptstyle \bullet$} [c] at 14 9
\put {$ \scriptstyle \bullet$} [c] at 16 9
\setlinear \plot 14 3 16 9 12 12 14 9 14 3 12 0 10 3 10 9 12 12 /
\setlinear \plot 10 9 14 3  /
\setlinear \plot 14 9 10 3  /
\put{$2{,}520$} [c] at 13 -2
\endpicture
\end{minipage}
$$
$$
\begin{minipage}{4cm}
\beginpicture
\setcoordinatesystem units   <1.5mm,2mm>
\setplotarea x from 0 to 16, y from -2 to 15
\put{181)} [l] at 2 12
\put {$ \scriptstyle \bullet$} [c] at 10 9
\put {$ \scriptstyle \bullet$} [c] at 11 0
\put {$ \scriptstyle \bullet$} [c] at 11 3
\put {$ \scriptstyle \bullet$} [c] at 11 6
\put {$ \scriptstyle \bullet$} [c] at 11 12
\put {$ \scriptstyle \bullet$} [c] at 12 9
\put {$ \scriptstyle \bullet$} [c] at 16 6
\setlinear \plot 11  0 11 6 10 9 11 12 12 9 11 6 /
\setlinear \plot 11 12 16 6 11 0 /
\put{$2{,}520$} [c] at 13 -2
\endpicture
\end{minipage}
\begin{minipage}{4cm}
\beginpicture
\setcoordinatesystem units   <1.5mm,2mm>
\setplotarea x from 0 to 16, y from -2 to 15
\put{182)} [l] at 2 12
\put {$ \scriptstyle \bullet$} [c] at 10 3
\put {$ \scriptstyle \bullet$} [c] at 11 12
\put {$ \scriptstyle \bullet$} [c] at 11 9
\put {$ \scriptstyle \bullet$} [c] at 11 6
\put {$ \scriptstyle \bullet$} [c] at 11 0
\put {$ \scriptstyle \bullet$} [c] at 12 3
\put {$ \scriptstyle \bullet$} [c] at 16 6
\setlinear \plot 11  12 11 6 10 3 11 0 12 3 11 6 /
\setlinear \plot 11 0 16 6 11 12 /
\put{$2{,}520$} [c] at 13 -2
\endpicture
\end{minipage}
\begin{minipage}{4cm}
\beginpicture
\setcoordinatesystem units   <1.5mm,2mm>
\setplotarea x from 0 to 16, y from -2 to 15
\put{183)} [l] at 2 12
\put {$ \scriptstyle \bullet$} [c] at 10 6
\put {$ \scriptstyle \bullet$} [c] at 12 0
\put {$ \scriptstyle \bullet$} [c] at 12 3
\put {$ \scriptstyle \bullet$} [c] at 12 9
\put {$ \scriptstyle \bullet$} [c] at 12 12
\put {$ \scriptstyle \bullet$} [c] at 14 6
\put {$ \scriptstyle \bullet$} [c] at 16 6
\setlinear \plot 12 12 16 6  12 0 12 3 10 6 12 9 12 12 /
\setlinear \plot 12 9 14 6  12 3 /
\put{$2{,}520$} [c] at 13 -2
\endpicture
\end{minipage}
\begin{minipage}{4cm}
\beginpicture
\setcoordinatesystem units   <1.5mm,2mm>
\setplotarea x from 0 to 16, y from -2 to 15
\put{184)} [l] at 2 12
\put {$ \scriptstyle \bullet$} [c] at 10 9
\put {$ \scriptstyle \bullet$} [c] at 11 0
\put {$ \scriptstyle \bullet$} [c] at 11 3
\put {$ \scriptstyle \bullet$} [c] at 11 6
\put {$ \scriptstyle \bullet$} [c] at 11 12
\put {$ \scriptstyle \bullet$} [c] at 12 9
\put {$ \scriptstyle \bullet$} [c] at 16 12
\setlinear \plot 11 0 11 6  10 9 11 12 12 9 11 6 /
\setlinear \plot 16 12  11 3 /
\put{$2{,}520$} [c] at 13 -2
\endpicture
\end{minipage}
\begin{minipage}{4cm}
\beginpicture
\setcoordinatesystem units   <1.5mm,2mm>
\setplotarea x from 0 to 16, y from -2 to 15
\put{185)} [l] at 2 12
\put {$ \scriptstyle \bullet$} [c] at 10 3
\put {$ \scriptstyle \bullet$} [c] at 11 0
\put {$ \scriptstyle \bullet$} [c] at 11 6
\put {$ \scriptstyle \bullet$} [c] at 11 9
\put {$ \scriptstyle \bullet$} [c] at 11 12
\put {$ \scriptstyle \bullet$} [c] at 12 3
\put {$ \scriptstyle \bullet$} [c] at 16 0
\setlinear \plot 11 12 11 6  10 3 11 0 12 3 11 6 /
\setlinear \plot 16 0  11 9 /
\put{$2{,}520$} [c] at 13 -2
\endpicture
\end{minipage}
\begin{minipage}{4cm}
\beginpicture
\setcoordinatesystem units   <1.5mm,2mm>
\setplotarea x from 0 to 16, y from -2 to 15
\put{186)} [l] at 2 12
\put {$ \scriptstyle \bullet$} [c] at 10 6
\put {$ \scriptstyle \bullet$} [c] at 11 0
\put {$ \scriptstyle \bullet$} [c] at 11 2
\put {$ \scriptstyle \bullet$} [c] at 11 10
\put {$ \scriptstyle \bullet$} [c] at 11 12
\put {$ \scriptstyle \bullet$} [c] at 12 6
\put {$ \scriptstyle \bullet$} [c] at 16 12
\setlinear \plot 11 12 11 10 10 6 11 2 11 0 /
\setlinear \plot 16 12 11 2 12 6 11 10 /
\put{$2{,}520$} [c] at 13 -2
\endpicture
\end{minipage}
$$
$$
\begin{minipage}{4cm}
\beginpicture
\setcoordinatesystem units   <1.5mm,2mm>
\setplotarea x from 0 to 16, y from -2 to 15
\put{187)} [l] at 2 12
\put {$ \scriptstyle \bullet$} [c] at 10 6
\put {$ \scriptstyle \bullet$} [c] at 11 0
\put {$ \scriptstyle \bullet$} [c] at 11 2
\put {$ \scriptstyle \bullet$} [c] at 11 10
\put {$ \scriptstyle \bullet$} [c] at 11 12
\put {$ \scriptstyle \bullet$} [c] at 12 6
\put {$ \scriptstyle \bullet$} [c] at 16 0
\setlinear \plot 11 12 11 10 10 6 11 2 11 0 /
\setlinear \plot 16 0 11 10 12 6 11 2 /
\put{$2{,}520$} [c] at 13 -2
\endpicture
\end{minipage}
\begin{minipage}{4cm}
\beginpicture
\setcoordinatesystem units   <1.5mm,2mm>
\setplotarea x from 0 to 16, y from -2 to 15
\put{188)} [l] at 2 12
\put {$ \scriptstyle \bullet$} [c] at 10 6
\put {$ \scriptstyle \bullet$} [c] at 10 12
\put {$ \scriptstyle \bullet$} [c] at 12 0
\put {$ \scriptstyle \bullet$} [c] at 12 3
\put {$ \scriptstyle \bullet$} [c] at 14 6
\put {$ \scriptstyle \bullet$} [c] at 14 12
\put {$ \scriptstyle \bullet$} [c] at 16 6
\setlinear \plot 12 0 12 3 16 6 14 12 14 6 10 12 10 6 14 12 /
\setlinear \plot 10 6 12 3 14 6 /
\put{$2{,}520$} [c] at 13 -2
\endpicture
\end{minipage}
\begin{minipage}{4cm}
\beginpicture
\setcoordinatesystem units   <1.5mm,2mm>
\setplotarea x from 0 to 16, y from -2 to 15
\put{189)} [l] at 2 12
\put {$ \scriptstyle \bullet$} [c] at 10 6
\put {$ \scriptstyle \bullet$} [c] at 10 0
\put {$ \scriptstyle \bullet$} [c] at 12 12
\put {$ \scriptstyle \bullet$} [c] at 12 9
\put {$ \scriptstyle \bullet$} [c] at 14 6
\put {$ \scriptstyle \bullet$} [c] at 14 0
\put {$ \scriptstyle \bullet$} [c] at 16 6
\setlinear \plot 12 12 12 9 16 6 14 0 14 6 10 0 10 6 14 0 /
\setlinear \plot 10 6 12 9 14 6 /
\put{$2{,}520$} [c] at 13 -2
\endpicture
\end{minipage}
\begin{minipage}{4cm}
\beginpicture
\setcoordinatesystem units   <1.5mm,2mm>
\setplotarea x from 0 to 16, y from -2 to 15
\put{190)} [l] at 2 12
\put {$ \scriptstyle \bullet$} [c] at 10 12
\put {$ \scriptstyle \bullet$} [c] at 13 12
\put {$ \scriptstyle \bullet$} [c] at 13 9
\put {$ \scriptstyle \bullet$} [c] at 13 6
\put {$ \scriptstyle \bullet$} [c] at 13 3
\put {$ \scriptstyle \bullet$} [c] at 13 0
\put {$ \scriptstyle \bullet$} [c] at 16 12
\setlinear \plot 13 12 13 0  /
\setlinear \plot 10 12  13 3   /
\setlinear \plot 16 12  13 9   /
\put{$2{,}520$} [c] at 13 -2
\endpicture
\end{minipage}
\begin{minipage}{4cm}
\beginpicture
\setcoordinatesystem units   <1.5mm,2mm>
\setplotarea x from 0 to 16, y from -2 to 15
\put{191)} [l] at 2 12
\put {$ \scriptstyle \bullet$} [c] at 10 0
\put {$ \scriptstyle \bullet$} [c] at 13 12
\put {$ \scriptstyle \bullet$} [c] at 13 9
\put {$ \scriptstyle \bullet$} [c] at 13 6
\put {$ \scriptstyle \bullet$} [c] at 13 3
\put {$ \scriptstyle \bullet$} [c] at 13 0
\put {$ \scriptstyle \bullet$} [c] at 16 0
\setlinear \plot 13 12 13 0  /
\setlinear \plot 10 0  13 9   /
\setlinear \plot 16 0  13 3   /
\put{$2{,}520$} [c] at 13 -2
\endpicture
\end{minipage}
\begin{minipage}{4cm}
\beginpicture
\setcoordinatesystem units   <1.5mm,2mm>
\setplotarea x from 0 to 16, y from -2 to 15
\put{192)} [l] at 2 12
\put {$ \scriptstyle \bullet$} [c] at 10 6
\put {$ \scriptstyle \bullet$} [c] at 10 12
\put {$ \scriptstyle \bullet$} [c] at 12 0
\put {$ \scriptstyle \bullet$} [c] at 12 3
\put {$ \scriptstyle \bullet$} [c] at 14 6
\put {$ \scriptstyle \bullet$} [c] at 14 12
\put {$ \scriptstyle \bullet$} [c] at 16 12
\setlinear \plot 16 12 14 6 12 3 12 0 /
\setlinear \plot 12 3 10 6 10 12 14 6 14 12 10 6 /
\put{$2{,}520$} [c] at 13 -2
\endpicture
\end{minipage}
$$
$$
\begin{minipage}{4cm}
\beginpicture
\setcoordinatesystem units   <1.5mm,2mm>
\setplotarea x from 0 to 16, y from -2 to 15
\put{193)} [l] at 2 12
\put {$ \scriptstyle \bullet$} [c] at 10 6
\put {$ \scriptstyle \bullet$} [c] at 10 0
\put {$ \scriptstyle \bullet$} [c] at 12 12
\put {$ \scriptstyle \bullet$} [c] at 12 9
\put {$ \scriptstyle \bullet$} [c] at 14 6
\put {$ \scriptstyle \bullet$} [c] at 14 0
\put {$ \scriptstyle \bullet$} [c] at 16 0
\setlinear \plot 16 0 14 6 12 9 12 12 /
\setlinear \plot 12 9 10 6 10 0 14 6 14 0 10 6 /
\put{$2{,}520$} [c] at 13 -2
\endpicture
\end{minipage}
\begin{minipage}{4cm}
\beginpicture
\setcoordinatesystem units   <1.5mm,2mm>
\setplotarea x from 0 to 16, y from -2 to 15
\put{194)} [l] at 2 12
\put {$ \scriptstyle \bullet$} [c] at  10 0
\put {$ \scriptstyle \bullet$} [c] at  10 3
\put {$ \scriptstyle \bullet$} [c] at  10 9
\put {$ \scriptstyle \bullet$} [c] at  10 12
\put {$ \scriptstyle \bullet$} [c] at  16 0
\put {$ \scriptstyle \bullet$} [c] at  16  6
\put {$ \scriptstyle \bullet$} [c] at  16 12
\setlinear \plot 10 12 10 0    /
\setlinear \plot  10 9 16 0 16 12     /
\setlinear \plot 10 3 16 6     /
\put{$2{,}520 $} [c] at 13 -2
\endpicture
\end{minipage}
\begin{minipage}{4cm}
\beginpicture
\setcoordinatesystem units   <1.5mm,2mm>
\setplotarea x from 0 to 16, y from -2 to 15
\put{195)} [l] at 2 12
\put {$ \scriptstyle \bullet$} [c] at  10 0
\put {$ \scriptstyle \bullet$} [c] at  10 3
\put {$ \scriptstyle \bullet$} [c] at  10 9
\put {$ \scriptstyle \bullet$} [c] at  10 12
\put {$ \scriptstyle \bullet$} [c] at  16 0
\put {$ \scriptstyle \bullet$} [c] at  16  6
\put {$ \scriptstyle \bullet$} [c] at  16 12
\setlinear \plot 10 12 10 0    /
\setlinear \plot  10 3 16 12 16 0     /
\setlinear \plot 10 9 16 6     /
\put{$2{,}520 $} [c] at 13 -2
\endpicture
\end{minipage}
\begin{minipage}{4cm}
\beginpicture
\setcoordinatesystem units   <1.5mm,2mm>
\setplotarea x from 0 to 16, y from -2 to 15
\put{196)} [l] at 2 12
\put {$ \scriptstyle \bullet$} [c] at  10 0
\put {$ \scriptstyle \bullet$} [c] at  10 6
\put {$ \scriptstyle \bullet$} [c] at  10 9
\put {$ \scriptstyle \bullet$} [c] at  10 12
\put {$ \scriptstyle \bullet$} [c] at  13 12
\put {$ \scriptstyle \bullet$} [c] at  16 0
\put {$ \scriptstyle \bullet$} [c] at  16 6
\setlinear \plot  10 12  10 6 16 0 16 6 13 12 10 6  10 0 16 6  /
\put{$2{,}520  $} [c] at 13 -2
\endpicture
\end{minipage}
\begin{minipage}{4cm}
\beginpicture
\setcoordinatesystem units   <1.5mm,2mm>
\setplotarea x from 0 to 16, y from -2 to 15
\put{197)} [l] at 2 12
\put {$ \scriptstyle \bullet$} [c] at  10 0
\put {$ \scriptstyle \bullet$} [c] at  10 6
\put {$ \scriptstyle \bullet$} [c] at  10 3
\put {$ \scriptstyle \bullet$} [c] at  10 12
\put {$ \scriptstyle \bullet$} [c] at  13 0
\put {$ \scriptstyle \bullet$} [c] at  16 12
\put {$ \scriptstyle \bullet$} [c] at  16 6
\setlinear \plot  10 0  10 6 16 12 16 6 13 0 10 6  10 12 16 6  /
\put{$2{,}520  $} [c] at 13 -2
\endpicture
\end{minipage}
\begin{minipage}{4cm}
\beginpicture
\setcoordinatesystem units   <1.5mm,2mm>
\setplotarea x from 0 to 16, y from -2 to 15
\put{198)} [l] at 2 12
\put {$ \scriptstyle \bullet$} [c] at  10 0
\put {$ \scriptstyle \bullet$} [c] at  10 6
\put {$ \scriptstyle \bullet$} [c] at  13 12
\put {$ \scriptstyle \bullet$} [c] at  16  0
\put {$ \scriptstyle \bullet$} [c] at  16 6
\put {$ \scriptstyle \bullet$} [c] at  16 12
\put {$ \scriptstyle \bullet$} [c] at  11.5 9
\setlinear \plot 16 12 16 0 10 6 13 12 16 6 10 0 10 6  /
\put{$2{,}520   $} [c] at 13 -2
\endpicture
\end{minipage}
$$
$$
\begin{minipage}{4cm}
\beginpicture
\setcoordinatesystem units   <1.5mm,2mm>
\setplotarea x from 0 to 16, y from -2 to 15
\put{199)} [l] at 2 12
\put {$ \scriptstyle \bullet$} [c] at  10 12
\put {$ \scriptstyle \bullet$} [c] at  10 6
\put {$ \scriptstyle \bullet$} [c] at  13 0
\put {$ \scriptstyle \bullet$} [c] at  16  0
\put {$ \scriptstyle \bullet$} [c] at  16 6
\put {$ \scriptstyle \bullet$} [c] at  16 12
\put {$ \scriptstyle \bullet$} [c] at  11.5 3
\setlinear \plot 16 0 16 12 10 6 13 0 16 6 10 12 10 6  /
\put{$2{,}520   $} [c] at 13 -2
\endpicture
\end{minipage}
\begin{minipage}{4cm}
\beginpicture
\setcoordinatesystem units   <1.5mm,2mm>
\setplotarea x from 0 to 16, y from -2 to 15
\put{200)} [l] at 2 12
\put {$ \scriptstyle \bullet$} [c] at 10 6
\put {$ \scriptstyle \bullet$} [c] at 10 12
\put {$ \scriptstyle \bullet$} [c] at 12 0
\put {$ \scriptstyle \bullet$} [c] at 12 3
\put {$ \scriptstyle \bullet$} [c] at 14 6
\put {$ \scriptstyle \bullet$} [c] at 14 12
\put {$ \scriptstyle \bullet$} [c] at 16 6
\setlinear \plot 16 6 12 0 12 3 10 6 10 12 16 6 14 12 14 6 12 3 /
\setlinear \plot 10 12 14  6   /
\setlinear \plot 10 6 14  12   /
\put{$1{,}260$} [c] at 13 -2
\endpicture
\end{minipage}
\begin{minipage}{4cm}
\beginpicture
\setcoordinatesystem units   <1.5mm,2mm>
\setplotarea x from 0 to 16, y from -2 to 15
\put{201)} [l] at 2 12
\put {$ \scriptstyle \bullet$} [c] at 10 6
\put {$ \scriptstyle \bullet$} [c] at 10 0
\put {$ \scriptstyle \bullet$} [c] at 12 12
\put {$ \scriptstyle \bullet$} [c] at 12 9
\put {$ \scriptstyle \bullet$} [c] at 14 6
\put {$ \scriptstyle \bullet$} [c] at 14 0
\put {$ \scriptstyle \bullet$} [c] at 16 6
\setlinear \plot 16 6 12 12 12 9 10 6 10 0 16 6 14 0 14 6 12 9 /
\setlinear \plot 10 0 14  6   /
\setlinear \plot 10 6 14  0   /
\put{$1{,}260$} [c] at 13 -2
\endpicture
\end{minipage}
\begin{minipage}{4cm}
\beginpicture
\setcoordinatesystem units   <1.5mm,2mm>
\setplotarea x from 0 to 16, y from -2 to 15
\put{202)} [l] at 2 12
\put {$ \scriptstyle \bullet$} [c] at 10 3
\put {$ \scriptstyle \bullet$} [c] at 10 9
\put {$ \scriptstyle \bullet$} [c] at 11 0
\put {$ \scriptstyle \bullet$} [c] at 11 12
\put {$ \scriptstyle \bullet$} [c] at 12 3
\put {$ \scriptstyle \bullet$} [c] at 12 9
\put {$ \scriptstyle \bullet$} [c] at 16 12
\setlinear \plot 11 0 10 3 10 9 11 12 12 9 12 3 11 0  /
\setlinear \plot 10 3 12 9  /
\setlinear \plot 10 9 12 3  /
\setlinear \plot 10 3 16 12 12 3  /
\put{$1{,}260$} [c] at 13 -2
\endpicture
\end{minipage}
\begin{minipage}{4cm}
\beginpicture
\setcoordinatesystem units   <1.5mm,2mm>
\setplotarea x from 0 to 16, y from -2 to 15
\put{203)} [l] at 2 12
\put {$ \scriptstyle \bullet$} [c] at 10 3
\put {$ \scriptstyle \bullet$} [c] at 10 9
\put {$ \scriptstyle \bullet$} [c] at 11 0
\put {$ \scriptstyle \bullet$} [c] at 11 12
\put {$ \scriptstyle \bullet$} [c] at 12 3
\put {$ \scriptstyle \bullet$} [c] at 12 9
\put {$ \scriptstyle \bullet$} [c] at 16 0
\setlinear \plot 11 0 10 3 10 9 11 12 12 9 12 3 11 0  /
\setlinear \plot 10 3 12 9  /
\setlinear \plot 10 9 12 3  /
\setlinear \plot 10 9 16 0 12 9  /
\put{$1{,}260$} [c] at 13 -2
\endpicture
\end{minipage}
\begin{minipage}{4cm}
\beginpicture
\setcoordinatesystem units   <1.5mm,2mm>
\setplotarea x from 0 to 16, y from -2 to 15
\put{204)} [l] at 2 12
\put {$ \scriptstyle \bullet$} [c] at 10 3
\put {$ \scriptstyle \bullet$} [c] at 10 6
\put {$ \scriptstyle \bullet$} [c] at 10 12
\put {$ \scriptstyle \bullet$} [c] at 13 0
\put {$ \scriptstyle \bullet$} [c] at 13 3
\put {$ \scriptstyle \bullet$} [c] at 16 3
\put {$ \scriptstyle \bullet$} [c] at 16 12
\setlinear \plot  10 12 10 3 13 0 16 3 16 12 10 6  /
\setlinear \plot  10 12  16 3   /
\setlinear \plot  10 6 13 3 13 0  /
\put{$1{,}260$} [c] at 13 -2
\endpicture
\end{minipage}
$$
$$
\begin{minipage}{4cm}
\beginpicture
\setcoordinatesystem units   <1.5mm,2mm>
\setplotarea x from 0 to 16, y from -2 to 15
\put{205)} [l] at 2 12
\put {$ \scriptstyle \bullet$} [c] at 10 0
\put {$ \scriptstyle \bullet$} [c] at 10 6
\put {$ \scriptstyle \bullet$} [c] at 10 9
\put {$ \scriptstyle \bullet$} [c] at 13 9
\put {$ \scriptstyle \bullet$} [c] at 13 12
\put {$ \scriptstyle \bullet$} [c] at 16 0
\put {$ \scriptstyle \bullet$} [c] at 16 9
\setlinear \plot  10 0 10 9 13 12 16 9 16 0 10 6  /
\setlinear \plot  10 0  16 9   /
\setlinear \plot  10 6 13 9 13 12  /
\put{$1{,}260$} [c] at 13 -2
\endpicture
\end{minipage}
\begin{minipage}{4cm}
\beginpicture
\setcoordinatesystem units   <1.5mm,2mm>
\setplotarea x from 0 to 16, y from -2 to 15
\put{206)} [l] at 2 12
\put {$ \scriptstyle \bullet$} [c] at 10 3
\put {$ \scriptstyle \bullet$} [c] at 10 6
\put {$ \scriptstyle \bullet$} [c] at 10 12
\put {$ \scriptstyle \bullet$} [c] at 13 12
\put {$ \scriptstyle \bullet$} [c] at 13 0
\put {$ \scriptstyle \bullet$} [c] at 16 3
\put {$ \scriptstyle \bullet$} [c] at 16 12
\setlinear \plot  10 12 10 3 13 0 16 3 16 12 10 3  /
\setlinear \plot  13 12 10 6  16 3   /
\put{$1{,}260$} [c] at 13 -2
\endpicture
\end{minipage}
\begin{minipage}{4cm}
\beginpicture
\setcoordinatesystem units   <1.5mm,2mm>
\setplotarea x from 0 to 16, y from -2 to 15
\put{207)} [l] at 2 12
\put {$ \scriptstyle \bullet$} [c] at 10 9
\put {$ \scriptstyle \bullet$} [c] at 10 6
\put {$ \scriptstyle \bullet$} [c] at 10 0
\put {$ \scriptstyle \bullet$} [c] at 13 12
\put {$ \scriptstyle \bullet$} [c] at 13 0
\put {$ \scriptstyle \bullet$} [c] at 16 9
\put {$ \scriptstyle \bullet$} [c] at 16 0
\setlinear \plot  10 0 10 9 13 12 16 9 16 0 10 9  /
\setlinear \plot  13 0 10 6  16 9   /
\put{$1{,}260$} [c] at 13 -2
\endpicture
\end{minipage}
\begin{minipage}{4cm}
\beginpicture
\setcoordinatesystem units   <1.5mm,2mm>
\setplotarea x from 0 to 16, y from -2 to 15
\put{208)} [l] at 2 12
\put {$ \scriptstyle \bullet$} [c] at  10 12
\put {$ \scriptstyle \bullet$} [c] at  10 3
\put {$ \scriptstyle \bullet$} [c] at  11 6
\put {$ \scriptstyle \bullet$} [c] at  11  0
\put {$ \scriptstyle \bullet$} [c] at  12 3
\put {$ \scriptstyle \bullet$} [c] at  12 12
\put {$ \scriptstyle \bullet$} [c] at  16 0
\setlinear \plot  10 12 11 6 10 3 11 0 12 3 11 6 12 12   /
\setlinear \plot  16 0 11  6  /
\put{$1{,}260  $} [c] at 13 -2
\endpicture
\end{minipage}
\begin{minipage}{4cm}
\beginpicture
\setcoordinatesystem units   <1.5mm,2mm>
\setplotarea x from 0 to 16, y from -2 to 15
\put{209)} [l] at 2 12
\put {$ \scriptstyle \bullet$} [c] at  10 0
\put {$ \scriptstyle \bullet$} [c] at  10 9
\put {$ \scriptstyle \bullet$} [c] at  11 6
\put {$ \scriptstyle \bullet$} [c] at  11  12
\put {$ \scriptstyle \bullet$} [c] at  12 0
\put {$ \scriptstyle \bullet$} [c] at  12 9
\put {$ \scriptstyle \bullet$} [c] at  16 12
\setlinear \plot  10 0 11 6 10 9 11 12 12 9 11 6 12 0   /
\setlinear \plot  16 12 11  6  /
\put{$1{,}260  $} [c] at 13 -2
\endpicture
\end{minipage}
\begin{minipage}{4cm}
\beginpicture
\setcoordinatesystem units   <1.5mm,2mm>
\setplotarea x from 0 to 16, y from -2 to 15
\put{210)} [l] at 2 12
\put {$ \scriptstyle \bullet$} [c] at  10 0
\put {$ \scriptstyle \bullet$} [c] at  10 9
\put {$ \scriptstyle \bullet$} [c] at  10 6
\put {$ \scriptstyle \bullet$} [c] at  10 12
\put {$ \scriptstyle \bullet$} [c] at  13 12
\put {$ \scriptstyle \bullet$} [c] at  16 0
\put {$ \scriptstyle \bullet$} [c] at  16 12
\setlinear \plot  10 0 10 12  /
\setlinear \plot  16 12 10 6 16 0  /
\setlinear \plot  10 9 13 12  /
\put{$1{,}260$} [c] at 13 -2
\endpicture
\end{minipage}
$$
$$
\begin{minipage}{4cm}
\beginpicture
\setcoordinatesystem units   <1.5mm,2mm>
\setplotarea x from 0 to 16, y from -2 to 15
\put{211)} [l] at 2 12
\put {$ \scriptstyle \bullet$} [c] at  10 0
\put {$ \scriptstyle \bullet$} [c] at  10 3
\put {$ \scriptstyle \bullet$} [c] at  10 6
\put {$ \scriptstyle \bullet$} [c] at  10 12
\put {$ \scriptstyle \bullet$} [c] at  13 0
\put {$ \scriptstyle \bullet$} [c] at  16 0
\put {$ \scriptstyle \bullet$} [c] at  16 12
\setlinear \plot  10 0 10 12  /
\setlinear \plot  16 12 10 6 16 0  /
\setlinear \plot  10 3 13 0  /
\put{$1{,}260$} [c] at 13 -2
\endpicture
\end{minipage}
\begin{minipage}{4cm}
\beginpicture
\setcoordinatesystem units   <1.5mm,2mm>
\setplotarea x from 0 to 16, y from -2 to 15
\put{212)} [l] at 2 12
\put {$ \scriptstyle \bullet$} [c] at 10 6
\put {$ \scriptstyle \bullet$} [c] at 10 12
\put {$ \scriptstyle \bullet$} [c] at 13 0
\put {$ \scriptstyle \bullet$} [c] at 13 6
\put {$ \scriptstyle \bullet$} [c] at 16 6
\put {$ \scriptstyle \bullet$} [c] at 16 9
\put {$ \scriptstyle \bullet$} [c] at 16 12
\setlinear \plot 16 6 10 12 10 6  13 0 16 6 16 12   /
\setlinear \plot 10 12 13 6   16 9  10 6  /
\setlinear \plot 13 0 13 6 /
\put{$840$} [c] at 13 -2
\endpicture
\end{minipage}
\begin{minipage}{4cm}
\beginpicture
\setcoordinatesystem units   <1.5mm,2mm>
\setplotarea x from 0 to 16, y from -2 to 15
\put{213)} [l] at 2 12
\put {$ \scriptstyle \bullet$} [c] at 10 6
\put {$ \scriptstyle \bullet$} [c] at 10 0
\put {$ \scriptstyle \bullet$} [c] at 13 12
\put {$ \scriptstyle \bullet$} [c] at 13 6
\put {$ \scriptstyle \bullet$} [c] at 16 6
\put {$ \scriptstyle \bullet$} [c] at 16 3
\put {$ \scriptstyle \bullet$} [c] at 16 0
\setlinear \plot 16 6 10 0 10 6  13 12 16 6 16 0   /
\setlinear \plot 10 0 13 6   16 3  10 6  /
\setlinear \plot 13 12 13 6 /
\put{$840$} [c] at 13 -2
\endpicture
\end{minipage}
\begin{minipage}{4cm}
\beginpicture
\setcoordinatesystem units   <1.5mm,2mm>
\setplotarea x from 0 to 16, y from -2 to 15
\put{214)} [l] at 2 12
\put {$ \scriptstyle \bullet$} [c] at 10 6
\put {$ \scriptstyle \bullet$} [c] at 10  12
\put {$ \scriptstyle \bullet$} [c] at 13 0
\put {$ \scriptstyle \bullet$} [c] at 13 12
\put {$ \scriptstyle \bullet$} [c] at 14.5 3
\put {$ \scriptstyle \bullet$} [c] at 16 6
\put {$ \scriptstyle \bullet$} [c] at 16 12
\setlinear \plot 13 0 10 6 10 12  16 6 14.5 3 13 0 10 6 16 12 16 6   /
\setlinear \plot 10 6 13 12  16 6   /
\put{$840$} [c] at 13 -2
\endpicture
\end{minipage}
\begin{minipage}{4cm}
\beginpicture
\setcoordinatesystem units   <1.5mm,2mm>
\setplotarea x from 0 to 16, y from -2 to 15
\put{215)} [l] at 2 12
\put {$ \scriptstyle \bullet$} [c] at 10 6
\put {$ \scriptstyle \bullet$} [c] at 10  0
\put {$ \scriptstyle \bullet$} [c] at 13 0
\put {$ \scriptstyle \bullet$} [c] at 13 12
\put {$ \scriptstyle \bullet$} [c] at 14.5 9
\put {$ \scriptstyle \bullet$} [c] at 16 6
\put {$ \scriptstyle \bullet$} [c] at 16 0
\setlinear \plot 13 12 10 6 10 0  16 6 14.5 9 13 12 10 6 16 0 16 6   /
\setlinear \plot 10 6 13 0  16 6   /
\put{$840$} [c] at 13 -2
\endpicture
\end{minipage}
\begin{minipage}{4cm}
\beginpicture
\setcoordinatesystem units   <1.5mm,2mm>
\setplotarea x from 0 to 16, y from -2 to 15
\put{216)} [l] at 2 12
\put {$ \scriptstyle \bullet$} [c] at  10  0
\put {$ \scriptstyle \bullet$} [c] at  10 6
\put {$ \scriptstyle \bullet$} [c] at  10 3
\put {$ \scriptstyle \bullet$} [c] at  10 12
\put {$ \scriptstyle \bullet$} [c] at  13 12
\put {$ \scriptstyle \bullet$} [c] at  16 0
\put {$ \scriptstyle \bullet$} [c] at  16 12
\setlinear \plot  10 0 10 12  /
\setlinear \plot  16 12 10 6 16 0 /
\setlinear \plot  10 6 13 12  /
\put{$840$} [c] at 13 -2
\endpicture
\end{minipage}
$$
$$
\begin{minipage}{4cm}
\beginpicture
\setcoordinatesystem units   <1.5mm,2mm>
\setplotarea x from 0 to 16, y from -2 to 15
\put{217)} [l] at 2 12
\put {$ \scriptstyle \bullet$} [c] at  10  0
\put {$ \scriptstyle \bullet$} [c] at  10 6
\put {$ \scriptstyle \bullet$} [c] at  10 9
\put {$ \scriptstyle \bullet$} [c] at  10 12
\put {$ \scriptstyle \bullet$} [c] at  13 0
\put {$ \scriptstyle \bullet$} [c] at  16 0
\put {$ \scriptstyle \bullet$} [c] at  16 12
\setlinear \plot  10 0 10 12  /
\setlinear \plot  16 12 10 6 16 0 /
\setlinear \plot  10 6 13 0  /
\put{$840$} [c] at 13 -2
\endpicture
\end{minipage}
\begin{minipage}{4cm}
\beginpicture
\setcoordinatesystem units   <1.5mm,2mm>
\setplotarea x from 0 to 16, y from -2 to 15
\put{218)} [l] at 2 12
\put {$ \scriptstyle \bullet$} [c] at  10 0
\put {$ \scriptstyle \bullet$} [c] at  10 6
\put {$ \scriptstyle \bullet$} [c] at  10 12
\put {$ \scriptstyle \bullet$} [c] at  13 12
\put {$ \scriptstyle \bullet$} [c] at  16 0
\put {$ \scriptstyle \bullet$} [c] at  16 6
\put {$ \scriptstyle \bullet$} [c] at  16 12
\setlinear \plot 16 6 10 12 10 0 16 6 16 12 10 6 16 0 16 6 13 12 10 6  /
\put{$210$} [c] at 13 -2
\endpicture
\end{minipage}
\begin{minipage}{4cm}
\beginpicture
\setcoordinatesystem units   <1.5mm,2mm>
\setplotarea x from 0 to 16, y from -2 to 15
\put{219)} [l] at 2 12
\put {$ \scriptstyle \bullet$} [c] at  10 0
\put {$ \scriptstyle \bullet$} [c] at  10 6
\put {$ \scriptstyle \bullet$} [c] at  10 12
\put {$ \scriptstyle \bullet$} [c] at  13 0
\put {$ \scriptstyle \bullet$} [c] at  16 0
\put {$ \scriptstyle \bullet$} [c] at  16 6
\put {$ \scriptstyle \bullet$} [c] at  16 12
\setlinear \plot 16 6 10 12 10 0 16 6 16 12 10 6 16 0 16 6 13 0 10 6  /
\put{$210$} [c] at 13 -2
\endpicture
\end{minipage}
\begin{minipage}{4cm}
\beginpicture
\setcoordinatesystem units   <1.5mm,2mm>
\setplotarea x from 0 to 16, y from -2 to 15
\put{220)} [l] at 2 12
\put {$ \scriptstyle \bullet$} [c] at  10 0
\put {$ \scriptstyle \bullet$} [c] at  10 6
\put {$ \scriptstyle \bullet$} [c] at  10 12
\put {$ \scriptstyle \bullet$} [c] at  13 0
\put {$ \scriptstyle \bullet$} [c] at  13 6
\put {$ \scriptstyle \bullet$} [c] at  13 12
\put {$ \scriptstyle \bullet$} [c] at  16 6
\setlinear \plot  10 12 10 0  16 6 13 0  13 12 16 6 10 12 13 6 10 0  /
\setlinear \plot 13 12 10 6 13 0 /
\put{$210$} [c] at 13 -2
\endpicture
\end{minipage}
\begin{minipage}{4cm}
\beginpicture
\setcoordinatesystem units   <1.5mm,2mm>
\setplotarea x from 0 to 16, y from -2 to 15
\put{221)} [l] at 2 12
\put {$ \scriptstyle \bullet$} [c] at 8 0
\put {$ \scriptstyle \bullet$} [c] at 8 2.3
\put {$ \scriptstyle \bullet$} [c] at 8 4.8
\put {$ \scriptstyle \bullet$} [c] at 8 7.2
\put {$ \scriptstyle \bullet$} [c] at 8 9.6
\put {$ \scriptstyle \bullet$} [c] at 8  12
\put {$ \scriptstyle \bullet$} [c] at 16  0
\setlinear \plot  8 0   8 12     /
\put{$5{,}040$} [c] at 12 -2
\endpicture
\end{minipage}
\begin{minipage}{4cm}
\beginpicture
\setcoordinatesystem units   <1.5mm,2mm>
\setplotarea x from 0 to 16, y from -2 to 15
\put{${\bf  15}$} [l] at 2 15
\put{222)} [l] at 2 12
\put {$ \scriptstyle \bullet$} [c] at 10 4
\put {$ \scriptstyle \bullet$} [c] at 10 8
\put {$ \scriptstyle \bullet$} [c] at 10 12
\put {$ \scriptstyle \bullet$} [c] at 13 0
\put {$ \scriptstyle \bullet$} [c] at 13 8
\put {$ \scriptstyle \bullet$} [c] at 13 12
\put {$ \scriptstyle \bullet$} [c] at 16 4
\setlinear \plot 10 12  10 4 13 0 16 4 13 8 10 4     /
\setlinear \plot 13 12  13 8     /
\put{$5{,}040$} [c] at 13 -2
\endpicture
\end{minipage}
$$
$$
\begin{minipage}{4cm}
\beginpicture
\setcoordinatesystem units   <1.5mm,2mm>
\setplotarea x from 0 to 16, y from -2 to 15
\put{223)} [l] at 2 12
\put {$ \scriptstyle \bullet$} [c] at 10 4
\put {$ \scriptstyle \bullet$} [c] at 10 8
\put {$ \scriptstyle \bullet$} [c] at 10 0
\put {$ \scriptstyle \bullet$} [c] at 13 0
\put {$ \scriptstyle \bullet$} [c] at 13 4
\put {$ \scriptstyle \bullet$} [c] at 13 12
\put {$ \scriptstyle \bullet$} [c] at 16 8
\setlinear \plot 10 0  10 8 13 12 16 8 13 4 10 8     /
\setlinear \plot 13 0  13 4     /
\put{$5{,}040$} [c] at 13 -2
\endpicture
\end{minipage}
\begin{minipage}{4cm}
\beginpicture
\setcoordinatesystem units   <1.5mm,2mm>
\setplotarea x from 0 to 16, y from -2 to 15
\put{224)} [l] at 2 12
\put {$ \scriptstyle \bullet$} [c] at 10 6
\put {$ \scriptstyle \bullet$} [c] at 11.5 9
\put {$ \scriptstyle \bullet$} [c] at 13 0
\put {$ \scriptstyle \bullet$} [c] at 13 12
\put {$ \scriptstyle \bullet$} [c] at 16 3
\put {$ \scriptstyle \bullet$} [c] at 16 6
\put {$ \scriptstyle \bullet$} [c] at 16 12
\setlinear \plot 13 0 10 6 13 12 16 6 16 3 13 0     /
\setlinear \plot 16 12  16 6     /
\put{$5{,}040$} [c] at 13 -2
\endpicture
\end{minipage}
\begin{minipage}{4cm}
\beginpicture
\setcoordinatesystem units   <1.5mm,2mm>
\setplotarea x from 0 to 16, y from -2 to 15
\put{225)} [l] at 2 12
\put {$ \scriptstyle \bullet$} [c] at 10 6
\put {$ \scriptstyle \bullet$} [c] at 11.5 3
\put {$ \scriptstyle \bullet$} [c] at 13 0
\put {$ \scriptstyle \bullet$} [c] at 13 12
\put {$ \scriptstyle \bullet$} [c] at 16 9
\put {$ \scriptstyle \bullet$} [c] at 16 6
\put {$ \scriptstyle \bullet$} [c] at 16 0
\setlinear \plot 13 12 10 6 13 0 16 6 16 9 13 12     /
\setlinear \plot 16 0  16 6     /
\put{$5{,}040$} [c] at 13 -2
\endpicture
\end{minipage}
\begin{minipage}{4cm}
\beginpicture
\setcoordinatesystem units   <1.5mm,2mm>
\setplotarea x from 0 to 16, y from -2 to 15
\put{226)} [l] at 2 12
\put {$ \scriptstyle \bullet$} [c] at 10 6
\put {$ \scriptstyle \bullet$} [c] at 10 8
\put {$ \scriptstyle \bullet$} [c] at 10 10
\put {$ \scriptstyle \bullet$} [c] at 10 12
\put {$ \scriptstyle \bullet$} [c] at 13 12
\put {$ \scriptstyle \bullet$} [c] at 13 0
\put {$ \scriptstyle \bullet$} [c] at 16 6
\setlinear \plot 10 12  10 6 13 0  16 6 13 12 10 6    /
\put{$5{,}040$} [c] at 13 -2
\endpicture
\end{minipage}
\begin{minipage}{4cm}
\beginpicture
\setcoordinatesystem units   <1.5mm,2mm>
\setplotarea x from 0 to 16, y from -2 to 15
\put{227)} [l] at 2 12
\put {$ \scriptstyle \bullet$} [c] at 10 6
\put {$ \scriptstyle \bullet$} [c] at 10 4
\put {$ \scriptstyle \bullet$} [c] at 10 2
\put {$ \scriptstyle \bullet$} [c] at 10 0
\put {$ \scriptstyle \bullet$} [c] at 13 12
\put {$ \scriptstyle \bullet$} [c] at 13 0
\put {$ \scriptstyle \bullet$} [c] at 16 6
\setlinear \plot 10 0  10 6 13 0  16 6 13 12 10 6    /
\put{$5{,}040$} [c] at 13 -2
\endpicture
\end{minipage}
\begin{minipage}{4cm}
\beginpicture
\setcoordinatesystem units   <1.5mm,2mm>
\setplotarea x from 0 to 16, y from -2 to 15
\put{228)} [l] at 2 12
\put {$ \scriptstyle \bullet$} [c] at 10 6
\put {$ \scriptstyle \bullet$} [c] at 10.9 7.7
\put {$ \scriptstyle \bullet$} [c] at 11.9 9.7
\put {$ \scriptstyle \bullet$} [c] at 13 12
\put {$ \scriptstyle \bullet$} [c] at 13 0
\put {$ \scriptstyle \bullet$} [c] at 16 6
\put {$ \scriptstyle \bullet$} [c] at 16 12
\setlinear \plot 16 12 16 6 13 0 10 6 13 12 16 6    /
\put{$5{,}040$} [c] at 13 -2
\endpicture
\end{minipage}
$$
$$
\begin{minipage}{4cm}
\beginpicture
\setcoordinatesystem units   <1.5mm,2mm>
\setplotarea x from 0 to 16, y from -2 to 15
\put{229)} [l] at 2 12
\put {$ \scriptstyle \bullet$} [c] at 10 6
\put {$ \scriptstyle \bullet$} [c] at 11.1 3.7
\put {$ \scriptstyle \bullet$} [c] at 12.1 1.7
\put {$ \scriptstyle \bullet$} [c] at 13 12
\put {$ \scriptstyle \bullet$} [c] at 13 0
\put {$ \scriptstyle \bullet$} [c] at 16 6
\put {$ \scriptstyle \bullet$} [c] at 16 0
\setlinear \plot 16 0 16 6 13 0 10 6 13  12 16 6    /
\put{$5{,}040$} [c] at 13 -2
\endpicture
\end{minipage}
\begin{minipage}{4cm}
\beginpicture
\setcoordinatesystem units   <1.5mm,2mm>
\setplotarea x from 0 to 16, y from -2 to 15
\put{230)} [l] at 2 12
\put {$ \scriptstyle \bullet$} [c] at 10 3
\put {$ \scriptstyle \bullet$} [c] at 10 6
\put {$ \scriptstyle \bullet$} [c] at 10 12
\put {$ \scriptstyle \bullet$} [c] at 11.5 9
\put {$ \scriptstyle \bullet$} [c] at 13 12
\put {$ \scriptstyle \bullet$} [c] at 13 0
\put {$ \scriptstyle \bullet$} [c] at 16 3
\setlinear \plot 10 12  10 3 13 0 16 3 13 12 10 6      /
\put{$5{,}040$} [c] at 13 -2
\endpicture
\end{minipage}
\begin{minipage}{4cm}
\beginpicture
\setcoordinatesystem units   <1.5mm,2mm>
\setplotarea x from 0 to 16, y from -2 to 15
\put{231)} [l] at 2 12
\put {$ \scriptstyle \bullet$} [c] at 10 9
\put {$ \scriptstyle \bullet$} [c] at 10 6
\put {$ \scriptstyle \bullet$} [c] at 10 0
\put {$ \scriptstyle \bullet$} [c] at 11.5 3
\put {$ \scriptstyle \bullet$} [c] at 13 12
\put {$ \scriptstyle \bullet$} [c] at 13 0
\put {$ \scriptstyle \bullet$} [c] at 16 9
\setlinear \plot 10 0  10 9 13 12 16 9 13 0 10 6      /
\put{$5{,}040$} [c] at 13 -2
\endpicture
\end{minipage}
\begin{minipage}{4cm}
\beginpicture
\setcoordinatesystem units   <1.5mm,2mm>
\setplotarea x from 0 to 16, y from -2 to 15
\put{232)} [l] at 2 12
\put {$ \scriptstyle \bullet$} [c] at 10 3
\put {$ \scriptstyle \bullet$} [c] at 13 0
\put {$ \scriptstyle \bullet$} [c] at 13 9
\put {$ \scriptstyle \bullet$} [c] at 13 12
\put {$ \scriptstyle \bullet$} [c] at 14.5 6
\put {$ \scriptstyle \bullet$} [c] at 16 3
\put {$ \scriptstyle \bullet$} [c] at 16 12
\setlinear \plot 16 12  16  3 13 0 10 3 13 9 16 3     /
\setlinear \plot 13 12 13 9     /
\put{$5{,}040$} [c] at 13 -2
\endpicture
\end{minipage}
\begin{minipage}{4cm}
\beginpicture
\setcoordinatesystem units   <1.5mm,2mm>
\setplotarea x from 0 to 16, y from -2 to 15
\put{233)} [l] at 2 12
\put {$ \scriptstyle \bullet$} [c] at 10 9
\put {$ \scriptstyle \bullet$} [c] at 13 0
\put {$ \scriptstyle \bullet$} [c] at 13 3
\put {$ \scriptstyle \bullet$} [c] at 13 12
\put {$ \scriptstyle \bullet$} [c] at 14.5 6
\put {$ \scriptstyle \bullet$} [c] at 16 9
\put {$ \scriptstyle \bullet$} [c] at 16 0
\setlinear \plot 16 0  16  9 13 12 10 9 13 3 16 9     /
\setlinear \plot 13 0 13 3     /
\put{$5{,}040$} [c] at 13 -2
\endpicture
\end{minipage}
\begin{minipage}{4cm}
\beginpicture
\setcoordinatesystem units   <1.5mm,2mm>
\setplotarea x from 0 to 16, y from -2 to 15
\put{234)} [l] at 2 12
\put {$ \scriptstyle \bullet$} [c] at 10 3
\put {$ \scriptstyle \bullet$} [c] at 10 9
\put {$ \scriptstyle \bullet$} [c] at 13 6
\put {$ \scriptstyle \bullet$} [c] at 13 12
\put {$ \scriptstyle \bullet$} [c] at 13 0
\put {$ \scriptstyle \bullet$} [c] at 16 6
\put {$ \scriptstyle \bullet$} [c] at 16 12
\setlinear \plot 13 6 16 12 16 6 13 12 13 6 10 3 10 9 13 12 /
\setlinear \plot 10 3    13 0 16 6      /
\put{$5{,}040$} [c] at 13 -2
\endpicture
\end{minipage}
$$
$$
\begin{minipage}{4cm}
\beginpicture
\setcoordinatesystem units   <1.5mm,2mm>
\setplotarea x from 0 to 16, y from -2 to 15
\put{235)} [l] at 2 12
\put {$ \scriptstyle \bullet$} [c] at 10 3
\put {$ \scriptstyle \bullet$} [c] at 10 9
\put {$ \scriptstyle \bullet$} [c] at 13 6
\put {$ \scriptstyle \bullet$} [c] at 13 12
\put {$ \scriptstyle \bullet$} [c] at 13 0
\put {$ \scriptstyle \bullet$} [c] at 16 6
\put {$ \scriptstyle \bullet$} [c] at 16 0
\setlinear \plot 13 6 16 0 16 6 13 0 13 6 10 9 10 3 13 0 /
\setlinear \plot 10 9    13 12 16 6      /
\put{$5{,}040$} [c] at 13 -2
\endpicture
\end{minipage}
\begin{minipage}{4cm}
\beginpicture
\setcoordinatesystem units   <1.5mm,2mm>
\setplotarea x from 0 to 16, y from -2 to 15
\put{236)} [l] at 2 12
\put {$ \scriptstyle \bullet$} [c] at 10 6
\put {$ \scriptstyle \bullet$} [c] at 12 12
\put {$ \scriptstyle \bullet$} [c] at 14 0
\put {$ \scriptstyle \bullet$} [c] at 16 6
\put {$ \scriptstyle \bullet$} [c] at 16 12
\put {$ \scriptstyle \bullet$} [c] at 11.4 4
\put {$ \scriptstyle \bullet$} [c] at 13.4 10
\setlinear \plot 14 0 10 6 12 12   16 6 14 0 /
\setlinear \plot 16 6 16 12 11.4 4 13.4 10 /
\put{$5{,}040$} [c] at 13 -2
\endpicture
\end{minipage}
\begin{minipage}{4cm}
\beginpicture
\setcoordinatesystem units   <1.5mm,2mm>
\setplotarea x from 0 to 16, y from -2 to 15
\put{237)} [l] at 2 12
\put {$ \scriptstyle \bullet$} [c] at 10 6
\put {$ \scriptstyle \bullet$} [c] at 12 0
\put {$ \scriptstyle \bullet$} [c] at 14 12
\put {$ \scriptstyle \bullet$} [c] at 16 6
\put {$ \scriptstyle \bullet$} [c] at 16 0
\put {$ \scriptstyle \bullet$} [c] at 11.4 8
\put {$ \scriptstyle \bullet$} [c] at 13.4 2
\setlinear \plot 14 12 10 6 12 0   16 6 14 12 /
\setlinear \plot 16 6 16 0 11.4 8 13.4 2 /
\put{$5{,}040$} [c] at 13 -2
\endpicture
\end{minipage}
\begin{minipage}{4cm}
\beginpicture
\setcoordinatesystem units   <1.5mm,2mm>
\setplotarea x from 0 to 16, y from -2 to 15
\put{238)} [l] at 2 12
\put {$ \scriptstyle \bullet$} [c] at  10 0
\put {$ \scriptstyle \bullet$} [c] at  10 3
\put {$ \scriptstyle \bullet$} [c] at  10 6
\put {$ \scriptstyle \bullet$} [c] at  10 9
\put {$ \scriptstyle \bullet$} [c] at  10 12
\put {$ \scriptstyle \bullet$} [c] at  16 0
\put {$ \scriptstyle \bullet$} [c] at  16 12
\setlinear \plot  10 0 10 12   /
\setlinear \plot  16 0 16 12 10  6 /
\put{$5{,}040  $} [c] at 13 -2
\endpicture
\end{minipage}
\begin{minipage}{4cm}
\beginpicture
\setcoordinatesystem units   <1.5mm,2mm>
\setplotarea x from 0 to 16, y from -2 to 15
\put{239)} [l] at 2 12
\put {$ \scriptstyle \bullet$} [c] at  10 0
\put {$ \scriptstyle \bullet$} [c] at  10 3
\put {$ \scriptstyle \bullet$} [c] at  10 6
\put {$ \scriptstyle \bullet$} [c] at  10 9
\put {$ \scriptstyle \bullet$} [c] at  10 12
\put {$ \scriptstyle \bullet$} [c] at  16 0
\put {$ \scriptstyle \bullet$} [c] at  16 12
\setlinear \plot  10 0 10 12   /
\setlinear \plot  16 12 16 0 10  6 /
\put{$5{,}040   $} [c] at 13 -2
\endpicture
\end{minipage}
\begin{minipage}{4cm}
\beginpicture
\setcoordinatesystem units   <1.5mm,2mm>
\setplotarea x from 0 to 16, y from -2 to 15
\put{240)} [l] at 2 12
\put {$ \scriptstyle \bullet$} [c] at  10 0
\put {$ \scriptstyle \bullet$} [c] at  10 3
\put {$ \scriptstyle \bullet$} [c] at  10 6
\put {$ \scriptstyle \bullet$} [c] at  10 9
\put {$ \scriptstyle \bullet$} [c] at  10 12
\put {$ \scriptstyle \bullet$} [c] at  16 0
\put {$ \scriptstyle \bullet$} [c] at  16 12
\setlinear \plot  10 0 10 12 16 0 16 12 10 3  /
\put{$5{,}040  $} [c] at 13 -2
\endpicture
\end{minipage}
$$
$$
\begin{minipage}{4cm}
\beginpicture
\setcoordinatesystem units   <1.5mm,2mm>
\setplotarea x from 0 to 16, y from -2 to 15
\put{241)} [l] at 2 12
\put {$ \scriptstyle \bullet$} [c] at  10 0
\put {$ \scriptstyle \bullet$} [c] at  10 3
\put {$ \scriptstyle \bullet$} [c] at  10 6
\put {$ \scriptstyle \bullet$} [c] at  10 9
\put {$ \scriptstyle \bullet$} [c] at  10 12
\put {$ \scriptstyle \bullet$} [c] at  16 0
\put {$ \scriptstyle \bullet$} [c] at  16 12
\setlinear \plot  10 12 10 0 16 12 16 0 10 9  /
\put{$5{,}040  $} [c] at 13 -2
\endpicture
\end{minipage}
\begin{minipage}{4cm}
\beginpicture
\setcoordinatesystem units   <1.5mm,2mm>
\setplotarea x from 0 to 16, y from -2 to 15
\put{242)} [l] at 2 12
\put {$ \scriptstyle \bullet$} [c] at  10 6
\put {$ \scriptstyle \bullet$} [c] at  10 12
\put {$ \scriptstyle \bullet$} [c] at  11.5 3
\put {$ \scriptstyle \bullet$} [c] at  13 0
\put {$ \scriptstyle \bullet$} [c] at  13 12
\put {$ \scriptstyle \bullet$} [c] at  16 0
\put {$ \scriptstyle \bullet$} [c] at  16 6
\setlinear \plot 16 6 16 0  10 12 10 6 13 0 16 6  13 12 10 6 /
\put{$5{,}040  $} [c] at 13 -2
\endpicture
\end{minipage}
\begin{minipage}{4cm}
\beginpicture
\setcoordinatesystem units   <1.5mm,2mm>
\setplotarea x from 0 to 16, y from -2 to 15
\put{243)} [l] at 2 12
\put {$ \scriptstyle \bullet$} [c] at  10 6
\put {$ \scriptstyle \bullet$} [c] at  10 0
\put {$ \scriptstyle \bullet$} [c] at  11.5 9
\put {$ \scriptstyle \bullet$} [c] at  13 0
\put {$ \scriptstyle \bullet$} [c] at  13 12
\put {$ \scriptstyle \bullet$} [c] at  16 12
\put {$ \scriptstyle \bullet$} [c] at  16 6
\setlinear \plot 16 6 16 12  10 0 10 6 13 0 16 6  13 12 10 6 /
\put{$5{,}040  $} [c] at 13 -2
\endpicture
\end{minipage}
\begin{minipage}{4cm}
\beginpicture
\setcoordinatesystem units   <1.5mm,2mm>
\setplotarea x from 0 to 16, y from -2 to 15
\put{244)} [l] at 2 12
\put {$ \scriptstyle \bullet$} [c] at  10   0
\put {$ \scriptstyle \bullet$} [c] at  10 6
\put {$ \scriptstyle \bullet$} [c] at  10 9
\put {$ \scriptstyle \bullet$} [c] at  10 12
\put {$ \scriptstyle \bullet$} [c] at  16 0
\put {$ \scriptstyle \bullet$} [c] at  16 6
\put {$ \scriptstyle \bullet$} [c] at  16 12
\setlinear \plot 10 12 10 0 16 12 16 0 /
\setlinear \plot 10 9 16 6 /
\put{$5{,}040  $} [c] at 13 -2
\endpicture
\end{minipage}
\begin{minipage}{4cm}
\beginpicture
\setcoordinatesystem units   <1.5mm,2mm>
\setplotarea x from 0 to 16, y from -2 to 15
\put{245)} [l] at 2 12
\put {$ \scriptstyle \bullet$} [c] at  10   0
\put {$ \scriptstyle \bullet$} [c] at  10 6
\put {$ \scriptstyle \bullet$} [c] at  10 3
\put {$ \scriptstyle \bullet$} [c] at  10 12
\put {$ \scriptstyle \bullet$} [c] at  16 0
\put {$ \scriptstyle \bullet$} [c] at  16 6
\put {$ \scriptstyle \bullet$} [c] at  16 12
\setlinear \plot 10 0 10 12 16 0 16 12 /
\setlinear \plot 10 3 16 6 /
\put{$5{,}040  $} [c] at 13 -2
\endpicture
\end{minipage}
\begin{minipage}{4cm}
\beginpicture
\setcoordinatesystem units   <1.5mm,2mm>
\setplotarea x from 0 to 16, y from -2 to 15
\put{246)} [l] at 2 12
\put {$ \scriptstyle \bullet$} [c] at  10 9
\put {$ \scriptstyle \bullet$} [c] at  10 12
\put {$ \scriptstyle \bullet$} [c] at  13 12
\put {$ \scriptstyle \bullet$} [c] at  13 6
\put {$ \scriptstyle \bullet$} [c] at  13 0
\put {$ \scriptstyle \bullet$} [c] at  16 9
\put {$ \scriptstyle \bullet$} [c] at  16 0
\setlinear \plot  13 0 13 6 10 9 10 12    /
\setlinear \plot  10 9 13 12 16 9 13 6     /
\setlinear \plot  16 0 16 9      /
\put{$5{,}040  $} [c] at 13 -2
\endpicture
\end{minipage}
$$
$$
\begin{minipage}{4cm}
\beginpicture
\setcoordinatesystem units   <1.5mm,2mm>
\setplotarea x from 0 to 16, y from -2 to 15
\put{247)} [l] at 2 12
\put {$ \scriptstyle \bullet$} [c] at  10 3
\put {$ \scriptstyle \bullet$} [c] at  10 0
\put {$ \scriptstyle \bullet$} [c] at  13 12
\put {$ \scriptstyle \bullet$} [c] at  13 6
\put {$ \scriptstyle \bullet$} [c] at  13 0
\put {$ \scriptstyle \bullet$} [c] at  16 3
\put {$ \scriptstyle \bullet$} [c] at  16 12
\setlinear \plot  13 12 13 6 10 3 10 0    /
\setlinear \plot  10 3 13 0 16 3 13 6     /
\setlinear \plot  16 3 16 12      /
\put{$5{,}040  $} [c] at 13 -2
\endpicture
\end{minipage}
\begin{minipage}{4cm}
\beginpicture
\setcoordinatesystem units   <1.5mm,2mm>
\setplotarea x from 0 to 16, y from -2 to 15
\put{248)} [l] at 2 12
\put {$ \scriptstyle \bullet$} [c] at  10 6
\put {$ \scriptstyle \bullet$} [c] at  10 9
\put {$ \scriptstyle \bullet$} [c] at  10 12
\put {$ \scriptstyle \bullet$} [c] at  13 0
\put {$ \scriptstyle \bullet$} [c] at  13 12
\put {$ \scriptstyle \bullet$} [c] at  16 6
\put {$ \scriptstyle \bullet$} [c] at  16 0
\setlinear \plot 10 12 10 6 13 0  16 6 13 12 10 6  /
\setlinear \plot 10 9 16 0 16 6 /
\put{$5{,}040   $} [c] at 13 -2
\endpicture
\end{minipage}
\begin{minipage}{4cm}
\beginpicture
\setcoordinatesystem units   <1.5mm,2mm>
\setplotarea x from 0 to 16, y from -2 to 15
\put{249)} [l] at 2 12
\put {$ \scriptstyle \bullet$} [c] at  10 6
\put {$ \scriptstyle \bullet$} [c] at  10 3
\put {$ \scriptstyle \bullet$} [c] at  10 0
\put {$ \scriptstyle \bullet$} [c] at  13 0
\put {$ \scriptstyle \bullet$} [c] at  13 12
\put {$ \scriptstyle \bullet$} [c] at  16 6
\put {$ \scriptstyle \bullet$} [c] at  16 12
\setlinear \plot 10 0 10 6 13 12  16 6 13 0 10 6  /
\setlinear \plot 10 3 16 12 16 6 /
\put{$5{,}040   $} [c] at 13 -2
\endpicture
\end{minipage}
\begin{minipage}{4cm}
\beginpicture
\setcoordinatesystem units   <1.5mm,2mm>
\setplotarea x from 0 to 16, y from -2 to 15
\put{250)} [l] at 2 12
\put {$ \scriptstyle \bullet$} [c] at  10 0
\put {$ \scriptstyle \bullet$} [c] at  10 12
\put {$ \scriptstyle \bullet$} [c] at  12 9
\put {$ \scriptstyle \bullet$} [c] at  14 6
\put {$ \scriptstyle \bullet$} [c] at  15 0
\put {$ \scriptstyle \bullet$} [c] at  15 12
\put {$ \scriptstyle \bullet$} [c] at  16 6
\setlinear \plot  10 12  14 6 15 0 16 6 15 12 14 6 10 0  /
\put{$5{,}040  $} [c] at 13 -2
\endpicture
\end{minipage}
\begin{minipage}{4cm}
\beginpicture
\setcoordinatesystem units   <1.5mm,2mm>
\setplotarea x from 0 to 16, y from -2 to 15
\put{251)} [l] at 2 12
\put {$ \scriptstyle \bullet$} [c] at  10 0
\put {$ \scriptstyle \bullet$} [c] at  10 12
\put {$ \scriptstyle \bullet$} [c] at  12 3
\put {$ \scriptstyle \bullet$} [c] at  14 6
\put {$ \scriptstyle \bullet$} [c] at  15 0
\put {$ \scriptstyle \bullet$} [c] at  15 12
\put {$ \scriptstyle \bullet$} [c] at  16 6
\setlinear \plot  10 12  14 6 15 0 16 6 15 12 14 6 10 0  /
\put{$5{,}040  $} [c] at 13 -2
\endpicture
\end{minipage}
\begin{minipage}{4cm}
\beginpicture
\setcoordinatesystem units   <1.5mm,2mm>
\setplotarea x from 0 to 16, y from -2 to 15
\put{252)} [l] at 2 12
\put {$ \scriptstyle \bullet$} [c] at  10 6
\put {$ \scriptstyle \bullet$} [c] at  10 9
\put {$ \scriptstyle \bullet$} [c] at  10 12
\put {$ \scriptstyle \bullet$} [c] at  13 0
\put {$ \scriptstyle \bullet$} [c] at  16 0
\put {$ \scriptstyle \bullet$} [c] at  16 6
\put {$ \scriptstyle \bullet$} [c] at  16 12
\setlinear \plot 10 12 10 6 13 0  16 6 16 12 10 6  /
\setlinear \plot 10 12 16 6 16 0 /
\put{$5{,}040   $} [c] at 13 -2
\endpicture
\end{minipage}
$$
$$
\begin{minipage}{4cm}
\beginpicture
\setcoordinatesystem units   <1.5mm,2mm>
\setplotarea x from 0 to 16, y from -2 to 15
\put{253)} [l] at 2 12
\put {$ \scriptstyle \bullet$} [c] at  10 6
\put {$ \scriptstyle \bullet$} [c] at  10 3
\put {$ \scriptstyle \bullet$} [c] at  10 0
\put {$ \scriptstyle \bullet$} [c] at  13 12
\put {$ \scriptstyle \bullet$} [c] at  16 0
\put {$ \scriptstyle \bullet$} [c] at  16 6
\put {$ \scriptstyle \bullet$} [c] at  16 12
\setlinear \plot 10 0 10 6 13 12  16 6 16 0 10 6  /
\setlinear \plot 10 0 16 6 16 12 /
\put{$5{,}040   $} [c] at 13 -2
\endpicture
\end{minipage}
\begin{minipage}{4cm}
\beginpicture
\setcoordinatesystem units   <1.5mm,2mm>
\setplotarea x from 0 to 16, y from -2 to 15
\put{254)} [l] at 2 12
\put {$ \scriptstyle \bullet$} [c] at 10 4
\put {$ \scriptstyle \bullet$} [c] at 10 8
\put {$ \scriptstyle \bullet$} [c] at 13 0
\put {$ \scriptstyle \bullet$} [c] at 13 8
\put {$ \scriptstyle \bullet$} [c] at 13 12
\put {$ \scriptstyle \bullet$} [c] at 16 4
\put {$ \scriptstyle \bullet$} [c] at 16 8
\setlinear \plot 16 4 13 0 10 4 10 8 13 12 16 8 16 4 13 8 13 12  /
\setlinear \plot 10 4 13 8   /
\put{$2{,}520$} [c] at 13 -2
\endpicture
\end{minipage}
\begin{minipage}{4cm}
\beginpicture
\setcoordinatesystem units   <1.5mm,2mm>
\setplotarea x from 0 to 16, y from -2 to 15
\put{255)} [l] at 2 12
\put {$ \scriptstyle \bullet$} [c] at 10 4
\put {$ \scriptstyle \bullet$} [c] at 10 8
\put {$ \scriptstyle \bullet$} [c] at 13 0
\put {$ \scriptstyle \bullet$} [c] at 13 4
\put {$ \scriptstyle \bullet$} [c] at 13 12
\put {$ \scriptstyle \bullet$} [c] at 16 4
\put {$ \scriptstyle \bullet$} [c] at 16 8
\setlinear \plot 16 8 13 12 10 8 10 4 13 0 16 4 16 8 13 4 13 0  /
\setlinear \plot 10 8 13 4   /
\put{$2{,}520$} [c] at 13 -2
\endpicture
\end{minipage}
\begin{minipage}{4cm}
\beginpicture
\setcoordinatesystem units   <1.5mm,2mm>
\setplotarea x from 0 to 16, y from -2 to 15
\put{256)} [l] at 2 12
\put {$ \scriptstyle \bullet$} [c] at 10 6
\put {$ \scriptstyle \bullet$} [c] at 10 10
\put {$ \scriptstyle \bullet$} [c] at 13 0
\put {$ \scriptstyle \bullet$} [c] at 13 8
\put {$ \scriptstyle \bullet$} [c] at 13 4
\put {$ \scriptstyle \bullet$} [c] at 13 12
\put {$ \scriptstyle \bullet$} [c] at 16 8
\setlinear \plot 13 0 13 12 16 8 13 4 10 6 10 10 13 12   /
\put{$2{,}520$} [c] at 13 -2
\endpicture
\end{minipage}
\begin{minipage}{4cm}
\beginpicture
\setcoordinatesystem units   <1.5mm,2mm>
\setplotarea x from 0 to 16, y from -2 to 15
\put{257)} [l] at 2 12
\put {$ \scriptstyle \bullet$} [c] at 10 6
\put {$ \scriptstyle \bullet$} [c] at 10 2
\put {$ \scriptstyle \bullet$} [c] at 13 0
\put {$ \scriptstyle \bullet$} [c] at 13 8
\put {$ \scriptstyle \bullet$} [c] at 13 4
\put {$ \scriptstyle \bullet$} [c] at 13 12
\put {$ \scriptstyle \bullet$} [c] at 16 4
\setlinear \plot 13 12 13 0 16 4 13 8 10 6 10 2 13 0   /
\put{$2{,}520$} [c] at 13 -2
\endpicture
\end{minipage}
\begin{minipage}{4cm}
\beginpicture
\setcoordinatesystem units   <1.5mm,2mm>
\setplotarea x from 0 to 16, y from -2 to 15
\put{258)} [l] at 2 12
\put {$ \scriptstyle \bullet$} [c] at 10 6
\put {$ \scriptstyle \bullet$} [c] at 10 12
\put {$ \scriptstyle \bullet$} [c] at 13 0
\put {$ \scriptstyle \bullet$} [c] at 13 3
\put {$ \scriptstyle \bullet$} [c] at 13 9
\put {$ \scriptstyle \bullet$} [c] at 16 6
\put {$ \scriptstyle \bullet$} [c] at 16 12
\setlinear \plot 13 0 13 3 10 6 10 12 13 9 13 3 16 6 16 12 13 9  /
\put{$2{,}520$} [c] at 13 -2
\endpicture
\end{minipage}
$$
$$
\begin{minipage}{4cm}
\beginpicture
\setcoordinatesystem units   <1.5mm,2mm>
\setplotarea x from 0 to 16, y from -2 to 15
\put{259)} [l] at 2 12
\put {$ \scriptstyle \bullet$} [c] at 10 6
\put {$ \scriptstyle \bullet$} [c] at 10 0
\put {$ \scriptstyle \bullet$} [c] at 13 12
\put {$ \scriptstyle \bullet$} [c] at 13 9
\put {$ \scriptstyle \bullet$} [c] at 13 3
\put {$ \scriptstyle \bullet$} [c] at 16 6
\put {$ \scriptstyle \bullet$} [c] at 16 0
\setlinear \plot 13 12 13 9 10 6 10 0 13 3 13 9 16 6 16 0 13 3  /
\put{$2{,}520$} [c] at 13 -2
\endpicture
\end{minipage}
\begin{minipage}{4cm}
\beginpicture
\setcoordinatesystem units   <1.5mm,2mm>
\setplotarea x from 0 to 16, y from -2 to 15
\put{260)} [l] at 2 12
\put {$ \scriptstyle \bullet$} [c] at 12 0
\put {$ \scriptstyle \bullet$} [c] at 12 3
\put {$ \scriptstyle \bullet$} [c] at 12 6
\put {$ \scriptstyle \bullet$} [c] at 10 9
\put {$ \scriptstyle \bullet$} [c] at 12 12
\put {$ \scriptstyle \bullet$} [c] at 14 9
\put {$ \scriptstyle \bullet$} [c] at 16 12
\setlinear \plot 16 12 12 0 12 6 10 9 12 12 14 9 12 6   /
\put{$2{,}520$} [c] at 13 -2
\endpicture
\end{minipage}
\begin{minipage}{4cm}
\beginpicture
\setcoordinatesystem units   <1.5mm,2mm>
\setplotarea x from 0 to 16, y from -2 to 15
\put{261)} [l] at 2 12
\put {$ \scriptstyle \bullet$} [c] at 12 12
\put {$ \scriptstyle \bullet$} [c] at 12 9
\put {$ \scriptstyle \bullet$} [c] at 12 6
\put {$ \scriptstyle \bullet$} [c] at 10 3
\put {$ \scriptstyle \bullet$} [c] at 12 0
\put {$ \scriptstyle \bullet$} [c] at 14 3
\put {$ \scriptstyle \bullet$} [c] at 16 0
\setlinear \plot 16 0 12 12 12 6 10 3 12 0 14 3 12 6   /
\put{$2{,}520$} [c] at 13 -2
\endpicture
\end{minipage}
\begin{minipage}{4cm}
\beginpicture
\setcoordinatesystem units   <1.5mm,2mm>
\setplotarea x from 0 to 16, y from -2 to 15
\put{262)} [l] at 2 12
\put {$ \scriptstyle \bullet$} [c] at 10 6
\put {$ \scriptstyle \bullet$} [c] at 11 0
\put {$ \scriptstyle \bullet$} [c] at 11 3
\put {$ \scriptstyle \bullet$} [c] at 11 9
\put {$ \scriptstyle \bullet$} [c] at 11 12
\put {$ \scriptstyle \bullet$} [c] at 12 6
\put {$ \scriptstyle \bullet$} [c] at 16 12
\setlinear \plot 16 12 11 0 11 3 10 6 11 9 11 12   /
\setlinear \plot  11 9 12 6 11 3   /
\put{$2{,}520$} [c] at 13 -2
\endpicture
\end{minipage}
\begin{minipage}{4cm}
\beginpicture
\setcoordinatesystem units   <1.5mm,2mm>
\setplotarea x from 0 to 16, y from -2 to 15
\put{263)} [l] at 2 12
\put {$ \scriptstyle \bullet$} [c] at 10 6
\put {$ \scriptstyle \bullet$} [c] at 11 0
\put {$ \scriptstyle \bullet$} [c] at 11 3
\put {$ \scriptstyle \bullet$} [c] at 11 9
\put {$ \scriptstyle \bullet$} [c] at 11 12
\put {$ \scriptstyle \bullet$} [c] at 12 6
\put {$ \scriptstyle \bullet$} [c] at 16 0
\setlinear \plot 16 0 11 12 11 9 10 6 11 3 11 0   /
\setlinear \plot  11 9 12 6 11 3   /
\put{$2{,}520$} [c] at 13 -2
\endpicture
\end{minipage}
\begin{minipage}{4cm}
\beginpicture
\setcoordinatesystem units   <1.5mm,2mm>
\setplotarea x from 0 to 16, y from -2 to 15
\put{264)} [l] at 2 12
\put {$ \scriptstyle \bullet$} [c] at 10 4
\put {$ \scriptstyle \bullet$} [c] at 11 0
\put {$ \scriptstyle \bullet$} [c] at 11 8
\put {$ \scriptstyle \bullet$} [c] at 11 10
\put {$ \scriptstyle \bullet$} [c] at 11 12
\put {$ \scriptstyle \bullet$} [c] at 12 4
\put {$ \scriptstyle \bullet$} [c] at 16 12
\setlinear \plot 16 12 11 0 10 4 11 8 11 12     /
\setlinear \plot 11 0 12 4 11 8    /
\put{$2{,}520$} [c] at 13 -2
\endpicture
\end{minipage}
$$
$$
\begin{minipage}{4cm}
\beginpicture
\setcoordinatesystem units   <1.5mm,2mm>
\setplotarea x from 0 to 16, y from -2 to 15
\put{265)} [l] at 2 12
\put {$ \scriptstyle \bullet$} [c] at 10 8
\put {$ \scriptstyle \bullet$} [c] at 11 12
\put {$ \scriptstyle \bullet$} [c] at 11 2
\put {$ \scriptstyle \bullet$} [c] at 11 4
\put {$ \scriptstyle \bullet$} [c] at 11 0
\put {$ \scriptstyle \bullet$} [c] at 12 8
\put {$ \scriptstyle \bullet$} [c] at 16 0
\setlinear \plot 16 0 11 12 10 8 11 4 11 0     /
\setlinear \plot 11 12 12 8 11 4    /
\put{$2{,}520$} [c] at 13 -2
\endpicture
\end{minipage}
\begin{minipage}{4cm}
\beginpicture
\setcoordinatesystem units   <1.5mm,2mm>
\setplotarea x from 0 to 16, y from -2 to 15
\put{266)} [l] at 2 12
\put {$ \scriptstyle \bullet$} [c] at 10 8
\put {$ \scriptstyle \bullet$} [c] at 10 12
\put {$ \scriptstyle \bullet$} [c] at 12 4
\put {$ \scriptstyle \bullet$} [c] at 12 8
\put {$ \scriptstyle \bullet$} [c] at 12 12
\put {$ \scriptstyle \bullet$} [c] at 14 0
\put {$ \scriptstyle \bullet$} [c] at 16 4
\setlinear \plot 14 0 12 4 12 12 16 4 14 0    /
\setlinear \plot 12 4 10 8 10 12 12 8    /
\setlinear \plot 12 12 10 8    /
\put{$2{,}520$} [c] at 13 -2
\endpicture
\end{minipage}
\begin{minipage}{4cm}
\beginpicture
\setcoordinatesystem units   <1.5mm,2mm>
\setplotarea x from 0 to 16, y from -2 to 15
\put{267)} [l] at 2 12
\put {$ \scriptstyle \bullet$} [c] at 10 4
\put {$ \scriptstyle \bullet$} [c] at 10 0
\put {$ \scriptstyle \bullet$} [c] at 12 4
\put {$ \scriptstyle \bullet$} [c] at 12 8
\put {$ \scriptstyle \bullet$} [c] at 12 0
\put {$ \scriptstyle \bullet$} [c] at 14 12
\put {$ \scriptstyle \bullet$} [c] at 16 8
\setlinear \plot 14 12 12 8 12 0 16 8 14 12    /
\setlinear \plot 12 8 10 4 10 0 12 4    /
\setlinear \plot 12 0 10 4    /
\put{$2{,}520$} [c] at 13 -2
\endpicture
\end{minipage}
\begin{minipage}{4cm}
\beginpicture
\setcoordinatesystem units   <1.5mm,2mm>
\setplotarea x from 0 to 16, y from -2 to 15
\put{268)} [l] at 2 12
\put {$ \scriptstyle \bullet$} [c] at 10 4
\put {$ \scriptstyle \bullet$} [c] at 10 8
\put {$ \scriptstyle \bullet$} [c] at 12 0
\put {$ \scriptstyle \bullet$} [c] at 12 12
\put {$ \scriptstyle \bullet$} [c] at 14 4
\put {$ \scriptstyle \bullet$} [c] at 14 8
\put {$ \scriptstyle \bullet$} [c] at 16 12
\setlinear \plot 16 12 14 4 12 0  10 4 10 8 12 12 14 8 10 4 10 8 14 4 14 8    /
\put{$2{,}520$} [c] at 13 -2
\endpicture
\end{minipage}
\begin{minipage}{4cm}
\beginpicture
\setcoordinatesystem units   <1.5mm,2mm>
\setplotarea x from 0 to 16, y from -2 to 15
\put{269)} [l] at 2 12
\put {$ \scriptstyle \bullet$} [c] at 10 4
\put {$ \scriptstyle \bullet$} [c] at 10 8
\put {$ \scriptstyle \bullet$} [c] at 12 0
\put {$ \scriptstyle \bullet$} [c] at 12 12
\put {$ \scriptstyle \bullet$} [c] at 14 4
\put {$ \scriptstyle \bullet$} [c] at 14 8
\put {$ \scriptstyle \bullet$} [c] at 16 0
\setlinear \plot 16 0 14 8 12 12  10 8 10 4 12 0 14 4 10 8 10 4 14 8  14 4  /
\put{$2{,}520$} [c] at 13 -2
\endpicture
\end{minipage}
\begin{minipage}{4cm}
\beginpicture
\setcoordinatesystem units   <1.5mm,2mm>
\setplotarea x from 0 to 16, y from -2 to 15
\put{270)} [l] at 2 12
\put {$ \scriptstyle \bullet$} [c] at 10 4
\put {$ \scriptstyle \bullet$} [c] at 10 12
\put {$ \scriptstyle \bullet$} [c] at 13 0
\put {$ \scriptstyle \bullet$} [c] at 13 4
\put {$ \scriptstyle \bullet$} [c] at 16 4
\put {$ \scriptstyle \bullet$} [c] at 16 8
\put {$ \scriptstyle \bullet$} [c] at 16 12
\setlinear \plot 13 0 10 4 10 12 13 4 13 0 16 4 16 12  10 4  /
\setlinear \plot 10 12 16 4      /
\setlinear \plot 13 4 16 8      /
\put{$2{,}520$} [c] at 13 -2
\endpicture
\end{minipage}
$$
$$
\begin{minipage}{4cm}
\beginpicture
\setcoordinatesystem units   <1.5mm,2mm>
\setplotarea x from 0 to 16, y from -2 to 15
\put{271)} [l] at 2 12
\put {$ \scriptstyle \bullet$} [c] at 10 8
\put {$ \scriptstyle \bullet$} [c] at 10 0
\put {$ \scriptstyle \bullet$} [c] at 13 12
\put {$ \scriptstyle \bullet$} [c] at 13 8
\put {$ \scriptstyle \bullet$} [c] at 16 4
\put {$ \scriptstyle \bullet$} [c] at 16 8
\put {$ \scriptstyle \bullet$} [c] at 16 0
\setlinear \plot 13 12 10 8 10 0 13 8 13 12 16 8 16 0  10 8  /
\setlinear \plot 10 0 16 8      /
\setlinear \plot 13 8 16 4      /
\put{$2{,}520$} [c] at 13 -2
\endpicture
\end{minipage}
\begin{minipage}{4cm}
\beginpicture
\setcoordinatesystem units   <1.5mm,2mm>
\setplotarea x from 0 to 16, y from -2 to 15
\put{272)} [l] at 2 12
\put {$ \scriptstyle \bullet$} [c] at 10 4
\put {$ \scriptstyle \bullet$} [c] at 10 12
\put {$ \scriptstyle \bullet$} [c] at 13 0
\put {$ \scriptstyle \bullet$} [c] at 13 4
\put {$ \scriptstyle \bullet$} [c] at 16 4
\put {$ \scriptstyle \bullet$} [c] at 16 8
\put {$ \scriptstyle \bullet$} [c] at 16 12
\setlinear \plot 13 0 13 4 10 12 10 4 16 8 13 4  16 8 10 4 /
\setlinear \plot 10 4 13 0 16 4 16 12     /
\put{$2{,}520$} [c] at 13 -2
\endpicture
\end{minipage}
\begin{minipage}{4cm}
\beginpicture
\setcoordinatesystem units   <1.5mm,2mm>
\setplotarea x from 0 to 16, y from -2 to 15
\put{273)} [l] at 2 12
\put {$ \scriptstyle \bullet$} [c] at 10 8
\put {$ \scriptstyle \bullet$} [c] at 10 0
\put {$ \scriptstyle \bullet$} [c] at 13 12
\put {$ \scriptstyle \bullet$} [c] at 13 8
\put {$ \scriptstyle \bullet$} [c] at 16 0
\put {$ \scriptstyle \bullet$} [c] at 16 4
\put {$ \scriptstyle \bullet$} [c] at 16 8
\setlinear \plot 13 12 13 8 10 0 10 8 16 4 13 8  16 4 10 8 /
\setlinear \plot 10 8 13 12 16 8 16 0     /
\put{$2{,}520$} [c] at 13 -2
\endpicture
\end{minipage}
\begin{minipage}{4cm}
\beginpicture
\setcoordinatesystem units   <1.5mm,2mm>
\setplotarea x from 0 to 16, y from -2 to 15
\put{274)} [l] at 2 12
\put {$ \scriptstyle \bullet$} [c] at 10 4
\put {$ \scriptstyle \bullet$} [c] at 10 8
\put {$ \scriptstyle \bullet$} [c] at 10 12
\put {$ \scriptstyle \bullet$} [c] at 13 0
\put {$ \scriptstyle \bullet$} [c] at 16 4
\put {$ \scriptstyle \bullet$} [c] at 16 8
\put {$ \scriptstyle \bullet$} [c] at 16 12
\setlinear \plot 16 4 13 0 10 4 10 12 16 4 16 12 10 4    /
\put{$2{,}520$} [c] at 13 -2
\endpicture
\end{minipage}
\begin{minipage}{4cm}
\beginpicture
\setcoordinatesystem units   <1.5mm,2mm>
\setplotarea x from 0 to 16, y from -2 to 15
\put{275)} [l] at 2 12
\put {$ \scriptstyle \bullet$} [c] at 10 8
\put {$ \scriptstyle \bullet$} [c] at 10 4
\put {$ \scriptstyle \bullet$} [c] at 10 0
\put {$ \scriptstyle \bullet$} [c] at 13 12
\put {$ \scriptstyle \bullet$} [c] at 16 8
\put {$ \scriptstyle \bullet$} [c] at 16 4
\put {$ \scriptstyle \bullet$} [c] at 16 0
\setlinear \plot 16 8 13 12 10 8 10 0 16 8 16 0 10 8    /
\put{$2{,}520$} [c] at 13 -2
\endpicture
\end{minipage}
\begin{minipage}{4cm}
\beginpicture
\setcoordinatesystem units   <1.5mm,2mm>
\setplotarea x from 0 to 16, y from -2 to 15
\put{276)} [l] at 2 12
\put {$ \scriptstyle \bullet$} [c] at 10 6
\put {$ \scriptstyle \bullet$} [c] at 13 0
\put {$ \scriptstyle \bullet$} [c] at 13 6
\put {$ \scriptstyle \bullet$} [c] at 13 12
\put {$ \scriptstyle \bullet$} [c] at 14.5 9
\put {$ \scriptstyle \bullet$} [c] at 16 6
\put {$ \scriptstyle \bullet$} [c] at 16 12
\setlinear \plot 13 0 13 6 14.5 9 16 12 /
\setlinear \plot 13 12 10 6 13 0 16 6 13 12   /
\put{$2{,}520$} [c] at 13 -2
\endpicture
\end{minipage}
$$
$$
\begin{minipage}{4cm}
\beginpicture
\setcoordinatesystem units   <1.5mm,2mm>
\setplotarea x from 0 to 16, y from -2 to 15
\put{277)} [l] at 2 12
\put {$ \scriptstyle \bullet$} [c] at 10 6
\put {$ \scriptstyle \bullet$} [c] at 13 0
\put {$ \scriptstyle \bullet$} [c] at 13 6
\put {$ \scriptstyle \bullet$} [c] at 13 12
\put {$ \scriptstyle \bullet$} [c] at 14.5 3
\put {$ \scriptstyle \bullet$} [c] at 16 6
\put {$ \scriptstyle \bullet$} [c] at 16 0
\setlinear \plot 13 12 13 6 14.5 3 16 0 /
\setlinear \plot  13 0 10 6 13 12 16 6  13 0  /
\put{$2{,}520$} [c] at 13 -2
\endpicture
\end{minipage}
\begin{minipage}{4cm}
\beginpicture
\setcoordinatesystem units   <1.5mm,2mm>
\setplotarea x from 0 to 16, y from -2 to 15
\put{278)} [l] at 2 12
\put {$ \scriptstyle \bullet$} [c] at 10 12
\put {$ \scriptstyle \bullet$} [c] at 12  0
\put {$ \scriptstyle \bullet$} [c] at 12  3
\put {$ \scriptstyle \bullet$} [c] at 12 6
\put {$ \scriptstyle \bullet$} [c] at 12 9
\put {$ \scriptstyle \bullet$} [c] at 14 12
\put {$ \scriptstyle \bullet$} [c] at 16 12
\setlinear \plot 16 12 12 0 12 9 14 12    /
\setlinear \plot 10 12 12 9    /
\put{$2{,}520$} [c] at 13 -2
\endpicture
\end{minipage}
\begin{minipage}{4cm}
\beginpicture
\setcoordinatesystem units   <1.5mm,2mm>
\setplotarea x from 0 to 16, y from -2 to 15
\put{279)} [l] at 2 12
\put {$ \scriptstyle \bullet$} [c] at 10 0
\put {$ \scriptstyle \bullet$} [c] at 12  12
\put {$ \scriptstyle \bullet$} [c] at 12  9
\put {$ \scriptstyle \bullet$} [c] at 12 6
\put {$ \scriptstyle \bullet$} [c] at 12 3
\put {$ \scriptstyle \bullet$} [c] at 14 0
\put {$ \scriptstyle \bullet$} [c] at 16 0
\setlinear \plot 16 0 12 12 12 3 14 0    /
\setlinear \plot 10 0 12 3    /
\put{$2{,}520$} [c] at 13 -2
\endpicture
\end{minipage}
\begin{minipage}{4cm}
\beginpicture
\setcoordinatesystem units   <1.5mm,2mm>
\setplotarea x from 0 to 16, y from -2 to 15
\put{280)} [l] at 2 12
\put {$ \scriptstyle \bullet$} [c] at 10 12
\put {$ \scriptstyle \bullet$} [c] at 13 0
\put {$ \scriptstyle \bullet$} [c] at 13 3
\put {$ \scriptstyle \bullet$} [c] at 13 6
\put {$ \scriptstyle \bullet$} [c] at 13 9
\put {$ \scriptstyle \bullet$} [c] at 13 12
\put {$ \scriptstyle \bullet$} [c] at 16 12
\setlinear \plot 13 0 13 12   /
\setlinear \plot 10 12 13 6 16  12  /
\put{$2{,}520$} [c] at 13 -2
\endpicture
\end{minipage}
\begin{minipage}{4cm}
\beginpicture
\setcoordinatesystem units   <1.5mm,2mm>
\setplotarea x from 0 to 16, y from -2 to 15
\put{281)} [l] at 2 12
\put {$ \scriptstyle \bullet$} [c] at 10 0
\put {$ \scriptstyle \bullet$} [c] at 13 0
\put {$ \scriptstyle \bullet$} [c] at 13 3
\put {$ \scriptstyle \bullet$} [c] at 13 6
\put {$ \scriptstyle \bullet$} [c] at 13 9
\put {$ \scriptstyle \bullet$} [c] at 13 12
\put {$ \scriptstyle \bullet$} [c] at 16 0
\setlinear \plot 13 0 13 12   /
\setlinear \plot 10 0 13 6 16  0  /
\put{$2{,}520$} [c] at 13 -2
\endpicture
\end{minipage}
\begin{minipage}{4cm}
\beginpicture
\setcoordinatesystem units   <1.5mm,2mm>
\setplotarea x from 0 to 16, y from -2 to 15
\put{282)} [l] at 2 12
\put {$ \scriptstyle \bullet$} [c] at 10 8
\put {$ \scriptstyle \bullet$} [c] at 10 12
\put {$ \scriptstyle \bullet$} [c] at 13  0
\put {$ \scriptstyle \bullet$} [c] at 13 4
\put {$ \scriptstyle \bullet$} [c] at 13 12
\put {$ \scriptstyle \bullet$} [c] at 16 12
\put {$ \scriptstyle \bullet$} [c] at 16 8
\setlinear \plot 10 12 10 8 13 4 13 0    /
\setlinear \plot 10 8 13 12 16 8  16 12    /
\setlinear \plot 13 4 16 8 /
\put{$2{,}520$} [c] at 13 -2
\endpicture
\end{minipage}
$$
$$
\begin{minipage}{4cm}
\beginpicture
\setcoordinatesystem units   <1.5mm,2mm>
\setplotarea x from 0 to 16, y from -2 to 15
\put{283)} [l] at 2 12
\put {$ \scriptstyle \bullet$} [c] at 10 4
\put {$ \scriptstyle \bullet$} [c] at 10 0
\put {$ \scriptstyle \bullet$} [c] at 13  0
\put {$ \scriptstyle \bullet$} [c] at 13 8
\put {$ \scriptstyle \bullet$} [c] at 13 12
\put {$ \scriptstyle \bullet$} [c] at 16 0
\put {$ \scriptstyle \bullet$} [c] at 16 4
\setlinear \plot 10 0 10 4 13 8 13 12    /
\setlinear \plot 10 4 13 0 16 4  16 0    /
\setlinear \plot 16 4 13 8 /
\put{$2{,}520$} [c] at 13 -2
\endpicture
\end{minipage}
\begin{minipage}{4cm}
\beginpicture
\setcoordinatesystem units   <1.5mm,2mm>
\setplotarea x from 0 to 16, y from -2 to 15
\put{284)} [l] at 2 12
\put {$ \scriptstyle \bullet$} [c] at 10 12
\put {$ \scriptstyle \bullet$} [c] at 12 4
\put {$ \scriptstyle \bullet$} [c] at 12 12
\put {$ \scriptstyle \bullet$} [c] at 14 0
\put {$ \scriptstyle \bullet$} [c] at 14 8
\put {$ \scriptstyle \bullet$} [c] at 16 4
\put {$ \scriptstyle \bullet$} [c] at 16 12
\setlinear \plot 10 12  12 4 14 0 16 4 14 8 16 12    /
\setlinear \plot 12 12 14 8  12 4   /
\put{$2{,}520$} [c] at 13 -2
\endpicture
\end{minipage}
\begin{minipage}{4cm}
\beginpicture
\setcoordinatesystem units   <1.5mm,2mm>
\setplotarea x from 0 to 16, y from -2 to 15
\put{285)} [l] at 2 12
\put {$ \scriptstyle \bullet$} [c] at 10 0
\put {$ \scriptstyle \bullet$} [c] at 12 8
\put {$ \scriptstyle \bullet$} [c] at 12 0
\put {$ \scriptstyle \bullet$} [c] at 14 12
\put {$ \scriptstyle \bullet$} [c] at 14 4
\put {$ \scriptstyle \bullet$} [c] at 16 8
\put {$ \scriptstyle \bullet$} [c] at 16 0
\setlinear \plot 10 0  12 8 14 12 16 8 14 4 16 0    /
\setlinear \plot 12 0 14 4  12 8   /
\put{$2{,}520$} [c] at 13 -2
\endpicture
\end{minipage}
\begin{minipage}{4cm}
\beginpicture
\setcoordinatesystem units   <1.5mm,2mm>
\setplotarea x from 0 to 16, y from -2 to 15
\put{286)} [l] at 2 12
\put {$ \scriptstyle \bullet$} [c] at 10 4
\put {$ \scriptstyle \bullet$} [c] at 10 12
\put {$ \scriptstyle \bullet$} [c] at 13 12
\put {$ \scriptstyle \bullet$} [c] at 13 0
\put {$ \scriptstyle \bullet$} [c] at 16 4
\put {$ \scriptstyle \bullet$} [c] at 16 8
\put {$ \scriptstyle \bullet$} [c] at 16 12
\setlinear \plot 16 4 13  0 10 4 10 12 16 4 16 12 10 4 13 12 16 8  /
\put{$2{,}520$} [c] at 13 -2
\endpicture
\end{minipage}
\begin{minipage}{4cm}
\beginpicture
\setcoordinatesystem units   <1.5mm,2mm>
\setplotarea x from 0 to 16, y from -2 to 15
\put{287)} [l] at 2 12
\put {$ \scriptstyle \bullet$} [c] at 10 8
\put {$ \scriptstyle \bullet$} [c] at 10 0
\put {$ \scriptstyle \bullet$} [c] at 13 12
\put {$ \scriptstyle \bullet$} [c] at 13 0
\put {$ \scriptstyle \bullet$} [c] at 16 4
\put {$ \scriptstyle \bullet$} [c] at 16 8
\put {$ \scriptstyle \bullet$} [c] at 16 0
\setlinear \plot 16 8 13  12 10 8 10 0 16 8 16 0 10 8 13 0 16 4  /
\put{$2{,}520$} [c] at 13 -2
\endpicture
\end{minipage}
\begin{minipage}{4cm}
\beginpicture
\setcoordinatesystem units   <1.5mm,2mm>
\setplotarea x from 0 to 16, y from -2 to 15
\put{288)} [l] at 2 12
\put {$ \scriptstyle \bullet$} [c] at 10 4
\put {$ \scriptstyle \bullet$} [c] at 10 12
\put {$ \scriptstyle \bullet$} [c] at 13 0
\put {$ \scriptstyle \bullet$} [c] at 13 12
\put {$ \scriptstyle \bullet$} [c] at 16 4
\put {$ \scriptstyle \bullet$} [c] at 16 8
\put {$ \scriptstyle \bullet$} [c] at 16 12
\setlinear \plot 13 0 10 4 10 12  16 8 16 4 13 0 /
\setlinear \plot 10 4 13 12 16 8 16 12     /
\put{$2{,}520$} [c] at 13 -2
\endpicture
\end{minipage}
$$
$$
\begin{minipage}{4cm}
\beginpicture
\setcoordinatesystem units   <1.5mm,2mm>
\setplotarea x from 0 to 16, y from -2 to 15
\put{289)} [l] at 2 12
\put {$ \scriptstyle \bullet$} [c] at 10 8
\put {$ \scriptstyle \bullet$} [c] at 10 0
\put {$ \scriptstyle \bullet$} [c] at 13 0
\put {$ \scriptstyle \bullet$} [c] at 13 12
\put {$ \scriptstyle \bullet$} [c] at 16 4
\put {$ \scriptstyle \bullet$} [c] at 16 8
\put {$ \scriptstyle \bullet$} [c] at 16 0
\setlinear \plot 13 12 10 8 10 0  16 4 16 8 13 12 /
\setlinear \plot 10 8 13 0 16 4 16 0     /
\put{$2{,}520$} [c] at 13 -2
\endpicture
\end{minipage}
\begin{minipage}{4cm}
\beginpicture
\setcoordinatesystem units   <1.5mm,2mm>
\setplotarea x from 0 to 16, y from -2 to 15
\put{290)} [l] at 2 12
\put {$ \scriptstyle \bullet$} [c] at  10 0
\put {$ \scriptstyle \bullet$} [c] at  10 6
\put {$ \scriptstyle \bullet$} [c] at  10 9
\put {$ \scriptstyle \bullet$} [c] at  10 12
\put {$ \scriptstyle \bullet$} [c] at  16 0
\put {$ \scriptstyle \bullet$} [c] at  16 6
\put {$ \scriptstyle \bullet$} [c] at  16 12
\setlinear \plot  10  12 10 0 16 6 16 12     /
\setlinear \plot  10 6 16 0 16 6      /
\put{$2{,}520 $} [c] at 13 -2
\endpicture
\end{minipage}
\begin{minipage}{4cm}
\beginpicture
\setcoordinatesystem units   <1.5mm,2mm>
\setplotarea x from 0 to 16, y from -2 to 15
\put{291)} [l] at 2 12
\put {$ \scriptstyle \bullet$} [c] at  10 0
\put {$ \scriptstyle \bullet$} [c] at  10 6
\put {$ \scriptstyle \bullet$} [c] at  10 3
\put {$ \scriptstyle \bullet$} [c] at  10 12
\put {$ \scriptstyle \bullet$} [c] at  16 0
\put {$ \scriptstyle \bullet$} [c] at  16 6
\put {$ \scriptstyle \bullet$} [c] at  16 12
\setlinear \plot  10  0 10 12 16 6 16 0     /
\setlinear \plot  10 6 16 12 16 6      /
\put{$2{,}520  $} [c] at 13 -2
\endpicture
\end{minipage}
\begin{minipage}{4cm}
\beginpicture
\setcoordinatesystem units   <1.5mm,2mm>
\setplotarea x from 0 to 16, y from -2 to 15
\put{292)} [l] at 2 12
\put {$ \scriptstyle \bullet$} [c] at  10 3
\put {$ \scriptstyle \bullet$} [c] at  10 9
\put {$ \scriptstyle \bullet$} [c] at  10 12
\put {$ \scriptstyle \bullet$} [c] at  11 0
\put {$ \scriptstyle \bullet$} [c] at  12 3
\put {$ \scriptstyle \bullet$} [c] at  12 12
\put {$ \scriptstyle \bullet$} [c] at  16 0
\setlinear \plot 10 12 10 3 11 0 12 3 12 12 16 0 10 9 12 3 /
\put{$2{,}520   $} [c] at 13 -2
\endpicture
\end{minipage}
\begin{minipage}{4cm}
\beginpicture
\setcoordinatesystem units   <1.5mm,2mm>
\setplotarea x from 0 to 16, y from -2 to 15
\put{293)} [l] at 2 12
\put {$ \scriptstyle \bullet$} [c] at  10 3
\put {$ \scriptstyle \bullet$} [c] at  10 9
\put {$ \scriptstyle \bullet$} [c] at  10 0
\put {$ \scriptstyle \bullet$} [c] at  11 12
\put {$ \scriptstyle \bullet$} [c] at  12 9
\put {$ \scriptstyle \bullet$} [c] at  12 0
\put {$ \scriptstyle \bullet$} [c] at  16 12
\setlinear \plot 10 0 10 9 11 12 12 9 12 0 16 12 10 3 12 9 /
\put{$2{,}520   $} [c] at 13 -2
\endpicture
\end{minipage}
\begin{minipage}{4cm}
\beginpicture
\setcoordinatesystem units   <1.5mm,2mm>
\setplotarea x from 0 to 16, y from -2 to 15
\put{294)} [l] at 2 12
\put {$ \scriptstyle \bullet$} [c] at  10 0
\put {$ \scriptstyle \bullet$} [c] at  10 6
\put {$ \scriptstyle \bullet$} [c] at  10 12
\put {$ \scriptstyle \bullet$} [c] at  11.5 3
\put {$ \scriptstyle \bullet$} [c] at  13 0
\put {$ \scriptstyle \bullet$} [c] at  16 6
\put {$ \scriptstyle \bullet$} [c] at  16 12
\setlinear \plot  10  0 10 12 16 6  16 12 10 6 13 0 16 6    /
\put{$2{,}520   $} [c] at 13 -2
\endpicture
\end{minipage}
$$
$$
\begin{minipage}{4cm}
\beginpicture
\setcoordinatesystem units   <1.5mm,2mm>
\setplotarea x from 0 to 16, y from -2 to 15
\put{295)} [l] at 2 12
\put {$ \scriptstyle \bullet$} [c] at  10 0
\put {$ \scriptstyle \bullet$} [c] at  10 6
\put {$ \scriptstyle \bullet$} [c] at  10 12
\put {$ \scriptstyle \bullet$} [c] at  11.5 9
\put {$ \scriptstyle \bullet$} [c] at  13 12
\put {$ \scriptstyle \bullet$} [c] at  16 6
\put {$ \scriptstyle \bullet$} [c] at  16 0
\setlinear \plot  10  12 10 0 16 6  16 0 10 6 13 12 16 6    /
\put{$2{,}520   $} [c] at 13 -2
\endpicture
\end{minipage}
\begin{minipage}{4cm}
\beginpicture
\setcoordinatesystem units   <1.5mm,2mm>
\setplotarea x from 0 to 16, y from -2 to 15
\put{296)} [l] at 2 12
\put {$ \scriptstyle \bullet$} [c] at  10 0
\put {$ \scriptstyle \bullet$} [c] at  10 4
\put {$ \scriptstyle \bullet$} [c] at  10 8
\put {$ \scriptstyle \bullet$} [c] at  10 12
\put {$ \scriptstyle \bullet$} [c] at  13 12
\put {$ \scriptstyle \bullet$} [c] at  16 0
\put {$ \scriptstyle \bullet$} [c] at  16  12
\setlinear \plot 10 0 10 4 16 12 16 0 10 8 13 12  /
\setlinear \plot  10 12  10 4     /
\put{$2{,}520$} [c] at 13 -2
\endpicture
\end{minipage}
\begin{minipage}{4cm}
\beginpicture
\setcoordinatesystem units   <1.5mm,2mm>
\setplotarea x from 0 to 16, y from -2 to 15
\put{297)} [l] at 2 12
\put {$ \scriptstyle \bullet$} [c] at  10 0
\put {$ \scriptstyle \bullet$} [c] at  10 4
\put {$ \scriptstyle \bullet$} [c] at  10 8
\put {$ \scriptstyle \bullet$} [c] at  10 12
\put {$ \scriptstyle \bullet$} [c] at  13 0
\put {$ \scriptstyle \bullet$} [c] at  16 0
\put {$ \scriptstyle \bullet$} [c] at  16  12
\setlinear \plot 10 12 10 8 16 0 16 12 10 4 13 0  /
\setlinear \plot  10 0  10 8     /
\put{$2{,}520$} [c] at 13 -2
\endpicture
\end{minipage}
\begin{minipage}{4cm}
\beginpicture
\setcoordinatesystem units   <1.5mm,2mm>
\setplotarea x from 0 to 16, y from -2 to 15
\put{298)} [l] at 2 12
\put {$ \scriptstyle \bullet$} [c] at  10 12
\put {$ \scriptstyle \bullet$} [c] at  10 8
\put {$ \scriptstyle \bullet$} [c] at  12 0
\put {$ \scriptstyle \bullet$} [c] at  12 4
\put {$ \scriptstyle \bullet$} [c] at  14 12
\put {$ \scriptstyle \bullet$} [c] at  14 8
\put {$ \scriptstyle \bullet$} [c] at  16 0
\setlinear \plot  12 0 12 4 10  8 10 12 16 0 14 12 14 8 12 4    /
\setlinear \plot  10 12 14 8     /
\setlinear \plot  10 8  14 12 /
\put{$1{,}260 $} [c] at 13 -2
\endpicture
\end{minipage}
\begin{minipage}{4cm}
\beginpicture
\setcoordinatesystem units   <1.5mm,2mm>
\setplotarea x from 0 to 16, y from -2 to 15
\put{299)} [l] at 2 12
\put {$ \scriptstyle \bullet$} [c] at  10 0
\put {$ \scriptstyle \bullet$} [c] at  10 4
\put {$ \scriptstyle \bullet$} [c] at  12 12
\put {$ \scriptstyle \bullet$} [c] at  12 8
\put {$ \scriptstyle \bullet$} [c] at  14 0
\put {$ \scriptstyle \bullet$} [c] at  14 4
\put {$ \scriptstyle \bullet$} [c] at  16 12
\setlinear \plot  12 12 12 8 10  4 10 0 16 12 14 0 14 4 12 8    /
\setlinear \plot  10 0 14 4     /
\setlinear \plot  10 4  14 0 /
\put{$1{,}260  $} [c] at 13 -2
\endpicture
\end{minipage}
\begin{minipage}{4cm}
\beginpicture
\setcoordinatesystem units   <1.5mm,2mm>
\setplotarea x from 0 to 16, y from -2 to 15
\put{300)} [l] at 2 12
\put {$ \scriptstyle \bullet$} [c] at  10 4
\put {$ \scriptstyle \bullet$} [c] at  11 0
\put {$ \scriptstyle \bullet$} [c] at  11 8
\put {$ \scriptstyle \bullet$} [c] at  11 12
\put {$ \scriptstyle \bullet$} [c] at  12 4
\put {$ \scriptstyle \bullet$} [c] at  16 0
\put {$ \scriptstyle \bullet$} [c] at  16 12
\setlinear \plot  11 0  10 4 11 8 12 4 11 0   /
\setlinear \plot 11 8 11 12 16 0 16 12 11 8  /
\put{$1{,}260  $} [c] at 13 -2
\endpicture
\end{minipage}
$$
$$
\begin{minipage}{4cm}
\beginpicture
\setcoordinatesystem units   <1.5mm,2mm>
\setplotarea x from 0 to 16, y from -2 to 15
\put{301)} [l] at 2 12
\put {$ \scriptstyle \bullet$} [c] at  10 8
\put {$ \scriptstyle \bullet$} [c] at  11 0
\put {$ \scriptstyle \bullet$} [c] at  11 4
\put {$ \scriptstyle \bullet$} [c] at  11 12
\put {$ \scriptstyle \bullet$} [c] at  12 8
\put {$ \scriptstyle \bullet$} [c] at  16 0
\put {$ \scriptstyle \bullet$} [c] at  16 12
\setlinear \plot  11 12  10 8 11 4 12 8 11 12   /
\setlinear \plot 11 4 11 0 16 12 16 0 11 4  /
\put{$1{,}260  $} [c] at 13 -2
\endpicture
\end{minipage}
\begin{minipage}{4cm}
\beginpicture
\setcoordinatesystem units   <1.5mm,2mm>
\setplotarea x from 0 to 16, y from -2 to 15
\put{302)} [l] at 2 12
\put {$ \scriptstyle \bullet$} [c] at  10 0
\put {$ \scriptstyle \bullet$} [c] at  10 6
\put {$ \scriptstyle \bullet$} [c] at  10 12
\put {$ \scriptstyle \bullet$} [c] at  13 0
\put {$ \scriptstyle \bullet$} [c] at  13 6
\put {$ \scriptstyle \bullet$} [c] at  13 12
\put {$ \scriptstyle \bullet$} [c] at  16 6
\setlinear \plot 10 0 10 12 16 6 13  0 10 6 13 12 16 6    /
\setlinear \plot 13 12 13  0  /
\setlinear \plot 10 12 13  6 10 0 /
\put{$1{,}260 $} [c] at 13 -2
\endpicture
\end{minipage}
\begin{minipage}{4cm}
\beginpicture
\setcoordinatesystem units   <1.5mm,2mm>
\setplotarea x from 0 to 16, y from -2 to 15
\put{303)} [l] at 2 12
\put {$ \scriptstyle \bullet$} [c] at  10 0
\put {$ \scriptstyle \bullet$} [c] at  10 6
\put {$ \scriptstyle \bullet$} [c] at  10 12
\put {$ \scriptstyle \bullet$} [c] at  13 0
\put {$ \scriptstyle \bullet$} [c] at  13 6
\put {$ \scriptstyle \bullet$} [c] at  13 12
\put {$ \scriptstyle \bullet$} [c] at  16 6
\setlinear \plot 10 12 10 0 16 6 13  12 10 6 13 0 16 6    /
\setlinear \plot 13 12 13  0  /
\setlinear \plot 10 0 13  6 10 12 /
\put{$1{,}260  $} [c] at 13 -2
\endpicture
\end{minipage}
\begin{minipage}{4cm}
\beginpicture
\setcoordinatesystem units   <1.5mm,2mm>
\setplotarea x from 0 to 16, y from -2 to 15
\put{304)} [l] at 2 12
\put {$ \scriptstyle \bullet$} [c] at  13 0
\put {$ \scriptstyle \bullet$} [c] at  13 6
\put {$ \scriptstyle \bullet$} [c] at  13 12
\put {$ \scriptstyle \bullet$} [c] at  10 12
\put {$ \scriptstyle \bullet$} [c] at  16  0
\put {$ \scriptstyle \bullet$} [c] at  16  6
\put {$ \scriptstyle \bullet$} [c] at  16 12
\setlinear \plot  10 12 13 6 16 0 16 12 13 6 13 12 16 6 13 0 13 6 16 12  /
\put{$1{,}260$} [c] at 13 -2
\endpicture
\end{minipage}
\begin{minipage}{4cm}
\beginpicture
\setcoordinatesystem units   <1.5mm,2mm>
\setplotarea x from 0 to 16, y from -2 to 15
\put{305)} [l] at 2 12
\put {$ \scriptstyle \bullet$} [c] at  13 0
\put {$ \scriptstyle \bullet$} [c] at  13 6
\put {$ \scriptstyle \bullet$} [c] at  13 12
\put {$ \scriptstyle \bullet$} [c] at  10 0
\put {$ \scriptstyle \bullet$} [c] at  16  0
\put {$ \scriptstyle \bullet$} [c] at  16  6
\put {$ \scriptstyle \bullet$} [c] at  16 12
\setlinear \plot  10 0 13 6 16 0 16 12 13 6 13 12 16 6 13 0 13 6 16 12  /
\put{$1{,}260$} [c] at 13 -2
\endpicture
\end{minipage}
\begin{minipage}{4cm}
\beginpicture
\setcoordinatesystem units   <1.5mm,2mm>
\setplotarea x from 0 to 16, y from -2 to 15
\put{306)} [l] at 2 12
\put {$ \scriptstyle \bullet$} [c] at  10 0
\put {$ \scriptstyle \bullet$} [c] at  10 4
\put {$ \scriptstyle \bullet$} [c] at  10 8
\put {$ \scriptstyle \bullet$} [c] at  10 12
\put {$ \scriptstyle \bullet$} [c] at  13 12
\put {$ \scriptstyle \bullet$} [c] at  16 0
\put {$ \scriptstyle \bullet$} [c] at  16  12
\setlinear \plot 10 12 10 0 16 12 16 0 10 4   /
\setlinear \plot  13 12  10 8     /
\put{$1{,}260$} [c] at 13 -2
\endpicture
\end{minipage}
$$
$$
\begin{minipage}{4cm}
\beginpicture
\setcoordinatesystem units   <1.5mm,2mm>
\setplotarea x from 0 to 16, y from -2 to 15
\put{307)} [l] at 2 12
\put {$ \scriptstyle \bullet$} [c] at  10 0
\put {$ \scriptstyle \bullet$} [c] at  10 4
\put {$ \scriptstyle \bullet$} [c] at  10 8
\put {$ \scriptstyle \bullet$} [c] at  10 12
\put {$ \scriptstyle \bullet$} [c] at  13 0
\put {$ \scriptstyle \bullet$} [c] at  16 0
\put {$ \scriptstyle \bullet$} [c] at  16  12
\setlinear \plot 10 0 10 12 16 0 16 12 10 8   /
\setlinear \plot  13 0  10 4     /
\put{$1{,}260$} [c] at 13 -2
\endpicture
\end{minipage}
\begin{minipage}{4cm}
\beginpicture
\setcoordinatesystem units   <1.5mm,2mm>
\setplotarea x from 0 to 16, y from -2 to 15
\put{308)} [l] at 2 12
\put {$ \scriptstyle \bullet$} [c] at  10 0
\put {$ \scriptstyle \bullet$} [c] at  10 4
\put {$ \scriptstyle \bullet$} [c] at  10 8
\put {$ \scriptstyle \bullet$} [c] at  10 12
\put {$ \scriptstyle \bullet$} [c] at  13 12
\put {$ \scriptstyle \bullet$} [c] at  16 0
\put {$ \scriptstyle \bullet$} [c] at  16  12
\setlinear \plot  10 0 10 12 16 0 16 12 10 8 13 12 16 0   /
\put{$840$} [c] at 13 -2
\endpicture
\end{minipage}
\begin{minipage}{4cm}
\beginpicture
\setcoordinatesystem units   <1.5mm,2mm>
\setplotarea x from 0 to 16, y from -2 to 15
\put{309)} [l] at 2 12
\put {$ \scriptstyle \bullet$} [c] at  10 0
\put {$ \scriptstyle \bullet$} [c] at  10 4
\put {$ \scriptstyle \bullet$} [c] at  10 8
\put {$ \scriptstyle \bullet$} [c] at  10 12
\put {$ \scriptstyle \bullet$} [c] at  13 0
\put {$ \scriptstyle \bullet$} [c] at  16 0
\put {$ \scriptstyle \bullet$} [c] at  16  12
\setlinear \plot  10 12 10 0 16 12 16 0 10 4 13 0 16 12   /
\put{$840$} [c] at 13 -2
\endpicture
\end{minipage}
\begin{minipage}{4cm}
\beginpicture
\setcoordinatesystem units   <1.5mm,2mm>
\setplotarea x from 0 to 16, y from -2 to 15
\put{310)} [l] at 2 12
\put {$ \scriptstyle \bullet$} [c] at  10 0
\put {$ \scriptstyle \bullet$} [c] at  10 6
\put {$ \scriptstyle \bullet$} [c] at  10 12
\put {$ \scriptstyle \bullet$} [c] at  13 0
\put {$ \scriptstyle \bullet$} [c] at  13 12
\put {$ \scriptstyle \bullet$} [c] at  16 6
\put {$ \scriptstyle \bullet$} [c] at  16 12
\setlinear \plot  10 0 10 12 16 6  16 12 10 6 13 12 16 6 13 0 10 6 /
\put{$840$} [c] at 13 -2
\endpicture
\end{minipage}
\begin{minipage}{4cm}
\beginpicture
\setcoordinatesystem units   <1.5mm,2mm>
\setplotarea x from 0 to 16, y from -2 to 15
\put{311)} [l] at 2 12
\put {$ \scriptstyle \bullet$} [c] at  10 0
\put {$ \scriptstyle \bullet$} [c] at  10 6
\put {$ \scriptstyle \bullet$} [c] at  10 12
\put {$ \scriptstyle \bullet$} [c] at  13 0
\put {$ \scriptstyle \bullet$} [c] at  13 12
\put {$ \scriptstyle \bullet$} [c] at  16 6
\put {$ \scriptstyle \bullet$} [c] at  16 0
\setlinear \plot  10 12 10 0 16 6  16 0 10 6 13 0 16 6 13 12 10 6 /
\put{$840$} [c] at 13 -2
\endpicture
\end{minipage}
\begin{minipage}{4cm}
\beginpicture
\setcoordinatesystem units   <1.5mm,2mm>
\setplotarea x from 0 to 16, y from -2 to 15
\put{${\bf  16}$} [l] at 2 15
\put{312)} [l] at 2 12
\put {$ \scriptstyle \bullet$} [c] at 10 6
\put {$ \scriptstyle \bullet$} [c] at 11.5 9
\put {$ \scriptstyle \bullet$} [c] at 11.5 3
\put {$ \scriptstyle \bullet$} [c] at 13 0
\put {$ \scriptstyle \bullet$} [c] at 13 6
\put {$ \scriptstyle \bullet$} [c] at 13 12
\put {$ \scriptstyle \bullet$} [c] at 16 6
\setlinear \plot 13 12 16 6 13 0 10 6 13 12 13 6 11.5 3 /
\put{$5{,}040$} [c] at 13 -2
\endpicture
\end{minipage}
$$
$$
\begin{minipage}{4cm}
\beginpicture
\setcoordinatesystem units   <1.5mm,2mm>
\setplotarea x from 0 to 16, y from -2 to 15
\put{313)} [l] at 2 12
\put {$ \scriptstyle \bullet$} [c] at 10 6
\put {$ \scriptstyle \bullet$} [c] at 11.5 9
\put {$ \scriptstyle \bullet$} [c] at 11.5 3
\put {$ \scriptstyle \bullet$} [c] at 13 0
\put {$ \scriptstyle \bullet$} [c] at 13 6
\put {$ \scriptstyle \bullet$} [c] at 13 12
\put {$ \scriptstyle \bullet$} [c] at 16 6
\setlinear \plot 13 0 16 6 13 12 10 6 13 0 13 6 11.5 9 /
\put{$5{,}040$} [c] at 13 -2
\endpicture
\end{minipage}
\begin{minipage}{4cm}
\beginpicture
\setcoordinatesystem units   <1.5mm,2mm>
\setplotarea x from 0 to 16, y from -2 to 15
\put{314)} [l] at 2 12
\put {$ \scriptstyle \bullet$} [c] at 10 8
\put {$ \scriptstyle \bullet$} [c] at 10.5 6
\put {$ \scriptstyle \bullet$} [c] at 11 0
\put {$ \scriptstyle \bullet$} [c] at 11 4
\put {$ \scriptstyle \bullet$} [c] at 11 12
\put {$ \scriptstyle \bullet$} [c] at 12 8
\put {$ \scriptstyle \bullet$} [c] at 16 12
\setlinear \plot 11 0 11 4 10 8 11 12 12 8 11 4 16 12  /
\put{$5{,}040$} [c] at 13 -2
\endpicture
\end{minipage}
\begin{minipage}{4cm}
\beginpicture
\setcoordinatesystem units   <1.5mm,2mm>
\setplotarea x from 0 to 16, y from -2 to 15
\put{315)} [l] at 2 12
\put {$ \scriptstyle \bullet$} [c] at 10 4
\put {$ \scriptstyle \bullet$} [c] at 10.5 6
\put {$ \scriptstyle \bullet$} [c] at 11 12
\put {$ \scriptstyle \bullet$} [c] at 11 8
\put {$ \scriptstyle \bullet$} [c] at 11 0
\put {$ \scriptstyle \bullet$} [c] at 12 4
\put {$ \scriptstyle \bullet$} [c] at 16 0
\setlinear \plot 11 12 11 8 10 4 11 0 12 4 11 8 16 0  /
\put{$5{,}040$} [c] at 13 -2
\endpicture
\end{minipage}
\begin{minipage}{4cm}
\beginpicture
\setcoordinatesystem units   <1.5mm,2mm>
\setplotarea x from 0 to 16, y from -2 to 15
\put{316)} [l] at 2 12
\put {$ \scriptstyle \bullet$} [c] at 10 8
\put {$ \scriptstyle \bullet$} [c] at 12 4
\put {$ \scriptstyle \bullet$} [c] at 12 12
\put {$ \scriptstyle \bullet$} [c] at 14 0
\put {$ \scriptstyle \bullet$} [c] at 14 8
\put {$ \scriptstyle \bullet$} [c] at 16 4
\put {$ \scriptstyle \bullet$} [c] at 16 12
\setlinear \plot 12 4 14 8 16 12 16 4 14 0 12 4 10 8 12 12 14 8  /
\put{$5{,}040$} [c] at 13 -2
\endpicture
\end{minipage}
\begin{minipage}{4cm}
\beginpicture
\setcoordinatesystem units   <1.5mm,2mm>
\setplotarea x from 0 to 16, y from -2 to 15
\put{317)} [l] at 2 12
\put {$ \scriptstyle \bullet$} [c] at 10 4
\put {$ \scriptstyle \bullet$} [c] at 12 8
\put {$ \scriptstyle \bullet$} [c] at 12 0
\put {$ \scriptstyle \bullet$} [c] at 14 12
\put {$ \scriptstyle \bullet$} [c] at 14 4
\put {$ \scriptstyle \bullet$} [c] at 16 8
\put {$ \scriptstyle \bullet$} [c] at 16 0
\setlinear \plot 12 8 14 4 16 0 16 8 14 12 12 8 10 4 12 0 14 4  /
\put{$5{,}040$} [c] at 13 -2
\endpicture
\end{minipage}
\begin{minipage}{4cm}
\beginpicture
\setcoordinatesystem units   <1.5mm,2mm>
\setplotarea x from 0 to 16, y from -2 to 15
\put{318)} [l] at 2 12
\put {$ \scriptstyle \bullet$} [c] at 13 0
\put {$ \scriptstyle \bullet$} [c] at 10 9
\put {$ \scriptstyle \bullet$} [c] at 13 12
\put {$ \scriptstyle \bullet$} [c] at 16 3
\put {$ \scriptstyle \bullet$} [c] at 16 12
\put {$ \scriptstyle \bullet$} [c] at 11.5 4.5
\put {$ \scriptstyle \bullet$} [c] at 14.5 7.5
\setlinear \plot 13 0 10 9 13 12 16 3 13 0    /
\setlinear \plot 16 12 16 3     /
\setlinear \plot 11.5 4.5 14.5 7.5     /
\put{$5{,}040$} [c] at 13 -2
\endpicture
\end{minipage}
$$
$$
\begin{minipage}{4cm}
\beginpicture
\setcoordinatesystem units   <1.5mm,2mm>
\setplotarea x from 0 to 16, y from -2 to 15
\put{319)} [l] at 2 12
\put {$ \scriptstyle \bullet$} [c] at 13 12
\put {$ \scriptstyle \bullet$} [c] at 10 3
\put {$ \scriptstyle \bullet$} [c] at 13 0
\put {$ \scriptstyle \bullet$} [c] at 16 9
\put {$ \scriptstyle \bullet$} [c] at 16 0
\put {$ \scriptstyle \bullet$} [c] at 11.5 7.5
\put {$ \scriptstyle \bullet$} [c] at 14.5 4.5
\setlinear \plot 13 12 10 3 13 0 16 9 13 12    /
\setlinear \plot 16 0 16 9     /
\setlinear \plot 11.5 7.5 14.5 4.5     /
\put{$5{,}040$} [c] at 13 -2
\endpicture
\end{minipage}
\begin{minipage}{4cm}
\beginpicture
\setcoordinatesystem units   <1.5mm,2mm>
\setplotarea x from 0 to 16, y from -2 to 15
\put{320)} [l] at 2 12
\put {$ \scriptstyle \bullet$} [c] at 10 4
\put {$ \scriptstyle \bullet$} [c] at 10 12
\put {$ \scriptstyle \bullet$} [c] at 13  0
\put {$ \scriptstyle \bullet$} [c] at 13 4
\put {$ \scriptstyle \bullet$} [c] at 16 4
\put {$ \scriptstyle \bullet$} [c] at 16 8
\put {$ \scriptstyle \bullet$} [c] at 16 12
\setlinear \plot 10 4 10 12 13 4 13 0 10 4 16 12 16 4 13 0   /
\setlinear \plot 13 4 16 8    /
\put{$5{,}040$} [c] at 13 -2
\endpicture
\end{minipage}
\begin{minipage}{4cm}
\beginpicture
\setcoordinatesystem units   <1.5mm,2mm>
\setplotarea x from 0 to 16, y from -2 to 15
\put{321)} [l] at 2 12
\put {$ \scriptstyle \bullet$} [c] at 10 8
\put {$ \scriptstyle \bullet$} [c] at 10 0
\put {$ \scriptstyle \bullet$} [c] at 13  12
\put {$ \scriptstyle \bullet$} [c] at 13 8
\put {$ \scriptstyle \bullet$} [c] at 16 4
\put {$ \scriptstyle \bullet$} [c] at 16 8
\put {$ \scriptstyle \bullet$} [c] at 16 0
\setlinear \plot 10 8 10 0 13 8 13 12 10 8 16 0 16 8 13 12   /
\setlinear \plot 13 8 16 4    /
\put{$5{,}040$} [c] at 13 -2
\endpicture
\end{minipage}
\begin{minipage}{4cm}
\beginpicture
\setcoordinatesystem units   <1.5mm,2mm>
\setplotarea x from 0 to 16, y from -2 to 15
\put{322)} [l] at 2 12
\put {$ \scriptstyle \bullet$} [c] at 10 4
\put {$ \scriptstyle \bullet$} [c] at 10 8
\put {$ \scriptstyle \bullet$} [c] at 10 12
\put {$ \scriptstyle \bullet$} [c] at 13 0
\put {$ \scriptstyle \bullet$} [c] at 16 4
\put {$ \scriptstyle \bullet$} [c] at 16  8
\put {$ \scriptstyle \bullet$} [c] at 16  12
\setlinear \plot 10 12 10 4 13 0 16 4 16 12 10 4    /
\put{$5{,}040$} [c] at 13 -2
\endpicture
\end{minipage}
\begin{minipage}{4cm}
\beginpicture
\setcoordinatesystem units   <1.5mm,2mm>
\setplotarea x from 0 to 16, y from -2 to 15
\put{323)} [l] at 2 12
\put {$ \scriptstyle \bullet$} [c] at 10 4
\put {$ \scriptstyle \bullet$} [c] at 10 8
\put {$ \scriptstyle \bullet$} [c] at 10 0
\put {$ \scriptstyle \bullet$} [c] at 13 12
\put {$ \scriptstyle \bullet$} [c] at 16 4
\put {$ \scriptstyle \bullet$} [c] at 16  8
\put {$ \scriptstyle \bullet$} [c] at 16  0
\setlinear \plot 10 0 10 8 13 12 16 8 16 0 10 8    /
\put{$5{,}040$} [c] at 13 -2
\endpicture
\end{minipage}
\begin{minipage}{4cm}
\beginpicture
\setcoordinatesystem units   <1.5mm,2mm>
\setplotarea x from 0 to 16, y from -2 to 15
\put{324)} [l] at 2 12
\put {$ \scriptstyle \bullet$} [c] at 10 3
\put {$ \scriptstyle \bullet$} [c] at 10 6
\put {$ \scriptstyle \bullet$} [c] at 10 9
\put {$ \scriptstyle \bullet$} [c] at 10 12
\put {$ \scriptstyle \bullet$} [c] at 13 0
\put {$ \scriptstyle \bullet$} [c] at 16 3
\put {$ \scriptstyle \bullet$} [c] at 16 12
\setlinear \plot 10 12 10 3 13 0 16 3 16 12    /
\put{$5{,}040$} [c] at 13 -2
\endpicture
\end{minipage}
$$
$$
\begin{minipage}{4cm}
\beginpicture
\setcoordinatesystem units   <1.5mm,2mm>
\setplotarea x from 0 to 16, y from -2 to 15
\put{325)} [l] at 2 12
\put {$ \scriptstyle \bullet$} [c] at 10 0
\put {$ \scriptstyle \bullet$} [c] at 10 3
\put {$ \scriptstyle \bullet$} [c] at 10 6
\put {$ \scriptstyle \bullet$} [c] at 10 9
\put {$ \scriptstyle \bullet$} [c] at 13 12
\put {$ \scriptstyle \bullet$} [c] at 16 0
\put {$ \scriptstyle \bullet$} [c] at 16 9
\setlinear \plot 10 0 10 9 13 12 16 9 16 0    /
\put{$5{,}040$} [c] at 13 -2
\endpicture
\end{minipage}
\begin{minipage}{4cm}
\beginpicture
\setcoordinatesystem units   <1.5mm,2mm>
\setplotarea x from 0 to 16, y from -2 to 15
\put{326)} [l] at 2 12
\put {$ \scriptstyle \bullet$} [c] at 10 12
\put {$ \scriptstyle \bullet$} [c] at 13  0
\put {$ \scriptstyle \bullet$} [c] at 13 3
\put {$ \scriptstyle \bullet$} [c] at 13 6
\put {$ \scriptstyle \bullet$} [c] at 13 9
\put {$ \scriptstyle \bullet$} [c] at 13 12
\put {$ \scriptstyle \bullet$} [c] at 16 12
\setlinear \plot 13 12  13 0   /
\setlinear \plot 10 12  13  3     /
\setlinear \plot 16 12  13  6     /
\put{$5{,}040$} [c] at 13 -2
\endpicture
\end{minipage}
\begin{minipage}{4cm}
\beginpicture
\setcoordinatesystem units   <1.5mm,2mm>
\setplotarea x from 0 to 16, y from -2 to 15
\put{327)} [l] at 2 12
\put {$ \scriptstyle \bullet$} [c] at 10 0
\put {$ \scriptstyle \bullet$} [c] at 13  0
\put {$ \scriptstyle \bullet$} [c] at 13 3
\put {$ \scriptstyle \bullet$} [c] at 13 6
\put {$ \scriptstyle \bullet$} [c] at 13 9
\put {$ \scriptstyle \bullet$} [c] at 13 12
\put {$ \scriptstyle \bullet$} [c] at 16 0
\setlinear \plot 13 12  13 0   /
\setlinear \plot 10 0  13  9     /
\setlinear \plot 16 0  13  6     /
\put{$5{,}040$} [c] at 13 -2
\endpicture
\end{minipage}
\begin{minipage}{4cm}
\beginpicture
\setcoordinatesystem units   <1.5mm,2mm>
\setplotarea x from 0 to 16, y from -2 to 15
\put{328)} [l] at 2 12
\put {$ \scriptstyle \bullet$} [c] at 10 4
\put {$ \scriptstyle \bullet$} [c] at 10 12
\put {$ \scriptstyle \bullet$} [c] at 13 12
\put {$ \scriptstyle \bullet$} [c] at 13 0
\put {$ \scriptstyle \bullet$} [c] at 16 4
\put {$ \scriptstyle \bullet$} [c] at 16 8
\put {$ \scriptstyle \bullet$} [c] at 16 12
\setlinear \plot 16 4 13 0 10 4 10 12  16 4 16 12   /
\setlinear \plot 10 12 16 4     /
\setlinear \plot 10 4 13 12 16 8     /
\put{$5{,}040$} [c] at 13 -2
\endpicture
\end{minipage}
\begin{minipage}{4cm}
\beginpicture
\setcoordinatesystem units   <1.5mm,2mm>
\setplotarea x from 0 to 16, y from -2 to 15
\put{329)} [l] at 2 12
\put {$ \scriptstyle \bullet$} [c] at 10 8
\put {$ \scriptstyle \bullet$} [c] at 10 0
\put {$ \scriptstyle \bullet$} [c] at 13 12
\put {$ \scriptstyle \bullet$} [c] at 13 0
\put {$ \scriptstyle \bullet$} [c] at 16 4
\put {$ \scriptstyle \bullet$} [c] at 16 8
\put {$ \scriptstyle \bullet$} [c] at 16 0
\setlinear \plot 16 8 13 12  10 8 10 0  16 8 16 0   /
\setlinear \plot 10 0 16 8     /
\setlinear \plot 10 8 13 0 16 4     /
\put{$5{,}040$} [c] at 13 -2
\endpicture
\end{minipage}
\begin{minipage}{4cm}
\beginpicture
\setcoordinatesystem units   <1.5mm,2mm>
\setplotarea x from 0 to 16, y from -2 to 15
\put{330)} [l] at 2 12
\put {$ \scriptstyle \bullet$} [c] at  10 0
\put {$ \scriptstyle \bullet$} [c] at  10 12
\put {$ \scriptstyle \bullet$} [c] at  11.5 9
\put {$ \scriptstyle \bullet$} [c] at  13 0
\put {$ \scriptstyle \bullet$} [c] at  13 3
\put {$ \scriptstyle \bullet$} [c] at  13 6
\put {$ \scriptstyle \bullet$} [c] at  16 12
\setlinear \plot  13 0 13 6 10 12 10 0  /
\setlinear \plot  13 6 16 12 /
\put{$5{,}040$} [c] at 13 -2
\endpicture
\end{minipage}
$$
$$
\begin{minipage}{4cm}
\beginpicture
\setcoordinatesystem units   <1.5mm,2mm>
\setplotarea x from 0 to 16, y from -2 to 15
\put{331)} [l] at 2 12
\put {$ \scriptstyle \bullet$} [c] at  10 0
\put {$ \scriptstyle \bullet$} [c] at  10 12
\put {$ \scriptstyle \bullet$} [c] at  11.5 3
\put {$ \scriptstyle \bullet$} [c] at  13 12
\put {$ \scriptstyle \bullet$} [c] at  13 9
\put {$ \scriptstyle \bullet$} [c] at  13 6
\put {$ \scriptstyle \bullet$} [c] at  16 0
\setlinear \plot  10 12 10 0 13 6 13 12   /
\setlinear \plot  16 0 13 6 /
\put{$5{,}040  $} [c] at 13 -2
\endpicture
\end{minipage}
\begin{minipage}{4cm}
\beginpicture
\setcoordinatesystem units   <1.5mm,2mm>
\setplotarea x from 0 to 16, y from -2 to 15
\put{332)} [l] at 2 12
\put {$ \scriptstyle \bullet$} [c] at  10 0
\put {$ \scriptstyle \bullet$} [c] at  10 3
\put {$ \scriptstyle \bullet$} [c] at  10 6
\put {$ \scriptstyle \bullet$} [c] at  10 9
\put {$ \scriptstyle \bullet$} [c] at  10 12
\put {$ \scriptstyle \bullet$} [c] at  16 0
\put {$ \scriptstyle \bullet$} [c] at  16 12
\setlinear \plot  10 0 10 12   /
\setlinear \plot  16 0 16 12  10 3   /
\put{$5{,}040  $} [c] at 13 -2
\endpicture
\end{minipage}
\begin{minipage}{4cm}
\beginpicture
\setcoordinatesystem units   <1.5mm,2mm>
\setplotarea x from 0 to 16, y from -2 to 15
\put{333)} [l] at 2 12
\put {$ \scriptstyle \bullet$} [c] at  10 0
\put {$ \scriptstyle \bullet$} [c] at  10 3
\put {$ \scriptstyle \bullet$} [c] at  10 6
\put {$ \scriptstyle \bullet$} [c] at  10 9
\put {$ \scriptstyle \bullet$} [c] at  10 12
\put {$ \scriptstyle \bullet$} [c] at  16 0
\put {$ \scriptstyle \bullet$} [c] at  16 12
\setlinear \plot  10 0 10 12   /
\setlinear \plot  16 12 16 0  10 9   /
\put{$5{,}040  $} [c] at 13 -2
\endpicture
\end{minipage}
\begin{minipage}{4cm}
\beginpicture
\setcoordinatesystem units   <1.5mm,2mm>
\setplotarea x from 0 to 16, y from -2 to 15
\put{334)} [l] at 2 12
\put {$ \scriptstyle \bullet$} [c] at  10 6
\put {$ \scriptstyle \bullet$} [c] at  10 9
\put {$ \scriptstyle \bullet$} [c] at  10 12
\put {$ \scriptstyle \bullet$} [c] at  12 0
\put {$ \scriptstyle \bullet$} [c] at  14 6
\put {$ \scriptstyle \bullet$} [c] at  14 12
\put {$ \scriptstyle \bullet$} [c] at  16 0
\setlinear \plot 10 9 16 0  14 12 14 6 10 12  10 6 14 12 /
\setlinear \plot 10 6 12  0  14 6  /
\put{$5{,}040  $} [c] at 13 -2
\endpicture
\end{minipage}
\begin{minipage}{4cm}
\beginpicture
\setcoordinatesystem units   <1.5mm,2mm>
\setplotarea x from 0 to 16, y from -2 to 15
\put{335)} [l] at 2 12
\put {$ \scriptstyle \bullet$} [c] at  10 6
\put {$ \scriptstyle \bullet$} [c] at  10 3
\put {$ \scriptstyle \bullet$} [c] at  10 0
\put {$ \scriptstyle \bullet$} [c] at  12 12
\put {$ \scriptstyle \bullet$} [c] at  14 6
\put {$ \scriptstyle \bullet$} [c] at  14 0
\put {$ \scriptstyle \bullet$} [c] at  16 12
\setlinear \plot 10 3 16 12  14 0 14 6 10 0  10 6 14 0 /
\setlinear \plot 10 6 12  12  14 6  /
\put{$5{,}040  $} [c] at 13 -2
\endpicture
\end{minipage}
\begin{minipage}{4cm}
\beginpicture
\setcoordinatesystem units   <1.5mm,2mm>
\setplotarea x from 0 to 16, y from -2 to 15
\put{336)} [l] at 2 12
\put {$ \scriptstyle \bullet$} [c] at  10 12
\put {$ \scriptstyle \bullet$} [c] at  10 9
\put {$ \scriptstyle \bullet$} [c] at  10 6
\put {$ \scriptstyle \bullet$} [c] at  10 0
\put {$ \scriptstyle \bullet$} [c] at  16 0
\put {$ \scriptstyle \bullet$} [c] at  16 6
\put {$ \scriptstyle \bullet$} [c] at  16 12
\setlinear \plot 10 12 10 0 16 12 16 0 10 6  /
\put{$5{,}040  $} [c] at 13 -2
\endpicture
\end{minipage}
$$
$$
\begin{minipage}{4cm}
\beginpicture
\setcoordinatesystem units   <1.5mm,2mm>
\setplotarea x from 0 to 16, y from -2 to 15
\put{337)} [l] at 2 12
\put {$ \scriptstyle \bullet$} [c] at  10 12
\put {$ \scriptstyle \bullet$} [c] at  10 3
\put {$ \scriptstyle \bullet$} [c] at  10 6
\put {$ \scriptstyle \bullet$} [c] at  10 0
\put {$ \scriptstyle \bullet$} [c] at  16 0
\put {$ \scriptstyle \bullet$} [c] at  16 6
\put {$ \scriptstyle \bullet$} [c] at  16 12
\setlinear \plot 10 0 10 12 16 0 16 12 10 6  /
\put{$5{,}040  $} [c] at 13 -2
\endpicture
\end{minipage}
\begin{minipage}{4cm}
\beginpicture
\setcoordinatesystem units   <1.5mm,2mm>
\setplotarea x from 0 to 16, y from -2 to 15
\put{338)} [l] at 2 12
\put {$ \scriptstyle \bullet$} [c] at  10 4
\put {$ \scriptstyle \bullet$} [c] at  11 0
\put {$ \scriptstyle \bullet$} [c] at  11 8
\put {$ \scriptstyle \bullet$} [c] at  11 12
\put {$ \scriptstyle \bullet$} [c] at  12 4
\put {$ \scriptstyle \bullet$} [c] at  16 12
\put {$ \scriptstyle \bullet$} [c] at  16 0
\setlinear \plot  11  12 11 8 10 4  11 0 12 4 11 8 16 0 16 12 12 4   /
\put{$5{,}040  $} [c] at 13 -2
\endpicture
\end{minipage}
\begin{minipage}{4cm}
\beginpicture
\setcoordinatesystem units   <1.5mm,2mm>
\setplotarea x from 0 to 16, y from -2 to 15
\put{339)} [l] at 2 12
\put {$ \scriptstyle \bullet$} [c] at  10 8
\put {$ \scriptstyle \bullet$} [c] at  11 0
\put {$ \scriptstyle \bullet$} [c] at  11 4
\put {$ \scriptstyle \bullet$} [c] at  11 12
\put {$ \scriptstyle \bullet$} [c] at  12 8
\put {$ \scriptstyle \bullet$} [c] at  16 12
\put {$ \scriptstyle \bullet$} [c] at  16 0
\setlinear \plot  11  0 11 4 10 8  11 12 12 8 11 4 16 12 16 0 12 8   /
\put{$5{,}040  $} [c] at 13 -2
\endpicture
\end{minipage}
\begin{minipage}{4cm}
\beginpicture
\setcoordinatesystem units   <1.5mm,2mm>
\setplotarea x from 0 to 16, y from -2 to 15
\put{340)} [l] at 2 12
\put {$ \scriptstyle \bullet$} [c] at  10 0
\put {$ \scriptstyle \bullet$} [c] at  10 6
\put {$ \scriptstyle \bullet$} [c] at  10 12
\put {$ \scriptstyle \bullet$} [c] at  11.5 9
\put {$ \scriptstyle \bullet$} [c] at  13 0
\put {$ \scriptstyle \bullet$} [c] at  13 12
\put {$ \scriptstyle \bullet$} [c] at  16 6
\setlinear \plot  10  0 10 12    /
\setlinear \plot  13 0 10 6 13 12 16 6 13 0   /
\put{$5{,}040  $} [c] at 13 -2
\endpicture
\end{minipage}
\begin{minipage}{4cm}
\beginpicture
\setcoordinatesystem units   <1.5mm,2mm>
\setplotarea x from 0 to 16, y from -2 to 15
\put{341)} [l] at 2 12
\put {$ \scriptstyle \bullet$} [c] at  10 0
\put {$ \scriptstyle \bullet$} [c] at  10 6
\put {$ \scriptstyle \bullet$} [c] at  10 12
\put {$ \scriptstyle \bullet$} [c] at  11.5 3
\put {$ \scriptstyle \bullet$} [c] at  13 0
\put {$ \scriptstyle \bullet$} [c] at  13 12
\put {$ \scriptstyle \bullet$} [c] at  16 6
\setlinear \plot  10  0 10 12    /
\setlinear \plot  13 0 10 6 13 12 16 6 13 0   /
\put{$5{,}040  $} [c] at 13 -2
\endpicture
\end{minipage}
\begin{minipage}{4cm}
\beginpicture
\setcoordinatesystem units   <1.5mm,2mm>
\setplotarea x from 0 to 16, y from -2 to 15
\put{342)} [l] at 2 12
\put {$ \scriptstyle \bullet$} [c] at  10 6
\put {$ \scriptstyle \bullet$} [c] at  10 12
\put {$ \scriptstyle \bullet$} [c] at  11.5 3
\put {$ \scriptstyle \bullet$} [c] at  13 0
\put {$ \scriptstyle \bullet$} [c] at  13 12
\put {$ \scriptstyle \bullet$} [c] at  16 6
\put {$ \scriptstyle \bullet$} [c] at  16 0
\setlinear \plot  10  12 10 6 13 0 16 6 13 12 10 6    /
\setlinear \plot  16 6 16 0   /
\put{$5{,}040  $} [c] at 13 -2
\endpicture
\end{minipage}
$$
$$
\begin{minipage}{4cm}
\beginpicture
\setcoordinatesystem units   <1.5mm,2mm>
\setplotarea x from 0 to 16, y from -2 to 15
\put{343)} [l] at 2 12
\put {$ \scriptstyle \bullet$} [c] at  10 6
\put {$ \scriptstyle \bullet$} [c] at  10 0
\put {$ \scriptstyle \bullet$} [c] at  11.5 9
\put {$ \scriptstyle \bullet$} [c] at  13 0
\put {$ \scriptstyle \bullet$} [c] at  13 12
\put {$ \scriptstyle \bullet$} [c] at  16 6
\put {$ \scriptstyle \bullet$} [c] at  16 12
\setlinear \plot  10  0 10 6 13 12 16 6 13 0 10 6    /
\setlinear \plot  16 6 16 12   /
\put{$5{,}040  $} [c] at 13 -2
\endpicture
\end{minipage}
\begin{minipage}{4cm}
\beginpicture
\setcoordinatesystem units   <1.5mm,2mm>
\setplotarea x from 0 to 16, y from -2 to 15
\put{344)} [l] at 2 12
\put {$ \scriptstyle \bullet$} [c] at  10 4
\put {$ \scriptstyle \bullet$} [c] at  11 0
\put {$ \scriptstyle \bullet$} [c] at  11 8
\put {$ \scriptstyle \bullet$} [c] at  11 12
\put {$ \scriptstyle \bullet$} [c] at  12  4
\put {$ \scriptstyle \bullet$} [c] at  16 12
\put {$ \scriptstyle \bullet$} [c] at  16 0
\setlinear \plot 11 12 11 8 10 4 11 0 12 4 11 8  /
\setlinear \plot 11  0  16 12 16 0 12 4   /
\put{$5{,}040  $} [c] at 13 -2
\endpicture
\end{minipage}
\begin{minipage}{4cm}
\beginpicture
\setcoordinatesystem units   <1.5mm,2mm>
\setplotarea x from 0 to 16, y from -2 to 15
\put{345)} [l] at 2 12
\put {$ \scriptstyle \bullet$} [c] at  10 8
\put {$ \scriptstyle \bullet$} [c] at  11 0
\put {$ \scriptstyle \bullet$} [c] at  11 4
\put {$ \scriptstyle \bullet$} [c] at  11 12
\put {$ \scriptstyle \bullet$} [c] at  12  8
\put {$ \scriptstyle \bullet$} [c] at  16 12
\put {$ \scriptstyle \bullet$} [c] at  16 0
\setlinear \plot 11 0 11 4 10 8 11 12 12 8 11 4  /
\setlinear \plot 11  12  16 0 16 12 12 8   /
\put{$5{,}040  $} [c] at 13 -2
\endpicture
\end{minipage}
\begin{minipage}{4cm}
\beginpicture
\setcoordinatesystem units   <1.5mm,2mm>
\setplotarea x from 0 to 16, y from -2 to 15
\put{346)} [l] at 2 12
\put {$ \scriptstyle \bullet$} [c] at  10 12
\put {$ \scriptstyle \bullet$} [c] at  10 8
\put {$ \scriptstyle \bullet$} [c] at  10 4
\put {$ \scriptstyle \bullet$} [c] at  10  0
\put {$ \scriptstyle \bullet$} [c] at  16 12
\put {$ \scriptstyle \bullet$} [c] at  16 8
\put {$ \scriptstyle \bullet$} [c] at  16  0
\setlinear \plot  16 12  16 0 10 8 10 0 16 8     /
\setlinear \plot  10 12 10 8  /
\put{$5{,}040  $} [c] at 13 -2
\endpicture
\end{minipage}
\begin{minipage}{4cm}
\beginpicture
\setcoordinatesystem units   <1.5mm,2mm>
\setplotarea x from 0 to 16, y from -2 to 15
\put{347)} [l] at 2 12
\put {$ \scriptstyle \bullet$} [c] at  10 12
\put {$ \scriptstyle \bullet$} [c] at  10 8
\put {$ \scriptstyle \bullet$} [c] at  10 4
\put {$ \scriptstyle \bullet$} [c] at  10  0
\put {$ \scriptstyle \bullet$} [c] at  16 12
\put {$ \scriptstyle \bullet$} [c] at  16 4
\put {$ \scriptstyle \bullet$} [c] at  16  0
\setlinear \plot  16 0  16 12 10 4 10 12 16 4     /
\setlinear \plot  10 0 10 4  /
\put{$5{,}040  $} [c] at 13 -2
\endpicture
\end{minipage}
\begin{minipage}{4cm}
\beginpicture
\setcoordinatesystem units   <1.5mm,2mm>
\setplotarea x from 0 to 16, y from -2 to 15
\put{348)} [l] at 2 12
\put {$ \scriptstyle \bullet$} [c] at  10 6
\put {$ \scriptstyle \bullet$} [c] at  10 9
\put {$ \scriptstyle \bullet$} [c] at  10 12
\put {$ \scriptstyle \bullet$} [c] at  13 0
\put {$ \scriptstyle \bullet$} [c] at  13 12
\put {$ \scriptstyle \bullet$} [c] at  16 6
\put {$ \scriptstyle \bullet$} [c] at  16 0
\setlinear \plot  10 12 10 6 13 12  16 6 13 0 10 6  /
\setlinear \plot  10 12 16 0 16 6    /
\put{$5{,}040   $} [c] at 13 -2
\endpicture
\end{minipage}
$$
$$
\begin{minipage}{4cm}
\beginpicture
\setcoordinatesystem units   <1.5mm,2mm>
\setplotarea x from 0 to 16, y from -2 to 15
\put{349)} [l] at 2 12
\put {$ \scriptstyle \bullet$} [c] at  10 6
\put {$ \scriptstyle \bullet$} [c] at  10 3
\put {$ \scriptstyle \bullet$} [c] at  10 0
\put {$ \scriptstyle \bullet$} [c] at  13 0
\put {$ \scriptstyle \bullet$} [c] at  13 12
\put {$ \scriptstyle \bullet$} [c] at  16 6
\put {$ \scriptstyle \bullet$} [c] at  16 12
\setlinear \plot  10 0 10 6 13 0  16 6 13 12 10 6  /
\setlinear \plot  10 0 16 12 16 6    /
\put{$5{,}040   $} [c] at 13 -2
\endpicture
\end{minipage}
\begin{minipage}{4cm}
\beginpicture
\setcoordinatesystem units   <1.5mm,2mm>
\setplotarea x from 0 to 16, y from -2 to 15
\put{350)} [l] at 2 12
\put {$ \scriptstyle \bullet$} [c] at  10 8
\put {$ \scriptstyle \bullet$} [c] at  10 12
\put {$ \scriptstyle \bullet$} [c] at  13 0
\put {$ \scriptstyle \bullet$} [c] at  13 4
\put {$ \scriptstyle \bullet$} [c] at  16 0
\put {$ \scriptstyle \bullet$} [c] at  16 8
\put {$ \scriptstyle \bullet$} [c] at  16 12
\setlinear \plot  10 12 10 8 13 4 13 0    /
\setlinear \plot  16 0  16 12      /
\setlinear \plot  13 4 16 8       /
\put{$5{,}040   $} [c] at 13 -2
\endpicture
\end{minipage}
\begin{minipage}{4cm}
\beginpicture
\setcoordinatesystem units   <1.5mm,2mm>
\setplotarea x from 0 to 16, y from -2 to 15
\put{351)} [l] at 2 12
\put {$ \scriptstyle \bullet$} [c] at  10 4
\put {$ \scriptstyle \bullet$} [c] at  10 0
\put {$ \scriptstyle \bullet$} [c] at  13 12
\put {$ \scriptstyle \bullet$} [c] at  13 8
\put {$ \scriptstyle \bullet$} [c] at  16 0
\put {$ \scriptstyle \bullet$} [c] at  16 4
\put {$ \scriptstyle \bullet$} [c] at  16 12
\setlinear \plot  10 0 10 4 13 8 13 12    /
\setlinear \plot  16 0  16 12      /
\setlinear \plot  13 8 16 4       /
\put{$5{,}040   $} [c] at 13 -2
\endpicture
\end{minipage}
\begin{minipage}{4cm}
\beginpicture
\setcoordinatesystem units   <1.5mm,2mm>
\setplotarea x from 0 to 16, y from -2 to 15
\put{352)} [l] at 2 12
\put {$ \scriptstyle \bullet$} [c] at  10 0
\put {$ \scriptstyle \bullet$} [c] at  10  3
\put {$ \scriptstyle \bullet$} [c] at  10 6
\put {$ \scriptstyle \bullet$} [c] at  10 9
\put {$ \scriptstyle \bullet$} [c] at  10 12
\put {$ \scriptstyle \bullet$} [c] at  16 0
\put {$ \scriptstyle \bullet$} [c] at  16 12
\setlinear \plot  10 0 10 12 16  0 16 12 10 0   /
\put{$5{,}040  $} [c] at 13 -2
\endpicture
\end{minipage}
\begin{minipage}{4cm}
\beginpicture
\setcoordinatesystem units   <1.5mm,2mm>
\setplotarea x from 0 to 16, y from -2 to 15
\put{353)} [l] at 2 12
\put {$ \scriptstyle \bullet$} [c] at  10 0
\put {$ \scriptstyle \bullet$} [c] at  10 9
\put {$ \scriptstyle \bullet$} [c] at  10 12
\put {$ \scriptstyle \bullet$} [c] at  13 6
\put {$ \scriptstyle \bullet$} [c] at  16 0
\put {$ \scriptstyle \bullet$} [c] at  16 3
\put {$ \scriptstyle \bullet$} [c] at  16 12
\setlinear \plot  10 12 10 0  /
\setlinear \plot  16 12 16 0 /
\setlinear \plot  10 9 16 3   /
\put{$5{,}040$} [c] at 13 -2
\endpicture
\end{minipage}
\begin{minipage}{4cm}
\beginpicture
\setcoordinatesystem units   <1.5mm,2mm>
\setplotarea x from 0 to 16, y from -2 to 15
\put{354)} [l] at 2 12
\put {$ \scriptstyle \bullet$} [c] at  10 0
\put {$ \scriptstyle \bullet$} [c] at  10 3
\put {$ \scriptstyle \bullet$} [c] at  10 9
\put {$ \scriptstyle \bullet$} [c] at  10 12
\put {$ \scriptstyle \bullet$} [c] at  16 0
\put {$ \scriptstyle \bullet$} [c] at  16 6
\put {$ \scriptstyle \bullet$} [c] at  16 12
\setlinear \plot 10 12 10 9  16 0 16  12  10 3 10 0  /
\setlinear \plot  10 3 10 9     /
\put{$5{,}040   $} [c] at 13 -2
\endpicture
\end{minipage}
$$
$$
\begin{minipage}{4cm}
\beginpicture
\setcoordinatesystem units   <1.5mm,2mm>
\setplotarea x from 0 to 16, y from -2 to 15
\put{355)} [l] at 2 12
\put {$ \scriptstyle \bullet$} [c] at  10 3
\put {$ \scriptstyle \bullet$} [c] at  10 9
\put {$ \scriptstyle \bullet$} [c] at  13 0
\put {$ \scriptstyle \bullet$} [c] at  13 12
\put {$ \scriptstyle \bullet$} [c] at  16 0
\put {$ \scriptstyle \bullet$} [c] at  16 6
\put {$ \scriptstyle \bullet$} [c] at  16  12
\setlinear \plot  13 0 10 3 10 9 13 12  16 6 13 0  /
\setlinear \plot  16 0 16  12     /
\put{$5{,}040   $} [c] at 13 -2
\endpicture
\end{minipage}
\begin{minipage}{4cm}
\beginpicture
\setcoordinatesystem units   <1.5mm,2mm>
\setplotarea x from 0 to 16, y from -2 to 15
\put{356)} [l] at 2 12
\put {$ \scriptstyle \bullet$} [c] at  10  12
\put {$ \scriptstyle \bullet$} [c] at  10  6
\put {$ \scriptstyle \bullet$} [c] at  13 0
\put {$ \scriptstyle \bullet$} [c] at  13 6
\put {$ \scriptstyle \bullet$} [c] at  13  12
\put {$ \scriptstyle \bullet$} [c] at  16 0
\put {$ \scriptstyle \bullet$} [c] at  16 6
\setlinear \plot 16 0 16 6 13 12 10 6 10 12 16 0   /
\setlinear \plot 10 6 13 0 16 6   /
\setlinear \plot 13 0 13 12   /
\put{$5{,}040$} [c] at 13 -2
\endpicture
\end{minipage}
\begin{minipage}{4cm}
\beginpicture
\setcoordinatesystem units   <1.5mm,2mm>
\setplotarea x from 0 to 16, y from -2 to 15
\put{357)} [l] at 2 12
\put {$ \scriptstyle \bullet$} [c] at 10 8
\put {$ \scriptstyle \bullet$} [c] at 12 4
\put {$ \scriptstyle \bullet$} [c] at 12 8
\put {$ \scriptstyle \bullet$} [c] at 12 12
\put {$ \scriptstyle \bullet$} [c] at 14 0
\put {$ \scriptstyle \bullet$} [c] at 14 8
\put {$ \scriptstyle \bullet$} [c] at 16 4
\setlinear \plot 14 0 10 8 12 12 16 4 14 0    /
\setlinear \plot 14 8 12 4 12 12    /
\put{$2{,}520$} [c] at 13 -2
\endpicture
\end{minipage}
\begin{minipage}{4cm}
\beginpicture
\setcoordinatesystem units   <1.5mm,2mm>
\setplotarea x from 0 to 16, y from -2 to 15
\put{358)} [l] at 2 12
\put {$ \scriptstyle \bullet$} [c] at 10 4
\put {$ \scriptstyle \bullet$} [c] at 12 0
\put {$ \scriptstyle \bullet$} [c] at 12 4
\put {$ \scriptstyle \bullet$} [c] at 12 8
\put {$ \scriptstyle \bullet$} [c] at 14 4
\put {$ \scriptstyle \bullet$} [c] at 14 12
\put {$ \scriptstyle \bullet$} [c] at 16 8
\setlinear \plot 12 0 10 4 14 12 16 8 12 0    /
\setlinear \plot 14 4 12 8 12 0    /
\put{$2{,}520$} [c] at 13 -2
\endpicture
\end{minipage}
\begin{minipage}{4cm}
\beginpicture
\setcoordinatesystem units   <1.5mm,2mm>
\setplotarea x from 0 to 16, y from -2 to 15
\put{359)} [l] at 2 12
\put {$ \scriptstyle \bullet$} [c] at 10 8
\put {$ \scriptstyle \bullet$} [c] at 13 0
\put {$ \scriptstyle \bullet$} [c] at 13 4
\put {$ \scriptstyle \bullet$} [c] at 13 8
\put {$ \scriptstyle \bullet$} [c] at 13 12
\put {$ \scriptstyle \bullet$} [c] at 16 8
\put {$ \scriptstyle \bullet$} [c] at 16 12
\setlinear \plot 16 12 16 8 13 12 10 8 13 4 16 8    /
\setlinear \plot 13 0 13  12      /
\put{$2{,}520$} [c] at 13 -2
\endpicture
\end{minipage}
\begin{minipage}{4cm}
\beginpicture
\setcoordinatesystem units   <1.5mm,2mm>
\setplotarea x from 0 to 16, y from -2 to 15
\put{360)} [l] at 2 12
\put {$ \scriptstyle \bullet$} [c] at 10 4
\put {$ \scriptstyle \bullet$} [c] at 13 0
\put {$ \scriptstyle \bullet$} [c] at 13 4
\put {$ \scriptstyle \bullet$} [c] at 13 8
\put {$ \scriptstyle \bullet$} [c] at 13 12
\put {$ \scriptstyle \bullet$} [c] at 16 4
\put {$ \scriptstyle \bullet$} [c] at 16 0
\setlinear \plot 16 0 16 4 13 0 10 4 13 8 16 4    /
\setlinear \plot 13 0 13  12      /
\put{$2{,}520$} [c] at 13 -2
\endpicture
\end{minipage}
$$
$$
\begin{minipage}{4cm}
\beginpicture
\setcoordinatesystem units   <1.5mm,2mm>
\setplotarea x from 0 to 16, y from -2 to 15
\put{361)} [l] at 2 12
\put {$ \scriptstyle \bullet$} [c] at 10 4
\put {$ \scriptstyle \bullet$} [c] at 10 8
\put {$ \scriptstyle \bullet$} [c] at 10 12
\put {$ \scriptstyle \bullet$} [c] at 11.5 6
\put {$ \scriptstyle \bullet$} [c] at 13 0
\put {$ \scriptstyle \bullet$} [c] at 16 4
\put {$ \scriptstyle \bullet$} [c] at 16 12
\setlinear \plot 16 4 13 0 10 4 10 12 16 4 16 12 10 4  /
\put{$2{,}520$} [c] at 13 -2
\endpicture
\end{minipage}
\begin{minipage}{4cm}
\beginpicture
\setcoordinatesystem units   <1.5mm,2mm>
\setplotarea x from 0 to 16, y from -2 to 15
\put{362)} [l] at 2 12
\put {$ \scriptstyle \bullet$} [c] at 10 4
\put {$ \scriptstyle \bullet$} [c] at 10 8
\put {$ \scriptstyle \bullet$} [c] at 10 0
\put {$ \scriptstyle \bullet$} [c] at 11.5 6
\put {$ \scriptstyle \bullet$} [c] at 13 12
\put {$ \scriptstyle \bullet$} [c] at 16 8
\put {$ \scriptstyle \bullet$} [c] at 16 0
\setlinear \plot 16 8 13 12 10 8 10 0 16 8 16 0 10 8  /
\put{$2{,}520$} [c] at 13 -2
\endpicture
\end{minipage}
\begin{minipage}{4cm}
\beginpicture
\setcoordinatesystem units   <1.5mm,2mm>
\setplotarea x from 0 to 16, y from -2 to 15
\put{363)} [l] at 2 12
\put {$ \scriptstyle \bullet$} [c] at 10 4
\put {$ \scriptstyle \bullet$} [c] at 10 8
\put {$ \scriptstyle \bullet$} [c] at 10 12
\put {$ \scriptstyle \bullet$} [c] at 11.5 10
\put {$ \scriptstyle \bullet$} [c] at 13 0
\put {$ \scriptstyle \bullet$} [c] at 16 4
\put {$ \scriptstyle \bullet$} [c] at 16 12
\setlinear \plot 16 4 13 0 10 4 10 12 16 4 16 12 10 4  /
\put{$2{,}520$} [c] at 13 -2
\endpicture
\end{minipage}
\begin{minipage}{4cm}
\beginpicture
\setcoordinatesystem units   <1.5mm,2mm>
\setplotarea x from 0 to 16, y from -2 to 15
\put{364)} [l] at 2 12
\put {$ \scriptstyle \bullet$} [c] at 10 4
\put {$ \scriptstyle \bullet$} [c] at 10 8
\put {$ \scriptstyle \bullet$} [c] at 10 0
\put {$ \scriptstyle \bullet$} [c] at 11.5 2
\put {$ \scriptstyle \bullet$} [c] at 13 12
\put {$ \scriptstyle \bullet$} [c] at 16 8
\put {$ \scriptstyle \bullet$} [c] at 16 0
\setlinear \plot 16 8 13 12 10 8 10 0 16 8 16 0 10 8  /
\put{$2{,}520$} [c] at 13 -2
\endpicture
\end{minipage}
\begin{minipage}{4cm}
\beginpicture
\setcoordinatesystem units   <1.5mm,2mm>
\setplotarea x from 0 to 16, y from -2 to 15
\put{365)} [l] at 2 12
\put {$ \scriptstyle \bullet$} [c] at 10 4
\put {$ \scriptstyle \bullet$} [c] at 10 12
\put {$ \scriptstyle \bullet$} [c] at 13 0
\put {$ \scriptstyle \bullet$} [c] at 13 4
\put {$ \scriptstyle \bullet$} [c] at 16 4
\put {$ \scriptstyle \bullet$} [c] at 16 8
\put {$ \scriptstyle \bullet$} [c] at 16 12
\setlinear \plot 13 0   10 4  10 12 13 4 13 0 16 4 10 12 /
\setlinear \plot 13 4 16 8 16 12     /
\setlinear \plot 16 4 16 8  /
\put{$2{,}520$} [c] at 13 -2
\endpicture
\end{minipage}
\begin{minipage}{4cm}
\beginpicture
\setcoordinatesystem units   <1.5mm,2mm>
\setplotarea x from 0 to 16, y from -2 to 15
\put{366)} [l] at 2 12
\put {$ \scriptstyle \bullet$} [c] at 10 8
\put {$ \scriptstyle \bullet$} [c] at 10 0
\put {$ \scriptstyle \bullet$} [c] at 13 12
\put {$ \scriptstyle \bullet$} [c] at 13 8
\put {$ \scriptstyle \bullet$} [c] at 16 4
\put {$ \scriptstyle \bullet$} [c] at 16 8
\put {$ \scriptstyle \bullet$} [c] at 16 0
\setlinear \plot 13 12   10 8  10 0 13 8 13 12 16 8 10 0 /
\setlinear \plot 13 8 16 4 16 0     /
\setlinear \plot 16 4 16 8  /
\put{$2{,}520$} [c] at 13 -2
\endpicture
\end{minipage}
$$
$$
\begin{minipage}{4cm}
\beginpicture
\setcoordinatesystem units   <1.5mm,2mm>
\setplotarea x from 0 to 16, y from -2 to 15
\put{367)} [l] at 2 12
\put {$ \scriptstyle \bullet$} [c] at 10 4
\put {$ \scriptstyle \bullet$} [c] at 10 12
\put {$ \scriptstyle \bullet$} [c] at 13 0
\put {$ \scriptstyle \bullet$} [c] at 13 12
\put {$ \scriptstyle \bullet$} [c] at 16 4
\put {$ \scriptstyle \bullet$} [c] at 16 8
\put {$ \scriptstyle \bullet$} [c] at 16 12
\setlinear \plot 10 12 10 4 13 0 16 4 16 12 10 4 13 12 16 8  /
\put{$2{,}520$} [c] at 13 -2
\endpicture
\end{minipage}
\begin{minipage}{4cm}
\beginpicture
\setcoordinatesystem units   <1.5mm,2mm>
\setplotarea x from 0 to 16, y from -2 to 15
\put{368)} [l] at 2 12
\put {$ \scriptstyle \bullet$} [c] at 10 8
\put {$ \scriptstyle \bullet$} [c] at 10 0
\put {$ \scriptstyle \bullet$} [c] at 13 0
\put {$ \scriptstyle \bullet$} [c] at 13 12
\put {$ \scriptstyle \bullet$} [c] at 16 4
\put {$ \scriptstyle \bullet$} [c] at 16 8
\put {$ \scriptstyle \bullet$} [c] at 16 0
\setlinear \plot 10 0 10 8 13 12 16 8 16 0 10 8 13 0 16 4  /
\put{$2{,}520$} [c] at 13 -2
\endpicture
\end{minipage}
\begin{minipage}{4cm}
\beginpicture
\setcoordinatesystem units   <1.5mm,2mm>
\setplotarea x from 0 to 16, y from -2 to 15
\put{369)} [l] at 2 12
\put {$ \scriptstyle \bullet$} [c] at 10 8
\put {$ \scriptstyle \bullet$} [c] at 12 0
\put {$ \scriptstyle \bullet$} [c] at 12 4
\put {$ \scriptstyle \bullet$} [c] at 12 12
\put {$ \scriptstyle \bullet$} [c] at 14 8
\put {$ \scriptstyle \bullet$} [c] at 14 12
\put {$ \scriptstyle \bullet$} [c] at 16 12
\setlinear \plot 16 12 14 8 12 4 12 0    /
\setlinear \plot 14 12 14 8 12 4 10 8 12 12 14 8    /
\put{$2{,}520$} [c] at 13 -2
\endpicture
\end{minipage}
\begin{minipage}{4cm}
\beginpicture
\setcoordinatesystem units   <1.5mm,2mm>
\setplotarea x from 0 to 16, y from -2 to 15
\put{370)} [l] at 2 12
\put {$ \scriptstyle \bullet$} [c] at 10 4
\put {$ \scriptstyle \bullet$} [c] at 12 12
\put {$ \scriptstyle \bullet$} [c] at 12 8
\put {$ \scriptstyle \bullet$} [c] at 12 0
\put {$ \scriptstyle \bullet$} [c] at 14 4
\put {$ \scriptstyle \bullet$} [c] at 14 0
\put {$ \scriptstyle \bullet$} [c] at 16 0
\setlinear \plot 16 0 14 4 12 8 12 12    /
\setlinear \plot 14 0 14 4 12 8 10 4 12 0 14 4    /
\put{$2{,}520$} [c] at 13 -2
\endpicture
\end{minipage}
\begin{minipage}{4cm}
\beginpicture
\setcoordinatesystem units   <1.5mm,2mm>
\setplotarea x from 0 to 16, y from -2 to 15
\put{371)} [l] at 2 12
\put {$ \scriptstyle \bullet$} [c] at  10 4
\put {$ \scriptstyle \bullet$} [c] at  10 8
\put {$ \scriptstyle \bullet$} [c] at  10 12
\put {$ \scriptstyle \bullet$} [c] at  11 0
\put {$ \scriptstyle \bullet$} [c] at  12 4
\put {$ \scriptstyle \bullet$} [c] at  12 12
\put {$ \scriptstyle \bullet$} [c] at  16 0
\setlinear \plot 10 12 10 4 12  12 12 4 11 0  10 4 /
\setlinear \plot 10 12 16 0 12 12  /
\setlinear \plot 10 8 12 4 /
\put{$2{,}520  $} [c] at 13 -2
\endpicture
\end{minipage}
\begin{minipage}{4cm}
\beginpicture
\setcoordinatesystem units   <1.5mm,2mm>
\setplotarea x from 0 to 16, y from -2 to 15
\put{372)} [l] at 2 12
\put {$ \scriptstyle \bullet$} [c] at  10 0
\put {$ \scriptstyle \bullet$} [c] at  10 4
\put {$ \scriptstyle \bullet$} [c] at  10 8
\put {$ \scriptstyle \bullet$} [c] at  11 12
\put {$ \scriptstyle \bullet$} [c] at  12 0
\put {$ \scriptstyle \bullet$} [c] at  12 8
\put {$ \scriptstyle \bullet$} [c] at  16 12
\setlinear \plot 10 0 10 8 12  0 12 8 11 12  10 8 /
\setlinear \plot 10 0 16 12 12 0  /
\setlinear \plot 10 4 12 8 /
\put{$2{,}520  $} [c] at 13 -2
\endpicture
\end{minipage}
$$
$$
\begin{minipage}{4cm}
\beginpicture
\setcoordinatesystem units   <1.5mm,2mm>
\setplotarea x from 0 to 16, y from -2 to 15
\put{373)} [l] at 2 12
\put {$ \scriptstyle \bullet$} [c] at  10 0
\put {$ \scriptstyle \bullet$} [c] at  12 6
\put {$ \scriptstyle \bullet$} [c] at  12 9
\put {$ \scriptstyle \bullet$} [c] at  12 12
\put {$ \scriptstyle \bullet$} [c] at  14 0
\put {$ \scriptstyle \bullet$} [c] at  16 6
\put {$ \scriptstyle \bullet$} [c] at  16 12
\setlinear \plot 10 0 12 9 12 6 16 12 16 6 14 0 12 6  /
\setlinear \plot 12 12 12 9 16 6 /
\put{$2{,}520$} [c] at 13 -2
\endpicture
\end{minipage}
\begin{minipage}{4cm}
\beginpicture
\setcoordinatesystem units   <1.5mm,2mm>
\setplotarea x from 0 to 16, y from -2 to 15
\put{374)} [l] at 2 12
\put {$ \scriptstyle \bullet$} [c] at  10 12
\put {$ \scriptstyle \bullet$} [c] at  12 0
\put {$ \scriptstyle \bullet$} [c] at  12 3
\put {$ \scriptstyle \bullet$} [c] at  12 6
\put {$ \scriptstyle \bullet$} [c] at  14 12
\put {$ \scriptstyle \bullet$} [c] at  16 6
\put {$ \scriptstyle \bullet$} [c] at  16 0
\setlinear \plot 10 12 12 3 12 6 16 0 16 6 14 12 12 6  /
\setlinear \plot 12 0 12 3 16 6 /
\put{$2{,}520 $} [c] at 13 -2
\endpicture
\end{minipage}
\begin{minipage}{4cm}
\beginpicture
\setcoordinatesystem units   <1.5mm,2mm>
\setplotarea x from 0 to 16, y from -2 to 15
\put{375)} [l] at 2 12
\put {$ \scriptstyle \bullet$} [c] at  10 8
\put {$ \scriptstyle \bullet$} [c] at  10 12
\put {$ \scriptstyle \bullet$} [c] at  11 0
\put {$ \scriptstyle \bullet$} [c] at  11 4
\put {$ \scriptstyle \bullet$} [c] at  12 8
\put {$ \scriptstyle \bullet$} [c] at  12 12
\put {$ \scriptstyle \bullet$} [c] at  16 0
\setlinear \plot 16 0 12 12  12 8  10 12 10  8 12 12    /
\setlinear \plot 10 8 11  4 12 8     /
\setlinear \plot 11 0 11  4   /
\put{$2{,}520   $} [c] at 13 -2
\endpicture
\end{minipage}
\begin{minipage}{4cm}
\beginpicture
\setcoordinatesystem units   <1.5mm,2mm>
\setplotarea x from 0 to 16, y from -2 to 15
\put{376)} [l] at 2 12
\put {$ \scriptstyle \bullet$} [c] at  10 4
\put {$ \scriptstyle \bullet$} [c] at  10 0
\put {$ \scriptstyle \bullet$} [c] at  11 12
\put {$ \scriptstyle \bullet$} [c] at  11 8
\put {$ \scriptstyle \bullet$} [c] at  12 4
\put {$ \scriptstyle \bullet$} [c] at  12 0
\put {$ \scriptstyle \bullet$} [c] at  16 12
\setlinear \plot 16 12 12 0  12 4  10 0 10  4 12 0    /
\setlinear \plot 10 4 11  8 12 4     /
\setlinear \plot 11 12 11  8   /
\put{$2{,}520  $} [c] at 13 -2
\endpicture
\end{minipage}
\begin{minipage}{4cm}
\beginpicture
\setcoordinatesystem units   <1.5mm,2mm>
\setplotarea x from 0 to 16, y from -2 to 15
\put{377)} [l] at 2 12
\put {$ \scriptstyle \bullet$} [c] at  10 4
\put {$ \scriptstyle \bullet$} [c] at  10 12
\put {$ \scriptstyle \bullet$} [c] at  11 0
\put {$ \scriptstyle \bullet$} [c] at  11 8
\put {$ \scriptstyle \bullet$} [c] at  12 4
\put {$ \scriptstyle \bullet$} [c] at  12 12
\put {$ \scriptstyle \bullet$} [c] at  16 0
\setlinear \plot 16 0 12 12 11 8 10 4 11 0 12 4  11 8 10 12  /
\put{$2{,}520  $} [c] at 13 -2
\endpicture
\end{minipage}
\begin{minipage}{4cm}
\beginpicture
\setcoordinatesystem units   <1.5mm,2mm>
\setplotarea x from 0 to 16, y from -2 to 15
\put{378)} [l] at 2 12
\put {$ \scriptstyle \bullet$} [c] at  10 8
\put {$ \scriptstyle \bullet$} [c] at  10 0
\put {$ \scriptstyle \bullet$} [c] at  11 12
\put {$ \scriptstyle \bullet$} [c] at  11 4
\put {$ \scriptstyle \bullet$} [c] at  12 8
\put {$ \scriptstyle \bullet$} [c] at  12 0
\put {$ \scriptstyle \bullet$} [c] at  16 12
\setlinear \plot 16 12 12 0 11 4 10 8 11 12 12 8  11 4 10 0  /
\put{$2{,}520   $} [c] at 13 -2
\endpicture
\end{minipage}
$$
$$
\begin{minipage}{4cm}
\beginpicture
\setcoordinatesystem units   <1.5mm,2mm>
\setplotarea x from 0 to 16, y from -2 to 15
\put{379)} [l] at 2 12
\put {$ \scriptstyle \bullet$} [c] at  10 0
\put {$ \scriptstyle \bullet$} [c] at  10 6
\put {$ \scriptstyle \bullet$} [c] at  10 12
\put {$ \scriptstyle \bullet$} [c] at  13 0
\put {$ \scriptstyle \bullet$} [c] at  13 6
\put {$ \scriptstyle \bullet$} [c] at  13 12
\put {$ \scriptstyle \bullet$} [c] at  16 6
\setlinear \plot 10 12 10 0  13 6 13 12 16 6 13 0 13 6 10 12 /
\setlinear \plot  13 0 10 6 13 12  /
\put{$2{,}520   $} [c] at 13 -2
\endpicture
\end{minipage}
\begin{minipage}{4cm}
\beginpicture
\setcoordinatesystem units   <1.5mm,2mm>
\setplotarea x from 0 to 16, y from -2 to 15
\put{380)} [l] at 2 12
\put {$ \scriptstyle \bullet$} [c] at  10 0
\put {$ \scriptstyle \bullet$} [c] at  10 4
\put {$ \scriptstyle \bullet$} [c] at  10 8
\put {$ \scriptstyle \bullet$} [c] at  10 12
\put {$ \scriptstyle \bullet$} [c] at  12 0
\put {$ \scriptstyle \bullet$} [c] at  12 12
\put {$ \scriptstyle \bullet$} [c] at  16 12
\setlinear \plot 10 0 10 12 12 0 12  12 10 4    /
\setlinear \plot  10 8 16 12 12 0     /
\put{$2{,}520$} [c] at 13 -2
\endpicture
\end{minipage}
\begin{minipage}{4cm}
\beginpicture
\setcoordinatesystem units   <1.5mm,2mm>
\setplotarea x from 0 to 16, y from -2 to 15
\put{381)} [l] at 2 12
\put {$ \scriptstyle \bullet$} [c] at  10 0
\put {$ \scriptstyle \bullet$} [c] at  10 4
\put {$ \scriptstyle \bullet$} [c] at  10 8
\put {$ \scriptstyle \bullet$} [c] at  10 12
\put {$ \scriptstyle \bullet$} [c] at  12 0
\put {$ \scriptstyle \bullet$} [c] at  12 12
\put {$ \scriptstyle \bullet$} [c] at  16 0
\setlinear \plot 10 12 10 0 12 12 12  0 10 8    /
\setlinear \plot  10 4 16 0 12 12     /
\put{$2{,}520$} [c] at 13 -2
\endpicture
\end{minipage}
\begin{minipage}{4cm}
\beginpicture
\setcoordinatesystem units   <1.5mm,2mm>
\setplotarea x from 0 to 16, y from -2 to 15
\put{382)} [l] at 2 12
\put {$ \scriptstyle \bullet$} [c] at  10 0
\put {$ \scriptstyle \bullet$} [c] at  10 12
\put {$ \scriptstyle \bullet$} [c] at  13 6
\put {$ \scriptstyle \bullet$} [c] at  13 12
\put {$ \scriptstyle \bullet$} [c] at  16 12
\put {$ \scriptstyle \bullet$} [c] at  16  6
\put {$ \scriptstyle \bullet$} [c] at  14.5 0
\setlinear \plot  16 12 16 6 14.5 0 13 6 13 12 16 6 /
\setlinear \plot  10 12 13 6 10 0  /
\setlinear \plot  16 12 13 6  /
\put{$2{,}520$} [c] at 13 -2
\endpicture
\end{minipage}
\begin{minipage}{4cm}
\beginpicture
\setcoordinatesystem units   <1.5mm,2mm>
\setplotarea x from 0 to 16, y from -2 to 15
\put{383)} [l] at 2 12
\put {$ \scriptstyle \bullet$} [c] at  10 0
\put {$ \scriptstyle \bullet$} [c] at  10 12
\put {$ \scriptstyle \bullet$} [c] at  13 6
\put {$ \scriptstyle \bullet$} [c] at  13 0
\put {$ \scriptstyle \bullet$} [c] at  16 0
\put {$ \scriptstyle \bullet$} [c] at  16  6
\put {$ \scriptstyle \bullet$} [c] at  14.5 12
\setlinear \plot  16 0 16 6 14.5 12 13 6 13 0 16 6 /
\setlinear \plot  10 12 13 6 10 0  /
\setlinear \plot  16 0 13 6  /
\put{$2{,}520$} [c] at 13 -2
\endpicture
\end{minipage}
\begin{minipage}{4cm}
\beginpicture
\setcoordinatesystem units   <1.5mm,2mm>
\setplotarea x from 0 to 16, y from -2 to 15
\put{384)} [l] at 2 12
\put {$ \scriptstyle \bullet$} [c] at  10 0
\put {$ \scriptstyle \bullet$} [c] at  10 4
\put {$ \scriptstyle \bullet$} [c] at  10 8
\put {$ \scriptstyle \bullet$} [c] at  10 12
\put {$ \scriptstyle \bullet$} [c] at  13 12
\put {$ \scriptstyle \bullet$} [c] at  16 0
\put {$ \scriptstyle \bullet$} [c] at  16 12
\setlinear \plot  10 0 10 12 /
\setlinear \plot  10 8 13 12 16 0 16 12 10 8 /
\put{$2{,}520$} [c] at 13 -2
\endpicture
\end{minipage}
$$
$$
\begin{minipage}{4cm}
\beginpicture
\setcoordinatesystem units   <1.5mm,2mm>
\setplotarea x from 0 to 16, y from -2 to 15
\put{385)} [l] at 2 12
\put {$ \scriptstyle \bullet$} [c] at  10 0
\put {$ \scriptstyle \bullet$} [c] at  10 4
\put {$ \scriptstyle \bullet$} [c] at  10 8
\put {$ \scriptstyle \bullet$} [c] at  10 12
\put {$ \scriptstyle \bullet$} [c] at  13 0
\put {$ \scriptstyle \bullet$} [c] at  16 0
\put {$ \scriptstyle \bullet$} [c] at  16 12
\setlinear \plot  10 0 10 12 /
\setlinear \plot  10 4 13 0 16 12 16 0 10 4 /
\put{$2{,}520$} [c] at 13 -2
\endpicture
\end{minipage}
\begin{minipage}{4cm}
\beginpicture
\setcoordinatesystem units   <1.5mm,2mm>
\setplotarea x from 0 to 16, y from -2 to 15
\put{386)} [l] at 2 12
\put {$ \scriptstyle \bullet$} [c] at  10 12
\put {$ \scriptstyle \bullet$} [c] at  13 0
\put {$ \scriptstyle \bullet$} [c] at  13 6
\put {$ \scriptstyle \bullet$} [c] at  13 12
\put {$ \scriptstyle \bullet$} [c] at  16 0
\put {$ \scriptstyle \bullet$} [c] at  16 6
\put {$ \scriptstyle \bullet$} [c] at  16 12
\setlinear \plot   13 12 13 0  16 6 16 0 13 12 /
\setlinear \plot   10 12 13 6  16 12 16 6 10 12 /
\put{$2{,}520$} [c] at 13 -2
\endpicture
\end{minipage}
\begin{minipage}{4cm}
\beginpicture
\setcoordinatesystem units   <1.5mm,2mm>
\setplotarea x from 0 to 16, y from -2 to 15
\put{387)} [l] at 2 12
\put {$ \scriptstyle \bullet$} [c] at  10 0
\put {$ \scriptstyle \bullet$} [c] at  13 0
\put {$ \scriptstyle \bullet$} [c] at  13 6
\put {$ \scriptstyle \bullet$} [c] at  13 12
\put {$ \scriptstyle \bullet$} [c] at  16 0
\put {$ \scriptstyle \bullet$} [c] at  16 6
\put {$ \scriptstyle \bullet$} [c] at  16 12
\setlinear \plot   13 0 13 12  16 6 16 12 13 0 /
\setlinear \plot   10 0 13 6  16 0 16 6 10 0 /
\put{$2{,}520$} [c] at 13 -2
\endpicture
\end{minipage}
\begin{minipage}{4cm}
\beginpicture
\setcoordinatesystem units   <1.5mm,2mm>
\setplotarea x from 0 to 16, y from -2 to 15
\put{388)} [l] at 2 12
\put {$ \scriptstyle \bullet$} [c] at  10 0
\put {$ \scriptstyle \bullet$} [c] at  10 4
\put {$ \scriptstyle \bullet$} [c] at  10 8
\put {$ \scriptstyle \bullet$} [c] at  10 12
\put {$ \scriptstyle \bullet$} [c] at  13 12
\put {$ \scriptstyle \bullet$} [c] at  16 0
\put {$ \scriptstyle \bullet$} [c] at  16 12
\setlinear \plot 10 0 10 12    /
\setlinear \plot   16  12 16 0 10 4     /
\setlinear \plot  13 12  10  8   /
\put{$2{,}520$} [c] at 13 -2
\endpicture
\end{minipage}
\begin{minipage}{4cm}
\beginpicture
\setcoordinatesystem units   <1.5mm,2mm>
\setplotarea x from 0 to 16, y from -2 to 15
\put{389)} [l] at 2 12
\put {$ \scriptstyle \bullet$} [c] at  10 0
\put {$ \scriptstyle \bullet$} [c] at  10 4
\put {$ \scriptstyle \bullet$} [c] at  10 8
\put {$ \scriptstyle \bullet$} [c] at  10 12
\put {$ \scriptstyle \bullet$} [c] at  13 0
\put {$ \scriptstyle \bullet$} [c] at  16 0
\put {$ \scriptstyle \bullet$} [c] at  16 12
\setlinear \plot 10 0 10 12    /
\setlinear \plot   16  0 16 12 10 8     /
\setlinear \plot  13 0  10  4   /
\put{$2{,}520$} [c] at 13 -2
\endpicture
\end{minipage}
\begin{minipage}{4cm}
\beginpicture
\setcoordinatesystem units   <1.5mm,2mm>
\setplotarea x from 0 to 16, y from -2 to 15
\put{390)} [l] at 2 12
\put {$ \scriptstyle \bullet$} [c] at  10 12
\put {$ \scriptstyle \bullet$} [c] at  13 0
\put {$ \scriptstyle \bullet$} [c] at  13 4
\put {$ \scriptstyle \bullet$} [c] at  13 8
\put {$ \scriptstyle \bullet$} [c] at  13 12
\put {$ \scriptstyle \bullet$} [c] at  16  0
\put {$ \scriptstyle \bullet$} [c] at  16 12
\setlinear \plot  13 0 13 12  /
\setlinear \plot  16 12 13 8 16 0  /
\setlinear \plot  10 12 13 4  /
\put{$2{,}520$} [c] at 13 -2
\endpicture
\end{minipage}
$$
$$
\begin{minipage}{4cm}
\beginpicture
\setcoordinatesystem units   <1.5mm,2mm>
\setplotarea x from 0 to 16, y from -2 to 15
\put{391)} [l] at 2 12
\put {$ \scriptstyle \bullet$} [c] at  10 0
\put {$ \scriptstyle \bullet$} [c] at  13 0
\put {$ \scriptstyle \bullet$} [c] at  13 4
\put {$ \scriptstyle \bullet$} [c] at  13 8
\put {$ \scriptstyle \bullet$} [c] at  13 12
\put {$ \scriptstyle \bullet$} [c] at  16  0
\put {$ \scriptstyle \bullet$} [c] at  16 12
\setlinear \plot  13 0 13 12  /
\setlinear \plot  16 0 13 4 16 12  /
\setlinear \plot  10 0 13 8  /
\put{$2{,}520$} [c] at 13 -2
\endpicture
\end{minipage}
\begin{minipage}{4cm}
\beginpicture
\setcoordinatesystem units   <1.5mm,2mm>
\setplotarea x from 0 to 16, y from -2 to 15
\put{392)} [l] at 2 12
\put {$ \scriptstyle \bullet$} [c] at  16 12
\put {$ \scriptstyle \bullet$} [c] at  16 0
\put {$ \scriptstyle \bullet$} [c] at  13 12
\put {$ \scriptstyle \bullet$} [c] at  10 0
\put {$ \scriptstyle \bullet$} [c] at  10 4
\put {$ \scriptstyle \bullet$} [c] at  10 8
\put {$ \scriptstyle \bullet$} [c] at  10 12
\setlinear \plot 10 12 10 0   16 12 16 0 10 8 13 12  /
\put{$2{,}520$} [c] at 13 -2
\endpicture
\end{minipage}
\begin{minipage}{4cm}
\beginpicture
\setcoordinatesystem units   <1.5mm,2mm>
\setplotarea x from 0 to 16, y from -2 to 15
\put{393)} [l] at 2 12
\put {$ \scriptstyle \bullet$} [c] at  16 12
\put {$ \scriptstyle \bullet$} [c] at  16 0
\put {$ \scriptstyle \bullet$} [c] at  13 0
\put {$ \scriptstyle \bullet$} [c] at  10 0
\put {$ \scriptstyle \bullet$} [c] at  10 4
\put {$ \scriptstyle \bullet$} [c] at  10 8
\put {$ \scriptstyle \bullet$} [c] at  10 12
\setlinear \plot 10 0 10 12   16 0 16 12 10 4 13 0  /
\put{$2{,}520$} [c] at 13 -2
\endpicture
\end{minipage}
\begin{minipage}{4cm}
\beginpicture
\setcoordinatesystem units   <1.5mm,2mm>
\setplotarea x from 0 to 16, y from -2 to 15
\put{394)} [l] at 2 12
\put {$ \scriptstyle \bullet$} [c] at 10 4
\put {$ \scriptstyle \bullet$} [c] at 10 8
\put {$ \scriptstyle \bullet$} [c] at 11.5 12
 \put {$ \scriptstyle \bullet$} [c] at 11.5 0
\put {$ \scriptstyle \bullet$} [c] at 13 4
\put {$ \scriptstyle \bullet$} [c] at 13 8
\put {$ \scriptstyle \bullet$} [c] at 16 6
\setlinear \plot 11.5 0 10 4 10 8 11.5 12  16 6 11.5 0 13 4 13 8  11.5 12  /
\setlinear \plot 10 4 13 8     /
\setlinear \plot 10 8 13 4     /
\put{$1{,}260$} [c] at 13 -2
\endpicture
\end{minipage}
\begin{minipage}{4cm}
\beginpicture
\setcoordinatesystem units   <1.5mm,2mm>
\setplotarea x from 0 to 16, y from -2 to 15
\put{395)} [l] at 2 12
\put {$ \scriptstyle \bullet$} [c] at 10 4
\put {$ \scriptstyle \bullet$} [c] at 10 12
\put {$ \scriptstyle \bullet$} [c] at 13 0
\put {$ \scriptstyle \bullet$} [c] at 13 4
\put {$ \scriptstyle \bullet$} [c] at 13 12
\put {$ \scriptstyle \bullet$} [c] at 16 4
\put {$ \scriptstyle \bullet$} [c] at 16 8
\setlinear \plot 13 0 10 4 10 12 13 4 13 0 16 4 16 8 13 12 10 4    /
\setlinear \plot 10 12  16 8  /
\setlinear \plot 13 12  13 4  /
\put{$1{,}260$} [c] at 13 -2
\endpicture
\end{minipage}
\begin{minipage}{4cm}
\beginpicture
\setcoordinatesystem units   <1.5mm,2mm>
\setplotarea x from 0 to 16, y from -2 to 15
\put{396)} [l] at 2 12
\put {$ \scriptstyle \bullet$} [c] at 10 8
\put {$ \scriptstyle \bullet$} [c] at 10 0
\put {$ \scriptstyle \bullet$} [c] at 13 12
\put {$ \scriptstyle \bullet$} [c] at 13 8
\put {$ \scriptstyle \bullet$} [c] at 13 0
\put {$ \scriptstyle \bullet$} [c] at 16 8
\put {$ \scriptstyle \bullet$} [c] at 16 4
\setlinear \plot 13 12 10 8 10 0 13 8 13 12 16 8 16 4 13 0 10 8    /
\setlinear \plot 10 0  16 4  /
\setlinear \plot 13 0  13 8  /
\put{$1{,}260$} [c] at 13 -2
\endpicture
\end{minipage}
$$
$$
\begin{minipage}{4cm}
\beginpicture
\setcoordinatesystem units   <1.5mm,2mm>
\setplotarea x from 0 to 16, y from -2 to 15
\put{397)} [l] at 2 12
\put {$ \scriptstyle \bullet$} [c] at 10 8
\put {$ \scriptstyle \bullet$} [c] at 10 12
\put {$ \scriptstyle \bullet$} [c] at 11 0
\put {$ \scriptstyle \bullet$} [c] at 12 8
\put {$ \scriptstyle \bullet$} [c] at 12 12
\put {$ \scriptstyle \bullet$} [c] at 11 4
\put {$ \scriptstyle \bullet$} [c] at 16 12
\setlinear \plot  12 8 11 4 10 8 10 12  12 8 12 12 10 8  /
\setlinear \plot 11  0 11 4 16 12  /
\put{$1{,}260$} [c] at 13 -2
\endpicture
\end{minipage}
\begin{minipage}{4cm}
\beginpicture
\setcoordinatesystem units   <1.5mm,2mm>
\setplotarea x from 0 to 16, y from -2 to 15
\put{398)} [l] at 2 12
\put {$ \scriptstyle \bullet$} [c] at 10 4
\put {$ \scriptstyle \bullet$} [c] at 10 0
\put {$ \scriptstyle \bullet$} [c] at 11 12
\put {$ \scriptstyle \bullet$} [c] at 12 4
\put {$ \scriptstyle \bullet$} [c] at 12 0
\put {$ \scriptstyle \bullet$} [c] at 11 8
\put {$ \scriptstyle \bullet$} [c] at 16 0
\setlinear \plot  12 4 11 8 10 4 10 0  12 4 12 0 10 4  /
\setlinear \plot 11  12 11 8 16 0  /
\put{$1{,}260$} [c] at 13 -2
\endpicture
\end{minipage}
\begin{minipage}{4cm}
\beginpicture
\setcoordinatesystem units   <1.5mm,2mm>
\setplotarea x from 0 to 16, y from -2 to 15
\put{399)} [l] at 2 12
\put {$ \scriptstyle \bullet$} [c] at 10 4
\put {$ \scriptstyle \bullet$} [c] at 10 12
\put {$ \scriptstyle \bullet$} [c] at 12 0
\put {$ \scriptstyle \bullet$} [c] at 14 4
\put {$ \scriptstyle \bullet$} [c] at 14 12
\put {$ \scriptstyle \bullet$} [c] at 16 8
\put {$ \scriptstyle \bullet$} [c] at 16 12
\setlinear \plot 10 12 10 4  12 0 14 4 16 8 16 12  /
\setlinear \plot 10 12 14 4 14 12 10 4 16 8 /
\put{$1{,}260$} [c] at 13 -2
\endpicture
\end{minipage}
\begin{minipage}{4cm}
\beginpicture
\setcoordinatesystem units   <1.5mm,2mm>
\setplotarea x from 0 to 16, y from -2 to 15
\put{400)} [l] at 2 12
\put {$ \scriptstyle \bullet$} [c] at 10 8
\put {$ \scriptstyle \bullet$} [c] at 10 0
\put {$ \scriptstyle \bullet$} [c] at 12 12
\put {$ \scriptstyle \bullet$} [c] at 14 8
\put {$ \scriptstyle \bullet$} [c] at 14 0
\put {$ \scriptstyle \bullet$} [c] at 16 0
\put {$ \scriptstyle \bullet$} [c] at 16 4
\setlinear \plot 10 0 10 8  12 12 14 8 16 4 16 0  /
\setlinear \plot 10 0 14 8 14 0 10 8 16 4 /
\put{$1{,}260$} [c] at 13 -2
\endpicture
\end{minipage}
\begin{minipage}{4cm}
\beginpicture
\setcoordinatesystem units   <1.5mm,2mm>
\setplotarea x from 0 to 16, y from -2 to 15
\put{401)} [l] at 2 12
\put {$ \scriptstyle \bullet$} [c] at  10 0
\put {$ \scriptstyle \bullet$} [c] at  10 6
\put {$ \scriptstyle \bullet$} [c] at  10 12
\put {$ \scriptstyle \bullet$} [c] at  13 3
\put {$ \scriptstyle \bullet$} [c] at  16 0
\put {$ \scriptstyle \bullet$} [c] at  16 6
\put {$ \scriptstyle \bullet$} [c] at  16 12
\setlinear \plot 10 12  10 0 13 3 16 12  16 0 13 3 10 12 16 6  /
\setlinear \plot  10 6 16 12  /
\put{$1{,}260$} [c] at 13 -2
\endpicture
\end{minipage}
\begin{minipage}{4cm}
\beginpicture
\setcoordinatesystem units   <1.5mm,2mm>
\setplotarea x from 0 to 16, y from -2 to 15
\put{402)} [l] at 2 12
\put {$ \scriptstyle \bullet$} [c] at  10 0
\put {$ \scriptstyle \bullet$} [c] at  10 6
\put {$ \scriptstyle \bullet$} [c] at  10 12
\put {$ \scriptstyle \bullet$} [c] at  13 9
\put {$ \scriptstyle \bullet$} [c] at  16 0
\put {$ \scriptstyle \bullet$} [c] at  16 6
\put {$ \scriptstyle \bullet$} [c] at  16 12
\setlinear \plot 10 0  10 12 13 9 16 0  16 12 13 9 10 0 16 6  /
\setlinear \plot  10 6 16 0  /
\put{$1{,}260$} [c] at 13 -2
\endpicture
\end{minipage}
$$
$$
\begin{minipage}{4cm}
\beginpicture
\setcoordinatesystem units   <1.5mm,2mm>
\setplotarea x from 0 to 16, y from -2 to 15
\put{403)} [l] at 1 12
\put {$ \scriptstyle \bullet$} [c] at  10 0
\put {$ \scriptstyle \bullet$} [c] at  10 6
\put {$ \scriptstyle \bullet$} [c] at  10 12
\put {$ \scriptstyle \bullet$} [c] at  13 12
\put {$ \scriptstyle \bullet$} [c] at  16 0
\put {$ \scriptstyle \bullet$} [c] at  16 6
\put {$ \scriptstyle \bullet$} [c] at  16 12
\setlinear \plot  10 12 10 0 16 6 16 0 10 6 13 12 16 6 16 12  /
\put{$1{,}260$} [c] at 13 -2
\endpicture
\end{minipage}
\begin{minipage}{4cm}
\beginpicture
\setcoordinatesystem units   <1.5mm,2mm>
\setplotarea x from 0 to 16, y from -2 to 15
\put{404)} [l] at 1 12
\put {$ \scriptstyle \bullet$} [c] at  10 0
\put {$ \scriptstyle \bullet$} [c] at  10 6
\put {$ \scriptstyle \bullet$} [c] at  10 12
\put {$ \scriptstyle \bullet$} [c] at  13 0
\put {$ \scriptstyle \bullet$} [c] at  16 0
\put {$ \scriptstyle \bullet$} [c] at  16 6
\put {$ \scriptstyle \bullet$} [c] at  16 12
\setlinear \plot  10 0 10 12 16 6 16 12 10 6 13 0 16 6 16 0  /
\put{$1{,}260$} [c] at 13 -2
\endpicture
\end{minipage}
\begin{minipage}{4cm}
\beginpicture
\setcoordinatesystem units   <1.5mm,2mm>
\setplotarea x from 0 to 16, y from -2 to 15
\put{405)} [l] at 2 12
\put {$ \scriptstyle \bullet$} [c] at  10 0
\put {$ \scriptstyle \bullet$} [c] at  10 4
\put {$ \scriptstyle \bullet$} [c] at  10 8
\put {$ \scriptstyle \bullet$} [c] at  10 12
\put {$ \scriptstyle \bullet$} [c] at  13 12
\put {$ \scriptstyle \bullet$} [c] at  16 0
\put {$ \scriptstyle \bullet$} [c] at  16 12
\setlinear \plot 10 0 10 12     /
\setlinear \plot   16  0  10 4 16 12    /
\setlinear \plot  13 12  10  4   /
\put{$1{,}260$} [c] at 13 -2
\endpicture
\end{minipage}
\begin{minipage}{4cm}
\beginpicture
\setcoordinatesystem units   <1.5mm,2mm>
\setplotarea x from 0 to 16, y from -2 to 15
\put{406)} [l] at 2 12
\put {$ \scriptstyle \bullet$} [c] at  10 0
\put {$ \scriptstyle \bullet$} [c] at  10 4
\put {$ \scriptstyle \bullet$} [c] at  10 8
\put {$ \scriptstyle \bullet$} [c] at  10 12
\put {$ \scriptstyle \bullet$} [c] at  13 0
\put {$ \scriptstyle \bullet$} [c] at  16 0
\put {$ \scriptstyle \bullet$} [c] at  16 12
\setlinear \plot 10 0 10 12     /
\setlinear \plot   16  0  10 8 16 12    /
\setlinear \plot  13 0  10  8   /
\put{$1{,}260$} [c] at 13 -2
\endpicture
\end{minipage}
\begin{minipage}{4cm}
\beginpicture
\setcoordinatesystem units   <1.5mm,2mm>
\setplotarea x from 0 to 16, y from -2 to 15
\put{407)} [l] at 2 12
\put {$ \scriptstyle \bullet$} [c] at  10 0
\put {$ \scriptstyle \bullet$} [c] at  10 6
\put {$ \scriptstyle \bullet$} [c] at  10 12
\put {$ \scriptstyle \bullet$} [c] at  13 12
\put {$ \scriptstyle \bullet$} [c] at  16 0
\put {$ \scriptstyle \bullet$} [c] at  16 6
\put {$ \scriptstyle \bullet$} [c] at  16 12
\setlinear \plot  16 0 16 12 10 6 10  0 /
\setlinear \plot  10 12 10 6 13 12 16 6 10 12 /
\put{$420$} [c] at 13 -2
\endpicture
\end{minipage}
\begin{minipage}{4cm}
\beginpicture
\setcoordinatesystem units   <1.5mm,2mm>
\setplotarea x from 0 to 16, y from -2 to 15
\put{408)} [l] at 2 12
\put {$ \scriptstyle \bullet$} [c] at  10 0
\put {$ \scriptstyle \bullet$} [c] at  10 6
\put {$ \scriptstyle \bullet$} [c] at  10 12
\put {$ \scriptstyle \bullet$} [c] at  13 0
\put {$ \scriptstyle \bullet$} [c] at  16 0
\put {$ \scriptstyle \bullet$} [c] at  16 6
\put {$ \scriptstyle \bullet$} [c] at  16 12
\setlinear \plot  16 12 16 0 10 6 10  12 /
\setlinear \plot  10 0 10 6 13 0 16 6 10 0 /
\put{$420$} [c] at 13 -2
\endpicture
\end{minipage}
$$
$$
\begin{minipage}{4cm}
\beginpicture
\setcoordinatesystem units   <1.5mm,2mm>
\setplotarea x from 0 to 16, y from -2 to 15
\put{409)} [l] at 2 12
\put {$ \scriptstyle \bullet$} [c] at 10 4
\put {$ \scriptstyle \bullet$} [c] at 10 12
\put {$ \scriptstyle \bullet$} [c] at 13 12
\put {$ \scriptstyle \bullet$} [c] at 13 4
\put {$ \scriptstyle \bullet$} [c] at 13 0
\put {$ \scriptstyle \bullet$} [c] at 16 4
\put {$ \scriptstyle \bullet$} [c] at 16 12
\setlinear \plot 13 0 10 4 10 12 13 4 13 0 16 4 16 12 10 4 13 12 13 4 16 12 10 4  /
\setlinear \plot 10 12   16 4 13 12     /
\put{$140$} [c] at 13 -2
\endpicture
\end{minipage}
\begin{minipage}{4cm}
\beginpicture
\setcoordinatesystem units   <1.5mm,2mm>
\setplotarea x from 0 to 16, y from -2 to 15
\put{410)} [l] at 2 12
\put {$ \scriptstyle \bullet$} [c] at 10 0
\put {$ \scriptstyle \bullet$} [c] at 10 8
\put {$ \scriptstyle \bullet$} [c] at 13 0
\put {$ \scriptstyle \bullet$} [c] at 13 8
\put {$ \scriptstyle \bullet$} [c] at 13 12
\put {$ \scriptstyle \bullet$} [c] at 16 8
\put {$ \scriptstyle \bullet$} [c] at 16 0
\setlinear \plot 13 12 10 8 10 0 13 8 13 12 16 8 16 0 10 8 13 0 13 8 16 0 10 8  /
\setlinear \plot 10 0  16 8 13 0     /
\put{$140$} [c] at 13 -2
\endpicture
\end{minipage}
\begin{minipage}{4cm}
\beginpicture
\setcoordinatesystem units   <1.5mm,2mm>
\setplotarea x from 0 to 16, y from -2 to 15
\put{411)} [l] at 2 12
\put {$ \scriptstyle \bullet$} [c] at 10 0
\put {$ \scriptstyle \bullet$} [c] at 13 0
\put {$ \scriptstyle \bullet$} [c] at 16 0
\put {$ \scriptstyle \bullet$} [c] at 13 6
\put {$ \scriptstyle \bullet$} [c] at 10 12
\put {$ \scriptstyle \bullet$} [c] at 13 12
\put {$ \scriptstyle \bullet$} [c] at 16 12
\setlinear \plot 10 0 13 6 16 12    /
\setlinear \plot 16 0 13 6 10 12    /
\setlinear \plot 13 0 13 12     /
\put{$140$} [c] at 13 -2
\endpicture
\end{minipage}
\begin{minipage}{4cm}
\beginpicture
\setcoordinatesystem units   <1.5mm,2mm>
\setplotarea x from 0 to 16, y from -2 to 15
\put{412)} [l] at 2 12
\put {$ \scriptstyle \bullet$} [c] at 12 12
\put {$ \scriptstyle \bullet$} [c] at 12 6
\put {$ \scriptstyle \bullet$} [c] at 12 3
\put {$ \scriptstyle \bullet$} [c] at 12 0
\put {$ \scriptstyle \bullet$} [c] at 10 9
\put {$ \scriptstyle \bullet$} [c] at 14 9
\put {$ \scriptstyle \bullet$} [c] at 16  0
\setlinear \plot  12 0 12 6  10 9 12 12 14 9 12  6     /
\put{$2{,}520$} [c] at 11 -2
\put{$\scriptstyle \bullet$} [c] at 16 0
 \endpicture
\end{minipage}
\begin{minipage}{4cm}
\beginpicture
\setcoordinatesystem units   <1.5mm,2mm>
\setplotarea x from 0 to 16, y from -2 to 15
\put{413)} [l] at 2 12
\put {$ \scriptstyle \bullet$} [c] at 12 12
\put {$ \scriptstyle \bullet$} [c] at 12 9
\put {$ \scriptstyle \bullet$} [c] at 12 6
\put {$ \scriptstyle \bullet$} [c] at 12 0
\put {$ \scriptstyle \bullet$} [c] at 10 3
\put {$ \scriptstyle \bullet$} [c] at 14 3
\put {$ \scriptstyle \bullet$} [c] at 16  0
\setlinear \plot  12 12 12 6  10 3 12 0 14 3 12  6     /
\put{$2{,}520$} [c] at 11 -2
\put{$\scriptstyle \bullet$} [c] at 16  0
\endpicture
\end{minipage}
\begin{minipage}{4cm}
\beginpicture
\setcoordinatesystem units   <1.5mm,2mm>
\setplotarea x from 0 to 16, y from -2 to 15
\put{414)} [l] at 2 12
\put {$ \scriptstyle \bullet$} [c] at 12 0
\put {$ \scriptstyle \bullet$} [c] at 12 2
\put {$ \scriptstyle \bullet$} [c] at 12 4
\put {$ \scriptstyle \bullet$} [c] at 12 6
\put {$ \scriptstyle \bullet$} [c] at 10 12
\put {$ \scriptstyle \bullet$} [c] at 14 12
\setlinear \plot  12 0 12 6  10 12      /
\setlinear \plot  12 6  14 12      /
\put{$2{,}520$} [c] at 11 -2
\put{$\scriptstyle \bullet$} [c] at 16  0
\endpicture
\end{minipage}
$$
$$
\begin{minipage}{4cm}
\beginpicture
\setcoordinatesystem units   <1.5mm,2mm>
\setplotarea x from 0 to 16, y from -2 to 15
\put{415)} [l] at 2 12
\put {$ \scriptstyle \bullet$} [c] at 12 12
\put {$ \scriptstyle \bullet$} [c] at 12 10
\put {$ \scriptstyle \bullet$} [c] at 12 8
\put {$ \scriptstyle \bullet$} [c] at 12 6
\put {$ \scriptstyle \bullet$} [c] at 10 0
\put {$ \scriptstyle \bullet$} [c] at 14 0
\setlinear \plot  12 12 12 6  10 0      /
\setlinear \plot  12 6  14 0      /
\put{$2{,}520$} [c] at 11 -2
\put{$\scriptstyle \bullet$} [c] at 16  0
\endpicture
\end{minipage}
\begin{minipage}{4cm}
\beginpicture
\setcoordinatesystem units   <1.5mm,2mm>
\setplotarea x from 0 to 16, y from -2 to 15
\put{416)} [l] at 2 12
\put {$ \scriptstyle \bullet$} [c] at 12 12
\put {$ \scriptstyle \bullet$} [c] at 12 10
\put {$ \scriptstyle \bullet$} [c] at 12 2
\put {$ \scriptstyle \bullet$} [c] at 12 0
\put {$ \scriptstyle \bullet$} [c] at 10  6
\put {$ \scriptstyle \bullet$} [c] at 14 6
\setlinear \plot  12 0 12 2 10 6 12  10  12 12     /
\setlinear \plot  12 10   14 6 12 2      /
\put{$2{,}520$} [c] at 11 -2
\put{$\scriptstyle \bullet$} [c] at 16  0
\endpicture
\end{minipage}
\begin{minipage}{4cm}
\beginpicture
\setcoordinatesystem units   <1.5mm,2mm>
\setplotarea x from 0 to 16, y from -2 to 15
\put{${\bf  17}$} [l] at 2 15
\put{417)} [l] at 2 12
\put {$ \scriptstyle \bullet$} [c] at  10 8
\put {$ \scriptstyle \bullet$} [c] at 12 8
\put {$ \scriptstyle \bullet$} [c] at 14 0
\put {$ \scriptstyle \bullet$} [c] at 14 4
\put {$ \scriptstyle \bullet$} [c] at 14 12
\put {$ \scriptstyle \bullet$} [c] at 16 8
\put {$ \scriptstyle \bullet$} [c] at 16 12
\setlinear \plot 16 12 16 8 14 4 14 0 10 8 14 12 12 8 14 4 16 8 14 12    /
\put{$5{,}040$} [c] at 13 -2
\endpicture
\end{minipage}
\begin{minipage}{4cm}
\beginpicture
\setcoordinatesystem units   <1.5mm,2mm>
\setplotarea x from 0 to 16, y from -2 to 15
\put{418)} [l] at 2 12
\put {$ \scriptstyle \bullet$} [c] at  10 4
\put {$ \scriptstyle \bullet$} [c] at 12 4
\put {$ \scriptstyle \bullet$} [c] at 14 12
\put {$ \scriptstyle \bullet$} [c] at 14 8
\put {$ \scriptstyle \bullet$} [c] at 14 0
\put {$ \scriptstyle \bullet$} [c] at 16 4
\put {$ \scriptstyle \bullet$} [c] at 16 0
\setlinear \plot 16 0 16 4 14 8 14 12 10 4 14 0 12 4 14 8 16 4 14 0    /
\put{$5{,}040$} [c] at 13 -2
\endpicture
\end{minipage}
\begin{minipage}{4cm}
\beginpicture
\setcoordinatesystem units   <1.5mm,2mm>
\setplotarea x from 0 to 16, y from -2 to 15
\put{419)} [l] at 2 12
\put {$ \scriptstyle \bullet$} [c] at 10 12
\put {$ \scriptstyle \bullet$} [c] at 12 0
\put {$ \scriptstyle \bullet$} [c] at 12 4
\put {$ \scriptstyle \bullet$} [c] at 12 8
\put {$ \scriptstyle \bullet$} [c] at 16 4
\put {$ \scriptstyle \bullet$} [c] at 16 8
\put {$ \scriptstyle \bullet$} [c] at 16 12
\setlinear \plot 10 12 12 4 12 0 16 4 16 12 12 8 12 4  16 8   /
\put{$5{,}040$} [c] at 13 -2
\endpicture
\end{minipage}
\begin{minipage}{4cm}
\beginpicture
\setcoordinatesystem units   <1.5mm,2mm>
\setplotarea x from 0 to 16, y from -2 to 15
\put{420)} [l] at 2 12
\put {$ \scriptstyle \bullet$} [c] at 10 0
\put {$ \scriptstyle \bullet$} [c] at 12 12
\put {$ \scriptstyle \bullet$} [c] at 12 4
\put {$ \scriptstyle \bullet$} [c] at 12 8
\put {$ \scriptstyle \bullet$} [c] at 16 4
\put {$ \scriptstyle \bullet$} [c] at 16 8
\put {$ \scriptstyle \bullet$} [c] at 16 0
\setlinear \plot 10 0 12 8 12 12 16 8 16 0 12 4 12 8  16 4   /
\put{$5{,}040$} [c] at 13 -2
\endpicture
\end{minipage}
$$

$$
\begin{minipage}{4cm}
\beginpicture
\setcoordinatesystem units   <1.5mm,2mm>
\setplotarea x from 0 to 16, y from -2 to 15
\put{421)} [l] at 2 12
\put {$ \scriptstyle \bullet$} [c] at 10 4
\put {$ \scriptstyle \bullet$} [c] at 10 12
\put {$ \scriptstyle \bullet$} [c] at 13 0
\put {$ \scriptstyle \bullet$} [c] at 13 4
\put {$ \scriptstyle \bullet$} [c] at 13 8
\put {$ \scriptstyle \bullet$} [c] at 13 12
\put {$ \scriptstyle \bullet$} [c] at 16 4
\setlinear \plot 13 8 10 12 10 4 13 0 13 12 16 4 13 0   /
\setlinear \plot 10 4 13 12   /
\put{$5{,}040$} [c] at 13 -2
\endpicture
\end{minipage}
\begin{minipage}{4cm}
\beginpicture
\setcoordinatesystem units   <1.5mm,2mm>
\setplotarea x from 0 to 16, y from -2 to 15
\put{422)} [l] at 2 12
\put {$ \scriptstyle \bullet$} [c] at 10 8
\put {$ \scriptstyle \bullet$} [c] at 10 0
\put {$ \scriptstyle \bullet$} [c] at 13 0
\put {$ \scriptstyle \bullet$} [c] at 13 4
\put {$ \scriptstyle \bullet$} [c] at 13 8
\put {$ \scriptstyle \bullet$} [c] at 13 12
\put {$ \scriptstyle \bullet$} [c] at 16 8
\setlinear \plot 13 4 10 0 10 8 13 12 13 0 16 8 13 12   /
\setlinear \plot 10 8 13 0   /
\put{$5{,}040$} [c] at 13 -2
\endpicture
\end{minipage}
\begin{minipage}{4cm}
\beginpicture
\setcoordinatesystem units   <1.5mm,2mm>
\setplotarea x from 0 to 16, y from -2 to 15
\put{423)} [l] at 2 12
\put {$ \scriptstyle \bullet$} [c] at 10 4
\put {$ \scriptstyle \bullet$} [c] at 10 12
\put {$ \scriptstyle \bullet$} [c] at 13 0
\put {$ \scriptstyle \bullet$} [c] at 13 4
\put {$ \scriptstyle \bullet$} [c] at 16 4
\put {$ \scriptstyle \bullet$} [c] at 16 8
\put {$ \scriptstyle \bullet$} [c] at 16 12
\setlinear \plot 13 4 13 0 10 4 10 12  13 4 16 8 16 12   /
\setlinear \plot 13 0 16 4 16 8  /
\put{$5{,}040$} [c] at 13 -2
\endpicture
\end{minipage}
\begin{minipage}{4cm}
\beginpicture
\setcoordinatesystem units   <1.5mm,2mm>
\setplotarea x from 0 to 16, y from -2 to 15
\put{424)} [l] at 2 12
\put {$ \scriptstyle \bullet$} [c] at 10 0
\put {$ \scriptstyle \bullet$} [c] at 10 8
\put {$ \scriptstyle \bullet$} [c] at 13 12
\put {$ \scriptstyle \bullet$} [c] at 13 8
\put {$ \scriptstyle \bullet$} [c] at 16 0
\put {$ \scriptstyle \bullet$} [c] at 16 4
\put {$ \scriptstyle \bullet$} [c] at 16 8
\setlinear \plot 13 8 13 12 10 8 10 0  13 8 16 4 16 0   /
\setlinear \plot 13 12 16 8 16 4  /
\put{$5{,}040$} [c] at 13 -2
\endpicture
\end{minipage}
\begin{minipage}{4cm}
\beginpicture
\setcoordinatesystem units   <1.5mm,2mm>
\setplotarea x from 0 to 16, y from -2 to 15
\put{425)} [l] at 2 12
\put {$ \scriptstyle \bullet$} [c] at 10 4
\put {$ \scriptstyle \bullet$} [c] at 10 12
\put {$ \scriptstyle \bullet$} [c] at 13 0
\put {$ \scriptstyle \bullet$} [c] at 16 4
\put {$ \scriptstyle \bullet$} [c] at 16  12
\put {$ \scriptstyle \bullet$} [c] at 11.5 10
\put {$ \scriptstyle \bullet$} [c] at 14.5 6
\setlinear \plot 16 12 16 4 13 0 10 4 10 12 16 4   /
\put{$5{,}040$} [c] at 13 -2
\endpicture
\end{minipage}
\begin{minipage}{4cm}
\beginpicture
\setcoordinatesystem units   <1.5mm,2mm>
\setplotarea x from 0 to 16, y from -2 to 15
\put{426)} [l] at 2 12
\put {$ \scriptstyle \bullet$} [c] at 10 0
\put {$ \scriptstyle \bullet$} [c] at 10 8
\put {$ \scriptstyle \bullet$} [c] at 13 12
\put {$ \scriptstyle \bullet$} [c] at 16 0
\put {$ \scriptstyle \bullet$} [c] at 16  8
\put {$ \scriptstyle \bullet$} [c] at 11.5 2
\put {$ \scriptstyle \bullet$} [c] at 14.5 6
\setlinear \plot 16 0 16 8 13 12 10 8 10 0 16 8   /
\put{$5{,}040$} [c] at 13 -2
\endpicture
\end{minipage}
$$
$$
\begin{minipage}{4cm}
\beginpicture
\setcoordinatesystem units   <1.5mm,2mm>
\setplotarea x from 0 to 16, y from -2 to 15
\put{427)} [l] at 2 12
\put {$ \scriptstyle \bullet$} [c] at 10 5
\put {$ \scriptstyle \bullet$} [c] at  10.5 7.5
\put {$ \scriptstyle \bullet$} [c] at 11 0
\put {$ \scriptstyle \bullet$} [c] at 11 10
\put {$ \scriptstyle \bullet$} [c] at 11 12
\put {$ \scriptstyle \bullet$} [c] at 12 5
\put {$ \scriptstyle \bullet$} [c] at 16 12
\setlinear \plot  16 12 11 0 10 5 11 10 11 12  /
\setlinear \plot  11 0 12 5 11 10  /
\put{$5{,}040$} [c] at 13 -2
\endpicture
\end{minipage}
\begin{minipage}{4cm}
\beginpicture
\setcoordinatesystem units   <1.5mm,2mm>
\setplotarea x from 0 to 16, y from -2 to 15
\put{428)} [l] at 2 12
\put {$ \scriptstyle \bullet$} [c] at 10 7
\put {$ \scriptstyle \bullet$} [c] at  10.5 4.5
\put {$ \scriptstyle \bullet$} [c] at 11 0
\put {$ \scriptstyle \bullet$} [c] at 11 2
\put {$ \scriptstyle \bullet$} [c] at 11 12
\put {$ \scriptstyle \bullet$} [c] at 12 7
\put {$ \scriptstyle \bullet$} [c] at 16 0
\setlinear \plot  16 0 11 12 10 7 11 2 11 0  /
\setlinear \plot  11 12 12 7 11 2  /
\put{$5{,}040$} [c] at 13 -2
\endpicture
\end{minipage}
\begin{minipage}{4cm}
\beginpicture
\setcoordinatesystem units   <1.5mm,2mm>
\setplotarea x from 0 to 16, y from -2 to 15
\put{429)} [l] at 2 12
\put {$ \scriptstyle \bullet$} [c] at 10 6
\put {$ \scriptstyle \bullet$} [c] at 10 10
\put {$ \scriptstyle \bullet$} [c] at 11 0
\put {$ \scriptstyle \bullet$} [c] at 11 4
\put {$ \scriptstyle \bullet$} [c] at 11 12
\put {$ \scriptstyle \bullet$} [c] at 12 8
\put {$ \scriptstyle \bullet$} [c] at 16 12
\setlinear \plot 16 12 11 0 11 4 10 6 10 10 11 12 12 8 11 4    /
\put{$5{,}040$} [c] at 13 -2
\endpicture
\end{minipage}
\begin{minipage}{4cm}
\beginpicture
\setcoordinatesystem units   <1.5mm,2mm>
\setplotarea x from 0 to 16, y from -2 to 15
\put{430)} [l] at 2 12
\put {$ \scriptstyle \bullet$} [c] at 10 2
\put {$ \scriptstyle \bullet$} [c] at 10 6
\put {$ \scriptstyle \bullet$} [c] at 11 12
\put {$ \scriptstyle \bullet$} [c] at 11 8
\put {$ \scriptstyle \bullet$} [c] at 11 0
\put {$ \scriptstyle \bullet$} [c] at 12 4
\put {$ \scriptstyle \bullet$} [c] at 16 0
\setlinear \plot 16 0 11 12 11 8 10 6 10 2 11 0 12 4 11 8    /
\put{$5{,}040$} [c] at 13 -2
\endpicture
\end{minipage}
\begin{minipage}{4cm}
\beginpicture
\setcoordinatesystem units   <1.5mm,2mm>
\setplotarea x from 0 to 16, y from -2 to 15
\put{431)} [l] at 2 12
\put {$ \scriptstyle \bullet$} [c] at 10 4
\put {$ \scriptstyle \bullet$} [c] at 10 12
\put {$ \scriptstyle \bullet$} [c] at 12 0
\put {$ \scriptstyle \bullet$} [c] at 14 4
\put {$ \scriptstyle \bullet$} [c] at 14 8
\put {$ \scriptstyle \bullet$} [c] at 14 12
\put {$ \scriptstyle \bullet$} [c] at 16 12
\setlinear \plot 16 12 14 4 12 0 10 4 10 12 14 4 14 12   /
\setlinear \plot 10 4 14 8   /
\put{$5{,}040$} [c] at 13 -2
\endpicture
\end{minipage}
\begin{minipage}{4cm}
\beginpicture
\setcoordinatesystem units   <1.5mm,2mm>
\setplotarea x from 0 to 16, y from -2 to 15
\put{432)} [l] at 2 12
\put {$ \scriptstyle \bullet$} [c] at 10 8
\put {$ \scriptstyle \bullet$} [c] at 10 0
\put {$ \scriptstyle \bullet$} [c] at 12 12
\put {$ \scriptstyle \bullet$} [c] at 14 4
\put {$ \scriptstyle \bullet$} [c] at 14 8
\put {$ \scriptstyle \bullet$} [c] at 14 0
\put {$ \scriptstyle \bullet$} [c] at 16 0
\setlinear \plot 16 0 14 8 12 12 10 8 10 0 14 8 14 0   /
\setlinear \plot 10 8 14 4   /
\put{$5{,}040$} [c] at 13 -2
\endpicture
\end{minipage}
$$

$$
\begin{minipage}{4cm}
\beginpicture
\setcoordinatesystem units   <1.5mm,2mm>
\setplotarea x from 0 to 16, y from -2 to 15
\put{433)} [l] at 2 12
\put {$ \scriptstyle \bullet$} [c] at 10 4
\put {$ \scriptstyle \bullet$} [c] at 10 12
\put {$ \scriptstyle \bullet$} [c] at 13 0
\put {$ \scriptstyle \bullet$} [c] at 13 12
\put {$ \scriptstyle \bullet$} [c] at 16 4
\put {$ \scriptstyle \bullet$} [c] at 16 8
\put {$ \scriptstyle \bullet$} [c] at 16 12
\setlinear \plot 10 12 10 4 13 0 16 4 16 12   /
\setlinear \plot 10 4 13 12 16 8  /
\put{$5{,}040$} [c] at 13 -2
\endpicture
\end{minipage}
\begin{minipage}{4cm}
\beginpicture
\setcoordinatesystem units   <1.5mm,2mm>
\setplotarea x from 0 to 16, y from -2 to 15
\put{434)} [l] at 2 12
\put {$ \scriptstyle \bullet$} [c] at 10 0
\put {$ \scriptstyle \bullet$} [c] at 10 8
\put {$ \scriptstyle \bullet$} [c] at 13 0
\put {$ \scriptstyle \bullet$} [c] at 13 12
\put {$ \scriptstyle \bullet$} [c] at 16 4
\put {$ \scriptstyle \bullet$} [c] at 16 8
\put {$ \scriptstyle \bullet$} [c] at 16 0
\setlinear \plot 10 0 10 8 13 12 16 8 16 0   /
\setlinear \plot 10 8 13 0 16 4  /
\put{$5{,}040$} [c] at 13 -2
\endpicture
\end{minipage}
\begin{minipage}{4cm}
\beginpicture
\setcoordinatesystem units   <1.5mm,2mm>
\setplotarea x from 0 to 16, y from -2 to 15
\put{435)} [l] at 2 12
\put {$ \scriptstyle \bullet$} [c] at 10 12
\put {$ \scriptstyle \bullet$} [c] at  10.5 8
\put {$ \scriptstyle \bullet$} [c] at 11 0
\put {$ \scriptstyle \bullet$} [c] at 11 2
\put {$ \scriptstyle \bullet$} [c] at 11 4
\put {$ \scriptstyle \bullet$} [c] at 12 12
\put {$ \scriptstyle \bullet$} [c] at 16 12
\setlinear \plot  16 12 11 0 11 4 10 12  /
\setlinear \plot  12 12 11 4  /
\put{$5{,}040$} [c] at 13 -2
\endpicture
\end{minipage}
\begin{minipage}{4cm}
\beginpicture
\setcoordinatesystem units   <1.5mm,2mm>
\setplotarea x from 0 to 16, y from -2 to 15
\put{436)} [l] at 2 12
\put {$ \scriptstyle \bullet$} [c] at 10 0
\put {$ \scriptstyle \bullet$} [c] at  10.5 4
\put {$ \scriptstyle \bullet$} [c] at 11 12
\put {$ \scriptstyle \bullet$} [c] at 11 10
\put {$ \scriptstyle \bullet$} [c] at 11 8
\put {$ \scriptstyle \bullet$} [c] at 12 0
\put {$ \scriptstyle \bullet$} [c] at 16 0
\setlinear \plot  16 0 11 12 11 8 10 0  /
\setlinear \plot  12 0 11 8  /
\put{$5{,}040$} [c] at 13 -2
\endpicture
\end{minipage}
\begin{minipage}{4cm}
\beginpicture
\setcoordinatesystem units   <1.5mm,2mm>
\setplotarea x from 0 to 16, y from -2 to 15
\put{437)} [l] at 2 12
\put {$ \scriptstyle \bullet$} [c] at  10 0
\put {$ \scriptstyle \bullet$} [c] at  10 4
\put {$ \scriptstyle \bullet$} [c] at  10 8
\put {$ \scriptstyle \bullet$} [c] at  10 12
\put {$ \scriptstyle \bullet$} [c] at  16 0
\put {$ \scriptstyle \bullet$} [c] at  16 4
\put {$ \scriptstyle \bullet$} [c] at  16 12
\setlinear \plot 16 12 16 0  10 12 10 0 16 4     /
\put{$5{,}040   $} [c] at 13 -2
\endpicture
\end{minipage}
\begin{minipage}{4cm}
\beginpicture
\setcoordinatesystem units   <1.5mm,2mm>
\setplotarea x from 0 to 16, y from -2 to 15
\put{438)} [l] at 2 12
\put {$ \scriptstyle \bullet$} [c] at  10 0
\put {$ \scriptstyle \bullet$} [c] at  10 4
\put {$ \scriptstyle \bullet$} [c] at  10 8
\put {$ \scriptstyle \bullet$} [c] at  10 12
\put {$ \scriptstyle \bullet$} [c] at  16 0
\put {$ \scriptstyle \bullet$} [c] at  16 8
\put {$ \scriptstyle \bullet$} [c] at  16 12
\setlinear \plot 16 0 16 12  10 0 10 12 16 8     /
\put{$5{,}040   $} [c] at 13 -2
\endpicture
\end{minipage}
$$
$$
\begin{minipage}{4cm}
\beginpicture
\setcoordinatesystem units   <1.5mm,2mm>
\setplotarea x from 0 to 16, y from -2 to 15
\put{439)} [l] at 2 12
\put {$ \scriptstyle \bullet$} [c] at  10 0
\put {$ \scriptstyle \bullet$} [c] at  12 4
\put {$ \scriptstyle \bullet$} [c] at  12 8
\put {$ \scriptstyle \bullet$} [c] at  12 12
\put {$ \scriptstyle \bullet$} [c] at  14 0
\put {$ \scriptstyle \bullet$} [c] at  16 4
\put {$ \scriptstyle \bullet$} [c] at  16 12
\setlinear \plot 10  0 12 8 12 4  14 0 16 4 16 12 12 4      /
\setlinear \plot 12 8 12 12 16 4  /
\put{$5{,}040   $} [c] at 13 -2
\endpicture
\end{minipage}
\begin{minipage}{4cm}
\beginpicture
\setcoordinatesystem units   <1.5mm,2mm>
\setplotarea x from 0 to 16, y from -2 to 15
\put{440)} [l] at 2 12
\put {$ \scriptstyle \bullet$} [c] at  10 12
\put {$ \scriptstyle \bullet$} [c] at  12 0
\put {$ \scriptstyle \bullet$} [c] at  12 4
\put {$ \scriptstyle \bullet$} [c] at  12 8
\put {$ \scriptstyle \bullet$} [c] at  14 12
\put {$ \scriptstyle \bullet$} [c] at  16 8
\put {$ \scriptstyle \bullet$} [c] at  16 0
\setlinear \plot 10  12 12 4 12 8  14 12 16 8 16 0 12 8      /
\setlinear \plot 12 4 12 0 16 8  /
\put{$5{,}040   $} [c] at 13 -2
\endpicture
\end{minipage}
\begin{minipage}{4cm}
\beginpicture
\setcoordinatesystem units   <1.5mm,2mm>
\setplotarea x from 0 to 16, y from -2 to 15
\put{441)} [l] at 2 12
\put {$ \scriptstyle \bullet$} [c] at  10 0
\put {$ \scriptstyle \bullet$} [c] at  10  3
\put {$ \scriptstyle \bullet$} [c] at  10 6
\put {$ \scriptstyle \bullet$} [c] at  10 9
\put {$ \scriptstyle \bullet$} [c] at  10 12
\put {$ \scriptstyle \bullet$} [c] at  16 0
\put {$ \scriptstyle \bullet$} [c] at  16 12
\setlinear \plot  10 12 10 0 16 12 16 0  /
\put{$5{,}040   $} [c] at 13 -2
\endpicture
\end{minipage}
\begin{minipage}{4cm}
\beginpicture
\setcoordinatesystem units   <1.5mm,2mm>
\setplotarea x from 0 to 16, y from -2 to 15
\put{442)} [l] at 2 12
\put {$ \scriptstyle \bullet$} [c] at  10 0
\put {$ \scriptstyle \bullet$} [c] at  10  3
\put {$ \scriptstyle \bullet$} [c] at  10 6
\put {$ \scriptstyle \bullet$} [c] at  10 9
\put {$ \scriptstyle \bullet$} [c] at  10 12
\put {$ \scriptstyle \bullet$} [c] at  16 0
\put {$ \scriptstyle \bullet$} [c] at  16 12
\setlinear \plot  10 0 10 12 16 0 16 12  /
\put{$5{,}040   $} [c] at 13 -2
\endpicture
\end{minipage}
\begin{minipage}{4cm}
\beginpicture
\setcoordinatesystem units   <1.5mm,2mm>
\setplotarea x from 0 to 16, y from -2 to 15
\put{443)} [l] at 2 12
\put {$ \scriptstyle \bullet$} [c] at  10 6
\put {$ \scriptstyle \bullet$} [c] at  10 9
\put {$ \scriptstyle \bullet$} [c] at  10 12
\put {$ \scriptstyle \bullet$} [c] at  13 0
\put {$ \scriptstyle \bullet$} [c] at  13 12
\put {$ \scriptstyle \bullet$} [c] at  16 6
\put {$ \scriptstyle \bullet$} [c] at  16 0
\setlinear \plot 10 12 10 6 13 12 16 6 16 0   /
\setlinear \plot 10 6 13 0 16 6 /
\put{$5{,}040   $} [c] at 13 -2
\endpicture
\end{minipage}
\begin{minipage}{4cm}
\beginpicture
\setcoordinatesystem units   <1.5mm,2mm>
\setplotarea x from 0 to 16, y from -2 to 15
\put{444)} [l] at 2 12
\put {$ \scriptstyle \bullet$} [c] at  10 6
\put {$ \scriptstyle \bullet$} [c] at  10 3
\put {$ \scriptstyle \bullet$} [c] at  10 0
\put {$ \scriptstyle \bullet$} [c] at  13 0
\put {$ \scriptstyle \bullet$} [c] at  13 12
\put {$ \scriptstyle \bullet$} [c] at  16 6
\put {$ \scriptstyle \bullet$} [c] at  16 12
\setlinear \plot 10 0 10 6 13 0 16 6 16 12   /
\setlinear \plot 10 6 13 12 16 6 /
\put{$5{,}040   $} [c] at 13 -2
\endpicture
\end{minipage}
$$

$$
\begin{minipage}{4cm}
\beginpicture
\setcoordinatesystem units   <1.5mm,2mm>
\setplotarea x from 0 to 16, y from -2 to 15
\put{445)} [l] at 2 12
\put {$ \scriptstyle \bullet$} [c] at  10 0
\put {$ \scriptstyle \bullet$} [c] at  10 12
\put {$ \scriptstyle \bullet$} [c] at  12 8
\put {$ \scriptstyle \bullet$} [c] at  14 0
\put {$ \scriptstyle \bullet$} [c] at  14 4
\put {$ \scriptstyle \bullet$} [c] at  14 12
\put {$ \scriptstyle \bullet$} [c] at  16 8
\setlinear \plot 10 0 10 12 12 8 14 4 14 0    /
\setlinear \plot 12 8 14 12  16 8 14 4 /
\put{$5{,}040   $} [c] at 13 -2
\endpicture
\end{minipage}
\begin{minipage}{4cm}
\beginpicture
\setcoordinatesystem units   <1.5mm,2mm>
\setplotarea x from 0 to 16, y from -2 to 15
\put{446)} [l] at 2 12
\put {$ \scriptstyle \bullet$} [c] at  10 0
\put {$ \scriptstyle \bullet$} [c] at  10 12
\put {$ \scriptstyle \bullet$} [c] at  12 4
\put {$ \scriptstyle \bullet$} [c] at  14 0
\put {$ \scriptstyle \bullet$} [c] at  14 8
\put {$ \scriptstyle \bullet$} [c] at  14 12
\put {$ \scriptstyle \bullet$} [c] at  16 4
\setlinear \plot 10 12 10 0 12 4 14 8 14 12    /
\setlinear \plot 12 4 14 0  16 4 14 8 /
\put{$5{,}040   $} [c] at 13 -2
\endpicture
\end{minipage}
\begin{minipage}{4cm}
\beginpicture
\setcoordinatesystem units   <1.5mm,2mm>
\setplotarea x from 0 to 16, y from -2 to 15
\put{447)} [l] at 2 12
\put {$ \scriptstyle \bullet$} [c] at  10 0
\put {$ \scriptstyle \bullet$} [c] at  10 4
\put {$ \scriptstyle \bullet$} [c] at  10 8
\put {$ \scriptstyle \bullet$} [c] at  10 12
\put {$ \scriptstyle \bullet$} [c] at  13 0
\put {$ \scriptstyle \bullet$} [c] at  16 4
\put {$ \scriptstyle \bullet$} [c] at  16 12
\setlinear \plot 10 12  10 4 13 0 16 4 16 12     /
\setlinear \plot 10 0 10 4   /
\put{$5{,}040   $} [c] at 13 -2
\endpicture
\end{minipage}
\begin{minipage}{4cm}
\beginpicture
\setcoordinatesystem units   <1.5mm,2mm>
\setplotarea x from 0 to 16, y from -2 to 15
\put{448)} [l] at 2 12
\put {$ \scriptstyle \bullet$} [c] at  10 0
\put {$ \scriptstyle \bullet$} [c] at  10 4
\put {$ \scriptstyle \bullet$} [c] at  10 8
\put {$ \scriptstyle \bullet$} [c] at  10 12
\put {$ \scriptstyle \bullet$} [c] at  13 12
\put {$ \scriptstyle \bullet$} [c] at  16 8
\put {$ \scriptstyle \bullet$} [c] at  16 0
\setlinear \plot 10 0  10 8 13 12 16 8 16 0     /
\setlinear \plot 10 12 10 8   /
\put{$5{,}040   $} [c] at 13 -2
\endpicture
\end{minipage}
\begin{minipage}{4cm}
\beginpicture
\setcoordinatesystem units   <1.5mm,2mm>
\setplotarea x from 0 to 16, y from -2 to 15
\put{449)} [l] at 2 12
\put {$ \scriptstyle \bullet$} [c] at  10 3
\put {$ \scriptstyle \bullet$} [c] at  12 0
\put {$ \scriptstyle \bullet$} [c] at  12 6
\put {$ \scriptstyle \bullet$} [c] at  12 12
\put {$ \scriptstyle \bullet$} [c] at  14 3
\put {$ \scriptstyle \bullet$} [c] at  16 0
\put {$ \scriptstyle \bullet$} [c] at  16 12
\setlinear \plot 16 0 12 12 12 6 10 3 12 0 14 3 12 6     /
\setlinear \plot 16 0 16 12 14 3    /
\put{$5{,}040   $} [c] at 13 -2
\endpicture
\end{minipage}
\begin{minipage}{4cm}
\beginpicture
\setcoordinatesystem units   <1.5mm,2mm>
\setplotarea x from 0 to 16, y from -2 to 15
\put{450)} [l] at 2 12
\put {$ \scriptstyle \bullet$} [c] at  10 9
\put {$ \scriptstyle \bullet$} [c] at  12 12
\put {$ \scriptstyle \bullet$} [c] at  12 6
\put {$ \scriptstyle \bullet$} [c] at  12 0
\put {$ \scriptstyle \bullet$} [c] at  14 9
\put {$ \scriptstyle \bullet$} [c] at  16 0
\put {$ \scriptstyle \bullet$} [c] at  16 12
\setlinear \plot 16 12 12 0 12 6 10 9 12 12 14 9 12 6     /
\setlinear \plot 16 12 16 0 14 9    /
\put{$5{,}040   $} [c] at 13 -2
\endpicture
\end{minipage}
$$
$$
\begin{minipage}{4cm}
\beginpicture
\setcoordinatesystem units   <1.5mm,2mm>
\setplotarea x from 0 to 16, y from -2 to 15
\put{451)} [l] at 2 12
\put {$ \scriptstyle \bullet$} [c] at  10 3
\put {$ \scriptstyle \bullet$} [c] at  10 6
\put {$ \scriptstyle \bullet$} [c] at  10 12
\put {$ \scriptstyle \bullet$} [c] at  12 0
\put {$ \scriptstyle \bullet$} [c] at  12 12
\put {$ \scriptstyle \bullet$} [c] at  14 3
\put {$ \scriptstyle \bullet$} [c] at  16 0
\setlinear \plot 10 12 10 3 12 0 14 3 12 12 10 3  /
\setlinear \plot 10 6 16 0 12 12 /

\put{$5{,}040$} [c] at 13 -2
\endpicture
\end{minipage}
\begin{minipage}{4cm}
\beginpicture
\setcoordinatesystem units   <1.5mm,2mm>
\setplotarea x from 0 to 16, y from -2 to 15
\put{452)} [l] at 2 12
\put {$ \scriptstyle \bullet$} [c] at  10 9
\put {$ \scriptstyle \bullet$} [c] at  10 6
\put {$ \scriptstyle \bullet$} [c] at  10 0
\put {$ \scriptstyle \bullet$} [c] at  12 0
\put {$ \scriptstyle \bullet$} [c] at  12 12
\put {$ \scriptstyle \bullet$} [c] at  14 9
\put {$ \scriptstyle \bullet$} [c] at  16 12
\setlinear \plot 10 0 10 9 12 12 14 9 12 0 10 9  /
\setlinear \plot 10 6 16 12 12 0 /
\put{$5{,}040  $} [c] at 13 -2
\endpicture
\end{minipage}
\begin{minipage}{4cm}
\beginpicture
\setcoordinatesystem units   <1.5mm,2mm>
\setplotarea x from 0 to 16, y from -2 to 15
\put{453)} [l] at 2 12
\put {$ \scriptstyle \bullet$} [c] at  10 6
\put {$ \scriptstyle \bullet$} [c] at  10 12
\put {$ \scriptstyle \bullet$} [c] at  10.8 7.5
\put {$ \scriptstyle \bullet$} [c] at  13 0
\put {$ \scriptstyle \bullet$} [c] at  13 12
\put {$ \scriptstyle \bullet$} [c] at  16 6
\put {$ \scriptstyle \bullet$} [c] at  16 0
\setlinear \plot 10 12 10 6  13 12 16 6 13 0 10 6   /
\setlinear \plot 10 12 16 0 16 6 /
\put{$5{,}040  $} [c] at 13 -2
\endpicture
\end{minipage}
\begin{minipage}{4cm}
\beginpicture
\setcoordinatesystem units   <1.5mm,2mm>
\setplotarea x from 0 to 16, y from -2 to 15
\put{454)} [l] at 2 12
\put {$ \scriptstyle \bullet$} [c] at  10 6
\put {$ \scriptstyle \bullet$} [c] at  10 0
\put {$ \scriptstyle \bullet$} [c] at  10.8 4.5
\put {$ \scriptstyle \bullet$} [c] at  13 0
\put {$ \scriptstyle \bullet$} [c] at  13 12
\put {$ \scriptstyle \bullet$} [c] at  16 6
\put {$ \scriptstyle \bullet$} [c] at  16 12
\setlinear \plot 10 0 10 6  13 0 16 6 13 12 10 6   /
\setlinear \plot 10 0 16 12 16 6 /
\put{$5{,}040  $} [c] at 13 -2
\endpicture
\end{minipage}
\begin{minipage}{4cm}
\beginpicture
\setcoordinatesystem units   <1.5mm,2mm>
\setplotarea x from 0 to 16, y from -2 to 15
\put{455)} [l] at 2 12
\put {$ \scriptstyle \bullet$} [c] at  10 0
\put {$ \scriptstyle \bullet$} [c] at  10 4
\put {$ \scriptstyle \bullet$} [c] at  10 8
\put {$ \scriptstyle \bullet$} [c] at  10 12
\put {$ \scriptstyle \bullet$} [c] at  16 0
\put {$ \scriptstyle \bullet$} [c] at  16 8
\put {$ \scriptstyle \bullet$} [c] at  16 12
\setlinear \plot  10 0  10 12 16 0 16 12 10 4 /
\put{$5{,}040   $} [c] at 13 -2
\endpicture
\end{minipage}
\begin{minipage}{4cm}
\beginpicture
\setcoordinatesystem units   <1.5mm,2mm>
\setplotarea x from 0 to 16, y from -2 to 15
\put{456)} [l] at 2 12
\put {$ \scriptstyle \bullet$} [c] at  10 0
\put {$ \scriptstyle \bullet$} [c] at  10 4
\put {$ \scriptstyle \bullet$} [c] at  10 8
\put {$ \scriptstyle \bullet$} [c] at  10 12
\put {$ \scriptstyle \bullet$} [c] at  16 0
\put {$ \scriptstyle \bullet$} [c] at  16 8
\put {$ \scriptstyle \bullet$} [c] at  16 12
\setlinear \plot  10 12  10 0 16 12 16 0 10 8 /
\put{$5{,}040   $} [c] at 13 -2
\endpicture
\end{minipage}
$$

$$
\begin{minipage}{4cm}
\beginpicture
\setcoordinatesystem units   <1.5mm,2mm>
\setplotarea x from 0 to 16, y from -2 to 15
\put{457)} [l] at 2 12
\put {$ \scriptstyle \bullet$} [c] at  10 4
\put {$ \scriptstyle \bullet$} [c] at  10 12
\put {$ \scriptstyle \bullet$} [c] at  13 0
\put {$ \scriptstyle \bullet$} [c] at  13 8
\put {$ \scriptstyle \bullet$} [c] at  16 0
\put {$ \scriptstyle \bullet$} [c] at  16 4
\put {$ \scriptstyle \bullet$} [c] at  16 12
\setlinear \plot 16 4 16 12  13 8 16 0 16 4 13 0 10 4 10 12 13 8 13 0   /
\put{$5{,}040   $} [c] at 13 -2
\endpicture
\end{minipage}
\begin{minipage}{4cm}
\beginpicture
\setcoordinatesystem units   <1.5mm,2mm>
\setplotarea x from 0 to 16, y from -2 to 15
\put{458)} [l] at 2 12
\put {$ \scriptstyle \bullet$} [c] at  10 8
\put {$ \scriptstyle \bullet$} [c] at  10 0
\put {$ \scriptstyle \bullet$} [c] at  13 12
\put {$ \scriptstyle \bullet$} [c] at  13 4
\put {$ \scriptstyle \bullet$} [c] at  16 0
\put {$ \scriptstyle \bullet$} [c] at  16 8
\put {$ \scriptstyle \bullet$} [c] at  16 12
\setlinear \plot 16 8 16 0  13 4 16 12 16 8 13 12 10 8 10 0 13 4 13 12   /
\put{$5{,}040  $} [c] at 13 -2
\endpicture
\end{minipage}
\begin{minipage}{4cm}
\beginpicture
\setcoordinatesystem units   <1.5mm,2mm>
\setplotarea x from 0 to 16, y from -2 to 15
\put{459)} [l] at 1 12
\put {$ \scriptstyle \bullet$} [c] at  10 0
\put {$ \scriptstyle \bullet$} [c] at  10 4
\put {$ \scriptstyle \bullet$} [c] at  10 8
\put {$ \scriptstyle \bullet$} [c] at  10 12
\put {$ \scriptstyle \bullet$} [c] at  13 12
\put {$ \scriptstyle \bullet$} [c] at  16 0
\put {$ \scriptstyle \bullet$} [c] at  16 12
\setlinear \plot 10 0 10 4 16 12 16 0 13 12 10 8 10 4 /
\setlinear \plot  10 12 10 8   /
\put{$5{,}040$} [c] at 13 -2
\endpicture
\end{minipage}
\begin{minipage}{4cm}
\beginpicture
\setcoordinatesystem units   <1.5mm,2mm>
\setplotarea x from 0 to 16, y from -2 to 15
\put{460)} [l] at 2 12
\put {$ \scriptstyle \bullet$} [c] at  10 0
\put {$ \scriptstyle \bullet$} [c] at  10 4
\put {$ \scriptstyle \bullet$} [c] at  10 8
\put {$ \scriptstyle \bullet$} [c] at  10 12
\put {$ \scriptstyle \bullet$} [c] at  13 0
\put {$ \scriptstyle \bullet$} [c] at  16 0
\put {$ \scriptstyle \bullet$} [c] at  16 12
\setlinear \plot 10 12 10 8 16 0 16 12 13 0 10 4 10 8 /
\setlinear \plot  10 0 10 4  /
\put{$5{,}040$} [c] at 13 -2
\endpicture
\end{minipage}
\begin{minipage}{4cm}
\beginpicture
\setcoordinatesystem units   <1.5mm,2mm>
\setplotarea x from 0 to 16, y from -2 to 15
\put{461)} [l] at 2 12
\put {$ \scriptstyle \bullet$} [c] at  10 0
\put {$ \scriptstyle \bullet$} [c] at  10 6
\put {$ \scriptstyle \bullet$} [c] at  10 12
\put {$ \scriptstyle \bullet$} [c] at  13 0
\put {$ \scriptstyle \bullet$} [c] at  13 6
\put {$ \scriptstyle \bullet$} [c] at  13 12
\put {$ \scriptstyle \bullet$} [c] at  16 12
\setlinear \plot 16 12 13 6 13 0 10 12 10 0 13 6 13 12 10 6 /
\put{$5{,}040$} [c] at 13 -2
\endpicture
\end{minipage}
\begin{minipage}{4cm}
\beginpicture
\setcoordinatesystem units   <1.5mm,2mm>
\setplotarea x from 0 to 16, y from -2 to 15
\put{462)} [l] at 2 12
\put {$ \scriptstyle \bullet$} [c] at  10 0
\put {$ \scriptstyle \bullet$} [c] at  10 6
\put {$ \scriptstyle \bullet$} [c] at  10 12
\put {$ \scriptstyle \bullet$} [c] at  13 0
\put {$ \scriptstyle \bullet$} [c] at  13 6
\put {$ \scriptstyle \bullet$} [c] at  13 12
\put {$ \scriptstyle \bullet$} [c] at  16 0
\setlinear \plot 16 0 13 6 13 12 10 0 10 12 13 6 13 0 10 6 /
\put{$5{,}040$} [c] at 13 -2
\endpicture
\end{minipage}
$$
$$
\begin{minipage}{4cm}
\beginpicture
\setcoordinatesystem units   <1.5mm,2mm>
\setplotarea x from 0 to 16, y from -2 to 15
\put{463)} [l] at 2 12
\put {$ \scriptstyle \bullet$} [c] at 10 4
\put {$ \scriptstyle \bullet$} [c] at 10 8
\put {$ \scriptstyle \bullet$} [c] at 13 0
\put {$ \scriptstyle \bullet$} [c] at 13 4
\put {$ \scriptstyle \bullet$} [c] at 13 12
\put {$ \scriptstyle \bullet$} [c] at 16 4
\put {$ \scriptstyle \bullet$} [c] at 16 8
\setlinear \plot 13 0 10 4 10 8 13 12 16 8  13 4 13 0 16 4 16 8     /
\put{$2{,}520$} [c] at 13 -2
\endpicture
\end{minipage}
\begin{minipage}{4cm}
\beginpicture
\setcoordinatesystem units   <1.5mm,2mm>
\setplotarea x from 0 to 16, y from -2 to 15
\put{464)} [l] at 2 12
\put {$ \scriptstyle \bullet$} [c] at 10 4
\put {$ \scriptstyle \bullet$} [c] at 10 8
\put {$ \scriptstyle \bullet$} [c] at 13 0
\put {$ \scriptstyle \bullet$} [c] at 13 8
\put {$ \scriptstyle \bullet$} [c] at 13 12
\put {$ \scriptstyle \bullet$} [c] at 16 4
\put {$ \scriptstyle \bullet$} [c] at 16 8
\setlinear \plot 13 12 10 8 10 4 13 0 16 4  13 8 13 12 16 8 16 4     /
\put{$2{,}520$} [c] at 13 -2
\endpicture
\end{minipage}
\begin{minipage}{4cm}
\beginpicture
\setcoordinatesystem units   <1.5mm,2mm>
\setplotarea x from 0 to 16, y from -2 to 15
\put{465)} [l] at 2 12
\put {$ \scriptstyle \bullet$} [c] at 10 8
\put {$ \scriptstyle \bullet$} [c] at 12 0
\put {$ \scriptstyle \bullet$} [c] at  12 4

\put {$ \scriptstyle \bullet$} [c] at 12 12
\put {$ \scriptstyle \bullet$} [c] at 14 8
\put {$ \scriptstyle \bullet$} [c] at 16 8
\put {$ \scriptstyle \bullet$} [c] at 16 12
\setlinear \plot 12 0 12 4 10 8 12 12 14 8 12 4 16 8 16 12    /
\put{$2{,}520$} [c] at 13 -2
\endpicture
\end{minipage}
\begin{minipage}{4cm}
\beginpicture
\setcoordinatesystem units   <1.5mm,2mm>
\setplotarea x from 0 to 16, y from -2 to 15
\put{466)} [l] at 2 12
\put {$ \scriptstyle \bullet$} [c] at 10 4
\put {$ \scriptstyle \bullet$} [c] at 12 0
\put {$ \scriptstyle \bullet$} [c] at  12 8
\put {$ \scriptstyle \bullet$} [c] at 12 12
\put {$ \scriptstyle \bullet$} [c] at 14 4
\put {$ \scriptstyle \bullet$} [c] at 16 4
\put {$ \scriptstyle \bullet$} [c] at 16 0
\setlinear \plot 12 12 12 8 10 4 12 0 14 4 12 8 16 4 16 0    /
\put{$2{,}520$} [c] at 13 -2
\endpicture
\end{minipage}
\begin{minipage}{4cm}
\beginpicture
\setcoordinatesystem units   <1.5mm,2mm>
\setplotarea x from 0 to 16, y from -2 to 15
\put{467)} [l] at 2 12
\put {$ \scriptstyle \bullet$} [c] at 10 4
\put {$ \scriptstyle \bullet$} [c] at 10 12
\put {$ \scriptstyle \bullet$} [c] at 13 0
\put {$ \scriptstyle \bullet$} [c] at 13 4
\put {$ \scriptstyle \bullet$} [c] at 16 4
\put {$ \scriptstyle \bullet$} [c] at 16 8
\put {$ \scriptstyle \bullet$} [c] at 16 12
\setlinear \plot 13 0 10 4 10 12 13 4 16 12 16 4  13 0 13 4  /
\setlinear \plot 10 12 16 4    /
\setlinear \plot 16 12 10 4    /
\put{$2{,}520$} [c] at 13 -2
\endpicture
\end{minipage}
\begin{minipage}{4cm}
\beginpicture
\setcoordinatesystem units   <1.5mm,2mm>
\setplotarea x from 0 to 16, y from -2 to 15
\put{468)} [l] at 2 12
\put {$ \scriptstyle \bullet$} [c] at 10 0
\put {$ \scriptstyle \bullet$} [c] at 10 8
\put {$ \scriptstyle \bullet$} [c] at 13 12
\put {$ \scriptstyle \bullet$} [c] at 13 8
\put {$ \scriptstyle \bullet$} [c] at 16 0
\put {$ \scriptstyle \bullet$} [c] at 16 4
\put {$ \scriptstyle \bullet$} [c] at 16 8
\setlinear \plot 13 12 10 8 10 0 13 8 16 0 16 8  13 12  13 8 /
\setlinear \plot 10 0 16 8    /
\setlinear \plot 16 0 10 8    /
\put{$2{,}520$} [c] at 13 -2
\endpicture
\end{minipage}
$$

$$
\begin{minipage}{4cm}
\beginpicture
\setcoordinatesystem units   <1.5mm,2mm>
\setplotarea x from 0 to 16, y from -2 to 15
\put{469)} [l] at 2 12
\put {$ \scriptstyle \bullet$} [c] at 10 4
\put {$ \scriptstyle \bullet$} [c] at 10 12
\put {$ \scriptstyle \bullet$} [c] at 13 0
\put {$ \scriptstyle \bullet$} [c] at 13 4
\put {$ \scriptstyle \bullet$} [c] at 13 8
\put {$ \scriptstyle \bullet$} [c] at 13 12
\put {$ \scriptstyle \bullet$} [c] at 16 4
\setlinear \plot 10 4 13 8 16 4 13 0 13 12   /
\setlinear \plot 10 12 10 4 13 0 /
\put{$2{,}520$} [c] at 13 -2
\endpicture
\end{minipage}
\begin{minipage}{4cm}
\beginpicture
\setcoordinatesystem units   <1.5mm,2mm>
\setplotarea x from 0 to 16, y from -2 to 15
\put{470)} [l] at 2 12
\put {$ \scriptstyle \bullet$} [c] at 10 8
\put {$ \scriptstyle \bullet$} [c] at 10 0
\put {$ \scriptstyle \bullet$} [c] at 13 0
\put {$ \scriptstyle \bullet$} [c] at 13 4
\put {$ \scriptstyle \bullet$} [c] at 13 8
\put {$ \scriptstyle \bullet$} [c] at 13 12
\put {$ \scriptstyle \bullet$} [c] at 16 8
\setlinear \plot 10 8 13 4 16 8 13 12 13 0   /
\setlinear \plot 10 0 10 8 13 12 /
\put{$2{,}520$} [c] at 13 -2
\endpicture
\end{minipage}
\begin{minipage}{4cm}
\beginpicture
\setcoordinatesystem units   <1.5mm,2mm>
\setplotarea x from 0 to 16, y from -2 to 15
\put{471)} [l] at 2 12
\put {$ \scriptstyle \bullet$} [c] at 10 4
\put {$ \scriptstyle \bullet$} [c] at 10 12
\put {$ \scriptstyle \bullet$} [c] at 13 0
\put {$ \scriptstyle \bullet$} [c] at 14 8
\put {$ \scriptstyle \bullet$} [c] at 15 4
\put {$ \scriptstyle \bullet$} [c] at 15 12
\put {$ \scriptstyle \bullet$} [c] at 16 8
\setlinear \plot 15 4 13 0 10 4 10 12 15 4 16 8 15 12 14 8 15 4 10 12 /
\setlinear \plot 10 4 15 12 /
\put{$2{,}520$} [c] at 13 -2
\endpicture
\end{minipage}
\begin{minipage}{4cm}
\beginpicture
\setcoordinatesystem units   <1.5mm,2mm>
\setplotarea x from 0 to 16, y from -2 to 15
\put{472)} [l] at 2 12
\put {$ \scriptstyle \bullet$} [c] at 10 8
\put {$ \scriptstyle \bullet$} [c] at 10 0
\put {$ \scriptstyle \bullet$} [c] at 13 12
\put {$ \scriptstyle \bullet$} [c] at 14 4
\put {$ \scriptstyle \bullet$} [c] at 15 0
\put {$ \scriptstyle \bullet$} [c] at 15 8
\put {$ \scriptstyle \bullet$} [c] at 16 4
\setlinear \plot 15 8 13 12 10 8 10 0 15 8 16 4 15 0 14 4 15 8 10 0 /
\setlinear \plot 10 8 15 0 /
\put{$2{,}520$} [c] at 13 -2
\endpicture
\end{minipage}
\begin{minipage}{4cm}
\beginpicture
\setcoordinatesystem units   <1.5mm,2mm>
\setplotarea x from 0 to 16, y from -2 to 15
\put{473)} [l] at 2 12
\put {$ \scriptstyle \bullet$} [c] at 10 4
\put {$ \scriptstyle \bullet$} [c] at 10 8
\put {$ \scriptstyle \bullet$} [c] at 10 12
\put {$ \scriptstyle \bullet$} [c] at 13 0
\put {$ \scriptstyle \bullet$} [c] at 16 4
\put {$ \scriptstyle \bullet$} [c] at 16 8
\put {$ \scriptstyle \bullet$} [c] at 16 12
\setlinear \plot 10 12  10 4  13 0  16 4 16 12   /
\put{$2{,}520$} [c] at 13 -2
\endpicture
\end{minipage}
\begin{minipage}{4cm}
\beginpicture
\setcoordinatesystem units   <1.5mm,2mm>
\setplotarea x from 0 to 16, y from -2 to 15
\put{474)} [l] at 2 12
\put {$ \scriptstyle \bullet$} [c] at 10 4
\put {$ \scriptstyle \bullet$} [c] at 10 8
\put {$ \scriptstyle \bullet$} [c] at 10 0
\put {$ \scriptstyle \bullet$} [c] at 13 12
\put {$ \scriptstyle \bullet$} [c] at 16 4
\put {$ \scriptstyle \bullet$} [c] at 16 8
\put {$ \scriptstyle \bullet$} [c] at 16 0
\setlinear \plot 10 0  10 8  13 12  16 8 16 0   /
\put{$2{,}520$} [c] at 13 -2
\endpicture
\end{minipage}
$$
$$
\begin{minipage}{4cm}
\beginpicture
\setcoordinatesystem units   <1.5mm,2mm>
\setplotarea x from 0 to 16, y from -2 to 15
\put{475)} [l] at 2 12
\put {$ \scriptstyle \bullet$} [c] at 10 4
\put {$ \scriptstyle \bullet$} [c] at 10 12
\put {$ \scriptstyle \bullet$} [c] at 13 0
\put {$ \scriptstyle \bullet$} [c] at 13 4
\put {$ \scriptstyle \bullet$} [c] at 13 8
\put {$ \scriptstyle \bullet$} [c] at 16 4
\put {$ \scriptstyle \bullet$} [c] at 16 12
\setlinear \plot 13 8 13 0 10  4 10 12 13 4  /
\setlinear \plot 13 0 16  4 16 12 13 8 10 4 /
\put{$2{,}520$} [c] at 13 -2
\endpicture
\end{minipage}
\begin{minipage}{4cm}
\beginpicture
\setcoordinatesystem units   <1.5mm,2mm>
\setplotarea x from 0 to 16, y from -2 to 15
\put{476)} [l] at 2 12
\put {$ \scriptstyle \bullet$} [c] at 10 8
\put {$ \scriptstyle \bullet$} [c] at 10 0
\put {$ \scriptstyle \bullet$} [c] at 13 12
\put {$ \scriptstyle \bullet$} [c] at 13 4
\put {$ \scriptstyle \bullet$} [c] at 13 8
\put {$ \scriptstyle \bullet$} [c] at 16 8
\put {$ \scriptstyle \bullet$} [c] at 16 0
\setlinear \plot 13 4 13 12 10  8 10 0 13 8  /
\setlinear \plot 13 12 16  8 16 0 13 4 10 8 /
\put{$2{,}520$} [c] at 13 -2
\endpicture
\end{minipage}
\begin{minipage}{4cm}
\beginpicture
\setcoordinatesystem units   <1.5mm,2mm>
\setplotarea x from 0 to 16, y from -2 to 15
\put{477)} [l] at 2 12
\put {$ \scriptstyle \bullet$} [c] at 10 4
\put {$ \scriptstyle \bullet$} [c] at 10 12
\put {$ \scriptstyle \bullet$} [c] at 13 0
\put {$ \scriptstyle \bullet$} [c] at 13 8
\put {$ \scriptstyle \bullet$} [c] at 13 12
\put {$ \scriptstyle \bullet$} [c] at 16 4
\put {$ \scriptstyle \bullet$} [c] at 16 12
\setlinear \plot 10 4 13 12 13 8 10 12 10 4 13 0 16 4 16 12   /
\setlinear \plot 16 4 13 8 /
\put{$2{,}520$} [c] at 13 -2
\endpicture
\end{minipage}
\begin{minipage}{4cm}
\beginpicture
\setcoordinatesystem units   <1.5mm,2mm>
\setplotarea x from 0 to 16, y from -2 to 15
\put{478)} [l] at 2 12
\put {$ \scriptstyle \bullet$} [c] at 10 8
\put {$ \scriptstyle \bullet$} [c] at 10 0
\put {$ \scriptstyle \bullet$} [c] at 13 0
\put {$ \scriptstyle \bullet$} [c] at 13 4
\put {$ \scriptstyle \bullet$} [c] at 13 12
\put {$ \scriptstyle \bullet$} [c] at 16 0
\put {$ \scriptstyle \bullet$} [c] at 16 8
\setlinear \plot 10 8 13 0 13 4 10 0 10 8 13 12 16 8 16 0   /
\setlinear \plot 16 8 13 4 /
\put{$2{,}520$} [c] at 13 -2
\endpicture
\end{minipage}
\begin{minipage}{4cm}
\beginpicture
\setcoordinatesystem units   <1.5mm,2mm>
\setplotarea x from 0 to 16, y from -2 to 15
\put{479)} [l] at 2 12
\put {$ \scriptstyle \bullet$} [c] at 10 8
\put {$ \scriptstyle \bullet$} [c] at 10 12
\put {$ \scriptstyle \bullet$} [c] at 13 0
\put {$ \scriptstyle \bullet$} [c] at 13 4
\put {$ \scriptstyle \bullet$} [c] at 13 12
\put {$ \scriptstyle \bullet$} [c] at 16 12
\put {$ \scriptstyle \bullet$} [c] at 16 8
\setlinear \plot 10 12 10 8 13 4 13 0     /
\setlinear \plot 13 4 16 8 16 12      /
\setlinear \plot 13 12 16 8       /
\put{$2{,}520$} [c] at 13 -2
\endpicture
\end{minipage}
\begin{minipage}{4cm}
\beginpicture
\setcoordinatesystem units   <1.5mm,2mm>
\setplotarea x from 0 to 16, y from -2 to 15
\put{480)} [l] at 2 12
\put {$ \scriptstyle \bullet$} [c] at 10 4
\put {$ \scriptstyle \bullet$} [c] at 10 0
\put {$ \scriptstyle \bullet$} [c] at 13 12
\put {$ \scriptstyle \bullet$} [c] at 13 8
\put {$ \scriptstyle \bullet$} [c] at 13 0
\put {$ \scriptstyle \bullet$} [c] at 16 0
\put {$ \scriptstyle \bullet$} [c] at 16 4
\setlinear \plot 10 0 10 4 13 8 13 12     /
\setlinear \plot 13 8 16 4 16 0      /
\setlinear \plot 13 0 16 4       /
\put{$2{,}520$} [c] at 13 -2
\endpicture
\end{minipage}
$$

$$
\begin{minipage}{4cm}
\beginpicture
\setcoordinatesystem units   <1.5mm,2mm>
\setplotarea x from 0 to 16, y from -2 to 15
\put{481)} [l] at 2 12
\put {$ \scriptstyle \bullet$} [c] at 10 4
\put {$ \scriptstyle \bullet$} [c] at 10 12
\put {$ \scriptstyle \bullet$} [c] at 13 0
\put {$ \scriptstyle \bullet$} [c] at 13 12
\put {$ \scriptstyle \bullet$} [c] at 16 4
\put {$ \scriptstyle \bullet$} [c] at 16 8
\put {$ \scriptstyle \bullet$} [c] at 16 12
\setlinear \plot 16 4  13 0 10 4 10  12  16 4 16 12 10 4 13 12 16 4   /
\put{$2{,}520$} [c] at 13 -2
\endpicture
\end{minipage}
\begin{minipage}{4cm}
\beginpicture
\setcoordinatesystem units   <1.5mm,2mm>
\setplotarea x from 0 to 16, y from -2 to 15
\put{482)} [l] at 2 12
\put {$ \scriptstyle \bullet$} [c] at 10 8
\put {$ \scriptstyle \bullet$} [c] at 10 0
\put {$ \scriptstyle \bullet$} [c] at 13 0
\put {$ \scriptstyle \bullet$} [c] at 13 12
\put {$ \scriptstyle \bullet$} [c] at 16 4
\put {$ \scriptstyle \bullet$} [c] at 16 8
\put {$ \scriptstyle \bullet$} [c] at 16 0
\setlinear \plot 16 8 13 12 10 8 10  0  16 8 16 0 10 8 13 0 16 8   /
\put{$2{,}520$} [c] at 13 -2
\endpicture
\end{minipage}
\begin{minipage}{4cm}
\beginpicture
\setcoordinatesystem units   <1.5mm,2mm>
\setplotarea x from 0 to 16, y from -2 to 15
\put{483)} [l] at 2 12
\put {$ \scriptstyle \bullet$} [c] at 10 4
\put {$ \scriptstyle \bullet$} [c] at 10 12
\put {$ \scriptstyle \bullet$} [c] at 13 0
\put {$ \scriptstyle \bullet$} [c] at 13 12
\put {$ \scriptstyle \bullet$} [c] at 16 4
\put {$ \scriptstyle \bullet$} [c] at 16 8
\put {$ \scriptstyle \bullet$} [c] at 16 12
\setlinear \plot 13 0 10 4 10 12  16 8 16 4 13 0   /
\setlinear \plot 13 12 16 8 16 12  /
\put{$2{,}520$} [c] at 13 -2
\endpicture
\end{minipage}
\begin{minipage}{4cm}
\beginpicture
\setcoordinatesystem units   <1.5mm,2mm>
\setplotarea x from 0 to 16, y from -2 to 15
\put{484)} [l] at 2 12
\put {$ \scriptstyle \bullet$} [c] at 10 8
\put {$ \scriptstyle \bullet$} [c] at 10 0
\put {$ \scriptstyle \bullet$} [c] at 13 0
\put {$ \scriptstyle \bullet$} [c] at 13 12
\put {$ \scriptstyle \bullet$} [c] at 16 4
\put {$ \scriptstyle \bullet$} [c] at 16 8
\put {$ \scriptstyle \bullet$} [c] at 16 0
\setlinear \plot 13 12 10 8 10 0  16 4 16 8 13 12   /
\setlinear \plot 13 0 16 4 16 0  /
\put{$2{,}520$} [c] at 13 -2
\endpicture
\end{minipage}
\begin{minipage}{4cm}
\beginpicture
\setcoordinatesystem units   <1.5mm,2mm>
\setplotarea x from 0 to 16, y from -2 to 15
\put{485)} [l] at 2 12
\put {$ \scriptstyle \bullet$} [c] at  10 6
\put {$ \scriptstyle \bullet$} [c] at  10.5 3
\put {$ \scriptstyle \bullet$} [c] at  11  0
\put {$ \scriptstyle \bullet$} [c] at  11 12
\put {$ \scriptstyle \bullet$} [c] at  12 6
\put {$ \scriptstyle \bullet$} [c] at  16 0
\put {$ \scriptstyle \bullet$} [c] at  16 12
\setlinear \plot 12 6 11 0 16 12 16 0 12 6 11 12 10 6 11 0 12 6    /
\setlinear \plot  16 0 10.5 3 /
\put{$2{,}520$} [c] at 13 -2
\endpicture
\end{minipage}
\begin{minipage}{4cm}
\beginpicture
\setcoordinatesystem units   <1.5mm,2mm>
\setplotarea x from 0 to 16, y from -2 to 15
\put{486)} [l] at 2 12
\put {$ \scriptstyle \bullet$} [c] at  10 6
\put {$ \scriptstyle \bullet$} [c] at  10.5 9
\put {$ \scriptstyle \bullet$} [c] at  11  0
\put {$ \scriptstyle \bullet$} [c] at  11 12
\put {$ \scriptstyle \bullet$} [c] at  12 6
\put {$ \scriptstyle \bullet$} [c] at  16 0
\put {$ \scriptstyle \bullet$} [c] at  16 12
\setlinear \plot 12 6 11 12 16 0 16 12 12 6 11 0 10 6 11 12 12 6    /
\setlinear \plot  16 12 10.5 9 /
\put{$2{,}520$} [c] at 13 -2
\endpicture
\end{minipage}
$$
$$
\begin{minipage}{4cm}
\beginpicture
\setcoordinatesystem units   <1.5mm,2mm>
\setplotarea x from 0 to 16, y from -2 to 15
\put{487)} [l] at 2 12
\put {$ \scriptstyle \bullet$} [c] at  10 8
\put {$ \scriptstyle \bullet$} [c] at  10 12
\put {$ \scriptstyle \bullet$} [c] at  12 0
\put {$ \scriptstyle \bullet$} [c] at  12 4
\put {$ \scriptstyle \bullet$} [c] at  14 8
\put {$ \scriptstyle \bullet$} [c] at  14 12
\put {$ \scriptstyle \bullet$} [c] at  16 0
\setlinear \plot  12 0 12 4 10 8 10 12 16 0 14 12 14 8 12 4 /
\put{$2{,}520   $} [c] at 13 -2
\endpicture
\end{minipage}
\begin{minipage}{4cm}
\beginpicture
\setcoordinatesystem units   <1.5mm,2mm>
\setplotarea x from 0 to 16, y from -2 to 15
\put{488)} [l] at 2 12
\put {$ \scriptstyle \bullet$} [c] at  10 4
\put {$ \scriptstyle \bullet$} [c] at  10 0
\put {$ \scriptstyle \bullet$} [c] at  12 12
\put {$ \scriptstyle \bullet$} [c] at  12 8
\put {$ \scriptstyle \bullet$} [c] at  14 0
\put {$ \scriptstyle \bullet$} [c] at  14 4
\put {$ \scriptstyle \bullet$} [c] at  16 12
\setlinear \plot  12 12 12 8 10 4 10 0 16 12 14 0 14 4 12 8 /
\put{$2{,}520   $} [c] at 13 -2
\endpicture
\end{minipage}
\begin{minipage}{4cm}
\beginpicture
\setcoordinatesystem units   <1.5mm,2mm>
\setplotarea x from 0 to 16, y from -2 to 15
\put{489)} [l] at 2 12
\put {$ \scriptstyle \bullet$} [c] at  10 4
\put {$ \scriptstyle \bullet$} [c] at  10 12
\put {$ \scriptstyle \bullet$} [c] at  13 0
\put {$ \scriptstyle \bullet$} [c] at  13 8
\put {$ \scriptstyle \bullet$} [c] at  16 0
\put {$ \scriptstyle \bullet$} [c] at  16 4
\put {$ \scriptstyle \bullet$} [c] at  16 12
\setlinear \plot 16 0 16 4 13 0 10 4 10 12 13 8 16 12 16 0 10 4   /
\setlinear \plot 13 8 13 0 /
\put{$2{,}520   $} [c] at 13 -2
\endpicture
\end{minipage}
\begin{minipage}{4cm}
\beginpicture
\setcoordinatesystem units   <1.5mm,2mm>
\setplotarea x from 0 to 16, y from -2 to 15
\put{490)} [l] at 2 12
\put {$ \scriptstyle \bullet$} [c] at  10 0
\put {$ \scriptstyle \bullet$} [c] at  10 8
\put {$ \scriptstyle \bullet$} [c] at  13 4
\put {$ \scriptstyle \bullet$} [c] at  13 12
\put {$ \scriptstyle \bullet$} [c] at  16 0
\put {$ \scriptstyle \bullet$} [c] at  16 8
\put {$ \scriptstyle \bullet$} [c] at  16 12
\setlinear \plot 16 12 16 8 13 12 10 8 10 0 13 4 16 0 16 12 10 8 /
\setlinear \plot 13 4 13 12 /
\put{$2{,}520   $} [c] at 13 -2
\endpicture
\end{minipage}
\begin{minipage}{4cm}
\beginpicture
\setcoordinatesystem units   <1.5mm,2mm>
\setplotarea x from 0 to 16, y from -2 to 15
\put{491)} [l] at 2 12
\put {$ \scriptstyle \bullet$} [c] at  10 4
\put {$ \scriptstyle \bullet$} [c] at  10 8
\put {$ \scriptstyle \bullet$} [c] at  10 12
\put {$ \scriptstyle \bullet$} [c] at  12 0
\put {$ \scriptstyle \bullet$} [c] at  14 4
\put {$ \scriptstyle \bullet$} [c] at  14 12
\put {$ \scriptstyle \bullet$} [c] at  16 0
\setlinear \plot 16 0 14 12 14 4 12 0 10 4 10 12  /
\setlinear \plot 10 4 14 12   /
\setlinear \plot  10 8  14 4 /
\put{$2{,}520$} [c] at 13 -2
\endpicture
\end{minipage}
\begin{minipage}{4cm}
\beginpicture
\setcoordinatesystem units   <1.5mm,2mm>
\setplotarea x from 0 to 16, y from -2 to 15
\put{492)} [l] at 2 12
\put {$ \scriptstyle \bullet$} [c] at  10 4
\put {$ \scriptstyle \bullet$} [c] at  10 8
\put {$ \scriptstyle \bullet$} [c] at  10 0
\put {$ \scriptstyle \bullet$} [c] at  12 12
\put {$ \scriptstyle \bullet$} [c] at  14 0
\put {$ \scriptstyle \bullet$} [c] at  14 8
\put {$ \scriptstyle \bullet$} [c] at  16 12
\setlinear \plot 16 12 14 0 14 8 12 12 10 8 10 0  /
\setlinear \plot 10 8 14 0   /
\setlinear \plot  10 4  14 8 /
\put{$2{,}520$} [c] at 13 -2
\endpicture
\end{minipage}
$$

$$
\begin{minipage}{4cm}
\beginpicture
\setcoordinatesystem units   <1.5mm,2mm>
\setplotarea x from 0 to 16, y from -2 to 15
\put{493)} [l] at 2 12
\put {$ \scriptstyle \bullet$} [c] at  10 0
\put {$ \scriptstyle \bullet$} [c] at  10 4
\put {$ \scriptstyle \bullet$} [c] at  10 8
\put {$ \scriptstyle \bullet$} [c] at  10 12
\put {$ \scriptstyle \bullet$} [c] at  13 0
\put {$ \scriptstyle \bullet$} [c] at  13 12
\put {$ \scriptstyle \bullet$} [c] at  16 12
\setlinear \plot 10 0 10 12  13 0  13 12 10 0   /
\setlinear \plot  10 8 16  12 13  0    /
\put{$2{,}520$} [c] at 13 -2
\endpicture
\end{minipage}
\begin{minipage}{4cm}
\beginpicture
\setcoordinatesystem units   <1.5mm,2mm>
\setplotarea x from 0 to 16, y from -2 to 15
\put{494)} [l] at 2 12
\put {$ \scriptstyle \bullet$} [c] at  10 0
\put {$ \scriptstyle \bullet$} [c] at  10 4
\put {$ \scriptstyle \bullet$} [c] at  10 8
\put {$ \scriptstyle \bullet$} [c] at  10 12
\put {$ \scriptstyle \bullet$} [c] at  13 0
\put {$ \scriptstyle \bullet$} [c] at  13 12
\put {$ \scriptstyle \bullet$} [c] at  16 0
\setlinear \plot 10 12 10 0  13 12  13 0 10 12   /
\setlinear \plot  10 4 16  0 13  12    /
\put{$2{,}520$} [c] at 13 -2
\endpicture
\end{minipage}
\begin{minipage}{4cm}
\beginpicture
\setcoordinatesystem units   <1.5mm,2mm>
\setplotarea x from 0 to 16, y from -2 to 15
\put{495)} [l] at 2 12
\put {$ \scriptstyle \bullet$} [c] at  10 0
\put {$ \scriptstyle \bullet$} [c] at  10 4
\put {$ \scriptstyle \bullet$} [c] at  10 8
\put {$ \scriptstyle \bullet$} [c] at  10 12
\put {$ \scriptstyle \bullet$} [c] at  13 0
\put {$ \scriptstyle \bullet$} [c] at  13 12
\put {$ \scriptstyle \bullet$} [c] at  16 12
\setlinear \plot  10 0 10  4 16 12 13 0 13 12 10 4 10 12   /
\setlinear \plot  10 8 13 0 /
\put{$2{,}520$} [c] at 13 -2
\endpicture
\end{minipage}
\begin{minipage}{4cm}
\beginpicture
\setcoordinatesystem units   <1.5mm,2mm>
\setplotarea x from 0 to 16, y from -2 to 15
\put{496)} [l] at 2 12
\put {$ \scriptstyle \bullet$} [c] at  10 0
\put {$ \scriptstyle \bullet$} [c] at  10 4
\put {$ \scriptstyle \bullet$} [c] at  10 8
\put {$ \scriptstyle \bullet$} [c] at  10 12
\put {$ \scriptstyle \bullet$} [c] at  13 0
\put {$ \scriptstyle \bullet$} [c] at  13 12
\put {$ \scriptstyle \bullet$} [c] at  16 0
\setlinear \plot  10 12 10  8 16 0 13 12 13 0 10 8 10 0   /
\setlinear \plot  10 4 13 12 /
\put{$2{,}520$} [c] at 13 -2
\endpicture
\end{minipage}
\begin{minipage}{4cm}
\beginpicture
\setcoordinatesystem units   <1.5mm,2mm>
\setplotarea x from 0 to 16, y from -2 to 15
\put{497)} [l] at 2 12
\put {$ \scriptstyle \bullet$} [c] at  10 0
\put {$ \scriptstyle \bullet$} [c] at  10 4
\put {$ \scriptstyle \bullet$} [c] at  10 8
\put {$ \scriptstyle \bullet$} [c] at  10 12
\put {$ \scriptstyle \bullet$} [c] at  13 12
\put {$ \scriptstyle \bullet$} [c] at  16 0
\put {$ \scriptstyle \bullet$} [c] at  16 12
\setlinear \plot 10  0 10 12     /
\setlinear \plot   16  12 16 0 10 8  13 12    /
\put{$2{,}520$} [c] at 13 -2
\endpicture
\end{minipage}
\begin{minipage}{4cm}
\beginpicture
\setcoordinatesystem units   <1.5mm,2mm>
\setplotarea x from 0 to 16, y from -2 to 15
\put{498)} [l] at 2 12
\put {$ \scriptstyle \bullet$} [c] at  10 0
\put {$ \scriptstyle \bullet$} [c] at  10 4
\put {$ \scriptstyle \bullet$} [c] at  10 8
\put {$ \scriptstyle \bullet$} [c] at  10 12
\put {$ \scriptstyle \bullet$} [c] at  13 0
\put {$ \scriptstyle \bullet$} [c] at  16 0
\put {$ \scriptstyle \bullet$} [c] at  16 12
\setlinear \plot 10  0 10 12     /
\setlinear \plot   16  0 16 12 10 4  13 0    /
\put{$2{,}520$} [c] at 13 -2
\endpicture
\end{minipage}
$$
$$
\begin{minipage}{4cm}
\beginpicture
\setcoordinatesystem units   <1.5mm,2mm>
\setplotarea x from 0 to 16, y from -2 to 15
\put{499)} [l] at 2 12
\put {$ \scriptstyle \bullet$} [c] at  10 12
\put {$ \scriptstyle \bullet$} [c] at  12 12
\put {$ \scriptstyle \bullet$} [c] at  12 6
\put {$ \scriptstyle \bullet$} [c] at  14 0
\put {$ \scriptstyle \bullet$} [c] at  16 0
\put {$ \scriptstyle \bullet$} [c] at  16 6
\put {$ \scriptstyle \bullet$} [c] at  16 12
\setlinear \plot 16 0 16 12 12 6 14 0 16 6 12  12 12 6 10 12  /
\put{$2{,}520$} [c] at 13 -2
\endpicture
\end{minipage}
\begin{minipage}{4cm}
\beginpicture
\setcoordinatesystem units   <1.5mm,2mm>
\setplotarea x from 0 to 16, y from -2 to 15
\put{500)} [l] at 2 12
\put {$ \scriptstyle \bullet$} [c] at  10 0
\put {$ \scriptstyle \bullet$} [c] at  12 0
\put {$ \scriptstyle \bullet$} [c] at  12 6
\put {$ \scriptstyle \bullet$} [c] at  14 12
\put {$ \scriptstyle \bullet$} [c] at  16 0
\put {$ \scriptstyle \bullet$} [c] at  16 6
\put {$ \scriptstyle \bullet$} [c] at  16 12
\setlinear \plot 16 12 16 0 12 6 14 12 16 6 12  0 12 6 10 0  /
\put{$2{,}520$} [c] at 13 -2
\endpicture
\end{minipage}
\begin{minipage}{4cm}
\beginpicture
\setcoordinatesystem units   <1.5mm,2mm>
\setplotarea x from 0 to 16, y from -2 to 15
\put{501)} [l] at 2 12
\put {$ \scriptstyle \bullet$} [c] at  16 0
\put {$ \scriptstyle \bullet$} [c] at  16 12
\put {$ \scriptstyle \bullet$} [c] at  13 12
\put {$ \scriptstyle \bullet$} [c] at  10 0
\put {$ \scriptstyle \bullet$} [c] at  10 4
\put {$ \scriptstyle \bullet$} [c] at  10 8
\put {$ \scriptstyle \bullet$} [c] at  10 12
\setlinear \plot 10 12 10 0 16 12 16 0 10 4 13 12    /
\put{$2{,}520$} [c] at 13 -2
\endpicture
\end{minipage}
\begin{minipage}{4cm}
\beginpicture
\setcoordinatesystem units   <1.5mm,2mm>
\setplotarea x from 0 to 16, y from -2 to 15
\put{502)} [l] at 2 12
\put {$ \scriptstyle \bullet$} [c] at  16 0
\put {$ \scriptstyle \bullet$} [c] at  16 12
\put {$ \scriptstyle \bullet$} [c] at  13 0
\put {$ \scriptstyle \bullet$} [c] at  10 0
\put {$ \scriptstyle \bullet$} [c] at  10 4
\put {$ \scriptstyle \bullet$} [c] at  10 8
\put {$ \scriptstyle \bullet$} [c] at  10 12
\setlinear \plot 10 0 10 12 16 0 16 12 10 8 13 0    /
\put{$2{,}520$} [c] at 13 -2
\endpicture
\end{minipage}
\begin{minipage}{4cm}
\beginpicture
\setcoordinatesystem units   <1.5mm,2mm>
\setplotarea x from 0 to 16, y from -2 to 15
\put{503)} [l] at 2 12
\put {$ \scriptstyle \bullet$} [c] at  10 0
\put {$ \scriptstyle \bullet$} [c] at  10 6
\put {$ \scriptstyle \bullet$} [c] at  10 12
\put {$ \scriptstyle \bullet$} [c] at  13 12
\put {$ \scriptstyle \bullet$} [c] at  16 0
\put {$ \scriptstyle \bullet$} [c] at  16 6
\put {$ \scriptstyle \bullet$} [c] at  16 12
\setlinear \plot  16 12 16 0 10 12 10 0 /
\setlinear \plot  16 6 13 12 10 6  16 12 /
\put{$2{,}520$} [c] at 13 -2
\endpicture
\end{minipage}
\begin{minipage}{4cm}
\beginpicture
\setcoordinatesystem units   <1.5mm,2mm>
\setplotarea x from 0 to 16, y from -2 to 15
\put{504)} [l] at 2 12
\put {$ \scriptstyle \bullet$} [c] at  10 0
\put {$ \scriptstyle \bullet$} [c] at  10 6
\put {$ \scriptstyle \bullet$} [c] at  10 12
\put {$ \scriptstyle \bullet$} [c] at  13 0
\put {$ \scriptstyle \bullet$} [c] at  16 0
\put {$ \scriptstyle \bullet$} [c] at  16 6
\put {$ \scriptstyle \bullet$} [c] at  16 12
\setlinear \plot  16 0 16 12 10 0 10 12 /
\setlinear \plot  16 6 13 0 10 6  16 0 /
\put{$2{,}520$} [c] at 13 -2
\endpicture
\end{minipage}
$$

$$
\begin{minipage}{4cm}
\beginpicture
\setcoordinatesystem units   <1.5mm,2mm>
\setplotarea x from 0 to 16, y from -2 to 15
\put{505)} [l] at 2 12
\put {$ \scriptstyle \bullet$} [c] at 10 4
\put {$ \scriptstyle \bullet$} [c] at 10 8
\put {$ \scriptstyle \bullet$} [c] at 11 0
\put {$ \scriptstyle \bullet$} [c] at 11 12
\put {$ \scriptstyle \bullet$} [c] at 12 4
\put {$ \scriptstyle \bullet$} [c] at 12 8
\put {$ \scriptstyle \bullet$} [c] at 16 12
\setlinear \plot  16 12 11 0 10 4 10 8 11 12 12 8 12 4 11 0  /
\setlinear \plot  10 4 12 8 /
\setlinear \plot  10 8 12 4 /
\put{$1{,}260$} [c] at 13 -2
\endpicture
\end{minipage}
\begin{minipage}{4cm}
\beginpicture
\setcoordinatesystem units   <1.5mm,2mm>
\setplotarea x from 0 to 16, y from -2 to 15
\put{506)} [l] at 2 12
\put {$ \scriptstyle \bullet$} [c] at 10 4
\put {$ \scriptstyle \bullet$} [c] at 10 8
\put {$ \scriptstyle \bullet$} [c] at 11 0
\put {$ \scriptstyle \bullet$} [c] at 11 12
\put {$ \scriptstyle \bullet$} [c] at 12 4
\put {$ \scriptstyle \bullet$} [c] at 12 8
\put {$ \scriptstyle \bullet$} [c] at 16 0
\setlinear \plot  16 0 11 12 10 8 10 4 11 0 12 4 12 8 11 12  /
\setlinear \plot  10 4 12 8 /
\setlinear \plot  10 8 12 4 /
\put{$1{,}260$} [c] at 13 -2
\endpicture
\end{minipage}
\begin{minipage}{4cm}
\beginpicture
\setcoordinatesystem units   <1.5mm,2mm>
\setplotarea x from 0 to 16, y from -2 to 15
\put{507)} [l] at 2 12
\put {$ \scriptstyle \bullet$} [c] at 13 0
\put {$ \scriptstyle \bullet$} [c] at 10 6
\put {$ \scriptstyle \bullet$} [c] at 13 6
\put {$ \scriptstyle \bullet$} [c] at 16 6
\put {$ \scriptstyle \bullet$} [c] at 10 12
\put {$ \scriptstyle \bullet$} [c] at 13 12
\put {$ \scriptstyle \bullet$} [c] at 16 12
\setlinear \plot  13  0 10 6 10 12 13 6 13 0 16 6  13 12  10 6 16 12 13 6 13 12 /
\setlinear \plot  10  12 16 6 /
\put{$1{,}260$} [c] at 13 -2
\endpicture
\end{minipage}
\begin{minipage}{4cm}
\beginpicture
\setcoordinatesystem units   <1.5mm,2mm>
\setplotarea x from 0 to 16, y from -2 to 15
\put{508)} [l] at 2 12
\put {$ \scriptstyle \bullet$} [c] at 13 12
\put {$ \scriptstyle \bullet$} [c] at 10 6
\put {$ \scriptstyle \bullet$} [c] at 13 6
\put {$ \scriptstyle \bullet$} [c] at 16 6
\put {$ \scriptstyle \bullet$} [c] at 10 0
\put {$ \scriptstyle \bullet$} [c] at 13 0
\put {$ \scriptstyle \bullet$} [c] at 16 0
\setlinear \plot  13  12 10 6 10 0 13 6 13 12 16 6  13 0  10 6 16 0 13 6 13 0 /
\setlinear \plot  10  0 16 6 /
\put{$1{,}260$} [c] at 13 -2
\endpicture
\end{minipage}
\begin{minipage}{4cm}
\beginpicture
\setcoordinatesystem units   <1.5mm,2mm>
\setplotarea x from 0 to 16, y from -2 to 15
\put{509)} [l] at 2 12
\put {$ \scriptstyle \bullet$} [c] at 10 8
\put {$ \scriptstyle \bullet$} [c] at 10 12
\put {$ \scriptstyle \bullet$} [c] at 11 0
\put {$ \scriptstyle \bullet$} [c] at 11 4
\put {$ \scriptstyle \bullet$} [c] at 12 8
\put {$ \scriptstyle \bullet$} [c] at 12 12
\put {$ \scriptstyle \bullet$} [c] at 16 12
\setlinear \plot  16 12 11 0 11 4 10 8 10 12 12 8 12 12 10 8 /
\setlinear \plot  11 4 12 8 /
\put{$1{,}260$} [c] at 13 -2
\endpicture
\end{minipage}
\begin{minipage}{4cm}
\beginpicture
\setcoordinatesystem units   <1.5mm,2mm>
\setplotarea x from 0 to 16, y from -2 to 15
\put{510)} [l] at 2 12
\put {$ \scriptstyle \bullet$} [c] at 10 4
\put {$ \scriptstyle \bullet$} [c] at 10 0
\put {$ \scriptstyle \bullet$} [c] at 11 12
\put {$ \scriptstyle \bullet$} [c] at 11 8
\put {$ \scriptstyle \bullet$} [c] at 12 4
\put {$ \scriptstyle \bullet$} [c] at 12 0
\put {$ \scriptstyle \bullet$} [c] at 16 0
\setlinear \plot  16 0 11 12 11 8  10 4 10 0 12 4 12 0 10 4 /
\setlinear \plot  11 8 12 4 /
\put{$1{,}260$} [c] at 13 -2
\endpicture
\end{minipage}
$$
$$
\begin{minipage}{4cm}
\beginpicture
\setcoordinatesystem units   <1.5mm,2mm>
\setplotarea x from 0 to 16, y from -2 to 15
\put{511)} [l] at 2 12
\put {$ \scriptstyle \bullet$} [c] at 10 4
\put {$ \scriptstyle \bullet$} [c] at 10 12
\put {$ \scriptstyle \bullet$} [c] at 11 0
\put {$ \scriptstyle \bullet$} [c] at 11 8
\put {$ \scriptstyle \bullet$} [c] at 12 4
\put {$ \scriptstyle \bullet$} [c] at 12 12
\put {$ \scriptstyle \bullet$} [c] at 16 12
\setlinear \plot  16 12 11 0 10 4 11 8 12 12  /
\setlinear \plot  10  12 11 8 12 4 11 0 /
\put{$1{,}260$} [c] at 13 -2
\endpicture
\end{minipage}
\begin{minipage}{4cm}
\beginpicture
\setcoordinatesystem units   <1.5mm,2mm>
\setplotarea x from 0 to 16, y from -2 to 15
\put{512)} [l] at 2 12
\put {$ \scriptstyle \bullet$} [c] at 10 8
\put {$ \scriptstyle \bullet$} [c] at 10 0
\put {$ \scriptstyle \bullet$} [c] at 11 12
\put {$ \scriptstyle \bullet$} [c] at 11 4
\put {$ \scriptstyle \bullet$} [c] at 12 8
\put {$ \scriptstyle \bullet$} [c] at 12 0
\put {$ \scriptstyle \bullet$} [c] at 16 0
\setlinear \plot  16 0 11 12 10 8 11 4 12 0  /
\setlinear \plot  10  0 11 4 12 8 11 12 /
\put{$1{,}260$} [c] at 13 -2
\endpicture
\end{minipage}
\begin{minipage}{4cm}
\beginpicture
\setcoordinatesystem units   <1.5mm,2mm>
\setplotarea x from 0 to 16, y from -2 to 15
\put{513)} [l] at 2 12
\put {$ \scriptstyle \bullet$} [c] at  10 6
\put {$ \scriptstyle \bullet$} [c] at  11.5 0
\put {$ \scriptstyle \bullet$} [c] at  13 6
\put {$ \scriptstyle \bullet$} [c] at  13 12
\put {$ \scriptstyle \bullet$} [c] at  14.5 0
\put {$ \scriptstyle \bullet$} [c] at  16 6
\put {$ \scriptstyle \bullet$} [c] at  16 12
\setlinear \plot 10 6 13 12 16 6 16 12    /
\setlinear \plot 16 6 11.5 0 10 6 14.5 0 16 6      /
\setlinear \plot 11.5 0 13 6 13 12 /
\setlinear \plot 14.5 0 13 6 /
\put{$1{,}260$} [c] at 13 -2
\endpicture
\end{minipage}
\begin{minipage}{4cm}
\beginpicture
\setcoordinatesystem units   <1.5mm,2mm>
\setplotarea x from 0 to 16, y from -2 to 15
\put{514)} [l] at 2 12
\put {$ \scriptstyle \bullet$} [c] at  10 6
\put {$ \scriptstyle \bullet$} [c] at  11.5 12
\put {$ \scriptstyle \bullet$} [c] at  13 6
\put {$ \scriptstyle \bullet$} [c] at  13 0
\put {$ \scriptstyle \bullet$} [c] at  14.5 12
\put {$ \scriptstyle \bullet$} [c] at  16 6
\put {$ \scriptstyle \bullet$} [c] at  16 0
\setlinear \plot 10 6 13 0 16 6 16 0    /
\setlinear \plot 16 6 11.5 12 10 6 14.5 12 16 6      /
\setlinear \plot 11.5 12 13 6 13 0 /
\setlinear \plot 14.5 12 13 6 /
\put{$1{,}260$} [c] at 13 -2
\endpicture
\end{minipage}
\begin{minipage}{4cm}
\beginpicture
\setcoordinatesystem units   <1.5mm,2mm>
\setplotarea x from 0 to 16, y from -2 to 15
\put{515)} [l] at 2 12
\put {$ \scriptstyle \bullet$} [c] at  10 12
\put {$ \scriptstyle \bullet$} [c] at  12 0
\put {$ \scriptstyle \bullet$} [c] at  12 6
\put {$ \scriptstyle \bullet$} [c] at  12 12
\put {$ \scriptstyle \bullet$} [c] at  14 12
\put {$ \scriptstyle \bullet$} [c] at  16 6
\put {$ \scriptstyle \bullet$} [c] at  16 0
\setlinear \plot  10 12 12 6 14 12 16 6 16 0 12 6 12 0 16 6 14 12 /
\setlinear \plot  12 12 12 6 /
\put{$1{,}260$} [c] at 13 -2
\endpicture
\end{minipage}
\begin{minipage}{4cm}
\beginpicture
\setcoordinatesystem units   <1.5mm,2mm>
\setplotarea x from 0 to 16, y from -2 to 15
\put{516)} [l] at 2 12
\put {$ \scriptstyle \bullet$} [c] at  10 0
\put {$ \scriptstyle \bullet$} [c] at  12 0
\put {$ \scriptstyle \bullet$} [c] at  12 6
\put {$ \scriptstyle \bullet$} [c] at  12 12
\put {$ \scriptstyle \bullet$} [c] at  14 0
\put {$ \scriptstyle \bullet$} [c] at  16 6
\put {$ \scriptstyle \bullet$} [c] at  16 12
\setlinear \plot  10 0 12 6 14 0 16 6 16 12 12 6 12 12 16 6 14 0 /
\setlinear \plot  12 0 12 6 /
\put{$1{,}260$} [c] at 13 -2
\endpicture
\end{minipage}
$$

$$
\begin{minipage}{4cm}
\beginpicture
\setcoordinatesystem units   <1.5mm,2mm>
\setplotarea x from 0 to 16, y from -2 to 15
\put{517)} [l] at 2 12
\put {$ \scriptstyle \bullet$} [c] at  10 12
\put {$ \scriptstyle \bullet$} [c] at  14 0
\put {$ \scriptstyle \bullet$} [c] at  14 6
\put {$ \scriptstyle \bullet$} [c] at  14 12
\put {$ \scriptstyle \bullet$} [c] at  16 0
\put {$ \scriptstyle \bullet$} [c] at  16 6
\put {$ \scriptstyle \bullet$} [c] at  16 12
\setlinear \plot 16 0 10 12 14 0 14 12  16 6 16 12 14 6 16 0 16 6 14 0  /
\put{$630$} [c] at 13 -2
\endpicture
\end{minipage}
\begin{minipage}{4cm}
\beginpicture
\setcoordinatesystem units   <1.5mm,2mm>
\setplotarea x from 0 to 16, y from -2 to 15
\put{518)} [l] at 2 12
\put {$ \scriptstyle \bullet$} [c] at  10 0
\put {$ \scriptstyle \bullet$} [c] at  14 0
\put {$ \scriptstyle \bullet$} [c] at  14 6
\put {$ \scriptstyle \bullet$} [c] at  14 12
\put {$ \scriptstyle \bullet$} [c] at  16 0
\put {$ \scriptstyle \bullet$} [c] at  16 6
\put {$ \scriptstyle \bullet$} [c] at  16 12
\setlinear \plot 16 12 10 0 14 12 14 0  16 6 16 0 14 6 16 12 16 6 14 12  /
\put{$630$} [c] at 13 -2
\endpicture
\end{minipage}
\begin{minipage}{4cm}
\beginpicture
\setcoordinatesystem units   <1.5mm,2mm>
\setplotarea x from 0 to 16, y from -2 to 15
\put{519)} [l] at 2 12
\put {$ \scriptstyle \bullet$} [c] at  10 12
\put {$ \scriptstyle \bullet$} [c] at  10 6
\put {$ \scriptstyle \bullet$} [c] at  11 0
\put {$ \scriptstyle \bullet$} [c] at  12  12
\put {$ \scriptstyle \bullet$} [c] at  12 6
\put {$ \scriptstyle \bullet$} [c] at  16 0
\put {$ \scriptstyle \bullet$} [c] at  16 12
\setlinear \plot  10 12 10 6 11 0 12 6 12 12 10 6 16 12 16 0 10 12 12 6 16 12  /
\setlinear \plot  12  12 16 0 /
\put{$420$} [c] at 13 -2
\endpicture
\end{minipage}
\begin{minipage}{4cm}
\beginpicture
\setcoordinatesystem units   <1.5mm,2mm>
\setplotarea x from 0 to 16, y from -2 to 15
\put{520)} [l] at 2 12
\put {$ \scriptstyle \bullet$} [c] at  10 0
\put {$ \scriptstyle \bullet$} [c] at  10 6
\put {$ \scriptstyle \bullet$} [c] at  11 12
\put {$ \scriptstyle \bullet$} [c] at  12  0
\put {$ \scriptstyle \bullet$} [c] at  12 6
\put {$ \scriptstyle \bullet$} [c] at  16 0
\put {$ \scriptstyle \bullet$} [c] at  16 12
\setlinear \plot  10 0 10 6 11 12 12 6 12 0 10 6 16 0 16 12 10 0 12 6 16 0  /
\setlinear \plot  12  0 16 12 /
\put{$420$} [c] at 13 -2
\endpicture
\end{minipage}
\begin{minipage}{4cm}
\beginpicture
\setcoordinatesystem units   <1.5mm,2mm>
\setplotarea x from 0 to 16, y from -2 to 15
\put{521)} [l] at 2 12
\put {$ \scriptstyle \bullet$} [c] at  10 0
\put {$ \scriptstyle \bullet$} [c] at  10 12
\put {$ \scriptstyle \bullet$} [c] at  13 0
\put {$ \scriptstyle \bullet$} [c] at  13 3
\put {$ \scriptstyle \bullet$} [c] at  13 12
\put {$ \scriptstyle \bullet$} [c] at  16  0
\put {$ \scriptstyle \bullet$} [c] at  16 12
\setlinear \plot 16 12 16 0  10 12 13 3 16 12 /
\setlinear \plot 16 0 13 12 13 0   /
\setlinear \plot 10 0 13 3    /
\put{$420$} [c] at 12 -2
\endpicture
\end{minipage}
\begin{minipage}{4cm}
\beginpicture
\setcoordinatesystem units   <1.5mm,2mm>
\setplotarea x from 0 to 16, y from -2 to 15
\put{522)} [l] at 2 12
\put {$ \scriptstyle \bullet$} [c] at  10 0
\put {$ \scriptstyle \bullet$} [c] at  10 12
\put {$ \scriptstyle \bullet$} [c] at  13 0
\put {$ \scriptstyle \bullet$} [c] at  13 9
\put {$ \scriptstyle \bullet$} [c] at  13 12
\put {$ \scriptstyle \bullet$} [c] at  16  0
\put {$ \scriptstyle \bullet$} [c] at  16 12
\setlinear \plot 16 0 16 12  10 0 13 9 16 0 /
\setlinear \plot 16 12 13 0 13 12   /
\setlinear \plot 10 12 13 9    /
\put{$420 $} [c] at 12 -2
\endpicture
\end{minipage}
$$

$$
\begin{minipage}{4cm}
\beginpicture
\setcoordinatesystem units   <1.5mm,2mm>
\setplotarea x from 0 to 16, y from -2 to 15
\put{${\bf  18}$} [l] at 2 15
\put{523)} [l] at 2 12
\put {$ \scriptstyle \bullet$} [c] at 10 6
\put {$ \scriptstyle \bullet$} [c] at 12 6
\put {$ \scriptstyle \bullet$} [c] at 16 6
\put {$ \scriptstyle \bullet$} [c] at 14 0
\put {$ \scriptstyle \bullet$} [c] at 13 12
\put {$ \scriptstyle \bullet$} [c] at 12.6 9
\put {$ \scriptstyle \bullet$} [c] at 15.2 3
\setlinear \plot 14  0 10 6 13 12 16 6 14 0 12 6  13 12   /
\setlinear \plot 12.6 9 15.2 3 /
\put{$5{,}040$} [c] at 13 -2
\endpicture
\end{minipage}
\begin{minipage}{4cm}
\beginpicture
\setcoordinatesystem units   <1.5mm,2mm>
\setplotarea x from 0 to 16, y from -2 to 15
\put{524)} [l] at 2 12
\put {$ \scriptstyle \bullet$} [c] at 10 4
\put {$ \scriptstyle \bullet$} [c] at 10 8
\put {$ \scriptstyle \bullet$} [c] at 10 12
\put {$ \scriptstyle \bullet$} [c] at 13 0
\put {$ \scriptstyle \bullet$} [c] at 13 8
\put {$ \scriptstyle \bullet$} [c] at 16 4
\put {$ \scriptstyle \bullet$} [c] at 16 12
\setlinear \plot  10 12 10 4 16 12 16 4 13 0 10 4   /
\put{$5{,}040$} [c] at 13 -2
\endpicture
\end{minipage}
\begin{minipage}{4cm}
\beginpicture
\setcoordinatesystem units   <1.5mm,2mm>
\setplotarea x from 0 to 16, y from -2 to 15
\put{525)} [l] at 2 12
\put {$ \scriptstyle \bullet$} [c] at 10 4
\put {$ \scriptstyle \bullet$} [c] at 10 8
\put {$ \scriptstyle \bullet$} [c] at 10 0
\put {$ \scriptstyle \bullet$} [c] at 13 12
\put {$ \scriptstyle \bullet$} [c] at 13 4
\put {$ \scriptstyle \bullet$} [c] at 16 8
\put {$ \scriptstyle \bullet$} [c] at 16 0
\setlinear \plot  10 0 10 8 16 0 16 8 13 12 10 8   /
\put{$5{,}040$} [c] at 13 -2
\endpicture
\end{minipage}
\begin{minipage}{4cm}
\beginpicture
\setcoordinatesystem units   <1.5mm,2mm>
\setplotarea x from 0 to 16, y from -2 to 15
\put{526)} [l] at 2 12
\put {$ \scriptstyle \bullet$} [c] at 10 4
\put {$ \scriptstyle \bullet$} [c] at 10 8
\put {$ \scriptstyle \bullet$} [c] at 10 12
\put {$ \scriptstyle \bullet$} [c] at 13 0
\put {$ \scriptstyle \bullet$} [c] at 13 8
\put {$ \scriptstyle \bullet$} [c] at 16 4
\put {$ \scriptstyle \bullet$} [c] at 16 12
\setlinear \plot  16 12 16 4 10 12 10 4 13 0 16 4   /
\put{$5{,}040$} [c] at 13 -2
\endpicture
\end{minipage}
\begin{minipage}{4cm}
\beginpicture
\setcoordinatesystem units   <1.5mm,2mm>
\setplotarea x from 0 to 16, y from -2 to 15
\put{527)} [l] at 2 12
\put {$ \scriptstyle \bullet$} [c] at 10 4
\put {$ \scriptstyle \bullet$} [c] at 10 8
\put {$ \scriptstyle \bullet$} [c] at 10 0
\put {$ \scriptstyle \bullet$} [c] at 13 12
\put {$ \scriptstyle \bullet$} [c] at 13 4
\put {$ \scriptstyle \bullet$} [c] at 16 8
\put {$ \scriptstyle \bullet$} [c] at 16 0
\setlinear \plot  16 0 16 8 10 0 10 8 13 12 16 8   /
\put{$5{,}040$} [c] at 13 -2
\endpicture
\end{minipage}
\begin{minipage}{4cm}
\beginpicture
\setcoordinatesystem units   <1.5mm,2mm>
\setplotarea x from 0 to 16, y from -2 to 15
\put{528)} [l] at 2 12
\put {$ \scriptstyle \bullet$} [c] at 10 4
\put {$ \scriptstyle \bullet$} [c] at 10 12
\put {$ \scriptstyle \bullet$} [c] at 13 0
\put {$ \scriptstyle \bullet$} [c] at 13 4
\put {$ \scriptstyle \bullet$} [c] at 16 4
\put {$ \scriptstyle \bullet$} [c] at 16 8
\put {$ \scriptstyle \bullet$} [c] at 16 12
\setlinear \plot 13 4 10 12  10 4 13 0 13 4 16 12 16 4 10 12 /
\setlinear \plot 13 0  16 4 /
\put{$5{,}040$} [c] at 13 -2
\endpicture
\end{minipage}
$$
$$
\begin{minipage}{4cm}
\beginpicture
\setcoordinatesystem units   <1.5mm,2mm>
\setplotarea x from 0 to 16, y from -2 to 15
\put{529)} [l] at 2 12
\put {$ \scriptstyle \bullet$} [c] at 10 0
\put {$ \scriptstyle \bullet$} [c] at 10 8
\put {$ \scriptstyle \bullet$} [c] at 13 12
\put {$ \scriptstyle \bullet$} [c] at 13 8
\put {$ \scriptstyle \bullet$} [c] at 16 4
\put {$ \scriptstyle \bullet$} [c] at 16 8
\put {$ \scriptstyle \bullet$} [c] at 16 0
\setlinear \plot 13 8 10 0  10 8 13 12 13 8 16 0 16 8 10 0 /
\setlinear \plot 13 12  16 8 /
\put{$5{,}040$} [c] at 13 -2
\endpicture
\end{minipage}
\begin{minipage}{4cm}
\beginpicture
\setcoordinatesystem units   <1.5mm,2mm>
\setplotarea x from 0 to 16, y from -2 to 15
\put{530)} [l] at 2 12
\put {$ \scriptstyle \bullet$} [c] at 10 12
\put {$ \scriptstyle \bullet$} [c] at 12 9
\put {$ \scriptstyle \bullet$} [c] at 13 0
\put {$ \scriptstyle \bullet$} [c] at 13 4
\put {$ \scriptstyle \bullet$} [c] at 13 12
\put {$ \scriptstyle \bullet$} [c] at 14 9
\put {$ \scriptstyle \bullet$} [c] at 16 12
\setlinear \plot 16 12 14 9 13 4 13 0  /
\setlinear \plot 10 12 13 4 12 9 13 12 14 9  /
\put{$5{,}040$} [c] at 13 -2
\endpicture
\end{minipage}
\begin{minipage}{4cm}
\beginpicture
\setcoordinatesystem units   <1.5mm,2mm>
\setplotarea x from 0 to 16, y from -2 to 15
\put{531)} [l] at 2 12
\put {$ \scriptstyle \bullet$} [c] at 10 0
\put {$ \scriptstyle \bullet$} [c] at 12 3
\put {$ \scriptstyle \bullet$} [c] at 13 0
\put {$ \scriptstyle \bullet$} [c] at 13 8
\put {$ \scriptstyle \bullet$} [c] at 13 12
\put {$ \scriptstyle \bullet$} [c] at 14 3
\put {$ \scriptstyle \bullet$} [c] at 16 0
\setlinear \plot 16 0 14 3 13 8 13 12  /
\setlinear \plot 10 0 13 8 12 3 13 0 14 3  /
\put{$5{,}040$} [c] at 13 -2
\endpicture
\end{minipage}
\begin{minipage}{4cm}
\beginpicture
\setcoordinatesystem units   <1.5mm,2mm>
\setplotarea x from 0 to 16, y from -2 to 15
\put{532)} [l] at 2 12
\put {$ \scriptstyle \bullet$} [c] at 10 4
\put {$ \scriptstyle \bullet$} [c] at 10 12
\put {$ \scriptstyle \bullet$} [c] at 13 0
\put {$ \scriptstyle \bullet$} [c] at 13 12
\put {$ \scriptstyle \bullet$} [c] at 16 4
\put {$ \scriptstyle \bullet$} [c] at 16 8
\put {$ \scriptstyle \bullet$} [c] at 16 12
\setlinear \plot 10 12 10 4  13 0 16 4 16 12 10 4 13 12 16 4 /
\put{$5{,}040$} [c] at 13 -2
\endpicture
\end{minipage}
\begin{minipage}{4cm}
\beginpicture
\setcoordinatesystem units   <1.5mm,2mm>
\setplotarea x from 0 to 16, y from -2 to 15
\put{533)} [l] at 2 12
\put {$ \scriptstyle \bullet$} [c] at 10 8
\put {$ \scriptstyle \bullet$} [c] at 10 0
\put {$ \scriptstyle \bullet$} [c] at 13 0
\put {$ \scriptstyle \bullet$} [c] at 13 12
\put {$ \scriptstyle \bullet$} [c] at 16 4
\put {$ \scriptstyle \bullet$} [c] at 16 8
\put {$ \scriptstyle \bullet$} [c] at 16 0
\setlinear \plot 10 0 10 8  13 12 16 8 16 0 10 8 13 0 16 8 /
\put{$5{,}040$} [c] at 13 -2
\endpicture
\end{minipage}
\begin{minipage}{4cm}
\beginpicture
\setcoordinatesystem units   <1.5mm,2mm>
\setplotarea x from 0 to 16, y from -2 to 15
\put{534)} [l] at 2 12
\put {$ \scriptstyle \bullet$} [c] at  10 4
\put {$ \scriptstyle \bullet$} [c] at  10 8
\put {$ \scriptstyle \bullet$} [c] at  10 12
\put {$ \scriptstyle \bullet$} [c] at  12 0
\put {$ \scriptstyle \bullet$} [c] at  14 4
\put {$ \scriptstyle \bullet$} [c] at  14 12
\put {$ \scriptstyle \bullet$} [c] at  16 0
\setlinear \plot 10 12 10 4 12 0 14 4 14 12 16 0 10 8    /
\put{$5{,}040   $} [c] at 13 -2
\endpicture
\end{minipage}
$$

$$
\begin{minipage}{4cm}
\beginpicture
\setcoordinatesystem units   <1.5mm,2mm>
\setplotarea x from 0 to 16, y from -2 to 15
\put{535)} [l] at 2 12
\put {$ \scriptstyle \bullet$} [c] at  10 4
\put {$ \scriptstyle \bullet$} [c] at  10 8
\put {$ \scriptstyle \bullet$} [c] at  10 0
\put {$ \scriptstyle \bullet$} [c] at  12 12
\put {$ \scriptstyle \bullet$} [c] at  14 8
\put {$ \scriptstyle \bullet$} [c] at  14 0
\put {$ \scriptstyle \bullet$} [c] at  16 12
\setlinear \plot 10 0 10 8 12 12 14 8 14 0 16 12 10 4    /
\put{$5{,}040   $} [c] at 13 -2
\endpicture
\end{minipage}
\begin{minipage}{4cm}
\beginpicture
\setcoordinatesystem units   <1.5mm,2mm>
\setplotarea x from 0 to 16, y from -2 to 15
\put{536)} [l] at 2 12
\put {$ \scriptstyle \bullet$} [c] at  10 4
\put {$ \scriptstyle \bullet$} [c] at  10 8
\put {$ \scriptstyle \bullet$} [c] at  10 12
\put {$ \scriptstyle \bullet$} [c] at  13 0
\put {$ \scriptstyle \bullet$} [c] at  16 0
\put {$ \scriptstyle \bullet$} [c] at  16 4
\put {$ \scriptstyle \bullet$} [c] at  16 12
\setlinear \plot 10 12 10 4 13 0 16  4 16 12   /
\setlinear \plot 16 0 16 4 /
\put{$5{,}040   $} [c] at 13 -2
\endpicture
\end{minipage}
\begin{minipage}{4cm}
\beginpicture
\setcoordinatesystem units   <1.5mm,2mm>
\setplotarea x from 0 to 16, y from -2 to 15
\put{537)} [l] at 2 12
\put {$ \scriptstyle \bullet$} [c] at  10 0
\put {$ \scriptstyle \bullet$} [c] at  10 4
\put {$ \scriptstyle \bullet$} [c] at  10 8
\put {$ \scriptstyle \bullet$} [c] at  13 12
\put {$ \scriptstyle \bullet$} [c] at  16 0
\put {$ \scriptstyle \bullet$} [c] at  16 8
\put {$ \scriptstyle \bullet$} [c] at  16 12
\setlinear \plot 10 0 10 8  13 12 16  8 16 0   /
\setlinear \plot 16 12 16 8 /
\put{$5{,}040   $} [c] at 13 -2
\endpicture
\end{minipage}
\begin{minipage}{4cm}
\beginpicture
\setcoordinatesystem units   <1.5mm,2mm>
\setplotarea x from 0 to 16, y from -2 to 15
\put{538)} [l] at 2 12
\put {$ \scriptstyle \bullet$} [c] at  10 0
\put {$ \scriptstyle \bullet$} [c] at  10 12
\put {$ \scriptstyle \bullet$} [c] at  11.1 9
\put {$ \scriptstyle \bullet$} [c] at  12.2 6.2
\put {$ \scriptstyle \bullet$} [c] at  13 4
\put {$ \scriptstyle \bullet$} [c] at  13 0
\put {$ \scriptstyle \bullet$} [c] at  16 12
\setlinear \plot  10   0 10  12  13 4 13 0    /
\setlinear \plot  13  4 16 12    /
\put{$5{,}040   $} [c] at 13 -2
\endpicture
\end{minipage}
\begin{minipage}{4cm}
\beginpicture
\setcoordinatesystem units   <1.5mm,2mm>
\setplotarea x from 0 to 16, y from -2 to 15
\put{539)} [l] at 2 12
\put {$ \scriptstyle \bullet$} [c] at  10 0
\put {$ \scriptstyle \bullet$} [c] at  10 12
\put {$ \scriptstyle \bullet$} [c] at  11.1 3
\put {$ \scriptstyle \bullet$} [c] at  12.2 5.8
\put {$ \scriptstyle \bullet$} [c] at  13 8
\put {$ \scriptstyle \bullet$} [c] at  13 12
\put {$ \scriptstyle \bullet$} [c] at  16 0
\setlinear \plot  10   12 10  0  13 8 13 12    /
\setlinear \plot  13  8 16 0    /
\put{$5{,}040   $} [c] at 13 -2
\endpicture
\end{minipage}
\begin{minipage}{4cm}
\beginpicture
\setcoordinatesystem units   <1.5mm,2mm>
\setplotarea x from 0 to 16, y from -2 to 15
\put{540)} [l] at 2 12
\put {$ \scriptstyle \bullet$} [c] at  10 12
\put {$ \scriptstyle \bullet$} [c] at  12 0
\put {$ \scriptstyle \bullet$} [c] at  12 8
\put {$ \scriptstyle \bullet$} [c] at  14 12
\put {$ \scriptstyle \bullet$} [c] at  16 0
\put {$ \scriptstyle \bullet$} [c] at  16 8
\put {$ \scriptstyle \bullet$} [c] at  13 10
\setlinear \plot  10 12 12 0 12 8  14 12  16 8 16 0  12 8    /
\setlinear \plot  12 0 16 8    /
\put{$5{,}040  $} [c] at 13 -2
\endpicture
\end{minipage}
$$
$$
\begin{minipage}{4cm}
\beginpicture
\setcoordinatesystem units   <1.5mm,2mm>
\setplotarea x from 0 to 16, y from -2 to 15
\put{541)} [l] at 2 12
\put {$ \scriptstyle \bullet$} [c] at  10 0
\put {$ \scriptstyle \bullet$} [c] at  12 12
\put {$ \scriptstyle \bullet$} [c] at  12 4
\put {$ \scriptstyle \bullet$} [c] at  14 0
\put {$ \scriptstyle \bullet$} [c] at  16 4
\put {$ \scriptstyle \bullet$} [c] at  16 12
\put {$ \scriptstyle \bullet$} [c] at  13 2
\setlinear \plot  10 0 12 12 12 4  14 0  16 4 16 12  12 4    /
\setlinear \plot  12 12 16 4    /
\put{$5{,}040  $} [c] at 13 -2
\endpicture
\end{minipage}
\begin{minipage}{4cm}
\beginpicture
\setcoordinatesystem units   <1.5mm,2mm>
\setplotarea x from 0 to 16, y from -2 to 15
\put{542)} [l] at 2 12
\put {$\scriptstyle  \bullet$} [c] at  10 12
\put {$\scriptstyle  \bullet$} [c] at  10 8
\put {$\scriptstyle  \bullet$} [c] at  11 0
\put {$\scriptstyle  \bullet$} [c] at  11 4
\put {$\scriptstyle  \bullet$} [c] at  11 12
\put {$\scriptstyle  \bullet$} [c] at  12 8
\put {$\scriptstyle  \bullet$} [c] at  16 0
\setlinear \plot 16 0  11 12 12 8 11 4 11 0    /
\setlinear \plot 10 12 10 8  11 4 /
\setlinear \plot 10 8  11 12 /
\put{$5{,}040  $} [c] at 13 -2
\endpicture
\end{minipage}
\begin{minipage}{4cm}
\beginpicture
\setcoordinatesystem units   <1.5mm,2mm>
\setplotarea x from 0 to 16, y from -2 to 15
\put{543)} [l] at 2 12
\put {$\scriptstyle  \bullet$} [c] at  10 0
\put {$\scriptstyle  \bullet$} [c] at  10 4
\put {$\scriptstyle  \bullet$} [c] at  11 12
\put {$\scriptstyle  \bullet$} [c] at  11 8
\put {$\scriptstyle  \bullet$} [c] at  11 0
\put {$\scriptstyle  \bullet$} [c] at  12 4
\put {$\scriptstyle  \bullet$} [c] at  16 12
\setlinear \plot 16 12  11 0 12 4 11 8 11 12    /
\setlinear \plot 10 0 10 4  11 8 /
\setlinear \plot 10 4  11 0 /

\put{$5{,}040 $} [c] at 13 -2
\endpicture
\end{minipage}
\begin{minipage}{4cm}
\beginpicture
\setcoordinatesystem units   <1.5mm,2mm>
\setplotarea x from 0 to 16, y from -2 to 15
\put{544)} [l] at 2 12
\put{$  $} [c] at 13 -2
\put {$ \scriptstyle \bullet$} [c] at  10 4
\put {$ \scriptstyle \bullet$} [c] at  12 0
\put {$ \scriptstyle \bullet$} [c] at  12 8
\put {$ \scriptstyle \bullet$} [c] at  12 12
\put {$ \scriptstyle \bullet$} [c] at  14 4
\put {$ \scriptstyle \bullet$} [c] at  16 12
\put {$ \scriptstyle \bullet$} [c] at  16 0
\setlinear \plot 12 12 12 8 10 4  12 0 14 4 12 8     /
\setlinear \plot 16 0  16 12 14 4   /
\put{$5{,}040   $} [c] at 13 -2
\endpicture
\end{minipage}
\begin{minipage}{4cm}
\beginpicture
\setcoordinatesystem units   <1.5mm,2mm>
\setplotarea x from 0 to 16, y from -2 to 15
\put{545)} [l] at 2 12
\put {$ \scriptstyle \bullet$} [c] at  10 8
\put {$ \scriptstyle \bullet$} [c] at  12 0
\put {$ \scriptstyle \bullet$} [c] at  12 4
\put {$ \scriptstyle \bullet$} [c] at  12 12
\put {$ \scriptstyle \bullet$} [c] at  14 8
\put {$ \scriptstyle \bullet$} [c] at  16 12
\put {$ \scriptstyle \bullet$} [c] at  16 0
\setlinear \plot 12 0 12 4 10 8  12 12 14 8 12 4     /
\setlinear \plot 16 12  16 0 14 8   /
\put{$5{,}040   $} [c] at 13 -2
\endpicture
\end{minipage}
\begin{minipage}{4cm}
\beginpicture
\setcoordinatesystem units   <1.5mm,2mm>
\setplotarea x from 0 to 16, y from -2 to 15
\put{546)} [l] at 2 12
\put {$ \scriptstyle \bullet$} [c] at  10 6
\put {$ \scriptstyle \bullet$} [c] at  10 12
\put {$ \scriptstyle \bullet$} [c] at  11.5 3
\put {$ \scriptstyle \bullet$} [c] at  13 0
\put {$ \scriptstyle \bullet$} [c] at  13 12
\put {$ \scriptstyle \bullet$} [c] at  16 6
\put {$ \scriptstyle \bullet$} [c] at  16 0
\setlinear \plot 10 6 10 12 16 0 13 12  16 6 13 0 10 6 13 12    /
\put{$5{,}040   $} [c] at 13 -2
\endpicture
\end{minipage}
$$

$$
\begin{minipage}{4cm}
\beginpicture
\setcoordinatesystem units   <1.5mm,2mm>
\setplotarea x from 0 to 16, y from -2 to 15
\put{547)} [l] at 2 12
\put {$ \scriptstyle \bullet$} [c] at  10 6
\put {$ \scriptstyle \bullet$} [c] at  10 0
\put {$ \scriptstyle \bullet$} [c] at  11.5 9
\put {$ \scriptstyle \bullet$} [c] at  13 0
\put {$ \scriptstyle \bullet$} [c] at  13 12
\put {$ \scriptstyle \bullet$} [c] at  16 6
\put {$ \scriptstyle \bullet$} [c] at  16 12
\setlinear \plot 10 6 10 0 16 12 13 0  16 6 13 12 10 6 13 0    /
\put{$5{,}040   $} [c] at 13 -2
\endpicture
\end{minipage}
\begin{minipage}{4cm}
\beginpicture
\setcoordinatesystem units   <1.5mm,2mm>
\setplotarea x from 0 to 16, y from -2 to 15
\put{548)} [l] at 2 12
\put {$\scriptstyle  \bullet$} [c] at  10 12
\put {$\scriptstyle  \bullet$} [c] at  10 4
\put {$\scriptstyle  \bullet$} [c] at  11 0
\put {$\scriptstyle  \bullet$} [c] at  11 8
\put {$\scriptstyle  \bullet$} [c] at  11 12
\put {$\scriptstyle  \bullet$} [c] at  12 4
\put {$\scriptstyle  \bullet$} [c] at  16 0
\setlinear \plot 10 12 10 4 11 8 12 4 11 0 10  4  /
\setlinear \plot 16 0 11 8 11 12 /
\put{$5{,}040  $} [c] at 13 -2
\endpicture
\end{minipage}
\begin{minipage}{4cm}
\beginpicture
\setcoordinatesystem units   <1.5mm,2mm>
\setplotarea x from 0 to 16, y from -2 to 15
\put{549)} [l] at 2 12
\put {$ \scriptstyle \bullet$} [c] at  10 0
\put {$ \scriptstyle \bullet$} [c] at  10 8
\put {$ \scriptstyle \bullet$} [c] at  11 0
\put {$ \scriptstyle \bullet$} [c] at  11 4
\put {$ \scriptstyle \bullet$} [c] at  11 12
\put {$ \scriptstyle \bullet$} [c] at  12 8
\put {$ \scriptstyle \bullet$} [c] at  16 12
\setlinear \plot 10 0 10 8 11 4 12 8 11 12 10  8  /
\setlinear \plot 16 12 11 4 11 0 /

\put{$5{,}040  $} [c] at 13 -2
\endpicture
\end{minipage}
\begin{minipage}{4cm}
\beginpicture
\setcoordinatesystem units   <1.5mm,2mm>
\setplotarea x from 0 to 16, y from -2 to 15
\put{550)} [l] at 2 12
\put {$ \scriptstyle \bullet$} [c] at  10 0
\put {$ \scriptstyle \bullet$} [c] at  10 6
\put {$ \scriptstyle \bullet$} [c] at  10 9
\put {$ \scriptstyle \bullet$} [c] at  10 12
\put {$ \scriptstyle \bullet$} [c] at  16 0
\put {$ \scriptstyle \bullet$} [c] at  16 6
\put {$ \scriptstyle \bullet$} [c] at  16 12
\setlinear \plot 16  0 16 12  10 6 10 12       /
\setlinear \plot 10 0 10 6  /
\put{$5{,}040   $} [c] at 13 -2
\endpicture
\end{minipage}
\begin{minipage}{4cm}
\beginpicture
\setcoordinatesystem units   <1.5mm,2mm>
\setplotarea x from 0 to 16, y from -2 to 15
\put{551)} [l] at 2 12
\put {$ \scriptstyle \bullet$} [c] at  10 0
\put {$ \scriptstyle \bullet$} [c] at  10 6
\put {$ \scriptstyle \bullet$} [c] at  10 3
\put {$ \scriptstyle \bullet$} [c] at  10 12
\put {$ \scriptstyle \bullet$} [c] at  16 0
\put {$ \scriptstyle \bullet$} [c] at  16 6
\put {$ \scriptstyle \bullet$} [c] at  16 12
\setlinear \plot 16  12 16 0  10 6 10 0       /
\setlinear \plot 10 12 10 6  /
\put{$5{,}040   $} [c] at 13 -2
\endpicture
\end{minipage}
\begin{minipage}{4cm}
\beginpicture
\setcoordinatesystem units   <1.5mm,2mm>
\setplotarea x from 0 to 16, y from -2 to 15
\put{552)} [l] at 2 12
\put {$ \scriptstyle \bullet$} [c] at  10 6
\put {$ \scriptstyle \bullet$} [c] at  10.5 3
\put {$ \scriptstyle \bullet$} [c] at  11 0
\put {$ \scriptstyle \bullet$} [c] at  11 12
\put {$ \scriptstyle \bullet$} [c] at  12 6
\put {$ \scriptstyle \bullet$} [c] at  16 12
\put {$ \scriptstyle \bullet$} [c] at  16 0
\setlinear \plot 16 0 16 12 11 0 12 6 11 12  10 6 11 0    /
\setlinear \plot 10 6 16 0 12 6  /
\put{$5{,}040  $} [c] at 13 -2
\endpicture
\end{minipage}
$$
$$
\begin{minipage}{4cm}
\beginpicture
\setcoordinatesystem units   <1.5mm,2mm>
\setplotarea x from 0 to 16, y from -2 to 15
\put{553)} [l] at 2 12
\put {$ \scriptstyle \bullet$} [c] at  10 6
\put {$ \scriptstyle \bullet$} [c] at  10.5 9
\put {$ \scriptstyle \bullet$} [c] at  11 0
\put {$ \scriptstyle \bullet$} [c] at  11 12
\put {$ \scriptstyle \bullet$} [c] at  12 6
\put {$ \scriptstyle \bullet$} [c] at  16 12
\put {$ \scriptstyle \bullet$} [c] at  16 0
\setlinear \plot 16 12 16 0 11 12 12 6 11 0  10 6 11 12    /
\setlinear \plot 10 6 16 12 12 6  /
\put{$5{,}040  $} [c] at 13 -2
\endpicture
\end{minipage}
\begin{minipage}{4cm}
\beginpicture
\setcoordinatesystem units   <1.5mm,2mm>
\setplotarea x from 0 to 16, y from -2 to 15
\put{554)} [l] at 2 12
\put {$ \scriptstyle \bullet$} [c] at  10 0
\put {$ \scriptstyle \bullet$} [c] at  10 9
\put {$ \scriptstyle \bullet$} [c] at  10 12
\put {$ \scriptstyle \bullet$} [c] at  12 6
\put {$ \scriptstyle \bullet$} [c] at  14 0
\put {$ \scriptstyle \bullet$} [c] at  14 12
\put {$ \scriptstyle \bullet$} [c] at  16 6
\setlinear \plot 10  0 10 12    /
\setlinear \plot 10 9 12 6 14 0 16 6 14 12 12 6  /
\put{$5{,}040   $} [c] at 13 -2
\endpicture
\end{minipage}
\begin{minipage}{4cm}
\beginpicture
\setcoordinatesystem units   <1.5mm,2mm>
\setplotarea x from 0 to 16, y from -2 to 15
\put{555)} [l] at 2 12
\put {$ \scriptstyle \bullet$} [c] at  10 0
\put {$ \scriptstyle \bullet$} [c] at  10 3
\put {$ \scriptstyle \bullet$} [c] at  10 12
\put {$ \scriptstyle \bullet$} [c] at  12 6
\put {$ \scriptstyle \bullet$} [c] at  14 0
\put {$ \scriptstyle \bullet$} [c] at  14 12
\put {$ \scriptstyle \bullet$} [c] at  16 6
\setlinear \plot 10  0 10 12    /
\setlinear \plot 10 3 12 6 14 12 16 6 14 0 12 6  /
\put{$5{,}040   $} [c] at 13 -2
\endpicture
\end{minipage}
\begin{minipage}{4cm}
\beginpicture
\setcoordinatesystem units   <1.5mm,2mm>
\setplotarea x from 0 to 16, y from -2 to 15
\put{556)} [l] at 2 12
\put {$ \scriptstyle \bullet$} [c] at  10 0
\put {$ \scriptstyle \bullet$} [c] at  10 6
\put {$ \scriptstyle \bullet$} [c] at  10 12
\put {$ \scriptstyle \bullet$} [c] at  13 0
\put {$ \scriptstyle \bullet$} [c] at  13 6
\put {$ \scriptstyle \bullet$} [c] at  13 12
\put {$ \scriptstyle \bullet$} [c] at  16 6
\setlinear \plot  13 0 10 6 13 12 13 0 16 6 13 12  /
\setlinear \plot  10 0 10 12  13 6   /
\put{$5{,}040 $} [c] at 13 -2
\endpicture
\end{minipage}
\begin{minipage}{4cm}
\beginpicture
\setcoordinatesystem units   <1.5mm,2mm>
\setplotarea x from 0 to 16, y from -2 to 15
\put{557)} [l] at 2 12
\put {$ \scriptstyle \bullet$} [c] at  10 0
\put {$ \scriptstyle \bullet$} [c] at  10 6
\put {$ \scriptstyle \bullet$} [c] at  10 12
\put {$ \scriptstyle \bullet$} [c] at  13 0
\put {$ \scriptstyle \bullet$} [c] at  13 6
\put {$ \scriptstyle \bullet$} [c] at  13 12
\put {$ \scriptstyle \bullet$} [c] at  16 6
\setlinear \plot  13 12 10 6 13 0 13 12 16 6 13 0  /
\setlinear \plot  10 12 10 0  13 6   /
\put{$5{,}040 $} [c] at 13 -2
\endpicture
\end{minipage}
\begin{minipage}{4cm}
\beginpicture
\setcoordinatesystem units   <1.5mm,2mm>
\setplotarea x from 0 to 16, y from -2 to 15
\put{558)} [l] at 2 12
\put {$ \scriptstyle \bullet$} [c] at  10 0
\put {$ \scriptstyle \bullet$} [c] at  10 12
\put {$ \scriptstyle \bullet$} [c] at  13 0
\put {$ \scriptstyle \bullet$} [c] at  13 4
\put {$ \scriptstyle \bullet$} [c] at  13 8
\put {$ \scriptstyle \bullet$} [c] at  13 12
\put {$ \scriptstyle \bullet$} [c] at  16 6
\setlinear \plot 13 0 13 12 16 6 13 0 /
\setlinear \plot 13 8 10 0 10 12 13 4 /
\put{$5{,}040   $} [c] at 13 -2
\endpicture
\end{minipage}
$$

$$
\begin{minipage}{4cm}
\beginpicture
\setcoordinatesystem units   <1.5mm,2mm>
\setplotarea x from 0 to 16, y from -2 to 15
\put{559)} [l] at 2 12
\put {$ \scriptstyle \bullet$} [c] at  10 0
\put {$ \scriptstyle \bullet$} [c] at  10  4
\put {$ \scriptstyle \bullet$} [c] at  10 8
\put {$ \scriptstyle \bullet$} [c] at  10 12
\put {$ \scriptstyle \bullet$} [c] at  16 0
\put {$ \scriptstyle \bullet$} [c] at  16 6
\put {$ \scriptstyle \bullet$} [c] at  16  12
\setlinear \plot 10  0 10 12  16 0  16 12 10 0    /
\put{$5{,}040$} [c] at 13 -2
\endpicture
\end{minipage}
\begin{minipage}{4cm}
\beginpicture
\setcoordinatesystem units   <1.5mm,2mm>
\setplotarea x from 0 to 16, y from -2 to 15
\put{560)} [l] at 2 12
\put {$ \scriptstyle \bullet$} [c] at  10 6
\put {$ \scriptstyle \bullet$} [c] at  10 12
\put {$ \scriptstyle \bullet$} [c] at  13 0
\put {$ \scriptstyle \bullet$} [c] at  13 6
\put {$ \scriptstyle \bullet$} [c] at  13 12
\put {$ \scriptstyle \bullet$} [c] at  16 0
\put {$ \scriptstyle \bullet$} [c] at  16 6
\setlinear \plot 16 0 16 6 13 12 10 6 10  12 16 0   /
\setlinear \plot 10 6 13 0 16 6   /
\put{$5{,}040 $} [c] at 13 -2
\endpicture
\end{minipage}
\begin{minipage}{4cm}
\beginpicture
\setcoordinatesystem units   <1.5mm,2mm>
\setplotarea x from 0 to 16, y from -2 to 15
\put{561)} [l] at 2 12
\put {$ \scriptstyle \bullet$} [c] at  10 0
\put {$ \scriptstyle \bullet$} [c] at  10 4
\put {$ \scriptstyle \bullet$} [c] at  10 8
\put {$ \scriptstyle \bullet$} [c] at  10 12
\put {$ \scriptstyle \bullet$} [c] at  13 12
\put {$ \scriptstyle \bullet$} [c] at  16 0
\put {$ \scriptstyle \bullet$} [c] at  16 12
\setlinear \plot 10 12 10 0 13 12 16 0 16 12 10 4 /
\setlinear \plot  10 8 16 0  /
\put{$5{,}040$} [c] at 13 -2
\endpicture
\end{minipage}
\begin{minipage}{4cm}
\beginpicture
\setcoordinatesystem units   <1.5mm,2mm>
\setplotarea x from 0 to 16, y from -2 to 15
\put{562)} [l] at 2 12
\put {$ \scriptstyle \bullet$} [c] at  10 0
\put {$ \scriptstyle \bullet$} [c] at  10 4
\put {$ \scriptstyle \bullet$} [c] at  10 8
\put {$ \scriptstyle \bullet$} [c] at  10 12
\put {$ \scriptstyle \bullet$} [c] at  13 0
\put {$ \scriptstyle \bullet$} [c] at  16 0
\put {$ \scriptstyle \bullet$} [c] at  16 12
\setlinear \plot 10 0 10 12 13 0 16 12 16 0 10 8 /
\setlinear \plot  10 4 16 12  /
\put{$5{,}040$} [c] at 13 -2
\endpicture
\end{minipage}
\begin{minipage}{4cm}
\beginpicture
\setcoordinatesystem units   <1.5mm,2mm>
\setplotarea x from 0 to 16, y from -2 to 15
\put{563)} [l] at 2 12
\put {$ \scriptstyle \bullet$} [c] at  10 0
\put {$ \scriptstyle \bullet$} [c] at  10 6
\put {$ \scriptstyle \bullet$} [c] at  10 12
\put {$ \scriptstyle \bullet$} [c] at  13 0
\put {$ \scriptstyle \bullet$} [c] at  13 12
\put {$ \scriptstyle \bullet$} [c] at  16 6
\put {$ \scriptstyle \bullet$} [c] at  16 12
\setlinear \plot  16 12 16 6  13  0 10 6 13 12 16 6  /
\setlinear \plot 10 12 10 0    /
\put{$5{,}040$} [c] at 13 -2
\endpicture
\end{minipage}
\begin{minipage}{4cm}
\beginpicture
\setcoordinatesystem units   <1.5mm,2mm>
\setplotarea x from 0 to 16, y from -2 to 15
\put{564)} [l] at 2 12
\put {$ \scriptstyle \bullet$} [c] at  10 0
\put {$ \scriptstyle \bullet$} [c] at  10 6
\put {$ \scriptstyle \bullet$} [c] at  10 12
\put {$ \scriptstyle \bullet$} [c] at  13 0
\put {$ \scriptstyle \bullet$} [c] at  13 12
\put {$ \scriptstyle \bullet$} [c] at  16 6
\put {$ \scriptstyle \bullet$} [c] at  16 0
\setlinear \plot  16 0 16 6  13  0 10 6 13 12 16 6  /
\setlinear \plot 10 12 10 0    /
\put{$5{,}040$} [c] at 13 -2
\endpicture
\end{minipage}
$$
$$
\begin{minipage}{4cm}
\beginpicture
\setcoordinatesystem units   <1.5mm,2mm>
\setplotarea x from 0 to 16, y from -2 to 15
\put{565)} [l] at 2 12
\put {$ \scriptstyle \bullet$} [c] at  10 0
\put {$ \scriptstyle \bullet$} [c] at  10 4
\put {$ \scriptstyle \bullet$} [c] at  10 8
\put {$ \scriptstyle \bullet$} [c] at  10 12
\put {$ \scriptstyle \bullet$} [c] at  13 12
\put {$ \scriptstyle \bullet$} [c] at  16 0
\put {$ \scriptstyle \bullet$} [c] at  16 12
\setlinear \plot 10 12 10 0 13 12 16 0 16 12  /
\setlinear \plot  10 8 16 12  /
\put{$5{,}040$} [c] at 13 -2
\endpicture
\end{minipage}
\begin{minipage}{4cm}
\beginpicture
\setcoordinatesystem units   <1.5mm,2mm>
\setplotarea x from 0 to 16, y from -2 to 15
\put{566)} [l] at 2 12
\put {$ \scriptstyle \bullet$} [c] at  10 0
\put {$ \scriptstyle \bullet$} [c] at  10 4
\put {$ \scriptstyle \bullet$} [c] at  10 8
\put {$ \scriptstyle \bullet$} [c] at  10 12
\put {$ \scriptstyle \bullet$} [c] at  13 0
\put {$ \scriptstyle \bullet$} [c] at  16 0
\put {$ \scriptstyle \bullet$} [c] at  16 12
\setlinear \plot 10 0 10 12 13 0 16 12 16 0  /
\setlinear \plot  10 4 16 0  /
\put{$5{,}040$} [c] at 13 -2
\endpicture
\end{minipage}
\begin{minipage}{4cm}
\beginpicture
\setcoordinatesystem units   <1.5mm,2mm>
\setplotarea x from 0 to 16, y from -2 to 15
\put{567)} [l] at 2 12
\put {$ \scriptstyle \bullet$} [c] at  10 12
\put {$ \scriptstyle \bullet$} [c] at  12 0
\put {$ \scriptstyle \bullet$} [c] at  13 9
\put {$ \scriptstyle \bullet$} [c] at  12 12
\put {$ \scriptstyle \bullet$} [c] at  14 6
\put {$ \scriptstyle \bullet$} [c] at  16 0
\put {$ \scriptstyle \bullet$} [c] at  16 12
\setlinear \plot 10 12 12 0 14 6 16 12  /
\setlinear \plot  16 0 14 6 12 12   /
\put{$5{,}040$} [c] at 13 -2
\endpicture
\end{minipage}
\begin{minipage}{4cm}
\beginpicture
\setcoordinatesystem units   <1.5mm,2mm>
\setplotarea x from 0 to 16, y from -2 to 15
\put{568)} [l] at 2 12
\put {$ \scriptstyle \bullet$} [c] at  10 0
\put {$ \scriptstyle \bullet$} [c] at  12 0
\put {$ \scriptstyle \bullet$} [c] at  13 3
\put {$ \scriptstyle \bullet$} [c] at  12 12
\put {$ \scriptstyle \bullet$} [c] at  14 6
\put {$ \scriptstyle \bullet$} [c] at  16 0
\put {$ \scriptstyle \bullet$} [c] at  16 12
\setlinear \plot 10 0 12 12 14 6 16 0  /
\setlinear \plot  16 12 14 6 12 0   /
\put{$5{,}040$} [c] at 13 -2
\endpicture
\end{minipage}
\begin{minipage}{4cm}
\beginpicture
\setcoordinatesystem units   <1.5mm,2mm>
\setplotarea x from 0 to 16, y from -2 to 15
\put{569)} [l] at 2 12
\put {$ \scriptstyle \bullet$} [c] at  10 12
\put {$ \scriptstyle \bullet$} [c] at  13 0
\put {$ \scriptstyle \bullet$} [c] at  13 4
\put {$ \scriptstyle \bullet$} [c] at  13 8
\put {$ \scriptstyle \bullet$} [c] at  13 12
\put {$ \scriptstyle \bullet$} [c] at  16 0
\put {$ \scriptstyle \bullet$} [c] at  16 12
\setlinear \plot  13 12 13   0 /
\setlinear \plot   10 12 13 4 16 12 16  0 13 8  /
\put{$5{,}040$} [c] at 13 -2
\endpicture
\end{minipage}
\begin{minipage}{4cm}
\beginpicture
\setcoordinatesystem units   <1.5mm,2mm>
\setplotarea x from 0 to 16, y from -2 to 15
\put{570)} [l] at 2 12
\put {$ \scriptstyle \bullet$} [c] at  10 0
\put {$ \scriptstyle \bullet$} [c] at  13 0
\put {$ \scriptstyle \bullet$} [c] at  13 4
\put {$ \scriptstyle \bullet$} [c] at  13 8
\put {$ \scriptstyle \bullet$} [c] at  13 12
\put {$ \scriptstyle \bullet$} [c] at  16 0
\put {$ \scriptstyle \bullet$} [c] at  16 12
\setlinear \plot  13 12 13   0 /
\setlinear \plot   10 0 13 8 16 0 16  12 13 4  /
\put{$5{,}040$} [c] at 13 -2
\endpicture
\end{minipage}
$$

$$
\begin{minipage}{4cm}
\beginpicture
\setcoordinatesystem units   <1.5mm,2mm>
\setplotarea x from 0 to 16, y from -2 to 15
\put{571)} [l] at 2 12
\put {$ \scriptstyle \bullet$} [c] at 10 6
\put {$ \scriptstyle \bullet$} [c] at 13 0
\put {$ \scriptstyle \bullet$} [c] at 13 3
\put {$ \scriptstyle \bullet$} [c] at 13 6
\put {$ \scriptstyle \bullet$} [c] at 13 9
\put {$ \scriptstyle \bullet$} [c] at 13 12
\put {$ \scriptstyle \bullet$} [c] at 16 6
\setlinear \plot 13  0 10 6 13 12 16 6 13 0 13 12   /
\put{$2{,}520$} [c] at 13 -2
\endpicture
\end{minipage}
\begin{minipage}{4cm}
\beginpicture
\setcoordinatesystem units   <1.5mm,2mm>
\setplotarea x from 0 to 16, y from -2 to 15
\put{572)} [l] at 2 12
\put {$ \scriptstyle \bullet$} [c] at 10 9
\put {$ \scriptstyle \bullet$} [c] at 11 0
\put {$ \scriptstyle \bullet$} [c] at 11 6
\put {$ \scriptstyle \bullet$} [c] at 11 12
\put {$ \scriptstyle \bullet$} [c] at 12 9
\put {$ \scriptstyle \bullet$} [c] at 16 6
\put {$ \scriptstyle \bullet$} [c] at 16 12
\setlinear \plot 11 6 16 12 16 6 11 0 11 6 10 9 11 12 12 9 11 6  /
\put{$2{,}520$} [c] at 13 -2
\endpicture
\end{minipage}
\begin{minipage}{4cm}
\beginpicture
\setcoordinatesystem units   <1.5mm,2mm>
\setplotarea x from 0 to 16, y from -2 to 15
\put{573)} [l] at 2 12
\put {$ \scriptstyle \bullet$} [c] at 10 3
\put {$ \scriptstyle \bullet$} [c] at 11 0
\put {$ \scriptstyle \bullet$} [c] at 11 6
\put {$ \scriptstyle \bullet$} [c] at 11 12
\put {$ \scriptstyle \bullet$} [c] at 12 3
\put {$ \scriptstyle \bullet$} [c] at 16 6
\put {$ \scriptstyle \bullet$} [c] at 16 0
\setlinear \plot 11 6 16 0 16 6 11 12 11 6 10 3 11 0 12 3 11 6  /
\put{$2{,}520$} [c] at 13 -2
\endpicture
\end{minipage}
\begin{minipage}{4cm}
\beginpicture
\setcoordinatesystem units   <1.5mm,2mm>
\setplotarea x from 0 to 16, y from -2 to 15
\put{574)} [l] at 2 12
\put {$ \scriptstyle \bullet$} [c] at 10 9
\put {$ \scriptstyle \bullet$} [c] at 11 0
\put {$ \scriptstyle \bullet$} [c] at 11 6
\put {$ \scriptstyle \bullet$} [c] at 11 12
\put {$ \scriptstyle \bullet$} [c] at 12 9
\put {$ \scriptstyle \bullet$} [c] at 16 6
\put {$ \scriptstyle \bullet$} [c] at 16 12
\setlinear \plot 11 6 11 0 16 6 16 12  /
\setlinear \plot 16 6 11 12 10 9 11 6 12 9 11 12  /
\put{$2{,}520$} [c] at 13 -2
\endpicture
\end{minipage}
\begin{minipage}{4cm}
\beginpicture
\setcoordinatesystem units   <1.5mm,2mm>
\setplotarea x from 0 to 16, y from -2 to 15
\put{575)} [l] at 2 12
\put {$ \scriptstyle \bullet$} [c] at 10 3
\put {$ \scriptstyle \bullet$} [c] at 11 0
\put {$ \scriptstyle \bullet$} [c] at 11 6
\put {$ \scriptstyle \bullet$} [c] at 11 12
\put {$ \scriptstyle \bullet$} [c] at 12 3
\put {$ \scriptstyle \bullet$} [c] at 16 6
\put {$ \scriptstyle \bullet$} [c] at 16 0
\setlinear \plot 11 6 11 12 16 6 16 0  /
\setlinear \plot 16 6 11 0 10 3 11 6 12 3 11 0  /
\put{$2{,}520$} [c] at 13 -2
\endpicture
\end{minipage}
\begin{minipage}{4cm}
\beginpicture
\setcoordinatesystem units   <1.5mm,2mm>
\setplotarea x from 0 to 16, y from -2 to 15
\put{576)} [l] at 2 12
\put {$ \scriptstyle \bullet$} [c] at 10 4
\put {$ \scriptstyle \bullet$} [c] at 10 12
\put {$ \scriptstyle \bullet$} [c] at 13 0
\put {$ \scriptstyle \bullet$} [c] at 13 4
\put {$ \scriptstyle \bullet$} [c] at 13 8
\put {$ \scriptstyle \bullet$} [c] at 16 4
\put {$ \scriptstyle \bullet$} [c] at 16 12
\setlinear \plot  13 0 10 4 10 12 13 8 13 0 16 4 16 12 13 8    /
\put{$2{,}520$} [c] at 13 -2
\endpicture
\end{minipage}
$$
$$
\begin{minipage}{4cm}
\beginpicture
\setcoordinatesystem units   <1.5mm,2mm>
\setplotarea x from 0 to 16, y from -2 to 15
\put{577)} [l] at 2 12
\put {$ \scriptstyle \bullet$} [c] at 10 0
\put {$ \scriptstyle \bullet$} [c] at 10 8
\put {$ \scriptstyle \bullet$} [c] at 13 12
\put {$ \scriptstyle \bullet$} [c] at 13 4
\put {$ \scriptstyle \bullet$} [c] at 13 8
\put {$ \scriptstyle \bullet$} [c] at 16 8
\put {$ \scriptstyle \bullet$} [c] at 16 0
\setlinear \plot  13 12 10 8 10 0 13 4 13 12 16 8 16 0 13 4    /
\put{$2{,}520$} [c] at 13 -2
\endpicture
\end{minipage}
\begin{minipage}{4cm}
\beginpicture
\setcoordinatesystem units   <1.5mm,2mm>
\setplotarea x from 0 to 16, y from -2 to 15
\put{578)} [l] at 2 12
\put {$ \scriptstyle \bullet$} [c] at 10 6
\put {$ \scriptstyle \bullet$} [c] at 10 12
\put {$ \scriptstyle \bullet$} [c] at 13 0
\put {$ \scriptstyle \bullet$} [c] at 13 6
\put {$ \scriptstyle \bullet$} [c] at 13 12
\put {$ \scriptstyle \bullet$} [c] at 14.5 9
\put {$ \scriptstyle \bullet$} [c] at 16 6
\setlinear \plot  13 0 13 12 16 6 13 0 10 6 13 12    /
\setlinear \plot  10 6 10 12 13 6    /
\put{$2{,}520$} [c] at 13 -2
\endpicture
\end{minipage}
\begin{minipage}{4cm}
\beginpicture
\setcoordinatesystem units   <1.5mm,2mm>
\setplotarea x from 0 to 16, y from -2 to 15
\put{579)} [l] at 2 12
\put {$ \scriptstyle \bullet$} [c] at 10 6
\put {$ \scriptstyle \bullet$} [c] at 10 0
\put {$ \scriptstyle \bullet$} [c] at 13 0
\put {$ \scriptstyle \bullet$} [c] at 13 6
\put {$ \scriptstyle \bullet$} [c] at 13 12
\put {$ \scriptstyle \bullet$} [c] at 14.5 3
\put {$ \scriptstyle \bullet$} [c] at 16 6
\setlinear \plot  13 12 13 0 16 6 13 12 10 6 13 0    /
\setlinear \plot  10 6 10 0 13 6    /
\put{$2{,}520$} [c] at 13 -2
\endpicture
\end{minipage}
\begin{minipage}{4cm}
\beginpicture
\setcoordinatesystem units   <1.5mm,2mm>
\setplotarea x from 0 to 16, y from -2 to 15
\put{580)} [l] at 2 12
\put {$ \scriptstyle \bullet$} [c] at 10 4
\put {$ \scriptstyle \bullet$} [c] at 10 12
\put {$ \scriptstyle \bullet$} [c] at 13 4
\put {$ \scriptstyle \bullet$} [c] at 14.5 0
\put {$ \scriptstyle \bullet$} [c] at 14.5 8
\put {$ \scriptstyle \bullet$} [c] at 14.5 12
\put {$ \scriptstyle \bullet$} [c] at 16 4
\setlinear \plot 10 12 10 4 14.5 0 13  4 14.5 8 14.5 12 10 4 /
\setlinear \plot 14.5 0 16 4 14.5 8  /
\put{$2{,}520$} [c] at 13 -2
\endpicture
\end{minipage}
\begin{minipage}{4cm}
\beginpicture
\setcoordinatesystem units   <1.5mm,2mm>
\setplotarea x from 0 to 16, y from -2 to 15
\put{581)} [l] at 2 12
\put {$ \scriptstyle \bullet$} [c] at 10 8
\put {$ \scriptstyle \bullet$} [c] at 10 0
\put {$ \scriptstyle \bullet$} [c] at 13 8
\put {$ \scriptstyle \bullet$} [c] at 14.5 12
\put {$ \scriptstyle \bullet$} [c] at 14.5 4
\put {$ \scriptstyle \bullet$} [c] at 14.5 0
\put {$ \scriptstyle \bullet$} [c] at 16 8
\setlinear \plot 10 0 10 8 14.5 12 13  8 14.5 4 14.5 0 10 8 /
\setlinear \plot 14.5 12 16 8 14.5 4  /
\put{$2{,}520$} [c] at 13 -2
\endpicture
\end{minipage}
\begin{minipage}{4cm}
\beginpicture
\setcoordinatesystem units   <1.5mm,2mm>
\setplotarea x from 0 to 16, y from -2 to 15
\put{582)} [l] at 2 12
\put {$ \scriptstyle \bullet$} [c] at 10 12
\put {$ \scriptstyle \bullet$} [c] at 13  0
\put {$ \scriptstyle \bullet$} [c] at 13 3
\put {$ \scriptstyle \bullet$} [c] at 13 6
\put {$ \scriptstyle \bullet$} [c] at 13 9
\put {$ \scriptstyle \bullet$} [c] at 13 12
\put {$ \scriptstyle \bullet$} [c] at 16 12
\setlinear \plot  13  12 13 0     /
\setlinear \plot  10  12 13 3  16 12    /
\put{$2{,}520$} [c] at 13 -2
\endpicture
\end{minipage}
$$

$$
\begin{minipage}{4cm}
\beginpicture
\setcoordinatesystem units   <1.5mm,2mm>
\setplotarea x from 0 to 16, y from -2 to 15
\put{583)} [l] at 2 12
\put {$ \scriptstyle \bullet$} [c] at 10 0
\put {$ \scriptstyle \bullet$} [c] at 13  0
\put {$ \scriptstyle \bullet$} [c] at 13 3
\put {$ \scriptstyle \bullet$} [c] at 13 6
\put {$ \scriptstyle \bullet$} [c] at 13 9
\put {$ \scriptstyle \bullet$} [c] at 13 12
\put {$ \scriptstyle \bullet$} [c] at 16 0
\setlinear \plot  13  12 13 0     /
\setlinear \plot  10  0 13 9  16 0    /
\put{$2{,}520$} [c] at 13 -2
\endpicture
\end{minipage}
\begin{minipage}{4cm}
\beginpicture
\setcoordinatesystem units   <1.5mm,2mm>
\setplotarea x from 0 to 16, y from -2 to 15
\put{584)} [l] at 2 12
\put {$ \scriptstyle \bullet$} [c] at 10 4
\put {$ \scriptstyle \bullet$} [c] at 10 12
\put {$ \scriptstyle \bullet$} [c] at 13 0
\put {$ \scriptstyle \bullet$} [c] at 13 8
\put {$ \scriptstyle \bullet$} [c] at 13 12
\put {$ \scriptstyle \bullet$} [c] at 16 4
\put {$ \scriptstyle \bullet$} [c] at 16 12
\setlinear \plot  10 12 10 4 13 8 13 12     /
\setlinear \plot  10 4 13  0  16 4 16 12    /
\setlinear \plot  13  8  16 4    /
\put{$2{,}520$} [c] at 13 -2
\endpicture
\end{minipage}
\begin{minipage}{4cm}
\beginpicture
\setcoordinatesystem units   <1.5mm,2mm>
\setplotarea x from 0 to 16, y from -2 to 15
\put{585)} [l] at 2 12
\put {$ \scriptstyle \bullet$} [c] at 10 8
\put {$ \scriptstyle \bullet$} [c] at 10 0
\put {$ \scriptstyle \bullet$} [c] at 13 0
\put {$ \scriptstyle \bullet$} [c] at 13 4
\put {$ \scriptstyle \bullet$} [c] at 13 12
\put {$ \scriptstyle \bullet$} [c] at 16 8
\put {$ \scriptstyle \bullet$} [c] at 16 0
\setlinear \plot  10 0 10 8 13 4 13 0     /
\setlinear \plot  10 8 13 12  16 8 16 0    /
\setlinear \plot  13  4  16 8    /
\put{$2{,}520$} [c] at 13 -2
\endpicture
\end{minipage}
\begin{minipage}{4cm}
\beginpicture
\setcoordinatesystem units   <1.5mm,2mm>
\setplotarea x from 0 to 16, y from -2 to 15
\put{586)} [l] at 2 12
\put {$ \scriptstyle \bullet$} [c] at 10 6
\put {$ \scriptstyle \bullet$} [c] at 10 12
\put {$ \scriptstyle \bullet$} [c] at 12 0
\put {$ \scriptstyle \bullet$} [c] at 14 6
\put {$ \scriptstyle \bullet$} [c] at 14 12
\put {$ \scriptstyle \bullet$} [c] at 15 9
\put {$ \scriptstyle \bullet$} [c] at 16 12
\setlinear \plot  16 12 14 6 12 0 10 6 10 12 14 6 14 12 10 6   /
\put{$2{,}520$} [c] at 13 -2
\endpicture
\end{minipage}
\begin{minipage}{4cm}
\beginpicture
\setcoordinatesystem units   <1.5mm,2mm>
\setplotarea x from 0 to 16, y from -2 to 15
\put{587)} [l] at 2 12
\put {$ \scriptstyle \bullet$} [c] at 10 6
\put {$ \scriptstyle \bullet$} [c] at 10 0
\put {$ \scriptstyle \bullet$} [c] at 12 12
\put {$ \scriptstyle \bullet$} [c] at 14 6
\put {$ \scriptstyle \bullet$} [c] at 14 0
\put {$ \scriptstyle \bullet$} [c] at 15 3
\put {$ \scriptstyle \bullet$} [c] at 16 0
\setlinear \plot  16 0 14 6 12 12 10 6 10 0 14 6 14 0 10 6   /
\put{$2{,}520$} [c] at 13 -2
\endpicture
\end{minipage}
\begin{minipage}{4cm}
\beginpicture
\setcoordinatesystem units   <1.5mm,2mm>
\setplotarea x from 0 to 16, y from -2 to 15
\put{588)} [l] at 2 12
\put {$ \scriptstyle \bullet$} [c] at 10 4
\put {$ \scriptstyle \bullet$} [c] at 10 12
\put {$ \scriptstyle \bullet$} [c] at 12 0
\put {$ \scriptstyle \bullet$} [c] at 12 12
\put {$ \scriptstyle \bullet$} [c] at 14 4
\put {$ \scriptstyle \bullet$} [c] at 14 8
\put {$ \scriptstyle \bullet$} [c] at 16 12
\setlinear \plot 16 12 14 8 14 4 12 0 10 4 10 12 14 4  /
\setlinear \plot 12 12 14 8  /
\put{$2{,}520$} [c] at 13 -2
\endpicture
\end{minipage}
$$
$$
\begin{minipage}{4cm}
\beginpicture
\setcoordinatesystem units   <1.5mm,2mm>
\setplotarea x from 0 to 16, y from -2 to 15
\put{589)} [l] at 2 12
\put {$ \scriptstyle \bullet$} [c] at 10 8
\put {$ \scriptstyle \bullet$} [c] at 10 0
\put {$ \scriptstyle \bullet$} [c] at 12 0
\put {$ \scriptstyle \bullet$} [c] at 12 12
\put {$ \scriptstyle \bullet$} [c] at 14 8
\put {$ \scriptstyle \bullet$} [c] at 14 4
\put {$ \scriptstyle \bullet$} [c] at 16 0
\setlinear \plot 16 0 14 4 14 8 12 12 10 8 10 0 14 8  /
\setlinear \plot 12 0 14 4  /
\put{$2{,}520$} [c] at 13 -2
\endpicture
\end{minipage}
\begin{minipage}{4cm}
\beginpicture
\setcoordinatesystem units   <1.5mm,2mm>
\setplotarea x from 0 to 16, y from -2 to 15
\put{590)} [l] at 2 12
\put {$ \scriptstyle \bullet$} [c] at 10 6
\put {$ \scriptstyle \bullet$} [c] at 10 12
\put {$ \scriptstyle \bullet$} [c] at 13 0
\put {$ \scriptstyle \bullet$} [c] at 13 6
\put {$ \scriptstyle \bullet$} [c] at 13 12
\put {$ \scriptstyle \bullet$} [c] at 16 6
\put {$ \scriptstyle \bullet$} [c] at 16 12
\setlinear \plot  13  0  10 6  10 12 13 6 13 0 16 6 16 12 13 6 13 12 /
\setlinear \plot  10 6 13 12 16 6    /
\put{$2{,}520$} [c] at 13 -2
\endpicture
\end{minipage}
\begin{minipage}{4cm}
\beginpicture
\setcoordinatesystem units   <1.5mm,2mm>
\setplotarea x from 0 to 16, y from -2 to 15
\put{591)} [l] at 2 12
\put {$ \scriptstyle \bullet$} [c] at 10 6
\put {$ \scriptstyle \bullet$} [c] at 10 0
\put {$ \scriptstyle \bullet$} [c] at 13 0
\put {$ \scriptstyle \bullet$} [c] at 13 6
\put {$ \scriptstyle \bullet$} [c] at 13 12
\put {$ \scriptstyle \bullet$} [c] at 16 6
\put {$ \scriptstyle \bullet$} [c] at 16 0
\setlinear \plot  13  12  10 6  10 0 13 6 13 12 16 6 16 0 13 6  13 0 /
\setlinear \plot  10 6 13 0 16 6    /
\put{$2{,}520$} [c] at 13 -2
\endpicture
\end{minipage}
\begin{minipage}{4cm}
\beginpicture
\setcoordinatesystem units   <1.5mm,2mm>
\setplotarea x from 0 to 16, y from -2 to 15
\put{592)} [l] at 2 12
\put {$ \scriptstyle \bullet$} [c] at  10 4
\put {$ \scriptstyle \bullet$} [c] at  11 0
\put {$ \scriptstyle \bullet$} [c] at  11 8
\put {$ \scriptstyle \bullet$} [c] at  11 12
\put {$ \scriptstyle \bullet$} [c] at  12 4
\put {$ \scriptstyle \bullet$} [c] at  16 0
\put {$ \scriptstyle \bullet$} [c] at  16 12
\setlinear \plot  11 12 11 8 16 0  16 12 11 0 12 4 11 8 10 4 11 0   /
\put{$2{,}520$} [c] at 13 -2
\endpicture
\end{minipage}
\begin{minipage}{4cm}
\beginpicture
\setcoordinatesystem units   <1.5mm,2mm>
\setplotarea x from 0 to 16, y from -2 to 15
\put{593)} [l] at 2 12
\put {$ \scriptstyle \bullet$} [c] at  10 8
\put {$ \scriptstyle \bullet$} [c] at  11 0
\put {$ \scriptstyle \bullet$} [c] at  11 4
\put {$ \scriptstyle \bullet$} [c] at  11 12
\put {$ \scriptstyle \bullet$} [c] at  12 8
\put {$ \scriptstyle \bullet$} [c] at  16 0
\put {$ \scriptstyle \bullet$} [c] at  16 12
\setlinear \plot  11 0 11 4 16 12  16 0 11 12 12 8 11 4 10 8 11 12   /
\put{$2{,}520$} [c] at 13 -2
\endpicture
\end{minipage}
\begin{minipage}{4cm}
\beginpicture
\setcoordinatesystem units   <1.5mm,2mm>
\setplotarea x from 0 to 16, y from -2 to 15
\put{594)} [l] at 2 12
\put {$ \scriptstyle \bullet$} [c] at  10 0
\put {$ \scriptstyle \bullet$} [c] at  12 12
\put {$ \scriptstyle \bullet$} [c] at  12 8
\put {$ \scriptstyle \bullet$} [c] at  12 4
\put {$ \scriptstyle \bullet$} [c] at  14 0
\put {$ \scriptstyle \bullet$} [c] at  16 4
\put {$ \scriptstyle \bullet$} [c] at  16 12
\setlinear \plot   10 0 12 12 12 4  16 12 16 4 12 8    /
\setlinear \plot  12 4 14 0  16 4  /
\put{$2{,}520   $} [c] at 13 -2
\endpicture
\end{minipage}
$$

$$
\begin{minipage}{4cm}
\beginpicture
\setcoordinatesystem units   <1.5mm,2mm>
\setplotarea x from 0 to 16, y from -2 to 15
\put{595)} [l] at 2 12
\put {$ \scriptstyle \bullet$} [c] at  10 12
\put {$ \scriptstyle \bullet$} [c] at  12 0
\put {$ \scriptstyle \bullet$} [c] at  12 8
\put {$ \scriptstyle \bullet$} [c] at  12 4
\put {$ \scriptstyle \bullet$} [c] at  14 12
\put {$ \scriptstyle \bullet$} [c] at  16 8
\put {$ \scriptstyle \bullet$} [c] at  16 0
\setlinear \plot   10 12 12 0 12 8  16 0 16 8 12 4    /
\setlinear \plot  12 8 14 12  16 8  /
\put{$2{,}520   $} [c] at 13 -2
\endpicture
\end{minipage}
\begin{minipage}{4cm}
\beginpicture
\setcoordinatesystem units   <1.5mm,2mm>
\setplotarea x from 0 to 16, y from -2 to 15
\put{596)} [l] at 2 12
\put {$ \scriptstyle \bullet$} [c] at  10 6
\put {$ \scriptstyle \bullet$} [c] at  10 12
\put {$ \scriptstyle \bullet$} [c] at  12 0
\put {$ \scriptstyle \bullet$} [c] at  12 6
\put {$ \scriptstyle \bullet$} [c] at  12 12
\put {$ \scriptstyle \bullet$} [c] at  16 0
\put {$ \scriptstyle \bullet$} [c] at  16 6
\setlinear \plot   10 6  12 12 12 0 16 6  12 12  /
\setlinear \plot  16 6 16 0 10 12 10 6 12  0   /
\setlinear \plot  10 12 12 6 /
\put{$2{,}520 $} [c] at 13 -2
\endpicture
\end{minipage}
\begin{minipage}{4cm}
\beginpicture
\setcoordinatesystem units   <1.5mm,2mm>
\setplotarea x from 0 to 16, y from -2 to 15
\put{597)} [l] at 2 12
\put {$ \scriptstyle \bullet$} [c] at  10 6
\put {$ \scriptstyle \bullet$} [c] at  10 0
\put {$ \scriptstyle \bullet$} [c] at  12 0
\put {$ \scriptstyle \bullet$} [c] at  12 6
\put {$ \scriptstyle \bullet$} [c] at  12 12
\put {$ \scriptstyle \bullet$} [c] at  16 12
\put {$ \scriptstyle \bullet$} [c] at  16 6
\setlinear \plot   10 6  12 0 12 12 16 6  12 0  /
\setlinear \plot  16 6 16 12 10 0 10 6 12  12   /
\setlinear \plot  10 0 12 6 /
\put{$2{,}520$} [c] at 13 -2
\endpicture
\end{minipage}
\begin{minipage}{4cm}
\beginpicture
\setcoordinatesystem units   <1.5mm,2mm>
\setplotarea x from 0 to 16, y from -2 to 15
\put{598)} [l] at 2 12
\put {$ \scriptstyle \bullet$} [c] at  10 0
\put {$ \scriptstyle \bullet$} [c] at  10 6
\put {$ \scriptstyle \bullet$} [c] at  10 12
\put {$ \scriptstyle \bullet$} [c] at  13 6
\put {$ \scriptstyle \bullet$} [c] at  16 0
\put {$ \scriptstyle \bullet$} [c] at  16 6
\put {$ \scriptstyle \bullet$} [c] at  16 12
\setlinear \plot  10 0  10 12  16 0 16 12  10 0  /
\put{$2{,}520  $} [c] at 13 -2
\endpicture
\end{minipage}
\begin{minipage}{4cm}
\beginpicture
\setcoordinatesystem units   <1.5mm,2mm>
\setplotarea x from 0 to 16, y from -2 to 15
\put{599)} [l] at 2 12
\put {$ \scriptstyle \bullet$} [c] at  10 0
\put {$ \scriptstyle \bullet$} [c] at  10 6
\put {$ \scriptstyle \bullet$} [c] at  10 12
\put {$ \scriptstyle \bullet$} [c] at  12 0
\put {$ \scriptstyle \bullet$} [c] at  12 6
\put {$ \scriptstyle \bullet$} [c] at  12 12
\put {$ \scriptstyle \bullet$} [c] at  16 12
\setlinear \plot 12 6 10 0 10 12  12 0 16 12  10 6 12 12 12 0 /
\put{$2{,}520$} [c] at 13 -2
\endpicture
\end{minipage}
\begin{minipage}{4cm}
\beginpicture
\setcoordinatesystem units   <1.5mm,2mm>
\setplotarea x from 0 to 16, y from -2 to 15
\put{600)} [l] at 2 12
\put {$ \scriptstyle \bullet$} [c] at  10 0
\put {$ \scriptstyle \bullet$} [c] at  10 6
\put {$ \scriptstyle \bullet$} [c] at  10 12
\put {$ \scriptstyle \bullet$} [c] at  12 0
\put {$ \scriptstyle \bullet$} [c] at  12 6
\put {$ \scriptstyle \bullet$} [c] at  12 12
\put {$ \scriptstyle \bullet$} [c] at  16 0
\setlinear \plot 12 6 10 12 10 0  12 12 16 0  10 6 12 0 12 12 /
\put{$2{,}520$} [c] at 13 -2
\endpicture
\end{minipage}
$$
$$
\begin{minipage}{4cm}
\beginpicture
\setcoordinatesystem units   <1.5mm,2mm>
\setplotarea x from 0 to 16, y from -2 to 15
\put{601)} [l] at 2 12
\put {$ \scriptstyle \bullet$} [c] at  10 0
\put {$ \scriptstyle \bullet$} [c] at  10 4
\put {$ \scriptstyle \bullet$} [c] at  10 8
\put {$ \scriptstyle \bullet$} [c] at  10 12
\put {$ \scriptstyle \bullet$} [c] at  13 12
\put {$ \scriptstyle \bullet$} [c] at  16 0
\put {$ \scriptstyle \bullet$} [c] at  16 12
\setlinear \plot  10 0 10 12  /
\setlinear \plot  16 0 16 12 10 8   /
\setlinear \plot  10 8 13 12  /
\put{$2{,}520$} [c] at 13 -2
\endpicture
\end{minipage}
\begin{minipage}{4cm}
\beginpicture
\setcoordinatesystem units   <1.5mm,2mm>
\setplotarea x from 0 to 16, y from -2 to 15
\put{602)} [l] at 2 12
\put {$ \scriptstyle \bullet$} [c] at  10 0
\put {$ \scriptstyle \bullet$} [c] at  10 4
\put {$ \scriptstyle \bullet$} [c] at  10 8
\put {$ \scriptstyle \bullet$} [c] at  10 12
\put {$ \scriptstyle \bullet$} [c] at  13 0
\put {$ \scriptstyle \bullet$} [c] at  16 0
\put {$ \scriptstyle \bullet$} [c] at  16 12
\setlinear \plot  10 0 10 12  /
\setlinear \plot  16 12 16 0 10 4   /
\setlinear \plot  10 4 13 0  /
\put{$2{,}520$} [c] at 13 -2
\endpicture
\end{minipage}
\begin{minipage}{4cm}
\beginpicture
\setcoordinatesystem units   <1.5mm,2mm>
\setplotarea x from 0 to 16, y from -2 to 15
\put{603)} [l] at 2 12
\put {$ \scriptstyle \bullet$} [c] at  10 0
\put {$ \scriptstyle \bullet$} [c] at  10 6
\put {$ \scriptstyle \bullet$} [c] at  10 12
\put {$ \scriptstyle \bullet$} [c] at  13 12
\put {$ \scriptstyle \bullet$} [c] at  16 0
\put {$ \scriptstyle \bullet$} [c] at  16 6
\put {$ \scriptstyle \bullet$} [c] at  16 12
\setlinear \plot  10 0 10 12 16 0 16 12 10 0 /
\setlinear \plot  10 6  13 12 16 6 /
\put{$2{,}520$} [c] at 13 -2
\endpicture
\end{minipage}
\begin{minipage}{4cm}
\beginpicture
\setcoordinatesystem units   <1.5mm,2mm>
\setplotarea x from 0 to 16, y from -2 to 15
\put{604)} [l] at 2 12
\put {$ \scriptstyle \bullet$} [c] at  10 0
\put {$ \scriptstyle \bullet$} [c] at  10 6
\put {$ \scriptstyle \bullet$} [c] at  10 12
\put {$ \scriptstyle \bullet$} [c] at  13 0
\put {$ \scriptstyle \bullet$} [c] at  16 0
\put {$ \scriptstyle \bullet$} [c] at  16 6
\put {$ \scriptstyle \bullet$} [c] at  16 12
\setlinear \plot  10 12 10 0 16 12 16 0 10 12 /
\setlinear \plot  10 6  13 0 16 6 /
\put{$2{,}520$} [c] at 13 -2
\endpicture
\end{minipage}
\begin{minipage}{4cm}
\beginpicture
\setcoordinatesystem units   <1.5mm,2mm>
\setplotarea x from 0 to 16, y from -2 to 15
\put{605)} [l] at 2 12
\put {$ \scriptstyle \bullet$} [c] at  10 0
\put {$ \scriptstyle \bullet$} [c] at  10 4
\put {$ \scriptstyle \bullet$} [c] at  10 8
\put {$ \scriptstyle \bullet$} [c] at  10 12
\put {$ \scriptstyle \bullet$} [c] at  13 12
\put {$ \scriptstyle \bullet$} [c] at  16 0
\put {$ \scriptstyle \bullet$} [c] at  16 12
\setlinear \plot  10 0 10 12 16 0 16 12 10 4 13 12 16 0 /
\put{$2{,}520$} [c] at 13 -2
\endpicture
\end{minipage}
\begin{minipage}{4cm}
\beginpicture
\setcoordinatesystem units   <1.5mm,2mm>
\setplotarea x from 0 to 16, y from -2 to 15
\put{606)} [l] at 2 12
\put {$ \scriptstyle \bullet$} [c] at  10 0
\put {$ \scriptstyle \bullet$} [c] at  10 4
\put {$ \scriptstyle \bullet$} [c] at  10 8
\put {$ \scriptstyle \bullet$} [c] at  10 12
\put {$ \scriptstyle \bullet$} [c] at  13 0
\put {$ \scriptstyle \bullet$} [c] at  16 0
\put {$ \scriptstyle \bullet$} [c] at  16 12
\setlinear \plot  10 12 10 0 16 12 16 0 10 8 13 0 16 12 /
\put{$2{,}520$} [c] at 13 -2
\endpicture
\end{minipage}
$$

$$
\begin{minipage}{4cm}
\beginpicture
\setcoordinatesystem units   <1.5mm,2mm>
\setplotarea x from 0 to 16, y from -2 to 15
\put{607)} [l] at 2 12
\put {$ \scriptstyle \bullet$} [c] at  10 12
\put {$ \scriptstyle \bullet$} [c] at  13 0
\put {$ \scriptstyle \bullet$} [c] at  13 6
\put {$ \scriptstyle \bullet$} [c] at  13 12
\put {$ \scriptstyle \bullet$} [c] at  16 0
\put {$ \scriptstyle \bullet$} [c] at  16 6
\put {$ \scriptstyle \bullet$} [c] at  16 12
\setlinear \plot  13 0 13 12 16 6 16 0 /
\setlinear \plot  10 12 13  6 16  12 16 6   /
\put{$2{,}520$} [c] at 13 -2
\endpicture
\end{minipage}
\begin{minipage}{4cm}
\beginpicture
\setcoordinatesystem units   <1.5mm,2mm>
\setplotarea x from 0 to 16, y from -2 to 15
\put{608)} [l] at 2 12
\put {$ \scriptstyle \bullet$} [c] at  10 0
\put {$ \scriptstyle \bullet$} [c] at  13 0
\put {$ \scriptstyle \bullet$} [c] at  13 6
\put {$ \scriptstyle \bullet$} [c] at  13 12
\put {$ \scriptstyle \bullet$} [c] at  16 0
\put {$ \scriptstyle \bullet$} [c] at  16 6
\put {$ \scriptstyle \bullet$} [c] at  16 12
\setlinear \plot  13 12 13 0 16 6 16 12 /
\setlinear \plot  10 0 13  6 16  0 16 6   /
\put{$2{,}520$} [c] at 13 -2
\endpicture
\end{minipage}
\begin{minipage}{4cm}
\beginpicture
\setcoordinatesystem units   <1.5mm,2mm>
\setplotarea x from 0 to 16, y from -2 to 15
\put{609)} [l] at 2 12
\put {$ \scriptstyle \bullet$} [c] at  10 0
\put {$ \scriptstyle \bullet$} [c] at  10 4
\put {$ \scriptstyle \bullet$} [c] at  10 8
\put {$ \scriptstyle \bullet$} [c] at  10 12
\put {$ \scriptstyle \bullet$} [c] at  13 12
\put {$ \scriptstyle \bullet$} [c] at  16 0
\put {$ \scriptstyle \bullet$} [c] at  16 12
\setlinear \plot  16 12 16 0 10 12 10 0  /
\setlinear \plot  10 8 13 12 16 0 /
\put{$2{,}520$} [c] at 13 -2
\endpicture
\end{minipage}
\begin{minipage}{4cm}
\beginpicture
\setcoordinatesystem units   <1.5mm,2mm>
\setplotarea x from 0 to 16, y from -2 to 15
\put{610)} [l] at 2 12
\put {$ \scriptstyle \bullet$} [c] at  10 0
\put {$ \scriptstyle \bullet$} [c] at  10 4
\put {$ \scriptstyle \bullet$} [c] at  10 8
\put {$ \scriptstyle \bullet$} [c] at  10 12
\put {$ \scriptstyle \bullet$} [c] at  13 0
\put {$ \scriptstyle \bullet$} [c] at  16 0
\put {$ \scriptstyle \bullet$} [c] at  16 12
\setlinear \plot  16 0 16 12 10 0 10 12  /
\setlinear \plot  10 4 13 0 16 12 /
\put{$2{,}520$} [c] at 13 -2
\endpicture
\end{minipage}
\begin{minipage}{4cm}
\beginpicture
\setcoordinatesystem units   <1.5mm,2mm>
\setplotarea x from 0 to 16, y from -2 to 15
\put{611)} [l] at 2 12
\put {$ \scriptstyle \bullet$} [c] at  10 12
\put {$ \scriptstyle \bullet$} [c] at  12 0
\put {$ \scriptstyle \bullet$} [c] at  12 6
\put {$ \scriptstyle \bullet$} [c] at  12 12
\put {$ \scriptstyle \bullet$} [c] at  14 0
\put {$ \scriptstyle \bullet$} [c] at  14 12
\put {$ \scriptstyle \bullet$} [c] at  16 6
\setlinear \plot 10 12 12 6 14 0 16 6 14 12 12 6 12 12   /
\setlinear \plot  12 0 12 6  /
\put{$2{,}520$} [c] at 13 -2
\endpicture
\end{minipage}
\begin{minipage}{4cm}
\beginpicture
\setcoordinatesystem units   <1.5mm,2mm>
\setplotarea x from 0 to 16, y from -2 to 15
\put{612)} [l] at 2 12
\put {$ \scriptstyle \bullet$} [c] at  10 0
\put {$ \scriptstyle \bullet$} [c] at  12 0
\put {$ \scriptstyle \bullet$} [c] at  12 6
\put {$ \scriptstyle \bullet$} [c] at  12 12
\put {$ \scriptstyle \bullet$} [c] at  14 0
\put {$ \scriptstyle \bullet$} [c] at  14 12
\put {$ \scriptstyle \bullet$} [c] at  16 6
\setlinear \plot 10 0 12 6 14 0 16 6 14 12 12 6 12 12   /
\setlinear \plot  12 0 12 6  /
\put{$2{,}520$} [c] at 13 -2
\endpicture
\end{minipage}
$$
$$
\begin{minipage}{4cm}
\beginpicture
\setcoordinatesystem units   <1.5mm,2mm>
\setplotarea x from 0 to 16, y from -2 to 15
\put{613)} [l] at 2 12
\put {$ \scriptstyle \bullet$} [c] at  10 12
\put {$ \scriptstyle \bullet$} [c] at  13 0
\put {$ \scriptstyle \bullet$} [c] at  13 4
\put {$ \scriptstyle \bullet$} [c] at  13 8
\put {$ \scriptstyle \bullet$} [c] at  13 12
\put {$ \scriptstyle \bullet$} [c] at  16 0
\put {$ \scriptstyle \bullet$} [c] at  16 12
\setlinear \plot  13 0 13 12 16  0 16 12 13 8   /
\setlinear \plot  10 12 13 4    /
\put{$2{,}520$} [c] at 13 -2
\endpicture
\end{minipage}
\begin{minipage}{4cm}
\beginpicture
\setcoordinatesystem units   <1.5mm,2mm>
\setplotarea x from 0 to 16, y from -2 to 15
\put{614)} [l] at 2 12
\put {$ \scriptstyle \bullet$} [c] at  10 0
\put {$ \scriptstyle \bullet$} [c] at  13 0
\put {$ \scriptstyle \bullet$} [c] at  13 4
\put {$ \scriptstyle \bullet$} [c] at  13 8
\put {$ \scriptstyle \bullet$} [c] at  13 12
\put {$ \scriptstyle \bullet$} [c] at  16 0
\put {$ \scriptstyle \bullet$} [c] at  16 12
\setlinear \plot  13 12 13 0 16 12 16 0 13 4   /
\setlinear \plot  10 0 13 8    /
\put{$2{,}520$} [c] at 13 -2
\endpicture
\end{minipage}
\begin{minipage}{4cm}
\beginpicture
\setcoordinatesystem units   <1.5mm,2mm>
\setplotarea x from 0 to 16, y from -2 to 15
\put{615)} [l] at 2 12
\put {$ \scriptstyle \bullet$} [c] at  10 12
\put {$ \scriptstyle \bullet$} [c] at  10 0
\put {$ \scriptstyle \bullet$} [c] at  13 6
\put {$ \scriptstyle \bullet$} [c] at  13  12
\put {$ \scriptstyle \bullet$} [c] at  14.5 3
\put {$ \scriptstyle \bullet$} [c] at  16 0
\put {$ \scriptstyle \bullet$} [c] at  16 12
\setlinear \plot 10 12   16 0 16 12   /
\setlinear \plot  10 0 13  6 13 12 /
\put{$2{,}520$} [c] at 13 -2
\endpicture
\end{minipage}
\begin{minipage}{4cm}
\beginpicture
\setcoordinatesystem units   <1.5mm,2mm>
\setplotarea x from 0 to 16, y from -2 to 15
\put{616)} [l] at 2 12
\put {$ \scriptstyle \bullet$} [c] at  10 12
\put {$ \scriptstyle \bullet$} [c] at  10 0
\put {$ \scriptstyle \bullet$} [c] at  13 6
\put {$ \scriptstyle \bullet$} [c] at  13  0
\put {$ \scriptstyle \bullet$} [c] at  14.5 9
\put {$ \scriptstyle \bullet$} [c] at  16 0
\put {$ \scriptstyle \bullet$} [c] at  16 12
\setlinear \plot 10 0   16 12 16 0   /
\setlinear \plot  10 12 13  6 13 0 /
\put{$2{,}520$} [c] at 13 -2
\endpicture
\end{minipage}
\begin{minipage}{4cm}
\beginpicture
\setcoordinatesystem units   <1.5mm,2mm>
\setplotarea x from 0 to 16, y from -2 to 15
\put{617)} [l] at 2 12
\put {$ \scriptstyle \bullet$} [c] at  10 0
\put {$ \scriptstyle \bullet$} [c] at  10 12
\put {$ \scriptstyle \bullet$} [c] at  14 6
\put {$ \scriptstyle \bullet$} [c] at  14 12
\put {$ \scriptstyle \bullet$} [c] at  15 0
\put {$ \scriptstyle \bullet$} [c] at  16 6
\put {$ \scriptstyle \bullet$} [c] at  16 12
\setlinear \plot 15 0 10 12 10 0 14 6 15 0  16 6 16 12 14 6 14 12 16 6  /
\put{$2{,}520$} [c] at 13 -2
\endpicture
\end{minipage}
\begin{minipage}{4cm}
\beginpicture
\setcoordinatesystem units   <1.5mm,2mm>
\setplotarea x from 0 to 16, y from -2 to 15
\put{618)} [l] at 2 12
\put {$ \scriptstyle \bullet$} [c] at  10 0
\put {$ \scriptstyle \bullet$} [c] at  10 12
\put {$ \scriptstyle \bullet$} [c] at  14 6
\put {$ \scriptstyle \bullet$} [c] at  14 0
\put {$ \scriptstyle \bullet$} [c] at  15 12
\put {$ \scriptstyle \bullet$} [c] at  16 6
\put {$ \scriptstyle \bullet$} [c] at  16 0
\setlinear \plot 15 12 10 0 10 12 14 6 15 12  16 6 16 0 14 6 14 0 16 6  /
\put{$2{,}520$} [c] at 13 -2
\endpicture
\end{minipage}
$$

$$
\begin{minipage}{4cm}
\beginpicture
\setcoordinatesystem units   <1.5mm,2mm>
\setplotarea x from 0 to 16, y from -2 to 15
\put{619)} [l] at 2 12
\put {$ \scriptstyle \bullet$} [c] at  10 0
\put {$ \scriptstyle \bullet$} [c] at  10 12
\put {$ \scriptstyle \bullet$} [c] at  14 6
\put {$ \scriptstyle \bullet$} [c] at  14 12
\put {$ \scriptstyle \bullet$} [c] at  15 0
\put {$ \scriptstyle \bullet$} [c] at  16 6
\put {$ \scriptstyle \bullet$} [c] at  16 12
\setlinear \plot 14 6 10 12 10 0 14 12 14 6 15 0  16 6 16 12 14 6 14 12 16 6  /
\setlinear \plot 14 12 10 0 16 12 /
\put{$2{,}520$} [c] at 13 -2
\endpicture
\end{minipage}
\begin{minipage}{4cm}
\beginpicture
\setcoordinatesystem units   <1.5mm,2mm>
\setplotarea x from 0 to 16, y from -2 to 15
\put{620)} [l] at 2 12
\put {$ \scriptstyle \bullet$} [c] at  10 0
\put {$ \scriptstyle \bullet$} [c] at  10 12
\put {$ \scriptstyle \bullet$} [c] at  14 6
\put {$ \scriptstyle \bullet$} [c] at  14 0
\put {$ \scriptstyle \bullet$} [c] at  15 12
\put {$ \scriptstyle \bullet$} [c] at  16 6
\put {$ \scriptstyle \bullet$} [c] at  16 0
\setlinear \plot 14 6 10 0 10 12 14 0 14 6 15 12  16 6 16 0 14 6 14 0 16 6  /
\setlinear \plot 14 0 10 12 16 0 /
\put{$2{,}520$} [c] at 13 -2
\endpicture
\end{minipage}
\begin{minipage}{4cm}
\beginpicture
\setcoordinatesystem units   <1.5mm,2mm>
\setplotarea x from 0 to 16, y from -2 to 15
\put{621)} [l] at 2 12
\put {$ \scriptstyle \bullet$} [c] at  10 0
\put {$ \scriptstyle \bullet$} [c] at  10 6
\put {$ \scriptstyle \bullet$} [c] at  12 12
\put {$ \scriptstyle \bullet$} [c] at  14 0
\put {$ \scriptstyle \bullet$} [c] at  14 6
\put {$ \scriptstyle \bullet$} [c] at  16 6
\put {$ \scriptstyle \bullet$} [c] at  16 12
\setlinear \plot 14 6 10 0 10  6 12 12 14 6  14 0  10 6  /
\setlinear \plot 16 12 16 6 14 0  /
\setlinear \plot 10 0 16 6  /
\put{$1{,}260    $} [c] at 13 -2
\endpicture
\end{minipage}
\begin{minipage}{4cm}
\beginpicture
\setcoordinatesystem units   <1.5mm,2mm>
\setplotarea x from 0 to 16, y from -2 to 15
\put{622)} [l] at 2 12
\put {$ \scriptstyle \bullet$} [c] at  10 12
\put {$ \scriptstyle \bullet$} [c] at  10 6
\put {$ \scriptstyle \bullet$} [c] at  12 0
\put {$ \scriptstyle \bullet$} [c] at  14 12
\put {$ \scriptstyle \bullet$} [c] at  14 6
\put {$ \scriptstyle \bullet$} [c] at  16 6
\put {$ \scriptstyle \bullet$} [c] at  16 0
\setlinear \plot 14 6 10 12 10  6 12 0 14 6  14 12  10 6  /
\setlinear \plot 16 0 16 6 14 12  /
\setlinear \plot 10 12 16 6  /
\put{$1{,}260    $} [c] at 13 -2
\endpicture
\end{minipage}
\begin{minipage}{4cm}
\beginpicture
\setcoordinatesystem units   <1.5mm,2mm>
\setplotarea x from 0 to 16, y from -2 to 15
\put{623)} [l] at 2 12
\put {$ \scriptstyle \bullet$} [c] at  10 0
\put {$ \scriptstyle \bullet$} [c] at  10 12
\put {$ \scriptstyle \bullet$} [c] at  14 6
\put {$ \scriptstyle \bullet$} [c] at  14 0
\put {$ \scriptstyle \bullet$} [c] at  15 12
\put {$ \scriptstyle \bullet$} [c] at  16 6
\put {$ \scriptstyle \bullet$} [c] at  16 0
\setlinear \plot 10 0 10 12 14 0 14 6 15 12  16 6 16 0 14 6 14 0 16 6  /
\setlinear \plot 14 0  16 6 10 0 14 6 /
\put{$1{,}260$} [c] at 13 -2
\endpicture
\end{minipage}
\begin{minipage}{4cm}
\beginpicture
\setcoordinatesystem units   <1.5mm,2mm>
\setplotarea x from 0 to 16, y from -2 to 15
\put{624)} [l] at 2 12
\put {$ \scriptstyle \bullet$} [c] at  10 0
\put {$ \scriptstyle \bullet$} [c] at  10 12
\put {$ \scriptstyle \bullet$} [c] at  14 6
\put {$ \scriptstyle \bullet$} [c] at  14 12
\put {$ \scriptstyle \bullet$} [c] at  15 0
\put {$ \scriptstyle \bullet$} [c] at  16 6
\put {$ \scriptstyle \bullet$} [c] at  16 12
\setlinear \plot 10 12 10 0 14 12 14 6 15 0  16 6 16 12 14 6 14 12 16 6  /
\setlinear \plot 14 12  16 6 10 12 14 6 /
\put{$1{,}260$} [c] at 13 -2
\endpicture
\end{minipage}
$$
$$
\begin{minipage}{4cm}
\beginpicture
\setcoordinatesystem units   <1.5mm,2mm>
\setplotarea x from 0 to 16, y from -2 to 15
\put{625)} [l] at 2 12
\put {$ \scriptstyle \bullet$} [c] at  10 0
\put {$ \scriptstyle \bullet$} [c] at  10 6
\put {$ \scriptstyle \bullet$} [c] at  10 12
\put {$ \scriptstyle \bullet$} [c] at  13 12
\put {$ \scriptstyle \bullet$} [c] at  16 0
\put {$ \scriptstyle \bullet$} [c] at  16 6
\put {$ \scriptstyle \bullet$} [c] at  16 12
\setlinear \plot  10 12 10 6 16 0 16 6 10 0 10 6 /
\setlinear \plot  13 12 16 6 16 12   /
\put{$1{,}260$} [c] at 13 -2
\endpicture
\end{minipage}
\begin{minipage}{4cm}
\beginpicture
\setcoordinatesystem units   <1.5mm,2mm>
\setplotarea x from 0 to 16, y from -2 to 15
\put{626)} [l] at 2 12
\put {$ \scriptstyle \bullet$} [c] at  10 0
\put {$ \scriptstyle \bullet$} [c] at  10 6
\put {$ \scriptstyle \bullet$} [c] at  10 12
\put {$ \scriptstyle \bullet$} [c] at  13 0
\put {$ \scriptstyle \bullet$} [c] at  16 0
\put {$ \scriptstyle \bullet$} [c] at  16 6
\put {$ \scriptstyle \bullet$} [c] at  16 12
\setlinear \plot  10 0 10 6 16 12 16 6 10 12 10 6 /
\setlinear \plot  13 0 16 6 16 0   /
\put{$1{,}260$} [c] at 13 -2
\endpicture
\end{minipage}
\begin{minipage}{4cm}
\beginpicture
\setcoordinatesystem units   <1.5mm,2mm>
\setplotarea x from 0 to 16, y from -2 to 15
\put{627)} [l] at 2 12
\put {$ \scriptstyle \bullet$} [c] at  10 0
\put {$ \scriptstyle \bullet$} [c] at  13 12
\put {$ \scriptstyle \bullet$} [c] at  13 6
\put {$ \scriptstyle \bullet$} [c] at  13 0
\put {$ \scriptstyle \bullet$} [c] at  16 0
\put {$ \scriptstyle \bullet$} [c] at  16 6
\put {$ \scriptstyle \bullet$} [c] at  16 12
\setlinear \plot 10 0 13 12  13 0 16 6 16 12 13 6 16 0 16 6 13 12 /
\put{$1{,}260$} [c] at 13 -2
\endpicture
\end{minipage}
\begin{minipage}{4cm}
\beginpicture
\setcoordinatesystem units   <1.5mm,2mm>
\setplotarea x from 0 to 16, y from -2 to 15
\put{628)} [l] at 2 12
\put {$ \scriptstyle \bullet$} [c] at  10 12
\put {$ \scriptstyle \bullet$} [c] at  13 12
\put {$ \scriptstyle \bullet$} [c] at  13 6
\put {$ \scriptstyle \bullet$} [c] at  13 0
\put {$ \scriptstyle \bullet$} [c] at  16 0
\put {$ \scriptstyle \bullet$} [c] at  16 6
\put {$ \scriptstyle \bullet$} [c] at  16 12
\setlinear \plot 10 12 13 0  13 12 16 6 16 0 13 6 16 12 16 6 13 0 /
\put{$1{,}260$} [c] at 13 -2
\endpicture
\end{minipage}
\begin{minipage}{4cm}
\beginpicture
\setcoordinatesystem units   <1.5mm,2mm>
\setplotarea x from 0 to 16, y from -2 to 15
\put{629)} [l] at 2 12
\put {$ \scriptstyle \bullet$} [c] at  10 0
\put {$ \scriptstyle \bullet$} [c] at  10 6
\put {$ \scriptstyle \bullet$} [c] at  10 12
\put {$ \scriptstyle \bullet$} [c] at  13 12
\put {$ \scriptstyle \bullet$} [c] at  13  0
\put {$ \scriptstyle \bullet$} [c] at  16  12
\put {$ \scriptstyle \bullet$} [c] at  16  0
\setlinear \plot  10  0 10  12 /
\setlinear \plot  16 0 16  12 13 0 10 6  13 12  /
\setlinear \plot  10 6 16  0 /
\put{$1{,}260$} [c] at 12 -2
\put{$$} [c] at 12 -2
\endpicture
\end{minipage}
\begin{minipage}{4cm}
\beginpicture
\setcoordinatesystem units   <1.5mm,2mm>
\setplotarea x from 0 to 16, y from -2 to 15
\put{630)} [l] at 2 12
\put {$ \scriptstyle \bullet$} [c] at  10 0
\put {$ \scriptstyle \bullet$} [c] at  10 6
\put {$ \scriptstyle \bullet$} [c] at  10 12
\put {$ \scriptstyle \bullet$} [c] at  13 12
\put {$ \scriptstyle \bullet$} [c] at  13  0
\put {$ \scriptstyle \bullet$} [c] at  16  12
\put {$ \scriptstyle \bullet$} [c] at  16  0
\setlinear \plot  10  0 10  12 /
\setlinear \plot  16 12 16  0 13 12 10 6  13 0  /
\setlinear \plot  10 6 16 12 /
\put{$1{,}260$} [c] at 12 -2
\endpicture
\end{minipage}
$$

$$
\begin{minipage}{4cm}
\beginpicture
\setcoordinatesystem units   <1.5mm,2mm>
\setplotarea x from 0 to 16, y from -2 to 15
\put{631)} [l] at 2 12
\put {$ \scriptstyle \bullet$} [c] at  10 0
\put {$ \scriptstyle \bullet$} [c] at  10 6
\put {$ \scriptstyle \bullet$} [c] at  10 12
\put {$ \scriptstyle \bullet$} [c] at  13 0
\put {$ \scriptstyle \bullet$} [c] at  13  12
\put {$ \scriptstyle \bullet$} [c] at  16  12
\put {$ \scriptstyle \bullet$} [c] at  16  0
\setlinear \plot  16 12 10  0 10 12 16 0  16 12 13 0 10 6 13 12 16 0 /
\put{$1{,}260 $} [c] at 12 -2
\endpicture
\end{minipage}
\begin{minipage}{4cm}
\beginpicture
\setcoordinatesystem units   <1.5mm,2mm>
\setplotarea x from 0 to 16, y from -2 to 15
\put{632)} [l] at 2 12
\put {$ \scriptstyle \bullet$} [c] at 10 3
\put {$ \scriptstyle \bullet$} [c] at 11 6
\put {$ \scriptstyle \bullet$} [c] at 12 0
\put {$ \scriptstyle \bullet$} [c] at 12 9
\put {$ \scriptstyle \bullet$} [c] at 12 12
\put {$ \scriptstyle \bullet$} [c] at 14 3
\put{$\scriptstyle \bullet$} [c] at 16  0
\setlinear \plot  12 12  12 9  10 3 12 0 14 3 12 9      /
\put{$5{,}040$} [c] at 13 -2
\endpicture
\end{minipage}
\begin{minipage}{4cm}
\beginpicture
\setcoordinatesystem units   <1.5mm,2mm>
\setplotarea x from 0 to 16, y from -2 to 15
\put{633)} [l] at 2 12
\put {$ \scriptstyle \bullet$} [c] at 10 9
\put {$ \scriptstyle \bullet$} [c] at 11 6
\put {$ \scriptstyle \bullet$} [c] at 12 0
\put {$ \scriptstyle \bullet$} [c] at 12 3
\put {$ \scriptstyle \bullet$} [c] at 12 12
\put {$ \scriptstyle \bullet$} [c] at 14 9
\put{$\scriptstyle \bullet$} [c] at 16  0
\setlinear \plot  12 0  12 3  10 9 12 12 14 9 12 3      /
\put{$5{,}040$} [c] at 13 -2
\endpicture
\end{minipage}
\begin{minipage}{4cm}
\beginpicture
\setcoordinatesystem units   <1.5mm,2mm>
\setplotarea x from 0 to 16, y from -2 to 15
\put{634)} [l] at 2 12
\put {$ \scriptstyle \bullet$} [c] at 10 12
\put {$ \scriptstyle \bullet$} [c] at 11 9
\put {$ \scriptstyle \bullet$} [c] at 12 0
\put {$ \scriptstyle \bullet$} [c] at 12 3
\put {$ \scriptstyle \bullet$} [c] at 12 6
\put {$ \scriptstyle \bullet$} [c] at 14 12
\put{$\scriptstyle \bullet$} [c] at 16  0
\setlinear \plot  12 0  12 6  10 12      /
\setlinear \plot  12 6  14 12      /
\put{$5{,}040$} [c] at 13 -2
 \endpicture
\end{minipage}
\begin{minipage}{4cm}
\beginpicture
\setcoordinatesystem units   <1.5mm,2mm>
\setplotarea x from 0 to 16, y from -2 to 15
\put{635)} [l] at 2 12
\put {$ \scriptstyle \bullet$} [c] at 10 0
\put {$ \scriptstyle \bullet$} [c] at 11 3
\put {$ \scriptstyle \bullet$} [c] at 12 6
\put {$ \scriptstyle \bullet$} [c] at 12 9
\put {$ \scriptstyle \bullet$} [c] at 12 12
\put {$ \scriptstyle \bullet$} [c] at 14 0
\put{$\scriptstyle \bullet$} [c] at 16  0
\setlinear \plot  12 12  12 6  10 0   /
\setlinear \plot  12 6  14 0   /
\put{$5{,}040$} [c] at 13 -2
\endpicture
\end{minipage}
\begin{minipage}{4cm}
\beginpicture
\setcoordinatesystem units   <1.5mm,2mm>
\setplotarea x from 0 to  16, y from -2 to 15
\put{636)} [l] at 2 12
\put {$ \scriptstyle \bullet$} [c] at 10 0
\put {$ \scriptstyle \bullet$} [c] at 10 3
\put {$ \scriptstyle \bullet$} [c] at 10 6
\put {$ \scriptstyle \bullet$} [c] at 10 9
\put {$ \scriptstyle \bullet$} [c] at 10 12
\put {$ \scriptstyle \bullet$} [c] at 16 0
\put {$ \scriptstyle \bullet$} [c] at 16 12
\setlinear \plot 10 0 10 12   /
\setlinear \plot 16 0 16 12   /
\put{$5{,}040$} [c] at 13 -2
\endpicture
\end{minipage}
$$
$$
\begin{minipage}{4cm}
\beginpicture
\setcoordinatesystem units   <1.5mm,2mm>
\setplotarea x from 0 to 16, y from -2 to 15
\put{637)} [l] at 2 12
\put {$ \scriptstyle \bullet$} [c] at 10 8
\put {$ \scriptstyle \bullet$} [c] at 10 12
\put {$ \scriptstyle \bullet$} [c] at  12 0
\put {$ \scriptstyle \bullet$} [c] at 12 4
\put {$ \scriptstyle \bullet$} [c] at 14 8
\put {$ \scriptstyle \bullet$} [c] at 14 12
\put{$\scriptstyle \bullet$} [c] at 16  0
\setlinear \plot   12 0 12 4 10 8 10 12 14 8 14 12 10 8  /
\setlinear \plot  12 4 14 8 /
\put{$1{,}260$} [c] at 13 -2
\endpicture
\end{minipage}
\begin{minipage}{4cm}
\beginpicture
\setcoordinatesystem units   <1.5mm,2mm>
\setplotarea x from 0 to 16, y from -2 to 15
\put{638)} [l] at 2 12
\put {$ \scriptstyle \bullet$} [c] at 10 4
\put {$ \scriptstyle \bullet$} [c] at 10 0
\put {$ \scriptstyle \bullet$} [c] at  12 12
\put {$ \scriptstyle \bullet$} [c] at 12 8
\put {$ \scriptstyle \bullet$} [c] at 14 4
\put {$ \scriptstyle \bullet$} [c] at 14 0
\put{$\scriptstyle \bullet$} [c] at 16  0
\setlinear \plot   12 12 12 8 10 4 10 0 14 4 14 0 10 4  /
\setlinear \plot  12 8 14 4 /
\put{$1{,}260$} [c] at 13 -2
\endpicture
\end{minipage}
\begin{minipage}{4cm}
\beginpicture
\setcoordinatesystem units   <1.5mm,2mm>
\setplotarea x from 0 to 16, y from -2 to 15
\put{639)} [l] at 2 12
\put {$ \scriptstyle \bullet$} [c] at 10 4
\put {$ \scriptstyle \bullet$} [c] at 10 12
\put {$ \scriptstyle \bullet$} [c] at 12 0
\put {$ \scriptstyle \bullet$} [c] at 12 8
\put {$ \scriptstyle \bullet$} [c] at 14 4
\put {$ \scriptstyle \bullet$} [c] at 14 12
\put{$\scriptstyle \bullet$} [c] at 16  0
\setlinear \plot  12 0  10 4 12 8  14  12      /
\setlinear \plot  12 0  14 4 12 8 10 12      /
\put{$1{,}260$} [c] at 13 -2
\endpicture
\end{minipage}
\begin{minipage}{4cm}
\beginpicture
\setcoordinatesystem units   <1.5mm,2mm>
\setplotarea x from 0 to 16, y from -2 to 15
\put{640)} [l] at 2 12
\put {$ \scriptstyle \bullet$} [c] at 10 8
\put {$ \scriptstyle \bullet$} [c] at 10 0
\put {$ \scriptstyle \bullet$} [c] at 12 12
\put {$ \scriptstyle \bullet$} [c] at 12 4
\put {$ \scriptstyle \bullet$} [c] at 14 8
\put {$ \scriptstyle \bullet$} [c] at 14 0
\put{$\scriptstyle \bullet$} [c] at 16  0
\setlinear \plot  12 12  10 8 12 4  14  0      /
\setlinear \plot  12 12  14 8 12 4 10 0      /
\put{$1{,}260$} [c] at 13 -2
\endpicture
\end{minipage}
\begin{minipage}{4cm}
\beginpicture
\setcoordinatesystem units   <1.5mm,2mm>
\setplotarea x from 0 to 16, y from -2 to 15
\put{641)} [l] at 2 12
\put {$ \scriptstyle \bullet$} [c] at 10 4
\put {$ \scriptstyle \bullet$} [c] at 10 8
\put {$ \scriptstyle \bullet$} [c] at 12 0
\put {$ \scriptstyle \bullet$} [c] at 12 12
\put {$ \scriptstyle \bullet$} [c] at 14 4
\put {$ \scriptstyle \bullet$} [c] at 14 8
 \put{$\scriptstyle \bullet$} [c] at 16  0
\setlinear \plot 10 8 14 4  12 0  10 4 14  8 12 12  10 8 10 4  14 8   14 4 /
\put{$1{,}260$} [c] at 13 -2
\endpicture
\end{minipage}
\begin{minipage}{4cm}
\beginpicture
\setcoordinatesystem units   <1.5mm,2mm>
\setplotarea x from 0 to 16, y from -2 to 15
\put{642)} [l] at 2 12
\put {$ \scriptstyle \bullet$} [c] at 10 0
\put {$ \scriptstyle \bullet$} [c] at 10 12
\put {$ \scriptstyle \bullet$} [c] at 12 4
\put {$ \scriptstyle \bullet$} [c] at 12 8
\put {$ \scriptstyle \bullet$} [c] at 14 0
\put {$ \scriptstyle \bullet$} [c] at 14 12
\put{$\scriptstyle \bullet$} [c] at 16  0
\setlinear \plot   10 0 12 4 12 8 10 12     /
\setlinear \plot   12 8 14 12     /
\setlinear \plot   14 0 12 4     /
\put{$1{,}260$} [c] at 13 -2
\endpicture
\end{minipage}
$$

$$
\begin{minipage}{4cm}
\beginpicture
\setcoordinatesystem units   <1.5mm,2mm>
\setplotarea x from 0 to 16, y from -2 to 15
\put{${\bf  19}$} [l] at 2 15
\put{643)} [l] at 2 12
\put {$ \scriptstyle \bullet$} [c] at 10 6
\put {$ \scriptstyle \bullet$} [c] at 10 12
\put {$ \scriptstyle \bullet$} [c] at 13 0
\put {$ \scriptstyle \bullet$} [c] at 13 6
\put {$ \scriptstyle \bullet$} [c] at 16 6
\put {$ \scriptstyle \bullet$} [c] at 16 9
\put {$ \scriptstyle \bullet$} [c] at 16 12
\setlinear \plot 13  0 10 6 10 12 13 6 13 0 16 6 16 12 13 6   /
\put{$5{,}040$} [c] at 13 -2
\endpicture
\end{minipage}
\begin{minipage}{4cm}
\beginpicture
\setcoordinatesystem units   <1.5mm,2mm>
\setplotarea x from 0 to 16, y from -2 to 15
\put{644)} [l] at 2 12
\put {$ \scriptstyle \bullet$} [c] at 10 6
\put {$ \scriptstyle \bullet$} [c] at 10 0
\put {$ \scriptstyle \bullet$} [c] at 13 12
\put {$ \scriptstyle \bullet$} [c] at 13 6
\put {$ \scriptstyle \bullet$} [c] at 16 6
\put {$ \scriptstyle \bullet$} [c] at 16 3
\put {$ \scriptstyle \bullet$} [c] at 16 0
\setlinear \plot 13  12 10 6 10 0 13 6 13 12 16 6 16 0 13 6   /
\put{$5{,}040$} [c] at 13 -2
\endpicture
\end{minipage}
\begin{minipage}{4cm}
\beginpicture
\setcoordinatesystem units   <1.5mm,2mm>
\setplotarea x from 0 to 16, y from -2 to 15
\put{645)} [l] at 2 12
\put {$ \scriptstyle \bullet$} [c] at 10 4
\put {$ \scriptstyle \bullet$} [c] at 10 12
\put {$ \scriptstyle \bullet$} [c] at 13 0
\put {$ \scriptstyle \bullet$} [c] at 13 4
\put {$ \scriptstyle \bullet$} [c] at 13 8
\put {$ \scriptstyle \bullet$} [c] at 16 4
\put {$ \scriptstyle \bullet$} [c] at 16 12
\setlinear \plot 10  12  10 4 13 0 16 4 16 12 10 4  /
\setlinear \plot 13 0 13 8  /
\put{$5{,}040$} [c] at 13 -2
\endpicture
\end{minipage}
\begin{minipage}{4cm}
\beginpicture
\setcoordinatesystem units   <1.5mm,2mm>
\setplotarea x from 0 to 16, y from -2 to 15
\put{646)} [l] at 2 12
\put {$ \scriptstyle \bullet$} [c] at 10 8
\put {$ \scriptstyle \bullet$} [c] at 10 0
\put {$ \scriptstyle \bullet$} [c] at 13 12
\put {$ \scriptstyle \bullet$} [c] at 13 8
\put {$ \scriptstyle \bullet$} [c] at 13 4
\put {$ \scriptstyle \bullet$} [c] at 16 8
\put {$ \scriptstyle \bullet$} [c] at 16 0
\setlinear \plot 10  0  10 8 13 12 16 8 16 0 10 8  /
\setlinear \plot 13 12 13 4  /
\put{$5{,}040$} [c] at 13 -2
\endpicture
\end{minipage}
\begin{minipage}{4cm}
\beginpicture
\setcoordinatesystem units   <1.5mm,2mm>
\setplotarea x from 0 to 16, y from -2 to 15
\put{647)} [l] at 2 12
\put {$ \scriptstyle \bullet$} [c] at 10 4
\put {$ \scriptstyle \bullet$} [c] at 11 0
\put {$ \scriptstyle \bullet$} [c] at 11 12
\put {$ \scriptstyle \bullet$} [c] at 12 4
\put {$ \scriptstyle \bullet$} [c] at 16 12
\put {$ \scriptstyle \bullet$} [c] at 10.3 6
\put {$ \scriptstyle \bullet$} [c] at 10.6 9
\setlinear \plot 16 12  11 0  10 4 11 12 12 4 11 0 /
\put{$5{,}040$} [c] at 13 -2
\endpicture
\end{minipage}
\begin{minipage}{4cm}
\beginpicture
\setcoordinatesystem units   <1.5mm,2mm>
\setplotarea x from 0 to 16, y from -2 to 15
\put{648)} [l] at 2 12
\put {$ \scriptstyle \bullet$} [c] at 10 8
\put {$ \scriptstyle \bullet$} [c] at 11 0
\put {$ \scriptstyle \bullet$} [c] at 11 12
\put {$ \scriptstyle \bullet$} [c] at 12 8
\put {$ \scriptstyle \bullet$} [c] at 16 0
\put {$ \scriptstyle \bullet$} [c] at 10.3 6
\put {$ \scriptstyle \bullet$} [c] at 10.6 3
\setlinear \plot 16 0  11 12  10 8 11 0 12 8 11 12 /
\put{$5{,}040$} [c] at 13 -2
\endpicture
\end{minipage}
$$
$$
\begin{minipage}{4cm}
\beginpicture
\setcoordinatesystem units   <1.5mm,2mm>
\setplotarea x from 0 to 16, y from -2 to 15
\put{649)} [l] at 2 12
\put {$ \scriptstyle \bullet$} [c] at 10 2
\put {$ \scriptstyle \bullet$} [c] at 11 7
\put {$ \scriptstyle \bullet$} [c] at 12 0
\put {$ \scriptstyle \bullet$} [c] at 13 5
\put {$ \scriptstyle \bullet$} [c] at 12 12
\put {$ \scriptstyle \bullet$} [c] at 14 10
\put {$ \scriptstyle \bullet$} [c] at 16 12
\setlinear \plot 16 12 12 0 14 10 12 12 10 2 12 0  /
\setlinear \plot 11 7 13 5 /
\put{$5{,}040$} [c] at 13 -2
\endpicture
\end{minipage}
\begin{minipage}{4cm}
\beginpicture
\setcoordinatesystem units   <1.5mm,2mm>
\setplotarea x from 0 to 16, y from -2 to 15
\put{650)} [l] at 2 12
\put {$ \scriptstyle \bullet$} [c] at 10 10
\put {$ \scriptstyle \bullet$} [c] at 11 5
\put {$ \scriptstyle \bullet$} [c] at 12 0
\put {$ \scriptstyle \bullet$} [c] at 13 7
\put {$ \scriptstyle \bullet$} [c] at 12 12
\put {$ \scriptstyle \bullet$} [c] at 14 2
\put {$ \scriptstyle \bullet$} [c] at 16 0
\setlinear \plot 16 0 12 12 14 2 12 0 10 10 12 12  /
\setlinear \plot 11 5 13 7 /
\put{$5{,}040$} [c] at 13 -2
\endpicture
\end{minipage}
\begin{minipage}{4cm}
\beginpicture
\setcoordinatesystem units   <1.5mm,2mm>
\setplotarea x from 0 to 16, y from -2 to 15
\put{651)} [l] at 2 12
\put {$ \scriptstyle \bullet$} [c] at 10 4
\put {$ \scriptstyle \bullet$} [c] at 10 12
\put {$ \scriptstyle \bullet$} [c] at 13 0
\put {$ \scriptstyle \bullet$} [c] at 13 4
\put {$ \scriptstyle \bullet$} [c] at 13 8
\put {$ \scriptstyle \bullet$} [c] at 13 12
\put {$ \scriptstyle \bullet$} [c] at 16 4
\setlinear \plot 13 12 16 4 13 0 13 12 10 4 13 0 /
\setlinear \plot 10 4 10 12 13 4 /
\put{$5{,}040$} [c] at 13 -2
\endpicture
\end{minipage}
\begin{minipage}{4cm}
\beginpicture
\setcoordinatesystem units   <1.5mm,2mm>
\setplotarea x from 0 to 16, y from -2 to 15
\put{652)} [l] at 2 12
\put {$ \scriptstyle \bullet$} [c] at 10 8
\put {$ \scriptstyle \bullet$} [c] at 10 0
\put {$ \scriptstyle \bullet$} [c] at 13 0
\put {$ \scriptstyle \bullet$} [c] at 13 4
\put {$ \scriptstyle \bullet$} [c] at 13 8
\put {$ \scriptstyle \bullet$} [c] at 13 12
\put {$ \scriptstyle \bullet$} [c] at 16 8
\setlinear \plot 13 0 16 8 13 12 13 0 10 8 13 12 /
\setlinear \plot 10 8 10 0 13 8 /
\put{$5{,}040$} [c] at 13 -2
\endpicture
\end{minipage}
\begin{minipage}{4cm}
\beginpicture
\setcoordinatesystem units   <1.5mm,2mm>
\setplotarea x from 0 to 16, y from -2 to 15
\put{653)} [l] at 2 12
\put {$ \scriptstyle \bullet$} [c] at 10 6
\put {$ \scriptstyle \bullet$} [c] at 10 12
\put {$ \scriptstyle \bullet$} [c] at 13  12
\put {$ \scriptstyle \bullet$} [c] at 13 0
\put {$ \scriptstyle \bullet$} [c] at 16 6
\put {$ \scriptstyle \bullet$} [c] at 16 9
\put {$ \scriptstyle \bullet$} [c] at 16 12
\setlinear \plot 10 12 10 6 13 0 16 6 13 12 10 6  /
\setlinear \plot 16 12 16 6 /
\put{$5{,}040$} [c] at 13 -2
\endpicture
\end{minipage}
\begin{minipage}{4cm}
\beginpicture
\setcoordinatesystem units   <1.5mm,2mm>
\setplotarea x from 0 to 16, y from -2 to 15
\put{654)} [l] at 2 12
\put {$ \scriptstyle \bullet$} [c] at 10 6
\put {$ \scriptstyle \bullet$} [c] at 10 0
\put {$ \scriptstyle \bullet$} [c] at 13  12
\put {$ \scriptstyle \bullet$} [c] at 13 0
\put {$ \scriptstyle \bullet$} [c] at 16 6
\put {$ \scriptstyle \bullet$} [c] at 16 3
\put {$ \scriptstyle \bullet$} [c] at 16 0
\setlinear \plot 10 0 10 6 13 12 16 6 13 0 10 6  /
\setlinear \plot 16 0 16 6 /
\put{$5{,}040$} [c] at 13 -2
\endpicture
\end{minipage}
$$

$$
\begin{minipage}{4cm}
\beginpicture
\setcoordinatesystem units   <1.5mm,2mm>
\setplotarea x from 0 to 16, y from -2 to 15
\put{655)} [l] at 2 12
\put {$ \scriptstyle \bullet$} [c] at 10 4
\put {$ \scriptstyle \bullet$} [c] at 10 12
\put {$ \scriptstyle \bullet$} [c] at 13 0
\put {$ \scriptstyle \bullet$} [c] at 13 8
\put {$ \scriptstyle \bullet$} [c] at 13 12
\put {$ \scriptstyle \bullet$} [c] at 16 4
\put {$ \scriptstyle \bullet$} [c] at 16 12
\setlinear \plot 16 12 16 4 13 0 10 4 10 12 16 4  /
\setlinear \plot 13 12 13 8  /
\put{$5{,}040$} [c] at 13 -2
\endpicture
\end{minipage}
\begin{minipage}{4cm}
\beginpicture
\setcoordinatesystem units   <1.5mm,2mm>
\setplotarea x from 0 to 16, y from -2 to 15
\put{656)} [l] at 2 12
\put {$ \scriptstyle \bullet$} [c] at 10 8
\put {$ \scriptstyle \bullet$} [c] at 10 0
\put {$ \scriptstyle \bullet$} [c] at 13 0
\put {$ \scriptstyle \bullet$} [c] at 13 4
\put {$ \scriptstyle \bullet$} [c] at 13 12
\put {$ \scriptstyle \bullet$} [c] at 16 8
\put {$ \scriptstyle \bullet$} [c] at 16 0
\setlinear \plot 16 0 16 8 13 12 10 8 10 0 16 8  /
\setlinear \plot 13 0 13 4  /
\put{$5{,}040$} [c] at 13 -2
\endpicture
\end{minipage}
\begin{minipage}{4cm}
\beginpicture
\setcoordinatesystem units   <1.5mm,2mm>
\setplotarea x from 0 to 16, y from -2 to 15
\put{657)} [l] at 2 12
\put {$ \scriptstyle \bullet$} [c] at 10 4
\put {$ \scriptstyle \bullet$} [c] at 10 12
\put {$ \scriptstyle \bullet$} [c] at 12 0
\put {$ \scriptstyle \bullet$} [c] at 14 4
\put {$ \scriptstyle \bullet$} [c] at 14 8
\put {$ \scriptstyle \bullet$} [c] at 14 12
\put {$ \scriptstyle \bullet$} [c] at 16 12
\setlinear \plot 16 12 14 4 12 0 10 4 10 12 14 4 14 12  10 4 /
\put{$5{,}040$} [c] at 13 -2
\endpicture
\end{minipage}
\begin{minipage}{4cm}
\beginpicture
\setcoordinatesystem units   <1.5mm,2mm>
\setplotarea x from 0 to 16, y from -2 to 15
\put{658)} [l] at 2 12
\put {$ \scriptstyle \bullet$} [c] at 10 8
\put {$ \scriptstyle \bullet$} [c] at 10 0
\put {$ \scriptstyle \bullet$} [c] at 12 12
\put {$ \scriptstyle \bullet$} [c] at 14 4
\put {$ \scriptstyle \bullet$} [c] at 14 8
\put {$ \scriptstyle \bullet$} [c] at 14 0
\put {$ \scriptstyle \bullet$} [c] at 16 0
\setlinear \plot 16 0 14 8 12 12 10 8 10 0 14 8 14 0  10 8 /
\put{$5{,}040$} [c] at 13 -2
\endpicture
\end{minipage}
\begin{minipage}{4cm}
\beginpicture
\setcoordinatesystem units   <1.5mm,2mm>
\setplotarea x from 0 to 16, y from -2 to 15
\put{659)} [l] at 2 12
\put {$ \scriptstyle \bullet$} [c] at 10 12
\put {$ \scriptstyle \bullet$} [c] at 13 0
\put {$ \scriptstyle \bullet$} [c] at 16 12
\put {$ \scriptstyle \bullet$} [c] at 14 12
\put {$ \scriptstyle \bullet$} [c] at 12 4
\put {$ \scriptstyle \bullet$} [c] at 10.7 9
\put {$ \scriptstyle \bullet$} [c] at 11.2 7
\setlinear \plot 10 12 13 0  16 12  /
\setlinear \plot 12 4 14 12  /
\put{$5{,}040$} [c] at 13 -2
\endpicture
\end{minipage}
\begin{minipage}{4cm}
\beginpicture
\setcoordinatesystem units   <1.5mm,2mm>
\setplotarea x from 0 to 16, y from -2 to 15
\put{660)} [l] at 2 12
\put {$ \scriptstyle \bullet$} [c] at 10 0
\put {$ \scriptstyle \bullet$} [c] at 13 12
\put {$ \scriptstyle \bullet$} [c] at 16 0
\put {$ \scriptstyle \bullet$} [c] at 14 0
\put {$ \scriptstyle \bullet$} [c] at 12 8
\put {$ \scriptstyle \bullet$} [c] at 10.7 3
\put {$ \scriptstyle \bullet$} [c] at 11.2 5
\setlinear \plot 10 0 13 12  16 0  /
\setlinear \plot 12 8 14 0  /
\put{$5{,}040$} [c] at 13 -2
\endpicture
\end{minipage}
$$
$$
\begin{minipage}{4cm}
\beginpicture
\setcoordinatesystem units   <1.5mm,2mm>
\setplotarea x from 0 to 16, y from -2 to 15
\put{661)} [l] at 2 12
\put {$ \scriptstyle \bullet$} [c] at 10 12
\put {$ \scriptstyle \bullet$} [c] at 12 8
\put {$ \scriptstyle \bullet$} [c] at 13 0
\put {$ \scriptstyle \bullet$} [c] at 13 4
\put {$ \scriptstyle \bullet$} [c] at 13 12
\put {$ \scriptstyle \bullet$} [c] at 14 8
\put {$ \scriptstyle \bullet$} [c] at 16 12
\setlinear \plot 16 12 13 0 13 4 12 8 10 12  /
\setlinear \plot 12 8 13 12 14 8 13 4 /
\put{$5{,}040$} [c] at 13 -2
\endpicture
\end{minipage}
\begin{minipage}{4cm}
\beginpicture
\setcoordinatesystem units   <1.5mm,2mm>
\setplotarea x from 0 to 16, y from -2 to 15
\put{662)} [l] at 2 12
\put {$ \scriptstyle \bullet$} [c] at 10 0
\put {$ \scriptstyle \bullet$} [c] at 12 4
\put {$ \scriptstyle \bullet$} [c] at 13 12
\put {$ \scriptstyle \bullet$} [c] at 13 8
\put {$ \scriptstyle \bullet$} [c] at 13 0
\put {$ \scriptstyle \bullet$} [c] at 14 4
\put {$ \scriptstyle \bullet$} [c] at 16 0
\setlinear \plot 16 0 13 12 13 8 12 4 10 0  /
\setlinear \plot 12 4 13 0 14 4 13 8 /
\put{$5{,}040$} [c] at 13 -2
\endpicture
\end{minipage}
\begin{minipage}{4cm}
\beginpicture
\setcoordinatesystem units   <1.5mm,2mm>
\setplotarea x from 0 to 16, y from -2 to 15
\put{663)} [l] at 2 12
\put {$ \scriptstyle \bullet$} [c] at  10 4
\put {$ \scriptstyle \bullet$} [c] at  10 8
\put {$ \scriptstyle \bullet$} [c] at  10 12
\put {$ \scriptstyle \bullet$} [c] at  14 4
\put {$ \scriptstyle \bullet$} [c] at  14 12
\put {$ \scriptstyle \bullet$} [c] at  12 0
\put {$ \scriptstyle \bullet$} [c] at  16 0
\setlinear \plot 16 0 14 12 14 4 12 0 10 4 10 12 14 4      /
\setlinear \plot 10 4 14 12 /
\put{$5{,}040   $} [c] at 13 -2
\endpicture
\end{minipage}
\begin{minipage}{4cm}
\beginpicture
\setcoordinatesystem units   <1.5mm,2mm>
\setplotarea x from 0 to 16, y from -2 to 15
\put{664)} [l] at 2 12
\put {$ \scriptstyle \bullet$} [c] at  10 4
\put {$ \scriptstyle \bullet$} [c] at  10 8
\put {$ \scriptstyle \bullet$} [c] at  10 0
\put {$ \scriptstyle \bullet$} [c] at  14 8
\put {$ \scriptstyle \bullet$} [c] at  14 0
\put {$ \scriptstyle \bullet$} [c] at  12 12
\put {$ \scriptstyle \bullet$} [c] at  16 12
\setlinear \plot 16 12 14 0 14 8 12 12 10 8 10 0 14 8      /
\setlinear \plot 10 8 14 0 /
\put{$5{,}040   $} [c] at 13 -2
\endpicture
\end{minipage}
\begin{minipage}{4cm}
\beginpicture
\setcoordinatesystem units   <1.5mm,2mm>
\setplotarea x from 0 to 16, y from -2 to 15
\put{665)} [l] at 2 12
\put {$ \scriptstyle \bullet$} [c] at  10 6
\put {$ \scriptstyle \bullet$} [c] at  10 12
\put {$ \scriptstyle \bullet$} [c] at  12 0
\put {$ \scriptstyle \bullet$} [c] at  12 6
\put {$ \scriptstyle \bullet$} [c] at  14 6
\put {$ \scriptstyle \bullet$} [c] at  14 12
\put {$ \scriptstyle \bullet$} [c] at  16 0
\setlinear \plot 10 6 16 0 14 12 14 6 12 0 12 6 14 12   /
\setlinear \plot 12 6 10 12 10 6 12 0 /
\put{$5{,}040   $} [c] at 13 -2
\endpicture
\end{minipage}
\begin{minipage}{4cm}
\beginpicture
\setcoordinatesystem units   <1.5mm,2mm>
\setplotarea x from 0 to 16, y from -2 to 15
\put{666)} [l] at 2 12
\put {$ \scriptstyle \bullet$} [c] at  10 6
\put {$ \scriptstyle \bullet$} [c] at  10 0
\put {$ \scriptstyle \bullet$} [c] at  12 12
\put {$ \scriptstyle \bullet$} [c] at  12 6
\put {$ \scriptstyle \bullet$} [c] at  14 6
\put {$ \scriptstyle \bullet$} [c] at  14 0
\put {$ \scriptstyle \bullet$} [c] at  16 12
\setlinear \plot 10 6 16 12 14 0 14 6 12 12 12 6 14 0   /
\setlinear \plot 12 6 10 0 10 6 12 12 /
\put{$5{,}040  $} [c] at 13 -2
\endpicture
\end{minipage}
$$

$$
\begin{minipage}{4cm}
\beginpicture
\setcoordinatesystem units   <1.5mm,2mm>
\setplotarea x from 0 to 16, y from -2 to 15
\put{667)} [l] at 2 12
\put {$ \scriptstyle \bullet$} [c] at  10 6
\put {$ \scriptstyle \bullet$} [c] at  10.5 9
\put {$ \scriptstyle \bullet$} [c] at  11 0
\put {$ \scriptstyle \bullet$} [c] at  11 12
\put {$ \scriptstyle \bullet$} [c] at  12 6
\put {$ \scriptstyle \bullet$} [c] at  16 0
\put {$ \scriptstyle \bullet$} [c] at  16 12
\setlinear \plot 12 6  16 12 16 0 11 12 12 6 11 0 10 6 11 12    /
\put{$5{,}040   $} [c] at 13 -2
\endpicture
\end{minipage}
\begin{minipage}{4cm}
\beginpicture
\setcoordinatesystem units   <1.5mm,2mm>
\setplotarea x from 0 to 16, y from -2 to 15
\put{668)} [l] at 2 12
\put {$ \scriptstyle \bullet$} [c] at  10 6
\put {$ \scriptstyle \bullet$} [c] at  10.5 3
\put {$ \scriptstyle \bullet$} [c] at  11 0
\put {$ \scriptstyle \bullet$} [c] at  11 12
\put {$ \scriptstyle \bullet$} [c] at  12 6
\put {$ \scriptstyle \bullet$} [c] at  16 0
\put {$ \scriptstyle \bullet$} [c] at  16 12
\setlinear \plot 12 6  16 0 16 12 11 0 12 6 11 12 10 6 11 0    /
\put{$5{,}040   $} [c] at 13 -2
\endpicture
\end{minipage}
\begin{minipage}{4cm}
\beginpicture
\setcoordinatesystem units   <1.5mm,2mm>
\setplotarea x from 0 to 16, y from -2 to 15
\put{669)} [l] at 2 12
\put {$ \scriptstyle \bullet$} [c] at  10 6
\put {$ \scriptstyle \bullet$} [c] at  10 12
\put {$ \scriptstyle \bullet$} [c] at  11.5 9
\put {$ \scriptstyle \bullet$} [c] at  13 0
\put {$ \scriptstyle \bullet$} [c] at  13 12
\put {$ \scriptstyle \bullet$} [c] at  16 6
\put {$ \scriptstyle \bullet$} [c] at  16 0
\setlinear \plot 10 12 10 6  13 0  16 6 13 12 10 6   /
\setlinear \plot 16 0 16 6 /
\put{$5{,}040   $} [c] at 13 -2
\endpicture
\end{minipage}
\begin{minipage}{4cm}
\beginpicture
\setcoordinatesystem units   <1.5mm,2mm>
\setplotarea x from 0 to 16, y from -2 to 15
\put{670)} [l] at 2 12
\put {$ \scriptstyle \bullet$} [c] at  10 6
\put {$ \scriptstyle \bullet$} [c] at  10 0
\put {$ \scriptstyle \bullet$} [c] at  11.5 3
\put {$ \scriptstyle \bullet$} [c] at  13 0
\put {$ \scriptstyle \bullet$} [c] at  13 12
\put {$ \scriptstyle \bullet$} [c] at  16 6
\put {$ \scriptstyle \bullet$} [c] at  16 12
\setlinear \plot 10 0 10 6  13 12  16 6 13 0 10 6   /
\setlinear \plot 16 12 16 6 /
\put{$5{,}040   $} [c] at 13 -2
\endpicture
\end{minipage}
\begin{minipage}{4cm}
\beginpicture
\setcoordinatesystem units   <1.5mm,2mm>
\setplotarea x from 0 to 16, y from -2 to 15
\put{671)} [l] at 2 12
\put {$ \scriptstyle \bullet$} [c] at  10 0
\put {$ \scriptstyle \bullet$} [c] at  10 12
\put {$ \scriptstyle \bullet$} [c] at  12 4
\put {$ \scriptstyle \bullet$} [c] at  12 8
\put {$ \scriptstyle \bullet$} [c] at  14 0
\put {$ \scriptstyle \bullet$} [c] at  14 12
\put {$ \scriptstyle \bullet$} [c] at  16 6
\setlinear \plot  10 0 10 12 12 8 12 4 14 0 16 6 14 12  12 8   /
\put{$5{,}040   $} [c] at 13 -2
\endpicture
\end{minipage}
\begin{minipage}{4cm}
\beginpicture
\setcoordinatesystem units   <1.5mm,2mm>
\setplotarea x from 0 to 16, y from -2 to 15
\put{672)} [l] at 2 12
\put {$ \scriptstyle \bullet$} [c] at  10 0
\put {$ \scriptstyle \bullet$} [c] at  10 12
\put {$ \scriptstyle \bullet$} [c] at  12 4
\put {$ \scriptstyle \bullet$} [c] at  12 8
\put {$ \scriptstyle \bullet$} [c] at  14 0
\put {$ \scriptstyle \bullet$} [c] at  14 12
\put {$ \scriptstyle \bullet$} [c] at  16 6
\setlinear \plot  10 12 10 0 12 4 12 8 14 12 16 6 14 0  12 4   /
\put{$5{,}040   $} [c] at 13 -2
\endpicture
\end{minipage}
$$
$$
\begin{minipage}{4cm}
\beginpicture
\setcoordinatesystem units   <1.5mm,2mm>
\setplotarea x from 0 to 16, y from -2 to 15
\put{673)} [l] at 2 12
\put {$ \scriptstyle \bullet$} [c] at  10 0
\put {$ \scriptstyle \bullet$} [c] at  10 4
\put {$ \scriptstyle \bullet$} [c] at  10 8
\put {$ \scriptstyle \bullet$} [c] at  10 12
\put {$ \scriptstyle \bullet$} [c] at  14 8
\put {$ \scriptstyle \bullet$} [c] at  14 12
\put {$ \scriptstyle \bullet$} [c] at  16 0
\setlinear \plot 10 4 10 12 16  0 14 12 14 8 10 0 10 4 14 12  /
\put{$5{,}040   $} [c] at 13 -2
\endpicture
\end{minipage}
\begin{minipage}{4cm}
\beginpicture
\setcoordinatesystem units   <1.5mm,2mm>
\setplotarea x from 0 to 16, y from -2 to 15
\put{674)} [l] at 2 12
\put {$ \scriptstyle \bullet$} [c] at  10 0
\put {$ \scriptstyle \bullet$} [c] at  10 4
\put {$ \scriptstyle \bullet$} [c] at  10 8
\put {$ \scriptstyle \bullet$} [c] at  10 12
\put {$ \scriptstyle \bullet$} [c] at  14 4
\put {$ \scriptstyle \bullet$} [c] at  14 0
\put {$ \scriptstyle \bullet$} [c] at  16 12
\setlinear \plot 10 8 10 0 16  12 14 0 14 4 10 12 10 8 14 0  /
\put{$5{,}040   $} [c] at 13 -2
\endpicture
\end{minipage}
\begin{minipage}{4cm}
\beginpicture
\setcoordinatesystem units   <1.5mm,2mm>
\setplotarea x from 0 to 16, y from -2 to 15
\put{675)} [l] at 2 12
\put {$ \scriptstyle \bullet$} [c] at  10 8
\put {$ \scriptstyle \bullet$} [c] at  10 12
\put {$ \scriptstyle \bullet$} [c] at  12 4
\put {$ \scriptstyle \bullet$} [c] at  12 0
\put {$ \scriptstyle \bullet$} [c] at  14 8
\put {$ \scriptstyle \bullet$} [c] at  14 12
\put {$ \scriptstyle \bullet$} [c] at  16 0
\setlinear \plot  10 12  10 8  12 4 14 8 14 12 16 0 /
\setlinear \plot  12 0 12 4 /
\put{$5{,}040  $} [c] at 13 -2
\endpicture
\end{minipage}
\begin{minipage}{4cm}
\beginpicture
\setcoordinatesystem units   <1.5mm,2mm>
\setplotarea x from 0 to 16, y from -2 to 15
\put{676)} [l] at 2 12
\put {$ \scriptstyle \bullet$} [c] at  10 4
\put {$ \scriptstyle \bullet$} [c] at  10 0
\put {$ \scriptstyle \bullet$} [c] at  12 12
\put {$ \scriptstyle \bullet$} [c] at  12 8
\put {$ \scriptstyle \bullet$} [c] at  14 4
\put {$ \scriptstyle \bullet$} [c] at  14 0
\put {$ \scriptstyle \bullet$} [c] at  16 12
\setlinear \plot  10 0  10 4  12 8 14 4 14 0 16 12 /
\setlinear \plot  12 12 12 8 /
\put{$5{,}040  $} [c] at 13 -2
\endpicture
\end{minipage}
\begin{minipage}{4cm}
\beginpicture
\setcoordinatesystem units   <1.5mm,2mm>
\setplotarea x from 0 to 16, y from -2 to 15
\put{677)} [l] at 2 12
\put {$ \scriptstyle \bullet$} [c] at  10  0
\put {$ \scriptstyle \bullet$} [c] at  10 4
\put {$ \scriptstyle \bullet$} [c] at  10 8
\put {$ \scriptstyle \bullet$} [c] at  10  12
\put {$ \scriptstyle \bullet$} [c] at  16 0
\put {$ \scriptstyle \bullet$} [c] at  16 6
\put {$ \scriptstyle \bullet$} [c] at  16 12
\setlinear \plot 10 12 10 0 16 12 16 0   /
\put{$5{,}040   $} [c] at 13 -2
\endpicture
\end{minipage}
\begin{minipage}{4cm}
\beginpicture
\setcoordinatesystem units   <1.5mm,2mm>
\setplotarea x from 0 to 16, y from -2 to 15
\put{678)} [l] at 2 12
\put {$ \scriptstyle \bullet$} [c] at  10  0
\put {$ \scriptstyle \bullet$} [c] at  10 4
\put {$ \scriptstyle \bullet$} [c] at  10 8
\put {$ \scriptstyle \bullet$} [c] at  10  12
\put {$ \scriptstyle \bullet$} [c] at  16 0
\put {$ \scriptstyle \bullet$} [c] at  16 6
\put {$ \scriptstyle \bullet$} [c] at  16 12
\setlinear \plot 10 0 10 12 16 0 16 12   /
\put{$5{,}040   $} [c] at 13 -2
\endpicture
\end{minipage}
$$

$$
\begin{minipage}{4cm}
\beginpicture
\setcoordinatesystem units   <1.5mm,2mm>
\setplotarea x from 0 to 16, y from -2 to 15
\put{679)} [l] at 2 12
\put {$ \scriptstyle \bullet$} [c] at  10 6
\put {$ \scriptstyle \bullet$} [c] at  12 0
\put {$ \scriptstyle \bullet$} [c] at  12 12
\put {$ \scriptstyle \bullet$} [c] at  14 6
\put {$ \scriptstyle \bullet$} [c] at  16 0
\put {$ \scriptstyle \bullet$} [c] at  16 6
\put {$ \scriptstyle \bullet$} [c] at  16 12
\setlinear \plot  16 0  16 12 14 6 12 12 10 6 12 0 14 6 /
\setlinear \plot  12 12 16 6 /
\put{$5{,}040  $} [c] at 13 -2
\endpicture
\end{minipage}
\begin{minipage}{4cm}
\beginpicture
\setcoordinatesystem units   <1.5mm,2mm>
\setplotarea x from 0 to 16, y from -2 to 15
\put{680)} [l] at 2 12
\put {$ \scriptstyle \bullet$} [c] at  10 6
\put {$ \scriptstyle \bullet$} [c] at  12 0
\put {$ \scriptstyle \bullet$} [c] at  12 12
\put {$ \scriptstyle \bullet$} [c] at  14 6
\put {$ \scriptstyle \bullet$} [c] at  16 0
\put {$ \scriptstyle \bullet$} [c] at  16 6
\put {$ \scriptstyle \bullet$} [c] at  16 12
\setlinear \plot  16 12  16 0 14 6 12 0 10 6 12 12 14 6 /
\setlinear \plot  12 0 16 6 /
\put{$5{,}040  $} [c] at 13 -2
\endpicture
\end{minipage}
\begin{minipage}{4cm}
\beginpicture
\setcoordinatesystem units   <1.5mm,2mm>
\setplotarea x from 0 to 16, y from -2 to 15
\put{681)} [l] at 2 12
\put {$ \scriptstyle \bullet$} [c] at  10 12
\put {$ \scriptstyle \bullet$} [c] at  10 8
\put {$ \scriptstyle \bullet$} [c] at  10 4
\put {$ \scriptstyle \bullet$} [c] at  10 0
\put {$ \scriptstyle \bullet$} [c] at  13 0
\put {$ \scriptstyle \bullet$} [c] at  13 12
\put {$ \scriptstyle \bullet$} [c] at  16 12
\setlinear \plot  10 12 10 0   /
\setlinear \plot   16 12 13 0 13 12  10 4 /
\setlinear \plot  10 8 13 0  /
\put{$5{,}040$} [c] at 13 -2
\endpicture
\end{minipage}
\begin{minipage}{4cm}
\beginpicture
\setcoordinatesystem units   <1.5mm,2mm>
\setplotarea x from 0 to 16, y from -2 to 15
\put{682)} [l] at 2 12
\put {$ \scriptstyle \bullet$} [c] at  10 12
\put {$ \scriptstyle \bullet$} [c] at  10 8
\put {$ \scriptstyle \bullet$} [c] at  10 4
\put {$ \scriptstyle \bullet$} [c] at  10 0
\put {$ \scriptstyle \bullet$} [c] at  13 0
\put {$ \scriptstyle \bullet$} [c] at  13 12
\put {$ \scriptstyle \bullet$} [c] at  16 0
\setlinear \plot  10 12 10 0   /
\setlinear \plot   16 0 13 12 13 0  10 8 /
\setlinear \plot  10 4 13 12  /
\put{$5{,}040$} [c] at 13 -2
\endpicture
\end{minipage}
\begin{minipage}{4cm}
\beginpicture
\setcoordinatesystem units   <1.5mm,2mm>
\setplotarea x from 0 to 16, y from -2 to 15
\put{683)} [l] at 2 12
\put {$ \scriptstyle \bullet$} [c] at  10 12
\put {$ \scriptstyle \bullet$} [c] at  13 0
\put {$ \scriptstyle \bullet$} [c] at  13 4
\put {$ \scriptstyle \bullet$} [c] at  13 8
\put {$ \scriptstyle \bullet$} [c] at  16 0
\put {$ \scriptstyle \bullet$} [c] at  16 12
\put {$ \scriptstyle \bullet$} [c] at  13 12
\setlinear \plot  10 12 13 8 /
\setlinear \plot  13 0 13 12 16 0 16 12 /
\put{$5{,}040$} [c] at 13 -2
\endpicture
\end{minipage}
\begin{minipage}{4cm}
\beginpicture
\setcoordinatesystem units   <1.5mm,2mm>
\setplotarea x from 0 to 16, y from -2 to 15
\put{684)} [l] at 2 12
\put {$ \scriptstyle \bullet$} [c] at  10 0
\put {$ \scriptstyle \bullet$} [c] at  13 0
\put {$ \scriptstyle \bullet$} [c] at  13 4
\put {$ \scriptstyle \bullet$} [c] at  13 8
\put {$ \scriptstyle \bullet$} [c] at  16 0
\put {$ \scriptstyle \bullet$} [c] at  16 12
\put {$ \scriptstyle \bullet$} [c] at  13 12
\setlinear \plot  10 0 13 4 /
\setlinear \plot  13 12 13 0 16 12 16 0 /
\put{$5{,}040$} [c] at 13 -2
\endpicture
\end{minipage}
$$
$$
\begin{minipage}{4cm}
\beginpicture
\setcoordinatesystem units   <1.5mm,2mm>
\setplotarea x from 0 to 16, y from -2 to 15
\put{685)} [l] at 2 12
\put {$ \scriptstyle \bullet$} [c] at  10 12
\put {$ \scriptstyle \bullet$} [c] at  13 0
\put {$ \scriptstyle \bullet$} [c] at  13 4
\put {$ \scriptstyle \bullet$} [c] at  13 8
\put {$ \scriptstyle \bullet$} [c] at  13 12
\put {$ \scriptstyle \bullet$} [c] at  16 0
\put {$ \scriptstyle \bullet$} [c] at  16 12
\setlinear \plot  13 12 13 0 16 12 16 0 13 8  /
\setlinear \plot  10 12 13 4  /
\put{$5{,}040$} [c] at 13 -2
\endpicture
\end{minipage}
\begin{minipage}{4cm}
\beginpicture
\setcoordinatesystem units   <1.5mm,2mm>
\setplotarea x from 0 to 16, y from -2 to 15
\put{686)} [l] at 2 12
\put {$ \scriptstyle \bullet$} [c] at  10 0
\put {$ \scriptstyle \bullet$} [c] at  13 0
\put {$ \scriptstyle \bullet$} [c] at  13 4
\put {$ \scriptstyle \bullet$} [c] at  13 8
\put {$ \scriptstyle \bullet$} [c] at  13 12
\put {$ \scriptstyle \bullet$} [c] at  16 0
\put {$ \scriptstyle \bullet$} [c] at  16 12
\setlinear \plot  13 0 13 12 16 0 16 12 13 4  /
\setlinear \plot  10 0 13 8  /
\put{$5{,}040$} [c] at 13 -2
\endpicture
\end{minipage}
\begin{minipage}{4cm}
\beginpicture
\setcoordinatesystem units   <1.5mm,2mm>
\setplotarea x from 0 to 16, y from -2 to 15
\put{687)} [l] at 2 12
\put {$ \scriptstyle \bullet$} [c] at  16 12
\put {$ \scriptstyle \bullet$} [c] at  16 0
\put {$ \scriptstyle \bullet$} [c] at  13 12
\put {$ \scriptstyle \bullet$} [c] at  10 0
\put {$ \scriptstyle \bullet$} [c] at  10 4
\put {$ \scriptstyle \bullet$} [c] at  10 8
\put {$ \scriptstyle \bullet$} [c] at  10 12
\setlinear \plot 16 0 10 12 10 0  16 12 16 0 13 12 10 4  /
\put{$5{,}040$} [c] at 13 -2
\endpicture
\end{minipage}
\begin{minipage}{4cm}
\beginpicture
\setcoordinatesystem units   <1.5mm,2mm>
\setplotarea x from 0 to 16, y from -2 to 15
\put{688)} [l] at 2 12
\put {$ \scriptstyle \bullet$} [c] at  16 12
\put {$ \scriptstyle \bullet$} [c] at  16 0
\put {$ \scriptstyle \bullet$} [c] at  13 0
\put {$ \scriptstyle \bullet$} [c] at  10 0
\put {$ \scriptstyle \bullet$} [c] at  10 4
\put {$ \scriptstyle \bullet$} [c] at  10 8
\put {$ \scriptstyle \bullet$} [c] at  10 12
\setlinear \plot 16 12 10 0 10 12  16 0 16 12 13 0 10 8  /
\put{$5{,}040$} [c] at 13 -2
\endpicture
\end{minipage}
\begin{minipage}{4cm}
\beginpicture
\setcoordinatesystem units   <1.5mm,2mm>
\setplotarea x from 0 to 16, y from -2 to 15
\put{689)} [l] at 2 12
\put {$ \scriptstyle \bullet$} [c] at  10 0
\put {$ \scriptstyle \bullet$} [c] at  10 6
\put {$ \scriptstyle \bullet$} [c] at  10 12
\put {$ \scriptstyle \bullet$} [c] at  13 12
\put {$ \scriptstyle \bullet$} [c] at  16 0
\put {$ \scriptstyle \bullet$} [c] at  16 6
\put {$ \scriptstyle \bullet$} [c] at  16 12
\setlinear \plot  10 0  10 12 16 0 16  12 /
\setlinear \plot  10 6 13 12 16 6 /
\put{$5{,}040$} [c] at 13 -2
\endpicture
\end{minipage}
\begin{minipage}{4cm}
\beginpicture
\setcoordinatesystem units   <1.5mm,2mm>
\setplotarea x from 0 to 16, y from -2 to 15
\put{690)} [l] at 2 12
\put {$ \scriptstyle \bullet$} [c] at  10 0
\put {$ \scriptstyle \bullet$} [c] at  10 6
\put {$ \scriptstyle \bullet$} [c] at  10 12
\put {$ \scriptstyle \bullet$} [c] at  13 0
\put {$ \scriptstyle \bullet$} [c] at  16 0
\put {$ \scriptstyle \bullet$} [c] at  16 6
\put {$ \scriptstyle \bullet$} [c] at  16 12
\setlinear \plot  10 12  10 0 16 12 16  0 /
\setlinear \plot  10 6 13 0 16 6 /
\put{$5{,}040$} [c] at 13 -2
\endpicture
\end{minipage}
$$

$$
\begin{minipage}{4cm}
\beginpicture
\setcoordinatesystem units   <1.5mm,2mm>
\setplotarea x from 0 to 16, y from -2 to 15
\put{691)} [l] at 2 12
\put {$ \scriptstyle \bullet$} [c] at  10 12
\put {$ \scriptstyle \bullet$} [c] at  13 0
\put {$ \scriptstyle \bullet$} [c] at  13 6
\put {$ \scriptstyle \bullet$} [c] at  13 12
\put {$ \scriptstyle \bullet$} [c] at  16 0
\put {$ \scriptstyle \bullet$} [c] at  16 6
\put {$ \scriptstyle \bullet$} [c] at  16 12
\setlinear \plot 10 12  13 6 13  12 16 0  16 12 13 6 13 0 16  6  /
\put{$5{,}040$} [c] at 13 -2
\endpicture
\end{minipage}
\begin{minipage}{4cm}
\beginpicture
\setcoordinatesystem units   <1.5mm,2mm>
\setplotarea x from 0 to 16, y from -2 to 15
\put{692)} [l] at 2 12
\put {$ \scriptstyle \bullet$} [c] at  10 0
\put {$ \scriptstyle \bullet$} [c] at  13 0
\put {$ \scriptstyle \bullet$} [c] at  13 6
\put {$ \scriptstyle \bullet$} [c] at  13 12
\put {$ \scriptstyle \bullet$} [c] at  16 0
\put {$ \scriptstyle \bullet$} [c] at  16 6
\put {$ \scriptstyle \bullet$} [c] at  16 12
\setlinear \plot 10 0  13 6 13  0 16 12  16 0 13 6 13 12 16  6  /
\put{$5{,}040$} [c] at 13 -2
\endpicture
\end{minipage}
\begin{minipage}{4cm}
\beginpicture
\setcoordinatesystem units   <1.5mm,2mm>
\setplotarea x from 0 to 16, y from -2 to 15
\put{693)} [l] at 2 12
\put {$ \scriptstyle \bullet$} [c] at  10 0
\put {$ \scriptstyle \bullet$} [c] at  10 12
\put {$ \scriptstyle \bullet$} [c] at  14 6
\put {$ \scriptstyle \bullet$} [c] at  14 12
\put {$ \scriptstyle \bullet$} [c] at  15 0
\put {$ \scriptstyle \bullet$} [c] at  16 6
\put {$ \scriptstyle \bullet$} [c] at  16 12
\setlinear \plot 16 12 10 0 10 12  14 6 16 12 16 6 15 0 14 6 14 12 16 6 /
\setlinear \plot  10 12 14 6  /
\put{$5{,}040$} [c] at 13 -2
\endpicture
\end{minipage}
\begin{minipage}{4cm}
\beginpicture
\setcoordinatesystem units   <1.5mm,2mm>
\setplotarea x from 0 to 16, y from -2 to 15
\put{694)} [l] at 2 12
\put {$ \scriptstyle \bullet$} [c] at  10 0
\put {$ \scriptstyle \bullet$} [c] at  10 12
\put {$ \scriptstyle \bullet$} [c] at  14 6
\put {$ \scriptstyle \bullet$} [c] at  14 0
\put {$ \scriptstyle \bullet$} [c] at  15 12
\put {$ \scriptstyle \bullet$} [c] at  16 6
\put {$ \scriptstyle \bullet$} [c] at  16 0
\setlinear \plot 16 0 10 12 10 0  14 6 16 0 16 6 15 12 14 6 14 0 16 6 /
\setlinear \plot  10 0 14 6  /
\put{$5{,}040$} [c] at 13 -2
\endpicture
\end{minipage}
\begin{minipage}{4cm}
\beginpicture
\setcoordinatesystem units   <1.5mm,2mm>
\setplotarea x from 0 to 16, y from -2 to 15
\put{695)} [l] at 2 12
\put {$ \scriptstyle \bullet$} [c] at 10 4
\put {$ \scriptstyle \bullet$} [c] at 12  0
\put {$ \scriptstyle \bullet$} [c] at 12 4
\put {$ \scriptstyle \bullet$} [c] at 12 12
\put {$ \scriptstyle \bullet$} [c] at 14 4
\put {$ \scriptstyle \bullet$} [c] at 14 8
\put {$ \scriptstyle \bullet$} [c] at 16 12
\setlinear \plot 12 0 10 4 12 12 12 0  14 4 14 8 12 12  /
\setlinear \plot 14 8 16 12  /
\put{$2{,}520$} [c] at 13 -2
\endpicture
\end{minipage}
\begin{minipage}{4cm}
\beginpicture
\setcoordinatesystem units   <1.5mm,2mm>
\setplotarea x from 0 to 16, y from -2 to 15
\put{696)} [l] at 2 12
\put {$ \scriptstyle \bullet$} [c] at 10 8
\put {$ \scriptstyle \bullet$} [c] at 12  0
\put {$ \scriptstyle \bullet$} [c] at 12 8
\put {$ \scriptstyle \bullet$} [c] at 12 12
\put {$ \scriptstyle \bullet$} [c] at 14 4
\put {$ \scriptstyle \bullet$} [c] at 14 8
\put {$ \scriptstyle \bullet$} [c] at 16 0
\setlinear \plot 12 12 10 8 12 0 12 12  14 8 14 4 12 0  /
\setlinear \plot 14 4 16 0  /
\put{$2{,}520$} [c] at 13 -2
\endpicture
\end{minipage}
$$
$$
\begin{minipage}{4cm}
\beginpicture
\setcoordinatesystem units   <1.5mm,2mm>
\setplotarea x from 0 to 16, y from -2 to 15
\put{697)} [l] at 2 12
\put {$ \scriptstyle \bullet$} [c] at 10 8
\put {$ \scriptstyle \bullet$} [c] at 12 0
\put {$ \scriptstyle \bullet$} [c] at 12 4
\put {$ \scriptstyle \bullet$} [c] at 12 12
\put {$ \scriptstyle \bullet$} [c] at 14 8
\put {$ \scriptstyle \bullet$} [c] at 16 8
\put {$ \scriptstyle \bullet$} [c] at 16 12
\setlinear \plot  12 4 10 8 12 12 14 8 12 4 12 0 16 8 16 12   /
\put{$2{,}520$} [c] at 13 -2
\endpicture
\end{minipage}
\begin{minipage}{4cm}
\beginpicture
\setcoordinatesystem units   <1.5mm,2mm>
\setplotarea x from 0 to 16, y from -2 to 15
\put{698)} [l] at 2 12
\put {$ \scriptstyle \bullet$} [c] at 10 4
\put {$ \scriptstyle \bullet$} [c] at 12 12
\put {$ \scriptstyle \bullet$} [c] at 12 8
\put {$ \scriptstyle \bullet$} [c] at 12 0
\put {$ \scriptstyle \bullet$} [c] at 14 4
\put {$ \scriptstyle \bullet$} [c] at 16 4
\put {$ \scriptstyle \bullet$} [c] at 16 0
\setlinear \plot  12 8 10 4 12 0 14 4 12 8 12 12 16 4 16 0   /
\put{$2{,}520$} [c] at 13 -2
\endpicture
\end{minipage}
\begin{minipage}{4cm}
\beginpicture
\setcoordinatesystem units   <1.5mm,2mm>
\setplotarea x from 0 to 16, y from -2 to 15
\put{699)} [l] at 2 12
\put {$ \scriptstyle \bullet$} [c] at 10 4
\put {$ \scriptstyle \bullet$} [c] at 12 0
\put {$ \scriptstyle \bullet$} [c] at 12 8
\put {$ \scriptstyle \bullet$} [c] at 12 12
\put {$ \scriptstyle \bullet$} [c] at 14 4
\put {$ \scriptstyle \bullet$} [c] at 16 4
\put {$ \scriptstyle \bullet$} [c] at 16 12
\setlinear \plot  16 12 16 4  12  0 10 4 12 8 12 12 /
\setlinear \plot  12 8 14 4 12 0  /
\put{$2{,}520$} [c] at 13 -2
\endpicture
\end{minipage}
\begin{minipage}{4cm}
\beginpicture
\setcoordinatesystem units   <1.5mm,2mm>
\setplotarea x from 0 to 16, y from -2 to 15
\put{700)} [l] at 2 12
\put {$ \scriptstyle \bullet$} [c] at 10 8
\put {$ \scriptstyle \bullet$} [c] at 12 0
\put {$ \scriptstyle \bullet$} [c] at 12 4
\put {$ \scriptstyle \bullet$} [c] at 12 12
\put {$ \scriptstyle \bullet$} [c] at 14 8
\put {$ \scriptstyle \bullet$} [c] at 16 8
\put {$ \scriptstyle \bullet$} [c] at 16 0
\setlinear \plot  16 0 16 8  12  12 10 8 12 4 12 0 /
\setlinear \plot  12 4 14 8 12 12  /
\put{$2{,}520$} [c] at 13 -2
\endpicture
\end{minipage}
\begin{minipage}{4cm}
\beginpicture
\setcoordinatesystem units   <1.5mm,2mm>
\setplotarea x from 0 to 16, y from -2 to 15
\put{701)} [l] at 2 12
\put {$ \scriptstyle \bullet$} [c] at 10 12
\put {$ \scriptstyle \bullet$} [c] at 12 4
\put {$ \scriptstyle \bullet$} [c] at 12 12
\put {$ \scriptstyle \bullet$} [c] at 14 0
\put {$ \scriptstyle \bullet$} [c] at 14 8
\put {$ \scriptstyle \bullet$} [c] at 14 12
\put {$ \scriptstyle \bullet$} [c] at 16 4
\setlinear \plot 10 12 12 4 14 0 16 4 14 8 14 12  /
\setlinear \plot 12 12 12 4 14 8  /
\put{$2{,}520$} [c] at 13 -2
\endpicture
\end{minipage}
\begin{minipage}{4cm}
\beginpicture
\setcoordinatesystem units   <1.5mm,2mm>
\setplotarea x from 0 to 16, y from -2 to 15
\put{702)} [l] at 2 12
\put {$ \scriptstyle \bullet$} [c] at 10 0
\put {$ \scriptstyle \bullet$} [c] at 12 8
\put {$ \scriptstyle \bullet$} [c] at 12 0
\put {$ \scriptstyle \bullet$} [c] at 14 0
\put {$ \scriptstyle \bullet$} [c] at 14 4
\put {$ \scriptstyle \bullet$} [c] at 14 12
\put {$ \scriptstyle \bullet$} [c] at 16 8
\setlinear \plot 10 0 12 8 14 12 16 8 14 4 14 0  /
\setlinear \plot 12 0 12 8 14 4  /
\put{$2{,}520$} [c] at 13 -2
\endpicture
\end{minipage}
$$

$$
\begin{minipage}{4cm}
\beginpicture
\setcoordinatesystem units   <1.5mm,2mm>
\setplotarea x from 0 to 16, y from -2 to 15
\put{703)} [l] at 2 12
\put {$ \scriptstyle \bullet$} [c] at 10 12
\put {$ \scriptstyle \bullet$} [c] at 12 4
\put {$ \scriptstyle \bullet$} [c] at 12 8
\put {$ \scriptstyle \bullet$} [c] at 14 0
\put {$ \scriptstyle \bullet$} [c] at 14 12
\put {$ \scriptstyle \bullet$} [c] at 16 4
\put {$ \scriptstyle \bullet$} [c] at 16 12
\setlinear \plot 16 12 16 4 14 0 12 4 12 8 14 12  /
\setlinear \plot 10 12 12 8 /
\put{$2{,}520$} [c] at 13 -2
\endpicture
\end{minipage}
\begin{minipage}{4cm}
\beginpicture
\setcoordinatesystem units   <1.5mm,2mm>
\setplotarea x from 0 to 16, y from -2 to 15
\put{704)} [l] at 2 12
\put {$ \scriptstyle \bullet$} [c] at 10 0
\put {$ \scriptstyle \bullet$} [c] at 12 4
\put {$ \scriptstyle \bullet$} [c] at 12 8
\put {$ \scriptstyle \bullet$} [c] at 14 12
\put {$ \scriptstyle \bullet$} [c] at 14 0
\put {$ \scriptstyle \bullet$} [c] at 16 8
\put {$ \scriptstyle \bullet$} [c] at 16 0
\setlinear \plot 16 0 16 8 14 12 12 8 12 4 14 0  /
\setlinear \plot 10 0 12 4 /
\put{$2{,}520$} [c] at 13 -2
\endpicture
\end{minipage}
\begin{minipage}{4cm}
\beginpicture
\setcoordinatesystem units   <1.5mm,2mm>
\setplotarea x from 0 to 16, y from -2 to 15
\put{705)} [l] at 2 12
\put {$ \scriptstyle \bullet$} [c] at 10 4
\put {$ \scriptstyle \bullet$} [c] at 10 8
\put {$ \scriptstyle \bullet$} [c] at 10  12
\put {$ \scriptstyle \bullet$} [c] at 13 0
\put {$ \scriptstyle \bullet$} [c] at 13 4
\put {$ \scriptstyle \bullet$} [c] at 13 12
\put {$ \scriptstyle \bullet$} [c] at 16 12
\setlinear \plot  16 12 13 0 10 4 10 12 13 4 13 12 10 8   /
\setlinear \plot   13 0 13 4 /
\put{$2{,}520$} [c] at 13 -2
\endpicture
\end{minipage}
\begin{minipage}{4cm}
\beginpicture
\setcoordinatesystem units   <1.5mm,2mm>
\setplotarea x from 0 to 16, y from -2 to 15
\put{706)} [l] at 2 12
\put {$ \scriptstyle \bullet$} [c] at 10 4
\put {$ \scriptstyle \bullet$} [c] at 10 8
\put {$ \scriptstyle \bullet$} [c] at 10  0
\put {$ \scriptstyle \bullet$} [c] at 13 0
\put {$ \scriptstyle \bullet$} [c] at 13 8
\put {$ \scriptstyle \bullet$} [c] at 13 12
\put {$ \scriptstyle \bullet$} [c] at 16 0
\setlinear \plot  16 0 13 12 10 8 10 0 13 8 13 0 10 4   /
\setlinear \plot   13 12 13 8 /
\put{$2{,}520$} [c] at 13 -2
\endpicture
\end{minipage}
\begin{minipage}{4cm}
\beginpicture
\setcoordinatesystem units   <1.5mm,2mm>
\setplotarea x from 0 to 16, y from -2 to 15
\put{707)} [l] at 2 12
\put {$ \scriptstyle \bullet$} [c] at 10 4
\put {$ \scriptstyle \bullet$} [c] at 10 8
\put {$ \scriptstyle \bullet$} [c] at 10  12
\put {$ \scriptstyle \bullet$} [c] at 13 0
\put {$ \scriptstyle \bullet$} [c] at 13 4
\put {$ \scriptstyle \bullet$} [c] at 13 12
\put {$ \scriptstyle \bullet$} [c] at 16 12
\setlinear \plot  16 12 13 0 10 4 10 12    /
\setlinear \plot   13 0 13 12 10 4 /
\setlinear \plot   10 8  13  4 /
\put{$2{,}520$} [c] at 13 -2
\endpicture
\end{minipage}
\begin{minipage}{4cm}
\beginpicture
\setcoordinatesystem units   <1.5mm,2mm>
\setplotarea x from 0 to 16, y from -2 to 15
\put{708)} [l] at 2 12
\put {$ \scriptstyle \bullet$} [c] at 10 4
\put {$ \scriptstyle \bullet$} [c] at 10 8
\put {$ \scriptstyle \bullet$} [c] at 10  0
\put {$ \scriptstyle \bullet$} [c] at 13 0
\put {$ \scriptstyle \bullet$} [c] at 13 8
\put {$ \scriptstyle \bullet$} [c] at 13 12
\put {$ \scriptstyle \bullet$} [c] at 16 0
\setlinear \plot  16 0 13 12 10 8 10 0    /
\setlinear \plot   13 12 13 0 10 8 /
\setlinear \plot   10 4  13  8 /
\put{$2{,}520$} [c] at 13 -2
\endpicture
\end{minipage}
$$
$$
\begin{minipage}{4cm}
\beginpicture
\setcoordinatesystem units   <1.5mm,2mm>
\setplotarea x from 0 to 16, y from -2 to 15
\put{709)} [l] at 2 12
\put {$ \scriptstyle \bullet$} [c] at  10 8
\put {$ \scriptstyle \bullet$} [c] at  11 0
\put {$ \scriptstyle \bullet$} [c] at  11 4
\put {$ \scriptstyle \bullet$} [c] at  11 12
\put {$ \scriptstyle \bullet$} [c] at  12 8
\put {$ \scriptstyle \bullet$} [c] at  16 0
\put {$ \scriptstyle \bullet$} [c] at  16 12
\setlinear \plot  11 0 11 4 10 8  11 12 12  8 11 4 16 12 16 0   /
\put{$2{,}520  $} [c] at 13 -2
\endpicture
\end{minipage}
\begin{minipage}{4cm}
\beginpicture
\setcoordinatesystem units   <1.5mm,2mm>
\setplotarea x from 0 to 16, y from -2 to 15
\put{710)} [l] at 2 12
\put {$ \scriptstyle \bullet$} [c] at  10 4
\put {$ \scriptstyle \bullet$} [c] at  11 12
\put {$ \scriptstyle \bullet$} [c] at  11 8
\put {$ \scriptstyle \bullet$} [c] at  11 0
\put {$ \scriptstyle \bullet$} [c] at  12 4
\put {$ \scriptstyle \bullet$} [c] at  16 0
\put {$ \scriptstyle \bullet$} [c] at  16 12
\setlinear \plot  11 12 11 8 10 4  11 0 12  4 11 8 16 0 16 12   /
\put{$2{,}520  $} [c] at 13 -2
\endpicture
\end{minipage}
\begin{minipage}{4cm}
\beginpicture
\setcoordinatesystem units   <1.5mm,2mm>
\setplotarea x from 0 to 16, y from -2 to 15
\put{711)} [l] at 2 12
\put {$ \scriptstyle \bullet$} [c] at  10 6
\put {$ \scriptstyle \bullet$} [c] at  10 12
\put {$ \scriptstyle \bullet$} [c] at  13 0
\put {$ \scriptstyle \bullet$} [c] at  13 6
\put {$ \scriptstyle \bullet$} [c] at  13 12
\put {$ \scriptstyle \bullet$} [c] at  16 0
\put {$ \scriptstyle \bullet$} [c] at  16 6
\setlinear \plot 16 0 16 6 13 0 10 6 13 12 16 6  /
\setlinear \plot 10 6 10 12 13 6 /
\setlinear \plot 13 0 13 12 /
\put{$2{,}520 $} [c] at 13 -2
\endpicture
\end{minipage}
\begin{minipage}{4cm}
\beginpicture
\setcoordinatesystem units   <1.5mm,2mm>
\setplotarea x from 0 to 16, y from -2 to 15
\put{712)} [l] at 2 12
\put {$ \scriptstyle \bullet$} [c] at  10 6
\put {$ \scriptstyle \bullet$} [c] at  10 0
\put {$ \scriptstyle \bullet$} [c] at  13 12
\put {$ \scriptstyle \bullet$} [c] at  13 6
\put {$ \scriptstyle \bullet$} [c] at  13 0
\put {$ \scriptstyle \bullet$} [c] at  16 12
\put {$ \scriptstyle \bullet$} [c] at  16 6
\setlinear \plot 16 12 16 6 13 12 10 6 13 0 16 6  /
\setlinear \plot 10 6 10 0 13 6 /
\setlinear \plot 13 0 13 12 /
\put{$2{,}520  $} [c] at 13 -2
\endpicture
\end{minipage}
\begin{minipage}{4cm}
\beginpicture
\setcoordinatesystem units   <1.5mm,2mm>
\setplotarea x from 0 to 16, y from -2 to 15
\put{713)} [l] at 2 12
\put {$ \scriptstyle \bullet$} [c] at  10 4
\put {$ \scriptstyle \bullet$} [c] at   11 0
\put {$ \scriptstyle \bullet$} [c] at  11 8
\put {$ \scriptstyle \bullet$} [c] at  11 12
\put {$ \scriptstyle \bullet$} [c] at  12 4
\put {$ \scriptstyle \bullet$} [c] at  16 0
\put {$ \scriptstyle \bullet$} [c] at  16 12
\setlinear \plot  16 0 16 12 11 0 12 4 11 8 11 12 16 0   /
\setlinear \plot  11 8 10 4  11 0    /
\put{$2{,}520  $} [c] at 13 -2
\endpicture
\end{minipage}
\begin{minipage}{4cm}
\beginpicture
\setcoordinatesystem units   <1.5mm,2mm>
\setplotarea x from 0 to 16, y from -2 to 15
\put{714)} [l] at 2 12
\put {$ \scriptstyle \bullet$} [c] at  10 8
\put {$ \scriptstyle \bullet$} [c] at   11 12
\put {$ \scriptstyle \bullet$} [c] at  11 4
\put {$ \scriptstyle \bullet$} [c] at  11 0
\put {$ \scriptstyle \bullet$} [c] at  12 8
\put {$ \scriptstyle \bullet$} [c] at  16 0
\put {$ \scriptstyle \bullet$} [c] at  16 12
\setlinear \plot  16 12 16 0 11 12 12 8 11 4 11 0 16 12   /
\setlinear \plot  11 4 10 8  11 12    /
\put{$2{,}520  $} [c] at 13 -2
\endpicture
\end{minipage}
$$

$$
\begin{minipage}{4cm}
\beginpicture
\setcoordinatesystem units   <1.5mm,2mm>
\setplotarea x from 0 to 16, y from -2 to 15
\put{715)} [l] at 2 12
\put {$ \scriptstyle \bullet$} [c] at  10 0
\put {$ \scriptstyle \bullet$} [c] at  10 6
\put {$ \scriptstyle \bullet$} [c] at  11 12
\put {$ \scriptstyle \bullet$} [c] at  12 0
\put {$ \scriptstyle \bullet$} [c] at  12 6
\put {$ \scriptstyle \bullet$} [c] at  13 6
\put {$ \scriptstyle \bullet$} [c] at  16 12
\setlinear \plot  10 0 10 6 11 12  12 6 12  0 16 12  10 0 12 6     /
\setlinear \plot  10 6 12 0 /
\put{$2{,}520  $} [c] at 13 -2
\endpicture
\end{minipage}
\begin{minipage}{4cm}
\beginpicture
\setcoordinatesystem units   <1.5mm,2mm>
\setplotarea x from 0 to 16, y from -2 to 15
\put{716)} [l] at 2 12
\put {$ \scriptstyle \bullet$} [c] at  10 12
\put {$ \scriptstyle \bullet$} [c] at  10 6
\put {$ \scriptstyle \bullet$} [c] at  11 0
\put {$ \scriptstyle \bullet$} [c] at  12 12
\put {$ \scriptstyle \bullet$} [c] at  12 6
\put {$ \scriptstyle \bullet$} [c] at  13 6
\put {$ \scriptstyle \bullet$} [c] at  16 0
\setlinear \plot  10 12 10 6 11 0  12 6 12  12 16 0  10 12 12 6     /
\setlinear \plot  10 6 12 12 /
\put{$2{,}520  $} [c] at 13 -2
\endpicture
\end{minipage}
\begin{minipage}{4cm}
\beginpicture
\setcoordinatesystem units   <1.5mm,2mm>
\setplotarea x from 0 to 16, y from -2 to 15
\put{717)} [l] at 2 12
\put {$ \scriptstyle \bullet$} [c] at  10 4
\put {$ \scriptstyle \bullet$} [c] at  10 12
\put {$ \scriptstyle \bullet$} [c] at  13 0
\put {$ \scriptstyle \bullet$} [c] at  13 8
\put {$ \scriptstyle \bullet$} [c] at  16 4
\put {$ \scriptstyle \bullet$} [c] at  16 12
\put {$ \scriptstyle \bullet$} [c] at  15 0
\setlinear \plot 13 0 10  4 10 12 13 8  13 0 16 4 16 12  13 8  /
\setlinear \plot 13 8 15 0 /
\put{$2{,}520   $} [c] at 13 -2
\endpicture
\end{minipage}
\begin{minipage}{4cm}
\beginpicture
\setcoordinatesystem units   <1.5mm,2mm>
\setplotarea x from 0 to 16, y from -2 to 15
\put{718)} [l] at 2 12
\put {$ \scriptstyle \bullet$} [c] at  10 8
\put {$ \scriptstyle \bullet$} [c] at  10 0
\put {$ \scriptstyle \bullet$} [c] at  13 12
\put {$ \scriptstyle \bullet$} [c] at  13 4
\put {$ \scriptstyle \bullet$} [c] at  16 8
\put {$ \scriptstyle \bullet$} [c] at  16 0
\put {$ \scriptstyle \bullet$} [c] at  15 12
\setlinear \plot 13 12 10  8 10 0 13 4  13 12 16 8 16 0  13 4  /
\setlinear \plot 13 4 15 12 /
\put{$2{,}520   $} [c] at 13 -2
\endpicture
\end{minipage}
\begin{minipage}{4cm}
\beginpicture
\setcoordinatesystem units   <1.5mm,2mm>
\setplotarea x from 0 to 16, y from -2 to 15
\put{719)} [l] at 2 12
\put {$ \scriptstyle \bullet$} [c] at  10 0
\put {$ \scriptstyle \bullet$} [c] at  10 6
\put {$ \scriptstyle \bullet$} [c] at  10 12
\put {$ \scriptstyle \bullet$} [c] at  16 12
\put {$ \scriptstyle \bullet$} [c] at  16 6
\put {$ \scriptstyle \bullet$} [c] at  16 0
\put {$ \scriptstyle \bullet$} [c] at  13 12
\setlinear \plot  10 12 10  0 13 12 16 6 16 0 10 6  /
\setlinear \plot  10 0 16 12 16 6  /
\put{$2{,}520$} [c] at 13 -2
\endpicture
\end{minipage}
\begin{minipage}{4cm}
\beginpicture
\setcoordinatesystem units   <1.5mm,2mm>
\setplotarea x from 0 to 16, y from -2 to 15
\put{720)} [l] at 2 12
\put {$ \scriptstyle \bullet$} [c] at  10 0
\put {$ \scriptstyle \bullet$} [c] at  10 6
\put {$ \scriptstyle \bullet$} [c] at  10 12
\put {$ \scriptstyle \bullet$} [c] at  16 12
\put {$ \scriptstyle \bullet$} [c] at  16 6
\put {$ \scriptstyle \bullet$} [c] at  16 0
\put {$ \scriptstyle \bullet$} [c] at  13 0
\setlinear \plot  10 0 10  12 13 0 16 6 16 12 10 6  /
\setlinear \plot  10 12 16 0 16 6  /
\put{$2{,}520$} [c] at 13 -2
\endpicture
\end{minipage}
$$
$$
\begin{minipage}{4cm}
\beginpicture
\setcoordinatesystem units   <1.5mm,2mm>
\setplotarea x from 0 to 16, y from -2 to 15
\put{721)} [l] at 2 12
\put {$ \scriptstyle \bullet$} [c] at  10 0
\put {$ \scriptstyle \bullet$} [c] at  10 4
\put {$ \scriptstyle \bullet$} [c] at  10 8
\put {$ \scriptstyle \bullet$} [c] at  10 12
\put {$ \scriptstyle \bullet$} [c] at  13 12
\put {$ \scriptstyle \bullet$} [c] at  16 0
\put {$ \scriptstyle \bullet$} [c] at  16 12
\setlinear \plot  10 0 10 12  /
\setlinear \plot  13 12 10 4 16 12 16 0  13 12 /
\put{$2{,}520$} [c] at 13 -2
\endpicture
\end{minipage}
\begin{minipage}{4cm}
\beginpicture
\setcoordinatesystem units   <1.5mm,2mm>
\setplotarea x from 0 to 16, y from -2 to 15
\put{722)} [l] at 2 12
\put {$ \scriptstyle \bullet$} [c] at  10 0
\put {$ \scriptstyle \bullet$} [c] at  10 4
\put {$ \scriptstyle \bullet$} [c] at  10 8
\put {$ \scriptstyle \bullet$} [c] at  10 12
\put {$ \scriptstyle \bullet$} [c] at  13 0
\put {$ \scriptstyle \bullet$} [c] at  16 0
\put {$ \scriptstyle \bullet$} [c] at  16 12
\setlinear \plot  10 0 10 12  /
\setlinear \plot  13 0 10 8 16 0 16 12  13 0 /
\put{$2{,}520$} [c] at 13 -2
\endpicture
\end{minipage}
\begin{minipage}{4cm}
\beginpicture
\setcoordinatesystem units   <1.5mm,2mm>
\setplotarea x from 0 to 16, y from -2 to 15
\put{723)} [l] at 2 12
\put {$ \scriptstyle \bullet$} [c] at  10 12
\put {$ \scriptstyle \bullet$} [c] at  10 6
\put {$ \scriptstyle \bullet$} [c] at  10 0
\put {$ \scriptstyle \bullet$} [c] at  13 12
\put {$ \scriptstyle \bullet$} [c] at  16 0
\put {$ \scriptstyle \bullet$} [c] at  16 6
\put {$ \scriptstyle \bullet$} [c] at  16 12
\setlinear \plot 13 12 10 6 10 0 16 12  16 0 10 6 10 12  /
\put{$2{,}520$} [c] at 13 -2
\endpicture
\end{minipage}
\begin{minipage}{4cm}
\beginpicture
\setcoordinatesystem units   <1.5mm,2mm>
\setplotarea x from 0 to 16, y from -2 to 15
\put{724)} [l] at 2 12
\put {$ \scriptstyle \bullet$} [c] at  10 12
\put {$ \scriptstyle \bullet$} [c] at  10 6
\put {$ \scriptstyle \bullet$} [c] at  10 0
\put {$ \scriptstyle \bullet$} [c] at  13 0
\put {$ \scriptstyle \bullet$} [c] at  16 0
\put {$ \scriptstyle \bullet$} [c] at  16 6
\put {$ \scriptstyle \bullet$} [c] at  16 12
\setlinear \plot 13 0 10 6 10 12 16 0  16 12 10 6 10 0  /
\put{$2{,}520$} [c] at 13 -2
\endpicture
\end{minipage}
\begin{minipage}{4cm}
\beginpicture
\setcoordinatesystem units   <1.5mm,2mm>
\setplotarea x from 0 to 16, y from -2 to 15
\put{725)} [l] at 2 12
\put {$ \scriptstyle \bullet$} [c] at  10  12
\put {$ \scriptstyle \bullet$} [c] at  11 0
\put {$ \scriptstyle \bullet$} [c] at  11 8
\put {$ \scriptstyle \bullet$} [c] at  12 12
\put {$ \scriptstyle \bullet$} [c] at  13 8
\put {$ \scriptstyle \bullet$} [c] at  13 0
\put {$ \scriptstyle \bullet$} [c] at  16 12
\setlinear \plot 10  12 11 8 11 0 13 8 13 0 11 8 12 12 13 8   /
\setlinear \plot  11 0 16 12  13 0  /
\put{$2{,}520$} [c] at 13 -2
\endpicture
\end{minipage}
\begin{minipage}{4cm}
\beginpicture
\setcoordinatesystem units   <1.5mm,2mm>
\setplotarea x from 0 to 16, y from -2 to 15
\put{726)} [l] at 2 12
\put {$ \scriptstyle \bullet$} [c] at  10  0
\put {$ \scriptstyle \bullet$} [c] at  11 12
\put {$ \scriptstyle \bullet$} [c] at  11 4
\put {$ \scriptstyle \bullet$} [c] at  12 0
\put {$ \scriptstyle \bullet$} [c] at  13 12
\put {$ \scriptstyle \bullet$} [c] at  13 4
\put {$ \scriptstyle \bullet$} [c] at  16 0
\setlinear \plot 10  0 11 4 11 12 13 4 13 12 11 4 12 0 13 4   /
\setlinear \plot  11 12 16 0  13 12  /
\put{$2{,}520$} [c] at 13 -2
\endpicture
\end{minipage}
$$

$$
\begin{minipage}{4cm}
\beginpicture
\setcoordinatesystem units   <1.5mm,2mm>
\setplotarea x from 0 to 16, y from -2 to 15
\put{727)} [l] at 2 12
\put {$ \scriptstyle \bullet$} [c] at  10 0
\put {$ \scriptstyle \bullet$} [c] at  10 6
\put {$ \scriptstyle \bullet$} [c] at  10 12
\put {$ \scriptstyle \bullet$} [c] at  13 12
\put {$ \scriptstyle \bullet$} [c] at  16 0
\put {$ \scriptstyle \bullet$} [c] at  16 6
\put {$ \scriptstyle \bullet$} [c] at  16 12
\setlinear \plot  10 0 10 12 16  0 16 12 10 6 13  12 16 0  /
\put{$2{,}520$} [c] at 13 -2
\endpicture
\end{minipage}
\begin{minipage}{4cm}
\beginpicture
\setcoordinatesystem units   <1.5mm,2mm>
\setplotarea x from 0 to 16, y from -2 to 15
\put{728)} [l] at 2 12
\put {$ \scriptstyle \bullet$} [c] at  10 0
\put {$ \scriptstyle \bullet$} [c] at  10 6
\put {$ \scriptstyle \bullet$} [c] at  10 12
\put {$ \scriptstyle \bullet$} [c] at  13 0
\put {$ \scriptstyle \bullet$} [c] at  16 0
\put {$ \scriptstyle \bullet$} [c] at  16 6
\put {$ \scriptstyle \bullet$} [c] at  16 12
\setlinear \plot  10 12 10 0 16  12 16 0 10 6 13  0 16 12  /
\put{$2{,}520$} [c] at 13 -2
\endpicture
\end{minipage}
\begin{minipage}{4cm}
\beginpicture
\setcoordinatesystem units   <1.5mm,2mm>
\setplotarea x from 0 to 16, y from -2 to 15
\put{729)} [l] at 2 12
\put {$ \scriptstyle \bullet$} [c] at  10 0
\put {$ \scriptstyle \bullet$} [c] at  10 12
\put {$ \scriptstyle \bullet$} [c] at  14 4
\put {$ \scriptstyle \bullet$} [c] at  14 12
\put {$ \scriptstyle \bullet$} [c] at  15 0
\put {$ \scriptstyle \bullet$} [c] at  16 4
\put {$ \scriptstyle \bullet$} [c] at  16 12
\setlinear \plot 14 4 10 12 10 0 16 12 16 4 15  0 14 4 14 12 16 4  /
\setlinear \plot  14 12 10 0 16 12  /
\put{$2{,}520$} [c] at 13 -2
\endpicture
\end{minipage}
\begin{minipage}{4cm}
\beginpicture
\setcoordinatesystem units   <1.5mm,2mm>
\setplotarea x from 0 to 16, y from -2 to 15
\put{730)} [l] at 2 12
\put {$ \scriptstyle \bullet$} [c] at  10 0
\put {$ \scriptstyle \bullet$} [c] at  10 12
\put {$ \scriptstyle \bullet$} [c] at  14 0
\put {$ \scriptstyle \bullet$} [c] at  14 8
\put {$ \scriptstyle \bullet$} [c] at  15 12
\put {$ \scriptstyle \bullet$} [c] at  16 0
\put {$ \scriptstyle \bullet$} [c] at  16 8
\setlinear \plot 14 8 10 0 10 12 16 0 16 8 15  12 14 8 14 0 16 8  /
\setlinear \plot  14 0 10 12 16 0  /
\put{$2{,}520$} [c] at 13 -2
\endpicture
\end{minipage}
\begin{minipage}{4cm}
\beginpicture
\setcoordinatesystem units   <1.5mm,2mm>
\setplotarea x from 0 to 16, y from -2 to 15
\put{731)} [l] at 2 12
\put {$ \scriptstyle \bullet$} [c] at  10 0
\put {$ \scriptstyle \bullet$} [c] at  10 4
\put {$ \scriptstyle \bullet$} [c] at  10 8
\put {$ \scriptstyle \bullet$} [c] at  10 12
\put {$ \scriptstyle \bullet$} [c] at  13 12
\put {$ \scriptstyle \bullet$} [c] at  16 0
\put {$ \scriptstyle \bullet$} [c] at  16 12
\setlinear \plot  10 0 10 12  /
\setlinear \plot  16 0 16 12 10 4  /
\setlinear \plot  10 8 13 12  /
\put{$2{,}520$} [c] at 13 -2
\endpicture
\end{minipage}
\begin{minipage}{4cm}
\beginpicture
\setcoordinatesystem units   <1.5mm,2mm>
\setplotarea x from 0 to 16, y from -2 to 15
\put{732)} [l] at 2 12
\put {$ \scriptstyle \bullet$} [c] at  10 0
\put {$ \scriptstyle \bullet$} [c] at  10 4
\put {$ \scriptstyle \bullet$} [c] at  10 8
\put {$ \scriptstyle \bullet$} [c] at  10 12
\put {$ \scriptstyle \bullet$} [c] at  13 0
\put {$ \scriptstyle \bullet$} [c] at  16 0
\put {$ \scriptstyle \bullet$} [c] at  16 12
\setlinear \plot  10 0 10 12  /
\setlinear \plot  16 12 16 0 10 8  /
\setlinear \plot  10 4 13 0  /
\put{$2{,}520$} [c] at 13 -2
\endpicture
\end{minipage}
$$
$$
\begin{minipage}{4cm}
\beginpicture
\setcoordinatesystem units   <1.5mm,2mm>
\setplotarea x from 0 to 16, y from -2 to 15
\put{733)} [l] at 2 12
\put {$ \scriptstyle \bullet$} [c] at  10 0
\put {$ \scriptstyle \bullet$} [c] at  10 12
\put {$ \scriptstyle \bullet$} [c] at  12 8
\put {$ \scriptstyle \bullet$} [c] at  12 12
\put {$ \scriptstyle \bullet$} [c] at  16 12
\put {$ \scriptstyle \bullet$} [c] at  16 8
\put {$ \scriptstyle \bullet$} [c] at  14 0
\setlinear \plot  16 8 14 0 12 8 12 12 16 8 16 12 12 8 10 0 10 12  /
\put{$2{,}520$} [c] at 13 -2
\endpicture
\end{minipage}
\begin{minipage}{4cm}
\beginpicture
\setcoordinatesystem units   <1.5mm,2mm>
\setplotarea x from 0 to 16, y from -2 to 15
\put{734)} [l] at 2 12
\put {$ \scriptstyle \bullet$} [c] at  10 0
\put {$ \scriptstyle \bullet$} [c] at  10 12
\put {$ \scriptstyle \bullet$} [c] at  12 4
\put {$ \scriptstyle \bullet$} [c] at  12 0
\put {$ \scriptstyle \bullet$} [c] at  16 0
\put {$ \scriptstyle \bullet$} [c] at  16 4
\put {$ \scriptstyle \bullet$} [c] at  14 12
\setlinear \plot  16 4 14 12 12 4 12 0 16 4 16 0 12 4 10 12 10 0  /
\put{$2{,}520$} [c] at 13 -2
\endpicture
\end{minipage}
\begin{minipage}{4cm}
\beginpicture
\setcoordinatesystem units   <1.5mm,2mm>
\setplotarea x from 0 to 16, y from -2 to 15
\put{735)} [l] at 2 12
\put {$ \scriptstyle \bullet$} [c] at  10 0
\put {$ \scriptstyle \bullet$} [c] at  10 12
\put {$ \scriptstyle \bullet$} [c] at  13 0
\put {$ \scriptstyle \bullet$} [c] at  13 12
\put {$ \scriptstyle \bullet$} [c] at  16 12
\put {$ \scriptstyle \bullet$} [c] at  16 6
\put {$ \scriptstyle \bullet$} [c] at  16 0
\setlinear \plot  10 0 10 12 16 0  16 12  /
\setlinear \plot  10 12 13 0 16 6 13 12 10 0 /
\put{$2{,}520$} [c] at 13 -2
\endpicture
\end{minipage}
\begin{minipage}{4cm}
\beginpicture
\setcoordinatesystem units   <1.5mm,2mm>
\setplotarea x from 0 to 16, y from -2 to 15
\put{736)} [l] at 2 12
\put {$ \scriptstyle \bullet$} [c] at  10 0
\put {$ \scriptstyle \bullet$} [c] at  10 12
\put {$ \scriptstyle \bullet$} [c] at  13 0
\put {$ \scriptstyle \bullet$} [c] at  13 12
\put {$ \scriptstyle \bullet$} [c] at  16 12
\put {$ \scriptstyle \bullet$} [c] at  16 6
\put {$ \scriptstyle \bullet$} [c] at  16 0
\setlinear \plot  10 12 10 0 16 12  16 0  /
\setlinear \plot  10 0 13 12 16 6 13 0 10 12 /
\put{$2{,}520$} [c] at 13 -2
\endpicture
\end{minipage}
\begin{minipage}{4cm}
\beginpicture
\setcoordinatesystem units   <1.5mm,2mm>
\setplotarea x from 0 to 16, y from -2 to 15
\put{737)} [l] at 2 12
\put {$ \scriptstyle \bullet$} [c] at 13 0
\put {$ \scriptstyle \bullet$} [c] at 10 6
\put {$ \scriptstyle \bullet$} [c] at 10 12
\put {$ \scriptstyle \bullet$} [c] at 13 6
\put {$ \scriptstyle \bullet$} [c] at 13 12
\put {$ \scriptstyle \bullet$} [c] at 16 6
\put {$ \scriptstyle \bullet$} [c] at 16 12
\setlinear \plot 16 12  16 6 13 0  10 6 10 12 16 6 13 12  10 6  /
\setlinear \plot 13 0 13 12   /
\setlinear \plot 10 12 13 6  /
\put{$1{,}260$} [c] at 13 -2
\endpicture
\end{minipage}
\begin{minipage}{4cm}
\beginpicture
\setcoordinatesystem units   <1.5mm,2mm>
\setplotarea x from 0 to 16, y from -2 to 15
\put{738)} [l] at 2 12
\put {$ \scriptstyle \bullet$} [c] at 13 12
\put {$ \scriptstyle \bullet$} [c] at 10 6
\put {$ \scriptstyle \bullet$} [c] at 10 0
\put {$ \scriptstyle \bullet$} [c] at 13 0
\put {$ \scriptstyle \bullet$} [c] at 13 6
\put {$ \scriptstyle \bullet$} [c] at 16 6
\put {$ \scriptstyle \bullet$} [c] at 16 0
\setlinear \plot 16 0  16 6 13 12  10 6 10 0 16 6 13 0  10 6  /
\setlinear \plot 13 0 13 12   /
\setlinear \plot 10 0 13 6  /
\put{$1{,}260$} [c] at 13 -2
\endpicture
\end{minipage}
$$

$$
\begin{minipage}{4cm}
\beginpicture
\setcoordinatesystem units   <1.5mm,2mm>
\setplotarea x from 0 to 16, y from -2 to 15
\put{739)} [l] at 2 12
\put {$ \scriptstyle \bullet$} [c] at 13 0
\put {$ \scriptstyle \bullet$} [c] at 10 6
\put {$ \scriptstyle \bullet$} [c] at 10 12
\put {$ \scriptstyle \bullet$} [c] at 13 6
\put {$ \scriptstyle \bullet$} [c] at 13 12
\put {$ \scriptstyle \bullet$} [c] at 16 6
\put {$ \scriptstyle \bullet$} [c] at 16 12
\setlinear \plot  13 6 13 0 16 6 16 12 13 6 13 12 10 6 10 12 13 6  /
\setlinear \plot 16 12 10 6 13 0   /
\put{$1{,}260$} [c] at 13 -2
\endpicture
\end{minipage}
\begin{minipage}{4cm}
\beginpicture
\setcoordinatesystem units   <1.5mm,2mm>
\setplotarea x from 0 to 16, y from -2 to 15
\put{740)} [l] at 2 12
\put {$ \scriptstyle \bullet$} [c] at 13 12
\put {$ \scriptstyle \bullet$} [c] at 10 6
\put {$ \scriptstyle \bullet$} [c] at 10 0
\put {$ \scriptstyle \bullet$} [c] at 13 0
\put {$ \scriptstyle \bullet$} [c] at 13 6
\put {$ \scriptstyle \bullet$} [c] at 16 6
\put {$ \scriptstyle \bullet$} [c] at 16 0
\setlinear \plot  13 6 13 12 16 6 16 0 13 6 13 0 10 6 10 0 13 6  /
\setlinear \plot 16 0 10 6 13 12   /
\put{$1{,}260$} [c] at 13 -2
\endpicture
\end{minipage}
\begin{minipage}{4cm}
\beginpicture
\setcoordinatesystem units   <1.5mm,2mm>
\setplotarea x from 0 to 16, y from -2 to 15
\put{741)} [l] at 2 12
\put {$ \scriptstyle \bullet$} [c] at  10 0
\put {$ \scriptstyle \bullet$} [c] at  10 4
\put {$ \scriptstyle \bullet$} [c] at  10 8
\put {$ \scriptstyle \bullet$} [c] at  10 12
\put {$ \scriptstyle \bullet$} [c] at  13 12
\put {$ \scriptstyle \bullet$} [c] at  16 0
\put {$ \scriptstyle \bullet$} [c] at  16 12
\setlinear \plot 10 12 10 0  16 12 16 0 13 12 10 0 /
\setlinear \plot  10 4 16 0    /
\put{$1{,}260$} [c] at 13 -2
\endpicture
\end{minipage}
\begin{minipage}{4cm}
\beginpicture
\setcoordinatesystem units   <1.5mm,2mm>
\setplotarea x from 0 to 16, y from -2 to 15
\put{742)} [l] at 2 12
\put {$ \scriptstyle \bullet$} [c] at  10 0
\put {$ \scriptstyle \bullet$} [c] at  10 4
\put {$ \scriptstyle \bullet$} [c] at  10 8
\put {$ \scriptstyle \bullet$} [c] at  10 12
\put {$ \scriptstyle \bullet$} [c] at  13 0
\put {$ \scriptstyle \bullet$} [c] at  16 0
\put {$ \scriptstyle \bullet$} [c] at  16 12
\setlinear \plot 10 0 10 12  16 0 16 12 13 0 10 12 /
\setlinear \plot  10 8 16 12    /
\put{$1{,}260$} [c] at 13 -2
\endpicture
\end{minipage}
\begin{minipage}{4cm}
\beginpicture
\setcoordinatesystem units   <1.5mm,2mm>
\setplotarea x from 0 to 16, y from -2 to 15
\put{743)} [l] at 2 12
\put {$ \scriptstyle \bullet$} [c] at  10 12
\put {$ \scriptstyle \bullet$} [c] at  13 0
\put {$ \scriptstyle \bullet$} [c] at  13 8
\put {$ \scriptstyle \bullet$} [c] at  13 12
\put {$ \scriptstyle \bullet$} [c] at  16 0
\put {$ \scriptstyle \bullet$} [c] at  16 8
\put {$ \scriptstyle \bullet$} [c] at  16 12
\setlinear \plot  16 8 13 12 13 8 13 0  10 12 16 0 16 12 13 8 13 0  /
\put{$1{,}260$} [c] at 13 -2
\endpicture
\end{minipage}
\begin{minipage}{4cm}
\beginpicture
\setcoordinatesystem units   <1.5mm,2mm>
\setplotarea x from 0 to 16, y from -2 to 15
\put{744)} [l] at 2 12
\put {$ \scriptstyle \bullet$} [c] at  10 0
\put {$ \scriptstyle \bullet$} [c] at  13 0
\put {$ \scriptstyle \bullet$} [c] at  13 4
\put {$ \scriptstyle \bullet$} [c] at  13 12
\put {$ \scriptstyle \bullet$} [c] at  16 0
\put {$ \scriptstyle \bullet$} [c] at  16 4
\put {$ \scriptstyle \bullet$} [c] at  16 12
\setlinear \plot  16 4 13 0 13 4 13 12  10 0 16 12 16 0 13 4 13 12  /
\put{$1{,}260$} [c] at 13 -2
\endpicture
\end{minipage}
$$
$$
\begin{minipage}{4cm}
\beginpicture
\setcoordinatesystem units   <1.5mm,2mm>
\setplotarea x from 0 to 16, y from -2 to 15
\put{745)} [l] at 2 12
\put {$ \scriptstyle \bullet$} [c] at 13 0
\put {$ \scriptstyle \bullet$} [c] at 10 6
\put {$ \scriptstyle \bullet$} [c] at 10 12
\put {$ \scriptstyle \bullet$} [c] at 13 6
\put {$ \scriptstyle \bullet$} [c] at 13 12
\put {$ \scriptstyle \bullet$} [c] at 16 6
\put {$ \scriptstyle \bullet$} [c] at 16 12
\setlinear \plot 16 12 16  6 13 0 10 6 10 12 13 6 16 12 /
\setlinear \plot 10 6  13 12 16 6  /
\setlinear \plot 13 0  13  6  /
\put{$840$} [c] at 13 -2
\endpicture
\end{minipage}
\begin{minipage}{4cm}
\beginpicture
\setcoordinatesystem units   <1.5mm,2mm>
\setplotarea x from 0 to 16, y from -2 to 15
\put{746)} [l] at 2 12
\put {$ \scriptstyle \bullet$} [c] at 13 12
\put {$ \scriptstyle \bullet$} [c] at 10 6
\put {$ \scriptstyle \bullet$} [c] at 10 0
\put {$ \scriptstyle \bullet$} [c] at 13 0
\put {$ \scriptstyle \bullet$} [c] at 13 6
\put {$ \scriptstyle \bullet$} [c] at 16 6
\put {$ \scriptstyle \bullet$} [c] at 16 0
\setlinear \plot 16 0 16  6 13 12 10 6 10 0 13 6 16 0 /
\setlinear \plot 10 6  13 0 16 6  /
\setlinear \plot 13 12  13  6  /
\put{$840$} [c] at 13 -2
\endpicture
\end{minipage}
\begin{minipage}{4cm}
\beginpicture
\setcoordinatesystem units   <1.5mm,2mm>
\setplotarea x from 0 to 16, y from -2 to 15
\put{747)} [l] at 2 12
\put {$ \scriptstyle \bullet$} [c] at  10 0
\put {$ \scriptstyle \bullet$} [c] at  10 12
\put {$ \scriptstyle \bullet$} [c] at  13 0
\put {$ \scriptstyle \bullet$} [c] at  13 12
\put {$ \scriptstyle \bullet$} [c] at  16 12
\put {$ \scriptstyle \bullet$} [c] at  16 6
\put {$ \scriptstyle \bullet$} [c] at  16 0
\setlinear \plot  16 0  16 12 10 0 10 12 13 0 16 12 /
\setlinear \plot  10 12  16 6 13 12 10 0 /
\setlinear \plot  13 12  13 0 /
\put{$420$} [c] at 13 -2

\endpicture
\end{minipage}
\begin{minipage}{4cm}
\beginpicture
\setcoordinatesystem units   <1.5mm,2mm>
\setplotarea x from 0 to 16, y from -2 to 15
\put{748)} [l] at 2 12
\put {$ \scriptstyle \bullet$} [c] at  10 0
\put {$ \scriptstyle \bullet$} [c] at  10 12
\put {$ \scriptstyle \bullet$} [c] at  13 0
\put {$ \scriptstyle \bullet$} [c] at  13 12
\put {$ \scriptstyle \bullet$} [c] at  16 12
\put {$ \scriptstyle \bullet$} [c] at  16 6
\put {$ \scriptstyle \bullet$} [c] at  16 0
\setlinear \plot  16 12  16 0 10 12 10 0 13 12 16 0 /
\setlinear \plot  10 0  16 6 13 0 10 12 /
\setlinear \plot  13 12  13 0 /
\put{$420$} [c] at 13 -2
\endpicture
\end{minipage}
\begin{minipage}{4cm}
\beginpicture
\setcoordinatesystem units   <1.5mm,2mm>
\setplotarea x from 0 to 16, y from -2 to 15
\put{749)} [l] at 2 12
\put {$ \scriptstyle \bullet$} [c] at 10 8
\put {$ \scriptstyle \bullet$} [c] at 12 8
\put {$ \scriptstyle \bullet$} [c] at 14 8
\put {$ \scriptstyle \bullet$} [c] at 16 8
\put {$ \scriptstyle \bullet$} [c] at 13 0
\put {$ \scriptstyle \bullet$} [c] at 13 12
\put {$ \scriptstyle \bullet$} [c] at 13 4
\setlinear \plot 13 0 13 4 10 8 13 12 16 8  13  4   /
\setlinear \plot  13 4 12 8 13 12 14 8 13 4   /
\put{$210$} [c] at 13 -2
\endpicture
\end{minipage}
\begin{minipage}{4cm}
\beginpicture
\setcoordinatesystem units   <1.5mm,2mm>
\setplotarea x from 0 to 16, y from -2 to 15
\put{750)} [l] at 2 12
\put {$ \scriptstyle \bullet$} [c] at 10 4
\put {$ \scriptstyle \bullet$} [c] at 12 4
\put {$ \scriptstyle \bullet$} [c] at 14 4
\put {$ \scriptstyle \bullet$} [c] at 16 4
\put {$ \scriptstyle \bullet$} [c] at 13 0
\put {$ \scriptstyle \bullet$} [c] at 13 12
\put {$ \scriptstyle \bullet$} [c] at 13 8
\setlinear \plot 13 12 13 8 10 4 13 0 16 4  13  8   /
\setlinear \plot  13 8 12 4 13 0 14 4 13 8   /
\put{$210$} [c] at 13 -2
\endpicture
\end{minipage}
$$

$$
\begin{minipage}{4cm}
\beginpicture
\setcoordinatesystem units   <1.5mm,2mm>
\setplotarea x from 0 to 16, y from -2 to 15
\put{751)} [l] at 2 12
\put {$ \scriptstyle \bullet$} [c] at 13 0
\put {$ \scriptstyle \bullet$} [c] at 13 4
\put {$ \scriptstyle \bullet$} [c] at 13 8
\put {$ \scriptstyle \bullet$} [c] at 10 12
\put {$ \scriptstyle \bullet$} [c] at 12 12
\put {$ \scriptstyle \bullet$} [c] at 14 12
\put {$ \scriptstyle \bullet$} [c] at 16 12
\setlinear \plot 13 0 13 8     /
\setlinear \plot 10 12 13 8 16 12    /
\setlinear \plot 12 12 13 8 14 12    /
\put{$210$} [c] at 13 -2
\endpicture
\end{minipage}
\begin{minipage}{4cm}
\beginpicture
\setcoordinatesystem units   <1.5mm,2mm>
\setplotarea x from 0 to 16, y from -2 to 15
\put{752)} [l] at 2 12
\put {$ \scriptstyle \bullet$} [c] at 13 12
\put {$ \scriptstyle \bullet$} [c] at 13 8
\put {$ \scriptstyle \bullet$} [c] at 13 4
\put {$ \scriptstyle \bullet$} [c] at 10 0
\put {$ \scriptstyle \bullet$} [c] at 12 0
\put {$ \scriptstyle \bullet$} [c] at 14 0
\put {$ \scriptstyle \bullet$} [c] at 16 0
\setlinear \plot 13 12 13 4     /
\setlinear \plot 10 0 13 4 16 0    /
\setlinear \plot 12 0 13 4 14 0    /
\put{$210$} [c] at 13 -2
\endpicture
\end{minipage}
\begin{minipage}{4cm}
\beginpicture
\setcoordinatesystem units   <1.5mm,2mm>
\setplotarea x from 0 to 16, y from -2 to 12
\put{${\bf  20}$} [l] at 2 15
\put{753)} [l] at 2 12
\put {$ \scriptstyle \bullet$} [c] at 10 6
\put {$ \scriptstyle \bullet$} [c] at 10 12
\put {$ \scriptstyle \bullet$} [c] at 13 0
\put {$ \scriptstyle \bullet$} [c] at 13 6
\put {$ \scriptstyle \bullet$} [c] at 16 6
\put {$ \scriptstyle \bullet$} [c] at 16 12
\put {$ \scriptstyle \bullet$} [c] at 14.5 9
\setlinear \plot 13 0 10  6  10 12 13 6 13 0 16 6 16 12 13 6 /
\put{$5{,}040$} [c] at 13 -2
\endpicture
\end{minipage}
\begin{minipage}{4cm}
\beginpicture
\setcoordinatesystem units   <1.5mm,2mm>
\setplotarea x from 0 to 16, y from -2 to 15
\put{754)} [l] at 2 12
\put {$ \scriptstyle \bullet$} [c] at 10 6
\put {$ \scriptstyle \bullet$} [c] at 10 0
\put {$ \scriptstyle \bullet$} [c] at 13 12
\put {$ \scriptstyle \bullet$} [c] at 13 6
\put {$ \scriptstyle \bullet$} [c] at 16 6
\put {$ \scriptstyle \bullet$} [c] at 16 0
\put {$ \scriptstyle \bullet$} [c] at 14.5 3
\setlinear \plot 13 12 10  6  10 0 13 6 13 12 16 6 16 0 13 6 /
\put{$5{,}040$} [c] at 13 -2
\endpicture
\end{minipage}
\begin{minipage}{4cm}
\beginpicture
\setcoordinatesystem units   <1.5mm,2mm>
\setplotarea x from 0 to 16, y from -2 to 15
\put{755)} [l] at 2 12
\put {$ \scriptstyle \bullet$} [c] at 10 6
\put {$ \scriptstyle \bullet$} [c] at 10 12
\put {$ \scriptstyle \bullet$} [c] at 13 0
\put {$ \scriptstyle \bullet$} [c] at 13 12
\put {$ \scriptstyle \bullet$} [c] at 14.5 9
\put {$ \scriptstyle \bullet$} [c] at 16 6
\put {$ \scriptstyle \bullet$} [c] at 16 12
\setlinear \plot 10  12  10 6  13 0 16  6 13  12 10 6 /
\setlinear \plot 16  12  16 6  /
\put{$5{,}040$} [c] at 13 -2
\endpicture
\end{minipage}
\begin{minipage}{4cm}
\beginpicture
\setcoordinatesystem units   <1.5mm,2mm>
\setplotarea x from 0 to 16, y from -2 to 15
\put{756)} [l] at 2 12
\put {$ \scriptstyle \bullet$} [c] at 10 6
\put {$ \scriptstyle \bullet$} [c] at 10 0
\put {$ \scriptstyle \bullet$} [c] at 13 0
\put {$ \scriptstyle \bullet$} [c] at 13 12
\put {$ \scriptstyle \bullet$} [c] at 14.5 3
\put {$ \scriptstyle \bullet$} [c] at 16 6
\put {$ \scriptstyle \bullet$} [c] at 16 0
\setlinear \plot 10  0  10 6  13 12 16  6 13  0 10 6 /
\setlinear \plot 16  0  16 6  /
\put{$5{,}040$} [c] at 13 -2
\endpicture
\end{minipage}
$$

$$
\begin{minipage}{4cm}
\beginpicture
\setcoordinatesystem units   <1.5mm,2mm>
\setplotarea x from 0 to 16, y from -2 to 15
\put{757)} [l] at 2 12
\put {$ \scriptstyle \bullet$} [c] at  10  4
\put {$ \scriptstyle \bullet$} [c] at  10 8
\put {$ \scriptstyle \bullet$} [c] at  10 12
\put {$ \scriptstyle \bullet$} [c] at  11.5 0
\put {$ \scriptstyle \bullet$} [c] at  13 4
\put {$ \scriptstyle \bullet$} [c] at  13 12
\put {$ \scriptstyle \bullet$} [c] at  16 0
\setlinear \plot 11.5 0 10 4 10 12 16 0 13 12 13 4 11.5 0      /
\put{$5{,}040   $} [c] at 13 -2
\endpicture
\end{minipage}
\begin{minipage}{4cm}
\beginpicture
\setcoordinatesystem units   <1.5mm,2mm>
\setplotarea x from 0 to 16, y from -2 to 15
\put{758)} [l] at 2 12
\put {$ \scriptstyle \bullet$} [c] at  10  4
\put {$ \scriptstyle \bullet$} [c] at  10 8
\put {$ \scriptstyle \bullet$} [c] at  10 0
\put {$ \scriptstyle \bullet$} [c] at  11.5 12
\put {$ \scriptstyle \bullet$} [c] at  13 8
\put {$ \scriptstyle \bullet$} [c] at  13 0
\put {$ \scriptstyle \bullet$} [c] at  16 12
\setlinear \plot 11.5 12 10 8 10 0 16 12 13 0 13 8 11.5 12      /
\put{$5{,}040   $} [c] at 13 -2
\endpicture
\end{minipage}
\begin{minipage}{4cm}
\beginpicture
\setcoordinatesystem units   <1.5mm,2mm>
\setplotarea x from 0 to 16, y from -2 to 15
\put{759)} [l] at 2 12
\put {$ \scriptstyle \bullet$} [c] at  10 6
\put {$ \scriptstyle \bullet$} [c] at  11  0
\put {$ \scriptstyle \bullet$} [c] at  11 12
\put {$ \scriptstyle \bullet$} [c] at  12 6
\put {$ \scriptstyle \bullet$} [c] at  16 0
\put {$ \scriptstyle \bullet$} [c] at  16  6
\put {$ \scriptstyle \bullet$} [c] at  16  12
\setlinear \plot 11 12 16 0 16 12 12 6 11 0 10 6 11 12  12 6  /
\put{$5{,}040$} [c] at 13 -2
\endpicture
\end{minipage}
\begin{minipage}{4cm}
\beginpicture
\setcoordinatesystem units   <1.5mm,2mm>
\setplotarea x from 0 to 16, y from -2 to 15
\put{760)} [l] at 2 12
\put {$ \scriptstyle \bullet$} [c] at  10 6
\put {$ \scriptstyle \bullet$} [c] at  11  0
\put {$ \scriptstyle \bullet$} [c] at  11 12
\put {$ \scriptstyle \bullet$} [c] at  12 6
\put {$ \scriptstyle \bullet$} [c] at  16 0
\put {$ \scriptstyle \bullet$} [c] at  16  6
\put {$ \scriptstyle \bullet$} [c] at  16  12
\setlinear \plot 11 0 16 12 16 0 12 6 11 12 10 6 11 0  12 6  /
\put{$5{,}040$} [c] at 13 -2
\endpicture
\end{minipage}
\begin{minipage}{4cm}
\beginpicture
\setcoordinatesystem units   <1.5mm,2mm>
\setplotarea x from 0 to 16, y from -2 to 15
\put{761)} [l] at 2 12
\put {$ \scriptstyle \bullet$} [c] at  10 0
\put {$ \scriptstyle \bullet$} [c] at  13 6
\put {$ \scriptstyle \bullet$} [c] at  13 12
\put {$ \scriptstyle \bullet$} [c] at  14.5 3
\put {$ \scriptstyle \bullet$} [c] at  16 0
\put {$ \scriptstyle \bullet$} [c] at  16 6
\put {$ \scriptstyle \bullet$} [c] at  16 12
\setlinear \plot 16 12 16 0 13 6 13  12  /
\setlinear \plot 10  0 13 6 /
\put{$5{,}040   $} [c] at 13 -2
\endpicture
\end{minipage}
\begin{minipage}{4cm}
\beginpicture
\setcoordinatesystem units   <1.5mm,2mm>
\setplotarea x from 0 to 16, y from -2 to 15
\put{762)} [l] at 2 12
\put {$ \scriptstyle \bullet$} [c] at  10 12
\put {$ \scriptstyle \bullet$} [c] at  13 6
\put {$ \scriptstyle \bullet$} [c] at  13 0
\put {$ \scriptstyle \bullet$} [c] at  14.5 9
\put {$ \scriptstyle \bullet$} [c] at  16 0
\put {$ \scriptstyle \bullet$} [c] at  16 6
\put {$ \scriptstyle \bullet$} [c] at  16 12
\setlinear \plot 16 0 16 12 13 6 13  0  /
\setlinear \plot 10  12 13 6 /
\put{$5{,}040   $} [c] at 13 -2
\endpicture
\end{minipage}
$$
$$
\begin{minipage}{4cm}
\beginpicture
\setcoordinatesystem units   <1.5mm,2mm>
\setplotarea x from 0 to 16, y from -2 to 15
\put{763)} [l] at 2 12
\put {$ \scriptstyle \bullet$} [c] at  10 6
\put {$ \scriptstyle \bullet$} [c] at  10 12
\put {$ \scriptstyle \bullet$} [c] at  10.5 3
\put {$ \scriptstyle \bullet$} [c] at  11 0
\put {$ \scriptstyle \bullet$} [c] at  11 12
\put {$ \scriptstyle \bullet$} [c] at  12 6
\put {$ \scriptstyle \bullet$} [c] at  16 0
\setlinear \plot  16 0 11 12  12 6 11 0 10 6 10 12    /
\setlinear \plot  10 6 11 12 /

\put{$5{,}040   $} [c] at 13 -2
\endpicture
\end{minipage}
\begin{minipage}{4cm}
\beginpicture
\setcoordinatesystem units   <1.5mm,2mm>
\setplotarea x from 0 to 16, y from -2 to 15
\put{764)} [l] at 2 12
\put {$ \scriptstyle \bullet$} [c] at  10 6
\put {$ \scriptstyle \bullet$} [c] at  10 0
\put {$ \scriptstyle \bullet$} [c] at  10.5 9
\put {$ \scriptstyle \bullet$} [c] at  11 0
\put {$ \scriptstyle \bullet$} [c] at  11 12
\put {$ \scriptstyle \bullet$} [c] at  12 6
\put {$ \scriptstyle \bullet$} [c] at  16 12
\setlinear \plot  16 12 11 0  12 6 11 12 10 6 10 0    /
\setlinear \plot  10 6 11 0 /

\put{$5{,}040   $} [c] at 13 -2
\endpicture
\end{minipage}
\begin{minipage}{4cm}
\beginpicture
\setcoordinatesystem units   <1.5mm,2mm>
\setplotarea x from 0 to 16, y from -2 to 15
\put{765)} [l] at 2 12
\put {$ \scriptstyle \bullet$} [c] at  10 6
\put {$ \scriptstyle \bullet$} [c] at  16 12
\put {$ \scriptstyle \bullet$} [c] at  10.5 3
\put {$ \scriptstyle \bullet$} [c] at  11 0
\put {$ \scriptstyle \bullet$} [c] at  11 12
\put {$ \scriptstyle \bullet$} [c] at  12 6
\put {$ \scriptstyle \bullet$} [c] at  16 0
\setlinear \plot  10 6 16 0 16 12 11 0 10 6 11 12 12 6 11 0   /
\put{$5{,}040   $} [c] at 13 -2
\endpicture
\end{minipage}
\begin{minipage}{4cm}
\beginpicture
\setcoordinatesystem units   <1.5mm,2mm>
\setplotarea x from 0 to 16, y from -2 to 15
\put{766)} [l] at 2 12
\put {$ \scriptstyle \bullet$} [c] at  10 6
\put {$ \scriptstyle \bullet$} [c] at  16 12
\put {$ \scriptstyle \bullet$} [c] at  10.5 9
\put {$ \scriptstyle \bullet$} [c] at  11 0
\put {$ \scriptstyle \bullet$} [c] at  11 12
\put {$ \scriptstyle \bullet$} [c] at  12 6
\put {$ \scriptstyle \bullet$} [c] at  16 0
\setlinear \plot  10 6 16 12 16 0 11 12 10 6 11 0 12 6 11 12   /
\put{$5{,}040   $} [c] at 13 -2
\endpicture
\end{minipage}
\begin{minipage}{4cm}
\beginpicture
\setcoordinatesystem units   <1.5mm,2mm>
\setplotarea x from 0 to 16, y from -2 to 15
\put{767)} [l] at 2 12
\put {$ \scriptstyle \bullet$} [c] at  10 6
\put {$ \scriptstyle \bullet$} [c] at  11 9
\put {$ \scriptstyle \bullet$} [c] at  12 0
\put {$ \scriptstyle \bullet$} [c] at  12 12
\put {$ \scriptstyle \bullet$} [c] at  14 6
\put {$ \scriptstyle \bullet$} [c] at  16 0
\put {$ \scriptstyle \bullet$} [c] at  16 12
\setlinear \plot 16 0 16 12 14 6 12 0 10 6 12 12 14 6 /
\put{$5{,}040   $} [c] at 13 -2
\endpicture
\end{minipage}
\begin{minipage}{4cm}
\beginpicture
\setcoordinatesystem units   <1.5mm,2mm>
\setplotarea x from 0 to 16, y from -2 to 15
\put{768)} [l] at 2 12
\put {$ \scriptstyle \bullet$} [c] at  10 6
\put {$ \scriptstyle \bullet$} [c] at  11 3
\put {$ \scriptstyle \bullet$} [c] at  12 0
\put {$ \scriptstyle \bullet$} [c] at  12 12
\put {$ \scriptstyle \bullet$} [c] at  14 6
\put {$ \scriptstyle \bullet$} [c] at  16 0
\put {$ \scriptstyle \bullet$} [c] at  16 12
\setlinear \plot 16 12 16 0 14 6 12 12 10 6 12 0 14 6 /
\put{$5{,}040   $} [c] at 13 -2
\endpicture
\end{minipage}
$$

$$
\begin{minipage}{4cm}
\beginpicture
\setcoordinatesystem units   <1.5mm,2mm>
\setplotarea x from 0 to 16, y from -2 to 15
\put{769)} [l] at 2 12
\put {$ \scriptstyle \bullet$} [c] at  10 0
\put {$ \scriptstyle \bullet$} [c] at  13 4
\put {$ \scriptstyle \bullet$} [c] at  14.5 0
\put {$ \scriptstyle \bullet$} [c] at  14.5 8
\put {$ \scriptstyle \bullet$} [c] at  14.5 12
\put {$ \scriptstyle \bullet$} [c] at  16 4
\put {$ \scriptstyle \bullet$} [c] at  16 12
\setlinear \plot  10 0 14.5 12 14.5 8 13 4 14.5 0 16 4 14.5 8  /
\setlinear \plot  16 4 16 12  /
\put{$5{,}040  $} [c] at 13 -2
\endpicture
\end{minipage}
\begin{minipage}{4cm}
\beginpicture
\setcoordinatesystem units   <1.5mm,2mm>
\setplotarea x from 0 to 16, y from -2 to 15
\put{770)} [l] at 2 12
\put {$ \scriptstyle \bullet$} [c] at  10 12
\put {$ \scriptstyle \bullet$} [c] at  13 8
\put {$ \scriptstyle \bullet$} [c] at  14.5 0
\put {$ \scriptstyle \bullet$} [c] at  14.5 4
\put {$ \scriptstyle \bullet$} [c] at  14.5 12
\put {$ \scriptstyle \bullet$} [c] at  16 8
\put {$ \scriptstyle \bullet$} [c] at  16 0
\setlinear \plot  10 12 14.5 0 14.5 4 16 8 14.5 12 13 8 14.5 4  /
\setlinear \plot  16 8 16 0  /
\put{$5{,}040  $} [c] at 13 -2
\endpicture
\end{minipage}
\begin{minipage}{4cm}
\beginpicture
\setcoordinatesystem units   <1.5mm,2mm>
\setplotarea x from 0 to 16, y from -2 to 15
\put{771)} [l] at 2 12
\put {$ \scriptstyle \bullet$} [c] at  10 0
\put {$ \scriptstyle \bullet$} [c] at  12 4
\put {$ \scriptstyle \bullet$} [c] at  12 8
\put {$ \scriptstyle \bullet$} [c] at  12 12
\put {$ \scriptstyle \bullet$} [c] at  14 0
\put {$ \scriptstyle \bullet$} [c] at  16 4
\put {$ \scriptstyle \bullet$} [c] at  16 12
\setlinear \plot 10  0 12 12 12 4  14 0 16 4 16 12 12 4      /
\setlinear \plot 12 4 12 12 16 4 /
\put{$5{,}040   $} [c] at 13 -2
\endpicture
\end{minipage}
\begin{minipage}{4cm}
\beginpicture
\setcoordinatesystem units   <1.5mm,2mm>
\setplotarea x from 0 to 16, y from -2 to 15
\put{772)} [l] at 2 12
\put {$ \scriptstyle \bullet$} [c] at  10 12
\put {$ \scriptstyle \bullet$} [c] at  12 0
\put {$ \scriptstyle \bullet$} [c] at  12 4
\put {$ \scriptstyle \bullet$} [c] at  12 8
\put {$ \scriptstyle \bullet$} [c] at  14 12
\put {$ \scriptstyle \bullet$} [c] at  16 0
\put {$ \scriptstyle \bullet$} [c] at  16 8
\setlinear \plot 10  12 12 0 12 8  14 12 16 8 16 0 12 8      /
\setlinear \plot 12 8 12 0 16 8 /
\put{$5{,}040   $} [c] at 13 -2
\endpicture
\end{minipage}
\begin{minipage}{4cm}
\beginpicture
\setcoordinatesystem units   <1.5mm,2mm>
\setplotarea x from 0 to 16, y from -2 to 15
\put{773)} [l] at 2 12
\put {$ \scriptstyle \bullet$} [c] at  10 0
\put {$ \scriptstyle \bullet$} [c] at  10 4
\put {$ \scriptstyle \bullet$} [c] at  10 12
\put {$ \scriptstyle \bullet$} [c] at  13 0
\put {$ \scriptstyle \bullet$} [c] at  13 8
\put {$ \scriptstyle \bullet$} [c] at  16 4
\put {$ \scriptstyle \bullet$} [c] at  16 12
\setlinear \plot 10   0 10 12  13 8 16  12 16 4  13 0 10 4 /
\setlinear \plot 13 0 13 8 /
\put{$5{,}040   $} [c] at 13 -2
\endpicture
\end{minipage}
\begin{minipage}{4cm}
\beginpicture
\setcoordinatesystem units   <1.5mm,2mm>
\setplotarea x from 0 to 16, y from -2 to 15
\put{774)} [l] at 2 12
\put {$ \scriptstyle \bullet$} [c] at  10 0
\put {$ \scriptstyle \bullet$} [c] at  10 8
\put {$ \scriptstyle \bullet$} [c] at  10 12
\put {$ \scriptstyle \bullet$} [c] at  13 12
\put {$ \scriptstyle \bullet$} [c] at  13 4
\put {$ \scriptstyle \bullet$} [c] at  16 8
\put {$ \scriptstyle \bullet$} [c] at  16 0
\setlinear \plot 10   12 10 0  13 4 16  0 16 8  13 12 10 8 /
\setlinear \plot 13 12 13 4 /
\put{$5{,}040$} [c] at 13 -2
\endpicture
\end{minipage}
$$
$$
\begin{minipage}{4cm}
\beginpicture
\setcoordinatesystem units   <1.5mm,2mm>
\setplotarea x from 0 to 16, y from -2 to 15
\put{775)} [l] at 2 12
\put {$ \scriptstyle \bullet$} [c] at  10 0
\put {$ \scriptstyle \bullet$} [c] at  10 12
\put {$ \scriptstyle \bullet$} [c] at  11.5 9
\put {$ \scriptstyle \bullet$} [c] at  13 6
\put {$ \scriptstyle \bullet$} [c] at  14.5 3
\put {$ \scriptstyle \bullet$} [c] at  16 0
\put {$ \scriptstyle \bullet$} [c] at  16 12
\setlinear \plot  10 0 10 12 16 0 16 12 /
\put{$5{,}040$} [c] at 13 -2
\endpicture
\end{minipage}
\begin{minipage}{4cm}
\beginpicture
\setcoordinatesystem units   <1.5mm,2mm>
\setplotarea x from 0 to 16, y from -2 to 15
\put{776)} [l] at 2 12
\put {$ \scriptstyle \bullet$} [c] at  10 6
\put {$ \scriptstyle \bullet$} [c] at  12 0
\put {$ \scriptstyle \bullet$} [c] at  12 6
\put {$ \scriptstyle \bullet$} [c] at  12 12
\put {$ \scriptstyle \bullet$} [c] at  14 6
\put {$ \scriptstyle \bullet$} [c] at  16 0
\put {$ \scriptstyle \bullet$} [c] at  16 12
\setlinear \plot  10 6  16 0 16 12 14 6 12 0 10 6  12 12 12 0  /
\setlinear \plot  14 6 12 12 /
\put{$5{,}040  $} [c] at 13 -2
\endpicture
\end{minipage}
\begin{minipage}{4cm}
\beginpicture
\setcoordinatesystem units   <1.5mm,2mm>
\setplotarea x from 0 to 16, y from -2 to 15
\put{777)} [l] at 2 12
\put {$ \scriptstyle \bullet$} [c] at  10 0
\put {$ \scriptstyle \bullet$} [c] at  10 12
\put {$ \scriptstyle \bullet$} [c] at  12 3
\put {$ \scriptstyle \bullet$} [c] at  12 9
\put {$ \scriptstyle \bullet$} [c] at  14 0
\put {$ \scriptstyle \bullet$} [c] at  14  12
\put {$ \scriptstyle \bullet$} [c] at  16 6
\setlinear \plot 10 0 12 9 14 12 16 6 14 0 12 3 12 9  /
\setlinear \plot 10 12 12 3 /
\put{$5{,}040   $} [c] at 13 -2
\endpicture
\end{minipage}
\begin{minipage}{4cm}
\beginpicture
\setcoordinatesystem units   <1.5mm,2mm>
\setplotarea x from 0 to 16, y from -2 to 15
\put{778)} [l] at 2 12
\put {$ \scriptstyle \bullet$} [c] at  10 12
\put {$ \scriptstyle \bullet$} [c] at  13 0
\put {$ \scriptstyle \bullet$} [c] at  13 4
\put {$ \scriptstyle \bullet$} [c] at  13 8
\put {$ \scriptstyle \bullet$} [c] at  13 12
\put {$ \scriptstyle \bullet$} [c] at  16 0
\put {$ \scriptstyle \bullet$} [c] at  16 12
\setlinear \plot  13 0 13 12 16 0 16 12 13 4 10 12 /
\put{$5{,}040$} [c] at 13 -2
\endpicture
\end{minipage}
\begin{minipage}{4cm}
\beginpicture
\setcoordinatesystem units   <1.5mm,2mm>
\setplotarea x from 0 to 16, y from -2 to 15
\put{779)} [l] at 2 12
\put {$ \scriptstyle \bullet$} [c] at  10 0
\put {$ \scriptstyle \bullet$} [c] at  13 0
\put {$ \scriptstyle \bullet$} [c] at  13 4
\put {$ \scriptstyle \bullet$} [c] at  13 8
\put {$ \scriptstyle \bullet$} [c] at  13 12
\put {$ \scriptstyle \bullet$} [c] at  16 0
\put {$ \scriptstyle \bullet$} [c] at  16 12
\setlinear \plot  13 12 13 0 16 12 16 0 13 8 10 0 /
\put{$5{,}040$} [c] at 13 -2
\endpicture
\end{minipage}
\begin{minipage}{4cm}
\beginpicture
\setcoordinatesystem units   <1.5mm,2mm>
\setplotarea x from 0 to 16, y from -2 to 15
\put{780)} [l] at 2 12
\put {$ \scriptstyle \bullet$} [c] at  16 0
\put {$ \scriptstyle \bullet$} [c] at  16 12
\put {$ \scriptstyle \bullet$} [c] at  13 12
\put {$ \scriptstyle \bullet$} [c] at  10 0
\put {$ \scriptstyle \bullet$} [c] at  10 4
\put {$ \scriptstyle \bullet$} [c] at  10 8
\put {$ \scriptstyle \bullet$} [c] at  10 12
\setlinear \plot  10 12 10 0  16 12 16 0 13 12 10 4 /
\put{$5{,}040$} [c] at 13 -2
\endpicture
\end{minipage}
$$

$$
\begin{minipage}{4cm}
\beginpicture
\setcoordinatesystem units   <1.5mm,2mm>
\setplotarea x from 0 to 16, y from -2 to 15
\put{781)} [l] at 2 12
\put {$ \scriptstyle \bullet$} [c] at  16 0
\put {$ \scriptstyle \bullet$} [c] at  16 12
\put {$ \scriptstyle \bullet$} [c] at  13 0
\put {$ \scriptstyle \bullet$} [c] at  10 0
\put {$ \scriptstyle \bullet$} [c] at  10 4
\put {$ \scriptstyle \bullet$} [c] at  10 8
\put {$ \scriptstyle \bullet$} [c] at  10 12
\setlinear \plot  10 0 10 12  16 0 16 12 13 0 10 8 /
\put{$5{,}040$} [c] at 13 -2
\endpicture
\end{minipage}
\begin{minipage}{4cm}
\beginpicture
\setcoordinatesystem units   <1.5mm,2mm>
\setplotarea x from 0 to 16, y from -2 to 15
\put{782)} [l] at 2 12
\put {$ \scriptstyle \bullet$} [c] at  10 0
\put {$ \scriptstyle \bullet$} [c] at  10 4
\put {$ \scriptstyle \bullet$} [c] at  10 8
\put {$ \scriptstyle \bullet$} [c] at  10 12
\put {$ \scriptstyle \bullet$} [c] at  13 12
\put {$ \scriptstyle \bullet$} [c] at  16 0
\put {$ \scriptstyle \bullet$} [c] at  16 12
\setlinear \plot  10 0 10 12 16  0 16 12 /
\setlinear \plot  10 4 13 12 16 0 /
\put{$5{,}040$} [c] at 13 -2
\endpicture
\end{minipage}
\begin{minipage}{4cm}
\beginpicture
\setcoordinatesystem units   <1.5mm,2mm>
\setplotarea x from 0 to 16, y from -2 to 15
\put{783)} [l] at 2 12
\put {$ \scriptstyle \bullet$} [c] at  10 0
\put {$ \scriptstyle \bullet$} [c] at  10 4
\put {$ \scriptstyle \bullet$} [c] at  10 8
\put {$ \scriptstyle \bullet$} [c] at  10 12
\put {$ \scriptstyle \bullet$} [c] at  13 0
\put {$ \scriptstyle \bullet$} [c] at  16 0
\put {$ \scriptstyle \bullet$} [c] at  16 12
\setlinear \plot  10 12 10 0 16  12 16 0 /
\setlinear \plot  10 8 13 0 16 12 /
\put{$5{,}040$} [c] at 13 -2
\endpicture
\end{minipage}
\begin{minipage}{4cm}
\beginpicture
\setcoordinatesystem units   <1.5mm,2mm>
\setplotarea x from 0 to 16, y from -2 to 15
\put{784)} [l] at 2 12
\put {$ \scriptstyle \bullet$} [c] at  13 12
\put {$ \scriptstyle \bullet$} [c] at  13 4
\put {$ \scriptstyle \bullet$} [c] at  13 8
\put {$ \scriptstyle \bullet$} [c] at  13 0
\put {$ \scriptstyle \bullet$} [c] at  16 0
\put {$ \scriptstyle \bullet$} [c] at  16 12
\put {$ \scriptstyle \bullet$} [c] at  10 12
\setlinear \plot  10 12 13 0  13 12   /
\setlinear \plot  13 4 16 12 16 0 13 8  /
\put{$5{,}040$} [c] at 13 -2
\endpicture
\end{minipage}
\begin{minipage}{4cm}
\beginpicture
\setcoordinatesystem units   <1.5mm,2mm>
\setplotarea x from 0 to 16, y from -2 to 15
\put{785)} [l] at 2 12
\put {$ \scriptstyle \bullet$} [c] at  13 0
\put {$ \scriptstyle \bullet$} [c] at  13 4
\put {$ \scriptstyle \bullet$} [c] at  13 8
\put {$ \scriptstyle \bullet$} [c] at  13 12
\put {$ \scriptstyle \bullet$} [c] at  16 0
\put {$ \scriptstyle \bullet$} [c] at  16 12
\put {$ \scriptstyle \bullet$} [c] at  10 0
\setlinear \plot  10 0 13 12  13 0   /
\setlinear \plot  13 8 16 0 16 12 13 4  /
\put{$5{,}040$} [c] at 13 -2
\endpicture
\end{minipage}
\begin{minipage}{4cm}
\beginpicture
\setcoordinatesystem units   <1.5mm,2mm>
\setplotarea x from 0 to 16, y from -2 to 15
\put{786)} [l] at 2 12
\put {$ \scriptstyle \bullet$} [c] at  10 0
\put {$ \scriptstyle \bullet$} [c] at  10 12
\put {$ \scriptstyle \bullet$} [c] at  13 12
\put {$ \scriptstyle \bullet$} [c] at  13 6
\put {$ \scriptstyle \bullet$} [c] at  14.5 0
\put {$ \scriptstyle \bullet$} [c] at  14.5 12
\put {$ \scriptstyle \bullet$} [c] at  16 6
\setlinear \plot 14.5 0 16 6 14.5 12 13 6 14.5 0 10 12 10 0 13 6 13 12  /
\put{$5{,}040$} [c] at 13 -2
\endpicture
\end{minipage}
$$
$$
\begin{minipage}{4cm}
\beginpicture
\setcoordinatesystem units   <1.5mm,2mm>
\setplotarea x from 0 to 16, y from -2 to 15
\put{787)} [l] at 2 12
\put {$ \scriptstyle \bullet$} [c] at  10 0
\put {$ \scriptstyle \bullet$} [c] at  10 12
\put {$ \scriptstyle \bullet$} [c] at  13 0
\put {$ \scriptstyle \bullet$} [c] at  13 6
\put {$ \scriptstyle \bullet$} [c] at  14.5 0
\put {$ \scriptstyle \bullet$} [c] at  14.5 12
\put {$ \scriptstyle \bullet$} [c] at  16 6
\setlinear \plot 14.5 12 16 6 14.5 0 13 6 14.5 12 10 0 10 12 13 6 13 0  /
\put{$5{,}040$} [c] at 13 -2
\endpicture
\end{minipage}
\begin{minipage}{4cm}
\beginpicture
\setcoordinatesystem units   <1.5mm,2mm>
\setplotarea x from 0 to 16, y from -2 to 15
\put{788)} [l] at 2 12
\put {$ \scriptstyle \bullet$} [c] at  10 0
\put {$ \scriptstyle \bullet$} [c] at  10 4
\put {$ \scriptstyle \bullet$} [c] at  10 12
\put {$ \scriptstyle \bullet$} [c] at  13 0
\put {$ \scriptstyle \bullet$} [c] at  13 8
\put {$ \scriptstyle \bullet$} [c] at  13 12
\put {$ \scriptstyle \bullet$} [c] at  16 12
\setlinear \plot  16 12 13 0  13 12  /
\setlinear \plot  10 12 10 0   /
\setlinear \plot  10 4 13  8    /
\put{$5{,}040$} [c] at 13 -2
\endpicture
\end{minipage}
\begin{minipage}{4cm}
\beginpicture
\setcoordinatesystem units   <1.5mm,2mm>
\setplotarea x from 0 to 16, y from -2 to 15
\put{789)} [l] at 2 12
\put {$ \scriptstyle \bullet$} [c] at  10 0
\put {$ \scriptstyle \bullet$} [c] at  10 8
\put {$ \scriptstyle \bullet$} [c] at  10 12
\put {$ \scriptstyle \bullet$} [c] at  13 0
\put {$ \scriptstyle \bullet$} [c] at  13 4
\put {$ \scriptstyle \bullet$} [c] at  13 12
\put {$ \scriptstyle \bullet$} [c] at  16 0
\setlinear \plot  16 0 13 12  13 0  /
\setlinear \plot  10 12 10 0   /
\setlinear \plot  10 8 13  4    /
\put{$5{,}040$} [c] at 13 -2
\endpicture
\end{minipage}
\begin{minipage}{4cm}
\beginpicture
\setcoordinatesystem units   <1.5mm,2mm>
\setplotarea x from 0 to 16, y from -2 to 15
\put{790)} [l] at 2 12
\put {$ \scriptstyle \bullet$} [c] at  10 6
\put {$ \scriptstyle \bullet$} [c] at  10 12
\put {$ \scriptstyle \bullet$} [c] at  10 0
\put {$ \scriptstyle \bullet$} [c] at  13 12
\put {$ \scriptstyle \bullet$} [c] at  16 0
\put {$ \scriptstyle \bullet$} [c] at  16 6
\put {$ \scriptstyle \bullet$} [c] at  16 12
\setlinear \plot  16 12 16 0 13 12 10 6 10 0 16 6 /
\setlinear \plot  10 12 10 6   /
\put{$5{,}040$} [c] at 13 -2
\endpicture
\end{minipage}
\begin{minipage}{4cm}
\beginpicture
\setcoordinatesystem units   <1.5mm,2mm>
\setplotarea x from 0 to 16, y from -2 to 15
\put{791)} [l] at 2 12
\put {$ \scriptstyle \bullet$} [c] at  10 6
\put {$ \scriptstyle \bullet$} [c] at  10 12
\put {$ \scriptstyle \bullet$} [c] at  10 0
\put {$ \scriptstyle \bullet$} [c] at  13 0
\put {$ \scriptstyle \bullet$} [c] at  16 0
\put {$ \scriptstyle \bullet$} [c] at  16 6
\put {$ \scriptstyle \bullet$} [c] at  16 12
\setlinear \plot  16 0 16 12 13 0 10 6 10 12 16 6 /
\setlinear \plot  10 0 10 6   /
\put{$5{,}040$} [c] at 13 -2
\endpicture
\end{minipage}
\begin{minipage}{4cm}
\beginpicture
\setcoordinatesystem units   <1.5mm,2mm>
\setplotarea x from 0 to 16, y from -2 to 15
\put{792)} [l] at 2 12
\put {$ \scriptstyle \bullet$} [c] at  13 0
\put {$ \scriptstyle \bullet$} [c] at  13 4
\put {$ \scriptstyle \bullet$} [c] at  13 8
\put {$ \scriptstyle \bullet$} [c] at  13 12
\put {$ \scriptstyle \bullet$} [c] at  16 0
\put {$ \scriptstyle \bullet$} [c] at  16 12
\put {$ \scriptstyle \bullet$} [c] at  10 12
\setlinear \plot  16 0 16  12  13 8 13 12   /
\setlinear \plot  13 0 13 4  10 12 /
\setlinear \plot  13 4 13 8      /
\put{$5{,}040$} [c] at 13 -2
\endpicture
\end{minipage}
$$

$$
\begin{minipage}{4cm}
\beginpicture
\setcoordinatesystem units   <1.5mm,2mm>
\setplotarea x from 0 to 16, y from -2 to 15
\put{793)} [l] at 2 12
\put {$ \scriptstyle \bullet$} [c] at  13 0
\put {$ \scriptstyle \bullet$} [c] at  13 4
\put {$ \scriptstyle \bullet$} [c] at  13 8
\put {$ \scriptstyle \bullet$} [c] at  13 12
\put {$ \scriptstyle \bullet$} [c] at  16 0
\put {$ \scriptstyle \bullet$} [c] at  16 12
\put {$ \scriptstyle \bullet$} [c] at  10 0
\setlinear \plot  16 12 16  0  13 4 13 0   /
\setlinear \plot  13 12 13 8  10 0 /
\setlinear \plot  13 4 13 8      /
\put{$5{,}040$} [c] at 13 -2
\endpicture
\end{minipage}
\begin{minipage}{4cm}
\beginpicture
\setcoordinatesystem units   <1.5mm,2mm>
\setplotarea x from 0 to 16, y from -2 to 15
\put{794)} [l] at 2 12
\put {$ \scriptstyle \bullet$} [c] at  10 0
\put {$ \scriptstyle \bullet$} [c] at  10 6
\put {$ \scriptstyle \bullet$} [c] at  10 12
\put {$ \scriptstyle \bullet$} [c] at  12 0
\put {$ \scriptstyle \bullet$} [c] at  12 6
\put {$ \scriptstyle \bullet$} [c] at  12 12
\put {$ \scriptstyle \bullet$} [c] at  16 12
\setlinear \plot 10 0 16 12 12 0 12 12 10 6 10 0 12 6     /
\setlinear \plot 10  6 10 12 12 0  /
\put{$5{,}040$} [c] at 13 -2
\endpicture
\end{minipage}
\begin{minipage}{4cm}
\beginpicture
\setcoordinatesystem units   <1.5mm,2mm>
\setplotarea x from 0 to 16, y from -2 to 15
\put{795)} [l] at 2 12
\put {$ \scriptstyle \bullet$} [c] at  10 0
\put {$ \scriptstyle \bullet$} [c] at  10 6
\put {$ \scriptstyle \bullet$} [c] at  10 12
\put {$ \scriptstyle \bullet$} [c] at  12 0
\put {$ \scriptstyle \bullet$} [c] at  12 6
\put {$ \scriptstyle \bullet$} [c] at  12 12
\put {$ \scriptstyle \bullet$} [c] at  16 0
\setlinear \plot 10 12 16 0 12 12 12 0 10 6 10 12 12 6     /
\setlinear \plot 10  6 10 0 12 12  /
\put{$5{,}040$} [c] at 13 -2
\endpicture
\end{minipage}
\begin{minipage}{4cm}
\beginpicture
\setcoordinatesystem units   <1.5mm,2mm>
\setplotarea x from 0 to 16, y from -2 to 15
\put{796)} [l] at 2 12
\put {$ \scriptstyle \bullet$} [c] at  10 0
\put {$ \scriptstyle \bullet$} [c] at  10 6
\put {$ \scriptstyle \bullet$} [c] at  10 12
\put {$ \scriptstyle \bullet$} [c] at  13 12
\put {$ \scriptstyle \bullet$} [c] at  16 0
\put {$ \scriptstyle \bullet$} [c] at  16 6
\put {$ \scriptstyle \bullet$} [c] at  16 12
\setlinear \plot  10 12 10 0 13 12 16 6 16  0 /
\setlinear \plot  10 12 16 6  16 12 /
\put{$5{,}040$} [c] at 13 -2
\endpicture
\end{minipage}
\begin{minipage}{4cm}
\beginpicture
\setcoordinatesystem units   <1.5mm,2mm>
\setplotarea x from 0 to 16, y from -2 to 15
\put{797)} [l] at 2 12
\put {$ \scriptstyle \bullet$} [c] at  10 0
\put {$ \scriptstyle \bullet$} [c] at  10 6
\put {$ \scriptstyle \bullet$} [c] at  10 12
\put {$ \scriptstyle \bullet$} [c] at  13 0
\put {$ \scriptstyle \bullet$} [c] at  16 0
\put {$ \scriptstyle \bullet$} [c] at  16 6
\put {$ \scriptstyle \bullet$} [c] at  16 12
\setlinear \plot  10 0 10 12 13 0 16 6 16  12 /
\setlinear \plot  10 0 16 6  16 0 /
\put{$5{,}040$} [c] at 13 -2
\endpicture
\end{minipage}
\begin{minipage}{4cm}
\beginpicture
\setcoordinatesystem units   <1.5mm,2mm>
\setplotarea x from 0 to 16, y from -2 to 15
\put{798)} [l] at 2 12
\put {$ \scriptstyle \bullet$} [c] at  10 12
\put {$ \scriptstyle \bullet$} [c] at  13 0
\put {$ \scriptstyle \bullet$} [c] at  13 6
\put {$ \scriptstyle \bullet$} [c] at  13 12
\put {$ \scriptstyle \bullet$} [c] at  16 0
\put {$ \scriptstyle \bullet$} [c] at  16 6
\put {$ \scriptstyle \bullet$} [c] at  16 12
\setlinear \plot  13 12 13 0 10 12 /
\setlinear \plot  13 0 16 6 16 0 13 6 16  12 16 6  /
\put{$5{,}040$} [c] at 13 -2
\endpicture
\end{minipage}
$$
$$
\begin{minipage}{4cm}
\beginpicture
\setcoordinatesystem units   <1.5mm,2mm>
\setplotarea x from 0 to 16, y from -2 to 15
\put{799)} [l] at 2 12
\put {$ \scriptstyle \bullet$} [c] at  10 0
\put {$ \scriptstyle \bullet$} [c] at  13 0
\put {$ \scriptstyle \bullet$} [c] at  13 6
\put {$ \scriptstyle \bullet$} [c] at  13 12
\put {$ \scriptstyle \bullet$} [c] at  16 0
\put {$ \scriptstyle \bullet$} [c] at  16 6
\put {$ \scriptstyle \bullet$} [c] at  16 12
\setlinear \plot  13 0 13 12 10 0 /
\setlinear \plot  13 12 16 6 16 0 13 6 16  12 16 6  /
\put{$5{,}040$} [c] at 13 -2
\endpicture
\end{minipage}
\begin{minipage}{4cm}
\beginpicture
\setcoordinatesystem units   <1.5mm,2mm>
\setplotarea x from 0 to 16, y from -2 to 15
\put{800)} [l] at 2 12
\put {$ \scriptstyle \bullet$} [c] at  13 0
\put {$ \scriptstyle \bullet$} [c] at  13 4
\put {$ \scriptstyle \bullet$} [c] at  13 8
\put {$ \scriptstyle \bullet$} [c] at  13 12
\put {$ \scriptstyle \bullet$} [c] at  16 0
\put {$ \scriptstyle \bullet$} [c] at  16 12
\put {$ \scriptstyle \bullet$} [c] at  10 12
\setlinear \plot 10 12 13 0 16 12 16 0 13 4 13 12  /
\setlinear \plot  13 0  13 4    /
\put{$5{,}040$} [c] at 13 -2
\endpicture
\end{minipage}
\begin{minipage}{4cm}
\beginpicture
\setcoordinatesystem units   <1.5mm,2mm>
\setplotarea x from 0 to 16, y from -2 to 15
\put{801)} [l] at 2 12
\put {$ \scriptstyle \bullet$} [c] at  13 0
\put {$ \scriptstyle \bullet$} [c] at  13 4
\put {$ \scriptstyle \bullet$} [c] at  13 8
\put {$ \scriptstyle \bullet$} [c] at  13 12
\put {$ \scriptstyle \bullet$} [c] at  16 0
\put {$ \scriptstyle \bullet$} [c] at  16 12
\put {$ \scriptstyle \bullet$} [c] at  10 0
\setlinear \plot 10 0 13 12 16 0 16 12 13 8 13 0  /
\setlinear \plot  13 12  13 8    /
\put{$5{,}040$} [c] at 13 -2
\endpicture
\end{minipage}
\begin{minipage}{4cm}
\beginpicture
\setcoordinatesystem units   <1.5mm,2mm>
\setplotarea x from 0 to 16, y from -2 to 15
\put{802)} [l] at 2 12
\put {$ \scriptstyle \bullet$} [c] at  10 0
\put {$ \scriptstyle \bullet$} [c] at  10 12
\put {$ \scriptstyle \bullet$} [c] at  13 12
\put {$ \scriptstyle \bullet$} [c] at  13 0
\put {$ \scriptstyle \bullet$} [c] at  16 0
\put {$ \scriptstyle \bullet$} [c] at  16  6
\put {$ \scriptstyle \bullet$} [c] at  16  12
\setlinear \plot  16  0 16  12     /
\setlinear \plot   16  6  13 12 10 0 10 12 13 0 16  6    /
\put{$5{,}040 $} [c] at 13 -2
\endpicture
\end{minipage}
\begin{minipage}{4cm}
\beginpicture
\setcoordinatesystem units   <1.5mm,2mm>
\setplotarea x from 0 to 16, y from -2 to 15
\put{803)} [l] at 2 12
\put {$ \scriptstyle \bullet$} [c] at 10 4
\put {$ \scriptstyle \bullet$} [c] at 10 8
\put {$ \scriptstyle \bullet$} [c] at  13 0
\put {$ \scriptstyle \bullet$} [c] at 13 6
\put {$ \scriptstyle \bullet$} [c] at 13 12
\put {$ \scriptstyle \bullet$} [c] at 16 4
\put {$ \scriptstyle \bullet$} [c] at 16 8
\setlinear \plot 13 0 10 4 10 8 13 12 16 8 16 4 13 0 13 12  /
\put{$2{,}520$} [c] at 13 -2
\endpicture
\end{minipage}
\begin{minipage}{4cm}
\beginpicture
\setcoordinatesystem units   <1.5mm,2mm>
\setplotarea x from 0 to 16, y from -2 to 15
\put{804)} [l] at 2 12
\put {$ \scriptstyle \bullet$} [c] at 10 4
\put {$ \scriptstyle \bullet$} [c] at 10 12
\put {$ \scriptstyle \bullet$} [c] at 13  0
\put {$ \scriptstyle \bullet$} [c] at 13 4
\put {$ \scriptstyle \bullet$} [c] at 16  4
\put {$ \scriptstyle \bullet$} [c] at 16 8
\put {$ \scriptstyle \bullet$} [c] at 16 12
\setlinear \plot 13 0 10 4 10 12 13 4 13 0 16 4 10 12 /
\setlinear \plot  16 12 16 4  /
\put{$2{,}520$} [c] at 13 -2
\endpicture
\end{minipage}
$$

$$
\begin{minipage}{4cm}
\beginpicture
\setcoordinatesystem units   <1.5mm,2mm>
\setplotarea x from 0 to 16, y from -2 to 15
\put{805)} [l] at 2 12
\put {$ \scriptstyle \bullet$} [c] at 10 8
\put {$ \scriptstyle \bullet$} [c] at 10 0
\put {$ \scriptstyle \bullet$} [c] at 13  12
\put {$ \scriptstyle \bullet$} [c] at 13 8
\put {$ \scriptstyle \bullet$} [c] at 16  8
\put {$ \scriptstyle \bullet$} [c] at 16 4
\put {$ \scriptstyle \bullet$} [c] at 16 0
\setlinear \plot 13 12 10 8 10 0 13 8 13 12 16 8 10 0 /
\setlinear \plot  16 0 16 8  /
\put{$2{,}520$} [c] at 13 -2
\endpicture
\end{minipage}
\begin{minipage}{4cm}
\beginpicture
\setcoordinatesystem units   <1.5mm,2mm>
\setplotarea x from 0 to 16, y from -2 to 15
\put{806)} [l] at 2 12
\put {$ \scriptstyle \bullet$} [c] at 10 6
\put {$ \scriptstyle \bullet$} [c] at 10 12
\put {$ \scriptstyle \bullet$} [c] at  13 0
\put {$ \scriptstyle \bullet$} [c] at 13 3
\put {$ \scriptstyle \bullet$} [c] at 13 12
\put {$ \scriptstyle \bullet$} [c] at 16 6
\put {$ \scriptstyle \bullet$} [c] at 16 12
\setlinear \plot 10 12 10 6 13 3 16 6 16 12  /
\setlinear \plot  13 0 13 12  /
\put{$2{,}520$} [c] at 13 -2
\endpicture
\end{minipage}
\begin{minipage}{4cm}
\beginpicture
\setcoordinatesystem units   <1.5mm,2mm>
\setplotarea x from 0 to 16, y from -2 to 15
\put{807)} [l] at 2 12
\put {$ \scriptstyle \bullet$} [c] at 10 6
\put {$ \scriptstyle \bullet$} [c] at 10 0
\put {$ \scriptstyle \bullet$} [c] at  13 0
\put {$ \scriptstyle \bullet$} [c] at 13 9
\put {$ \scriptstyle \bullet$} [c] at 13 12
\put {$ \scriptstyle \bullet$} [c] at 16 6
\put {$ \scriptstyle \bullet$} [c] at 16 0
\setlinear \plot 10 0 10 6 13 9 16 6 16 0  /
\setlinear \plot  13 0 13 12  /
\put{$2{,}520$} [c] at 13 -2
\endpicture
\end{minipage}
\begin{minipage}{4cm}
\beginpicture
\setcoordinatesystem units   <1.5mm,2mm>
\setplotarea x from 0 to 16, y from -2 to 15
\put{808)} [l] at 2 12
\put {$ \scriptstyle \bullet$} [c] at 10 6
\put {$ \scriptstyle \bullet$} [c] at 10 12
\put {$ \scriptstyle \bullet$} [c] at  13 0
\put {$ \scriptstyle \bullet$} [c] at 13 12
\put {$ \scriptstyle \bullet$} [c] at 16 6
\put {$ \scriptstyle \bullet$} [c] at 16 9
\put {$ \scriptstyle \bullet$} [c] at 16 12
\setlinear \plot 10 12 10 6 13 0 16 6 16 12 10 6 /
\setlinear \plot  13 12 10 6   /
\put{$2{,}520$} [c] at 13 -2
\endpicture
\end{minipage}
\begin{minipage}{4cm}
\beginpicture
\setcoordinatesystem units   <1.5mm,2mm>
\setplotarea x from 0 to 16, y from -2 to 15
\put{809)} [l] at 2 12
\put {$ \scriptstyle \bullet$} [c] at 10 6
\put {$ \scriptstyle \bullet$} [c] at 10 0
\put {$ \scriptstyle \bullet$} [c] at  13 0
\put {$ \scriptstyle \bullet$} [c] at 13 12
\put {$ \scriptstyle \bullet$} [c] at 16 0
\put {$ \scriptstyle \bullet$} [c] at 16 3
\put {$ \scriptstyle \bullet$} [c] at 16 6
\setlinear \plot 10 0 10 6 13 12 16 6 16 0 10 6 /
\setlinear \plot  13 0 10 6   /
\put{$2{,}520$} [c] at 13 -2
\endpicture
\end{minipage}
\begin{minipage}{4cm}
\beginpicture
\setcoordinatesystem units   <1.5mm,2mm>
\setplotarea x from 0 to 16, y from -2 to 15
\put{810)} [l] at 2 12
\put {$ \scriptstyle \bullet$} [c] at 10 6
\put {$ \scriptstyle \bullet$} [c] at 10 12
\put {$ \scriptstyle \bullet$} [c] at 13 0
\put {$ \scriptstyle \bullet$} [c] at 13 6
\put {$ \scriptstyle \bullet$} [c] at 13 12
\put {$ \scriptstyle \bullet$} [c] at 16 6
\put {$ \scriptstyle \bullet$} [c] at 16 12
\setlinear \plot 10 6 13  0 16 6 16 12 13 6 10 12 10 6 13 12 13 0  /
\put{$2{,}520$} [c] at 13 -2
\endpicture
\end{minipage}
$$
$$
\begin{minipage}{4cm}
\beginpicture
\setcoordinatesystem units   <1.5mm,2mm>
\setplotarea x from 0 to 16, y from -2 to 15
\put{811)} [l] at 2 12
\put {$ \scriptstyle \bullet$} [c] at 10 6
\put {$ \scriptstyle \bullet$} [c] at 10 0
\put {$ \scriptstyle \bullet$} [c] at 13 0
\put {$ \scriptstyle \bullet$} [c] at 13 6
\put {$ \scriptstyle \bullet$} [c] at 13 12
\put {$ \scriptstyle \bullet$} [c] at 16 6
\put {$ \scriptstyle \bullet$} [c] at 16 0
\setlinear \plot 10 6 13  12 16 6 16 0 13 6 10 0 10 6 13 0 13 12  /
\put{$2{,}520$} [c] at 13 -2
\endpicture
\end{minipage}
\begin{minipage}{4cm}
\beginpicture
\setcoordinatesystem units   <1.5mm,2mm>
\setplotarea x from 0 to 16, y from -2 to 15
\put{812)} [l] at 2 12
\put {$ \scriptstyle \bullet$} [c] at 10 6
\put {$ \scriptstyle \bullet$} [c] at 10 12
\put {$ \scriptstyle \bullet$} [c] at 13 0
\put {$ \scriptstyle \bullet$} [c] at 13 6
\put {$ \scriptstyle \bullet$} [c] at 13 12
\put {$ \scriptstyle \bullet$} [c] at 16 6
\put {$ \scriptstyle \bullet$} [c] at 16 12
\setlinear \plot 16 12 16 6 13  0 10 6 10 12 13 6 13 12 16 6   /
\setlinear \plot 10 6 13 12 /
\setlinear \plot 13 6 13 0 /
\put{$2{,}520$} [c] at 13 -2
\endpicture
\end{minipage}
\begin{minipage}{4cm}
\beginpicture
\setcoordinatesystem units   <1.5mm,2mm>
\setplotarea x from 0 to 16, y from -2 to 15
\put{813)} [l] at 2 12
\put {$ \scriptstyle \bullet$} [c] at 10 6
\put {$ \scriptstyle \bullet$} [c] at 10 0
\put {$ \scriptstyle \bullet$} [c] at 13 0
\put {$ \scriptstyle \bullet$} [c] at 13 6
\put {$ \scriptstyle \bullet$} [c] at 13 12
\put {$ \scriptstyle \bullet$} [c] at 16 6
\put {$ \scriptstyle \bullet$} [c] at 16 0
\setlinear \plot 16 0 16 6 13 12 10 6 10 0 13 6 13 0 16 6   /
\setlinear \plot 10 6 13 0 /
\setlinear \plot 13 6 13 12 /
\put{$2{,}520$} [c] at 13 -2
\endpicture
\end{minipage}
\begin{minipage}{4cm}
\beginpicture
\setcoordinatesystem units   <1.5mm,2mm>
\setplotarea x from 0 to 16, y from -2 to 15
\put{814)} [l] at 2 12
\put {$ \scriptstyle \bullet$} [c] at  10 0
\put {$ \scriptstyle \bullet$} [c] at  10 6
\put {$ \scriptstyle \bullet$} [c] at  10 12
\put {$ \scriptstyle \bullet$} [c] at  13 6
\put {$ \scriptstyle \bullet$} [c] at  16 0
\put {$ \scriptstyle \bullet$} [c] at  16 6
\put {$ \scriptstyle \bullet$} [c] at  16 12
\setlinear \plot 16 0 16 12  10  0 10 12 16 6      /
\put{$2{,}520   $} [c] at 13 -2
\endpicture
\end{minipage}
\begin{minipage}{4cm}
\beginpicture
\setcoordinatesystem units   <1.5mm,2mm>
\setplotarea x from 0 to 16, y from -2 to 15
\put{815)} [l] at 2 12
\put {$ \scriptstyle \bullet$} [c] at  10 0
\put {$ \scriptstyle \bullet$} [c] at  10 6
\put {$ \scriptstyle \bullet$} [c] at  10 12
\put {$ \scriptstyle \bullet$} [c] at  13 6
\put {$ \scriptstyle \bullet$} [c] at  16 0
\put {$ \scriptstyle \bullet$} [c] at  16 6
\put {$ \scriptstyle \bullet$} [c] at  16 12
\setlinear \plot 16 12 16 0  10  12 10 0 16 6      /
\put{$2{,}520   $} [c] at 13 -2
\endpicture
\end{minipage}
\begin{minipage}{4cm}
\beginpicture
\setcoordinatesystem units   <1.5mm,2mm>
\setplotarea x from 0 to 16, y from -2 to 15
\put{816)} [l] at 2 12
\put {$ \scriptstyle \bullet$} [c] at  10 6
\put {$ \scriptstyle \bullet$} [c] at  10 12
\put {$ \scriptstyle \bullet$} [c] at  11.5 0
\put {$ \scriptstyle \bullet$} [c] at  13 6
\put {$ \scriptstyle \bullet$} [c] at  13 12
\put {$ \scriptstyle \bullet$} [c] at  14.5 6
\put {$ \scriptstyle \bullet$} [c] at  16 0
\setlinear \plot 13 6 10 12 10 6 11.5 0  13 6 13 12 10 6      /
\setlinear \plot 16 0 13 12  /
\put{$2{,}520   $} [c] at 13 -2
\endpicture
\end{minipage}
$$

$$
\begin{minipage}{4cm}
\beginpicture
\setcoordinatesystem units   <1.5mm,2mm>
\setplotarea x from 0 to 16, y from -2 to 15
\put{817)} [l] at 2 12
\put {$ \scriptstyle \bullet$} [c] at  10 6
\put {$ \scriptstyle \bullet$} [c] at  10 0
\put {$ \scriptstyle \bullet$} [c] at  11.5 12
\put {$ \scriptstyle \bullet$} [c] at  13 6
\put {$ \scriptstyle \bullet$} [c] at  13 0
\put {$ \scriptstyle \bullet$} [c] at  14.5 6
\put {$ \scriptstyle \bullet$} [c] at  16 12
\setlinear \plot 13 6 10 0 10 6 11.5 12  13 6 13 0 10 6      /
\setlinear \plot 16 12 13 0  /
\put{$2{,}520   $} [c] at 13 -2
\endpicture
\end{minipage}
\begin{minipage}{4cm}
\beginpicture
\setcoordinatesystem units   <1.5mm,2mm>
\setplotarea x from 0 to 16, y from -2 to 15
\put{818)} [l] at 2 12
\put {$ \scriptstyle \bullet$} [c] at  10 8
\put {$ \scriptstyle \bullet$} [c] at  11 0
\put {$ \scriptstyle \bullet$} [c] at  11 4
\put {$ \scriptstyle \bullet$} [c] at  11  12
\put {$ \scriptstyle \bullet$} [c] at  12 8
\put {$ \scriptstyle \bullet$} [c] at  16 0
\put {$ \scriptstyle \bullet$} [c] at  16 12
\setlinear \plot 16 0 16 12 11 0 11 4 10 8 11 12 12 8 11 4 /
\put{$2{,}520  $} [c] at 13 -2
\endpicture
\end{minipage}
\begin{minipage}{4cm}
\beginpicture
\setcoordinatesystem units   <1.5mm,2mm>
\setplotarea x from 0 to 16, y from -2 to 15
\put{819)} [l] at 2 12
\put {$ \scriptstyle \bullet$} [c] at  10 4
\put {$ \scriptstyle \bullet$} [c] at  11 12
\put {$ \scriptstyle \bullet$} [c] at  11 8
\put {$ \scriptstyle \bullet$} [c] at  11  0
\put {$ \scriptstyle \bullet$} [c] at  12 4
\put {$ \scriptstyle \bullet$} [c] at  16 0
\put {$ \scriptstyle \bullet$} [c] at  16 12
\setlinear \plot 16 12 16 0 11 12 11 8 10 4 11 0 12 4 11 8 /
\put{$2{,}520  $} [c] at 13 -2
\endpicture
\end{minipage}
\begin{minipage}{4cm}
\beginpicture
\setcoordinatesystem units   <1.5mm,2mm>
\setplotarea x from 0 to 16, y from -2 to 15
\put{820)} [l] at 2 12
\put {$ \scriptstyle \bullet$} [c] at  10 4
\put {$ \scriptstyle \bullet$} [c] at  11 0
\put {$ \scriptstyle \bullet$} [c] at  11 8
\put {$ \scriptstyle \bullet$} [c] at  11 12
\put {$ \scriptstyle \bullet$} [c] at  12  4
\put {$ \scriptstyle \bullet$} [c] at  16  0
\put {$ \scriptstyle \bullet$} [c] at  16  12
\setlinear \plot  16 0 16 12 11 0  10 4 11 8 11 12    /
\setlinear \plot  11 0 12 4  11 8     /
\put{$2{,}520$} [c] at 13 -2
\endpicture
\end{minipage}
\begin{minipage}{4cm}
\beginpicture
\setcoordinatesystem units   <1.5mm,2mm>
\setplotarea x from 0 to 16, y from -2 to 15
\put{821)} [l] at 2 12
\put {$ \scriptstyle \bullet$} [c] at  10 8
\put {$ \scriptstyle \bullet$} [c] at  11 0
\put {$ \scriptstyle \bullet$} [c] at  11 4
\put {$ \scriptstyle \bullet$} [c] at  11 12
\put {$ \scriptstyle \bullet$} [c] at  12  8
\put {$ \scriptstyle \bullet$} [c] at  16  0
\put {$ \scriptstyle \bullet$} [c] at  16  12
\setlinear \plot  16 12 16 0 11 12 10 8 11 4 11 0    /
\setlinear \plot  11 4 12 8 11 12   /
\put{$2{,}520$} [c] at 13 -2
\endpicture
\end{minipage}
\begin{minipage}{4cm}
\beginpicture
\setcoordinatesystem units   <1.5mm,2mm>
\setplotarea x from 0 to 16, y from -2 to 15
\put{822)} [l] at 2 12
\put {$ \scriptstyle \bullet$} [c] at  10 6
\put {$ \scriptstyle \bullet$} [c] at  13 0
\put {$ \scriptstyle \bullet$} [c] at  13 6
\put {$ \scriptstyle \bullet$} [c] at  13 12
\put {$ \scriptstyle \bullet$} [c] at  16 0
\put {$ \scriptstyle \bullet$} [c] at  16 6
\put {$ \scriptstyle \bullet$} [c] at  16 12
\setlinear \plot  13 0 10 6 13 12  13 0 16 6 13 12    /
\setlinear \plot  16  12  16 0    /
\put{$2{,}520  $} [c] at 13 -2
\endpicture
\end{minipage}
$$
$$
\begin{minipage}{4cm}
\beginpicture
\setcoordinatesystem units   <1.5mm,2mm>
\setplotarea x from 0 to 16, y from -2 to 15
\put{823)} [l] at 2 12
\put {$ \scriptstyle \bullet$} [c] at  16 12
\put {$ \scriptstyle \bullet$} [c] at  16 0
\put {$ \scriptstyle \bullet$} [c] at  13 12
\put {$ \scriptstyle \bullet$} [c] at  10 0
\put {$ \scriptstyle \bullet$} [c] at  10 4
\put {$ \scriptstyle \bullet$} [c] at  10 8
\put {$ \scriptstyle \bullet$} [c] at  10 12
\setlinear \plot 10 8 16 0 16 12 10 0 10 12  /
\setlinear \plot  16 8 16 0 13 12 10 0      /
\put{$2{,}520$} [c] at 13 -2
\endpicture
\end{minipage}
\begin{minipage}{4cm}
\beginpicture
\setcoordinatesystem units   <1.5mm,2mm>
\setplotarea x from 0 to 16, y from -2 to 15
\put{824)} [l] at 2 12
\put {$ \scriptstyle \bullet$} [c] at  16 12
\put {$ \scriptstyle \bullet$} [c] at  16 0
\put {$ \scriptstyle \bullet$} [c] at  13 0
\put {$ \scriptstyle \bullet$} [c] at  10 0
\put {$ \scriptstyle \bullet$} [c] at  10 4
\put {$ \scriptstyle \bullet$} [c] at  10 8
\put {$ \scriptstyle \bullet$} [c] at  10 12
\setlinear \plot 10 4 16 12 16 0 10 12 10 0  /
\setlinear \plot  16 4 16 12 13 0 10 12      /
\put{$2{,}520$} [c] at 13 -2
\endpicture
\end{minipage}
\begin{minipage}{4cm}
\beginpicture
\setcoordinatesystem units   <1.5mm,2mm>
\setplotarea x from 0 to 16, y from -2 to 15
\put{825)} [l] at 2 12
\put {$ \scriptstyle \bullet$} [c] at  10 0
\put {$ \scriptstyle \bullet$} [c] at  10 6
\put {$ \scriptstyle \bullet$} [c] at  10 12
\put {$ \scriptstyle \bullet$} [c] at  13 12
\put {$ \scriptstyle \bullet$} [c] at  16 0
\put {$ \scriptstyle \bullet$} [c] at  16 6
\put {$ \scriptstyle \bullet$} [c] at  16  12
\setlinear \plot  10  12 10 0  /
\setlinear \plot  16  12  16 0  /
\setlinear \plot  10 6 13 12 16 6  /
\put{$2{,}520$} [c] at 13 -2
\endpicture
\end{minipage}
\begin{minipage}{4cm}
\beginpicture
\setcoordinatesystem units   <1.5mm,2mm>
\setplotarea x from 0 to 16, y from -2 to 15
\put{826)} [l] at 2 12
\put {$ \scriptstyle \bullet$} [c] at  10 0
\put {$ \scriptstyle \bullet$} [c] at  10 6
\put {$ \scriptstyle \bullet$} [c] at  10 12
\put {$ \scriptstyle \bullet$} [c] at  13 0
\put {$ \scriptstyle \bullet$} [c] at  16 0
\put {$ \scriptstyle \bullet$} [c] at  16 6
\put {$ \scriptstyle \bullet$} [c] at  16  12
\setlinear \plot  10  12 10 0  /
\setlinear \plot  16  12  16 0  /
\setlinear \plot  10 6 13 0 16 6  /
\put{$2{,}520$} [c] at 13 -2
\endpicture
\end{minipage}
\begin{minipage}{4cm}
\beginpicture
\setcoordinatesystem units   <1.5mm,2mm>
\setplotarea x from 0 to 16, y from -2 to 15
\put{827)} [l] at 2 12
\put {$ \scriptstyle \bullet$} [c] at  10 12
\put {$ \scriptstyle \bullet$} [c] at  13 0
\put {$ \scriptstyle \bullet$} [c] at  13 4
\put {$ \scriptstyle \bullet$} [c] at  13 8
\put {$ \scriptstyle \bullet$} [c] at  13 12
\put {$ \scriptstyle \bullet$} [c] at  16 0
\put {$ \scriptstyle \bullet$} [c] at  16 12
\setlinear \plot 10 12 13 0  13 12 16 0 16 12 13 8  /
 \put{$2{,}520$} [c] at 13 -2
\endpicture
\end{minipage}
\begin{minipage}{4cm}
\beginpicture
\setcoordinatesystem units   <1.5mm,2mm>
\setplotarea x from 0 to 16, y from -2 to 15
\put{828)} [l] at 2 12
\put {$ \scriptstyle \bullet$} [c] at  10 0
\put {$ \scriptstyle \bullet$} [c] at  13 0
\put {$ \scriptstyle \bullet$} [c] at  13 4
\put {$ \scriptstyle \bullet$} [c] at  13 8
\put {$ \scriptstyle \bullet$} [c] at  13 12
\put {$ \scriptstyle \bullet$} [c] at  16 0
\put {$ \scriptstyle \bullet$} [c] at  16 12
\setlinear \plot 10 0 13 12  13 0 16 12 16 0 13 4  /
 \put{$2{,}520$} [c] at 13 -2
\endpicture
\end{minipage}
$$

$$
\begin{minipage}{4cm}
\beginpicture
\setcoordinatesystem units   <1.5mm,2mm>
\setplotarea x from 0 to 16, y from -2 to 15
\put{829)} [l] at 2 12
\put {$ \scriptstyle \bullet$} [c] at  10 12
\put {$ \scriptstyle \bullet$} [c] at  10 4
\put {$ \scriptstyle \bullet$} [c] at  10 8
\put {$ \scriptstyle \bullet$} [c] at  10 0
\put {$ \scriptstyle \bullet$} [c] at  13 12
\put {$ \scriptstyle \bullet$} [c] at  16 0
\put {$ \scriptstyle \bullet$} [c] at  16 12
\setlinear \plot 10 12 10 0 16 12  16 0  /
\setlinear \plot  13 12 10  8     /
\put{$2{,}520$} [c] at 13 -2
\endpicture
\end{minipage}
\begin{minipage}{4cm}
\beginpicture
\setcoordinatesystem units   <1.5mm,2mm>
\setplotarea x from 0 to 16, y from -2 to 15
\put{830)} [l] at 2 12
\put {$ \scriptstyle \bullet$} [c] at  10 0
\put {$ \scriptstyle \bullet$} [c] at  10 4
\put {$ \scriptstyle \bullet$} [c] at  10 8
\put {$ \scriptstyle \bullet$} [c] at  10 12
\put {$ \scriptstyle \bullet$} [c] at  13 0
\put {$ \scriptstyle \bullet$} [c] at  16 0
\put {$ \scriptstyle \bullet$} [c] at  16 12
\setlinear \plot 10 0 10 12 16 0  16 12  /
\setlinear \plot  13 0 10  4     /
\put{$2{,}520$} [c] at 13 -2
\endpicture
\end{minipage}
\begin{minipage}{4cm}
\beginpicture
\setcoordinatesystem units   <1.5mm,2mm>
\setplotarea x from 0 to 16, y from -2 to 15
\put{831)} [l] at 2 12
\put {$ \scriptstyle \bullet$} [c] at  10 12
\put {$ \scriptstyle \bullet$} [c] at  13 0
\put {$ \scriptstyle \bullet$} [c] at  13 6
\put {$ \scriptstyle \bullet$} [c] at  13 12
\put {$ \scriptstyle \bullet$} [c] at  16 0
\put {$ \scriptstyle \bullet$} [c] at  16 9
\put {$ \scriptstyle \bullet$} [c] at  16  12
\setlinear \plot  10 12 13 6 16 9 16 12   /
\setlinear \plot  13 0 13 12  /
\setlinear \plot  16 9 16 0 /
\put{$2{,}520$} [c] at 13 -2
\endpicture
\end{minipage}
\begin{minipage}{4cm}
\beginpicture
\setcoordinatesystem units   <1.5mm,2mm>
\setplotarea x from 0 to 16, y from -2 to 15
\put{832)} [l] at 2 12
\put {$ \scriptstyle \bullet$} [c] at  10 0
\put {$ \scriptstyle \bullet$} [c] at  13 0
\put {$ \scriptstyle \bullet$} [c] at  13 9
\put {$ \scriptstyle \bullet$} [c] at  13 12
\put {$ \scriptstyle \bullet$} [c] at  16 0
\put {$ \scriptstyle \bullet$} [c] at  16 6
\put {$ \scriptstyle \bullet$} [c] at  16  12
\setlinear \plot  10 0 13 9 16 6 16 0   /
\setlinear \plot  13 0 13 12  /
\setlinear \plot  16 6 16 12 /
\put{$2{,}520$} [c] at 13 -2
\endpicture
\end{minipage}
\begin{minipage}{4cm}
\beginpicture
\setcoordinatesystem units   <1.5mm,2mm>
\setplotarea x from 0 to 16, y from -2 to 15
\put{833)} [l] at 2 12
\put {$ \scriptstyle \bullet$} [c] at  10 6
\put {$ \scriptstyle \bullet$} [c] at  11 0
\put {$ \scriptstyle \bullet$} [c] at  11 12
\put {$ \scriptstyle \bullet$} [c] at  12 6
\put {$ \scriptstyle \bullet$} [c] at  14 12
\put {$ \scriptstyle \bullet$} [c] at  16 0
\put {$ \scriptstyle \bullet$} [c] at  16 12
\setlinear \plot 11 12  10 6 11  0 12 6 11  12 16 0 16 12 12  6   /
\setlinear \plot 12 6 14 12 16 0   /
\put{$2{,}520$} [c] at 13 -2
\endpicture
\end{minipage}
\begin{minipage}{4cm}
\beginpicture
\setcoordinatesystem units   <1.5mm,2mm>
\setplotarea x from 0 to 16, y from -2 to 15
\put{834)} [l] at 2 12
\put {$ \scriptstyle \bullet$} [c] at  10 6
\put {$ \scriptstyle \bullet$} [c] at  11 0
\put {$ \scriptstyle \bullet$} [c] at  11 12
\put {$ \scriptstyle \bullet$} [c] at  12 6
\put {$ \scriptstyle \bullet$} [c] at  14 0
\put {$ \scriptstyle \bullet$} [c] at  16 0
\put {$ \scriptstyle \bullet$} [c] at  16 12
\setlinear \plot 11 0  10 6 11  12 12 6 11  0 16 12 16 0 12  6   /
\setlinear \plot 12 6 14 0 16 12   /
\put{$2{,}520$} [c] at 13 -2
\endpicture
\end{minipage}
$$
$$
\begin{minipage}{4cm}
\beginpicture
\setcoordinatesystem units   <1.5mm,2mm>
\setplotarea x from 0 to 16, y from -2 to 15
\put{835)} [l] at 2 12
\put {$ \scriptstyle \bullet$} [c] at  10 6
\put {$ \scriptstyle \bullet$} [c] at  10 12
\put {$ \scriptstyle \bullet$} [c] at  13 12
\put {$ \scriptstyle \bullet$} [c] at  13 0
\put {$ \scriptstyle \bullet$} [c] at  16 6
\put {$ \scriptstyle \bullet$} [c] at  16 12
\put {$ \scriptstyle \bullet$} [c] at  16 0
\setlinear \plot 10 6  10 12 16 0 16  12 16 6 13 12 10 6 13 0 16 6 13 12   /
\put{$2{,}520$} [c] at 13 -2
\endpicture
\end{minipage}
\begin{minipage}{4cm}
\beginpicture
\setcoordinatesystem units   <1.5mm,2mm>
\setplotarea x from 0 to 16, y from -2 to 15
\put{836)} [l] at 2 12
\put {$ \scriptstyle \bullet$} [c] at  10 0
\put {$ \scriptstyle \bullet$} [c] at  10 6
\put {$ \scriptstyle \bullet$} [c] at  13 12
\put {$ \scriptstyle \bullet$} [c] at  13 0
\put {$ \scriptstyle \bullet$} [c] at  16 6
\put {$ \scriptstyle \bullet$} [c] at  16 12
\put {$ \scriptstyle \bullet$} [c] at  16 0
\setlinear \plot 10 6  10 0 16 12 16  0 16 6 13 0 10 6 13 12 16 6 13 0   /
\put{$2{,}520$} [c] at 13 -2
\endpicture
\end{minipage}
\begin{minipage}{4cm}
\beginpicture
\setcoordinatesystem units   <1.5mm,2mm>
\setplotarea x from 0 to 16, y from -2 to 15
\put{837)} [l] at 2 12
\put {$ \scriptstyle \bullet$} [c] at  10 12
\put {$ \scriptstyle \bullet$} [c] at  11.5 12
\put {$ \scriptstyle \bullet$} [c] at  11.5 4
\put {$ \scriptstyle \bullet$} [c] at  12.5 0
\put {$ \scriptstyle \bullet$} [c] at  13.5 12
\put {$ \scriptstyle \bullet$} [c] at  13.5 4
\put {$ \scriptstyle \bullet$} [c] at  16 0
\setlinear \plot 10 12 11.5 4 12.5 0 13.5 4  13.5 12  11.5 4 11.5 12 13.5 4   /
\setlinear \plot 11.5 12 16 0 13.5  12   /
\put{$2{,}520$} [c] at 13 -2
\endpicture
\end{minipage}
\begin{minipage}{4cm}
\beginpicture
\setcoordinatesystem units   <1.5mm,2mm>
\setplotarea x from 0 to 16, y from -2 to 15
\put{838)} [l] at 2 12
\put {$ \scriptstyle \bullet$} [c] at  10 0
\put {$ \scriptstyle \bullet$} [c] at  11.5 0
\put {$ \scriptstyle \bullet$} [c] at  11.5 8
\put {$ \scriptstyle \bullet$} [c] at  12.5 12
\put {$ \scriptstyle \bullet$} [c] at  13.5 0
\put {$ \scriptstyle \bullet$} [c] at  13.5 8
\put {$ \scriptstyle \bullet$} [c] at  16 12
\setlinear \plot 10 0 11.5 8 12.5 12 13.5 8  13.5 0  11.5 8 11.5 0 13.5 8   /
\setlinear \plot 11.5 0 16 12 13.5  0   /
\put{$2{,}520$} [c] at 13 -2
\endpicture
\end{minipage}
\begin{minipage}{4cm}
\beginpicture
\setcoordinatesystem units   <1.5mm,2mm>
\setplotarea x from 0 to 16, y from -2 to 15
\put{839)} [l] at 2 12
\put {$ \scriptstyle \bullet$} [c] at  10 12
\put {$ \scriptstyle \bullet$} [c] at  13 12
\put {$ \scriptstyle \bullet$} [c] at  13 6
\put {$ \scriptstyle \bullet$} [c] at  14.5 0
\put {$ \scriptstyle \bullet$} [c] at  16 12
\put {$ \scriptstyle \bullet$} [c] at  16 6
\put {$ \scriptstyle \bullet$} [c] at  16 0
\setlinear \plot 10 12 14.5 0 13 6 13 12 16 6 16 12 13 6 /
\setlinear \plot  16 12 16 0 /
\setlinear \plot  14.5 0 16 6 /
\put{$2{,}520$} [c] at 13 -2
\endpicture
\end{minipage}
\begin{minipage}{4cm}
\beginpicture
\setcoordinatesystem units   <1.5mm,2mm>
\setplotarea x from 0 to 16, y from -2 to 15
\put{840)} [l] at 2 12
\put {$ \scriptstyle \bullet$} [c] at  10 0
\put {$ \scriptstyle \bullet$} [c] at  13 0
\put {$ \scriptstyle \bullet$} [c] at  13 6
\put {$ \scriptstyle \bullet$} [c] at  14.5 12
\put {$ \scriptstyle \bullet$} [c] at  16 12
\put {$ \scriptstyle \bullet$} [c] at  16 6
\put {$ \scriptstyle \bullet$} [c] at  16 0
\setlinear \plot 10 0 14.5 12 13 6 13 0 16 6 16 0 13 6 /
\setlinear \plot  16 12 16 0 /
\setlinear \plot  14.5 12 16 6 /
\put{$2{,}520$} [c] at 13 -2
\endpicture
\end{minipage}
$$

$$
\begin{minipage}{4cm}
\beginpicture
\setcoordinatesystem units   <1.5mm,2mm>
\setplotarea x from 0 to 16, y from -2 to 15
\put{841)} [l] at 2 12
\put {$ \scriptstyle \bullet$} [c] at  10 0
\put {$ \scriptstyle \bullet$} [c] at  10 6
\put {$ \scriptstyle \bullet$} [c] at  10 12
\put {$ \scriptstyle \bullet$} [c] at  13 12
\put {$ \scriptstyle \bullet$} [c] at  16 0
\put {$ \scriptstyle \bullet$} [c] at  16 6
\put {$ \scriptstyle \bullet$} [c] at  16 12
\setlinear \plot  10 12 10 0 16  12 16 0 10 12 /
\setlinear \plot  16 6  13 12 10 0 16 12    /
\put{$2{,}520$} [c] at 13 -2
\endpicture
\end{minipage}
\begin{minipage}{4cm}
\beginpicture
\setcoordinatesystem units   <1.5mm,2mm>
\setplotarea x from 0 to 16, y from -2 to 15
\put{842)} [l] at 2 12
\put {$ \scriptstyle \bullet$} [c] at  10 0
\put {$ \scriptstyle \bullet$} [c] at  10 6
\put {$ \scriptstyle \bullet$} [c] at  10 12
\put {$ \scriptstyle \bullet$} [c] at  13 0
\put {$ \scriptstyle \bullet$} [c] at  16 0
\put {$ \scriptstyle \bullet$} [c] at  16 6
\put {$ \scriptstyle \bullet$} [c] at  16 12
\setlinear \plot  10 0 10 12 16  0 16 12 10 0 /
\setlinear \plot  16 6  13 0 10 12 16 0    /
\put{$2{,}520$} [c] at 13 -2
\endpicture
\end{minipage}
\begin{minipage}{4cm}
\beginpicture
\setcoordinatesystem units   <1.5mm,2mm>
\setplotarea x from 0 to 16, y from -2 to 15
\put{843)} [l] at 2 12
\put {$ \scriptstyle \bullet$} [c] at  10 0
\put {$ \scriptstyle \bullet$} [c] at  10 12
\put {$ \scriptstyle \bullet$} [c] at  13 12
\put {$ \scriptstyle \bullet$} [c] at  13 6
\put {$ \scriptstyle \bullet$} [c] at  16 0
\put {$ \scriptstyle \bullet$} [c] at  16 6
\put {$ \scriptstyle \bullet$} [c] at  16 12
\setlinear \plot 10 0 13 6  10 12  /
\setlinear \plot 13  12  13 6 16 0 16 12   /
\put{$2{,}520$} [c] at 13 -2
\endpicture
\end{minipage}
\begin{minipage}{4cm}
\beginpicture
\setcoordinatesystem units   <1.5mm,2mm>
\setplotarea x from 0 to 16, y from -2 to 15
\put{844)} [l] at 2 12
\put {$ \scriptstyle \bullet$} [c] at  10 0
\put {$ \scriptstyle \bullet$} [c] at  10 12
\put {$ \scriptstyle \bullet$} [c] at  13 0
\put {$ \scriptstyle \bullet$} [c] at  13 6
\put {$ \scriptstyle \bullet$} [c] at  16 0
\put {$ \scriptstyle \bullet$} [c] at  16 6
\put {$ \scriptstyle \bullet$} [c] at  16 12
\setlinear \plot 10 12 13 6  10 0  /
\setlinear \plot 13  0  13 6 16 12 16 0   /
\put{$2{,}520$} [c] at 13 -2
\endpicture
\end{minipage}
\begin{minipage}{4cm}
\beginpicture
\setcoordinatesystem units   <1.5mm,2mm>
\setplotarea x from 0 to 16, y from -2 to 15
\put{845)} [l] at 2 12
\put {$ \scriptstyle \bullet$} [c] at  10 0
\put {$ \scriptstyle \bullet$} [c] at  10 12
\put {$ \scriptstyle \bullet$} [c] at  13 12
\put {$ \scriptstyle \bullet$} [c] at  13 0
\put {$ \scriptstyle \bullet$} [c] at  16 0
\put {$ \scriptstyle \bullet$} [c] at  16  6
\put {$ \scriptstyle \bullet$} [c] at  16  12
\setlinear \plot 10 0 10 12  13 0 16 6 16 12   /
\setlinear \plot 10 12 16 0 16 6  /
\setlinear \plot 10 0 13 12 16 0 /
\setlinear \plot 10 0 16 6 /
\put{$2{,}520$} [c] at 13 -2
\endpicture
\end{minipage}
\begin{minipage}{4cm}
\beginpicture
\setcoordinatesystem units   <1.5mm,2mm>
\setplotarea x from 0 to 16, y from -2 to 15
\put{846)} [l] at 2 12
\put {$ \scriptstyle \bullet$} [c] at  10 0
\put {$ \scriptstyle \bullet$} [c] at  10 12
\put {$ \scriptstyle \bullet$} [c] at  13 12
\put {$ \scriptstyle \bullet$} [c] at  13 0
\put {$ \scriptstyle \bullet$} [c] at  16 0
\put {$ \scriptstyle \bullet$} [c] at  16  6
\put {$ \scriptstyle \bullet$} [c] at  16  12
\setlinear \plot 10 0 10 12  13 0 16 12 16 6 16 0  /
\setlinear \plot 10 0 13 12 16 6  /
\setlinear \plot 10 12 16 6 /
\setlinear \plot 10 0 16 12 /
\put{$2{,}520$} [c] at 13 -2
\endpicture
\end{minipage}
$$
$$
\begin{minipage}{4cm}
\beginpicture
\setcoordinatesystem units   <1.5mm,2mm>
\setplotarea x from 0 to 16, y from -2 to 15
\put{847)} [l] at 2 12
\put {$ \scriptstyle \bullet$} [c] at  10 0
\put {$ \scriptstyle \bullet$} [c] at  10 12
\put {$ \scriptstyle \bullet$} [c] at  13 12
\put {$ \scriptstyle \bullet$} [c] at  13 9
\put {$ \scriptstyle \bullet$} [c] at  16 12
\put {$ \scriptstyle \bullet$} [c] at  16 9
\put {$ \scriptstyle \bullet$} [c] at  16 0
\setlinear \plot 13 12 10 0 10 12 16 0 16 9 16 12 10 0 /
\setlinear \plot  16 0 13  9 13 12 16 9  /
\setlinear \plot 16 12 13 9    /
\put{$1{,}260$} [c] at 13 -2
\endpicture
\end{minipage}
\begin{minipage}{4cm}
\beginpicture
\setcoordinatesystem units   <1.5mm,2mm>
\setplotarea x from 0 to 16, y from -2 to 15
\put{848)} [l] at 2 12
\put {$ \scriptstyle \bullet$} [c] at  10 0
\put {$ \scriptstyle \bullet$} [c] at  10 12
\put {$ \scriptstyle \bullet$} [c] at  13 0
\put {$ \scriptstyle \bullet$} [c] at  13 3
\put {$ \scriptstyle \bullet$} [c] at  16 12
\put {$ \scriptstyle \bullet$} [c] at  16 3
\put {$ \scriptstyle \bullet$} [c] at  16 0
\setlinear \plot 13 0 10 12 10 0 16 12 16 3 16 0 10 12 /
\setlinear \plot  16 12 13  3 13 0 16 3  /
\setlinear \plot 16 0 13 3    /
\put{$1{,}260$} [c] at 13 -2
\endpicture
\end{minipage}
\begin{minipage}{4cm}
\beginpicture
\setcoordinatesystem units   <1.5mm,2mm>
\setplotarea x from 0 to 16, y from -2 to 15
\put{849)} [l] at 2 12
\put {$ \scriptstyle \bullet$} [c] at  10 12
\put {$ \scriptstyle \bullet$} [c] at  10 6
\put {$ \scriptstyle \bullet$} [c] at  12 12
\put {$ \scriptstyle \bullet$} [c] at  12 0
\put {$ \scriptstyle \bullet$} [c] at  14 12
\put {$ \scriptstyle \bullet$} [c] at  14 6
\put {$ \scriptstyle \bullet$} [c] at  16 0
\setlinear \plot  16 0 14 12 14 6 12 12 10 6 12 0 14 6 10 12 10 6 14 12 /
\put{$1{,}260$} [c] at 13 -2
\endpicture
\end{minipage}
\begin{minipage}{4cm}
\beginpicture
\setcoordinatesystem units   <1.5mm,2mm>
\setplotarea x from 0 to 16, y from -2 to 15
\put{850)} [l] at 2 12
\put {$ \scriptstyle \bullet$} [c] at  10 0
\put {$ \scriptstyle \bullet$} [c] at  10 6
\put {$ \scriptstyle \bullet$} [c] at  12 12
\put {$ \scriptstyle \bullet$} [c] at  12 0
\put {$ \scriptstyle \bullet$} [c] at  14 0
\put {$ \scriptstyle \bullet$} [c] at  14 6
\put {$ \scriptstyle \bullet$} [c] at  16 12
\setlinear \plot  16 12 14 0 14 6 12 0 10 6 12 12 14 6 10 0 10 6 14 0 /
\put{$1{,}260$} [c] at 13 -2
\endpicture
\end{minipage}
\begin{minipage}{4cm}
\beginpicture
\setcoordinatesystem units   <1.5mm,2mm>
\setplotarea x from 0 to 16, y from -2 to 15
\put{851)} [l] at 2 12
\put {$ \scriptstyle \bullet$} [c] at  10 0
\put {$ \scriptstyle \bullet$} [c] at  10 12
\put {$ \scriptstyle \bullet$} [c] at  13 12
\put {$ \scriptstyle \bullet$} [c] at  13 0
\put {$ \scriptstyle \bullet$} [c] at  16 0
\put {$ \scriptstyle \bullet$} [c] at  16  6
\put {$ \scriptstyle \bullet$} [c] at  16  12
\setlinear \plot 10 0 10 12  16 0 16 12 13 0 13 12 10 0  /
\setlinear \plot 10 12 13 0  /
\setlinear \plot 10 0 16 12 /
\setlinear \plot 13 12 16 6 /
\put{$1{,}260 $} [c] at 13 -2
\endpicture
\end{minipage}
\begin{minipage}{4cm}
\beginpicture
\setcoordinatesystem units   <1.5mm,2mm>
\setplotarea x from 0 to 16, y from -2 to 15
\put{852)} [l] at 2 12
\put {$ \scriptstyle \bullet$} [c] at  10 0
\put {$ \scriptstyle \bullet$} [c] at  10 12
\put {$ \scriptstyle \bullet$} [c] at  13 12
\put {$ \scriptstyle \bullet$} [c] at  13 0
\put {$ \scriptstyle \bullet$} [c] at  16 0
\put {$ \scriptstyle \bullet$} [c] at  16  6
\put {$ \scriptstyle \bullet$} [c] at  16  12
\setlinear \plot 10 12 10 0  16 12 16 0 13 12 13 0 10 12  /
\setlinear \plot 10 0 13 12  /
\setlinear \plot 10 12 16 0 /
\setlinear \plot 13 0 16 6 /
\put{$1{,}260 $} [c] at 13 -2
\endpicture
\end{minipage}
$$

$$
\begin{minipage}{4cm}
\beginpicture
\setcoordinatesystem units   <1.5mm,2mm>
\setplotarea x from 0 to 16, y from -2 to 15
\put{853)} [l] at 2 12
\put {$ \scriptstyle \bullet$} [c] at  10  0
\put {$ \scriptstyle \bullet$} [c] at  10 12
\put {$ \scriptstyle \bullet$} [c] at  13 0
\put {$ \scriptstyle \bullet$} [c] at  13 12
\put {$ \scriptstyle \bullet$} [c] at  13 6
\put {$ \scriptstyle \bullet$} [c] at  16  0
\put {$ \scriptstyle \bullet$} [c] at  16 12
\setlinear \plot  10 12  13 6 10 0    /
\setlinear \plot  13 0 13 12   /
\setlinear \plot   16 0 16 12 13 6  /
\put{$1{,}260$} [c] at 13 -2
\endpicture
\end{minipage}
\begin{minipage}{4cm}
\beginpicture
\setcoordinatesystem units   <1.5mm,2mm>
\setplotarea x from 0 to 16, y from -2 to 15
\put{854)} [l] at 2 12
\put {$ \scriptstyle \bullet$} [c] at  10  0
\put {$ \scriptstyle \bullet$} [c] at  10 12
\put {$ \scriptstyle \bullet$} [c] at  13 0
\put {$ \scriptstyle \bullet$} [c] at  13 12
\put {$ \scriptstyle \bullet$} [c] at  13 6
\put {$ \scriptstyle \bullet$} [c] at  16  0
\put {$ \scriptstyle \bullet$} [c] at  16 12
\setlinear \plot  10 0  13 6 10 12    /
\setlinear \plot  13 0 13 12   /
\setlinear \plot   16 12 16 0 13 6  /
\put{$1{,}260$} [c] at 13 -2
\endpicture
\end{minipage}
\begin{minipage}{4cm}
\beginpicture
\setcoordinatesystem units   <1.5mm,2mm>
\setplotarea x from 0 to 16, y from -2 to 15
\put{855)} [l] at 2 12
\put {$ \scriptstyle \bullet$} [c] at  10 0
\put {$ \scriptstyle \bullet$} [c] at  10 12
\put {$ \scriptstyle \bullet$} [c] at  13 0
\put {$ \scriptstyle \bullet$} [c] at  13 12
\put {$ \scriptstyle \bullet$} [c] at  13 6
\put {$ \scriptstyle \bullet$} [c] at  16  0
\put {$ \scriptstyle \bullet$} [c] at  16  12
\setlinear \plot  10  0 10  12  16 0 16  12 13 0 13 12  10 0   /
\put{$1{,}260 $} [c] at 13 -2
\endpicture
\end{minipage}
\begin{minipage}{4cm}
\beginpicture
\setcoordinatesystem units   <1.5mm,2mm>
\setplotarea x from 0 to 16, y from -2 to 15
\put{856)} [l] at 2 12
\put {$ \scriptstyle \bullet$} [c] at 10 8
\put {$ \scriptstyle \bullet$} [c] at 13 0
\put {$ \scriptstyle \bullet$} [c] at 13 4
\put {$ \scriptstyle \bullet$} [c] at 13 12
\put {$ \scriptstyle \bullet$} [c] at 12 8
\put {$ \scriptstyle \bullet$} [c] at 14 8
\put {$ \scriptstyle \bullet$} [c] at 16 8
\setlinear \plot 13 0 13 4 10 8 13 12 16 8 13 0 13 4 12 8 13 12 14 8 13 4 /
\put{$840$} [c] at 13 -2
\endpicture
\end{minipage}
\begin{minipage}{4cm}
\beginpicture
\setcoordinatesystem units   <1.5mm,2mm>
\setplotarea x from 0 to 16, y from -2 to 15
\put{857)} [l] at 2 12
\put {$ \scriptstyle \bullet$} [c] at 10 4
\put {$ \scriptstyle \bullet$} [c] at 13 12
\put {$ \scriptstyle \bullet$} [c] at 13 8
\put {$ \scriptstyle \bullet$} [c] at 13 0
\put {$ \scriptstyle \bullet$} [c] at 12 4
\put {$ \scriptstyle \bullet$} [c] at 14 4
\put {$ \scriptstyle \bullet$} [c] at 16 4
\setlinear \plot 13 12 13 8 10 4 13 0 16 4 13 12 13 8 12 4 13 0 14 4 13 8 /
\put{$840$} [c] at 13 -2
\endpicture
\end{minipage}
\begin{minipage}{4cm}
\beginpicture
\setcoordinatesystem units   <1.5mm,2mm>
\setplotarea x from 0 to 16, y from -2 to 15
\put{858)} [l] at 2 12
\put {$ \scriptstyle \bullet$} [c] at 10 8
\put {$ \scriptstyle \bullet$} [c] at 11.5 0
\put {$ \scriptstyle \bullet$} [c] at 11.5 4
\put {$ \scriptstyle \bullet$} [c] at 11.5 8
\put {$ \scriptstyle \bullet$} [c] at 11.5 12
\put {$ \scriptstyle \bullet$} [c] at 13 8
\put {$ \scriptstyle \bullet$} [c] at 16 12
\setlinear \plot 11.5 4 10 8 11.5 12 13 8 11.5 4 11.5 12 /
\setlinear \plot 16 12 11.5 4 11.5 0 /
\put{$840$} [c] at 13 -2
\endpicture
\end{minipage}
$$
$$
\begin{minipage}{4cm}
\beginpicture
\setcoordinatesystem units   <1.5mm,2mm>
\setplotarea x from 0 to 16, y from -2 to 15
\put{859)} [l] at 2 12
\put {$ \scriptstyle \bullet$} [c] at 10 4
\put {$ \scriptstyle \bullet$} [c] at 11.5 0
\put {$ \scriptstyle \bullet$} [c] at 11.5 4
\put {$ \scriptstyle \bullet$} [c] at 11.5 8
\put {$ \scriptstyle \bullet$} [c] at 11.5 12
\put {$ \scriptstyle \bullet$} [c] at 13 4
\put {$ \scriptstyle \bullet$} [c] at 16 0
\setlinear \plot 11.5 8 10 4 11.5 0 13 4 11.5 8 11.5 0 /
\setlinear \plot 16 0 11.5 8 11.5 12 /
\put{$840$} [c] at 13 -2
\endpicture
\end{minipage}
\begin{minipage}{4cm}
\beginpicture
\setcoordinatesystem units   <1.5mm,2mm>
\setplotarea x from 0 to 16, y from -2 to 15
\put{860)} [l] at 2 12
\put {$ \scriptstyle \bullet$} [c] at 10 12
\put {$ \scriptstyle \bullet$} [c] at 14 0
\put {$ \scriptstyle \bullet$} [c] at 14 4
\put {$ \scriptstyle \bullet$} [c] at 14 8
\put {$ \scriptstyle \bullet$} [c] at 14 12
\put {$ \scriptstyle \bullet$} [c] at 12 12
\put {$ \scriptstyle \bullet$} [c] at 16 12
\setlinear \plot 14 0  14 12 /
\setlinear \plot 12 12  14 8  16 12 /
\setlinear \plot 10 12  14 4 /
\put{$840$} [c] at 13 -2
\endpicture
\end{minipage}
\begin{minipage}{4cm}
\beginpicture
\setcoordinatesystem units   <1.5mm,2mm>
\setplotarea x from 0 to 16, y from -2 to 15
\put{861)} [l] at 2 12
\put {$ \scriptstyle \bullet$} [c] at 10 0
\put {$ \scriptstyle \bullet$} [c] at 14 0
\put {$ \scriptstyle \bullet$} [c] at 14 4
\put {$ \scriptstyle \bullet$} [c] at 14 8
\put {$ \scriptstyle \bullet$} [c] at 14 12
\put {$ \scriptstyle \bullet$} [c] at 12 0
\put {$ \scriptstyle \bullet$} [c] at 16 0
\setlinear \plot 14 0  14 12 /
\setlinear \plot 12 0  14 4  16 0 /
\setlinear \plot 10 0  14 8 /
\put{$840$} [c] at 13 -2
\endpicture
\end{minipage}
\begin{minipage}{4cm}
\beginpicture
\setcoordinatesystem units   <1.5mm,2mm>
\setplotarea x from 0 to 16, y from -2 to 15
\put{862)} [l] at 2 12
\put {$ \scriptstyle \bullet$} [c] at 13 0
\put {$ \scriptstyle \bullet$} [c] at 10 4
\put {$ \scriptstyle \bullet$} [c] at 12 4
\put {$ \scriptstyle \bullet$} [c] at 14 4
\put {$ \scriptstyle \bullet$} [c] at 16 4
\put {$ \scriptstyle \bullet$} [c] at 10 12
\put {$ \scriptstyle \bullet$} [c] at 16 12
\setlinear \plot 13 0 10 4 10 12 12 4 13 0 14 4 16 12 16 4 13 0  /
\setlinear \plot  10 4 16 12 12 4 /
\setlinear \plot 16 4 10 12 14 4  /
\put{$105$} [c] at 13 -2
\endpicture
\end{minipage}
\begin{minipage}{4cm}
\beginpicture
\setcoordinatesystem units   <1.5mm,2mm>
\setplotarea x from 0 to 16, y from -2 to 15
\put{863)} [l] at 2 12
\put {$ \scriptstyle \bullet$} [c] at 13 12
\put {$ \scriptstyle \bullet$} [c] at 10 8
\put {$ \scriptstyle \bullet$} [c] at 12 8
\put {$ \scriptstyle \bullet$} [c] at 14 8
\put {$ \scriptstyle \bullet$} [c] at 16 8
\put {$ \scriptstyle \bullet$} [c] at 10 0
\put {$ \scriptstyle \bullet$} [c] at 16 0
\setlinear \plot 13 12 10 8 10 0 12 8 13 12 14 8 16 0 16 8 13 12  /
\setlinear \plot  10 8 16 0 12 8 /
\setlinear \plot 16 8 10 0 14 8  /
\put{$105$} [c] at 13 -2
\endpicture
\end{minipage}
\begin{minipage}{4cm}
\beginpicture
\setcoordinatesystem units   <1.5mm,2mm>
\setplotarea x from 0 to 16, y from -2 to 15
\put{864)} [l] at 2 12
\put {$ \scriptstyle \bullet$} [c] at 10 6
\put {$ \scriptstyle \bullet$} [c] at 10  12
\put {$ \scriptstyle \bullet$} [c] at 12 12
\put {$ \scriptstyle \bullet$} [c] at 14 12
\put {$ \scriptstyle \bullet$} [c] at 16 12
\put {$ \scriptstyle \bullet$} [c] at 16 6
\put {$ \scriptstyle \bullet$} [c] at 13 0
\setlinear \plot 16 6 13 0 10 6 10 12  16 6  14 12 10 6 12 12 16 6 16 12 10 6 /
\put{$105$} [c] at 13 -2
\endpicture
\end{minipage}
$$

$$
\begin{minipage}{4cm}
\beginpicture
\setcoordinatesystem units   <1.5mm,2mm>
\setplotarea x from 0 to 16, y from -2 to 15
\put{865)} [l] at 2 12
\put {$ \scriptstyle \bullet$} [c] at 10 6
\put {$ \scriptstyle \bullet$} [c] at 10  0
\put {$ \scriptstyle \bullet$} [c] at 12 0
\put {$ \scriptstyle \bullet$} [c] at 14 0
\put {$ \scriptstyle \bullet$} [c] at 16 0
\put {$ \scriptstyle \bullet$} [c] at 16 6
\put {$ \scriptstyle \bullet$} [c] at 13 12
\setlinear \plot 16 6 13 12 10 6 10 0  16 6  14 0 10 6 12 0 16 6 16 0 10 6 /
\put{$105$} [c] at 13 -2
\endpicture
\end{minipage}
\begin{minipage}{4cm}
\beginpicture
\setcoordinatesystem units   <1.5mm,2mm>
\setplotarea x from 0 to 16, y from -2 to 15
\put {866)} [l] at  2 12
\put {$ \scriptstyle \bullet$} [c] at  10 0
\put {$ \scriptstyle \bullet$} [c] at  16 0
\put {$ \scriptstyle \bullet$} [c] at  10 12
\put {$ \scriptstyle \bullet$} [c] at  12 12
\put {$ \scriptstyle \bullet$} [c] at  14 12
\put {$ \scriptstyle \bullet$} [c] at  16 12
\put {$ \scriptstyle \bullet$} [c] at  13 6
\setlinear \plot   10 0 13 6 16 12  /
\setlinear \plot  16 0 13 6 10 12 /
\setlinear \plot  12 12 13 6 14 12 /
\put{$105$} [c] at 13 -2
\endpicture
\end{minipage}
\begin{minipage}{4cm}
\beginpicture
\setcoordinatesystem units   <1.5mm,2mm>
\setplotarea x from 0 to 16, y from -2 to 15
\put {867)} [l] at  2 12
\put {$ \scriptstyle \bullet$} [c] at  10 12
\put {$ \scriptstyle \bullet$} [c] at  16 12
\put {$ \scriptstyle \bullet$} [c] at  10 0
\put {$ \scriptstyle \bullet$} [c] at  12 0
\put {$ \scriptstyle \bullet$} [c] at  14 0
\put {$ \scriptstyle \bullet$} [c] at  16 0
\put {$ \scriptstyle \bullet$} [c] at  13 6
\setlinear \plot   10 12 13 6 16 0  /
\setlinear \plot  16 12 13 6 10 0 /
\setlinear \plot  12 0 13 6 14 0 /
\put{$105$} [c] at 13 -2
\endpicture
\end{minipage}
\begin{minipage}{4cm}
\beginpicture
\setcoordinatesystem units   <1.5mm,2mm>
\setplotarea x from 0 to 16, y from -2 to 15
\put{868)} [l] at 2 12
\put {$ \scriptstyle \bullet$} [c] at 10 12
\put {$ \scriptstyle \bullet$} [c] at 12 0
\put {$ \scriptstyle \bullet$} [c] at 12 4
\put {$ \scriptstyle \bullet$} [c] at 10.8 9
\put {$ \scriptstyle \bullet$} [c] at 11.5 6.2
\put {$ \scriptstyle \bullet$} [c] at 14 12
\put{$\scriptstyle \bullet$} [c] at 16  0
\setlinear \plot    10 12 12 4 12 0  /
\setlinear \plot    12 4 14 12 /
\put{$5{,}040$} [c] at 13 -2
\endpicture
\end{minipage}
\begin{minipage}{4cm}
\beginpicture
\setcoordinatesystem units   <1.5mm,2mm>
\setplotarea x from 0 to 16, y from -2 to 15
\put{869)} [l] at 2 12
\put {$ \scriptstyle \bullet$} [c] at 10 0
\put {$ \scriptstyle \bullet$} [c] at 12 12
\put {$ \scriptstyle \bullet$} [c] at 12 8
\put {$ \scriptstyle \bullet$} [c] at 10.8 3
\put {$ \scriptstyle \bullet$} [c] at 11.5 5.8
\put {$ \scriptstyle \bullet$} [c] at 14 0
\put{$\scriptstyle \bullet$} [c] at 16  0
\setlinear \plot    10 0 12 8 12 12  /
\setlinear \plot    12 8 14 0 /
\put{$5{,}040$} [c] at 13 -2
\endpicture
\end{minipage}
\begin{minipage}{4cm}
\beginpicture
\setcoordinatesystem units   <1.5mm,2mm>
\setplotarea x from 0 to 16, y from -2 to 15
\put{870)} [l] at 2 12
\put {$ \scriptstyle \bullet$} [c] at 10 8
\put {$ \scriptstyle \bullet$} [c] at 10 12
\put {$ \scriptstyle \bullet$} [c] at 12 0
\put {$ \scriptstyle \bullet$} [c] at 12 4
\put {$ \scriptstyle \bullet$} [c] at 12 12
\put {$ \scriptstyle \bullet$} [c] at 14 8
\put{$\scriptstyle \bullet$} [c] at 16  0
\setlinear \plot  12 0 12 4 10 8 12 12 14 8  12 4     /
\setlinear \plot  10 8 10 12 /
\put{$5{,}040$} [c] at 13 -2
 \endpicture
\end{minipage}
$$
$$
\begin{minipage}{4cm}
\beginpicture
\setcoordinatesystem units   <1.5mm,2mm>
\setplotarea x from 0 to 16, y from -2 to 15
\put{871)} [l] at 2 12
\put {$ \scriptstyle \bullet$} [c] at 10 4
\put {$ \scriptstyle \bullet$} [c] at 10 0
\put {$ \scriptstyle \bullet$} [c] at 12 12
\put {$ \scriptstyle \bullet$} [c] at 12 8
\put {$ \scriptstyle \bullet$} [c] at 12 0
\put {$ \scriptstyle \bullet$} [c] at 14 4
\put{$\scriptstyle \bullet$} [c] at 16  0
\setlinear \plot  12 12 12 8 10 4 12 0 14 4  12 8     /
\setlinear \plot  10 4 10 0 /
\put{$5{,}040$} [c] at 13 -2
\endpicture
\end{minipage}
\begin{minipage}{4cm}
\beginpicture
\setcoordinatesystem units   <1.5mm,2mm>
\setplotarea x from 0 to 16, y from -2 to 15
\put{872)} [l] at 2 12
\put {$ \scriptstyle \bullet$} [c] at 10 6
\put {$ \scriptstyle \bullet$} [c] at 12 0
\put {$ \scriptstyle \bullet$} [c] at 12 12
\put {$ \scriptstyle \bullet$} [c] at 14 6
\put {$ \scriptstyle \bullet$} [c] at 11 3
\put {$ \scriptstyle \bullet$} [c] at 13 9
\put{$\scriptstyle \bullet$} [c] at 16  0
\setlinear \plot  12 0  10 6 12 12  14  6 12 0     /
\setlinear \plot  11 3 13 9 /
\put{$5{,}040$} [c] at 13 -2
 \endpicture
\end{minipage}
\begin{minipage}{4cm}
\beginpicture
\setcoordinatesystem units   <1.5mm,2mm>
\setplotarea x from 0 to 16, y from -2 to 15
\put{873)} [l] at 2 12
\put {$ \scriptstyle \bullet$} [c] at 10 4
\put {$ \scriptstyle \bullet$} [c] at 10 6
\put {$ \scriptstyle \bullet$} [c] at 10 8
\put {$ \scriptstyle \bullet$} [c] at 12 0
\put {$ \scriptstyle \bullet$} [c] at 12 12
\put {$ \scriptstyle \bullet$} [c] at 14 6
\put{$\scriptstyle \bullet$} [c] at 16  0
\setlinear \plot    12 0 10 4 10 8 12 12 14 6 12 0  /
\put{$5{,}040$} [c] at 13 -2
 \endpicture
\end{minipage}
\begin{minipage}{4cm}
\beginpicture
\setcoordinatesystem units   <1.5mm,2mm>
\setplotarea x from 0 to 16, y from -2 to 15
\put{874)} [l] at 2 12
\put {$ \scriptstyle \bullet$} [c] at 10 0
\put {$ \scriptstyle \bullet$} [c] at 10 4
\put {$ \scriptstyle \bullet$} [c] at 10 8
\put {$ \scriptstyle \bullet$} [c] at 10 12
\put {$ \scriptstyle \bullet$} [c] at 16 0
\put {$ \scriptstyle \bullet$} [c] at 16 6
\put{$\scriptstyle \bullet$} [c] at 16  12
\setlinear \plot    10 12 10 0   /
\setlinear \plot    16 0 16 12 /
\put{$5{,}040$} [c] at 13 -2
 \endpicture
\end{minipage}
\begin{minipage}{4cm}
\beginpicture
\setcoordinatesystem units   <1.5mm,2mm>
\setplotarea x from 0 to 16, y from -2 to 15
\put{875)} [l] at 2 12
\put {$ \scriptstyle \bullet$} [c] at 10 4
\put {$ \scriptstyle \bullet$} [c] at 10 8
\put {$ \scriptstyle \bullet$} [c] at 10 12
\put {$ \scriptstyle \bullet$} [c] at 12 0
\put {$ \scriptstyle \bullet$} [c] at 14 4
\put {$ \scriptstyle \bullet$} [c] at 14 12
\put{$\scriptstyle \bullet$} [c] at 16  0
\setlinear \plot    10 12 10 4 12 0 14 4 14 12 10 4 /
\setlinear \plot    10 8 14 4 /
\put{$2{,}520$} [c] at 13 -2
 \endpicture
\end{minipage}
\begin{minipage}{4cm}
\beginpicture
\setcoordinatesystem units   <1.5mm,2mm>
\setplotarea x from 0 to 16, y from -2 to 15
\put{876)} [l] at 2 12
\put {$ \scriptstyle \bullet$} [c] at 10 4
\put {$ \scriptstyle \bullet$} [c] at 10 8
\put {$ \scriptstyle \bullet$} [c] at 10 0
\put {$ \scriptstyle \bullet$} [c] at 12 12
\put {$ \scriptstyle \bullet$} [c] at 14 8
\put {$ \scriptstyle \bullet$} [c] at 14 0
\put{$\scriptstyle \bullet$} [c] at 16  0
\setlinear \plot    10 0 10 8 12 12 14 8 14 0 10 8 /
\setlinear \plot    10 4 14 8 /
\put{$2{,}520$} [c] at 13 -2
\endpicture
\end{minipage}
$$

$$
\begin{minipage}{4cm}
\beginpicture
\setcoordinatesystem units   <1.5mm,2mm>
\setplotarea x from 0 to 16, y from -2 to 15
\put{877)} [l] at 2 12
\put {$ \scriptstyle \bullet$} [c] at 10 6
\put {$ \scriptstyle \bullet$} [c] at 10 12
\put {$ \scriptstyle \bullet$} [c] at 14 12
\put {$ \scriptstyle \bullet$} [c] at 14 6
\put {$ \scriptstyle \bullet$} [c] at 12 0
\put {$ \scriptstyle \bullet$} [c] at 11 3
\put{$\scriptstyle \bullet$} [c] at 16  0
\setlinear \plot    12 0  10 6  10 12 14 6 14 12 10 6  /
\setlinear \plot    12 0 14 6 /
\put{$2{,}520$} [c] at 13 -2
\endpicture
\end{minipage}
\begin{minipage}{4cm}
\beginpicture
\setcoordinatesystem units   <1.5mm,2mm>
\setplotarea x from 0 to 16, y from -2 to 15
\put{878)} [l] at 2 12
\put {$ \scriptstyle \bullet$} [c] at 10 6
\put {$ \scriptstyle \bullet$} [c] at 10 0
\put {$ \scriptstyle \bullet$} [c] at 14 0
\put {$ \scriptstyle \bullet$} [c] at 14 6
\put {$ \scriptstyle \bullet$} [c] at 12 12
\put {$ \scriptstyle \bullet$} [c] at 11 9
\put{$\scriptstyle \bullet$} [c] at 16  0
\setlinear \plot    12 12  10 6 10 0 14 6 14 0 10 6  /
\setlinear \plot    12 12 14 6 /
\put{$2{,}520$} [c] at 13 -2
\endpicture
\end{minipage}
\begin{minipage}{4cm}
\beginpicture
\setcoordinatesystem units   <1.5mm,2mm>
\setplotarea x from 0 to 16, y from -2 to 15
\put{879)} [l] at 2 12
\put {$ \scriptstyle \bullet$} [c] at 10 0
\put {$ \scriptstyle \bullet$} [c] at 10 12
\put {$ \scriptstyle \bullet$} [c] at 12 6
\put {$ \scriptstyle \bullet$} [c] at 11 9
\put {$ \scriptstyle \bullet$} [c] at 14 0
\put {$ \scriptstyle \bullet$} [c] at 14 12
\put{$\scriptstyle \bullet$} [c] at 16  0
\setlinear \plot    10 0 14  12 /
\setlinear \plot    10 12 14 0 /
\put{$2{,}520$} [c] at 13 -2
 \endpicture
\end{minipage}
\begin{minipage}{4cm}
\beginpicture
\setcoordinatesystem units   <1.5mm,2mm>
\setplotarea x from 0 to 16, y from -2 to 15
\put{880)} [l] at 2 12
\put {$ \scriptstyle \bullet$} [c] at 10 0
\put {$ \scriptstyle \bullet$} [c] at 10 12
\put {$ \scriptstyle \bullet$} [c] at 12 6
\put {$ \scriptstyle \bullet$} [c] at 11 3
\put {$ \scriptstyle \bullet$} [c] at 14 0
\put {$ \scriptstyle \bullet$} [c] at 14 12
\put{$\scriptstyle \bullet$} [c] at 16  0
\setlinear \plot    10 0 14  12 /
\setlinear \plot    10 12 14 0 /
\put{$2{,}520$} [c] at 13 -2
\endpicture
\end{minipage}
\begin{minipage}{4cm}
\beginpicture
\setcoordinatesystem units   <1.5mm,2mm>
\setplotarea x from 0 to 16, y from -2 to 15
\put{881)} [l] at 2 12
\put {$ \scriptstyle \bullet$} [c] at 10 0
\put {$ \scriptstyle \bullet$} [c] at 10 6
\put {$ \scriptstyle \bullet$} [c] at 10 12
\put {$ \scriptstyle \bullet$} [c] at 14 0
\put {$ \scriptstyle \bullet$} [c] at 14 6
\put {$ \scriptstyle \bullet$} [c] at 14 12
\put{$\scriptstyle \bullet$} [c] at 16  0
\setlinear \plot    10 12 10 0 14 6  10 12  /
\setlinear \plot    14 12 14 0 10 6 14 12 /
\put{$630$} [c] at 13 -2
 \endpicture
\end{minipage}
\begin{minipage}{4cm}
\beginpicture
\setcoordinatesystem units   <1.5mm,2mm>
\setplotarea x from 0 to 16, y from -2 to 15
\put{${\bf  21}$} [l] at 2 15
\put{882)} [l] at 2 12
\put {$ \scriptstyle \bullet$} [c] at 10 4
\put {$ \scriptstyle \bullet$} [c] at 10 12
\put {$ \scriptstyle \bullet$} [c] at 13 0
\put {$ \scriptstyle \bullet$} [c] at 13 4
\put {$ \scriptstyle \bullet$} [c] at 13 12
\put {$ \scriptstyle \bullet$} [c] at 14.5 8
\put {$ \scriptstyle \bullet$} [c] at 16 4
\setlinear \plot 10 12 10 4 13 0 16 4 13 12 10 4 /
\setlinear \plot 13 0 13 12 /
\put{$5{,}040$} [c] at 13 -2
\endpicture
\end{minipage}
$$
$$
\begin{minipage}{4cm}
\beginpicture
\setcoordinatesystem units   <1.5mm,2mm>
\setplotarea x from 0 to 16, y from -2 to 15
\put{883)} [l] at 2 12
\put {$ \scriptstyle \bullet$} [c] at 10 8
\put {$ \scriptstyle \bullet$} [c] at 10 0
\put {$ \scriptstyle \bullet$} [c] at 13 0
\put {$ \scriptstyle \bullet$} [c] at 13 8
\put {$ \scriptstyle \bullet$} [c] at 13 12
\put {$ \scriptstyle \bullet$} [c] at 14.5 4
\put {$ \scriptstyle \bullet$} [c] at 16 8
\setlinear \plot 10 0 10 8 13 12 16 8 13 0 10 8 /
\setlinear \plot 13 0 13 12 /
\put{$5{,}040$} [c] at 13 -2
\endpicture
\end{minipage}
\begin{minipage}{4cm}
\beginpicture
\setcoordinatesystem units   <1.5mm,2mm>
\setplotarea x from 0 to 16, y from -2 to 15
\put{884)} [l] at 2 12
\put {$ \scriptstyle \bullet$} [c] at 10 4
\put {$ \scriptstyle \bullet$} [c] at 10 12
\put {$ \scriptstyle \bullet$} [c] at 13 0
\put {$ \scriptstyle \bullet$} [c] at 13 12
\put {$ \scriptstyle \bullet$} [c] at 16 4
\put {$ \scriptstyle \bullet$} [c] at 16 8
\put {$ \scriptstyle \bullet$} [c] at 16 12
\setlinear \plot 10 12 10 4 13 0 16 4 16 12 /
\setlinear \plot 10 12 16 4 13 12 /
\put{$5{,}040$} [c] at 13 -2
\endpicture
\end{minipage}
\begin{minipage}{4cm}
\beginpicture
\setcoordinatesystem units   <1.5mm,2mm>
\setplotarea x from 0 to 16, y from -2 to 15
\put{885)} [l] at 2 12
\put {$ \scriptstyle \bullet$} [c] at 10 0
\put {$ \scriptstyle \bullet$} [c] at 10 8
\put {$ \scriptstyle \bullet$} [c] at 13 0
\put {$ \scriptstyle \bullet$} [c] at 13 12
\put {$ \scriptstyle \bullet$} [c] at 16 4
\put {$ \scriptstyle \bullet$} [c] at 16 8
\put {$ \scriptstyle \bullet$} [c] at 16 0
\setlinear \plot 10 0 10 8 13 12 16 8 16 0 /
\setlinear \plot 10 0 16 8 13 0 /
\put{$5{,}040$} [c] at 13 -2
\endpicture
\end{minipage}
\begin{minipage}{4cm}
\beginpicture
\setcoordinatesystem units   <1.5mm,2mm>
\setplotarea x from 0 to 16, y from -2 to 15
\put{886)} [l] at 2 12
\put {$ \scriptstyle \bullet$} [c] at 10 6
\put {$ \scriptstyle \bullet$} [c] at 12 6
\put {$ \scriptstyle \bullet$} [c] at 12 12
\put {$ \scriptstyle \bullet$} [c] at 14 6
\put {$ \scriptstyle \bullet$} [c] at 14 12
\put {$ \scriptstyle \bullet$} [c] at 16 12
\put {$ \scriptstyle \bullet$} [c] at 12 0
\setlinear \plot 12 12 10 6 12 0 12 6 12 12 14 6 14 12 12 6  /
\setlinear \plot  16 12 14 6 12 0  /
\put{$5{,}040$} [c] at 13 -2
\endpicture
\end{minipage}
\begin{minipage}{4cm}
\beginpicture
\setcoordinatesystem units   <1.5mm,2mm>
\setplotarea x from 0 to 16, y from -2 to 15
\put{887)} [l] at 2 12
\put {$ \scriptstyle \bullet$} [c] at 10 6
\put {$ \scriptstyle \bullet$} [c] at 12 6
\put {$ \scriptstyle \bullet$} [c] at 12 12
\put {$ \scriptstyle \bullet$} [c] at 14 6
\put {$ \scriptstyle \bullet$} [c] at 14 0
\put {$ \scriptstyle \bullet$} [c] at 16 0
\put {$ \scriptstyle \bullet$} [c] at 12 0
\setlinear \plot 12 0 10 6 12 12 12 6 12 0 14 6 14 0 12 6  /
\setlinear \plot  16 0 14 6 12 12  /
\put{$5{,}040$} [c] at 13 -2
\endpicture
\end{minipage}
\begin{minipage}{4cm}
\beginpicture
\setcoordinatesystem units   <1.5mm,2mm>
\setplotarea x from 0 to 16, y from -2 to 15
\put{888)} [l] at 2 12
\put {$ \scriptstyle \bullet$} [c] at 10 6
\put {$ \scriptstyle \bullet$} [c] at 10 12
\put {$ \scriptstyle \bullet$} [c] at 11.5 12
\put {$ \scriptstyle \bullet$} [c] at 11.5 10
\put {$ \scriptstyle \bullet$} [c] at 11.5 0
\put {$ \scriptstyle \bullet$} [c] at 13 6
\put {$ \scriptstyle \bullet$} [c] at 16 12
\setlinear \plot 16 12 11.5 0 13 6 11.5 10 11.5 12  /
\setlinear \plot 10 12 10 6 11.5 10  /
\setlinear \plot 11.5 0 10 6   /
\put{$5{,}040$} [c] at 13 -2
\endpicture
\end{minipage}
$$
$$
\begin{minipage}{4cm}
\beginpicture
\setcoordinatesystem units   <1.5mm,2mm>
\setplotarea x from 0 to 16, y from -2 to 15
\put{889)} [l] at 2 12
\put {$ \scriptstyle \bullet$} [c] at 10 6
\put {$ \scriptstyle \bullet$} [c] at 10 0
\put {$ \scriptstyle \bullet$} [c] at 11.5 12
\put {$ \scriptstyle \bullet$} [c] at 11.5 2
\put {$ \scriptstyle \bullet$} [c] at 11.5 0
\put {$ \scriptstyle \bullet$} [c] at 13 6
\put {$ \scriptstyle \bullet$} [c] at 16 0
\setlinear \plot 16 0 11.5 12 13 6 11.5 2 11.5 0  /
\setlinear \plot 10 0 10 6 11.5 2  /
\setlinear \plot 11.5 12 10 6   /
\put{$5{,}040$} [c] at 13 -2
\endpicture
\end{minipage}
\begin{minipage}{4cm}
\beginpicture
\setcoordinatesystem units   <1.5mm,2mm>
\setplotarea x from 0 to 16, y from -2 to 15
\put{890)} [l] at 2 12
\put {$ \scriptstyle \bullet$} [c] at 10 4
\put {$ \scriptstyle \bullet$} [c] at 10 12
\put {$ \scriptstyle \bullet$} [c] at 12 0
\put {$ \scriptstyle \bullet$} [c] at 12 4
\put {$ \scriptstyle \bullet$} [c] at 12  8
\put {$ \scriptstyle \bullet$} [c] at 14 12
\put {$ \scriptstyle \bullet$} [c] at 16 12
\setlinear \plot 16 12 12 0 12 8 14 12  /
\setlinear \plot 12 8 10 12 10 4 12 0 /
\put{$5{,}040$} [c] at 13 -2
\endpicture
\end{minipage}
\begin{minipage}{4cm}
\beginpicture
\setcoordinatesystem units   <1.5mm,2mm>
\setplotarea x from 0 to 16, y from -2 to 15
\put{891)} [l] at 2 12
\put {$ \scriptstyle \bullet$} [c] at 10 8
\put {$ \scriptstyle \bullet$} [c] at 10 0
\put {$ \scriptstyle \bullet$} [c] at 12 12
\put {$ \scriptstyle \bullet$} [c] at 12 4
\put {$ \scriptstyle \bullet$} [c] at 12  8
\put {$ \scriptstyle \bullet$} [c] at 14 0
\put {$ \scriptstyle \bullet$} [c] at 16 0
\setlinear \plot 16 0 12 12 12 4 14 0  /
\setlinear \plot 12 4 10 0 10 8 12 12 /
\put{$5{,}040$} [c] at 13 -2
\endpicture
\end{minipage}
\begin{minipage}{4cm}
\beginpicture
\setcoordinatesystem units   <1.5mm,2mm>
\setplotarea x from 0 to 16, y from -2 to 15
\put{892)} [l] at 2 12
\put {$ \scriptstyle \bullet$} [c] at 10 4
\put {$ \scriptstyle \bullet$} [c] at 10 8
\put {$ \scriptstyle \bullet$} [c] at 10 12
\put {$ \scriptstyle \bullet$} [c] at 13 0
\put {$ \scriptstyle \bullet$} [c] at 13 4
\put {$ \scriptstyle \bullet$} [c] at 13 12
\put {$ \scriptstyle \bullet$} [c] at 16 12
\setlinear \plot 16 12 13 0 13 12 10 4 10 12 13 4 /
\setlinear \plot 10 4 13 0  /
\put{$5{,}040$} [c] at 13 -2
\endpicture
\end{minipage}
\begin{minipage}{4cm}
\beginpicture
\setcoordinatesystem units   <1.5mm,2mm>
\setplotarea x from 0 to 16, y from -2 to 15
\put{893)} [l] at 2 12
\put {$ \scriptstyle \bullet$} [c] at 10 4
\put {$ \scriptstyle \bullet$} [c] at 10 8
\put {$ \scriptstyle \bullet$} [c] at 10 0
\put {$ \scriptstyle \bullet$} [c] at 13 0
\put {$ \scriptstyle \bullet$} [c] at 13 8
\put {$ \scriptstyle \bullet$} [c] at 13 12
\put {$ \scriptstyle \bullet$} [c] at 16 0
\setlinear \plot 16 0 13 12 13 0 10 8 10 0 13 8 /
\setlinear \plot 10 8 13 12  /
\put{$5{,}040$} [c] at 13 -2
\endpicture
\end{minipage}
\begin{minipage}{4cm}
\beginpicture
\setcoordinatesystem units   <1.5mm,2mm>
\setplotarea x from 0 to 16, y from -2 to 15
\put{894)} [l] at 2 12
\put {$ \scriptstyle \bullet$} [c] at  10 6
\put {$ \scriptstyle \bullet$} [c] at  10 9
\put {$ \scriptstyle \bullet$} [c] at  10 12
\put {$ \scriptstyle \bullet$} [c] at  11.5 0
\put {$ \scriptstyle \bullet$} [c] at  11.5 12
\put {$ \scriptstyle \bullet$} [c] at  13 6
\put {$ \scriptstyle \bullet$} [c] at  16 0
\setlinear \plot 16 0 11.5 12 13 6 11.5 0 10 6 10 12        /
\setlinear \plot  10 6  11.5 12 /
\put{$5{,}040   $} [c] at 13 -2
\endpicture
\end{minipage}
$$
$$
\begin{minipage}{4cm}
\beginpicture
\setcoordinatesystem units   <1.5mm,2mm>
\setplotarea x from 0 to 16, y from -2 to 15
\put{895)} [l] at 2 12
\put {$ \scriptstyle \bullet$} [c] at  10 6
\put {$ \scriptstyle \bullet$} [c] at  10 3
\put {$ \scriptstyle \bullet$} [c] at  10 0
\put {$ \scriptstyle \bullet$} [c] at  11.5 0
\put {$ \scriptstyle \bullet$} [c] at  11.5 12
\put {$ \scriptstyle \bullet$} [c] at  13 6
\put {$ \scriptstyle \bullet$} [c] at  16 12
\setlinear \plot 16 12 11.5 0 13 6 11.5 12 10 6 10 0        /
\setlinear \plot  10 6  11.5 0 /
\put{$5{,}040   $} [c] at 13 -2
\endpicture
\end{minipage}
\begin{minipage}{4cm}
\beginpicture
\setcoordinatesystem units   <1.5mm,2mm>
\setplotarea x from 0 to 16, y from -2 to 15
\put{896)} [l] at 2 12
\put {$ \scriptstyle \bullet$} [c] at  10 0
\put {$ \scriptstyle \bullet$} [c] at  13 6
\put {$ \scriptstyle \bullet$} [c] at  13 9
\put {$ \scriptstyle \bullet$} [c] at  13 12
\put {$ \scriptstyle \bullet$} [c] at  14.5 0
\put {$ \scriptstyle \bullet$} [c] at  14.5 12
\put {$ \scriptstyle \bullet$} [c] at  16 6
\setlinear \plot 10 0 13 12 13 6 14.5 12 16 6 14.5 0 13 6        /
\put{$5{,}040   $} [c] at 13 -2
\endpicture
\end{minipage}
\begin{minipage}{4cm}
\beginpicture
\setcoordinatesystem units   <1.5mm,2mm>
\setplotarea x from 0 to 16, y from -2 to 15
\put{897)} [l] at 2 12
\put {$ \scriptstyle \bullet$} [c] at  10 12
\put {$ \scriptstyle \bullet$} [c] at  13 6
\put {$ \scriptstyle \bullet$} [c] at  13 3
\put {$ \scriptstyle \bullet$} [c] at  13 0
\put {$ \scriptstyle \bullet$} [c] at  14.5 0
\put {$ \scriptstyle \bullet$} [c] at  14.5 12
\put {$ \scriptstyle \bullet$} [c] at  16 6
\setlinear \plot 10 12 13 0 13 6 14.5 0 16 6 14.5 12 13 6        /
\put{$5{,}040   $} [c] at 13 -2
\endpicture
\end{minipage}
\begin{minipage}{4cm}
\beginpicture
\setcoordinatesystem units   <1.5mm,2mm>
\setplotarea x from 0 to 16, y from -2 to 15
\put{898)} [l] at 2 12
\put {$ \scriptstyle \bullet$} [c] at  10 0
\put {$ \scriptstyle \bullet$} [c] at  10 12
\put {$ \scriptstyle \bullet$} [c] at  12 6
\put {$ \scriptstyle \bullet$} [c] at  13 9
\put {$ \scriptstyle \bullet$} [c] at  14 0
\put {$ \scriptstyle \bullet$} [c] at  14 12
\put {$ \scriptstyle \bullet$} [c] at  16 6
\setlinear \plot 10 0 10 12 12 6 14 12 16 6 14 0 12 6        /
\put{$5{,}040   $} [c] at 13 -2
\endpicture
\end{minipage}
\begin{minipage}{4cm}
\beginpicture
\setcoordinatesystem units   <1.5mm,2mm>
\setplotarea x from 0 to 16, y from -2 to 15
\put{899)} [l] at 2 12
\put {$ \scriptstyle \bullet$} [c] at  10 0
\put {$ \scriptstyle \bullet$} [c] at  10 12
\put {$ \scriptstyle \bullet$} [c] at  12 6
\put {$ \scriptstyle \bullet$} [c] at  13 3
\put {$ \scriptstyle \bullet$} [c] at  14 0
\put {$ \scriptstyle \bullet$} [c] at  14 12
\put {$ \scriptstyle \bullet$} [c] at  16 6
\setlinear \plot 10 12 10 0 12 6 14 0 16 6 14 12 12 6        /
\put{$5{,}040   $} [c] at 13 -2
\endpicture
\end{minipage}
\begin{minipage}{4cm}
\beginpicture
\setcoordinatesystem units   <1.5mm,2mm>
\setplotarea x from 0 to 16, y from -2 to 15
\put{900)} [l] at 2 12
\put {$ \scriptstyle \bullet$} [c] at  10 6
\put {$ \scriptstyle \bullet$} [c] at  12 0
\put {$ \scriptstyle \bullet$} [c] at  12 12
\put {$ \scriptstyle \bullet$} [c] at  14 6
\put {$ \scriptstyle \bullet$} [c] at  16 6
\put {$ \scriptstyle \bullet$} [c] at  16 0
\put {$ \scriptstyle \bullet$} [c] at  16 12
\setlinear \plot 14 6 12 0 10 6 12 12 14 6 16 12 16 0       /
\put{$5{,}040  $} [c] at 13 -2
\endpicture
\end{minipage}
$$
$$
\begin{minipage}{4cm}
\beginpicture
\setcoordinatesystem units   <1.5mm,2mm>
\setplotarea x from 0 to 16, y from -2 to 15
\put{901)} [l] at 2 12
\put {$ \scriptstyle \bullet$} [c] at  10 6
\put {$ \scriptstyle \bullet$} [c] at  12 0
\put {$ \scriptstyle \bullet$} [c] at  12 12
\put {$ \scriptstyle \bullet$} [c] at  14 6
\put {$ \scriptstyle \bullet$} [c] at  16 6
\put {$ \scriptstyle \bullet$} [c] at  16 0
\put {$ \scriptstyle \bullet$} [c] at  16 12
\setlinear \plot 14 6 12 0 10 6 12 12 14 6 16 0 16 12       /
\put{$5{,}040  $} [c] at 13 -2
\endpicture
\end{minipage}
\begin{minipage}{4cm}
\beginpicture
\setcoordinatesystem units   <1.5mm,2mm>
\setplotarea x from 0 to 16, y from -2 to 15
\put{902)} [l] at  2 12
\put {$ \scriptstyle \bullet$} [c] at  10 6
\put {$ \scriptstyle \bullet$} [c] at  11 0
\put {$ \scriptstyle \bullet$} [c] at  11 12
\put {$ \scriptstyle \bullet$} [c] at  12 6
\put {$ \scriptstyle \bullet$} [c] at  14 6
\put {$ \scriptstyle \bullet$} [c] at  16  0
\put {$ \scriptstyle \bullet$} [c] at  16  12
\setlinear \plot 12 6 16 0 16 12 14 6 11 0 10 6 11 12 12 6 16 0    /
\setlinear \plot 11 0 12 6  /
\put{$5{,}040 $} [c] at 13 -2
\endpicture
\end{minipage}
\begin{minipage}{4cm}
\beginpicture
\setcoordinatesystem units   <1.5mm,2mm>
\setplotarea x from 0 to 16, y from -2 to 15
\put{903)} [l] at 2 12
\put {$ \scriptstyle \bullet$} [c] at  10 6
\put {$ \scriptstyle \bullet$} [c] at  11 0
\put {$ \scriptstyle \bullet$} [c] at  11 12
\put {$ \scriptstyle \bullet$} [c] at  12 6
\put {$ \scriptstyle \bullet$} [c] at  14 6
\put {$ \scriptstyle \bullet$} [c] at  16  0
\put {$ \scriptstyle \bullet$} [c] at  16  12
\setlinear \plot 12 6 16 12 16 0 14 6 11 12 10 6 11 0 12 6 16 12    /
\setlinear \plot 11 12 12 6  /
\put{$5{,}040$} [c] at 13 -2
\endpicture
\end{minipage}
\begin{minipage}{4cm}
\beginpicture
\setcoordinatesystem units   <1.5mm,2mm>
\setplotarea x from 0 to 16, y from -2 to 15
\put{904)} [l] at 2 12
\put {$ \scriptstyle \bullet$} [c] at  10 12
\put {$ \scriptstyle \bullet$} [c] at  13 0
\put {$ \scriptstyle \bullet$} [c] at  13  4
\put {$ \scriptstyle \bullet$} [c] at  13 8
\put {$ \scriptstyle \bullet$} [c] at  13 12
\put {$ \scriptstyle \bullet$} [c] at  16 0
\put {$ \scriptstyle \bullet$} [c] at  16 12
\setlinear \plot  13  12 13  0  16 12 16 0 13 12 /
\setlinear \plot  10 12 13 4  /
\put{$5{,}040$} [c] at 13 -2
\endpicture
\end{minipage}
\begin{minipage}{4cm}
\beginpicture
\setcoordinatesystem units   <1.5mm,2mm>
\setplotarea x from 0 to 16, y from -2 to 15
\put{905)} [l] at 2 12
\put {$ \scriptstyle \bullet$} [c] at  10 0
\put {$ \scriptstyle \bullet$} [c] at  13 0
\put {$ \scriptstyle \bullet$} [c] at  13  4
\put {$ \scriptstyle \bullet$} [c] at  13 8
\put {$ \scriptstyle \bullet$} [c] at  13 12
\put {$ \scriptstyle \bullet$} [c] at  16 0
\put {$ \scriptstyle \bullet$} [c] at  16 12
\setlinear \plot  13  12 13  0  16 12 16 0 13 12 /
\setlinear \plot  10 0 13 8  /
\put{$5{,}040$} [c] at 13 -2
\endpicture
\end{minipage}
\begin{minipage}{4cm}
\beginpicture
\setcoordinatesystem units   <1.5mm,2mm>
\setplotarea x from 0 to 16, y from -2 to 15
\put{906)} [l] at 2 12
\put {$ \scriptstyle \bullet$} [c] at  10 12
\put {$ \scriptstyle \bullet$} [c] at  10 6
\put {$ \scriptstyle \bullet$} [c] at  11 0
\put {$ \scriptstyle \bullet$} [c] at  11 12
\put {$ \scriptstyle \bullet$} [c] at  12 6
\put {$ \scriptstyle \bullet$} [c] at  16 0
\put {$ \scriptstyle \bullet$} [c] at  16 12
\setlinear \plot  10 12 10 6 11 0 16 12 16 0 12 6 11 12 10 6    /
\setlinear \plot  11 0  12 6   /
\put{$5{,}040$} [c] at 13 -2
\endpicture
\end{minipage}
$$
$$
\begin{minipage}{4cm}
\beginpicture
\setcoordinatesystem units   <1.5mm,2mm>
\setplotarea x from 0 to 16, y from -2 to 15
\put{907)} [l] at 2 12
\put {$ \scriptstyle \bullet$} [c] at  10 0
\put {$ \scriptstyle \bullet$} [c] at  10 6
\put {$ \scriptstyle \bullet$} [c] at  11 0
\put {$ \scriptstyle \bullet$} [c] at  11 12
\put {$ \scriptstyle \bullet$} [c] at  12 6
\put {$ \scriptstyle \bullet$} [c] at  16 0
\put {$ \scriptstyle \bullet$} [c] at  16 12
\setlinear \plot  10 0 10 6 11 12 16 0 16 12 12 6 11 0 10 6    /
\setlinear \plot  11 12  12 6   /
\put{$5{,}040$} [c] at 13 -2
\endpicture
\end{minipage}
\begin{minipage}{4cm}
\beginpicture
\setcoordinatesystem units   <1.5mm,2mm>
\setplotarea x from 0 to 16, y from -2 to 15
\put{908)} [l] at 2 12
\put {$ \scriptstyle \bullet$} [c] at  10 0
\put {$ \scriptstyle \bullet$} [c] at  10 4
\put {$ \scriptstyle \bullet$} [c] at  10 8
\put {$ \scriptstyle \bullet$} [c] at  10 12
\put {$ \scriptstyle \bullet$} [c] at  13 12
\put {$ \scriptstyle \bullet$} [c] at  16 0
\put {$ \scriptstyle \bullet$} [c] at  16 12
\setlinear \plot  16 12 16 0 13 12 10 4 10 12 /
\setlinear \plot  10 4 10 0 /
\put{$5{,}040$} [c] at 13 -2
\endpicture
\end{minipage}
\begin{minipage}{4cm}
\beginpicture
\setcoordinatesystem units   <1.5mm,2mm>
\setplotarea x from 0 to 16, y from -2 to 15
\put{909)} [l] at 2 12
\put {$ \scriptstyle \bullet$} [c] at  10 0
\put {$ \scriptstyle \bullet$} [c] at  10 4
\put {$ \scriptstyle \bullet$} [c] at  10 8
\put {$ \scriptstyle \bullet$} [c] at  10 12
\put {$ \scriptstyle \bullet$} [c] at  13 0
\put {$ \scriptstyle \bullet$} [c] at  16 0
\put {$ \scriptstyle \bullet$} [c] at  16 12
\setlinear \plot  16 0 16 12 13 0 10 8 10 0 /
\setlinear \plot  10 8 10 12 /
\put{$5{,}040$} [c] at 13 -2
\endpicture
\end{minipage}
\begin{minipage}{4cm}
\beginpicture
\setcoordinatesystem units   <1.5mm,2mm>
\setplotarea x from 0 to 16, y from -2 to 15
\put{910)} [l] at 2 12
\put {$ \scriptstyle \bullet$} [c] at  10 0
\put {$ \scriptstyle \bullet$} [c] at  10 12
\put {$ \scriptstyle \bullet$} [c] at  12 12
\put {$ \scriptstyle \bullet$} [c] at  12 6
\put {$ \scriptstyle \bullet$} [c] at  14 12
\put {$ \scriptstyle \bullet$} [c] at  14 0
\put {$ \scriptstyle \bullet$} [c] at  16 6
\setlinear \plot  10 12 10 0 12 6 12 12 /
\setlinear \plot  12 6 14 12 16 6 14 0 12 6 /
\put{$5{,}040$} [c] at 13 -2
\endpicture
\end{minipage}
\begin{minipage}{4cm}
\beginpicture
\setcoordinatesystem units   <1.5mm,2mm>
\setplotarea x from 0 to 16, y from -2 to 15
\put{911)} [l] at 2 12
\put {$ \scriptstyle \bullet$} [c] at  10 0
\put {$ \scriptstyle \bullet$} [c] at  10 12
\put {$ \scriptstyle \bullet$} [c] at  12 0
\put {$ \scriptstyle \bullet$} [c] at  12 6
\put {$ \scriptstyle \bullet$} [c] at  14 12
\put {$ \scriptstyle \bullet$} [c] at  14 0
\put {$ \scriptstyle \bullet$} [c] at  16 6
\setlinear \plot  10 0 10 12 12 6 12 0 /
\setlinear \plot  12 6 14 12 16 6 14 0 12 6 /
\put{$5{,}040$} [c] at 13 -2
\endpicture
\end{minipage}
\begin{minipage}{4cm}
\beginpicture
\setcoordinatesystem units   <1.5mm,2mm>
\setplotarea x from 0 to 16, y from -2 to 15
\put{912)} [l] at 2 12
\put {$ \scriptstyle \bullet$} [c] at  16 0
\put {$ \scriptstyle \bullet$} [c] at  16 12
\put {$ \scriptstyle \bullet$} [c] at  13 12
\put {$ \scriptstyle \bullet$} [c] at  10 0
\put {$ \scriptstyle \bullet$} [c] at  10 4
\put {$ \scriptstyle \bullet$} [c] at  10 8
\put {$ \scriptstyle \bullet$} [c] at  10 12
\setlinear \plot  10 12  10 0 13 12  16 0 10 8  /
\setlinear \plot  16 0  16 12 /
\put{$5{,}040$} [c] at 13 -2
\endpicture
\end{minipage}
$$
$$
\begin{minipage}{4cm}
\beginpicture
\setcoordinatesystem units   <1.5mm,2mm>
\setplotarea x from 0 to 16, y from -2 to 15
\put{913)} [l] at 2 12
\put {$ \scriptstyle \bullet$} [c] at  16 0
\put {$ \scriptstyle \bullet$} [c] at  16 12
\put {$ \scriptstyle \bullet$} [c] at  13 0
\put {$ \scriptstyle \bullet$} [c] at  10 0
\put {$ \scriptstyle \bullet$} [c] at  10 4
\put {$ \scriptstyle \bullet$} [c] at  10 8
\put {$ \scriptstyle \bullet$} [c] at  10 12
\setlinear \plot  10 0  10 12 13 0  16 12 10 4  /
\setlinear \plot  16 0  16 12 /
\put{$5{,}040$} [c] at 13 -2
\endpicture
\end{minipage}
\begin{minipage}{4cm}
\beginpicture
\setcoordinatesystem units   <1.5mm,2mm>
\setplotarea x from 0 to 16, y from -2 to 15
\put{914)} [l] at 2 12
\put {$ \scriptstyle \bullet$} [c] at  10 12
\put {$ \scriptstyle \bullet$} [c] at  13 0
\put {$ \scriptstyle \bullet$} [c] at  13 8
\put {$ \scriptstyle \bullet$} [c] at  13 12
\put {$ \scriptstyle \bullet$} [c] at  16 0
\put {$ \scriptstyle \bullet$} [c] at  16 8
\put {$ \scriptstyle \bullet$} [c] at  16 12
\setlinear \plot 16 0 10 12 13 0 13 12 16 0 16 12 13 8  /
\put{$5{,}040$} [c] at 13 -2
\endpicture
\end{minipage}
\begin{minipage}{4cm}
\beginpicture
\setcoordinatesystem units   <1.5mm,2mm>
\setplotarea x from 0 to 16, y from -2 to 15
\put{915)} [l] at 2 12
\put {$ \scriptstyle \bullet$} [c] at  10 0
\put {$ \scriptstyle \bullet$} [c] at  13 0
\put {$ \scriptstyle \bullet$} [c] at  13 4
\put {$ \scriptstyle \bullet$} [c] at  13 12
\put {$ \scriptstyle \bullet$} [c] at  16 0
\put {$ \scriptstyle \bullet$} [c] at  16 4
\put {$ \scriptstyle \bullet$} [c] at  16 12
\setlinear \plot 16 12 10 0 13 12 13 0 16 12 16 0 13 4  /
\put{$5{,}040$} [c] at 13 -2
\endpicture
\end{minipage}
\begin{minipage}{4cm}
\beginpicture
\setcoordinatesystem units   <1.5mm,2mm>
\setplotarea x from 0 to 16, y from -2 to 15
\put{916)} [l] at 2 12
\put {$ \scriptstyle \bullet$} [c] at  16 12
\put {$ \scriptstyle \bullet$} [c] at  13 12
\put {$ \scriptstyle \bullet$} [c] at  13 6
\put {$ \scriptstyle \bullet$} [c] at  13 0
\put {$ \scriptstyle \bullet$} [c] at  10 0
\put {$ \scriptstyle \bullet$} [c] at  10 6
\put {$ \scriptstyle \bullet$} [c] at  10 12
\setlinear \plot  16 12 13 0 10 12  10 0 13 6 13 12 10 6  /
\setlinear \plot  13 6  13 0  /
\put{$5{,}040$} [c] at 13 -2
\endpicture
\end{minipage}
\begin{minipage}{4cm}
\beginpicture
\setcoordinatesystem units   <1.5mm,2mm>
\setplotarea x from 0 to 16, y from -2 to 15
\put{917)} [l] at 2 12
\put {$ \scriptstyle \bullet$} [c] at  16 0
\put {$ \scriptstyle \bullet$} [c] at  13 12
\put {$ \scriptstyle \bullet$} [c] at  13 6
\put {$ \scriptstyle \bullet$} [c] at  13 0
\put {$ \scriptstyle \bullet$} [c] at  10 0
\put {$ \scriptstyle \bullet$} [c] at  10 6
\put {$ \scriptstyle \bullet$} [c] at  10 12
\setlinear \plot  16 0 13 12 10 0  10 12 13 6 13 0 10 6  /
\setlinear \plot  13 6  13 12  /
\put{$5{,}040$} [c] at 13 -2
\endpicture
\end{minipage}
\begin{minipage}{4cm}
\beginpicture
\setcoordinatesystem units   <1.5mm,2mm>
\setplotarea x from 0 to 16, y from -2 to 15
\put{918)} [l] at 2 12
\put {$ \scriptstyle \bullet$} [c] at  16 12
\put {$ \scriptstyle \bullet$} [c] at  13 12
\put {$ \scriptstyle \bullet$} [c] at  13 6
\put {$ \scriptstyle \bullet$} [c] at  13 0
\put {$ \scriptstyle \bullet$} [c] at  10 0
\put {$ \scriptstyle \bullet$} [c] at  10 6
\put {$ \scriptstyle \bullet$} [c] at  10 12
\setlinear \plot  13 0 13 12 10 0  10 12 13 0  /
\setlinear \plot  16 12 13 6 /
\put{$5{,}040$} [c] at 13 -2
\endpicture
\end{minipage}
$$
$$
\begin{minipage}{4cm}
\beginpicture
\setcoordinatesystem units   <1.5mm,2mm>
\setplotarea x from 0 to 16, y from -2 to 15
\put{919)} [l] at 2 12
\put {$ \scriptstyle \bullet$} [c] at  16 0
\put {$ \scriptstyle \bullet$} [c] at  13 12
\put {$ \scriptstyle \bullet$} [c] at  13 6
\put {$ \scriptstyle \bullet$} [c] at  13 0
\put {$ \scriptstyle \bullet$} [c] at  10 0
\put {$ \scriptstyle \bullet$} [c] at  10 6
\put {$ \scriptstyle \bullet$} [c] at  10 12
\setlinear \plot  13 0 13 12 10 0  10 12 13 0  /
\setlinear \plot  16 0 13 6 /
\put{$5{,}040$} [c] at 13 -2
\endpicture
\end{minipage}
\begin{minipage}{4cm}
\beginpicture
\setcoordinatesystem units   <1.5mm,2mm>
\setplotarea x from 0 to 16, y from -2 to 15
\put{920)} [l] at 2 12
\put {$ \scriptstyle \bullet$} [c] at  10 0
\put {$ \scriptstyle \bullet$} [c] at  10 12
\put {$ \scriptstyle \bullet$} [c] at  14 6
\put {$ \scriptstyle \bullet$} [c] at  15 0
\put {$ \scriptstyle \bullet$} [c] at  15 12
\put {$ \scriptstyle \bullet$} [c] at  16 6
\put {$ \scriptstyle \bullet$} [c] at  16 12
\setlinear \plot 16 12 16 6 15 0  14 6 15 12 16 6 /
\setlinear \plot 14 6 10 12 10 0 15  12  /
\put{$5{,}040$} [c] at 13 -2
\endpicture
\end{minipage}
\begin{minipage}{4cm}
\beginpicture
\setcoordinatesystem units   <1.5mm,2mm>
\setplotarea x from 0 to 16, y from -2 to 15
\put{921)} [l] at 2 12
\put {$ \scriptstyle \bullet$} [c] at  10 0
\put {$ \scriptstyle \bullet$} [c] at  10 12
\put {$ \scriptstyle \bullet$} [c] at  14 6
\put {$ \scriptstyle \bullet$} [c] at  15 0
\put {$ \scriptstyle \bullet$} [c] at  15 12
\put {$ \scriptstyle \bullet$} [c] at  16 6
\put {$ \scriptstyle \bullet$} [c] at  16 0
\setlinear \plot 16 0 16 6 15 12  14 6 15 0 16 6 /
\setlinear \plot 14 6 10  0 10 12 15  0  /
\put{$5{,}040$} [c] at 13 -2
\endpicture
\end{minipage}
\begin{minipage}{4cm}
\beginpicture
\setcoordinatesystem units   <1.5mm,2mm>
\setplotarea x from 0 to 16, y from -2 to 15
\put{922)} [l] at 2 12
\put {$ \scriptstyle \bullet$} [c] at  10 0
\put {$ \scriptstyle \bullet$} [c] at  10 12
\put {$ \scriptstyle \bullet$} [c] at  13 0
\put {$ \scriptstyle \bullet$} [c] at  13 3
\put {$ \scriptstyle \bullet$} [c] at  13 12
\put {$ \scriptstyle \bullet$} [c] at  16 12
\put {$ \scriptstyle \bullet$} [c] at  16  0
\setlinear \plot 10 0 10 12 13 0  16 12 16 0 13 12 13 0 /
\setlinear \plot  10 0 13 3 /
\setlinear \plot  10 12 16 0 /
\put{$5{,}040 $} [c] at 13 -2
\endpicture
\end{minipage}
\begin{minipage}{4cm}
\beginpicture
\setcoordinatesystem units   <1.5mm,2mm>
\setplotarea x from 0 to 16, y from -2 to 15
\put{923)} [l] at 2 12
\put {$ \scriptstyle \bullet$} [c] at  10 0
\put {$ \scriptstyle \bullet$} [c] at  10 12
\put {$ \scriptstyle \bullet$} [c] at  13 0
\put {$ \scriptstyle \bullet$} [c] at  13 9
\put {$ \scriptstyle \bullet$} [c] at  13 12
\put {$ \scriptstyle \bullet$} [c] at  16 12
\put {$ \scriptstyle \bullet$} [c] at  16  0
\setlinear \plot 10 12 10 0 13 12  16 0 16 12 13 0 13 12 /
\setlinear \plot  10 12 13 9 /
\setlinear \plot  10 0 16 12 /
\put{$5{,}040 $} [c] at 13 -2
\endpicture
\end{minipage}
\begin{minipage}{4cm}
\beginpicture
\setcoordinatesystem units   <1.5mm,2mm>
\setplotarea x from 0 to 16, y from -2 to 15
\put{924)} [l] at 2 12
\put {$ \scriptstyle \bullet$} [c] at 10 12
\put {$ \scriptstyle \bullet$} [c] at 13 4
\put {$ \scriptstyle \bullet$} [c] at 13 12
\put {$ \scriptstyle \bullet$} [c] at 12.5 8
\put {$ \scriptstyle \bullet$} [c] at 13.5 8
\put {$ \scriptstyle \bullet$} [c] at 14.5 0
\put {$ \scriptstyle \bullet$} [c] at 16 4
\setlinear \plot 10 12 13 4 14.5 0 16 4 13 12 12.5 8 13 4 13.5 8 13 12 /
\put{$2{,}520$} [c] at 13 -2
\endpicture
\end{minipage}
$$
$$
\begin{minipage}{4cm}
\beginpicture
\setcoordinatesystem units   <1.5mm,2mm>
\setplotarea x from 0 to 16, y from -2 to 15
\put{925)} [l] at 2 12
\put {$ \scriptstyle \bullet$} [c] at 10 0
\put {$ \scriptstyle \bullet$} [c] at 13 8
\put {$ \scriptstyle \bullet$} [c] at 13 0
\put {$ \scriptstyle \bullet$} [c] at 12.5 4
\put {$ \scriptstyle \bullet$} [c] at 13.5 4
\put {$ \scriptstyle \bullet$} [c] at 14.5 12
\put {$ \scriptstyle \bullet$} [c] at 16 8
\setlinear \plot 10 0 13 8 14.5 12 16 8 13 0 12.5 4 13 8 13.5 4 13 0 /
\put{$2{,}520$} [c] at 13 -2
\endpicture
\end{minipage}
\begin{minipage}{4cm}
\beginpicture
\setcoordinatesystem units   <1.5mm,2mm>
\setplotarea x from 0 to 16, y from -2 to 15
\put{926)} [l] at 2 12
\put {$ \scriptstyle \bullet$} [c] at 10 12
\put {$ \scriptstyle \bullet$} [c] at 12 4
\put {$ \scriptstyle \bullet$} [c] at 14 12
\put {$ \scriptstyle \bullet$} [c] at 14 0
\put {$ \scriptstyle \bullet$} [c] at 16 4
\put {$ \scriptstyle \bullet$} [c] at 16 8
\put {$ \scriptstyle \bullet$} [c] at 16 12
\setlinear \plot 16 12 16 4 14 0 12 4 14 12 /
\setlinear \plot 10 12 12 4 /
\put{$2{,}520$} [c] at 13 -2
\endpicture
\end{minipage}
\begin{minipage}{4cm}
\beginpicture
\setcoordinatesystem units   <1.5mm,2mm>
\setplotarea x from 0 to 16, y from -2 to 15
\put{927)} [l] at 2 12
\put {$ \scriptstyle \bullet$} [c] at 10 0
\put {$ \scriptstyle \bullet$} [c] at 12 8
\put {$ \scriptstyle \bullet$} [c] at 14 0
\put {$ \scriptstyle \bullet$} [c] at 14 12
\put {$ \scriptstyle \bullet$} [c] at 16 4
\put {$ \scriptstyle \bullet$} [c] at 16 8
\put {$ \scriptstyle \bullet$} [c] at 16 0
\setlinear \plot 16 0 16 8 14 12 12 8 14 0 /
\setlinear \plot 10 0 12 8 /
\put{$2{,}520$} [c] at 13 -2
\endpicture
\end{minipage}
\begin{minipage}{4cm}
\beginpicture
\setcoordinatesystem units   <1.5mm,2mm>
\setplotarea x from 0 to 16, y from -2 to 15
\put{928)} [l] at 2 12
\put {$ \scriptstyle \bullet$} [c] at 10 8
\put {$ \scriptstyle \bullet$} [c] at 10 12
\put {$ \scriptstyle \bullet$} [c] at 12 4
\put {$ \scriptstyle \bullet$} [c] at 14 8
\put {$ \scriptstyle \bullet$} [c] at 14 12
\put {$ \scriptstyle \bullet$} [c] at 16 0
\put {$ \scriptstyle \bullet$} [c] at 16 12
\setlinear \plot 16 12 16  0 12 4 14 8 14 12  /
\setlinear \plot 10 12  10 8  12 4   /
\put{$2{,}520$} [c] at 13 -2
\endpicture
\end{minipage}
\begin{minipage}{4cm}
\beginpicture
\setcoordinatesystem units   <1.5mm,2mm>
\setplotarea x from 0 to 16, y from -2 to 15
\put{929)} [l] at 2 12
\put {$ \scriptstyle \bullet$} [c] at 10 4
\put {$ \scriptstyle \bullet$} [c] at 10 0
\put {$ \scriptstyle \bullet$} [c] at 12 8
\put {$ \scriptstyle \bullet$} [c] at 14 4
\put {$ \scriptstyle \bullet$} [c] at 16 12
\put {$ \scriptstyle \bullet$} [c] at 14 0
\put {$ \scriptstyle \bullet$} [c] at 16 0
\setlinear \plot 16 0 16  12 12 8 14 4 14 0  /
\setlinear \plot 10 0  10 4  12 8   /
\put{$2{,}520$} [c] at 13 -2
\endpicture
\end{minipage}
\begin{minipage}{4cm}
\beginpicture
\setcoordinatesystem units   <1.5mm,2mm>
\setplotarea x from 0 to 16, y from -2 to 15
\put{930)} [l] at 2 12
\put {$ \scriptstyle \bullet$} [c] at 10 4
\put {$ \scriptstyle \bullet$} [c] at 12 12
\put {$ \scriptstyle \bullet$} [c] at 12 0
\put {$ \scriptstyle \bullet$} [c] at 14 4
\put {$ \scriptstyle \bullet$} [c] at 16 4
\put {$ \scriptstyle \bullet$} [c] at 16 8
\put {$ \scriptstyle \bullet$} [c] at 16 12
\setlinear \plot 16 12 16 4 12 0 14  4 12 12 10 4 12 0 /
\put{$2{,}520$} [c] at 13 -2
\endpicture
\end{minipage}
$$
$$
\begin{minipage}{4cm}
\beginpicture
\setcoordinatesystem units   <1.5mm,2mm>
\setplotarea x from 0 to 16, y from -2 to 15
\put{931)} [l] at 2 12
\put {$ \scriptstyle \bullet$} [c] at 10 8
\put {$ \scriptstyle \bullet$} [c] at 12 12
\put {$ \scriptstyle \bullet$} [c] at 12 0
\put {$ \scriptstyle \bullet$} [c] at 14 8
\put {$ \scriptstyle \bullet$} [c] at 16 8
\put {$ \scriptstyle \bullet$} [c] at 16 4
\put {$ \scriptstyle \bullet$} [c] at 16 0
\setlinear \plot 16 0 16 8 12 12 14  8 12 0 10 8 12 12 /
\put{$2{,}520$} [c] at 13 -2
\endpicture
\end{minipage}
\begin{minipage}{4cm}
\beginpicture
\setcoordinatesystem units   <1.5mm,2mm>
\setplotarea x from 0 to 16, y from -2 to 15
\put{932)} [l] at 2 12
\put {$ \scriptstyle \bullet$} [c] at 10 4
\put {$ \scriptstyle \bullet$} [c] at 10 8
\put {$ \scriptstyle \bullet$} [c] at 11 0
\put {$ \scriptstyle \bullet$} [c] at 11 12
\put {$ \scriptstyle \bullet$} [c] at 12 4
\put {$ \scriptstyle \bullet$} [c] at 12 8
\put {$ \scriptstyle \bullet$} [c] at 16 12
\setlinear \plot 16 12 11 0 10 4 10 8 11 12 12 8 12 4 11 0  /
\put{$2{,}520$} [c] at 13 -2
\endpicture
\end{minipage}
\begin{minipage}{4cm}
\beginpicture
\setcoordinatesystem units   <1.5mm,2mm>
\setplotarea x from 0 to 16, y from -2 to 15
\put{933)} [l] at 2 12
\put {$ \scriptstyle \bullet$} [c] at 10 4
\put {$ \scriptstyle \bullet$} [c] at 10 8
\put {$ \scriptstyle \bullet$} [c] at 11 0
\put {$ \scriptstyle \bullet$} [c] at 11 12
\put {$ \scriptstyle \bullet$} [c] at 12 4
\put {$ \scriptstyle \bullet$} [c] at 12 8
\put {$ \scriptstyle \bullet$} [c] at 16 0
\setlinear \plot 16 0 11 12 10 8 10 4 11 0 12 4 12 8 11 12  /
\put{$2{,}520$} [c] at 13 -2
\endpicture
\end{minipage}
\begin{minipage}{4cm}
\beginpicture
\setcoordinatesystem units   <1.5mm,2mm>
\setplotarea x from 0 to 16, y from -2 to 15
\put{934)} [l] at 2 12
\put {$ \scriptstyle \bullet$} [c] at 10 3
\put {$ \scriptstyle \bullet$} [c] at 10 6
\put {$ \scriptstyle \bullet$} [c] at 10 9
\put {$ \scriptstyle \bullet$} [c] at 10 12
\put {$ \scriptstyle \bullet$} [c] at 13 0
\put {$ \scriptstyle \bullet$} [c] at 13 12
\put {$ \scriptstyle \bullet$} [c] at 16 12
\setlinear \plot 10 12 10 3 13 0 16 12 /
\setlinear \plot 13 0 13 12 /
\put{$2{,}520$} [c] at 13 -2
\endpicture
\end{minipage}
\begin{minipage}{4cm}
\beginpicture
\setcoordinatesystem units   <1.5mm,2mm>
\setplotarea x from 0 to 16, y from -2 to 15
\put{935)} [l] at 2 12
\put {$ \scriptstyle \bullet$} [c] at 10 3
\put {$ \scriptstyle \bullet$} [c] at 10 6
\put {$ \scriptstyle \bullet$} [c] at 10 9
\put {$ \scriptstyle \bullet$} [c] at 10 0
\put {$ \scriptstyle \bullet$} [c] at 13 0
\put {$ \scriptstyle \bullet$} [c] at 13 12
\put {$ \scriptstyle \bullet$} [c] at 16 0
\setlinear \plot 10 0 10 9 13 12 16 0 /
\setlinear \plot 13 0 13 12 /
\put{$2{,}520$} [c] at 13 -2
\endpicture
\end{minipage}
\begin{minipage}{4cm}
\beginpicture
\setcoordinatesystem units   <1.5mm,2mm>
\setplotarea x from 0 to 16, y from -2 to 15
\put{936)} [l] at 2 12
\put {$ \scriptstyle \bullet$} [c] at  10 0
\put {$ \scriptstyle \bullet$} [c] at  10 6
\put {$ \scriptstyle \bullet$} [c] at  10 12
\put {$ \scriptstyle \bullet$} [c] at  14 6
\put {$ \scriptstyle \bullet$} [c] at  15 0
\put {$ \scriptstyle \bullet$} [c] at  15 12
\put {$ \scriptstyle \bullet$} [c] at  16 6
\setlinear \plot 10 12 10 0 15 12 16 6 15 0 14 6 15 12       /
\setlinear \plot 10 6 15 0 /
\put{$2{,}520  $} [c] at 13 -2
\endpicture
\end{minipage}
$$
$$
\begin{minipage}{4cm}
\beginpicture
\setcoordinatesystem units   <1.5mm,2mm>
\setplotarea x from 0 to 16, y from -2 to 15
\put{937)} [l] at 2 12
\put {$ \scriptstyle \bullet$} [c] at  10 0
\put {$ \scriptstyle \bullet$} [c] at  10 6
\put {$ \scriptstyle \bullet$} [c] at  10 12
\put {$ \scriptstyle \bullet$} [c] at  14 6
\put {$ \scriptstyle \bullet$} [c] at  15 0
\put {$ \scriptstyle \bullet$} [c] at  15 12
\put {$ \scriptstyle \bullet$} [c] at  16 6
\setlinear \plot 10 0 10 12 15 0 16 6 15 12 14 6 15 0       /
\setlinear \plot 10 6 15 12 /
\put{$2{,}520  $} [c] at 13 -2
\endpicture
\end{minipage}
\begin{minipage}{4cm}
\beginpicture
\setcoordinatesystem units   <1.5mm,2mm>
\setplotarea x from 0 to 16, y from -2 to 15
\put{938)} [l] at 2 12
\put {$ \scriptstyle \bullet$} [c] at  10 12
\put {$ \scriptstyle \bullet$} [c] at  13 0
\put {$ \scriptstyle \bullet$} [c] at  13 4
\put {$ \scriptstyle \bullet$} [c] at  13 8
\put {$ \scriptstyle \bullet$} [c] at  13 12
\put {$ \scriptstyle \bullet$} [c] at  16 0
\put {$ \scriptstyle \bullet$} [c] at  16 12
\setlinear \plot  10 12  13 0 13 12  /
\setlinear \plot  16 12 16 0 13 4   /
\put{$2{,}520$} [c] at 13 -2
\endpicture
\end{minipage}
\begin{minipage}{4cm}
\beginpicture
\setcoordinatesystem units   <1.5mm,2mm>
\setplotarea x from 0 to 16, y from -2 to 15
\put{939)} [l] at 2 12
\put {$ \scriptstyle \bullet$} [c] at  10 0
\put {$ \scriptstyle \bullet$} [c] at  13 0
\put {$ \scriptstyle \bullet$} [c] at  13 4
\put {$ \scriptstyle \bullet$} [c] at  13 8
\put {$ \scriptstyle \bullet$} [c] at  13 12
\put {$ \scriptstyle \bullet$} [c] at  16 0
\put {$ \scriptstyle \bullet$} [c] at  16 12
\setlinear \plot  10 0  13 12 13 0  /
\setlinear \plot  16 0 16 12 13 8   /
\put{$2{,}520$} [c] at 13 -2
\endpicture
\end{minipage}
\begin{minipage}{4cm}
\beginpicture
\setcoordinatesystem units   <1.5mm,2mm>
\setplotarea x from 0 to 16, y from -2 to 15
\put{940)} [l] at 2 12
\put {$ \scriptstyle \bullet$} [c] at  10 12
\put {$ \scriptstyle \bullet$} [c] at  10 0
\put {$ \scriptstyle \bullet$} [c] at  13 12
\put {$ \scriptstyle \bullet$} [c] at  13 6
\put {$ \scriptstyle \bullet$} [c] at  14.5 0
\put {$ \scriptstyle \bullet$} [c] at  16 6
\put {$ \scriptstyle \bullet$} [c] at  16 12
\setlinear \plot  10 0 10 12 13 6 13 12 16 6 16 12 13 6 14.5 0  16 6  /
\put{$2{,}520$} [c] at 13 -2
\endpicture
\end{minipage}
\begin{minipage}{4cm}
\beginpicture
\setcoordinatesystem units   <1.5mm,2mm>
\setplotarea x from 0 to 16, y from -2 to 15
\put{941)} [l] at 2 12
\put {$ \scriptstyle \bullet$} [c] at  10 12
\put {$ \scriptstyle \bullet$} [c] at  10 0
\put {$ \scriptstyle \bullet$} [c] at  13 0
\put {$ \scriptstyle \bullet$} [c] at  13 6
\put {$ \scriptstyle \bullet$} [c] at  14.5 12
\put {$ \scriptstyle \bullet$} [c] at  16 6
\put {$ \scriptstyle \bullet$} [c] at  16 0
\setlinear \plot  10 12 10 0 13 6 13 0 16 6 16 0 13 6 14.5 12  16 6  /
\put{$2{,}520$} [c] at 13 -2
\endpicture
\end{minipage}
\begin{minipage}{4cm}
\beginpicture
\setcoordinatesystem units   <1.5mm,2mm>
\setplotarea x from 0 to 16, y from -2 to 15
\put{942)} [l] at 2 12
\put {$ \scriptstyle \bullet$} [c] at  10 0
\put {$ \scriptstyle \bullet$} [c] at  10 6
\put {$ \scriptstyle \bullet$} [c] at  10 12
\put {$ \scriptstyle \bullet$} [c] at  13 12
\put {$ \scriptstyle \bullet$} [c] at  16 0
\put {$ \scriptstyle \bullet$} [c] at  16 6
\put {$ \scriptstyle \bullet$} [c] at  16 12
\setlinear \plot  16 12 16 0 13 12 10 6 10 0   /
\setlinear \plot  10 6 10 12 16 0  /
\put{$2{,}520$} [c] at 13 -2
\endpicture
\end{minipage}
$$
$$
\begin{minipage}{4cm}
\beginpicture
\setcoordinatesystem units   <1.5mm,2mm>
\setplotarea x from 0 to 16, y from -2 to 15
\put{943)} [l] at 2 12
\put {$ \scriptstyle \bullet$} [c] at  10 0
\put {$ \scriptstyle \bullet$} [c] at  10 6
\put {$ \scriptstyle \bullet$} [c] at  10 12
\put {$ \scriptstyle \bullet$} [c] at  13 0
\put {$ \scriptstyle \bullet$} [c] at  16 0
\put {$ \scriptstyle \bullet$} [c] at  16 6
\put {$ \scriptstyle \bullet$} [c] at  16 12
\setlinear \plot  16 0 16 12 13 0 10 6 10 12   /
\setlinear \plot  10 6 10 0 16 12  /
\put{$2{,}520$} [c] at 13 -2
\endpicture
\end{minipage}
\begin{minipage}{4cm}
\beginpicture
\setcoordinatesystem units   <1.5mm,2mm>
\setplotarea x from 0 to 16, y from -2 to 15
\put{944)} [l] at 2 12
\put {$ \scriptstyle \bullet$} [c] at  10 12
\put {$ \scriptstyle \bullet$} [c] at  10 8
\put {$ \scriptstyle \bullet$} [c] at  10 4
\put {$ \scriptstyle \bullet$} [c] at  10 0
\put {$ \scriptstyle \bullet$} [c] at  13 12
\put {$ \scriptstyle \bullet$} [c] at  16 12
\put {$ \scriptstyle \bullet$} [c] at  16 0
\setlinear \plot 16 0 10 12 10 0  13 12 16 0 16 12 10 0  /
\put{$2{,}520$} [c] at 13 -2
\endpicture
\end{minipage}
\begin{minipage}{4cm}
\beginpicture
\setcoordinatesystem units   <1.5mm,2mm>
\setplotarea x from 0 to 16, y from -2 to 15
\put{945)} [l] at 2 12
\put {$ \scriptstyle \bullet$} [c] at  10 12
\put {$ \scriptstyle \bullet$} [c] at  10 8
\put {$ \scriptstyle \bullet$} [c] at  10 4
\put {$ \scriptstyle \bullet$} [c] at  10 0
\put {$ \scriptstyle \bullet$} [c] at  13 0
\put {$ \scriptstyle \bullet$} [c] at  16 12
\put {$ \scriptstyle \bullet$} [c] at  16 0
\setlinear \plot 16 12 10 0 10 12  13 0 16 12 16 0 10 12  /
\put{$2{,}520$} [c] at 13 -2
\endpicture
\end{minipage}
\begin{minipage}{4cm}
\beginpicture
\setcoordinatesystem units   <1.5mm,2mm>
\setplotarea x from 0 to 16, y from -2 to 15
\put{946)} [l] at 2 12
\put {$ \scriptstyle \bullet$} [c] at  10 0
\put {$ \scriptstyle \bullet$} [c] at  14 12
\put {$ \scriptstyle \bullet$} [c] at  14 4
\put {$ \scriptstyle \bullet$} [c] at  15 0
\put {$ \scriptstyle \bullet$} [c] at  15 12
\put {$ \scriptstyle \bullet$} [c] at  16 12
\put {$ \scriptstyle \bullet$} [c] at  16 4
\setlinear \plot 10 0 14 12 14 4 15 0 16 4 16 12 10 0 15 12   16 4  /
\put{$2{,}520$} [c] at 13 -2
\endpicture
\end{minipage}
\begin{minipage}{4cm}
\beginpicture
\setcoordinatesystem units   <1.5mm,2mm>
\setplotarea x from 0 to 16, y from -2 to 15
\put{947)} [l] at 2 12
\put {$ \scriptstyle \bullet$} [c] at  10 12
\put {$ \scriptstyle \bullet$} [c] at  14 0
\put {$ \scriptstyle \bullet$} [c] at  14 8
\put {$ \scriptstyle \bullet$} [c] at  15 0
\put {$ \scriptstyle \bullet$} [c] at  15 12
\put {$ \scriptstyle \bullet$} [c] at  16 0
\put {$ \scriptstyle \bullet$} [c] at  16 8
\setlinear \plot 10 12 14 0 14 8 15 12 16 8 16 0 10 12 15 0   16 8  /
\put{$2{,}520$} [c] at 13 -2
\endpicture
\end{minipage}
\begin{minipage}{4cm}
\beginpicture
\setcoordinatesystem units   <1.5mm,2mm>
\setplotarea x from 0 to 16, y from -2 to 15
\put{948)} [l] at 2 12
\put {$ \scriptstyle \bullet$} [c] at  10 0
\put {$ \scriptstyle \bullet$} [c] at  10 12
\put {$ \scriptstyle \bullet$} [c] at  13 12
\put {$ \scriptstyle \bullet$} [c] at  13 6
\put {$ \scriptstyle \bullet$} [c] at  14.5 12
\put {$ \scriptstyle \bullet$} [c] at  14.5 0
\put {$ \scriptstyle \bullet$} [c] at  16 6
\setlinear \plot 10 0 10 12 13 6 13 12 10 0  /
\setlinear \plot  13 6 14.5 12 16 6 14.5 0 13 6  /
\put{$2{,}520$} [c] at 13 -2
\endpicture
\end{minipage}
$$
$$
\begin{minipage}{4cm}
\beginpicture
\setcoordinatesystem units   <1.5mm,2mm>
\setplotarea x from 0 to 16, y from -2 to 15
\put{949)} [l] at 2 12
\put {$ \scriptstyle \bullet$} [c] at  10 0
\put {$ \scriptstyle \bullet$} [c] at  10 12
\put {$ \scriptstyle \bullet$} [c] at  13 0
\put {$ \scriptstyle \bullet$} [c] at  13 6
\put {$ \scriptstyle \bullet$} [c] at  14.5 12
\put {$ \scriptstyle \bullet$} [c] at  14.5 0
\put {$ \scriptstyle \bullet$} [c] at  16 6
\setlinear \plot 10 12 10 0 13 6 13 0 10 12  /
\setlinear \plot  13 6 14.5 0 16 6 14.5 12 13 6  /
\put{$2{,}520$} [c] at 13 -2
\endpicture
\end{minipage}
\begin{minipage}{4cm}
\beginpicture
\setcoordinatesystem units   <1.5mm,2mm>
\setplotarea x from 0 to 16, y from -2 to 15
\put{950)} [l] at 2 12
\put {$ \scriptstyle \bullet$} [c] at  10 0
\put {$ \scriptstyle \bullet$} [c] at  10 12
\put {$ \scriptstyle \bullet$} [c] at  14 12
\put {$ \scriptstyle \bullet$} [c] at  14 6
\put {$ \scriptstyle \bullet$} [c] at  15 0
\put {$ \scriptstyle \bullet$} [c] at  16 6
\put {$ \scriptstyle \bullet$} [c] at  16 12
\setlinear \plot 10 0 10 12 15  0 16 6 16 12 14  6 14 12 16 6 /
\setlinear \plot  10 0 16 12  /
\setlinear \plot 14 12 14 6 15  0    /
\put{$2{,}520$} [c] at 13 -2
\endpicture
\end{minipage}
\begin{minipage}{4cm}
\beginpicture
\setcoordinatesystem units   <1.5mm,2mm>
\setplotarea x from 0 to 16, y from -2 to 15
\put{951)} [l] at 2 12
\put {$ \scriptstyle \bullet$} [c] at  10 0
\put {$ \scriptstyle \bullet$} [c] at  10 12
\put {$ \scriptstyle \bullet$} [c] at  14 0
\put {$ \scriptstyle \bullet$} [c] at  14 6
\put {$ \scriptstyle \bullet$} [c] at  15 12
\put {$ \scriptstyle \bullet$} [c] at  16 6
\put {$ \scriptstyle \bullet$} [c] at  16 0
\setlinear \plot 10 12 10 0 15  12 16 6 16 0 14  6 14 0 16 6 /
\setlinear \plot  10 12 16 0  /
\setlinear \plot 14 0 14 6 15  12    /
\put{$2{,}520$} [c] at 13 -2
\endpicture
\end{minipage}
\begin{minipage}{4cm}
\beginpicture
\setcoordinatesystem units   <1.5mm,2mm>
\setplotarea x from 0 to 16, y from -2 to 15
\put{952)} [l] at 2 12
\put {$ \scriptstyle \bullet$} [c] at  10 12
\put {$ \scriptstyle \bullet$} [c] at  12 12
\put {$ \scriptstyle \bullet$} [c] at  12 6
\put {$ \scriptstyle \bullet$} [c] at  14 0
\put {$ \scriptstyle \bullet$} [c] at  14 12
\put {$ \scriptstyle \bullet$} [c] at  16 6
\put {$ \scriptstyle \bullet$} [c] at  16 0
\setlinear \plot  10 12 12 6 14 0 16 6 16 0   /
\setlinear \plot  12 12  12 6 14 12 16 6  /
\put{$2{,}520$} [c] at 13 -2
\endpicture
\end{minipage}
\begin{minipage}{4cm}
\beginpicture
\setcoordinatesystem units   <1.5mm,2mm>
\setplotarea x from 0 to 16, y from -2 to 15
\put{953)} [l] at 2 12
\put {$ \scriptstyle \bullet$} [c] at  10 0
\put {$ \scriptstyle \bullet$} [c] at  12 0
\put {$ \scriptstyle \bullet$} [c] at  12 6
\put {$ \scriptstyle \bullet$} [c] at  14 0
\put {$ \scriptstyle \bullet$} [c] at  14 12
\put {$ \scriptstyle \bullet$} [c] at  16 6
\put {$ \scriptstyle \bullet$} [c] at  16 12
\setlinear \plot  10 0 12 6 14 12 16 6 16 12   /
\setlinear \plot  12 0  12 6 14 0 16 6  /
\put{$2{,}520$} [c] at 13 -2
\endpicture
\end{minipage}
\begin{minipage}{4cm}
\beginpicture
\setcoordinatesystem units   <1.5mm,2mm>
\setplotarea x from 0 to 16, y from -2 to 15
\put{954)} [l] at 2 12
\put {$ \scriptstyle \bullet$} [c] at  16 0
\put {$ \scriptstyle \bullet$} [c] at  16 6
\put {$ \scriptstyle \bullet$} [c] at  16 12
\put {$ \scriptstyle \bullet$} [c] at  13 6
\put {$ \scriptstyle \bullet$} [c] at  13 0
\put {$ \scriptstyle \bullet$} [c] at  13 12
\put {$ \scriptstyle \bullet$} [c] at  10 12
\setlinear \plot  10 12 13 0 13 12 16 6 16 0 /
\setlinear \plot  13  6 16 12 16  6  /
\put{$2{,}520$} [c] at 13 -2
\endpicture
\end{minipage}
$$
$$
\begin{minipage}{4cm}
\beginpicture
\setcoordinatesystem units   <1.5mm,2mm>
\setplotarea x from 0 to 16, y from -2 to 15
\put{955)} [l] at 2 12
\put {$ \scriptstyle \bullet$} [c] at  16 0
\put {$ \scriptstyle \bullet$} [c] at  16 6
\put {$ \scriptstyle \bullet$} [c] at  16 12
\put {$ \scriptstyle \bullet$} [c] at  13 6
\put {$ \scriptstyle \bullet$} [c] at  13 0
\put {$ \scriptstyle \bullet$} [c] at  13 12
\put {$ \scriptstyle \bullet$} [c] at  10 0
\setlinear \plot  10 0 13 12 13 0 16 6 16 12 /
\setlinear \plot  13  6 16 0 16  6  /
\put{$2{,}520$} [c] at 13 -2
\endpicture
\end{minipage}
\begin{minipage}{4cm}
\beginpicture
\setcoordinatesystem units   <1.5mm,2mm>
\setplotarea x from 0 to 16, y from -2 to 15
\put{956)} [l] at 2 12
\put {$ \scriptstyle \bullet$} [c] at  10 0
\put {$ \scriptstyle \bullet$} [c] at  10 12
\put {$ \scriptstyle \bullet$} [c] at  13 12
\put {$ \scriptstyle \bullet$} [c] at  13 0
\put {$ \scriptstyle \bullet$} [c] at  16 0
\put {$ \scriptstyle \bullet$} [c] at  16  6
\put {$ \scriptstyle \bullet$} [c] at  16  12
\setlinear \plot 10 0 10 12  13 0 16 12 16 0  /
\setlinear \plot 10 0 13 12 16 6  /
\setlinear \plot 10 12 16 6 /
\put{$2{,}520$} [c] at 13 -2
\endpicture
\end{minipage}
\begin{minipage}{4cm}
\beginpicture
\setcoordinatesystem units   <1.5mm,2mm>
\setplotarea x from 0 to 16, y from -2 to 15
\put{957)} [l] at 2 12
\put {$ \scriptstyle \bullet$} [c] at  10 0
\put {$ \scriptstyle \bullet$} [c] at  10 12
\put {$ \scriptstyle \bullet$} [c] at  13 12
\put {$ \scriptstyle \bullet$} [c] at  13 0
\put {$ \scriptstyle \bullet$} [c] at  16 0
\put {$ \scriptstyle \bullet$} [c] at  16  6
\put {$ \scriptstyle \bullet$} [c] at  16  12
\setlinear \plot 10 12 10 0  13 12 16 0 16 12  /
\setlinear \plot 10 12 13 0 16 6  /
\setlinear \plot 10 0 16 6 /
\put{$2{,}520$} [c] at 13 -2
\endpicture
\end{minipage}
\begin{minipage}{4cm}
\beginpicture
\setcoordinatesystem units   <1.5mm,2mm>
\setplotarea x from 0 to 16, y from -2 to 15
\put{958)} [l] at 2 12
\put {$ \scriptstyle \bullet$} [c] at  10 12
\put {$ \scriptstyle \bullet$} [c] at  10 6
\put {$ \scriptstyle \bullet$} [c] at  10 0
\put {$ \scriptstyle \bullet$} [c] at  13 12
\put {$ \scriptstyle \bullet$} [c] at  16 0
\put {$ \scriptstyle \bullet$} [c] at  16 6
\put {$ \scriptstyle \bullet$} [c] at  16 12
\setlinear \plot  10 12  10 0   13 12 16 0 16 12   /
\setlinear \plot  10 6  16 0   /
\setlinear \plot  10 0 16 6   /
\put{$1{,}260$} [c] at 13 -2
\endpicture
\end{minipage}
\begin{minipage}{4cm}
\beginpicture
\setcoordinatesystem units   <1.5mm,2mm>
\setplotarea x from 0 to 16, y from -2 to 15
\put{959)} [l] at 2 12
\put {$ \scriptstyle \bullet$} [c] at  10 12
\put {$ \scriptstyle \bullet$} [c] at  10 6
\put {$ \scriptstyle \bullet$} [c] at  10 0
\put {$ \scriptstyle \bullet$} [c] at  13 0
\put {$ \scriptstyle \bullet$} [c] at  16 0
\put {$ \scriptstyle \bullet$} [c] at  16 6
\put {$ \scriptstyle \bullet$} [c] at  16 12
\setlinear \plot  10 0  10 12   13 0 16 12 16 0   /
\setlinear \plot  10 6  16 12   /
\setlinear \plot  10 12 16 6   /
\put{$1{,}260$} [c] at 13 -2
\endpicture
\end{minipage}
\begin{minipage}{4cm}
\beginpicture
\setcoordinatesystem units   <1.5mm,2mm>
\setplotarea x from 0 to 16, y from -2 to 15
\put{960)} [l] at 2 12
\put {$ \scriptstyle \bullet$} [c] at  10 12
\put {$ \scriptstyle \bullet$} [c] at  10 0
\put {$ \scriptstyle \bullet$} [c] at  14 12
\put {$ \scriptstyle \bullet$} [c] at  14 6
\put {$ \scriptstyle \bullet$} [c] at  15 0
\put {$ \scriptstyle \bullet$} [c] at  16 12
\put {$ \scriptstyle \bullet$} [c] at  16 6
\setlinear \plot  10 12  10 0 16 12 16 6 14 12 14 6 15 0  16 6   /
\setlinear \plot  10 0 14 12  /
\setlinear \plot  14 6  16 12  /
\put{$1{,}260$} [c] at 13 -2
\endpicture
\end{minipage}
$$
$$
\begin{minipage}{4cm}
\beginpicture
\setcoordinatesystem units   <1.5mm,2mm>
\setplotarea x from 0 to 16, y from -2 to 15
\put{961)} [l] at 2 12
\put {$ \scriptstyle \bullet$} [c] at  10 12
\put {$ \scriptstyle \bullet$} [c] at  10 0
\put {$ \scriptstyle \bullet$} [c] at  14 0
\put {$ \scriptstyle \bullet$} [c] at  14 6
\put {$ \scriptstyle \bullet$} [c] at  15 12
\put {$ \scriptstyle \bullet$} [c] at  16 0
\put {$ \scriptstyle \bullet$} [c] at  16 6
\setlinear \plot  10 0  10 12 16 0 16 6 14 0 14 6 15 12  16 6   /
\setlinear \plot  10 12 14 0  /
\setlinear \plot  14 6  16 0  /
\put{$1{,}260$} [c] at 13 -2
\endpicture
\end{minipage}
\begin{minipage}{4cm}
\beginpicture
\setcoordinatesystem units   <1.5mm,2mm>
\setplotarea x from 0 to 16, y from -2 to 15
\put{962)} [l] at 2 12
\put {$ \scriptstyle \bullet$} [c] at  10 0
\put {$ \scriptstyle \bullet$} [c] at  10 12
\put {$ \scriptstyle \bullet$} [c] at  13 0
\put {$ \scriptstyle \bullet$} [c] at  13 12
\put {$ \scriptstyle \bullet$} [c] at  16  0
\put {$ \scriptstyle \bullet$} [c] at  16  6
\put {$ \scriptstyle \bullet$} [c] at  16  12
\setlinear \plot  16 0  16 12 13 0  13 12  10  0 10 12 13 0 /
\setlinear \plot  10 0 16  12   /
\setlinear \plot  13 12 16  6   /
\put{$1{,}260 $} [c] at 13 -2
\endpicture
\end{minipage}
\begin{minipage}{4cm}
\beginpicture
\setcoordinatesystem units   <1.5mm,2mm>
\setplotarea x from 0 to 16, y from -2 to 15
\put{963)} [l] at 2 12
\put {$ \scriptstyle \bullet$} [c] at  10 0
\put {$ \scriptstyle \bullet$} [c] at  10 12
\put {$ \scriptstyle \bullet$} [c] at  13 0
\put {$ \scriptstyle \bullet$} [c] at  13 12
\put {$ \scriptstyle \bullet$} [c] at  16  0
\put {$ \scriptstyle \bullet$} [c] at  16  6
\put {$ \scriptstyle \bullet$} [c] at  16  12
\setlinear \plot  16 12  16 0 13 12  13 0  10  12 10 0 13 12 /
\setlinear \plot  10 12 16  0   /
\setlinear \plot  13 0 16  6   /
\put{$1{,}260 $} [c] at 13 -2
\endpicture
\end{minipage}
\begin{minipage}{4cm}
\beginpicture
\setcoordinatesystem units   <1.5mm,2mm>
\setplotarea x from 0 to 16, y from -2 to 15
\put{964)} [l] at 2 12
\put {$ \scriptstyle \bullet$} [c] at  10 0
\put {$ \scriptstyle \bullet$} [c] at  10 6
\put {$ \scriptstyle \bullet$} [c] at  10 12
\put {$ \scriptstyle \bullet$} [c] at  13 12
\put {$ \scriptstyle \bullet$} [c] at  13  0
\put {$ \scriptstyle \bullet$} [c] at  16  12
\put {$ \scriptstyle \bullet$} [c] at  16  0
\setlinear \plot  10 12 10 0  16 12  16  0  /
\setlinear \plot  16 12 13 0 10 6 13 12   /
\put{$1{,}260$} [c] at 13 -2
\endpicture
\end{minipage}
\begin{minipage}{4cm}
\beginpicture
\setcoordinatesystem units   <1.5mm,2mm>
\setplotarea x from 0 to 16, y from -2 to 15
\put{965)} [l] at 2 12
\put {$ \scriptstyle \bullet$} [c] at  10 0
\put {$ \scriptstyle \bullet$} [c] at  10 6
\put {$ \scriptstyle \bullet$} [c] at  10 12
\put {$ \scriptstyle \bullet$} [c] at  13 12
\put {$ \scriptstyle \bullet$} [c] at  13  0
\put {$ \scriptstyle \bullet$} [c] at  16  12
\put {$ \scriptstyle \bullet$} [c] at  16  0
\setlinear \plot  10 0 10 12  16 0  16  12  /
\setlinear \plot  16 0 13 12 10 6 13 0   /
\put{$1{,}260$} [c] at 13 -2
\endpicture
\end{minipage}
\begin{minipage}{4cm}
\beginpicture
\setcoordinatesystem units   <1.5mm,2mm>
\setplotarea x from 0 to 16, y from -2 to 15
\put{966)} [l] at 2 12
\put {$ \scriptstyle \bullet$} [c] at  10 6
\put {$ \scriptstyle \bullet$} [c] at  11 0
\put {$ \scriptstyle \bullet$} [c] at 11 6
\put {$ \scriptstyle \bullet$} [c] at 11 9
\put {$ \scriptstyle \bullet$} [c] at 11 12
\put {$ \scriptstyle \bullet$} [c] at 12 6
\put {$ \scriptstyle \bullet$} [c] at 16 12
\setlinear \plot 11 12 11 9 10 6 11 0 11 9 12 6 11 0 16 12  /
\put{$840$} [c] at 13 -2
\endpicture
\end{minipage}
$$
$$
\begin{minipage}{4cm}
\beginpicture
\setcoordinatesystem units   <1.5mm,2mm>
\setplotarea x from 0 to 16, y from -2 to 15
\put{967)} [l] at 2 12
\put {$ \scriptstyle \bullet$} [c] at  10 8
\put {$ \scriptstyle \bullet$} [c] at  11 0
\put {$ \scriptstyle \bullet$} [c] at 11 4
\put {$ \scriptstyle \bullet$} [c] at 11 8
\put {$ \scriptstyle \bullet$} [c] at 11 12
\put {$ \scriptstyle \bullet$} [c] at 12 8
\put {$ \scriptstyle \bullet$} [c] at 16 0
\setlinear \plot 11 0 11 4 10 8 11 12 11 4 12 8 11 12 16 0  /
\put{$840$} [c] at 13 -2
\endpicture
\end{minipage}
\begin{minipage}{4cm}
\beginpicture
\setcoordinatesystem units   <1.5mm,2mm>
\setplotarea x from 0 to 16, y from -2 to 15
\put{968)} [l] at 2 12
\put {$ \scriptstyle \bullet$} [c] at 10 8
\put {$ \scriptstyle \bullet$} [c] at 12 4
\put {$ \scriptstyle \bullet$} [c] at 12 8
\put {$ \scriptstyle \bullet$} [c] at 12 12
\put {$ \scriptstyle \bullet$} [c] at 14 8
\put {$ \scriptstyle \bullet$} [c] at 16 0
\put {$ \scriptstyle \bullet$} [c] at 16 12
\setlinear \plot 16 12 16 0 12 4 14 8 12 12 12 4 10 8 12 12  /
\put{$840$} [c] at 13 -2
\endpicture
\end{minipage}
\begin{minipage}{4cm}
\beginpicture
\setcoordinatesystem units   <1.5mm,2mm>
\setplotarea x from 0 to 16, y from -2 to 15
\put{969)} [l] at 2 12
\put {$ \scriptstyle \bullet$} [c] at 10 4
\put {$ \scriptstyle \bullet$} [c] at 12 0
\put {$ \scriptstyle \bullet$} [c] at 12 4
\put {$ \scriptstyle \bullet$} [c] at 12 8
\put {$ \scriptstyle \bullet$} [c] at 14 4
\put {$ \scriptstyle \bullet$} [c] at 16 0
\put {$ \scriptstyle \bullet$} [c] at 16 12
\setlinear \plot 16 0 16 12 12 8 14 4 12 0 12 8 10 4 12 0  /
\put{$840$} [c] at 13 -2
\endpicture
\end{minipage}
\begin{minipage}{4cm}
\beginpicture
\setcoordinatesystem units   <1.5mm,2mm>
\setplotarea x from 0 to 16, y from -2 to 15
\put{970)} [l] at 2 12
\put {$ \scriptstyle \bullet$} [c] at 13 0
\put {$ \scriptstyle \bullet$} [c] at 10 6
\put {$ \scriptstyle \bullet$} [c] at 12 12
\put {$ \scriptstyle \bullet$} [c] at 12 6
\put {$ \scriptstyle \bullet$} [c] at 14 6
\put {$ \scriptstyle \bullet$} [c] at 16 6
\put {$ \scriptstyle \bullet$} [c] at 16 12
\setlinear \plot 13 0 10 6 12 12 12 6 13  0 14 6 16 12  16 6 13 0 /
\setlinear \plot 14 6 12 12 16 6    /
\setlinear \plot 12 6 16 12     /
\put{$840$} [c] at 13 -2
\endpicture
\end{minipage}
\begin{minipage}{4cm}
\beginpicture
\setcoordinatesystem units   <1.5mm,2mm>
\setplotarea x from 0 to 16, y from -2 to 15
\put{971)} [l] at 2 12
\put {$ \scriptstyle \bullet$} [c] at 13 12
\put {$ \scriptstyle \bullet$} [c] at 10 6
\put {$ \scriptstyle \bullet$} [c] at 12 0
\put {$ \scriptstyle \bullet$} [c] at 12 6
\put {$ \scriptstyle \bullet$} [c] at 14 6
\put {$ \scriptstyle \bullet$} [c] at 16 6
\put {$ \scriptstyle \bullet$} [c] at 16 0
\setlinear \plot 13 12 10 6 12 0 12 6 13  12 14 6 16 0  16 6 13 12 /
\setlinear \plot 14 6 12 0 16 6    /
\setlinear \plot 12 6 16 0     /
\put{$840$} [c] at 13 -2
\endpicture
\end{minipage}
\begin{minipage}{4cm}
\beginpicture
\setcoordinatesystem units   <1.5mm,2mm>
\setplotarea x from 0 to 16, y from -2 to 15
\put{972)} [l] at 2 12
\put {$ \scriptstyle \bullet$} [c] at 13 0
\put {$ \scriptstyle \bullet$} [c] at 10 6
\put {$ \scriptstyle \bullet$} [c] at 10 12
\put {$ \scriptstyle \bullet$} [c] at 12 12
\put {$ \scriptstyle \bullet$} [c] at 14 12
\put {$ \scriptstyle \bullet$} [c] at 16 12
\put {$ \scriptstyle \bullet$} [c] at 16 6
\setlinear \plot 10 12  10 6 13  0 16 6 16 12  10 6 12 12  16 6 14 12 10 6 /
\setlinear \plot 14 12 16 6    /
\put{$840$} [c] at 13 -2
\endpicture
\end{minipage}
$$
$$
\begin{minipage}{4cm}
\beginpicture
\setcoordinatesystem units   <1.5mm,2mm>
\setplotarea x from 0 to 16, y from -2 to 15
\put{973)} [l] at 2 12
\put {$ \scriptstyle \bullet$} [c] at 13 12
\put {$ \scriptstyle \bullet$} [c] at 10 6
\put {$ \scriptstyle \bullet$} [c] at 10 0
\put {$ \scriptstyle \bullet$} [c] at 12 0
\put {$ \scriptstyle \bullet$} [c] at 14 0
\put {$ \scriptstyle \bullet$} [c] at 16 0
\put {$ \scriptstyle \bullet$} [c] at 16 6
\setlinear \plot 10 0  10 6 13  12 16 6 16 0  10 6 12 0  16 6 14 0 10 6 /
\setlinear \plot 14 0 16 6    /
\put{$840$} [c] at 13 -2
\endpicture
\end{minipage}
\begin{minipage}{4cm}
\beginpicture
\setcoordinatesystem units   <1.5mm,2mm>
\setplotarea x from 0 to 16, y from -2 to 15
\put{974)} [l] at 2 12
\put {$ \scriptstyle \bullet$} [c] at 10 12
\put {$ \scriptstyle \bullet$} [c] at 12 12
\put {$ \scriptstyle \bullet$} [c] at 14 12
\put {$ \scriptstyle \bullet$} [c] at 12 8
\put {$ \scriptstyle \bullet$} [c] at 12 4
\put {$ \scriptstyle \bullet$} [c] at 16 0
\put {$ \scriptstyle \bullet$} [c] at 16 12
\setlinear \plot 16 12  16 0 12 4 12 8 14 12  /
\setlinear \plot 10 12 12 8 12 12   /
\put{$840$} [c] at 13 -2
\endpicture
\end{minipage}
\begin{minipage}{4cm}
\beginpicture
\setcoordinatesystem units   <1.5mm,2mm>
\setplotarea x from 0 to 16, y from -2 to 15
\put{975)} [l] at 2 12
\put {$ \scriptstyle \bullet$} [c] at 10 0
\put {$ \scriptstyle \bullet$} [c] at 12 0
\put {$ \scriptstyle \bullet$} [c] at 14 0
\put {$ \scriptstyle \bullet$} [c] at 12 8
\put {$ \scriptstyle \bullet$} [c] at 12 4
\put {$ \scriptstyle \bullet$} [c] at 16 0
\put {$ \scriptstyle \bullet$} [c] at 16 12
\setlinear \plot 16 0  16 12 12 8 12 4 14 0  /
\setlinear \plot 10 0 12 4 12 0   /
\put{$840$} [c] at 13 -2
\endpicture
\end{minipage}
\begin{minipage}{4cm}
\beginpicture
\setcoordinatesystem units   <1.5mm,2mm>
\setplotarea x from 0 to 16, y from -2 to 15
\put{976)} [l] at 2 12
\put {$ \scriptstyle \bullet$} [c] at  10 6
\put {$ \scriptstyle \bullet$} [c] at  11 0
\put {$ \scriptstyle \bullet$} [c] at  11 6
\put {$ \scriptstyle \bullet$} [c] at  11 12
\put {$ \scriptstyle \bullet$} [c] at  12 6
\put {$ \scriptstyle \bullet$} [c] at  16 0
\put {$ \scriptstyle \bullet$} [c] at  16 12
\setlinear \plot 11 0 10 6 11 12 11 0 12 6 11 12        /
\setlinear \plot 11 0 16 12 16 0 12 6      /
\setlinear \plot 10 6 16 0  11 6        /
\put{$420  $} [c] at 13 -2
\endpicture
\end{minipage}
\begin{minipage}{4cm}
\beginpicture
\setcoordinatesystem units   <1.5mm,2mm>
\setplotarea x from 0 to 16, y from -2 to 15
\put{977)} [l] at 2 12
\put {$ \scriptstyle \bullet$} [c] at  10 6
\put {$ \scriptstyle \bullet$} [c] at  11 0
\put {$ \scriptstyle \bullet$} [c] at  11 6
\put {$ \scriptstyle \bullet$} [c] at  11 12
\put {$ \scriptstyle \bullet$} [c] at  12 6
\put {$ \scriptstyle \bullet$} [c] at  16 0
\put {$ \scriptstyle \bullet$} [c] at  16 12
\setlinear \plot 11 12 10 6 11 0 11 12 12 6 11 0        /
\setlinear \plot 11 12 16 0 16 12 12 6      /
\setlinear \plot 10 6 16 12  11 6        /
\put{$420  $} [c] at 13 -2
\endpicture
\end{minipage}
\begin{minipage}{4cm}
\beginpicture
\setcoordinatesystem units   <1.5mm,2mm>
\setplotarea x from 0 to 16, y from -2 to 15
\put {978)} [l] at 2 12
\put {$ \scriptstyle \bullet$} [c] at  10 12
\put {$ \scriptstyle \bullet$} [c] at  12  12
\put {$ \scriptstyle \bullet$} [c] at  14 12
\put {$ \scriptstyle \bullet$} [c] at  14 6
\put {$ \scriptstyle \bullet$} [c] at  14 0
\put {$ \scriptstyle \bullet$} [c] at  16  0
\put {$ \scriptstyle \bullet$} [c] at  16  12
\setlinear \plot  10  12 14 6 16 0 16 12 14  0 14 12 /
\setlinear \plot  12 12 14 6 /
\put{$420$}[c] at 13 -2
\endpicture
\end{minipage}
$$
$$
\begin{minipage}{4cm}
\beginpicture
\setcoordinatesystem units   <1.5mm,2mm>
\setplotarea x from 0 to 16, y from -2 to 15
\put {979)} [l] at 2 12
\put {$ \scriptstyle \bullet$} [c] at  10 0
\put {$ \scriptstyle \bullet$} [c] at  12 0
\put {$ \scriptstyle \bullet$} [c] at  14 12
\put {$ \scriptstyle \bullet$} [c] at  14 6
\put {$ \scriptstyle \bullet$} [c] at  14 0
\put {$ \scriptstyle \bullet$} [c] at  16  0
\put {$ \scriptstyle \bullet$} [c] at  16  12
\setlinear \plot  10  0 14 6 16 12 16 0 14  12 14 0 /
\setlinear \plot  12 0 14 6 /
\put{$420$}[c] at 13 -2
\endpicture
\end{minipage}
\begin{minipage}{4cm}
\beginpicture
\setcoordinatesystem units   <1.5mm,2mm>
\setplotarea x from 0 to 16, y from -2 to 15
\put {980)} [l] at 2 12
\put {$ \scriptstyle \bullet$} [c] at  10 0
\put {$ \scriptstyle \bullet$} [c] at  10 12
\put {$ \scriptstyle \bullet$} [c] at  12 12
\put {$ \scriptstyle \bullet$} [c] at  14 12
\put {$ \scriptstyle \bullet$} [c] at  16 0
\put {$ \scriptstyle \bullet$} [c] at  16  6
\put {$ \scriptstyle \bullet$} [c] at  16  12
\setlinear \plot  10  12 10 0 16 12 16 0  /
\setlinear \plot  10 12 16 6 14 12 10 0 12 12 16 6 /
\put{$210$}[c] at 13 -2
\endpicture
\end{minipage}
\begin{minipage}{4cm}
\beginpicture
\setcoordinatesystem units   <1.5mm,2mm>
\setplotarea x from 0 to 16, y from -2 to 15
\put {981)} [l] at 2 12
\put {$ \scriptstyle \bullet$} [c] at  10 0
\put {$ \scriptstyle \bullet$} [c] at  10 12
\put {$ \scriptstyle \bullet$} [c] at  12 0
\put {$ \scriptstyle \bullet$} [c] at  14 0
\put {$ \scriptstyle \bullet$} [c] at  16 0
\put {$ \scriptstyle \bullet$} [c] at  16  6
\put {$ \scriptstyle \bullet$} [c] at  16  12
\setlinear \plot  10  0 10 12 16 0 16 12  /
\setlinear \plot  10 0 16 6 14 0 10 12 12 0 16 6 /
\put{$210$}[c] at 13 -2
\endpicture
\end{minipage}
\begin{minipage}{4cm}
\beginpicture
\setcoordinatesystem units   <1.5mm,2mm>
\setplotarea x from 0 to  16, y from -2 to 15
\put{982)} [l] at 2 12
\put {$ \scriptstyle \bullet$} [c] at 10 8
\put {$ \scriptstyle \bullet$} [c] at 12 0
\put {$ \scriptstyle \bullet$} [c] at 12 4
\put {$ \scriptstyle \bullet$} [c] at 12 12
\put {$ \scriptstyle \bullet$} [c] at 14 8
\put {$ \scriptstyle \bullet$} [c] at 16 0
\put {$ \scriptstyle \bullet$} [c] at 16 12
\setlinear \plot 12 0 12 4 10  8 12 12 14 8  12 4 /
\setlinear \plot 16 0 16  12 /
\put{$2{,}520$} [c] at 13 -2
\endpicture
\end{minipage}
\begin{minipage}{4cm}
\beginpicture
\setcoordinatesystem units   <1.5mm,2mm>
\setplotarea x from 0 to  16, y from -2 to 15
\put{983)} [l] at 2 12
\put {$ \scriptstyle \bullet$} [c] at 10 4
\put {$ \scriptstyle \bullet$} [c] at 12 12
\put {$ \scriptstyle \bullet$} [c] at 12 8
\put {$ \scriptstyle \bullet$} [c] at 12 0
\put {$ \scriptstyle \bullet$} [c] at 14 4
\put {$ \scriptstyle \bullet$} [c] at 16 0
\put {$ \scriptstyle \bullet$} [c] at 16 12
\setlinear \plot 12 12 12 8 10  4 12 0 14 4  12 8 /
\setlinear \plot 16 0 16  12 /
\put{$2{,}520$} [c] at 13 -2
\endpicture
\end{minipage}
\begin{minipage}{4cm}
\beginpicture
\setcoordinatesystem units   <1.5mm,2mm>
\setplotarea x from 0 to  16, y from -2 to 15
\put{984)} [l] at 2 12
\put {$ \scriptstyle \bullet$} [c] at 10 12
\put {$ \scriptstyle \bullet$} [c] at 12 8
\put {$ \scriptstyle \bullet$} [c] at 12  4
\put {$ \scriptstyle \bullet$} [c] at 12 0
\put {$ \scriptstyle \bullet$} [c] at 14 12
\put {$ \scriptstyle \bullet$} [c] at 16 0
\put {$ \scriptstyle \bullet$} [c] at 16 12
\setlinear \plot 12 0  12 8  10 12   /
\setlinear \plot 12 8  14 12   /
\setlinear \plot 16 0 16  12 /
\put{$2{,}520$} [c] at 13 -2
\endpicture
\end{minipage}
$$
$$
\begin{minipage}{4cm}
\beginpicture
\setcoordinatesystem units   <1.5mm,2mm>
\setplotarea x from 0 to  16, y from -2 to 15
\put{985)} [l] at 2 12
\put {$ \scriptstyle \bullet$} [c] at 10 0
\put {$ \scriptstyle \bullet$} [c] at 12 8
\put {$ \scriptstyle \bullet$} [c] at 12  4
\put {$ \scriptstyle \bullet$} [c] at 12 12
\put {$ \scriptstyle \bullet$} [c] at 14 0
\put {$ \scriptstyle \bullet$} [c] at 16 0
\put {$ \scriptstyle \bullet$} [c] at 16 12
\setlinear \plot 12 12  12 4  10 0   /
\setlinear \plot 12 4  14 0   /
\setlinear \plot 16 0 16  12 /
\put{$2{,}520$} [c] at 13 -2
\endpicture
\end{minipage}
\begin{minipage}{4cm}
\beginpicture
\setcoordinatesystem units   <1.5mm,2mm>
\setplotarea x from 0 to 16, y from -2 to 15
\put{${\bf  22}$} [l] at 2 15
\put{986)} [l] at 2 12
\put {$ \scriptstyle \bullet$} [c] at 11.5 0
\put {$ \scriptstyle \bullet$} [c] at 10 6
\put {$ \scriptstyle \bullet$} [c] at 13 6
\put {$ \scriptstyle \bullet$} [c] at 16 6
\put {$ \scriptstyle \bullet$} [c] at 16 12
\put {$ \scriptstyle \bullet$} [c] at 11.5 12
\put {$ \scriptstyle \bullet$} [c] at 10.8 9
\setlinear \plot 16 12 16 6 11.5 0 10 6 11.5 12 13 6  11.5 0 /
\put{$5{,}040$} [c] at  13 -2
\endpicture
\end{minipage}
\begin{minipage}{4cm}
\beginpicture
\setcoordinatesystem units   <1.5mm,2mm>
\setplotarea x from 0 to 16, y from -2 to 15
\put{987)} [l] at 2 12
\put {$ \scriptstyle \bullet$} [c] at 11.5 0
\put {$ \scriptstyle \bullet$} [c] at 10 6
\put {$ \scriptstyle \bullet$} [c] at 13 6
\put {$ \scriptstyle \bullet$} [c] at 16 6
\put {$ \scriptstyle \bullet$} [c] at 16 0
\put {$ \scriptstyle \bullet$} [c] at 11.5 12
\put {$ \scriptstyle \bullet$} [c] at 10.8 3
\setlinear \plot 16 0 16 6 11.5 12 10 6 11.5 0 13 6  11.5 12 /
\put{$5{,}040$} [c] at  13 -2
\endpicture
\end{minipage}
\begin{minipage}{4cm}
\beginpicture
\setcoordinatesystem units   <1.5mm,2mm>
\setplotarea x from 0 to 16, y from -2 to 15
\put{988)} [l] at 2 12
\put {$ \scriptstyle \bullet$} [c] at 13 0
\put {$ \scriptstyle \bullet$} [c] at 10 6
\put {$ \scriptstyle \bullet$} [c] at 13 6
\put {$ \scriptstyle \bullet$} [c] at 16 6
\put {$ \scriptstyle \bullet$} [c] at 16 12
\put {$ \scriptstyle \bullet$} [c] at 10 12
\put {$ \scriptstyle \bullet$} [c] at 13 12
\setlinear \plot 13 0 10 6 10 12 13 6  13 12 16 6 16  12 /
\setlinear \plot 13 6 13 0 16 6 /
\put{$5{,}040$} [c] at  13 -2
\endpicture
\end{minipage}
\begin{minipage}{4cm}
\beginpicture
\setcoordinatesystem units   <1.5mm,2mm>
\setplotarea x from 0 to 16, y from -2 to 15
\put{989)} [l] at 2 12
\put {$ \scriptstyle \bullet$} [c] at 13 0
\put {$ \scriptstyle \bullet$} [c] at 10 6
\put {$ \scriptstyle \bullet$} [c] at 13 6
\put {$ \scriptstyle \bullet$} [c] at 16 6
\put {$ \scriptstyle \bullet$} [c] at 16 0
\put {$ \scriptstyle \bullet$} [c] at 10 0
\put {$ \scriptstyle \bullet$} [c] at 13 12
\setlinear \plot 13 12 10 6 10 0 13 6  13 0 16 6 16  0 /
\setlinear \plot 13 6 13 12 16 6 /
\put{$5{,}040$} [c] at  13 -2
\endpicture
\end{minipage}
\begin{minipage}{4cm}
\beginpicture
\setcoordinatesystem units   <1.5mm,2mm>
\setplotarea x from 0 to 16, y from -2 to 15
\put{990)} [l] at 2 12
\put {$ \scriptstyle \bullet$} [c] at 10 12
\put {$ \scriptstyle \bullet$} [c] at 11 9
\put {$ \scriptstyle \bullet$} [c] at 16 6
\put {$ \scriptstyle \bullet$} [c] at 16 12
\put {$ \scriptstyle \bullet$} [c] at 12 6
\put {$ \scriptstyle \bullet$} [c] at 14 12
\put {$ \scriptstyle \bullet$} [c] at 14 0
\setlinear \plot 10 12  12 6 14 0 16 6 16 12   /
\setlinear \plot  14 12 12 6   /
\put{$5{,}040$} [c] at  13 -2
\endpicture
\end{minipage}
$$
$$
\begin{minipage}{4cm}
\beginpicture
\setcoordinatesystem units   <1.5mm,2mm>
\setplotarea x from 0 to 16, y from -2 to 15
\put{991)} [l] at 2 12
\put {$ \scriptstyle \bullet$} [c] at 10 0
\put {$ \scriptstyle \bullet$} [c] at 11 3
\put {$ \scriptstyle \bullet$} [c] at 16 6
\put {$ \scriptstyle \bullet$} [c] at 16 0
\put {$ \scriptstyle \bullet$} [c] at 12 6
\put {$ \scriptstyle \bullet$} [c] at 14 0
\put {$ \scriptstyle \bullet$} [c] at 14 12
\setlinear \plot 10 0  12 6 14 12 16 6 16 0   /
\setlinear \plot  14 0 12 6   /
\put{$5{,}040$} [c] at  13 -2
\endpicture
\end{minipage}
\begin{minipage}{4cm}
\beginpicture
\setcoordinatesystem units   <1.5mm,2mm>
\setplotarea x from 0 to 16, y from -2 to 15
\put{992)} [l] at 2 12
\put {$ \scriptstyle \bullet$} [c] at  10 4
\put {$ \scriptstyle \bullet$} [c] at  10 8
\put {$ \scriptstyle \bullet$} [c] at  10 12
\put {$ \scriptstyle \bullet$} [c] at  12 0
\put {$ \scriptstyle \bullet$} [c] at  14 4
\put {$ \scriptstyle \bullet$} [c] at  14 12
\put {$ \scriptstyle \bullet$} [c] at  16 0
\setlinear \plot 10 12 10 4 12 0 14 4 14 12 16 0    /
\put{$5{,}040   $} [c] at 13 -2
\endpicture
\end{minipage}
\begin{minipage}{4cm}
\beginpicture
\setcoordinatesystem units   <1.5mm,2mm>
\setplotarea x from 0 to 16, y from -2 to 15
\put{993)} [l] at 2 12
\put {$ \scriptstyle \bullet$} [c] at  10 4
\put {$ \scriptstyle \bullet$} [c] at  10 8
\put {$ \scriptstyle \bullet$} [c] at  10 0
\put {$ \scriptstyle \bullet$} [c] at  12 12
\put {$ \scriptstyle \bullet$} [c] at  14 0
\put {$ \scriptstyle \bullet$} [c] at  14 8
\put {$ \scriptstyle \bullet$} [c] at  16 12
\setlinear \plot 10 0 10 8 12 12 14 8 14 0 16 12    /
\put{$5{,}040   $} [c] at 13 -2
\endpicture
\end{minipage}
\begin{minipage}{4cm}
\beginpicture
\setcoordinatesystem units   <1.5mm,2mm>
\setplotarea x from 0 to 16, y from -2 to 15
\put{994)} [l] at 2 12
\put {$ \scriptstyle \bullet$} [c] at  10 0
\put {$ \scriptstyle \bullet$} [c] at  14 6
\put {$ \scriptstyle \bullet$} [c] at  14.5 9
\put {$ \scriptstyle \bullet$} [c] at  15 0
\put {$ \scriptstyle \bullet$} [c] at  15 12
\put {$ \scriptstyle \bullet$} [c] at  16 6
\put {$ \scriptstyle \bullet$} [c] at  16 12
\setlinear \plot 10 0 15 12 14 6 15 0 16 6 16 12      /
\setlinear \plot 16 6 15 12  /
\put{$5{,}040  $} [c] at 13 -2
\endpicture
\end{minipage}
\begin{minipage}{4cm}
\beginpicture
\setcoordinatesystem units   <1.5mm,2mm>
\setplotarea x from 0 to 16, y from -2 to 15
\put{995)} [l] at 2 12
\put {$ \scriptstyle \bullet$} [c] at  10 12
\put {$ \scriptstyle \bullet$} [c] at  14 6
\put {$ \scriptstyle \bullet$} [c] at  14.5 3
\put {$ \scriptstyle \bullet$} [c] at  15 0
\put {$ \scriptstyle \bullet$} [c] at  15 12
\put {$ \scriptstyle \bullet$} [c] at  16 6
\put {$ \scriptstyle \bullet$} [c] at  16 0
\setlinear \plot 10 12 15 0 14 6 15 12 16 6 16 0      /
\setlinear \plot 16 6 15 0  /
\put{$5{,}040   $} [c] at 13 -2
\endpicture
\end{minipage}
\begin{minipage}{4cm}
\beginpicture
\setcoordinatesystem units   <1.5mm,2mm>
\setplotarea x from 0 to 16, y from -2 to 15
\put{996)} [l] at 2 12
\put {$ \scriptstyle \bullet$} [c] at  10 4
\put {$ \scriptstyle \bullet$} [c] at  10 8
\put {$ \scriptstyle \bullet$} [c] at  11 0
\put {$ \scriptstyle \bullet$} [c] at  11 12
\put {$ \scriptstyle \bullet$} [c] at  12 6
\put {$ \scriptstyle \bullet$} [c] at  16 12
\put {$ \scriptstyle \bullet$} [c] at  16 0
\setlinear \plot 11 12 16 0 16 12 11 0 12 6 11 12  10 8 10 4 11 0    /
\put{$5{,}040  $} [c] at 13 -2
\endpicture
\end{minipage}
$$
$$
\begin{minipage}{4cm}
\beginpicture
\setcoordinatesystem units   <1.5mm,2mm>
\setplotarea x from 0 to 16, y from -2 to 15
\put{997)} [l] at 2 12
\put {$ \scriptstyle \bullet$} [c] at  10  0
\put {$ \scriptstyle \bullet$} [c] at  10 6
\put {$ \scriptstyle \bullet$} [c] at  10 12
\put {$ \scriptstyle \bullet$} [c] at  14.5 3
\put {$ \scriptstyle \bullet$} [c] at  16 0
\put {$ \scriptstyle \bullet$} [c] at  16 6
\put {$ \scriptstyle \bullet$} [c] at  16 12
\setlinear \plot  10 0 10 12 16 0 16 12 10  0    /
\put{$5{,}040$} [c] at 13 -2
\endpicture
\end{minipage}
\begin{minipage}{4cm}
\beginpicture
\setcoordinatesystem units   <1.5mm,2mm>
\setplotarea x from 0 to 16, y from -2 to 15
\put{998)} [l] at 2 12
\put {$ \scriptstyle \bullet$} [c] at  10 0
\put {$ \scriptstyle \bullet$} [c] at  10 6
\put {$ \scriptstyle \bullet$} [c] at  13 0
\put {$ \scriptstyle \bullet$} [c] at  13 6
\put {$ \scriptstyle \bullet$} [c] at  13 12
\put {$ \scriptstyle \bullet$} [c] at  16 6
\put {$ \scriptstyle \bullet$} [c] at  16 12
\setlinear \plot 16 12 16 6 13 12 13 0 16 6     /
\setlinear \plot  10 0 10 6 13 12  /
\setlinear \plot  13 0 10 6   /
\put{$5{,}040   $} [c] at 13 -2
\endpicture
\end{minipage}
\begin{minipage}{4cm}
\beginpicture
\setcoordinatesystem units   <1.5mm,2mm>
\setplotarea x from 0 to 16, y from -2 to 15
\put{999)} [l] at 2 12
\put {$ \scriptstyle \bullet$} [c] at  10 0
\put {$ \scriptstyle \bullet$} [c] at  10 12
\put {$ \scriptstyle \bullet$} [c] at  12.3 6
\put {$ \scriptstyle \bullet$} [c] at  13 0
\put {$ \scriptstyle \bullet$} [c] at  13 12
\put {$ \scriptstyle \bullet$} [c] at  13.6 6
\put {$ \scriptstyle \bullet$} [c] at  16 12
\setlinear \plot 13 0 16 12 10 0 10 12 12.3 6 13 0  13.6 6  13 12  12.3 6     /
\setlinear \plot  10 0  13 12   /
\put{$5{,}040$} [c] at 13 -2
\endpicture
\end{minipage}
\begin{minipage}{4cm}
\beginpicture
\setcoordinatesystem units   <1.5mm,2mm>
\setplotarea x from 0 to 16, y from -2 to 15
\put{1.000)} [l] at 2 12
\put {$ \scriptstyle \bullet$} [c] at  10 0
\put {$ \scriptstyle \bullet$} [c] at  10 12
\put {$ \scriptstyle \bullet$} [c] at  12.3 6
\put {$ \scriptstyle \bullet$} [c] at  13 0
\put {$ \scriptstyle \bullet$} [c] at  13 12
\put {$ \scriptstyle \bullet$} [c] at  13.6 6
\put {$ \scriptstyle \bullet$} [c] at  16 0
\setlinear \plot 13 12 16 0 10 12 10 0 12.3 6 13 12  13.6 6  13 0  12.3 6     /
\setlinear \plot  10 12  13 0   /
\put{$5{,}040$} [c] at 13 -2
\endpicture
\end{minipage}
\begin{minipage}{4cm}
\beginpicture
\setcoordinatesystem units   <1.5mm,2mm>
\setplotarea x from 0 to 16, y from -2 to 15
\put{1.001)} [l] at 2 12
\put {$ \scriptstyle \bullet$} [c] at  10 0
\put {$ \scriptstyle \bullet$} [c] at  10 6
\put {$ \scriptstyle \bullet$} [c] at  10 12
\put {$ \scriptstyle \bullet$} [c] at  11 0
\put {$ \scriptstyle \bullet$} [c] at  11 12
\put {$ \scriptstyle \bullet$} [c] at  12 6
\put {$ \scriptstyle \bullet$} [c] at  16 12
\setlinear \plot 16 12 11 0  10 6 11 12 12 6 11 0   /
\setlinear \plot  10 12  10 0  /
\put{$5{,}040$} [c] at 13 -2
\endpicture
\end{minipage}
\begin{minipage}{4cm}
\beginpicture
\setcoordinatesystem units   <1.5mm,2mm>
\setplotarea x from 0 to 16, y from -2 to 15
\put{1.002)} [l] at 2 12
\put {$ \scriptstyle \bullet$} [c] at  10 0
\put {$ \scriptstyle \bullet$} [c] at  10 6
\put {$ \scriptstyle \bullet$} [c] at  10 12
\put {$ \scriptstyle \bullet$} [c] at  11 0
\put {$ \scriptstyle \bullet$} [c] at  11 12
\put {$ \scriptstyle \bullet$} [c] at  12 6
\put {$ \scriptstyle \bullet$} [c] at  16 0
\setlinear \plot 16 0 11 12  10 6 11 0 12 6 11 12   /
\setlinear \plot  10 12  10 0  /
\put{$5{,}040$} [c] at 13 -2
\endpicture
\end{minipage}
$$
$$
\begin{minipage}{4cm}
\beginpicture
\setcoordinatesystem units   <1.5mm,2mm>
\setplotarea x from 0 to 16, y from -2 to 15
\put{1.003)} [l] at 2 12
\put {$ \scriptstyle \bullet$} [c] at  10 12
\put {$ \scriptstyle \bullet$} [c] at  12 0
\put {$ \scriptstyle \bullet$} [c] at  12 6
\put {$ \scriptstyle \bullet$} [c] at  12 12
\put {$ \scriptstyle \bullet$} [c] at  16 0
\put {$ \scriptstyle \bullet$} [c] at  16 6
\put {$ \scriptstyle \bullet$} [c] at  16 12
\setlinear \plot 10 12 12 0 12 12 16 0 16  12 12 6   /
\setlinear \plot 12 0 16 6    /
\put{$5{,}040$} [c] at 13 -2
\endpicture
\end{minipage}
\begin{minipage}{4cm}
\beginpicture
\setcoordinatesystem units   <1.5mm,2mm>
\setplotarea x from 0 to 16, y from -2 to 15
\put{1.004)} [l] at 2 12
\put {$ \scriptstyle \bullet$} [c] at  10 0
\put {$ \scriptstyle \bullet$} [c] at  12 0
\put {$ \scriptstyle \bullet$} [c] at  12 6
\put {$ \scriptstyle \bullet$} [c] at  12 12
\put {$ \scriptstyle \bullet$} [c] at  16 0
\put {$ \scriptstyle \bullet$} [c] at  16 6
\put {$ \scriptstyle \bullet$} [c] at  16 12
\setlinear \plot 10 0 12 12 12 0 16 12 16  0 12 6   /
\setlinear \plot 12 12 16 6    /
\put{$5{,}040$} [c] at 13 -2
\endpicture
\end{minipage}
\begin{minipage}{4cm}
\beginpicture
\setcoordinatesystem units   <1.5mm,2mm>
\setplotarea x from 0 to 16, y from -2 to 15
\put{1.005)} [l] at 2 12
\put {$ \scriptstyle \bullet$} [c] at  10 12
\put {$ \scriptstyle \bullet$} [c] at  12 0
\put {$ \scriptstyle \bullet$} [c] at  12 4
\put {$ \scriptstyle \bullet$} [c] at  12 8
\put {$ \scriptstyle \bullet$} [c] at  12 12
\put {$ \scriptstyle \bullet$} [c] at  16 12
\put {$ \scriptstyle \bullet$} [c] at  16  0
\setlinear \plot 10 12 12 0 12 12  16 0 16 12 12 4 /
\put{$5{,}040$} [c] at 13 -2
\endpicture
\end{minipage}
\begin{minipage}{4cm}
\beginpicture
\setcoordinatesystem units   <1.5mm,2mm>
\setplotarea x from 0 to 16, y from -2 to 15
\put{1.006)} [l] at 2 12
\put {$ \scriptstyle \bullet$} [c] at  10 0
\put {$ \scriptstyle \bullet$} [c] at  12 0
\put {$ \scriptstyle \bullet$} [c] at  12 4
\put {$ \scriptstyle \bullet$} [c] at  12 8
\put {$ \scriptstyle \bullet$} [c] at  12 12
\put {$ \scriptstyle \bullet$} [c] at  16 12
\put {$ \scriptstyle \bullet$} [c] at  16  0
\setlinear \plot 10 0 12 12 12 0  16 12 16 0 12 8  /
\put{$5{,}040$} [c] at 13 -2
\endpicture
\end{minipage}
\begin{minipage}{4cm}
\beginpicture
\setcoordinatesystem units   <1.5mm,2mm>
\setplotarea x from 0 to 16, y from -2 to 15
\put{1.007)} [l] at 2 12
\put{$ \scriptstyle \bullet$} [c] at  10 0
\put {$ \scriptstyle \bullet$} [c] at  10 6
\put {$ \scriptstyle \bullet$} [c] at  10 12
\put {$ \scriptstyle \bullet$} [c] at  13 12
\put {$ \scriptstyle \bullet$} [c] at  13 3
\put {$ \scriptstyle \bullet$} [c] at  16 0
\put {$ \scriptstyle \bullet$} [c] at  16 12
\setlinear \plot 10  12 10 0   /
\setlinear \plot  16 12 16 0 10 6  /
\setlinear \plot  13 3 13 12   /
\put{$5{,}040$} [c] at 13 -2
\endpicture
\end{minipage}
\begin{minipage}{4cm}
\beginpicture
\setcoordinatesystem units   <1.5mm,2mm>
\setplotarea x from 0 to 16, y from -2 to 15
\put{1.008)} [l] at 2 12
\put{$ \scriptstyle \bullet$} [c] at  10 0
\put {$ \scriptstyle \bullet$} [c] at  10 6
\put {$ \scriptstyle \bullet$} [c] at  10 12
\put {$ \scriptstyle \bullet$} [c] at  13 0
\put {$ \scriptstyle \bullet$} [c] at  13 9
\put {$ \scriptstyle \bullet$} [c] at  16 0
\put {$ \scriptstyle \bullet$} [c] at  16 12
\setlinear \plot 10  12 10 0   /
\setlinear \plot  16 0 16 12 10 6  /
\setlinear \plot  13 9 13 0   /
\put{$5{,}040$} [c] at 13 -2
\endpicture
\end{minipage}
$$
$$
\begin{minipage}{4cm}
\beginpicture
\setcoordinatesystem units   <1.5mm,2mm>
\setplotarea x from 0 to 16, y from -2 to 15
\put{1.009)} [l] at 2 12
\put {$ \scriptstyle \bullet$} [c] at  10 12
\put {$ \scriptstyle \bullet$} [c] at  13 0
\put {$ \scriptstyle \bullet$} [c] at  13 6
\put {$ \scriptstyle \bullet$} [c] at  13 12
\put {$ \scriptstyle \bullet$} [c] at  16 0
\put {$ \scriptstyle \bullet$} [c] at  16 6
\put {$ \scriptstyle \bullet$} [c] at  16 12
\setlinear \plot  13 0  13 12  16 0 16 12 /
\setlinear \plot  10 12 13 6  /
\put{$5{,}040$} [c] at 13 -2
\endpicture
\end{minipage}
\begin{minipage}{4cm}
\beginpicture
\setcoordinatesystem units   <1.5mm,2mm>
\setplotarea x from 0 to 16, y from -2 to 15
\put{1.010)} [l] at 2 12
\put {$ \scriptstyle \bullet$} [c] at  10 0
\put {$ \scriptstyle \bullet$} [c] at  13 0
\put {$ \scriptstyle \bullet$} [c] at  13 6
\put {$ \scriptstyle \bullet$} [c] at  13 12
\put {$ \scriptstyle \bullet$} [c] at  16 0
\put {$ \scriptstyle \bullet$} [c] at  16 6
\put {$ \scriptstyle \bullet$} [c] at  16 12
\setlinear \plot  16 0  16 12  13 0 13 12 /
\setlinear \plot  10 0 13 6  /
\put{$5{,}040$} [c] at 13 -2
\endpicture
\end{minipage}
\begin{minipage}{4cm}
\beginpicture
\setcoordinatesystem units   <1.5mm,2mm>
\setplotarea x from 0 to 16, y from -2 to 15
\put{1.011)} [l] at 2 12
\put {$ \scriptstyle \bullet$} [c] at  10 0
\put {$ \scriptstyle \bullet$} [c] at  10 4
\put {$ \scriptstyle \bullet$} [c] at  10 8
\put {$ \scriptstyle \bullet$} [c] at  10 12
\put {$ \scriptstyle \bullet$} [c] at  13 0
\put {$ \scriptstyle \bullet$} [c] at  13 12
\put {$ \scriptstyle \bullet$} [c] at  16 12
\setlinear \plot  16 12  13  0 13 12 10 0 10 12 13 0 /
\put{$5{,}040$} [c] at 13 -2
\endpicture
\end{minipage}
\begin{minipage}{4cm}
\beginpicture
\setcoordinatesystem units   <1.5mm,2mm>
\setplotarea x from 0 to 16, y from -2 to 15
\put{1.012)} [l] at 2 12
\put {$ \scriptstyle \bullet$} [c] at  10 0
\put {$ \scriptstyle \bullet$} [c] at  10 4
\put {$ \scriptstyle \bullet$} [c] at  10 8
\put {$ \scriptstyle \bullet$} [c] at  10 12
\put {$ \scriptstyle \bullet$} [c] at  13 0
\put {$ \scriptstyle \bullet$} [c] at  13 12
\put {$ \scriptstyle \bullet$} [c] at  16 0
\setlinear \plot  16 0  13  12 13 0 10 12 10 0 13 12 /
\put{$5{,}040$} [c] at 13 -2
\endpicture
\end{minipage}
\begin{minipage}{4cm}
\beginpicture
\setcoordinatesystem units   <1.5mm,2mm>
\setplotarea x from 0 to 16, y from -2 to 15
\put{1.013)} [l] at 2 12
\put {$ \scriptstyle \bullet$} [c] at  10 6
\put {$ \scriptstyle \bullet$} [c] at  10 12
\put {$ \scriptstyle \bullet$} [c] at  12 0
\put {$ \scriptstyle \bullet$} [c] at  12 12
\put {$ \scriptstyle \bullet$} [c] at  14 6
\put {$ \scriptstyle \bullet$} [c] at  16 0
\put {$ \scriptstyle \bullet$} [c] at  16 12
\setlinear \plot 16 12 16 0  14  6 12  0  10 6 12 12 14 6 /
\setlinear \plot 10 6  10 12 /
\put{$5{,}040$} [c] at 13 -2
\endpicture
\end{minipage}
\begin{minipage}{4cm}
\beginpicture
\setcoordinatesystem units   <1.5mm,2mm>
\setplotarea x from 0 to 16, y from -2 to 15
\put{1.014)} [l] at 2 12
\put {$ \scriptstyle \bullet$} [c] at  10 6
\put {$ \scriptstyle \bullet$} [c] at  10 0
\put {$ \scriptstyle \bullet$} [c] at  12 0
\put {$ \scriptstyle \bullet$} [c] at  12 12
\put {$ \scriptstyle \bullet$} [c] at  14 6
\put {$ \scriptstyle \bullet$} [c] at  16 0
\put {$ \scriptstyle \bullet$} [c] at  16 12
\setlinear \plot 16 0 16 12  14  6 12  0  10 6 12 12 14 6 /
\setlinear \plot 10 6  10 0 /
\put{$5{,}040$} [c] at 13 -2
\endpicture
\end{minipage}
$$
$$
\begin{minipage}{4cm}
\beginpicture
\setcoordinatesystem units   <1.5mm,2mm>
\setplotarea x from 0 to 16, y from -2 to 15
\put{1.015)} [l] at 2 12
\put {$ \scriptstyle \bullet$} [c] at  13 0
\put {$ \scriptstyle \bullet$} [c] at  13 4
\put {$ \scriptstyle \bullet$} [c] at  13 8
\put {$ \scriptstyle \bullet$} [c] at  13 12
\put {$ \scriptstyle \bullet$} [c] at  16 0
\put {$ \scriptstyle \bullet$} [c] at  16 12
\put {$ \scriptstyle \bullet$} [c] at  10 12
\setlinear \plot 10 12 13  0  16 12 16 0 13 8 13 0   /
\setlinear \plot  13 8 13 12     /
\put{$5{,}040$} [c] at 13 -2
\endpicture
\end{minipage}
\begin{minipage}{4cm}
\beginpicture
\setcoordinatesystem units   <1.5mm,2mm>
\setplotarea x from 0 to 16, y from -2 to 15
\put{1.016)} [l] at 2 12
\put {$ \scriptstyle \bullet$} [c] at  13 0
\put {$ \scriptstyle \bullet$} [c] at  13 4
\put {$ \scriptstyle \bullet$} [c] at  13 8
\put {$ \scriptstyle \bullet$} [c] at  13 12
\put {$ \scriptstyle \bullet$} [c] at  16 0
\put {$ \scriptstyle \bullet$} [c] at  16 12
\put {$ \scriptstyle \bullet$} [c] at  10 0
\setlinear \plot 10 0 13  12  16 0 16 12 13 4 13 0   /
\setlinear \plot  13 4 13 12     /
\put{$5{,}040$} [c] at 13 -2
\endpicture
\end{minipage}
\begin{minipage}{4cm}
\beginpicture
\setcoordinatesystem units   <1.5mm,2mm>
\setplotarea x from 0 to 16, y from -2 to 15
\put{1.017)} [l] at 2 12
\put {$ \scriptstyle \bullet$} [c] at  13 0
\put {$ \scriptstyle \bullet$} [c] at  13 4
\put {$ \scriptstyle \bullet$} [c] at  13 8
\put {$ \scriptstyle \bullet$} [c] at  13 12
\put {$ \scriptstyle \bullet$} [c] at  16 0
\put {$ \scriptstyle \bullet$} [c] at  16 12
\put {$ \scriptstyle \bullet$} [c] at  10 12
\setlinear \plot 16 0 16 12  13 8 13 0 10 12   /
\setlinear \plot  13 12 13 8    /
\put{$5{,}040$} [c] at 13 -2
\endpicture
\end{minipage}
\begin{minipage}{4cm}
\beginpicture
\setcoordinatesystem units   <1.5mm,2mm>
\setplotarea x from 0 to 16, y from -2 to 15
\put{1.018)} [l] at 2 12
\put {$ \scriptstyle \bullet$} [c] at  13 0
\put {$ \scriptstyle \bullet$} [c] at  13 4
\put {$ \scriptstyle \bullet$} [c] at  13 8
\put {$ \scriptstyle \bullet$} [c] at  13 12
\put {$ \scriptstyle \bullet$} [c] at  16 0
\put {$ \scriptstyle \bullet$} [c] at  16 12
\put {$ \scriptstyle \bullet$} [c] at  10 0
\setlinear \plot 16 12 16 0  13 4 13 12 10 0   /
\setlinear \plot  13 0 13 4    /
\put{$5{,}040$} [c] at 13 -2
\endpicture
\end{minipage}
\begin{minipage}{4cm}
\beginpicture
\setcoordinatesystem units   <1.5mm,2mm>
\setplotarea x from 0 to 16, y from -2 to 15
\put{1.019)} [l] at 2 12
\put {$ \scriptstyle \bullet$} [c] at  10 12
\put {$ \scriptstyle \bullet$} [c] at  13 12
\put {$ \scriptstyle \bullet$} [c] at  13 8
\put {$ \scriptstyle \bullet$} [c] at  13 4
\put {$ \scriptstyle \bullet$} [c] at  13 0
\put {$ \scriptstyle \bullet$} [c] at  16 0
\put {$ \scriptstyle \bullet$} [c] at  16 12
\setlinear \plot  16 12  16 0 13 12 13 0  /
\setlinear \plot  13 4 10 12  /
\put{$5{,}040$} [c] at 13 -2
\endpicture
\end{minipage}
\begin{minipage}{4cm}
\beginpicture
\setcoordinatesystem units   <1.5mm,2mm>
\setplotarea x from 0 to 16, y from -2 to 15
\put{1.020)} [l] at 2 12
\put {$ \scriptstyle \bullet$} [c] at  10 0
\put {$ \scriptstyle \bullet$} [c] at  13 12
\put {$ \scriptstyle \bullet$} [c] at  13 8
\put {$ \scriptstyle \bullet$} [c] at  13 4
\put {$ \scriptstyle \bullet$} [c] at  13 0
\put {$ \scriptstyle \bullet$} [c] at  16 0
\put {$ \scriptstyle \bullet$} [c] at  16 12
\setlinear \plot  16 0  16 12 13 0 13 12  /
\setlinear \plot  13 8 10 0  /
\put{$5{,}040$} [c] at 13 -2
\endpicture
\end{minipage}
$$
$$
\begin{minipage}{4cm}
\beginpicture
\setcoordinatesystem units   <1.5mm,2mm>
\setplotarea x from 0 to 16, y from -2 to 15
\put{1.021)} [l] at 2 12
\put {$ \scriptstyle \bullet$} [c] at  10 0
\put {$ \scriptstyle \bullet$} [c] at  10 4
\put {$ \scriptstyle \bullet$} [c] at  10 8
\put {$ \scriptstyle \bullet$} [c] at  10 12
\put {$ \scriptstyle \bullet$} [c] at  13 12
\put {$ \scriptstyle \bullet$} [c] at  16 0
\put {$ \scriptstyle \bullet$} [c] at  16 12
\setlinear \plot  10 0 10 12  /
\setlinear \plot  16 0 16 12 10 4 13 12  /
\put{$5{,}040$} [c] at 13 -2
\endpicture
\end{minipage}
\begin{minipage}{4cm}
\beginpicture
\setcoordinatesystem units   <1.5mm,2mm>
\setplotarea x from 0 to 16, y from -2 to 15
\put{1.022)} [l] at 2 12
\put {$ \scriptstyle \bullet$} [c] at  10 0
\put {$ \scriptstyle \bullet$} [c] at  10 4
\put {$ \scriptstyle \bullet$} [c] at  10 8
\put {$ \scriptstyle \bullet$} [c] at  10 12
\put {$ \scriptstyle \bullet$} [c] at  13 0
\put {$ \scriptstyle \bullet$} [c] at  16 0
\put {$ \scriptstyle \bullet$} [c] at  16 12
\setlinear \plot  10 0 10 12  /
\setlinear \plot  16 12 16 0 10 8 13 0  /
\put{$5{,}040$} [c] at 13 -2
\endpicture
\end{minipage}
\begin{minipage}{4cm}
\beginpicture
\setcoordinatesystem units   <1.5mm,2mm>
\setplotarea x from 0 to 16, y from -2 to 15
\put{1.023)} [l] at 2 12
\put {$ \scriptstyle \bullet$} [c] at  10 12
\put {$ \scriptstyle \bullet$} [c] at  12 6
\put {$ \scriptstyle \bullet$} [c] at  12 12
\put {$ \scriptstyle \bullet$} [c] at  13 0
\put {$ \scriptstyle \bullet$} [c] at  14 6
\put {$ \scriptstyle \bullet$} [c] at  14 12
\put {$ \scriptstyle \bullet$} [c] at  16 0
\setlinear \plot  10 12 12 6 12 12 14 6 14 12 16 0 /
\setlinear \plot 14 12 12 6 13 0 14 6  /
\put{$5{,}040$} [c] at 13 -2
\endpicture
\end{minipage}
\begin{minipage}{4cm}
\beginpicture
\setcoordinatesystem units   <1.5mm,2mm>
\setplotarea x from 0 to 16, y from -2 to 15
\put{1.024)} [l] at 2 12
\put {$ \scriptstyle \bullet$} [c] at  10 0
\put {$ \scriptstyle \bullet$} [c] at  12 6
\put {$ \scriptstyle \bullet$} [c] at  12 0
\put {$ \scriptstyle \bullet$} [c] at  13 12
\put {$ \scriptstyle \bullet$} [c] at  14 6
\put {$ \scriptstyle \bullet$} [c] at  14 0
\put {$ \scriptstyle \bullet$} [c] at  16 12
\setlinear \plot  10 0 12 6 12 0 14 6 14 0 16 12 /
\setlinear \plot 14 0 12 6 13 12 14 6  /
\put{$5{,}040$} [c] at 13 -2
\endpicture
\end{minipage}
\begin{minipage}{4cm}
\beginpicture
\setcoordinatesystem units   <1.5mm,2mm>
\setplotarea x from 0 to 16, y from -2 to 15
\put{1.025)} [l] at 2 12
\put {$ \scriptstyle \bullet$} [c] at  10 0
\put {$ \scriptstyle \bullet$} [c] at  10 12
\put {$ \scriptstyle \bullet$} [c] at  14 12
\put {$ \scriptstyle \bullet$} [c] at  14 6
\put {$ \scriptstyle \bullet$} [c] at  15 0
\put {$ \scriptstyle \bullet$} [c] at  15 12
\put {$ \scriptstyle \bullet$} [c] at  16 6
\setlinear \plot 14 12 14 6  10 12 10 0 15 12 16 6  15 0 14 6 15 12   /
\put{$5{,}040$} [c] at 13 -2
\endpicture
\end{minipage}
\begin{minipage}{4cm}
\beginpicture
\setcoordinatesystem units   <1.5mm,2mm>
\setplotarea x from 0 to 16, y from -2 to 15
\put{1.026)} [l] at 2 12
\put {$ \scriptstyle \bullet$} [c] at  10 0
\put {$ \scriptstyle \bullet$} [c] at  10 12
\put {$ \scriptstyle \bullet$} [c] at  14 0
\put {$ \scriptstyle \bullet$} [c] at  14 6
\put {$ \scriptstyle \bullet$} [c] at  15 0
\put {$ \scriptstyle \bullet$} [c] at  15 12
\put {$ \scriptstyle \bullet$} [c] at  16 6
\setlinear \plot 14 0 14 6  10 0 10 12 15 0 16 6  15 12 14 6 15 0   /
\put{$5{,}040$} [c] at 13 -2
\endpicture
\end{minipage}
$$
$$
\begin{minipage}{4cm}
\beginpicture
\setcoordinatesystem units   <1.5mm,2mm>
\setplotarea x from 0 to 16, y from -2 to 15
\put{1.027)} [l] at 2 12
\put {$ \scriptstyle \bullet$} [c] at  10 0
\put {$ \scriptstyle \bullet$} [c] at  10 6
\put {$ \scriptstyle \bullet$} [c] at  10 12
\put {$ \scriptstyle \bullet$} [c] at  13 12
\put {$ \scriptstyle \bullet$} [c] at  16 0
\put {$ \scriptstyle \bullet$} [c] at  16 6
\put {$ \scriptstyle \bullet$} [c] at  16 12
\setlinear \plot  10 12  10 0 13 12 16 0 16 12   /
\setlinear \plot   10 12 16  6   /
\put{$5{,}040$} [c] at 13 -2
\endpicture
\end{minipage}
\begin{minipage}{4cm}
\beginpicture
\setcoordinatesystem units   <1.5mm,2mm>
\setplotarea x from 0 to 16, y from -2 to 15
\put{1.028)} [l] at 2 12
\put {$ \scriptstyle \bullet$} [c] at  10 0
\put {$ \scriptstyle \bullet$} [c] at  10 6
\put {$ \scriptstyle \bullet$} [c] at  10 12
\put {$ \scriptstyle \bullet$} [c] at  13 0
\put {$ \scriptstyle \bullet$} [c] at  16 0
\put {$ \scriptstyle \bullet$} [c] at  16 6
\put {$ \scriptstyle \bullet$} [c] at  16 12
\setlinear \plot  10 0  10 12 13 0 16 12 16 0   /
\setlinear \plot   10 0 16  6   /
\put{$5{,}040$} [c] at 13 -2
\endpicture
\end{minipage}
\begin{minipage}{4cm}
\beginpicture
\setcoordinatesystem units   <1.5mm,2mm>
\setplotarea x from 0 to 16, y from -2 to 15
\put{1.029)} [l] at  2 12
\put {$ \scriptstyle \bullet$} [c] at  10 0
\put {$ \scriptstyle \bullet$} [c] at  10 6
\put {$ \scriptstyle \bullet$} [c] at  10 12
\put {$ \scriptstyle \bullet$} [c] at  13 12
\put {$ \scriptstyle \bullet$} [c] at  16 0
\put {$ \scriptstyle \bullet$} [c] at  16 6
\put {$ \scriptstyle \bullet$} [c] at  16 12
\setlinear \plot 10 12 10 0 13 12 16 0 16 12 10 0  /
\setlinear \plot 10 6 16 0   /
\put{$5{,}040$} [c] at 13 -2
\endpicture
\end{minipage}
\begin{minipage}{4cm}
\beginpicture
\setcoordinatesystem units   <1.5mm,2mm>
\setplotarea x from 0 to 16, y from -2 to 15
\put{1.030)} [l] at  2 12
\put {$ \scriptstyle \bullet$} [c] at  10 0
\put {$ \scriptstyle \bullet$} [c] at  10 6
\put {$ \scriptstyle \bullet$} [c] at  10 12
\put {$ \scriptstyle \bullet$} [c] at  13 0
\put {$ \scriptstyle \bullet$} [c] at  16 0
\put {$ \scriptstyle \bullet$} [c] at  16 6
\put {$ \scriptstyle \bullet$} [c] at  16 12
\setlinear \plot 10 0 10 12 13 0 16 12 16 0 10 12  /
\setlinear \plot 10 6 16 12   /
\put{$5{,}040$} [c] at 13 -2
\endpicture
\end{minipage}
\begin{minipage}{4cm}
\beginpicture
\setcoordinatesystem units   <1.5mm,2mm>
\setplotarea x from 0 to 16, y from -2 to 15
\put{1.031)} [l] at 2 12
\put {$ \scriptstyle \bullet$} [c] at  10 0
\put {$ \scriptstyle \bullet$} [c] at  10 12
\put {$ \scriptstyle \bullet$} [c] at  13 12
\put {$ \scriptstyle \bullet$} [c] at  13 0
\put {$ \scriptstyle \bullet$} [c] at  16 0
\put {$ \scriptstyle \bullet$} [c] at  16 6
\put {$ \scriptstyle \bullet$} [c] at  16  12
\setlinear \plot  10 0 10 12 13 0 13 12 10 0   /
\setlinear \plot  13 12 16 6 16 0  /
\setlinear \plot 13 0 16 12 16 6 /
\put{$5{,}040 $} [c] at 13 -2
\endpicture
\end{minipage}
\begin{minipage}{4cm}
\beginpicture
\setcoordinatesystem units   <1.5mm,2mm>
\setplotarea x from 0 to 16, y from -2 to 15
\put{1.032)} [l] at 2 12
\put {$ \scriptstyle \bullet$} [c] at  10 0
\put {$ \scriptstyle \bullet$} [c] at  10 12
\put {$ \scriptstyle \bullet$} [c] at  13 12
\put {$ \scriptstyle \bullet$} [c] at  13 0
\put {$ \scriptstyle \bullet$} [c] at  16 0
\put {$ \scriptstyle \bullet$} [c] at  16 6
\put {$ \scriptstyle \bullet$} [c] at  16  12
\setlinear \plot  10 12 10 0 13 12 13 0 10 12   /
\setlinear \plot  13 0 16 6 16 12  /
\setlinear \plot 13 12 16 0 16 6 /
\put{$5{,}040 $} [c] at 13 -2
\endpicture
\end{minipage}
$$
$$
\begin{minipage}{4cm}
\beginpicture
\setcoordinatesystem units   <1.5mm,2mm>
\setplotarea x from 0 to 16, y from -2 to 15
\put{1.033)} [l] at 2 12
\put {$ \scriptstyle \bullet$} [c] at  10 6
\put {$ \scriptstyle \bullet$} [c] at  12 0
\put {$ \scriptstyle \bullet$} [c] at  12 12
\put {$ \scriptstyle \bullet$} [c] at  14 6
\put {$ \scriptstyle \bullet$} [c] at  16 0
\put {$ \scriptstyle \bullet$} [c] at  16 6
\put {$ \scriptstyle \bullet$} [c] at  16 12
\setlinear \plot 12 0 10 6 12 12  14 6 12 0 16 6 16 12      /
\setlinear \plot  16 0 16 6  /
\put{$2{,}520  $} [c] at 13 -2
\endpicture
\end{minipage}
\begin{minipage}{4cm}
\beginpicture
\setcoordinatesystem units   <1.5mm,2mm>
\setplotarea x from 0 to 16, y from -2 to 15
\put{1.034)} [l] at 2 12
\put {$ \scriptstyle \bullet$} [c] at  10 6
\put {$ \scriptstyle \bullet$} [c] at  12 0
\put {$ \scriptstyle \bullet$} [c] at  12 12
\put {$ \scriptstyle \bullet$} [c] at  14 6
\put {$ \scriptstyle \bullet$} [c] at  16 0
\put {$ \scriptstyle \bullet$} [c] at  16 6
\put {$ \scriptstyle \bullet$} [c] at  16 12
\setlinear \plot 12 12 10 6 12 0  14 6 12 12 16 6 16 0      /
\setlinear \plot  16 12 16 6  /
\put{$2{,}520  $} [c] at 13 -2
\endpicture
\end{minipage}
\begin{minipage}{4cm}
\beginpicture
\setcoordinatesystem units   <1.5mm,2mm>
\setplotarea x from 0 to 16, y from -2 to 15
\put{1.035)} [l] at 2 12
\put {$ \scriptstyle \bullet$} [c] at  10 6
\put {$ \scriptstyle \bullet$} [c] at  11 0
\put {$ \scriptstyle \bullet$} [c] at  11 6
\put {$ \scriptstyle \bullet$} [c] at  11 12
\put {$ \scriptstyle \bullet$} [c] at  12 6
\put {$ \scriptstyle \bullet$} [c] at  16 0
\put {$ \scriptstyle \bullet$} [c] at  16 12
\setlinear \plot 11 12 11 0  10 6  11 12 12 6  11 0 16 12 16 0 12 6      /
\setlinear \plot   11 6 16 0   /
\put{$2{,}520  $} [c] at 13 -2
\endpicture
\end{minipage}
\begin{minipage}{4cm}
\beginpicture
\setcoordinatesystem units   <1.5mm,2mm>
\setplotarea x from 0 to 16, y from -2 to 15
\put{1.036)} [l] at 2 12
\put {$ \scriptstyle \bullet$} [c] at  10 6
\put {$ \scriptstyle \bullet$} [c] at  11 0
\put {$ \scriptstyle \bullet$} [c] at  11 6
\put {$ \scriptstyle \bullet$} [c] at  11 12
\put {$ \scriptstyle \bullet$} [c] at  12 6
\put {$ \scriptstyle \bullet$} [c] at  16 0
\put {$ \scriptstyle \bullet$} [c] at  16 12
\setlinear \plot 11 0 11 12  10 6  11 0 12 6  11 12 16 0 16 12 12 6      /
\setlinear \plot   11 6 16 12   /
\put{$2{,}520 $} [c] at 13 -2
\endpicture
\end{minipage}
\begin{minipage}{4cm}
\beginpicture
\setcoordinatesystem units   <1.5mm,2mm>
\setplotarea x from 0 to 16, y from -2 to 15
\put{1.037)} [l] at 2 12
\put {$ \scriptstyle \bullet$} [c] at  10 0
\put {$ \scriptstyle \bullet$} [c] at  12.5 4
\put {$ \scriptstyle \bullet$} [c] at  13 0
\put {$ \scriptstyle \bullet$} [c] at  13 10
\put {$ \scriptstyle \bullet$} [c] at  13 12
\put {$ \scriptstyle \bullet$} [c] at  13.5 4
\put {$ \scriptstyle \bullet$} [c] at  16 12
\setlinear \plot  10 0 13  10 13 12   /
\setlinear \plot  16 12 13 0 13.5 4 13 10 12.5 4 13 0   /
\put{$2{,}520$} [c] at 13 -2
\endpicture
\end{minipage}
\begin{minipage}{4cm}
\beginpicture
\setcoordinatesystem units   <1.5mm,2mm>
\setplotarea x from 0 to 16, y from -2 to 15
\put{1.038)} [l] at 2 12
\put {$ \scriptstyle \bullet$} [c] at  10 12
\put {$ \scriptstyle \bullet$} [c] at  12.5 8
\put {$ \scriptstyle \bullet$} [c] at  13 0
\put {$ \scriptstyle \bullet$} [c] at  13 2
\put {$ \scriptstyle \bullet$} [c] at  13 12
\put {$ \scriptstyle \bullet$} [c] at  13.5 8
\put {$ \scriptstyle \bullet$} [c] at  16 0
\setlinear \plot  10 12 13  2 13 0   /
\setlinear \plot  16 0 13 12 12.5 8 13 2 13.5 8 13 12   /
\put{$2{,}520$} [c] at 13 -2
\endpicture
\end{minipage}
$$
$$
\begin{minipage}{4cm}
\beginpicture
\setcoordinatesystem units   <1.5mm,2mm>
\setplotarea x from 0 to 16, y from -2 to 15
\put{1.039)} [l] at 2 12
\put {$ \scriptstyle \bullet$} [c] at  10 6
\put {$ \scriptstyle \bullet$} [c] at  11 0
\put {$ \scriptstyle \bullet$} [c] at  11 12
\put {$ \scriptstyle \bullet$} [c] at  12 6
\put {$ \scriptstyle \bullet$} [c] at  16 0
\put {$ \scriptstyle \bullet$} [c] at  16 6
\put {$ \scriptstyle \bullet$} [c] at  16  12
\setlinear \plot  11 12 16 0 16 12 11 0  12 6 11 12 10 6 11 0    /
\put{$2{,}520$} [c] at 13 -2
\endpicture
\end{minipage}
\begin{minipage}{4cm}
\beginpicture
\setcoordinatesystem units   <1.5mm,2mm>
\setplotarea x from 0 to 16, y from -2 to 15
\put{1.040)} [l] at 2 12
\put {$ \scriptstyle \bullet$} [c] at  10 12
\put {$ \scriptstyle \bullet$} [c] at  10 8
\put {$ \scriptstyle \bullet$} [c] at  10 4
\put {$ \scriptstyle \bullet$} [c] at  10 0
\put {$ \scriptstyle \bullet$} [c] at  13 12
\put {$ \scriptstyle \bullet$} [c] at  16 12
\put {$ \scriptstyle \bullet$} [c] at  16 0
\setlinear \plot  10 12 10 0  13 12 16 0 16 12 10 0  /
\put{$2{,}520$} [c] at 13 -2
\endpicture
\end{minipage}
\begin{minipage}{4cm}
\beginpicture
\setcoordinatesystem units   <1.5mm,2mm>
\setplotarea x from 0 to 16, y from -2 to 15
\put{1.041)} [l] at 2 12
\put {$ \scriptstyle \bullet$} [c] at  10 12
\put {$ \scriptstyle \bullet$} [c] at  10 8
\put {$ \scriptstyle \bullet$} [c] at  10 4
\put {$ \scriptstyle \bullet$} [c] at  10 0
\put {$ \scriptstyle \bullet$} [c] at  13 0
\put {$ \scriptstyle \bullet$} [c] at  16 12
\put {$ \scriptstyle \bullet$} [c] at  16 0
\setlinear \plot  10 0 10 12  13 0 16 12 16 0 10 12  /
\put{$2{,}520$} [c] at 13 -2
\endpicture
\end{minipage}
\begin{minipage}{4cm}
\beginpicture
\setcoordinatesystem units   <1.5mm,2mm>
\setplotarea x from 0 to 16, y from -2 to 15
\put{1.042)} [l] at 2 12
\put {$ \scriptstyle \bullet$} [c] at  10 0
\put {$ \scriptstyle \bullet$} [c] at  10 6
\put {$ \scriptstyle \bullet$} [c] at  10 12
\put {$ \scriptstyle \bullet$} [c] at  13 12
\put {$ \scriptstyle \bullet$} [c] at  16 0
\put {$ \scriptstyle \bullet$} [c] at  16 6
\put {$ \scriptstyle \bullet$} [c] at  16 12
\setlinear \plot  10 0 10 12  /
\setlinear \plot  16 0 16 12  /
\setlinear \plot  10 0  16 6   /
\setlinear \plot 10 6  13 12    /
\put{$2{,}520$} [c] at 13 -2
\endpicture
\end{minipage}
\begin{minipage}{4cm}
\beginpicture
\setcoordinatesystem units   <1.5mm,2mm>
\setplotarea x from 0 to 16, y from -2 to 15
\put{1.043)} [l] at 2 12
\put {$ \scriptstyle \bullet$} [c] at  10 0
\put {$ \scriptstyle \bullet$} [c] at  10 6
\put {$ \scriptstyle \bullet$} [c] at  10 12
\put {$ \scriptstyle \bullet$} [c] at  13 0
\put {$ \scriptstyle \bullet$} [c] at  16 0
\put {$ \scriptstyle \bullet$} [c] at  16 6
\put {$ \scriptstyle \bullet$} [c] at  16 12
\setlinear \plot  10 0 10 12  /
\setlinear \plot  16 0 16 12  /
\setlinear \plot  10 12  16 6   /
\setlinear \plot 10 6  13 0    /
\put{$2{,}520$} [c] at 13 -2
\endpicture
\end{minipage}
\begin{minipage}{4cm}
\beginpicture
\setcoordinatesystem units   <1.5mm,2mm>
\setplotarea x from 0 to 16, y from -2 to 15
\put{1.044)} [l] at 2 12
\put {$ \scriptstyle \bullet$} [c] at  10 0
\put {$ \scriptstyle \bullet$} [c] at  10 4
\put {$ \scriptstyle \bullet$} [c] at  10 8
\put {$ \scriptstyle \bullet$} [c] at  10 12
\put {$ \scriptstyle \bullet$} [c] at  13 0
\put {$ \scriptstyle \bullet$} [c] at  13 12
\put {$ \scriptstyle \bullet$} [c] at  16 12
\setlinear \plot 16 12 13 0 13 12     /
\setlinear \plot  10 0 10 12     /
\setlinear \plot  13 0   10 4     /
\put{$2{,}520$} [c] at 13 -2
\endpicture
\end{minipage}
$$
$$
\begin{minipage}{4cm}
\beginpicture
\setcoordinatesystem units   <1.5mm,2mm>
\setplotarea x from 0 to 16, y from -2 to 15
\put{1.045)} [l] at 2 12
\put {$ \scriptstyle \bullet$} [c] at  10 0
\put {$ \scriptstyle \bullet$} [c] at  10 4
\put {$ \scriptstyle \bullet$} [c] at  10 8
\put {$ \scriptstyle \bullet$} [c] at  10 12
\put {$ \scriptstyle \bullet$} [c] at  13 0
\put {$ \scriptstyle \bullet$} [c] at  13 12
\put {$ \scriptstyle \bullet$} [c] at  16 0
\setlinear \plot 16 0 13 12 13 0     /
\setlinear \plot  10 0 10 12     /
\setlinear \plot  13 12   10 8 /
\put{$2{,}520$} [c] at 13 -2
\endpicture
\end{minipage}
\begin{minipage}{4cm}
\beginpicture
\setcoordinatesystem units   <1.5mm,2mm>
\setplotarea x from 0 to 16, y from -2 to 15
\put{1.046)} [l] at 2 12
\put {$ \scriptstyle \bullet$} [c] at  10 0
\put {$ \scriptstyle \bullet$} [c] at  10 6
\put {$ \scriptstyle \bullet$} [c] at  10 12
\put {$ \scriptstyle \bullet$} [c] at  13 12
\put {$ \scriptstyle \bullet$} [c] at  16 0
\put {$ \scriptstyle \bullet$} [c] at  16 6
\put {$ \scriptstyle \bullet$} [c] at  16 12
\setlinear \plot  13 12 16 6 10 12 10 0   /
\setlinear \plot  16 0 16 12   /
\put{$2{,}520$} [c] at 13 -2
\endpicture
\end{minipage}
\begin{minipage}{4cm}
\beginpicture
\setcoordinatesystem units   <1.5mm,2mm>
\setplotarea x from 0 to 16, y from -2 to 15
\put{1.047)} [l] at 2 12
\put {$ \scriptstyle \bullet$} [c] at  10 0
\put {$ \scriptstyle \bullet$} [c] at  10 6
\put {$ \scriptstyle \bullet$} [c] at  10 12
\put {$ \scriptstyle \bullet$} [c] at  13 0
\put {$ \scriptstyle \bullet$} [c] at  16 0
\put {$ \scriptstyle \bullet$} [c] at  16 6
\put {$ \scriptstyle \bullet$} [c] at  16 12
\setlinear \plot  13 0 16 6 10 0 10 12   /
\setlinear \plot  16 0 16 12   /
\put{$2{,}520$} [c] at 13 -2
\endpicture
\end{minipage}
\begin{minipage}{4cm}
\beginpicture
\setcoordinatesystem units   <1.5mm,2mm>
\setplotarea x from 0 to 16, y from -2 to 15
\put{1.048)} [l] at 2 12
\put {$ \scriptstyle \bullet$} [c] at  10 12
\put {$ \scriptstyle \bullet$} [c] at  10 0
\put {$ \scriptstyle \bullet$} [c] at  13 12
\put {$ \scriptstyle \bullet$} [c] at  13 6
\put {$ \scriptstyle \bullet$} [c] at  14.5 0
\put {$ \scriptstyle \bullet$} [c] at  16 6
\put {$ \scriptstyle \bullet$} [c] at  16 12
\setlinear \plot 10 12 10 0 13 12 13 6 14.5 0 16  6 16 12  13 6 /
\setlinear \plot 13 12  16 6 /
\put{$2{,}520$} [c] at 13 -2
\endpicture
\end{minipage}
\begin{minipage}{4cm}
\beginpicture
\setcoordinatesystem units   <1.5mm,2mm>
\setplotarea x from 0 to 16, y from -2 to 15
\put{1.049)} [l] at 2 12
\put {$ \scriptstyle \bullet$} [c] at  10 12
\put {$ \scriptstyle \bullet$} [c] at  10 0
\put {$ \scriptstyle \bullet$} [c] at  13 0
\put {$ \scriptstyle \bullet$} [c] at  13 6
\put {$ \scriptstyle \bullet$} [c] at  14.5 12
\put {$ \scriptstyle \bullet$} [c] at  16 6
\put {$ \scriptstyle \bullet$} [c] at  16 0
\setlinear \plot 10 0 10 12 13 0 13 6 14.5 12 16  6 16 0  13 6 /
\setlinear \plot 13 0  16 6 /
\put{$2{,}520$} [c] at 13 -2
\endpicture
\end{minipage}
\begin{minipage}{4cm}
\beginpicture
\setcoordinatesystem units   <1.5mm,2mm>
\setplotarea x from 0 to 16, y from -2 to 15
\put{1.050)} [l] at 2 12
\put {$ \scriptstyle \bullet$} [c] at  10 0
\put {$ \scriptstyle \bullet$} [c] at  10 6
\put {$ \scriptstyle \bullet$} [c] at  10 12
\put {$ \scriptstyle \bullet$} [c] at  13 0
\put {$ \scriptstyle \bullet$} [c] at  13 6
\put {$ \scriptstyle \bullet$} [c] at  13 12
\put {$ \scriptstyle \bullet$} [c] at  16 12
\setlinear \plot  13  12  13 0 16 12    /
\setlinear \plot  10 12  10 0 13 6   /
\setlinear \plot  10 6 13 0  /
\put{$2{,}520$} [c] at 13 -2
\endpicture
\end{minipage}
$$
$$
\begin{minipage}{4cm}
\beginpicture
\setcoordinatesystem units   <1.5mm,2mm>
\setplotarea x from 0 to 16, y from -2 to 15
\put{1.051)} [l] at 2 12
\put {$ \scriptstyle \bullet$} [c] at  10 0
\put {$ \scriptstyle \bullet$} [c] at  10 6
\put {$ \scriptstyle \bullet$} [c] at  10 12
\put {$ \scriptstyle \bullet$} [c] at  13 0
\put {$ \scriptstyle \bullet$} [c] at  13 6
\put {$ \scriptstyle \bullet$} [c] at  13 12
\put {$ \scriptstyle \bullet$} [c] at  16 0
\setlinear \plot  13  0  13 12 16 0    /
\setlinear \plot  10 0  10 12 13 6   /
\setlinear \plot  10 6 13 12  /
\put{$2{,}520$} [c] at 13 -2
\endpicture
\end{minipage}
\begin{minipage}{4cm}
\beginpicture
\setcoordinatesystem units   <1.5mm,2mm>
\setplotarea x from 0 to 16, y from -2 to 15
\put{1.052)} [l] at 2 12
\put {$ \scriptstyle \bullet$} [c] at  10 0
\put {$ \scriptstyle \bullet$} [c] at  10 12
\put {$ \scriptstyle \bullet$} [c] at  13 12
\put {$ \scriptstyle \bullet$} [c] at  13 0
\put {$ \scriptstyle \bullet$} [c] at  16 0
\put {$ \scriptstyle \bullet$} [c] at  16 6
\put {$ \scriptstyle \bullet$} [c] at  16  12
\setlinear \plot  10  12 10 0 16 12 16 0 13 12 13 0 10  12 16 6   /
\setlinear \plot  13 0 16 12  /
\put{$2{,}520$} [c] at 13 -2
\endpicture
\end{minipage}
\begin{minipage}{4cm}
\beginpicture
\setcoordinatesystem units   <1.5mm,2mm>
\setplotarea x from 0 to 16, y from -2 to 15
\put{1.053)} [l] at 2 12
\put {$ \scriptstyle \bullet$} [c] at  10 0
\put {$ \scriptstyle \bullet$} [c] at  10 12
\put {$ \scriptstyle \bullet$} [c] at  13 12
\put {$ \scriptstyle \bullet$} [c] at  13 0
\put {$ \scriptstyle \bullet$} [c] at  16 0
\put {$ \scriptstyle \bullet$} [c] at  16 6
\put {$ \scriptstyle \bullet$} [c] at  16  12
\setlinear \plot  10  0 10 12 16 0 16 12 13 0 13 12 10  0 16 6   /
\setlinear \plot  13 12 16 0  /
\put{$2{,}520$} [c] at 13 -2
\endpicture
\end{minipage}
\begin{minipage}{4cm}
\beginpicture
\setcoordinatesystem units   <1.5mm,2mm>
\setplotarea x from 0 to 16, y from -2 to 15
\put{1.054)} [l] at 2 12
\put {$ \scriptstyle \bullet$} [c] at  10 0
\put {$ \scriptstyle \bullet$} [c] at  10 12
\put {$ \scriptstyle \bullet$} [c] at  13 12
\put {$ \scriptstyle \bullet$} [c] at  13 6
\put {$ \scriptstyle \bullet$} [c] at  13 0
\put {$ \scriptstyle \bullet$} [c] at  16  0
\put {$ \scriptstyle \bullet$} [c] at  16 12
\setlinear \plot  13 0 13 12 10 0 10 12 16  0  16 12    /
\put{$2{,}520 $} [c] at 13 -2
\endpicture
\end{minipage}
\begin{minipage}{4cm}
\beginpicture
\setcoordinatesystem units   <1.5mm,2mm>
\setplotarea x from 0 to 16, y from -2 to 15
\put{1.055)} [l] at 2 12
\put {$ \scriptstyle \bullet$} [c] at  10 0
\put {$ \scriptstyle \bullet$} [c] at  10 12
\put {$ \scriptstyle \bullet$} [c] at  13 12
\put {$ \scriptstyle \bullet$} [c] at  13 6
\put {$ \scriptstyle \bullet$} [c] at  13 0
\put {$ \scriptstyle \bullet$} [c] at  16  0
\put {$ \scriptstyle \bullet$} [c] at  16 12
\setlinear \plot  13 12 13 0 10 12 10 0 16  12  16 0    /
\put{$2{,}520 $} [c] at 13 -2
\endpicture
\end{minipage}
\begin{minipage}{4cm}
\beginpicture
\setcoordinatesystem units   <1.5mm,2mm>
\setplotarea x from 0 to 16, y from -2 to 15
\put{1.056)} [l] at 2 12
\put {$ \scriptstyle \bullet$} [c] at  10 0
\put {$ \scriptstyle \bullet$} [c] at  10 12
\put {$ \scriptstyle \bullet$} [c] at  13 0
\put {$ \scriptstyle \bullet$} [c] at  13 12
\put {$ \scriptstyle \bullet$} [c] at  16 0
\put {$ \scriptstyle \bullet$} [c] at  16 6
\put {$ \scriptstyle \bullet$} [c] at  16  12
\setlinear \plot  10  0 10 12 13 0 16  12 16 0    /
\setlinear \plot  10 12 16 6 13 12 13 0     /
\put{$2{,}520$} [c] at 13 -2
\endpicture
\end{minipage}
$$
$$
\begin{minipage}{4cm}
\beginpicture
\setcoordinatesystem units   <1.5mm,2mm>
\setplotarea x from 0 to 16, y from -2 to 15
\put{1.057)} [l] at 2 12
\put {$ \scriptstyle \bullet$} [c] at  10 0
\put {$ \scriptstyle \bullet$} [c] at  10 12
\put {$ \scriptstyle \bullet$} [c] at  13 0
\put {$ \scriptstyle \bullet$} [c] at  13 12
\put {$ \scriptstyle \bullet$} [c] at  16 0
\put {$ \scriptstyle \bullet$} [c] at  16 6
\put {$ \scriptstyle \bullet$} [c] at  16  12
\setlinear \plot  10  12 10 0 13 12 16  0 16 12    /
\setlinear \plot  10 0 16 6 13 0 13 12     /
\put{$2{,}520$} [c] at 13 -2
\endpicture
\end{minipage}
\begin{minipage}{4cm}
\beginpicture
\setcoordinatesystem units   <1.5mm,2mm>
\setplotarea x from 0 to 16, y from -2 to 15
\put{1.058)} [l] at 2 12
\put {$ \scriptstyle \bullet$} [c] at  10 0
\put {$ \scriptstyle \bullet$} [c] at  10 12
\put {$ \scriptstyle \bullet$} [c] at  13 0
\put {$ \scriptstyle \bullet$} [c] at  13 4
\put {$ \scriptstyle \bullet$} [c] at  13  12
\put {$ \scriptstyle \bullet$} [c] at  16  12
\put {$ \scriptstyle \bullet$} [c] at  16  0
\setlinear \plot 13 12 16 0 10 12 10 0 13 4 13 12 16 0   /
\setlinear \plot 13 4 13 0  16 12 16 0  /
\put{$2{,}520 $} [c] at 13 -2
\endpicture
\end{minipage}
\begin{minipage}{4cm}
\beginpicture
\setcoordinatesystem units   <1.5mm,2mm>
\setplotarea x from 0 to 16, y from -2 to 15
\put{1.059)} [l] at 2 12
\put {$ \scriptstyle \bullet$} [c] at  10 0
\put {$ \scriptstyle \bullet$} [c] at  10 12
\put {$ \scriptstyle \bullet$} [c] at  13 0
\put {$ \scriptstyle \bullet$} [c] at  13 8
\put {$ \scriptstyle \bullet$} [c] at  13  12
\put {$ \scriptstyle \bullet$} [c] at  16  12
\put {$ \scriptstyle \bullet$} [c] at  16  0
\setlinear \plot 13 0 16 12 10 0 10 12 13 8 13 0 16 12   /
\setlinear \plot 13 8 13 12  16 0 16 12  /
\put{$2{,}520 $} [c] at 13 -2
\endpicture
\end{minipage}
\begin{minipage}{4cm}
\beginpicture
\setcoordinatesystem units   <1.5mm,2mm>
\setplotarea x from 0 to 16, y from -2 to 15
\put{1.060)} [l] at 2 12
\put {$ \scriptstyle \bullet$} [c] at  10 0
\put {$ \scriptstyle \bullet$} [c] at  10 12
\put {$ \scriptstyle \bullet$} [c] at  13 12
\put {$ \scriptstyle \bullet$} [c] at  13 0
\put {$ \scriptstyle \bullet$} [c] at  13 6
\put {$ \scriptstyle \bullet$} [c] at  16 12
\put {$ \scriptstyle \bullet$} [c] at  16 0
\setlinear \plot  10  12 10 0  13 12 13 0  /
\setlinear \plot  16 12 13 6  16 0 /
\put{$2{,}520$} [c] at 13 -2
\endpicture
\end{minipage}
\begin{minipage}{4cm}
\beginpicture
\setcoordinatesystem units   <1.5mm,2mm>
\setplotarea x from 0 to 16, y from -2 to 15
\put{1.061)} [l] at 2 12
\put {$ \scriptstyle \bullet$} [c] at  10 0
\put {$ \scriptstyle \bullet$} [c] at  10 12
\put {$ \scriptstyle \bullet$} [c] at  13 12
\put {$ \scriptstyle \bullet$} [c] at  13 0
\put {$ \scriptstyle \bullet$} [c] at  13 6
\put {$ \scriptstyle \bullet$} [c] at  16 12
\put {$ \scriptstyle \bullet$} [c] at  16 0
\setlinear \plot  10  0 10 12  13 0 13 12  /
\setlinear \plot  16 0 13 6  16 12 /
\put{$2{,}520 $} [c] at 13 -2
\endpicture
\end{minipage}
\begin{minipage}{4cm}
\beginpicture
\setcoordinatesystem units   <1.5mm,2mm>
\setplotarea x from 0 to 16, y from -2 to 15
\put{1.062)} [l] at 2 12
\put {$ \scriptstyle \bullet$} [c] at 10 3
\put {$ \scriptstyle \bullet$} [c] at 12 6
\put {$ \scriptstyle \bullet$} [c] at 13 0
\put {$ \scriptstyle \bullet$} [c] at 13 3
\put {$ \scriptstyle \bullet$} [c] at 13 12
\put {$ \scriptstyle \bullet$} [c] at 14 6
\put {$ \scriptstyle \bullet$} [c] at 16 3
\setlinear \plot 10 3 13 12 16 3 13 0 10 3 /
\setlinear \plot 13 0 13 3 14 6 13 12 12 6 13 3 /
\put{$1{,}260$} [c] at 13 -2
\endpicture
\end{minipage}
$$
$$
\begin{minipage}{4cm}
\beginpicture
\setcoordinatesystem units   <1.5mm,2mm>
\setplotarea x from 0 to 16, y from -2 to 15
\put{1.063)} [l] at 2 12
\put {$ \scriptstyle \bullet$} [c] at 10 9
\put {$ \scriptstyle \bullet$} [c] at 12 6
\put {$ \scriptstyle \bullet$} [c] at 13 12
\put {$ \scriptstyle \bullet$} [c] at 13 9
\put {$ \scriptstyle \bullet$} [c] at 13 0
\put {$ \scriptstyle \bullet$} [c] at 14 6
\put {$ \scriptstyle \bullet$} [c] at 16 9
\setlinear \plot 10 9 13 12 16 9 13 0 10 9 /
\setlinear \plot 13 12 13 9 14 6 13 0 12 6 13 9 /
\put{$1{,}260$} [c] at 13 -2
\endpicture
\end{minipage}
\begin{minipage}{4cm}
\beginpicture
\setcoordinatesystem units   <1.5mm,2mm>
\setplotarea x from 0 to 16, y from -2 to 15
\put{1.064)} [l] at 2 12
\put {$ \scriptstyle \bullet$} [c] at 13 0
\put {$ \scriptstyle \bullet$} [c] at 10 4
\put {$ \scriptstyle \bullet$} [c] at 12 4
\put {$ \scriptstyle \bullet$} [c] at 14 4
\put {$ \scriptstyle \bullet$} [c] at 16 4
\put {$ \scriptstyle \bullet$} [c] at 12 12
\put {$ \scriptstyle \bullet$} [c] at 14 12
\setlinear \plot 13 0 10 4 12 12 12 4 13 0 16 4 14 12 14 4 13 0 /
\setlinear \plot 12 12 14  4    /
\setlinear \plot 14 12 12  4    /
\put{$1{,}260$} [c] at  13 -2
\endpicture
\end{minipage}
\begin{minipage}{4cm}
\beginpicture
\setcoordinatesystem units   <1.5mm,2mm>
\setplotarea x from 0 to 16, y from -2 to 15
\put{1.065)} [l] at 2 12
\put {$ \scriptstyle \bullet$} [c] at 13 12
\put {$ \scriptstyle \bullet$} [c] at 10 8
\put {$ \scriptstyle \bullet$} [c] at 12 8
\put {$ \scriptstyle \bullet$} [c] at 14 8
\put {$ \scriptstyle \bullet$} [c] at 16 8
\put {$ \scriptstyle \bullet$} [c] at 12 0
\put {$ \scriptstyle \bullet$} [c] at 14 0
\setlinear \plot 13 12 10 8 12 0 12 8 13 12 16 8 14 0 14 8 13 12 /
\setlinear \plot 12 0 14  8    /
\setlinear \plot 14 0 12  8    /
\put{$1{,}260$} [c] at  13 -2
\endpicture
\end{minipage}
\begin{minipage}{4cm}
\beginpicture
\setcoordinatesystem units   <1.5mm,2mm>
\setplotarea x from 0 to 16, y from -2 to 15
\put{1.066)} [l] at 2 12
\put {$ \scriptstyle \bullet$} [c] at 10 12
\put {$ \scriptstyle \bullet$} [c] at  12 9
\put {$ \scriptstyle \bullet$} [c] at 13 12
\put {$ \scriptstyle \bullet$} [c] at 13 0
\put {$ \scriptstyle \bullet$} [c] at 13 3
\put {$ \scriptstyle \bullet$} [c] at 14 9
\put {$ \scriptstyle \bullet$} [c] at 16 12
\setlinear \plot 10 12 13 3 12 9 13 12 14 9 13 3 13 0 /
\setlinear \plot 16 12 13 3  /
\put{$1{,}260$} [c] at 13 -2
\endpicture
\end{minipage}
\begin{minipage}{4cm}
\beginpicture
\setcoordinatesystem units   <1.5mm,2mm>
\setplotarea x from 0 to 16, y from -2 to 15
\put{1.067)} [l] at 2 12
\put {$ \scriptstyle \bullet$} [c] at 10 0
\put {$ \scriptstyle \bullet$} [c] at  12 3
\put {$ \scriptstyle \bullet$} [c] at 13 12
\put {$ \scriptstyle \bullet$} [c] at 13 0
\put {$ \scriptstyle \bullet$} [c] at 13 9
\put {$ \scriptstyle \bullet$} [c] at 14 3
\put {$ \scriptstyle \bullet$} [c] at 16 0
\setlinear \plot 10 0 13 9 12 3 13 0 14 3 13 9 13 12 /
\setlinear \plot 16 0 13 9  /
\put{$1{,}260$} [c] at 13 -2
\endpicture
\end{minipage}
\begin{minipage}{4cm}
\beginpicture
\setcoordinatesystem units   <1.5mm,2mm>
\setplotarea x from 0 to 16, y from -2 to 15
\put{1.068)} [l] at 2 12
\put {$ \scriptstyle \bullet$} [c] at 10 6
\put {$ \scriptstyle \bullet$} [c] at 10 0
\put {$ \scriptstyle \bullet$} [c] at 13 0
\put {$ \scriptstyle \bullet$} [c] at 13 6
\put {$ \scriptstyle \bullet$} [c] at 13 12
\put {$ \scriptstyle \bullet$} [c] at 16 0
\put {$ \scriptstyle \bullet$} [c] at 16 6
\setlinear \plot 16 0 16 6 13  12 13  0 10 6 10  0 13 6 /
\setlinear \plot 10 6 13 12     /
\put{$1{,}260$} [c] at  13 -2
\endpicture
\end{minipage}
$$
$$
\begin{minipage}{4cm}
\beginpicture
\setcoordinatesystem units   <1.5mm,2mm>
\setplotarea x from 0 to 16, y from -2 to 15
\put{1.069)} [l] at 2 12
\put {$ \scriptstyle \bullet$} [c] at 10 6
\put {$ \scriptstyle \bullet$} [c] at 10 12
\put {$ \scriptstyle \bullet$} [c] at 13 0
\put {$ \scriptstyle \bullet$} [c] at 13 6
\put {$ \scriptstyle \bullet$} [c] at 13 12
\put {$ \scriptstyle \bullet$} [c] at 16 12
\put {$ \scriptstyle \bullet$} [c] at 16 6
\setlinear \plot 16 12 16 6 13  0 13  12 10 6 10  12 13 6 /
\setlinear \plot 10 6 13 0     /
\put{$1{,}260$} [c] at  13 -2
\endpicture
\end{minipage}
\begin{minipage}{4cm}
\beginpicture
\setcoordinatesystem units   <1.5mm,2mm>
\setplotarea x from 0 to 16, y from -2 to 15
\put{1.070)} [l] at 2 12
\put {$ \scriptstyle \bullet$} [c] at 13 12
\put {$ \scriptstyle \bullet$} [c] at 14 6
\put {$ \scriptstyle \bullet$} [c] at 12 6
\put {$ \scriptstyle \bullet$} [c] at 10 0
\put {$ \scriptstyle \bullet$} [c] at 12 0
\put {$ \scriptstyle \bullet$} [c] at 14 0
\put {$ \scriptstyle \bullet$} [c] at 16 0
\setlinear \plot 10 0 12 6 13 12 14 6 16 0 /
\setlinear \plot 12 6 12 0 14 6 14 0 12 6    /
\put{$1{,}260$} [c] at  13 -2
\endpicture
\end{minipage}
\begin{minipage}{4cm}
\beginpicture
\setcoordinatesystem units   <1.5mm,2mm>
\setplotarea x from 0 to 16, y from -2 to 15
\put{1.071)} [l] at 2 12
\put {$ \scriptstyle \bullet$} [c] at 13 0
\put {$ \scriptstyle \bullet$} [c] at 14 6
\put {$ \scriptstyle \bullet$} [c] at 12 6
\put {$ \scriptstyle \bullet$} [c] at 10 12
\put {$ \scriptstyle \bullet$} [c] at 12 12
\put {$ \scriptstyle \bullet$} [c] at 14 12
\put {$ \scriptstyle \bullet$} [c] at 16 12
\setlinear \plot 10 12 12 6 13 0 14 6 16 12 /
\setlinear \plot 12 6 12 12 14 6 14 12 12 6    /
\put{$1{,}260$} [c] at  13 -2
\endpicture
\end{minipage}
\begin{minipage}{4cm}
\beginpicture
\setcoordinatesystem units   <1.5mm,2mm>
\setplotarea x from 0 to 16, y from -2 to 15
\put{1.072)} [l] at 2 12
\put {$ \scriptstyle \bullet$} [c] at 10 0
\put {$ \scriptstyle \bullet$} [c] at 12 0
\put {$ \scriptstyle \bullet$} [c] at 14 0
\put {$ \scriptstyle \bullet$} [c] at 16 0
\put {$ \scriptstyle \bullet$} [c] at 11 8
\put {$ \scriptstyle \bullet$} [c] at 11 12
\put {$ \scriptstyle \bullet$} [c] at 15 4
\setlinear \plot 10 0 11 8  12 0   /
\setlinear \plot 14 0 15 4 16 0   /
\setlinear \plot 11 12 11 8 15 4   /
\put{$1{,}260$} [c] at  13 -2
\endpicture
\end{minipage}
\begin{minipage}{4cm}
\beginpicture
\setcoordinatesystem units   <1.5mm,2mm>
\setplotarea x from 0 to 16, y from -2 to 15
\put{1.073)} [l] at 2 12
\put {$ \scriptstyle \bullet$} [c] at 10 12
\put {$ \scriptstyle \bullet$} [c] at 12 12
\put {$ \scriptstyle \bullet$} [c] at 14 12
\put {$ \scriptstyle \bullet$} [c] at 16 12
\put {$ \scriptstyle \bullet$} [c] at 11 4
\put {$ \scriptstyle \bullet$} [c] at 11 0
\put {$ \scriptstyle \bullet$} [c] at 15 8
\setlinear \plot 10 12 11 4  12 12   /
\setlinear \plot 14 12 15 8 16 12   /
\setlinear \plot 11 0 11 4 15 8   /
\put{$1{,}260$} [c] at  13 -2
\endpicture
\end{minipage}
\begin{minipage}{4cm}
\beginpicture
\setcoordinatesystem units   <1.5mm,2mm>
\setplotarea x from 0 to 16, y from -2 to 15
\put{1.074)} [l] at 2 12
\put {$ \scriptstyle \bullet$} [c] at  10 0
\put {$ \scriptstyle \bullet$} [c] at  10  12
\put {$ \scriptstyle \bullet$} [c] at  13 0
\put {$ \scriptstyle \bullet$} [c] at  13 12
\put {$ \scriptstyle \bullet$} [c] at  16 6
\put {$ \scriptstyle \bullet$} [c] at  16 12
\put {$ \scriptstyle \bullet$} [c] at  16  0
\setlinear \plot  10  0 10 12  13 0 13 12  10 0 16 6 13 0    /
\setlinear \plot   16  0  16 12    /
\put{$1{,}260$} [c] at 13 -2
\endpicture
\end{minipage}
$$
$$
\begin{minipage}{4cm}
\beginpicture
\setcoordinatesystem units   <1.5mm,2mm>
\setplotarea x from 0 to 16, y from -2 to 15
\put{1.075)} [l] at 2 12
\put {$ \scriptstyle \bullet$} [c] at  10 0
\put {$ \scriptstyle \bullet$} [c] at  10  12
\put {$ \scriptstyle \bullet$} [c] at  13 0
\put {$ \scriptstyle \bullet$} [c] at  13 12
\put {$ \scriptstyle \bullet$} [c] at  16 6
\put {$ \scriptstyle \bullet$} [c] at  16 12
\put {$ \scriptstyle \bullet$} [c] at  16  0
\setlinear \plot  10  12 10 0  13 12 13 0  10 12 16 6 13 12    /
\setlinear \plot   16  0  16 12    /
\put{$1{,}260$} [c] at 13 -2
\endpicture
\end{minipage}
\begin{minipage}{4cm}
\beginpicture
\setcoordinatesystem units   <1.5mm,2mm>
\setplotarea x from 0 to 16, y from -2 to 15
\put{1.076)} [l] at 2 12
\put {$ \scriptstyle \bullet$} [c] at  10 0
\put {$ \scriptstyle \bullet$} [c] at  10 12
\put {$ \scriptstyle \bullet$} [c] at  13 12
\put {$ \scriptstyle \bullet$} [c] at  13 0
\put {$ \scriptstyle \bullet$} [c] at  16 0
\put {$ \scriptstyle \bullet$} [c] at  16  6
\put {$ \scriptstyle \bullet$} [c] at  16 12
\setlinear \plot 10 0 10 12  13 0 13  12 10 0 16 12  13 0   /
\setlinear \plot 10 12 16 0 16 12 /
\setlinear \plot 16 0 13 12  /
\put{$1{,}260 $} [c] at 13 -2
\endpicture
\end{minipage}
\begin{minipage}{4cm}
\beginpicture
\setcoordinatesystem units   <1.5mm,2mm>
\setplotarea x from 0 to 16, y from -2 to 15
\put{1.077)} [l] at 2 12
\put {$ \scriptstyle \bullet$} [c] at  12.1 6
\put {$ \scriptstyle \bullet$} [c] at  13 0
\put {$ \scriptstyle \bullet$} [c] at  13 6
\put {$ \scriptstyle \bullet$} [c] at  13 12
\put {$ \scriptstyle \bullet$} [c] at  13.9 6
\put {$ \scriptstyle \bullet$} [c] at  16 12
\put {$ \scriptstyle \bullet$} [c] at  10 0
\setlinear \plot 16 12 13 0 13 12 12.1 6 13 0 13.9 6 13 12    /
\setlinear \plot  13 6 10 0 12 6    /
\setlinear \plot   10 0  13.9  6    /
\put{$840 $} [c] at 13 -2
\endpicture
\end{minipage}
\begin{minipage}{4cm}
\beginpicture
\setcoordinatesystem units   <1.5mm,2mm>
\setplotarea x from 0 to 16, y from -2 to 15
\put{1.078)} [l] at 2 12
\put {$ \scriptstyle \bullet$} [c] at  12.1 6
\put {$ \scriptstyle \bullet$} [c] at  13 0
\put {$ \scriptstyle \bullet$} [c] at  13 6
\put {$ \scriptstyle \bullet$} [c] at  13 12
\put {$ \scriptstyle \bullet$} [c] at  13.9 6
\put {$ \scriptstyle \bullet$} [c] at  16 0
\put {$ \scriptstyle \bullet$} [c] at  10 12
\setlinear \plot 16 0 13 12 13 0 12.1 6 13 12 13.9 6 13 0    /
\setlinear \plot  13 6 10 12 12 6    /
\setlinear \plot   10 12  13.9  6    /
\put{$840 $} [c] at 13 -2
\endpicture
\end{minipage}
\begin{minipage}{4cm}
\beginpicture
\setcoordinatesystem units   <1.5mm,2mm>
\setplotarea x from 0 to 16, y from -2 to 15
\put {1.079)} [l] at 2 12
\put {$ \scriptstyle \bullet$} [c] at  10  0
\put {$ \scriptstyle \bullet$} [c] at  10 12
\put {$ \scriptstyle \bullet$} [c] at  12  6
\put {$ \scriptstyle \bullet$} [c] at  12 12
\put {$ \scriptstyle \bullet$} [c] at  14 12
\put {$ \scriptstyle \bullet$} [c] at  14  0
\put {$ \scriptstyle \bullet$} [c] at  16  12
\setlinear \plot   16 12 14 0 12 6  14 12  /
\setlinear \plot  10 0 12 6 12 12 /
\setlinear \plot  10 12 12 6 /
\put{$840$}[c] at 13 -2
\endpicture
\end{minipage}
\begin{minipage}{4cm}
\beginpicture
\setcoordinatesystem units   <1.5mm,2mm>
\setplotarea x from 0 to 16, y from -2 to 15
\put {1.080)} [l] at 2 12
\put {$ \scriptstyle \bullet$} [c] at  10  0
\put {$ \scriptstyle \bullet$} [c] at  10 12
\put {$ \scriptstyle \bullet$} [c] at  12  6
\put {$ \scriptstyle \bullet$} [c] at  12 0
\put {$ \scriptstyle \bullet$} [c] at  14 12
\put {$ \scriptstyle \bullet$} [c] at  14  0
\put {$ \scriptstyle \bullet$} [c] at  16  0
\setlinear \plot   16 0 14 12 12 6  14 0  /
\setlinear \plot  10 12 12 6 12 0 /
\setlinear \plot  10 0 12 6 /
\put{$840$}[c] at 13 -2
\endpicture
\end{minipage}
$$
$$
\begin{minipage}{4cm}
\beginpicture
\setcoordinatesystem units   <1.5mm,2mm>
\setplotarea x from 0 to 16, y from -2 to 15
\put {1.081)} [l] at 2 12
\put {$ \scriptstyle \bullet$} [c] at  10 0
\put {$ \scriptstyle \bullet$} [c] at  10 12
\put {$ \scriptstyle \bullet$} [c] at  12 12
\put {$ \scriptstyle \bullet$} [c] at  14 12
\put {$ \scriptstyle \bullet$} [c] at  14 6
\put {$ \scriptstyle \bullet$} [c] at  14  0
\put {$ \scriptstyle \bullet$} [c] at  16  12
\setlinear \plot  10 12 10 0 12 12 14 6 14 0 /
\setlinear \plot  10 0 14 12 14 6 /
\setlinear  \plot 10 12 14 6 16 12 /
\put{$840$}[c] at 13 -2
\endpicture
\end{minipage}
\begin{minipage}{4cm}
\beginpicture
\setcoordinatesystem units   <1.5mm,2mm>
\setplotarea x from 0 to 16, y from -2 to 15
\put {1.082)} [l] at 2 12
\put {$ \scriptstyle \bullet$} [c] at  10 0
\put {$ \scriptstyle \bullet$} [c] at  10 12
\put {$ \scriptstyle \bullet$} [c] at  12 0
\put {$ \scriptstyle \bullet$} [c] at  14 0
\put {$ \scriptstyle \bullet$} [c] at  14 6
\put {$ \scriptstyle \bullet$} [c] at  14  12
\put {$ \scriptstyle \bullet$} [c] at  16  0
\setlinear \plot  10 0 10 12 12 0 14 6 14 12 /
\setlinear \plot  10 12 14 0 14 6 /
\setlinear  \plot 10 0 14 6 16 0 /
\put{$840$}[c] at 13 -2
\endpicture
\end{minipage}
\begin{minipage}{4cm}
\beginpicture
\setcoordinatesystem units   <1.5mm,2mm>
\setplotarea x from 0 to 16, y from -2 to 15
\put {1.083)} [l] at 2 12
\put {$ \scriptstyle \bullet$} [c] at  10 0
\put {$ \scriptstyle \bullet$} [c] at  10 12
\put {$ \scriptstyle \bullet$} [c] at  12 12
\put {$ \scriptstyle \bullet$} [c] at  14 12
\put {$ \scriptstyle \bullet$} [c] at  14 4
\put {$ \scriptstyle \bullet$} [c] at  14 0
\put {$ \scriptstyle \bullet$} [c] at  16 12
\setlinear \plot  10 12 10 0 16 12 14 0 14  12 10 0 12 12  14 4  12 12 /
\setlinear \plot 10 12 14 4 /
\put{$840$}[c] at 13 -2
\endpicture
\end{minipage}
\begin{minipage}{4cm}
\beginpicture
\setcoordinatesystem units   <1.5mm,2mm>
\setplotarea x from 0 to 16, y from -2 to 15
\put {1.084)} [l] at 2 12
\put {$ \scriptstyle \bullet$} [c] at  10 0
\put {$ \scriptstyle \bullet$} [c] at  10 12
\put {$ \scriptstyle \bullet$} [c] at  12 0
\put {$ \scriptstyle \bullet$} [c] at  14 0
\put {$ \scriptstyle \bullet$} [c] at  14 8
\put {$ \scriptstyle \bullet$} [c] at  14 12
\put {$ \scriptstyle \bullet$} [c] at  16 0
\setlinear \plot  10 0 10 12 16 0 14 12 14  0 10 12 12 0  14 8  12 0 /
\setlinear \plot 10 0 14 8 /
\put{$840$}[c] at 13 -2
\endpicture
\end{minipage}
\begin{minipage}{4cm}
\beginpicture
\setcoordinatesystem units   <1.5mm,2mm>
\setplotarea x from 0 to 16, y from -2 to 15
\put{1.085)} [l] at 2 12
\put {$ \scriptstyle \bullet$} [c] at 10 12
\put {$ \scriptstyle \bullet$} [c] at 14 12
\put {$ \scriptstyle \bullet$} [c] at 12 0
\put {$ \scriptstyle \bullet$} [c] at 10.5 9
\put {$ \scriptstyle \bullet$} [c] at 11 6
\put {$ \scriptstyle \bullet$} [c] at  11.5 3
\put{$\scriptstyle \bullet$} [c] at 16  0
\setlinear \plot  10 12  12 0 14 12      /
\put{$5{,}040$} [c] at 13 -2
\endpicture
\end{minipage}
\begin{minipage}{4cm}
\beginpicture
\setcoordinatesystem units   <1.5mm,2mm>
\setplotarea x from 0 to 16, y from -2 to 15
\put{1.086)} [l] at 2 12
\put {$ \scriptstyle \bullet$} [c] at 10 0
\put {$ \scriptstyle \bullet$} [c] at 14 0
\put {$ \scriptstyle \bullet$} [c] at 12 12
\put {$ \scriptstyle \bullet$} [c] at 10.5 3
\put {$ \scriptstyle \bullet$} [c] at 11 6
\put {$ \scriptstyle \bullet$} [c] at  11.5 9
\put{$\scriptstyle \bullet$} [c] at 16  0
\setlinear \plot  10 0  12 12 14 0      /
\put{$5{,}040$} [c] at 13 -2
\endpicture
\end{minipage}
$$
$$
\begin{minipage}{4cm}
\beginpicture
\setcoordinatesystem units   <1.5mm,2mm>
\setplotarea x from 0 to 16, y from -2 to 15
\put{1.087)} [l] at 2 12
\put {$ \scriptstyle \bullet$} [c] at 10 6
\put {$ \scriptstyle \bullet$} [c] at 10 12
\put {$ \scriptstyle \bullet$} [c] at 12 0
\put {$ \scriptstyle \bullet$} [c] at 12 12
\put {$ \scriptstyle \bullet$} [c] at 14 6
\put {$ \scriptstyle \bullet$} [c] at 11 3
\put{$\scriptstyle \bullet$} [c] at 16  0
\setlinear \plot  10 12  10 6 12 0 14 6 12 12 10 6      /
\put{$5{,}040$} [c] at 13 -2
 \endpicture
\end{minipage}
\begin{minipage}{4cm}
\beginpicture
\setcoordinatesystem units   <1.5mm,2mm>
\setplotarea x from 0 to 16, y from -2 to 15
\put{1.088)} [l] at 2 12
\put {$ \scriptstyle \bullet$} [c] at 10 6
\put {$ \scriptstyle \bullet$} [c] at 10 0
\put {$ \scriptstyle \bullet$} [c] at 12 0
\put {$ \scriptstyle \bullet$} [c] at 12 12
\put {$ \scriptstyle \bullet$} [c] at 14 6
\put {$ \scriptstyle \bullet$} [c] at 11 9
\put{$\scriptstyle \bullet$} [c] at 16  0
\setlinear \plot  10 0  10 6 12 12 14 6 12 0 10 6      /
\put{$5{,}040$} [c] at 13 -2
\put{$\scriptstyle \bullet$} [c] at 16  0 \endpicture
\end{minipage}
\begin{minipage}{4cm}
\beginpicture
\setcoordinatesystem units   <1.5mm,2mm>
\setplotarea x from 0 to 16, y from -2 to 15
\put{1.089)} [l] at 2 12
\put {$ \scriptstyle \bullet$} [c] at 10 4
\put {$ \scriptstyle \bullet$} [c] at 12 0
\put {$ \scriptstyle \bullet$} [c] at 12 8
\put {$ \scriptstyle \bullet$} [c] at 12 12
\put {$ \scriptstyle \bullet$} [c] at 14  4
\put {$ \scriptstyle \bullet$} [c] at 14 12
\put{$\scriptstyle \bullet$} [c] at 16  0
\setlinear \plot  12 12 12 8  10 4 12 0 14 4 12 8    /
\setlinear \plot   14 4 14 12  /
\put{$5{,}040$} [c] at 13 -2
 \endpicture
\end{minipage}
\begin{minipage}{4cm}
\beginpicture
\setcoordinatesystem units   <1.5mm,2mm>
\setplotarea x from 0 to 16, y from -2 to 15
\put{1.090)} [l] at 2 12
\put {$ \scriptstyle \bullet$} [c] at 10 8
\put {$ \scriptstyle \bullet$} [c] at 12 0
\put {$ \scriptstyle \bullet$} [c] at 12 4
\put {$ \scriptstyle \bullet$} [c] at 12 12
\put {$ \scriptstyle \bullet$} [c] at 14  8
\put {$ \scriptstyle \bullet$} [c] at 14 0
\put{$\scriptstyle \bullet$} [c] at 16  0
\setlinear \plot  12 0 12 4  10 8 12 12 14 8 12 4    /
\setlinear \plot   14 8 14 0  /
\put{$5{,}040$} [c] at 13 -2
\endpicture
\end{minipage}
\begin{minipage}{4cm}
\beginpicture
\setcoordinatesystem units   <1.5mm,2mm>
\setplotarea x from 0 to 16, y from -2 to 15
\put{1.091)} [l] at 2 12
\put {$ \scriptstyle \bullet$} [c] at 10 4
\put {$ \scriptstyle \bullet$} [c] at 10 8
\put {$ \scriptstyle \bullet$} [c] at 10 12
\put {$ \scriptstyle \bullet$} [c] at 12 0
\put {$ \scriptstyle \bullet$} [c] at 14 4
\put {$ \scriptstyle \bullet$} [c] at 14 12
\put{$\scriptstyle \bullet$} [c] at 16  0
\setlinear \plot 10 4  10 12  14 4 14 12 10 4 12 0 14 4 10 12     /
\put{$5{,}040$} [c] at 13 -2
 \endpicture
\end{minipage}
\begin{minipage}{4cm}
\beginpicture
\setcoordinatesystem units   <1.5mm,2mm>
\setplotarea x from 0 to 16, y from -2 to 15
\put{1.092)} [l] at 2 12
\put {$ \scriptstyle \bullet$} [c] at 10 4
\put {$ \scriptstyle \bullet$} [c] at 10 8
\put {$ \scriptstyle \bullet$} [c] at 10 0
\put {$ \scriptstyle \bullet$} [c] at 12 12
\put {$ \scriptstyle \bullet$} [c] at 14 8
\put {$ \scriptstyle \bullet$} [c] at 14 0
\put{$\scriptstyle \bullet$} [c] at 16  0
\setlinear \plot 10 8  10 0  14 8 14 0 10 8 12 12 14 8 10 0     /
\put{$5{,}040$} [c] at 13 -2
\endpicture
\end{minipage}
$$
$$
\begin{minipage}{4cm}
\beginpicture
\setcoordinatesystem units   <1.5mm,2mm>
\setplotarea x from 0 to 16, y from -2 to 15
\put{1.093)} [l] at 2 12
\put {$ \scriptstyle \bullet$} [c] at 10 0
\put {$ \scriptstyle \bullet$} [c] at 10 4
\put {$ \scriptstyle \bullet$} [c] at 10 8
\put {$ \scriptstyle \bullet$} [c] at 10 12
\put {$ \scriptstyle \bullet$} [c] at 14 0
\put {$ \scriptstyle \bullet$} [c] at 14 12
\put{$\scriptstyle \bullet$} [c] at 16  0
\setlinear \plot    10 0 10 12   /
\setlinear \plot    10 4 14 12 14 0 10 8   /
\put{$5{,}040$} [c] at 13 -2
\endpicture
\end{minipage}
\begin{minipage}{4cm}
\beginpicture
\setcoordinatesystem units   <1.5mm,2mm>
\setplotarea x from 0 to 16, y from -2 to 15
\put{1.094)} [l] at 2 12
\put {$ \scriptstyle \bullet$} [c] at 10 12
\put {$ \scriptstyle \bullet$} [c] at 10 8
\put {$ \scriptstyle \bullet$} [c] at 12 4
\put {$ \scriptstyle \bullet$} [c] at 12 0
\put {$ \scriptstyle \bullet$} [c] at 14 12
\put {$ \scriptstyle \bullet$} [c] at 14 8
\put{$\scriptstyle \bullet$} [c] at 16  0
\setlinear \plot    10 12 10 8 12 4 12 0    /
\setlinear \plot  12 4 14 8 14 12 /
\put{$2{,}520$} [c] at 13 -2
\endpicture
\end{minipage}
\begin{minipage}{4cm}
\beginpicture
\setcoordinatesystem units   <1.5mm,2mm>
\setplotarea x from 0 to 16, y from -2 to 15
\put{1.095)} [l] at 2 12
\put {$ \scriptstyle \bullet$} [c] at 10 0
\put {$ \scriptstyle \bullet$} [c] at 10 4
\put {$ \scriptstyle \bullet$} [c] at 12 8
\put {$ \scriptstyle \bullet$} [c] at 12 12
\put {$ \scriptstyle \bullet$} [c] at 14 0
\put {$ \scriptstyle \bullet$} [c] at 14 4
\put{$\scriptstyle \bullet$} [c] at 16  0
\setlinear \plot    10 0 10 4 12 8 12 12    /
\setlinear \plot  12 8 14 4 14 0 /
\put{$2{,}520$} [c] at 13 -2
\endpicture
\end{minipage}
\begin{minipage}{4cm}
\beginpicture
\setcoordinatesystem units   <1.5mm,2mm>
\setplotarea x from 0 to 16, y from -2 to 15
\put{1.096)} [l] at 2 12
\put {$ \scriptstyle \bullet$} [c] at 10 0
\put {$ \scriptstyle \bullet$} [c] at 10 4
\put {$ \scriptstyle \bullet$} [c] at 10 8
\put {$ \scriptstyle \bullet$} [c] at 10 12
\put {$ \scriptstyle \bullet$} [c] at 14 0
\put {$ \scriptstyle \bullet$} [c] at 14 12
\put{$\scriptstyle \bullet$} [c] at 16  0
\setlinear \plot 10 12 10 0 14 12 14 0 10 4    /
\put{$2{,}520$} [c] at 13 -2
\endpicture
\end{minipage}
\begin{minipage}{4cm}
\beginpicture
\setcoordinatesystem units   <1.5mm,2mm>
\setplotarea x from 0 to 16, y from -2 to 15
\put{1.097)} [l] at 2 12
\put {$ \scriptstyle \bullet$} [c] at 10 0
\put {$ \scriptstyle \bullet$} [c] at 10 4
\put {$ \scriptstyle \bullet$} [c] at 10 8
\put {$ \scriptstyle \bullet$} [c] at 10 12
\put {$ \scriptstyle \bullet$} [c] at 14 0
\put {$ \scriptstyle \bullet$} [c] at 14 12
\put{$\scriptstyle \bullet$} [c] at 16  0
\setlinear \plot 10 0 10 12 14 0 14 12 10 8    /
\put{$2{,}520$} [c] at 13 -2
\endpicture
\end{minipage}
\begin{minipage}{4cm}
\beginpicture
\setcoordinatesystem units   <1.5mm,2mm>
\setplotarea x from 0 to 16, y from -2 to 15
\put{1.098)} [l] at 2 12
\put {$ \scriptstyle \bullet$} [c] at 10 6
\put {$ \scriptstyle \bullet$} [c] at 10 12
\put {$ \scriptstyle \bullet$} [c] at 12 0
\put {$ \scriptstyle \bullet$} [c] at 14 0
\put {$ \scriptstyle \bullet$} [c] at 14 6
\put {$ \scriptstyle \bullet$} [c] at 14 12
\put{$\scriptstyle \bullet$} [c] at 16  0
\setlinear \plot 14 0 14 12 10 6 10 12 14 6 12 0 10 6    /
\put{$2{,}520$} [c] at 13 -2
 \endpicture
\end{minipage}
$$
$$
\begin{minipage}{4cm}
\beginpicture
\setcoordinatesystem units   <1.5mm,2mm>
\setplotarea x from 0 to 16, y from -2 to 15
\put{1.099)} [l] at 2 12
\put {$ \scriptstyle \bullet$} [c] at 10 6
\put {$ \scriptstyle \bullet$} [c] at 10 0
\put {$ \scriptstyle \bullet$} [c] at 12 12
\put {$ \scriptstyle \bullet$} [c] at 14 0
\put {$ \scriptstyle \bullet$} [c] at 14 6
\put {$ \scriptstyle \bullet$} [c] at 14 12
\put{$\scriptstyle \bullet$} [c] at 16  0
\setlinear \plot 14 12 14 0 10 6 10 0 14 6 12 12 10 6    /
\put{$2{,}520$} [c] at 13 -2
\endpicture
\end{minipage}
\begin{minipage}{4cm}
\beginpicture
\setcoordinatesystem units   <1.5mm,2mm>
\setplotarea x from 0 to 16, y from -2 to 15
\put{1.100)} [l] at 2 12
\put {$ \scriptstyle \bullet$} [c] at 10 4
\put {$ \scriptstyle \bullet$} [c] at 10 8
\put {$ \scriptstyle \bullet$} [c] at 12 0
\put {$ \scriptstyle \bullet$} [c] at 12 12
\put {$ \scriptstyle \bullet$} [c] at 14 4
\put {$ \scriptstyle \bullet$} [c] at 14 8
\put{$\scriptstyle \bullet$} [c] at 16  0
\setlinear \plot  12 0  10 4 10 8  12 12  14 8 14 4  12 0     /
\put{$2{,}520$} [c] at 13 -2
 \endpicture
\end{minipage}
\begin{minipage}{4cm}
\beginpicture
\setcoordinatesystem units   <1.5mm,2mm>
\setplotarea x from 0 to 16, y from -2 to 15
\put{1.101)} [l] at 2 12
\put {$ \scriptstyle \bullet$} [c] at 10 4
\put {$ \scriptstyle \bullet$} [c] at 12 0
\put {$ \scriptstyle \bullet$} [c] at 12 4
\put {$ \scriptstyle \bullet$} [c] at 12 8
\put {$ \scriptstyle \bullet$} [c] at 12 12
\put {$ \scriptstyle \bullet$} [c] at 14 4
\put{$\scriptstyle \bullet$} [c] at 16  0
\setlinear \plot    12 0 10 4 12 8 14 4 12 0   /
\setlinear \plot   12 0 12 12  /
\put{$840$} [c] at 13 -2
\endpicture
\end{minipage}
\begin{minipage}{4cm}
\beginpicture
\setcoordinatesystem units   <1.5mm,2mm>
\setplotarea x from 0 to 16, y from -2 to 15
\put{1.102)} [l] at 2 12
\put {$ \scriptstyle \bullet$} [c] at 10 8
\put {$ \scriptstyle \bullet$} [c] at 12 0
\put {$ \scriptstyle \bullet$} [c] at 12 4
\put {$ \scriptstyle \bullet$} [c] at 12 8
\put {$ \scriptstyle \bullet$} [c] at 12 12
\put {$ \scriptstyle \bullet$} [c] at 14 8
\put{$\scriptstyle \bullet$} [c] at 16  0
\setlinear \plot    12 12 10 8 12 4 14 8 12 12   /
\setlinear \plot   12 0 12 12  /
\put{$840$} [c] at 13 -2
\endpicture
\end{minipage}
\begin{minipage}{4cm}
\beginpicture
\setcoordinatesystem units   <1.5mm,2mm>
\setplotarea x from 0 to 16, y from -2 to 15
\put{1.103)} [l] at 2 12
\put {$ \scriptstyle \bullet$} [c] at 12 0
\put {$ \scriptstyle \bullet$} [c] at 12 4
\put {$ \scriptstyle \bullet$} [c] at 12 8
\put {$ \scriptstyle \bullet$} [c] at 12 12
\put {$ \scriptstyle \bullet$} [c] at 10 12
\put {$ \scriptstyle \bullet$} [c] at 14 12
\put{$\scriptstyle \bullet$} [c] at 16  0
\setlinear \plot  12 0 12 12     /
\setlinear \plot  10 12 12 8 14 12  /
\put{$840$} [c] at 13 -2
 \endpicture
\end{minipage}
\begin{minipage}{4cm}
\beginpicture
\setcoordinatesystem units   <1.5mm,2mm>
\setplotarea x from 0 to 16, y from -2 to 15
\put{1.104)} [l] at 2 12
\put {$ \scriptstyle \bullet$} [c] at 12 0
\put {$ \scriptstyle \bullet$} [c] at 12 4
\put {$ \scriptstyle \bullet$} [c] at 12 8
\put {$ \scriptstyle \bullet$} [c] at 12 12
\put {$ \scriptstyle \bullet$} [c] at 10 0
\put {$ \scriptstyle \bullet$} [c] at 14 0
\put{$\scriptstyle \bullet$} [c] at 16  0
\setlinear \plot  12 0 12 12     /
\setlinear \plot  10 0 12 4 14 0  /
\put{$840$} [c] at 13 -2
\endpicture
\end{minipage}
$$

$$
\begin{minipage}{4cm}
\beginpicture
\setcoordinatesystem units    <1.5mm,2mm>
\setplotarea x from 0 to 16, y from -2 to 15
\put{${\bf  23}$} [l] at 2 15
\put{1.105)} [l] at 2 12
\put {$ \scriptstyle \bullet$} [c] at 10 4
\put {$ \scriptstyle \bullet$} [c] at 10 8
\put {$ \scriptstyle \bullet$} [c] at 10 12
\put {$ \scriptstyle \bullet$} [c] at 13 0
\put {$ \scriptstyle \bullet$} [c] at 13 4
\put {$ \scriptstyle \bullet$} [c] at 13 12
\put {$ \scriptstyle \bullet$} [c] at 16 12
\setlinear \plot 16 12 13 0  13 12 10 4 13 0 /
\setlinear \plot 10 4 10 12 /
\put{$5{,}040$} [c] at 13 -2
\endpicture
\end{minipage}
\begin{minipage}{4cm}
\beginpicture
\setcoordinatesystem units    <1.5mm,2mm>
\setplotarea x from 0 to 16, y from -2 to 15
\put{1.106)} [l] at 2 12
\put {$ \scriptstyle \bullet$} [c] at 10 4
\put {$ \scriptstyle \bullet$} [c] at 10 8
\put {$ \scriptstyle \bullet$} [c] at 10 0
\put {$ \scriptstyle \bullet$} [c] at 13 0
\put {$ \scriptstyle \bullet$} [c] at 13 8
\put {$ \scriptstyle \bullet$} [c] at 13 12
\put {$ \scriptstyle \bullet$} [c] at 16 0
\setlinear \plot 16 0 13 12  13 0 10 8 13 12 /
\setlinear \plot 10 8 10 0 /
\put{$5{,}040$} [c] at 13 -2
\endpicture
\end{minipage}
\begin{minipage}{4cm}
\beginpicture
\setcoordinatesystem units    <1.5mm,2mm>
\setplotarea x from 0 to 16, y from -2 to 15
\put{1.107)} [l] at 2 12
\put {$ \scriptstyle \bullet$} [c] at 10 4
\put {$ \scriptstyle \bullet$} [c] at 10 8
\put {$ \scriptstyle \bullet$} [c] at 10 12
\put {$ \scriptstyle \bullet$} [c] at 13 0
\put {$ \scriptstyle \bullet$} [c] at 13 4
\put {$ \scriptstyle \bullet$} [c] at 13 12
\put {$ \scriptstyle \bullet$} [c] at 16 12
\setlinear \plot 13 12 13 0  10 4 10  12 13 4 /
\setlinear \plot 16 12 13 0 /
\put{$5{,}040$} [c] at 13 -2
\endpicture
\end{minipage}
\begin{minipage}{4cm}
\beginpicture
\setcoordinatesystem units    <1.5mm,2mm>
\setplotarea x from 0 to 16, y from -2 to 15
\put{1.108)} [l] at 2 12
\put {$ \scriptstyle \bullet$} [c] at 10 4
\put {$ \scriptstyle \bullet$} [c] at 10 8
\put {$ \scriptstyle \bullet$} [c] at 10 0
\put {$ \scriptstyle \bullet$} [c] at 13 0
\put {$ \scriptstyle \bullet$} [c] at 13 8
\put {$ \scriptstyle \bullet$} [c] at 13 12
\put {$ \scriptstyle \bullet$} [c] at 16 0
\setlinear \plot 13 0 13 12  10 8 10 0  13 8 /
\setlinear \plot 16 0 13 12 /
\put{$5{,}040$} [c] at 13 -2
\endpicture
\end{minipage}
\begin{minipage}{4cm}
\beginpicture
\setcoordinatesystem units   <1.5mm,2mm>
\setplotarea x from 0 to 16, y from -2 to 15
\put{1.109)} [l] at 2 12
\put {$ \scriptstyle \bullet$} [c] at  10 0
\put {$ \scriptstyle \bullet$} [c] at  10 6
\put {$ \scriptstyle \bullet$} [c] at  10 12
\put {$ \scriptstyle \bullet$} [c] at  16 0
\put {$ \scriptstyle \bullet$} [c] at  16 12
\put {$ \scriptstyle \bullet$} [c] at  12 4
\put {$ \scriptstyle \bullet$} [c] at  14 8
\setlinear \plot  10 12 10 0 16 12 16 0  /
\put{$5{,}040$} [c] at 13 -2
\endpicture
\end{minipage}
\begin{minipage}{4cm}
\beginpicture
\setcoordinatesystem units   <1.5mm,2mm>
\setplotarea x from 0 to 16, y from -2 to 15
\put{1.110)} [l] at 2 12
\put {$ \scriptstyle \bullet$} [c] at  10 0
\put {$ \scriptstyle \bullet$} [c] at  10 6
\put {$ \scriptstyle \bullet$} [c] at  10 12
\put {$ \scriptstyle \bullet$} [c] at  16 0
\put {$ \scriptstyle \bullet$} [c] at  16 12
\put {$ \scriptstyle \bullet$} [c] at  12 8
\put {$ \scriptstyle \bullet$} [c] at  14 4
\setlinear \plot  10 0 10 12 16 0 16 12  /
\put{$5{,}040$} [c] at 13 -2
\endpicture
\end{minipage}
$$
$$
\begin{minipage}{4cm}
\beginpicture
\setcoordinatesystem units   <1.5mm,2mm>
\setplotarea x from 0 to 16, y from -2 to 15
\put{1.111)} [l] at 2 12
\put {$ \scriptstyle \bullet$} [c] at  16 12
\put {$ \scriptstyle \bullet$} [c] at  16 0
\put {$ \scriptstyle \bullet$} [c] at  12 6
\put {$ \scriptstyle \bullet$} [c] at  11 0
\put {$ \scriptstyle \bullet$} [c] at  11 12
\put {$ \scriptstyle \bullet$} [c] at  10 6
\put {$ \scriptstyle \bullet$} [c] at  11.5 9
\setlinear \plot  16 0 16  12  11 0 12 6 11 12 10 6 11  0  /
\put{$5{,}040$} [c] at 13 -2
\endpicture
\end{minipage}
\begin{minipage}{4cm}
\beginpicture
\setcoordinatesystem units   <1.5mm,2mm>
\setplotarea x from 0 to 16, y from -2 to 15
\put{1.112)} [l] at 2 12
\put {$ \scriptstyle \bullet$} [c] at  16 12
\put {$ \scriptstyle \bullet$} [c] at  16 0
\put {$ \scriptstyle \bullet$} [c] at  12 6
\put {$ \scriptstyle \bullet$} [c] at  11 0
\put {$ \scriptstyle \bullet$} [c] at  11 12
\put {$ \scriptstyle \bullet$} [c] at  10 6
\put {$ \scriptstyle \bullet$} [c] at  11.5 3
\setlinear \plot  16 12 16  0  11 12 12 6 11 0 10 6 11  12  /
\put{$5{,}040$} [c] at 13 -2
\endpicture
\end{minipage}
\begin{minipage}{4cm}
\beginpicture
\setcoordinatesystem units   <1.5mm,2mm>
\setplotarea x from 0 to 16, y from -2 to 15
\put{1.113)} [l] at 2 12
\put {$ \scriptstyle \bullet$} [c] at  10 0
\put {$ \scriptstyle \bullet$} [c] at  10 6
\put {$ \scriptstyle \bullet$} [c] at  12 0
\put {$ \scriptstyle \bullet$} [c] at  12 12
\put {$ \scriptstyle \bullet$} [c] at  14 6
\put {$ \scriptstyle \bullet$} [c] at  16 6
\put {$ \scriptstyle \bullet$} [c] at  16 12
\setlinear \plot  16 12 16 6 12 0 14 6 12 12 10 6 12  0  /
\setlinear \plot  10 6 10 0 /
\put{$5{,}040$} [c] at 13 -2
\endpicture
\end{minipage}
\begin{minipage}{4cm}
\beginpicture
\setcoordinatesystem units   <1.5mm,2mm>
\setplotarea x from 0 to 16, y from -2 to 15
\put{1.114)} [l] at 2 12
\put {$ \scriptstyle \bullet$} [c] at  10 12
\put {$ \scriptstyle \bullet$} [c] at  10 6
\put {$ \scriptstyle \bullet$} [c] at  12 0
\put {$ \scriptstyle \bullet$} [c] at  12 12
\put {$ \scriptstyle \bullet$} [c] at  14 6
\put {$ \scriptstyle \bullet$} [c] at  16 6
\put {$ \scriptstyle \bullet$} [c] at  16 0
\setlinear \plot  16 0 16 6 12 12 14 6 12 0 10 6 12  12  /
\setlinear \plot  10 6 10 12 /
\put{$5{,}040$} [c] at 13 -2
\endpicture
\end{minipage}
\begin{minipage}{4cm}
\beginpicture
\setcoordinatesystem units   <1.5mm,2mm>
\setplotarea x from 0 to 16, y from -2 to 15
\put{1.115)} [l] at 2 12
\put {$ \scriptstyle \bullet$} [c] at  16 0
\put {$ \scriptstyle \bullet$} [c] at  16 6
\put {$ \scriptstyle \bullet$} [c] at  16 12
\put {$ \scriptstyle \bullet$} [c] at  12 0
\put {$ \scriptstyle \bullet$} [c] at  12 6
\put {$ \scriptstyle \bullet$} [c] at  12 12
\put {$ \scriptstyle \bullet$} [c] at  10 12
\setlinear \plot  10 12 12 0 16 12 16 0  /
\setlinear \plot  12 0 12 12 16 6 /
\put{$5{,}040$} [c] at 13 -2
\endpicture
\end{minipage}
\begin{minipage}{4cm}
\beginpicture
\setcoordinatesystem units   <1.5mm,2mm>
\setplotarea x from 0 to 16, y from -2 to 15
\put{1.116)} [l] at 2 12
\put {$ \scriptstyle \bullet$} [c] at  16 0
\put {$ \scriptstyle \bullet$} [c] at  16 6
\put {$ \scriptstyle \bullet$} [c] at  16 12
\put {$ \scriptstyle \bullet$} [c] at  12 0
\put {$ \scriptstyle \bullet$} [c] at  12 6
\put {$ \scriptstyle \bullet$} [c] at  12 12
\put {$ \scriptstyle \bullet$} [c] at  10 0
\setlinear \plot  10 0 12 12 16 0 16 12  /
\setlinear \plot  12 12 12 0 16 6 /
\put{$5{,}040$} [c] at 13 -2
\endpicture
\end{minipage}
$$

$$
\begin{minipage}{4cm}
\beginpicture
\setcoordinatesystem units   <1.5mm,2mm>
\setplotarea x from 0 to 16, y from -2 to 15
\put{1.117)} [l] at 2 12
\put {$ \scriptstyle \bullet$} [c] at  16 12
\put {$ \scriptstyle \bullet$} [c] at  14 12
\put {$ \scriptstyle \bullet$} [c] at  14 0
\put {$ \scriptstyle \bullet$} [c] at  10 0
\put {$ \scriptstyle \bullet$} [c] at  10 4
\put {$ \scriptstyle \bullet$} [c] at  10 8
\put {$ \scriptstyle \bullet$} [c] at  10 12
\setlinear \plot  10 12  10 0 14 12  14 0 16 12  /
\put{$5{,}040$} [c] at 13 -2
\endpicture
\end{minipage}
\begin{minipage}{4cm}
\beginpicture
\setcoordinatesystem units   <1.5mm,2mm>
\setplotarea x from 0 to 16, y from -2 to 15
\put{1.118)} [l] at 2 12
\put {$ \scriptstyle \bullet$} [c] at  16 0
\put {$ \scriptstyle \bullet$} [c] at  14 12
\put {$ \scriptstyle \bullet$} [c] at  14 0
\put {$ \scriptstyle \bullet$} [c] at  10 0
\put {$ \scriptstyle \bullet$} [c] at  10 4
\put {$ \scriptstyle \bullet$} [c] at  10 8
\put {$ \scriptstyle \bullet$} [c] at  10 12
\setlinear \plot  10 0  10 12 14 0  14 12 16 0  /
\put{$5{,}040$} [c] at 13 -2
\endpicture
\end{minipage}
\begin{minipage}{4cm}
\beginpicture
\setcoordinatesystem units   <1.5mm,2mm>
\setplotarea x from 0 to 16, y from -2 to 15
\put{1.119)} [l] at 2 12
\put {$ \scriptstyle \bullet$} [c] at  10 12
\put {$ \scriptstyle \bullet$} [c] at  12 0
\put {$ \scriptstyle \bullet$} [c] at  12 4
\put {$ \scriptstyle \bullet$} [c] at  12 8
\put {$ \scriptstyle \bullet$} [c] at  12 12
\put {$ \scriptstyle \bullet$} [c] at  16 0
\put {$ \scriptstyle \bullet$} [c] at  16 12
\setlinear \plot  10 12  12 0 12 12   /
\setlinear \plot  16 12 16 0  12 8  /
\put{$5{,}040$} [c] at 13 -2
\endpicture
\end{minipage}
\begin{minipage}{4cm}
\beginpicture
\setcoordinatesystem units   <1.5mm,2mm>
\setplotarea x from 0 to 16, y from -2 to 15
\put{1.120)} [l] at 2 12
\put {$ \scriptstyle \bullet$} [c] at  10 0
\put {$ \scriptstyle \bullet$} [c] at  12 0
\put {$ \scriptstyle \bullet$} [c] at  12 4
\put {$ \scriptstyle \bullet$} [c] at  12 8
\put {$ \scriptstyle \bullet$} [c] at  12 12
\put {$ \scriptstyle \bullet$} [c] at  16 0
\put {$ \scriptstyle \bullet$} [c] at  16 12
\setlinear \plot  10 0  12 12 12 0   /
\setlinear \plot  16 0 16 12  12 4  /
\put{$5{,}040$} [c] at 13 -2
\endpicture
\end{minipage}
\begin{minipage}{4cm}
\beginpicture
\setcoordinatesystem units   <1.5mm,2mm>
\setplotarea x from 0 to 16, y from -2 to 15
\put{1.121)} [l] at 2 12
\put {$ \scriptstyle \bullet$} [c] at  10 0
\put {$ \scriptstyle \bullet$} [c] at  10 6
\put {$ \scriptstyle \bullet$} [c] at  10 12
\put {$ \scriptstyle \bullet$} [c] at  14  0
\put {$ \scriptstyle \bullet$} [c] at  14 6
\put {$ \scriptstyle \bullet$} [c] at  14 12
\put {$ \scriptstyle \bullet$} [c] at  16 12
\setlinear \plot  10 12 10 0  14 12 14 6  10 12   /
\setlinear \plot  16 12 14 0 14 6 /
\put{$5{,}040$} [c] at 13 -2
\endpicture
\end{minipage}
\begin{minipage}{4cm}
\beginpicture
\setcoordinatesystem units   <1.5mm,2mm>
\setplotarea x from 0 to 16, y from -2 to 15
\put{1.122)} [l] at 2 12
\put {$ \scriptstyle \bullet$} [c] at  10 0
\put {$ \scriptstyle \bullet$} [c] at  10 6
\put {$ \scriptstyle \bullet$} [c] at  10 12
\put {$ \scriptstyle \bullet$} [c] at  14  0
\put {$ \scriptstyle \bullet$} [c] at  14 6
\put {$ \scriptstyle \bullet$} [c] at  14 12
\put {$ \scriptstyle \bullet$} [c] at  16 0
\setlinear \plot  10 0 10 12  14 0 14 6  10 0   /
\setlinear \plot  16 0 14 12 14 6 /
\put{$5{,}040$} [c] at 13 -2
\endpicture
\end{minipage}
$$
$$
\begin{minipage}{4cm}
\beginpicture
\setcoordinatesystem units   <1.5mm,2mm>
\setplotarea x from 0 to 16, y from -2 to 15
\put{1.123)} [l] at 2 12
\put {$ \scriptstyle \bullet$} [c] at  10 12
\put {$ \scriptstyle \bullet$} [c] at  10 6
\put {$ \scriptstyle \bullet$} [c] at  12 12
\put {$ \scriptstyle \bullet$} [c] at  12 0
\put {$ \scriptstyle \bullet$} [c] at  14 6
\put {$ \scriptstyle \bullet$} [c] at  16 0
\put {$ \scriptstyle \bullet$} [c] at  16 12
\setlinear \plot  10 12 10 6 12 0 14 6 12 12 10 6   /
\setlinear \plot  16 0 16  12 14 6  /
\put{$5{,}040$} [c] at 13 -2
\endpicture
\end{minipage}
\begin{minipage}{4cm}
\beginpicture
\setcoordinatesystem units   <1.5mm,2mm>
\setplotarea x from 0 to 16, y from -2 to 15
\put{1.124)} [l] at 2 12
\put {$ \scriptstyle \bullet$} [c] at  10 0
\put {$ \scriptstyle \bullet$} [c] at  10 6
\put {$ \scriptstyle \bullet$} [c] at  12 12
\put {$ \scriptstyle \bullet$} [c] at  12 0
\put {$ \scriptstyle \bullet$} [c] at  14 6
\put {$ \scriptstyle \bullet$} [c] at  16 0
\put {$ \scriptstyle \bullet$} [c] at  16 12
\setlinear \plot  10 0 10 6 12 12 14 6 12 0 10 6   /
\setlinear \plot  16 12 16  0 14 6  /
\put{$5{,}040$} [c] at 13 -2
\endpicture
\end{minipage}
\begin{minipage}{4cm}
\beginpicture
\setcoordinatesystem units   <1.5mm,2mm>
\setplotarea x from 0 to 16, y from -2 to 15
\put{1.125)} [l] at 2 12
\put {$ \scriptstyle \bullet$} [c] at  10 6
\put {$ \scriptstyle \bullet$} [c] at  11.5 0
\put {$ \scriptstyle \bullet$} [c] at  11.5 12
\put {$ \scriptstyle \bullet$} [c] at  13 6
\put {$ \scriptstyle \bullet$} [c] at  13 12
\put {$ \scriptstyle \bullet$} [c] at  16 0
\put {$ \scriptstyle \bullet$} [c] at  16 12
\setlinear \plot  11.5 0 16 12 16 0 13 12   13 6 11.5 0 10 6 11.5 12 13 6   /
\put{$5{,}040$} [c] at 13 -2
\endpicture
\end{minipage}
\begin{minipage}{4cm}
\beginpicture
\setcoordinatesystem units   <1.5mm,2mm>
\setplotarea x from 0 to 16, y from -2 to 15
\put{1.126)} [l] at 2 12
\put {$ \scriptstyle \bullet$} [c] at  10 6
\put {$ \scriptstyle \bullet$} [c] at  11.5 0
\put {$ \scriptstyle \bullet$} [c] at  11.5 12
\put {$ \scriptstyle \bullet$} [c] at  13 6
\put {$ \scriptstyle \bullet$} [c] at  13 0
\put {$ \scriptstyle \bullet$} [c] at  16 0
\put {$ \scriptstyle \bullet$} [c] at  16 12
\setlinear \plot  11.5 12 16 0 16 12 13 0   13 6 11.5 12 10 6 11.5 0 13 6   /
\put{$5{,}040$} [c] at 13 -2
\endpicture
\end{minipage}
\begin{minipage}{4cm}
\beginpicture
\setcoordinatesystem units   <1.5mm,2mm>
\setplotarea x from 0 to 16, y from -2 to 15
\put{1.127)} [l] at 2 12
\put {$ \scriptstyle \bullet$} [c] at  10 12
\put {$ \scriptstyle \bullet$} [c] at  12 0
\put {$ \scriptstyle \bullet$} [c] at  12 6
\put {$ \scriptstyle \bullet$} [c] at  12 12
\put {$ \scriptstyle \bullet$} [c] at  16 0
\put {$ \scriptstyle \bullet$} [c] at  16 6
\put {$ \scriptstyle \bullet$} [c] at  16 12
\setlinear \plot  10 12  12 0 16 12  16 0 12 6 12 12   /
\setlinear \plot  12 6 12 0  /
\put{$5{,}040$} [c] at 13 -2
\endpicture
\end{minipage}
\begin{minipage}{4cm}
\beginpicture
\setcoordinatesystem units   <1.5mm,2mm>
\setplotarea x from 0 to 16, y from -2 to 15
\put{1.128)} [l] at 2 12
\put {$ \scriptstyle \bullet$} [c] at  10 0
\put {$ \scriptstyle \bullet$} [c] at  12 0
\put {$ \scriptstyle \bullet$} [c] at  12 6
\put {$ \scriptstyle \bullet$} [c] at  12 12
\put {$ \scriptstyle \bullet$} [c] at  16 0
\put {$ \scriptstyle \bullet$} [c] at  16 6
\put {$ \scriptstyle \bullet$} [c] at  16 12
\setlinear \plot  10 0  12 12 16  0  16 12 12 6 12 0   /
\setlinear \plot  12 6 12 12  /
\put{$5{,}040$} [c] at 13 -2
\endpicture
\end{minipage}
$$

$$
\begin{minipage}{4cm}
\beginpicture
\setcoordinatesystem units   <1.5mm,2mm>
\setplotarea x from 0 to 16, y from -2 to 15
\put{1.129)} [l] at 2 12
\put {$ \scriptstyle \bullet$} [c] at  10 0
\put {$ \scriptstyle \bullet$} [c] at  10 6
\put {$ \scriptstyle \bullet$} [c] at  10 12
\put {$ \scriptstyle \bullet$} [c] at  13 12
\put {$ \scriptstyle \bullet$} [c] at  16 0
\put {$ \scriptstyle \bullet$} [c] at  16 6
\put {$ \scriptstyle \bullet$} [c] at  16 12
\setlinear \plot  10 12  10 0 13  12  16 0 16 12    /
\setlinear \plot  10 0  16 6  /
\put{$5{,}040$} [c] at 13 -2
\endpicture
\end{minipage}
\begin{minipage}{4cm}
\beginpicture
\setcoordinatesystem units   <1.5mm,2mm>
\setplotarea x from 0 to 16, y from -2 to 15
\put{1.130)} [l] at 2 12
\put {$ \scriptstyle \bullet$} [c] at  10 0
\put {$ \scriptstyle \bullet$} [c] at  10 6
\put {$ \scriptstyle \bullet$} [c] at  10 12
\put {$ \scriptstyle \bullet$} [c] at  13 0
\put {$ \scriptstyle \bullet$} [c] at  16 0
\put {$ \scriptstyle \bullet$} [c] at  16 6
\put {$ \scriptstyle \bullet$} [c] at  16 12
\setlinear \plot  10 0  10 12 13  0  16 12 16 0    /
\setlinear \plot  10 12  16 6  /
\put{$5{,}040$} [c] at 13 -2
\endpicture
\end{minipage}
\begin{minipage}{4cm}
\beginpicture
\setcoordinatesystem units   <1.5mm,2mm>
\setplotarea x from 0 to 16, y from -2 to 15
\put{1.131)} [l] at 2 12
\put {$ \scriptstyle \bullet$} [c] at  10 12
\put {$ \scriptstyle \bullet$} [c] at  11 0
\put {$ \scriptstyle \bullet$} [c] at  11 12
\put {$ \scriptstyle \bullet$} [c] at  14  6
\put {$ \scriptstyle \bullet$} [c] at  15 0
\put {$ \scriptstyle \bullet$} [c] at  15 12
\put {$ \scriptstyle \bullet$} [c] at  16 6
\setlinear \plot  10 12  11 0 11 12 14 6 15 0  16 6 15 12 14 6   /
\setlinear \plot  11 0 15 12 /
\put{$5{,}040$} [c] at 13 -2
\endpicture
\end{minipage}
\begin{minipage}{4cm}
\beginpicture
\setcoordinatesystem units   <1.5mm,2mm>
\setplotarea x from 0 to 16, y from -2 to 15
\put{1.132)} [l] at 2 12
\put {$ \scriptstyle \bullet$} [c] at  10 0
\put {$ \scriptstyle \bullet$} [c] at  11 0
\put {$ \scriptstyle \bullet$} [c] at  11 12
\put {$ \scriptstyle \bullet$} [c] at  14  6
\put {$ \scriptstyle \bullet$} [c] at  15 0
\put {$ \scriptstyle \bullet$} [c] at  15 12
\put {$ \scriptstyle \bullet$} [c] at  16 6
\setlinear \plot  10 0  11 12 11 0 14 6 15 12  16 6 15 0 14 6   /
\setlinear \plot  11 12 15 0 /
\put{$5{,}040$} [c] at 13 -2
\endpicture
\end{minipage}
\begin{minipage}{4cm}
\beginpicture
\setcoordinatesystem units   <1.5mm,2mm>
\setplotarea x from 0 to 16, y from -2 to 15
\put{1.133)} [l] at 2 12
\put {$ \scriptstyle \bullet$} [c] at  10 12
\put {$ \scriptstyle \bullet$} [c] at  12 0
\put {$ \scriptstyle \bullet$} [c] at  12 6
\put {$ \scriptstyle \bullet$} [c] at  12 12
\put {$ \scriptstyle \bullet$} [c] at  13 3
\put {$ \scriptstyle \bullet$} [c] at  16 0
\put {$ \scriptstyle \bullet$} [c] at  16 12
\setlinear \plot  12 12 12 0 16 12 16 0 12 12  /
\setlinear \plot   10 12 12 6  /
\put{$5{,}040$} [c] at 13 -2
\endpicture
\end{minipage}
\begin{minipage}{4cm}
\beginpicture
\setcoordinatesystem units   <1.5mm,2mm>
\setplotarea x from 0 to 16, y from -2 to 15
\put{1.134)} [l] at 2 12
\put {$ \scriptstyle \bullet$} [c] at  10 0
\put {$ \scriptstyle \bullet$} [c] at  12 0
\put {$ \scriptstyle \bullet$} [c] at  12 6
\put {$ \scriptstyle \bullet$} [c] at  12 12
\put {$ \scriptstyle \bullet$} [c] at  13 9
\put {$ \scriptstyle \bullet$} [c] at  16 0
\put {$ \scriptstyle \bullet$} [c] at  16 12
\setlinear \plot  12 0 12 12 16 0 16 12 12 0  /
\setlinear \plot   10 0 12 6  /
\put{$5{,}040$} [c] at 13 -2
\endpicture
\end{minipage}
$$
$$
\begin{minipage}{4cm}
\beginpicture
\setcoordinatesystem units   <1.5mm,2mm>
\setplotarea x from 0 to 16, y from -2 to 15
\put{1.135)} [l] at 2 12
\put {$ \scriptstyle \bullet$} [c] at  10 0
\put {$ \scriptstyle \bullet$} [c] at  10 4
\put {$ \scriptstyle \bullet$} [c] at  10 8
\put {$ \scriptstyle \bullet$} [c] at  10 12
\put {$ \scriptstyle \bullet$} [c] at  13 12
\put {$ \scriptstyle \bullet$} [c] at  16 0
\put {$ \scriptstyle \bullet$} [c] at  16 12
\setlinear \plot  10 12 10 0 16 12 16 0  /
\setlinear \plot  10 4 13 12  /
\put{$5{,}040$} [c] at 13 -2
\endpicture
\end{minipage}
\begin{minipage}{4cm}
\beginpicture
\setcoordinatesystem units   <1.5mm,2mm>
\setplotarea x from 0 to 16, y from -2 to 15
\put{1.136)} [l] at 2 12
\put {$ \scriptstyle \bullet$} [c] at  10 0
\put {$ \scriptstyle \bullet$} [c] at  10 4
\put {$ \scriptstyle \bullet$} [c] at  10 8
\put {$ \scriptstyle \bullet$} [c] at  10 12
\put {$ \scriptstyle \bullet$} [c] at  13 0
\put {$ \scriptstyle \bullet$} [c] at  16 0
\put {$ \scriptstyle \bullet$} [c] at  16 12
\setlinear \plot  10 0 10 12 16 0 16 12  /
\setlinear \plot  10 8 13 0  /
\put{$5{,}040$} [c] at 13 -2
\endpicture
\end{minipage}
\begin{minipage}{4cm}
\beginpicture
\setcoordinatesystem units   <1.5mm,2mm>
\setplotarea x from 0 to 16, y from -2 to 15
\put{1.137)} [l] at 2 12
\put {$ \scriptstyle \bullet$} [c] at  10  0
\put {$ \scriptstyle \bullet$} [c] at  10 12
\put {$ \scriptstyle \bullet$} [c] at  13 0
\put {$ \scriptstyle \bullet$} [c] at  13 12
\put {$ \scriptstyle \bullet$} [c] at  13 6
\put {$ \scriptstyle \bullet$} [c] at  16  0
\put {$ \scriptstyle \bullet$} [c] at  16 12
\setlinear \plot  10 0  10 12 13 0 13 12 16 0 16 12 13 6   /
\setlinear \plot  10 0 13 12   /
\put{$5{,}040$} [c] at 13 -2
\endpicture
\end{minipage}
\begin{minipage}{4cm}
\beginpicture
\setcoordinatesystem units   <1.5mm,2mm>
\setplotarea x from 0 to 16, y from -2 to 15
\put{1.138)} [l] at 2 12
\put {$ \scriptstyle \bullet$} [c] at  10  0
\put {$ \scriptstyle \bullet$} [c] at  10 12
\put {$ \scriptstyle \bullet$} [c] at  13 0
\put {$ \scriptstyle \bullet$} [c] at  13 12
\put {$ \scriptstyle \bullet$} [c] at  13 6
\put {$ \scriptstyle \bullet$} [c] at  16  0
\put {$ \scriptstyle \bullet$} [c] at  16 12
\setlinear \plot  10 12  10 0 13 12 13 0 16 12 16 0 13 6   /
\setlinear \plot  10 12 13 0   /
\put{$5{,}040$} [c] at 13 -2
\endpicture
\end{minipage}
\begin{minipage}{4cm}
\beginpicture
\setcoordinatesystem units    <1.5mm,2mm>
\setplotarea x from 0 to 16, y from -2 to 15
\put{1.139)} [l] at 2 12
\put {$ \scriptstyle \bullet$} [c] at 10 6
\put {$ \scriptstyle \bullet$} [c] at 13 6
\put {$ \scriptstyle \bullet$} [c] at 16 6
\put {$ \scriptstyle \bullet$} [c] at 13 0
\put {$ \scriptstyle \bullet$} [c] at 13 12
\put {$ \scriptstyle \bullet$} [c] at 14.5 9
\put {$ \scriptstyle \bullet$} [c] at 16 12
\setlinear \plot 16 12 16 6 13 0 10 6 13 12 16  6   /
\setlinear \plot 13 0 13 12    /
\put{$2{,}520$} [c] at 13 -2
\endpicture
\end{minipage}
\begin{minipage}{4cm}
\beginpicture
\setcoordinatesystem units    <1.5mm,2mm>
\setplotarea x from 0 to 16, y from -2 to 15
\put{1.140)} [l] at 2 12
\put {$ \scriptstyle \bullet$} [c] at 10 6
\put {$ \scriptstyle \bullet$} [c] at 13 6
\put {$ \scriptstyle \bullet$} [c] at 16 6
\put {$ \scriptstyle \bullet$} [c] at 13 0
\put {$ \scriptstyle \bullet$} [c] at 13 12
\put {$ \scriptstyle \bullet$} [c] at 14.5 3
\put {$ \scriptstyle \bullet$} [c] at 16 0
\setlinear \plot 16 0 16 6 13 12 10 6 13 0 16  6   /
\setlinear \plot 13 0 13 12    /
\put{$2{,}520$} [c] at 13 -2
\endpicture
\end{minipage}
$$

$$
\begin{minipage}{4cm}
\beginpicture
\setcoordinatesystem units    <1.5mm,2mm>
\setplotarea x from 0 to 16, y from -2 to 15
\put{1.141)} [l] at 2 12
\put {$ \scriptstyle \bullet$} [c] at 12 0
\put {$ \scriptstyle \bullet$} [c] at 10 4
\put {$ \scriptstyle \bullet$} [c] at 11.5 4
\put {$ \scriptstyle \bullet$} [c] at 12.5 4
\put {$ \scriptstyle \bullet$} [c] at 12 8
\put {$ \scriptstyle \bullet$} [c] at 12 12
\put {$ \scriptstyle \bullet$} [c] at 16 12
\setlinear \plot 16 12 12 0 10  4 12 12  12 8 11.5 4 12 0 12.5 4 12 8 /
\put{$2{,}520$} [c] at 13 -2
\endpicture
\end{minipage}
\begin{minipage}{4cm}
\beginpicture
\setcoordinatesystem units    <1.5mm,2mm>
\setplotarea x from 0 to 16, y from -2 to 15
\put{1.142)} [l] at 2 12
\put {$ \scriptstyle \bullet$} [c] at 12 0
\put {$ \scriptstyle \bullet$} [c] at 10 8
\put {$ \scriptstyle \bullet$} [c] at 11.5 8
\put {$ \scriptstyle \bullet$} [c] at 12.5 8
\put {$ \scriptstyle \bullet$} [c] at 12 4
\put {$ \scriptstyle \bullet$} [c] at 12 12
\put {$ \scriptstyle \bullet$} [c] at 16 0
\setlinear \plot 16 0 12 12 10  8 12 0  12 4 11.5 8 12 12 12.5 8 12 4 /
\put{$2{,}520$} [c] at 13 -2
\endpicture
\end{minipage}
\begin{minipage}{4cm}
\beginpicture
\setcoordinatesystem units    <1.5mm,2mm>
\setplotarea x from 0 to 16, y from -2 to 15
\put{1.143)} [l] at 2 12
\put {$ \scriptstyle \bullet$} [c] at 12 0
\put {$ \scriptstyle \bullet$} [c] at 10 8
\put {$ \scriptstyle \bullet$} [c] at 11.5 8
\put {$ \scriptstyle \bullet$} [c] at 12.5 8
\put {$ \scriptstyle \bullet$} [c] at 12 4
\put {$ \scriptstyle \bullet$} [c] at 12 12
\put {$ \scriptstyle \bullet$} [c] at 16 12
\setlinear \plot 16 12 12 0 12  4 12.5 8 12 12 11.5 8 12 4 /
\setlinear \plot  12 0 10  8 12 12 /
\put{$2{,}520$} [c] at 13 -2
\endpicture
\end{minipage}
\begin{minipage}{4cm}
\beginpicture
\setcoordinatesystem units    <1.5mm,2mm>
\setplotarea x from 0 to 16, y from -2 to 15
\put{1.144)} [l] at 2 12
\put {$ \scriptstyle \bullet$} [c] at 12 0
\put {$ \scriptstyle \bullet$} [c] at 10 4
\put {$ \scriptstyle \bullet$} [c] at 11.5 4
\put {$ \scriptstyle \bullet$} [c] at 12.5 4
\put {$ \scriptstyle \bullet$} [c] at 12 8
\put {$ \scriptstyle \bullet$} [c] at 12 12
\put {$ \scriptstyle \bullet$} [c] at 16 0
\setlinear \plot 16 0 12 12 12  8 12.5 4 12 0 11.5 4 12 8 /
\setlinear \plot  12 12 10  4 12 0 /
\put{$2{,}520$} [c] at 13 -2
\endpicture
\end{minipage}
\begin{minipage}{4cm}
\beginpicture
\setcoordinatesystem units    <1.5mm,2mm>
\setplotarea x from 0 to 16, y from -2 to 15
\put{1.145)} [l] at 2 12
\put {$ \scriptstyle \bullet$} [c] at 10 12
\put {$ \scriptstyle \bullet$} [c] at 13 4
\put {$ \scriptstyle \bullet$} [c] at 13 0
\put {$ \scriptstyle \bullet$} [c] at 13 12
\put {$ \scriptstyle \bullet$} [c] at 14.5 8
\put {$ \scriptstyle \bullet$} [c] at 16 4
\put {$ \scriptstyle \bullet$} [c] at 16 12
\setlinear \plot  13  4 16 12 16 4 13 0 13 12  /
\setlinear \plot 10 12 13 4   /
\put{$2{,}520$} [c] at 13 -2
\endpicture
\end{minipage}
\begin{minipage}{4cm}
\beginpicture
\setcoordinatesystem units    <1.5mm,2mm>
\setplotarea x from 0 to 16, y from -2 to 15
\put{1.146)} [l] at 2 12
\put {$ \scriptstyle \bullet$} [c] at 10 0
\put {$ \scriptstyle \bullet$} [c] at 13 8
\put {$ \scriptstyle \bullet$} [c] at 13 0
\put {$ \scriptstyle \bullet$} [c] at 13 12
\put {$ \scriptstyle \bullet$} [c] at 14.5 4
\put {$ \scriptstyle \bullet$} [c] at 16 8
\put {$ \scriptstyle \bullet$} [c] at 16 0
\setlinear \plot  13  8 16 0 16 8 13 12 13 0  /
\setlinear \plot 10 0 13 8   /
\put{$2{,}520$} [c] at 13 -2
\endpicture
\end{minipage}
$$
$$
\begin{minipage}{4cm}
\beginpicture
\setcoordinatesystem units    <1.5mm,2mm>
\setplotarea x from 0 to 16, y from -2 to 15
\put{1.147)} [l] at 2 12
\put {$ \scriptstyle \bullet$} [c] at 13 0
\put {$ \scriptstyle \bullet$} [c] at 13 4
\put {$ \scriptstyle \bullet$} [c] at 13 12
\put {$ \scriptstyle \bullet$} [c] at 10 4
\put {$ \scriptstyle \bullet$} [c] at 10  12
\put {$ \scriptstyle \bullet$} [c] at 16 12
\put {$ \scriptstyle \bullet$} [c] at 16 4
\setlinear \plot 13 12 13  0  10 4 10 12 13 4 16 12 16 4 13 0  /
\put{$2{,}520$} [c] at 13 -2
\endpicture
\end{minipage}
\begin{minipage}{4cm}
\beginpicture
\setcoordinatesystem units    <1.5mm,2mm>
\setplotarea x from 0 to 16, y from -2 to 15
\put{1.148)} [l] at 2 12
\put {$ \scriptstyle \bullet$} [c] at 13 0
\put {$ \scriptstyle \bullet$} [c] at 13 8
\put {$ \scriptstyle \bullet$} [c] at 13 12
\put {$ \scriptstyle \bullet$} [c] at 10 8
\put {$ \scriptstyle \bullet$} [c] at 10  0
\put {$ \scriptstyle \bullet$} [c] at 16 0
\put {$ \scriptstyle \bullet$} [c] at 16 8
\setlinear \plot 13 0 13  12 10 8 10 0 13 8 16 0 16 8 13 12  /
\put{$2{,}520$} [c] at 13 -2
\endpicture
\end{minipage}
\begin{minipage}{4cm}
\beginpicture
\setcoordinatesystem units    <1.5mm,2mm>
\setplotarea x from 0 to 16, y from -2 to 15
\put{1.149)} [l] at 2 12
\put {$ \scriptstyle \bullet$} [c] at 13 0
\put {$ \scriptstyle \bullet$} [c] at 13 4
\put {$ \scriptstyle \bullet$} [c] at 13 12
\put {$ \scriptstyle \bullet$} [c] at 10 4
\put {$ \scriptstyle \bullet$} [c] at 10  12
\put {$ \scriptstyle \bullet$} [c] at 16 12
\put {$ \scriptstyle \bullet$} [c] at 16 4
\setlinear \plot 10 12 10 4  13 12 13 0 10 4  /
\setlinear \plot 13 12 16 4  13 0  /
\setlinear \plot 16 12 16 4   /
\put{$2{,}520$} [c] at 13 -2
\endpicture
\end{minipage}
\begin{minipage}{4cm}
\beginpicture
\setcoordinatesystem units    <1.5mm,2mm>
\setplotarea x from 0 to 16, y from -2 to 15
\put{1.150)} [l] at 2 12
\put {$ \scriptstyle \bullet$} [c] at 13 0
\put {$ \scriptstyle \bullet$} [c] at 13 8
\put {$ \scriptstyle \bullet$} [c] at 13 12
\put {$ \scriptstyle \bullet$} [c] at 10 8
\put {$ \scriptstyle \bullet$} [c] at 10  0
\put {$ \scriptstyle \bullet$} [c] at 16 0
\put {$ \scriptstyle \bullet$} [c] at 16 8
\setlinear \plot 10 0 10 8  13 0 13 12 10 8  /
\setlinear \plot 13 0 16 8  13 12  /
\setlinear \plot 16 0 16 8   /
\put{$2{,}520$} [c] at 13 -2
\endpicture
\end{minipage}
\begin{minipage}{4cm}
\beginpicture
\setcoordinatesystem units    <1.5mm,2mm>
\setplotarea x from 0 to 16, y from -2 to 15
\put{1.151)} [l] at 2 12
\put {$ \scriptstyle \bullet$} [c] at 10 12
\put {$ \scriptstyle \bullet$} [c] at 13 12
\put {$ \scriptstyle \bullet$} [c] at 12.5 8
\put {$ \scriptstyle \bullet$} [c] at 13.5 8
\put {$ \scriptstyle \bullet$} [c] at 13 4
\put {$ \scriptstyle \bullet$} [c] at 13 0
\put {$ \scriptstyle \bullet$} [c] at 16 12
\setlinear \plot 16 12 13  0 13 4 13.5 8 13 12 12.5 8 13 4  10 12  /
\put{$2{,}520$} [c] at 13 -2
\endpicture
\end{minipage}
\begin{minipage}{4cm}
\beginpicture
\setcoordinatesystem units    <1.5mm,2mm>
\setplotarea x from 0 to 16, y from -2 to 15
\put{1.152)} [l] at 2 12
\put {$ \scriptstyle \bullet$} [c] at 10 0
\put {$ \scriptstyle \bullet$} [c] at 13 12
\put {$ \scriptstyle \bullet$} [c] at 12.5 4
\put {$ \scriptstyle \bullet$} [c] at 13.5 4
\put {$ \scriptstyle \bullet$} [c] at 13 8
\put {$ \scriptstyle \bullet$} [c] at 13 0
\put {$ \scriptstyle \bullet$} [c] at 16 0
\setlinear \plot 16 0 13  12 13 8 13.5 4 13 0 12.5 4 13 8  10 0  /
\put{$2{,}520$} [c] at 13 -2
\endpicture
\end{minipage}
$$

$$
\begin{minipage}{4cm}
\beginpicture
\setcoordinatesystem units    <1.5mm,2mm>
\setplotarea x from 0 to 16, y from -2 to 15
\put{1.153)} [l] at 2 12
\put {$ \scriptstyle \bullet$} [c] at 10 12
\put {$ \scriptstyle \bullet$} [c] at 12 12
\put {$ \scriptstyle \bullet$} [c] at 14 12
\put {$ \scriptstyle \bullet$} [c] at 13 8
\put {$ \scriptstyle \bullet$} [c] at 13 4
\put {$ \scriptstyle \bullet$} [c] at 13 0
\put {$ \scriptstyle \bullet$} [c] at 16 12
\setlinear \plot 16 12 13  0 13 8 14 12 /
\setlinear \plot 10 12  13 4   /
\setlinear \plot 12 12  13 8   /
\put{$2{,}520$} [c] at 13 -2
\endpicture
\end{minipage}
\begin{minipage}{4cm}
\beginpicture
\setcoordinatesystem units    <1.5mm,2mm>
\setplotarea x from 0 to 16, y from -2 to 15
\put{1.154)} [l] at 2 12
\put {$ \scriptstyle \bullet$} [c] at 10 0
\put {$ \scriptstyle \bullet$} [c] at 12 0
\put {$ \scriptstyle \bullet$} [c] at 14 0
\put {$ \scriptstyle \bullet$} [c] at 13 12
\put {$ \scriptstyle \bullet$} [c] at 13 8
\put {$ \scriptstyle \bullet$} [c] at 13 4
\put {$ \scriptstyle \bullet$} [c] at 16 0
\setlinear \plot 16 0 13  12 13 4 14 0 /
\setlinear \plot 10 0  13 8   /
\setlinear \plot 12 0  13 4   /
\put{$2{,}520$} [c] at 13 -2
\endpicture
\end{minipage}
\begin{minipage}{4cm}
\beginpicture
\setcoordinatesystem units    <1.5mm,2mm>
\setplotarea x from 0 to 16, y from -2 to 15
\put{1.155)} [l] at 2 12
\put {$ \scriptstyle \bullet$} [c] at 10 0
\put {$ \scriptstyle \bullet$} [c] at 12 6
\put {$ \scriptstyle \bullet$} [c] at 12 12
\put {$ \scriptstyle \bullet$} [c] at 14 0
\put {$ \scriptstyle \bullet$} [c] at 14 6
\put {$ \scriptstyle \bullet$} [c] at 14 12
\put {$ \scriptstyle \bullet$} [c] at 16 6
\setlinear \plot 12 12  14 6   /
\setlinear \plot 14 12  12 6   /
\setlinear \plot 14 0  14 12  16 6 14 0 12 6 12 12 10 0  /
\put{$2{,}520$} [c] at 13 -2
\endpicture
\end{minipage}
\begin{minipage}{4cm}
\beginpicture
\setcoordinatesystem units    <1.5mm,2mm>
\setplotarea x from 0 to 16, y from -2 to 15
\put{1.156)} [l] at 2 12
\put {$ \scriptstyle \bullet$} [c] at 10 12
\put {$ \scriptstyle \bullet$} [c] at 12 6
\put {$ \scriptstyle \bullet$} [c] at 12 0
\put {$ \scriptstyle \bullet$} [c] at 14 0
\put {$ \scriptstyle \bullet$} [c] at 14 6
\put {$ \scriptstyle \bullet$} [c] at 14 12
\put {$ \scriptstyle \bullet$} [c] at 16 6
\setlinear \plot 12 0  14 6   /
\setlinear \plot 14 0  12 6   /
\setlinear \plot 14 12  14 0  16 6 14 12 12 6 12 0 10 12  /
\put{$2{,}520$} [c] at 13 -2
\endpicture
\end{minipage}
\begin{minipage}{4cm}
\beginpicture
\setcoordinatesystem units    <1.5mm,2mm>
\setplotarea x from 0 to 16, y from -2 to 15
\put{1.157)} [l] at 2 12
\put {$ \scriptstyle \bullet$} [c] at 10 12
\put {$ \scriptstyle \bullet$} [c] at 10 6
\put {$ \scriptstyle \bullet$} [c] at 12 0
\put {$ \scriptstyle \bullet$} [c] at 12 6
\put {$ \scriptstyle \bullet$} [c] at 14 6
\put {$ \scriptstyle \bullet$} [c] at 14 12
\put {$ \scriptstyle \bullet$} [c] at 16 0
\setlinear \plot 10 12 16 0 14  12 14 6 12 0 10 6 10 12 12 6 12 0  /
\setlinear \plot 14 12 12 6  /
\put{$2{,}520$} [c] at 13 -2
\endpicture
\end{minipage}
\begin{minipage}{4cm}
\beginpicture
\setcoordinatesystem units    <1.5mm,2mm>
\setplotarea x from 0 to 16, y from -2 to 15
\put{1.158)} [l] at 2 12
\put {$ \scriptstyle \bullet$} [c] at 10 0
\put {$ \scriptstyle \bullet$} [c] at 10 6
\put {$ \scriptstyle \bullet$} [c] at 12 12
\put {$ \scriptstyle \bullet$} [c] at 12 6
\put {$ \scriptstyle \bullet$} [c] at 14 6
\put {$ \scriptstyle \bullet$} [c] at 14 0
\put {$ \scriptstyle \bullet$} [c] at 16 12
\setlinear \plot 10 0 16 12 14  0 14 6 12 12 10 6 10 0 12 6 12 12  /
\setlinear \plot 14 0  12 6  /
\put{$2{,}520$} [c] at 13 -2
\endpicture
\end{minipage}
$$
$$
\begin{minipage}{4cm}
\beginpicture
\setcoordinatesystem units    <1.5mm,2mm>
\setplotarea x from 0 to 16, y from -2 to 15
\put{1.159)} [l] at 2 12
\put {$ \scriptstyle \bullet$} [c] at 10 6
\put {$ \scriptstyle \bullet$} [c] at 11 0
\put {$ \scriptstyle \bullet$} [c] at 11 12
\put {$ \scriptstyle \bullet$} [c] at 12 6
\put {$ \scriptstyle \bullet$} [c] at 16 6
\put {$ \scriptstyle \bullet$} [c] at 16 12
\put {$ \scriptstyle \bullet$} [c] at 16 0
\setlinear \plot 16 0 16 12  11 0 12 6 11 12 10 6 11 0  /
\put{$2{,}520$} [c] at 13 -2
\endpicture
\end{minipage}
\begin{minipage}{4cm}
\beginpicture
\setcoordinatesystem units    <1.5mm,2mm>
\setplotarea x from 0 to 16, y from -2 to 15
\put{1.160)} [l] at 2 12
\put {$ \scriptstyle \bullet$} [c] at 10 6
\put {$ \scriptstyle \bullet$} [c] at 11 0
\put {$ \scriptstyle \bullet$} [c] at 11 12
\put {$ \scriptstyle \bullet$} [c] at 12 6
\put {$ \scriptstyle \bullet$} [c] at 16 6
\put {$ \scriptstyle \bullet$} [c] at 16 12
\put {$ \scriptstyle \bullet$} [c] at 16 0
\setlinear \plot 16 12 16 0  11 12 12 6 11 0 10 6 11 12  /
\put{$2{,}520$} [c] at 13 -2
\endpicture
\end{minipage}
\begin{minipage}{4cm}
\beginpicture
\setcoordinatesystem units   <1.5mm,2mm>
\setplotarea x from 0 to 16, y from -2 to 15
\put{1.161)} [l] at 2 12
\put {$ \scriptstyle \bullet$} [c] at  10 0
\put {$ \scriptstyle \bullet$} [c] at  13 3
\put {$ \scriptstyle \bullet$} [c] at  13 12
\put {$ \scriptstyle \bullet$} [c] at  14.5 12
\put {$ \scriptstyle \bullet$} [c] at  14.5 0
\put {$ \scriptstyle \bullet$} [c] at  16 3
\put {$ \scriptstyle \bullet$} [c] at  16 12
\setlinear \plot 10 0 14.5 12 16 3 16 12 10 0     /
\setlinear \plot 13 12 13 3 14.5 0 16 3    /
\put{$2{,}520$} [c] at 13 -2
\endpicture
\end{minipage}
\begin{minipage}{4cm}
\beginpicture
\setcoordinatesystem units   <1.5mm,2mm>
\setplotarea x from 0 to 16, y from -2 to 15
\put{1.162)} [l] at 2 12
\put {$ \scriptstyle \bullet$} [c] at  10 12
\put {$ \scriptstyle \bullet$} [c] at  13 0
\put {$ \scriptstyle \bullet$} [c] at  13 9
\put {$ \scriptstyle \bullet$} [c] at  14.5 12
\put {$ \scriptstyle \bullet$} [c] at  14.5 0
\put {$ \scriptstyle \bullet$} [c] at  16 9
\put {$ \scriptstyle \bullet$} [c] at  16 0
\setlinear \plot 10 12 14.5 0 16 9 16 0 10 12     /
\setlinear \plot 13 0 13 9 14.5 12 16 9    /
\put{$2{,}520$} [c] at 13 -2
\endpicture
\end{minipage}
\begin{minipage}{4cm}
\beginpicture
\setcoordinatesystem units   <1.5mm,2mm>
\setplotarea x from 0 to 16, y from -2 to 15
\put{1.163)} [l] at 2 12
\put {$ \scriptstyle \bullet$} [c] at  10 0
\put {$ \scriptstyle \bullet$} [c] at  10 4
\put {$ \scriptstyle \bullet$} [c] at  10 8
\put {$ \scriptstyle \bullet$} [c] at  10 12
\put {$ \scriptstyle \bullet$} [c] at  13 12
\put {$ \scriptstyle \bullet$} [c] at  16 0
\put {$ \scriptstyle \bullet$} [c] at  16 12
\setlinear \plot  10 0 10 12  /
\setlinear \plot  16 12 16 0 10 8 /
\setlinear \plot  13 12 16 0 /
\put{$2{,}520$} [c] at 13 -2
\endpicture
\end{minipage}
\begin{minipage}{4cm}
\beginpicture
\setcoordinatesystem units   <1.5mm,2mm>
\setplotarea x from 0 to 16, y from -2 to 15
\put{1.164)} [l] at 2 12
\put {$ \scriptstyle \bullet$} [c] at  10 0
\put {$ \scriptstyle \bullet$} [c] at  10 4
\put {$ \scriptstyle \bullet$} [c] at  10 8
\put {$ \scriptstyle \bullet$} [c] at  10 12
\put {$ \scriptstyle \bullet$} [c] at  13 0
\put {$ \scriptstyle \bullet$} [c] at  16 0
\put {$ \scriptstyle \bullet$} [c] at  16 12
\setlinear \plot  10 0 10 12  /
\setlinear \plot  16 0 16 12 10 4 /
\setlinear \plot  13 0 16 12 /
\put{$2{,}520$} [c] at 13 -2
\endpicture
\end{minipage}
$$

$$
\begin{minipage}{4cm}
\beginpicture
\setcoordinatesystem units   <1.5mm,2mm>
\setplotarea x from 0 to 16, y from -2 to 15
\put{1.165)} [l] at 2 12
\put {$ \scriptstyle \bullet$} [c] at  10 0
\put {$ \scriptstyle \bullet$} [c] at  10 6
\put {$ \scriptstyle \bullet$} [c] at  10 12
\put {$ \scriptstyle \bullet$} [c] at  13 12
\put {$ \scriptstyle \bullet$} [c] at  16 0
\put {$ \scriptstyle \bullet$} [c] at  16 6
\put {$ \scriptstyle \bullet$} [c] at  16 12
\setlinear \plot  10 0 10 12 16 0 16 12 /
\setlinear \plot  16 6 13 12 /
\put{$2{,}520$} [c] at 13 -2
\endpicture
\end{minipage}
\begin{minipage}{4cm}
\beginpicture
\setcoordinatesystem units   <1.5mm,2mm>
\setplotarea x from 0 to 16, y from -2 to 15
\put{1.166)} [l] at 2 12
\put {$ \scriptstyle \bullet$} [c] at  10 0
\put {$ \scriptstyle \bullet$} [c] at  10 6
\put {$ \scriptstyle \bullet$} [c] at  10 12
\put {$ \scriptstyle \bullet$} [c] at  13 0
\put {$ \scriptstyle \bullet$} [c] at  16 0
\put {$ \scriptstyle \bullet$} [c] at  16 6
\put {$ \scriptstyle \bullet$} [c] at  16 12
\setlinear \plot  10 12 10 0 16 12 16 0 /
\setlinear \plot  16 6 13 0 /
\put{$2{,}520$} [c] at 13 -2
\endpicture
\end{minipage}
\begin{minipage}{4cm}
\beginpicture
\setcoordinatesystem units   <1.5mm,2mm>
\setplotarea x from 0 to 16, y from -2 to 15
\put{1.167)} [l] at 2 12
\put {$ \scriptstyle \bullet$} [c] at  10 0
\put {$ \scriptstyle \bullet$} [c] at  10 6
\put {$ \scriptstyle \bullet$} [c] at  10 12
\put {$ \scriptstyle \bullet$} [c] at  13 12
\put {$ \scriptstyle \bullet$} [c] at  16 0
\put {$ \scriptstyle \bullet$} [c] at  16 6
\put {$ \scriptstyle \bullet$} [c] at  16 12
\setlinear \plot 16 0 13 12 10 0  10 12  16 0 16 12  10 0 /
\put{$2{,}520$} [c] at 13 -2
\endpicture
\end{minipage}
\begin{minipage}{4cm}
\beginpicture
\setcoordinatesystem units   <1.5mm,2mm>
\setplotarea x from 0 to 16, y from -2 to 15
\put{1.168)} [l] at 2 12
\put {$ \scriptstyle \bullet$} [c] at  10 0
\put {$ \scriptstyle \bullet$} [c] at  10 6
\put {$ \scriptstyle \bullet$} [c] at  10 12
\put {$ \scriptstyle \bullet$} [c] at  13 0
\put {$ \scriptstyle \bullet$} [c] at  16 0
\put {$ \scriptstyle \bullet$} [c] at  16 6
\put {$ \scriptstyle \bullet$} [c] at  16 12
\setlinear \plot 16 12 13 0 10 12  10 0  16 12 16 0  10 12 /
\put{$2{,}520$} [c] at 13 -2
\endpicture
\end{minipage}
\begin{minipage}{4cm}
\beginpicture
\setcoordinatesystem units   <1.5mm,2mm>
\setplotarea x from 0  to 16, y from -2 to 15
\put {1.169)} [l] at 2 12
\put {$ \scriptstyle \bullet$} [c] at  10 0
\put {$ \scriptstyle \bullet$} [c] at  10 12
\put {$ \scriptstyle \bullet$} [c] at  12 12
\put {$ \scriptstyle \bullet$} [c] at  14 0
\put {$ \scriptstyle \bullet$} [c] at  14 6
\put {$ \scriptstyle \bullet$} [c] at  14  12
\put {$ \scriptstyle \bullet$} [c] at  16  12
\setlinear \plot  10 12 10 0 14 12 14 0   10 12 /
\setlinear \plot  10 0 12 12 14 6 /
\setlinear \plot 16 12 14 6 /
\put{$2{,}520$}[c] at 13 -2
\endpicture
\end{minipage}
\begin{minipage}{4cm}
\beginpicture
\setcoordinatesystem units   <1.5mm,2mm>
\setplotarea x from 0  to 16, y from -2 to 15
\put {1.170)} [l] at 2 12
\put {$ \scriptstyle \bullet$} [c] at  10 0
\put {$ \scriptstyle \bullet$} [c] at  10 12
\put {$ \scriptstyle \bullet$} [c] at  12 0
\put {$ \scriptstyle \bullet$} [c] at  14 0
\put {$ \scriptstyle \bullet$} [c] at  14 6
\put {$ \scriptstyle \bullet$} [c] at  14  12
\put {$ \scriptstyle \bullet$} [c] at  16  0
\setlinear \plot  10 0 10 12 14 0 14 12   10 0 /
\setlinear \plot  10 12 12 0 14 6 /
\setlinear \plot 16 0 14 6 /
\put{$2{,}520$}[c] at 13 -2
\endpicture
\end{minipage}
$$
$$
\begin{minipage}{4cm}
\beginpicture
\setcoordinatesystem units   <1.5mm,2mm>
\setplotarea x from 0 to 16, y from -2 to 15
\put{1.171)} [l] at 2 12
\put {$ \scriptstyle \bullet$} [c] at  10 0
\put {$ \scriptstyle \bullet$} [c] at  10 12
\put {$ \scriptstyle \bullet$} [c] at  13 0
\put {$ \scriptstyle \bullet$} [c] at  13 12
\put {$ \scriptstyle \bullet$} [c] at  16  0
\put {$ \scriptstyle \bullet$} [c] at  16  6
\put {$ \scriptstyle \bullet$} [c] at  16  12
\setlinear \plot  10 0  10 12 13 0  16 6  16 12 /
\setlinear \plot  10  0  13 12  16 0 16 6 /
\put{$2{,}520$} [c] at 13 -2
\endpicture
\end{minipage}
\begin{minipage}{4cm}
\beginpicture
\setcoordinatesystem units   <1.5mm,2mm>
\setplotarea x from 0 to 16, y from -2 to 15
\put{1.172)} [l] at 2 12
\put {$ \scriptstyle \bullet$} [c] at  10 0
\put {$ \scriptstyle \bullet$} [c] at  10 12
\put {$ \scriptstyle \bullet$} [c] at  13 0
\put {$ \scriptstyle \bullet$} [c] at  13 12
\put {$ \scriptstyle \bullet$} [c] at  16  0
\put {$ \scriptstyle \bullet$} [c] at  16  6
\put {$ \scriptstyle \bullet$} [c] at  16  12
\setlinear \plot  10 12  10 0 13 12  16 6  16 0 /
\setlinear \plot  10  12  13 0  16 12 16 6 /
\put{$2{,}520$} [c] at 13 -2
\endpicture
\end{minipage}
\begin{minipage}{4cm}
\beginpicture
\setcoordinatesystem units   <1.5mm,2mm>
\setplotarea x from 0 to 16, y from -2 to 15
\put{1.173)} [l] at 2 12
\put {$ \scriptstyle \bullet$} [c] at  10 12
\put {$ \scriptstyle \bullet$} [c] at  10 0
\put {$ \scriptstyle \bullet$} [c] at  13 12
\put {$ \scriptstyle \bullet$} [c] at  13 0
\put {$ \scriptstyle \bullet$} [c] at  13 6
\put {$ \scriptstyle \bullet$} [c] at  16 0
\put {$ \scriptstyle \bullet$} [c] at  16 12
\setlinear \plot  10 0 10 12 13 6  13 12  10 0  /
\setlinear \plot  16 0 16 12 13 6 13 0   /
\put{$2{,}520$} [c] at 13 -2
\endpicture
\end{minipage}
\begin{minipage}{4cm}
\beginpicture
\setcoordinatesystem units   <1.5mm,2mm>
\setplotarea x from 0 to 16, y from -2 to 15
\put{1.174)} [l] at 2 12
\put {$ \scriptstyle \bullet$} [c] at  10 12
\put {$ \scriptstyle \bullet$} [c] at  10 0
\put {$ \scriptstyle \bullet$} [c] at  13 12
\put {$ \scriptstyle \bullet$} [c] at  13 0
\put {$ \scriptstyle \bullet$} [c] at  13 6
\put {$ \scriptstyle \bullet$} [c] at  16 0
\put {$ \scriptstyle \bullet$} [c] at  16 12
\setlinear \plot  10 12 10 0 13 6  13 0  10 12  /
\setlinear \plot  16 12 16 0 13 6 13 12   /
\put{$2{,}520$} [c] at 13 -2
\endpicture
\end{minipage}
\begin{minipage}{4cm}
\beginpicture
\setcoordinatesystem units   <1.5mm,2mm>
\setplotarea x from 0 to 16, y from -2 to 15
\put{1.175)} [l] at 2 12
\put {$ \scriptstyle \bullet$} [c] at  10 12
\put {$ \scriptstyle \bullet$} [c] at  10 0
\put {$ \scriptstyle \bullet$} [c] at  13 12
\put {$ \scriptstyle \bullet$} [c] at  13 0
\put {$ \scriptstyle \bullet$} [c] at  16 12
\put {$ \scriptstyle \bullet$} [c] at  16 6
\put {$ \scriptstyle \bullet$} [c] at  16 0
\setlinear \plot 10 0 10 12 13  0 13 12 10 0 16 6 16 12  13 0 /
\setlinear \plot 16 0 16 6   /
\setlinear \plot 10 12 13 0   /
\put{$2{,}520$} [c] at 13 -2
\endpicture
\end{minipage}
\begin{minipage}{4cm}
\beginpicture
\setcoordinatesystem units   <1.5mm,2mm>
\setplotarea x from 0 to 16, y from -2 to 15
\put{1.176)} [l] at 2 12
\put {$ \scriptstyle \bullet$} [c] at  10 12
\put {$ \scriptstyle \bullet$} [c] at  10 0
\put {$ \scriptstyle \bullet$} [c] at  13 12
\put {$ \scriptstyle \bullet$} [c] at  13 0
\put {$ \scriptstyle \bullet$} [c] at  16 12
\put {$ \scriptstyle \bullet$} [c] at  16 6
\put {$ \scriptstyle \bullet$} [c] at  16 0
\setlinear \plot 10 12 10 0 13  12 13 0 10 12 16 6 16 0  13 12 /
\setlinear \plot 16 12 16 6   /
\setlinear \plot 10 0 13 12   /
\put{$2{,}520 $} [c] at 13 -2
\endpicture
\end{minipage}
$$

$$
\begin{minipage}{4cm}
\beginpicture
\setcoordinatesystem units   <1.5mm,2mm>
\setplotarea x from 0 to 16, y from -2 to 15
\put{1.177)} [l] at 2 12
\put {$ \scriptstyle \bullet$} [c] at  10 0
\put {$ \scriptstyle \bullet$} [c] at  10  12
\put {$ \scriptstyle \bullet$} [c] at  13 0
\put {$ \scriptstyle \bullet$} [c] at  13 12
\put {$ \scriptstyle \bullet$} [c] at  16 0
\put {$ \scriptstyle \bullet$} [c] at  16  6
\put {$ \scriptstyle \bullet$} [c] at  16  12
\setlinear \plot 10 0 10 12  16 0  16  12 13 0 13 12 16 6 /
\setlinear \plot 10 12 13 0  /
\put{$2{,}520$} [c] at 13 -2
\endpicture
\end{minipage}
\begin{minipage}{4cm}
\beginpicture
\setcoordinatesystem units   <1.5mm,2mm>
\setplotarea x from 0 to 16, y from -2 to 15
\put{1.178)} [l] at 2 12
\put {$ \scriptstyle \bullet$} [c] at  10 0
\put {$ \scriptstyle \bullet$} [c] at  10  12
\put {$ \scriptstyle \bullet$} [c] at  13 0
\put {$ \scriptstyle \bullet$} [c] at  13 12
\put {$ \scriptstyle \bullet$} [c] at  16 0
\put {$ \scriptstyle \bullet$} [c] at  16  6
\put {$ \scriptstyle \bullet$} [c] at  16  12
\setlinear \plot 10 12 10 0  16 12  16  0 13 12 13 0 16 6 /
\setlinear \plot 10 0 13 12  /
\put{$2{,}520$} [c] at 13 -2
\endpicture
\end{minipage}
\begin{minipage}{4cm}
\beginpicture
\setcoordinatesystem units   <1.5mm,2mm>
\setplotarea x from 0 to 16, y from -2 to 15
\put{1.179)} [l] at 2 12
\put {$ \scriptstyle \bullet$} [c] at  10  0
\put {$ \scriptstyle \bullet$} [c] at  10 12
\put {$ \scriptstyle \bullet$} [c] at  13 12
\put {$ \scriptstyle \bullet$} [c] at  13 0
\put {$ \scriptstyle \bullet$} [c] at  16 0
\put {$ \scriptstyle \bullet$} [c] at  16 6
\put {$ \scriptstyle \bullet$} [c] at  16 12
\setlinear \plot 13 0 10 12 10  0 13 12  13 0  16 12 16 0 13 12  /
\setlinear \plot  10 0 16 12  /
\put{$2{,}520 $} [c] at 13 -2
\endpicture
\end{minipage}
\begin{minipage}{4cm}
\beginpicture
\setcoordinatesystem units   <1.5mm,2mm>
\setplotarea x from 0 to 16, y from -2 to 15
\put{1.180)} [l] at 2 12
\put {$ \scriptstyle \bullet$} [c] at  10  0
\put {$ \scriptstyle \bullet$} [c] at  10 12
\put {$ \scriptstyle \bullet$} [c] at  13 12
\put {$ \scriptstyle \bullet$} [c] at  13 0
\put {$ \scriptstyle \bullet$} [c] at  16 0
\put {$ \scriptstyle \bullet$} [c] at  16 6
\put {$ \scriptstyle \bullet$} [c] at  16 12
\setlinear \plot 13 12 10 0 10  12 13 0  13 12  16 0 16 12 13 0  /
\setlinear \plot  10 12 16 0  /
\put{$2{,}520 $} [c] at 13 -2
\endpicture
\end{minipage}
\begin{minipage}{4cm}
\beginpicture
\setcoordinatesystem units    <1.5mm,2mm>
\setplotarea x from 0 to 16, y from -2 to 15
\put{1.181)} [l] at 2 12
\put {$ \scriptstyle \bullet$} [c] at 13 0
\put {$ \scriptstyle \bullet$} [c] at 10 4
\put {$ \scriptstyle \bullet$} [c] at 12 4
\put {$ \scriptstyle \bullet$} [c] at 14 4
\put {$ \scriptstyle \bullet$} [c] at 16 4
\put {$ \scriptstyle \bullet$} [c] at 13 12
\put {$ \scriptstyle \bullet$} [c] at 16 12
\setlinear \plot  13 12 10 4  13 0 16 4 16 12 14 4 13 0 12 4  13 12 16 4 /
\setlinear \plot  13 12 14 4  /
\put{$1{,}260$} [c] at 13 -2
\endpicture
\end{minipage}
\begin{minipage}{4cm}
\beginpicture
\setcoordinatesystem units    <1.5mm,2mm>
\setplotarea x from 0 to 16, y from -2 to 15
\put{1.182)} [l] at 2 12
\put {$ \scriptstyle \bullet$} [c] at 13 0
\put {$ \scriptstyle \bullet$} [c] at 10 8
\put {$ \scriptstyle \bullet$} [c] at 12 8
\put {$ \scriptstyle \bullet$} [c] at 14 8
\put {$ \scriptstyle \bullet$} [c] at 16 8
\put {$ \scriptstyle \bullet$} [c] at 13 12
\put {$ \scriptstyle \bullet$} [c] at 16 0
\setlinear \plot  13 0 10 8  13 12 16 8 16 0 14 8 13 12 12 8  13 0 16 8 /
\setlinear \plot  13 0 14 8  /
\put{$1{,}260$} [c] at 13 -2
\endpicture
\end{minipage}
$$
$$
\begin{minipage}{4cm}
\beginpicture
\setcoordinatesystem units    <1.5mm,2mm>
\setplotarea x from 0 to 16, y from -2 to 15
\put{1.183)} [l] at 2 12
\put {$ \scriptstyle \bullet$} [c] at 10 12
\put {$ \scriptstyle \bullet$} [c] at 12 12
\put {$ \scriptstyle \bullet$} [c] at 14 12
\put {$ \scriptstyle \bullet$} [c] at 16 12
\put {$ \scriptstyle \bullet$} [c] at 14 6
\put {$ \scriptstyle \bullet$} [c] at 15 0
\put {$ \scriptstyle \bullet$} [c] at 16 6
\setlinear \plot  14 6 15 0 16 6 16 12 14 6 14 12 16 6 /
\setlinear \plot  10 12 14 6 12  12  /
\put{$1{,}260$} [c] at 13 -2
\endpicture
\end{minipage}
\begin{minipage}{4cm}
\beginpicture
\setcoordinatesystem units    <1.5mm,2mm>
\setplotarea x from 0 to 16, y from -2 to 15
\put{1.184)} [l] at 2 12
\put {$ \scriptstyle \bullet$} [c] at 10 0
\put {$ \scriptstyle \bullet$} [c] at 12 0
\put {$ \scriptstyle \bullet$} [c] at 14 0
\put {$ \scriptstyle \bullet$} [c] at 16 0
\put {$ \scriptstyle \bullet$} [c] at 14 6
\put {$ \scriptstyle \bullet$} [c] at 15 12
\put {$ \scriptstyle \bullet$} [c] at 16 6
\setlinear \plot  14 6 15 12 16 6 16 0 14 6 14 0 16 6 /
\setlinear \plot  10 0 14 6 12  0  /
\put{$1{,}260$} [c] at 13 -2
\endpicture
\end{minipage}
\begin{minipage}{4cm}
\beginpicture
\setcoordinatesystem units   <1.5mm,2mm>
\setplotarea x from 0 to 16, y from -2 to 15
\put{1.185)} [l] at 2 12
\put {$ \scriptstyle \bullet$} [c] at  10 6
\put {$ \scriptstyle \bullet$} [c] at  10 12
\put {$ \scriptstyle \bullet$} [c] at  11 0
\put {$ \scriptstyle \bullet$} [c] at  12 6
\put {$ \scriptstyle \bullet$} [c] at  12 12
\put {$ \scriptstyle \bullet$} [c] at  16 0
\put {$ \scriptstyle \bullet$} [c] at  16 12
\setlinear \plot 10 6 12 12 12 6 10 12 10 6 11 0 16 12 16 0 /
\setlinear \plot 11 0 12 6 /
\put{$1{,}260$} [c] at 13 -2
\endpicture
\end{minipage}
\begin{minipage}{4cm}
\beginpicture
\setcoordinatesystem units   <1.5mm,2mm>
\setplotarea x from 0 to 16, y from -2 to 15
\put{1.186)} [l] at 2 12
\put {$ \scriptstyle \bullet$} [c] at  10 6
\put {$ \scriptstyle \bullet$} [c] at  10 0
\put {$ \scriptstyle \bullet$} [c] at  11 12
\put {$ \scriptstyle \bullet$} [c] at  12 6
\put {$ \scriptstyle \bullet$} [c] at  12 0
\put {$ \scriptstyle \bullet$} [c] at  16 0
\put {$ \scriptstyle \bullet$} [c] at  16 12
\setlinear \plot 10 6 12 0 12 6 10 0 10 6 11 12 16 0 16 12 /
\setlinear \plot 11 12 12 6 /
\put{$1{,}260$} [c] at 13 -2
\endpicture
\end{minipage}
\begin{minipage}{4cm}
\beginpicture
\setcoordinatesystem units   <1.5mm,2mm>
\setplotarea x from 0  to 16, y from -2 to 15
\put {1.187)} [l] at 2 12
\put {$ \scriptstyle \bullet$} [c] at  10 12
\put {$ \scriptstyle \bullet$} [c] at  12 0
\put {$ \scriptstyle \bullet$} [c] at  12 12
\put {$ \scriptstyle \bullet$} [c] at  14 12
\put {$ \scriptstyle \bullet$} [c] at  16 0
\put {$ \scriptstyle \bullet$} [c] at  16  6
\put {$ \scriptstyle \bullet$} [c] at  16  12
\setlinear \plot  10 12 12 0  16 12 16 0 /
\setlinear \plot  12 0  12 12 16 6 14 12 12 0 /
\put{$840$}[c] at 13 -2
\endpicture
\end{minipage}
\begin{minipage}{4cm}
\beginpicture
\setcoordinatesystem units   <1.5mm,2mm>
\setplotarea x from 0  to 16, y from -2 to 15
\put {1.188)} [l] at 2 12
\put {$ \scriptstyle \bullet$} [c] at  10 0
\put {$ \scriptstyle \bullet$} [c] at  12 0
\put {$ \scriptstyle \bullet$} [c] at  12 12
\put {$ \scriptstyle \bullet$} [c] at  14 0
\put {$ \scriptstyle \bullet$} [c] at  16 0
\put {$ \scriptstyle \bullet$} [c] at  16  6
\put {$ \scriptstyle \bullet$} [c] at  16  12
\setlinear \plot  10 0 12 12  16 0 16 12 /
\setlinear \plot  12 12  12 0 16 6 14 0 12 12 /
\put{$840$}[c] at 13 -2
\endpicture
\end{minipage}
$$

$$
\begin{minipage}{4cm}
\beginpicture
\setcoordinatesystem units   <1.5mm,2mm>
\setplotarea x from 0 to 16, y from -2 to 15
\put{1.189)} [l] at 2 12
\put {$ \scriptstyle \bullet$} [c] at  10 12
\put {$ \scriptstyle \bullet$} [c] at  12 9
\put {$ \scriptstyle \bullet$} [c] at  12  0
\put {$ \scriptstyle \bullet$} [c] at  13 12
\put {$ \scriptstyle \bullet$} [c] at  14 9
\put {$ \scriptstyle \bullet$} [c] at  14 0
\put {$ \scriptstyle \bullet$} [c] at  16 12
\setlinear \plot 12 9 14 0 10 12 12 0 12 9 13 12 14 9 14 0 16 12 12 0 14 9 /
\put{$630$} [c] at 13 -2
\endpicture
\end{minipage}
\begin{minipage}{4cm}
\beginpicture
\setcoordinatesystem units   <1.5mm,2mm>
\setplotarea x from 0 to 16, y from -2 to 15
\put{1.190)} [l] at 2 12
\put {$ \scriptstyle \bullet$} [c] at  10 0
\put {$ \scriptstyle \bullet$} [c] at  12 3
\put {$ \scriptstyle \bullet$} [c] at  12  12
\put {$ \scriptstyle \bullet$} [c] at  13 0
\put {$ \scriptstyle \bullet$} [c] at  14 3
\put {$ \scriptstyle \bullet$} [c] at  14 12
\put {$ \scriptstyle \bullet$} [c] at  16 0
\setlinear \plot 12 3 14 12 10 0 12 12 12 3 13 0 14 3 14 12 16 0 12 12 14 3 /
\put{$630$} [c] at 13 -2
\endpicture
\end{minipage}
\begin{minipage}{4cm}
\beginpicture
\setcoordinatesystem units   <1.5mm,2mm>
\setplotarea x from 0  to 16, y from -2 to 15
\put {1.191)} [l] at 2 12
\put {$ \scriptstyle \bullet$} [c] at  10 0
\put {$ \scriptstyle \bullet$} [c] at  12 0
\put {$ \scriptstyle \bullet$} [c] at  14 12
\put {$ \scriptstyle \bullet$} [c] at  16 12
\put {$ \scriptstyle \bullet$} [c] at  10 12
\put {$ \scriptstyle \bullet$} [c] at  12  12
\put {$ \scriptstyle \bullet$} [c] at  15 6
\setlinear \plot  10 0 10 12 12 0 12  12 10 0  /
\setlinear \plot  16 12 15 6  12  0   /
\setlinear \plot  14 12  15 6 10 0  /
\put{$630$} [c] at 13 -2
\endpicture
\end{minipage}
\begin{minipage}{4cm}
\beginpicture
\setcoordinatesystem units   <1.5mm,2mm>
\setplotarea x from 0  to 16, y from -2 to 15
\put {1.192)} [l] at 2 12
\put {$ \scriptstyle \bullet$} [c] at  10 0
\put {$ \scriptstyle \bullet$} [c] at  12 0
\put {$ \scriptstyle \bullet$} [c] at  14 0
\put {$ \scriptstyle \bullet$} [c] at  16 0
\put {$ \scriptstyle \bullet$} [c] at  10 12
\put {$ \scriptstyle \bullet$} [c] at  12  12
\put {$ \scriptstyle \bullet$} [c] at  15 6
\setlinear \plot  10 0 10 12 12 0 12  12 10 0  /
\setlinear \plot  16 0 15 6  12  12   /
\setlinear \plot  14 0  15 6 10 12  /
\put{$630$} [c] at 13 -2
\endpicture
\end{minipage}
\begin{minipage}{4cm}
\beginpicture
\setcoordinatesystem units    <1.5mm,2mm>
\setplotarea x from 0 to 16, y from -2 to 15
\put{1.193)} [l] at 2 12
\put {$ \scriptstyle \bullet$} [c] at 10 4
\put {$ \scriptstyle \bullet$} [c] at 11.5 4
\put {$ \scriptstyle \bullet$} [c] at 13 4
\put {$ \scriptstyle \bullet$} [c] at 10 12
\put {$ \scriptstyle \bullet$} [c] at 13 12
\put {$ \scriptstyle \bullet$} [c] at 13 0
\put {$ \scriptstyle \bullet$} [c] at 16 12
\setlinear \plot 10 12 10 4 13 0 13 12 11.5 4 10 12 /
\setlinear \plot 16 12 13 0  11.5 4   /
\setlinear \plot 10 12 13 4     /
\setlinear \plot 13 12 10 4     /
\put{$420$} [c] at 13 -2
\endpicture
\end{minipage}
\begin{minipage}{4cm}
\beginpicture
\setcoordinatesystem units    <1.5mm,2mm>
\setplotarea x from 0 to 16, y from -2 to 15
\put{1.194)} [l] at 2 12
\put {$ \scriptstyle \bullet$} [c] at 10 8
\put {$ \scriptstyle \bullet$} [c] at 11.5 8
\put {$ \scriptstyle \bullet$} [c] at 13 8
\put {$ \scriptstyle \bullet$} [c] at 10 0
\put {$ \scriptstyle \bullet$} [c] at 13 0
\put {$ \scriptstyle \bullet$} [c] at 13 12
\put {$ \scriptstyle \bullet$} [c] at 16 0
\setlinear \plot 10 0 10 8 13 12 13 0 11.5 8 10 0 /
\setlinear \plot 16 0 13 12  11.5 8   /
\setlinear \plot 10 0 13 8     /
\setlinear \plot 13 0 10 8     /
\put{$420$} [c] at 13 -2
\endpicture
\end{minipage}
$$
$$
\begin{minipage}{4cm}
\beginpicture
\setcoordinatesystem units    <1.5mm,2mm>
\setplotarea x from 0 to 16, y from -2 to 15
\put{1.195)} [l] at 2 12
\put {$ \scriptstyle \bullet$} [c] at 10 12
\put {$ \scriptstyle \bullet$} [c] at 11.5 12
\put {$ \scriptstyle \bullet$} [c] at 13 12
\put {$ \scriptstyle \bullet$} [c] at 10 6
\put {$ \scriptstyle \bullet$} [c] at 13 6
\put {$ \scriptstyle \bullet$} [c] at 11.5 0
\put {$ \scriptstyle \bullet$} [c] at 16 12
\setlinear \plot 16 12 11.5 0 10 6 10 12 13 6 13 12 10 6  /
\setlinear \plot 10 6 11.5 12 13 6 11.5 0    /
\put{$420$} [c] at 13 -2
\endpicture
\end{minipage}
\begin{minipage}{4cm}
\beginpicture
\setcoordinatesystem units    <1.5mm,2mm>
\setplotarea x from 0 to 16, y from -2 to 15
\put{1.196)} [l] at 2 12
\put {$ \scriptstyle \bullet$} [c] at 10 0
\put {$ \scriptstyle \bullet$} [c] at 11.5 0
\put {$ \scriptstyle \bullet$} [c] at 13 0
\put {$ \scriptstyle \bullet$} [c] at 10 6
\put {$ \scriptstyle \bullet$} [c] at 13 6
\put {$ \scriptstyle \bullet$} [c] at 11.5 12
\put {$ \scriptstyle \bullet$} [c] at 16 0
\setlinear \plot 16 0 11.5 12 10 6 10 0 13 6 13 0 10 6  /
\setlinear \plot 10 6 11.5 0 13 6 11.5 12    /
\put{$420$} [c] at 13 -2
\endpicture
\end{minipage}
\begin{minipage}{4cm}
\beginpicture
\setcoordinatesystem units   <1.5mm,2mm>
\setplotarea x from 0  to 16, y from -2 to 15
\put {1.197)} [l] at  2 12
\put {$ \scriptstyle \bullet$} [c] at  10 12
\put {$ \scriptstyle \bullet$} [c] at  12 12
\put {$ \scriptstyle \bullet$} [c] at  14 12
\put {$ \scriptstyle \bullet$} [c] at  16  12
\put {$ \scriptstyle \bullet$} [c] at  10 0
\put {$ \scriptstyle \bullet$} [c] at  13 0
\put {$ \scriptstyle \bullet$} [c] at  16 0
\setlinear \plot 16 12  10 0 10 12 13  0 16 12 16 0 14 12 13 0 12 12 10 0 14 12  /
\setlinear \plot  10 12  16 0 12 12 /
\put{$35$} [c] at 13 -2
\endpicture
\end{minipage}
\begin{minipage}{4cm}
\beginpicture
\setcoordinatesystem units   <1.5mm,2mm>
\setplotarea x from 0  to 16, y from -2 to 15
\put {1.198)} [l] at  2 12
\put {$ \scriptstyle \bullet$} [c] at  10 0
\put {$ \scriptstyle \bullet$} [c] at  12 0
\put {$ \scriptstyle \bullet$} [c] at  14 0
\put {$ \scriptstyle \bullet$} [c] at  16  0
\put {$ \scriptstyle \bullet$} [c] at  10 12
\put {$ \scriptstyle \bullet$} [c] at  13 12
\put {$ \scriptstyle \bullet$} [c] at  16 12
\setlinear \plot 16 0  10 12 10 0 13  12 16 0 16 12 14 0 13 12 12 0 10 12 14 0  /
\setlinear \plot  10 0  16 12 12 0 /
\put{$35$} [c] at 13 -2
\endpicture
\end{minipage}
\begin{minipage}{4cm}
\beginpicture
\setcoordinatesystem units   <1.5mm,2mm>
\setplotarea x from 0 to 16, y from -2 to 15
\put{${\bf  24}$} [l] at 2 15
\put{1.199)} [l] at 2 12
\put {$ \scriptstyle \bullet$} [c] at  10 12
\put {$ \scriptstyle \bullet$} [c] at  10 6
\put {$ \scriptstyle \bullet$} [c] at  11 12
\put {$ \scriptstyle \bullet$} [c] at  11 0
\put {$ \scriptstyle \bullet$} [c] at  12 6
\put {$ \scriptstyle \bullet$} [c] at  16 0
\put {$ \scriptstyle \bullet$} [c] at  10.5 9
\setlinear \plot  10 12 10 6 11 12 12 6 11 0 10 6   /
\setlinear \plot  16 0 11 12   /
\put{$5{,}040$} [c] at 13 -2
\endpicture
\end{minipage}
\begin{minipage}{4cm}
\beginpicture
\setcoordinatesystem units   <1.5mm,2mm>
\setplotarea x from 0 to 16, y from -2 to 15
\put{1.200)} [l] at 2 12
\put {$ \scriptstyle \bullet$} [c] at  10 0
\put {$ \scriptstyle \bullet$} [c] at  10 6
\put {$ \scriptstyle \bullet$} [c] at  11 12
\put {$ \scriptstyle \bullet$} [c] at  11 0
\put {$ \scriptstyle \bullet$} [c] at  12 6
\put {$ \scriptstyle \bullet$} [c] at  16 12
\put {$ \scriptstyle \bullet$} [c] at  10.5 3
\setlinear \plot  10 0 10 6 11 0 12 6 11 12 10 6   /
\setlinear \plot  16 12 11 0   /
\put{$5{,}040$} [c] at 13 -2
\endpicture
\end{minipage}
$$
$$
\begin{minipage}{4cm}
\beginpicture
\setcoordinatesystem units   <1.5mm,2mm>
\setplotarea x from 0 to 16, y from -2 to 15
\put{1.201)} [l] at 2 12
\put {$ \scriptstyle \bullet$} [c] at  10 0
\put {$ \scriptstyle \bullet$} [c] at  10 6
\put {$ \scriptstyle \bullet$} [c] at  10 12
\put {$ \scriptstyle \bullet$} [c] at  13 6
\put {$ \scriptstyle \bullet$} [c] at  16 0
\put {$ \scriptstyle \bullet$} [c] at  16 6
\put {$ \scriptstyle \bullet$} [c] at  16 12
\setlinear \plot  10 12 10 0 16 12 16 0   /
\put{$5{,}040$} [c] at 13 -2
\endpicture
\end{minipage}
\begin{minipage}{4cm}
\beginpicture
\setcoordinatesystem units   <1.5mm,2mm>
\setplotarea x from 0 to 16, y from -2 to 15
\put{1.202)} [l] at 2 12
\put {$ \scriptstyle \bullet$} [c] at  13 0
\put {$ \scriptstyle \bullet$} [c] at  13 4
\put {$ \scriptstyle \bullet$} [c] at  13 8
\put {$ \scriptstyle \bullet$} [c] at  13 12
\put {$ \scriptstyle \bullet$} [c] at  16 0
\put {$ \scriptstyle \bullet$} [c] at  16 12
\put {$ \scriptstyle \bullet$} [c] at  10 12
\setlinear \plot  16 0 16 12  13 4 13 0 10 12   /
\setlinear \plot  13 12 13 4  /
\put{$5{,}040$} [c] at 13 -2
\endpicture
\end{minipage}
\begin{minipage}{4cm}
\beginpicture
\setcoordinatesystem units   <1.5mm,2mm>
\setplotarea x from 0 to 16, y from -2 to 15
\put{1.203)} [l] at 2 12
\put {$ \scriptstyle \bullet$} [c] at  13 0
\put {$ \scriptstyle \bullet$} [c] at  13 4
\put {$ \scriptstyle \bullet$} [c] at  13 8
\put {$ \scriptstyle \bullet$} [c] at  13 12
\put {$ \scriptstyle \bullet$} [c] at  16 0
\put {$ \scriptstyle \bullet$} [c] at  16 12
\put {$ \scriptstyle \bullet$} [c] at  10 0
\setlinear \plot  16 12 16 0  13 8 13 12 10 0   /
\setlinear \plot  13 0 13 8  /
\put{$5{,}040$} [c] at 13 -2
\endpicture
\end{minipage}
\begin{minipage}{4cm}
\beginpicture
\setcoordinatesystem units   <1.5mm,2mm>
\setplotarea x from 0 to 16, y from -2 to 15
\put{1.204)} [l] at 2 12
\put {$ \scriptstyle \bullet$} [c] at  10 12
\put {$ \scriptstyle \bullet$} [c] at  13 0
\put {$ \scriptstyle \bullet$} [c] at  13 4
\put {$ \scriptstyle \bullet$} [c] at  13 8
\put {$ \scriptstyle \bullet$} [c] at  13 12
\put {$ \scriptstyle \bullet$} [c] at  16 0
\put {$ \scriptstyle \bullet$} [c] at  16 12
\setlinear \plot 10 12 13 0  13 12  16 0 16 12 13 0 /
\put{$5{,}040$} [c] at 13 -2
\endpicture
\end{minipage}
\begin{minipage}{4cm}
\beginpicture
\setcoordinatesystem units   <1.5mm,2mm>
\setplotarea x from 0 to 16, y from -2 to 15
\put{1.205)} [l] at 2 12
\put {$ \scriptstyle \bullet$} [c] at  10 0
\put {$ \scriptstyle \bullet$} [c] at  13 0
\put {$ \scriptstyle \bullet$} [c] at  13 4
\put {$ \scriptstyle \bullet$} [c] at  13 8
\put {$ \scriptstyle \bullet$} [c] at  13 12
\put {$ \scriptstyle \bullet$} [c] at  16 0
\put {$ \scriptstyle \bullet$} [c] at  16 12
\setlinear \plot 10 0 13 12  13 0  16 12 16 0 13 12 /
\put{$5{,}040$} [c] at 13 -2
\endpicture
\end{minipage}
\begin{minipage}{4cm}
\beginpicture
\setcoordinatesystem units   <1.5mm,2mm>
\setplotarea x from 0 to 16, y from -2 to 15
\put{1.206)} [l] at 2 12
\put {$ \scriptstyle \bullet$} [c] at  10 12
\put {$ \scriptstyle \bullet$} [c] at  14 12
\put {$ \scriptstyle \bullet$} [c] at  14 6
\put {$ \scriptstyle \bullet$} [c] at  15 0
\put {$ \scriptstyle \bullet$} [c] at  15 12
\put {$ \scriptstyle \bullet$} [c] at  16 6
\put {$ \scriptstyle \bullet$} [c] at  16 0
\setlinear \plot 10  12 15 0 14 6 15 12  16 6  16 0  /
\setlinear \plot  14 12 14 6   /
\setlinear \plot 16 6 15 0     /
\put{$5{,}040$} [c] at 13 -2
\endpicture
\end{minipage}
$$

$$
\begin{minipage}{4cm}
\beginpicture
\setcoordinatesystem units   <1.5mm,2mm>
\setplotarea x from 0 to 16, y from -2 to 15
\put{1.207)} [l] at 2 12
\put {$ \scriptstyle \bullet$} [c] at  10 0
\put {$ \scriptstyle \bullet$} [c] at  14 0
\put {$ \scriptstyle \bullet$} [c] at  14 6
\put {$ \scriptstyle \bullet$} [c] at  15 0
\put {$ \scriptstyle \bullet$} [c] at  15 12
\put {$ \scriptstyle \bullet$} [c] at  16 6
\put {$ \scriptstyle \bullet$} [c] at  16 12
\setlinear \plot 10  0 15 12 14 6 15 0  16 6  16 12  /
\setlinear \plot  14 0 14 6   /
\setlinear \plot 16 6 15 12     /
\put{$5{,}040$} [c] at 13 -2
\endpicture
\end{minipage}
\begin{minipage}{4cm}
\beginpicture
\setcoordinatesystem units   <1.5mm,2mm>
\setplotarea x from 0 to 16, y from -2 to 15
\put{1.208)} [l] at 2 12
\put {$ \scriptstyle \bullet$} [c] at  10 12
\put {$ \scriptstyle \bullet$} [c] at  13 0
\put {$ \scriptstyle \bullet$} [c] at  13 6
\put {$ \scriptstyle \bullet$} [c] at  13 12
\put {$ \scriptstyle \bullet$} [c] at  14.5  6
\put {$ \scriptstyle \bullet$} [c] at  16 12
\put {$ \scriptstyle \bullet$} [c] at  16 0
\setlinear \plot  13  0  13 12 16 0 16 12 /
\setlinear \plot  13 6 10 12  /
\put{$5{,}040$} [c] at 13 -2
\endpicture
\end{minipage}
\begin{minipage}{4cm}
\beginpicture
\setcoordinatesystem units   <1.5mm,2mm>
\setplotarea x from 0 to 16, y from -2 to 15
\put{1.209)} [l] at 2 12
\put {$ \scriptstyle \bullet$} [c] at  10 0
\put {$ \scriptstyle \bullet$} [c] at  13 0
\put {$ \scriptstyle \bullet$} [c] at  13 6
\put {$ \scriptstyle \bullet$} [c] at  13 12
\put {$ \scriptstyle \bullet$} [c] at  14.5  6
\put {$ \scriptstyle \bullet$} [c] at  16 12
\put {$ \scriptstyle \bullet$} [c] at  16 0
\setlinear \plot  13  12  13 0 16 12 16 0 /
\setlinear \plot  13 6 10 0  /
\put{$5{,}040$} [c] at 13 -2
\endpicture
\end{minipage}
\begin{minipage}{4cm}
\beginpicture
\setcoordinatesystem units   <1.5mm,2mm>
\setplotarea x from 0 to 16, y from -2 to 15
\put{1.210)} [l] at 2 12
\put {$ \scriptstyle \bullet$} [c] at  10 0
\put {$ \scriptstyle \bullet$} [c] at  10 12
\put {$ \scriptstyle \bullet$} [c] at  13 6
\put {$ \scriptstyle \bullet$} [c] at  13 12
\put {$ \scriptstyle \bullet$} [c] at  16 6
\put {$ \scriptstyle \bullet$} [c] at  16 12
\put {$ \scriptstyle \bullet$} [c] at  16 0
\setlinear \plot  10 12   10 0  13 6 16 0 16 12 /
\setlinear \plot  13 12 13 6 /
\put{$5{,}040$} [c] at 13 -2
\endpicture
\end{minipage}
\begin{minipage}{4cm}
\beginpicture
\setcoordinatesystem units   <1.5mm,2mm>
\setplotarea x from 0 to 16, y from -2 to 15
\put{1.211)} [l] at 2 12
\put {$ \scriptstyle \bullet$} [c] at  10 0
\put {$ \scriptstyle \bullet$} [c] at  10 12
\put {$ \scriptstyle \bullet$} [c] at  13 6
\put {$ \scriptstyle \bullet$} [c] at  13 0
\put {$ \scriptstyle \bullet$} [c] at  16 6
\put {$ \scriptstyle \bullet$} [c] at  16 12
\put {$ \scriptstyle \bullet$} [c] at  16 0
\setlinear \plot  10 0   10 12  13 6 16 12 16 0 /
\setlinear \plot  13 0 13 6 /
\put{$5{,}040$} [c] at 13 -2
\endpicture
\end{minipage}
\begin{minipage}{4cm}
\beginpicture
\setcoordinatesystem units   <1.5mm,2mm>
\setplotarea x from 0 to 16, y from -2 to 15
\put{1.212)} [l] at 2 12
\put {$ \scriptstyle \bullet$} [c] at  10 0
\put {$ \scriptstyle \bullet$} [c] at  10 6
\put {$ \scriptstyle \bullet$} [c] at  10 12
\put {$ \scriptstyle \bullet$} [c] at  13  12
\put {$ \scriptstyle \bullet$} [c] at  16 0
\put {$ \scriptstyle \bullet$} [c] at  16 6
\put {$ \scriptstyle \bullet$} [c] at  16 12
\setlinear \plot  16 12  16  0 10  12 10 0  13  12 16 0    /
\put{$5{,}040$} [c] at 13 -2
\endpicture
\end{minipage}
$$
$$
\begin{minipage}{4cm}
\beginpicture
\setcoordinatesystem units   <1.5mm,2mm>
\setplotarea x from 0 to 16, y from -2 to 15
\put{1.213)} [l] at 2 12
\put {$ \scriptstyle \bullet$} [c] at  10 0
\put {$ \scriptstyle \bullet$} [c] at  10 6
\put {$ \scriptstyle \bullet$} [c] at  10 12
\put {$ \scriptstyle \bullet$} [c] at  13  0
\put {$ \scriptstyle \bullet$} [c] at  16 0
\put {$ \scriptstyle \bullet$} [c] at  16 6
\put {$ \scriptstyle \bullet$} [c] at  16 12
\setlinear \plot  16 0  16  12 10  0 10 12  13  0 16 12    /
\put{$5{,}040$} [c] at 13 -2
\endpicture
\end{minipage}
\begin{minipage}{4cm}
\beginpicture
\setcoordinatesystem units   <1.5mm,2mm>
\setplotarea x from 0 to 16, y from -2 to 15
\put{1.214)} [l] at 2 12
\put {$ \scriptstyle \bullet$} [c] at  10 6
\put {$ \scriptstyle \bullet$} [c] at  12 0
\put {$ \scriptstyle \bullet$} [c] at  12 12
\put {$ \scriptstyle \bullet$} [c] at  14 6
\put {$ \scriptstyle \bullet$} [c] at  14 12
\put {$ \scriptstyle \bullet$} [c] at  16 0
\put {$ \scriptstyle \bullet$} [c] at  16 12
\setlinear \plot 16 0 16 12 14 6 12 0 10 6 12 12 14 6 /
\setlinear \plot  14 12 14 6  /
\put{$5{,}040$} [c] at 13 -2
\endpicture
\end{minipage}
\begin{minipage}{4cm}
\beginpicture
\setcoordinatesystem units   <1.5mm,2mm>
\setplotarea x from 0 to 16, y from -2 to 15
\put{1.215)} [l] at 2 12
\put {$ \scriptstyle \bullet$} [c] at  10 6
\put {$ \scriptstyle \bullet$} [c] at  12 0
\put {$ \scriptstyle \bullet$} [c] at  12 12
\put {$ \scriptstyle \bullet$} [c] at  14 6
\put {$ \scriptstyle \bullet$} [c] at  14 0
\put {$ \scriptstyle \bullet$} [c] at  16 0
\put {$ \scriptstyle \bullet$} [c] at  16 12
\setlinear \plot 16 12 16 0 14 6 12 12 10 6 12 0 14 6 /
\setlinear \plot  14 0 14 6  /
\put{$5{,}040$} [c] at 13 -2
\endpicture
\end{minipage}
\begin{minipage}{4cm}
\beginpicture
\setcoordinatesystem units   <1.5mm,2mm>
\setplotarea x from 0 to 16, y from -2 to 15
\put{1.216)} [l] at 2 12
\put {$ \scriptstyle \bullet$} [c] at  10 6
\put {$ \scriptstyle \bullet$} [c] at  10 12
\put {$ \scriptstyle \bullet$} [c] at  11 0
\put {$ \scriptstyle \bullet$} [c] at  11 12
\put {$ \scriptstyle \bullet$} [c] at  12 6
\put {$ \scriptstyle \bullet$} [c] at  16 0
\put {$ \scriptstyle \bullet$} [c] at  16 12
\setlinear \plot  10 12  10 6 11 0 12 6 11 12 10 6 /
\setlinear \plot  11 0  16 12 16 0 11 12  /
\put{$5{,}040$} [c] at 13 -2
\endpicture
\end{minipage}
\begin{minipage}{4cm}
\beginpicture
\setcoordinatesystem units   <1.5mm,2mm>
\setplotarea x from 0 to 16, y from -2 to 15
\put{1.217)} [l] at 2 12
\put {$ \scriptstyle \bullet$} [c] at  10 0
\put {$ \scriptstyle \bullet$} [c] at  10 6
\put {$ \scriptstyle \bullet$} [c] at  11 0
\put {$ \scriptstyle \bullet$} [c] at  11 12
\put {$ \scriptstyle \bullet$} [c] at  12 6
\put {$ \scriptstyle \bullet$} [c] at  16 0
\put {$ \scriptstyle \bullet$} [c] at  16 12
\setlinear \plot  10 0  10 6 11 12 12 6 11 0 10 6 /
\setlinear \plot  11 12  16 0 16 12 11 0  /
\put{$5{,}040$} [c] at 13 -2
\endpicture
\end{minipage}
\begin{minipage}{4cm}
\beginpicture
\setcoordinatesystem units   <1.5mm,2mm>
\setplotarea x from 0 to 16, y from -2 to 15
\put{1.218)} [l] at 2 12
\put {$ \scriptstyle \bullet$} [c] at  10 0
\put {$ \scriptstyle \bullet$} [c] at  10 12
\put {$ \scriptstyle \bullet$} [c] at  12 6
\put {$ \scriptstyle \bullet$} [c] at  12 12
\put {$ \scriptstyle \bullet$} [c] at  14 12
\put {$ \scriptstyle \bullet$} [c] at  14 0
\put {$ \scriptstyle \bullet$} [c] at  16 6
\setlinear \plot  10 12   10 0  12 12 12 6 14 12 16 6 14  0  12 6 /
\put{$5{,}040$} [c] at 13 -2
\endpicture
\end{minipage}
$$

$$
\begin{minipage}{4cm}
\beginpicture
\setcoordinatesystem units   <1.5mm,2mm>
\setplotarea x from 0 to 16, y from -2 to 15
\put{1.219)} [l] at 2 12
\put {$ \scriptstyle \bullet$} [c] at  10 0
\put {$ \scriptstyle \bullet$} [c] at  10 12
\put {$ \scriptstyle \bullet$} [c] at  12 0
\put {$ \scriptstyle \bullet$} [c] at  12 6
\put {$ \scriptstyle \bullet$} [c] at  14 12
\put {$ \scriptstyle \bullet$} [c] at  14 0
\put {$ \scriptstyle \bullet$} [c] at  16 6
\setlinear \plot  10 0   10 12  12 0 12 6 14 0 16 6 14  12  12 6 /
\put{$5{,}040$} [c] at 13 -2
\endpicture
\end{minipage}
\begin{minipage}{4cm}
\beginpicture
\setcoordinatesystem units   <1.5mm,2mm>
\setplotarea x from 0 to 16, y from -2 to 15
\put{1.220)} [l] at 2 12
\put {$ \scriptstyle \bullet$} [c] at  10 12
\put {$ \scriptstyle \bullet$} [c] at  10 0
\put {$ \scriptstyle \bullet$} [c] at  13 0
\put {$ \scriptstyle \bullet$} [c] at  13 6
\put {$ \scriptstyle \bullet$} [c] at  13 12
\put {$ \scriptstyle \bullet$} [c] at  11.5 6
\put {$ \scriptstyle \bullet$} [c] at  16 12
\setlinear \plot  16 12 13 0 10 12 10 0 13 6 13 0  /
\setlinear \plot 13  12 13 6    /
\put{$5{,}040$} [c] at 13 -2
\endpicture
\end{minipage}
\begin{minipage}{4cm}
\beginpicture
\setcoordinatesystem units   <1.5mm,2mm>
\setplotarea x from 0 to 16, y from -2 to 15
\put{1.221)} [l] at 2 12
\put {$ \scriptstyle \bullet$} [c] at  10 12
\put {$ \scriptstyle \bullet$} [c] at  10 0
\put {$ \scriptstyle \bullet$} [c] at  13 0
\put {$ \scriptstyle \bullet$} [c] at  13 6
\put {$ \scriptstyle \bullet$} [c] at  13 12
\put {$ \scriptstyle \bullet$} [c] at  11.5 6
\put {$ \scriptstyle \bullet$} [c] at  16 0
\setlinear \plot  16 0 13 12 10 0 10 12 13 6 13 12  /
\setlinear \plot 13  0 13 6    /
\put{$5{,}040$} [c] at 13 -2
\endpicture
\end{minipage}
\begin{minipage}{4cm}
\beginpicture
\setcoordinatesystem units   <1.5mm,2mm>
\setplotarea x from 0 to 16, y from -2 to 15
\put{1.222)} [l] at 2 12
\put {$ \scriptstyle \bullet$} [c] at  10 12
\put {$ \scriptstyle \bullet$} [c] at  10  0
\put {$ \scriptstyle \bullet$} [c] at  13 0
\put {$ \scriptstyle \bullet$} [c] at  13 12
\put {$ \scriptstyle \bullet$} [c] at  16 0
\put {$ \scriptstyle \bullet$} [c] at  16  12
\put {$ \scriptstyle \bullet$} [c] at  16  6
\setlinear \plot  10  0 10 12 13 0 16  12 16  0  13 12 13 0 /
\setlinear \plot  10 12 16 6 /
\put{$5{,}040$} [c] at 13 -2
\endpicture
\end{minipage}
\begin{minipage}{4cm}
\beginpicture
\setcoordinatesystem units   <1.5mm,2mm>
\setplotarea x from 0 to 16, y from -2 to 15
\put{1.223)} [l] at 2 12
\put {$ \scriptstyle \bullet$} [c] at  10 12
\put {$ \scriptstyle \bullet$} [c] at  10  0
\put {$ \scriptstyle \bullet$} [c] at  13 0
\put {$ \scriptstyle \bullet$} [c] at  13 12
\put {$ \scriptstyle \bullet$} [c] at  16 0
\put {$ \scriptstyle \bullet$} [c] at  16  12
\put {$ \scriptstyle \bullet$} [c] at  16  6
\setlinear \plot  10  12 10 0 13 12 16   0 16  12  13 0 13 12 /
\setlinear \plot  10 0 16 6 /
\put{$5{,}040$} [c] at 13 -2
\endpicture
\end{minipage}
\begin{minipage}{4cm}
\beginpicture
\setcoordinatesystem units   <1.5mm,2mm>
\setplotarea x from 0 to 16, y from -2 to 15
\put{1.224)} [l] at 2 12
\put {$ \scriptstyle \bullet$} [c] at  10 0
\put {$ \scriptstyle \bullet$} [c] at  10 12
\put {$ \scriptstyle \bullet$} [c] at  13 12
\put {$ \scriptstyle \bullet$} [c] at  13 0
\put {$ \scriptstyle \bullet$} [c] at  16 12
\put {$ \scriptstyle \bullet$} [c] at  16 6
\put {$ \scriptstyle \bullet$} [c] at  16 0
\setlinear \plot  10 0 10  12 13 0 13 12 16 6 10 12   /
\setlinear \plot  16 0 16 12   /
\put{$5{,}040$} [c] at 13 -2
\endpicture
\end{minipage}
$$
$$
\begin{minipage}{4cm}
\beginpicture
\setcoordinatesystem units   <1.5mm,2mm>
\setplotarea x from 0 to 16, y from -2 to 15
\put{1.225)} [l] at 2 12
\put {$ \scriptstyle \bullet$} [c] at  10 0
\put {$ \scriptstyle \bullet$} [c] at  10 12
\put {$ \scriptstyle \bullet$} [c] at  13 12
\put {$ \scriptstyle \bullet$} [c] at  13 0
\put {$ \scriptstyle \bullet$} [c] at  16 12
\put {$ \scriptstyle \bullet$} [c] at  16 6
\put {$ \scriptstyle \bullet$} [c] at  16 0
\setlinear \plot  10 12 10  0 13 12 13 0 16 6 10 0   /
\setlinear \plot  16 0 16 12   /
\put{$5{,}040$} [c] at 13 -2
\endpicture
\end{minipage}
\begin{minipage}{4cm}
\beginpicture
\setcoordinatesystem units   <1.5mm,2mm>
\setplotarea x from 0 to 16, y from -2 to 15
\put{1.226)} [l] at 2 12
\put {$ \scriptstyle \bullet$} [c] at  10 0
\put {$ \scriptstyle \bullet$} [c] at  10 12
\put {$ \scriptstyle \bullet$} [c] at  13 0
\put {$ \scriptstyle \bullet$} [c] at  13 4
\put {$ \scriptstyle \bullet$} [c] at  13 12
\put {$ \scriptstyle \bullet$} [c] at  16  12
\put {$ \scriptstyle \bullet$} [c] at  16  0
\setlinear \plot 16 0 10 12 10 0  13 12 16 0 16 12 13 0 13 12  /
\setlinear \plot  10 12 13 0   /
\put{$5{,}040$} [c] at 13 -2
\endpicture
\end{minipage}
\begin{minipage}{4cm}
\beginpicture
\setcoordinatesystem units   <1.5mm,2mm>
\setplotarea x from 0 to 16, y from -2 to 15
\put{1.227)} [l] at 2 12
\put {$ \scriptstyle \bullet$} [c] at  10 0
\put {$ \scriptstyle \bullet$} [c] at  10 12
\put {$ \scriptstyle \bullet$} [c] at  13 0
\put {$ \scriptstyle \bullet$} [c] at  13 12
\put {$ \scriptstyle \bullet$} [c] at  16 6
\put {$ \scriptstyle \bullet$} [c] at  16 0
\put {$ \scriptstyle \bullet$} [c] at  16  12
\setlinear \plot 13 12 10  0 10 12 13 0 13 12 16 0 16 12 13 0   /
\put{$5{,}040 $} [c] at 13 -2
\endpicture
\end{minipage}
\begin{minipage}{4cm}
\beginpicture
\setcoordinatesystem units   <1.5mm,2mm>
\setplotarea x from 0 to 16, y from -2 to 15
\put{1.228)} [l] at 2 12
\put {$ \scriptstyle \bullet$} [c] at  10 0
\put {$ \scriptstyle \bullet$} [c] at  10 12
\put {$ \scriptstyle \bullet$} [c] at  13 0
\put {$ \scriptstyle \bullet$} [c] at  13 6
\put {$ \scriptstyle \bullet$} [c] at  13  12
\put {$ \scriptstyle \bullet$} [c] at  16 12
\put {$ \scriptstyle \bullet$} [c] at  16  0
\setlinear \plot  10  12 10 0 16 12 16 0    /
\setlinear \plot  13 12 13 0  /
\put{$5{,}040$} [c] at 13 -2
\endpicture
\end{minipage}
\begin{minipage}{4cm}
\beginpicture
\setcoordinatesystem units <1.5mm, 2mm>
\setplotarea x from 0 to 16, y from -2 to 15
\put{1.229)} [l] at 2 12
\put {$ \scriptstyle \bullet$} [c] at 10 6
\put {$ \scriptstyle \bullet$} [c] at 12  6
\put {$ \scriptstyle \bullet$} [c] at 14 6
\put {$ \scriptstyle \bullet$} [c] at 12 12
\put {$ \scriptstyle \bullet$} [c] at 16 12
\put {$ \scriptstyle \bullet$} [c] at 13 0
\put {$ \scriptstyle \bullet$} [c] at 16 6
\setlinear \plot 12 12 10 6 13 0 14 6 16 12 16 6 13 0 12 6 12 12 14 6 /
\put{$2{,}520$} [c] at 13 -2
\endpicture
\end{minipage}
\begin{minipage}{4cm}
\beginpicture
\setcoordinatesystem units <1.5mm, 2mm>
\setplotarea x from 0 to 16, y from -2 to 15
\put{1.230)} [l] at 2 12
\put {$ \scriptstyle \bullet$} [c] at 10 6
\put {$ \scriptstyle \bullet$} [c] at 12  6
\put {$ \scriptstyle \bullet$} [c] at 14 6
\put {$ \scriptstyle \bullet$} [c] at 12 0
\put {$ \scriptstyle \bullet$} [c] at 16 0
\put {$ \scriptstyle \bullet$} [c] at 13 12
\put {$ \scriptstyle \bullet$} [c] at 16 6
\setlinear \plot 12 0 10 6 13 12 14 6 16 0 16 6 13 12 12 6 12 0 14 6 /
\put{$2{,}520$} [c] at 13 -2
\endpicture
\end{minipage}
$$

$$
\begin{minipage}{4cm}
\beginpicture
\setcoordinatesystem units <1.5mm, 2mm>
\setplotarea x from 0 to 16, y from -2 to 15
\put{1.231)} [l] at 2 12
\put {$ \scriptstyle \bullet$} [c] at 10 12
\put {$ \scriptstyle \bullet$} [c] at 11.5 12
\put {$ \scriptstyle \bullet$} [c] at 13 12
\put {$ \scriptstyle \bullet$} [c] at 16 12
\put {$ \scriptstyle \bullet$} [c] at 12 6
\put {$ \scriptstyle \bullet$} [c] at 14 6
\put {$ \scriptstyle \bullet$} [c] at 13 0
\setlinear \plot 10 12 12 6 13 0 14 6  16 12  /
\setlinear \plot 11.5 12 12 6 13 12 14 6    /
\put{$2{,}520$} [c] at 13 -2
\endpicture
\end{minipage}
\begin{minipage}{4cm}
\beginpicture
\setcoordinatesystem units <1.5mm, 2mm>
\setplotarea x from 0 to 16, y from -2 to 15
\put{1.232)} [l] at 2 12
\put {$ \scriptstyle \bullet$} [c] at 10 0
\put {$ \scriptstyle \bullet$} [c] at 11.5 0
\put {$ \scriptstyle \bullet$} [c] at 13 0
\put {$ \scriptstyle \bullet$} [c] at 16 0
\put {$ \scriptstyle \bullet$} [c] at 12 6
\put {$ \scriptstyle \bullet$} [c] at 14 6
\put {$ \scriptstyle \bullet$} [c] at 13 12
\setlinear \plot 10 0 12 6 13 12 14 6  16 0  /
\setlinear \plot 11.5 0 12 6 13 0 14 6    /
\put{$2{,}520$} [c] at 13 -2
\endpicture
\end{minipage}
\begin{minipage}{4cm}
\beginpicture
\setcoordinatesystem units <1.5mm, 2mm>
\setplotarea x from 0 to 16, y from -2 to 15
\put{1.233)} [l] at 2 12
\put {$ \scriptstyle \bullet$} [c] at 10 12
\put {$ \scriptstyle \bullet$} [c] at 11.5 7.5
\put {$ \scriptstyle \bullet$} [c] at 12.5 7.5
\put {$ \scriptstyle \bullet$} [c] at 12 0
\put {$ \scriptstyle \bullet$} [c] at 12 3
\put {$ \scriptstyle \bullet$} [c] at 12 12
\put {$ \scriptstyle \bullet$} [c] at 16 0
\setlinear \plot 10 12 12 0 12 3 11.5 7.5 12 12 12.5 7.5  12 3   /
\setlinear \plot 16 0 12 12    /
\put{$2{,}520$} [c] at 13 -2
\endpicture
\end{minipage}
\begin{minipage}{4cm}
\beginpicture
\setcoordinatesystem units <1.5mm, 2mm>
\setplotarea x from 0 to 16, y from -2 to 15
\put{1.234)} [l] at 2 12
\put {$ \scriptstyle \bullet$} [c] at 10 0
\put {$ \scriptstyle \bullet$} [c] at 11.5 4.5
\put {$ \scriptstyle \bullet$} [c] at 12.5 4.5
\put {$ \scriptstyle \bullet$} [c] at 12 12
\put {$ \scriptstyle \bullet$} [c] at 12 9
\put {$ \scriptstyle \bullet$} [c] at 12 0
\put {$ \scriptstyle \bullet$} [c] at 16 12
\setlinear \plot 10 0 12 12 12 9 11.5 4.5 12 0 12.5 4.5  12 9   /
\setlinear \plot 16 12 12 0    /
\put{$2{,}520$} [c] at 13 -2
\endpicture
\end{minipage}
\begin{minipage}{4cm}
\beginpicture
\setcoordinatesystem units <1.5mm, 2mm>
\setplotarea x from 0 to 16, y from -2 to 15
\put{1.235)} [l] at 2 12
\put {$ \scriptstyle \bullet$} [c] at 10 12
\put {$ \scriptstyle \bullet$} [c] at 12 3
\put {$ \scriptstyle \bullet$} [c] at 13 0
\put {$ \scriptstyle \bullet$} [c] at 13 3
\put {$ \scriptstyle \bullet$} [c] at 13 12
\put {$ \scriptstyle \bullet$} [c] at 14 3
\put {$ \scriptstyle \bullet$} [c] at 16 0
\setlinear \plot 13 3 10 12 12 3 13 0 14 3  13 12 12 3 /
\setlinear \plot 16 0 13 12 13 0    /
\put{$2{,}520$} [c] at 13 -2
\endpicture
\end{minipage}
\begin{minipage}{4cm}
\beginpicture
\setcoordinatesystem units <1.5mm, 2mm>
\setplotarea x from 0 to 16, y from -2 to 15
\put{1.236)} [l] at 2 12
\put {$ \scriptstyle \bullet$} [c] at 10 0
\put {$ \scriptstyle \bullet$} [c] at 12 9
\put {$ \scriptstyle \bullet$} [c] at 13 0
\put {$ \scriptstyle \bullet$} [c] at 13 9
\put {$ \scriptstyle \bullet$} [c] at 13 12
\put {$ \scriptstyle \bullet$} [c] at 14 9
\put {$ \scriptstyle \bullet$} [c] at 16 12
\setlinear \plot 13 9 10 0 12 9 13 12 14 9  13 0 12 9 /
\setlinear \plot 16 12 13 0 13 12    /
\put{$2{,}520$} [c] at 13 -2
\endpicture
\end{minipage}
$$
$$
\begin{minipage}{4cm}
\beginpicture
\setcoordinatesystem units <1.5mm, 2mm>
\setplotarea x from 0 to 16, y from -2 to 15
\put{1.237)} [l] at 2 12
\put {$ \scriptstyle \bullet$} [c] at 10 6
\put {$ \scriptstyle \bullet$} [c] at 11 0
\put {$ \scriptstyle \bullet$} [c] at 11 6
\put {$ \scriptstyle \bullet$} [c] at 11 12
\put {$ \scriptstyle \bullet$} [c] at 12 6
\put {$ \scriptstyle \bullet$} [c] at 16 0
\put {$ \scriptstyle \bullet$} [c] at 16 12
\setlinear \plot 11 12 10 6 11 0 12 6 11 12 16 0 16 12 12 6 /
\setlinear \plot 11 12 11 0    /
\put{$2{,}520$} [c] at 13 -2
\endpicture
\end{minipage}
\begin{minipage}{4cm}
\beginpicture
\setcoordinatesystem units <1.5mm, 2mm>
\setplotarea x from 0 to 16, y from -2 to 15
\put{1.238)} [l] at 2 12
\put {$ \scriptstyle \bullet$} [c] at 10 6
\put {$ \scriptstyle \bullet$} [c] at 11 0
\put {$ \scriptstyle \bullet$} [c] at 11 6
\put {$ \scriptstyle \bullet$} [c] at 11 12
\put {$ \scriptstyle \bullet$} [c] at 12 6
\put {$ \scriptstyle \bullet$} [c] at 16 0
\put {$ \scriptstyle \bullet$} [c] at 16 12
\setlinear \plot 11 0 10 6 11 12 12 6 11 0 16 12 16 0 12 6 /
\setlinear \plot 11 12 11 0    /
\put{$2{,}520$} [c] at 13 -2
\endpicture
\end{minipage}
\begin{minipage}{4cm}
\beginpicture
\setcoordinatesystem units   <1.5mm,2mm>
\setplotarea x from 0 to 16, y from -2 to 15
\put{1.239)} [l] at 2 12
\put {$ \scriptstyle \bullet$} [c] at  10 0
\put {$ \scriptstyle \bullet$} [c] at  10 6
\put {$ \scriptstyle \bullet$} [c] at  10 12
\put {$ \scriptstyle \bullet$} [c] at  11.7 7
\put {$ \scriptstyle \bullet$} [c] at  13  12
\put {$ \scriptstyle \bullet$} [c] at  16 0
\put {$ \scriptstyle \bullet$} [c] at  16 12
\setlinear \plot  10 0 10 12 16 0 16 12  10 0 13 12 16 0 /
\put{$2{,}520$} [c] at 13 -2
\endpicture
\end{minipage}
\begin{minipage}{4cm}
\beginpicture
\setcoordinatesystem units   <1.5mm,2mm>
\setplotarea x from 0 to 16, y from -2 to 15
\put{1.240)} [l] at 2 12
\put {$ \scriptstyle \bullet$} [c] at  10 0
\put {$ \scriptstyle \bullet$} [c] at  10 6
\put {$ \scriptstyle \bullet$} [c] at  10 12
\put {$ \scriptstyle \bullet$} [c] at  11.2 7
\put {$ \scriptstyle \bullet$} [c] at  13  0
\put {$ \scriptstyle \bullet$} [c] at  16 0
\put {$ \scriptstyle \bullet$} [c] at  16 12
\setlinear \plot  10 12 10 0 16 12 16 0  10 12 13 0 16 12 /
\put{$2{,}520$} [c] at 13 -2
\endpicture
\end{minipage}
\begin{minipage}{4cm}
\beginpicture
\setcoordinatesystem units   <1.5mm,2mm>
\setplotarea x from 0 to 16, y from -2 to 15
\put{1.241)} [l] at 2 12
\put {$ \scriptstyle \bullet$} [c] at  10 12
\put {$ \scriptstyle \bullet$} [c] at  12 6.5
\put {$ \scriptstyle \bullet$} [c] at  13 0
\put {$ \scriptstyle \bullet$} [c] at  13 12
\put {$ \scriptstyle \bullet$} [c] at  14 6.5
\put {$ \scriptstyle \bullet$} [c] at  16 0
\put {$ \scriptstyle \bullet$} [c] at  16 12
\setlinear \plot 16 0 14 6.5 13 12 12 6.5 13 0  10  12 16 0 16 12  13 0  /
\setlinear \plot 13 0 14 6.5  /
\put{$2{,}520$} [c] at 13 -2
\endpicture
\end{minipage}
\begin{minipage}{4cm}
\beginpicture
\setcoordinatesystem units   <1.5mm,2mm>
\setplotarea x from 0 to 16, y from -2 to 15
\put{1.242)} [l] at 2 12
\put {$ \scriptstyle \bullet$} [c] at  10 0
\put {$ \scriptstyle \bullet$} [c] at  12 5.5
\put {$ \scriptstyle \bullet$} [c] at  13 0
\put {$ \scriptstyle \bullet$} [c] at  13 12
\put {$ \scriptstyle \bullet$} [c] at  14 5.5
\put {$ \scriptstyle \bullet$} [c] at  16 0
\put {$ \scriptstyle \bullet$} [c] at  16 12
\setlinear \plot 16 12 14 5.5 13 0 12 5.5 13 12  10  0 16 12 16 0  13 12  /
\setlinear \plot 13 12 14 5.5  /
\put{$2{,}520$} [c] at 13 -2
\endpicture
\end{minipage}
$$

$$
\begin{minipage}{4cm}
\beginpicture
\setcoordinatesystem units   <1.5mm,2mm>
\setplotarea x from 0 to 16, y from -2 to 15
\put{1.243)} [l] at 2 12
\put {$ \scriptstyle \bullet$} [c] at  16 12
\put {$ \scriptstyle \bullet$} [c] at  13 0
\put {$ \scriptstyle \bullet$} [c] at  13 12
\put {$ \scriptstyle \bullet$} [c] at  10 0
\put {$ \scriptstyle \bullet$} [c] at  10 4
\put {$ \scriptstyle \bullet$} [c] at  10 8
\put {$ \scriptstyle \bullet$} [c] at  10 12
\setlinear \plot  10 0 10 12 13 0 16 12 /
\setlinear \plot  13 12 13 0 /
\put{$2{,}520$} [c] at 13 -2
\endpicture
\end{minipage}
\begin{minipage}{4cm}
\beginpicture
\setcoordinatesystem units   <1.5mm,2mm>
\setplotarea x from 0 to 16, y from -2 to 15
\put{1.244)} [l] at 2 12
\put {$ \scriptstyle \bullet$} [c] at  16 0
\put {$ \scriptstyle \bullet$} [c] at  13 0
\put {$ \scriptstyle \bullet$} [c] at  13 12
\put {$ \scriptstyle \bullet$} [c] at  10 0
\put {$ \scriptstyle \bullet$} [c] at  10 4
\put {$ \scriptstyle \bullet$} [c] at  10 8
\put {$ \scriptstyle \bullet$} [c] at  10 12
\setlinear \plot  10 12 10 0 13 12 16 0 /
\setlinear \plot  13 12 13 0 /
\put{$2{,}520$} [c] at 13 -2
\endpicture
\end{minipage}
\begin{minipage}{4cm}
\beginpicture
\setcoordinatesystem units   <1.5mm,2mm>
\setplotarea x from 0 to 16, y from -2 to 15
\put{1.245)} [l] at 2 12
\put {$ \scriptstyle \bullet$} [c] at  10 0
\put {$ \scriptstyle \bullet$} [c] at  14 6
\put {$ \scriptstyle \bullet$} [c] at  14 12
\put {$ \scriptstyle \bullet$} [c] at  15 0
\put {$ \scriptstyle \bullet$} [c] at  15 12
\put {$ \scriptstyle \bullet$} [c] at  16 6
\put {$ \scriptstyle \bullet$} [c] at  16 12
\setlinear \plot  10 0 15 12 14 6 15 0 16 6  16 12   /
\setlinear \plot  14 12 14  6  /
\setlinear \plot  16 12 16 6  15 12   /
\put{$2{,}520$} [c] at 13 -2
\endpicture
\end{minipage}
\begin{minipage}{4cm}
\beginpicture
\setcoordinatesystem units   <1.5mm,2mm>
\setplotarea x from 0 to 16, y from -2 to 15
\put{1.246)} [l] at 2 12
\put {$ \scriptstyle \bullet$} [c] at  10 12
\put {$ \scriptstyle \bullet$} [c] at  14 6
\put {$ \scriptstyle \bullet$} [c] at  14 0
\put {$ \scriptstyle \bullet$} [c] at  15 0
\put {$ \scriptstyle \bullet$} [c] at  15 12
\put {$ \scriptstyle \bullet$} [c] at  16 6
\put {$ \scriptstyle \bullet$} [c] at  16 0
\setlinear \plot  10 12 15 0 14 6 15 12 16 6  16 0   /
\setlinear \plot  14 0 14  6  /
\setlinear \plot  16 0 16 6  15 0   /
\put{$2{,}520$} [c] at 13 -2
\endpicture
\end{minipage}
\begin{minipage}{4cm}
\beginpicture
\setcoordinatesystem units   <1.5mm,2mm>
\setplotarea x from 0 to 16, y from -2 to 15
\put{1.247)} [l] at 2 12
\put {$ \scriptstyle \bullet$} [c] at  10 12
\put {$ \scriptstyle \bullet$} [c] at  12  0
\put {$ \scriptstyle \bullet$} [c] at  12 6
\put {$ \scriptstyle \bullet$} [c] at  12 12
\put {$ \scriptstyle \bullet$} [c] at  14 9
\put {$ \scriptstyle \bullet$} [c] at  16 0
\put {$ \scriptstyle \bullet$} [c] at  16 12
\setlinear \plot  10 12 12 6 16 12 16 0 /
\setlinear \plot  12 12 12 0  /
\put{$2{,}520$} [c] at 13 -2
\endpicture
\end{minipage}
\begin{minipage}{4cm}
\beginpicture
\setcoordinatesystem units   <1.5mm,2mm>
\setplotarea x from 0 to 16, y from -2 to 15
\put{1.248)} [l] at 2 12
\put {$ \scriptstyle \bullet$} [c] at  10 0
\put {$ \scriptstyle \bullet$} [c] at  12  0
\put {$ \scriptstyle \bullet$} [c] at  12 6
\put {$ \scriptstyle \bullet$} [c] at  12 12
\put {$ \scriptstyle \bullet$} [c] at  14 3
\put {$ \scriptstyle \bullet$} [c] at  16 0
\put {$ \scriptstyle \bullet$} [c] at  16 12
\setlinear \plot  10 0 12 6 16 0 16 12 /
\setlinear \plot  12 12 12 0  /
\put{$2{,}520$} [c] at 13 -2
\endpicture
\end{minipage}
$$
$$
\begin{minipage}{4cm}
\beginpicture
\setcoordinatesystem units   <1.5mm,2mm>
\setplotarea x from 0 to 16, y from -2 to 15
\put{1.249)} [l] at 2 12
\put {$ \scriptstyle \bullet$} [c] at  10 12
\put {$ \scriptstyle \bullet$} [c] at  12 0
\put {$ \scriptstyle \bullet$} [c] at  12 3
\put {$ \scriptstyle \bullet$} [c] at  13 12
\put {$ \scriptstyle \bullet$} [c] at  14 0
\put {$ \scriptstyle \bullet$} [c] at  14 3
\put {$ \scriptstyle \bullet$} [c] at  16 12
\setlinear \plot  10 12   12 0  12 3 13 12 14 3 14 0 12 3 12 0  14 3 /
\setlinear \plot  12 0 16 12  14 0  /
\put{$2{,}520$} [c] at 13 -2
\endpicture
\end{minipage}
\begin{minipage}{4cm}
\beginpicture
\setcoordinatesystem units   <1.5mm,2mm>
\setplotarea x from 0 to 16, y from -2 to 15
\put{1.250)} [l] at 2 12
\put {$ \scriptstyle \bullet$} [c] at  10 0
\put {$ \scriptstyle \bullet$} [c] at  12 12
\put {$ \scriptstyle \bullet$} [c] at  12 9
\put {$ \scriptstyle \bullet$} [c] at  13 0
\put {$ \scriptstyle \bullet$} [c] at  14 12
\put {$ \scriptstyle \bullet$} [c] at  14 9
\put {$ \scriptstyle \bullet$} [c] at  16 0
\setlinear \plot  10 0   12 12  12 9 13 0 14 9 14 12 12 9 12 12  14 9 /
\setlinear \plot  12 12 16 0  14 12  /
\put{$2{,}520$} [c] at 13 -2
\endpicture
\end{minipage}
\begin{minipage}{4cm}
\beginpicture
\setcoordinatesystem units   <1.5mm,2mm>
\setplotarea x from 0 to 16, y from -2 to 15
\put {1.251)} [l] at 2 12
\put {$ \scriptstyle \bullet$} [c] at  10 0
\put {$ \scriptstyle \bullet$} [c] at  10 12
\put {$ \scriptstyle \bullet$} [c] at  12 12
\put {$ \scriptstyle \bullet$} [c] at  14 0
\put {$ \scriptstyle \bullet$} [c] at  14 6
\put {$ \scriptstyle \bullet$} [c] at  14  12
\put {$ \scriptstyle \bullet$} [c] at  16  12
\setlinear \plot  10 12 10 0  14 12 14 0 /
\setlinear \plot  10 0 12 12 14 6 /
\setlinear \plot 16 12 14 6 /
\put{$2{,}520$}[c] at 13 -2
\endpicture
\end{minipage}
\begin{minipage}{4cm}
\beginpicture
\setcoordinatesystem units   <1.5mm,2mm>
\setplotarea x from 0 to 16, y from -2 to 15
\put {1.252)} [l] at 2 12
\put {$ \scriptstyle \bullet$} [c] at  10 0
\put {$ \scriptstyle \bullet$} [c] at  10 12
\put {$ \scriptstyle \bullet$} [c] at  12 0
\put {$ \scriptstyle \bullet$} [c] at  14 0
\put {$ \scriptstyle \bullet$} [c] at  14 6
\put {$ \scriptstyle \bullet$} [c] at  14  12
\put {$ \scriptstyle \bullet$} [c] at  16  0
\setlinear \plot  10 0 10 12  14 0 14 12 /
\setlinear \plot  10 12 12 0 14 6 /
\setlinear \plot 16 0 14 6 /
\put{$2{,}520$}[c] at 13 -2
\endpicture
\end{minipage}
\begin{minipage}{4cm}
\beginpicture
\setcoordinatesystem units   <1.5mm,2mm>
\setplotarea x from 0 to 16, y from -2 to 15
\put {1.253)} [l] at 2 12
\put {$ \scriptstyle \bullet$} [c] at  10 0
\put {$ \scriptstyle \bullet$} [c] at  10 6
\put {$ \scriptstyle \bullet$} [c] at  10 12
\put {$ \scriptstyle \bullet$} [c] at  12 12
\put {$ \scriptstyle \bullet$} [c] at  14 0
\put {$ \scriptstyle \bullet$} [c] at  14 12
\put {$ \scriptstyle \bullet$} [c] at  16 12
\setlinear \plot  10 12 10 0  14 12 14 0  10 6 12 12 /
\setlinear \plot 16 12 14 0  /
\put{$2{,}520$}[c] at 13 -2
\endpicture
\end{minipage}
\begin{minipage}{4cm}
\beginpicture
\setcoordinatesystem units   <1.5mm,2mm>
\setplotarea x from 0 to 16, y from -2 to 15
\put {1.254)} [l] at 2 12
\put {$ \scriptstyle \bullet$} [c] at  10 0
\put {$ \scriptstyle \bullet$} [c] at  10 6
\put {$ \scriptstyle \bullet$} [c] at  10 12
\put {$ \scriptstyle \bullet$} [c] at  12 0
\put {$ \scriptstyle \bullet$} [c] at  14 0
\put {$ \scriptstyle \bullet$} [c] at  14 12
\put {$ \scriptstyle \bullet$} [c] at  16 0
\setlinear \plot  10 0 10 12  14 0 14 12  10 6 12 0 /
\setlinear \plot 16 0 14 12  /
\put{$2{,}520$}[c] at 13 -2
\endpicture
\end{minipage}
$$

$$
\begin{minipage}{4cm}
\beginpicture
\setcoordinatesystem units   <1.5mm,2mm>
\setplotarea x from 0 to 16, y from -2 to 15
\put{1.255)} [l] at 2 12
\put {$ \scriptstyle \bullet$} [c] at  10 0
\put {$ \scriptstyle \bullet$} [c] at  10 12
\put {$ \scriptstyle \bullet$} [c] at  13 12
\put {$ \scriptstyle \bullet$} [c] at  13 0
\put {$ \scriptstyle \bullet$} [c] at  16 0
\put {$ \scriptstyle \bullet$} [c] at  16  6
\put {$ \scriptstyle \bullet$} [c] at  16  12
\setlinear \plot  16 12 16  6 13 0 13  12  10 0 10  12  13 0  /
\setlinear \plot   16  0  16 6    /
\put{$2{,}520 $} [c] at 13 -2
\endpicture
\end{minipage}
\begin{minipage}{4cm}
\beginpicture
\setcoordinatesystem units   <1.5mm,2mm>
\setplotarea x from 0 to 16, y from -2 to 15
\put{1.256)} [l] at 2 12
\put {$ \scriptstyle \bullet$} [c] at  10 0
\put {$ \scriptstyle \bullet$} [c] at  10 12
\put {$ \scriptstyle \bullet$} [c] at  13 12
\put {$ \scriptstyle \bullet$} [c] at  13 0
\put {$ \scriptstyle \bullet$} [c] at  16 0
\put {$ \scriptstyle \bullet$} [c] at  16  6
\put {$ \scriptstyle \bullet$} [c] at  16  12
\setlinear \plot  16 0 16  6 13 12 13 0  10 12 10  0  13 12  /
\setlinear \plot   16  12  16 6    /
\put{$2{,}520 $} [c] at 13 -2
\endpicture
\end{minipage}
\begin{minipage}{4cm}
\beginpicture
\setcoordinatesystem units   <1.5mm,2mm>
\setplotarea x from 0 to 16, y from -2 to 15
\put{1.257)} [l] at 2 12
\put {$ \scriptstyle \bullet$} [c] at  10 0
\put {$ \scriptstyle \bullet$} [c] at  10 12
\put {$ \scriptstyle \bullet$} [c] at  13 0
\put {$ \scriptstyle \bullet$} [c] at  13 12
\put {$ \scriptstyle \bullet$} [c] at  16 0
\put {$ \scriptstyle \bullet$} [c] at  16  6
\put {$ \scriptstyle \bullet$} [c] at  16  12
\setlinear \plot 16 12 16 0 10 12  13 0 13 12 16 0    /
\setlinear \plot   10  0  10 12    /
\setlinear \plot   13  0  16 6    /
\put{$2{,}520$} [c] at 13 -2
\endpicture
\end{minipage}
\begin{minipage}{4cm}
\beginpicture
\setcoordinatesystem units   <1.5mm,2mm>
\setplotarea x from 0 to 16, y from -2 to 15
\put{1.258)} [l] at 2 12
\put {$ \scriptstyle \bullet$} [c] at  10 0
\put {$ \scriptstyle \bullet$} [c] at  10 12
\put {$ \scriptstyle \bullet$} [c] at  13 0
\put {$ \scriptstyle \bullet$} [c] at  13 12
\put {$ \scriptstyle \bullet$} [c] at  16 0
\put {$ \scriptstyle \bullet$} [c] at  16  6
\put {$ \scriptstyle \bullet$} [c] at  16  12
\setlinear \plot 16 0 16 12 10 0  13 12 13 0 16 12    /
\setlinear \plot   10  0  10 12    /
\setlinear \plot   13  12  16 6    /
\put{$2{,}520$} [c] at 13 -2
\endpicture
\end{minipage}
\begin{minipage}{4cm}
\beginpicture
\setcoordinatesystem units   <1.5mm,2mm>
\setplotarea x from 0 to 16, y from -2 to 15
\put{1.259)} [l] at 2 12
\put {$ \scriptstyle \bullet$} [c] at  10 0
\put {$ \scriptstyle \bullet$} [c] at  10 6
\put {$ \scriptstyle \bullet$} [c] at  10 12
\put {$ \scriptstyle \bullet$} [c] at  13 12
\put {$ \scriptstyle \bullet$} [c] at  13 0
\put {$ \scriptstyle \bullet$} [c] at  16  0
\put {$ \scriptstyle \bullet$} [c] at  16  12
\setlinear \plot 10 12 10 0  13 12 16 0 16 12    /
\setlinear \plot   10 6  13 0 13 12    /
\put{$2{,}520$} [c] at 13 -2
\endpicture
\end{minipage}
\begin{minipage}{4cm}
\beginpicture
\setcoordinatesystem units   <1.5mm,2mm>
\setplotarea x from 0 to 16, y from -2 to 15
\put{1.260)} [l] at 2 12
\put {$ \scriptstyle \bullet$} [c] at  10 0
\put {$ \scriptstyle \bullet$} [c] at  10 6
\put {$ \scriptstyle \bullet$} [c] at  10 12
\put {$ \scriptstyle \bullet$} [c] at  13 12
\put {$ \scriptstyle \bullet$} [c] at  13 0
\put {$ \scriptstyle \bullet$} [c] at  16  0
\put {$ \scriptstyle \bullet$} [c] at  16  12
\setlinear \plot 10 0 10 12  13 0 16 12 16 0    /
\setlinear \plot   10 6  13 12 13 0    /
\put{$2{,}520$} [c] at 13 -2
\endpicture
\end{minipage}
$$
$$
\begin{minipage}{4cm}
\beginpicture
\setcoordinatesystem units   <1.5mm,2mm>
\setplotarea x from 0 to 16, y from -2 to 15
\put{1.261)} [l] at 2 12
\put {$ \scriptstyle \bullet$} [c] at  16 12
\put {$ \scriptstyle \bullet$} [c] at  10 0
\put {$ \scriptstyle \bullet$} [c] at  13 0
\put {$ \scriptstyle \bullet$} [c] at  13.5 5
\put {$ \scriptstyle \bullet$} [c] at  12.5 5
\put {$ \scriptstyle \bullet$} [c] at  13.5 12
\put {$ \scriptstyle \bullet$} [c] at  12.5 12
\setlinear \plot 10 0 12.5 12 12.5 5 13 0 13.5 5 13.5 12 12.5 5 /
\setlinear \plot  10 0 13.5 12  /
\setlinear \plot 16 12 13 0    /
\setlinear \plot 12.5 12  13.5 5    /
\put{$1{,}260$} [c] at 13 -2
\endpicture
\end{minipage}
\begin{minipage}{4cm}
\beginpicture
\setcoordinatesystem units   <1.5mm,2mm>
\setplotarea x from 0 to 16, y from -2 to 15
\put{1.262)} [l] at 2 12
\put {$ \scriptstyle \bullet$} [c] at  16 0
\put {$ \scriptstyle \bullet$} [c] at  10 12
\put {$ \scriptstyle \bullet$} [c] at  13 12
\put {$ \scriptstyle \bullet$} [c] at  13.5 7
\put {$ \scriptstyle \bullet$} [c] at  12.5 7
\put {$ \scriptstyle \bullet$} [c] at  13.5 0
\put {$ \scriptstyle \bullet$} [c] at  12.5 0
\setlinear \plot 10 12 12.5 0 12.5 7 13 12 13.5 7 13.5 0 12.5 7 /
\setlinear \plot  10 12 13.5 0  /
\setlinear \plot 16 0 13 12    /
\setlinear \plot 12.5 0  13.5 7    /
\put{$1{,}260$} [c] at 13 -2
\endpicture
\end{minipage}
\begin{minipage}{4cm}
\beginpicture
\setcoordinatesystem units   <1.5mm,2mm>
\setplotarea x from 0 to 16, y from -2 to 15
\put {1.263)} [l] at 2 12
\put {$ \scriptstyle \bullet$} [c] at  10 0
\put {$ \scriptstyle \bullet$} [c] at  10 12
\put {$ \scriptstyle \bullet$} [c] at  12 12
\put {$ \scriptstyle \bullet$} [c] at  14 0
\put {$ \scriptstyle \bullet$} [c] at  14 6
\put {$ \scriptstyle \bullet$} [c] at  14  12
\put {$ \scriptstyle \bullet$} [c] at  16  12
\setlinear \plot  10 12 10 0  14 12 14 0   10 12 /
\setlinear \plot  10 0 12 12 14 6 /
\setlinear \plot 10 0 16 12 14 0 /
\put{$1{,}260$}[c] at 13 -2
\endpicture
\end{minipage}
\begin{minipage}{4cm}
\beginpicture
\setcoordinatesystem units   <1.5mm,2mm>
\setplotarea x from 0 to 16, y from -2 to 15
\put {1.264)} [l] at 2 12
\put {$ \scriptstyle \bullet$} [c] at  10 0
\put {$ \scriptstyle \bullet$} [c] at  10 12
\put {$ \scriptstyle \bullet$} [c] at  12 0
\put {$ \scriptstyle \bullet$} [c] at  14 0
\put {$ \scriptstyle \bullet$} [c] at  14 6
\put {$ \scriptstyle \bullet$} [c] at  14  12
\put {$ \scriptstyle \bullet$} [c] at  16  0
\setlinear \plot  10 0 10 12  14 0 14 12   10 0 /
\setlinear \plot  10 12 12 0 14 6 /
\setlinear \plot 10 12 16 0 14 12 /
\put{$1{,}260$}[c] at 13 -2
\endpicture
\end{minipage}
\begin{minipage}{4cm}
\beginpicture
\setcoordinatesystem units   <1.5mm,2mm>
\setplotarea x from 0 to 16, y from -2 to 15
\put {1.265)} [l] at 2 12
\put {$ \scriptstyle \bullet$} [c] at  10 0
\put {$ \scriptstyle \bullet$} [c] at  12 12
\put {$ \scriptstyle \bullet$} [c] at  14 0
\put {$ \scriptstyle \bullet$} [c] at  16 12
\put {$ \scriptstyle \bullet$} [c] at  10 12
\put {$ \scriptstyle \bullet$} [c] at  14  6
\put {$ \scriptstyle \bullet$} [c] at  14 12
\setlinear \plot  10 12 10 0  12 12 14 6  16 12 /
\setlinear \plot  14 12 14 0   /
\setlinear \plot 10 12 14 6 /
\put{$1{,}260$}[c] at 13 -2
\endpicture
\end{minipage}
\begin{minipage}{4cm}
\beginpicture
\setcoordinatesystem units   <1.5mm,2mm>
\setplotarea x from 0 to 16, y from -2 to 15
\put {1.266)} [l] at 2 12
\put {$ \scriptstyle \bullet$} [c] at  10 0
\put {$ \scriptstyle \bullet$} [c] at  12 0
\put {$ \scriptstyle \bullet$} [c] at  14 0
\put {$ \scriptstyle \bullet$} [c] at  16 0
\put {$ \scriptstyle \bullet$} [c] at  10 12
\put {$ \scriptstyle \bullet$} [c] at  14  6
\put {$ \scriptstyle \bullet$} [c] at  14 12
\setlinear \plot  10 0 10 12  12 0 14 6  16 0 /
\setlinear \plot  14 0 14 12   /
\setlinear \plot 10 0 14 6 /
\put{$1{,}260$}[c] at 13 -2
\endpicture
\end{minipage}
$$

$$
\begin{minipage}{4cm}
\beginpicture
\setcoordinatesystem units   <1.5mm,2mm>
\setplotarea x from 0 to 16, y from -2 to 15
\put {1.267)} [l] at 2 12
\put {$ \scriptstyle \bullet$} [c] at  10 0
\put {$ \scriptstyle \bullet$} [c] at  10 12
\put {$ \scriptstyle \bullet$} [c] at  12 12
\put {$ \scriptstyle \bullet$} [c] at  14 0
\put {$ \scriptstyle \bullet$} [c] at  14 6
\put {$ \scriptstyle \bullet$} [c] at  14 12
\put {$ \scriptstyle \bullet$} [c] at  16  12
\setlinear \plot  10 12 10 0  14 12 14 0 /
\setlinear \plot  10 12 14 6 12 12 10 0 /
\setlinear \plot 16 12 14 0 /
\put{$840$}[c] at 13 -2
\endpicture
\end{minipage}
\begin{minipage}{4cm}
\beginpicture
\setcoordinatesystem units   <1.5mm,2mm>
\setplotarea x from 0 to 16, y from -2 to 15
\put {1.268)} [l] at 2 12
\put {$ \scriptstyle \bullet$} [c] at  10 0
\put {$ \scriptstyle \bullet$} [c] at  10 12
\put {$ \scriptstyle \bullet$} [c] at  12 0
\put {$ \scriptstyle \bullet$} [c] at  14 0
\put {$ \scriptstyle \bullet$} [c] at  14 6
\put {$ \scriptstyle \bullet$} [c] at  14 12
\put {$ \scriptstyle \bullet$} [c] at  16  0
\setlinear \plot  10 0 10 12  14 0 14 12 /
\setlinear \plot  10 0 14 6 12 0 10 12 /
\setlinear \plot 16 0 14 12 /
\put{$840$}[c] at 13 -2
\endpicture
\end{minipage}
\begin{minipage}{4cm}
\beginpicture
\setcoordinatesystem units   <1.5mm,2mm>
\setplotarea x from 0 to 16, y from -2 to 15
\put{1.269)} [l] at 2 12
\put {$ \scriptstyle \bullet$} [c] at 10 12
\put {$ \scriptstyle \bullet$} [c] at 12  12
\put {$ \scriptstyle \bullet$} [c] at 14 12
\put {$ \scriptstyle \bullet$} [c] at 16 12
\put {$ \scriptstyle \bullet$} [c] at 10 0
\put {$ \scriptstyle \bullet$} [c] at 16 0
\put{$\scriptstyle \bullet$} [c] at 13 0
\setlinear \plot 16 12 13 0  10 12 10 0 16 12 16 0 12 12 10 0 /
\setlinear \plot 10 0  14 12 16 0  /
\setlinear \plot 12 12 13 0  14 12  /
\put{$420$} [c] at 13 -2
\put{$\scriptstyle \bullet$} [c] at 16  0 \endpicture
\end{minipage}
\begin{minipage}{4cm}
\beginpicture
\setcoordinatesystem units   <1.5mm,2mm>
\setplotarea x from 0 to 16, y from -2 to 15
\put{1.270)} [l] at 2 12
\put {$ \scriptstyle \bullet$} [c] at 10 0
\put {$ \scriptstyle \bullet$} [c] at 12  0
\put {$ \scriptstyle \bullet$} [c] at 14 0
\put {$ \scriptstyle \bullet$} [c] at 16 0
\put {$ \scriptstyle \bullet$} [c] at 10 12
\put {$ \scriptstyle \bullet$} [c] at 16 12
\put{$\scriptstyle \bullet$} [c] at 13 12
\setlinear \plot 16 0 13 12  10 0 10 12 16 0 16 12 12 0 10 12 /
\setlinear \plot 10 12  14 0 16 12  /
\setlinear \plot 12 0 13 12  14 0  /
\put{$420$} [c] at 13 -2
 \endpicture
\end{minipage}
\begin{minipage}{4cm}
\beginpicture
\setcoordinatesystem units   <1.5mm,2mm>
\setplotarea x from 0 to 16, y from -2 to 15
\put{1.271)} [l] at 2 12
\put {$ \scriptstyle \bullet$} [c] at 10 6
\put {$ \scriptstyle \bullet$} [c] at 11 9
\put {$ \scriptstyle \bullet$} [c] at 12 0
\put {$ \scriptstyle \bullet$} [c] at 12 12
\put {$ \scriptstyle \bullet$} [c] at 14 6
\put {$ \scriptstyle \bullet$} [c] at 14 12
\put{$\scriptstyle \bullet$} [c] at 16  0
\setlinear \plot 14 6 12 12  10 6 12  0  14 6 14 12 /
\put{$5{,}040$} [c] at 13 -2
 \endpicture
\end{minipage}
\begin{minipage}{4cm}
\beginpicture
\setcoordinatesystem units   <1.5mm,2mm>
\setplotarea x from 0 to 16, y from -2 to 15
\put{1.272)} [l] at 2 12
\put {$ \scriptstyle \bullet$} [c] at 10 6
\put {$ \scriptstyle \bullet$} [c] at 11 9
\put {$ \scriptstyle \bullet$} [c] at 12 0
\put {$ \scriptstyle \bullet$} [c] at 12 12
\put {$ \scriptstyle \bullet$} [c] at 14 6
\put {$ \scriptstyle \bullet$} [c] at 14 0
\put{$\scriptstyle \bullet$} [c] at 16  0
\setlinear \plot 14 6 12 12  10 6 12  0  14 6 14 0 /
\put{$5{,}040$} [c] at 13 -2
\endpicture
\end{minipage}
$$
$$
\begin{minipage}{4cm}
\beginpicture
\setcoordinatesystem units   <1.5mm,2mm>
\setplotarea x from 0 to 16, y from -2 to 15
\put{1.273)} [l] at 2 12
\put {$ \scriptstyle \bullet$} [c] at 10 12
\put {$ \scriptstyle \bullet$} [c] at 10 9
\put {$ \scriptstyle \bullet$} [c] at 10 6
\put {$ \scriptstyle \bullet$} [c] at 12 0
\put {$ \scriptstyle \bullet$} [c] at 12 12
\put {$ \scriptstyle \bullet$} [c] at 14 6
\put{$\scriptstyle \bullet$} [c] at 16  0
\setlinear \plot 10 12 10 6 12 0 14  6 12  12 10 6  /
\put{$5{,}040$} [c] at 13 -2
\endpicture
\end{minipage}
\begin{minipage}{4cm}
\beginpicture
\setcoordinatesystem units   <1.5mm,2mm>
\setplotarea x from 0 to 16, y from -2 to 15
\put{1.274)} [l]  at 2 12
\put {$ \scriptstyle \bullet$} [c] at 10 0
\put {$ \scriptstyle \bullet$} [c] at 10 3
\put {$ \scriptstyle \bullet$} [c] at 10 6
\put {$ \scriptstyle \bullet$} [c] at 12 0
\put {$ \scriptstyle \bullet$} [c] at 12 12
\put {$ \scriptstyle \bullet$} [c] at 14 6
\put{$\scriptstyle \bullet$} [c] at 16  0
\setlinear \plot 10 0 10 6 12 0 14  6 12  12 10 6  /
\put{$5{,}040$} [c] at 13 -2
\endpicture
\end{minipage}
\begin{minipage}{4cm}
\beginpicture
\setcoordinatesystem units   <1.5mm,2mm>
\setplotarea x from 0 to 16, y from -2 to 15
\put{1.275)} [l] at 2 12
\put {$ \scriptstyle \bullet$} [c] at 10 0
\put {$ \scriptstyle \bullet$} [c] at 10 4
\put {$ \scriptstyle \bullet$} [c] at 10 8
\put {$ \scriptstyle \bullet$} [c] at 10 12
\put {$ \scriptstyle \bullet$} [c] at 14 0
\put {$ \scriptstyle \bullet$} [c] at 14 12
\setlinear \plot 10 12 10 0  /
\setlinear \plot 10 8 14 12 14 0  /
\put{$5{,}040$} [c] at 13 -2
\put{$\scriptstyle \bullet$} [c] at 16  0 \endpicture
\end{minipage}
\begin{minipage}{4cm}
\beginpicture
\setcoordinatesystem units   <1.5mm,2mm>
\setplotarea x from 0 to 16, y from -2 to 15
\put{1.276)} [l] at 2 12
\put {$ \scriptstyle \bullet$} [c] at 10 0
\put {$ \scriptstyle \bullet$} [c] at 10 4
\put {$ \scriptstyle \bullet$} [c] at 10 8
\put {$ \scriptstyle \bullet$} [c] at 10 12
\put {$ \scriptstyle \bullet$} [c] at 14 0
\put {$ \scriptstyle \bullet$} [c] at 14 12
\setlinear \plot 10 12 10 0  /
\setlinear \plot 10 4 14 0 14 12  /
\put{$5{,}040$} [c] at 13 -2
\put{$\scriptstyle \bullet$} [c] at 16  0 \endpicture
\end{minipage}
\begin{minipage}{4cm}
\beginpicture
\setcoordinatesystem units   <1.5mm,2mm>
\setplotarea x from 0 to 16, y from -2 to 15
\put{1.277)} [l] at 2 12
\put {$ \scriptstyle \bullet$} [c] at 10 0
\put {$ \scriptstyle \bullet$} [c] at 10 4
\put {$ \scriptstyle \bullet$} [c] at 10 8
\put {$ \scriptstyle \bullet$} [c] at 10 12
\put {$ \scriptstyle \bullet$} [c] at 14 0
\put {$ \scriptstyle \bullet$} [c] at 14 12
\setlinear \plot 10 12 10 0  14 12 14 0 10 8 /
\put{$5{,}040$} [c] at 13 -2
\put{$\scriptstyle \bullet$} [c] at 16  0 \endpicture
\end{minipage}
\begin{minipage}{4cm}
\beginpicture
\setcoordinatesystem units   <1.5mm,2mm>
\setplotarea x from 0 to 16, y from -2 to 15
\put{1.278)} [l] at 2 12
\put {$ \scriptstyle \bullet$} [c] at 10 0
\put {$ \scriptstyle \bullet$} [c] at 10 4
\put {$ \scriptstyle \bullet$} [c] at 10 8
\put {$ \scriptstyle \bullet$} [c] at 10 12
\put {$ \scriptstyle \bullet$} [c] at 14 0
\put {$ \scriptstyle \bullet$} [c] at 14 12
\setlinear \plot 10 0 10 12  14 0 14 12 10 4 /
\put{$5{,}040$} [c] at 13 -2
\put{$\scriptstyle \bullet$} [c] at 16  0 \endpicture
\end{minipage}
$$

$$
\begin{minipage}{4cm}
\beginpicture
\setcoordinatesystem units   <1.5mm,2mm>
\setplotarea x from 0 to 16, y from -2 to 15
\put{1.279)} [l] at 2 12
\put {$ \scriptstyle \bullet$} [c] at 10 12
\put {$ \scriptstyle \bullet$} [c] at 11 9
\put {$ \scriptstyle \bullet$} [c] at 12 6
\put {$ \scriptstyle \bullet$} [c] at 12 0
\put {$ \scriptstyle \bullet$} [c] at 14 12
\setlinear \plot 10 12 12 6  12 0   /
\setlinear \plot 12 6  14 12   /
\put {$ \scriptstyle \bullet$} [c] at 16 0
\put {$ \scriptstyle \bullet$} [c] at 16 12
\setlinear \plot 16 0 16 12   /
\put{$5{,}040$} [c] at 13 -2
\endpicture
\end{minipage}
\begin{minipage}{4cm}
\beginpicture
\setcoordinatesystem units   <1.5mm,2mm>
\setplotarea x from 0 to 16, y from -2 to 15
\put{1.280)} [l] at 2 12
\put {$ \scriptstyle \bullet$} [c] at 10 0
\put {$ \scriptstyle \bullet$} [c] at 11 3
\put {$ \scriptstyle \bullet$} [c] at 12 6
\put {$ \scriptstyle \bullet$} [c] at 12 12
\put {$ \scriptstyle \bullet$} [c] at 14 0
\setlinear \plot 10 0 12 6  12 12   /
\setlinear \plot 12 6  14 0   /
\put {$ \scriptstyle \bullet$} [c] at 16 0
\put {$ \scriptstyle \bullet$} [c] at 16 12
\setlinear \plot 16 0 16 12   /
\put{$5{,}040$} [c] at 13 -2
\endpicture
\end{minipage}
\begin{minipage}{4cm}
\beginpicture
\setcoordinatesystem units   <1.5mm,2mm>
\setplotarea x from 0 to 16, y from -2 to 15
\put{1.281)} [l] at 2 12
\put {$ \scriptstyle \bullet$} [c] at 10 6
\put {$ \scriptstyle \bullet$} [c] at 11.5 0
\put {$ \scriptstyle \bullet$} [c] at 11.5 12
\put {$ \scriptstyle \bullet$} [c] at 13 6
\put {$ \scriptstyle \bullet$} [c] at 14 0
\put {$ \scriptstyle \bullet$} [c] at 14 12
\setlinear \plot 11.5 0 10 6 11.5 12 13 6 11.5 0 /
\setlinear \plot 14 0 13 6  14 12 /
\put{$5{,}040$} [c] at 13 -2
\put{$\scriptstyle \bullet$} [c] at 16  0
\endpicture
\end{minipage}
\begin{minipage}{4cm}
\beginpicture
\setcoordinatesystem units   <1.5mm,2mm>
\setplotarea x from 0 to 16, y from -2 to 15
\put{1.282)} [l] at 2 12
\put {$ \scriptstyle \bullet$} [c] at 10 0
\put {$ \scriptstyle \bullet$} [c] at 10 9
\put {$ \scriptstyle \bullet$} [c] at 10 12
\put {$ \scriptstyle \bullet$} [c] at 14 0
\put {$ \scriptstyle \bullet$} [c] at 14 3
\put {$ \scriptstyle \bullet$} [c] at 14 12
\put{$\scriptstyle \bullet$} [c] at 16  0
\setlinear \plot 10 12 10 0   /
\setlinear \plot 14 12 14 0 /
\setlinear \plot 10 9 14 3 /
\put{$5{,}040$} [c] at 13 -2
 \endpicture
\end{minipage}
\begin{minipage}{4cm}
\beginpicture
\setcoordinatesystem units   <1.5mm,2mm>
\setplotarea x from 0 to 16, y from -2 to 15
\put{1.283)} [l] at 2 12
\put {$ \scriptstyle \bullet$} [c] at 10 0
\put {$ \scriptstyle \bullet$} [c] at 10 6
\put {$ \scriptstyle \bullet$} [c] at 10 12
\put {$ \scriptstyle \bullet$} [c] at 14 0
\put {$ \scriptstyle \bullet$} [c] at 14 6
\put {$ \scriptstyle \bullet$} [c] at 14 12
\setlinear \plot 10 6 10 0 14 12 14 0 10 6 10 12 14 6 /
\put{$5{,}040$} [c] at 13 -2
\put{$\scriptstyle \bullet$} [c] at 16  0 \endpicture
\end{minipage}
\begin{minipage}{4cm}
\beginpicture
\setcoordinatesystem units   <1.5mm,2mm>
\setplotarea x from 0 to 16, y from -2 to 15
\put{1.284)} [l] at 2 12
\put {$ \scriptstyle \bullet$} [c] at 10 3
\put {$ \scriptstyle \bullet$} [c] at 10 9
\put {$ \scriptstyle \bullet$} [c] at 12 0
\put {$ \scriptstyle \bullet$} [c] at 12 12
\put {$ \scriptstyle \bullet$} [c] at 14 6
\setlinear \plot 12 0 10 3 10 9 12 12 14 6  12 0 /
\put {$ \scriptstyle \bullet$} [c] at 16 0
\put {$ \scriptstyle \bullet$} [c] at 16 12
\setlinear \plot 16 0 16 12   /

\put{$5{,}040$} [c] at 13 -2
\endpicture
\end{minipage}
$$
$$
\begin{minipage}{4cm}
\beginpicture
\setcoordinatesystem units   <1.5mm,2mm>
\setplotarea x from 0 to 16, y from -2 to 15
\put{1.285)} [l] at 2 12
\put {$ \scriptstyle \bullet$} [c] at 10 6
\put {$ \scriptstyle \bullet$} [c] at 11 0
\put {$ \scriptstyle \bullet$} [c] at 11 3
\put {$ \scriptstyle \bullet$} [c] at 11 12
\put {$ \scriptstyle \bullet$} [c] at 12 6
\put {$ \scriptstyle \bullet$} [c] at 14 6
\setlinear \plot 10 6 11 3 12 6 11 12  10 6 /
\setlinear \plot  11 12 14 6 11 0 11 3  /
\put{$2{,}520$} [c] at 13 -2
\put{$\scriptstyle \bullet$} [c] at 16  0 \endpicture
\end{minipage}
\begin{minipage}{4cm}
\beginpicture
\setcoordinatesystem units   <1.5mm,2mm>
\setplotarea x from 0 to 16, y from -2 to 15
\put{1.286)} [l] at 2 12
\put {$ \scriptstyle \bullet$} [c] at 10 6
\put {$ \scriptstyle \bullet$} [c] at 11 0
\put {$ \scriptstyle \bullet$} [c] at 11 9
\put {$ \scriptstyle \bullet$} [c] at 11 12
\put {$ \scriptstyle \bullet$} [c] at 12 6
\put {$ \scriptstyle \bullet$} [c] at 14 6
\setlinear \plot 10 6 11 9 12 6 11 0  10 6 /
\setlinear \plot  11 0 14 6 11 12 11 9  /
\put{$2{,}520$} [c] at 13 -2
\put{$\scriptstyle \bullet$} [c] at 16  0 \endpicture
\end{minipage}
\begin{minipage}{4cm}
\beginpicture
\setcoordinatesystem units   <1.5mm,2mm>
\setplotarea x from 0 to 16, y from -2 to 15
\put{1.287)} [l] at 2 12
\put {$ \scriptstyle \bullet$} [c] at 10 9
\put {$ \scriptstyle \bullet$} [c] at 11 0
\put {$ \scriptstyle \bullet$} [c] at 11 3
\put {$ \scriptstyle \bullet$} [c] at 11 12
\put {$ \scriptstyle \bullet$} [c] at 12 9
\put {$ \scriptstyle \bullet$} [c] at 14 12
\setlinear \plot 10 9 11 12 12 9 11 3 11 0 /
\setlinear \plot  10 9 11 3 14 12  /
\put{$2{,}520$} [c] at 13 -2
\put{$\scriptstyle \bullet$} [c] at 16  0 \endpicture
\end{minipage}
\begin{minipage}{4cm}
\beginpicture
\setcoordinatesystem units   <1.5mm,2mm>
\setplotarea x from 0 to 16, y from -2 to 15
\put{1.288)} [l] at 2 12
\put {$ \scriptstyle \bullet$} [c] at 10 3
\put {$ \scriptstyle \bullet$} [c] at 11 0
\put {$ \scriptstyle \bullet$} [c] at 11 9
\put {$ \scriptstyle \bullet$} [c] at 11 12
\put {$ \scriptstyle \bullet$} [c] at 12 3
\put {$ \scriptstyle \bullet$} [c] at 14 0
\setlinear \plot 10 3 11 0 12 3 11 9 11 12 /
\setlinear \plot  10 3 11 9 14 0  /
\put{$2{,}520$} [c] at 13 -2
\put{$\scriptstyle \bullet$} [c] at 16  0 \endpicture
\end{minipage}
\begin{minipage}{4cm}
\beginpicture
\setcoordinatesystem units   <1.5mm,2mm>
\setplotarea x from 0 to 16, y from -2 to 15
\put{1.289)} [l] at 2 12
\put {$ \scriptstyle \bullet$} [c] at 10 12
\put {$ \scriptstyle \bullet$} [c] at 12 0
\put {$ \scriptstyle \bullet$} [c] at 12 4
\put {$ \scriptstyle \bullet$} [c] at 12 8
\put {$ \scriptstyle \bullet$} [c] at 12 12
\put {$ \scriptstyle \bullet$} [c] at 14 12
\setlinear \plot 12 0 12 12  /
\setlinear \plot 10 12 12 4 /
\setlinear \plot 14 12 12 8 /
\put{$2{,}520$} [c] at 13 -2
\put{$\scriptstyle \bullet$} [c] at 16  0 \endpicture
\end{minipage}
\begin{minipage}{4cm}
\beginpicture
\setcoordinatesystem units   <1.5mm,2mm>
\setplotarea x from 0 to 16, y from -2 to 15
\put{1.290)} [l] at 2 12
\put {$ \scriptstyle \bullet$} [c] at 10 0
\put {$ \scriptstyle \bullet$} [c] at 12 0
\put {$ \scriptstyle \bullet$} [c] at 12 4
\put {$ \scriptstyle \bullet$} [c] at 12 8
\put {$ \scriptstyle \bullet$} [c] at 12 12
\put {$ \scriptstyle \bullet$} [c] at 14 0
\setlinear \plot 12 0 12 12  /
\setlinear \plot 10 0 12 8 /
\setlinear \plot 14 0 12 4 /
\put{$2{,}520$} [c] at 13 -2
\put{$\scriptstyle \bullet$} [c] at 16  0
\endpicture
\end{minipage}
$$
$$
\begin{minipage}{4cm}
\beginpicture
\setcoordinatesystem units   <1.5mm,2mm>
\setplotarea x from 0 to 16, y from -2 to 15
\put{1.291)} [l] at 2 12
\put {$ \scriptstyle \bullet$} [c] at 10 0
\put {$ \scriptstyle \bullet$} [c] at 12 6
\put {$ \scriptstyle \bullet$} [c] at 12 12
\put {$ \scriptstyle \bullet$} [c] at 14 0
\setlinear \plot 10 0 12 6 12  12   /
\setlinear \plot 14 0 12 6  /
\put {$ \scriptstyle \bullet$} [c] at 16 0
\put {$ \scriptstyle \bullet$} [c] at 16 6
\put {$ \scriptstyle \bullet$} [c] at 16 12
\setlinear \plot 16 0 16 12   /
\put{$2{,}520$} [c] at 13 -2
\endpicture
\end{minipage}
\begin{minipage}{4cm}
\beginpicture
\setcoordinatesystem units   <1.5mm,2mm>
\setplotarea x from 0 to 16, y from -2 to 15
\put{1.292)} [l] at 2 12
\put {$ \scriptstyle \bullet$} [c] at 10 12
\put {$ \scriptstyle \bullet$} [c] at 12 0
\put {$ \scriptstyle \bullet$} [c] at 12 6
\put {$ \scriptstyle \bullet$} [c] at 14 12
\setlinear \plot 12 0 12 6 10  12    /
\setlinear \plot 12 6 14  12    /
\put {$ \scriptstyle \bullet$} [c] at 16 0
\put {$ \scriptstyle \bullet$} [c] at 16 6
\put {$ \scriptstyle \bullet$} [c] at 16 12
\setlinear \plot 16 0 16 12   /
\put{$2{,}520$} [c] at 13 -2
\endpicture
\end{minipage}
\begin{minipage}{4cm}
\beginpicture
\setcoordinatesystem units   <1.5mm,2mm>
\setplotarea x from 0 to 16, y from -2 to 15
\put{1.293)} [l] at 2 12
\put {$ \scriptstyle \bullet$} [c] at 10 0
\put {$ \scriptstyle \bullet$} [c] at 10 3
\put {$ \scriptstyle \bullet$} [c] at 10 6
\put {$ \scriptstyle \bullet$} [c] at 10 9
\put {$ \scriptstyle \bullet$} [c] at 10 12
\put {$ \scriptstyle \bullet$} [c] at 13 0
\put{$\scriptstyle \bullet$} [c] at 16  0
\setlinear \plot 10 12 10 0  /
\put{$2{,}520$} [c] at 13 -2
\endpicture
\end{minipage}
\begin{minipage}{4cm}
\beginpicture
\setcoordinatesystem units   <1.5mm,2mm>
\setplotarea x from 0 to 16, y from -2 to 15
\put{1.294)} [l] at 2 12
\put {$ \scriptstyle \bullet$} [c] at 10 6
\put {$ \scriptstyle \bullet$} [c] at 12 0
\put {$ \scriptstyle \bullet$} [c] at 12 12
\put {$ \scriptstyle \bullet$} [c] at 14 6
\setlinear \plot 12 0 10 6 12  12  14 6 12 0   /
\put {$ \scriptstyle \bullet$} [c] at 16 0
\put {$ \scriptstyle \bullet$} [c] at 16 6
\put {$ \scriptstyle \bullet$} [c] at 16 12
\setlinear \plot 16 0 16 12   /
\put{$2{,}520$} [c] at 13 -2
\endpicture
\end{minipage}
\begin{minipage}{4cm}
\beginpicture
\setcoordinatesystem units   <1.5mm,2mm>
\setplotarea x from 0 to 16, y from -2 to 15
\put{1.295)} [l] at 2 12
\put {$ \scriptstyle \bullet$} [c] at 10 0
\put {$ \scriptstyle \bullet$} [c] at 10 6
\put {$ \scriptstyle \bullet$} [c] at 10 12
\put {$ \scriptstyle \bullet$} [c] at 14 0
\put {$ \scriptstyle \bullet$} [c] at 14 6
\put {$ \scriptstyle \bullet$} [c] at 14 12
\setlinear \plot 10 0 10 12 14 6 14 12 10  6  /
\setlinear \plot 14 0 14 6 /
\put{$1{,}260$} [c] at 13 -2
\put{$\scriptstyle \bullet$} [c] at 16  0 \endpicture
\end{minipage}
\begin{minipage}{4cm}
\beginpicture
\setcoordinatesystem units   <1.5mm,2mm>
\setplotarea x from 0 to 16, y from -2 to 15
\put{1.296)} [l] at 2 12
\put {$ \scriptstyle \bullet$} [c] at 10 0
\put {$ \scriptstyle \bullet$} [c] at 10 6
\put {$ \scriptstyle \bullet$} [c] at 10 12
\put {$ \scriptstyle \bullet$} [c] at 14 0
\put {$ \scriptstyle \bullet$} [c] at 14 6
\put {$ \scriptstyle \bullet$} [c] at 14 12
\setlinear \plot 10 12 10 0 14 6 14 0 10  6  /
\setlinear \plot 14 12 14 6 /
\put{$1{,}260$} [c] at 13 -2
\put{$\scriptstyle \bullet$} [c] at 16  0 \endpicture
\end{minipage}
$$
$$
\begin{minipage}{4cm}
\beginpicture
\setcoordinatesystem units   <1.5mm,2mm>
\setplotarea x from 0 to 16, y from -2 to 15
\put{1.297)} [l] at 2 12
\put {$ \scriptstyle \bullet$} [c] at 10 6
\put {$ \scriptstyle \bullet$} [c] at 10 12
\put {$ \scriptstyle \bullet$} [c] at 12   0
\put {$ \scriptstyle \bullet$} [c] at 14 6
\put {$ \scriptstyle \bullet$} [c] at 14 12
\setlinear \plot 10 6  12 0  14 6 14 12 10 6 10 12 14 6   /
\put {$ \scriptstyle \bullet$} [c] at 16 0
\put {$ \scriptstyle \bullet$} [c] at 16 12
\setlinear \plot 16 0 16 12   /
\put{$1{,}260$} [c] at 13 -2
\endpicture
\end{minipage}
\begin{minipage}{4cm}
\beginpicture
\setcoordinatesystem units   <1.5mm,2mm>
\setplotarea x from 0 to 16, y from -2 to 15
\put{1.298)} [l] at 2 12
\put {$ \scriptstyle \bullet$} [c] at 10 6
\put {$ \scriptstyle \bullet$} [c] at 10 0
\put {$ \scriptstyle \bullet$} [c] at 12  12
\put {$ \scriptstyle \bullet$} [c] at 14 6
\put {$ \scriptstyle \bullet$} [c] at 14 0
\setlinear \plot 10 6  12 12  14 6 14 0 10 6 10 0 14 6   /
\put {$ \scriptstyle \bullet$} [c] at 16 0
\put {$ \scriptstyle \bullet$} [c] at 16 12
\setlinear \plot 16 0 16 12   /
\put{$1{,}260$} [c] at 13 -2
\endpicture
\end{minipage}
\begin{minipage}{4cm}
\beginpicture
\setcoordinatesystem units   <1.5mm,2mm>
\setplotarea x from 0 to 16, y from -2 to 15
\put{1.299)} [l] at 2 12
\put {$ \scriptstyle \bullet$} [c] at 10 0
\put {$ \scriptstyle \bullet$} [c] at 10 12
\put {$ \scriptstyle \bullet$} [c] at 12 6
\put {$ \scriptstyle \bullet$} [c] at 14 0
\put {$ \scriptstyle \bullet$} [c] at 14 12
\setlinear \plot 10 0 14 12   /
\setlinear \plot 10 12  14 0  /
\put {$ \scriptstyle \bullet$} [c] at 16 0
\put {$ \scriptstyle \bullet$} [c] at 16 12
\setlinear \plot 16 0 16 12   /

\put{$1{,}260$} [c] at 13 -2
\endpicture
\end{minipage}
\begin{minipage}{4cm}
\beginpicture
\setcoordinatesystem units   <1.5mm,2mm>
\setplotarea x from 0 to 16, y from -2 to 15
\put{1.300)} [l] at 2 12
\put {$ \scriptstyle \bullet$} [c] at 10 6
\put {$ \scriptstyle \bullet$} [c] at  10 12
\put {$ \scriptstyle \bullet$} [c] at 12 0
\put {$ \scriptstyle \bullet$} [c] at 12 12
\put {$ \scriptstyle \bullet$} [c] at 14 6
\put {$ \scriptstyle \bullet$} [c] at 14 12
\setlinear \plot 10 12 10 6  12 0 14  6 14  12 10 6  12 12 14 6 10 12 /
\put{$420$} [c] at 13 -2
\put{$\scriptstyle \bullet$} [c] at 16  0 \endpicture
\end{minipage}
\begin{minipage}{4cm}
\beginpicture
\setcoordinatesystem units   <1.5mm,2mm>
\setplotarea x from 0 to 16, y from -2 to 15
\put{1.301)} [l] at 2 12
\put {$ \scriptstyle \bullet$} [c] at 10 6
\put {$ \scriptstyle \bullet$} [c] at  10 0
\put {$ \scriptstyle \bullet$} [c] at 12 0
\put {$ \scriptstyle \bullet$} [c] at 12 12
\put {$ \scriptstyle \bullet$} [c] at 14 6
\put {$ \scriptstyle \bullet$} [c] at 14 0
\setlinear \plot 10 0 10 6  12 12 14  6 14  0 10 6  12 0 14 6 10 0 /
\put{$420$} [c] at 13 -2
\put{$\scriptstyle \bullet$} [c] at 16  0 \endpicture
\end{minipage}
\begin{minipage}{4cm}
\beginpicture
\setcoordinatesystem units   <1.5mm,2mm>
\setplotarea x from 0 to 16, y from -2 to 15
\put{1.302)} [l] at 2 12
\put {$ \scriptstyle \bullet$} [c] at 10 12
\put {$ \scriptstyle \bullet$} [c] at 10 0
\put {$ \scriptstyle \bullet$} [c] at 12 12
\put {$ \scriptstyle \bullet$} [c] at 12 6
\put {$ \scriptstyle \bullet$} [c] at 14 0
\put {$ \scriptstyle \bullet$} [c] at 14 12
\setlinear \plot 10 12 14 0  /
\setlinear \plot 10 0 14 12   /
\setlinear \plot 12 6 12 12   /
\put{$420$} [c] at 13 -2
\put{$\scriptstyle \bullet$} [c] at 16  0 \endpicture
\end{minipage}
$$

$$
\begin{minipage}{4cm}
\beginpicture
\setcoordinatesystem units   <1.5mm,2mm>
\setplotarea x from 0 to 16, y from -2 to 15
\put{1.303)} [l] at 2 12
\put {$ \scriptstyle \bullet$} [c] at 10 12
\put {$ \scriptstyle \bullet$} [c] at 10 0
\put {$ \scriptstyle \bullet$} [c] at 12 0
\put {$ \scriptstyle \bullet$} [c] at 12 6
\put {$ \scriptstyle \bullet$} [c] at 14 0
\put {$ \scriptstyle \bullet$} [c] at 14 12
\setlinear \plot 10 12 14 0  /
\setlinear \plot 10 0 14 12   /
\setlinear \plot 12 6 12 0   /
\put{$420$} [c] at 13 -2
\put{$\scriptstyle \bullet$} [c] at 16  0 \endpicture
\end{minipage}
\begin{minipage}{4cm}
\beginpicture
\setcoordinatesystem units   <1.5mm,2mm>
\setplotarea x from 0 to 16, y from -2 to 15
\put{1.304)} [l] at 2 12
\put {$ \scriptstyle \bullet$} [c] at 10 6
\put {$ \scriptstyle \bullet$} [c] at 10 12
\put {$ \scriptstyle \bullet$} [c] at 11 0
\put {$ \scriptstyle \bullet$} [c] at 12 12
\put {$ \scriptstyle \bullet$} [c] at 12 6
\put {$ \scriptstyle \bullet$} [c] at 14 6
\setlinear \plot 10 12 10 6 11 0 12 6 12  12 14 6 10 12 12 6  /
\setlinear \plot  10 6 12 12  /
\setlinear \plot  11 0 14  6   /
\put{$420$} [c] at 13 -2
\put{$\scriptstyle \bullet$} [c] at 16  0 \endpicture
\end{minipage}
\begin{minipage}{4cm}
\beginpicture
\setcoordinatesystem units   <1.5mm,2mm>
\setplotarea x from 0 to 16, y from -2 to 15
\put{1.305)} [l] at 2 12
\put {$ \scriptstyle \bullet$} [c] at 10 6
\put {$ \scriptstyle \bullet$} [c] at 10 0
\put {$ \scriptstyle \bullet$} [c] at 11 12
\put {$ \scriptstyle \bullet$} [c] at 12 0
\put {$ \scriptstyle \bullet$} [c] at 12 6
\put {$ \scriptstyle \bullet$} [c] at 14 6
\setlinear \plot 10 0 10 6 11 12 12 6 12  0 14 6 10 0 12 6  /
\setlinear \plot  10 6 12 0  /
\setlinear \plot  11 12 14  6   /
\put{$420$} [c] at 13 -2
\put{$\scriptstyle \bullet$} [c] at 16  0 \endpicture
\end{minipage}
\begin{minipage}{4cm}
\beginpicture
\setcoordinatesystem units    <1.5mm,2mm>
\setplotarea x from 0 to 16, y from -2 to 15
\put{${\bf  25}$} [l]  at 2 15

\put{1.306)} [l]  at 2 12
\put {$ \scriptstyle \bullet$} [c] at 10 12
\put {$ \scriptstyle \bullet$} [c] at 10 6
\put {$ \scriptstyle \bullet$} [c] at 11.5 9
\put {$ \scriptstyle \bullet$} [c] at 13 0
\put {$ \scriptstyle \bullet$} [c] at 13 6
\put {$ \scriptstyle \bullet$} [c] at 13 12
\put {$ \scriptstyle \bullet$} [c] at 16 12
\setlinear \plot 16 12 13  0 13 12 10 6 13 0    /
\setlinear \plot 10 6 10  12 /
\put{$5{,}040$} [c]  at 13 -2
\endpicture
\end{minipage}
\begin{minipage}{4cm}
\beginpicture
\setcoordinatesystem units    <1.5mm,2mm>
\setplotarea x from 0 to 16, y from -2 to 15
\put{1.307)} [l]  at 2 12
\put {$ \scriptstyle \bullet$} [c] at 10 0
\put {$ \scriptstyle \bullet$} [c] at 10 6
\put {$ \scriptstyle \bullet$} [c] at 11.5 3
\put {$ \scriptstyle \bullet$} [c] at 13 0
\put {$ \scriptstyle \bullet$} [c] at 13 6
\put {$ \scriptstyle \bullet$} [c] at 13 12
\put {$ \scriptstyle \bullet$} [c] at 16 0
\setlinear \plot 16 0 13  12 13 0 10 6 13 12    /
\setlinear \plot 10 6 10  0 /
\put{$5{,}040$} [c]  at 13 -2
\endpicture
\end{minipage}
\begin{minipage}{4cm}
\beginpicture
\setcoordinatesystem units    <1.5mm,2mm>
\setplotarea x from 0 to 16, y from -2 to 15
\put{1.308)} [l]  at 2 12
\put {$ \scriptstyle \bullet$} [c] at 10 12
\put {$ \scriptstyle \bullet$} [c] at 10 6
\put {$ \scriptstyle \bullet$} [c] at 10 9
\put {$ \scriptstyle \bullet$} [c] at 13 0
\put {$ \scriptstyle \bullet$} [c] at 13 6
\put {$ \scriptstyle \bullet$} [c] at 13 12
\put {$ \scriptstyle \bullet$} [c] at 16 12
\setlinear \plot 10 12 10 6 13 0 16 12     /
\setlinear \plot 13 0 13 12   /
\put{$5{,}040$} [c]  at 13 -2
\endpicture
\end{minipage}
$$
$$
\begin{minipage}{4cm}
\beginpicture
\setcoordinatesystem units    <1.5mm,2mm>
\setplotarea x from 0 to 16, y from -2 to 15
\put{1.309)} [l]  at 2 12
\put {$ \scriptstyle \bullet$} [c] at 10 0
\put {$ \scriptstyle \bullet$} [c] at 10 6
\put {$ \scriptstyle \bullet$} [c] at 10 3
\put {$ \scriptstyle \bullet$} [c] at 13 0
\put {$ \scriptstyle \bullet$} [c] at 13 6
\put {$ \scriptstyle \bullet$} [c] at 13 12
\put {$ \scriptstyle \bullet$} [c] at 16 0
\setlinear \plot 10 0 10 6 13 12 16 0     /
\setlinear \plot 13 0 13 12   /
\put{$5{,}040$} [c]  at 13 -2
\endpicture
\end{minipage}
\begin{minipage}{4cm}
\beginpicture
\setcoordinatesystem units    <1.5mm,2mm>
\setplotarea x from 0 to 16, y from -2 to 15
\put{1.310)} [l]  at 2 12
\put {$ \scriptstyle \bullet$} [c] at 10 6
\put {$ \scriptstyle \bullet$} [c] at 10  12
\put {$ \scriptstyle \bullet$} [c] at 13 0
\put {$ \scriptstyle \bullet$} [c] at 13 6
\put {$ \scriptstyle \bullet$} [c] at 13 12
\put {$ \scriptstyle \bullet$} [c] at 16 6
\put {$ \scriptstyle \bullet$} [c] at 16 12
\setlinear \plot 16 12 16 6 13 0 13 12 10 6 13  0    /
\setlinear \plot 10  6 10 12      /
\put{$5{,}040$} [c]  at 13 -2
\endpicture
\end{minipage}
\begin{minipage}{4cm}
\beginpicture
\setcoordinatesystem units    <1.5mm,2mm>
\setplotarea x from 0 to 16, y from -2 to 15
\put{1.311)} [l]  at 2 12
\put {$ \scriptstyle \bullet$} [c] at 10 6
\put {$ \scriptstyle \bullet$} [c] at 10  0
\put {$ \scriptstyle \bullet$} [c] at 13 0
\put {$ \scriptstyle \bullet$} [c] at 13 6
\put {$ \scriptstyle \bullet$} [c] at 13 12
\put {$ \scriptstyle \bullet$} [c] at 16 6
\put {$ \scriptstyle \bullet$} [c] at 16 0
\setlinear \plot 16 0 16 6 13 12 13 0 10 6 13  12    /
\setlinear \plot 10  6 10 0      /
\put{$5{,}040$} [c]  at 13 -2
\endpicture
\end{minipage}
\begin{minipage}{4cm}
\beginpicture
\setcoordinatesystem units    <1.5mm,2mm>
\setplotarea x  from 0 to 16, y  from -2 to 15
\put{1.312)} [l]  at 2 12
\put {$ \scriptstyle \bullet$} [c] at  10 0
\put {$ \scriptstyle \bullet$} [c] at  12 4
\put {$ \scriptstyle \bullet$} [c] at  12 12
\put {$ \scriptstyle \bullet$} [c] at  14 0
\put {$ \scriptstyle \bullet$} [c] at  14 8
\put {$ \scriptstyle \bullet$} [c] at  16 4
\put {$ \scriptstyle \bullet$} [c] at  16 12
\setlinear \plot 10 0 12 12 12 4 14 0 16 4 16 12 14 8 12 12  /
\setlinear \plot 14 0 14 8 /
\put{$5{,}040  $} [c] at 13 -2
\endpicture
\end{minipage}
\begin{minipage}{4cm}
\beginpicture
\setcoordinatesystem units    <1.5mm,2mm>
\setplotarea x  from 0 to 16, y  from -2 to 15
\put{1.313)} [l]  at 2 12
\put {$ \scriptstyle \bullet$} [c] at  10 12
\put {$ \scriptstyle \bullet$} [c] at  12 8
\put {$ \scriptstyle \bullet$} [c] at  12 0
\put {$ \scriptstyle \bullet$} [c] at  14 12
\put {$ \scriptstyle \bullet$} [c] at  14 4
\put {$ \scriptstyle \bullet$} [c] at  16 8
\put {$ \scriptstyle \bullet$} [c] at  16 0
\setlinear \plot 10 12 12 0 12 8 14 12 16 8 16 0 14 4 12 0  /
\setlinear \plot 14 12 14 4 /
\put{$5{,}040  $} [c] at 13 -2
\endpicture
\end{minipage}
\begin{minipage}{4cm}
\beginpicture
\setcoordinatesystem units    <1.5mm,2mm>
\setplotarea x  from 0 to 16, y from -2 to 15
\put{1.314)} [l] at 2 12
\put {$ \scriptstyle \bullet$} [c] at  10 12
\put {$ \scriptstyle \bullet$} [c] at  11.5 0
\put {$ \scriptstyle \bullet$} [c] at  11.5 12
\put {$ \scriptstyle \bullet$} [c] at  14 6
\put {$ \scriptstyle \bullet$} [c] at  15 0
\put {$ \scriptstyle \bullet$} [c] at  15 12
\put {$ \scriptstyle \bullet$} [c] at  16 6
\setlinear \plot  10 12 11.5 0 11.5 12 15 0 16 6 15 12 14 6 15 0  /
\setlinear \plot  14 6 11.5 0  /
\put{$5{,}040$} [c] at 13 -2
\endpicture
\end{minipage}
$$
$$
\begin{minipage}{4cm}
\beginpicture
\setcoordinatesystem units    <1.5mm,2mm>
\setplotarea x  from 0 to 16, y from -2 to 15
\put{1.315)} [l] at 2 12
\put {$ \scriptstyle \bullet$} [c] at  10 0
\put {$ \scriptstyle \bullet$} [c] at  11.5 0
\put {$ \scriptstyle \bullet$} [c] at  11.5 12
\put {$ \scriptstyle \bullet$} [c] at  14 6
\put {$ \scriptstyle \bullet$} [c] at  15 0
\put {$ \scriptstyle \bullet$} [c] at  15 12
\put {$ \scriptstyle \bullet$} [c] at  16 6
\setlinear \plot  10 0 11.5 12 11.5 0 15 12 16 6 15 0 14 6 15 12  /
\setlinear \plot  14 6 11.5 12  /
\put{$5{,}040$} [c] at 13 -2
\endpicture
\end{minipage}
\begin{minipage}{4cm}
\beginpicture
\setcoordinatesystem units    <1.5mm,2mm>
\setplotarea x  from 0 to 16, y from -2 to 15
\put{1.316)} [l] at 2 12
\put {$ \scriptstyle \bullet$} [c] at  10 12
\put {$ \scriptstyle \bullet$} [c] at  13 12
\put {$ \scriptstyle \bullet$} [c] at  13 0
\put {$ \scriptstyle \bullet$} [c] at  13 4
\put {$ \scriptstyle \bullet$} [c] at  13 8
\put {$ \scriptstyle \bullet$} [c] at  16 0
\put {$ \scriptstyle \bullet$} [c] at  16 12
\setlinear \plot  10 12  13 0 13 12 16 0 16 12 /
\put{$5{,}040$} [c] at 13 -2
\endpicture
\end{minipage}
\begin{minipage}{4cm}
\beginpicture
\setcoordinatesystem units    <1.5mm,2mm>
\setplotarea x  from 0 to 16, y from -2 to 15
\put{1.317)} [l] at 2 12
\put {$ \scriptstyle \bullet$} [c] at  10 0
\put {$ \scriptstyle \bullet$} [c] at  13 12
\put {$ \scriptstyle \bullet$} [c] at  13 0
\put {$ \scriptstyle \bullet$} [c] at  13 4
\put {$ \scriptstyle \bullet$} [c] at  13 8
\put {$ \scriptstyle \bullet$} [c] at  16 0
\put {$ \scriptstyle \bullet$} [c] at  16 12
\setlinear \plot  10 0  13 12 13 0 16 12 16 0 /
\put{$5{,}040$} [c] at 13 -2
\endpicture
\end{minipage}
\begin{minipage}{4cm}
\beginpicture
\setcoordinatesystem units    <1.5mm,2mm>
\setplotarea x  from 0 to 16, y from -2 to 15
\put{1.318)} [l] at 2 12
\put {$ \scriptstyle \bullet$} [c] at  10 12
\put {$ \scriptstyle \bullet$} [c] at  10 6
\put {$ \scriptstyle \bullet$} [c] at  11 0
\put {$ \scriptstyle \bullet$} [c] at  11 12
\put {$ \scriptstyle \bullet$} [c] at  12 6
\put {$ \scriptstyle \bullet$} [c] at  16 0
\put {$ \scriptstyle \bullet$} [c] at  16 12
\setlinear \plot  10 12 10 6 11 0 12 6 11 12 10 6 /
\setlinear \plot 11 12 16 0 16 12  /
\put{$5{,}040$} [c] at 13 -2
\endpicture
\end{minipage}
\begin{minipage}{4cm}
\beginpicture
\setcoordinatesystem units    <1.5mm,2mm>
\setplotarea x  from 0 to 16, y from -2 to 15
\put{1.319)} [l] at 2 12
\put {$ \scriptstyle \bullet$} [c] at  10 0
\put {$ \scriptstyle \bullet$} [c] at  10 6
\put {$ \scriptstyle \bullet$} [c] at  11 0
\put {$ \scriptstyle \bullet$} [c] at  11 12
\put {$ \scriptstyle \bullet$} [c] at  12 6
\put {$ \scriptstyle \bullet$} [c] at  16 0
\put {$ \scriptstyle \bullet$} [c] at  16 12
\setlinear \plot  10 0 10 6 11 12 12 6 11 0 10 6 /
\setlinear \plot 11 0 16 12 16 0  /
\put{$5{,}040$} [c] at 13 -2
\endpicture
\end{minipage}
\begin{minipage}{4cm}
\beginpicture
\setcoordinatesystem units    <1.5mm,2mm>
\setplotarea x  from 0 to 16, y from -2 to 15
\put{1.320)} [l] at 2 12
\put {$ \scriptstyle \bullet$} [c] at  10 12
\put {$ \scriptstyle \bullet$} [c] at  10 6
\put {$ \scriptstyle \bullet$} [c] at  10 0
\put {$ \scriptstyle \bullet$} [c] at  14 12
\put {$ \scriptstyle \bullet$} [c] at  14 6
\put {$ \scriptstyle \bullet$} [c] at  14 0
\put {$ \scriptstyle \bullet$} [c] at  16 12
\setlinear \plot  16 12  14 0 14  12 10 0 10 12  14 0 /
\put{$5{,}040$} [c] at 13 -2
\endpicture
\end{minipage}
$$
$$
\begin{minipage}{4cm}
\beginpicture
\setcoordinatesystem units    <1.5mm,2mm>
\setplotarea x  from 0 to 16, y from -2 to 15
\put{1.321)} [l] at 2 12
\put {$ \scriptstyle \bullet$} [c] at  10 12
\put {$ \scriptstyle \bullet$} [c] at  10 6
\put {$ \scriptstyle \bullet$} [c] at  10 0
\put {$ \scriptstyle \bullet$} [c] at  14 12
\put {$ \scriptstyle \bullet$} [c] at  14 6
\put {$ \scriptstyle \bullet$} [c] at  14 0
\put {$ \scriptstyle \bullet$} [c] at  16 0
\setlinear \plot  16 0  14 12 14  0 10 12 10 0  14 12 /
\put{$5{,}040$} [c] at 13 -2
\endpicture
\end{minipage}
\begin{minipage}{4cm}
\beginpicture
\setcoordinatesystem units    <1.5mm,2mm>
\setplotarea x  from 0 to 16, y from -2 to 15
\put{1.322)} [l] at 2 12
\put {$ \scriptstyle \bullet$} [c] at  10 0
\put {$ \scriptstyle \bullet$} [c] at  10 6
\put {$ \scriptstyle \bullet$} [c] at 10 12
\put {$ \scriptstyle \bullet$} [c] at  13 12
\put {$ \scriptstyle \bullet$} [c] at  13  6
\put {$ \scriptstyle \bullet$} [c] at  16 0
\put {$ \scriptstyle \bullet$} [c] at  16 12
\setlinear \plot  10 0 10 12 16 0 16 12   /
\setlinear \plot 13 6 13 12  /
\put{$5{,}040$} [c] at 13 -2
\endpicture
\end{minipage}
\begin{minipage}{4cm}
\beginpicture
\setcoordinatesystem units    <1.5mm,2mm>
\setplotarea x  from 0 to 16, y from -2 to 15
\put{1.323)} [l] at 2 12
\put {$ \scriptstyle \bullet$} [c] at  10 0
\put {$ \scriptstyle \bullet$} [c] at  10 6
\put {$ \scriptstyle \bullet$} [c] at 10 12
\put {$ \scriptstyle \bullet$} [c] at  13 0
\put {$ \scriptstyle \bullet$} [c] at  13  6
\put {$ \scriptstyle \bullet$} [c] at  16 0
\put {$ \scriptstyle \bullet$} [c] at  16 12
\setlinear \plot  10 12 10 0 16 12 16 0   /
\setlinear \plot 13 6 13 0  /
\put{$5{,}040$} [c] at 13 -2
\endpicture
\end{minipage}
\begin{minipage}{4cm}
\beginpicture
\setcoordinatesystem units    <1.5mm,2mm>
\setplotarea x  from 0 to 16, y  from -2 to 15
\put{1.324)} [l]  at 2 12
\put {$ \scriptstyle \bullet$} [c] at  10 0
\put {$ \scriptstyle \bullet$} [c] at  10 12
\put {$ \scriptstyle \bullet$} [c] at  13 12
\put {$ \scriptstyle \bullet$} [c] at  13 0
\put {$ \scriptstyle \bullet$} [c] at  13 6
\put {$ \scriptstyle \bullet$} [c] at  16 12
\put {$ \scriptstyle \bullet$} [c] at  16  0
\setlinear \plot  10  0 10 12 13 0 13 12 10 0   /
\setlinear \plot  16 0 16 12  13 6  /
\put{$5{,}040 $} [c]  at 13 -2

\endpicture
\end{minipage}
\begin{minipage}{4cm}
\beginpicture
\setcoordinatesystem units    <1.5mm,2mm>
\setplotarea x  from 0 to 16, y  from -2 to 15
\put{1.325)} [l]  at 2 12
\put {$ \scriptstyle \bullet$} [c] at  10 0
\put {$ \scriptstyle \bullet$} [c] at  10 12
\put {$ \scriptstyle \bullet$} [c] at  13 12
\put {$ \scriptstyle \bullet$} [c] at  13 0
\put {$ \scriptstyle \bullet$} [c] at  13 6
\put {$ \scriptstyle \bullet$} [c] at  16 12
\put {$ \scriptstyle \bullet$} [c] at  16  0
\setlinear \plot  10  12 10 0 13 12 13 0 10 12   /
\setlinear \plot  16 12 16 0  13 6  /
\put{$5{,}040 $} [c]  at 13 -2
\endpicture
\end{minipage}
\begin{minipage}{4cm}
\beginpicture
\setcoordinatesystem units    <1.5mm,2mm>
\setplotarea x  from 0 to 16, y  from -2 to 15
\put{1.326)} [l]  at 2 12
\put {$ \scriptstyle \bullet$} [c] at  10 0
\put {$ \scriptstyle \bullet$} [c] at  10 12
\put {$ \scriptstyle \bullet$} [c] at  13 0
\put {$ \scriptstyle \bullet$} [c] at  13 12
\put {$ \scriptstyle \bullet$} [c] at  16 6
\put {$ \scriptstyle \bullet$} [c] at  16 0
\put {$ \scriptstyle \bullet$} [c] at  16 12
\setlinear \plot  10 0 10  12  13 0 13 12  16 6 16 12 /
\setlinear \plot  10 12 16 0 16 6   /
\put{$5{,}040$} [c]  at 13 -2
\endpicture
\end{minipage}
$$
$$
\begin{minipage}{4cm}
\beginpicture
\setcoordinatesystem units    <1.5mm,2mm>
\setplotarea x  from 0 to 16, y  from -2 to 15
\put{1.327)} [l]  at 2 12
\put {$ \scriptstyle \bullet$} [c] at  10 0
\put {$ \scriptstyle \bullet$} [c] at  10 12
\put {$ \scriptstyle \bullet$} [c] at  13 0
\put {$ \scriptstyle \bullet$} [c] at  13 12
\put {$ \scriptstyle \bullet$} [c] at  16 6
\put {$ \scriptstyle \bullet$} [c] at  16 0
\put {$ \scriptstyle \bullet$} [c] at  16 12
\setlinear \plot  10 12 10  0  13 12 13 0  16 6 16 0 /
\setlinear \plot  10 0 16 12 16 6   /
\put{$5{,}040$} [c]  at 13 -2
\endpicture
\end{minipage}
\begin{minipage}{4cm}
\beginpicture
\setcoordinatesystem units    <1.5mm,2mm>
\setplotarea x  from 0 to 16, y  from -2 to 15
\put{1.328)} [l]  at 2 12
\put {$ \scriptstyle \bullet$} [c] at  10 0
\put {$ \scriptstyle \bullet$} [c] at  10 12
\put {$ \scriptstyle \bullet$} [c] at  13 0
\put {$ \scriptstyle \bullet$} [c] at  13 4
\put {$ \scriptstyle \bullet$} [c] at  13 12
\put {$ \scriptstyle \bullet$} [c] at  16 0
\put {$ \scriptstyle \bullet$} [c] at  16 12
\setlinear \plot  10 0 10  12  13 0 16 12  16 0 10 12 /
\setlinear \plot  10 0 13 12 13 0   /

\put{$5{,}040 $} [c]  at 13 -2
\endpicture
\end{minipage}
\begin{minipage}{4cm}
\beginpicture
\setcoordinatesystem units    <1.5mm,2mm>
\setplotarea x  from 0 to 16, y  from -2 to 15
\put{1.329)} [l]  at 2 12
\put {$ \scriptstyle \bullet$} [c] at  10 0
\put {$ \scriptstyle \bullet$} [c] at  10 12
\put {$ \scriptstyle \bullet$} [c] at  13 0
\put {$ \scriptstyle \bullet$} [c] at  13 8
\put {$ \scriptstyle \bullet$} [c] at  13 12
\put {$ \scriptstyle \bullet$} [c] at  16 0
\put {$ \scriptstyle \bullet$} [c] at  16 12
\setlinear \plot  10 12 10  0  13 12 16 0  16 12 10 0 /
\setlinear \plot  10 12 13 0 13 12   /
\put{$5{,}040 $} [c]  at 13 -2
\endpicture
\end{minipage}
\begin{minipage}{4cm}
\beginpicture
\setcoordinatesystem units    <1.5mm,2mm>
\setplotarea x  from 0 to 16, y  from -2 to 15
\put{1.330)} [l]  at 2 12
\put {$ \scriptstyle \bullet$} [c] at  10 0
\put {$ \scriptstyle \bullet$} [c] at  10 12
\put {$ \scriptstyle \bullet$} [c] at  13 0
\put {$ \scriptstyle \bullet$} [c] at  13 6
\put {$ \scriptstyle \bullet$} [c] at  13 12
\put {$ \scriptstyle \bullet$} [c] at  16  12
\put {$ \scriptstyle \bullet$} [c] at  16  0
\setlinear \plot  16 0 16  12 13  0 10 12  10 0 13 6 /
\setlinear \plot  13 0 13 12  /
\put{$5{,}040$} [c]  at 13 -2
\endpicture
\end{minipage}
\begin{minipage}{4cm}
\beginpicture
\setcoordinatesystem units    <1.5mm,2mm>
\setplotarea x  from 0 to 16, y  from -2 to 15
\put{1.331)} [l]  at 2 12
\put {$ \scriptstyle \bullet$} [c] at  10 0
\put {$ \scriptstyle \bullet$} [c] at  10 12
\put {$ \scriptstyle \bullet$} [c] at  13 0
\put {$ \scriptstyle \bullet$} [c] at  13 6
\put {$ \scriptstyle \bullet$} [c] at  13 12
\put {$ \scriptstyle \bullet$} [c] at  16  12
\put {$ \scriptstyle \bullet$} [c] at  16  0
\setlinear \plot  16 12 16  0 13 12 10 0  10 12 13 6 /
\setlinear \plot  13 0 13 12  /
\put{$5{,}040$} [c]  at 13 -2
\endpicture
\end{minipage}
\begin{minipage}{4cm}
\beginpicture
\setcoordinatesystem units    <1.5mm,2mm>
\setplotarea x from 0 to 16, y from -2 to 15
\put{1.332)} [l]  at 2 12
\put {$ \scriptstyle \bullet$} [c] at 10 6
\put {$ \scriptstyle \bullet$} [c] at 11 6
\put {$ \scriptstyle \bullet$} [c] at 11 12
\put {$ \scriptstyle \bullet$} [c] at 12 0
\put {$ \scriptstyle \bullet$} [c] at 13 12
\put {$ \scriptstyle \bullet$} [c] at 13 6
\put {$ \scriptstyle \bullet$} [c] at 16 12
\setlinear \plot 16 12 12 0 10 6  11 12 11 6 12 0 13 6 13 12 11 6    /
\setlinear \plot 11 12 13 6 /
\put{$2{,}520$} [c]  at 13 -2
\endpicture
\end{minipage}
$$
$$
\begin{minipage}{4cm}
\beginpicture
\setcoordinatesystem units    <1.5mm,2mm>
\setplotarea x from 0 to 16, y from -2 to 15
\put{1.333)} [l]  at 2 12
\put {$ \scriptstyle \bullet$} [c] at 10 6
\put {$ \scriptstyle \bullet$} [c] at 11 6
\put {$ \scriptstyle \bullet$} [c] at 11 0
\put {$ \scriptstyle \bullet$} [c] at 12 12
\put {$ \scriptstyle \bullet$} [c] at 13 0
\put {$ \scriptstyle \bullet$} [c] at 13 6
\put {$ \scriptstyle \bullet$} [c] at 16 0
\setlinear \plot 16 0 12 12 10 6  11 0 11 6 12 12 13 6 13 0 11 6    /
\setlinear \plot 11 0 13 6 /
\put{$2{,}520$} [c]  at 13 -2
\endpicture
\end{minipage}
\begin{minipage}{4cm}
\beginpicture
\setcoordinatesystem units    <1.5mm,2mm>
\setplotarea x from 0 to 16, y from -2 to 15
\put{1.334)} [l]  at 2 12
\put {$ \scriptstyle \bullet$} [c] at 10 12
\put {$ \scriptstyle \bullet$} [c] at 11 6
\put {$ \scriptstyle \bullet$} [c] at 11 12
\put {$ \scriptstyle \bullet$} [c] at 13 0
\put {$ \scriptstyle \bullet$} [c] at 13 12
\put {$ \scriptstyle \bullet$} [c] at 13 6
\put {$ \scriptstyle \bullet$} [c] at 16 12
\setlinear \plot 16 12 13 0 13 12  11 6  11 12 13  6  /
\setlinear \plot 10 12 11 6 13 0 /
\put{$2{,}520$} [c]  at 13 -2
\endpicture
\end{minipage}
\begin{minipage}{4cm}
\beginpicture
\setcoordinatesystem units    <1.5mm,2mm>
\setplotarea x from 0 to 16, y from -2 to 15
\put{1.335)} [l]  at 2 12
\put {$ \scriptstyle \bullet$} [c] at 10 0
\put {$ \scriptstyle \bullet$} [c] at 11 6
\put {$ \scriptstyle \bullet$} [c] at 11 0
\put {$ \scriptstyle \bullet$} [c] at 13 12
\put {$ \scriptstyle \bullet$} [c] at 13 0
\put {$ \scriptstyle \bullet$} [c] at 13 6
\put {$ \scriptstyle \bullet$} [c] at 16 0
\setlinear \plot 16 0 13 12 13 0  11 6  11 0 13  6  /
\setlinear \plot 10 0 11 6 13 12 /
\put{$2{,}520$} [c]  at 13 -2
\endpicture
\end{minipage}
\begin{minipage}{4cm}
\beginpicture
\setcoordinatesystem units    <1.5mm,2mm>
\setplotarea x  from 0 to 16, y  from -2 to 15
\put{1.336)} [l]  at 2 12
\put {$ \scriptstyle \bullet$} [c] at  10 0
\put {$ \scriptstyle \bullet$} [c] at  12 6
\put {$ \scriptstyle \bullet$} [c] at  12 12
\put {$ \scriptstyle \bullet$} [c] at  14 0
\put {$ \scriptstyle \bullet$} [c] at  14 6
\put {$ \scriptstyle \bullet$} [c] at  14 12
\put {$ \scriptstyle \bullet$} [c] at  16 6
\setlinear \plot 10 0 12 12 12 6 14 12 14 0 12 6   /
\setlinear \plot 14 0 16 6 14 12    /
\put{$2{,}520  $} [c] at 13 -2
\endpicture
\end{minipage}
\begin{minipage}{4cm}
\beginpicture
\setcoordinatesystem units    <1.5mm,2mm>
\setplotarea x  from 0 to 16, y  from -2 to 15
\put{1.337)} [l]  at 2 12
\put {$ \scriptstyle \bullet$} [c] at  10 12
\put {$ \scriptstyle \bullet$} [c] at  12 6
\put {$ \scriptstyle \bullet$} [c] at  12 0
\put {$ \scriptstyle \bullet$} [c] at  14 0
\put {$ \scriptstyle \bullet$} [c] at  14 6
\put {$ \scriptstyle \bullet$} [c] at  14 12
\put {$ \scriptstyle \bullet$} [c] at  16 6
\setlinear \plot 10 12 12 0 12 6 14 0 14 12 12 6   /
\setlinear \plot 14 12 16 6 14 0    /
\put{$2{,}520  $} [c] at 13 -2
\endpicture
\end{minipage}
\begin{minipage}{4cm}
\beginpicture
\setcoordinatesystem units    <1.5mm,2mm>
\setplotarea x  from 0 to 16, y  from -2 to 15
\put{1.338)} [l]  at 2 12
\put {$ \scriptstyle \bullet$} [c] at  10 0
\put {$ \scriptstyle \bullet$} [c] at  10 12
\put {$ \scriptstyle \bullet$} [c] at  11.2 9
\put {$ \scriptstyle \bullet$} [c] at  14 6
\put {$ \scriptstyle \bullet$} [c] at  15 12
\put {$ \scriptstyle \bullet$} [c] at  15 0
\put {$ \scriptstyle \bullet$} [c] at  16 6
\setlinear \plot 15 0 10 12 10 0 15 12 14 6 15 0 16 6 15 12  /
\put{$2{,}520  $} [c] at 13 -2
\endpicture
\end{minipage}
$$
$$
\begin{minipage}{4cm}
\beginpicture
\setcoordinatesystem units    <1.5mm,2mm>
\setplotarea x  from 0 to 16, y  from -2 to 15
\put{1.339)} [l]  at 2 12
\put {$ \scriptstyle \bullet$} [c] at  10 0
\put {$ \scriptstyle \bullet$} [c] at  10 12
\put {$ \scriptstyle \bullet$} [c] at  11.2 3
\put {$ \scriptstyle \bullet$} [c] at  14 6
\put {$ \scriptstyle \bullet$} [c] at  15 12
\put {$ \scriptstyle \bullet$} [c] at  15 0
\put {$ \scriptstyle \bullet$} [c] at  16 6
\setlinear \plot 15 12 10 0 10 12 15 0 14 6 15 12 16 6 15 0  /
\put{$2{,}520  $} [c] at 13 -2
\endpicture
\end{minipage}
\begin{minipage}{4cm}
\beginpicture
\setcoordinatesystem units    <1.5mm,2mm>
\setplotarea x  from 0 to 16, y from -2 to 15
\put{1.340)} [l] at 2 12
\put {$ \scriptstyle \bullet$} [c] at  10 0
\put {$ \scriptstyle \bullet$} [c] at  10 6
\put {$ \scriptstyle \bullet$} [c] at  10 12
\put {$ \scriptstyle \bullet$} [c] at  11.2 9
\put {$ \scriptstyle \bullet$} [c] at  15 0
\put {$ \scriptstyle \bullet$} [c] at  15 12
\put {$ \scriptstyle \bullet$} [c] at  16 0
\setlinear \plot  16 0 15 12 15 0 10 12 10 0 15 12 /
\put{$2{,}520$} [c] at 13 -2
\endpicture
\end{minipage}
\begin{minipage}{4cm}
\beginpicture
\setcoordinatesystem units    <1.5mm,2mm>
\setplotarea x  from 0 to 16, y from -2 to 15
\put{1.341)} [l] at 2 12
\put {$ \scriptstyle \bullet$} [c] at  10 0
\put {$ \scriptstyle \bullet$} [c] at  10 6
\put {$ \scriptstyle \bullet$} [c] at  10 12
\put {$ \scriptstyle \bullet$} [c] at  11.2 3
\put {$ \scriptstyle \bullet$} [c] at  15 0
\put {$ \scriptstyle \bullet$} [c] at  15 12
\put {$ \scriptstyle \bullet$} [c] at  16 12
\setlinear \plot  16 12 15 0 15 12 10 0 10 12 15 0 /
\put{$2{,}520$} [c] at 13 -2
\endpicture
\end{minipage}
\begin{minipage}{4cm}
\beginpicture
\setcoordinatesystem units <1.5mm, 2mm>
\setplotarea x  from 0 to 16, y from -2 to 15
\put{1.342)} [l] at 2 12
\put {$ \scriptstyle \bullet$} [c] at  10 0
\put {$ \scriptstyle \bullet$} [c] at  10 6
\put {$ \scriptstyle \bullet$} [c] at  10 12
\put {$ \scriptstyle \bullet$} [c] at  16 0
\put {$ \scriptstyle \bullet$} [c] at  16 6
\put {$ \scriptstyle \bullet$} [c] at  16 12
\put {$ \scriptstyle \bullet$} [c] at  13 12
\setlinear \plot  10 12 10 0 13 12 16 0 16 12  /
\put{$2{,}520$} [c] at 13 -2
\endpicture
\end{minipage}
\begin{minipage}{4cm}
\beginpicture
\setcoordinatesystem units <1.5mm, 2mm>
\setplotarea x  from 0 to 16, y from -2 to 15
\put{1.343)} [l] at 2 12
\put {$ \scriptstyle \bullet$} [c] at  10 0
\put {$ \scriptstyle \bullet$} [c] at  10 6
\put {$ \scriptstyle \bullet$} [c] at  10 12
\put {$ \scriptstyle \bullet$} [c] at  16 0
\put {$ \scriptstyle \bullet$} [c] at  16 6
\put {$ \scriptstyle \bullet$} [c] at  16 12
\put {$ \scriptstyle \bullet$} [c] at  13 0
\setlinear \plot  10 0 10 12 13 0 16 12 16 0  /
\put{$2{,}520$} [c] at 13 -2
\endpicture
\end{minipage}
\begin{minipage}{4cm}
\beginpicture
\setcoordinatesystem units    <1.5mm,2mm>
\setplotarea x from 0 to 16, y from -2 to 15
\put {1.344)} [l] at 2 12
\put {$ \scriptstyle \bullet$} [c] at  10 0
\put {$ \scriptstyle \bullet$} [c] at  10 12
\put {$ \scriptstyle \bullet$} [c] at  12 12
\put {$ \scriptstyle \bullet$} [c] at  14 0
\put {$ \scriptstyle \bullet$} [c] at  14  6
\put {$ \scriptstyle \bullet$} [c] at  14 12
\put {$ \scriptstyle \bullet$} [c] at  16  12
\setlinear \plot  10 12 10 0  14 12 14 0 /
\setlinear \plot  10 0 12 12 14 6 /
\plot 16 12 10 0  /
\plot 16 12 14 0  /
\put{$2{,}520$}[c]  at 13 -2
\endpicture
\end{minipage}
$$
$$
\begin{minipage}{4cm}
\beginpicture
\setcoordinatesystem units    <1.5mm,2mm>
\setplotarea x from 0 to 16, y from -2 to 15
\put {1.345)} [l] at 2 12
\put {$ \scriptstyle \bullet$} [c] at  10 0
\put {$ \scriptstyle \bullet$} [c] at  10 12
\put {$ \scriptstyle \bullet$} [c] at  12 0
\put {$ \scriptstyle \bullet$} [c] at  14 0
\put {$ \scriptstyle \bullet$} [c] at  14  6
\put {$ \scriptstyle \bullet$} [c] at  14 12
\put {$ \scriptstyle \bullet$} [c] at  16  0
\setlinear \plot  10 0 10 12  14 0 14 12 /
\setlinear \plot  10 12 12 0 14 6 /
\plot 16 0 10 12  /
\plot 16 0 14 12  /
\put{$2{,}520$}[c]  at 13 -2
\endpicture
\end{minipage}
\begin{minipage}{4cm}
\beginpicture
\setcoordinatesystem units    <1.5mm,2mm>
\setplotarea x from 0 to 16, y from -2 to 15
\put {1.346)} [l] at 2 12
\put {$ \scriptstyle \bullet$} [c] at  10 0
\put {$ \scriptstyle \bullet$} [c] at  10 12
\put {$ \scriptstyle \bullet$} [c] at  12 12
\put {$ \scriptstyle \bullet$} [c] at  14 0
\put {$ \scriptstyle \bullet$} [c] at  14 6
\put {$ \scriptstyle \bullet$} [c] at  14 12
\put {$ \scriptstyle \bullet$} [c] at  16 12
\setlinear \plot  10 12 10 0 14 12 14 0 12 12 10 0 /
\setlinear \plot  10 12 14 0 /
\setlinear \plot 16 12 14 6 /
\put{$2{,}520$}[c]  at 13 -2
\endpicture
\end{minipage}
\begin{minipage}{4cm}
\beginpicture
\setcoordinatesystem units    <1.5mm,2mm>
\setplotarea x from 0 to 16, y from -2 to 15
\put {1.347)} [l] at 2 12
\put {$ \scriptstyle \bullet$} [c] at  10 0
\put {$ \scriptstyle \bullet$} [c] at  10 12
\put {$ \scriptstyle \bullet$} [c] at  12 0
\put {$ \scriptstyle \bullet$} [c] at  14 0
\put {$ \scriptstyle \bullet$} [c] at  14 6
\put {$ \scriptstyle \bullet$} [c] at  14 12
\put {$ \scriptstyle \bullet$} [c] at  16 0
\setlinear \plot  10 0 10 12 14 0 14 12 12 0 10 12 /
\setlinear \plot  10 0 14 12 /
\setlinear \plot 16 0 14 6 /
\put{$2{,}520$}[c]  at 13 -2
\endpicture
\end{minipage}
\begin{minipage}{4cm}
\beginpicture
\setcoordinatesystem units    <1.5mm,2mm>
\setplotarea x from 0 to 16, y from -2 to 15
\put {1.348)} [l] at 2 12
\put {$ \scriptstyle \bullet$} [c] at  10 0
\put {$ \scriptstyle \bullet$} [c] at  10 12
\put {$ \scriptstyle \bullet$} [c] at  12 12
\put {$ \scriptstyle \bullet$} [c] at  14 0
\put {$ \scriptstyle \bullet$} [c] at  14 6
\put {$ \scriptstyle \bullet$} [c] at  14  12
\put {$ \scriptstyle \bullet$} [c] at  16 12
\setlinear \plot  10 12 10 0  12 12 14 6  16 12  /
\setlinear \plot  10 12 14 0 14 12 /
\put{$2{,}520$}[c]  at 13 -2
\endpicture
\end{minipage}
\begin{minipage}{4cm}
\beginpicture
\setcoordinatesystem units    <1.5mm,2mm>
\setplotarea x from 0 to 16, y from -2 to 15
\put {1.349)} [l] at 2 12
\put {$ \scriptstyle \bullet$} [c] at  10 0
\put {$ \scriptstyle \bullet$} [c] at  10 12
\put {$ \scriptstyle \bullet$} [c] at  12 0
\put {$ \scriptstyle \bullet$} [c] at  14 0
\put {$ \scriptstyle \bullet$} [c] at  14 6
\put {$ \scriptstyle \bullet$} [c] at  14  12
\put {$ \scriptstyle \bullet$} [c] at  16 0
\setlinear \plot  10 0 10 12  12 0 14 6  16 0  /
\setlinear \plot  10 0 14 12 14 0 /
\put{$2{,}520$}[c]  at 13 -2
\endpicture
\end{minipage}
\begin{minipage}{4cm}
\beginpicture
\setcoordinatesystem units    <1.5mm,2mm>
\setplotarea x  from 0 to 16, y  from -2 to 15
\put{1.350)} [l]  at 2 12
\put {$ \scriptstyle \bullet$} [c] at  10 0
\put {$ \scriptstyle \bullet$} [c] at  10 12
\put {$ \scriptstyle \bullet$} [c] at  13 12
\put {$ \scriptstyle \bullet$} [c] at  13 6
\put {$ \scriptstyle \bullet$} [c] at  13 0
\put {$ \scriptstyle \bullet$} [c] at  16 12
\put {$ \scriptstyle \bullet$} [c] at  16  0
\setlinear \plot  10 0 10 12 13 6 13 0  16 12 16 0   /
\setlinear \plot  10 0 13  12  13 6  /
\put{$2{,}520 $} [c]  at 13 -2
\endpicture
\end{minipage}
$$
$$
\begin{minipage}{4cm}
\beginpicture
\setcoordinatesystem units    <1.5mm,2mm>
\setplotarea x  from 0 to 16, y  from -2 to 15
\put{1.351)} [l]  at 2 12
\put {$ \scriptstyle \bullet$} [c] at  10 0
\put {$ \scriptstyle \bullet$} [c] at  10 12
\put {$ \scriptstyle \bullet$} [c] at  13 12
\put {$ \scriptstyle \bullet$} [c] at  13 6
\put {$ \scriptstyle \bullet$} [c] at  13 0
\put {$ \scriptstyle \bullet$} [c] at  16 12
\put {$ \scriptstyle \bullet$} [c] at  16  0
\setlinear \plot  10 12 10 0 13 6 13 12  16 0 16 12   /
\setlinear \plot  10 12 13  0  13 6  /
\put{$2{,}520 $} [c]  at 13 -2
\endpicture
\end{minipage}
\begin{minipage}{4cm}
\beginpicture
\setcoordinatesystem units    <1.5mm,2mm>
\setplotarea x from 0 to 16, y from -2 to 15
\put{1.352)} [l]  at 2 12
\put {$ \scriptstyle \bullet$} [c] at 10 12
\put {$ \scriptstyle \bullet$} [c] at 13 12
\put {$ \scriptstyle \bullet$} [c] at 13 6
\put {$ \scriptstyle \bullet$} [c] at 14.5 0
\put {$ \scriptstyle \bullet$} [c] at 14.5 6
\put {$ \scriptstyle \bullet$} [c] at 14.5 12
\put {$ \scriptstyle \bullet$} [c] at 16 6
\setlinear \plot 13 12 13 6 14.5 0 14.5 12 16 6 14.5 0  /
\setlinear \plot 10 12 13 6  14.5 12  /
\put{$1{,}260$} [c]  at 13 -2
\endpicture
\end{minipage}
\begin{minipage}{4cm}
\beginpicture
\setcoordinatesystem units    <1.5mm,2mm>
\setplotarea x from 0 to 16, y from -2 to 15
\put{1.353)} [l]  at 2 12
\put {$ \scriptstyle \bullet$} [c] at 10 0
\put {$ \scriptstyle \bullet$} [c] at 13 0
\put {$ \scriptstyle \bullet$} [c] at 13 6
\put {$ \scriptstyle \bullet$} [c] at 14.5 0
\put {$ \scriptstyle \bullet$} [c] at 14.5 6
\put {$ \scriptstyle \bullet$} [c] at 14.5 12
\put {$ \scriptstyle \bullet$} [c] at 16 6
\setlinear \plot 13 0 13 6 14.5 12 14.5 0 16 6 14.5 12  /
\setlinear \plot 10 0 13 6  14.5 0  /
\put{$1{,}260$} [c]  at 13 -2
\endpicture
\end{minipage}
\begin{minipage}{4cm}
\beginpicture
\setcoordinatesystem units    <1.5mm,2mm>
\setplotarea x from 0 to 16, y from -2 to 15
\put{1.354)} [l]  at 2 12
\put {$ \scriptstyle \bullet$} [c] at 10 9
\put {$ \scriptstyle \bullet$} [c] at 11 12
\put {$ \scriptstyle \bullet$} [c] at 11 6
\put {$ \scriptstyle \bullet$} [c] at 12 9
\put {$ \scriptstyle \bullet$} [c] at 13 0
\put {$ \scriptstyle \bullet$} [c] at 14 12
\put {$ \scriptstyle \bullet$} [c] at 16 12

\setlinear \plot 16 12 13 0 11 6 10 9 11 12 12 9 11 6   /
\setlinear \plot 14 12 13 0 /
\put{$1{,}260$} [c]  at 13 -2
\endpicture
\end{minipage}
\begin{minipage}{4cm}
\beginpicture
\setcoordinatesystem units    <1.5mm,2mm>
\setplotarea x from 0 to 16, y from -2 to 15
\put{1.355)} [l]  at 2 12
\put {$ \scriptstyle \bullet$} [c] at 10 3
\put {$ \scriptstyle \bullet$} [c] at 11 0
\put {$ \scriptstyle \bullet$} [c] at 11 6
\put {$ \scriptstyle \bullet$} [c] at 12 3
\put {$ \scriptstyle \bullet$} [c] at 13 12
\put {$ \scriptstyle \bullet$} [c] at 14 0
\put {$ \scriptstyle \bullet$} [c] at 16 0
\setlinear \plot 16 0 13 12 11 6 10 3 11 0 12 3 11 6   /
\setlinear \plot 14 0 13 12 /
\put{$1{,}260$} [c]  at 13 -2
\endpicture
\end{minipage}
\begin{minipage}{4cm}
\beginpicture
\setcoordinatesystem units    <1.5mm,2mm>
\setplotarea x from 0 to 16, y from -2 to 15
\put{1.356)} [l]  at 2 12
\put {$ \scriptstyle \bullet$} [c] at 10 12
\put {$ \scriptstyle \bullet$} [c] at 11  9
\put {$ \scriptstyle \bullet$} [c] at 12 12
\put {$ \scriptstyle \bullet$} [c] at 11 6
\put {$ \scriptstyle \bullet$} [c] at 13 0
\put {$ \scriptstyle \bullet$} [c] at 14 12
\put {$ \scriptstyle \bullet$} [c] at 16 12
\setlinear \plot 16 12  13 0 11 6 11 9 10 12    /
\setlinear \plot 12 12 11 9      /
\setlinear \plot 14 12 13 0      /

\put{$1{,}260$} [c]  at 13 -2
\endpicture
\end{minipage}
$$
$$
\begin{minipage}{4cm}
\beginpicture
\setcoordinatesystem units    <1.5mm,2mm>
\setplotarea x from 0 to 16, y from -2 to 15
\put{1.357)} [l]  at 2 12
\put {$ \scriptstyle \bullet$} [c] at 10 0
\put {$ \scriptstyle \bullet$} [c] at 11  3
\put {$ \scriptstyle \bullet$} [c] at 12 0
\put {$ \scriptstyle \bullet$} [c] at 11 6
\put {$ \scriptstyle \bullet$} [c] at 13 12
\put {$ \scriptstyle \bullet$} [c] at 14 0
\put {$ \scriptstyle \bullet$} [c] at 16 0
\setlinear \plot 16 0  13 12 11 6 11 3 10 0    /
\setlinear \plot 12 0 11 3      /
\setlinear \plot 14 0 13 12      /
\put{$1{,}260$} [c]  at 13 -2
\endpicture
\end{minipage}
\begin{minipage}{4cm}
\beginpicture
\setcoordinatesystem units    <1.5mm,2mm>
\setplotarea x from 0 to 16, y from -2 to 15
\put{1.358)} [l]  at 2 12
\put {$ \scriptstyle \bullet$} [c] at 11 12
\put {$ \scriptstyle \bullet$} [c] at 11 10
\put {$ \scriptstyle \bullet$} [c] at 10 7
\put {$ \scriptstyle \bullet$} [c] at 12 7
\put {$ \scriptstyle \bullet$} [c] at 11 0
\put {$ \scriptstyle \bullet$} [c] at 14 12
\put {$ \scriptstyle \bullet$} [c] at 16 12

\setlinear \plot 14 12 11 0 16 12    /
\setlinear \plot 11 12 11 10 10 7  11 0 12 7 11 10      /
\put{$1{,}260$} [c]  at 13 -2
\endpicture
\end{minipage}
\begin{minipage}{4cm}
\beginpicture
\setcoordinatesystem units    <1.5mm,2mm>
\setplotarea x from 0 to 16, y from -2 to 15
\put{1.359)} [l]  at 2 12
\put {$ \scriptstyle \bullet$} [c] at 11 0
\put {$ \scriptstyle \bullet$} [c] at 11 2
\put {$ \scriptstyle \bullet$} [c] at 10 5
\put {$ \scriptstyle \bullet$} [c] at 12 5
\put {$ \scriptstyle \bullet$} [c] at 11 12
\put {$ \scriptstyle \bullet$} [c] at 14 0
\put {$ \scriptstyle \bullet$} [c] at 16 0
\setlinear \plot 14 0 11 12 16 0    /
\setlinear \plot 11 0 11 2 10 5  11 12 12 5 11 2      /
\put{$1{,}260$} [c]  at 13 -2
\endpicture
\end{minipage}
\begin{minipage}{4cm}
\beginpicture
\setcoordinatesystem units    <1.5mm,2mm>
\setplotarea x  from 0 to 16, y from -2 to 15
\put{1.360)} [l] at 2 12
\put {$ \scriptstyle \bullet$} [c] at  10 0
\put {$ \scriptstyle \bullet$} [c] at  10 6
\put {$ \scriptstyle \bullet$} [c] at  10 12
\put {$ \scriptstyle \bullet$} [c] at  11.5 9
\put {$ \scriptstyle \bullet$} [c] at  13 12
\put {$ \scriptstyle \bullet$} [c] at  16 0
\put {$ \scriptstyle \bullet$} [c] at  16 12
\setlinear \plot  10 0 10 12 16 0 16  12   10  0 13  12 16 0  /
\put{$1{,}260$} [c] at 13 -2
\endpicture
\end{minipage}
\begin{minipage}{4cm}
\beginpicture
\setcoordinatesystem units    <1.5mm,2mm>
\setplotarea x  from 0 to 16, y from -2 to 15
\put{1.361)} [l] at 2 12
\put {$ \scriptstyle \bullet$} [c] at  10 0
\put {$ \scriptstyle \bullet$} [c] at  10 6
\put {$ \scriptstyle \bullet$} [c] at  10 12
\put {$ \scriptstyle \bullet$} [c] at  11.5 3
\put {$ \scriptstyle \bullet$} [c] at  13 0
\put {$ \scriptstyle \bullet$} [c] at  16 0
\put {$ \scriptstyle \bullet$} [c] at  16 12
\setlinear \plot  10 12 10 0 16 12 16  0   10  12 13  0 16 12  /
\put{$1{,}260$} [c] at 13 -2
\endpicture
\end{minipage}
\begin{minipage}{4cm}
\beginpicture
\setcoordinatesystem units    <1.5mm,2mm>
\setplotarea x  from 0 to 16, y from -2 to 15
\put{1.362)} [l] at 2 12
\put {$ \scriptstyle \bullet$} [c] at  10 12
\put {$ \scriptstyle \bullet$} [c] at  12 6
\put {$ \scriptstyle \bullet$} [c] at  12 0
\put {$ \scriptstyle \bullet$} [c] at  13 12
\put {$ \scriptstyle \bullet$} [c] at  14 6
\put {$ \scriptstyle \bullet$} [c] at  14 0
\put {$ \scriptstyle \bullet$} [c] at  16 12
\setlinear \plot  10 12  12 0 12 6 13 12 14 6 12 0 /
\setlinear \plot  12 6 14  0 16 12 /
\setlinear \plot  14 6 14  0 /
\put{$1{,}260$} [c] at 13 -2
\endpicture
\end{minipage}
$$
$$
\begin{minipage}{4cm}
\beginpicture
\setcoordinatesystem units    <1.5mm,2mm>
\setplotarea x  from 0 to 16, y from -2 to 15
\put{1.363)} [l] at 2 12
\put {$ \scriptstyle \bullet$} [c] at  10 0
\put {$ \scriptstyle \bullet$} [c] at  12 6
\put {$ \scriptstyle \bullet$} [c] at  12 12
\put {$ \scriptstyle \bullet$} [c] at  13 0
\put {$ \scriptstyle \bullet$} [c] at  14 6
\put {$ \scriptstyle \bullet$} [c] at  14 12
\put {$ \scriptstyle \bullet$} [c] at  16 0
\setlinear \plot  10 0  12 12 12 6 13 0 14 6 12 12 /
\setlinear \plot  12 6 14  12 16 0 /
\setlinear \plot  14 6 14  12 /
\put{$1{,}260$} [c] at 13 -2
\endpicture
\end{minipage}
\begin{minipage}{4cm}
\beginpicture
\setcoordinatesystem units    <1.5mm,2mm>
\setplotarea x from 0 to 16, y from -2 to 15
\put {1.364)} [l] at 2 12
\put {$ \scriptstyle \bullet$} [c] at  10 0
\put {$ \scriptstyle \bullet$} [c] at  14 12
\put {$ \scriptstyle \bullet$} [c] at  12 12
\put {$ \scriptstyle \bullet$} [c] at  16 0
\put {$ \scriptstyle \bullet$} [c] at  10 12
\put {$ \scriptstyle \bullet$} [c] at  13 6
\put {$ \scriptstyle \bullet$} [c] at  16 12
\setlinear \plot  10 12 10 0 13 6 16 0 16 12  /
\setlinear \plot  12 12 13 6  14  12   /
\put{$1{,}260$}[c]  at 13 -2
\endpicture
\end{minipage}
\begin{minipage}{4cm}
\beginpicture
\setcoordinatesystem units    <1.5mm,2mm>
\setplotarea x from 0 to 16, y from -2 to 15
\put {1.365)} [l] at 2 12
\put {$ \scriptstyle \bullet$} [c] at  10 0
\put {$ \scriptstyle \bullet$} [c] at  14 0
\put {$ \scriptstyle \bullet$} [c] at  12 0
\put {$ \scriptstyle \bullet$} [c] at  16 0
\put {$ \scriptstyle \bullet$} [c] at  10 12
\put {$ \scriptstyle \bullet$} [c] at  13 6
\put {$ \scriptstyle \bullet$} [c] at  16 12
\setlinear \plot  10 0 10 12 13 6 16 12 16 0  /
\setlinear \plot  12 0 13 6  14  0   /
\put{$1{,}260$}[c]  at 13 -2
\endpicture
\end{minipage}
\begin{minipage}{4cm}
\beginpicture
\setcoordinatesystem units    <1.5mm,2mm>
\setplotarea x  from 0 to 16, y  from -2 to 15
\put{1.366)} [l]  at 2 12
\put {$ \scriptstyle \bullet$} [c] at  10 0
\put {$ \scriptstyle \bullet$} [c] at  10 12
\put {$ \scriptstyle \bullet$} [c] at  13 12
\put {$ \scriptstyle \bullet$} [c] at  13 0
\put {$ \scriptstyle \bullet$} [c] at  16  0
\put {$ \scriptstyle \bullet$} [c] at  16  6
\put {$ \scriptstyle \bullet$} [c] at  16  12
\setlinear \plot  10 0 10 12 13  0 13 12  10 0 16 12 13 0 /
\setlinear \plot  16 12  16  0   /
\put{$1{,}260$} [c]  at 13 -2
\endpicture
\end{minipage}
\begin{minipage}{4cm}
\beginpicture
\setcoordinatesystem units    <1.5mm,2mm>
\setplotarea x  from 0 to 16, y  from -2 to 15
\put{1.367)} [l]  at 2 12
\put {$ \scriptstyle \bullet$} [c] at  10 0
\put {$ \scriptstyle \bullet$} [c] at  10 12
\put {$ \scriptstyle \bullet$} [c] at  13 12
\put {$ \scriptstyle \bullet$} [c] at  13 0
\put {$ \scriptstyle \bullet$} [c] at  16  0
\put {$ \scriptstyle \bullet$} [c] at  16  6
\put {$ \scriptstyle \bullet$} [c] at  16  12
\setlinear \plot  10 12 10 0 13  12 13 0  10 12 16 0 13 12 /
\setlinear \plot  16 12  16  0   /
\put{$1{,}260$} [c]  at 13 -2
\endpicture
\end{minipage}
\begin{minipage}{4cm}
\beginpicture
\setcoordinatesystem units    <1.5mm,2mm>
\setplotarea x from 0 to 16, y  from -2 to 15
\put {1.368)} [l] at 2 12
\put {$ \scriptstyle \bullet$} [c] at  10 12
\put {$ \scriptstyle \bullet$} [c] at  12 12
\put {$ \scriptstyle \bullet$} [c] at  14 12
\put {$ \scriptstyle \bullet$} [c] at  16 12
\put {$ \scriptstyle \bullet$} [c] at  10 0
\put {$ \scriptstyle \bullet$} [c] at  13 0
\put {$ \scriptstyle \bullet$} [c] at  16 0
\setlinear \plot   10  12 10 0   12 12 13 0 14  12 16 0  16 12  10 0   14 12 /
\setlinear \plot  10 12  13 0  /
\setlinear \plot  12 12  16 0  /
\put{$1{,}260$} [c]  at 13 -2
\endpicture
\end{minipage}
$$
$$
\begin{minipage}{4cm}
\beginpicture
\setcoordinatesystem units    <1.5mm,2mm>
\setplotarea x from 0 to 16, y  from -2 to 15
\put {1.369)} [l] at 2 12
\put {$ \scriptstyle \bullet$} [c] at  10 0
\put {$ \scriptstyle \bullet$} [c] at  12 0
\put {$ \scriptstyle \bullet$} [c] at  14 0
\put {$ \scriptstyle \bullet$} [c] at  16 0
\put {$ \scriptstyle \bullet$} [c] at  10 12
\put {$ \scriptstyle \bullet$} [c] at  13 12
\put {$ \scriptstyle \bullet$} [c] at  16 12
\setlinear \plot   10  0 10 12   12 0 13 12 14  0 16 12  16 0  10 12   14 0 /
\setlinear \plot  10 0  13 12  /
\setlinear \plot  12 0  16 12  /
\put{$1{,}260$} [c]  at 13 -2
\endpicture
\end{minipage}
\begin{minipage}{4cm}
\beginpicture
\setcoordinatesystem units    <1.5mm,2mm>
\setplotarea x from 0 to  16, y  from -2 to 15
\put{1.370)} [l]  at 2 12
\put {$ \scriptstyle \bullet$} [c] at 10 0
\put {$ \scriptstyle \bullet$} [c] at  10 4
\put {$ \scriptstyle \bullet$} [c] at 10 8
\put {$ \scriptstyle \bullet$} [c] at 10 12
\setlinear \plot 10 0  10  12   /
\put {$ \scriptstyle \bullet$} [c] at 12 12
\put {$ \scriptstyle \bullet$} [c] at 14 0
\put {$ \scriptstyle \bullet$} [c] at 16 12
\setlinear \plot 12 12 14  0 16 12   /
\put{$2{,}520$} [c] at 13 -2
\endpicture
\end{minipage}
\begin{minipage}{4cm}
\beginpicture
\setcoordinatesystem units    <1.5mm,2mm>
\setplotarea x from 0 to  16, y  from -2 to 15
\put{1.371)} [l]  at 2 12
\put {$ \scriptstyle \bullet$} [c] at 10 0
\put {$ \scriptstyle \bullet$} [c] at  10 4
\put {$ \scriptstyle \bullet$} [c] at 10 8
\put {$ \scriptstyle \bullet$} [c] at 10 12
\setlinear \plot 10 0  10  12   /
\put {$ \scriptstyle \bullet$} [c] at 12 0
\put {$ \scriptstyle \bullet$} [c] at 14 12
\put {$ \scriptstyle \bullet$} [c] at 16 0
\setlinear \plot 12 0 14  12 16 0   /
\put{$2{,}520$} [c] at 13 -2
\endpicture
\end{minipage}
\begin{minipage}{4cm}
\beginpicture
\setcoordinatesystem units    <1.5mm,2mm>
\setplotarea x from 0 to 16, y from -2 to 15
\put{${\bf  26}$} [l] at 2 15

\put{1.372)} [l] at 2 12
\put {$ \scriptstyle \bullet$} [c] at  10 12
\put {$ \scriptstyle \bullet$} [c] at  10 0
\put {$ \scriptstyle \bullet$} [c] at  12.5 6
\put {$ \scriptstyle \bullet$} [c] at  13 12
\put {$ \scriptstyle \bullet$} [c] at  13 0
\put {$ \scriptstyle \bullet$} [c] at  13.5 6
\put {$ \scriptstyle \bullet$} [c] at  16 12
\setlinear \plot 10 0 10 12 12.5 6 13 0 13.5 6 13 12 12.5 6  /
\setlinear \plot  10 0 13 12  /
\setlinear \plot 13 0 16 12    /
\put{$5{,}040$} [c] at 13 -2
\endpicture
\end{minipage}
\begin{minipage}{4cm}
\beginpicture
\setcoordinatesystem units    <1.5mm,2mm>
\setplotarea x from 0 to 16, y from -2 to 15
\put{1.373)} [l] at 2 12
\put {$ \scriptstyle \bullet$} [c] at  10 12
\put {$ \scriptstyle \bullet$} [c] at  10 0
\put {$ \scriptstyle \bullet$} [c] at  12.5 6
\put {$ \scriptstyle \bullet$} [c] at  13 12
\put {$ \scriptstyle \bullet$} [c] at  13 0
\put {$ \scriptstyle \bullet$} [c] at  13.5 6
\put {$ \scriptstyle \bullet$} [c] at  16 0
\setlinear \plot 10 12 10 0 12.5 6 13 12 13.5 6 13 0 12.5 6  /
\setlinear \plot  10 12 13 0  /
\setlinear \plot 13 12 16 0    /
\put{$5{,}040$} [c] at 13 -2
\endpicture
\end{minipage}
\begin{minipage}{4cm}
\beginpicture
\setcoordinatesystem units    <1.5mm,2mm>
\setplotarea x from 0 to 16, y from -2 to 15
\put{1.374)} [l] at 2 12
\put {$ \scriptstyle \bullet$} [c] at  10 12
\put {$ \scriptstyle \bullet$} [c] at  12.5 6
\put {$ \scriptstyle \bullet$} [c] at  13 0
\put {$ \scriptstyle \bullet$} [c] at  13 12
\put {$ \scriptstyle \bullet$} [c] at  13.5 6
\put {$ \scriptstyle \bullet$} [c] at  16 0
\put {$ \scriptstyle \bullet$} [c] at  16 12
\setlinear \plot 16 0 16 12 13 0 12.5 6 13 12 13.5 6 13 0  /
\setlinear \plot  12.5 6 10  12   /
\put{$5{,}040$} [c] at 13 -2
\endpicture
\end{minipage}
$$
$$
\begin{minipage}{4cm}
\beginpicture
\setcoordinatesystem units    <1.5mm,2mm>
\setplotarea x from 0 to 16, y from -2 to 15
\put{1.375)} [l] at 2 12
\put {$ \scriptstyle \bullet$} [c] at  10 0
\put {$ \scriptstyle \bullet$} [c] at  12.5 6
\put {$ \scriptstyle \bullet$} [c] at  13 0
\put {$ \scriptstyle \bullet$} [c] at  13 12
\put {$ \scriptstyle \bullet$} [c] at  13.5 6
\put {$ \scriptstyle \bullet$} [c] at  16 0
\put {$ \scriptstyle \bullet$} [c] at  16 12
\setlinear \plot 16 12 16 0 13 12 12.5 6 13 0 13.5 6 13 12  /
\setlinear \plot  12.5 6 10  0   /
\put{$5{,}040$} [c] at 13 -2
\endpicture
\end{minipage}
\begin{minipage}{4cm}
\beginpicture
\setcoordinatesystem units    <1.5mm,2mm>
\setplotarea x from 0 to 16, y from -2 to 15
\put{1.376)} [l] at 2 12
\put {$ \scriptstyle \bullet$} [c] at  10 12
\put {$ \scriptstyle \bullet$} [c] at  12 12
\put {$ \scriptstyle \bullet$} [c] at  12 6
\put {$ \scriptstyle \bullet$} [c] at  12 0
\put {$ \scriptstyle \bullet$} [c] at  14 9
\put {$ \scriptstyle \bullet$} [c] at  16 12
\put {$ \scriptstyle \bullet$} [c] at  16 0
\setlinear \plot  16 0 16 12 12 6 12 12 /
\setlinear \plot  10 12 12 0 12 6 /
\put{$5{,}040$} [c] at 13 -2
\endpicture
\end{minipage}
\begin{minipage}{4cm}
\beginpicture
\setcoordinatesystem units    <1.5mm,2mm>
\setplotarea x from 0 to 16, y from -2 to 15
\put{1.377)} [l] at 2 12
\put {$ \scriptstyle \bullet$} [c] at  10 0
\put {$ \scriptstyle \bullet$} [c] at  12 12
\put {$ \scriptstyle \bullet$} [c] at  12 6
\put {$ \scriptstyle \bullet$} [c] at  12 0
\put {$ \scriptstyle \bullet$} [c] at  14 3
\put {$ \scriptstyle \bullet$} [c] at  16 12
\put {$ \scriptstyle \bullet$} [c] at  16 0
\setlinear \plot  16 12 16 0 12 6 12 0 /
\setlinear \plot  10 0 12 12 12 6 /
\put{$5{,}040$} [c] at 13 -2
\endpicture
\end{minipage}
\begin{minipage}{4cm}
\beginpicture
\setcoordinatesystem units    <1.5mm,2mm>
\setplotarea x from 0 to 16, y from -2 to 15
\put{1.378)} [l] at 2 12
\put {$ \scriptstyle \bullet$} [c] at  10 0
\put {$ \scriptstyle \bullet$} [c] at  10 4
\put {$ \scriptstyle \bullet$} [c] at  10 8
\put {$ \scriptstyle \bullet$} [c] at  10 12
\put {$ \scriptstyle \bullet$} [c] at  13 12
\put {$ \scriptstyle \bullet$} [c] at  16 0
\put {$ \scriptstyle \bullet$} [c] at  16 12
\setlinear \plot  10 12 10 0  16 12 16 0 /
\setlinear \plot  10 0 13 12   /
\put{$5{,}040$} [c] at 13 -2
\endpicture
\end{minipage}
\begin{minipage}{4cm}
\beginpicture
\setcoordinatesystem units    <1.5mm,2mm>
\setplotarea x from 0 to 16, y from -2 to 15
\put{1.379)} [l] at 2 12
\put {$ \scriptstyle \bullet$} [c] at  10 0
\put {$ \scriptstyle \bullet$} [c] at  10 4
\put {$ \scriptstyle \bullet$} [c] at  10 8
\put {$ \scriptstyle \bullet$} [c] at  10 12
\put {$ \scriptstyle \bullet$} [c] at  13 0
\put {$ \scriptstyle \bullet$} [c] at  16 0
\put {$ \scriptstyle \bullet$} [c] at  16 12
\setlinear \plot  10 0 10 12  16 0 16 12 /
\setlinear \plot  10 12 13 0   /
\put{$5{,}040$} [c] at 13 -2
\endpicture
\end{minipage}
\begin{minipage}{4cm}
\beginpicture
\setcoordinatesystem units    <1.5mm,2mm>
\setplotarea x from 0 to 16, y from -2 to 15
\put{1.380)} [l] at 2 12
\put {$ \scriptstyle \bullet$} [c] at  16 12
\put {$ \scriptstyle \bullet$} [c] at  13.5 6
\put {$ \scriptstyle \bullet$} [c] at  13 0
\put {$ \scriptstyle \bullet$} [c] at  13 12
\put {$ \scriptstyle \bullet$} [c] at  12.5 6
\put {$ \scriptstyle \bullet$} [c] at  10 12
\put {$ \scriptstyle \bullet$} [c] at  10 0
\setlinear \plot 16 12 13 0 10 12 10 0 13.5 6 13 0 12.5 6 13 12 13.5 6 /
\put{$5{,}040$} [c] at 13 -2
\endpicture
\end{minipage}
$$
$$
\begin{minipage}{4cm}
\beginpicture
\setcoordinatesystem units    <1.5mm,2mm>
\setplotarea x from 0 to 16, y from -2 to 15
\put{1.381)} [l] at 2 12
\put {$ \scriptstyle \bullet$} [c] at  16 0
\put {$ \scriptstyle \bullet$} [c] at  13.5 6
\put {$ \scriptstyle \bullet$} [c] at  13 0
\put {$ \scriptstyle \bullet$} [c] at  13 12
\put {$ \scriptstyle \bullet$} [c] at  12.5 6
\put {$ \scriptstyle \bullet$} [c] at  10 12
\put {$ \scriptstyle \bullet$} [c] at  10 0
\setlinear \plot 16 0 13 12 10 0 10 12 13.5 6 13 12 12.5 6 13 0 13.5 6 /
\put{$5{,}040$} [c] at 13 -2
\endpicture
\end{minipage}
\begin{minipage}{4cm}
\beginpicture
\setcoordinatesystem units    <1.5mm,2mm>
\setplotarea x from 0 to 16, y from -2 to 15
\put{1.382)} [l] at 2 12
\put {$ \scriptstyle \bullet$} [c] at  10 12
\put {$ \scriptstyle \bullet$} [c] at  12 0
\put {$ \scriptstyle \bullet$} [c] at  12 6
\put {$ \scriptstyle \bullet$} [c] at  12 12
\put {$ \scriptstyle \bullet$} [c] at  16 6
\put {$ \scriptstyle \bullet$} [c] at  16 0
\put {$ \scriptstyle \bullet$} [c] at  16 12
\setlinear \plot  10 12 12 0 12 12 16 0 16 12  /
\put{$5{,}040$} [c] at 13 -2
\endpicture
\end{minipage}
\begin{minipage}{4cm}
\beginpicture
\setcoordinatesystem units    <1.5mm,2mm>
\setplotarea x from 0 to 16, y from -2 to 15
\put{1.383)} [l] at 2 12
\put {$ \scriptstyle \bullet$} [c] at  10 0
\put {$ \scriptstyle \bullet$} [c] at  12 0
\put {$ \scriptstyle \bullet$} [c] at  12 6
\put {$ \scriptstyle \bullet$} [c] at  12 12
\put {$ \scriptstyle \bullet$} [c] at  16 6
\put {$ \scriptstyle \bullet$} [c] at  16 0
\put {$ \scriptstyle \bullet$} [c] at  16 12
\setlinear \plot  10 0 12  12 12 0 16 12 16 0  /
\put{$5{,}040$} [c] at 13 -2
\endpicture
\end{minipage}
\begin{minipage}{4cm}
\beginpicture
\setcoordinatesystem units    <1.5mm,2mm>
\setplotarea x from 0 to 16, y from -2 to 15
\put{1.384)} [l] at 2 12
\put {$ \scriptstyle \bullet$} [c] at  10 0
\put {$ \scriptstyle \bullet$} [c] at  12 12
\put {$ \scriptstyle \bullet$} [c] at  12 3
\put {$ \scriptstyle \bullet$} [c] at  14 12
\put {$ \scriptstyle \bullet$} [c] at  14 0
\put {$ \scriptstyle \bullet$} [c] at  16 12
\put {$ \scriptstyle \bullet$} [c] at  16 3
\setlinear \plot  10 0 12 12 12 3 14 0 14 12  10 0   /
\setlinear \plot  16 12 16 3 14 0   /
\put{$5{,}040$} [c] at 13 -2
\endpicture
\end{minipage}
\begin{minipage}{4cm}
\beginpicture
\setcoordinatesystem units    <1.5mm,2mm>
\setplotarea x from 0 to 16, y from -2 to 15
\put{1.385)} [l] at 2 12
\put {$ \scriptstyle \bullet$} [c] at  10 12
\put {$ \scriptstyle \bullet$} [c] at  12 0
\put {$ \scriptstyle \bullet$} [c] at  12 9
\put {$ \scriptstyle \bullet$} [c] at  14 12
\put {$ \scriptstyle \bullet$} [c] at  14 0
\put {$ \scriptstyle \bullet$} [c] at  16 0
\put {$ \scriptstyle \bullet$} [c] at  16 9
\setlinear \plot  10 12 12 0 12 9 14 12 14 0  10 12   /
\setlinear \plot  16 0 16 9 14 12   /
\put{$5{,}040$} [c] at 13 -2
\endpicture
\end{minipage}
\begin{minipage}{4cm}
\beginpicture
\setcoordinatesystem units    <1.5mm,2mm>
\setplotarea x from 0 to 16, y from -2 to 15
\put{1.386)} [l] at 2 12
\put {$ \scriptstyle \bullet$} [c] at  10 12
\put {$ \scriptstyle \bullet$} [c] at  12 0
\put {$ \scriptstyle \bullet$} [c] at  12 6
\put {$ \scriptstyle \bullet$} [c] at  12 12
\put {$ \scriptstyle \bullet$} [c] at  16 0
\put {$ \scriptstyle \bullet$} [c] at  16 6
\put {$ \scriptstyle \bullet$} [c] at  16 12
\setlinear \plot  10 12 12  0 12 12   /
\setlinear \plot  16 0 16 12   /
\setlinear \plot  12 0 16  6      /
\put{$5{,}040$} [c] at 13 -2
\endpicture
\end{minipage}
$$
$$
\begin{minipage}{4cm}
\beginpicture
\setcoordinatesystem units    <1.5mm,2mm>
\setplotarea x from 0 to 16, y from -2 to 15
\put{1.387)} [l] at 2 12
\put {$ \scriptstyle \bullet$} [c] at  10 0
\put {$ \scriptstyle \bullet$} [c] at  12 0
\put {$ \scriptstyle \bullet$} [c] at  12 6
\put {$ \scriptstyle \bullet$} [c] at  12 12
\put {$ \scriptstyle \bullet$} [c] at  16 0
\put {$ \scriptstyle \bullet$} [c] at  16 6
\put {$ \scriptstyle \bullet$} [c] at  16 12
\setlinear \plot  10 0 12  12 12 0   /
\setlinear \plot  16 0 16 12   /
\setlinear \plot  12 12 16  6      /
\put{$5{,}040$} [c] at 13 -2
\endpicture
\end{minipage}
\begin{minipage}{4cm}
\beginpicture
\setcoordinatesystem units    <1.5mm,2mm>
\setplotarea x from 0 to 16, y from -2 to 15
\put{1.388)} [l] at 2 12
\put {$ \scriptstyle \bullet$} [c] at  10 0
\put {$ \scriptstyle \bullet$} [c] at  10 6
\put {$ \scriptstyle \bullet$} [c] at  10 12
\put {$ \scriptstyle \bullet$} [c] at  13 6
\put {$ \scriptstyle \bullet$} [c] at  13 12
\put {$ \scriptstyle \bullet$} [c] at  16 0
\put {$ \scriptstyle \bullet$} [c] at  16 12
\setlinear \plot 10 12 10 0 16 12 16 0  /
\setlinear \plot  13 12 13 6    /
\put{$5{,}040$} [c] at 13 -2
\endpicture
\end{minipage}
\begin{minipage}{4cm}
\beginpicture
\setcoordinatesystem units    <1.5mm,2mm>
\setplotarea x from 0 to 16, y from -2 to 15
\put{1.389)} [l] at 2 12
\put {$ \scriptstyle \bullet$} [c] at  10 0
\put {$ \scriptstyle \bullet$} [c] at  10 6
\put {$ \scriptstyle \bullet$} [c] at  10 12
\put {$ \scriptstyle \bullet$} [c] at  13 6
\put {$ \scriptstyle \bullet$} [c] at  13 0
\put {$ \scriptstyle \bullet$} [c] at  16 0
\put {$ \scriptstyle \bullet$} [c] at  16 12
\setlinear \plot 10 0 10 12 16 0 16 12  /
\setlinear \plot  13 0 13 6    /
\put{$5{,}040$} [c] at 13 -2
\endpicture
\end{minipage}
\begin{minipage}{4cm}
\beginpicture
\setcoordinatesystem units    <1.5mm,2mm>
\setplotarea x from 0 to 16, y from -2 to 15
\put {1.390)} [l] at 2 12
\put {$ \scriptstyle \bullet$} [c] at  10 12
\put {$ \scriptstyle \bullet$} [c] at  12 12
\put {$ \scriptstyle \bullet$} [c] at  14 12
\put {$ \scriptstyle \bullet$} [c] at  16 12
\put {$ \scriptstyle \bullet$} [c] at  10 0
\put {$ \scriptstyle \bullet$} [c] at  14 6
\put {$ \scriptstyle \bullet$} [c] at  14 0
\setlinear \plot  10 12 10 0  14 12 14 0  12 12 10 0 /
\setlinear \plot 16 12 14 6 /
\put{$5{,}040$}[c] at 13 -2
\endpicture
\end{minipage}
\begin{minipage}{4cm}
\beginpicture
\setcoordinatesystem units    <1.5mm,2mm>
\setplotarea x from 0 to 16, y from -2 to 15
\put {1.391)} [l] at 2 12
\put {$ \scriptstyle \bullet$} [c] at  10 0
\put {$ \scriptstyle \bullet$} [c] at  12 0
\put {$ \scriptstyle \bullet$} [c] at  14 0
\put {$ \scriptstyle \bullet$} [c] at  16 0
\put {$ \scriptstyle \bullet$} [c] at  10 12
\put {$ \scriptstyle \bullet$} [c] at  14 6
\put {$ \scriptstyle \bullet$} [c] at  14 12
\setlinear \plot  10 0 10 12  14 0 14 12  12 0 10 12 /
\setlinear \plot 16 0 14 6 /
\put{$5{,}040$}[c] at 13 -2
\endpicture
\end{minipage}
\begin{minipage}{4cm}
\beginpicture
\setcoordinatesystem units    <1.5mm,2mm>
\setplotarea x from 0 to 16, y from -2 to 15
\put{1.392)} [l] at 2 12
\put {$ \scriptstyle \bullet$} [c] at  10 0
\put {$ \scriptstyle \bullet$} [c] at  10 12
\put {$ \scriptstyle \bullet$} [c] at  13 0
\put {$ \scriptstyle \bullet$} [c] at  13 6
\put {$ \scriptstyle \bullet$} [c] at  13 12
\put {$ \scriptstyle \bullet$} [c] at  16  12
\put {$ \scriptstyle \bullet$} [c] at  16  0
\setlinear \plot  16 12 13 0  13 12  10 0 10 12 13 0   /
\setlinear \plot   16  0  13 6    /
\put{$5{,}040 $} [c] at 12 -2
\endpicture
\end{minipage}
$$
$$
\begin{minipage}{4cm}
\beginpicture
\setcoordinatesystem units    <1.5mm,2mm>
\setplotarea x from 0 to 16, y from -2 to 15
\put{1.393)} [l] at 2 12
\put {$ \scriptstyle \bullet$} [c] at  10 0
\put {$ \scriptstyle \bullet$} [c] at  10 12
\put {$ \scriptstyle \bullet$} [c] at  13 0
\put {$ \scriptstyle \bullet$} [c] at  13 6
\put {$ \scriptstyle \bullet$} [c] at  13 12
\put {$ \scriptstyle \bullet$} [c] at  16  12
\put {$ \scriptstyle \bullet$} [c] at  16  0
\setlinear \plot  16 0 13 12  13 0  10 12 10 0 13 12   /
\setlinear \plot   16  12  13 6    /
\put{$5{,}040 $} [c] at 12 -2
\endpicture
\end{minipage}
\begin{minipage}{4cm}
\beginpicture
\setcoordinatesystem units    <1.5mm,2mm>
\setplotarea x from 0 to 16, y from -2 to 15
\put{1.394)} [l] at 2 12
\put {$ \scriptstyle \bullet$} [c] at  10 0
\put {$ \scriptstyle \bullet$} [c] at  10 12
\put {$ \scriptstyle \bullet$} [c] at  13 12
\put {$ \scriptstyle \bullet$} [c] at  13 0
\put {$ \scriptstyle \bullet$} [c] at  16 6
\put {$ \scriptstyle \bullet$} [c] at  16 12
\put {$ \scriptstyle \bullet$} [c] at  16  0
\setlinear \plot 10 0 10 12 16 0 16 12 13 0 13 12 16 0 /
\setlinear \plot 10 12 13 0 /
\put{$5{,}040$} [c] at 12 -2
\endpicture
\end{minipage}
\begin{minipage}{4cm}
\beginpicture
\setcoordinatesystem units    <1.5mm,2mm>
\setplotarea x from 0 to 16, y from -2 to 15
\put{1.395)} [l] at 2 12
\put {$ \scriptstyle \bullet$} [c] at  10 0
\put {$ \scriptstyle \bullet$} [c] at  10 12
\put {$ \scriptstyle \bullet$} [c] at  13 12
\put {$ \scriptstyle \bullet$} [c] at  13 0
\put {$ \scriptstyle \bullet$} [c] at  16 6
\put {$ \scriptstyle \bullet$} [c] at  16 12
\put {$ \scriptstyle \bullet$} [c] at  16  0
\setlinear \plot 10 12 10 0 16 12 16 0 13 12 13 0 16 12 /
\setlinear \plot 10 0 13 12 /
\put{$5{,}040$} [c] at 12 -2
\endpicture
\end{minipage}
\begin{minipage}{4cm}
\beginpicture
\setcoordinatesystem units    <1.5mm,2mm>
\setplotarea x from 0 to 16, y from -2 to 15
\put{1.396)} [l] at 2 12
\put {$ \scriptstyle \bullet$} [c] at  10 0
\put {$ \scriptstyle \bullet$} [c] at  10 12
\put {$ \scriptstyle \bullet$} [c] at  13 12
\put {$ \scriptstyle \bullet$} [c] at  13 0
\put {$ \scriptstyle \bullet$} [c] at  16 0
\put {$ \scriptstyle \bullet$} [c] at  16 6
\put {$ \scriptstyle \bullet$} [c] at  16 12
\setlinear \plot  10  12 10 0  13 12 13 0  16 6   /
\setlinear \plot  16   0    16 12   /
\put{$5{,}040$} [c] at 13 -2
\endpicture
\end{minipage}
\begin{minipage}{4cm}
\beginpicture
\setcoordinatesystem units    <1.5mm,2mm>
\setplotarea x from 0 to 16, y from -2 to 15
\put{1.397)} [l] at 2 12
\put {$ \scriptstyle \bullet$} [c] at  10 0
\put {$ \scriptstyle \bullet$} [c] at  10 12
\put {$ \scriptstyle \bullet$} [c] at  13 12
\put {$ \scriptstyle \bullet$} [c] at  13 0
\put {$ \scriptstyle \bullet$} [c] at  16 0
\put {$ \scriptstyle \bullet$} [c] at  16 6
\put {$ \scriptstyle \bullet$} [c] at  16 12
\setlinear \plot  10  0 10 12  13 0 13 12  16 6   /
\setlinear \plot  16   0    16 12   /
\put{$5{,}040$} [c] at 13 -2
\endpicture
\end{minipage}
\begin{minipage}{4cm}
\beginpicture
\setcoordinatesystem units    <1.5mm,2mm>
\setplotarea x from 0 to 16, y from -2 to 15
\put{1.398)} [l] at 2 12
\put {$ \scriptstyle \bullet$} [c] at  10 0
\put {$ \scriptstyle \bullet$} [c] at  10 12
\put {$ \scriptstyle \bullet$} [c] at  13 12
\put {$ \scriptstyle \bullet$} [c] at  13 0
\put {$ \scriptstyle \bullet$} [c] at  16 6
\put {$ \scriptstyle \bullet$} [c] at  16  12
\put {$ \scriptstyle \bullet$} [c] at  16  0
\setlinear \plot  10 12  10 0  16 12  16 0  /
\setlinear \plot  10 0  13 12  16 6  /
\setlinear \plot  13 0 13  12    /
\put{$5{,}040 $} [c] at 13 -2
\endpicture
\end{minipage}
$$
$$
\begin{minipage}{4cm}
\beginpicture
\setcoordinatesystem units    <1.5mm,2mm>
\setplotarea x from 0 to 16, y from -2 to 15
\put{1.399)} [l] at 2 12
\put {$ \scriptstyle \bullet$} [c] at  10 0
\put {$ \scriptstyle \bullet$} [c] at  10 12
\put {$ \scriptstyle \bullet$} [c] at  13 12
\put {$ \scriptstyle \bullet$} [c] at  13 0
\put {$ \scriptstyle \bullet$} [c] at  16 6
\put {$ \scriptstyle \bullet$} [c] at  16  12
\put {$ \scriptstyle \bullet$} [c] at  16  0
\setlinear \plot  10 0  10 12  16 0  16 12  /
\setlinear \plot  10 12  13 0  16 6  /
\setlinear \plot  13 0 13  12    /
\put{$5{,}040 $} [c] at 13 -2
\endpicture
\end{minipage}
\begin{minipage}{4cm}
\beginpicture
\setcoordinatesystem units    <1.5mm,2mm>
\setplotarea x from 0 to 16, y from -2 to 15
\put{1.400)} [l] at 2 12
\put {$ \scriptstyle \bullet$} [c] at  10 0
\put {$ \scriptstyle \bullet$} [c] at  10 12
\put {$ \scriptstyle \bullet$} [c] at  13 12
\put {$ \scriptstyle \bullet$} [c] at  13 0
\put {$ \scriptstyle \bullet$} [c] at  16  6
\put {$ \scriptstyle \bullet$} [c] at  16  12
\put {$ \scriptstyle \bullet$} [c] at  16  0
\setlinear \plot 13 0 10 12 10 0 13  12 16  0 16 12 13 0  /
\put{$5{,}040$} [c] at 13 -2
\endpicture
\end{minipage}
\begin{minipage}{4cm}
\beginpicture
\setcoordinatesystem units    <1.5mm,2mm>
\setplotarea x from 0 to 16, y from -2 to 15
\put{1.401)} [l] at 2 12
\put {$ \scriptstyle \bullet$} [c] at  10 0
\put {$ \scriptstyle \bullet$} [c] at  10 6
\put {$ \scriptstyle \bullet$} [c] at  10 12
\put {$ \scriptstyle \bullet$} [c] at  12 3
\put {$ \scriptstyle \bullet$} [c] at  14 0
\put {$ \scriptstyle \bullet$} [c] at  14 12
\put {$ \scriptstyle \bullet$} [c] at  16 12
\setlinear \plot  10 0 10 12    /
\setlinear \plot 10 6 14 0 14 12  /
\setlinear \plot  14 0 16 12      /
\put{$2{,}520$} [c] at 13 -2
\endpicture
\end{minipage}
\begin{minipage}{4cm}
\beginpicture
\setcoordinatesystem units    <1.5mm,2mm>
\setplotarea x from 0 to 16, y from -2 to 15
\put{1.402)} [l] at 2 12
\put {$ \scriptstyle \bullet$} [c] at  10 0
\put {$ \scriptstyle \bullet$} [c] at  10 6
\put {$ \scriptstyle \bullet$} [c] at  10 12
\put {$ \scriptstyle \bullet$} [c] at  12 9
\put {$ \scriptstyle \bullet$} [c] at  14 0
\put {$ \scriptstyle \bullet$} [c] at  14 12
\put {$ \scriptstyle \bullet$} [c] at  16 0
\setlinear \plot  10 0 10 12    /
\setlinear \plot 10 6 14 12 14 0  /
\setlinear \plot  14 12 16 0      /
\put{$2{,}520$} [c] at 13 -2
\endpicture
\end{minipage}
\begin{minipage}{4cm}
\beginpicture
\setcoordinatesystem units    <1.5mm,2mm>
\setplotarea x from 0 to 16, y from -2 to 15
\put{1.403)} [l] at 2 12
\put {$ \scriptstyle \bullet$} [c] at  10 12
\put {$ \scriptstyle \bullet$} [c] at  12.5 6
\put {$ \scriptstyle \bullet$} [c] at  11 12
\put {$ \scriptstyle \bullet$} [c] at  13 0
\put {$ \scriptstyle \bullet$} [c] at  13 12
\put {$ \scriptstyle \bullet$} [c] at  13.5 6
\put {$ \scriptstyle \bullet$} [c] at  16 0
\setlinear \plot 16 0 13 12 13.5 6 13 0  12.5 6 13 12   /
\setlinear \plot  10 12 12.5 6 11 12  /
\put{$2{,}520$} [c] at 13 -2
\endpicture
\end{minipage}
\begin{minipage}{4cm}
\beginpicture
\setcoordinatesystem units    <1.5mm,2mm>
\setplotarea x from 0 to 16, y from -2 to 15
\put{1.404)} [l] at 2 12
\put {$ \scriptstyle \bullet$} [c] at  10 0
\put {$ \scriptstyle \bullet$} [c] at  12.5 6
\put {$ \scriptstyle \bullet$} [c] at  11 0
\put {$ \scriptstyle \bullet$} [c] at  13 0
\put {$ \scriptstyle \bullet$} [c] at  13 12
\put {$ \scriptstyle \bullet$} [c] at  13.5 6
\put {$ \scriptstyle \bullet$} [c] at  16 12
\setlinear \plot 16 12 13 0 13.5 6 13 12  12.5 6 13 0   /
\setlinear \plot  10 0 12.5 6 11 0  /

\put{$2{,}520$} [c] at 13 -2
\endpicture
\end{minipage}
$$
$$
\begin{minipage}{4cm}
\beginpicture
\setcoordinatesystem units    <1.5mm,2mm>
\setplotarea x from 0 to 16, y from -2 to 15
\put{1.405)} [l] at 2 12
\put {$ \scriptstyle \bullet$} [c] at  10 0
\put {$ \scriptstyle \bullet$} [c] at  12.5 12
\put {$ \scriptstyle \bullet$} [c] at  12.5 6
\put {$ \scriptstyle \bullet$} [c] at  13 0
\put {$ \scriptstyle \bullet$} [c] at  13.5 12
\put {$ \scriptstyle \bullet$} [c] at  13.5 6
\put {$ \scriptstyle \bullet$} [c] at  16 12
\setlinear \plot 16 12 13 0 13.5  6  13.5 12  /
\setlinear \plot  13 0 12.5 6 12.5 12 13.5 6 13.5 12 12.5 6 /
\setlinear \plot  10 0 12.5 12 /
\put{$2{,}520$} [c] at 13 -2
\endpicture
\end{minipage}
\begin{minipage}{4cm}
\beginpicture
\setcoordinatesystem units    <1.5mm,2mm>
\setplotarea x from 0 to 16, y from -2 to 15
\put{1.406)} [l] at 2 12
\put {$ \scriptstyle \bullet$} [c] at  10 12
\put {$ \scriptstyle \bullet$} [c] at  12.5 0
\put {$ \scriptstyle \bullet$} [c] at  12.5 6
\put {$ \scriptstyle \bullet$} [c] at  13 12
\put {$ \scriptstyle \bullet$} [c] at  13.5 0
\put {$ \scriptstyle \bullet$} [c] at  13.5 6
\put {$ \scriptstyle \bullet$} [c] at  16 0
\setlinear \plot 16 0 13 12 13.5  6  13.5 0  /
\setlinear \plot  13 12 12.5 6 12.5 0 13.5 6 13.5 0 12.5 6 /
\setlinear \plot  10 12 12.5 0 /
\put{$2{,}520$} [c] at 13 -2
\endpicture
\end{minipage}
\begin{minipage}{4cm}
\beginpicture
\setcoordinatesystem units    <1.5mm,2mm>
\setplotarea x from 0 to 16, y from -2 to 15
\put {1.407)} [l] at 2 12
\put {$ \scriptstyle \bullet$} [c] at  10 0
\put {$ \scriptstyle \bullet$} [c] at  10 12
\put {$ \scriptstyle \bullet$} [c] at  12 12
\put {$ \scriptstyle \bullet$} [c] at  14 12
\put {$ \scriptstyle \bullet$} [c] at  16 12
\put {$ \scriptstyle \bullet$} [c] at  14  6
\put {$ \scriptstyle \bullet$} [c] at  14  0
\setlinear \plot  10  0 10 12 14 0 16 12  /
\setlinear \plot  10 0 12 12 14 6 14 12 10 0 /
\setlinear \plot 14 0 14 6 /
\put{$2{,}520$}[c] at 13 -2
\endpicture
\end{minipage}
\begin{minipage}{4cm}
\beginpicture
\setcoordinatesystem units    <1.5mm,2mm>
\setplotarea x from 0 to 16, y from -2 to 15
\put {1.408)} [l] at 2 12
\put {$ \scriptstyle \bullet$} [c] at  10 0
\put {$ \scriptstyle \bullet$} [c] at  10 12
\put {$ \scriptstyle \bullet$} [c] at  12 0
\put {$ \scriptstyle \bullet$} [c] at  14 12
\put {$ \scriptstyle \bullet$} [c] at  16 0
\put {$ \scriptstyle \bullet$} [c] at  14  6
\put {$ \scriptstyle \bullet$} [c] at  14  0
\setlinear \plot  10  12 10 0 14 12 16 0  /
\setlinear \plot  10 12 12 0 14 6 14 0 10 12 /
\setlinear \plot 14 12 14 6 /
\put{$2{,}520$}[c] at 13 -2
\endpicture
\end{minipage}
\begin{minipage}{4cm}
\beginpicture
\setcoordinatesystem units    <1.5mm,2mm>
\setplotarea x from 0 to 16, y from -2 to 15
\put {1.409)} [l] at 2 12
\put {$ \scriptstyle \bullet$} [c] at  10 0
\put {$ \scriptstyle \bullet$} [c] at  10 12
\put {$ \scriptstyle \bullet$} [c] at  12 12
\put {$ \scriptstyle \bullet$} [c] at  14 0
\put {$ \scriptstyle \bullet$} [c] at  14 6
\put {$ \scriptstyle \bullet$} [c] at  14 12
\put {$ \scriptstyle \bullet$} [c] at  16 12
\setlinear \plot  10 12 10 0 12 12 14 6 16 12 /
\setlinear \plot  14 0 14 12 /
\put{$2{,}520$}[c] at 13 -2
\endpicture
\end{minipage}
\begin{minipage}{4cm}
\beginpicture
\setcoordinatesystem units    <1.5mm,2mm>
\setplotarea x from 0 to 16, y from -2 to 15
\put {1.410)} [l] at 2 12
\put {$ \scriptstyle \bullet$} [c] at  10 0
\put {$ \scriptstyle \bullet$} [c] at  10 12
\put {$ \scriptstyle \bullet$} [c] at  12 0
\put {$ \scriptstyle \bullet$} [c] at  14 0
\put {$ \scriptstyle \bullet$} [c] at  14 6
\put {$ \scriptstyle \bullet$} [c] at  14 12
\put {$ \scriptstyle \bullet$} [c] at  16 0
\setlinear \plot  10 0 10 12 12 0 14 6 16 0 /
\setlinear \plot  14 0 14 12 /
\put{$2{,}520$}[c] at 13 -2
\endpicture
\end{minipage}
$$
$$
\begin{minipage}{4cm}
\beginpicture
\setcoordinatesystem units    <1.5mm,2mm>
\setplotarea x from 0 to 16, y from -2 to 15
\put {1.411)} [l] at 2 12
\put {$ \scriptstyle \bullet$} [c] at  10 0
\put {$ \scriptstyle \bullet$} [c] at  10 12
\put {$ \scriptstyle \bullet$} [c] at  12 12
\put {$ \scriptstyle \bullet$} [c] at  14 0
\put {$ \scriptstyle \bullet$} [c] at  14 6
\put {$ \scriptstyle \bullet$} [c] at  14 12
\put {$ \scriptstyle \bullet$} [c] at  16 12
\setlinear \plot  10 0 10 12 14 6 12 12 10 0  /
\setlinear \plot  14 12 14 0 16 12   /
\put{$2{,}520$}[c] at 13 -2
\endpicture
\end{minipage}
\begin{minipage}{4cm}
\beginpicture
\setcoordinatesystem units    <1.5mm,2mm>
\setplotarea x from 0 to 16, y from -2 to 15
\put {1.412)} [l] at 2 12
\put {$ \scriptstyle \bullet$} [c] at  10 0
\put {$ \scriptstyle \bullet$} [c] at  10 12
\put {$ \scriptstyle \bullet$} [c] at  12 0
\put {$ \scriptstyle \bullet$} [c] at  14 0
\put {$ \scriptstyle \bullet$} [c] at  14 6
\put {$ \scriptstyle \bullet$} [c] at  14 12
\put {$ \scriptstyle \bullet$} [c] at  16 0
\setlinear \plot  10 12 10 0 14 6 12 0 10 12  /
\setlinear \plot  14 0 14 12 16 0   /
\put{$2{,}520$}[c] at 13 -2
\endpicture
\end{minipage}
\begin{minipage}{4cm}
\beginpicture
\setcoordinatesystem units    <1.5mm,2mm>
\setplotarea x from 0 to 16, y from -2 to 15
\put{1.413)} [l] at 2 12
\put {$ \scriptstyle \bullet$} [c] at  10 12
\put {$ \scriptstyle \bullet$} [c] at  10 0
\put {$ \scriptstyle \bullet$} [c] at  13 0
\put {$ \scriptstyle \bullet$} [c] at  13 6
\put {$ \scriptstyle \bullet$} [c] at  13 12
\put {$ \scriptstyle \bullet$} [c] at  16 12
\put {$ \scriptstyle \bullet$} [c] at  16 0
\setlinear \plot  10  12 10 0  13 6 16 0  16 12  /
\setlinear \plot  13 0 13 12    /
\put{$2{,}520$} [c] at 13 -2
\endpicture
\end{minipage}
\begin{minipage}{4cm}
\beginpicture
\setcoordinatesystem units    <1.5mm,2mm>
\setplotarea x from 0 to 16, y from -2 to 15
\put{1.414)} [l] at 2 12
\put {$ \scriptstyle \bullet$} [c] at  10 12
\put {$ \scriptstyle \bullet$} [c] at  10 0
\put {$ \scriptstyle \bullet$} [c] at  13 0
\put {$ \scriptstyle \bullet$} [c] at  13 6
\put {$ \scriptstyle \bullet$} [c] at  13 12
\put {$ \scriptstyle \bullet$} [c] at  16 12
\put {$ \scriptstyle \bullet$} [c] at  16 0
\setlinear \plot  10  0 10 12  13 6 16 12  16 0  /
\setlinear \plot  13 0 13 12    /
\put{$2{,}520$} [c] at 13 -2
\endpicture
\end{minipage}
\begin{minipage}{4cm}
\beginpicture
\setcoordinatesystem units    <1.5mm,2mm>
\setplotarea x from 0 to 16, y from -2 to 15
\put{1.415)} [l] at 2 12
\put {$ \scriptstyle \bullet$} [c] at  10 0
\put {$ \scriptstyle \bullet$} [c] at  10 12
\put {$ \scriptstyle \bullet$} [c] at  13 0
\put {$ \scriptstyle \bullet$} [c] at  13 12
\put {$ \scriptstyle \bullet$} [c] at  16 0
\put {$ \scriptstyle \bullet$} [c] at  16  6
\put {$ \scriptstyle \bullet$} [c] at  16  12
\setlinear \plot 10 12 10 0  13 12 16 0 16 12    /
\setlinear \plot   10 12 13 0  13 12    /
\put{$2{,}520$} [c] at 13 -2
\endpicture
\end{minipage}
\begin{minipage}{4cm}
\beginpicture
\setcoordinatesystem units    <1.5mm,2mm>
\setplotarea x from 0 to 16, y from -2 to 15
\put{1.416)} [l] at 2 12
\put {$ \scriptstyle \bullet$} [c] at  10 0
\put {$ \scriptstyle \bullet$} [c] at  10 12
\put {$ \scriptstyle \bullet$} [c] at  13 0
\put {$ \scriptstyle \bullet$} [c] at  13 12
\put {$ \scriptstyle \bullet$} [c] at  16 0
\put {$ \scriptstyle \bullet$} [c] at  16  6
\put {$ \scriptstyle \bullet$} [c] at  16  12
\setlinear \plot 10 0 10 12  13 0 16 12 16 0    /
\setlinear \plot   10 0 13 12  13 0    /
\put{$2{,}520$} [c] at 13 -2
\endpicture
\end{minipage}
$$
$$
\begin{minipage}{4cm}
\beginpicture
\setcoordinatesystem units    <1.5mm,2mm>
\setplotarea x from 0 to 16, y from -2 to 15
\put{1.417)} [l] at 2 12
\put {$ \scriptstyle \bullet$} [c] at  10 0
\put {$ \scriptstyle \bullet$} [c] at  10 12
\put {$ \scriptstyle \bullet$} [c] at  13 0
\put {$ \scriptstyle \bullet$} [c] at  13 6
\put {$ \scriptstyle \bullet$} [c] at  13 12
\put {$ \scriptstyle \bullet$} [c] at  16  0
\put {$ \scriptstyle \bullet$} [c] at  16  12
\setlinear \plot  10 12 10 0 13 12 16 0 16 12 13 0 10 12  /
\setlinear \plot  13 0 13 12 /
\put{$2{,}520 $} [c] at 13 -2
\endpicture
\end{minipage}
\begin{minipage}{4cm}
\beginpicture
\setcoordinatesystem units    <1.5mm,2mm>
\setplotarea x from 0 to 16, y from -2 to 15
\put{1.418)} [l] at 2 12
\put {$ \scriptstyle \bullet$} [c] at 11  0
\put {$ \scriptstyle \bullet$} [c] at 10 6
\put {$ \scriptstyle \bullet$} [c] at 12 6
\put {$ \scriptstyle \bullet$} [c] at 11 12
\put {$ \scriptstyle \bullet$} [c] at 15 6
\put {$ \scriptstyle \bullet$} [c] at 14 12
\put {$ \scriptstyle \bullet$} [c] at 16 12
\setlinear \plot 16 12 15 6 11 0 10 6 11 12 12 6 11 0    /
\setlinear \plot 14 12 15 6     /
\put{$1{,}260$} [c] at 13 -2
\endpicture
\end{minipage}
\begin{minipage}{4cm}
\beginpicture
\setcoordinatesystem units    <1.5mm,2mm>
\setplotarea x from 0 to 16, y from -2 to 15
\put{1.419)} [l] at 2 12
\put {$ \scriptstyle \bullet$} [c] at 11  0
\put {$ \scriptstyle \bullet$} [c] at 10 6
\put {$ \scriptstyle \bullet$} [c] at 12 6
\put {$ \scriptstyle \bullet$} [c] at 11 12
\put {$ \scriptstyle \bullet$} [c] at 15 6
\put {$ \scriptstyle \bullet$} [c] at 14 0
\put {$ \scriptstyle \bullet$} [c] at 16 0
\setlinear \plot 16 0 15 6 11 12 10 6 11 0 12 6 11 12    /
\setlinear \plot 14 0 15 6     /
\put{$1{,}260$} [c] at 13 -2
\endpicture
\end{minipage}
\begin{minipage}{4cm}
\beginpicture
\setcoordinatesystem units    <1.5mm,2mm>
\setplotarea x from 0 to 16, y from -2 to 15
\put{1.420)} [l] at 2 12
\put {$ \scriptstyle \bullet$} [c] at  16 12
\put {$ \scriptstyle \bullet$} [c] at  14 12
\put {$ \scriptstyle \bullet$} [c] at  12 6
\put {$ \scriptstyle \bullet$} [c] at  12 0
\put {$ \scriptstyle \bullet$} [c] at  11 12
\put {$ \scriptstyle \bullet$} [c] at  10 6
\put {$ \scriptstyle \bullet$} [c] at  10 0
\setlinear \plot 12 0 10 6 10 0 12 6 12 0 10 6 11 12 12 6   /
\setlinear \plot 16 12 12 0 14 12 /
\put{$1{,}260$} [c] at 13 -2
\endpicture
\end{minipage}
\begin{minipage}{4cm}
\beginpicture
\setcoordinatesystem units    <1.5mm,2mm>
\setplotarea x from 0 to 16, y from -2 to 15
\put{1.421)} [l] at 2 12
\put {$ \scriptstyle \bullet$} [c] at  16 0
\put {$ \scriptstyle \bullet$} [c] at  14 0
\put {$ \scriptstyle \bullet$} [c] at  12 6
\put {$ \scriptstyle \bullet$} [c] at  12 12
\put {$ \scriptstyle \bullet$} [c] at  11 0
\put {$ \scriptstyle \bullet$} [c] at  10 6
\put {$ \scriptstyle \bullet$} [c] at  10 12
\setlinear \plot 12 12 10 6 10 12 12 6 12 12 10 6 11 0 12 6   /
\setlinear \plot 16 0 12 12 14 0 /
\put{$1{,}260$} [c] at 13 -2
\endpicture
\end{minipage}
\begin{minipage}{4cm}
\beginpicture
\setcoordinatesystem units    <1.5mm,2mm>
\setplotarea x from 0 to 16, y from -2 to 15
\put{1.422)} [l] at 2 12
\put {$ \scriptstyle \bullet$} [c] at  16 12
\put {$ \scriptstyle \bullet$} [c] at  12.3 6
\put {$ \scriptstyle \bullet$} [c] at  13 0
\put {$ \scriptstyle \bullet$} [c] at  13 12
\put {$ \scriptstyle \bullet$} [c] at  13.7 6
\put {$ \scriptstyle \bullet$} [c] at  10 0
\put {$ \scriptstyle \bullet$} [c] at  10 12
\setlinear \plot 16 12 10 0 10 12 13 0 12.3 6 13 12 13.7 6 13 0 16 12 /
\setlinear \plot  10 0 13 12 /
\put{$1{,}260$} [c] at 13 -2
\endpicture
\end{minipage}
$$
$$
\begin{minipage}{4cm}
\beginpicture
\setcoordinatesystem units    <1.5mm,2mm>
\setplotarea x from 0 to 16, y from -2 to 15
\put{1.423)} [l] at 2 12
\put {$ \scriptstyle \bullet$} [c] at  16 0
\put {$ \scriptstyle \bullet$} [c] at  12.3 6
\put {$ \scriptstyle \bullet$} [c] at  13 0
\put {$ \scriptstyle \bullet$} [c] at  13 12
\put {$ \scriptstyle \bullet$} [c] at  13.7 6
\put {$ \scriptstyle \bullet$} [c] at  10 0
\put {$ \scriptstyle \bullet$} [c] at  10 12
\setlinear \plot 16 0 10 12 10 0 13 12 12.3 6 13 0 13.7 6 13 12 16 0 /
\setlinear \plot  10 12 13 0 /
\put{$1{,}260$} [c] at 13 -2
\endpicture
\end{minipage}
\begin{minipage}{4cm}
\beginpicture
\setcoordinatesystem units    <1.5mm,2mm>
\setplotarea x from 0 to 16, y from -2 to 15
\put {1.424)} [l] at 2 12
\put {$ \scriptstyle \bullet$} [c] at  10 12
\put {$ \scriptstyle \bullet$} [c] at  12 12
\put {$ \scriptstyle \bullet$} [c] at  11 0
\put {$ \scriptstyle \bullet$} [c] at  14 12
\put {$ \scriptstyle \bullet$} [c] at  16 12
\put {$ \scriptstyle \bullet$} [c] at  15  6
\put {$ \scriptstyle \bullet$} [c] at  15  0
\setlinear \plot  10 12 11 0 15 6 15 0 /
\setlinear \plot  12 12 11 0 /
\setlinear \plot  14 12 15 6 16 12 /
\put{$1{,}260$}[c] at 13 -2
\endpicture
\end{minipage}
\begin{minipage}{4cm}
\beginpicture
\setcoordinatesystem units    <1.5mm,2mm>
\setplotarea x from 0 to 16, y from -2 to 15
\put {1.425)} [l] at 2 12
\put {$ \scriptstyle \bullet$} [c] at  10 0
\put {$ \scriptstyle \bullet$} [c] at  12 0
\put {$ \scriptstyle \bullet$} [c] at  11 12
\put {$ \scriptstyle \bullet$} [c] at  14 0
\put {$ \scriptstyle \bullet$} [c] at  16 0
\put {$ \scriptstyle \bullet$} [c] at  15  6
\put {$ \scriptstyle \bullet$} [c] at  15  12
\setlinear \plot  10 0 11 12 15 6 15 12 /
\setlinear \plot  12 0 11 12 /
\setlinear \plot  14 0 15 6 16 0 /
\put{$1{,}260$}[c] at 13 -2
\endpicture
\end{minipage}
\begin{minipage}{4cm}
\beginpicture
\setcoordinatesystem units    <1.5mm,2mm>
\setplotarea x from 0 to 16, y from -2 to 15
\put{1.426)} [l] at 2 12
\put {$ \scriptstyle \bullet$} [c] at  10 0
\put {$ \scriptstyle \bullet$} [c] at  10 12
\put {$ \scriptstyle \bullet$} [c] at  13 12
\put {$ \scriptstyle \bullet$} [c] at  13 0
\put {$ \scriptstyle \bullet$} [c] at  14.5 6
\put {$ \scriptstyle \bullet$} [c] at  16  12
\put {$ \scriptstyle \bullet$} [c] at  16  0
\setlinear \plot  10 0 10 12 13 0 13 12  10 0  14.5 6  13 0 /
\setlinear \plot  16 0 16 12 14.5 6 /
\put{$1{,}260 $} [c] at 13 -2
\endpicture
\end{minipage}
\begin{minipage}{4cm}
\beginpicture
\setcoordinatesystem units    <1.5mm,2mm>
\setplotarea x from 0 to 16, y from -2 to 15
\put{1.427)} [l] at 2 12
\put {$ \scriptstyle \bullet$} [c] at  10 0
\put {$ \scriptstyle \bullet$} [c] at  10 12
\put {$ \scriptstyle \bullet$} [c] at  13 12
\put {$ \scriptstyle \bullet$} [c] at  13 0
\put {$ \scriptstyle \bullet$} [c] at  14.5 6
\put {$ \scriptstyle \bullet$} [c] at  16  12
\put {$ \scriptstyle \bullet$} [c] at  16  0
\setlinear \plot  10 12 10 0 13 12 13 0  10 12  14.5 6  13 12 /
\setlinear \plot  16 12 16 0 14.5 6 /
\put{$1{,}260 $} [c] at 13 -2
\endpicture
\end{minipage}
\begin{minipage}{4cm}
\beginpicture
\setcoordinatesystem units    <1.5mm,2mm>
\setplotarea x from 0 to 16, y from -2 to 15
\put{1.428)} [l] at 2 12
\put {$ \scriptstyle \bullet$} [c] at 13 0
\put {$ \scriptstyle \bullet$} [c] at 10 6
\put {$ \scriptstyle \bullet$} [c] at 12 6
\put {$ \scriptstyle \bullet$} [c] at 14 6
\put {$ \scriptstyle \bullet$} [c] at 16 8
\put {$ \scriptstyle \bullet$} [c] at 13 12
\put {$ \scriptstyle \bullet$} [c] at 16 4
\setlinear \plot 13 0 10 6 13 12 16 8 16 4 13 0   /
\setlinear \plot 13 0 12 6 13 12 14 6 13 0   /
\put{$840$} [c] at 13 -2
\endpicture
\end{minipage}
$$
$$
\begin{minipage}{4cm}
\beginpicture
\setcoordinatesystem units    <1.5mm,2mm>
\setplotarea x from 0 to 16, y from -2 to 15
\put{1.429)} [l] at 2 12
\put {$ \scriptstyle \bullet$} [c] at 10 12
\put {$ \scriptstyle \bullet$} [c] at 12 12
\put {$ \scriptstyle \bullet$} [c] at 14 12
\put {$ \scriptstyle \bullet$} [c] at 16 12
\put {$ \scriptstyle \bullet$} [c] at 13 6
\put {$ \scriptstyle \bullet$} [c] at 13 0
\put {$ \scriptstyle \bullet$} [c] at 15 10
\setlinear \plot 13 0 13 6 10 12   /
\setlinear \plot 12 12 13 6 14 12   /
\setlinear \plot 16 12 13 6  /
\put{$840$} [c] at 13 -2
\endpicture
\end{minipage}
\begin{minipage}{4cm}
\beginpicture
\setcoordinatesystem units    <1.5mm,2mm>
\setplotarea x from 0 to 16, y from -2 to 15
\put{1.430)} [l] at 2 12
\put {$ \scriptstyle \bullet$} [c] at 10 0
\put {$ \scriptstyle \bullet$} [c] at 12 0
\put {$ \scriptstyle \bullet$} [c] at 14 0
\put {$ \scriptstyle \bullet$} [c] at 16 0
\put {$ \scriptstyle \bullet$} [c] at 13 6
\put {$ \scriptstyle \bullet$} [c] at 13 12
\put {$ \scriptstyle \bullet$} [c] at 15 2
\setlinear \plot 13 12 13 6 10 0   /
\setlinear \plot 12 0 13 6 14 0   /
\setlinear \plot 16 0 13 6  /
\put{$840$} [c] at 13 -2
\endpicture
\end{minipage}
\begin{minipage}{4cm}
\beginpicture
\setcoordinatesystem units    <1.5mm,2mm>
\setplotarea x from 0 to 16, y from -2 to 15
\put {1.431)} [l] at  2 12
\put {$ \scriptstyle \bullet$} [c] at  10 12
\put {$ \scriptstyle \bullet$} [c] at  12 12
\put {$ \scriptstyle \bullet$} [c] at  14 12
\put {$ \scriptstyle \bullet$} [c] at  16 12
\put {$ \scriptstyle \bullet$} [c] at  10 0
\put {$ \scriptstyle \bullet$} [c] at  13 0
\put {$ \scriptstyle \bullet$} [c] at  16 0
\setlinear \plot   10 12 10 0  12 12 16 0 16 12  13 0   10 12 /
\setlinear \plot   13 0 14  12 16 0 /
\setlinear \plot   10  0 16 12 /
\put{$840$} [c] at 13 -2
\endpicture
\end{minipage}
\begin{minipage}{4cm}
\beginpicture
\setcoordinatesystem units    <1.5mm,2mm>
\setplotarea x from 0 to 16, y from -2 to 15
\put {1.432)} [l] at  2 12
\put {$ \scriptstyle \bullet$} [c] at  10 0
\put {$ \scriptstyle \bullet$} [c] at  10 12
\put {$ \scriptstyle \bullet$} [c] at  14 0
\put {$ \scriptstyle \bullet$} [c] at  16 0
\put {$ \scriptstyle \bullet$} [c] at  13 12
\put {$ \scriptstyle \bullet$} [c] at  16 12
\put {$ \scriptstyle \bullet$} [c] at  12 0
\setlinear \plot   10 0 10 12  12 0 16 12 16 0  13 12   10 0 /
\setlinear \plot   13 12 14  0 16 12 /
\setlinear \plot   10  12 16  0 /
\put{$840$} [c] at 13 -2
\endpicture
\end{minipage}
\begin{minipage}{4cm}
\beginpicture
\setcoordinatesystem units    <1.5mm,2mm>
\setplotarea x from 0 to 16, y from -2 to 15
\put{1.433)} [l] at 2 12
\put {$ \scriptstyle \bullet$} [c] at 13 0
\put {$ \scriptstyle \bullet$} [c] at 10 6
\put {$ \scriptstyle \bullet$} [c] at 12 6
\put {$ \scriptstyle \bullet$} [c] at 14 6
\put {$ \scriptstyle \bullet$} [c] at 16 6
\put {$ \scriptstyle \bullet$} [c] at 11 12
\put {$ \scriptstyle \bullet$} [c] at 15 12
\setlinear \plot 13  0 10 6 11 12 12 6 13 0     /
\setlinear \plot 13  0 14 6 15 12 16 6 13 0     /
\put{$630$} [c] at 13 -2
\endpicture
\end{minipage}
\begin{minipage}{4cm}
\beginpicture
\setcoordinatesystem units    <1.5mm,2mm>
\setplotarea x from 0 to 16, y from -2 to 15
\put{1.434)} [l] at 2 12
\put {$ \scriptstyle \bullet$} [c] at 13 12
\put {$ \scriptstyle \bullet$} [c] at 10 6
\put {$ \scriptstyle \bullet$} [c] at 12 6
\put {$ \scriptstyle \bullet$} [c] at 14 6
\put {$ \scriptstyle \bullet$} [c] at 16 6
\put {$ \scriptstyle \bullet$} [c] at 11 0
\put {$ \scriptstyle \bullet$} [c] at 15 0
\setlinear \plot 13  12 10 6 11 0 12 6 13 12     /
\setlinear \plot 13  12 14 6 15 0 16 6 13 12     /
\put{$630$} [c] at 13 -2
\endpicture
\end{minipage}
$$
$$
\begin{minipage}{4cm}
\beginpicture
\setcoordinatesystem units    <1.5mm,2mm>
\setplotarea x from 0 to 16, y from -2 to 15
\put{1.435)} [l] at 2 12
\put {$ \scriptstyle \bullet$} [c] at 10 12
\put {$ \scriptstyle \bullet$} [c] at 12 12
\put {$ \scriptstyle \bullet$} [c] at 14 12
\put {$ \scriptstyle \bullet$} [c] at 16 12
\put {$ \scriptstyle \bullet$} [c] at 11 6
\put {$ \scriptstyle \bullet$} [c] at 15 6
\put {$ \scriptstyle \bullet$} [c] at 13 0
\setlinear \plot 10 12 11 6 13 0 15 6 16 12    /
\setlinear \plot 11 6 12  12      /
\setlinear \plot 15 6 14  12      /
\put{$630$} [c] at 13 -2
\endpicture
\end{minipage}
\begin{minipage}{4cm}
\beginpicture
\setcoordinatesystem units    <1.5mm,2mm>
\setplotarea x from 0 to 16, y from -2 to 15
\put{1.436)} [l] at 2 12
\put {$ \scriptstyle \bullet$} [c] at 10 0
\put {$ \scriptstyle \bullet$} [c] at 12 0
\put {$ \scriptstyle \bullet$} [c] at 14 0
\put {$ \scriptstyle \bullet$} [c] at 16 0
\put {$ \scriptstyle \bullet$} [c] at 11 6
\put {$ \scriptstyle \bullet$} [c] at 15 6
\put {$ \scriptstyle \bullet$} [c] at 13 12
\setlinear \plot 10 0 11 6 13 12 15 6 16 0    /
\setlinear \plot 11 6 12  0      /
\setlinear \plot 15 6 14  0      /
\put{$630$} [c] at 13 -2
\endpicture
\end{minipage}
\begin{minipage}{4cm}
\beginpicture
\setcoordinatesystem units    <1.5mm,2mm>
\setplotarea x from 0 to 16, y from -2 to 15
\put {1.437)} [l] at 2 12
\put {$ \scriptstyle \bullet$} [c] at  10 12
\put {$ \scriptstyle \bullet$} [c] at  12 12
\put {$ \scriptstyle \bullet$} [c] at  14 12
\put {$ \scriptstyle \bullet$} [c] at  16 12
\put {$ \scriptstyle \bullet$} [c] at  10 0
\put {$ \scriptstyle \bullet$} [c] at  13  0
\put {$ \scriptstyle \bullet$} [c] at  16 0
\setlinear \plot 16 12  10 0 10 12 13 0  16 12 16 0 14 12 10 0 12 12 13 0 14 12 /
\put{$630$} [c] at 13 -2
\endpicture
\end{minipage}
\begin{minipage}{4cm}
\beginpicture
\setcoordinatesystem units    <1.5mm,2mm>
\setplotarea x from 0 to 16, y from -2 to 15
\put {1.438)} [l] at 2 12
\put {$ \scriptstyle \bullet$} [c] at  10 0
\put {$ \scriptstyle \bullet$} [c] at  12 0
\put {$ \scriptstyle \bullet$} [c] at  14 0
\put {$ \scriptstyle \bullet$} [c] at  16 0
\put {$ \scriptstyle \bullet$} [c] at  10 12
\put {$ \scriptstyle \bullet$} [c] at  13  12
\put {$ \scriptstyle \bullet$} [c] at  16 12
\setlinear \plot 16 0  10 12 10 0 13 12  16 0 16 12 14 0 10 12 12 0 13 12 14 0 /
\put{$630$} [c] at 13 -2
\endpicture
\end{minipage}
\begin{minipage}{4cm}
\beginpicture
\setcoordinatesystem units    <1.5mm,2mm>
\setplotarea x from 0 to 16, y from -2 to 15
\put {1.439)} [l] at 2 12
\put {$ \scriptstyle \bullet$} [c] at  10 0
\put {$ \scriptstyle \bullet$} [c] at  10 12
\put {$ \scriptstyle \bullet$} [c] at  12 12
\put {$ \scriptstyle \bullet$} [c] at  13  0
\put {$ \scriptstyle \bullet$} [c] at  14 12
\put {$ \scriptstyle \bullet$} [c] at  16  0
\put {$ \scriptstyle \bullet$} [c] at  16  12
\setlinear \plot   10 12 10 0  16 12 16 0 14 12 13 0 16  12  /
\setlinear \plot   14 12 10  0  12 12 16 0 /
\setlinear \plot   13 0 12  12 /
\put{$420$} [c] at 13 -2
\endpicture
\end{minipage}
\begin{minipage}{4cm}
\beginpicture
\setcoordinatesystem units    <1.5mm,2mm>
\setplotarea x from 0 to 16, y from -2 to 15
\put {1.440)} [l] at 2 12
\put {$ \scriptstyle \bullet$} [c] at  10 0
\put {$ \scriptstyle \bullet$} [c] at  10 12
\put {$ \scriptstyle \bullet$} [c] at  12 0
\put {$ \scriptstyle \bullet$} [c] at  13  12
\put {$ \scriptstyle \bullet$} [c] at  14 0
\put {$ \scriptstyle \bullet$} [c] at  16  0
\put {$ \scriptstyle \bullet$} [c] at  16  12
\setlinear \plot   10 0 10 12  16 0 16 12 14 0 13 12 16  0  /
\setlinear \plot   14 0 10  12  12 0 16 12 /
\setlinear \plot   13 12 12  0 /
\put{$420$} [c] at 13 -2
\endpicture
\end{minipage}
$$
$$
\begin{minipage}{4cm}
\beginpicture
\setcoordinatesystem units    <1.5mm,2mm>
\setplotarea x from 0 to 16, y from -2 to 15
\put{1.441)} [l] at 2 12
\put {$ \scriptstyle \bullet$} [c] at 10 4
\put {$ \scriptstyle \bullet$} [c] at 10 8
\put {$ \scriptstyle \bullet$} [c] at 10 12
\put {$ \scriptstyle \bullet$} [c] at 12 0
\put {$ \scriptstyle \bullet$} [c] at 14 4
\put {$ \scriptstyle \bullet$} [c] at 14 12
\setlinear \plot 10 12 10 4 12 0 14 4 14 12  /
\put{$5{,}040$} [c] at 11 -2
\put{$\scriptstyle \bullet$} [c] at 16  0
\endpicture
\end{minipage}
\begin{minipage}{4cm}
\beginpicture
\setcoordinatesystem units    <1.5mm,2mm>
\setplotarea x from 0 to 16, y from -2 to 15
\put{1.442)} [l] at 2 12
\put {$ \scriptstyle \bullet$} [c] at 10 0
\put {$ \scriptstyle \bullet$} [c] at 10 4
\put {$ \scriptstyle \bullet$} [c] at 10 8
\put {$ \scriptstyle \bullet$} [c] at 12 12
\put {$ \scriptstyle \bullet$} [c] at 14 0
\put {$ \scriptstyle \bullet$} [c] at 14 8
\setlinear \plot 10 0 10 8 12 12 14 8 14 0  /
\put{$5{,}040$} [c] at 11 -2
\put{$\scriptstyle \bullet$} [c] at 16  0
\endpicture
\end{minipage}
\begin{minipage}{4cm}
\beginpicture
\setcoordinatesystem units    <1.5mm,2mm>
\setplotarea x from 0 to 16, y from -2 to 15
\put{1.443)} [l] at 2 12
\put {$ \scriptstyle \bullet$} [c] at 10 6
\put {$ \scriptstyle \bullet$} [c] at 10 12
\put {$ \scriptstyle \bullet$} [c] at 11 9
\put {$ \scriptstyle \bullet$} [c] at 12 0
\put {$ \scriptstyle \bullet$} [c] at 12 12
\put {$ \scriptstyle \bullet$} [c] at 14 6
\setlinear \plot 10 12 10 6 12 0 14 6 12 12 10 6   /
\put{$5{,}040$} [c] at 11 -2
\put{$\scriptstyle \bullet$} [c] at 16  0
\endpicture
\end{minipage}
\begin{minipage}{4cm}
\beginpicture
\setcoordinatesystem units    <1.5mm,2mm>
\setplotarea x from 0 to 16, y from -2 to 15
\put{1.444)} [l] at 2 12
\put {$ \scriptstyle \bullet$} [c] at 10 6
\put {$ \scriptstyle \bullet$} [c] at 10 0
\put {$ \scriptstyle \bullet$} [c] at 11 3
\put {$ \scriptstyle \bullet$} [c] at 12 0
\put {$ \scriptstyle \bullet$} [c] at 12 12
\put {$ \scriptstyle \bullet$} [c] at 14 6
\setlinear \plot 10 0 10 6 12 12 14 6 12 0 10 6   /
\put{$5{,}040$} [c] at 11 -2
\put{$\scriptstyle \bullet$} [c] at 16  0
\endpicture
\end{minipage}
\begin{minipage}{4cm}
\beginpicture
\setcoordinatesystem units    <1.5mm,2mm>
\setplotarea x from 0 to 16, y from -2 to 15
\put{1.445)} [l] at 2 12
\put {$ \scriptstyle \bullet$} [c] at 10 0
\put {$ \scriptstyle \bullet$} [c] at 10 6
\put {$ \scriptstyle \bullet$} [c] at 10 12
\put {$ \scriptstyle \bullet$} [c] at 14 0
\put {$ \scriptstyle \bullet$} [c] at 14 6
\put {$ \scriptstyle \bullet$} [c] at 14 12
\setlinear \plot 10 12 10 0 14 12 14 0 10 6  /
\put{$5{,}040$} [c] at 11 -2
\put{$\scriptstyle \bullet$} [c] at 16  0
\endpicture
\end{minipage}
\begin{minipage}{4cm}
\beginpicture
\setcoordinatesystem units    <1.5mm,2mm>
\setplotarea x from 0 to 16, y from -2 to 15
\put{1.446)} [l] at 2 12
\put {$ \scriptstyle \bullet$} [c] at 10 0
\put {$ \scriptstyle \bullet$} [c] at 10 6
\put {$ \scriptstyle \bullet$} [c] at 10 12
\put {$ \scriptstyle \bullet$} [c] at 14 0
\put {$ \scriptstyle \bullet$} [c] at 14 6
\put {$ \scriptstyle \bullet$} [c] at 14 12
\setlinear \plot 10 0 10 12 14 0 14 12 10 6  /
\put{$5{,}040$} [c] at 11 -2
\put{$\scriptstyle \bullet$} [c] at 16  0
\endpicture
\end{minipage}
$$
$$
\begin{minipage}{4cm}
\beginpicture
\setcoordinatesystem units    <1.5mm,2mm>
\setplotarea x from 0 to 16, y from -2 to 15
\put{1.447)} [l] at 2 12
\put {$ \scriptstyle \bullet$} [c] at 10 0
\put {$ \scriptstyle \bullet$} [c] at 10 4
\put {$ \scriptstyle \bullet$} [c] at 10 8
\put {$ \scriptstyle \bullet$} [c] at 10 12
\put {$ \scriptstyle \bullet$} [c] at 14 0
\put {$ \scriptstyle \bullet$} [c] at 14 12
\setlinear \plot 10 0 10 12  /
\setlinear \plot 14  0 14  12 10 4  /
\put{$5{,}040$} [c] at 11 -2
\put{$\scriptstyle \bullet$} [c] at 16  0
\endpicture
\end{minipage}
\begin{minipage}{4cm}
\beginpicture
\setcoordinatesystem units    <1.5mm,2mm>
\setplotarea x from 0 to 16, y from -2 to 15
\put{1.448)} [l] at 2 12
\put {$ \scriptstyle \bullet$} [c] at 10 0
\put {$ \scriptstyle \bullet$} [c] at 10 4
\put {$ \scriptstyle \bullet$} [c] at 10 8
\put {$ \scriptstyle \bullet$} [c] at 10 12
\put {$ \scriptstyle \bullet$} [c] at 14 0
\put {$ \scriptstyle \bullet$} [c] at 14 12
\setlinear \plot 10 0 10 12  /
\setlinear \plot 14  12 14  0 10 8  /
\put{$5{,}040$} [c] at 11 -2
\put{$\scriptstyle \bullet$} [c] at 16  0
\endpicture
\end{minipage}
\begin{minipage}{4cm}
\beginpicture
\setcoordinatesystem units    <1.5mm,2mm>
\setplotarea x from 0 to 16, y from -2 to 15
\put{1.449)} [l] at 2 12
\put {$ \scriptstyle \bullet$} [c] at 10 0
\put {$ \scriptstyle \bullet$} [c] at 10 4
\put {$ \scriptstyle \bullet$} [c] at 10 8
\put {$ \scriptstyle \bullet$} [c] at 10 12
\put {$ \scriptstyle \bullet$} [c] at 14 0
\put {$ \scriptstyle \bullet$} [c] at 14 12
\setlinear \plot 10 12 10 0 14 12 14 0 10 12  /
\put{$5{,}040$} [c] at 11 -2
\put{$\scriptstyle \bullet$} [c] at 16  0
\endpicture
\end{minipage}
\begin{minipage}{4cm}
\beginpicture
\setcoordinatesystem units    <1.5mm,2mm>
\setplotarea x from 0 to 16, y from -2 to 15
\put{1.450)} [l] at 2 12
\put {$ \scriptstyle \bullet$} [c] at 10 0
\put {$ \scriptstyle \bullet$} [c] at 10 6
\put {$ \scriptstyle \bullet$} [c] at 12 0
\put {$ \scriptstyle \bullet$} [c] at 12 12
\put {$ \scriptstyle \bullet$} [c] at 14 6
\put {$ \scriptstyle \bullet$} [c] at 14 12
\setlinear \plot  10 0 10 6 12 12 14 6 12 0 10 6  /
\setlinear \plot 14 6  14 12   /
\put{$5{,}040$} [c] at 11 -2
\put{$\scriptstyle \bullet$} [c] at 16  0
\endpicture
\end{minipage}
\begin{minipage}{4cm}
\beginpicture
\setcoordinatesystem units    <1.5mm,2mm>
\setplotarea x from 0 to 16, y from -2 to 15
\put{1.451)} [l] at 2 12
\put {$ \scriptstyle \bullet$} [c] at 10 9
\put {$ \scriptstyle \bullet$} [c] at 11 0
\put {$ \scriptstyle \bullet$} [c] at 11 3
\put {$ \scriptstyle \bullet$} [c] at 11 12
\put {$ \scriptstyle \bullet$} [c] at 12 9
\put {$ \scriptstyle \bullet$} [c] at 14 12
\setlinear \plot 14  12 11 0 11 3 10 9  11 12 12 9 11 3 /
\put{$2{,}520$} [c] at 11 -2
\put{$\scriptstyle \bullet$} [c] at 16  0
\endpicture
\end{minipage}
\begin{minipage}{4cm}
\beginpicture
\setcoordinatesystem units    <1.5mm,2mm>
\setplotarea x from 0 to 16, y from -2 to 15
\put{1.452)} [l] at 2 12
\put {$ \scriptstyle \bullet$} [c] at 10 3
\put {$ \scriptstyle \bullet$} [c] at 11 0
\put {$ \scriptstyle \bullet$} [c] at 11 9
\put {$ \scriptstyle \bullet$} [c] at 11 12
\put {$ \scriptstyle \bullet$} [c] at 12 3
\put {$ \scriptstyle \bullet$} [c] at 14 0
\setlinear \plot 14  0 11 12 11 9 10 3  11 0 12 3 11 9 /
\put{$2{,}520$} [c] at 11 -2
\put{$\scriptstyle \bullet$} [c] at 16  0
\endpicture
\end{minipage}
$$
$$
\begin{minipage}{4cm}
\beginpicture
\setcoordinatesystem units    <1.5mm,2mm>
\setplotarea x from 0 to 16, y from -2 to 15
\put{1.453)} [l] at 2 12
\put {$ \scriptstyle \bullet$} [c] at 10 8
\put {$ \scriptstyle \bullet$} [c] at 11 0
\put {$ \scriptstyle \bullet$} [c] at 11 10
\put {$ \scriptstyle \bullet$} [c] at 11 12
\put {$ \scriptstyle \bullet$} [c] at 12 8
\put {$ \scriptstyle \bullet$} [c] at 14 12
\setlinear \plot 14  12 11 0 10 8 11 10 11 12 /
\setlinear \plot 11 10 12 8  11 0 /
\put{$2{,}520$} [c] at 11 -2
\put{$\scriptstyle \bullet$} [c] at 16  0
\endpicture
\end{minipage}
\begin{minipage}{4cm}
\beginpicture
\setcoordinatesystem units    <1.5mm,2mm>
\setplotarea x from 0 to 16, y from -2 to 15
\put{1.454)} [l] at 2 12
\put {$ \scriptstyle \bullet$} [c] at 10 4
\put {$ \scriptstyle \bullet$} [c] at 11 0
\put {$ \scriptstyle \bullet$} [c] at 11 2
\put {$ \scriptstyle \bullet$} [c] at 11 12
\put {$ \scriptstyle \bullet$} [c] at 12 4
\put {$ \scriptstyle \bullet$} [c] at 14 0
\setlinear \plot 14  0 11 12 10 4 11 2 11 0 /
\setlinear \plot 11 2 12 4  11 12 /
\put{$2{,}520$} [c] at 11 -2
\put{$\scriptstyle \bullet$} [c] at 16  0
\endpicture
\end{minipage}
\begin{minipage}{4cm}
\beginpicture
\setcoordinatesystem units    <1.5mm,2mm>
\setplotarea x from 0 to 16, y from -2 to 15
\put{1.455)} [l] at 2 12
\put {$ \scriptstyle \bullet$} [c] at 10 12
\put {$ \scriptstyle \bullet$} [c] at 12 12
\put {$ \scriptstyle \bullet$} [c] at 12 8
\put {$ \scriptstyle \bullet$} [c] at 12 4
\put {$ \scriptstyle \bullet$} [c] at 12 0
\put {$ \scriptstyle \bullet$} [c] at 14 12
\setlinear \plot 10 12 12 0 12 12 /
\setlinear \plot 12 8  14 12   /
\put{$2{,}520$} [c] at 11 -2
\put{$\scriptstyle \bullet$} [c] at 16  0
\endpicture
\end{minipage}
\begin{minipage}{4cm}
\beginpicture
\setcoordinatesystem units    <1.5mm,2mm>
\setplotarea x from 0 to 16, y from -2 to 15
\put{1.456)} [l] at 2 12
\put {$ \scriptstyle \bullet$} [c] at 10 0
\put {$ \scriptstyle \bullet$} [c] at 12 12
\put {$ \scriptstyle \bullet$} [c] at 12 8
\put {$ \scriptstyle \bullet$} [c] at 12 4
\put {$ \scriptstyle \bullet$} [c] at 12 0
\put {$ \scriptstyle \bullet$} [c] at 14 0
\setlinear \plot 10 0 12 12 12 0 /
\setlinear \plot 12 4  14 0   /
\put{$2{,}520$} [c] at 11 -2
\put{$\scriptstyle \bullet$} [c] at 16  0
\endpicture
\end{minipage}
\begin{minipage}{4cm}
\beginpicture
\setcoordinatesystem units    <1.5mm,2mm>
\setplotarea x from 0 to 16, y from -2 to 15
\put{1.457)} [l] at 2 12
\put {$ \scriptstyle \bullet$} [c] at 10 6
\put {$ \scriptstyle \bullet$} [c] at 10 12
\put {$ \scriptstyle \bullet$} [c] at 11 0
\put {$ \scriptstyle \bullet$} [c] at 12 6
\put {$ \scriptstyle \bullet$} [c] at 12 12
\put {$ \scriptstyle \bullet$} [c] at 14 12
\setlinear \plot 14 12 12 6 11 0 10 6 10 12 12 6 12 12 10 6 /
\put{$2{,}520$} [c] at 11 -2
\put{$\scriptstyle \bullet$} [c] at 16  0
\endpicture
\end{minipage}
\begin{minipage}{4cm}
\beginpicture
\setcoordinatesystem units    <1.5mm,2mm>
\setplotarea x from 0 to 16, y from -2 to 15
\put{1.458)} [l] at 2 12
\put {$ \scriptstyle \bullet$} [c] at 10 6
\put {$ \scriptstyle \bullet$} [c] at 10 0
\put {$ \scriptstyle \bullet$} [c] at 11 12
\put {$ \scriptstyle \bullet$} [c] at 12 6
\put {$ \scriptstyle \bullet$} [c] at 12 0
\put {$ \scriptstyle \bullet$} [c] at 14 0
\setlinear \plot 14 0 12 6 11 12 10 6 10 0 12 6 12 0 10 6 /
\put{$2{,}520$} [c] at 11 -2
\put{$\scriptstyle \bullet$} [c] at 16  0
\endpicture
\end{minipage}
$$
$$
\begin{minipage}{4cm}
\beginpicture
\setcoordinatesystem units    <1.5mm,2mm>
\setplotarea x from 0 to 16, y from -2 to 15
\put{1.459)} [l] at 2 12
\put {$ \scriptstyle \bullet$} [c] at 10 6
\put {$ \scriptstyle \bullet$} [c] at 10 12
\put {$ \scriptstyle \bullet$} [c] at 12 0
\put {$ \scriptstyle \bullet$} [c] at 12 6
\put {$ \scriptstyle \bullet$} [c] at 12 12
\put {$ \scriptstyle \bullet$} [c] at 14 6
\setlinear \plot 10 6 10 12 12 6 12 12  14 6 12 0 10 6 12 12 /
\setlinear \plot 12 0 12 6 /
\put{$2{,}520$} [c] at 11 -2
\put{$\scriptstyle \bullet$} [c] at 16  0
\endpicture
\end{minipage}
\begin{minipage}{4cm}
\beginpicture
\setcoordinatesystem units    <1.5mm,2mm>
\setplotarea x from 0 to 16, y from -2 to 15
\put{1.460)} [l] at 2 12
\put {$ \scriptstyle \bullet$} [c] at 10 6
\put {$ \scriptstyle \bullet$} [c] at 10 0
\put {$ \scriptstyle \bullet$} [c] at 12 0
\put {$ \scriptstyle \bullet$} [c] at 12 6
\put {$ \scriptstyle \bullet$} [c] at 12 12
\put {$ \scriptstyle \bullet$} [c] at 14 6
\setlinear \plot 10 6 10 0 12 6 12 0  14 6 12 12 10 6 12 0 /
\setlinear \plot 12 12 12 6 /
\put{$2{,}520$} [c] at 11 -2
\put{$\scriptstyle \bullet$} [c] at 16  0
\endpicture
\end{minipage}
\begin{minipage}{4cm}
\beginpicture
\setcoordinatesystem units    <1.5mm,2mm>
\setplotarea x from 0 to 16, y from -2 to 15
\put{1.461)} [l] at 2 12
\put {$ \scriptstyle \bullet$} [c] at 10 2
\put {$ \scriptstyle \bullet$} [c] at 10 12
\put {$ \scriptstyle \bullet$} [c] at 11  0
\put {$ \scriptstyle \bullet$} [c] at 12 2
\put {$ \scriptstyle \bullet$} [c] at 12 12
\put {$ \scriptstyle \bullet$} [c] at 14 0
\setlinear \plot 10 2 11 0 12 2 12 12 10 2 10 12 12  2  /
\setlinear \plot 10 12 14  0 12 12  /
\put{$1{,}260$} [c] at 11 -2
\put{$\scriptstyle \bullet$} [c] at 16  0
\endpicture
\end{minipage}
\begin{minipage}{4cm}
\beginpicture
\setcoordinatesystem units    <1.5mm,2mm>
\setplotarea x from 0 to 16, y from -2 to 15
\put{1.462)} [l] at 2 12
\put {$ \scriptstyle \bullet$} [c] at 10 0
\put {$ \scriptstyle \bullet$} [c] at 10 10
\put {$ \scriptstyle \bullet$} [c] at 11  12
\put {$ \scriptstyle \bullet$} [c] at 12 0
\put {$ \scriptstyle \bullet$} [c] at 12 10
\put {$ \scriptstyle \bullet$} [c] at 14 12
\setlinear \plot 10 10 11 12 12 10 12 0 10 10 10 0 12  10  /
\setlinear \plot 10 0 14  12 12 0  /
\put{$1{,}260$} [c] at 11 -2
\put{$\scriptstyle \bullet$} [c] at 16  0
\endpicture
\end{minipage}
\begin{minipage}{4cm}
\beginpicture
\setcoordinatesystem units    <1.5mm,2mm>
\setplotarea x from 0 to 16, y from -2 to 15
\put{1.463)} [l] at 2 12
\put {$ \scriptstyle \bullet$} [c] at 10 0
\put {$ \scriptstyle \bullet$} [c] at 10 6
\put {$ \scriptstyle \bullet$} [c] at 10 12
\put {$ \scriptstyle \bullet$} [c] at 12 12
\put {$ \scriptstyle \bullet$} [c] at 14 0
\put {$ \scriptstyle \bullet$} [c] at 14 12
\setlinear \plot  10 12 10 0 14 12 14  0 10 6 12 12  /
\put{$1{,}260$} [c] at 11 -2
\put{$\scriptstyle \bullet$} [c] at 16  0
\endpicture
\end{minipage}
\begin{minipage}{4cm}
\beginpicture
\setcoordinatesystem units    <1.5mm,2mm>
\setplotarea x from 0 to 16, y from -2 to 15
\put{1.464)} [l] at 2 12
\put {$ \scriptstyle \bullet$} [c] at 10 0
\put {$ \scriptstyle \bullet$} [c] at 10 6
\put {$ \scriptstyle \bullet$} [c] at 10 12
\put {$ \scriptstyle \bullet$} [c] at 12 0
\put {$ \scriptstyle \bullet$} [c] at 14 0
\put {$ \scriptstyle \bullet$} [c] at 14 12
\setlinear \plot  10 0 10 12 14 0 14  12 10 6 12 0  /
\put{$1{,}260$} [c] at 11 -2
\put{$\scriptstyle \bullet$} [c] at 16  0
\endpicture
\end{minipage}
$$
$$
\begin{minipage}{4cm}
\beginpicture
\setcoordinatesystem units    <1.5mm,2mm>
\setplotarea x from 0 to 16, y from -2 to 15
\put{1.465)} [l] at 2 12
\put {$ \scriptstyle \bullet$} [c] at 10 0
\put {$ \scriptstyle \bullet$} [c] at 10 12
\put {$ \scriptstyle \bullet$} [c] at 12 12
\put {$ \scriptstyle \bullet$} [c] at 14 0
\put {$ \scriptstyle \bullet$} [c] at 14 6
\put {$ \scriptstyle \bullet$} [c] at 14 12
\setlinear \plot 14 0 14 12  10 0  10 12 14 6 12 12 10 0   /
\put{$840$} [c] at 13 -2
\put{$\scriptstyle \bullet$} [c] at 16  0
\endpicture
\end{minipage}
\begin{minipage}{4cm}
\beginpicture
\setcoordinatesystem units    <1.5mm,2mm>
\setplotarea x from 0 to 16, y from -2 to 15
\put{1.466)} [l] at 2 12
\put {$ \scriptstyle \bullet$} [c] at 10 0
\put {$ \scriptstyle \bullet$} [c] at 10 12
\put {$ \scriptstyle \bullet$} [c] at 12 0
\put {$ \scriptstyle \bullet$} [c] at 14 0
\put {$ \scriptstyle \bullet$} [c] at 14 6
\put {$ \scriptstyle \bullet$} [c] at 14 12
\setlinear \plot 14 12 14 0  10 12  10 0 14 6 12 0 10 12   /
\put{$840$} [c] at 13 -2
\put{$\scriptstyle \bullet$} [c] at 16  0
\endpicture
\end{minipage}
\begin{minipage}{4cm}
\beginpicture
\setcoordinatesystem units    <1.5mm,2mm>
\setplotarea x from 0 to 16, y from -2 to 15
\put{${\bf  27}$} [l] at 2 15
\put{1.467)} [l] at 2 12
\put {$ \scriptstyle \bullet$} [c] at  10 0
\put {$ \scriptstyle \bullet$} [c] at  10 6
\put {$ \scriptstyle \bullet$} [c] at  10 12
\put {$ \scriptstyle \bullet$} [c] at  12 6
\put {$ \scriptstyle \bullet$} [c] at  14 0
\put {$ \scriptstyle \bullet$} [c] at  14 12
\put {$ \scriptstyle \bullet$} [c] at  16 12
\setlinear \plot  10 12 10 0 14 12 14 0 16 12  /
\put{$5{,}040$} [c] at 13 -2
\endpicture
\end{minipage}
\begin{minipage}{4cm}
\beginpicture
\setcoordinatesystem units    <1.5mm,2mm>
\setplotarea x from 0 to 16, y from -2 to 15
\put{1.468)} [l] at 2 12
\put {$ \scriptstyle \bullet$} [c] at  10 0
\put {$ \scriptstyle \bullet$} [c] at  10 6
\put {$ \scriptstyle \bullet$} [c] at  10 12
\put {$ \scriptstyle \bullet$} [c] at  12 6
\put {$ \scriptstyle \bullet$} [c] at  14 0
\put {$ \scriptstyle \bullet$} [c] at  14 12
\put {$ \scriptstyle \bullet$} [c] at  16 0
\setlinear \plot  10 0 10 12 14 0 14 12 16 0  /
\put{$5{,}040$} [c] at 13 -2
\endpicture
\end{minipage}
\begin{minipage}{4cm}
\beginpicture
\setcoordinatesystem units    <1.5mm,2mm>
\setplotarea x from 0 to 16, y from -2 to 15
\put{1.469)} [l] at 2 12
\put {$ \scriptstyle \bullet$} [c] at  10 0
\put {$ \scriptstyle \bullet$} [c] at  10 6
\put {$ \scriptstyle \bullet$} [c] at  10 12
\put {$ \scriptstyle \bullet$} [c] at  13 0
\put {$ \scriptstyle \bullet$} [c] at  13 6
\put {$ \scriptstyle \bullet$} [c] at  13 12
\put {$ \scriptstyle \bullet$} [c] at  16 12
\setlinear \plot  10 0 10 12 13 0 13 12  /
\setlinear \plot   13 0 16 12  /
\put{$5{,}040$} [c] at 13 -2
\endpicture
\end{minipage}
\begin{minipage}{4cm}
\beginpicture
\setcoordinatesystem units    <1.5mm,2mm>
\setplotarea x from 0 to 16, y from -2 to 15
\put{1.470)} [l] at 2 12
\put {$ \scriptstyle \bullet$} [c] at  10 0
\put {$ \scriptstyle \bullet$} [c] at  10 6
\put {$ \scriptstyle \bullet$} [c] at  10 12
\put {$ \scriptstyle \bullet$} [c] at  13 0
\put {$ \scriptstyle \bullet$} [c] at  13 6
\put {$ \scriptstyle \bullet$} [c] at  13 12
\put {$ \scriptstyle \bullet$} [c] at  16 0
\setlinear \plot  10 12 10 0 13 12 13 0  /
\setlinear \plot   13 12 16 0  /
\put{$5{,}040$} [c] at 13 -2
\endpicture
\end{minipage}
$$

$$
\begin{minipage}{4cm}
\beginpicture
\setcoordinatesystem units    <1.5mm,2mm>
\setplotarea x from 0 to 16, y from -2 to 15
\put{1.471)} [l] at 2 12
\put {$ \scriptstyle \bullet$} [c] at  10 0
\put {$ \scriptstyle \bullet$} [c] at  10 12
\put {$ \scriptstyle \bullet$} [c] at  12.5 6
\put {$ \scriptstyle \bullet$} [c] at  13 0
\put {$ \scriptstyle \bullet$} [c] at  13 12
\put {$ \scriptstyle \bullet$} [c] at  13.5 6
\put {$ \scriptstyle \bullet$} [c] at  16 12
\setlinear \plot  10 12 10 0 12.5 6 13 12 13.5 6 13 0 12.5 6  /
\setlinear \plot  13 0 16 12 /
\put{$5{,}040$} [c] at 13 -2
\endpicture
\end{minipage}
\begin{minipage}{4cm}
\beginpicture
\setcoordinatesystem units    <1.5mm,2mm>
\setplotarea x from 0 to 16, y from -2 to 15
\put{1.472)} [l] at 2 12
\put {$ \scriptstyle \bullet$} [c] at  10 0
\put {$ \scriptstyle \bullet$} [c] at  10 12
\put {$ \scriptstyle \bullet$} [c] at  12.5 6
\put {$ \scriptstyle \bullet$} [c] at  13 0
\put {$ \scriptstyle \bullet$} [c] at  13 12
\put {$ \scriptstyle \bullet$} [c] at  13.5 6
\put {$ \scriptstyle \bullet$} [c] at  16 0
\setlinear \plot  10 0 10 12 12.5 6 13 0 13.5 6 13 12 12.5 6  /
\setlinear \plot  13 12 16 0 /
\put{$5{,}040$} [c] at 13 -2
\endpicture
\end{minipage}
\begin{minipage}{4cm}
\beginpicture
\setcoordinatesystem units    <1.5mm,2mm>
\setplotarea x from 0 to 16, y from -2 to 15
\put{1.473)} [l] at 2 12
\put {$ \scriptstyle \bullet$} [c] at  10 0
\put {$ \scriptstyle \bullet$} [c] at  10 6
\put {$ \scriptstyle \bullet$} [c] at  10 12
\put {$ \scriptstyle \bullet$} [c] at  11 9
\put {$ \scriptstyle \bullet$} [c] at  14 0
\put {$ \scriptstyle \bullet$} [c] at  14 12
\put {$ \scriptstyle \bullet$} [c] at  16 12
\setlinear \plot  16 12 14 0 14 12 10 0 10 12 14 0    /
\put{$5{,}040$} [c] at 13 -2
\endpicture
\end{minipage}
\begin{minipage}{4cm}
\beginpicture
\setcoordinatesystem units    <1.5mm,2mm>
\setplotarea x from 0 to 16, y from -2 to 15
\put{1.474)} [l] at 2 12
\put {$ \scriptstyle \bullet$} [c] at  10 0
\put {$ \scriptstyle \bullet$} [c] at  10 6
\put {$ \scriptstyle \bullet$} [c] at  10 12
\put {$ \scriptstyle \bullet$} [c] at  11 3
\put {$ \scriptstyle \bullet$} [c] at  14 0
\put {$ \scriptstyle \bullet$} [c] at  14 12
\put {$ \scriptstyle \bullet$} [c] at  16 0
\setlinear \plot  16 0 14 12 14 0 10 12 10 0 14 12    /
\put{$5{,}040$} [c] at 13 -2
\endpicture
\end{minipage}
\begin{minipage}{4cm}
\beginpicture
\setcoordinatesystem units    <1.5mm,2mm>
\setplotarea x from 0 to 16, y from  -2 to  15
\put{1.475)} [l]  at 2 12
\put {$ \scriptstyle \bullet$} [c] at  10 0
\put {$ \scriptstyle \bullet$} [c] at  10 12
\put {$ \scriptstyle \bullet$} [c] at  13 12
\put {$ \scriptstyle \bullet$} [c] at  13 6
\put {$ \scriptstyle \bullet$} [c] at  13  0
\put {$ \scriptstyle \bullet$} [c] at  16  12
\put {$ \scriptstyle \bullet$} [c] at  16  0
\setlinear \plot  10  12 10 0 13 6 13 12  /
\setlinear \plot  16 0 16 12 13 0 13 6   /
\put{$5{,}040$} [c] at 13 -2
\endpicture
\end{minipage}
\begin{minipage}{4cm}
\beginpicture
\setcoordinatesystem units    <1.5mm,2mm>
\setplotarea x from 0 to 16, y from  -2 to  15
\put{1.476)} [l]  at 2 12
\put {$ \scriptstyle \bullet$} [c] at  10 0
\put {$ \scriptstyle \bullet$} [c] at  10 12
\put {$ \scriptstyle \bullet$} [c] at  13 12
\put {$ \scriptstyle \bullet$} [c] at  13 6
\put {$ \scriptstyle \bullet$} [c] at  13  0
\put {$ \scriptstyle \bullet$} [c] at  16  12
\put {$ \scriptstyle \bullet$} [c] at  16  0
\setlinear \plot  10  0 10 12 13 6 13 0  /
\setlinear \plot  16 12 16 0 13 12 13 6   /
\put{$5{,}040 $} [c] at 13 -2
\endpicture
\end{minipage}
$$
$$
\begin{minipage}{4cm}
\beginpicture
\setcoordinatesystem units    <1.5mm,2mm>
\setplotarea x from 0 to 16, y from  -2 to  15
\put{1.477)} [l]  at 2 12
\put {$ \scriptstyle \bullet$} [c] at  10 0
\put {$ \scriptstyle \bullet$} [c] at  10 12
\put {$ \scriptstyle \bullet$} [c] at  13 0
\put {$ \scriptstyle \bullet$} [c] at  13 12
\put {$ \scriptstyle \bullet$} [c] at  16 0
\put {$ \scriptstyle \bullet$} [c] at  16  6
\put {$ \scriptstyle \bullet$} [c] at  16  12
\setlinear \plot 10 12 13 0 13 12 16 0 16 12  13 0   /
\setlinear \plot 10 0 10 12 /
\put{$5{,}040 $} [c] at 13 -2
\endpicture
\end{minipage}
\begin{minipage}{4cm}
\beginpicture
\setcoordinatesystem units    <1.5mm,2mm>
\setplotarea x from 0 to 16, y from  -2 to  15
\put{1.478)} [l]  at 2 12
\put {$ \scriptstyle \bullet$} [c] at  10 0
\put {$ \scriptstyle \bullet$} [c] at  10 12
\put {$ \scriptstyle \bullet$} [c] at  13 0
\put {$ \scriptstyle \bullet$} [c] at  13 12
\put {$ \scriptstyle \bullet$} [c] at  16 0
\put {$ \scriptstyle \bullet$} [c] at  16  6
\put {$ \scriptstyle \bullet$} [c] at  16  12
\setlinear \plot 10 0 13 12 13 0 16 12 16 0  13 12   /
\setlinear \plot 10 0 10 12 /
\put{$5{,}040 $} [c] at 13 -2
\endpicture
\end{minipage}
\begin{minipage}{4cm}
\beginpicture
\setcoordinatesystem units    <1.5mm,2mm>
\setplotarea x from 0 to 16, y from -2 to 15
\put{1.479)} [l]  at 2 12
\put {$ \scriptstyle \bullet$} [c] at 10 6
\put {$ \scriptstyle \bullet$} [c] at 11.5 0
\put {$ \scriptstyle \bullet$} [c] at 11.5 6
\put {$ \scriptstyle \bullet$} [c] at 11.5 8
\put {$ \scriptstyle \bullet$} [c] at 11.5 12
\put {$ \scriptstyle \bullet$} [c] at 13 6
\put {$ \scriptstyle \bullet$} [c] at 16 12
\setlinear \plot 16 12 11.5 0 11.5 12     /
\setlinear \plot 10 6 11.5 12 13 6 11.5 0 10 6   /
\put{$2{,}520$} [c]  at 13 -2
\endpicture
\end{minipage}
\begin{minipage}{4cm}
\beginpicture
\setcoordinatesystem units    <1.5mm,2mm>
\setplotarea x from 0 to 16, y from -2 to 15
\put{1.480)} [l]  at 2 12
\put {$ \scriptstyle \bullet$} [c] at 10 6
\put {$ \scriptstyle \bullet$} [c] at 11.5 0
\put {$ \scriptstyle \bullet$} [c] at 11.5 6
\put {$ \scriptstyle \bullet$} [c] at 11.5 4
\put {$ \scriptstyle \bullet$} [c] at 11.5 12
\put {$ \scriptstyle \bullet$} [c] at 13 6
\put {$ \scriptstyle \bullet$} [c] at 16 0
\setlinear \plot 16 0 11.5 12 11.5 0     /
\setlinear \plot 10 6 11.5 0 13 6 11.5 12 10 6   /
\put{$2{,}520$} [c]  at 13 -2
\endpicture
\end{minipage}
\begin{minipage}{4cm}
\beginpicture
\setcoordinatesystem units    <1.5mm,2mm>
\setplotarea x from 0 to 16, y from -2 to 15
\put{1.481)} [l]  at 2 12
\put {$ \scriptstyle \bullet$} [c] at 10 6
\put {$ \scriptstyle \bullet$} [c] at 10 12
\put {$ \scriptstyle \bullet$} [c] at 11.5 0
\put {$ \scriptstyle \bullet$} [c] at 11.5 6
\put {$ \scriptstyle \bullet$} [c] at 13 12
\put {$ \scriptstyle \bullet$} [c] at 13 6
\put {$ \scriptstyle \bullet$} [c] at 16 12
\setlinear \plot 16 12 11.5 0 10 6 10 12 11.5 6 11.5 0     /
\setlinear \plot  11.5 0 13 6 13 12 11.5  6   /
\put{$2{,}520$} [c]  at 13 -2
\endpicture
\end{minipage}
\begin{minipage}{4cm}
\beginpicture
\setcoordinatesystem units    <1.5mm,2mm>
\setplotarea x from 0 to 16, y from -2 to 15
\put{1.482)} [l]  at 2 12
\put {$ \scriptstyle \bullet$} [c] at 10 6
\put {$ \scriptstyle \bullet$} [c] at 10 0
\put {$ \scriptstyle \bullet$} [c] at 11.5 12
\put {$ \scriptstyle \bullet$} [c] at 11.5 6
\put {$ \scriptstyle \bullet$} [c] at 13 0
\put {$ \scriptstyle \bullet$} [c] at 13 6
\put {$ \scriptstyle \bullet$} [c] at 16 0
\setlinear \plot 16 0 11.5 12 10 6 10 0 11.5 6 11.5 12     /
\setlinear \plot  11.5 12 13 6 13 0 11.5  6   /
\put{$2{,}520$} [c]  at 13 -2
\endpicture
\end{minipage}
$$

$$
\begin{minipage}{4cm}
\beginpicture
\setcoordinatesystem units    <1.5mm,2mm>
\setplotarea x from 0 to 16, y from -2 to 15
\put{1.483)} [l]  at 2 12
\put {$ \scriptstyle \bullet$} [c] at 10 12
\put {$ \scriptstyle \bullet$} [c] at 12 12
\put {$ \scriptstyle \bullet$} [c] at 14 12
\put {$ \scriptstyle \bullet$} [c] at 12 9
\put {$ \scriptstyle \bullet$} [c] at 12 6
\put {$ \scriptstyle \bullet$} [c] at 12 0
\put {$ \scriptstyle \bullet$} [c] at 16 12
\setlinear \plot 16 12 12 0 12 12     /
\setlinear \plot 10 12 12 6  14 12     /
\put{$2{,}520$} [c]  at 13 -2
\endpicture
\end{minipage}
\begin{minipage}{4cm}
\beginpicture
\setcoordinatesystem units    <1.5mm,2mm>
\setplotarea x from 0 to 16, y from -2 to 15
\put{1.484)} [l]  at 2 12
\put {$ \scriptstyle \bullet$} [c] at 10 0
\put {$ \scriptstyle \bullet$} [c] at 12 0
\put {$ \scriptstyle \bullet$} [c] at 14 0
\put {$ \scriptstyle \bullet$} [c] at 12 3
\put {$ \scriptstyle \bullet$} [c] at 12 6
\put {$ \scriptstyle \bullet$} [c] at 12 12
\put {$ \scriptstyle \bullet$} [c] at 16 0
\setlinear \plot 16 0 12 12 12 0     /
\setlinear \plot 10 0 12 6  14 0     /
\put{$2{,}520$} [c]  at 13 -2
\endpicture
\end{minipage}
\begin{minipage}{4cm}
\beginpicture
\setcoordinatesystem units    <1.5mm,2mm>
\setplotarea x from 0 to 16, y from -2 to 15
\put{1.485)} [l]  at 2 12
\put {$ \scriptstyle \bullet$} [c] at 10 6
\put {$ \scriptstyle \bullet$} [c] at 10 12
\put {$ \scriptstyle \bullet$} [c] at 11.5 12
\put {$ \scriptstyle \bullet$} [c] at 11.5 0
\put {$ \scriptstyle \bullet$} [c] at 13 6
\put {$ \scriptstyle \bullet$} [c] at 13 12
\put {$ \scriptstyle \bullet$} [c] at 16 12
\setlinear \plot 10  12 10 6 11.5 12 13 6 11.5 0 10 6    /
\setlinear \plot 13  6 13  12 /
\setlinear \plot 11.5 0 16  12 /
\put{$2{,}520$} [c]  at 13 -2
\endpicture
\end{minipage}
\begin{minipage}{4cm}
\beginpicture
\setcoordinatesystem units    <1.5mm,2mm>
\setplotarea x from 0 to 16, y from -2 to 15
\put{1.486)} [l]  at 2 12
\put {$ \scriptstyle \bullet$} [c] at 10 6
\put {$ \scriptstyle \bullet$} [c] at 10 0
\put {$ \scriptstyle \bullet$} [c] at 11.5 12
\put {$ \scriptstyle \bullet$} [c] at 11.5 0
\put {$ \scriptstyle \bullet$} [c] at 13 6
\put {$ \scriptstyle \bullet$} [c] at 13 0
\put {$ \scriptstyle \bullet$} [c] at 16 0
\setlinear \plot 10  0 10 6 11.5 0 13 6 11.5 12 10 6    /
\setlinear \plot 13  6 13  0 /
\setlinear \plot 11.5 12 16  0 /
\put{$2{,}520$} [c]  at 13 -2
\endpicture
\end{minipage}
\begin{minipage}{4cm}
\beginpicture
\setcoordinatesystem units    <1.5mm,2mm>
\setplotarea x from 0 to 16, y from  -2 to  15
\put{1.487)} [l]  at 2 12
\put {$ \scriptstyle \bullet$} [c] at  10 6
\put {$ \scriptstyle \bullet$} [c] at  11 0
\put {$ \scriptstyle \bullet$} [c] at  11 12
\put {$ \scriptstyle \bullet$} [c] at  12 6
\put {$ \scriptstyle \bullet$} [c] at  13.5 6
\put {$ \scriptstyle \bullet$} [c] at  16 12
\put {$ \scriptstyle \bullet$} [c] at  16 0
\setlinear \plot 16 0 16 12  11  0 10 6 11 12 12 6 11 0   /
\put{$2{,}520  $} [c] at 13 -2
\endpicture
\end{minipage}
\begin{minipage}{4cm}
\beginpicture
\setcoordinatesystem units    <1.5mm,2mm>
\setplotarea x from 0 to 16, y from  -2 to  15
\put{1.488)} [l]  at 2 12
\put {$ \scriptstyle \bullet$} [c] at  10 6
\put {$ \scriptstyle \bullet$} [c] at  11 0
\put {$ \scriptstyle \bullet$} [c] at  11 12
\put {$ \scriptstyle \bullet$} [c] at  12 6
\put {$ \scriptstyle \bullet$} [c] at  13.5 6
\put {$ \scriptstyle \bullet$} [c] at  16 12
\put {$ \scriptstyle \bullet$} [c] at  16 0
\setlinear \plot 16 12 16 0  11  12 10 6 11 0 12 6 11 12   /
\put{$2{,}520  $} [c] at 13 -2
\endpicture
\end{minipage}
$$
$$
\begin{minipage}{4cm}
\beginpicture
\setcoordinatesystem units    <1.5mm,2mm>
\setplotarea x from 0 to 16, y from -2 to 15
\put{1.489)} [l] at 2 12
\put {$ \scriptstyle \bullet$} [c] at  10 0
\put {$ \scriptstyle \bullet$} [c] at  10 6
\put {$ \scriptstyle \bullet$} [c] at  10 12
\put {$ \scriptstyle \bullet$} [c] at  13 12
\put {$ \scriptstyle \bullet$} [c] at  13 6
\put {$ \scriptstyle \bullet$} [c] at  16 0
\put {$ \scriptstyle \bullet$} [c] at  16 12
\setlinear \plot  10 12 10 0 16 12 16 0  /
\setlinear \plot  13 12 10 6  /
\put{$2{,}520$} [c] at 13 -2
\endpicture
\end{minipage}
\begin{minipage}{4cm}
\beginpicture
\setcoordinatesystem units    <1.5mm,2mm>
\setplotarea x from 0 to 16, y from -2 to 15
\put{1.490)} [l] at 2 12
\put {$ \scriptstyle \bullet$} [c] at  10 0
\put {$ \scriptstyle \bullet$} [c] at  10 6
\put {$ \scriptstyle \bullet$} [c] at  10 12
\put {$ \scriptstyle \bullet$} [c] at  13 0
\put {$ \scriptstyle \bullet$} [c] at  13 6
\put {$ \scriptstyle \bullet$} [c] at  16 0
\put {$ \scriptstyle \bullet$} [c] at  16 12
\setlinear \plot  10 0 10 12 16 0 16 12  /
\setlinear \plot  13 0 10 6  /
\put{$2{,}520$} [c] at 13 -2
\endpicture
\end{minipage}
\begin{minipage}{4cm}
\beginpicture
\setcoordinatesystem units    <1.5mm,2mm>
\setplotarea x from 0 to 16, y from -2 to 15
\put{1.491)} [l] at 2 12
\put {$ \scriptstyle \bullet$} [c] at  16 12
\put {$ \scriptstyle \bullet$} [c] at  14 0
\put {$ \scriptstyle \bullet$} [c] at  14 12
\put {$ \scriptstyle \bullet$} [c] at  10 6
\put {$ \scriptstyle \bullet$} [c] at  10.5 0
\put {$ \scriptstyle \bullet$} [c] at  10.5 12
\put {$ \scriptstyle \bullet$} [c] at  11 6
\setlinear \plot 10.5 0 10 6 10.5 12 11 6 10.5 0 14 12 14 0 16 12  /
\setlinear \plot  10.5 12 14 0 /
\put{$2{,}520$} [c] at 13 -2
\endpicture
\end{minipage}
\begin{minipage}{4cm}
\beginpicture
\setcoordinatesystem units    <1.5mm,2mm>
\setplotarea x from 0 to 16, y from -2 to 15
\put{1.492)} [l] at 2 12
\put {$ \scriptstyle \bullet$} [c] at  16 0
\put {$ \scriptstyle \bullet$} [c] at  14 0
\put {$ \scriptstyle \bullet$} [c] at  14 12
\put {$ \scriptstyle \bullet$} [c] at  10 6
\put {$ \scriptstyle \bullet$} [c] at  10.5 0
\put {$ \scriptstyle \bullet$} [c] at  10.5 12
\put {$ \scriptstyle \bullet$} [c] at  11 6
\setlinear \plot 10.5 12 10 6 10.5 0 11 6 10.5 12 14 0 14 12 16 0  /
\setlinear \plot  10.5 0 14 12 /
\put{$2{,}520$} [c] at 13 -2
\endpicture
\end{minipage}
\begin{minipage}{4cm}
\beginpicture
\setcoordinatesystem units    <1.5mm,2mm>
\setplotarea x from 0 to 16, y from -2 to 15
\put{1.493)} [l] at 2 12
\put {$ \scriptstyle \bullet$} [c] at  10 0
\put {$ \scriptstyle \bullet$} [c] at  10 12
\put {$ \scriptstyle \bullet$} [c] at  12 12
\put {$ \scriptstyle \bullet$} [c] at  14 6
\put {$ \scriptstyle \bullet$} [c] at  15 0
\put {$ \scriptstyle \bullet$} [c] at  15 12
\put {$ \scriptstyle \bullet$} [c] at  16 6
\setlinear \plot  10 12 10 0 14 6 15 0  16 6 15 12 14 6 /
\setlinear \plot  10 0 12 12 /
\put{$2{,}520$} [c] at 13 -2
\endpicture
\end{minipage}
\begin{minipage}{4cm}
\beginpicture
\setcoordinatesystem units    <1.5mm,2mm>
\setplotarea x from 0 to 16, y from -2 to 15
\put{1.494)} [l] at 2 12
\put {$ \scriptstyle \bullet$} [c] at  10 0
\put {$ \scriptstyle \bullet$} [c] at  10 12
\put {$ \scriptstyle \bullet$} [c] at  12 0
\put {$ \scriptstyle \bullet$} [c] at  14 6
\put {$ \scriptstyle \bullet$} [c] at  15 0
\put {$ \scriptstyle \bullet$} [c] at  15 12
\put {$ \scriptstyle \bullet$} [c] at  16 6
\setlinear \plot  10 0 10 12 14 6 15 12  16 6 15 0 14 6 /
\setlinear \plot  10 12 12 0 /
\put{$2{,}520$} [c] at 13 -2
\endpicture
\end{minipage}
$$

$$
\begin{minipage}{4cm}
\beginpicture
\setcoordinatesystem units    <1.5mm,2mm>
\setplotarea x from 0 to 16, y from -2 to 15
\put{1.495)} [l] at 2 12
\put {$ \scriptstyle \bullet$} [c] at  10 0
\put {$ \scriptstyle \bullet$} [c] at  12 12
\put {$ \scriptstyle \bullet$} [c] at  14 0
\put {$ \scriptstyle \bullet$} [c] at  14 6
\put {$ \scriptstyle \bullet$} [c] at  14 12
\put {$ \scriptstyle \bullet$} [c] at  10 12
\put {$ \scriptstyle \bullet$} [c] at  16 12
\setlinear \plot  10 12 10 0 14 12 14 0  16 12 /
\setlinear \plot  10 0 12 12 14 6 /
\put{$2{,}520$} [c] at 13 -2
\endpicture
\end{minipage}
\begin{minipage}{4cm}
\beginpicture
\setcoordinatesystem units    <1.5mm,2mm>
\setplotarea x from 0 to 16, y from -2 to 15
\put{1.496)} [l] at 2 12
\put {$ \scriptstyle \bullet$} [c] at  10 0
\put {$ \scriptstyle \bullet$} [c] at  12 0
\put {$ \scriptstyle \bullet$} [c] at  14 0
\put {$ \scriptstyle \bullet$} [c] at  14 6
\put {$ \scriptstyle \bullet$} [c] at  14 12
\put {$ \scriptstyle \bullet$} [c] at  10 12
\put {$ \scriptstyle \bullet$} [c] at  16 0
\setlinear \plot  10 0 10 12 14 0 14 12  16 0 /
\setlinear \plot  10 12 12 0 14 6 /
\put{$2{,}520$} [c] at 13 -2
\endpicture
\end{minipage}
\begin{minipage}{4cm}
\beginpicture
\setcoordinatesystem units    <1.5mm,2mm>
\setplotarea x from 0 to 16, y from -2 to 15
\put{1.497)} [l] at 2 12
\put {$ \scriptstyle \bullet$} [c] at  10 12
\put {$ \scriptstyle \bullet$} [c] at  12 12
\put {$ \scriptstyle \bullet$} [c] at  14 12
\put {$ \scriptstyle \bullet$} [c] at  16 12
\put {$ \scriptstyle \bullet$} [c] at  10 0
\put {$ \scriptstyle \bullet$} [c] at  13 0
\put {$ \scriptstyle \bullet$} [c] at  16 0
\setlinear \plot  10 12 10 0 12 12  16 0  16 12  13 0 10 12 /
\setlinear \plot  10 12 13 0 14 12  16 0 /
\setlinear \plot  13 0  12 12 /
\put{$2{,}520$} [c] at 13 -2
\endpicture
\end{minipage}
\begin{minipage}{4cm}
\beginpicture
\setcoordinatesystem units    <1.5mm,2mm>
\setplotarea x from 0 to 16, y from -2 to 15
\put{1.498)} [l] at 2 12
\put {$ \scriptstyle \bullet$} [c] at  10 0
\put {$ \scriptstyle \bullet$} [c] at  12 0
\put {$ \scriptstyle \bullet$} [c] at  14 0
\put {$ \scriptstyle \bullet$} [c] at  16 0
\put {$ \scriptstyle \bullet$} [c] at  10 12
\put {$ \scriptstyle \bullet$} [c] at  13 12
\put {$ \scriptstyle \bullet$} [c] at  16 12
\setlinear \plot  10 0 10 12 12 0  16 12  16 0  13 12 10 0 /
\setlinear \plot  10 0 13 12  14 0  16 12 /
\setlinear \plot  13 12 12 0  /

\put{$2{,}520$} [c] at 13 -2
\endpicture
\end{minipage}
\begin{minipage}{4cm}
\beginpicture
\setcoordinatesystem units    <1.5mm,2mm>
\setplotarea x from 0 to 16, y from -2 to 15
\put{1.499)} [l] at 2 12
\put {$ \scriptstyle \bullet$} [c] at  10 0
\put {$ \scriptstyle \bullet$} [c] at  10 12
\put {$ \scriptstyle \bullet$} [c] at  12.3 6
\put {$ \scriptstyle \bullet$} [c] at  13 12
\put {$ \scriptstyle \bullet$} [c] at  13 0
\put {$ \scriptstyle \bullet$} [c] at  13.7 6
\put {$ \scriptstyle \bullet$} [c] at  16 12
\setlinear \plot 16 12 10 0  10 12 13 0 16 12 /
\setlinear \plot  13 0 12.3 6 13 12 13.7 6 13 0   /
\put{$1{,}260$} [c] at 13 -2
\endpicture
\end{minipage}
\begin{minipage}{4cm}
\beginpicture
\setcoordinatesystem units    <1.5mm,2mm>
\setplotarea x from 0 to 16, y from -2 to 15
\put{1.500)} [l] at 2 12
\put {$ \scriptstyle \bullet$} [c] at  10 0
\put {$ \scriptstyle \bullet$} [c] at  10 12
\put {$ \scriptstyle \bullet$} [c] at  12.3 6
\put {$ \scriptstyle \bullet$} [c] at  13 12
\put {$ \scriptstyle \bullet$} [c] at  13 0
\put {$ \scriptstyle \bullet$} [c] at  13.7 6
\put {$ \scriptstyle \bullet$} [c] at  16 0
\setlinear \plot 16 0 10 12  10 0 13 12 16 0 /
\setlinear \plot  13 0 12.3 6 13 12 13.7 6 13 0   /
\put{$1{,}260$} [c] at 13 -2
\endpicture
\end{minipage}
$$
$$
\begin{minipage}{4cm}
\beginpicture
\setcoordinatesystem units    <1.5mm,2mm>
\setplotarea x from 0 to 16, y from -2 to 15
\put{1.501)} [l] at 2 12
\put {$ \scriptstyle \bullet$} [c] at  10 0
\put {$ \scriptstyle \bullet$} [c] at  10 12
\put {$ \scriptstyle \bullet$} [c] at  10 6
\put {$ \scriptstyle \bullet$} [c] at  12 12
\put {$ \scriptstyle \bullet$} [c] at  14 0
\put {$ \scriptstyle \bullet$} [c] at  14 12
\put {$ \scriptstyle \bullet$} [c] at  16 12
\setlinear \plot 10 0 10 12 14 0  16 12  /
\setlinear \plot  10 6 12 12 14 0 14 12   /
\put{$1{,}260$} [c] at 13 -2
\endpicture
\end{minipage}
\begin{minipage}{4cm}
\beginpicture
\setcoordinatesystem units    <1.5mm,2mm>
\setplotarea x from 0 to 16, y from -2 to 15
\put{1.502)} [l] at 2 12
\put {$ \scriptstyle \bullet$} [c] at  10 0
\put {$ \scriptstyle \bullet$} [c] at  10 12
\put {$ \scriptstyle \bullet$} [c] at  10 6
\put {$ \scriptstyle \bullet$} [c] at  12 0
\put {$ \scriptstyle \bullet$} [c] at  14 12
\put {$ \scriptstyle \bullet$} [c] at  14 0
\put {$ \scriptstyle \bullet$} [c] at  16 0
\setlinear \plot 10 12 10 0 14 12  16 0  /
\setlinear \plot  10 6 12 0 14 12 14 0   /
\put{$1{,}260$} [c] at 13 -2
\endpicture
\end{minipage}
\begin{minipage}{4cm}
\beginpicture
\setcoordinatesystem units    <1.5mm,2mm>
\setplotarea x from 0 to 16, y from -2 to 15
\put{1.503)} [l] at 2 12
\put {$ \scriptstyle \bullet$} [c] at  10 12
\put {$ \scriptstyle \bullet$} [c] at  12 12
\put {$ \scriptstyle \bullet$} [c] at  14 12
\put {$ \scriptstyle \bullet$} [c] at  16 12
\put {$ \scriptstyle \bullet$} [c] at  10 0
\put {$ \scriptstyle \bullet$} [c] at  15 0
\put {$ \scriptstyle \bullet$} [c] at  15 6
\setlinear \plot 12 12 10 0 10 12 15 0  12 12  /
\setlinear \plot  14 12 15 6 15 0   /
\setlinear \plot  16 12 15 6  /
\put{$1{,}260$} [c] at 13 -2
\endpicture
\end{minipage}
\begin{minipage}{4cm}
\beginpicture
\setcoordinatesystem units    <1.5mm,2mm>
\setplotarea x from 0 to 16, y from -2 to 15
\put{1.504)} [l] at 2 12
\put {$ \scriptstyle \bullet$} [c] at  10 0
\put {$ \scriptstyle \bullet$} [c] at  12 0
\put {$ \scriptstyle \bullet$} [c] at  14 0
\put {$ \scriptstyle \bullet$} [c] at  16 0
\put {$ \scriptstyle \bullet$} [c] at  10 12
\put {$ \scriptstyle \bullet$} [c] at  15 12
\put {$ \scriptstyle \bullet$} [c] at  15 6
\setlinear \plot 12 0 10 12 10 0 15 12  12 0  /
\setlinear \plot  14 0 15 6 15  12   /
\setlinear \plot  16 0 15 6  /
\put{$1{,}260$} [c] at 13 -2
\endpicture
\end{minipage}
\begin{minipage}{4cm}
\beginpicture
\setcoordinatesystem units    <1.5mm,2mm>
\setplotarea x from 0 to 16, y from -2 to 15
\put{1.505)} [l] at 2 12
\put {$ \scriptstyle \bullet$} [c] at  10 12
\put {$ \scriptstyle \bullet$} [c] at  12 12
\put {$ \scriptstyle \bullet$} [c] at  14 12
\put {$ \scriptstyle \bullet$} [c] at  16 12
\put {$ \scriptstyle \bullet$} [c] at  10 0
\put {$ \scriptstyle \bullet$} [c] at  13 0
\put {$ \scriptstyle \bullet$} [c] at  16 0
\setlinear \plot  10 0 10 12 13 0  14 12 16 0 16 12 /
\setlinear \plot 12 12 10 0 14 12   /
\setlinear \plot  13 0 12 12 16 0  /
\put{$1{,}260$} [c] at 13 -2
\endpicture
\end{minipage}
\begin{minipage}{4cm}
\beginpicture
\setcoordinatesystem units    <1.5mm,2mm>
\setplotarea x from 0 to 16, y from -2 to 15
\put{1.506)} [l] at 2 12
\put {$ \scriptstyle \bullet$} [c] at  10 0
\put {$ \scriptstyle \bullet$} [c] at  12 0
\put {$ \scriptstyle \bullet$} [c] at  14 0
\put {$ \scriptstyle \bullet$} [c] at  16 0
\put {$ \scriptstyle \bullet$} [c] at  10 12
\put {$ \scriptstyle \bullet$} [c] at  13 12
\put {$ \scriptstyle \bullet$} [c] at  16 12
\setlinear \plot  10 12 10 0 13 12  14 0 16 12 16 0 /
\setlinear \plot 12 0 10 12 14 0   /
\setlinear \plot  13 12 12 0 16 12  /
\put{$1{,}260$} [c] at 13 -2
\endpicture
\end{minipage}
$$

$$
\begin{minipage}{4cm}
\beginpicture
\setcoordinatesystem units    <1.5mm,2mm>
\setplotarea x from 0 to 16, y from -2 to 15
\put{1.507)} [l]  at 2 12
\put {$ \scriptstyle \bullet$} [c] at 10 6
\put {$ \scriptstyle \bullet$} [c] at 12 6
\put {$ \scriptstyle \bullet$} [c] at 14 6
\put {$ \scriptstyle \bullet$} [c] at 16 6
\put {$ \scriptstyle \bullet$} [c] at 13 0
\put {$ \scriptstyle \bullet$} [c] at 13 12
\put {$ \scriptstyle \bullet$} [c] at 16 12
\setlinear \plot 16 12 16 6 13 0 10 6 13 12 16 6  /
\setlinear \plot 13 0 12 6 13 12 14 6 13 0  /
\put{$840$} [c]  at 13 -2
\endpicture
\end{minipage}
\begin{minipage}{4cm}
\beginpicture
\setcoordinatesystem units    <1.5mm,2mm>
\setplotarea x from 0 to 16, y from -2 to 15
\put{1.508)} [l]  at 2 12
\put {$ \scriptstyle \bullet$} [c] at 10 6
\put {$ \scriptstyle \bullet$} [c] at 12 6
\put {$ \scriptstyle \bullet$} [c] at 14 6
\put {$ \scriptstyle \bullet$} [c] at 16 6
\put {$ \scriptstyle \bullet$} [c] at 13 0
\put {$ \scriptstyle \bullet$} [c] at 13 12
\put {$ \scriptstyle \bullet$} [c] at 16 0
\setlinear \plot 16 0 16 6 13 12 10 6 13 0 16 6  /
\setlinear \plot 13 12 12 6 13 0 14 6 13 12  /
\put{$840$} [c]  at 13 -2
\endpicture
\end{minipage}
\begin{minipage}{4cm}
\beginpicture
\setcoordinatesystem units    <1.5mm,2mm>
\setplotarea x from 0 to 16, y from -2 to 15
\put{1.509)} [l]  at 2 12
\put {$ \scriptstyle \bullet$} [c] at 10 12
\put {$ \scriptstyle \bullet$} [c] at 11 12
\put {$ \scriptstyle \bullet$} [c] at 12 12
\put {$ \scriptstyle \bullet$} [c] at 12 6
\put {$ \scriptstyle \bullet$} [c] at 14 0
\put {$ \scriptstyle \bullet$} [c] at 14 12
\put {$ \scriptstyle \bullet$} [c] at 16 6
\setlinear \plot 12 6 14 12 16 6 14 0 12 6 12 12  /
\setlinear \plot 10 12 12 6 11 12   /
\put{$840$} [c]  at 13 -2
\endpicture
\end{minipage}
\begin{minipage}{4cm}
\beginpicture
\setcoordinatesystem units    <1.5mm,2mm>
\setplotarea x from 0 to 16, y from -2 to 15
\put{1.510)} [l]  at 2 12
\put {$ \scriptstyle \bullet$} [c] at 10 0
\put {$ \scriptstyle \bullet$} [c] at 11 0
\put {$ \scriptstyle \bullet$} [c] at 12 0
\put {$ \scriptstyle \bullet$} [c] at 12 6
\put {$ \scriptstyle \bullet$} [c] at 14 0
\put {$ \scriptstyle \bullet$} [c] at 14 12
\put {$ \scriptstyle \bullet$} [c] at 16 6
\setlinear \plot 12 6 14 0 16 6 14 12 12 6 12 0  /
\setlinear \plot 10 0 12 6 11 0   /
\put{$840$} [c]  at 13 -2
\endpicture
\end{minipage}
\begin{minipage}{4cm}
\beginpicture
\setcoordinatesystem units    <1.5mm,2mm>
\setplotarea x from 0 to 16, y from -2 to 15
\put{1.511)} [l]  at 2 12
\put {$ \scriptstyle \bullet$} [c] at 10 12
\put {$ \scriptstyle \bullet$} [c] at 12 0
\put {$ \scriptstyle \bullet$} [c] at 12 12
\put {$ \scriptstyle \bullet$} [c] at 14 0
\put {$ \scriptstyle \bullet$} [c] at 14 12
\put {$ \scriptstyle \bullet$} [c] at 16 12
\put {$ \scriptstyle \bullet$} [c] at 16 6
\setlinear \plot  12 0 10 12 14 0  16 6 16 12   /
\setlinear \plot 16 6 12 0 12 12 14 0 14 12 12 0  /
\put{$420$} [c]  at 13 -2
\endpicture
\end{minipage}
\begin{minipage}{4cm}
\beginpicture
\setcoordinatesystem units    <1.5mm,2mm>
\setplotarea x from 0 to 16, y from -2 to 15
\put{1.512)} [l]  at 2 12
\put {$ \scriptstyle \bullet$} [c] at 10 0
\put {$ \scriptstyle \bullet$} [c] at 12 0
\put {$ \scriptstyle \bullet$} [c] at 12 12
\put {$ \scriptstyle \bullet$} [c] at 14 0
\put {$ \scriptstyle \bullet$} [c] at 14 12
\put {$ \scriptstyle \bullet$} [c] at 16 0
\put {$ \scriptstyle \bullet$} [c] at 16 6
\setlinear \plot  12 12 10 0 14 12  16 6 16 0   /
\setlinear \plot 16 6 12 12 12 0 14 12 14 0 12 12  /
\put{$420$} [c]  at 13 -2
\endpicture
\end{minipage}
$$
$$
\begin{minipage}{4cm}
\beginpicture
\setcoordinatesystem units    <1.5mm,2mm>
\setplotarea x from 0 to 16, y  from -2 to 15
\put{1.513)} [l]  at 2 12
\put {$ \scriptstyle \bullet$} [c] at 10 12
\put {$ \scriptstyle \bullet$} [c] at 14 12
\put {$ \scriptstyle \bullet$} [c] at 12 0
\put {$ \scriptstyle \bullet$} [c] at 16 12
\put {$ \scriptstyle \bullet$} [c] at 16 0
\put {$ \scriptstyle \bullet$} [c] at 10.5 9
\put {$ \scriptstyle \bullet$} [c] at 11.1 5
\setlinear \plot  10 12  12 0   14 12   /
\setlinear \plot 16 0 16 12   /
\put{$5{,}040$} [c] at 13 -2
\endpicture
\end{minipage}
\begin{minipage}{4cm}
\beginpicture
\setcoordinatesystem units    <1.5mm,2mm>
\setplotarea x from 0 to 16, y  from -2 to 15
\put{1.514)} [l]  at 2 12
\put {$ \scriptstyle \bullet$} [c] at 10 0
\put {$ \scriptstyle \bullet$} [c] at 14 0
\put {$ \scriptstyle \bullet$} [c] at 12 12
\put {$ \scriptstyle \bullet$} [c] at 16 12
\put {$ \scriptstyle \bullet$} [c] at 16 0
\put {$ \scriptstyle \bullet$} [c] at 10.5 3
\put {$ \scriptstyle \bullet$} [c] at 11.1 7
\setlinear \plot  10 0  12 12   14 0   /
\setlinear \plot 16 0 16 12   /
\put{$5{,}040$} [c] at 13 -2
\endpicture
\end{minipage}
\begin{minipage}{4cm}
\beginpicture
\setcoordinatesystem units    <1.5mm,2mm>
\setplotarea x from 0 to 16, y  from -2 to 15
\put{1.515)} [l]  at 2 12
\put {$ \scriptstyle \bullet$} [c] at 10 12
\put {$ \scriptstyle \bullet$} [c] at 10 6
\put {$ \scriptstyle \bullet$} [c] at 12 0
\put {$ \scriptstyle \bullet$} [c] at 12 12
\put {$ \scriptstyle \bullet$} [c] at 14 6
\put {$ \scriptstyle \bullet$} [c] at 16 0
\put {$ \scriptstyle \bullet$} [c] at 16 12
\setlinear \plot 10 12 10 6   12 0 14 6 12 12 10 6  /
\setlinear \plot 16 0 16 12   /
\put{$5{,}040$} [c] at 13 -2
\endpicture
\end{minipage}
\begin{minipage}{4cm}
\beginpicture
\setcoordinatesystem units    <1.5mm,2mm>
\setplotarea x from 0 to 16, y  from -2 to 15
\put{1.516)} [l]  at 2 12
\put {$ \scriptstyle \bullet$} [c] at 10 0
\put {$ \scriptstyle \bullet$} [c] at 10 6
\put {$ \scriptstyle \bullet$} [c] at 12 0
\put {$ \scriptstyle \bullet$} [c] at 12 12
\put {$ \scriptstyle \bullet$} [c] at 14 6
\put {$ \scriptstyle \bullet$} [c] at 16 0
\put {$ \scriptstyle \bullet$} [c] at 16 12
\setlinear \plot 10 0 10 6   12 0 14 6 12 12 10 6  /
\setlinear \plot 16 0 16 12   /
\put{$5{,}040$} [c] at 13 -2
\endpicture
\end{minipage}
\begin{minipage}{4cm}
\beginpicture
\setcoordinatesystem units    <1.5mm,2mm>
\setplotarea x from 0 to 16, y  from -2 to 15
\put{1.517)} [l]  at 2 12
\put {$ \scriptstyle \bullet$} [c] at 10 0
\put {$ \scriptstyle \bullet$} [c] at 10 6
\put {$ \scriptstyle \bullet$} [c] at 10 12
\put {$ \scriptstyle \bullet$} [c] at 14 0
\put {$ \scriptstyle \bullet$} [c] at 14 12
\setlinear \plot  10 12 10  0 14 12 14 0 10 6   /
\put {$ \scriptstyle \bullet$} [c] at 16 0
\put {$ \scriptstyle \bullet$} [c] at 16 12
\setlinear \plot 16 0 16 12   /
\put{$2{,}520$} [c] at 13 -2
\endpicture
\end{minipage}
\begin{minipage}{4cm}
\beginpicture
\setcoordinatesystem units    <1.5mm,2mm>
\setplotarea x from 0 to 16, y  from -2 to 15
\put{1.518)} [l]  at 2 12
\put {$ \scriptstyle \bullet$} [c] at 10 0
\put {$ \scriptstyle \bullet$} [c] at 10 6
\put {$ \scriptstyle \bullet$} [c] at 10 12
\put {$ \scriptstyle \bullet$} [c] at 14 0
\put {$ \scriptstyle \bullet$} [c] at 14 12
\setlinear \plot  10 0 10  12 14 0 14 12 10 6   /
\put {$ \scriptstyle \bullet$} [c] at 16 0
\put {$ \scriptstyle \bullet$} [c] at 16 12
\setlinear \plot 16 0 16 12   /
\put{$2{,}520$} [c] at 13 -2
\endpicture
\end{minipage}
$$

$$
\begin{minipage}{4cm}
\beginpicture
\setcoordinatesystem units    <1.5mm,2mm>
\setplotarea x from 0 to 16, y from -2 to 15
\put{${\bf  28}$} [l] at 2 15

\put{1.519)} [l] at 2 12
\put {$ \scriptstyle \bullet$} [c] at  10 0
\put {$ \scriptstyle \bullet$} [c] at  12.5 4.5
\put {$ \scriptstyle \bullet$} [c] at  12.5 7.5
\put {$ \scriptstyle \bullet$} [c] at  13 0
\put {$ \scriptstyle \bullet$} [c] at  13 12
\put {$ \scriptstyle \bullet$} [c] at  13.5 6
\put {$ \scriptstyle \bullet$} [c] at  16 12
\setlinear \plot 10  0 13 12 12.5 7.5 12.5 4.5  13 0 13.5 6 13 12  /
\setlinear \plot 13 0 16 12 /
\put{$5{,}040$} [c] at 13 -2
\endpicture
\end{minipage}
\begin{minipage}{4cm}
\beginpicture
\setcoordinatesystem units    <1.5mm,2mm>
\setplotarea x from 0 to 16, y from -2 to 15
\put{1.520)} [l] at 2 12
\put {$ \scriptstyle \bullet$} [c] at  10 0
\put {$ \scriptstyle \bullet$} [c] at  10 12
\put {$ \scriptstyle \bullet$} [c] at  11 6
\put {$ \scriptstyle \bullet$} [c] at  12 0
\put {$ \scriptstyle \bullet$} [c] at  12 12
\put {$ \scriptstyle \bullet$} [c] at  13 6
\put {$ \scriptstyle \bullet$} [c] at  16 12
\setlinear \plot 10 0 10 12 11 6 12 12 13 6 12 0 11 6 /
\setlinear \plot  12 0 16 12  /
\put{$5{,}040$} [c] at 13 -2
\endpicture
\end{minipage}
\begin{minipage}{4cm}
\beginpicture
\setcoordinatesystem units    <1.5mm,2mm>
\setplotarea x from 0 to 16, y from -2 to 15
\put{1.521)} [l] at 2 12
\put {$ \scriptstyle \bullet$} [c] at  10 0
\put {$ \scriptstyle \bullet$} [c] at  10 12
\put {$ \scriptstyle \bullet$} [c] at  11 6
\put {$ \scriptstyle \bullet$} [c] at  12 0
\put {$ \scriptstyle \bullet$} [c] at  12 12
\put {$ \scriptstyle \bullet$} [c] at  13 6
\put {$ \scriptstyle \bullet$} [c] at  16 0
\setlinear \plot 10 12 10 0 11 6 12 0 13 6 12 12 11 6 /
\setlinear \plot  12 12 16 0  /
\put{$5{,}040$} [c] at 13 -2
\endpicture
\end{minipage}
\begin{minipage}{4cm}
\beginpicture
\setcoordinatesystem units    <1.5mm,2mm>
\setplotarea x from 0 to 16, y from -2 to 15
\put {1.522)} [l] at 2 12
\put {$ \scriptstyle \bullet$} [c] at  10 0
\put {$ \scriptstyle \bullet$} [c] at  10 12
\put {$ \scriptstyle \bullet$} [c] at   12 12
\put {$ \scriptstyle \bullet$} [c] at  14 0
\put {$ \scriptstyle \bullet$} [c] at  14 6
\put {$ \scriptstyle \bullet$} [c] at  14  12
\put {$ \scriptstyle \bullet$} [c] at  16 12
\setlinear \plot  16 12 14 0 10 12 10 0 12 12 14 6 14 0 /
\setlinear \plot  14 6 14 12  /
\put{$5{,}040$} [c] at 13 -2
\endpicture
\end{minipage}
\begin{minipage}{4cm}
\beginpicture
\setcoordinatesystem units    <1.5mm,2mm>
\setplotarea x from 0 to 16, y from -2 to 15
\put {1.523)} [l] at 2 12
\put {$ \scriptstyle \bullet$} [c] at  10 0
\put {$ \scriptstyle \bullet$} [c] at  10 12
\put {$ \scriptstyle \bullet$} [c] at   12 0
\put {$ \scriptstyle \bullet$} [c] at  14 0
\put {$ \scriptstyle \bullet$} [c] at  14 6
\put {$ \scriptstyle \bullet$} [c] at  14  12
\put {$ \scriptstyle \bullet$} [c] at  16 0
\setlinear \plot  16 0 14 12 10 0 10 12 12 0 14 6 14 12 /
\setlinear \plot  14 6 14 0  /
\put{$5{,}040$} [c] at 13 -2
\endpicture
\end{minipage}
\begin{minipage}{4cm}
\beginpicture
\setcoordinatesystem units    <1.5mm,2mm>
\setplotarea x from 0 to 16, y from -2 to 15
\put{1.524)} [l] at 2 12
\put {$ \scriptstyle \bullet$} [c] at  10 0
\put {$ \scriptstyle \bullet$} [c] at  10 12
\put {$ \scriptstyle \bullet$} [c] at  13 0
\put {$ \scriptstyle \bullet$} [c] at  13 6
\put {$ \scriptstyle \bullet$} [c] at  13  12
\put {$ \scriptstyle \bullet$} [c] at  16  12
\put {$ \scriptstyle \bullet$} [c] at  16  0
\setlinear \plot 16 12 13 0 13 12  10 0  10 12 13 0 /
\setlinear \plot 16 0 16 12  /
\put{$5{,}040$} [c]  at 13 -2
\endpicture
\end{minipage}
$$
$$
\begin{minipage}{4cm}
\beginpicture
\setcoordinatesystem units    <1.5mm,2mm>
\setplotarea x from 0 to 16, y from -2 to 15
\put{1.525)} [l] at 2 12
\put {$ \scriptstyle \bullet$} [c] at  10 0
\put {$ \scriptstyle \bullet$} [c] at  10 12
\put {$ \scriptstyle \bullet$} [c] at  13 0
\put {$ \scriptstyle \bullet$} [c] at  13 6
\put {$ \scriptstyle \bullet$} [c] at  13  12
\put {$ \scriptstyle \bullet$} [c] at  16  12
\put {$ \scriptstyle \bullet$} [c] at  16  0
\setlinear \plot 16 0 13 12 13 0  10 12  10 0 13 12 /
\setlinear \plot 16 0 16 12  /
\put{$5{,}040$} [c]  at 13 -2
\endpicture
\end{minipage}
\begin{minipage}{4cm}
\beginpicture
\setcoordinatesystem units    <1.5mm,2mm>
\setplotarea x from 0 to 16, y from -2 to 15
\put{1.526)} [l] at 2 12
\put {$ \scriptstyle \bullet$} [c] at  10 0
\put {$ \scriptstyle \bullet$} [c] at  10 12
\put {$ \scriptstyle \bullet$} [c] at  13 0
\put {$ \scriptstyle \bullet$} [c] at  13 12
\put {$ \scriptstyle \bullet$} [c] at  16 0
\put {$ \scriptstyle \bullet$} [c] at  16  6
\put {$ \scriptstyle \bullet$} [c] at  16  12
\setlinear \plot  10 12 10 0 13 12 13 0   /
\setlinear \plot 16 0 16 12 /
\setlinear \plot 13 12 16 6  /
\put{$5{,}040 $} [c]  at 13 -2
\endpicture
\end{minipage}
\begin{minipage}{4cm}
\beginpicture
\setcoordinatesystem units    <1.5mm,2mm>
\setplotarea x from 0 to 16, y from -2 to 15
\put{1.527)} [l] at 2 12
\put {$ \scriptstyle \bullet$} [c] at  10 0
\put {$ \scriptstyle \bullet$} [c] at  10 12
\put {$ \scriptstyle \bullet$} [c] at  13 0
\put {$ \scriptstyle \bullet$} [c] at  13 12
\put {$ \scriptstyle \bullet$} [c] at  16 0
\put {$ \scriptstyle \bullet$} [c] at  16  6
\put {$ \scriptstyle \bullet$} [c] at  16  12
\setlinear \plot  10 0 10 12 13 0 13 12   /
\setlinear \plot 16 0 16 12 /
\setlinear \plot 13 0 16 6  /
\put{$5{,}040 $} [c]  at 13 -2
\endpicture
\end{minipage}
\begin{minipage}{4cm}
\beginpicture
\setcoordinatesystem units    <1.5mm,2mm>
\setplotarea x from 0 to 16, y from -2 to 15
\put{1.528)} [l] at 2 12
\put {$ \scriptstyle \bullet$} [c] at  10 0
\put {$ \scriptstyle \bullet$} [c] at  10 12
\put {$ \scriptstyle \bullet$} [c] at  13 0
\put {$ \scriptstyle \bullet$} [c] at  13 12
\put {$ \scriptstyle \bullet$} [c] at  13 6
\put {$ \scriptstyle \bullet$} [c] at  16  12
\put {$ \scriptstyle \bullet$} [c] at  16  0
\setlinear \plot  10 12 10 0 13 6 /
\setlinear \plot 16 12 16 0 13 12 13 0 /
\put{$5{,}040 $} [c]  at 13 -2
\endpicture
\end{minipage}
\begin{minipage}{4cm}
\beginpicture
\setcoordinatesystem units    <1.5mm,2mm>
\setplotarea x from 0 to 16, y from -2 to 15
\put{1.529)} [l] at 2 12
\put {$ \scriptstyle \bullet$} [c] at  10 0
\put {$ \scriptstyle \bullet$} [c] at  10 12
\put {$ \scriptstyle \bullet$} [c] at  13 0
\put {$ \scriptstyle \bullet$} [c] at  13 12
\put {$ \scriptstyle \bullet$} [c] at  13 6
\put {$ \scriptstyle \bullet$} [c] at  16  12
\put {$ \scriptstyle \bullet$} [c] at  16  0
\setlinear \plot  10 0 10 12 13 6 /
\setlinear \plot 16 0 16 12 13 0 13 12 /
\put{$5{,}040 $} [c]  at 13 -2
\endpicture
\end{minipage}
\begin{minipage}{4cm}
\beginpicture
\setcoordinatesystem units    <1.5mm,2mm>
\setplotarea x from 0 to 16, y from -2 to 15
\put{1.530)} [l] at 2 12
\put {$ \scriptstyle \bullet$} [c] at  10 0
\put {$ \scriptstyle \bullet$} [c] at  10 12
\put {$ \scriptstyle \bullet$} [c] at  13 12
\put {$ \scriptstyle \bullet$} [c] at  13 6
\put {$ \scriptstyle \bullet$} [c] at  13 0
\put {$ \scriptstyle \bullet$} [c] at  16 12
\put {$ \scriptstyle \bullet$} [c] at  16  0
\setlinear \plot  10  0 10 12 13 6 13 0 16 12    /
\setlinear \plot  10 0 13 12  13 6  /
\setlinear \plot  13 12  16 0  /
\put{$5{,}040 $} [c]  at 13 -2
\endpicture
\end{minipage}
$$

$$
\begin{minipage}{4cm}
\beginpicture
\setcoordinatesystem units    <1.5mm,2mm>
\setplotarea x from 0 to 16, y from -2 to 15
\put{1.531)} [l] at 2 12
\put {$ \scriptstyle \bullet$} [c] at  10 0
\put {$ \scriptstyle \bullet$} [c] at  10 12
\put {$ \scriptstyle \bullet$} [c] at  13 12
\put {$ \scriptstyle \bullet$} [c] at  13 6
\put {$ \scriptstyle \bullet$} [c] at  13 0
\put {$ \scriptstyle \bullet$} [c] at  16 12
\put {$ \scriptstyle \bullet$} [c] at  16  0
\setlinear \plot  10  12 10 0 13 6 13 12 16 0    /
\setlinear \plot  10 12 13 0  13 6  /
\setlinear \plot  13 0  16 12  /
\put{$5{,}040 $} [c]  at 13 -2
\endpicture
\end{minipage}
\begin{minipage}{4cm}
\beginpicture
\setcoordinatesystem units    <1.5mm,2mm>
\setplotarea x from 0 to 16, y from -2 to 15
\put{1.532)} [l] at 2 12
\put {$ \scriptstyle \bullet$} [c] at  10 0
\put {$ \scriptstyle \bullet$} [c] at  10 12
\put {$ \scriptstyle \bullet$} [c] at  13 12
\put {$ \scriptstyle \bullet$} [c] at  13 0
\put {$ \scriptstyle \bullet$} [c] at  16  0
\put {$ \scriptstyle \bullet$} [c] at  16 6
\put {$ \scriptstyle \bullet$} [c] at  16  12
\setlinear \plot 10 0 10 12 16 0 16 12  /
\setlinear \plot 13 0 10 12 16 0 13 12 13 0  /
\put{$5{,}040$} [c]  at 13 -2
\endpicture
\end{minipage}
\begin{minipage}{4cm}
\beginpicture
\setcoordinatesystem units    <1.5mm,2mm>
\setplotarea x from 0 to 16, y from -2 to 15
\put{1.533)} [l] at 2 12
\put {$ \scriptstyle \bullet$} [c] at  10 0
\put {$ \scriptstyle \bullet$} [c] at  10 12
\put {$ \scriptstyle \bullet$} [c] at  13 12
\put {$ \scriptstyle \bullet$} [c] at  13 0
\put {$ \scriptstyle \bullet$} [c] at  16  0
\put {$ \scriptstyle \bullet$} [c] at  16 6
\put {$ \scriptstyle \bullet$} [c] at  16  12
\setlinear \plot 10 12 10 0 16 12 16 0  /
\setlinear \plot 13 12 10 0 16 12 13 0 13 12  /
\put{$5{,}040$} [c]  at 13 -2
\endpicture
\end{minipage}
\begin{minipage}{4cm}
\beginpicture
\setcoordinatesystem units    <1.5mm,2mm>
\setplotarea x from 0 to 16, y from -2 to 15
\put{1.534)} [l] at 2 12
\put {$ \scriptstyle \bullet$} [c] at  10 0
\put {$ \scriptstyle \bullet$} [c] at  10 12
\put {$ \scriptstyle \bullet$} [c] at  13 0
\put {$ \scriptstyle \bullet$} [c] at  13 12
\put {$ \scriptstyle \bullet$} [c] at  16 0
\put {$ \scriptstyle \bullet$} [c] at  16 6
\put {$ \scriptstyle \bullet$} [c] at  16  12
\setlinear \plot 10 12 13 0 16 12  16 0 13 12 10 0 /
\setlinear \plot 13 0 13 12  /
\put{$5{,}040$} [c]  at 13 -2
\endpicture
\end{minipage}
\begin{minipage}{4cm}
\beginpicture
\setcoordinatesystem units    <1.5mm,2mm>
\setplotarea x from 0 to 16, y from -2 to 15
\put{1.535)} [l] at 2 12
\put {$ \scriptstyle \bullet$} [c] at  10 6
\put {$ \scriptstyle \bullet$} [c] at  11 0
\put {$ \scriptstyle \bullet$} [c] at  11 6
\put {$ \scriptstyle \bullet$} [c] at  11 12
\put {$ \scriptstyle \bullet$} [c] at  12 6
\put {$ \scriptstyle \bullet$} [c] at  10 0
\put {$ \scriptstyle \bullet$} [c] at  16 12
\setlinear \plot 10  0  10 6 11 12 11 0 12 6 11 12   /
\setlinear \plot 10 6 11 0 16 12   /
\put{$2{,}520 $} [c] at 13 -2
\endpicture
\end{minipage}
\begin{minipage}{4cm}
\beginpicture
\setcoordinatesystem units    <1.5mm,2mm>
\setplotarea x from 0 to 16, y from -2 to 15
\put{1.536)} [l] at 2 12
\put {$ \scriptstyle \bullet$} [c] at  10 6
\put {$ \scriptstyle \bullet$} [c] at  11 0
\put {$ \scriptstyle \bullet$} [c] at  11 6
\put {$ \scriptstyle \bullet$} [c] at  11 12
\put {$ \scriptstyle \bullet$} [c] at  12 6
\put {$ \scriptstyle \bullet$} [c] at  10 12
\put {$ \scriptstyle \bullet$} [c] at  16 0
\setlinear \plot 10  12  10 6 11 0 11 12 12 6 11 0   /
\setlinear \plot 10 6 11 12 16 0   /
\put{$2{,}520  $} [c] at 13 -2
\endpicture
\end{minipage}
$$
$$
\begin{minipage}{4cm}
\beginpicture
\setcoordinatesystem units    <1.5mm,2mm>
\setplotarea x from 0 to 16, y from -2 to 15
\put{1.537)} [l] at 2 12
\put {$ \scriptstyle \bullet$} [c] at  10 12
\put {$ \scriptstyle \bullet$} [c] at   13 0
\put {$ \scriptstyle \bullet$} [c] at  13 6
\put {$ \scriptstyle \bullet$} [c] at  13 12
\put {$ \scriptstyle \bullet$} [c] at  13.8 3
\put {$ \scriptstyle \bullet$} [c] at  16 0
\put {$ \scriptstyle \bullet$} [c] at  16 12
\setlinear \plot  10 12 13  0 16 12 16 0 13 12 13 0    /
\put{$2{,}520$} [c] at 13 -2
\endpicture
\end{minipage}
\begin{minipage}{4cm}
\beginpicture
\setcoordinatesystem units    <1.5mm,2mm>
\setplotarea x from 0 to 16, y from -2 to 15
\put{1.538)} [l] at 2 12
\put {$ \scriptstyle \bullet$} [c] at  10 0
\put {$ \scriptstyle \bullet$} [c] at   13 0
\put {$ \scriptstyle \bullet$} [c] at  13 6
\put {$ \scriptstyle \bullet$} [c] at  13 12
\put {$ \scriptstyle \bullet$} [c] at  13.8 9
\put {$ \scriptstyle \bullet$} [c] at  16 0
\put {$ \scriptstyle \bullet$} [c] at  16 12
\setlinear \plot  10 0 13  12 16 0 16 12 13 0 13 12    /
\put{$2{,}520$} [c] at 13 -2
\endpicture
\end{minipage}
\begin{minipage}{4cm}
\beginpicture
\setcoordinatesystem units    <1.5mm,2mm>
\setplotarea x from 0 to 16, y from -2 to 15
\put{1.539)} [l] at 2 12
\put {$ \scriptstyle \bullet$} [c] at 10 6
\put {$ \scriptstyle \bullet$} [c] at  10.5 0
\put {$ \scriptstyle \bullet$} [c] at  10.5 12
\put {$ \scriptstyle \bullet$} [c] at  11 6
\put {$ \scriptstyle \bullet$} [c] at  14 12
\put {$ \scriptstyle \bullet$} [c] at  16 12
\put {$ \scriptstyle \bullet$} [c] at  16 0
\setlinear \plot 10.5 0 10 6 10.5 12 11 6 10.5 0 14 12 16 0 16 12 /
\put{$2{,}520$} [c] at 13 -2
\endpicture
\end{minipage}
\begin{minipage}{4cm}
\beginpicture
\setcoordinatesystem units    <1.5mm,2mm>
\setplotarea x from 0 to 16, y from -2 to 15
\put{1.540)} [l] at 2 12
\put {$ \scriptstyle \bullet$} [c] at 10 6
\put {$ \scriptstyle \bullet$} [c] at  10.5 0
\put {$ \scriptstyle \bullet$} [c] at  10.5 12
\put {$ \scriptstyle \bullet$} [c] at  11 6
\put {$ \scriptstyle \bullet$} [c] at  14 0
\put {$ \scriptstyle \bullet$} [c] at  16 12
\put {$ \scriptstyle \bullet$} [c] at  16 0
\setlinear \plot 10.5 12 10 6 10.5 0 11 6 10.5 12 14 0 16 12 16 0 /
\put{$2{,}520$} [c] at 13 -2
\endpicture
\end{minipage}
\begin{minipage}{4cm}
\beginpicture
\setcoordinatesystem units    <1.5mm,2mm>
\setplotarea x from 0 to 16, y from -2 to 15
\put {1.541)} [l] at 2 12
\put {$ \scriptstyle \bullet$} [c] at  10 12
\put {$ \scriptstyle \bullet$} [c] at  12  0
\put {$ \scriptstyle \bullet$} [c] at  12 12
\put {$ \scriptstyle \bullet$} [c] at  14 0
\put {$ \scriptstyle \bullet$} [c] at  14 6
\put {$ \scriptstyle \bullet$} [c] at  14 12
\put {$ \scriptstyle \bullet$} [c] at  16  12
\setlinear \plot  10 12 12 0  14 12 14 0 /
\setlinear \plot  12 12 12 0 /
\setlinear \plot 16 12 14 6  /
\put{$2{,}520$} [c] at 13 -2
\endpicture
\end{minipage}
\begin{minipage}{4cm}
\beginpicture
\setcoordinatesystem units    <1.5mm,2mm>
\setplotarea x from 0 to 16, y from -2 to 15
\put {1.542)} [l] at 2 12
\put {$ \scriptstyle \bullet$} [c] at  10 0
\put {$ \scriptstyle \bullet$} [c] at  12  0
\put {$ \scriptstyle \bullet$} [c] at  12 12
\put {$ \scriptstyle \bullet$} [c] at  14 0
\put {$ \scriptstyle \bullet$} [c] at  14 6
\put {$ \scriptstyle \bullet$} [c] at  14 12
\put {$ \scriptstyle \bullet$} [c] at  16  0
\setlinear \plot  10 0 12 12  14 0 14 12 /
\setlinear \plot  12 0 12 12 /
\setlinear \plot 16 0 14 6  /
\put{$2{,}520$} [c] at 13 -2
\endpicture
\end{minipage}
$$

$$
\begin{minipage}{4cm}
\beginpicture
\setcoordinatesystem units    <1.5mm,2mm>
\setplotarea x from 0 to 16, y from -2 to 15
\put {1.543)} [l] at 2 12
\put {$ \scriptstyle \bullet$} [c] at  10 0
\put {$ \scriptstyle \bullet$} [c] at  10 12
\put {$ \scriptstyle \bullet$} [c] at  12 12
\put {$ \scriptstyle \bullet$} [c] at  14 0
\put {$ \scriptstyle \bullet$} [c] at  14 6
\put {$ \scriptstyle \bullet$} [c] at  14  12
\put {$ \scriptstyle \bullet$} [c] at  16  12
\setlinear \plot  10 12 10 0  12 12 14 0   14 12 /
\setlinear \plot 16 12 14 6 /
\put{$2{,}520$} [c] at 13 -2
\endpicture
\end{minipage}
\begin{minipage}{4cm}
\beginpicture
\setcoordinatesystem units    <1.5mm,2mm>
\setplotarea x from 0 to 16, y from -2 to 15
\put {1.544)} [l] at 2 12
\put {$ \scriptstyle \bullet$} [c] at  10 0
\put {$ \scriptstyle \bullet$} [c] at  10 12
\put {$ \scriptstyle \bullet$} [c] at  12 0
\put {$ \scriptstyle \bullet$} [c] at  14 0
\put {$ \scriptstyle \bullet$} [c] at  14 6
\put {$ \scriptstyle \bullet$} [c] at  14  12
\put {$ \scriptstyle \bullet$} [c] at  16  0
\setlinear \plot  10 0 10 12  12 0 14 12   14 0 /
\setlinear \plot 16 0 14 6 /
\put{$2{,}520$} [c] at 13 -2
\endpicture
\end{minipage}
\begin{minipage}{4cm}
\beginpicture
\setcoordinatesystem units    <1.5mm,2mm>
\setplotarea x from 0 to 16, y from -2 to 15
\put {1.545)} [l] at 2 12
\put {$ \scriptstyle \bullet$} [c] at  10 12
\put {$ \scriptstyle \bullet$} [c] at  12 0
\put {$ \scriptstyle \bullet$} [c] at  12 12
\put {$ \scriptstyle \bullet$} [c] at  14 0
\put {$ \scriptstyle \bullet$} [c] at  14  12
\put {$ \scriptstyle \bullet$} [c] at  16 12
\put {$ \scriptstyle \bullet$} [c] at  16 6
\setlinear \plot  10 12 12 0 12  12 14 0 16 6 16 12  /
\setlinear \plot  16 6 12 0 14 12 14 0   /
\put{$2{,}520$} [c] at 13 -2
\endpicture
\end{minipage}
\begin{minipage}{4cm}
\beginpicture
\setcoordinatesystem units    <1.5mm,2mm>
\setplotarea x from 0 to 16, y from -2 to 15
\put {1.546)} [l] at 2 12
\put {$ \scriptstyle \bullet$} [c] at  10 0
\put {$ \scriptstyle \bullet$} [c] at  12 0
\put {$ \scriptstyle \bullet$} [c] at  12 12
\put {$ \scriptstyle \bullet$} [c] at  14 0
\put {$ \scriptstyle \bullet$} [c] at  14  12
\put {$ \scriptstyle \bullet$} [c] at  16 0
\put {$ \scriptstyle \bullet$} [c] at  16 6
\setlinear \plot  10 0 12 12 12  0 14 12 16 6 16 0  /
\setlinear \plot  16 6 12 12 14 0 14 12   /
\put{$2{,}520$} [c] at 13 -2
\endpicture
\end{minipage}
\begin{minipage}{4cm}
\beginpicture
\setcoordinatesystem units    <1.5mm,2mm>
\setplotarea x from 0 to 16, y from -2 to 15
\put{1.547)} [l] at 2 12
\put {$ \scriptstyle \bullet$} [c] at  10 0
\put {$ \scriptstyle \bullet$} [c] at  10 12
\put {$ \scriptstyle \bullet$} [c] at  13 0
\put {$ \scriptstyle \bullet$} [c] at  13 12
\put {$ \scriptstyle \bullet$} [c] at  14.5 6
\put {$ \scriptstyle \bullet$} [c] at  16  12
\put {$ \scriptstyle \bullet$} [c] at  16  0
\setlinear \plot  10 0 10 12 13 0 13 12  10 0  16 12  13 0 /
\setlinear \plot 16 12  16 0 /
\put{$2{,}520 $} [c]  at 13 -2
\endpicture
\end{minipage}
\begin{minipage}{4cm}
\beginpicture
\setcoordinatesystem units    <1.5mm,2mm>
\setplotarea x from 0 to 16, y from -2 to 15
\put{1.548)} [l] at 2 12
\put {$ \scriptstyle \bullet$} [c] at  10 0
\put {$ \scriptstyle \bullet$} [c] at  10 12
\put {$ \scriptstyle \bullet$} [c] at  13 0
\put {$ \scriptstyle \bullet$} [c] at  13 12
\put {$ \scriptstyle \bullet$} [c] at  14.5 6
\put {$ \scriptstyle \bullet$} [c] at  16  12
\put {$ \scriptstyle \bullet$} [c] at  16  0
\setlinear \plot  10 12 10 0 13 12 13 0  10 12  16 0  13 12 /
\setlinear \plot 16 12  16 0 /
\put{$2{,}520 $} [c]  at 13 -2
\endpicture
\end{minipage}
$$
$$
\begin{minipage}{4cm}
\beginpicture
\setcoordinatesystem units    <1.5mm,2mm>
\setplotarea x from 0 to 16, y from -2 to 15
\put {1.549)} [l] at  2 12
\put {$ \scriptstyle \bullet$} [c] at  10 12
\put {$ \scriptstyle \bullet$} [c] at  12 12
\put {$ \scriptstyle \bullet$} [c] at  14  12
\put {$ \scriptstyle \bullet$} [c] at  16 12
\put {$ \scriptstyle \bullet$} [c] at  10 0
\put {$ \scriptstyle \bullet$} [c] at   13 0
\put {$ \scriptstyle \bullet$} [c] at  16 0
\setlinear \plot   10 12 10  0  12 12 13 0 16 12  16 0  14 12 10  0 16 12 13 0 14 12 /
\setlinear \plot   13 0 16  12  /
\put{$2{,}520$} [c] at 13 -2
\endpicture
\end{minipage}
\begin{minipage}{4cm}
\beginpicture
\setcoordinatesystem units    <1.5mm,2mm>
\setplotarea x from 0 to 16, y from -2 to 15
\put {1.550)} [l] at  2 12
\put {$ \scriptstyle \bullet$} [c] at  10 0
\put {$ \scriptstyle \bullet$} [c] at  12 0
\put {$ \scriptstyle \bullet$} [c] at  14  0
\put {$ \scriptstyle \bullet$} [c] at  16 0
\put {$ \scriptstyle \bullet$} [c] at  10 12
\put {$ \scriptstyle \bullet$} [c] at   13 12
\put {$ \scriptstyle \bullet$} [c] at  16 12
\setlinear \plot   10 0 10  12  12 0 13 12 16 0  16 12  14 0 10  12 16 0 13 12 14 0 /
\setlinear \plot   13 12 16  0  /
\put{$2{,}520$} [c] at 13 -2
\endpicture
\end{minipage}
\begin{minipage}{4cm}
\beginpicture
\setcoordinatesystem units    <1.5mm,2mm>
\setplotarea x from 0 to 16, y from -2 to 15
\put{1.551)} [l] at 2 12
\put {$ \scriptstyle \bullet$} [c] at  10 12
\put {$ \scriptstyle \bullet$} [c] at  10  0
\put {$ \scriptstyle \bullet$} [c] at  13 12
\put {$ \scriptstyle \bullet$} [c] at  13 0
\put {$ \scriptstyle \bullet$} [c] at  16 12
\put {$ \scriptstyle \bullet$} [c] at  16 6
\put {$ \scriptstyle \bullet$} [c] at  16 0
\setlinear \plot  10 12 10 0 13 12   /
\setlinear \plot  10 0 16 6 16 0  /
\setlinear \plot  13 0 16 6 16 12  /
\put{$1{,}260$} [c]  at 13 -2
\endpicture
\end{minipage}
\begin{minipage}{4cm}
\beginpicture
\setcoordinatesystem units    <1.5mm,2mm>
\setplotarea x from 0 to 16, y from -2 to 15
\put{1.552)} [l] at 2 12
\put {$ \scriptstyle \bullet$} [c] at  10 12
\put {$ \scriptstyle \bullet$} [c] at  10  0
\put {$ \scriptstyle \bullet$} [c] at  13 12
\put {$ \scriptstyle \bullet$} [c] at  13 0
\put {$ \scriptstyle \bullet$} [c] at  16 12
\put {$ \scriptstyle \bullet$} [c] at  16 6
\put {$ \scriptstyle \bullet$} [c] at  16 0
\setlinear \plot  10 0 10 12 13 0   /
\setlinear \plot  10 12 16 6 16 12  /
\setlinear \plot  13 12 16 6 16 0  /
\put{$1{,}260$} [c]  at 13 -2
\endpicture
\end{minipage}
\begin{minipage}{4cm}
\beginpicture
\setcoordinatesystem units    <1.5mm,2mm>
\setplotarea x from 0 to 16, y from -2 to 15
\put {1.553)} [l] at  2 12
\put {$ \scriptstyle \bullet$} [c] at  10 12
\put {$ \scriptstyle \bullet$} [c] at  12 12
\put {$ \scriptstyle \bullet$} [c] at  14  12
\put {$ \scriptstyle \bullet$} [c] at  16 12
\put {$ \scriptstyle \bullet$} [c] at  10 0
\put {$ \scriptstyle \bullet$} [c] at   13 0
\put {$ \scriptstyle \bullet$} [c] at  16 0
\setlinear \plot   10 0 10 12 13 0 16 12 16 0 14 12 13 0 /
\setlinear \plot   10 0 12 12  16 0  /
\put{$1{,}260$} [c] at 13 -2
\endpicture
\end{minipage}
\begin{minipage}{4cm}
\beginpicture
\setcoordinatesystem units    <1.5mm,2mm>
\setplotarea x from 0 to 16, y from -2 to 15
\put {1.554)} [l] at  2 12
\put {$ \scriptstyle \bullet$} [c] at  10 0
\put {$ \scriptstyle \bullet$} [c] at  12 0
\put {$ \scriptstyle \bullet$} [c] at  14  0
\put {$ \scriptstyle \bullet$} [c] at  16 0
\put {$ \scriptstyle \bullet$} [c] at  10 12
\put {$ \scriptstyle \bullet$} [c] at   13 12
\put {$ \scriptstyle \bullet$} [c] at  16 12
\setlinear \plot   10 12 10 0 13 12 16 0 16 12 14 0 13 12 /
\setlinear \plot   10 12 12 0  16 12  /
\put{$1{,}260$} [c] at 13 -2
\endpicture
\end{minipage}
$$

$$
\begin{minipage}{4cm}
\beginpicture
\setcoordinatesystem units    <1.5mm,2mm>
\setplotarea x from 0 to 16, y from -2 to 15
\put{1.555)} [l] at 2 12
\put {$ \scriptstyle \bullet$} [c] at 10 6
\put {$ \scriptstyle \bullet$} [c] at 12 0
\put {$ \scriptstyle \bullet$} [c] at 12 6
\put {$ \scriptstyle \bullet$} [c] at 12 12
\put {$ \scriptstyle \bullet$} [c] at 14 6
\put {$ \scriptstyle \bullet$} [c] at 16 6
\put {$ \scriptstyle \bullet$} [c] at 16 12
\setlinear \plot 16 12 16 6 12 0 10 6 12 12 14 6 12 0    /
\setlinear \plot 12  0 12 12     /
\put{$840$} [c] at 13 -2
\endpicture
\end{minipage}
\begin{minipage}{4cm}
\beginpicture
\setcoordinatesystem units    <1.5mm,2mm>
\setplotarea x from 0 to 16, y from -2 to 15
\put{1.556)} [l] at 2 12
\put {$ \scriptstyle \bullet$} [c] at 10 6
\put {$ \scriptstyle \bullet$} [c] at 12 0
\put {$ \scriptstyle \bullet$} [c] at 12 6
\put {$ \scriptstyle \bullet$} [c] at 12 12
\put {$ \scriptstyle \bullet$} [c] at 14 6
\put {$ \scriptstyle \bullet$} [c] at 16 6
\put {$ \scriptstyle \bullet$} [c] at 16 0
\setlinear \plot 16 0 16 6 12 12 10 6 12 0 14 6 12 12    /
\setlinear \plot 12  0 12 12     /
\put{$840$} [c] at 13 -2
\endpicture
\end{minipage}
\begin{minipage}{4cm}
\beginpicture
\setcoordinatesystem units    <1.5mm,2mm>
\setplotarea x from 0 to 16, y from -2 to 15
\put{1.557)} [l] at 2 12
\put {$ \scriptstyle \bullet$} [c] at 10 6
\put {$ \scriptstyle \bullet$} [c] at 10 12
\put {$ \scriptstyle \bullet$} [c] at 13 0
\put {$ \scriptstyle \bullet$} [c] at 13 6
\put {$ \scriptstyle \bullet$} [c] at 13 12
\put {$ \scriptstyle \bullet$} [c] at 16 6
\put {$ \scriptstyle \bullet$} [c] at 16 12
\setlinear \plot 10 12 10 6 13 0 16 6 16 12   /
\setlinear \plot 13 0 13 12   /
\put{$840$} [c] at 13 -2
\endpicture
\end{minipage}
\begin{minipage}{4cm}
\beginpicture
\setcoordinatesystem units    <1.5mm,2mm>
\setplotarea x from 0 to 16, y from -2 to 15
\put{1.558)} [l] at 2 12
\put {$ \scriptstyle \bullet$} [c] at 10 6
\put {$ \scriptstyle \bullet$} [c] at 10 0
\put {$ \scriptstyle \bullet$} [c] at 13 0
\put {$ \scriptstyle \bullet$} [c] at 13 6
\put {$ \scriptstyle \bullet$} [c] at 13 12
\put {$ \scriptstyle \bullet$} [c] at 16 6
\put {$ \scriptstyle \bullet$} [c] at 16 0
\setlinear \plot 10 0 10 6 13 12 16 6 16 0   /
\setlinear \plot 13 0 13 12   /
\put{$840$} [c] at 13 -2
\endpicture
\end{minipage}
\begin{minipage}{4cm}
\beginpicture
\setcoordinatesystem units    <1.5mm,2mm>
\setplotarea x from 0 to 16, y from -2 to 15
\put{1.559)} [l] at 2 12
\put {$ \scriptstyle \bullet$} [c] at 10 12
\put {$ \scriptstyle \bullet$} [c] at  12 12
\put {$ \scriptstyle \bullet$} [c] at 14 12
\put {$ \scriptstyle \bullet$} [c] at 12 6
\put {$ \scriptstyle \bullet$} [c] at 12  0
\put {$ \scriptstyle \bullet$} [c] at 16 6
\put {$ \scriptstyle \bullet$} [c] at 16 12

\setlinear \plot 16 12 16 6 12 0 12 12    /
\setlinear \plot 10 12 12 6 14 12     /
\put{$840$} [c] at 13 -2
\endpicture
\end{minipage}
\begin{minipage}{4cm}
\beginpicture
\setcoordinatesystem units    <1.5mm,2mm>
\setplotarea x from 0 to 16, y from -2 to 15
\put{1.560)} [l] at 2 12
\put {$ \scriptstyle \bullet$} [c] at 10 0
\put {$ \scriptstyle \bullet$} [c] at  12 0
\put {$ \scriptstyle \bullet$} [c] at 14 0
\put {$ \scriptstyle \bullet$} [c] at 12 6
\put {$ \scriptstyle \bullet$} [c] at 12  12
\put {$ \scriptstyle \bullet$} [c] at 16 6
\put {$ \scriptstyle \bullet$} [c] at 16 0
\setlinear \plot 16 0 16 6 12 12 12 0    /
\setlinear \plot 10 0 12 6 14 0     /
\put{$840$} [c] at 13 -2
\endpicture
\end{minipage}
$$
$$
\begin{minipage}{4cm}
\beginpicture
\setcoordinatesystem units    <1.5mm,2mm>
\setplotarea x from 0 to 16, y from -2 to 15
\put{1.561)} [l] at 2 12
\put {$ \scriptstyle \bullet$} [c] at  10 6
\put {$ \scriptstyle \bullet$} [c] at  11.5 0
\put {$ \scriptstyle \bullet$} [c] at  11.5 6
\put {$ \scriptstyle \bullet$} [c] at  11.5 12
\put {$ \scriptstyle \bullet$} [c] at  13 6
\put {$ \scriptstyle \bullet$} [c] at  16 0
\put {$ \scriptstyle \bullet$} [c] at  16 12
\setlinear \plot 11.5 0 10 6 11.5 12 11.5 0 13 6 11.5 12 16 0 16 12 11.5 0  /
\put{$840$} [c] at 13 -2
\endpicture
\end{minipage}
\begin{minipage}{4cm}
\beginpicture
\setcoordinatesystem units    <1.5mm,2mm>
\setplotarea x from 0 to 16, y from -2 to 15
\put {1.562)} [l] at 2 12
\put {$ \scriptstyle \bullet$} [c] at  10 0
\put {$ \scriptstyle \bullet$} [c] at  10 12
\put {$ \scriptstyle \bullet$} [c] at  12 12
\put {$ \scriptstyle \bullet$} [c] at  14 0
\put {$ \scriptstyle \bullet$} [c] at  14 6
\put {$ \scriptstyle \bullet$} [c] at  14 12
\put {$ \scriptstyle \bullet$} [c] at  16  12
\setlinear \plot  10 0 10 12  14 6 12 12 /
\setlinear \plot  14 0 14 12  /
\setlinear \plot 16 12 14 6 /
\put{$840$} [c] at 13 -2
\endpicture
\end{minipage}
\begin{minipage}{4cm}
\beginpicture
\setcoordinatesystem units    <1.5mm,2mm>
\setplotarea x from 0 to 16, y from -2 to 15
\put {1.563)} [l] at 2 12
\put {$ \scriptstyle \bullet$} [c] at  10 0
\put {$ \scriptstyle \bullet$} [c] at  10 12
\put {$ \scriptstyle \bullet$} [c] at  12 0
\put {$ \scriptstyle \bullet$} [c] at  14 0
\put {$ \scriptstyle \bullet$} [c] at  14 6
\put {$ \scriptstyle \bullet$} [c] at  14 12
\put {$ \scriptstyle \bullet$} [c] at  16  0
\setlinear \plot  10 12 10 0  14 6 12 0 /
\setlinear \plot  14 0 14 12  /
\setlinear \plot 16 0 14 6 /
\put{$840$} [c] at 13 -2
\endpicture
\end{minipage}
\begin{minipage}{4cm}
\beginpicture
\setcoordinatesystem units    <1.5mm,2mm>
\setplotarea x from 0 to 16, y from -2 to 15
\put {1.564)} [l] at 2 12
\put {$ \scriptstyle \bullet$} [c] at  10 0
\put {$ \scriptstyle \bullet$} [c] at  10 12
\put {$ \scriptstyle \bullet$} [c] at  12 12
\put {$ \scriptstyle \bullet$} [c] at  14 0
\put {$ \scriptstyle \bullet$} [c] at  14 6
\put {$ \scriptstyle \bullet$} [c] at  14 12
\put {$ \scriptstyle \bullet$} [c] at  16  12
\setlinear \plot  10 12 10 0  14 12 14 0  12 12 10 0 /
\setlinear \plot  10 12 14 0  16 12 10 0 /
\put{$840$} [c] at 13 -2
\endpicture
\end{minipage}
\begin{minipage}{4cm}
\beginpicture
\setcoordinatesystem units    <1.5mm,2mm>
\setplotarea x from 0 to 16, y from -2 to 15
\put {1.565)} [l] at 2 12
\put {$ \scriptstyle \bullet$} [c] at  10 0
\put {$ \scriptstyle \bullet$} [c] at  10 12
\put {$ \scriptstyle \bullet$} [c] at  12 0
\put {$ \scriptstyle \bullet$} [c] at  14 0
\put {$ \scriptstyle \bullet$} [c] at  14 6
\put {$ \scriptstyle \bullet$} [c] at  14 12
\put {$ \scriptstyle \bullet$} [c] at  16  0
\setlinear \plot  10 0 10 12  14 0 14 12  12 0 10 12 /
\setlinear \plot  10 0 14 12  16 0 10 12 /
\put{$840$} [c] at 13 -2
\endpicture
\end{minipage}
\begin{minipage}{4cm}
\beginpicture
\setcoordinatesystem units    <1.5mm,2mm>
\setplotarea x from 0 to 16, y from -2 to 15
\put{1.566)} [l] at 2 12
\put {$ \scriptstyle \bullet$} [c] at 10 6
\put {$ \scriptstyle \bullet$} [c] at 10 0
\put {$ \scriptstyle \bullet$} [c] at 10 12
\put {$ \scriptstyle \bullet$} [c] at 14 6
\put {$ \scriptstyle \bullet$} [c] at 14 12
\put {$ \scriptstyle \bullet$} [c] at 14 0
\setlinear \plot 10 12  10 0    /
\setlinear \plot 14 12  14 0     /
\setlinear \plot 10 0 14 6     /
\put{$5{,}040$} [c] at 11 -2
\put{$\scriptstyle \bullet$} [c] at 16  0 \endpicture
\end{minipage}
$$

$$
\begin{minipage}{4cm}
\beginpicture
\setcoordinatesystem units    <1.5mm,2mm>
\setplotarea x from 0 to 16, y from -2 to 15
\put{1.567)} [l] at 2 12
\put {$ \scriptstyle \bullet$} [c] at 10 6
\put {$ \scriptstyle \bullet$} [c] at 10 0
\put {$ \scriptstyle \bullet$} [c] at 10 12
\put {$ \scriptstyle \bullet$} [c] at 14 6
\put {$ \scriptstyle \bullet$} [c] at 14 12
\put {$ \scriptstyle \bullet$} [c] at 14 0
\setlinear \plot 10 12  10 0    /
\setlinear \plot 14 12  14 0     /
\setlinear \plot 10 12  14 6     /
\put{$5{,}040$} [c] at 13 -2
\put{$\scriptstyle \bullet$} [c] at 16  0 \endpicture
\end{minipage}
\begin{minipage}{4cm}
\beginpicture
\setcoordinatesystem units    <1.5mm,2mm>
\setplotarea x from 0 to 16, y from -2 to 15
\put{1.568)} [l] at 2 12
\put {$ \scriptstyle \bullet$} [c] at 10 0
\put {$ \scriptstyle \bullet$} [c] at 10 4
\put {$ \scriptstyle \bullet$} [c] at 10 8
\put {$ \scriptstyle \bullet$} [c] at 10 12
\put {$ \scriptstyle \bullet$} [c] at 14 12
\put {$ \scriptstyle \bullet$} [c] at 14 0
\setlinear \plot 10 12  10 0 14 12 14 0    /
\put{$5{,}040$} [c] at 13 -2
\put{$\scriptstyle \bullet$} [c] at 16  0 \endpicture
\end{minipage}
\begin{minipage}{4cm}
\beginpicture
\setcoordinatesystem units    <1.5mm,2mm>
\setplotarea x from 0 to 16, y from -2 to 15
\put{1.569)} [l] at 2 12
\put {$ \scriptstyle \bullet$} [c] at 10 0
\put {$ \scriptstyle \bullet$} [c] at 10 4
\put {$ \scriptstyle \bullet$} [c] at 10 8
\put {$ \scriptstyle \bullet$} [c] at 10 12
\put {$ \scriptstyle \bullet$} [c] at 14 12
\put {$ \scriptstyle \bullet$} [c] at 14 0
\setlinear \plot 10 0 10 12  14 0 14 12    /
\put{$5{,}040$} [c] at 13 -2
\put{$\scriptstyle \bullet$} [c] at 16  0 \endpicture
\end{minipage}
\begin{minipage}{4cm}
\beginpicture
\setcoordinatesystem units    <1.5mm,2mm>
\setplotarea x from 0 to 16, y from -2 to 15
\put{1.570)} [l] at 2 12
\put {$ \scriptstyle \bullet$} [c] at 10 6
\put {$ \scriptstyle \bullet$} [c] at 10 0
\put {$ \scriptstyle \bullet$} [c] at 10 12
\put {$ \scriptstyle \bullet$} [c] at 12 9
\put {$ \scriptstyle \bullet$} [c] at 14 12
\put {$ \scriptstyle \bullet$} [c] at 14 0
\setlinear \plot 10 12  10 0    /
\setlinear \plot 10 6  14 12 14 0     /
\put{$5{,}040$} [c] at 13 -2
\put{$\scriptstyle \bullet$} [c] at 16  0 \endpicture
\end{minipage}
\begin{minipage}{4cm}
\beginpicture
\setcoordinatesystem units    <1.5mm,2mm>
\setplotarea x from 0 to 16, y from -2 to 15
\put{1.571)} [l] at 2 12
\put {$ \scriptstyle \bullet$} [c] at 10 6
\put {$ \scriptstyle \bullet$} [c] at 10 0
\put {$ \scriptstyle \bullet$} [c] at 10 12
\put {$ \scriptstyle \bullet$} [c] at 12 3
\put {$ \scriptstyle \bullet$} [c] at 14 12
\put {$ \scriptstyle \bullet$} [c] at 14 0
\setlinear \plot 10 12  10 0    /
\setlinear \plot 10 6  14 0 14 12     /
\put{$5{,}040$} [c] at 13 -2
\put{$\scriptstyle \bullet$} [c] at 16  0 \endpicture
\end{minipage}
\begin{minipage}{4cm}
\beginpicture
\setcoordinatesystem units    <1.5mm,2mm>
\setplotarea x from 0 to 16, y from -2 to 15
\put{1.572)} [l] at 2 12
\put {$ \scriptstyle \bullet$} [c] at 10 6
\put {$ \scriptstyle \bullet$} [c] at 10.5 0
\put {$ \scriptstyle \bullet$} [c] at 10.5 12
\put {$ \scriptstyle \bullet$} [c] at 11 6
\put {$ \scriptstyle \bullet$} [c] at 14 12
\put {$ \scriptstyle \bullet$} [c] at 14 0
\setlinear \plot 10.5 12   14 0 14 12 11 6 10.5 0 10 6 10.5 12 11 6     /
\put{$5{,}040$} [c] at 13 -2
\put{$\scriptstyle \bullet$} [c] at 16  0 \endpicture
\end{minipage}
$$
$$
\begin{minipage}{4cm}
\beginpicture
\setcoordinatesystem units    <1.5mm,2mm>
\setplotarea x from 0 to 16, y from -2 to 15
\put{1.573)} [l] at 2 12
\put {$ \scriptstyle \bullet$} [c] at 10 6
\put {$ \scriptstyle \bullet$} [c] at 10.5 0
\put {$ \scriptstyle \bullet$} [c] at 10.5 12
\put {$ \scriptstyle \bullet$} [c] at 11 6
\put {$ \scriptstyle \bullet$} [c] at 14 12
\put {$ \scriptstyle \bullet$} [c] at 14 0
\setlinear \plot 10.5 0   14 12 14 0 11 6 10.5 12 10 6 10.5 0 11 6     /
\put{$5{,}040$} [c] at 13 -2
\put{$\scriptstyle \bullet$} [c] at 16  0 \endpicture
\end{minipage}
\begin{minipage}{4cm}
\beginpicture
\setcoordinatesystem units    <1.5mm,2mm>
\setplotarea x from 0 to 16, y from -2 to 15
\put{1.574)} [l] at 2 12
\put {$ \scriptstyle \bullet$} [c] at 10 12
\put {$ \scriptstyle \bullet$} [c] at 13 12
\put {$ \scriptstyle \bullet$} [c] at  11.5 0
\put {$ \scriptstyle \bullet$} [c] at 10.8 6
\put {$ \scriptstyle \bullet$} [c] at 14 0
\put {$ \scriptstyle \bullet$} [c] at 14 6
\put {$ \scriptstyle \bullet$} [c] at 14 12
\put {$ \scriptstyle \bullet$} [c] at 14 0
\setlinear \plot 10 12 11.5 0 13 12   /
\setlinear \plot 14 12 14 0 /
\put{$5{,}040$} [c] at 13 -2
\endpicture
\end{minipage}
\begin{minipage}{4cm}
\beginpicture
\setcoordinatesystem units    <1.5mm,2mm>
\setplotarea x from 0 to 16, y from -2 to 15
\put{1.575)} [l] at 2 12
\put {$ \scriptstyle \bullet$} [c] at 10 0
\put {$ \scriptstyle \bullet$} [c] at 13 0
\put {$ \scriptstyle \bullet$} [c] at  11.5 12
\put {$ \scriptstyle \bullet$} [c] at 10.8 6
\put {$ \scriptstyle \bullet$} [c] at 14 0
\put {$ \scriptstyle \bullet$} [c] at 14 6
\put {$ \scriptstyle \bullet$} [c] at 14 12
\put {$ \scriptstyle \bullet$} [c] at 14 0
\setlinear \plot 10 0 11.5 12 13 0   /
\setlinear \plot 14 12 14 0 /
\put{$5{,}040$} [c] at 13 -2
\endpicture
\end{minipage}
\begin{minipage}{4cm}
\beginpicture
\setcoordinatesystem units    <1.5mm,2mm>
\setplotarea x from 0 to 16, y from -2 to 15
\put{1.576)} [l] at 2 12
\put {$ \scriptstyle \bullet$} [c] at 10 6
\put {$ \scriptstyle \bullet$} [c] at 10 12
\put {$ \scriptstyle \bullet$} [c] at 12 0
\put {$ \scriptstyle \bullet$} [c] at 12 6
\put {$ \scriptstyle \bullet$} [c] at 14 6
\put {$ \scriptstyle \bullet$} [c] at 14 12
\setlinear \plot 10 12 10 6 12 0 14 6 14 12 12 6  10 12  /
\setlinear \plot 12 0 12 6  /
\put{$2{,}520$} [c] at 13 -2
\put{$\scriptstyle \bullet$} [c] at 16  0 \endpicture
\end{minipage}
\begin{minipage}{4cm}
\beginpicture
\setcoordinatesystem units    <1.5mm,2mm>
\setplotarea x from 0 to 16, y from -2 to 15
\put{1.577)} [l] at 2 12
\put {$ \scriptstyle \bullet$} [c] at 10 6
\put {$ \scriptstyle \bullet$} [c] at 10 0
\put {$ \scriptstyle \bullet$} [c] at 12 12
\put {$ \scriptstyle \bullet$} [c] at 12 6
\put {$ \scriptstyle \bullet$} [c] at 14 6
\put {$ \scriptstyle \bullet$} [c] at 14 0
\setlinear \plot 10 0 10 6 12 12 14 6 14 0 12 6  10 0  /
\setlinear \plot 12 12 12 6  /
\put{$2{,}520$} [c] at 13 -2
\put{$\scriptstyle \bullet$} [c] at 16  0 \endpicture
\end{minipage}
\begin{minipage}{4cm}
\beginpicture
\setcoordinatesystem units    <1.5mm,2mm>
\setplotarea x from 0 to 16, y from -2 to 15
\put{1.578)} [l] at 2 12
\put {$ \scriptstyle \bullet$} [c] at 12 12
\put {$ \scriptstyle \bullet$} [c] at 12 0
\put {$ \scriptstyle \bullet$} [c] at 12 4
\put {$ \scriptstyle \bullet$} [c] at 12 8
\put {$ \scriptstyle \bullet$} [c] at 10 12
\put {$ \scriptstyle \bullet$} [c] at 14 12
\setlinear \plot 10 12 12 4 14 12   /
\setlinear \plot 12 0 12 12  /
\put{$2{,}520$} [c] at 13 -2
\put{$\scriptstyle \bullet$} [c] at 16  0
\endpicture
\end{minipage}
$$

$$
\begin{minipage}{4cm}
\beginpicture
\setcoordinatesystem units    <1.5mm,2mm>
\setplotarea x from 0 to 16, y from -2 to 15
\put{1.579)} [l] at 2 12
\put {$ \scriptstyle \bullet$} [c] at 12 0
\put {$ \scriptstyle \bullet$} [c] at 12 12
\put {$ \scriptstyle \bullet$} [c] at 12 4
\put {$ \scriptstyle \bullet$} [c] at 12 8
\put {$ \scriptstyle \bullet$} [c] at 10 0
\put {$ \scriptstyle \bullet$} [c] at 14 0
\setlinear \plot 10 0 12 8 14 0   /
\setlinear \plot 12 0 12 12  /
\put{$2{,}520$} [c] at 13 -2
\put{$\scriptstyle \bullet$} [c] at 16  0 \endpicture
\end{minipage}
\begin{minipage}{4cm}
\beginpicture
\setcoordinatesystem units    <1.5mm,2mm>
\setplotarea x from 0 to 16, y from -2 to 15
\put{1.580)} [l] at 2 12
\put {$ \scriptstyle \bullet$} [c] at 10 6
\put {$ \scriptstyle \bullet$} [c] at 12 12
\put {$ \scriptstyle \bullet$} [c] at 12 0
\put {$ \scriptstyle \bullet$} [c] at 10 12
\put {$ \scriptstyle \bullet$} [c] at 14 6
\put {$ \scriptstyle \bullet$} [c] at 14 12
\setlinear \plot 10 12 10 6 12 0 14 6 14 12    /
\setlinear \plot 10 6 12 12 14 6    /
\put{$2{,}520$} [c] at 13 -2
\put{$\scriptstyle \bullet$} [c] at 16  0 \endpicture
\end{minipage}
\begin{minipage}{4cm}
\beginpicture
\setcoordinatesystem units    <1.5mm,2mm>
\setplotarea x from 0 to 16, y from -2 to 15
\put{1.581)} [l] at 2 12
\put {$ \scriptstyle \bullet$} [c] at 10 6
\put {$ \scriptstyle \bullet$} [c] at 12 12
\put {$ \scriptstyle \bullet$} [c] at 12 0
\put {$ \scriptstyle \bullet$} [c] at 10 0
\put {$ \scriptstyle \bullet$} [c] at 14 6
\put {$ \scriptstyle \bullet$} [c] at 14 0
\setlinear \plot 10 0 10 6 12 12 14 6 14 0    /
\setlinear \plot 10 6 12 0 14 6    /
\put{$2{,}520$} [c] at 13 -2
\put{$\scriptstyle \bullet$} [c] at 16  0 \endpicture
\end{minipage}
\begin{minipage}{4cm}
\beginpicture
\setcoordinatesystem units    <1.5mm,2mm>
\setplotarea x from 0 to 16, y from -2 to 15
\put{1.582)} [l] at 2 12
\put {$ \scriptstyle \bullet$} [c] at 10 6
\put {$ \scriptstyle \bullet$} [c] at 10 12
\put {$ \scriptstyle \bullet$} [c] at 11 0
\put {$ \scriptstyle \bullet$} [c] at 12 12
\put {$ \scriptstyle \bullet$} [c] at 12 6
\put {$ \scriptstyle \bullet$} [c] at 14 0
\put{$\scriptstyle \bullet$} [c] at 16  0
\setlinear \plot 14 0 12 12 10 6 10 12  12  6 11 0 10 6    /
\setlinear \plot 12 6 12 12 /
\put{$2{,}520$} [c] at 13 -2
 \endpicture
\end{minipage}
\begin{minipage}{4cm}
\beginpicture
\setcoordinatesystem units    <1.5mm,2mm>
\setplotarea x from 0 to 16, y from -2 to 15
\put{1.583)} [l] at 2 12
\put {$ \scriptstyle \bullet$} [c] at 10 6
\put {$ \scriptstyle \bullet$} [c] at 10 0
\put {$ \scriptstyle \bullet$} [c] at 11 12
\put {$ \scriptstyle \bullet$} [c] at 12 0
\put {$ \scriptstyle \bullet$} [c] at 12 6
\put {$ \scriptstyle \bullet$} [c] at 14 12
\put{$\scriptstyle \bullet$} [c] at 16  0
\setlinear \plot 14 12 12 0 10 6 10 0  12  6 11 12 10 6    /
\setlinear \plot 12 6 12 0 /
\put{$2{,}520$} [c] at 13 -2
\put{$\scriptstyle \bullet$} [c] at 16  0 \endpicture
\end{minipage}
\begin{minipage}{4cm}
\beginpicture
\setcoordinatesystem units    <1.5mm,2mm>
\setplotarea x from 0 to 16, y from -2 to 15
\put{1.584)} [l] at 2 12
\put {$ \scriptstyle \bullet$} [c] at 10 12
\put {$ \scriptstyle \bullet$} [c] at 12  0
\put {$ \scriptstyle \bullet$} [c] at 12 12
\put {$ \scriptstyle \bullet$} [c] at 13 6
\put {$ \scriptstyle \bullet$} [c] at 14 12
\put {$ \scriptstyle \bullet$} [c] at 14 0
\setlinear \plot 10 12  12 0 13 6 14 12    /
\setlinear \plot 12 12 13 6  14 0     /
\put{$2{,}520$} [c] at 13 -2
\put{$\scriptstyle \bullet$} [c] at 16  0 \endpicture
\end{minipage}
$$
$$
\begin{minipage}{4cm}
\beginpicture
\setcoordinatesystem units    <1.5mm,2mm>
\setplotarea x from 0 to 16, y from -2 to 15
\put{1.585)} [l] at 2 12
\put {$ \scriptstyle \bullet$} [c] at 10 0
\put {$ \scriptstyle \bullet$} [c] at 12  0
\put {$ \scriptstyle \bullet$} [c] at 12 12
\put {$ \scriptstyle \bullet$} [c] at 13 6
\put {$ \scriptstyle \bullet$} [c] at 14 12
\put {$ \scriptstyle \bullet$} [c] at 14 0
\setlinear \plot 10 0  12 12 13 6 14 0    /
\setlinear \plot 12 0 13 6  14 12     /
\put{$2{,}520$} [c] at 13 -2
\put{$\scriptstyle \bullet$} [c] at 16  0 \endpicture
\end{minipage}
\begin{minipage}{4cm}
\beginpicture
\setcoordinatesystem units    <1.5mm,2mm>
\setplotarea x from 0 to 16, y from -2 to 15
\put{1.586)} [l] at 2 12
\put {$ \scriptstyle \bullet$} [c] at 10 0
\put {$ \scriptstyle \bullet$} [c] at 10 12
\put {$ \scriptstyle \bullet$} [c] at 14 12
\put {$ \scriptstyle \bullet$} [c] at 14 6
\put {$ \scriptstyle \bullet$} [c] at 14 0
\put {$ \scriptstyle \bullet$} [c] at 12 12
\setlinear \plot 10 0 10 12  14 0 14 12 10 0 14 12  /
\setlinear \plot 10 0 12  12  14  6    /
\put{$2{,}520$} [c] at 13 -2
\put{$\scriptstyle \bullet$} [c] at 16  0 \endpicture
\end{minipage}
\begin{minipage}{4cm}
\beginpicture
\setcoordinatesystem units    <1.5mm,2mm>
\setplotarea x from 0 to 16, y from -2 to 15
\put{1.587)} [l] at 2 12
\put {$ \scriptstyle \bullet$} [c] at 10 0
\put {$ \scriptstyle \bullet$} [c] at 10 12
\put {$ \scriptstyle \bullet$} [c] at 14 12
\put {$ \scriptstyle \bullet$} [c] at 14 6
\put {$ \scriptstyle \bullet$} [c] at 14 0
\put {$ \scriptstyle \bullet$} [c] at 12 0
\setlinear \plot 10 12 10 0  14 12 14 0 10 12 14 0  /
\setlinear \plot 10 12 12  0  14  6    /
\put{$2{,}520$} [c] at 13 -2
\put{$\scriptstyle \bullet$} [c] at 16  0 \endpicture
\end{minipage}
\begin{minipage}{4cm}
\beginpicture
\setcoordinatesystem units    <1.5mm,2mm>
\setplotarea x from 0 to 16, y from -2 to 15
\put{1.588)} [l] at 2 12
\put {$ \scriptstyle \bullet$} [c] at 10 0
\put {$ \scriptstyle \bullet$} [c] at 10 12
\put {$ \scriptstyle \bullet$} [c] at 12 12
\put {$ \scriptstyle \bullet$} [c] at 14 6
\put {$ \scriptstyle \bullet$} [c] at 14 12
\put {$ \scriptstyle \bullet$} [c] at 14 0
\setlinear \plot 10 0 10 12  14 6 12 12 10 0 /
\setlinear \plot 14 0 14  12     /
\put{$2{,}520$} [c] at 13 -2
\put{$\scriptstyle \bullet$} [c] at 16  0 \endpicture
\end{minipage}
\begin{minipage}{4cm}
\beginpicture
\setcoordinatesystem units    <1.5mm,2mm>
\setplotarea x from 0 to 16, y from -2 to 15
\put{1.589)} [l] at 2 12
\put {$ \scriptstyle \bullet$} [c] at 10 0
\put {$ \scriptstyle \bullet$} [c] at 10 12
\put {$ \scriptstyle \bullet$} [c] at 12 0
\put {$ \scriptstyle \bullet$} [c] at 14 6
\put {$ \scriptstyle \bullet$} [c] at 14 12
\put {$ \scriptstyle \bullet$} [c] at 14 0
\setlinear \plot 10 12 10 0  14 6 12 0 10 12 /
\setlinear \plot 14 0 14  12     /
\put{$2{,}520$} [c] at 13 -2
\put{$\scriptstyle \bullet$} [c] at 16  0 \endpicture
\end{minipage}
\begin{minipage}{4cm}
\beginpicture
\setcoordinatesystem units    <1.5mm,2mm>
\setplotarea x from 0 to 16, y from -2 to 15
\put{1.590)} [l] at 2 12
\put {$ \scriptstyle \bullet$} [c] at 10 4
\put {$ \scriptstyle \bullet$} [c] at 10 8
\put {$ \scriptstyle \bullet$} [c] at 12 0
\put {$ \scriptstyle \bullet$} [c] at 12 6
\put {$ \scriptstyle \bullet$} [c] at 12 12
\put {$ \scriptstyle \bullet$} [c] at 14 6
\setlinear \plot 12 0 10 4 10 8 12 12 14 6 12 0   /
\setlinear \plot 12 0 12 12   /
\put{$2{,}520$} [c] at 13 -2
\put{$\scriptstyle \bullet$} [c] at 16  0 \endpicture
\end{minipage}
$$

$$
\begin{minipage}{4cm}
\beginpicture
\setcoordinatesystem units    <1.5mm,2mm>
\setplotarea x from 0 to 16, y from -2 to 15
\put{1.591)} [l] at 2 12
\put {$ \scriptstyle \bullet$} [c] at 10 6
\put {$ \scriptstyle \bullet$} [c] at 10 0
\put {$ \scriptstyle \bullet$} [c] at 10 12
\put {$ \scriptstyle \bullet$} [c] at 14 12
\put {$ \scriptstyle \bullet$} [c] at 14 6
\put {$ \scriptstyle \bullet$} [c] at 14 0
\setlinear \plot 10 12  10 0 14 12 14 0 10 12    /
\put{$2{,}520$} [c] at 13 -2
\put{$\scriptstyle \bullet$} [c] at 16  0 \endpicture
\end{minipage}
\begin{minipage}{4cm}
\beginpicture
\setcoordinatesystem units    <1.5mm,2mm>
\setplotarea x from 0 to 16, y from -2 to 15
\put{1.592)} [l] at 2 12
\put {$ \scriptstyle \bullet$} [c] at 10 8
\put {$ \scriptstyle \bullet$} [c] at 11 0
\put {$ \scriptstyle \bullet$} [c] at 11 4
\put {$ \scriptstyle \bullet$} [c] at 11 12
\put {$ \scriptstyle \bullet$} [c] at 12 8
\put {$ \scriptstyle \bullet$} [c] at 14 0
\put{$\scriptstyle \bullet$} [c] at 16  0
\setlinear \plot 11 0 11 4 10 8 11 12 12 8 11 4   /
\put{$1{,}260$} [c] at 13 -2
\endpicture
\end{minipage}
\begin{minipage}{4cm}
\beginpicture
\setcoordinatesystem units    <1.5mm,2mm>
\setplotarea x from 0 to 16, y from -2 to 15
\put{1.593)} [l] at 2 12
\put {$ \scriptstyle \bullet$} [c] at 10 4
\put {$ \scriptstyle \bullet$} [c] at 11 0
\put {$ \scriptstyle \bullet$} [c] at 11 8
\put {$ \scriptstyle \bullet$} [c] at 11 12
\put {$ \scriptstyle \bullet$} [c] at 12 4
\put {$ \scriptstyle \bullet$} [c] at 14 0
\put{$\scriptstyle \bullet$} [c] at 16  0
\setlinear \plot 11 12 11 8 10 4 11 0 12 4 11 8   /
\put{$1{,}260$} [c] at 13 -2
\endpicture
\end{minipage}
\begin{minipage}{4cm}
\beginpicture
\setcoordinatesystem units    <1.5mm,2mm>
\setplotarea x from 0 to 16, y from -2 to 15
\put{1.594)} [l] at 2 12
\put {$ \scriptstyle \bullet$} [c] at 10 12
\put {$ \scriptstyle \bullet$} [c] at 11 0
\put {$ \scriptstyle \bullet$} [c] at 11 4
\put {$ \scriptstyle \bullet$} [c] at 11 8
\put {$ \scriptstyle \bullet$} [c] at 12 12
\put {$ \scriptstyle \bullet$} [c] at 14 0
\put{$\scriptstyle \bullet$} [c] at 16  0
\setlinear \plot 11 0 11 8 10  12    /
\setlinear \plot 11 8 12 12    /
\put{$1{,}260$} [c] at 13 -2
\endpicture
\end{minipage}
\begin{minipage}{4cm}
\beginpicture
\setcoordinatesystem units    <1.5mm,2mm>
\setplotarea x from 0 to 16, y from -2 to 15
\put{1.595)} [l] at 2 12
\put {$ \scriptstyle \bullet$} [c] at 10 0
\put {$ \scriptstyle \bullet$} [c] at 11 12
\put {$ \scriptstyle \bullet$} [c] at 11 8
\put {$ \scriptstyle \bullet$} [c] at 11 4
\put {$ \scriptstyle \bullet$} [c] at 12 0
\put {$ \scriptstyle \bullet$} [c] at 14 0
\put{$\scriptstyle \bullet$} [c] at 16  0
\setlinear \plot 11 12 11 4 10  0    /
\setlinear \plot 12 0 11 4   /
\put{$1{,}260$} [c] at 13 -2
\endpicture
\end{minipage}
\begin{minipage}{4cm}
\beginpicture
\setcoordinatesystem units    <1.5mm,2mm>
\setplotarea x from 0 to 16, y from -2 to 15
\put{1.596)} [l] at 2 12
\put {$ \scriptstyle \bullet$} [c] at 10 12
\put {$ \scriptstyle \bullet$} [c] at 10 0
\put {$ \scriptstyle \bullet$} [c] at 12 0
\put {$ \scriptstyle \bullet$} [c] at 12 12
\put {$ \scriptstyle \bullet$} [c] at 14 0
\put {$ \scriptstyle \bullet$} [c] at 14 6
\put {$ \scriptstyle \bullet$} [c] at 14 12
\setlinear \plot 10 0 10 12 12 0  12 12 10 0  /
\setlinear \plot 14 0 14  12    /
\put{$1{,}260$} [c] at 13 -2
\endpicture
\end{minipage}
$$

$$
\begin{minipage}{4cm}
\beginpicture
\setcoordinatesystem units    <1.5mm,2mm>
\setplotarea x from 0 to 16, y from -2 to 15
\put{${\bf  29}$} [l] at 2 15
\put{1.597)} [l] at 2 12

\put {$ \scriptstyle \bullet$} [c] at  10 0
\put {$ \scriptstyle \bullet$} [c] at  10 12
\put {$ \scriptstyle \bullet$} [c] at  12 12
\put {$ \scriptstyle \bullet$} [c] at  14 0
\put {$ \scriptstyle \bullet$} [c] at  14 6
\put {$ \scriptstyle \bullet$} [c] at  14 12
\put {$ \scriptstyle \bullet$} [c] at  16  12
\setlinear \plot  10 12 10 0  12 12 14 6   14 0 /
\setlinear \plot  14 12 14 6 /
\setlinear \plot 16 12 14 0 /
\put{$5{,}040$} [c] at 13 -2
\endpicture
\end{minipage}
\begin{minipage}{4cm}
\beginpicture
\setcoordinatesystem units    <1.5mm,2mm>
\setplotarea x from 0 to 16, y from -2 to 15
\put {1.598)} [l] at 2 12
\put {$ \scriptstyle \bullet$} [c] at  10 0
\put {$ \scriptstyle \bullet$} [c] at  10 12
\put {$ \scriptstyle \bullet$} [c] at  12 0
\put {$ \scriptstyle \bullet$} [c] at  14 0
\put {$ \scriptstyle \bullet$} [c] at  14 6
\put {$ \scriptstyle \bullet$} [c] at  14 12
\put {$ \scriptstyle \bullet$} [c] at  16  0
\setlinear \plot  10 0 10 12  12 0 14 6   14 12 /
\setlinear \plot  14 0 14 6 /
\setlinear \plot 16 0 14 12 /
\put{$5{,}040$} [c] at 13 -2
\endpicture
\end{minipage}
\begin{minipage}{4cm}
\beginpicture
\setcoordinatesystem units    <1.5mm,2mm>
\setplotarea x from  0 to 16, y from -2 to 15
\put{1.599)} [l] at 2 12
\put {$ \scriptstyle \bullet$} [c] at  10 0
\put {$ \scriptstyle \bullet$} [c] at  10 12
\put {$ \scriptstyle \bullet$} [c] at  13 0
\put {$ \scriptstyle \bullet$} [c] at  13 12
\put {$ \scriptstyle \bullet$} [c] at  16  0
\put {$ \scriptstyle \bullet$} [c] at  16  6
\put {$ \scriptstyle \bullet$} [c] at  16  12
\setlinear \plot  10 0  10 12 13 0  13 12 16 0  16 12 /
\put{$5{,}040$} [c] at 13 -2
\endpicture
\end{minipage}
\begin{minipage}{4cm}
\beginpicture
\setcoordinatesystem units    <1.5mm,2mm>
\setplotarea x from  0 to 16, y from -2 to 15
\put{1.600)} [l] at 2 12
\put {$ \scriptstyle \bullet$} [c] at  10 0
\put {$ \scriptstyle \bullet$} [c] at  10 12
\put {$ \scriptstyle \bullet$} [c] at  13 0
\put {$ \scriptstyle \bullet$} [c] at  13 12
\put {$ \scriptstyle \bullet$} [c] at  16  0
\put {$ \scriptstyle \bullet$} [c] at  16  6
\put {$ \scriptstyle \bullet$} [c] at  16  12
\setlinear \plot  10 12  10 0 13 12  13 0 16 12  16 0 /
\put{$5{,}040$} [c] at 13 -2
\endpicture
\end{minipage}
\begin{minipage}{4cm}
\beginpicture
\setcoordinatesystem units    <1.5mm,2mm>
\setplotarea x from 0 to 16, y from -2 to 15
\put {1.601)} [l] at 2 12
\put {$ \scriptstyle \bullet$} [c] at  10 12
\put {$ \scriptstyle \bullet$} [c] at  12 12
\put {$ \scriptstyle \bullet$} [c] at  14 12
\put {$ \scriptstyle \bullet$} [c] at  16 12
\put {$ \scriptstyle \bullet$} [c] at  10 0
\put {$ \scriptstyle \bullet$} [c] at  13 0
\put {$ \scriptstyle \bullet$} [c] at  16 0
\setlinear \plot  10 12 10 0  16 12 16 0   14 12 13 0 12 12 10 0  /
\setlinear \plot  13 0 16 12 /
\put{$5{,}040$} [c] at 13 -2
\endpicture
\end{minipage}
\begin{minipage}{4cm}
\beginpicture
\setcoordinatesystem units    <1.5mm,2mm>
\setplotarea x from 0 to 16, y from -2 to 15
\put {1.602)} [l] at 2 12
\put {$ \scriptstyle \bullet$} [c] at  10 0
\put {$ \scriptstyle \bullet$} [c] at  12 0
\put {$ \scriptstyle \bullet$} [c] at  14 0
\put {$ \scriptstyle \bullet$} [c] at  16 0
\put {$ \scriptstyle \bullet$} [c] at  10 12
\put {$ \scriptstyle \bullet$} [c] at  13 12
\put {$ \scriptstyle \bullet$} [c] at  16 12
\setlinear \plot  10 0 10 12  16 0 16 12   14 0 13 12 12 0 10 12  /
\setlinear \plot  13 12 16 0 /
\put{$5{,}040$} [c] at 13 -2
\endpicture
\end{minipage}
$$
$$
\begin{minipage}{4cm}
\beginpicture
\setcoordinatesystem units    <1.5mm,2mm>
\setplotarea x from 0 to 16, y from -2 to 15
\put{1.603)} [l] at 2 12
\put {$ \scriptstyle \bullet$} [c] at 10 6
\put {$ \scriptstyle \bullet$} [c] at 11 0
\put {$ \scriptstyle \bullet$} [c] at 11 6
\put {$ \scriptstyle \bullet$} [c] at 11 12
\put {$ \scriptstyle \bullet$} [c] at 12 6
\put {$ \scriptstyle \bullet$} [c] at 12 12
\put {$ \scriptstyle \bullet$} [c] at 16 12
\setlinear \plot 16 12 11 0 10 6 11 12 11 0 12 6 11 12   /
\setlinear \plot 12 6 12  12   /
\put{$2{,}520$} [c] at 13 -2
\endpicture
\end{minipage}
\begin{minipage}{4cm}
\beginpicture
\setcoordinatesystem units    <1.5mm,2mm>
\setplotarea x from 0 to 16, y from -2 to 15
\put{1.604)} [l] at 2 12
\put {$ \scriptstyle \bullet$} [c] at 10 6
\put {$ \scriptstyle \bullet$} [c] at 11 0
\put {$ \scriptstyle \bullet$} [c] at 11 6
\put {$ \scriptstyle \bullet$} [c] at 11 12
\put {$ \scriptstyle \bullet$} [c] at 12 6
\put {$ \scriptstyle \bullet$} [c] at 12 0
\put {$ \scriptstyle \bullet$} [c] at 16 0
\setlinear \plot 16 0 11 12 10 6 11 0 11 12 12 6 11 0   /
\setlinear \plot 12 6 12  0   /
\put{$2{,}520$} [c] at 13 -2
\endpicture
\end{minipage}
\begin{minipage}{4cm}
\beginpicture
\setcoordinatesystem units    <1.5mm,2mm>
\setplotarea x from 0 to 16, y from -2 to 15
\put{1.605)} [l] at 2 12
\put {$ \scriptstyle \bullet$} [c] at 10 6
\put {$ \scriptstyle \bullet$} [c] at 10.5 9
\put {$ \scriptstyle \bullet$} [c] at 11 0
\put {$ \scriptstyle \bullet$} [c] at 11 12
\put {$ \scriptstyle \bullet$} [c] at 12 6
\put {$ \scriptstyle \bullet$} [c] at 14.5 12
\put {$ \scriptstyle \bullet$} [c] at 16 12
\setlinear \plot 16 12 11 0 14.5 12    /
\setlinear \plot 10 6 11 0 12 6 11 12 10 6    /
\put{$2{,}520$} [c] at 13 -2
\endpicture
\end{minipage}
\begin{minipage}{4cm}
\beginpicture
\setcoordinatesystem units    <1.5mm,2mm>
\setplotarea x from 0 to 16, y from -2 to 15
\put{1.606)} [l] at 2 12
\put {$ \scriptstyle \bullet$} [c] at 10 6
\put {$ \scriptstyle \bullet$} [c] at 10.5 3
\put {$ \scriptstyle \bullet$} [c] at 11 0
\put {$ \scriptstyle \bullet$} [c] at 11 12
\put {$ \scriptstyle \bullet$} [c] at 12 6
\put {$ \scriptstyle \bullet$} [c] at 14.5 0
\put {$ \scriptstyle \bullet$} [c] at 16 0
\setlinear \plot 16 0 11 12 14.5 0    /
\setlinear \plot 10 6 11 0 12 6 11 12 10 6    /
\put{$2{,}520$} [c] at 13 -2
\endpicture
\end{minipage}
\begin{minipage}{4cm}
\beginpicture
\setcoordinatesystem units    <1.5mm,2mm>
\setplotarea x from 0 to 16, y from -2 to 15
\put{1.607)} [l] at 2 12
\put {$ \scriptstyle \bullet$} [c] at 10 12
\put {$ \scriptstyle \bullet$} [c] at 11 12
\put {$ \scriptstyle \bullet$} [c] at 11 6
\put {$ \scriptstyle \bullet$} [c] at 12 12
\put {$ \scriptstyle \bullet$} [c] at 12 0
\put {$ \scriptstyle \bullet$} [c] at 13 6
\put {$ \scriptstyle \bullet$} [c] at 16 12
\setlinear \plot 16 12 12 0 13 6 12 12 11 6 12 0   /
\setlinear \plot 10 12 11 6 11 12     /
\put{$2{,}520$} [c] at 13 -2
\endpicture
\end{minipage}
\begin{minipage}{4cm}
\beginpicture
\setcoordinatesystem units    <1.5mm,2mm>
\setplotarea x from 0 to 16, y from -2 to 15
\put{1.608)} [l] at 2 12
\put {$ \scriptstyle \bullet$} [c] at 10 0
\put {$ \scriptstyle \bullet$} [c] at 11 0
\put {$ \scriptstyle \bullet$} [c] at 11 6
\put {$ \scriptstyle \bullet$} [c] at 12 12
\put {$ \scriptstyle \bullet$} [c] at 12 0
\put {$ \scriptstyle \bullet$} [c] at 13 6
\put {$ \scriptstyle \bullet$} [c] at 16 0
\setlinear \plot 16 0 12 12 13 6 12 0 11 6 12 12   /
\setlinear \plot 10 0 11 6 11 0     /
\put{$2{,}520$} [c] at 13 -2
\endpicture
\end{minipage}
$$

$$
\begin{minipage}{4cm}
\beginpicture
\setcoordinatesystem units    <1.5mm,2mm>
\setplotarea x from 0 to 16, y from -2 to 15
\put{1.609)} [l] at 2 12
\put {$ \scriptstyle \bullet$} [c] at 10 12
\put {$ \scriptstyle \bullet$} [c] at 13 12
\put {$ \scriptstyle \bullet$} [c] at 14  12
\put {$ \scriptstyle \bullet$} [c] at 16 12
\put {$ \scriptstyle \bullet$} [c] at 11.5 6
\put {$ \scriptstyle \bullet$} [c] at 11.5 0
\put {$ \scriptstyle \bullet$} [c] at 10.7 9
\setlinear \plot 13 12 11.5 6 11.5 0 14 12    /
\setlinear \plot 11.5 0 16 12    /
\setlinear \plot 10 12 11.5 6   /
\put{$2{,}520$} [c] at 13 -2
\endpicture
\end{minipage}
\begin{minipage}{4cm}
\beginpicture
\setcoordinatesystem units    <1.5mm,2mm>
\setplotarea x from 0 to 16, y from -2 to 15
\put{1.610)} [l] at 2 12
\put {$ \scriptstyle \bullet$} [c] at 10 0
\put {$ \scriptstyle \bullet$} [c] at 13 0
\put {$ \scriptstyle \bullet$} [c] at 14  0
\put {$ \scriptstyle \bullet$} [c] at 16 0
\put {$ \scriptstyle \bullet$} [c] at 11.5 6
\put {$ \scriptstyle \bullet$} [c] at 11.5 12
\put {$ \scriptstyle \bullet$} [c] at 10.7 3
\setlinear \plot 13 0 11.5 6 11.5 12 14 0    /
\setlinear \plot 11.5 12 16 0    /
\setlinear \plot 10 0 11.5 6   /
\put{$2{,}520$} [c] at 13 -2
\endpicture
\end{minipage}
\begin{minipage}{4cm}
\beginpicture
\setcoordinatesystem units    <1.5mm,2mm>
\setplotarea x from  0 to 16, y from -2 to 15
\put{1.611)} [l] at 2 12
\put {$ \scriptstyle \bullet$} [c] at  10 0
\put {$ \scriptstyle \bullet$} [c] at  12.5 6
\put {$ \scriptstyle \bullet$} [c] at  13 0
\put {$ \scriptstyle \bullet$} [c] at  13 12
\put {$ \scriptstyle \bullet$} [c] at  13.5 6
\put {$ \scriptstyle \bullet$} [c] at  16 6
\put {$ \scriptstyle \bullet$} [c] at  16 12
\setlinear \plot 10 0 13 12 12.5 6 13 0 16 6 16 12  /
\setlinear \plot 13 0 13.5 6 13 12  /
\put{$2{,}520  $} [c] at 13 -2
\endpicture
\end{minipage}
\begin{minipage}{4cm}
\beginpicture
\setcoordinatesystem units    <1.5mm,2mm>
\setplotarea x from  0 to 16, y from -2 to 15
\put{1.612)} [l] at 2 12
\put {$ \scriptstyle \bullet$} [c] at  10 12
\put {$ \scriptstyle \bullet$} [c] at  12.5 6
\put {$ \scriptstyle \bullet$} [c] at  13 0
\put {$ \scriptstyle \bullet$} [c] at  13 12
\put {$ \scriptstyle \bullet$} [c] at  13.5 6
\put {$ \scriptstyle \bullet$} [c] at  16 6
\put {$ \scriptstyle \bullet$} [c] at  16 0
\setlinear \plot 10 12 13 0 12.5 6 13 12 16 6 16 0  /
\setlinear \plot 13 12 13.5 6 13 0  /
\put{$2{,}520  $} [c] at 13 -2
\endpicture
\end{minipage}
\begin{minipage}{4cm}
\beginpicture
\setcoordinatesystem units    <1.5mm,2mm>
\setplotarea x from 0 to 16, y from -2 to 15
\put{1.613)} [l] at 2 12
\put {$ \scriptstyle \bullet$} [c] at  10 0
\put {$ \scriptstyle \bullet$} [c] at  10 12
\put {$ \scriptstyle \bullet$} [c] at  11.5 6
\put {$ \scriptstyle \bullet$} [c] at  13 12
\put {$ \scriptstyle \bullet$} [c] at  14.5 6
\put {$ \scriptstyle \bullet$} [c] at  16 12
\put {$ \scriptstyle \bullet$} [c] at  16 0
\setlinear \plot 10 12 10 0  13 12 16 0 16  12   /
\put{$2{,}520$} [c] at 13 -2
\endpicture
\end{minipage}
\begin{minipage}{4cm}
\beginpicture
\setcoordinatesystem units    <1.5mm,2mm>
\setplotarea x from 0 to 16, y from -2 to 15
\put{1.614)} [l] at 2 12
\put {$ \scriptstyle \bullet$} [c] at  10 0
\put {$ \scriptstyle \bullet$} [c] at  10 12
\put {$ \scriptstyle \bullet$} [c] at  11.5 6
\put {$ \scriptstyle \bullet$} [c] at  13 0
\put {$ \scriptstyle \bullet$} [c] at  14.5 6
\put {$ \scriptstyle \bullet$} [c] at  16 12
\put {$ \scriptstyle \bullet$} [c] at  16 0
\setlinear \plot 10 0 10 12  13 0 16 12 16  0   /
\put{$2{,}520$} [c] at 13 -2
\endpicture
\end{minipage}
$$
$$
\begin{minipage}{4cm}
\beginpicture
\setcoordinatesystem units    <1.5mm,2mm>
\setplotarea x from  0 to 16, y from -2 to 15
\put{1.615)} [l] at 2 12
\put {$ \scriptstyle \bullet$} [c] at  10 12
\put {$ \scriptstyle \bullet$} [c] at  10 6
\put {$ \scriptstyle \bullet$} [c] at  11.5 0
\put {$ \scriptstyle \bullet$} [c] at  13 12
\put {$ \scriptstyle \bullet$} [c] at  13 6
\put {$ \scriptstyle \bullet$} [c] at  16 0
\put {$ \scriptstyle \bullet$} [c] at  16 12
\setlinear \plot  10 12 10 6 11.5 0 13 6 13 12    /
\setlinear \plot  11.5 0 16 12 16 0   /
\put{$2{,}520$} [c] at 13 -2
\endpicture
\end{minipage}
\begin{minipage}{4cm}
\beginpicture
\setcoordinatesystem units    <1.5mm,2mm>
\setplotarea x from  0 to 16, y from -2 to 15
\put{1.616)} [l] at 2 12
\put {$ \scriptstyle \bullet$} [c] at  10 0
\put {$ \scriptstyle \bullet$} [c] at  10 6
\put {$ \scriptstyle \bullet$} [c] at  11.5 12
\put {$ \scriptstyle \bullet$} [c] at  13 0
\put {$ \scriptstyle \bullet$} [c] at  13 6
\put {$ \scriptstyle \bullet$} [c] at  16 0
\put {$ \scriptstyle \bullet$} [c] at  16 12
\setlinear \plot  10 0 10 6 11.5  12 13 6 13 0    /
\setlinear \plot  11.5 12 16 0 16 12   /
\put{$2{,}520$} [c] at 13 -2
\endpicture
\end{minipage}
\begin{minipage}{4cm}
\beginpicture
\setcoordinatesystem units    <1.5mm,2mm>
\setplotarea x from 0 to 16, y from -2 to 15
\put {1.617)} [l] at 2 12
\put {$\scriptstyle \bullet$} [c] at  10 12
\put {$ \scriptstyle \bullet$} [c] at  12 12
\put {$ \scriptstyle \bullet$} [c] at  14 12
\put {$ \scriptstyle \bullet$} [c] at  16 12
\put {$ \scriptstyle \bullet$} [c] at  12 0
\put {$ \scriptstyle \bullet$} [c] at  14 6
\put {$ \scriptstyle \bullet$} [c] at  14 0
\setlinear \plot  10 12 12 0  12 12 14 0  16 12 /
\setlinear \plot  12 0  14  6 /
\setlinear \plot  14 12  14  0 /
\put{$2{,}520$} [c] at 13 -2
\endpicture
\end{minipage}
\begin{minipage}{4cm}
\beginpicture
\setcoordinatesystem units    <1.5mm,2mm>
\setplotarea x from 0 to 16, y from -2 to 15
\put {1.618)} [l] at 2 12
\put {$ \scriptstyle \bullet$} [c] at  10 0
\put {$ \scriptstyle \bullet$} [c] at  12 0
\put {$ \scriptstyle \bullet$} [c] at  14 0
\put {$ \scriptstyle \bullet$} [c] at  16 0
\put {$ \scriptstyle \bullet$} [c] at  12 12
\put {$ \scriptstyle \bullet$} [c] at  14 6
\put {$ \scriptstyle \bullet$} [c] at  14 12
\setlinear \plot  10 0 12 12  12 0 14 12  16 0 /
\setlinear \plot  12 12  14  6 /
\setlinear \plot  14 12  14  0 /
\put{$2{,}520$} [c] at 13 -2
\endpicture
\end{minipage}
\begin{minipage}{4cm}
\beginpicture
\setcoordinatesystem units    <1.5mm,2mm>
\setplotarea x from 0 to 16, y from -2 to 15
\put {1.619)} [l] at 2 12
\put {$ \scriptstyle \bullet$} [c] at  10 12
\put {$ \scriptstyle \bullet$} [c] at  12 12
\put {$ \scriptstyle \bullet$} [c] at  14 12
\put {$ \scriptstyle \bullet$} [c] at  16 12
\put {$ \scriptstyle \bullet$} [c] at  10 0
\put {$ \scriptstyle \bullet$} [c] at  14 6
\put {$ \scriptstyle \bullet$} [c] at  14  0
\setlinear \plot  10 12 10 0 16 12 14 0 14 12 10 0 12 12 14 0 /
\put{$2{,}520$} [c] at 13 -2
\endpicture
\end{minipage}
\begin{minipage}{4cm}
\beginpicture
\setcoordinatesystem units    <1.5mm,2mm>
\setplotarea x from 0 to 16, y from -2 to 15
\put {1.620)} [l] at 2 12
\put {$ \scriptstyle \bullet$} [c] at  10 0
\put {$ \scriptstyle \bullet$} [c] at  12 0
\put {$ \scriptstyle \bullet$} [c] at  14 0
\put {$ \scriptstyle \bullet$} [c] at  16 0
\put {$ \scriptstyle \bullet$} [c] at  10 12
\put {$ \scriptstyle \bullet$} [c] at  14 6
\put {$ \scriptstyle \bullet$} [c] at  14  12
\setlinear \plot  10 0 10 12 16 0 14 12 14 0 10 12 12 0 14 12 /
\put{$2{,}520$} [c] at 13 -2
\endpicture
\end{minipage}
$$

$$
\begin{minipage}{4cm}
\beginpicture
\setcoordinatesystem units    <1.5mm,2mm>
\setplotarea x from  0 to 16, y from -2 to 15
\put{1.621)} [l] at 2 12
\put {$ \scriptstyle \bullet$} [c] at  10 0
\put {$ \scriptstyle \bullet$} [c] at  10 12
\put {$ \scriptstyle \bullet$} [c] at  13 12
\put {$ \scriptstyle \bullet$} [c] at  13 0
\put {$ \scriptstyle \bullet$} [c] at  14.5  6
\put {$ \scriptstyle \bullet$} [c] at  16  12
\put {$ \scriptstyle \bullet$} [c] at  16  0
\setlinear \plot 13 0 10 12 10 0 13  12 13  0 16 12 16 0  /
\put{$2{,}520$} [c] at 13 -2
\endpicture
\end{minipage}
\begin{minipage}{4cm}
\beginpicture
\setcoordinatesystem units    <1.5mm,2mm>
\setplotarea x from  0 to 16, y from -2 to 15
\put{1.622)} [l] at 2 12
\put {$ \scriptstyle \bullet$} [c] at  10 0
\put {$ \scriptstyle \bullet$} [c] at  10 12
\put {$ \scriptstyle \bullet$} [c] at  13 12
\put {$ \scriptstyle \bullet$} [c] at  13 0
\put {$ \scriptstyle \bullet$} [c] at  14.5  6
\put {$ \scriptstyle \bullet$} [c] at  16  12
\put {$ \scriptstyle \bullet$} [c] at  16  0
\setlinear \plot 13 12 10 0 10 12 13  0 13  12 16 0 16 12  /
\put{$2{,}520$} [c] at 13 -2
\endpicture
\end{minipage}
\begin{minipage}{4cm}
\beginpicture
\setcoordinatesystem units    <1.5mm,2mm>
\setplotarea x from  0 to 16, y from -2 to 15
\put{1.623)} [l] at 2 12
\put {$ \scriptstyle \bullet$} [c] at  10 6
\put {$ \scriptstyle \bullet$} [c] at  11 0
\put {$ \scriptstyle \bullet$} [c] at  11 12
\put {$ \scriptstyle \bullet$} [c] at  12 6
\put {$ \scriptstyle \bullet$} [c] at  16 0
\put {$ \scriptstyle \bullet$} [c] at  16 12
\put {$ \scriptstyle \bullet$} [c] at  14 12
\setlinear \plot 16 12 16 0  11 12 10 6 11 0 12 6 11  12 /
\setlinear \plot  16 0 14 12  /
\put{$1{,}260$} [c] at 13 -2
\endpicture
\end{minipage}
\begin{minipage}{4cm}
\beginpicture
\setcoordinatesystem units    <1.5mm,2mm>
\setplotarea x from  0 to 16, y from -2 to 15
\put{1.624)} [l] at 2 12
\put {$ \scriptstyle \bullet$} [c] at  10 6
\put {$ \scriptstyle \bullet$} [c] at  11 0
\put {$ \scriptstyle \bullet$} [c] at  11 12
\put {$ \scriptstyle \bullet$} [c] at  12 6
\put {$ \scriptstyle \bullet$} [c] at  16 0
\put {$ \scriptstyle \bullet$} [c] at  16 12
\put {$ \scriptstyle \bullet$} [c] at  14 0
\setlinear \plot 16 0 16 12  11 0 10 6 11 12 12 6 11  0 /
\setlinear \plot  16 12 14 0  /
\put{$1{,}260$} [c] at 13 -2
\endpicture
\end{minipage}
\begin{minipage}{4cm}
\beginpicture
\setcoordinatesystem units    <1.5mm,2mm>
\setplotarea x from  0 to 16, y from -2 to 15
\put{1.625)} [l] at 2 12
\put {$ \scriptstyle \bullet$} [c] at  10 0
\put {$ \scriptstyle \bullet$} [c] at  10 12
\put {$ \scriptstyle \bullet$} [c] at  13 0
\put {$ \scriptstyle \bullet$} [c] at  13 12
\put {$ \scriptstyle \bullet$} [c] at  16 0
\put {$ \scriptstyle \bullet$} [c] at  16  6
\put {$ \scriptstyle \bullet$} [c] at  16  12
\setlinear \plot  10  12 10 0 16 12 16 0  /
\setlinear \plot  10 0 13  12      /
\setlinear \plot  13 0 16  6      /
\put{$1{,}260$} [c] at 13 -2
\endpicture
\end{minipage}
\begin{minipage}{4cm}
\beginpicture
\setcoordinatesystem units    <1.5mm,2mm>
\setplotarea x from  0 to 16, y from -2 to 15
\put{1.626)} [l] at 2 12
\put {$ \scriptstyle \bullet$} [c] at  10 0
\put {$ \scriptstyle \bullet$} [c] at  10 12
\put {$ \scriptstyle \bullet$} [c] at  13 0
\put {$ \scriptstyle \bullet$} [c] at  13 12
\put {$ \scriptstyle \bullet$} [c] at  16 0
\put {$ \scriptstyle \bullet$} [c] at  16  6
\put {$ \scriptstyle \bullet$} [c] at  16  12
\setlinear \plot  10  0 10 12 16 0 16 12  /
\setlinear \plot  10 12 13  0      /
\setlinear \plot  13 12 16  6      /
\put{$1{,}260$} [c] at 13 -2
\endpicture
\end{minipage}
$$
$$
\begin{minipage}{4cm}
\beginpicture
\setcoordinatesystem units    <1.5mm,2mm>
\setplotarea x from  0 to 16, y from -2 to 15
\put{1.627)} [l] at 2 12
\put {$ \scriptstyle \bullet$} [c] at  10 12
\put {$ \scriptstyle \bullet$} [c] at  12 12
\put {$ \scriptstyle \bullet$} [c] at  14 12
\put {$ \scriptstyle \bullet$} [c] at  16 12
\put {$ \scriptstyle \bullet$} [c] at  10 0
\put {$ \scriptstyle \bullet$} [c] at  13 0
\put {$ \scriptstyle \bullet$} [c] at  16  0
\setlinear \plot  10  12 10 0 14 12 16 0 16 12 /
\setlinear \plot  10 0 12  12 13 0 14 12     /
\setlinear \plot  10 12 13 0      /
\put{$1{,}260$} [c] at 13 -2
\endpicture
\end{minipage}
\begin{minipage}{4cm}
\beginpicture
\setcoordinatesystem units    <1.5mm,2mm>
\setplotarea x from  0 to 16, y from -2 to 15
\put{1.628)} [l] at 2 12
\put {$ \scriptstyle \bullet$} [c] at  10 0
\put {$ \scriptstyle \bullet$} [c] at  12 0
\put {$ \scriptstyle \bullet$} [c] at  14 0
\put {$ \scriptstyle \bullet$} [c] at  16 0
\put {$ \scriptstyle \bullet$} [c] at  10 12
\put {$ \scriptstyle \bullet$} [c] at  13 12
\put {$ \scriptstyle \bullet$} [c] at  16  12
\setlinear \plot  10  0 10 12 14 0 16 12 16 0 /
\setlinear \plot  10 12 12  0 13 12 14 0     /
\setlinear \plot  10 0 13 12      /
\put{$1{,}260$} [c] at 13 -2
\endpicture
\end{minipage}
\begin{minipage}{4cm}
\beginpicture
\setcoordinatesystem units    <1.5mm,2mm>
\setplotarea x from  0 to 16, y from -2 to 15
\put{1.629)} [l] at 2 12
\put {$ \scriptstyle \bullet$} [c] at  10 6
\put {$ \scriptstyle \bullet$} [c] at  11 0
\put {$ \scriptstyle \bullet$} [c] at  11 6
\put {$ \scriptstyle \bullet$} [c] at  11 12
\put {$ \scriptstyle \bullet$} [c] at  12 6
\put {$ \scriptstyle \bullet$} [c] at  16 0
\put {$ \scriptstyle \bullet$} [c] at  16 12
\setlinear \plot 11 0 10 6 11 12 11 0 12 6 11 12   /
\setlinear \plot 11 0 16 12  16 0 /
\put{$840$} [c] at 13 -2
\endpicture
\end{minipage}
\begin{minipage}{4cm}
\beginpicture
\setcoordinatesystem units    <1.5mm,2mm>
\setplotarea x from  0 to 16, y from -2 to 15
\put{1.630)} [l] at 2 12
\put {$ \scriptstyle \bullet$} [c] at  10 6
\put {$ \scriptstyle \bullet$} [c] at  11 0
\put {$ \scriptstyle \bullet$} [c] at  11 6
\put {$ \scriptstyle \bullet$} [c] at  11 12
\put {$ \scriptstyle \bullet$} [c] at  12 6
\put {$ \scriptstyle \bullet$} [c] at  16 0
\put {$ \scriptstyle \bullet$} [c] at  16 12
\setlinear \plot 11 0 10 6 11 12 11 0 12 6 11 12   /
\setlinear \plot 11 12 16 0  16 12 /
\put{$840 $} [c] at 13 -2
\endpicture
\end{minipage}
\begin{minipage}{4cm}
\beginpicture
\setcoordinatesystem units    <1.5mm,2mm>
\setplotarea x from 0 to 16, y from -2 to 15
\put {1.631)} [l] at 2 12
\put {$ \scriptstyle \bullet$} [c] at  10 12
\put {$ \scriptstyle \bullet$} [c] at  12 12
\put {$ \scriptstyle \bullet$} [c] at  14 12
\put {$ \scriptstyle \bullet$} [c] at  16 12
\put {$ \scriptstyle \bullet$} [c] at  10 0
\put {$ \scriptstyle \bullet$} [c] at  14 6
\put {$ \scriptstyle \bullet$} [c] at  14 0
\setlinear \plot  10 0 10 12  14 0 14 12 /
\setlinear \plot  12 12 14 6 16 12  /
\put{$840$} [c] at 13 -2
\endpicture
\end{minipage}
\begin{minipage}{4cm}
\beginpicture
\setcoordinatesystem units    <1.5mm,2mm>
\setplotarea x from 0 to 16, y from -2 to 15
\put {1.632)} [l] at 2 12
\put {$ \scriptstyle \bullet$} [c] at  10 0
\put {$ \scriptstyle \bullet$} [c] at  12 0
\put {$ \scriptstyle \bullet$} [c] at  14 0
\put {$ \scriptstyle \bullet$} [c] at  16 0
\put {$ \scriptstyle \bullet$} [c] at  10 12
\put {$ \scriptstyle \bullet$} [c] at  14 6
\put {$ \scriptstyle \bullet$} [c] at  14 12
\setlinear \plot  10 12 10 0  14 12 14 0 /
\setlinear \plot  12 0 14 6 16 0  /
\put{$840$} [c] at 13 -2
\endpicture
\end{minipage}
$$

$$
\begin{minipage}{4cm}
\beginpicture
\setcoordinatesystem units    <1.5mm,2mm>
\setplotarea x from 0 to 16, y from -2 to 15
\put {1.633)} [l] at 2 12
\put {$ \scriptstyle \bullet$} [c] at  10 12
\put {$ \scriptstyle \bullet$} [c] at  12 12
\put {$ \scriptstyle \bullet$} [c] at  14 12
\put {$ \scriptstyle \bullet$} [c] at  16 12
\put {$ \scriptstyle \bullet$} [c] at  12 0
\put {$ \scriptstyle \bullet$} [c] at  14 0
\put {$ \scriptstyle \bullet$} [c] at  11 6
\setlinear \plot  10 12  12 0 16 12  14 0 14 12 12 0   /
\setlinear \plot  12 0  12 12 14 0  /
\put{$840$} [c] at 13 -2
\endpicture
\end{minipage}
\begin{minipage}{4cm}
\beginpicture
\setcoordinatesystem units    <1.5mm,2mm>
\setplotarea x from 0 to 16, y from -2 to 15
\put {1.634)} [l] at 2 12
\put {$ \scriptstyle \bullet$} [c] at  10 0
\put {$ \scriptstyle \bullet$} [c] at  12 0
\put {$ \scriptstyle \bullet$} [c] at  14 0
\put {$ \scriptstyle \bullet$} [c] at  16 0
\put {$ \scriptstyle \bullet$} [c] at  12 12
\put {$ \scriptstyle \bullet$} [c] at  14 12
\put {$ \scriptstyle \bullet$} [c] at  11 6
\setlinear \plot  10 0  12 12 16 0  14 12 14 0 12 12   /
\setlinear \plot  12 12  12 0 14 12  /
\put{$840$} [c] at 13 -2
\endpicture
\end{minipage}
\begin{minipage}{4cm}
\beginpicture
\setcoordinatesystem units    <1.5mm,2mm>
\setplotarea x from 0 to 16, y from -2 to 15
\put{1.635)} [l] at 2 12
\put {$ \scriptstyle \bullet$} [c] at 10 6
\put {$ \scriptstyle \bullet$} [c] at 10 12
\put {$ \scriptstyle \bullet$} [c] at 13 0
\put {$ \scriptstyle \bullet$} [c] at 12 6
\put {$ \scriptstyle \bullet$} [c] at 12 12
\put {$ \scriptstyle \bullet$} [c] at 14 12
\put {$ \scriptstyle \bullet$} [c] at  16 12
\setlinear \plot 12 6 13 0 10 6 10 12 12 6 12 12 10 6   /
\setlinear \plot 14 12 13 0 16 12 /
\put{$630$} [c] at 13 -2
\endpicture
\end{minipage}
\begin{minipage}{4cm}
\beginpicture
\setcoordinatesystem units    <1.5mm,2mm>
\setplotarea x from 0 to 16, y from -2 to 15
\put{1.636)} [l] at 2 12
\put {$ \scriptstyle \bullet$} [c] at 10 6
\put {$ \scriptstyle \bullet$} [c] at 10 0
\put {$ \scriptstyle \bullet$} [c] at 13 12
\put {$ \scriptstyle \bullet$} [c] at 12 6
\put {$ \scriptstyle \bullet$} [c] at 12 0
\put {$ \scriptstyle \bullet$} [c] at 14 0
\put {$ \scriptstyle \bullet$} [c] at  16 0
\setlinear \plot 12 6 13 12 10 6 10 0 12 6 12 0 10 6   /
\setlinear \plot 14 0 13 12 16 0 /
\put{$630$} [c] at 13 -2
\endpicture
\end{minipage}
\begin{minipage}{4cm}
\beginpicture
\setcoordinatesystem units    <1.5mm,2mm>
\setplotarea x from 0 to 16, y from -2 to 15
\put{1.637)} [l] at 2 12
\put {$ \scriptstyle \bullet$} [c] at 10 12
\put {$ \scriptstyle \bullet$} [c] at 12 12
\put {$ \scriptstyle \bullet$} [c] at 14 12
\put {$ \scriptstyle \bullet$} [c] at 16 12
\put {$ \scriptstyle \bullet$} [c] at 10 0
\put {$ \scriptstyle \bullet$} [c] at 13 0
\put {$ \scriptstyle \bullet$} [c] at  16 0
\setlinear \plot 10 0 10 12 13  0 14 12 16 0 16 12 13 0 12 12 10 0   /
\put{$630$} [c] at 13 -2
\endpicture
\end{minipage}
\begin{minipage}{4cm}
\beginpicture
\setcoordinatesystem units    <1.5mm,2mm>
\setplotarea x from 0 to 16, y from -2 to 15
\put{1.638)} [l] at 2 12
\put {$ \scriptstyle \bullet$} [c] at 10 0
\put {$ \scriptstyle \bullet$} [c] at 12 0
\put {$ \scriptstyle \bullet$} [c] at 14 0
\put {$ \scriptstyle \bullet$} [c] at 16 0
\put {$ \scriptstyle \bullet$} [c] at 10 12
\put {$ \scriptstyle \bullet$} [c] at 13 12
\put {$ \scriptstyle \bullet$} [c] at  16 12
\setlinear \plot 10 12 10 0 13  12 14 0 16 12 16 0 13 12 12 0 10 12   /
\put{$630$} [c] at 13 -2
\endpicture
\end{minipage}
$$

$$
\begin{minipage}{4cm}
\beginpicture
\setcoordinatesystem units    <1.5mm,2mm>
\setplotarea x from 0 to 16, y from -2 to 15
\put{${\bf  30}$} [l] at 2 12

\put{1.639)} [l] at 2 12
\put {$ \scriptstyle \bullet$} [c] at  10 12
\put {$ \scriptstyle \bullet$} [c] at  11 12
\put {$ \scriptstyle \bullet$} [c] at  12.5 6
\put {$ \scriptstyle \bullet$} [c] at  13 0
\put {$ \scriptstyle \bullet$} [c] at  13 12
\put {$ \scriptstyle \bullet$} [c] at  13.5 6
\put {$ \scriptstyle \bullet$} [c] at  16 0
\setlinear \plot  16 0 13  12  12.5 6  13 0 13.5 6 13 12 /
\setlinear \plot  10 12 13 0  /
\setlinear \plot  11 12 12.5 6   /
\put{$5{,}040$} [c] at 13 -2
\endpicture
\end{minipage}
\begin{minipage}{4cm}
\beginpicture
\setcoordinatesystem units    <1.5mm,2mm>
\setplotarea x from 0 to 16, y from -2 to 15
\put{1.640)} [l] at 2 12
\put {$ \scriptstyle \bullet$} [c] at  10 0
\put {$ \scriptstyle \bullet$} [c] at  11 0
\put {$ \scriptstyle \bullet$} [c] at  12.5 6
\put {$ \scriptstyle \bullet$} [c] at  13 0
\put {$ \scriptstyle \bullet$} [c] at  13 12
\put {$ \scriptstyle \bullet$} [c] at  13.5 6
\put {$ \scriptstyle \bullet$} [c] at  16 12
\setlinear \plot  16 12 13  0  12.5 6  13 12 13.5 6 13 0 /
\setlinear \plot  10 0 13 12  /
\setlinear \plot  11 0 12.5 6   /
\put{$5{,}040$} [c] at 13 -2
\endpicture
\end{minipage}
\begin{minipage}{4cm}
\beginpicture
\setcoordinatesystem units    <1.5mm,2mm>
\setplotarea x from 0 to 16, y from -2 to 15
\put{1.641)} [l] at 2 12
\put {$ \scriptstyle \bullet$} [c] at  10 12
\put {$ \scriptstyle \bullet$} [c] at  10 0
\put {$ \scriptstyle \bullet$} [c] at  13 12
\put {$ \scriptstyle \bullet$} [c] at  13 0
\put {$ \scriptstyle \bullet$} [c] at  16 6
\put {$ \scriptstyle \bullet$} [c] at  16 12
\put {$ \scriptstyle \bullet$} [c] at  16 0
\setlinear \plot  10 12 10 0 13 12 16 0 16 12 /
\setlinear \plot  13 0 13 12 /
\put{$5{,}040 $} [c] at 13 -2
\endpicture
\end{minipage}
\begin{minipage}{4cm}
\beginpicture
\setcoordinatesystem units    <1.5mm,2mm>
\setplotarea x from 0 to 16, y from -2 to 15
\put{1.642)} [l] at 2 12
\put {$ \scriptstyle \bullet$} [c] at  10 12
\put {$ \scriptstyle \bullet$} [c] at  10 0
\put {$ \scriptstyle \bullet$} [c] at  13 12
\put {$ \scriptstyle \bullet$} [c] at  13 0
\put {$ \scriptstyle \bullet$} [c] at  16 6
\put {$ \scriptstyle \bullet$} [c] at  16 12
\put {$ \scriptstyle \bullet$} [c] at  16 0
\setlinear \plot  10 0 10 12 13 0 16 12 16 0 /
\setlinear \plot  13 0 13 12 /
\put{$5{,}040 $} [c] at 13 -2
\endpicture
\end{minipage}
\begin{minipage}{4cm}
\beginpicture
\setcoordinatesystem units    <1.5mm,2mm>
\setplotarea x from 0 to 16, y from -2 to 15
\put{1.643)} [l] at 2 12
\put {$ \scriptstyle \bullet$} [c] at  10 0
\put {$ \scriptstyle \bullet$} [c] at  10 12
\put {$ \scriptstyle \bullet$} [c] at  13 0
\put {$ \scriptstyle \bullet$} [c] at  13 12
\put {$ \scriptstyle \bullet$} [c] at  16 0
\put {$ \scriptstyle \bullet$} [c] at  16 12
\put {$ \scriptstyle \bullet$} [c] at  14.5 6
\setlinear \plot  10 0 10  12  13 0 13 12 16 0 16 12 /
\setlinear \plot  10 12 16 0  /
\put{$5{,}040$} [c] at 13 -2
\endpicture
\end{minipage}
\begin{minipage}{4cm}
\beginpicture
\setcoordinatesystem units    <1.5mm,2mm>
\setplotarea x from 0 to 16, y from -2 to 15
\put{1.644)} [l] at 2 12
\put {$ \scriptstyle \bullet$} [c] at  10 0
\put {$ \scriptstyle \bullet$} [c] at  10 12
\put {$ \scriptstyle \bullet$} [c] at  13 0
\put {$ \scriptstyle \bullet$} [c] at  13 12
\put {$ \scriptstyle \bullet$} [c] at  16 0
\put {$ \scriptstyle \bullet$} [c] at  16 12
\put {$ \scriptstyle \bullet$} [c] at  14.5 6
\setlinear \plot  10 12 10  0  13 12 13 0 16 12 16 0 /
\setlinear \plot  10 0 16 12  /
\put{$5{,}040$} [c] at 13 -2
\endpicture
\end{minipage}
$$
$$
\begin{minipage}{4cm}
\beginpicture
\setcoordinatesystem units    <1.5mm,2mm>
\setplotarea x from 0 to 16, y from -2 to 15
\put{1.645)} [l] at 2 12
\put {$ \scriptstyle \bullet$} [c] at  10 0
\put {$ \scriptstyle \bullet$} [c] at  10 12
\put {$ \scriptstyle \bullet$} [c] at  13 12
\put {$ \scriptstyle \bullet$} [c] at  13 6
\put {$ \scriptstyle \bullet$} [c] at  13 0
\put {$ \scriptstyle \bullet$} [c] at  16  12
\put {$ \scriptstyle \bullet$} [c] at  16  0
\setlinear \plot  10 12 10 0  13 12 13 0  16 12  16 0   /
\put{$5{,}040$} [c] at 12 -2
\endpicture
\end{minipage}
\begin{minipage}{4cm}
\beginpicture
\setcoordinatesystem units    <1.5mm,2mm>
\setplotarea x from 0 to 16, y from -2 to 15
\put{1.646)} [l] at 2 12
\put {$ \scriptstyle \bullet$} [c] at  10 0
\put {$ \scriptstyle \bullet$} [c] at  10 6
\put {$ \scriptstyle \bullet$} [c] at  11 0
\put {$ \scriptstyle \bullet$} [c] at  11 12
\put {$ \scriptstyle \bullet$} [c] at  12 6
\put {$ \scriptstyle \bullet$} [c] at  15 12
\put {$ \scriptstyle \bullet$} [c] at  16 12
\setlinear \plot  10 0 10 6 11 12 12 6  11 0 10 6  /
\setlinear \plot  15 12 11 0 16 12 /
\put{$2{,}520$} [c] at 13 -2
\endpicture
\end{minipage}
\begin{minipage}{4cm}
\beginpicture
\setcoordinatesystem units    <1.5mm,2mm>
\setplotarea x from 0 to 16, y from -2 to 15
\put{1.647)} [l] at 2 12
\put {$ \scriptstyle \bullet$} [c] at  10 12
\put {$ \scriptstyle \bullet$} [c] at  10 6
\put {$ \scriptstyle \bullet$} [c] at  11 0
\put {$ \scriptstyle \bullet$} [c] at  11 12
\put {$ \scriptstyle \bullet$} [c] at  12 6
\put {$ \scriptstyle \bullet$} [c] at  15 0
\put {$ \scriptstyle \bullet$} [c] at  16 0
\setlinear \plot  10 12 10 6 11 0 12 6  11 12 10 6  /
\setlinear \plot  15 0 11 12 16 0 /
\put{$2{,}520$} [c] at 13 -2
\endpicture
\end{minipage}
\begin{minipage}{4cm}
\beginpicture
\setcoordinatesystem units    <1.5mm,2mm>
\setplotarea x from 0 to 16, y from -2 to 15
\put{1.648)} [l] at 2 12
\put {$ \scriptstyle \bullet$} [c] at  10 12
\put {$ \scriptstyle \bullet$} [c] at  10 0
\put {$ \scriptstyle \bullet$} [c] at  12.5 6
\put {$ \scriptstyle \bullet$} [c] at  13 12
\put {$ \scriptstyle \bullet$} [c] at  13 0
\put {$ \scriptstyle \bullet$} [c] at  13.5 6
\put {$ \scriptstyle \bullet$} [c] at  16 12
\setlinear \plot  13 12 10 0 10 12 13 0 12.5 6 13 12 13.5 6 13 0 16 12 /
\put{$2{,}520$} [c] at 13 -2
\endpicture
\end{minipage}
\begin{minipage}{4cm}
\beginpicture
\setcoordinatesystem units    <1.5mm,2mm>
\setplotarea x from 0 to 16, y from -2 to 15
\put{1.649)} [l] at 2 12
\put {$ \scriptstyle \bullet$} [c] at  10 12
\put {$ \scriptstyle \bullet$} [c] at  10 0
\put {$ \scriptstyle \bullet$} [c] at  12.5 6
\put {$ \scriptstyle \bullet$} [c] at  13 12
\put {$ \scriptstyle \bullet$} [c] at  13 0
\put {$ \scriptstyle \bullet$} [c] at  13.5 6
\put {$ \scriptstyle \bullet$} [c] at  16 0
\setlinear \plot  13 0 10 12 10 0 13 12 12.5 6 13 0 13.5 6 13 12 16 0 /
\put{$2{,}520$} [c] at 13 -2
\endpicture
\end{minipage}
\begin{minipage}{4cm}
\beginpicture
\setcoordinatesystem units    <1.5mm,2mm>
\setplotarea x from 0 to 16, y from -2 to 15
\put{1.650)} [l] at 2 12
\put {$ \scriptstyle \bullet$} [c] at  10 12
\put {$ \scriptstyle \bullet$} [c] at  10 0
\put {$ \scriptstyle \bullet$} [c] at  14 0
\put {$ \scriptstyle \bullet$} [c] at  14 12
\put {$ \scriptstyle \bullet$} [c] at  16 12
\put {$ \scriptstyle \bullet$} [c] at  11.5 7.6
\put {$ \scriptstyle \bullet$} [c] at  12.5 4.4
\setlinear \plot 10  0 10 12 14 0  16 12   /
\setlinear \plot  14 0 14 12    /
\put{$2{,}520$} [c] at 13 -2
\endpicture
\end{minipage}
$$

$$
\begin{minipage}{4cm}
\beginpicture
\setcoordinatesystem units    <1.5mm,2mm>
\setplotarea x from 0 to 16, y from -2 to 15
\put{1.651)} [l] at 2 12
\put {$ \scriptstyle \bullet$} [c] at  10 12
\put {$ \scriptstyle \bullet$} [c] at  10 0
\put {$ \scriptstyle \bullet$} [c] at  14 0
\put {$ \scriptstyle \bullet$} [c] at  14 12
\put {$ \scriptstyle \bullet$} [c] at  16 0
\put {$ \scriptstyle \bullet$} [c] at  11.5 4.4
\put {$ \scriptstyle \bullet$} [c] at  12.5 7.6
\setlinear \plot 10  12 10 0 14 12  16 0   /
\setlinear \plot  14 0 14 12    /
\put{$2{,}520$} [c] at 13 -2
\endpicture
\end{minipage}
\begin{minipage}{4cm}
\beginpicture
\setcoordinatesystem units    <1.5mm,2mm>
\setplotarea x from 0 to 16, y from -2 to 15
\put {1.652)} [l] at 2 12
\put {$ \scriptstyle \bullet$} [c] at  10 12
\put {$ \scriptstyle \bullet$} [c] at  12 12
\put {$ \scriptstyle \bullet$} [c] at  14 12
\put {$ \scriptstyle \bullet$} [c] at  16 12
\put {$ \scriptstyle \bullet$} [c] at  10 0
\put {$ \scriptstyle \bullet$} [c] at  12 6
\put {$ \scriptstyle \bullet$} [c] at  12  0
\setlinear \plot 12 6 10 0 10 12 12 0 16 12 /
\setlinear \plot  14 12 12 0 12  12 /
\put{$2{,}520$} [c] at 13 -2
\endpicture
\end{minipage}
\begin{minipage}{4cm}
\beginpicture
\setcoordinatesystem units    <1.5mm,2mm>
\setplotarea x from 0 to 16, y from -2 to 15
\put {1.653)} [l] at 2 12
\put {$ \scriptstyle \bullet$} [c] at  10 0
\put {$ \scriptstyle \bullet$} [c] at  12 0
\put {$ \scriptstyle \bullet$} [c] at  14 0
\put {$ \scriptstyle \bullet$} [c] at  16 0
\put {$ \scriptstyle \bullet$} [c] at  10 12
\put {$ \scriptstyle \bullet$} [c] at  12 6
\put {$ \scriptstyle \bullet$} [c] at  12  12
\setlinear \plot 12 6 10 12 10 0 12 12 16 0 /
\setlinear \plot  14 0 12 12 12  0 /
\put{$2{,}520$} [c] at 13 -2
\endpicture
\end{minipage}
\begin{minipage}{4cm}
\beginpicture
\setcoordinatesystem units    <1.5mm,2mm>
\setplotarea x from 0 to 16, y from -2 to 15
\put {1.654)} [l] at 2 12
\put {$ \scriptstyle \bullet$} [c] at  10 12
\put {$ \scriptstyle \bullet$} [c] at  12 12
\put {$ \scriptstyle \bullet$} [c] at  14 12
\put {$ \scriptstyle \bullet$} [c] at  16 12
\put {$ \scriptstyle \bullet$} [c] at  10 0
\put {$ \scriptstyle \bullet$} [c] at  14 0
\put {$ \scriptstyle \bullet$} [c] at  15  6
\setlinear \plot  10 12 10 0 12 12 14 0 16 12 /
\setlinear  \plot 10 0 14 12 14 0 /
\put{$2{,}520$} [c] at 13 -2
\endpicture
\end{minipage}
\begin{minipage}{4cm}
\beginpicture
\setcoordinatesystem units    <1.5mm,2mm>
\setplotarea x from 0 to 16, y from -2 to 15
\put {1.655)} [l] at 2 12
\put {$ \scriptstyle \bullet$} [c] at  10 0
\put {$ \scriptstyle \bullet$} [c] at  12 0
\put {$ \scriptstyle \bullet$} [c] at  14 0
\put {$ \scriptstyle \bullet$} [c] at  16 0
\put {$ \scriptstyle \bullet$} [c] at  10 12
\put {$ \scriptstyle \bullet$} [c] at  14 12
\put {$ \scriptstyle \bullet$} [c] at  15  6
\setlinear \plot  10 0 10 12 12 0 14 12 16 0 /
\setlinear  \plot 10 12 14 0 14 12 /
\put{$2{,}520$} [c] at 13 -2
\endpicture
\end{minipage}
\begin{minipage}{4cm}
\beginpicture
\setcoordinatesystem units    <1.5mm,2mm>
\setplotarea x from 0 to 16, y from -2 to 15
\put {1.656)} [l] at 2 12
\put {$ \scriptstyle \bullet$} [c] at  10 12
\put {$ \scriptstyle \bullet$} [c] at  12 12
\put {$ \scriptstyle \bullet$} [c] at  14 12
\put {$ \scriptstyle \bullet$} [c] at  16 12
\put {$ \scriptstyle \bullet$} [c] at  14 6
\put {$ \scriptstyle \bullet$} [c] at  14 0
\put {$ \scriptstyle \bullet$} [c] at  10 0
\setlinear \plot  10 0 10 12  14 6 14 0 16 12 /
\setlinear \plot  12 12 14 6 14 12  /
\put{$2{,}520$} [c] at 13 -2
\endpicture
\end{minipage}
$$
$$
\begin{minipage}{4cm}
\beginpicture
\setcoordinatesystem units    <1.5mm,2mm>
\setplotarea x from 0 to 16, y from -2 to 15
\put {1.657)} [l] at 2 12
\put {$ \scriptstyle \bullet$} [c] at  10 0
\put {$ \scriptstyle \bullet$} [c] at  12 0
\put {$ \scriptstyle \bullet$} [c] at  14 0
\put {$ \scriptstyle \bullet$} [c] at  16 0
\put {$ \scriptstyle \bullet$} [c] at  14 6
\put {$ \scriptstyle \bullet$} [c] at  14 12
\put {$ \scriptstyle \bullet$} [c] at  10 12
\setlinear \plot  10 12 10 0  14 6 14 12 16 0 /
\setlinear \plot  12 0 14 6 14 0  /
\put{$2{,}520$} [c] at 13 -2
\endpicture
\end{minipage}
\begin{minipage}{4cm}
\beginpicture
\setcoordinatesystem units    <1.5mm,2mm>
\setplotarea x from 0 to 16, y from -2 to 15
\put {1.658)} [l] at 2 12
\put {$ \scriptstyle \bullet$} [c] at  10 12
\put {$ \scriptstyle \bullet$} [c] at  12 12
\put {$ \scriptstyle \bullet$} [c] at  14 12
\put {$ \scriptstyle \bullet$} [c] at  16 12
\put {$ \scriptstyle \bullet$} [c] at  14 6
\put {$ \scriptstyle \bullet$} [c] at  14 0
\put {$ \scriptstyle \bullet$} [c] at  10 0
\setlinear \plot  10 0 10 12  14 0  16 12 /
\setlinear \plot  10 0 12 12 14 0 14 12 10 0  /
\put{$2{,}520$} [c] at 13 -2
\endpicture
\end{minipage}
\begin{minipage}{4cm}
\beginpicture
\setcoordinatesystem units    <1.5mm,2mm>
\setplotarea x from 0 to 16, y from -2 to 15
\put {1.659)} [l] at 2 12
\put {$ \scriptstyle \bullet$} [c] at  10 0
\put {$ \scriptstyle \bullet$} [c] at  12 0
\put {$ \scriptstyle \bullet$} [c] at  14 0
\put {$ \scriptstyle \bullet$} [c] at  16 0
\put {$ \scriptstyle \bullet$} [c] at  14 6
\put {$ \scriptstyle \bullet$} [c] at  14 12
\put {$ \scriptstyle \bullet$} [c] at  10 12
\setlinear \plot  10 12 10 0  14 12  16 0 /
\setlinear \plot  10 12 12 0 14 12 14 0 10 12  /
\put{$2{,}520$} [c] at 13 -2
\endpicture
\end{minipage}
\begin{minipage}{4cm}
\beginpicture
\setcoordinatesystem units    <1.5mm,2mm>
\setplotarea x from 0 to 16, y from -2 to 15
\put{1.660)} [l] at 2 12
\put {$ \scriptstyle \bullet$} [c] at  10 0
\put {$ \scriptstyle \bullet$} [c] at  10 12
\put {$ \scriptstyle \bullet$} [c] at  13 0
\put {$ \scriptstyle \bullet$} [c] at  13 12
\put {$ \scriptstyle \bullet$} [c] at  13 6
\put {$ \scriptstyle \bullet$} [c] at  16  12
\put {$ \scriptstyle \bullet$} [c] at  16  0
\setlinear \plot  10  12 10 0 13 6 13 12 16 0  /
\setlinear \plot  16 12 13 0 13 6     /
\put{$2{,}520$} [c] at 13 -2
\endpicture
\end{minipage}
\begin{minipage}{4cm}
\beginpicture
\setcoordinatesystem units    <1.5mm,2mm>
\setplotarea x from 0 to 16, y from -2 to 15
\put{1.661)} [l] at 2 12
\put {$ \scriptstyle \bullet$} [c] at  10 0
\put {$ \scriptstyle \bullet$} [c] at  10 12
\put {$ \scriptstyle \bullet$} [c] at  13 0
\put {$ \scriptstyle \bullet$} [c] at  13 12
\put {$ \scriptstyle \bullet$} [c] at  13 6
\put {$ \scriptstyle \bullet$} [c] at  16  12
\put {$ \scriptstyle \bullet$} [c] at  16  0
\setlinear \plot  10  0 10 12 13 6 13 0 16 12  /
\setlinear \plot  16 0 13 12 13 6     /
\put{$2{,}520$} [c] at 13 -2
\endpicture
\end{minipage}
\begin{minipage}{4cm}
\beginpicture
\setcoordinatesystem units    <1.5mm,2mm>
\setplotarea x from 0 to 16, y from -2 to 15
\put {1.662)} [l] at 2 12
\put {$ \scriptstyle \bullet$} [c] at  10 12
\put {$ \scriptstyle \bullet$} [c] at  12 12
\put {$ \scriptstyle \bullet$} [c] at  14 12
\put {$ \scriptstyle \bullet$} [c] at  16  12
\put {$ \scriptstyle \bullet$} [c] at  10 0
\put {$ \scriptstyle \bullet$} [c] at  13 0
\put {$ \scriptstyle \bullet$} [c] at  16 0
\setlinear \plot   10 12 10 0  16 12 16 0 14  12  13 0 12 12  /
\setlinear \plot  10 12 13 0   16  12  /
\put{$2{,}520$} [c] at 13 -2
\endpicture
\end{minipage}
$$

$$
\begin{minipage}{4cm}
\beginpicture
\setcoordinatesystem units    <1.5mm,2mm>
\setplotarea x from 0 to 16, y from -2 to 15
\put {1.663)} [l] at 2 12
\put {$ \scriptstyle \bullet$} [c] at  10 0
\put {$ \scriptstyle \bullet$} [c] at  12 0
\put {$ \scriptstyle \bullet$} [c] at  14 0
\put {$ \scriptstyle \bullet$} [c] at  16  0
\put {$ \scriptstyle \bullet$} [c] at  10 12
\put {$ \scriptstyle \bullet$} [c] at  13 12
\put {$ \scriptstyle \bullet$} [c] at  16 12
\setlinear \plot   10 0 10 12  16 0 16 12 14  0  13 12 12 0  /
\setlinear \plot  10 0 13 12   16  0  /
\put{$2{,}520$} [c] at 13 -2
\endpicture
\end{minipage}
\begin{minipage}{4cm}
\beginpicture
\setcoordinatesystem units    <1.5mm,2mm>
\setplotarea x from 0 to 16, y from -2 to 15
\put {1.664)} [l] at 2 12
\put {$ \scriptstyle \bullet$} [c] at  10 0
\put {$ \scriptstyle \bullet$} [c] at  10 12
\put {$ \scriptstyle \bullet$} [c] at  12 12
\put {$ \scriptstyle \bullet$} [c] at  12 6
\put {$ \scriptstyle \bullet$} [c] at  14 0
\put {$ \scriptstyle \bullet$} [c] at  14 12
\put {$ \scriptstyle \bullet$} [c] at  16  12
\setlinear \plot  10 0 10 12 12 6 14 0 16  12  /
\setlinear \plot  12 6 12 12 10 0 /
\setlinear \plot  14  0 14 12 /
\put{$1{,}260$} [c] at 13 -2
\endpicture
\end{minipage}
\begin{minipage}{4cm}
\beginpicture
\setcoordinatesystem units    <1.5mm,2mm>
\setplotarea x from 0 to 16, y from -2 to 15
\put {1.665)} [l] at 2 12
\put {$ \scriptstyle \bullet$} [c] at  10 0
\put {$ \scriptstyle \bullet$} [c] at  10 12
\put {$ \scriptstyle \bullet$} [c] at  12 0
\put {$ \scriptstyle \bullet$} [c] at  12 6
\put {$ \scriptstyle \bullet$} [c] at  14 0
\put {$ \scriptstyle \bullet$} [c] at  14 12
\put {$ \scriptstyle \bullet$} [c] at  16  0
\setlinear \plot  10 12 10 0 12 6 14 12 16  0  /
\setlinear \plot  12 6 12 0 10 12 /
\setlinear \plot  14  0 14 12 /
\put{$1{,}260$} [c] at 13 -2
\endpicture
\end{minipage}
\begin{minipage}{4cm}
\beginpicture
\setcoordinatesystem units    <1.5mm,2mm>
\setplotarea x from 0 to 16, y from -2 to 15
\put {1.666)} [l] at 2 12
\put {$ \scriptstyle \bullet$} [c] at  10 0
\put {$ \scriptstyle \bullet$} [c] at  10 12
\put {$ \scriptstyle \bullet$} [c] at  12 12
\put {$ \scriptstyle \bullet$} [c] at  14 0
\put {$ \scriptstyle \bullet$} [c] at  14 12
\put {$ \scriptstyle \bullet$} [c] at  16 12
\put {$ \scriptstyle \bullet$} [c] at  16  0
\setlinear \plot   10 12 10 0 16 12 16 0  12 12 /
\setlinear \plot 10  0 12 12 14 0 16  12 /
\setlinear \plot 14  12 14 0  /
\put{$1{,}260$} [c] at 13 -2
\endpicture
\end{minipage}
\begin{minipage}{4cm}
\beginpicture
\setcoordinatesystem units    <1.5mm,2mm>
\setplotarea x from 0 to 16, y from -2 to 15
\put {1.667)} [l] at 2 12
\put {$ \scriptstyle \bullet$} [c] at  10 0
\put {$ \scriptstyle \bullet$} [c] at  10 12
\put {$ \scriptstyle \bullet$} [c] at  12 0
\put {$ \scriptstyle \bullet$} [c] at  14 0
\put {$ \scriptstyle \bullet$} [c] at  14 12
\put {$ \scriptstyle \bullet$} [c] at  16 12
\put {$ \scriptstyle \bullet$} [c] at  16  0
\setlinear \plot   10 0 10 12 16 0 16 12 12 0 /
\setlinear \plot 10  12 12 0 14 12 16  0 /
\setlinear \plot 14  12 14 0  /
\put{$1{,}260$} [c] at 13 -2
\endpicture
\end{minipage}
\begin{minipage}{4cm}
\beginpicture
\setcoordinatesystem units    <1.5mm,2mm>
\setplotarea x from 0 to 16, y from -2 to 15
\put {1.668)} [l] at 2 12
\put {$ \scriptstyle \bullet$} [c] at  10 0
\put {$ \scriptstyle \bullet$} [c] at  10 12
\put {$ \scriptstyle \bullet$} [c] at  12 12
\put {$ \scriptstyle \bullet$} [c] at  14 12
\put {$ \scriptstyle \bullet$} [c] at  14 0
\put {$ \scriptstyle \bullet$} [c] at  16 12
\put {$ \scriptstyle \bullet$} [c] at  16  0
\setlinear \plot 16 12 14 0  10 12 10 0  16 12 16 0 /
\setlinear \plot  10 0  12 12 14 0 14 12 10 0  /
\put{$420$} [c] at 13 -2
\endpicture
\end{minipage}
$$
$$
\begin{minipage}{4cm}
\beginpicture
\setcoordinatesystem units    <1.5mm,2mm>
\setplotarea x from 0 to 16, y from -2 to 15
\put {1.669)} [l] at 2 12
\put {$ \scriptstyle \bullet$} [c] at  10 0
\put {$ \scriptstyle \bullet$} [c] at  10 12
\put {$ \scriptstyle \bullet$} [c] at  12 0
\put {$ \scriptstyle \bullet$} [c] at  14 12
\put {$ \scriptstyle \bullet$} [c] at  14 0
\put {$ \scriptstyle \bullet$} [c] at  16 12
\put {$ \scriptstyle \bullet$} [c] at  16  0
\setlinear \plot 16 0 14 12  10 0 10 12  16 0 16 12 /
\setlinear \plot  10 12  12 0 14 12 14 0 10 12  /
\put{$420$} [c] at 13 -2
\endpicture
\end{minipage}
\begin{minipage}{4cm}
\beginpicture
\setcoordinatesystem units    <1.5mm,2mm>
\setplotarea x from 0 to 16, y from -2 to 15
\put{1.670)} [l] at 2 12
\put {$ \scriptstyle \bullet$} [c] at 10 6
\put {$ \scriptstyle \bullet$} [c] at 11 12
\put {$ \scriptstyle \bullet$} [c] at 11 0
\put {$ \scriptstyle \bullet$} [c] at 12 6
\put {$ \scriptstyle \bullet$} [c] at 10.5 3
\put {$ \scriptstyle \bullet$} [c] at 15 12
\put{$\scriptstyle \bullet$} [c] at 16  0
\setlinear \plot 11 0 10 6  11 12 12 6 11 0 15 12 /
\put{$5{,}040$} [c] at 13 -2
 \endpicture
\end{minipage}
\begin{minipage}{4cm}
\beginpicture
\setcoordinatesystem units    <1.5mm,2mm>
\setplotarea x from 0 to 16, y from -2 to 15
\put{1.671)} [l] at 2 12
\put {$ \scriptstyle \bullet$} [c] at 10 6
\put {$ \scriptstyle \bullet$} [c] at 11 12
\put {$ \scriptstyle \bullet$} [c] at 11 0
\put {$ \scriptstyle \bullet$} [c] at 12 6
\put {$ \scriptstyle \bullet$} [c] at 10.5 3
\put {$ \scriptstyle \bullet$} [c] at 15 0
\put{$\scriptstyle \bullet$} [c] at 16  0
\setlinear \plot 11 12 10 6  11 0 12 6 11 12 15 0 /
\put{$5{,}040$} [c] at 13 -2
\endpicture
\end{minipage}
\begin{minipage}{4cm}
\beginpicture
\setcoordinatesystem units    <1.5mm,2mm>
\setplotarea x from 0 to 16, y from -2 to 15
\put{1.672)} [l] at 2 12
\put {$ \scriptstyle \bullet$} [c] at 10 0
\put {$ \scriptstyle \bullet$} [c] at 12  0
\put {$ \scriptstyle \bullet$} [c] at 12 4
\put {$ \scriptstyle \bullet$} [c] at 12 8
\put {$ \scriptstyle \bullet$} [c] at 12 12
\put {$ \scriptstyle \bullet$} [c] at 14 0
\setlinear \plot 10 0  12 12 12 0  /
\setlinear \plot 12 8  14  0     /
\put{$5{,}040$} [c] at 13 -2
\put{$\scriptstyle \bullet$} [c] at 16  0
\endpicture
\end{minipage}
\begin{minipage}{4cm}
\beginpicture
\setcoordinatesystem units    <1.5mm,2mm>
\setplotarea x from 0 to 16, y from -2 to 15
\put{1.673)} [l] at 2 12
\put {$ \scriptstyle \bullet$} [c] at 10 12
\put {$ \scriptstyle \bullet$} [c] at 12  0
\put {$ \scriptstyle \bullet$} [c] at 12 4
\put {$ \scriptstyle \bullet$} [c] at 12 8
\put {$ \scriptstyle \bullet$} [c] at 12 12
\put {$ \scriptstyle \bullet$} [c] at 14 12
\setlinear \plot 10 12  12 0 12 12  /
\setlinear \plot 12 4  14  12     /
\put{$5{,}040$} [c] at 13 -2
\put{$\scriptstyle \bullet$} [c] at 16  0
\endpicture
\end{minipage}
\begin{minipage}{4cm}
\beginpicture
\setcoordinatesystem units    <1.5mm,2mm>
\setplotarea x from 0 to 16, y from -2 to 15
\put{1.674)} [l] at 2 12
\put {$ \scriptstyle \bullet$} [c] at 10 6
\put {$ \scriptstyle \bullet$} [c] at 11 0
\put {$ \scriptstyle \bullet$} [c] at 11 12
\put {$ \scriptstyle \bullet$} [c] at 12 6
\put {$ \scriptstyle \bullet$} [c] at 14 12
\put {$ \scriptstyle \bullet$} [c] at 14 0
\setlinear \plot 12 6 11 0  10 6 11 12 12 6 14 12 14 0  /
\put{$5{,}040$} [c] at 13 -2
\put{$\scriptstyle \bullet$} [c] at 16  0
\endpicture
\end{minipage}
$$

$$
\begin{minipage}{4cm}
\beginpicture
\setcoordinatesystem units    <1.5mm,2mm>
\setplotarea x from 0 to 16, y from -2 to 15
\put{1.675)} [l] at 2 12
\put {$ \scriptstyle \bullet$} [c] at 10 6
\put {$ \scriptstyle \bullet$} [c] at 11 0
\put {$ \scriptstyle \bullet$} [c] at 11 12
\put {$ \scriptstyle \bullet$} [c] at 12 6
\put {$ \scriptstyle \bullet$} [c] at 14 12
\put {$ \scriptstyle \bullet$} [c] at 14 0
\setlinear \plot 12 6 11 0  10 6 11 12 12 6 14 0 14 12  /
\put{$5{,}040$} [c] at 13 -2
\put{$\scriptstyle \bullet$} [c] at 16  0
\endpicture
\end{minipage}
\begin{minipage}{4cm}
\beginpicture
\setcoordinatesystem units    <1.5mm,2mm>
\setplotarea x from 0 to 16, y from -2 to 15
\put{1.676)} [l] at 2 12
\put {$ \scriptstyle \bullet$} [c] at 10 0
\put {$ \scriptstyle \bullet$} [c] at 10 12
\put {$ \scriptstyle \bullet$} [c] at 12 0
\put {$ \scriptstyle \bullet$} [c] at 14 0
\put {$ \scriptstyle \bullet$} [c] at 14 6
\put {$ \scriptstyle \bullet$} [c] at 14 12
\setlinear \plot 10 12 10  0 14 12 14 0   /
\setlinear \plot 10 12 12 0 14 6     /
\put{$5{,}040$} [c] at 13 -2
\put{$\scriptstyle \bullet$} [c] at 16  0
\endpicture
\end{minipage}
\begin{minipage}{4cm}
\beginpicture
\setcoordinatesystem units    <1.5mm,2mm>
\setplotarea x from 0 to 16, y from -2 to 15
\put{1.677)} [l] at 2 12
\put {$ \scriptstyle \bullet$} [c] at 10 0
\put {$ \scriptstyle \bullet$} [c] at 10 12
\put {$ \scriptstyle \bullet$} [c] at 12 12
\put {$ \scriptstyle \bullet$} [c] at 14 0
\put {$ \scriptstyle \bullet$} [c] at 14 6
\put {$ \scriptstyle \bullet$} [c] at 14 12
\setlinear \plot 10 0 10  12 14 0 14 12   /
\setlinear \plot 10 0 12 12 14 6     /
\put{$5{,}040$} [c] at 13 -2
\put{$\scriptstyle \bullet$} [c] at 16  0
\endpicture
\end{minipage}
\begin{minipage}{4cm}
\beginpicture
\setcoordinatesystem units    <1.5mm,2mm>
\setplotarea x from 0 to  16, y from -2 to 15
\put{1.678)} [l] at 2 12
\put {$ \scriptstyle \bullet$} [c] at 10 0
\put {$ \scriptstyle \bullet$} [c] at 10 12
\put {$ \scriptstyle \bullet$} [c] at 12 0
\put {$ \scriptstyle \bullet$} [c] at 12 6
\put {$ \scriptstyle \bullet$} [c] at 14 12
\put {$ \scriptstyle \bullet$} [c] at 16 0
\put {$ \scriptstyle \bullet$} [c] at 16 12
\setlinear \plot 10 0 10 12    12 6 12 0    /
\setlinear \plot  12 6 14 12    /
\setlinear \plot 16 0 16 12   /
\put{$5{,}040$} [c] at 13 -2
\endpicture
\end{minipage}
\begin{minipage}{4cm}
\beginpicture
\setcoordinatesystem units    <1.5mm,2mm>
\setplotarea x from 0 to  16, y from -2 to 15
\put{1.679)} [l] at 2 12
\put {$ \scriptstyle \bullet$} [c] at 10 0
\put {$ \scriptstyle \bullet$} [c] at 10 12
\put {$ \scriptstyle \bullet$} [c] at 12 12
\put {$ \scriptstyle \bullet$} [c] at 12 6
\put {$ \scriptstyle \bullet$} [c] at 14 0
\put {$ \scriptstyle \bullet$} [c] at 16 0
\put {$ \scriptstyle \bullet$} [c] at 16 12
\setlinear \plot 10 12 10 0    12 6 12 12    /
\setlinear \plot  12 6 14 0    /
\setlinear \plot 16 0 16 12   /
\put{$5{,}040$} [c] at 13 -2
\endpicture
\end{minipage}
\begin{minipage}{4cm}
\beginpicture
\setcoordinatesystem units    <1.5mm,2mm>
\setplotarea x from 0 to  16, y from -2 to 15
\put{1.680)} [l] at 2 12
\put {$ \scriptstyle \bullet$} [c] at 10 0
\put {$ \scriptstyle \bullet$} [c] at 10 6
\put {$ \scriptstyle \bullet$} [c] at 10 12
\put {$ \scriptstyle \bullet$} [c] at 14 0
\put {$ \scriptstyle \bullet$} [c] at 14 12
\put {$ \scriptstyle \bullet$} [c] at 16 12
\setlinear \plot 16 0 16 12   /
\setlinear \plot  10 0  10 12   14 0 14 12 10 0   /
\put {$ \scriptstyle \bullet$} [c] at 16 0
\put{$5{,}040$} [c] at 13 -2
\endpicture
\end{minipage}
$$
$$
\begin{minipage}{4cm}
\beginpicture
\setcoordinatesystem units    <1.5mm,2mm>
\setplotarea x from 0 to 16, y from -2 to 15
\put{1.681)} [l] at 2 12
\put {$ \scriptstyle \bullet$} [c] at 10 0
\put {$ \scriptstyle \bullet$} [c] at 10 6
\put {$ \scriptstyle \bullet$} [c] at 10 12
\put {$ \scriptstyle \bullet$} [c] at 14 0
\put {$ \scriptstyle \bullet$} [c] at 14 6
\put {$ \scriptstyle \bullet$} [c] at 14 12
\setlinear \plot 10 12 10 0 14 12 14 0   /
\put{$5{,}040$} [c] at 13 -2
\put{$\scriptstyle \bullet$} [c] at 16  0
\endpicture
\end{minipage}
\begin{minipage}{4cm}
\beginpicture
\setcoordinatesystem units    <1.5mm,2mm>
\setplotarea x from 0 to  16, y from -2 to 15
\put{1.682)} [l] at 2 12
\put {$ \scriptstyle \bullet$} [c] at 10 0
\put {$ \scriptstyle \bullet$} [c] at 10 4
\put {$ \scriptstyle \bullet$} [c] at 10 8
\put {$ \scriptstyle \bullet$} [c] at 10 12
\setlinear \plot 10 0  10 12    /
\put {$ \scriptstyle \bullet$} [c] at 13 0
\put {$ \scriptstyle \bullet$} [c] at 13 12
\setlinear \plot 13 0 13 12   /
\put {$ \scriptstyle \bullet$} [c] at 16 0
\put{$5{,}040$} [c] at 13 -2
\endpicture
\end{minipage}
\begin{minipage}{4cm}
\beginpicture
\setcoordinatesystem units    <1.5mm,2mm>
\setplotarea x from 0 to 16, y from -2 to 15
\put{1.683)} [l] at 2 12
\put {$ \scriptstyle \bullet$} [c] at 10 6
\put {$ \scriptstyle \bullet$} [c] at 12 0
\put {$ \scriptstyle \bullet$} [c] at 12 6
\put {$ \scriptstyle \bullet$} [c] at 12 12
\put {$ \scriptstyle \bullet$} [c] at 14 6
\put {$ \scriptstyle \bullet$} [c] at 14 12
\setlinear \plot 14 12 14 6 12 12 10 6 12  0 14 6    /
\setlinear \plot 12 12 12  0     /
\put{$2{,}520$} [c] at 13 -2
\put{$\scriptstyle \bullet$} [c] at 16  0
\endpicture
\end{minipage}
\begin{minipage}{4cm}
\beginpicture
\setcoordinatesystem units    <1.5mm,2mm>
\setplotarea x from 0 to 16, y from -2 to 15
\put{1.684)} [l] at 2 12
\put {$ \scriptstyle \bullet$} [c] at 10 6
\put {$ \scriptstyle \bullet$} [c] at 12 0
\put {$ \scriptstyle \bullet$} [c] at 12 6
\put {$ \scriptstyle \bullet$} [c] at 12 12
\put {$ \scriptstyle \bullet$} [c] at 14 6
\put {$ \scriptstyle \bullet$} [c] at 14 0
\setlinear \plot 14 0 14 6 12 12 10 6 12  0 14 6    /
\setlinear \plot 12 12 12  0     /
\put{$2{,}520$} [c] at 13 -2
\put{$\scriptstyle \bullet$} [c] at 16  0
\endpicture
\end{minipage}
\begin{minipage}{4cm}
\beginpicture
\setcoordinatesystem units    <1.5mm,2mm>
\setplotarea x from 0 to 16, y from -2 to 15
\put{1.685)} [l] at 2 12
\put {$ \scriptstyle \bullet$} [c] at 10 6
\put {$ \scriptstyle \bullet$} [c] at 11 0
\put {$ \scriptstyle \bullet$} [c] at 11 12
\put {$ \scriptstyle \bullet$} [c] at 12 6
\put {$ \scriptstyle \bullet$} [c] at 13 12
\put {$ \scriptstyle \bullet$} [c] at 14 12
\setlinear \plot 13 12 12 6 11 12  10 6 11 0 12 6 14 12    /
\put{$2{,}520$} [c] at 13 -2
\put{$\scriptstyle \bullet$} [c] at 16  0
\endpicture
\end{minipage}
\begin{minipage}{4cm}
\beginpicture
\setcoordinatesystem units    <1.5mm,2mm>
\setplotarea x from 0 to 16, y from -2 to 15
\put{1.686)} [l] at 2 12
\put {$ \scriptstyle \bullet$} [c] at 10 6
\put {$ \scriptstyle \bullet$} [c] at 11 0
\put {$ \scriptstyle \bullet$} [c] at 11 12
\put {$ \scriptstyle \bullet$} [c] at 12 6
\put {$ \scriptstyle \bullet$} [c] at 13 0
\put {$ \scriptstyle \bullet$} [c] at 14 0
\setlinear \plot 13 0 12 6 11 12  10 6 11 0 12 6 14 0    /
\put{$2{,}520$} [c] at 13 -2
\put{$\scriptstyle \bullet$} [c] at 16  0
 \endpicture
\end{minipage}
$$

$$
\begin{minipage}{4cm}
\beginpicture
\setcoordinatesystem units    <1.5mm,2mm>
\setplotarea x from 0 to 16, y from -2 to 15
\put{1.687)} [l] at 2 12
\put {$ \scriptstyle \bullet$} [c] at 10 0
\put {$ \scriptstyle \bullet$} [c] at 10 6
\put {$ \scriptstyle \bullet$} [c] at 10 12
\put {$ \scriptstyle \bullet$} [c] at 11 3
\put {$ \scriptstyle \bullet$} [c] at 14 0
\put {$ \scriptstyle \bullet$} [c] at 14 12
\setlinear \plot 10 12 10 0 14 12 14 0 10 12  /
\put{$2{,}520$} [c] at 13 -2
\put{$\scriptstyle \bullet$} [c] at 16  0 \endpicture
\end{minipage}
\begin{minipage}{4cm}
\beginpicture
\setcoordinatesystem units    <1.5mm,2mm>
\setplotarea x from 0 to 16, y from -2 to 15
\put{1.688)} [l] at 2 12
\put {$ \scriptstyle \bullet$} [c] at 10 0
\put {$ \scriptstyle \bullet$} [c] at 10 6
\put {$ \scriptstyle \bullet$} [c] at 10 12
\put {$ \scriptstyle \bullet$} [c] at 11 9
\put {$ \scriptstyle \bullet$} [c] at 14 0
\put {$ \scriptstyle \bullet$} [c] at 14 12
\setlinear \plot 10 12 10 0 14 12 14 0 10 12  /
\put{$2{,}520$} [c] at 13 -2
\put{$\scriptstyle \bullet$} [c] at 16  0 \endpicture
\end{minipage}
\begin{minipage}{4cm}
\beginpicture
\setcoordinatesystem units    <1.5mm,2mm>
\setplotarea x from 0 to 16, y from -2 to 15
\put{1.689)} [l] at 2 12
\put {$ \scriptstyle \bullet$} [c] at 10 0
\put {$ \scriptstyle \bullet$} [c] at 10 12
\put {$ \scriptstyle \bullet$} [c] at 12 12
\put {$ \scriptstyle \bullet$} [c] at 14 0
\put {$ \scriptstyle \bullet$} [c] at 14 6
\put {$ \scriptstyle \bullet$} [c] at 14 12
\setlinear \plot 10 12 10 0  12 12 14 6 14  12 10 0 /
\setlinear \plot 14 0 14  6     /
\put{$2{,}520$} [c] at 13 -2
\put{$\scriptstyle \bullet$} [c] at 16  0 \endpicture
\end{minipage}
\begin{minipage}{4cm}
\beginpicture
\setcoordinatesystem units    <1.5mm,2mm>
\setplotarea x from 0 to 16, y from -2 to 15
\put{1.690)} [l] at 2 12
\put {$ \scriptstyle \bullet$} [c] at 10 0
\put {$ \scriptstyle \bullet$} [c] at 10 12
\put {$ \scriptstyle \bullet$} [c] at 12 0
\put {$ \scriptstyle \bullet$} [c] at 14 0
\put {$ \scriptstyle \bullet$} [c] at 14 6
\put {$ \scriptstyle \bullet$} [c] at 14 12
\setlinear \plot 10 0 10 12  12 0 14 6 14  0 10 12 /
\setlinear \plot 14 12 14  6     /
\put{$2{,}520$} [c] at 13 -2
\put{$\scriptstyle \bullet$} [c] at 16  0 \endpicture
\end{minipage}
\begin{minipage}{4cm}
\beginpicture
\setcoordinatesystem units    <1.5mm,2mm>
\setplotarea x from 0 to  16, y from -2 to 15
\put{1.691)} [l] at 2 12
\put {$ \scriptstyle \bullet$} [c] at 10 12
\put {$ \scriptstyle \bullet$} [c] at 10 6
\put {$ \scriptstyle \bullet$} [c] at 12 0
\put {$ \scriptstyle \bullet$} [c] at 14 12
\put {$ \scriptstyle \bullet$} [c] at 14 6
\setlinear \plot 10 12 10 6   12 0 14 6 14 12   /
\put {$ \scriptstyle \bullet$} [c] at 16 0
\put {$ \scriptstyle \bullet$} [c] at 16 12
\setlinear \plot 16 0 16 12   /
\put{$2{,}520$} [c] at 13 -2
\endpicture
\end{minipage}
\begin{minipage}{4cm}
\beginpicture
\setcoordinatesystem units    <1.5mm,2mm>
\setplotarea x from 0 to  16, y from -2 to 15
\put{1.692)} [l] at 2 12
\put {$ \scriptstyle \bullet$} [c] at 10 0
\put {$ \scriptstyle \bullet$} [c] at 10 6
\put {$ \scriptstyle \bullet$} [c] at 12 12
\put {$ \scriptstyle \bullet$} [c] at 14 0
\put {$ \scriptstyle \bullet$} [c] at 14 6
\setlinear \plot 10 0 10 6   12 12 14 6 14 0  /
\put {$ \scriptstyle \bullet$} [c] at 16 0
\put {$ \scriptstyle \bullet$} [c] at 16 12
\setlinear \plot 16 0 16 12   /
\put{$2{,}520$} [c] at 13 -2
\endpicture
\end{minipage}
$$
$$
\begin{minipage}{4cm}
\beginpicture
\setcoordinatesystem units    <1.5mm,2mm>
\setplotarea x from 0 to 16, y from -2 to 15
\put{1.693)} [l] at 2 12
\put {$ \scriptstyle \bullet$} [c] at 10 6
\put {$ \scriptstyle \bullet$} [c] at 10 12
\put {$ \scriptstyle \bullet$} [c] at 11 0
\put {$ \scriptstyle \bullet$} [c] at 12 6
\put {$ \scriptstyle \bullet$} [c] at 12 12
\put {$ \scriptstyle \bullet$} [c] at 15 12
\setlinear \plot 15 12 11 0  10 6 10 12 12 6 12 12 10 6  /
\setlinear \plot 11 0 12  6     /
\put{$1{,}260$} [c] at 13 -2
\put{$\scriptstyle \bullet$} [c] at 16  0 \endpicture
\end{minipage}
\begin{minipage}{4cm}
\beginpicture
\setcoordinatesystem units    <1.5mm,2mm>
\setplotarea x from 0 to 16, y from -2 to 15
\put{1.694)} [l] at 2 12
\put {$ \scriptstyle \bullet$} [c] at 10 6
\put {$ \scriptstyle \bullet$} [c] at 10 0
\put {$ \scriptstyle \bullet$} [c] at 11 12
\put {$ \scriptstyle \bullet$} [c] at 12 6
\put {$ \scriptstyle \bullet$} [c] at 12 0
\put {$ \scriptstyle \bullet$} [c] at 15 0
\setlinear \plot 15 0 11 12  10 6 10 0 12 6 12 0 10 6  /
\setlinear \plot 11 12 12  6     /
\put{$1{,}260$} [c] at 13 -2
\put{$\scriptstyle \bullet$} [c] at 16  0 \endpicture
\end{minipage}
\begin{minipage}{4cm}
\beginpicture
\setcoordinatesystem units    <1.5mm,2mm>
\setplotarea x from 0 to 16, y from -2 to 15
\put{1.695)} [l] at 2 12
\put {$ \scriptstyle \bullet$} [c] at 10 0
\put {$ \scriptstyle \bullet$} [c] at 10 12
\put {$ \scriptstyle \bullet$} [c] at 12 0
\put {$ \scriptstyle \bullet$} [c] at 12 12
\put {$ \scriptstyle \bullet$} [c] at 14 6
\put {$ \scriptstyle \bullet$} [c] at 14 12
\setlinear \plot 14 12 14 6 12 0 12 12 10 0 10 12 12 0  /
\setlinear \plot 10 0  14 6 /
\put{$1{,}260$} [c] at 13 -2
\put{$\scriptstyle \bullet$} [c] at 16  0 \endpicture
\end{minipage}
\begin{minipage}{4cm}
\beginpicture
\setcoordinatesystem units    <1.5mm,2mm>
\setplotarea x from 0 to 16, y from -2 to 15
\put{1.696)} [l] at 2 12
\put {$ \scriptstyle \bullet$} [c] at 10 0
\put {$ \scriptstyle \bullet$} [c] at 10 12
\put {$ \scriptstyle \bullet$} [c] at 12 0
\put {$ \scriptstyle \bullet$} [c] at 12 12
\put {$ \scriptstyle \bullet$} [c] at 14 6
\put {$ \scriptstyle \bullet$} [c] at 14 0
\setlinear \plot 14 0 14 6 12 12 12 0 10 12 10 0 12 12  /
\setlinear \plot 10 12  14 6 /
\put{$1{,}260$} [c] at 13 -2
\put{$\scriptstyle \bullet$} [c] at 16  0 \endpicture
\end{minipage}
\begin{minipage}{4cm}
\beginpicture
\setcoordinatesystem units    <1.5mm,2mm>
\setplotarea x from 0 to 16, y from -2 to 15
\put{1.697)} [l] at 2 12
\put {$ \scriptstyle \bullet$} [c] at 10 6
\put {$ \scriptstyle \bullet$} [c] at 11 0
\put {$ \scriptstyle \bullet$} [c] at 11 12
\put {$ \scriptstyle \bullet$} [c] at 12 6
\put {$ \scriptstyle \bullet$} [c] at 14 0
\put {$ \scriptstyle \bullet$} [c] at 15 12
\put{$\scriptstyle \bullet$} [c] at 16  0
\setlinear \plot 10 6  11 12 12 6 11 0 10 6 /
\setlinear \plot 14 0 15  12 16 0    /
\put{$1{,}260$} [c] at 13 -2
\endpicture
\end{minipage}
\begin{minipage}{4cm}
\beginpicture
\setcoordinatesystem units    <1.5mm,2mm>
\setplotarea x from 0 to 16, y from -2 to 15
\put{1.698)} [l] at 2 12
\put {$ \scriptstyle \bullet$} [c] at 10 6
\put {$ \scriptstyle \bullet$} [c] at 11 0
\put {$ \scriptstyle \bullet$} [c] at 11 12
\put {$ \scriptstyle \bullet$} [c] at 12 6
\put {$ \scriptstyle \bullet$} [c] at 14 12
\put {$ \scriptstyle \bullet$} [c] at 15 0
\put{$\scriptstyle \bullet$} [c] at 16  12
\setlinear \plot 10 6  11 12 12 6 11 0 10 6 /
\setlinear \plot 14 12 15  0 16 12    /
\put{$1{,}260$} [c] at 13 -2
\endpicture
\end{minipage}
$$

$$
\begin{minipage}{4cm}
\beginpicture
\setcoordinatesystem units    <1.5mm,2mm>
\setplotarea x from 0 to 16, y from -2 to 15
\put{1.699)} [l] at 2 12
\put {$ \scriptstyle \bullet$} [c] at 10 12
\put {$ \scriptstyle \bullet$} [c] at 12 12
\put {$ \scriptstyle \bullet$} [c] at 11 6
\put {$ \scriptstyle \bullet$} [c] at 11 0
\put {$ \scriptstyle \bullet$} [c] at 14 0
\put {$ \scriptstyle \bullet$} [c] at 15 12
\put{$\scriptstyle \bullet$} [c] at 16  0
\setlinear \plot 10 12  11 6 11 0 /
\setlinear \plot 11 6 12 12 /
\setlinear \plot 14 0 15  12 16 0    /
\put{$1{,}260$} [c] at 13 -2
\endpicture
\end{minipage}
\begin{minipage}{4cm}
\beginpicture
\setcoordinatesystem units    <1.5mm,2mm>
\setplotarea x from 0 to 16, y from -2 to 15
\put{1.700)} [l] at 2 12
\put {$ \scriptstyle \bullet$} [c] at 10 0
\put {$ \scriptstyle \bullet$} [c] at 12 0
\put {$ \scriptstyle \bullet$} [c] at 11 6
\put {$ \scriptstyle \bullet$} [c] at 11 12
\put {$ \scriptstyle \bullet$} [c] at 14 12
\put {$ \scriptstyle \bullet$} [c] at 15 0
\put{$\scriptstyle \bullet$} [c] at 16  12
\setlinear \plot 10 0  11 6 11 12 /
\setlinear \plot 11 6 12 0 /
\setlinear \plot 14 12 15  0 16 12    /
\put{$1{,}260$} [c] at 13 -2
\endpicture
\end{minipage}
\begin{minipage}{4cm}
\beginpicture
\setcoordinatesystem units    <1.5mm,2mm>
\setplotarea x from 0 to 16, y from -2 to 15
\put{1.701)} [l] at 2 12
\put {$ \scriptstyle \bullet$} [c] at 10 12
\put {$ \scriptstyle \bullet$} [c] at 12 12
\put {$ \scriptstyle \bullet$} [c] at 11 6
\put {$ \scriptstyle \bullet$} [c] at 11 0
\put {$ \scriptstyle \bullet$} [c] at 14 12
\put {$ \scriptstyle \bullet$} [c] at 15 0
\put{$\scriptstyle \bullet$} [c] at 16  12
\setlinear \plot 10 12  11 6 11 0 /
\setlinear \plot 11 6 12 12 /
\setlinear \plot 14 12 15  0 16 12    /
\put{$1{,}260$} [c] at 13 -2
\endpicture
\end{minipage}
\begin{minipage}{4cm}
\beginpicture
\setcoordinatesystem units    <1.5mm,2mm>
\setplotarea x from 0 to 16, y from -2 to 15
\put{1.702)} [l] at 2 12
\put {$ \scriptstyle \bullet$} [c] at 10 0
\put {$ \scriptstyle \bullet$} [c] at 12 0
\put {$ \scriptstyle \bullet$} [c] at 11 6
\put {$ \scriptstyle \bullet$} [c] at 11 12
\put {$ \scriptstyle \bullet$} [c] at 14 0
\put {$ \scriptstyle \bullet$} [c] at 15 12
\put{$\scriptstyle \bullet$} [c] at 16  0
\setlinear \plot 10 0  11 6 11 12 /
\setlinear \plot 11 6 12 0 /
\setlinear \plot 14 0 15  12 16 0    /
\put{$1{,}260$} [c] at 13 -2
\endpicture
\end{minipage}
\begin{minipage}{4cm}
\beginpicture
\setcoordinatesystem units    <1.5mm,2mm>
\setplotarea x from 0 to  16, y from -2 to 15
\put{1.703)} [l] at 2 12
\put {$ \scriptstyle \bullet$} [c] at 10 12
\put {$ \scriptstyle \bullet$} [c] at 11 12
\put {$ \scriptstyle \bullet$} [c] at 12 12
\put {$ \scriptstyle \bullet$} [c] at 11 6
\put {$ \scriptstyle \bullet$} [c] at 11 0
\put {$ \scriptstyle \bullet$} [c] at 16 0
\put {$ \scriptstyle \bullet$} [c] at 16 12
\setlinear \plot  11 0  11 12     /
\setlinear \plot  10 12 11 6  12 12     /
\setlinear \plot 16 0 16 12   /
\put{$840$} [c] at 13 -2
\endpicture
\end{minipage}
\begin{minipage}{4cm}
\beginpicture
\setcoordinatesystem units    <1.5mm,2mm>
\setplotarea x from 0 to  16, y from -2 to 15
\put{1.704)} [l] at 2 12
\put {$ \scriptstyle \bullet$} [c] at 10 0
\put {$ \scriptstyle \bullet$} [c] at 11 0
\put {$ \scriptstyle \bullet$} [c] at 12 0
\put {$ \scriptstyle \bullet$} [c] at 11 6
\put {$ \scriptstyle \bullet$} [c] at 11 12
\put {$ \scriptstyle \bullet$} [c] at 16 0
\put {$ \scriptstyle \bullet$} [c] at 16 12
\setlinear \plot  11 0  11 12     /
\setlinear \plot  10 0 11 6  12 0     /
\setlinear \plot 16 0 16 12   /
\put{$840$} [c] at 13 -2
\endpicture
\end{minipage}
$$
$$
\begin{minipage}{4cm}
\beginpicture
\setcoordinatesystem units    <1.5mm,2mm>
\setplotarea x from 0 to  16, y from -2 to 15
\put{1.705)} [l] at 2 12
\put {$ \scriptstyle \bullet$} [c] at 10 6
\put {$ \scriptstyle \bullet$} [c] at 12 0
\put {$ \scriptstyle \bullet$} [c] at 12 6
\put {$ \scriptstyle \bullet$} [c] at 12 12
\put {$ \scriptstyle \bullet$} [c] at 14 6
\setlinear \plot  10 6   12 12   14 6 12 0 10 6    /
\setlinear \plot  12 0 12 12    /
\put {$ \scriptstyle \bullet$} [c] at 16 0
\put {$ \scriptstyle \bullet$} [c] at 16 12
\setlinear \plot 16 0 16 12   /
\put{$840$} [c] at 13 -2
\endpicture
\end{minipage}
\begin{minipage}{4cm}
\beginpicture
\setcoordinatesystem units    <1.5mm,2mm>
\setplotarea x from 0 to 16, y from -2 to 15
\put{1.706)} [l] at 2 12
\put {$ \scriptstyle \bullet$} [c] at 10 0
\put {$ \scriptstyle \bullet$} [c] at 10 12
\put {$ \scriptstyle \bullet$} [c] at 12 12
\put {$ \scriptstyle \bullet$} [c] at 12 0
\put {$ \scriptstyle \bullet$} [c] at 14 0
\put {$ \scriptstyle \bullet$} [c] at 14 12
\setlinear \plot 10 12 10 0  12 12 14 0 14  12 12 0 10 12 14 0 /
\setlinear \plot 10 0 14  12     /
\setlinear \plot 12 0 12  12     /
\put{$140$} [c] at 13 -2
\put{$\scriptstyle \bullet$} [c] at 16  0
\endpicture
\end{minipage}
\begin{minipage}{4cm}
\beginpicture
\setcoordinatesystem units   <1.5mm,2mm>
\setplotarea x from 0 to 16, y from -2 to 15
\put {${\bf  31}$} [l] at 2 15

\put {1.707)} [l] at 2 12
\put {$ \scriptstyle \bullet$} [c] at  10 0
\put {$ \scriptstyle \bullet$} [c] at  10 12
\put {$ \scriptstyle \bullet$} [c] at  12 12
\put {$ \scriptstyle \bullet$} [c] at  14 0
\put {$ \scriptstyle \bullet$} [c] at  14 6
\put {$ \scriptstyle \bullet$} [c] at  14 12
\put {$ \scriptstyle \bullet$} [c] at  16  12
\setlinear \plot  10 12 10 0  14 12 14 0  12 12 10 0 /
\setlinear \plot 16 12 14 0 /
\put{$5{,}040$} [c] at 13 -2
\endpicture
\end{minipage}
\begin{minipage}{4cm}
\beginpicture
\setcoordinatesystem units   <1.5mm,2mm>
\setplotarea x from 0 to 16, y from -2 to 15
\put {1.708)} [l] at 2 12
\put {$ \scriptstyle \bullet$} [c] at  10 0
\put {$ \scriptstyle \bullet$} [c] at  10 12
\put {$ \scriptstyle \bullet$} [c] at  12 0
\put {$ \scriptstyle \bullet$} [c] at  14 0
\put {$ \scriptstyle \bullet$} [c] at  14 6
\put {$ \scriptstyle \bullet$} [c] at  14 12
\put {$ \scriptstyle \bullet$} [c] at  16  0
\setlinear \plot  10 0 10 12  14 0 14 12  12 0 10 12 /
\setlinear \plot 16 0 14 12 /
\put{$5{,}040$} [c] at 13 -2
\endpicture
\end{minipage}
\begin{minipage}{4cm}
\beginpicture
\setcoordinatesystem units   <1.5mm,2mm>
\setplotarea x  from 0 to 16, y from -2 to 15
\put{1.709)} [l] at 2 12
\put {$ \scriptstyle \bullet$} [c] at  10 0
\put {$ \scriptstyle \bullet$} [c] at  10 12
\put {$ \scriptstyle \bullet$} [c] at  13 0
\put {$ \scriptstyle \bullet$} [c] at  13 12
\put {$ \scriptstyle \bullet$} [c] at  14.5 6
\put {$ \scriptstyle \bullet$} [c] at  16 0
\put {$ \scriptstyle \bullet$} [c] at  16 12
\setlinear \plot  10 0 10 12 13  0 13 12 16 0 16 12   /
\put{$5{,}040$} [c] at 13 -2
\endpicture
\end{minipage}
\begin{minipage}{4cm}
\beginpicture
\setcoordinatesystem units   <1.5mm,2mm>
\setplotarea x  from 0 to 16, y from -2 to 15
\put{1.710)} [l] at 2 12
\put {$ \scriptstyle \bullet$} [c] at  10 0
\put {$ \scriptstyle \bullet$} [c] at  10 12
\put {$ \scriptstyle \bullet$} [c] at  13 0
\put {$ \scriptstyle \bullet$} [c] at  13 12
\put {$ \scriptstyle \bullet$} [c] at  14.5 6
\put {$ \scriptstyle \bullet$} [c] at  16 0
\put {$ \scriptstyle \bullet$} [c] at  16 12
\setlinear \plot  10 12 10 0 13  12 13 0 16 12 16 0   /
\put{$5{,}040$} [c] at 13 -2
\endpicture
\end{minipage}
$$

$$
\begin{minipage}{4cm}
\beginpicture
\setcoordinatesystem units   <1.5mm,2mm>
\setplotarea x from 0 to 16, y from -2 to 15
\put{1.711)} [l] at 2 12
\put {$ \scriptstyle \bullet$} [c] at 10 6
\put {$ \scriptstyle \bullet$} [c] at 10.5 0
\put {$ \scriptstyle \bullet$} [c] at 10.5 12
\put {$ \scriptstyle \bullet$} [c] at 11 6
\put {$ \scriptstyle \bullet$} [c] at 12 6
\put {$ \scriptstyle \bullet$} [c] at 12 12
\put {$ \scriptstyle \bullet$} [c] at 16 12
\setlinear \plot 16 12 10.5 0 10 6 10.5 12 11 6 10.5 0 12 6 12 12   /
\put{$2{,}520$} [c] at 13 -2
\endpicture
\end{minipage}
\begin{minipage}{4cm}
\beginpicture
\setcoordinatesystem units   <1.5mm,2mm>
\setplotarea x from 0 to 16, y from -2 to 15
\put{1.712)} [l] at 2 12
\put {$ \scriptstyle \bullet$} [c] at 10 6
\put {$ \scriptstyle \bullet$} [c] at 10.5 0
\put {$ \scriptstyle \bullet$} [c] at 10.5 12
\put {$ \scriptstyle \bullet$} [c] at 11 6
\put {$ \scriptstyle \bullet$} [c] at 12 6
\put {$ \scriptstyle \bullet$} [c] at 12 0
\put {$ \scriptstyle \bullet$} [c] at 16 0
\setlinear \plot 16 0 10.5 12 10 6 10.5 0 11 6 10.5 12 12 6 12 0   /
\put{$2{,}520$} [c] at 13 -2
\endpicture
\end{minipage}
\begin{minipage}{4cm}
\beginpicture
\setcoordinatesystem units   <1.5mm,2mm>
\setplotarea x from 0 to 16, y from -2 to 15
\put{1.713)} [l] at 2 12
\put {$ \scriptstyle \bullet$} [c] at 10 12
\put {$ \scriptstyle \bullet$} [c] at 12 12
\put {$ \scriptstyle \bullet$} [c] at 13 12
\put {$ \scriptstyle \bullet$} [c] at 16 12
\put {$ \scriptstyle \bullet$} [c] at 11 6
\put {$ \scriptstyle \bullet$} [c] at 13 6
\put {$ \scriptstyle \bullet$} [c] at 13 0
\setlinear \plot 16 12 13 0 13 12    /
\setlinear \plot 13  0 11 6 10 12     /
\setlinear \plot 12  12 11 6    /
\put{$2{,}520$} [c] at 13 -2
\endpicture
\end{minipage}
\begin{minipage}{4cm}
\beginpicture
\setcoordinatesystem units   <1.5mm,2mm>
\setplotarea x from 0 to 16, y from -2 to 15
\put{1.714)} [l] at 2 12
\put {$ \scriptstyle \bullet$} [c] at 10 0
\put {$ \scriptstyle \bullet$} [c] at 12 0
\put {$ \scriptstyle \bullet$} [c] at 13 0
\put {$ \scriptstyle \bullet$} [c] at 16 0
\put {$ \scriptstyle \bullet$} [c] at 11 6
\put {$ \scriptstyle \bullet$} [c] at 13 6
\put {$ \scriptstyle \bullet$} [c] at 13 12
\setlinear \plot 16 0 13 12 13 0    /
\setlinear \plot 13  12 11 6 10 0     /
\setlinear \plot 12  0 11 6    /
\put{$2{,}520$} [c] at 13 -2
\endpicture
\end{minipage}
\begin{minipage}{4cm}
\beginpicture
\setcoordinatesystem units   <1.5mm,2mm>
\setplotarea x  from 0 to 16, y from -2 to 15
\put{1.715)} [l] at 2 12
\put {$ \scriptstyle \bullet$} [c] at  10 12
\put {$ \scriptstyle \bullet$} [c] at  10 0
\put {$ \scriptstyle \bullet$} [c] at  10 6
\put {$ \scriptstyle \bullet$} [c] at  11.5 6
\put {$ \scriptstyle \bullet$} [c] at  13 0
\put {$ \scriptstyle \bullet$} [c] at  13 12
\put {$ \scriptstyle \bullet$} [c] at  16 12
\setlinear \plot 10 0 10 12 13 0 13 12  /
\setlinear \plot  16 12  13  0   /
\put{$2{,}520$} [c] at 13 -2
\endpicture
\end{minipage}
\begin{minipage}{4cm}
\beginpicture
\setcoordinatesystem units   <1.5mm,2mm>
\setplotarea x  from 0 to 16, y from -2 to 15
\put{1.716)} [l] at 2 12
\put {$ \scriptstyle \bullet$} [c] at  10 12
\put {$ \scriptstyle \bullet$} [c] at  10 0
\put {$ \scriptstyle \bullet$} [c] at  10 6
\put {$ \scriptstyle \bullet$} [c] at  11.5 6
\put {$ \scriptstyle \bullet$} [c] at  13 0
\put {$ \scriptstyle \bullet$} [c] at  13 12
\put {$ \scriptstyle \bullet$} [c] at  16 0
\setlinear \plot 10 12 10 0 13 12 13 0  /
\setlinear \plot  16 0  13  12   /
\put{$2{,}520$} [c] at 13 -2
\endpicture
\end{minipage}
$$
$$
\begin{minipage}{4cm}
\beginpicture
\setcoordinatesystem units   <1.5mm,2mm>
\setplotarea x  from 0 to 16, y from -2 to 15
\put{1.717)} [l] at 2 12
\put {$ \scriptstyle \bullet$} [c] at  10 12
\put {$ \scriptstyle \bullet$} [c] at  12.5 6
\put {$ \scriptstyle \bullet$} [c] at  13.5 6
\put {$ \scriptstyle \bullet$} [c] at  13 0
\put {$ \scriptstyle \bullet$} [c] at  13 12
\put {$ \scriptstyle \bullet$} [c] at  16 0
\put {$ \scriptstyle \bullet$} [c] at  16 12
\setlinear \plot  10 12 13 0 12.5 6 13 12 13.5 6 13 0   /
\setlinear \plot  13 12 16 0 16 12  /
\put{$2{,}520$} [c] at 13 -2
\endpicture
\end{minipage}
\begin{minipage}{4cm}
\beginpicture
\setcoordinatesystem units   <1.5mm,2mm>
\setplotarea x  from 0 to 16, y from -2 to 15
\put{1.718)} [l] at 2 12
\put {$ \scriptstyle \bullet$} [c] at  10 0
\put {$ \scriptstyle \bullet$} [c] at  12.5 6
\put {$ \scriptstyle \bullet$} [c] at  13.5 6
\put {$ \scriptstyle \bullet$} [c] at  13 0
\put {$ \scriptstyle \bullet$} [c] at  13 12
\put {$ \scriptstyle \bullet$} [c] at  16 0
\put {$ \scriptstyle \bullet$} [c] at  16 12
\setlinear \plot  10 0 13 12 12.5 6 13 0 13.5 6 13 12   /
\setlinear \plot  13 0 16 12 16 0  /
\put{$2{,}520$} [c] at 13 -2
\endpicture
\end{minipage}
\begin{minipage}{4cm}
\beginpicture
\setcoordinatesystem units   <1.5mm,2mm>
\setplotarea x from 0 to 16, y from -2 to 15
\put {1.719)} [l] at 2 12
\put {$ \scriptstyle \bullet$} [c] at  10 12
\put {$ \scriptstyle \bullet$} [c] at  12 12
\put {$ \scriptstyle \bullet$} [c] at  14 12
\put {$ \scriptstyle \bullet$} [c] at  16 12
\put {$ \scriptstyle \bullet$} [c] at  12 6
\put {$ \scriptstyle \bullet$} [c] at  14 0
\put {$ \scriptstyle \bullet$} [c] at  10 0
\setlinear \plot  10 12 10 0 12 6 12  12  /
\setlinear \plot  14 12 14 0  12 6   /
\setlinear \plot 16 12 14 0  /
\put{$2{,}520$} [c] at 13 -2
\endpicture
\end{minipage}
\begin{minipage}{4cm}
\beginpicture
\setcoordinatesystem units   <1.5mm,2mm>
\setplotarea x from 0 to 16, y from -2 to 15
\put {1.720)} [l] at 2 12
\put {$ \scriptstyle \bullet$} [c] at  10 0
\put {$ \scriptstyle \bullet$} [c] at  12 0
\put {$ \scriptstyle \bullet$} [c] at  14 0
\put {$ \scriptstyle \bullet$} [c] at  16 0
\put {$ \scriptstyle \bullet$} [c] at  12 6
\put {$ \scriptstyle \bullet$} [c] at  14 12
\put {$ \scriptstyle \bullet$} [c] at  10 12
\setlinear \plot  10 0 10 12 12 6 12  0  /
\setlinear \plot  14 0 14 12  12 6   /
\setlinear \plot 16 0 14 12  /
\put{$2{,}520$} [c] at 13 -2
\endpicture
\end{minipage}
\begin{minipage}{4cm}
\beginpicture
\setcoordinatesystem units   <1.5mm,2mm>
\setplotarea x from 0 to 16, y from -2 to 15
\put {1.721)} [l] at 2 12
\put {$ \scriptstyle \bullet$} [c] at  10 12
\put {$ \scriptstyle \bullet$} [c] at  12 12
\put {$ \scriptstyle \bullet$} [c] at  14 12
\put {$ \scriptstyle \bullet$} [c] at  16 12
\put {$ \scriptstyle \bullet$} [c] at  10 0
\put {$ \scriptstyle \bullet$} [c] at  15  6
\put {$ \scriptstyle \bullet$} [c] at  16 0
\setlinear \plot  10 12 10 0 16 12 16 0 14 12  10 0 12 12 /
\put{$2{,}520$} [c] at 13 -2
\endpicture
\end{minipage}
\begin{minipage}{4cm}
\beginpicture
\setcoordinatesystem units   <1.5mm,2mm>
\setplotarea x from 0 to 16, y from -2 to 15
\put {1.722)} [l] at 2 12
\put {$ \scriptstyle \bullet$} [c] at  10 0
\put {$ \scriptstyle \bullet$} [c] at  12 0
\put {$ \scriptstyle \bullet$} [c] at  14 0
\put {$ \scriptstyle \bullet$} [c] at  16 0
\put {$ \scriptstyle \bullet$} [c] at  10 12
\put {$ \scriptstyle \bullet$} [c] at  15  6
\put {$ \scriptstyle \bullet$} [c] at  16 12
\setlinear \plot  10 0 10 12 16 0 16 12 14 0  10 12 12 0 /
\put{$2{,}520$} [c] at 13 -2
\endpicture
\end{minipage}
$$

$$
\begin{minipage}{4cm}
\beginpicture
\setcoordinatesystem units   <1.5mm,2mm>
\setplotarea x  from 0 to 16, y from -2 to 15
\put{1.723)} [l] at 2 12
\put {$ \scriptstyle \bullet$} [c] at  10 0
\put {$ \scriptstyle \bullet$} [c] at  10 12
\put {$ \scriptstyle \bullet$} [c] at  16 12
\put {$ \scriptstyle \bullet$} [c] at  16  0
\put {$ \scriptstyle \bullet$} [c] at  13 6
\put {$ \scriptstyle \bullet$} [c] at  13 12
\put {$ \scriptstyle \bullet$} [c] at  13  0
\setlinear \plot 10 0 10 12 13 0 13 12 /
\setlinear \plot  13 0 16 12 16 0      /
\put{$2{,}520$} [c] at 13 -2
\endpicture
\end{minipage}
\begin{minipage}{4cm}
\beginpicture
\setcoordinatesystem units   <1.5mm,2mm>
\setplotarea x  from 0 to 16, y from -2 to 15
\put{1.724)} [l] at 2 12
\put {$ \scriptstyle \bullet$} [c] at  10 0
\put {$ \scriptstyle \bullet$} [c] at  10 12
\put {$ \scriptstyle \bullet$} [c] at  16 12
\put {$ \scriptstyle \bullet$} [c] at  16  0
\put {$ \scriptstyle \bullet$} [c] at  13 6
\put {$ \scriptstyle \bullet$} [c] at  13 12
\put {$ \scriptstyle \bullet$} [c] at  13  0
\setlinear \plot 10 12 10 0 13 12 13 0 /
\setlinear \plot  13 12 16 0 16 12      /
\put{$2{,}520$} [c] at 13 -2
\endpicture
\end{minipage}
\begin{minipage}{4cm}
\beginpicture
\setcoordinatesystem units   <1.5mm,2mm>
\setplotarea x from 0 to 16, y from -2 to 15
\put {1.725)} [l] at 2 12
\put {$ \scriptstyle \bullet$} [c] at  10 12
\put {$ \scriptstyle \bullet$} [c] at  12 12
\put {$ \scriptstyle \bullet$} [c] at  14 12
\put {$ \scriptstyle \bullet$} [c] at  16 12
\put {$ \scriptstyle \bullet$} [c] at  10 0
\put {$ \scriptstyle \bullet$} [c] at  13 0
\put {$ \scriptstyle \bullet$} [c] at  16 0
\setlinear \plot   10 12 10  0  12 12 13 0 14  12  16 0 16  12 13 0 /
\put{$2{,}520$} [c] at 13 -2
\endpicture
\end{minipage}
\begin{minipage}{4cm}
\beginpicture
\setcoordinatesystem units   <1.5mm,2mm>
\setplotarea x from 0 to 16, y from -2 to 15
\put {1.726)} [l] at 2 12
\put {$ \scriptstyle \bullet$} [c] at  10 0
\put {$ \scriptstyle \bullet$} [c] at  12 0
\put {$ \scriptstyle \bullet$} [c] at  14 0
\put {$ \scriptstyle \bullet$} [c] at  16 0
\put {$ \scriptstyle \bullet$} [c] at  10 12
\put {$ \scriptstyle \bullet$} [c] at  13 12
\put {$ \scriptstyle \bullet$} [c] at  16 12
\setlinear \plot   10 0 10  12  12 0 13 12 14  0  16 12 16  0 13 12 /
\put{$2{,}520$} [c] at 13 -2
\endpicture
\end{minipage}
\begin{minipage}{4cm}
\beginpicture
\setcoordinatesystem units   <1.5mm,2mm>
\setplotarea x from 0 to 16, y from -2 to 15
\put {1.727)} [l] at 2 12
\put {$ \scriptstyle \bullet$} [c] at  10 12
\put {$ \scriptstyle \bullet$} [c] at  12 12
\put {$ \scriptstyle \bullet$} [c] at  14 12
\put {$ \scriptstyle \bullet$} [c] at  16 12
\put {$ \scriptstyle \bullet$} [c] at  10 0
\put {$ \scriptstyle \bullet$} [c] at  13 0
\put {$ \scriptstyle \bullet$} [c] at  16 0
\setlinear \plot   10 12 10 0  12 12 13 0 14 12 16 0  16 12 10 0  /
\put{$2{,}520$} [c] at 13 -2
\endpicture
\end{minipage}
\begin{minipage}{4cm}
\beginpicture
\setcoordinatesystem units   <1.5mm,2mm>
\setplotarea x from 0 to 16, y from -2 to 15
\put {1.728)} [l] at 2 12
\put {$ \scriptstyle \bullet$} [c] at  10 0
\put {$ \scriptstyle \bullet$} [c] at  12 0
\put {$ \scriptstyle \bullet$} [c] at  14 0
\put {$ \scriptstyle \bullet$} [c] at  16 0
\put {$ \scriptstyle \bullet$} [c] at  10 12
\put {$ \scriptstyle \bullet$} [c] at  13 12
\put {$ \scriptstyle \bullet$} [c] at  16 12
\setlinear \plot   10 0 10 12  12 0 13 12 14 0 16 12  16 0 10 12  /
\put{$2{,}520$} [c] at 13 -2
\endpicture
\end{minipage}
$$
$$
\begin{minipage}{4cm}
\beginpicture
\setcoordinatesystem units   <1.5mm,2mm>
\setplotarea x from 0 to 16, y from -2 to 15
\put {1.729)} [l] at 2 12
\put {$ \scriptstyle \bullet$} [c] at  10 12
\put {$ \scriptstyle \bullet$} [c] at  12 12
\put {$ \scriptstyle \bullet$} [c] at  14 12
\put {$ \scriptstyle \bullet$} [c] at  16 12
\put {$ \scriptstyle \bullet$} [c] at  10  0
\put {$ \scriptstyle \bullet$} [c] at  13 0
\put {$ \scriptstyle \bullet$} [c] at  16 0
\setlinear \plot   10 12 10  0  12 12 13 0 16  12  10 0   /
\setlinear \plot   10 12 13 0 14 12 16 0  /
\put{$840$} [c] at 13 -2
\endpicture
\end{minipage}
\begin{minipage}{4cm}
\beginpicture
\setcoordinatesystem units   <1.5mm,2mm>
\setplotarea x from 0 to 16, y from -2 to 15
\put {1.730)} [l] at 2 12
\put {$ \scriptstyle \bullet$} [c] at  10 0
\put {$ \scriptstyle \bullet$} [c] at  12 0
\put {$ \scriptstyle \bullet$} [c] at  14 0
\put {$ \scriptstyle \bullet$} [c] at  16 0
\put {$ \scriptstyle \bullet$} [c] at  10  12
\put {$ \scriptstyle \bullet$} [c] at  13 12
\put {$ \scriptstyle \bullet$} [c] at  16 12
\setlinear \plot   10 0 10  12  12 0 13 12 16  0  10 12   /
\setlinear \plot   10 0 13 12 14 0 16 12  /
\put{$840$} [c] at 13 -2
\endpicture
\end{minipage}
\begin{minipage}{4cm}
\beginpicture
\setcoordinatesystem units    <1.5mm,2mm>
\setplotarea x from  0 to 16, y from -2 to 15
\put{${\bf  32}$} [l] at 2 15

\put{1.731)} [l] at 2 12
\put {$ \scriptstyle \bullet$} [c] at  10 0
\put {$ \scriptstyle \bullet$} [c] at  10 6
\put {$ \scriptstyle \bullet$} [c] at  10 12
\put {$ \scriptstyle \bullet$} [c] at  13 12
\put {$ \scriptstyle \bullet$} [c] at  13 6
\put {$ \scriptstyle \bullet$} [c] at  16 0
\put {$ \scriptstyle \bullet$} [c] at  16 12
\setlinear \plot  10 12 10 0 16 12 16 0    /
\setlinear \plot  13 12 10 0  /
\put{$5{,}040$} [c] at 13 -2
\endpicture
\end{minipage}
\begin{minipage}{4cm}
\beginpicture
\setcoordinatesystem units    <1.5mm,2mm>
\setplotarea x from  0 to 16, y from -2 to 15
\put{1.732)} [l] at 2 12
\put {$ \scriptstyle \bullet$} [c] at  10 0
\put {$ \scriptstyle \bullet$} [c] at  10 6
\put {$ \scriptstyle \bullet$} [c] at  10 12
\put {$ \scriptstyle \bullet$} [c] at  13 0
\put {$ \scriptstyle \bullet$} [c] at  13 6
\put {$ \scriptstyle \bullet$} [c] at  16 0
\put {$ \scriptstyle \bullet$} [c] at  16 12
\setlinear \plot  10 0 10 12 16 0 16 12    /
\setlinear \plot  13 0 10 12  /
\put{$5{,}040$} [c] at 13 -2
\endpicture
\end{minipage}
\begin{minipage}{4cm}
\beginpicture
\setcoordinatesystem units    <1.5mm,2mm>
\setplotarea x from 0 to 16, y from -2 to 15
\put {1.733)} [l]  at 2 12
\put {$ \scriptstyle \bullet$} [c] at  10 12
\put {$ \scriptstyle \bullet$} [c] at  12 12
\put {$ \scriptstyle \bullet$} [c] at  14 12
\put {$ \scriptstyle \bullet$} [c] at  16 12
\put {$ \scriptstyle \bullet$} [c] at  10 0
\put {$ \scriptstyle \bullet$} [c] at  13 0
\put {$ \scriptstyle \bullet$} [c] at  16 0
\setlinear \plot   10 12 10 0  12 12 13 0 14  12  16 0 16 12  /
\setlinear \plot   10 0  14  12  /
\put{$5{,}040$} [c] at 13 -2
\endpicture
\end{minipage}
\begin{minipage}{4cm}
\beginpicture
\setcoordinatesystem units    <1.5mm,2mm>
\setplotarea x from 0 to 16, y from -2 to 15
\put {1.734)} [l]  at 2 12
\put {$ \scriptstyle \bullet$} [c] at  10 0
\put {$ \scriptstyle \bullet$} [c] at  12 0
\put {$ \scriptstyle \bullet$} [c] at  14 0
\put {$ \scriptstyle \bullet$} [c] at  16 0
\put {$ \scriptstyle \bullet$} [c] at  10 12
\put {$ \scriptstyle \bullet$} [c] at  13 12
\put {$ \scriptstyle \bullet$} [c] at  16 12
\setlinear \plot   10 0 10 12  12 0 13 12 14  0  16 12 16 0  /
\setlinear \plot   10 12  14  0  /
\put{$5{,}040$} [c] at 13 -2
\endpicture
\end{minipage}
$$

$$
\begin{minipage}{4cm}
\beginpicture
\setcoordinatesystem units    <1.5mm,2mm>
\setplotarea x from  0 to 16, y from -2 to 15
\put{1.735)} [l] at 2 12
\put {$ \scriptstyle \bullet$} [c] at  10 0
\put {$ \scriptstyle \bullet$} [c] at  10 12
\put {$ \scriptstyle \bullet$} [c] at  13 0
\put {$ \scriptstyle \bullet$} [c] at  13 6
\put {$ \scriptstyle \bullet$} [c] at  13  12
\put {$ \scriptstyle \bullet$} [c] at  16  12
\put {$ \scriptstyle \bullet$} [c] at  16  0
\setlinear \plot 16 12 13 0 13 12  10 0  10 12 13 0  /
\setlinear \plot 16 0 13 12  /
\put{$5{,}040$} [c] at 12 -2
\endpicture
\end{minipage}
\begin{minipage}{4cm}
\beginpicture
\setcoordinatesystem units    <1.5mm,2mm>
\setplotarea x from 0 to 16, y from -2 to 15
\put {1.736)} [l]  at 2 12
\put {$ \scriptstyle \bullet$} [c] at  10 12
\put {$ \scriptstyle \bullet$} [c] at  12 12
\put {$ \scriptstyle \bullet$} [c] at  14 12
\put {$ \scriptstyle \bullet$} [c] at  16 12
\put {$ \scriptstyle \bullet$} [c] at  10 0
\put {$ \scriptstyle \bullet$} [c] at  13 0
\put {$ \scriptstyle \bullet$} [c] at  16 0
\setlinear \plot   10 12 10  0  12 12 13 0 16  12  10 0 14 12 13 0  /
\setlinear \plot   14 12  16  0  /
\put{$2{,}520$} [c] at 13 -2
\endpicture
\end{minipage}
\begin{minipage}{4cm}
\beginpicture
\setcoordinatesystem units    <1.5mm,2mm>
\setplotarea x from 0 to 16, y from -2 to 15
\put {1.737)} [l]  at 2 12
\put {$ \scriptstyle \bullet$} [c] at  10 0
\put {$ \scriptstyle \bullet$} [c] at  12 0
\put {$ \scriptstyle \bullet$} [c] at  14 0
\put {$ \scriptstyle \bullet$} [c] at  16 0
\put {$ \scriptstyle \bullet$} [c] at  10 12
\put {$ \scriptstyle \bullet$} [c] at  13 12
\put {$ \scriptstyle \bullet$} [c] at  16 12
\setlinear \plot   10 0 10  12  12 0 13 12 16  0  10 12 14 0 13 12  /
\setlinear \plot   14 0  16  12  /
\put{$2{,}520$} [c] at 13 -2
\endpicture
\end{minipage}
\begin{minipage}{4cm}
\beginpicture
\setcoordinatesystem units    <1.5mm,2mm>
\setplotarea x from  0 to 16, y from -2 to 15
\put{1.738)} [l] at 2 12
\put {$ \scriptstyle \bullet$} [c] at  10 12
\put {$ \scriptstyle \bullet$} [c] at  12 6
\put {$ \scriptstyle \bullet$} [c] at  12.5 0
\put {$ \scriptstyle \bullet$} [c] at  12.5 12
\put {$ \scriptstyle \bullet$} [c] at  13 6
\put {$ \scriptstyle \bullet$} [c] at  16 12
\put {$ \scriptstyle \bullet$} [c] at  16 0
\setlinear \plot  10 12 12.5 0 12 6 12.5 12 13 6 12.5 0 16 12 16 0  /
\put{$2{,}520$} [c] at 13 -2
\endpicture
\end{minipage}
\begin{minipage}{4cm}
\beginpicture
\setcoordinatesystem units    <1.5mm,2mm>
\setplotarea x from  0 to 16, y from -2 to 15
\put{1.739)} [l] at 2 12
\put {$ \scriptstyle \bullet$} [c] at  10 0
\put {$ \scriptstyle \bullet$} [c] at  12 6
\put {$ \scriptstyle \bullet$} [c] at  12.5 0
\put {$ \scriptstyle \bullet$} [c] at  12.5 12
\put {$ \scriptstyle \bullet$} [c] at  13 6
\put {$ \scriptstyle \bullet$} [c] at  16 12
\put {$ \scriptstyle \bullet$} [c] at  16 0
\setlinear \plot  10 0 12.5 12 12 6 12.5 0 13 6 12.5 12 16 0 16 12  /
\put{$2{,}520$} [c] at 13 -2
\endpicture
\end{minipage}
\begin{minipage}{4cm}
\beginpicture
\setcoordinatesystem units    <1.5mm,2mm>
\setplotarea x from  0 to 16, y from -2 to 15
\put{1.740)} [l] at 2 12
\put {$ \scriptstyle \bullet$} [c] at  10 0
\put {$ \scriptstyle \bullet$} [c] at  10 12
\put {$ \scriptstyle \bullet$} [c] at  13 0
\put {$ \scriptstyle \bullet$} [c] at  13 12
\put {$ \scriptstyle \bullet$} [c] at  13 6
\put {$ \scriptstyle \bullet$} [c] at  16 12
\put {$ \scriptstyle \bullet$} [c] at  16 0
\setlinear \plot  10 12  10  0  16 12 16 0   /
\setlinear \plot  13  12 10  0   /
\setlinear \plot  13 0 13 6   /
\put{$2{,}520$} [c] at 12 -2
\endpicture
\end{minipage}
$$
$$
\begin{minipage}{4cm}
\beginpicture
\setcoordinatesystem units    <1.5mm,2mm>
\setplotarea x from  0 to 16, y from -2 to 15
\put{1.741)} [l] at 2 12
\put {$ \scriptstyle \bullet$} [c] at  10 0
\put {$ \scriptstyle \bullet$} [c] at  10 12
\put {$ \scriptstyle \bullet$} [c] at  13 0
\put {$ \scriptstyle \bullet$} [c] at  13 12
\put {$ \scriptstyle \bullet$} [c] at  13 6
\put {$ \scriptstyle \bullet$} [c] at  16 12
\put {$ \scriptstyle \bullet$} [c] at  16 0
\setlinear \plot  10 0  10  12  16 0 16 12   /
\setlinear \plot  13  0 10  12   /
\setlinear \plot  13 12 13 6   /
\put{$2{,}520$} [c] at 12 -2
\endpicture
\end{minipage}
\begin{minipage}{4cm}
\beginpicture
\setcoordinatesystem units    <1.5mm,2mm>
\setplotarea x from 0 to 16, y from -2 to 15
\put {1.742)} [l] at 2 12
\put {$ \scriptstyle \bullet$} [c] at  10 12
\put {$ \scriptstyle \bullet$} [c] at  12 12
\put {$ \scriptstyle \bullet$} [c] at  14 12
\put {$ \scriptstyle \bullet$} [c] at  16 0
\put {$ \scriptstyle \bullet$} [c] at  12 0
\put {$ \scriptstyle \bullet$} [c] at  16  6
\put {$ \scriptstyle \bullet$} [c] at  16  12
\setlinear \plot  10 12 12  0  14 12   16 0 16  12 /
\setlinear \plot 12 0 12 12 /
\put{$2{,}520$} [c] at 13 -2
\endpicture
\end{minipage}
\begin{minipage}{4cm}
\beginpicture
\setcoordinatesystem units    <1.5mm,2mm>
\setplotarea x from 0 to 16, y from -2 to 15
\put {1.743)} [l] at 2 12
\put {$ \scriptstyle \bullet$} [c] at  10 0
\put {$ \scriptstyle \bullet$} [c] at  12 0
\put {$ \scriptstyle \bullet$} [c] at  14 0
\put {$ \scriptstyle \bullet$} [c] at  16 0
\put {$ \scriptstyle \bullet$} [c] at  12 12
\put {$ \scriptstyle \bullet$} [c] at  16  6
\put {$ \scriptstyle \bullet$} [c] at  16  12
\setlinear \plot  10 0 12  12  14 0   16 12 16  0 /
\setlinear \plot 12 0 12 12 /
\put{$2{,}520$} [c] at 13 -2
\endpicture
\end{minipage}
\begin{minipage}{4cm}
\beginpicture
\setcoordinatesystem units    <1.5mm,2mm>
\setplotarea x from 0 to 16, y from -2 to 15
\put {1.744)} [l] at 2 12
\put {$ \scriptstyle \bullet$} [c] at  10 12
\put {$ \scriptstyle \bullet$} [c] at  12 12
\put {$ \scriptstyle \bullet$} [c] at  14 12
\put {$ \scriptstyle \bullet$} [c] at  16 12
\put {$ \scriptstyle \bullet$} [c] at  14 6
\put {$ \scriptstyle \bullet$} [c] at  14  0
\put {$ \scriptstyle \bullet$} [c] at  10 0
\setlinear \plot  10 12 10 0  12 12 14 0 14 12 /
\setlinear \plot  10 12 14 0   /
\setlinear \plot 16 12 14 0 /
\put{$2{,}520$} [c] at 13 -2
\endpicture
\end{minipage}
\begin{minipage}{4cm}
\beginpicture
\setcoordinatesystem units    <1.5mm,2mm>
\setplotarea x from 0 to 16, y from -2 to 15
\put {1.745)} [l] at 2 12
\put {$ \scriptstyle \bullet$} [c] at  10 0
\put {$ \scriptstyle \bullet$} [c] at  12 0
\put {$ \scriptstyle \bullet$} [c] at  14 0
\put {$ \scriptstyle \bullet$} [c] at  16 0
\put {$ \scriptstyle \bullet$} [c] at  14 6
\put {$ \scriptstyle \bullet$} [c] at  14  12
\put {$ \scriptstyle \bullet$} [c] at  10 12
\setlinear \plot  10 0 10 12  12 0 14 12 14 0 /
\setlinear \plot  10 0 14 12   /
\setlinear \plot 16 0 14 12 /
\put{$2{,}520$} [c] at 13 -2
\endpicture
\end{minipage}
\begin{minipage}{4cm}
\beginpicture
\setcoordinatesystem units    <1.5mm,2mm>
\setplotarea x from 0 to 16, y from -2 to 15
\put {1.746)} [l] at 2 12
\put {$ \scriptstyle \bullet$} [c] at  10 12
\put {$ \scriptstyle \bullet$} [c] at  12 12
\put {$ \scriptstyle \bullet$} [c] at  14 12
\put {$ \scriptstyle \bullet$} [c] at  16 12
\put {$ \scriptstyle \bullet$} [c] at  14 6
\put {$ \scriptstyle \bullet$} [c] at  14  0
\put {$ \scriptstyle \bullet$} [c] at  10 0
\setlinear \plot  10 0 10 12  14 0 14 12 /
\setlinear \plot  12 12 14 0  /
\setlinear \plot 16 12 14 6 /
\put{$2{,}520$} [c] at 13 -2
\endpicture
\end{minipage}
$$

$$
\begin{minipage}{4cm}
\beginpicture
\setcoordinatesystem units    <1.5mm,2mm>
\setplotarea x from 0 to 16, y from -2 to 15
\put {1.747)} [l] at 2 12
\put {$ \scriptstyle \bullet$} [c] at  10 0
\put {$ \scriptstyle \bullet$} [c] at  12 0
\put {$ \scriptstyle \bullet$} [c] at  14 0
\put {$ \scriptstyle \bullet$} [c] at  16 0
\put {$ \scriptstyle \bullet$} [c] at  14 6
\put {$ \scriptstyle \bullet$} [c] at  14  12
\put {$ \scriptstyle \bullet$} [c] at  10 12
\setlinear \plot  10 12 10 0  14 12 14 0 /
\setlinear \plot  12 0 14 12  /
\setlinear \plot 16 0 14 6 /
\put{$2{,}520$} [c] at 13 -2
\endpicture
\end{minipage}
\begin{minipage}{4cm}
\beginpicture
\setcoordinatesystem units    <1.5mm,2mm>
\setplotarea x from 0 to 16, y from -2 to 15
\put {1.748)} [l]  at 2 12
\put {$ \scriptstyle \bullet$} [c] at  10 0
\put {$ \scriptstyle \bullet$} [c] at  10 12
\put {$ \scriptstyle \bullet$} [c] at  11.5 0
\put {$ \scriptstyle \bullet$} [c] at  13   0
\put {$ \scriptstyle \bullet$} [c] at  13 12
\put {$ \scriptstyle \bullet$} [c] at  14  12
\put {$ \scriptstyle \bullet$} [c] at  16  12
\setlinear \plot   10 12 10 0  13 12  13 0  16 12  /
\setlinear \plot   14 12 13 0 10 12 11.5 0 13 12 /
\put{$630$} [c] at 13 -2
\endpicture
\end{minipage}
\begin{minipage}{4cm}
\beginpicture
\setcoordinatesystem units    <1.5mm,2mm>
\setplotarea x from 0 to 16, y from -2 to 15
\put {1.749)} [l]  at 2 12
\put {$ \scriptstyle \bullet$} [c] at  10 0
\put {$ \scriptstyle \bullet$} [c] at  10 12
\put {$ \scriptstyle \bullet$} [c] at  11.5 12
\put {$ \scriptstyle \bullet$} [c] at  13   0
\put {$ \scriptstyle \bullet$} [c] at  13 12
\put {$ \scriptstyle \bullet$} [c] at  14  0
\put {$ \scriptstyle \bullet$} [c] at  16  0
\setlinear \plot   10 0 10 12  13  0  13 12  16 0  /
\setlinear \plot   14 0 13 12 10 0 11.5 12 13 0 /
\put{$630$} [c] at 13 -2
\endpicture
\end{minipage}
\begin{minipage}{4cm}
\beginpicture
\setcoordinatesystem units    <1.5mm,2mm>
\setplotarea x from  0 to 16, y from -2 to 15
\put{1.750)} [l] at 2 12
\put {$ \scriptstyle \bullet$} [c] at 10 0
\put {$ \scriptstyle \bullet$} [c] at 12.5 6
\put {$ \scriptstyle \bullet$} [c] at 13 12
\put {$ \scriptstyle \bullet$} [c] at 13 0
\put {$ \scriptstyle \bullet$} [c] at 13.5 6
\put {$ \scriptstyle \bullet$} [c] at 14 12
\setlinear \plot 10 0 13 12 12.5 6 13 0 13.5 6 14 12     /
\setlinear \plot 13.5 6  13 12 /
\put{$5{,}040$} [c] at 13 -2
\put{$\scriptstyle \bullet$} [c] at 16  0 \endpicture
\end{minipage}
\begin{minipage}{4cm}
\beginpicture
\setcoordinatesystem units    <1.5mm,2mm>
\setplotarea x from  0 to 16, y from -2 to 15
\put{1.751)} [l] at 2 12
\put {$ \scriptstyle \bullet$} [c] at 10 12
\put {$ \scriptstyle \bullet$} [c] at 12.5 6
\put {$ \scriptstyle \bullet$} [c] at 13 12
\put {$ \scriptstyle \bullet$} [c] at 13 0
\put {$ \scriptstyle \bullet$} [c] at 13.5 6
\put {$ \scriptstyle \bullet$} [c] at 14 0
\setlinear \plot 10 12 13 0 12.5 6 13 12 13.5 6 14 0     /
\setlinear \plot 13.5 6  13 0 /
\put{$5{,}040$} [c] at 13 -2
\put{$\scriptstyle \bullet$} [c] at 16  0 \endpicture
\end{minipage}
\begin{minipage}{4cm}
\beginpicture
\setcoordinatesystem units    <1.5mm,2mm>
\setplotarea x from  0 to 16, y from -2 to 15
\put{1.752)} [l] at 2 12
\put {$ \scriptstyle \bullet$} [c] at 10 0
\put {$ \scriptstyle \bullet$} [c] at 10 12
\put {$ \scriptstyle \bullet$} [c] at 12 0
\put {$ \scriptstyle \bullet$} [c] at 12 6
\put {$ \scriptstyle \bullet$} [c] at 12 12
\put {$ \scriptstyle \bullet$} [c] at 14 12
\setlinear \plot 12 6 10 0 10 12 12 0 12 12   /
\setlinear \plot 12 0 14 12   /
\put{$5{,}040$} [c] at 13 -2
\put{$\scriptstyle \bullet$} [c] at 16  0 \endpicture
\end{minipage}
$$
$$
\begin{minipage}{4cm}
\beginpicture
\setcoordinatesystem units    <1.5mm,2mm>
\setplotarea x from  0 to 16, y from -2 to 15
\put{1.753)} [l] at 2 12
\put {$ \scriptstyle \bullet$} [c] at 10 0
\put {$ \scriptstyle \bullet$} [c] at 10 12
\put {$ \scriptstyle \bullet$} [c] at 12 0
\put {$ \scriptstyle \bullet$} [c] at 12 6
\put {$ \scriptstyle \bullet$} [c] at 12 12
\put {$ \scriptstyle \bullet$} [c] at 14 0
\setlinear \plot 12 6 10 12 10 0 12 12 12 0   /
\setlinear \plot 12 12 14 0   /
\put{$5{,}040$} [c] at 13 -2
\put{$\scriptstyle \bullet$} [c] at 16  0 \endpicture
\end{minipage}
\begin{minipage}{4cm}
\beginpicture
\setcoordinatesystem units    <1.5mm,2mm>
\setplotarea x from  0 to 16, y from -2 to 15
\put{1.754)} [l] at 2 12
\put {$ \scriptstyle \bullet$} [c] at 10 12
\put {$ \scriptstyle \bullet$} [c] at 10 0
\put {$ \scriptstyle \bullet$} [c] at 12 12
\put {$ \scriptstyle \bullet$} [c] at 14 0
\put {$ \scriptstyle \bullet$} [c] at 14  6
\put {$ \scriptstyle \bullet$} [c] at 14 12
\setlinear \plot 10 12 10 0 12 12 14 6 14 12   /
\setlinear \plot 14 0 14 6  /
\put{$5{,}040$} [c] at 13 -2
\put{$\scriptstyle \bullet$} [c] at 16  0 \endpicture
\end{minipage}
\begin{minipage}{4cm}
\beginpicture
\setcoordinatesystem units    <1.5mm,2mm>
\setplotarea x from  0 to 16, y from -2 to 15
\put{1.755)} [l] at 2 12
\put {$ \scriptstyle \bullet$} [c] at 10 12
\put {$ \scriptstyle \bullet$} [c] at 10 0
\put {$ \scriptstyle \bullet$} [c] at 12 0
\put {$ \scriptstyle \bullet$} [c] at 14 0
\put {$ \scriptstyle \bullet$} [c] at 14  6
\put {$ \scriptstyle \bullet$} [c] at 14 12
\setlinear \plot 10 0 10 12 12 0 14 6 14 0   /
\setlinear \plot 14 12 14 6  /
\put{$5{,}040$} [c] at 13 -2
\put{$\scriptstyle \bullet$} [c] at 16  0 \endpicture
\end{minipage}
\begin{minipage}{4cm}
\beginpicture
\setcoordinatesystem units    <1.5mm,2mm>
\setplotarea x from 0 to 16, y from -2 to 15
\put{1.756)} [l] at 2 12
\put {$ \scriptstyle \bullet$} [c] at 10 0
\put {$ \scriptstyle \bullet$} [c] at 10 12
\put {$ \scriptstyle \bullet$} [c] at 11 9
\put {$ \scriptstyle \bullet$} [c] at 13 3
\put {$ \scriptstyle \bullet$} [c] at 14  12
\put {$ \scriptstyle \bullet$} [c] at 14 0
\setlinear \plot 10 0 10 12  14 0 14 12   /
\put{$5{,}040$} [c] at 13 -2
\put{$\scriptstyle \bullet$} [c] at 16  0 \endpicture
\end{minipage}
\begin{minipage}{4cm}
\beginpicture
\setcoordinatesystem units    <1.5mm,2mm>
\setplotarea x from 0 to 16, y from -2 to 15
\put{1.757)} [l] at 2 12
\put {$ \scriptstyle \bullet$} [c] at 10 0
\put {$ \scriptstyle \bullet$} [c] at 10 12
\put {$ \scriptstyle \bullet$} [c] at 14 0
\put {$ \scriptstyle \bullet$} [c] at 14 12
\put {$ \scriptstyle \bullet$} [c] at 16 0
\put {$ \scriptstyle \bullet$} [c] at 16 6
\put{$\scriptstyle \bullet$} [c] at 16  12
\setlinear \plot  10 12 10 0 14 12 14 0 /
\setlinear \plot  16 0 16 12   /
\put{$5{,}040$} [c] at 13 -2
 \endpicture
\end{minipage}
\begin{minipage}{4cm}
\beginpicture
\setcoordinatesystem units    <1.5mm,2mm>
\setplotarea x from  0 to 16, y from -2 to 15
\put{1.758)} [l] at 2 12
\put {$ \scriptstyle \bullet$} [c] at 10 6
\put {$ \scriptstyle \bullet$} [c] at 11 0
\put {$ \scriptstyle \bullet$} [c] at 11 12
\put {$ \scriptstyle \bullet$} [c] at 12 6
\put {$ \scriptstyle \bullet$} [c] at 14 6
\put {$ \scriptstyle \bullet$} [c] at 14 12
\setlinear \plot 11  0 12 6 11 12 10  6 11 0 14 6 14 12    /
\put{$2{,}520$} [c] at 13 -2
\put{$\scriptstyle \bullet$} [c] at 16  0 \endpicture
\end{minipage}
$$

$$
\begin{minipage}{4cm}
\beginpicture
\setcoordinatesystem units    <1.5mm,2mm>
\setplotarea x from  0 to 16, y from -2 to 15
\put{1.759)} [l] at 2 12
\put {$ \scriptstyle \bullet$} [c] at 10 6
\put {$ \scriptstyle \bullet$} [c] at 11 0
\put {$ \scriptstyle \bullet$} [c] at 11 12
\put {$ \scriptstyle \bullet$} [c] at 12 6
\put {$ \scriptstyle \bullet$} [c] at 14 6
\put {$ \scriptstyle \bullet$} [c] at 14 0
\setlinear \plot 11  12 12 6 11 0 10  6 11 12 14 6 14 0    /
\put{$2{,}520$} [c] at 13 -2
\put{$\scriptstyle \bullet$} [c] at 16  0 \endpicture
\end{minipage}
\begin{minipage}{4cm}
\beginpicture
\setcoordinatesystem units    <1.5mm,2mm>
\setplotarea x from  0 to 16, y from -2 to 15
\put{1.760)} [l] at 2 12
\put {$ \scriptstyle \bullet$} [c] at 10 6
\put {$ \scriptstyle \bullet$} [c] at 10 12
\put {$ \scriptstyle \bullet$} [c] at 12 0
\put {$ \scriptstyle \bullet$} [c] at 12 12
\put {$ \scriptstyle \bullet$} [c] at 14 6
\put {$ \scriptstyle \bullet$} [c] at 14 12
\setlinear \plot 10 12 10 6 12 0 14 6 14 12     /
\setlinear \plot 12 12 14 6      /
\put{$2{,}520$} [c] at 13 -2
\put{$\scriptstyle \bullet$} [c] at 16  0 \endpicture
\end{minipage}
\begin{minipage}{4cm}
\beginpicture
\setcoordinatesystem units    <1.5mm,2mm>
\setplotarea x from  0 to 16, y from -2 to 15
\put{1.761)} [l] at 2 12
\put {$ \scriptstyle \bullet$} [c] at 10 6
\put {$ \scriptstyle \bullet$} [c] at 10 0
\put {$ \scriptstyle \bullet$} [c] at 12 0
\put {$ \scriptstyle \bullet$} [c] at 12 12
\put {$ \scriptstyle \bullet$} [c] at 14 6
\put {$ \scriptstyle \bullet$} [c] at 14 0
\setlinear \plot 10 0 10 6 12 12 14 6 14 0     /
\setlinear \plot 12 0 14 6      /
\put{$2{,}520$} [c] at 13 -2
\put{$\scriptstyle \bullet$} [c] at 16  0 \endpicture
\end{minipage}
\begin{minipage}{4cm}
\beginpicture
\setcoordinatesystem units    <1.5mm,2mm>
\setplotarea x from  0 to 16, y from -2 to 15
\put{1.762)} [l] at 2 12
\put {$ \scriptstyle \bullet$} [c] at 10 12
\put {$ \scriptstyle \bullet$} [c] at 10 0
\put {$ \scriptstyle \bullet$} [c] at 12 12
\put {$ \scriptstyle \bullet$} [c] at 14 0
\put {$ \scriptstyle \bullet$} [c] at 14 6
\put {$ \scriptstyle \bullet$} [c] at 14 12
\setlinear \plot 14 0 10 12 10 0 14 12 14 0 12 12 10 0   /
\put{$2{,}520$} [c] at 13 -2
\put{$\scriptstyle \bullet$} [c] at 16  0 \endpicture
\end{minipage}
\begin{minipage}{4cm}
\beginpicture
\setcoordinatesystem units    <1.5mm,2mm>
\setplotarea x from  0 to 16, y from -2 to 15
\put{1.763)} [l] at 2 12
\put {$ \scriptstyle \bullet$} [c] at 10 12
\put {$ \scriptstyle \bullet$} [c] at 10 0
\put {$ \scriptstyle \bullet$} [c] at 12 0
\put {$ \scriptstyle \bullet$} [c] at 14 0
\put {$ \scriptstyle \bullet$} [c] at 14 6
\put {$ \scriptstyle \bullet$} [c] at 14 12
\setlinear \plot 14 12 10 0 10 12 14 0 14 12 12 0 10 12   /
\put{$2{,}520$} [c] at 13 -2
\put{$\scriptstyle \bullet$} [c] at 16  0 \endpicture
\end{minipage}
\begin{minipage}{4cm}
\beginpicture
\setcoordinatesystem units    <1.5mm,2mm>
\setplotarea x from  0 to 16, y from -2 to 15
\put{1.764)} [l] at 2 12
\put {$ \scriptstyle \bullet$} [c] at 10 12
\put {$ \scriptstyle \bullet$} [c] at 10 0
\put {$ \scriptstyle \bullet$} [c] at 12 12
\put {$ \scriptstyle \bullet$} [c] at 14 0
\put {$ \scriptstyle \bullet$} [c] at 14  6
\put {$ \scriptstyle \bullet$} [c] at 14 12
\setlinear \plot 10 0 10 12 14 6 12 12   /
\setlinear \plot 14 0 14 12  /
\put{$2{,}520$} [c] at 13 -2
\put{$\scriptstyle \bullet$} [c] at 16  0 \endpicture
\end{minipage}
$$
$$
\begin{minipage}{4cm}
\beginpicture
\setcoordinatesystem units    <1.5mm,2mm>
\setplotarea x from 0 to 16, y from -2 to 15
\put{1.765)} [l] at 2 12
\put {$ \scriptstyle \bullet$} [c] at 10 12
\put {$ \scriptstyle \bullet$} [c] at 10 0
\put {$ \scriptstyle \bullet$} [c] at 12 0
\put {$ \scriptstyle \bullet$} [c] at 14 0
\put {$ \scriptstyle \bullet$} [c] at 14  6
\put {$ \scriptstyle \bullet$} [c] at 14 12
\setlinear \plot 10 12 10 0 14 6 12 0   /
\setlinear \plot 14 0 14 12  /
\put{$2{,}520$} [c] at 13 -2
\put{$\scriptstyle \bullet$} [c] at 16  0 \endpicture
\end{minipage}
\begin{minipage}{4cm}
\beginpicture
\setcoordinatesystem units    <1.5mm,2mm>
\setplotarea x from 0 to 16, y from -2 to 15
\put{1.766)} [l] at 2 12
\put {$ \scriptstyle \bullet$} [c] at 10 12
\put {$ \scriptstyle \bullet$} [c] at 10 0
\put {$ \scriptstyle \bullet$} [c] at 12 12
\put {$ \scriptstyle \bullet$} [c] at 14 0
\put {$ \scriptstyle \bullet$} [c] at 12  6
\put {$ \scriptstyle \bullet$} [c] at 14 12
\setlinear \plot 14 12 14  0 10 12 10 0 12 12 12 6    /
\put{$2{,}520$} [c] at 13 -2
\put{$\scriptstyle \bullet$} [c] at 16  0 \endpicture
\end{minipage}
\begin{minipage}{4cm}
\beginpicture
\setcoordinatesystem units    <1.5mm,2mm>
\setplotarea x from 0 to 16, y from -2 to 15
\put{1.767)} [l] at 2 12
\put {$ \scriptstyle \bullet$} [c] at 10 12
\put {$ \scriptstyle \bullet$} [c] at 10 0
\put {$ \scriptstyle \bullet$} [c] at 12 0
\put {$ \scriptstyle \bullet$} [c] at 14 0
\put {$ \scriptstyle \bullet$} [c] at 12  6
\put {$ \scriptstyle \bullet$} [c] at 14 12
\setlinear \plot 14 0 14  12 10 0 10 12 12 0 12 6    /
\put{$2{,}520$} [c] at 13 -2
\put{$\scriptstyle \bullet$} [c] at 16  0 \endpicture
\end{minipage}
\begin{minipage}{4cm}
\beginpicture
\setcoordinatesystem units    <1.5mm,2mm>
\setplotarea x from 0 to 16, y from -2 to 15
\put{1.768)} [l] at 2 12
\put {$ \scriptstyle \bullet$} [c] at 10 12
\put {$ \scriptstyle \bullet$} [c] at 12 12
\put {$ \scriptstyle \bullet$} [c] at  11 6
\put {$ \scriptstyle \bullet$} [c] at  11 0
\put {$ \scriptstyle \bullet$} [c] at 10.5 9
\put {$ \scriptstyle \bullet$} [c] at 14 0
\setlinear \plot  11 0 11 6  10 12  /
\setlinear \plot  11 6 12 12  /
\put{$2{,}520$} [c] at 13 -2
\put{$\scriptstyle \bullet$} [c] at 16  0 \endpicture
\end{minipage}
\begin{minipage}{4cm}
\beginpicture
\setcoordinatesystem units    <1.5mm,2mm>
\setplotarea x from 0 to 16, y from -2 to 15
\put{1.769)} [l] at 2 12
\put {$ \scriptstyle \bullet$} [c] at 10 0
\put {$ \scriptstyle \bullet$} [c] at 12 0
\put {$ \scriptstyle \bullet$} [c] at  11 6
\put {$ \scriptstyle \bullet$} [c] at  11 12
\put {$ \scriptstyle \bullet$} [c] at 10.5 3
\put {$ \scriptstyle \bullet$} [c] at 14 0
\setlinear \plot  11 12 11 6  10 0   /
\setlinear \plot  11 6 12 0  /
\put{$2{,}520$} [c] at 13 -2
\put{$\scriptstyle \bullet$} [c] at 16  0 \endpicture
\end{minipage}
\begin{minipage}{4cm}
\beginpicture
\setcoordinatesystem units    <1.5mm,2mm>
\setplotarea x from 0 to 16, y from -2 to 15
\put{1.770)} [l] at 2 12
\put {$ \scriptstyle \bullet$} [c] at 10 8
\put {$ \scriptstyle \bullet$} [c] at 10 4
\put {$ \scriptstyle \bullet$} [c] at 11 12
\put {$ \scriptstyle \bullet$} [c] at 11 0
\put {$ \scriptstyle \bullet$} [c] at 12  6
\put {$ \scriptstyle \bullet$} [c] at 14 0
\setlinear \plot 11 0 10 4  10 8  11 12 12 6  11 0   /
\put{$2{,}520$} [c] at 13 -2
\put{$\scriptstyle \bullet$} [c] at 16  0 \endpicture
\end{minipage}
$$
$$
\begin{minipage}{4cm}
\beginpicture
\setcoordinatesystem units    <1.5mm,2mm>
\setplotarea x from 0 to 16, y from -2 to 15
\put{1.771)} [l] at 2 12
\put {$ \scriptstyle \bullet$} [c] at 10 0
\put {$ \scriptstyle \bullet$} [c] at 10 6
\put {$ \scriptstyle \bullet$} [c] at 10 12
\put {$ \scriptstyle \bullet$} [c] at 13 6
\put {$ \scriptstyle \bullet$} [c] at 13 0
\put {$ \scriptstyle \bullet$} [c] at 13 12
\put{$\scriptstyle \bullet$} [c] at 16  0
\setlinear \plot  10 0 10 12  /
\setlinear \plot  13 0 13 12  /
\put{$2{,}520$} [c] at 13 -2
\endpicture
\end{minipage}
\begin{minipage}{4cm}
\beginpicture
\setcoordinatesystem units    <1.5mm,2mm>
\setplotarea x from 0 to 16, y from -2 to 15
\put{1.772)} [l] at 2 12
\put {$ \scriptstyle \bullet$} [c] at 10 12
\put {$ \scriptstyle \bullet$} [c] at 10 0
\put {$ \scriptstyle \bullet$} [c] at 13 6
\put {$ \scriptstyle \bullet$} [c] at 13.5 12
\put {$ \scriptstyle \bullet$} [c] at 13.5 0
\put {$ \scriptstyle \bullet$} [c] at 14 6
\setlinear \plot 10 0 10 12  13.5 0 13 6 13.5 12 14 6 13.5 0 10 12 /
\setlinear \plot 10 0 13.5 12 /
\put{$2{,}520$} [c] at 13 -2
\put{$\scriptstyle \bullet$} [c] at 16  0 \endpicture
\end{minipage}
\begin{minipage}{4cm}
\beginpicture
\setcoordinatesystem units    <1.5mm,2mm>
\setplotarea x from 0 to 16, y from -2 to 15
\put{1.773)} [l] at 2 12
\put {$ \scriptstyle \bullet$} [c] at 10 12
\put {$ \scriptstyle \bullet$} [c] at 10 0
\put {$ \scriptstyle \bullet$} [c] at 12 12
\put {$ \scriptstyle \bullet$} [c] at 12 0
\put {$ \scriptstyle \bullet$} [c] at 14  0
\put {$ \scriptstyle \bullet$} [c] at 14 12
\setlinear \plot  10 0 10 12 14 0 12 12 10 0 14 12 12 0 12  12 /
\setlinear \plot  10 12 12 0 /
\put{$1{,}260$} [c] at 13 -2
\put{$\scriptstyle \bullet$} [c] at 16  0 \endpicture
\end{minipage}
\begin{minipage}{4cm}
\beginpicture
\setcoordinatesystem units    <1.5mm,2mm>
\setplotarea x from 0 to 16, y from -2 to 15
\put{1.774)} [l] at 2 12
\put {$ \scriptstyle \bullet$} [c] at 10 12
\put {$ \scriptstyle \bullet$} [c] at 10 6
\put {$ \scriptstyle \bullet$} [c] at 11 0
\put {$ \scriptstyle \bullet$} [c] at 12 6
\put {$ \scriptstyle \bullet$} [c] at 12 12
\put {$ \scriptstyle \bullet$} [c] at 14 0
\setlinear \plot  11 0 10 6 10 12 12 6 12 12 10 6 11 0 12 6 /
\put{$630$} [c] at 13 -2
\put{$\scriptstyle \bullet$} [c] at 16  0 \endpicture
\end{minipage}
\begin{minipage}{4cm}
\beginpicture
\setcoordinatesystem units    <1.5mm,2mm>
\setplotarea x from 0 to 16, y from -2 to 15
\put{1.775)} [l] at 2 12
\put {$ \scriptstyle \bullet$} [c] at 10 0
\put {$ \scriptstyle \bullet$} [c] at 10 6
\put {$ \scriptstyle \bullet$} [c] at 11 12
\put {$ \scriptstyle \bullet$} [c] at 12 6
\put {$ \scriptstyle \bullet$} [c] at 12 0
\put {$ \scriptstyle \bullet$} [c] at 14 0
\setlinear \plot  11 12 10 6 10 0 12 6 12 0 10 6 11 12 12 6 /
\put{$630$} [c] at 13 -2
\put{$\scriptstyle \bullet$} [c] at 16  0 \endpicture
\end{minipage}
\begin{minipage}{4cm}
\beginpicture
\setcoordinatesystem units    <1.5mm,2mm>
\setplotarea x from 0 to 16, y from -2 to 15
\put{1.776)} [l] at 2 12
\put {$ \scriptstyle \bullet$} [c] at 10 0
\put {$ \scriptstyle \bullet$} [c] at 10 12
\put {$ \scriptstyle \bullet$} [c] at 11 6
\put {$ \scriptstyle \bullet$} [c] at 12 12
\put {$ \scriptstyle \bullet$} [c] at 12 0
\put {$ \scriptstyle \bullet$} [c] at 14 0
\put{$\scriptstyle \bullet$} [c] at 16  0
\setlinear \plot  10 0 12 12  /
\setlinear \plot  10 12 12 0 /
\put{$630$} [c] at 13 -2
 \endpicture
\end{minipage}
$$
$$
\begin{minipage}{4cm}
\beginpicture
\setcoordinatesystem units    <1.5mm,2mm>
\setplotarea x from 0 to 16, y from -2 to 15
\put {${\bf  33}$} [l] at 2 15

\put {1.777)} [l] at 2 12
\put {$ \scriptstyle \bullet$} [c] at  10 12
\put {$ \scriptstyle \bullet$} [c] at  12 12
\put {$ \scriptstyle \bullet$} [c] at  14 12
\put {$ \scriptstyle \bullet$} [c] at  16 12
\put {$ \scriptstyle \bullet$} [c] at  14 6
\put {$ \scriptstyle \bullet$} [c] at  14  0
\put {$ \scriptstyle \bullet$} [c] at  10  0
\setlinear \plot  10 12 10 0 12 12 14 0   14 12 /
\setlinear  \plot 16 12 14 0 /
\put{$5{,}040$} [c] at 13 -2
\endpicture
\end{minipage}
\begin{minipage}{4cm}
\beginpicture
\setcoordinatesystem units    <1.5mm,2mm>
\setplotarea x from 0 to 16, y from -2 to 15
\put {1.778)} [l] at 2 12
\put {$ \scriptstyle \bullet$} [c] at  10 0
\put {$ \scriptstyle \bullet$} [c] at  12 0
\put {$ \scriptstyle \bullet$} [c] at  14 0
\put {$ \scriptstyle \bullet$} [c] at  16 0
\put {$ \scriptstyle \bullet$} [c] at  14 6
\put {$ \scriptstyle \bullet$} [c] at  14  12
\put {$ \scriptstyle \bullet$} [c] at  10  12
\setlinear \plot  10 0 10 12 12 0 14 12   14 0 /
\setlinear  \plot 16 0 14 12 /
\put{$5{,}040$} [c] at 13 -2
\endpicture
\end{minipage}
\begin{minipage}{4cm}
\beginpicture
\setcoordinatesystem units    <1.5mm,2mm>
\setplotarea x from 0 to 16, y from -2 to 15
\put {1.779)} [l] at 2 12
\put {$ \scriptstyle \bullet$} [c] at  10 12
\put {$ \scriptstyle \bullet$} [c] at  12 12
\put {$ \scriptstyle \bullet$} [c] at  14 12
\put {$ \scriptstyle \bullet$} [c] at  16 12
\put {$ \scriptstyle \bullet$} [c] at  12 0
\put {$ \scriptstyle \bullet$} [c] at  14 6
\put {$ \scriptstyle \bullet$} [c] at  14  0
\setlinear \plot  10 12 12 0  14 12 14 0 16 12 /
\setlinear \plot  12 0 12 12 /
\put{$2{,}520$} [c] at 13 -2
\endpicture
\end{minipage}
\begin{minipage}{4cm}
\beginpicture
\setcoordinatesystem units    <1.5mm,2mm>
\setplotarea x from 0 to 16, y from -2 to 15
\put {1.780)} [l] at 2 12
\put {$ \scriptstyle \bullet$} [c] at  10 0
\put {$ \scriptstyle \bullet$} [c] at  12 0
\put {$ \scriptstyle \bullet$} [c] at  14 0
\put {$ \scriptstyle \bullet$} [c] at  16 0
\put {$ \scriptstyle \bullet$} [c] at  12 12
\put {$ \scriptstyle \bullet$} [c] at  14 6
\put {$ \scriptstyle \bullet$} [c] at  14  12
\setlinear \plot  10 0 12 12  14 0 14 12 16 0 /
\setlinear \plot  12 0 12 12 /
\put{$2{,}520$} [c] at 13 -2
\endpicture
\end{minipage}
\begin{minipage}{4cm}
\beginpicture
\setcoordinatesystem units    <1.5mm,2mm>
\setplotarea x from 0 to 16, y from -2 to 15
\put {1.781)}  [l] at 2 12
\put {$ \scriptstyle \bullet$} [c] at  10 12
\put {$ \scriptstyle \bullet$} [c] at  12 12
\put {$ \scriptstyle \bullet$} [c] at  14 12
\put {$ \scriptstyle \bullet$} [c] at  16  12
\put {$ \scriptstyle \bullet$} [c] at  10 0
\put {$ \scriptstyle \bullet$} [c] at  13  0
\put {$ \scriptstyle \bullet$} [c] at  16  0
\setlinear \plot   10 12 10 0  12 12 13 0 16  12    /
\setlinear \plot   10 12  13 0    /
\setlinear \plot   10 0  14 12 16 0    /
\put{$2{,}520$} [c] at 13 -2
\endpicture
\end{minipage}
\begin{minipage}{4cm}
\beginpicture
\setcoordinatesystem units    <1.5mm,2mm>
\setplotarea x from 0 to 16, y from -2 to 15
\put {1.782)}  [l] at 2 12
\put {$ \scriptstyle \bullet$} [c] at  10 0
\put {$ \scriptstyle \bullet$} [c] at  12 0
\put {$ \scriptstyle \bullet$} [c] at  14 0
\put {$ \scriptstyle \bullet$} [c] at  16  0
\put {$ \scriptstyle \bullet$} [c] at  10 12
\put {$ \scriptstyle \bullet$} [c] at  13  12
\put {$ \scriptstyle \bullet$} [c] at  16  12
\setlinear \plot   10 0 10 12  12 0 13 12 16  0    /
\setlinear \plot   10 0  13 12    /
\setlinear \plot   10 12  14 0 16 12    /
\put{$2{,}520$} [c] at 13 -2
\endpicture
\end{minipage}
$$
$$
\begin{minipage}{4cm}
\beginpicture
\setcoordinatesystem units    <1.5mm,2mm>
\setplotarea x from 0 to 16, y from -2 to 15
\put{1.783)} [l] at 2 12
\put {$ \scriptstyle \bullet$} [c] at 10 6
\put {$ \scriptstyle \bullet$} [c] at 10 12
\put {$ \scriptstyle \bullet$} [c] at 11 12
\put {$ \scriptstyle \bullet$} [c] at 11 0
\put {$ \scriptstyle \bullet$} [c] at 12 6
\put {$ \scriptstyle \bullet$} [c] at 15 12
\put {$ \scriptstyle \bullet$} [c] at 16 12
\setlinear \plot 10 12 10 6  11 0 12 6  11 12 10 6  /
\setlinear \plot 15 12 11  0 16 12   /
\put{$2{,}520$} [c] at 13 -2
\endpicture
\end{minipage}
\begin{minipage}{4cm}
\beginpicture
\setcoordinatesystem units    <1.5mm,2mm>
\setplotarea x from 0 to 16, y from -2 to 15
\put{1.784)} [l] at 2 12
\put {$ \scriptstyle \bullet$} [c] at 10 6
\put {$ \scriptstyle \bullet$} [c] at 10 0
\put {$ \scriptstyle \bullet$} [c] at 11 12
\put {$ \scriptstyle \bullet$} [c] at 11 0
\put {$ \scriptstyle \bullet$} [c] at 12 6
\put {$ \scriptstyle \bullet$} [c] at 15 0
\put {$ \scriptstyle \bullet$} [c] at 16 0
\setlinear \plot 10 0 10 6  11 12 12 6  11 0 10 6  /
\setlinear \plot 15 0 11  12 16 0   /
\put{$2{,}520$} [c] at 13 -2
\endpicture
\end{minipage}
\begin{minipage}{4cm}
\beginpicture
\setcoordinatesystem units    <1.5mm,2mm>
\setplotarea x from 0 to 16, y from -2 to 15
\put {1.785)}  [l] at 2 12
\put {$ \scriptstyle \bullet$} [c] at  10 12
\put {$ \scriptstyle \bullet$} [c] at  12 12
\put {$ \scriptstyle \bullet$} [c] at  14 12
\put {$ \scriptstyle \bullet$} [c] at  16 12
\put {$ \scriptstyle \bullet$} [c] at  10 0
\put {$ \scriptstyle \bullet$} [c] at  12 0
\put {$ \scriptstyle \bullet$} [c] at  16 0
\setlinear \plot   10 12 10 0  12 12 12 0 10 12    /
\setlinear \plot   12 12 16 0 16 12    /
\setlinear \plot   14 12 16 0    /
\put{$1{,}260$} [c] at 13 -2
\endpicture
\end{minipage}
\begin{minipage}{4cm}
\beginpicture
\setcoordinatesystem units    <1.5mm,2mm>
\setplotarea x from 0 to 16, y from -2 to 15
\put {1.786)}  [l] at 2 12
\put {$ \scriptstyle \bullet$} [c] at  10 0
\put {$ \scriptstyle \bullet$} [c] at  12 0
\put {$ \scriptstyle \bullet$} [c] at  14 0
\put {$ \scriptstyle \bullet$} [c] at  16 0
\put {$ \scriptstyle \bullet$} [c] at  10 12
\put {$ \scriptstyle \bullet$} [c] at  12 12
\put {$ \scriptstyle \bullet$} [c] at  16 12
\setlinear \plot   10 0 10 12  12  0 12 12 10 0    /
\setlinear \plot   12 0 16 12 16 0    /
\setlinear \plot   14 0 16 12    /
\put{$1{,}260$} [c] at 13 -2
\endpicture
\end{minipage}
\begin{minipage}{4cm}
\beginpicture
\setcoordinatesystem units    <1.5mm,2mm>
\setplotarea x from 0 to 16, y from -2 to 15
\put{1.787)} [l] at 2 12
\put {$ \scriptstyle \bullet$} [c] at 10 6
\put {$ \scriptstyle \bullet$} [c] at 10 9
\put {$ \scriptstyle \bullet$} [c] at 10 12
\put {$ \scriptstyle \bullet$} [c] at 13 0
\put {$ \scriptstyle \bullet$} [c] at 12 12
\put {$ \scriptstyle \bullet$} [c] at 14 12
\put {$ \scriptstyle \bullet$} [c] at 16 12
\setlinear \plot 10 12 10 6 13 0 12 12      /
\setlinear \plot 14 12 13 0 16 12      /
\put{$840$} [c] at 13 -2
\endpicture
\end{minipage}
\begin{minipage}{4cm}
\beginpicture
\setcoordinatesystem units    <1.5mm,2mm>
\setplotarea x from 0 to 16, y from -2 to 15
\put{1.788)} [l] at 2 12
\put {$ \scriptstyle \bullet$} [c] at 10 6
\put {$ \scriptstyle \bullet$} [c] at 10 3
\put {$ \scriptstyle \bullet$} [c] at 10 0
\put {$ \scriptstyle \bullet$} [c] at 13 12
\put {$ \scriptstyle \bullet$} [c] at 12 0
\put {$ \scriptstyle \bullet$} [c] at 14 0
\put {$ \scriptstyle \bullet$} [c] at 16 0
\setlinear \plot 10 0 10 6 13 12 12 0      /
\setlinear \plot 14 0 13 12 16 0      /
\put{$840$} [c] at 13 -2
\endpicture
\end{minipage}
$$
$$
\begin{minipage}{4cm}
\beginpicture
\setcoordinatesystem units    <1.5mm,2mm>
\setplotarea x from 0 to 16, y from -2 to 15
\put{1.789)} [l] at 2 12
\put {$ \scriptstyle \bullet$} [c] at 10 0
\put {$ \scriptstyle \bullet$} [c] at 10 6
\put {$ \scriptstyle \bullet$} [c] at 10 12
\put {$ \scriptstyle \bullet$} [c] at 14 0
\put {$ \scriptstyle \bullet$} [c] at 14 12
\setlinear \plot  10 12  10 0   14 12 14  0   /
\put {$ \scriptstyle \bullet$} [c] at 16 0
\put {$ \scriptstyle \bullet$} [c] at 16 12
\setlinear \plot 16 0 16 12   /
\put{$5{,}040$} [c] at 13 -2
\endpicture
\end{minipage}
\begin{minipage}{4cm}
\beginpicture
\setcoordinatesystem units    <1.5mm,2mm>
\setplotarea x from 0 to 16, y from -2 to 15
\put{1.790)} [l] at 2 12
\put {$ \scriptstyle \bullet$} [c] at 10 0
\put {$ \scriptstyle \bullet$} [c] at 10 6
\put {$ \scriptstyle \bullet$} [c] at 10 12
\put {$ \scriptstyle \bullet$} [c] at 14 0
\put {$ \scriptstyle \bullet$} [c] at 14 12
\setlinear \plot  10 0  10 12   14 0 14  12   /
\put {$ \scriptstyle \bullet$} [c] at 16 0
\put {$ \scriptstyle \bullet$} [c] at 16 12
\setlinear \plot 16 0 16 12   /
\put{$5{,}040$} [c] at 13 -2
\endpicture
\end{minipage}
\begin{minipage}{4cm}
\beginpicture
\setcoordinatesystem units    <1.5mm,2mm>
\setplotarea x from 0 to 16, y from -2 to 15
\put{1.791)} [l] at 2 12
\put {$ \scriptstyle \bullet$} [c] at 10 12
\put {$ \scriptstyle \bullet$} [c] at 12.5 6
\put {$ \scriptstyle \bullet$} [c] at 13 0
\put {$ \scriptstyle \bullet$} [c] at 13 12
\put {$ \scriptstyle \bullet$} [c] at 13.5 6
\setlinear \plot  10 12   13 0 12.5 6 13 12 13.5 6 13 0   /
\put {$ \scriptstyle \bullet$} [c] at 16 0
\put {$ \scriptstyle \bullet$} [c] at 16 12
\setlinear \plot 16 0 16 12   /
\put{$2{,}520$} [c] at 13 -2
\endpicture
\end{minipage}
\begin{minipage}{4cm}
\beginpicture
\setcoordinatesystem units    <1.5mm,2mm>
\setplotarea x from 0 to 16, y from -2 to 15
\put{1.792)} [l] at 2 12
\put {$ \scriptstyle \bullet$} [c] at 10 0
\put {$ \scriptstyle \bullet$} [c] at 12.5 6
\put {$ \scriptstyle \bullet$} [c] at 13 0
\put {$ \scriptstyle \bullet$} [c] at 13 12
\put {$ \scriptstyle \bullet$} [c] at 13.5 6
\setlinear \plot  10 0   13 12 12.5 6 13 0 13.5 6 13 12   /
\put {$ \scriptstyle \bullet$} [c] at 16 0
\put {$ \scriptstyle \bullet$} [c] at 16 12
\setlinear \plot 16 0 16 12   /
\put{$2{,}520$} [c] at 13 -2
\endpicture
\end{minipage}
\begin{minipage}{4cm}
\beginpicture
\setcoordinatesystem units    <1.5mm,2mm>
\setplotarea x from 0 to 16, y from -2 to 15
\put{1.793)} [l] at 2 12
\put {$ \scriptstyle \bullet$} [c] at 10 12
\put {$ \scriptstyle \bullet$} [c] at 12 12
\put {$ \scriptstyle \bullet$} [c] at 14 12
\put {$ \scriptstyle \bullet$} [c] at 13 0
\put {$ \scriptstyle \bullet$} [c] at 13 6
\setlinear \plot 10 12 13 0 13  6  12 12    /
\setlinear \plot 13  6   14 12    /
\put {$ \scriptstyle \bullet$} [c] at 16 0
\put {$ \scriptstyle \bullet$} [c] at 16 12
\setlinear \plot 16 0 16 12   /
\put{$2{,}520$} [c] at 13 -2
\endpicture
\end{minipage}
\begin{minipage}{4cm}
\beginpicture
\setcoordinatesystem units    <1.5mm,2mm>
\setplotarea x from 0 to 16, y from -2 to 15
\put{1.794)} [l] at 2 12
\put {$ \scriptstyle \bullet$} [c] at 10 0
\put {$ \scriptstyle \bullet$} [c] at 12 0
\put {$ \scriptstyle \bullet$} [c] at 14 0
\put {$ \scriptstyle \bullet$} [c] at 13 12
\put {$ \scriptstyle \bullet$} [c] at 13 6
\setlinear \plot 10 0 13 12 13  6  12 0    /
\setlinear \plot 13  6   14 0    /
\put {$ \scriptstyle \bullet$} [c] at 16 0
\put {$ \scriptstyle \bullet$} [c] at 16 12
\setlinear \plot 16 0 16 12   /
\put{$2{,}520$} [c] at 13 -2
\endpicture
\end{minipage}
$$
$$
\begin{minipage}{4cm}
\beginpicture
\setcoordinatesystem units    <1.5mm,2mm>
\setplotarea x from 0 to 16, y from -2 to 15
\put{1.795)} [l] at 2 12
\put {$ \scriptstyle \bullet$} [c] at 10 12
\put {$ \scriptstyle \bullet$} [c] at 12 12
\put {$ \scriptstyle \bullet$} [c] at 14 12
\put {$ \scriptstyle \bullet$} [c] at 10 0
\put {$ \scriptstyle \bullet$} [c] at 14 0
\setlinear \plot 14 0 10 12 10 0 12 12 14 0 14 12 10 0  /
\put {$ \scriptstyle \bullet$} [c] at 16 0
\put {$ \scriptstyle \bullet$} [c] at 16 12
\setlinear \plot 16 0 16 12   /
\put{$420$} [c] at 13 -2
\endpicture
\end{minipage}
\begin{minipage}{4cm}
\beginpicture
\setcoordinatesystem units    <1.5mm,2mm>
\setplotarea x from 0 to 16, y from -2 to 15
\put{1.796)} [l] at 2 12
\put {$ \scriptstyle \bullet$} [c] at 10 12
\put {$ \scriptstyle \bullet$} [c] at 12 0
\put {$ \scriptstyle \bullet$} [c] at 14 12
\put {$ \scriptstyle \bullet$} [c] at 10 0
\put {$ \scriptstyle \bullet$} [c] at 14 0
\setlinear \plot 14 12 10 0 10 12 12 0 14 12 14 0 10 12  /
\put {$ \scriptstyle \bullet$} [c] at 16 0
\put {$ \scriptstyle \bullet$} [c] at 16 12
\setlinear \plot 16 0 16 12   /
\put{$420$} [c] at 13 -2
\endpicture
\end{minipage}
\begin{minipage}{4cm}
\beginpicture
\setcoordinatesystem units    <1.5mm,2mm>
\setplotarea x  from 0 to 16, y from -2 to 15
\put{${\bf  34}$} [l] at 2 15

\put{1.797)} [l] at 2 12
\put {$ \scriptstyle \bullet$} [c] at  10 0
\put {$ \scriptstyle \bullet$} [c] at  10 12
\put {$ \scriptstyle \bullet$} [c] at  13 12
\put {$ \scriptstyle \bullet$} [c] at  13 0
\put {$ \scriptstyle \bullet$} [c] at  14.5 6
\put {$ \scriptstyle \bullet$} [c] at  16 12
\put {$ \scriptstyle \bullet$} [c] at  16 0
\setlinear \plot  10 12 10 0 13 12 13 0   /
\setlinear \plot  13 12  16 0  16 12 /
\put{$5{,}040$} [c] at 13 -2
\endpicture
\end{minipage}
\begin{minipage}{4cm}
\beginpicture
\setcoordinatesystem units    <1.5mm,2mm>
\setplotarea x  from 0 to 16, y from -2 to 15
\put{1.798)} [l] at 2 12
\put {$ \scriptstyle \bullet$} [c] at  10 0
\put {$ \scriptstyle \bullet$} [c] at  10 12
\put {$ \scriptstyle \bullet$} [c] at  13 12
\put {$ \scriptstyle \bullet$} [c] at  13 0
\put {$ \scriptstyle \bullet$} [c] at  14.5 6
\put {$ \scriptstyle \bullet$} [c] at  16 12
\put {$ \scriptstyle \bullet$} [c] at  16 0
\setlinear \plot  10 0 10 12 13 0 13 12   /
\setlinear \plot  13 0  16 12  16 0 /
\put{$5{,}040$} [c] at 13 -2
\endpicture
\end{minipage}
\begin{minipage}{4cm}
\beginpicture
\setcoordinatesystem units    <1.5mm,2mm>
\setplotarea x  from 0 to 16, y from -2 to 15
\put {1.799)} [l]  at 2 12
\put {$ \scriptstyle \bullet$} [c] at  10 0
\put {$ \scriptstyle \bullet$} [c] at  10 12
\put {$ \scriptstyle \bullet$} [c] at  12 12
\put {$ \scriptstyle \bullet$} [c] at  12 6
\put {$ \scriptstyle \bullet$} [c] at  14 0
\put {$ \scriptstyle \bullet$} [c] at  14 12
\put {$ \scriptstyle \bullet$} [c] at  16  12
\setlinear \plot  10 0 10 12 14 0 16 12 /
\setlinear \plot  12 12 12 6 /
\setlinear \plot 14 0 14 12 /
\put{$2{,}520$} [c]  at 13 -2
\endpicture
\end{minipage}
\begin{minipage}{4cm}
\beginpicture
\setcoordinatesystem units    <1.5mm,2mm>
\setplotarea x  from 0 to 16, y from -2 to 15
\put {1.800)} [l]  at 2 12
\put {$ \scriptstyle \bullet$} [c] at  10 0
\put {$ \scriptstyle \bullet$} [c] at  10 12
\put {$ \scriptstyle \bullet$} [c] at  12 0
\put {$ \scriptstyle \bullet$} [c] at  12 6
\put {$ \scriptstyle \bullet$} [c] at  14 0
\put {$ \scriptstyle \bullet$} [c] at  14 12
\put {$ \scriptstyle \bullet$} [c] at  16  0
\setlinear \plot  10 12 10 0 14 12 16 0 /
\setlinear \plot  12 0 12 6 /
\setlinear \plot 14 0 14 12 /
\put{$2{,}520$} [c]  at 13 -2
\endpicture
\end{minipage}
$$
$$
\begin{minipage}{4cm}
\beginpicture
\setcoordinatesystem units    <1.5mm,2mm>
\setplotarea x  from 0 to 16, y from -2 to 15
\put {1.801)} [l]  at 2 12
\put {$ \scriptstyle \bullet$} [c] at  10 12
\put {$ \scriptstyle \bullet$} [c] at  12 12
\put {$ \scriptstyle \bullet$} [c] at  14 12
\put {$ \scriptstyle \bullet$} [c] at  16 12
\put {$ \scriptstyle \bullet$} [c] at  12 0
\put {$ \scriptstyle \bullet$} [c] at  12  6
\put {$ \scriptstyle \bullet$} [c] at  10  0
\setlinear \plot  10 12 10 0 12 12 12 0 10 12 /
\setlinear \plot  14 12 12 0 16 12 /
\put{$2{,}520$} [c]  at 13 -2
\endpicture
\end{minipage}
\begin{minipage}{4cm}
\beginpicture
\setcoordinatesystem units    <1.5mm,2mm>
\setplotarea x  from 0 to 16, y from -2 to 15
\put {1.802)} [l]  at 2 12
\put {$ \scriptstyle \bullet$} [c] at  10 0
\put {$ \scriptstyle \bullet$} [c] at  12 0
\put {$ \scriptstyle \bullet$} [c] at  14 0
\put {$ \scriptstyle \bullet$} [c] at  16 0
\put {$ \scriptstyle \bullet$} [c] at  12 12
\put {$ \scriptstyle \bullet$} [c] at  12  6
\put {$ \scriptstyle \bullet$} [c] at  10  12
\setlinear \plot  10 0 10 12 12 0 12 12 10 0 /
\setlinear \plot  14 0 12 12 16 0 /
\put{$2{,}520$} [c]  at 13 -2
\endpicture
\end{minipage}
\begin{minipage}{4cm}
\beginpicture
\setcoordinatesystem units    <1.5mm,2mm>
\setplotarea x  from 0 to 16, y from -2 to 15
\put{1.803)} [l] at 2 12
\put {$ \scriptstyle \bullet$} [c] at  10 0
\put {$ \scriptstyle \bullet$} [c] at  10 12
\put {$ \scriptstyle \bullet$} [c] at  13 0
\put {$ \scriptstyle \bullet$} [c] at  13 12
\put {$ \scriptstyle \bullet$} [c] at  16 0
\put {$ \scriptstyle \bullet$} [c] at  16 6
\put {$ \scriptstyle \bullet$} [c] at  16 12
\setlinear \plot  10 0 10  12 16 0 16 12   /
\setlinear \plot  16 0 13 12   /
\setlinear \plot  10 12 13 0   /
\put{$2{,}520$} [c] at 13 -2
\endpicture
\end{minipage}
\begin{minipage}{4cm}
\beginpicture
\setcoordinatesystem units    <1.5mm,2mm>
\setplotarea x  from 0 to 16, y from -2 to 15
\put{1.804)} [l] at 2 12
\put {$ \scriptstyle \bullet$} [c] at  10 0
\put {$ \scriptstyle \bullet$} [c] at  10 12
\put {$ \scriptstyle \bullet$} [c] at  13 0
\put {$ \scriptstyle \bullet$} [c] at  13 12
\put {$ \scriptstyle \bullet$} [c] at  16 0
\put {$ \scriptstyle \bullet$} [c] at  16 6
\put {$ \scriptstyle \bullet$} [c] at  16 12
\setlinear \plot  10 12 10  0 16 12 16 0   /
\setlinear \plot  16 12 13 0   /
\setlinear \plot  10 0 13 12   /
\put{$2{,}520$} [c] at 13 -2
\endpicture
\end{minipage}
\begin{minipage}{4cm}
\beginpicture
\setcoordinatesystem units    <1.5mm,2mm>
\setplotarea x  from 0 to 16, y from -2 to 15
\put {1.805)} [l] at 2 12
\put {$ \scriptstyle \bullet$} [c] at  10 0
\put {$ \scriptstyle \bullet$} [c] at  12 12
\put {$ \scriptstyle \bullet$} [c] at  14 12
\put {$ \scriptstyle \bullet$} [c] at  16 0
\put {$ \scriptstyle \bullet$} [c] at  10 12
\put {$ \scriptstyle \bullet$} [c] at  13 0
\put {$ \scriptstyle \bullet$} [c] at  16  12
\setlinear \plot   10 12 10 0  16 12 16 0  12 12 10 0   /
\setlinear \plot   12 12 13 0   /
\setlinear \plot   14 12 16 0    /
\put{$2{,}520$} [c]  at 13 -2
\endpicture
\end{minipage}
\begin{minipage}{4cm}
\beginpicture
\setcoordinatesystem units    <1.5mm,2mm>
\setplotarea x  from 0 to 16, y from -2 to 15
\put {1.806)} [l] at 2 12
\put {$ \scriptstyle \bullet$} [c] at  10 0
\put {$ \scriptstyle \bullet$} [c] at  12 0
\put {$ \scriptstyle \bullet$} [c] at  14 0
\put {$ \scriptstyle \bullet$} [c] at  16 0
\put {$ \scriptstyle \bullet$} [c] at  10 12
\put {$ \scriptstyle \bullet$} [c] at  13 12
\put {$ \scriptstyle \bullet$} [c] at  16  12
\setlinear \plot   10 0 10 12  16 0 16 12  12 0 10 12   /
\setlinear \plot   12 0 13 12   /
\setlinear \plot   14 0 16 12    /
\put{$2{,}520$} [c]  at 13 -2
\endpicture
\end{minipage}
$$
$$
\begin{minipage}{4cm}
\beginpicture
\setcoordinatesystem units    <1.5mm,2mm>
\setplotarea x  from 0 to 16, y from -2 to 15
\put {1.807)} [l] at 2 12
\put {$ \scriptstyle \bullet$} [c] at  10 0
\put {$ \scriptstyle \bullet$} [c] at  12 0
\put {$ \scriptstyle \bullet$} [c] at  14 0
\put {$ \scriptstyle \bullet$} [c] at  16 0
\put {$ \scriptstyle \bullet$} [c] at  10 12
\put {$ \scriptstyle \bullet$} [c] at  13 12
\put {$ \scriptstyle \bullet$} [c] at  16 12
\setlinear \plot   10 12 10  0 13 12 16 0 16 12  14 0 13 12 12 0   /
\put{$2{,}520$} [c]  at 13 -2
\endpicture
\end{minipage}
\begin{minipage}{4cm}
\beginpicture
\setcoordinatesystem units    <1.5mm,2mm>
\setplotarea x  from 0 to 16, y from -2 to 15
\put {1.808)} [l] at 2 12
\put {$ \scriptstyle \bullet$} [c] at  10 12
\put {$ \scriptstyle \bullet$} [c] at  12 12
\put {$ \scriptstyle \bullet$} [c] at  14 12
\put {$ \scriptstyle \bullet$} [c] at  16 12
\put {$ \scriptstyle \bullet$} [c] at  10 0
\put {$ \scriptstyle \bullet$} [c] at  13 0
\put {$ \scriptstyle \bullet$} [c] at  16 0
\setlinear \plot   10 0 10  12 13 0 16 12 16 0  14 12 13 0 12 12   /
\put{$2{,}520$} [c]  at 13 -2
\endpicture
\end{minipage}
\begin{minipage}{4cm}
\beginpicture
\setcoordinatesystem units    <1.5mm,2mm>
\setplotarea x  from 0 to 16, y from -2 to 15
\put {1.809)} [l] at 2 12
\put {$ \scriptstyle \bullet$} [c] at  10 0
\put {$ \scriptstyle \bullet$} [c] at  12 12
\put {$ \scriptstyle \bullet$} [c] at  14 12
\put {$ \scriptstyle \bullet$} [c] at  16 0
\put {$ \scriptstyle \bullet$} [c] at  10 12
\put {$ \scriptstyle \bullet$} [c] at  13 0
\put {$ \scriptstyle \bullet$} [c] at  16 12
\setlinear \plot   10 12 10  0 12 12 13 0 14 12  16 0 16 12   /
\put{$2{,}520$} [c]  at 13 -2
\endpicture
\end{minipage}
\begin{minipage}{4cm}
\beginpicture
\setcoordinatesystem units    <1.5mm,2mm>
\setplotarea x  from 0 to 16, y from -2 to 15
\put {1.810)}  [l] at 2 12
\put {$ \scriptstyle \bullet$} [c] at  10 12
\put {$ \scriptstyle \bullet$} [c] at  12 0
\put {$ \scriptstyle \bullet$} [c] at  14 0
\put {$ \scriptstyle \bullet$} [c] at  16 12
\put {$ \scriptstyle \bullet$} [c] at  10 0
\put {$ \scriptstyle \bullet$} [c] at  13 12
\put {$ \scriptstyle \bullet$} [c] at  16 0
\setlinear \plot   10 0 10  12 12 0 13 12 14 0  16 12 16 0   /
\put{$2{,}520$} [c]  at 13 -2
\endpicture
\end{minipage}
\begin{minipage}{4cm}
\beginpicture
\setcoordinatesystem units    <1.5mm,2mm>
\setplotarea x  from 0 to 16, y from -2 to 15
\put {1.811)} [l]  at 2 12
\put {$ \scriptstyle \bullet$} [c] at  10 12
\put {$ \scriptstyle \bullet$} [c] at  12 12
\put {$ \scriptstyle \bullet$} [c] at  14 12
\put {$ \scriptstyle \bullet$} [c] at  16 0
\put {$ \scriptstyle \bullet$} [c] at  12  0
\put {$ \scriptstyle \bullet$} [c] at  16  6
\put {$ \scriptstyle \bullet$} [c] at  16 12
\setlinear \plot  10 12 12 0  16 6 /
\setlinear \plot  12 12 12 0 14 12 /
\setlinear \plot  16 0 16  12 /
\put{$840$} [c]  at 13 -2
\endpicture
\end{minipage}
\begin{minipage}{4cm}
\beginpicture
\setcoordinatesystem units    <1.5mm,2mm>
\setplotarea x  from 0 to 16, y from -2 to 15
\put {1.812)} [l]  at 2 12
\put {$ \scriptstyle \bullet$} [c] at  10 0
\put {$ \scriptstyle \bullet$} [c] at  12 0
\put {$ \scriptstyle \bullet$} [c] at  14 0
\put {$ \scriptstyle \bullet$} [c] at  16 0
\put {$ \scriptstyle \bullet$} [c] at  12  12
\put {$ \scriptstyle \bullet$} [c] at  16  6
\put {$ \scriptstyle \bullet$} [c] at  16 12
\setlinear \plot  10 0 12 12  16 6 /
\setlinear \plot  12 0 12 12 14 0 /
\setlinear \plot  16 0 16  12 /
\put{$840$} [c]  at 13 -2
\endpicture
\end{minipage}
$$
$$
\begin{minipage}{4cm}
\beginpicture
\setcoordinatesystem units    <1.5mm,2mm>
\setplotarea x  from 0 to 16, y from -2 to 15
\put{1.813)} [l]  at 2 12
\put {$ \scriptstyle \bullet$} [c] at 10 6
\put {$ \scriptstyle \bullet$} [c] at 12 6
\put {$ \scriptstyle \bullet$} [c] at 13 6
\put {$ \scriptstyle \bullet$} [c] at 14 6
\put {$ \scriptstyle \bullet$} [c] at 16 6
\put {$ \scriptstyle \bullet$} [c] at 13 12
\put {$ \scriptstyle \bullet$} [c] at 13 0
\setlinear \plot 13 0 10 6 13 12 16 6 13 0     /
\setlinear \plot 13 0 12 6 13 12 14 6 13 0     /
\setlinear \plot 13  0 13 12     /
\put{$42$} [c] at 13 -2
\endpicture
\end{minipage}
\begin{minipage}{4cm}
\beginpicture
\setcoordinatesystem units    <1.5mm,2mm>
\setplotarea x  from 0 to 16, y from -2 to 15
\put{1.814)} [l]  at 2 12
\put {$ \scriptstyle \bullet$} [c] at 10 12
\put {$ \scriptstyle \bullet$} [c] at 12 12
\put {$ \scriptstyle \bullet$} [c] at 13 12
\put {$ \scriptstyle \bullet$} [c] at 14 12
\put {$ \scriptstyle \bullet$} [c] at 16 12
\put {$ \scriptstyle \bullet$} [c] at 13 6
\put {$ \scriptstyle \bullet$} [c] at 13 0
\setlinear \plot 13 0 13  12     /
\setlinear \plot 10 12 13 6 16 12      /
\setlinear \plot 12 12 13 6 14 12      /
\put{$42$} [c] at 13 -2
\endpicture
\end{minipage}
\begin{minipage}{4cm}
\beginpicture
\setcoordinatesystem units    <1.5mm,2mm>
\setplotarea x  from 0 to 16, y from -2 to 15
\put{1.815)} [l]  at 2 12
\put {$ \scriptstyle \bullet$} [c] at 10 0
\put {$ \scriptstyle \bullet$} [c] at 12 0
\put {$ \scriptstyle \bullet$} [c] at 13 0
\put {$ \scriptstyle \bullet$} [c] at 14 0
\put {$ \scriptstyle \bullet$} [c] at 16 0
\put {$ \scriptstyle \bullet$} [c] at 13 6
\put {$ \scriptstyle \bullet$} [c] at 13 12
\setlinear \plot 13 0 13  12     /
\setlinear \plot 10 0 13 6 16 0      /
\setlinear \plot 12 0 13 6 14 0      /
\put{$42$} [c] at 13 -2
\endpicture
\end{minipage}
\begin{minipage}{4cm}
\beginpicture
\setcoordinatesystem units    <1.5mm,2mm>
\setplotarea x  from 0 to 16, y  from -2 to 15
\put{1.816)} [l]  at 2 12
\put {$ \scriptstyle \bullet$} [c] at 10 12
\put {$ \scriptstyle \bullet$} [c] at 11.5  8
\put {$ \scriptstyle \bullet$} [c] at 12 0
\put {$ \scriptstyle \bullet$} [c] at 12 12
\put {$ \scriptstyle \bullet$} [c] at 12.5 8
\put {$ \scriptstyle \bullet$} [c] at 14 12
\setlinear \plot 10 12 12 0  11.5  8 12 12 12.5 8 12 0  /
\setlinear \plot 12.5 8 14 12 /
\put{$5{,}040$} [c] at 13 -2
\put{$\scriptstyle \bullet$} [c] at 16  0 \endpicture
\end{minipage}
\begin{minipage}{4cm}
\beginpicture
\setcoordinatesystem units    <1.5mm,2mm>
\setplotarea x  from 0 to 16, y  from -2 to 15
\put{1.817)} [l]  at 2 12
\put {$ \scriptstyle \bullet$} [c] at 10 0
\put {$ \scriptstyle \bullet$} [c] at 11.5  4
\put {$ \scriptstyle \bullet$} [c] at 12 0
\put {$ \scriptstyle \bullet$} [c] at 12 12
\put {$ \scriptstyle \bullet$} [c] at 12.5 4
\put {$ \scriptstyle \bullet$} [c] at 14 0
\setlinear \plot 10 0 12 12  11.5  4 12 0 12.5 4 12 12  /
\setlinear \plot 12.5 4 14 0 /
\put{$5{,}040$} [c] at 13 -2
\put{$\scriptstyle \bullet$} [c] at 16  0 \endpicture
\end{minipage}
\begin{minipage}{4cm}
\beginpicture
\setcoordinatesystem units    <1.5mm,2mm>
\setplotarea x  from 0 to 16, y  from -2 to 15
\put{1.818)} [l]  at 2 12
\put {$ \scriptstyle \bullet$} [c] at 10 0
\put {$ \scriptstyle \bullet$} [c] at 10 6
\put {$ \scriptstyle \bullet$} [c] at 10 12
\put {$ \scriptstyle \bullet$} [c] at 12 6
\put {$ \scriptstyle \bullet$} [c] at 14 12
\put {$ \scriptstyle \bullet$} [c] at 14 0
\setlinear \plot  10 12 10 0 14 12 14 0   /
\put{$5{,}040$} [c] at 13 -2
\put{$\scriptstyle \bullet$} [c] at 16  0 \endpicture
\end{minipage}
$$
$$
\begin{minipage}{4cm}
\beginpicture
\setcoordinatesystem units    <1.5mm,2mm>
\setplotarea x  from 0 to 16, y  from -2 to 15
\put{1.819)} [l]  at 2 12
\put {$ \scriptstyle \bullet$} [c] at 10 0
\put {$ \scriptstyle \bullet$} [c] at 10 6
\put {$ \scriptstyle \bullet$} [c] at 10 12
\put {$ \scriptstyle \bullet$} [c] at 12 6
\put {$ \scriptstyle \bullet$} [c] at 14 12
\put {$ \scriptstyle \bullet$} [c] at 14 0
\setlinear \plot  10 0 10 12 14 0 14 12   /
\put{$5{,}040$} [c] at 13 -2
\put{$\scriptstyle \bullet$} [c] at 16  0 \endpicture
\end{minipage}
\begin{minipage}{4cm}
\beginpicture
\setcoordinatesystem units    <1.5mm,2mm>
\setplotarea x  from 0 to 16, y  from -2 to 15
\put{1.820)} [l]  at 2 12
\put {$ \scriptstyle \bullet$} [c] at 10 0
\put {$ \scriptstyle \bullet$} [c] at 10 12
\put {$ \scriptstyle \bullet$} [c] at 12 12
\put {$ \scriptstyle \bullet$} [c] at 14 0
\put {$ \scriptstyle \bullet$} [c] at 14 6
\put {$ \scriptstyle \bullet$} [c] at 14 12
\setlinear \plot 10 12 10 0  12 12 14 0 14 12 10 0 /
\put{$5{,}040$} [c] at 13 -2
\put{$\scriptstyle \bullet$} [c] at 16  0 \endpicture
\end{minipage}
\begin{minipage}{4cm}
\beginpicture
\setcoordinatesystem units    <1.5mm,2mm>
\setplotarea x  from 0 to 16, y  from -2 to 15
\put{1.821)} [l]  at 2 12
\put {$ \scriptstyle \bullet$} [c] at 10 0
\put {$ \scriptstyle \bullet$} [c] at 10 12
\put {$ \scriptstyle \bullet$} [c] at 12 0
\put {$ \scriptstyle \bullet$} [c] at 14 0
\put {$ \scriptstyle \bullet$} [c] at 14 6
\put {$ \scriptstyle \bullet$} [c] at 14 12
\setlinear \plot 10 0 10 12  12 0 14 12 14 0 10 12 /
\put{$5{,}040$} [c] at 13 -2
\put{$\scriptstyle \bullet$} [c] at 16  0 \endpicture
\end{minipage}
\begin{minipage}{4cm}
\beginpicture
\setcoordinatesystem units    <1.5mm,2mm>
\setplotarea x  from 0 to 16, y  from -2 to 15
\put{1.822)} [l]  at 2 12
\put {$ \scriptstyle \bullet$} [c] at 10 12
\put {$ \scriptstyle \bullet$} [c] at 12 0
\put {$ \scriptstyle \bullet$} [c] at 12 4
\put {$ \scriptstyle \bullet$} [c] at 12 8
\put {$ \scriptstyle \bullet$} [c] at 12 12
\put {$ \scriptstyle \bullet$} [c] at 14 12
\setlinear \plot  10  12 12 0 14 12    /
\setlinear \plot  12  0  12 12    /
\put{$2{,}520$} [c] at 13 -2
\put{$\scriptstyle \bullet$} [c] at 16  0 \endpicture
\end{minipage}
\begin{minipage}{4cm}
\beginpicture
\setcoordinatesystem units    <1.5mm,2mm>
\setplotarea x  from 0 to 16, y  from -2 to 15
\put{1.823)} [l]  at 2 12
\put {$ \scriptstyle \bullet$} [c] at 10 0
\put {$ \scriptstyle \bullet$} [c] at 12 0
\put {$ \scriptstyle \bullet$} [c] at 12 4
\put {$ \scriptstyle \bullet$} [c] at 12 8
\put {$ \scriptstyle \bullet$} [c] at 12 12
\put {$ \scriptstyle \bullet$} [c] at 14 0
\setlinear \plot  10  0 12 12 14 0    /
\setlinear \plot  12  0  12 12    /
\put{$2{,}520$} [c] at 13 -2
\put{$\scriptstyle \bullet$} [c] at 16  0 \endpicture
\end{minipage}
\begin{minipage}{4cm}
\beginpicture
\setcoordinatesystem units    <1.5mm,2mm>
\setplotarea x  from 0 to 16, y  from -2 to 15
\put{1.824)} [l]  at 2 12
\put {$ \scriptstyle \bullet$} [c] at 10 6
\put {$ \scriptstyle \bullet$} [c] at 10.5  0
\put {$ \scriptstyle \bullet$} [c] at 10.5  12
\put {$ \scriptstyle \bullet$} [c] at 11 6
\put {$ \scriptstyle \bullet$} [c] at 14 12
\put {$ \scriptstyle \bullet$} [c] at 14 0
\setlinear \plot 14 0 14  12 10.5 0  10 6 10.5 12 11 6 10.5 0  /
\put{$2{,}520$} [c] at 13 -2
\put{$\scriptstyle \bullet$} [c] at 16  0 \endpicture
\end{minipage}
$$
$$
\begin{minipage}{4cm}
\beginpicture
\setcoordinatesystem units    <1.5mm,2mm>
\setplotarea x  from 0 to 16, y  from -2 to 15
\put{1.825)} [l]  at 2 12
\put {$ \scriptstyle \bullet$} [c] at 10 6
\put {$ \scriptstyle \bullet$} [c] at 10.5  0
\put {$ \scriptstyle \bullet$} [c] at 10.5  12
\put {$ \scriptstyle \bullet$} [c] at 11 6
\put {$ \scriptstyle \bullet$} [c] at 14 12
\put {$ \scriptstyle \bullet$} [c] at 14 0
\setlinear \plot 14 12 14  0 10.5 12  10 6 10.5 0 11 6 10.5 12  /
\put{$2{,}520$} [c] at 13 -2
\put{$\scriptstyle \bullet$} [c] at 16  0 \endpicture
\end{minipage}
\begin{minipage}{4cm}
\beginpicture
\setcoordinatesystem units    <1.5mm,2mm>
\setplotarea x  from 0 to 16, y  from -2 to 15
\put{1.826)} [l]  at 2 12
\put {$ \scriptstyle \bullet$} [c] at 10 0
\put {$ \scriptstyle \bullet$} [c] at 10 12
\put {$ \scriptstyle \bullet$} [c] at 12 6
\put {$ \scriptstyle \bullet$} [c] at 12 12
\put {$ \scriptstyle \bullet$} [c] at 14 12
\put {$ \scriptstyle \bullet$} [c] at 14 0
\setlinear \plot  10 12 10  0  12 6 12  12    /
\setlinear \plot  12 6 14  0  14 12    /
\put{$2{,}520$} [c] at 13 -2
\put{$\scriptstyle \bullet$} [c] at 16  0 \endpicture
\end{minipage}
\begin{minipage}{4cm}
\beginpicture
\setcoordinatesystem units    <1.5mm,2mm>
\setplotarea x  from 0 to 16, y  from -2 to 15
\put{1.827)} [l]  at 2 12
\put {$ \scriptstyle \bullet$} [c] at 10 0
\put {$ \scriptstyle \bullet$} [c] at 10 12
\put {$ \scriptstyle \bullet$} [c] at 12 6
\put {$ \scriptstyle \bullet$} [c] at 12 0
\put {$ \scriptstyle \bullet$} [c] at 14 12
\put {$ \scriptstyle \bullet$} [c] at 14 0
\setlinear \plot  10 0 10 12  12 6 12  0    /
\setlinear \plot  12 6 14  12  14 0    /
\put{$2{,}520$} [c] at 13 -2
\put{$\scriptstyle \bullet$} [c] at 16  0 \endpicture
\end{minipage}
\begin{minipage}{4cm}
\beginpicture
\setcoordinatesystem units    <1.5mm,2mm>
\setplotarea x  from 0 to 16, y  from -2 to 15
\put{1.828)} [l]  at 2 12
\put {$ \scriptstyle \bullet$} [c] at 10 0
\put {$ \scriptstyle \bullet$} [c] at 10 12
\put {$ \scriptstyle \bullet$} [c] at 12 0
\put {$ \scriptstyle \bullet$} [c] at 12 12
\put {$ \scriptstyle \bullet$} [c] at 13 6
\put {$ \scriptstyle \bullet$} [c] at 14 12
\setlinear \plot  14 12 12  0 12 12 10 0  10 12 12 0   /
\put{$2{,}520$} [c] at 13 -2
\put{$\scriptstyle \bullet$} [c] at 16  0 \endpicture
\end{minipage}
\begin{minipage}{4cm}
\beginpicture
\setcoordinatesystem units    <1.5mm,2mm>
\setplotarea x  from 0 to 16, y  from -2 to 15
\put{1.829)} [l]  at 2 12
\put {$ \scriptstyle \bullet$} [c] at 10 0
\put {$ \scriptstyle \bullet$} [c] at 10 12
\put {$ \scriptstyle \bullet$} [c] at 12 0
\put {$ \scriptstyle \bullet$} [c] at 12 12
\put {$ \scriptstyle \bullet$} [c] at 13 6
\put {$ \scriptstyle \bullet$} [c] at 14 0
\setlinear \plot  14 0 12  12 12 0 10 12  10 0 12 12   /
\put{$2{,}520$} [c] at 13 -2
\put{$\scriptstyle \bullet$} [c] at 16  0 \endpicture
\end{minipage}
\begin{minipage}{4cm}
\beginpicture
\setcoordinatesystem units    <1.5mm,2mm>
\setplotarea x  from 0 to 16, y  from -2 to 15
\put{1.830)} [l]  at 2 12
\put {$ \scriptstyle \bullet$} [c] at 10 0
\put {$ \scriptstyle \bullet$} [c] at 10 12
\put {$ \scriptstyle \bullet$} [c] at 13 6
\put {$ \scriptstyle \bullet$} [c] at 13 0
\put {$ \scriptstyle \bullet$} [c] at 12 12
\put {$ \scriptstyle \bullet$} [c] at 14 12
\setlinear \plot 10 0 10 12 13 0  13 6 12 12  /
\setlinear \plot 13 6 14 12  /
\put{$2{,}520$} [c] at 13 -2
\put{$\scriptstyle \bullet$} [c] at 16  0 \endpicture
\end{minipage}
$$
$$
\begin{minipage}{4cm}
\beginpicture
\setcoordinatesystem units    <1.5mm,2mm>
\setplotarea x  from 0 to 16, y  from -2 to 15
\put{1.831)} [l]  at 2 12
\put {$ \scriptstyle \bullet$} [c] at 10 0
\put {$ \scriptstyle \bullet$} [c] at 10 12
\put {$ \scriptstyle \bullet$} [c] at 13 6
\put {$ \scriptstyle \bullet$} [c] at 13 12
\put {$ \scriptstyle \bullet$} [c] at 12 0
\put {$ \scriptstyle \bullet$} [c] at 14 0
\setlinear \plot 10 12 10 0 13 12  13 6 12 0  /
\setlinear \plot 13 6 14 0  /
\put{$2{,}520$} [c] at 13 -2
\put{$\scriptstyle \bullet$} [c] at 16  0 \endpicture
\end{minipage}
\begin{minipage}{4cm}
\beginpicture
\setcoordinatesystem units    <1.5mm,2mm>
\setplotarea x  from 0 to 16, y  from -2 to 15
\put{1.832)} [l]  at 2 12
\put {$ \scriptstyle \bullet$} [c] at 10 0
\put {$ \scriptstyle \bullet$} [c] at 10 12
\put {$ \scriptstyle \bullet$} [c] at 12 12
\put {$ \scriptstyle \bullet$} [c] at 12 0
\put {$ \scriptstyle \bullet$} [c] at 14 12
\put {$ \scriptstyle \bullet$} [c] at 14 0
\setlinear \plot  10  0  10  12 12 0 14 12 14 0  12 12  10 0  /
\setlinear \plot  12  0  12 12    /
\put{$2{,}520$} [c] at 13 -2
\put{$\scriptstyle \bullet$} [c] at 16  0 \endpicture
\end{minipage}
\begin{minipage}{4cm}
\beginpicture
\setcoordinatesystem units    <1.5mm,2mm>
\setplotarea x  from 0 to 16, y from -2 to 15
\put {${\bf  35}$} [l] at 2 15

\put {1.833)} [l] at 2 12
\put {$ \scriptstyle \bullet$} [c] at  10 0
\put {$ \scriptstyle \bullet$} [c] at  10 12
\put {$ \scriptstyle \bullet$} [c] at  12 12
\put {$ \scriptstyle \bullet$} [c] at  13 6
\put {$ \scriptstyle \bullet$} [c] at  14 12
\put {$ \scriptstyle \bullet$} [c] at  14 0
\put {$ \scriptstyle \bullet$} [c] at  16  12
\setlinear \plot  10 12 10 0 12 12 14 0 16 12 /
\setlinear \plot 14 0 14 12 /
\put{$2{,}520$} [c]  at 13 -2
\endpicture
\end{minipage}
\begin{minipage}{4cm}
\beginpicture
\setcoordinatesystem units    <1.5mm,2mm>
\setplotarea x  from 0 to 16, y from -2 to 15
\put {1.834)} [l] at 2 12
\put {$ \scriptstyle \bullet$} [c] at  10 0
\put {$ \scriptstyle \bullet$} [c] at  10 12
\put {$ \scriptstyle \bullet$} [c] at  12 0
\put {$ \scriptstyle \bullet$} [c] at  13 6
\put {$ \scriptstyle \bullet$} [c] at  14 12
\put {$ \scriptstyle \bullet$} [c] at  14 0
\put {$ \scriptstyle \bullet$} [c] at  16  0
\setlinear \plot  10 0 10 12 12 0 14 12 16 0 /
\setlinear \plot 14 0 14 12 /
\put{$2{,}520$} [c]  at 13 -2
\endpicture
\end{minipage}
\begin{minipage}{4cm}
\beginpicture
\setcoordinatesystem units    <1.5mm,2mm>
\setplotarea x from 0 to 16, y from -2 to 15
\put {1.835)} [l] at 2 12
\put {$ \scriptstyle \bullet$} [c] at  10 12
\put {$ \scriptstyle \bullet$} [c] at  12 12
\put {$ \scriptstyle \bullet$} [c] at  14 12
\put {$ \scriptstyle \bullet$} [c] at  16  12
\put {$ \scriptstyle \bullet$} [c] at  10 0
\put {$ \scriptstyle \bullet$} [c] at  13 0
\put {$ \scriptstyle \bullet$} [c] at  16 0
\setlinear \plot   10 12 10 0  12 12 13 0 14 12    /
\setlinear \plot   12 12 16 0 16 12    /
\put{$840$} [c] at 13 -2
\endpicture
\end{minipage}
\begin{minipage}{4cm}
\beginpicture
\setcoordinatesystem units    <1.5mm,2mm>
\setplotarea x from 0 to 16, y from -2 to 15
\put {1.836)} [l] at 2 12
\put {$ \scriptstyle \bullet$} [c] at  10 0
\put {$ \scriptstyle \bullet$} [c] at  12 0
\put {$ \scriptstyle \bullet$} [c] at  14 0
\put {$ \scriptstyle \bullet$} [c] at  16  0
\put {$ \scriptstyle \bullet$} [c] at  10 12
\put {$ \scriptstyle \bullet$} [c] at  13 12
\put {$ \scriptstyle \bullet$} [c] at  16 12
\setlinear \plot   10 0 10 12  12 0 13 12 14 0    /
\setlinear \plot   12 0 16 12 16 0    /
\put{$840$} [c] at 13 -2
\endpicture
\end{minipage}
$$
$$
\begin{minipage}{4cm}
\beginpicture
\setcoordinatesystem units    <1.5mm,2mm>
\setplotarea x  from 0 to 16, y from -2 to 15
\put {1.837)} [l] at 2 12
\put {$ \scriptstyle \bullet$} [c] at  10 12
\put {$ \scriptstyle \bullet$} [c] at  12 12
\put {$ \scriptstyle \bullet$} [c] at  14 12
\put {$ \scriptstyle \bullet$} [c] at  16 12
\put {$ \scriptstyle \bullet$} [c] at  16 6
\put {$ \scriptstyle \bullet$} [c] at  16 0
\put {$ \scriptstyle \bullet$} [c] at  13 0
\setlinear \plot  10 12 13  0  16 12 16 0 /
\setlinear \plot  12 12 13 0 14  12 /
\put{$840$} [c]  at 13 -2
\endpicture
\end{minipage}
\begin{minipage}{4cm}
\beginpicture
\setcoordinatesystem units    <1.5mm,2mm>
\setplotarea x  from 0 to 16, y from -2 to 15
\put {1.838)} [l] at 2 12
\put {$ \scriptstyle \bullet$} [c] at  10 0
\put {$ \scriptstyle \bullet$} [c] at  12 0
\put {$ \scriptstyle \bullet$} [c] at  14 0
\put {$ \scriptstyle \bullet$} [c] at  16 0
\put {$ \scriptstyle \bullet$} [c] at  16 6
\put {$ \scriptstyle \bullet$} [c] at  16 12
\put {$ \scriptstyle \bullet$} [c] at  13 12
\setlinear \plot  10 0 13  12  16 0 16 12 /
\setlinear \plot  12 0 13 12 14  0 /
\put{$840$} [c]  at 13 -2
\endpicture
\end{minipage}
\begin{minipage}{4cm}
\beginpicture
\setcoordinatesystem units    <1.5mm,2mm>
\setplotarea x from 0 to 16, y from -2 to 15
\put{1.839)} [l] at 2 12
\put {$ \scriptstyle \bullet$} [c] at 10 12
\put {$ \scriptstyle \bullet$} [c] at 11.5 12
\put {$ \scriptstyle \bullet$} [c] at 12.5 12
\put {$ \scriptstyle \bullet$} [c] at 14 12
\put {$ \scriptstyle \bullet$} [c] at 16 12
\put {$ \scriptstyle \bullet$} [c] at 12 6
\put {$ \scriptstyle \bullet$} [c] at 12  0
\setlinear \plot 10 12 12 6 12 0 16  12     /
\setlinear \plot 11.5 12 12 6 14  12     /
\setlinear \plot  12 6 12.5  12     /
\put{$210$} [c] at 13 -2
\endpicture
\end{minipage}
\begin{minipage}{4cm}
\beginpicture
\setcoordinatesystem units    <1.5mm,2mm>
\setplotarea x from 0 to 16, y from -2 to 15
\put{1.840)} [l] at 2 12
\put {$ \scriptstyle \bullet$} [c] at 10 0
\put {$ \scriptstyle \bullet$} [c] at 11.5 0
\put {$ \scriptstyle \bullet$} [c] at 12.5 0
\put {$ \scriptstyle \bullet$} [c] at 14 0
\put {$ \scriptstyle \bullet$} [c] at 16 0
\put {$ \scriptstyle \bullet$} [c] at 12 6
\put {$ \scriptstyle \bullet$} [c] at 12  12
\setlinear \plot 10 0 12 6 12 12 16  0     /
\setlinear \plot 11.5 0 12 6 14  0     /
\setlinear \plot  12 6 12.5  0     /
\put{$210$} [c] at 13 -2
\endpicture
\end{minipage}
\begin{minipage}{4cm}
\beginpicture
\setcoordinatesystem units    <1.5mm,2mm>
\setplotarea x from 0 to 16, y from -2 to 15
\put{1.841)} [l] at 2 12
\put {$ \scriptstyle \bullet$} [c] at 10 6
\put {$ \scriptstyle \bullet$} [c] at 11 6
\put {$ \scriptstyle \bullet$} [c] at 12 6
\put {$ \scriptstyle \bullet$} [c] at 13 6
\put {$ \scriptstyle \bullet$} [c] at 16 12
\put {$ \scriptstyle \bullet$} [c] at 11.5 12
\put {$ \scriptstyle \bullet$} [c] at 11.5  0
\setlinear \plot 16 12 11.5 0 10 6 11.5 12 11 6 11.5 0 12 6 11.5 12 13 6 11.5 0      /
\put{$210$} [c] at 13 -2
\endpicture
\end{minipage}
\begin{minipage}{4cm}
\beginpicture
\setcoordinatesystem units    <1.5mm,2mm>
\setplotarea x from 0 to 16, y from -2 to 15
\put{1.842)} [l] at 2 12
\put {$ \scriptstyle \bullet$} [c] at 10 6
\put {$ \scriptstyle \bullet$} [c] at 11 6
\put {$ \scriptstyle \bullet$} [c] at 12 6
\put {$ \scriptstyle \bullet$} [c] at 13 6
\put {$ \scriptstyle \bullet$} [c] at 16 0
\put {$ \scriptstyle \bullet$} [c] at 11.5 12
\put {$ \scriptstyle \bullet$} [c] at 11.5  0
\setlinear \plot 16 0 11.5 12 10 6 11.5 0 11 6 11.5 12 12 6 11.5 0 13 6 11.5 12      /
\put{$210$} [c] at 13 -2
\endpicture
\end{minipage}
$$
$$
\begin{minipage}{4cm}
\beginpicture
\setcoordinatesystem units    <1.5mm,2mm>
\setplotarea x  from 0 to 16, y from -2 to 15
\put {1.843)}  [l] at 2 12
\put {$ \scriptstyle \bullet$} [c] at  10 12
\put {$ \scriptstyle \bullet$} [c] at  11.5 12
\put {$ \scriptstyle \bullet$} [c] at  13 12
\put {$ \scriptstyle \bullet$} [c] at  14.5 12
\put {$ \scriptstyle \bullet$} [c] at  16  12
\put {$ \scriptstyle \bullet$} [c] at  11.5  0
\put {$ \scriptstyle \bullet$} [c] at  14.5 0
\setlinear \plot  10 12 11.5 0 13 12 14.5 0 16 12 11.5 0 11.5  12 14.5 0 14.5 12 11.5 0  /
\setlinear \plot  10 12 14.5 0 /
\put{$21$} [c] at 13 -2
\endpicture
\end{minipage}
\begin{minipage}{4cm}
\beginpicture
\setcoordinatesystem units    <1.5mm,2mm>
\setplotarea x  from 0 to 16, y from -2 to 15
\put {1.844)}  [l] at 2 12
\put {$ \scriptstyle \bullet$} [c] at  10 0
\put {$ \scriptstyle \bullet$} [c] at  11.5 0
\put {$ \scriptstyle \bullet$} [c] at  13 0
\put {$ \scriptstyle \bullet$} [c] at  14.5 0
\put {$ \scriptstyle \bullet$} [c] at  16  0
\put {$ \scriptstyle \bullet$} [c] at  11.5  12
\put {$ \scriptstyle \bullet$} [c] at  14.5 12
\setlinear \plot  10 0 11.5 12 13 0 14.5 12 16 0 11.5 12 11.5  0 14.5 12 14.5 0 11.5 12  /
\setlinear \plot  10 0 14.5 12 /
\put{$21$} [c] at 13 -2
\endpicture
\end{minipage}
\begin{minipage}{4cm}
\beginpicture
\setcoordinatesystem units    <1.5mm,2mm>
\setplotarea x  from 0 to 16, y from -2 to 15
\put{1.845)} [l]  at 2 12
\put {$ \scriptstyle \bullet$} [c] at 10 12
\put {$ \scriptstyle \bullet$} [c] at 10.5 6
\put {$ \scriptstyle \bullet$} [c] at 11 0
\put {$ \scriptstyle \bullet$} [c] at 12 12
\put {$ \scriptstyle \bullet$} [c] at 14 12
\put {$ \scriptstyle \bullet$} [c] at 15 0
\put{$\scriptstyle \bullet$} [c] at 16  12
\setlinear \plot   10  12 11 0 12 12   /
\setlinear \plot  14  12  15 0 16 12    /
\put{$2{,}520$} [c]  at 13 -2
\endpicture
\end{minipage}
\begin{minipage}{4cm}
\beginpicture
\setcoordinatesystem units    <1.5mm,2mm>
\setplotarea x  from 0 to 16, y from -2 to 15
\put{1.846)} [l]  at 2 12
\put {$ \scriptstyle \bullet$} [c] at 10 0
\put {$ \scriptstyle \bullet$} [c] at 10.5 6
\put {$ \scriptstyle \bullet$} [c] at 11 12
\put {$ \scriptstyle \bullet$} [c] at 12 0
\put {$ \scriptstyle \bullet$} [c] at 14 0
\put {$ \scriptstyle \bullet$} [c] at 15 12
\put{$\scriptstyle \bullet$} [c] at 16  0
\setlinear \plot   10  0 11 12 12 0   /
\setlinear \plot  14  0  15 12 16 0    /
\put{$2{,}520$} [c]  at 13 -2
\endpicture
\end{minipage}
\begin{minipage}{4cm}
\beginpicture
\setcoordinatesystem units    <1.5mm,2mm>
\setplotarea x  from 0 to 16, y from -2 to 15
\put{1.847)} [l]  at 2 12
\put {$ \scriptstyle \bullet$} [c] at 10 0
\put {$ \scriptstyle \bullet$} [c] at 10.5 6
\put {$ \scriptstyle \bullet$} [c] at 11 12
\put {$ \scriptstyle \bullet$} [c] at 12 0
\put {$ \scriptstyle \bullet$} [c] at 14 12
\put {$ \scriptstyle \bullet$} [c] at 15 0
\put{$\scriptstyle \bullet$} [c] at 16  12
\setlinear \plot   10  0 11 12 12 0   /
\setlinear \plot  14  12  15 0 16 12    /
\put{$2{,}520$} [c]  at 13 -2
\endpicture
\end{minipage}
\begin{minipage}{4cm}
\beginpicture
\setcoordinatesystem units    <1.5mm,2mm>
\setplotarea x  from 0 to 16, y from -2 to 15
\put{1.848)} [l]  at 2 12
\put {$ \scriptstyle \bullet$} [c] at 10 12
\put {$ \scriptstyle \bullet$} [c] at 10.5 6
\put {$ \scriptstyle \bullet$} [c] at 11 0
\put {$ \scriptstyle \bullet$} [c] at 12 12
\put {$ \scriptstyle \bullet$} [c] at 14 0
\put {$ \scriptstyle \bullet$} [c] at 15 12
\put{$\scriptstyle \bullet$} [c] at 16  0
\setlinear \plot   10  12 11 0 12 12   /
\setlinear \plot  14  0  15 12 16 0    /
\put{$2{,}520$} [c]  at 13 -2
\endpicture
\end{minipage}
$$
$$
\begin{minipage}{4cm}
\beginpicture
\setcoordinatesystem units    <1.5mm,2mm>
\setplotarea x  from 0 to 16, y from -2 to 15
\put{1.849)} [l]  at 2 12
\put {$ \scriptstyle \bullet$} [c] at 10 0
\put {$ \scriptstyle \bullet$} [c] at 10 12
\put {$ \scriptstyle \bullet$} [c] at 12 12
\put {$ \scriptstyle \bullet$} [c] at 12 0
\put {$ \scriptstyle \bullet$} [c] at 14 12
\put {$ \scriptstyle \bullet$} [c] at 15 0
\put{$\scriptstyle \bullet$} [c] at 16  12
\setlinear \plot   10  0 10 12 12 0 12 12 10 0   /
\setlinear \plot  14  12  15 0 16 12    /
\put{$630$} [c]  at 13 -2
\endpicture
\end{minipage}
\begin{minipage}{4cm}
\beginpicture
\setcoordinatesystem units    <1.5mm,2mm>
\setplotarea x  from 0 to 16, y from -2 to 15
\put{1.850)} [l]  at 2 12
\put {$ \scriptstyle \bullet$} [c] at 10 0
\put {$ \scriptstyle \bullet$} [c] at 10 12
\put {$ \scriptstyle \bullet$} [c] at 12 12
\put {$ \scriptstyle \bullet$} [c] at 12 0
\put {$ \scriptstyle \bullet$} [c] at 14 0
\put {$ \scriptstyle \bullet$} [c] at 15 12
\put{$\scriptstyle \bullet$} [c] at 16  0
\setlinear \plot   10  0 10 12 12 0 12 12 10 0   /
\setlinear \plot  14  0  15 12 16 0    /
\put{$630$} [c]  at 13 -2
\endpicture
\end{minipage}
\begin{minipage}{4cm}
\beginpicture
\setcoordinatesystem units <1.5mm, 2mm>
\setplotarea x from 0 to 16, y from -2 to 15

\put{${\bf  36}$} [l] at 2 15
\put {1.851)} [l] at 2 12
\put {$ \scriptstyle \bullet$} [c] at  10 12
\put {$ \scriptstyle \bullet$} [c] at  12 12
\put {$ \scriptstyle \bullet$} [c] at  14 12
\put {$ \scriptstyle \bullet$} [c] at  16  12
\put {$ \scriptstyle \bullet$} [c] at  10 0
\put {$ \scriptstyle \bullet$} [c] at  13  0
\put {$ \scriptstyle \bullet$} [c] at  16  0
\setlinear \plot   10 0 10  12  13 0 14   12 16 0 16 12   /
\setlinear \plot   12 12 13 0     /
\put{$5{,}040$} [c] at 13 -2
\endpicture
\end{minipage}
\begin{minipage}{4cm}
\beginpicture
\setcoordinatesystem units <1.5mm, 2mm>
\setplotarea x from 0 to 16, y from -2 to 15
\put {1.852)} [l] at 2 12
\put {$ \scriptstyle \bullet$} [c] at  10 0
\put {$ \scriptstyle \bullet$} [c] at  12 0
\put {$ \scriptstyle \bullet$} [c] at  14 0
\put {$ \scriptstyle \bullet$} [c] at  16   0
\put {$ \scriptstyle \bullet$} [c] at  10 12
\put {$ \scriptstyle \bullet$} [c] at  13  12
\put {$ \scriptstyle \bullet$} [c] at  16  12
\setlinear \plot   10 12 10  0  13 12 14   0 16 12 16 0   /
\setlinear \plot   12 0 13 12     /
\put{$5{,}040$} [c] at 13 -2
\endpicture
\end{minipage}
\begin{minipage}{4cm}
\beginpicture
\setcoordinatesystem units <1.5mm, 2mm>
\setplotarea x from 0 to 16, y from -2 to 15
\put {1.853)} [l] at 2 12
\put {$ \scriptstyle \bullet$} [c] at  10 12
\put {$ \scriptstyle \bullet$} [c] at  10 0
\put {$ \scriptstyle \bullet$} [c] at  11.5 0
\put {$ \scriptstyle \bullet$} [c] at  13   12
\put {$ \scriptstyle \bullet$} [c] at  13 0
\put {$ \scriptstyle \bullet$} [c] at  14  12
\put {$ \scriptstyle \bullet$} [c] at  16  12
\setlinear \plot   10 0 10  12   13 0 13 12 11.5 0 10 12  /
\setlinear \plot   14 12 13 0 16 12    /
\put{$2{,}520$} [c] at 12 -3
\endpicture
\end{minipage}
\begin{minipage}{4cm}
\beginpicture
\setcoordinatesystem units <1.5mm, 2mm>
\setplotarea x from 0 to 16, y from -2 to 15
\put {1.854)} [l] at 2 12
\put {$ \scriptstyle \bullet$} [c] at  10 0
\put {$ \scriptstyle \bullet$} [c] at  10 12
\put {$ \scriptstyle \bullet$} [c] at  11.5 12
\put {$ \scriptstyle \bullet$} [c] at  13   0
\put {$ \scriptstyle \bullet$} [c] at  13 12
\put {$ \scriptstyle \bullet$} [c] at  14  0
\put {$ \scriptstyle \bullet$} [c] at  16  0
\setlinear \plot   10 12 10  0   13 12 13  0 11.5 12 10 0  /
\setlinear \plot   14 0 13 12 16 0    /
\put{$2{,}520$} [c] at 13 -2
\endpicture
\end{minipage}
$$

$$
\begin{minipage}{4cm}
\beginpicture
\setcoordinatesystem units <1.5mm, 2mm>
\setplotarea x from 0 to 16, y from -2 to 15
\put{1.855)} [l] at 2 12
\put {$ \scriptstyle \bullet$} [c] at  10 0
\put {$ \scriptstyle \bullet$} [c] at  12.2 6
\put {$ \scriptstyle \bullet$} [c] at  13 0
\put {$ \scriptstyle \bullet$} [c] at  13 6
\put {$ \scriptstyle \bullet$} [c] at  13 12
\put {$ \scriptstyle \bullet$} [c] at  13.8 6
\put {$ \scriptstyle \bullet$} [c] at  16 12
\setlinear \plot 10 0 13 12 12.2 6 13 0 13 12 13.8  6 13  0  16 12 /
\put{$840$} [c] at 13 -2
\endpicture
\end{minipage}
\begin{minipage}{4cm}
\beginpicture
\setcoordinatesystem units <1.5mm, 2mm>
\setplotarea x from 0 to 16, y from -2 to 15
\put {1.856)} [l] at  2 12
\put {$ \scriptstyle \bullet$} [c] at  10 12
\put {$ \scriptstyle \bullet$} [c] at  12 12
\put {$ \scriptstyle \bullet$} [c] at  13 12
\put {$ \scriptstyle \bullet$} [c] at  14 12
\put {$ \scriptstyle \bullet$} [c] at  16  12
\put {$ \scriptstyle \bullet$} [c] at  10 0
\put {$ \scriptstyle \bullet$} [c] at  14  0
\setlinear \plot  10 12 10 0 12 12 14 0 14 12 10 0 12 12 14 0 /
\setlinear \plot 16 12 14 0 /
\setlinear \plot 10 12 14 0 13 12 /
\setlinear \plot 10 0 13 12 /
\put{ $210$} [c] at 13 -2
\endpicture
\end{minipage}
\begin{minipage}{4cm}
\beginpicture
\setcoordinatesystem units <1.5mm, 2mm>
\setplotarea x from 0 to 16, y from -2 to 15
\put {1.857)} [l] at  2 12
\put {$ \scriptstyle \bullet$} [c] at  10 0
\put {$ \scriptstyle \bullet$} [c] at  12 0
\put {$ \scriptstyle \bullet$} [c] at  13 0
\put {$ \scriptstyle \bullet$} [c] at  14 0
\put {$ \scriptstyle \bullet$} [c] at  16  0
\put {$ \scriptstyle \bullet$} [c] at  10 12
\put {$ \scriptstyle \bullet$} [c] at  14  12
\setlinear \plot  10 0 10 12 12 0 14 12 14 0 10 12 12 0 14 12 /
\setlinear \plot 16 0 14 12 /
\setlinear \plot 10 0 14 12 13 0 /
\setlinear \plot 10 12 13 0 /
\put{ $210$} [c] at 13 -2
\endpicture
\end{minipage}
\begin{minipage}{4cm}
\beginpicture
\setcoordinatesystem units   <1.5mm,2mm>
\setplotarea x from 0 to 16, y from -2 to 15
\put{1.858)} [l] at 2 12
\put {$ \scriptstyle \bullet$} [c] at 10 0
\put {$ \scriptstyle \bullet$} [c] at 10 12
\put {$ \scriptstyle \bullet$} [c] at 12 12
\put {$ \scriptstyle \bullet$} [c] at 14 12
\put {$ \scriptstyle \bullet$} [c] at 14 0
\put {$ \scriptstyle \bullet$} [c] at 14 6
\setlinear \plot 10 12 10 0 12 12 14 0 14 12   /
\put{$5{,}040$} [c]  at 13 -2
\put{$\scriptstyle \bullet$} [c] at 16  0 \endpicture
\end{minipage}
\begin{minipage}{4cm}
\beginpicture
\setcoordinatesystem units   <1.5mm,2mm>
\setplotarea x from 0 to 16, y from -2 to 15
\put{1.859)}  [l] at 2 12
\put {$ \scriptstyle \bullet$} [c] at 10 0
\put {$ \scriptstyle \bullet$} [c] at 10 12
\put {$ \scriptstyle \bullet$} [c] at 12 0
\put {$ \scriptstyle \bullet$} [c] at 14 12
\put {$ \scriptstyle \bullet$} [c] at 14 0
\put {$ \scriptstyle \bullet$} [c] at 14 6
\setlinear \plot 10 0 10 12 12 0 14 12 14 0   /
\put{$5{,}040$} [c]  at 13 -2
\put{$\scriptstyle \bullet$} [c] at 16  0 \endpicture
\end{minipage}
\begin{minipage}{4cm}
\beginpicture
\setcoordinatesystem units   <1.5mm,2mm>
\setplotarea x from 0 to 16, y from -2 to 15
\put{1.860)}  [l] at 2 12
\put {$ \scriptstyle \bullet$} [c] at 10 0
\put {$ \scriptstyle \bullet$} [c] at 10 12
\put {$ \scriptstyle \bullet$} [c] at 12 6
\put {$ \scriptstyle \bullet$} [c] at 12 12
\put {$ \scriptstyle \bullet$} [c] at 14 0
\put {$ \scriptstyle \bullet$} [c] at 14 12
\setlinear \plot 10 0  10 12 14 0 14 12  /
\setlinear \plot 12 6 12  12  /
\put{$5{,}040$} [c]  at 13 -2
\put{$\scriptstyle \bullet$} [c] at 16  0 \endpicture
\end{minipage}
$$
$$
\begin{minipage}{4cm}
\beginpicture
\setcoordinatesystem units   <1.5mm,2mm>
\setplotarea x from 0 to 16, y from -2 to 15
\put{1.861)}  [l] at 2 12
\put {$ \scriptstyle \bullet$} [c] at 10 0
\put {$ \scriptstyle \bullet$} [c] at 10 12
\put {$ \scriptstyle \bullet$} [c] at 12 6
\put {$ \scriptstyle \bullet$} [c] at 12 0
\put {$ \scriptstyle \bullet$} [c] at 14 0
\put {$ \scriptstyle \bullet$} [c] at 14 12
\setlinear \plot 10 12  10 0 14 12 14 0  /
\setlinear \plot 12 6 12  0  /
\put{$5{,}040$} [c]  at 13 -2
\put{$\scriptstyle \bullet$} [c] at 16  0 \endpicture
\end{minipage}
\begin{minipage}{4cm}
\beginpicture
\setcoordinatesystem units   <1.5mm,2mm>
\setplotarea x from 0 to 16, y from -2 to 15
\put{1.862)}  [l] at 2 12
\put {$ \scriptstyle \bullet$} [c] at 10 0
\put {$ \scriptstyle \bullet$} [c] at 10 12
\put {$ \scriptstyle \bullet$} [c] at 12 0
\put {$ \scriptstyle \bullet$} [c] at 12 6
\put {$ \scriptstyle \bullet$} [c] at 12 12
\put {$ \scriptstyle \bullet$} [c] at 14 12
\setlinear \plot 14 12 12 0 12  12 10 0 10 12 12 0 /
\put{$5{,}040$} [c]  at 13 -2
\put{$\scriptstyle \bullet$} [c] at 16  0 \endpicture
\end{minipage}
\begin{minipage}{4cm}
\beginpicture
\setcoordinatesystem units   <1.5mm,2mm>
\setplotarea x from 0 to 16, y from -2 to 15
\put{1.863)}  [l] at 2 12
\put {$ \scriptstyle \bullet$} [c] at 10 0
\put {$ \scriptstyle \bullet$} [c] at 10 12
\put {$ \scriptstyle \bullet$} [c] at 12 0
\put {$ \scriptstyle \bullet$} [c] at 12 6
\put {$ \scriptstyle \bullet$} [c] at 12 12
\put {$ \scriptstyle \bullet$} [c] at 14 0
\setlinear \plot 14 0 12 12 12  0 10 12 10 0 12 12 /
\put{$5{,}040$} [c]  at 13 -2
\put{$\scriptstyle \bullet$} [c] at 16  0 \endpicture
\end{minipage}
\begin{minipage}{4cm}
\beginpicture
\setcoordinatesystem units   <1.5mm,2mm>
\setplotarea x from 0 to  16, y from -2 to 15
\put{1.864)}  [l] at 2 12
\put {$ \scriptstyle \bullet$} [c] at 10 0
\put {$ \scriptstyle \bullet$} [c] at 10 12
\put {$ \scriptstyle \bullet$} [c] at 12 6
\put {$ \scriptstyle \bullet$} [c] at 14 0
\put {$ \scriptstyle \bullet$} [c] at 14 12
\setlinear \plot  10 12  10 0  14 12 14 0   /
\put {$ \scriptstyle \bullet$} [c] at 16 0
\put {$ \scriptstyle \bullet$} [c] at 16 12
\setlinear \plot 16 0 16 12   /
\put{$5{,}040$} [c] at 13 -2
\endpicture
\end{minipage}
\begin{minipage}{4cm}
\beginpicture
\setcoordinatesystem units   <1.5mm,2mm>
\setplotarea x from 0 to 16, y from -2 to 15
\put{1.865)}  [l] at 2 12
\put {$ \scriptstyle \bullet$} [c] at 10 12
\put {$ \scriptstyle \bullet$} [c] at 12 0
\put {$ \scriptstyle \bullet$} [c] at 12 12
\put {$ \scriptstyle \bullet$} [c] at 14 6
\put {$ \scriptstyle \bullet$} [c] at 14 12
\put {$ \scriptstyle \bullet$} [c] at 14 0
\setlinear \plot 10 12 12 0 14 6 14 12  /
\setlinear \plot 12 0  12 12  /
\setlinear \plot 14 0 14 6 /
\put{$2{,}520$} [c]  at 13 -2
\put{$\scriptstyle \bullet$} [c] at 16  0 \endpicture
\end{minipage}
\begin{minipage}{4cm}
\beginpicture
\setcoordinatesystem units   <1.5mm,2mm>
\setplotarea x from 0 to 16, y from -2 to 15
\put{1.866)}  [l] at 2 12
\put {$ \scriptstyle \bullet$} [c] at 10 0
\put {$ \scriptstyle \bullet$} [c] at 12 0
\put {$ \scriptstyle \bullet$} [c] at 12 12
\put {$ \scriptstyle \bullet$} [c] at 14 6
\put {$ \scriptstyle \bullet$} [c] at 14 12
\put {$ \scriptstyle \bullet$} [c] at 14 0
\setlinear \plot 10 0 12 12 14 6 14 0  /
\setlinear \plot 12 0  12 12  /
\setlinear \plot 14 12 14 6 /
\put{$2{,}520$} [c]  at 13 -2
\put{$\scriptstyle \bullet$} [c] at 16  0 \endpicture
\end{minipage}
$$

$$
\begin{minipage}{4cm}
\beginpicture
\setcoordinatesystem units   <1.5mm,2mm>
\setplotarea x from 0 to  16, y from -2 to 15
\put{1.867)}  [l] at 2 12
\put {$ \scriptstyle \bullet$} [c] at 10 12
\put {$ \scriptstyle \bullet$} [c] at 12 0
\put {$ \scriptstyle \bullet$} [c] at 12 12
\put {$ \scriptstyle \bullet$} [c] at 14 0
\put {$ \scriptstyle \bullet$} [c] at 14 12
\setlinear \plot  10 12   12 0 12 12 14 0 14 12 12 0 /
\put {$ \scriptstyle \bullet$} [c] at 16 0
\put {$ \scriptstyle \bullet$} [c] at 16 12
\setlinear \plot 16 0 16 12   /
\put{$2{,}520$} [c] at 13 -2
\endpicture
\end{minipage}
\begin{minipage}{4cm}
\beginpicture
\setcoordinatesystem units   <1.5mm,2mm>
\setplotarea x from 0 to  16, y from -2 to 15
\put{1.868)}  [l] at 2 12
\put {$ \scriptstyle \bullet$} [c] at 10 0
\put {$ \scriptstyle \bullet$} [c] at 12 0
\put {$ \scriptstyle \bullet$} [c] at 12 12
\put {$ \scriptstyle \bullet$} [c] at 14 0
\put {$ \scriptstyle \bullet$} [c] at 14 12
\setlinear \plot  10 0   12 12 12 0 14 12 14 0 12 12 /
\put {$ \scriptstyle \bullet$} [c] at 16 0
\put {$ \scriptstyle \bullet$} [c] at 16 12
\setlinear \plot 16 0 16 12   /
\put{$2{,}520$} [c] at 13 -2
\endpicture
\end{minipage}
\begin{minipage}{4cm}
\beginpicture
\setcoordinatesystem units   <1.5mm,2mm>
\setplotarea x from 0 to 16, y from -2 to 15
\put{1.869)}  [l] at 2 12
\put {$ \scriptstyle \bullet$} [c] at 10 12
\put {$ \scriptstyle \bullet$} [c] at 10.45 8
\put {$ \scriptstyle \bullet$} [c] at 11 3.9
\put {$ \scriptstyle \bullet$} [c] at 11.5 0
\put {$ \scriptstyle \bullet$} [c] at 13 12
\put {$ \scriptstyle \bullet$} [c] at 14 0
\setlinear \plot 10 12 11.5 0  13 12  /
\put{$2{,}520$} [c]  at 13 -2
\put{$\scriptstyle \bullet$} [c] at 16  0 \endpicture
\end{minipage}
\begin{minipage}{4cm}
\beginpicture
\setcoordinatesystem units   <1.5mm,2mm>
\setplotarea x from 0 to 16, y from -2 to 15
\put{1.870)}  [l] at 2 12
\put {$ \scriptstyle \bullet$} [c] at 10 0
\put {$ \scriptstyle \bullet$} [c] at 10.45 3.9
\put {$ \scriptstyle \bullet$} [c] at 11 8
\put {$ \scriptstyle \bullet$} [c] at 11.5 12
\put {$ \scriptstyle \bullet$} [c] at 13 0
\put {$ \scriptstyle \bullet$} [c] at 14 0
\setlinear \plot 10 0 11.5 12  13 0  /
\put{$2{,}520$} [c]  at 13 -2
\put{$\scriptstyle \bullet$} [c] at 16  0 \endpicture
\end{minipage}
\begin{minipage}{4cm}
\beginpicture
\setcoordinatesystem units   <1.5mm,2mm>
\setplotarea x from 0 to 16, y from -2 to 15
\put{1.871)}  [l] at 2 12
\put {$ \scriptstyle \bullet$} [c] at 10 0
\put {$ \scriptstyle \bullet$} [c] at 10 6
\put {$ \scriptstyle \bullet$} [c] at 11 0
\put {$ \scriptstyle \bullet$} [c] at 11 12
\put {$ \scriptstyle \bullet$} [c] at 12 6
\put {$ \scriptstyle \bullet$} [c] at 14 0
\setlinear \plot 10 0 10 6 11 12 12 6 11  0 10 6  /
\put{$2{,}520$} [c]  at 13 -2
\put{$\scriptstyle \bullet$} [c] at 16  0 \endpicture
\end{minipage}
\begin{minipage}{4cm}
\beginpicture
\setcoordinatesystem units   <1.5mm,2mm>
\setplotarea x from 0 to 16, y from -2 to 15
\put{1.872)}  [l] at 2 12
\put {$ \scriptstyle \bullet$} [c] at 10 6
\put {$ \scriptstyle \bullet$} [c] at 10 12
\put {$ \scriptstyle \bullet$} [c] at 11 0
\put {$ \scriptstyle \bullet$} [c] at 11 12
\put {$ \scriptstyle \bullet$} [c] at 12 6
\put {$ \scriptstyle \bullet$} [c] at 14 0
\setlinear \plot 10 12 10 6 11 12 12 6 11  0 10 6  /
\put{$2{,}520$} [c]  at 13 -2
\put{$\scriptstyle \bullet$} [c] at 16  0
\endpicture
\end{minipage}
$$
$$
\begin{minipage}{4cm}
\beginpicture
\setcoordinatesystem units   <1.5mm,2mm>
\setplotarea x from 0 to 16, y from -2 to 15
\put{1.873)}  [l] at 2 12
\put {$ \scriptstyle \bullet$} [c] at 10 12
\put {$ \scriptstyle \bullet$} [c] at 11 0
\put {$ \scriptstyle \bullet$} [c] at 11 6
\put {$ \scriptstyle \bullet$} [c] at 12 12
\put {$ \scriptstyle \bullet$} [c] at 14 12
\put {$ \scriptstyle \bullet$} [c] at 14 0
\setlinear \plot 10 12 11 6 11 0   /
\setlinear \plot 12 12 11 6    /
\setlinear \plot 14  0 14 12   /
\put{$2{,}520$} [c]  at 13 -2
\put{$\scriptstyle \bullet$} [c] at 16  0
\endpicture
\end{minipage}
\begin{minipage}{4cm}
\beginpicture
\setcoordinatesystem units   <1.5mm,2mm>
\setplotarea x from 0 to 16, y from -2 to 15
\put{1.874)}  [l] at 2 12
\put {$ \scriptstyle \bullet$} [c] at 10 0
\put {$ \scriptstyle \bullet$} [c] at 11 12
\put {$ \scriptstyle \bullet$} [c] at 11 6
\put {$ \scriptstyle \bullet$} [c] at 12 0
\put {$ \scriptstyle \bullet$} [c] at 14 12
\put {$ \scriptstyle \bullet$} [c] at 14 0
\setlinear \plot 10 0 11 6 11 12   /
\setlinear \plot 12 0 11 6    /
\setlinear \plot 14  0 14 12   /
\put{$2{,}520$} [c]  at 13 -2
\put{$\scriptstyle \bullet$} [c] at 16  0
\endpicture
\end{minipage}
\begin{minipage}{4cm}
\beginpicture
\setcoordinatesystem units   <1.5mm,2mm>
\setplotarea x from 0 to 16, y from -2 to 15
\put{1.875)}  [l] at 2 12
\put {$ \scriptstyle \bullet$} [c] at 10 6
\put {$ \scriptstyle \bullet$} [c] at 11 0
\put {$ \scriptstyle \bullet$} [c] at 11 12
\put {$ \scriptstyle \bullet$} [c] at 12 6
\put {$ \scriptstyle \bullet$} [c] at 14 12
\put {$ \scriptstyle \bullet$} [c] at 14 0
\setlinear \plot  11 0 10  6  11  12 12 6 11 0   /
\setlinear \plot  14  0  14 12    /
\put{$2{,}520$} [c]  at 13 -2
\put{$\scriptstyle \bullet$} [c] at 16  0
\endpicture
\end{minipage}
\begin{minipage}{4cm}
\beginpicture
\setcoordinatesystem units   <1.5mm,2mm>
\setplotarea x from 0 to 16, y from -2 to 15
\put{1.876)}  [l] at 2 12
\put {$ \scriptstyle \bullet$} [c] at 10 0
\put {$ \scriptstyle \bullet$} [c] at 10 6
\put {$ \scriptstyle \bullet$} [c] at 10 12
\put {$ \scriptstyle \bullet$} [c] at 13 0
\put {$ \scriptstyle \bullet$} [c] at 13 12
\put {$ \scriptstyle \bullet$} [c] at 16 0
\put{$\scriptstyle \bullet$} [c] at 16  12
\setlinear \plot  10 0 10  12    /
\setlinear \plot  13  0  13 12    /
\setlinear \plot  16  0  16 12    /
\put{$2{,}520$} [c]  at 13 -2
\endpicture
\end{minipage}
\begin{minipage}{4cm}
\beginpicture
\setcoordinatesystem units   <1.5mm,2mm>
\setplotarea x from 0 to 16, y from -2 to 15
\put{1.877)}  [l] at 2 12
\put {$ \scriptstyle \bullet$} [c] at 10 0
\put {$ \scriptstyle \bullet$} [c] at 10 12
\put {$ \scriptstyle \bullet$} [c] at 12 0
\put {$ \scriptstyle \bullet$} [c] at 12 12
\put {$ \scriptstyle \bullet$} [c] at 14 0
\put {$ \scriptstyle \bullet$} [c] at 14 12
\setlinear \plot 14 0 12 12 12 0  10 12 10  0 12 12  /
\setlinear \plot 10 12 14  0  14 12   /
\put{$1{,}260$} [c]  at 13 -2
\put{$\scriptstyle \bullet$} [c] at 16  0 \endpicture
\end{minipage}
\begin{minipage}{4cm}
\beginpicture
\setcoordinatesystem units   <1.5mm,2mm>
\setplotarea x from 0 to 16, y from -2 to 15
\put{1.878)}  [l] at 2 12
\put {$ \scriptstyle \bullet$} [c] at 10 0
\put {$ \scriptstyle \bullet$} [c] at 10 12
\put {$ \scriptstyle \bullet$} [c] at 12 0
\put {$ \scriptstyle \bullet$} [c] at 12 12
\put {$ \scriptstyle \bullet$} [c] at 14 0
\put {$ \scriptstyle \bullet$} [c] at 14 12
\setlinear \plot 14 12 12 0 12 12  10 0 10  12 12 0  /
\setlinear \plot 10 0 14  12  14 0   /
\put{$1{,}260$} [c]  at 13 -2
\put{$\scriptstyle \bullet$} [c] at 16  0 \endpicture
\end{minipage}
$$
$$
\begin{minipage}{4cm}
\beginpicture
\setcoordinatesystem units   <1.5mm,2mm>
\setplotarea x from 0 to 16, y from -2 to 15
\put{1.879)}  [l] at 2 12
\put {$ \scriptstyle \bullet$} [c] at 10 0
\put {$ \scriptstyle \bullet$} [c] at 10 6
\put {$ \scriptstyle \bullet$} [c] at 10  12
\put {$ \scriptstyle \bullet$} [c] at 12 0
\put {$ \scriptstyle \bullet$} [c] at 12 12
\put {$ \scriptstyle \bullet$} [c] at 14 0
\setlinear \plot  10 0  10 12 12 0 12 12 10 6   /
\put{$1{,}260$} [c]  at 13 -2
\put{$\scriptstyle \bullet$} [c] at 16  0
\endpicture
\end{minipage}
\begin{minipage}{4cm}
\beginpicture
\setcoordinatesystem units   <1.5mm,2mm>
\setplotarea x from 0 to 16, y from -2 to 15
\put{1.880)}  [l] at 2 12
\put {$ \scriptstyle \bullet$} [c] at 10 0
\put {$ \scriptstyle \bullet$} [c] at 10 6
\put {$ \scriptstyle \bullet$} [c] at 10  12
\put {$ \scriptstyle \bullet$} [c] at 12 0
\put {$ \scriptstyle \bullet$} [c] at 12 12
\put {$ \scriptstyle \bullet$} [c] at 14 0
\setlinear \plot  10 12  10 0 12 12 12 0 10 6   /
\put{$1{,}260$} [c]  at 13 -2
\put{$\scriptstyle \bullet$} [c] at 16  0
\endpicture
\end{minipage}
\begin{minipage}{4cm}
\beginpicture
\setcoordinatesystem units   <1.5mm,2mm>
\setplotarea x from 0 to  16, y from -2 to 15
\put{1.881)}  [l] at 2 12
\put {$ \scriptstyle \bullet$} [c] at 10 12
\put {$ \scriptstyle \bullet$} [c] at 12 12
\put {$ \scriptstyle \bullet$} [c] at 14 12
\put {$ \scriptstyle \bullet$} [c] at 12 0
\put {$ \scriptstyle \bullet$} [c] at 16 0
\put {$ \scriptstyle \bullet$} [c] at 16 6
\put {$ \scriptstyle \bullet$} [c] at 16 12
\setlinear \plot  10 12   12 0  14 12  /
\setlinear \plot 12 0 12 12   /
\setlinear \plot 16 0 16 12   /
\put{$840$} [c] at 13 -2
\endpicture
\end{minipage}
\begin{minipage}{4cm}
\beginpicture
\setcoordinatesystem units   <1.5mm,2mm>
\setplotarea x from 0 to  16, y from -2 to 15
\put{1.882)}  [l] at 2 12
\put {$ \scriptstyle \bullet$} [c] at 10 0
\put {$ \scriptstyle \bullet$} [c] at 12 0
\put {$ \scriptstyle \bullet$} [c] at 14 0
\put {$ \scriptstyle \bullet$} [c] at 12 12
\put {$ \scriptstyle \bullet$} [c] at 16 0
\put {$ \scriptstyle \bullet$} [c] at 16 6
\put {$ \scriptstyle \bullet$} [c] at 16 12
\setlinear \plot  10 0   12 12  14 0  /
\setlinear \plot 12 0 12 12   /
\setlinear \plot 16 0 16 12   /
\put{$840$} [c] at 13 -2
\endpicture
\end{minipage}
\begin{minipage}{4cm}
\beginpicture
\setcoordinatesystem units   <1.5mm,2mm>
\setplotarea x from 0 to 16, y from -2 to 15
\put{1.883)}  [l] at 2 12
\put {$ \scriptstyle \bullet$} [c] at 10 0
\put {$ \scriptstyle \bullet$} [c] at 10 12
\put {$ \scriptstyle \bullet$} [c] at 12 0
\put {$ \scriptstyle \bullet$} [c] at 12 12
\put {$ \scriptstyle \bullet$} [c] at 14 0
\put {$ \scriptstyle \bullet$} [c] at 14 12
\setlinear \plot 10 0 10 12  12 0 14 12 14 0 12 12 10 0 /
\put{$840$} [c]  at 13 -2
\put{$\scriptstyle \bullet$} [c] at 16  0 \endpicture
\end{minipage}
\begin{minipage}{4cm}
\beginpicture
\setcoordinatesystem units   <1.5mm,2mm>
\setplotarea x from 0 to 16, y from -2 to 15
\put{1.884)}  [l] at 2 12
\put {$ \scriptstyle \bullet$} [c] at 10 12
\put {$ \scriptstyle \bullet$} [c] at 11.5 12
\put {$ \scriptstyle \bullet$} [c] at 12.5 12
\put {$ \scriptstyle \bullet$} [c] at 14 12
\put {$ \scriptstyle \bullet$} [c] at 12 6
\put {$ \scriptstyle \bullet$} [c] at 12 0
\setlinear \plot 10 12 12 6 14 12    /
\setlinear \plot 11.5 12 12 6 12.5 12    /
\setlinear \plot 12 0 12 6    /
\put{$210$} [c]  at 13 -2
\put{$\scriptstyle \bullet$} [c] at 16  0 \endpicture
\end{minipage}
$$
$$
\begin{minipage}{4cm}
\beginpicture
\setcoordinatesystem units   <1.5mm,2mm>
\setplotarea x from 0 to 16, y from -2 to 15
\put{1.885)}  [l] at 2 12
\put {$ \scriptstyle \bullet$} [c] at 10 0
\put {$ \scriptstyle \bullet$} [c] at 11.5 0
\put {$ \scriptstyle \bullet$} [c] at 12.5 0
\put {$ \scriptstyle \bullet$} [c] at 14 0
\put {$ \scriptstyle \bullet$} [c] at 12 6
\put {$ \scriptstyle \bullet$} [c] at 12 12
\setlinear \plot 10 0 12 6 14 0     /
\setlinear \plot 11.5 0 12 6 12.5 0    /
\setlinear \plot 12 12 12 6    /
\put{$210$} [c]  at 13 -2
\put{$\scriptstyle \bullet$} [c] at 16  0 \endpicture
\end{minipage}
\begin{minipage}{4cm}
\beginpicture
\setcoordinatesystem units   <1.5mm,2mm>
\setplotarea x from 0 to 16, y from -2 to 15
\put{1.886)}  [l] at 2 12
\put {$ \scriptstyle \bullet$} [c] at 10 6
\put {$ \scriptstyle \bullet$} [c] at 11.5 6
\put {$ \scriptstyle \bullet$} [c] at 13.5 6
\put {$ \scriptstyle \bullet$} [c] at 15 6
\put {$ \scriptstyle \bullet$} [c] at 12.5 0
\put {$ \scriptstyle \bullet$} [c] at 12.5 12
\setlinear \plot  12.5 0 10  6  12.5 12 15 6 12.5  0   /
\setlinear \plot  12.5 0 11.5  6  12.5 12 13.5 6 12.5  0   /
\put{$210$} [c]  at 13 -2
\put{$\scriptstyle \bullet$} [c] at 16  0 \endpicture
\end{minipage}
\begin{minipage}{4cm}
\beginpicture
\setcoordinatesystem units <1.5mm,2mm>
\setplotarea x from 0 to 16, y from -2 to 15
\put {${\bf  37}$} [l] at 2 15

\put {1.887)} [l] at 2 12
\put {$ \scriptstyle \bullet$} [c] at  10 0
\put {$ \scriptstyle \bullet$} [c] at  10 12
\put {$ \scriptstyle \bullet$} [c] at  12 12
\put {$ \scriptstyle \bullet$} [c] at  14 0
\put {$ \scriptstyle \bullet$} [c] at  14 12
\put {$ \scriptstyle \bullet$} [c] at  16  0
\put {$ \scriptstyle \bullet$} [c] at  16  12
\setlinear \plot   10 12 10  0  14 12 14 0 16  12 16 0   /
\setlinear \plot   12 12 10 0    /
\put{$2{,}520$} [c] at 13 -2
\endpicture
\end{minipage}
\begin{minipage}{4cm}
\beginpicture
\setcoordinatesystem units <1.5mm,2mm>
\setplotarea x from 0 to 16, y from -2 to 15
\put {1.888)} [l] at 2 12
\put {$ \scriptstyle \bullet$} [c] at  10 0
\put {$ \scriptstyle \bullet$} [c] at  10 12
\put {$ \scriptstyle \bullet$} [c] at  12 0
\put {$ \scriptstyle \bullet$} [c] at  14 0
\put {$ \scriptstyle \bullet$} [c] at  14 12
\put {$ \scriptstyle \bullet$} [c] at  16  0
\put {$ \scriptstyle \bullet$} [c] at  16  12
\setlinear \plot   10 0 10  12  14 0 14 12 16  0 16 12   /
\setlinear \plot   12 0 10 12    /
\put{$2{,}520$} [c] at 13 -2
\endpicture
\end{minipage}
\begin{minipage}{4cm}
\beginpicture
\setcoordinatesystem units <1.5mm,2mm>
\setplotarea x from 0 to 16, y from -2 to 15
\put{1.889)} [l] at 2 12
\put {$ \scriptstyle \bullet$} [c] at 10 6
\put {$ \scriptstyle \bullet$} [c] at 10 12
\put {$ \scriptstyle \bullet$} [c] at 12 6
\put {$ \scriptstyle \bullet$} [c] at 12 12
\put {$ \scriptstyle \bullet$} [c] at 14 12
\put {$ \scriptstyle \bullet$} [c] at 16 12
\put {$ \scriptstyle \bullet$} [c] at 12 0
\setlinear \plot 10 12 10 6 12  0 12 12    /
\setlinear \plot 14 12 12 0 16 12     /
\put{$1{,}260$} [c] at 13 -2
\endpicture
\end{minipage}
\begin{minipage}{4cm}
\beginpicture
\setcoordinatesystem units <1.5mm,2mm>
\setplotarea x from 0 to 16, y from -2 to 15
\put{1.890)} [l] at 2 12
\put {$ \scriptstyle \bullet$} [c] at 10 6
\put {$ \scriptstyle \bullet$} [c] at 10 0
\put {$ \scriptstyle \bullet$} [c] at 12 6
\put {$ \scriptstyle \bullet$} [c] at 12 0
\put {$ \scriptstyle \bullet$} [c] at 14 0
\put {$ \scriptstyle \bullet$} [c] at 16 0
\put {$ \scriptstyle \bullet$} [c] at 12 12
\setlinear \plot 10 0 10 6 12  12 12 0    /
\setlinear \plot 14 0 12 12 16 0     /
\put{$1{,}260$} [c] at 13 -2
\endpicture
\end{minipage}
$$

$$
\begin{minipage}{4cm}
\beginpicture
\setcoordinatesystem units <1.5mm,2mm>
\setplotarea x from 0 to 16, y from -2 to 15
\put{1.891)} [l] at 2 12
\put {$ \scriptstyle \bullet$} [c] at 10 6
\put {$ \scriptstyle \bullet$} [c] at 11 6
\put {$ \scriptstyle \bullet$} [c] at 11 0
\put {$ \scriptstyle \bullet$} [c] at 11 12
\put {$ \scriptstyle \bullet$} [c] at 12 6
\put {$ \scriptstyle \bullet$} [c] at 14.5 12
\put {$ \scriptstyle \bullet$} [c] at 16 12
\setlinear \plot 10 6 11 12 11 0 10 6     /
\setlinear \plot 11 12 12 6 11 0 14.5 12      /
\setlinear \plot 16 12 11 0       /
\put{$420$} [c] at 13 -2
\endpicture
\end{minipage}
\begin{minipage}{4cm}
\beginpicture
\setcoordinatesystem units <1.5mm,2mm>
\setplotarea x from 0 to 16, y from -2 to 15
\put{1.892)} [l] at 2 12
\put {$ \scriptstyle \bullet$} [c] at 10 6
\put {$ \scriptstyle \bullet$} [c] at 11 6
\put {$ \scriptstyle \bullet$} [c] at 11 0
\put {$ \scriptstyle \bullet$} [c] at 11 12
\put {$ \scriptstyle \bullet$} [c] at 12 6
\put {$ \scriptstyle \bullet$} [c] at 14.5 0
\put {$ \scriptstyle \bullet$} [c] at 16 0
\setlinear \plot 10 6 11 12 11 0 10 6     /
\setlinear \plot 11 0 12 6 11 12 14.5 0      /
\setlinear \plot 16 0 11 12       /
\put{$420$} [c] at 13 -2
\endpicture
\end{minipage}
\begin{minipage}{4cm}
\beginpicture
\setcoordinatesystem units <1.5mm,2mm>
\setplotarea x from 0 to 16, y from -2 to 15
\put{1.893)} [l] at 2 12
\put {$ \scriptstyle \bullet$} [c] at 10 12
\put {$ \scriptstyle \bullet$} [c] at 11.5 12
\put {$ \scriptstyle \bullet$} [c] at 11.5 6
\put {$ \scriptstyle \bullet$} [c] at 13 12
\put {$ \scriptstyle \bullet$} [c] at 11.5 0
\put {$ \scriptstyle \bullet$} [c] at 14.5 12
\put {$ \scriptstyle \bullet$} [c] at 16  12
\setlinear \plot 10 12 11.5 6 13 12     /
\setlinear \plot  11.5 12 11.5 0 14.5 12     /
\setlinear \plot  16 12 11.5 0 /
\put{$420$} [c] at 13 -2
\endpicture
\end{minipage}
\begin{minipage}{4cm}
\beginpicture
\setcoordinatesystem units <1.5mm,2mm>
\setplotarea x from 0 to 16, y from -2 to 15
\put{1.894)} [l] at 2 12
\put {$ \scriptstyle \bullet$} [c] at 10 0
\put {$ \scriptstyle \bullet$} [c] at 11.5 12
\put {$ \scriptstyle \bullet$} [c] at 11.5 6
\put {$ \scriptstyle \bullet$} [c] at 13 0
\put {$ \scriptstyle \bullet$} [c] at 11.5 0
\put {$ \scriptstyle \bullet$} [c] at 14.5 0
\put {$ \scriptstyle \bullet$} [c] at 16  0
\setlinear \plot 10 0 11.5 6 13 0     /
\setlinear \plot  11.5 0 11.5 12 14.5 0     /
\setlinear \plot  16 0 11.5 12 /
\put{$420$} [c] at 13 -2
\endpicture
\end{minipage}
\begin{minipage}{4cm}
\beginpicture
\setcoordinatesystem units <1.5mm,2mm>
\setplotarea x  from 0 to 16, y from -2 to 15
\put {1.895)} [l] at  2 12
\put {$ \scriptstyle \bullet$} [c] at  10 12
\put {$ \scriptstyle \bullet$} [c] at  11.5 0
\put {$ \scriptstyle \bullet$} [c] at  13 12
\put {$ \scriptstyle \bullet$} [c] at  14.5 0
\put {$ \scriptstyle \bullet$} [c] at  16  12
\put {$ \scriptstyle \bullet$} [c] at  11.5  12
\put {$ \scriptstyle \bullet$} [c] at  14.5  12
\setlinear \plot  10 12 11.5 0 14.5 12 14.5 0  11.5 12 11.5  0  /
\setlinear  \plot  11.5 0 13 12 14.5 0 16 12 /
\put{$420$} [c] at 13 -2
\endpicture
\end{minipage}
\begin{minipage}{4cm}
\beginpicture
\setcoordinatesystem units <1.5mm,2mm>
\setplotarea x  from 0 to 16, y from -2 to 15
\put {1.896)} [l] at  2 12
\put {$ \scriptstyle \bullet$} [c] at  10 0
\put {$ \scriptstyle \bullet$} [c] at  11.5 0
\put {$ \scriptstyle \bullet$} [c] at  13 0
\put {$ \scriptstyle \bullet$} [c] at  14.5 0
\put {$ \scriptstyle \bullet$} [c] at  16  0
\put {$ \scriptstyle \bullet$} [c] at  11.5  12
\put {$ \scriptstyle \bullet$} [c] at  14.5  12
\setlinear \plot  10 0 11.5 12 14.5 0 14.5 12  11.5 0 11.5  12  /
\setlinear  \plot  11.5 12 13 0 14.5 12 16 0 /
\put{$420$} [c] at 13 -2
\endpicture
\end{minipage}
$$
$$
\begin{minipage}{4cm}
\beginpicture
\setcoordinatesystem units <1.5mm,2mm>
\setplotarea x from 0 to 16, y from -2 to 15
\put {${\bf  38}$} [l] at 2 15

\put {1.897)} [l] at 2 12
\put {$ \scriptstyle \bullet$} [c] at  10 0
\put {$ \scriptstyle \bullet$} [c] at  10 12
\put {$ \scriptstyle \bullet$} [c] at  12 0
\put {$ \scriptstyle \bullet$} [c] at  12 12
\put {$ \scriptstyle \bullet$} [c] at  14 12
\put {$ \scriptstyle \bullet$} [c] at  16  0
\put {$ \scriptstyle \bullet$} [c] at  16 12
\setlinear \plot   10 12 10 0  12 12 16 0 16  12    /
\setlinear \plot   12 0 12 12     /
\setlinear \plot   14 12 16 0    /
\put{$2{,}520$} [c] at 13 -2
\endpicture
\end{minipage}
\begin{minipage}{4cm}
\beginpicture
\setcoordinatesystem units <1.5mm,2mm>
\setplotarea x from 0 to 16, y from -2 to 15
\put {1.898)} [l] at 2 12
\put {$ \scriptstyle \bullet$} [c] at  10 0
\put {$ \scriptstyle \bullet$} [c] at  10 12
\put {$ \scriptstyle \bullet$} [c] at  12 0
\put {$ \scriptstyle \bullet$} [c] at  12 12
\put {$ \scriptstyle \bullet$} [c] at  14 0
\put {$ \scriptstyle \bullet$} [c] at  16  0
\put {$ \scriptstyle \bullet$} [c] at  16 12
\setlinear \plot   10 0 10 12  12 0 16 12 16  0    /
\setlinear \plot   12 0 12 12     /
\setlinear \plot   14 0 16 12    /
\put{$2{,}520$} [c] at 13 -2
\endpicture
\end{minipage}
\begin{minipage}{4cm}
\beginpicture
\setcoordinatesystem units <1.5mm,2mm>
\setplotarea x from 0 to 16, y from -2 to 15
\put {1.899)} [l] at 2 12
\put {$ \scriptstyle \bullet$} [c] at  10 0
\put {$ \scriptstyle \bullet$} [c] at  10 12
\put {$ \scriptstyle \bullet$} [c] at  12 12
\put {$ \scriptstyle \bullet$} [c] at  14 0
\put {$ \scriptstyle \bullet$} [c] at  14 6
\put {$ \scriptstyle \bullet$} [c] at  14 12
\put {$ \scriptstyle \bullet$} [c] at  16  12
\setlinear \plot  10 0 10 12 14 0 16 12 /
\setlinear \plot  12 12 14 0  /
\setlinear \plot 14 0 14 12 /
\put{$2{,}520$} [c] at 13  -2
\endpicture
\end{minipage}
\begin{minipage}{4cm}
\beginpicture
\setcoordinatesystem units <1.5mm,2mm>
\setplotarea x from 0 to 16, y from -2 to 15
\put {1.900)} [l] at 2 12
\put {$ \scriptstyle \bullet$} [c] at  10 0
\put {$ \scriptstyle \bullet$} [c] at  10 12
\put {$ \scriptstyle \bullet$} [c] at  12 0
\put {$ \scriptstyle \bullet$} [c] at  14 0
\put {$ \scriptstyle \bullet$} [c] at  14 6
\put {$ \scriptstyle \bullet$} [c] at  14 12
\put {$ \scriptstyle \bullet$} [c] at  16  0
\setlinear \plot  10 12 10 0 14 12 16 0 /
\setlinear \plot  12 0 14 12  /
\setlinear \plot 14 0 14 12 /
\put{$2{,}520$} [c] at 13  -2
\endpicture
\end{minipage}
\begin{minipage}{4cm}
\beginpicture
\setcoordinatesystem units <1.5mm,2mm>
\setplotarea x from 0 to 16, y from -2 to 15
\put{1.901)} [l] at 2 12
\put {$ \scriptstyle \bullet$} [c] at  10 12
\put {$ \scriptstyle \bullet$} [c] at  10.9 12
\put {$ \scriptstyle \bullet$} [c] at  12.5 6
\put {$ \scriptstyle \bullet$} [c] at  13 0
\put {$ \scriptstyle \bullet$} [c] at  13 12
\put {$ \scriptstyle \bullet$} [c] at  13.5 6
\put {$ \scriptstyle \bullet$} [c] at  16 0
\setlinear \plot  16 0 13  12  12.5 6  13 0 13.5 6 13 12 /
\setlinear \plot 10  12  13 0 10.9 12  /
\put{$1{,}260$} [c] at 13 -2
\endpicture
\end{minipage}
\begin{minipage}{4cm}
\beginpicture
\setcoordinatesystem units <1.5mm,2mm>
\setplotarea x from 0 to 16, y from -2 to 15
\put{1.902)} [l] at 2 12
\put {$ \scriptstyle \bullet$} [c] at  10 0
\put {$ \scriptstyle \bullet$} [c] at  10.9 0
\put {$ \scriptstyle \bullet$} [c] at  12.5 6
\put {$ \scriptstyle \bullet$} [c] at  13 0
\put {$ \scriptstyle \bullet$} [c] at  13 12
\put {$ \scriptstyle \bullet$} [c] at  13.5 6
\put {$ \scriptstyle \bullet$} [c] at  16 12
\setlinear \plot  16 12 13  0  12.5 6  13 12 13.5 6 13 0 /
\setlinear \plot 10  0  13 12 10.9 0  /
\put{$1{,}260$} [c] at 13 -2
\endpicture
\end{minipage}
$$
$$
\begin{minipage}{4cm}
\beginpicture
\setcoordinatesystem units <1.5mm,2mm>
\setplotarea x from 0 to 16, y from -2 to 15
\put {1.903)} [l] at 2 12
\put {$ \scriptstyle \bullet$} [c] at  10 12
\put {$ \scriptstyle \bullet$} [c] at  11 12
\put {$ \scriptstyle \bullet$} [c] at  12 0
\put {$ \scriptstyle \bullet$} [c] at  14 0
\put {$ \scriptstyle \bullet$} [c] at  16  12
\put {$ \scriptstyle \bullet$} [c] at  12 12
\put {$ \scriptstyle \bullet$} [c] at  14 12
\setlinear \plot  10 12 12 0 11 12  /
\setlinear \plot  12 12 12 0 14 12 14 0 12 12  /
\setlinear  \plot  12 0 16 12 14 0  /
\put{$420$} [c] at 13 -2
\endpicture
\end{minipage}
\begin{minipage}{4cm}
\beginpicture
\setcoordinatesystem units <1.5mm,2mm>
\setplotarea x from 0 to 16, y from -2 to 15
\put {1.904)} [l] at 2 12
\put {$ \scriptstyle \bullet$} [c] at  10 0
\put {$ \scriptstyle \bullet$} [c] at  11 0
\put {$ \scriptstyle \bullet$} [c] at  12 0
\put {$ \scriptstyle \bullet$} [c] at  14 0
\put {$ \scriptstyle \bullet$} [c] at  16  0
\put {$ \scriptstyle \bullet$} [c] at  12 12
\put {$ \scriptstyle \bullet$} [c] at  14 12
\setlinear \plot  10 0 12 12 11 0  /
\setlinear \plot  12 0 12 12 14 0 14 12 12 0  /
\setlinear  \plot  12 12 16 0 14 12  /
\put{$420$} [c] at 13 -2
\endpicture
\end{minipage}
\begin{minipage}{4cm}
\beginpicture
\setcoordinatesystem units <1.5mm,2mm>
\setplotarea x from 0 to 16, y from -2 to 15
\put{1.905)} [l] at 2 12
\put {$ \scriptstyle \bullet$} [c] at 10 0
\put {$ \scriptstyle \bullet$} [c] at 10 12
\put {$ \scriptstyle \bullet$} [c] at 12 12
\put {$ \scriptstyle \bullet$} [c] at 13 6
\put {$ \scriptstyle \bullet$} [c] at 14 0
\put {$ \scriptstyle \bullet$} [c] at 14 12
\put{$\scriptstyle \bullet$} [c] at 16  0
\setlinear \plot 10 12  10 0 12 12 14 0  14 12 /
\put{$5{,}040$} [c] at 13 -2
\endpicture
\end{minipage}
\begin{minipage}{4cm}
\beginpicture
\setcoordinatesystem units <1.5mm,2mm>
\setplotarea x from 0 to 16, y from -2 to 15
\put{1.906)} [l] at 2 12
\put {$ \scriptstyle \bullet$} [c] at 10 0
\put {$ \scriptstyle \bullet$} [c] at 10 12
\put {$ \scriptstyle \bullet$} [c] at 12 0
\put {$ \scriptstyle \bullet$} [c] at 13 6
\put {$ \scriptstyle \bullet$} [c] at 14 0
\put {$ \scriptstyle \bullet$} [c] at 14 12
\put{$\scriptstyle \bullet$} [c] at 16  0
\setlinear \plot 10 0  10 12 12 0 14 12  14 0 /
\put{$5{,}040$} [c] at 13 -2
\endpicture
\end{minipage}
\begin{minipage}{4cm}
\beginpicture
\setcoordinatesystem units <1.5mm,2mm>
\setplotarea x from 0 to 16, y from -2 to 15
\put{1.907)} [l] at 2 12
\put {$ \scriptstyle \bullet$} [c] at 10 6
\put {$ \scriptstyle \bullet$} [c] at 10 12
\put {$ \scriptstyle \bullet$} [c] at  12 12
\put {$ \scriptstyle \bullet$} [c] at 12 0
\put {$ \scriptstyle \bullet$} [c] at 14 6
\put {$ \scriptstyle \bullet$} [c] at 14 12
\put{$\scriptstyle \bullet$} [c] at 16  0
\setlinear \plot 10 12 10 6 12 0 14 6 14 12  /
\setlinear \plot  12 0 12 12  /
\put{$2{,}520$} [c] at 13 -2
 \endpicture
\end{minipage}
\begin{minipage}{4cm}
\beginpicture
\setcoordinatesystem units <1.5mm,2mm>
\setplotarea x from 0 to 16, y from -2 to 15
\put{1.908)} [l] at 2 12
\put {$ \scriptstyle \bullet$} [c] at 10 6
\put {$ \scriptstyle \bullet$} [c] at 10 0
\put {$ \scriptstyle \bullet$} [c] at  12 12
\put {$ \scriptstyle \bullet$} [c] at 12 0
\put {$ \scriptstyle \bullet$} [c] at 14 6
\put {$ \scriptstyle \bullet$} [c] at 14 0
\put{$\scriptstyle \bullet$} [c] at 16  0
\setlinear \plot 10 0 10 6 12 12 14 6 14 0  /
\setlinear \plot  12 0 12 12  /
\put{$2{,}520$} [c] at 13 -2
\endpicture
\end{minipage}
$$

$$
\begin{minipage}{4cm}
\beginpicture
\setcoordinatesystem units <1.5mm,2mm>
\setplotarea x from 0 to 16, y from -2 to 15
\put{1.909)} [l] at 2 12
\put {$ \scriptstyle \bullet$} [c] at 10 0
\put {$ \scriptstyle \bullet$} [c] at 10 12
\put {$ \scriptstyle \bullet$} [c] at 12 12
\put {$ \scriptstyle \bullet$} [c] at 14 0
\put {$ \scriptstyle \bullet$} [c] at 14 6
\put {$ \scriptstyle \bullet$} [c] at 14 12
\setlinear \plot 10 12 10 0 14 12 14 0 /
\setlinear \plot  10 0 12 12  /
\put{$2{,}520$} [c] at 13 -2
\put{$\scriptstyle \bullet$} [c] at 16  0 \endpicture
\end{minipage}
\begin{minipage}{4cm}
\beginpicture
\setcoordinatesystem units <1.5mm,2mm>
\setplotarea x from 0 to 16, y from -2 to 15
\put{1.910)} [l] at 2 12
\put {$ \scriptstyle \bullet$} [c] at 10 0
\put {$ \scriptstyle \bullet$} [c] at 10 12
\put {$ \scriptstyle \bullet$} [c] at 12 0
\put {$ \scriptstyle \bullet$} [c] at 14 0
\put {$ \scriptstyle \bullet$} [c] at 14 6
\put {$ \scriptstyle \bullet$} [c] at 14 12
\setlinear \plot 10 0 10 12 14 0 14 12 /
\setlinear \plot  10 12 12 0  /
\put{$2{,}520$} [c] at 13 -2
\put{$\scriptstyle \bullet$} [c] at 16  0 \endpicture
\end{minipage}
\begin{minipage}{4cm}
\beginpicture
\setcoordinatesystem units <1.5mm,2mm>
\setplotarea x from 0 to 16, y from -2 to 15
\put{1.911)} [l] at 2 12
\put {$ \scriptstyle \bullet$} [c] at 10 0
\put {$ \scriptstyle \bullet$} [c] at 10 12
\put {$ \scriptstyle \bullet$} [c] at 12 0
\put {$ \scriptstyle \bullet$} [c] at 12 12
\put {$ \scriptstyle \bullet$} [c] at 14 0
\put {$ \scriptstyle \bullet$} [c] at 14 12
\setlinear \plot 14  0 14  12 12 0  12 12 10 0 10 12 12 0 /
\put{$2{,}520$} [c] at 13 -2
\put{$\scriptstyle \bullet$} [c] at 16  0 \endpicture
\end{minipage}
\begin{minipage}{4cm}
\beginpicture
\setcoordinatesystem units <1.5mm,2mm>
\setplotarea x from 0 to 16, y from -2 to 15
\put{1.912)} [l] at 2 12
\put {$ \scriptstyle \bullet$} [c] at 10 0
\put {$ \scriptstyle \bullet$} [c] at 10 12
\put {$ \scriptstyle \bullet$} [c] at 12 0
\put {$ \scriptstyle \bullet$} [c] at 12 12
\put {$ \scriptstyle \bullet$} [c] at 14 0
\put {$ \scriptstyle \bullet$} [c] at 14 12
\setlinear \plot 14  12 14  0 12 12  12 0 10 12 10 0 12 12 /
\put{$2{,}520$} [c] at 13 -2
\put{$\scriptstyle \bullet$} [c] at 16  0 \endpicture
\end{minipage}
\begin{minipage}{4cm}
\beginpicture
\setcoordinatesystem units <1.5mm,2mm>
\setplotarea x from 0 to 16, y from -2 to 15
\put{1.913)} [l] at 2 12
\put {$ \scriptstyle \bullet$} [c] at 10 6
\put {$ \scriptstyle \bullet$} [c] at 11 0
\put {$ \scriptstyle \bullet$} [c] at  11 6
\put {$ \scriptstyle \bullet$} [c] at 11 12
\put {$ \scriptstyle \bullet$} [c] at 12 6
\put {$ \scriptstyle \bullet$} [c] at 15 12
\setlinear \plot 15 12 11 0 10 6 11 12 12 6 11 0 11 12 /
\put{$840$} [c] at 13 -2
\put{$\scriptstyle \bullet$} [c] at 16  0 \endpicture
\end{minipage}
\begin{minipage}{4cm}
\beginpicture
\setcoordinatesystem units <1.5mm,2mm>
\setplotarea x from 0 to 16, y from -2 to 15
\put{1.914)} [l] at 2 12
\put {$ \scriptstyle \bullet$} [c] at 10 6
\put {$ \scriptstyle \bullet$} [c] at 11 0
\put {$ \scriptstyle \bullet$} [c] at  11 6
\put {$ \scriptstyle \bullet$} [c] at 11 12
\put {$ \scriptstyle \bullet$} [c] at 12 6
\put {$ \scriptstyle \bullet$} [c] at 15 0
\setlinear \plot 15 0 11 12 10 6 11 0 12 6 11 12 11 0 /
\put{$840$} [c] at 13 -2
\put{$\scriptstyle \bullet$} [c] at 16  0 \endpicture
\end{minipage}
$$
$$
\begin{minipage}{4cm}
\beginpicture
\setcoordinatesystem units <1.5mm,2mm>
\setplotarea x from 0 to 16, y from -2 to 15
\put{1.915)} [l] at 2 12
\put {$ \scriptstyle \bullet$} [c] at 10 12
\put {$ \scriptstyle \bullet$} [c] at 12 12
\put {$ \scriptstyle \bullet$} [c] at 13 12
\put {$ \scriptstyle \bullet$} [c] at 13 6
\put {$ \scriptstyle \bullet$} [c] at 13 0
\put {$ \scriptstyle \bullet$} [c] at 14 12
\setlinear \plot 10 12  13 0  13 12 /
\setlinear \plot 12 12  13 6  14 12 /
\put{$840$} [c] at 13 -2
\put{$\scriptstyle \bullet$} [c] at 16  0 \endpicture
\end{minipage}
\begin{minipage}{4cm}
\beginpicture
\setcoordinatesystem units <1.5mm,2mm>
\setplotarea x from 0 to 16, y from -2 to 15
\put{1.916)} [l] at 2 12
\put {$ \scriptstyle \bullet$} [c] at 10 0
\put {$ \scriptstyle \bullet$} [c] at 12 0
\put {$ \scriptstyle \bullet$} [c] at 13 12
\put {$ \scriptstyle \bullet$} [c] at 13 6
\put {$ \scriptstyle \bullet$} [c] at 13 0
\put {$ \scriptstyle \bullet$} [c] at 14 0
\setlinear \plot 10 0  13 12  13 0 /
\setlinear \plot 12 0  13 6  14 0 /
\put{$840$} [c] at 13 -2
\put{$\scriptstyle \bullet$} [c] at 16  0 \endpicture
\end{minipage}
\begin{minipage}{4cm}
\beginpicture
\setcoordinatesystem units <1.5mm,2mm>
\setplotarea x from 0 to 16, y from -2 to 15
\put{1.917)} [l] at 2 12
\put {$ \scriptstyle \bullet$} [c] at 10 0
\put {$ \scriptstyle \bullet$} [c] at 10 12
\put {$ \scriptstyle \bullet$} [c] at 11.5 12
\put {$ \scriptstyle \bullet$} [c] at 12.5 12
\put {$ \scriptstyle \bullet$} [c] at 14 0
\put {$ \scriptstyle \bullet$} [c] at 14 12
\setlinear \plot 10 12 10 0 11.5 12 14 0 14 12 10 0 12.5 12 14 0 /
\setlinear \plot  10 12 14 0  /
\put{$105$} [c] at 13 -2
\put{$\scriptstyle \bullet$} [c] at 16  0 \endpicture
\end{minipage}
\begin{minipage}{4cm}
\beginpicture
\setcoordinatesystem units <1.5mm,2mm>
\setplotarea x from 0 to 16, y from -2 to 15
\put{1.918)} [l] at 2 12
\put {$ \scriptstyle \bullet$} [c] at 10 0
\put {$ \scriptstyle \bullet$} [c] at 10 12
\put {$ \scriptstyle \bullet$} [c] at 11.5 0
\put {$ \scriptstyle \bullet$} [c] at 12.5 0
\put {$ \scriptstyle \bullet$} [c] at 14 0
\put {$ \scriptstyle \bullet$} [c] at 14 12
\setlinear \plot 10 0 10 12 11.5 0 14 12 14 0 10 12 12.5 0 14 12 /
\setlinear \plot  10 0 14 12  /
\put{$105$} [c] at 13 -2
\put{$\scriptstyle \bullet$} [c] at 16  0
\endpicture
\end{minipage}
\begin{minipage}{4cm}
\beginpicture
\setcoordinatesystem units <1.5mm,2mm>
\setplotarea x from 0 to 16, y from -2 to 15
\put {${\bf  39}$} [l] at 2 15

\put {1.919)} [l] at 2 12
\put {$ \scriptstyle \bullet$} [c] at  10 12
\put {$ \scriptstyle \bullet$} [c] at  11 12
\put {$ \scriptstyle \bullet$} [c] at  12 0
\put {$ \scriptstyle \bullet$} [c] at  14 0
\put {$ \scriptstyle \bullet$} [c] at  16  12
\put {$ \scriptstyle \bullet$} [c] at  12 12
\put {$ \scriptstyle \bullet$} [c] at  14 12
\setlinear \plot  10 12 12 0 11 12  /
\setlinear \plot  12 12 12 0 14 12 14 0 12 12  /
\setlinear  \plot  14 0 16 12  /
\put{$1{,}260$} [c] at 13 -2
\endpicture
\end{minipage}
\begin{minipage}{4cm}
\beginpicture
\setcoordinatesystem units <1.5mm,2mm>
\setplotarea x from 0 to 16, y from -2 to 15
\put {1.920)} [l] at 2 12
\put {$ \scriptstyle \bullet$} [c] at  10 0
\put {$ \scriptstyle \bullet$} [c] at  11 0
\put {$ \scriptstyle \bullet$} [c] at  12 12
\put {$ \scriptstyle \bullet$} [c] at  14 12
\put {$ \scriptstyle \bullet$} [c] at  16  0
\put {$ \scriptstyle \bullet$} [c] at  12 0
\put {$ \scriptstyle \bullet$} [c] at  14 0
\setlinear \plot  10 0 12 12 11 0  /
\setlinear \plot  12 0 12 12 14 0 14 12 12 0  /
\setlinear  \plot  14 12 16 0  /
\put{$1{,}260$} [c] at 13 -2
\endpicture
\end{minipage}
$$
$$
\begin{minipage}{4cm}
\beginpicture
\setcoordinatesystem units <1.5mm,2mm>
\setplotarea x from 0 to  16, y from -2 to 15
\put{1.921)} [l] at 2 12
\put {$ \scriptstyle \bullet$} [c] at 10 12
\put {$ \scriptstyle \bullet$} [c] at 12 12
\put {$ \scriptstyle \bullet$} [c] at 12 6
\put {$ \scriptstyle \bullet$} [c] at 12 0
\put {$ \scriptstyle \bullet$} [c] at 14 12
\put {$ \scriptstyle \bullet$} [c] at 16 0
\put {$ \scriptstyle \bullet$} [c] at 16 12
\setlinear \plot  10 12   12 0 14 12    /
\setlinear \plot  12 0  12 12    /
\setlinear \plot 16 0 16 12   /
\put{$2{,}520$} [c] at 13 -2
\endpicture
\end{minipage}
\begin{minipage}{4cm}
\beginpicture
\setcoordinatesystem units <1.5mm,2mm>
\setplotarea x from 0 to  16, y from -2 to 15
\put{1.922)} [l] at 2 12
\put {$ \scriptstyle \bullet$} [c] at 10 0
\put {$ \scriptstyle \bullet$} [c] at 12 12
\put {$ \scriptstyle \bullet$} [c] at 12 6
\put {$ \scriptstyle \bullet$} [c] at 12 0
\put {$ \scriptstyle \bullet$} [c] at 14 0
\put {$ \scriptstyle \bullet$} [c] at 16 0
\put {$ \scriptstyle \bullet$} [c] at 16 12
\setlinear \plot  10 0   12 12 14 0    /
\setlinear \plot  12 0  12 12    /
\setlinear \plot 16 0 16 12   /
\put{$2{,}520$} [c] at 13 -2
\endpicture
\end{minipage}
\begin{minipage}{4cm}
\beginpicture
\setcoordinatesystem units <1.5mm,2mm>
\setplotarea x from 0 to  16, y from -2 to 15
\put{1.923)} [l] at 2 12
\put {$ \scriptstyle \bullet$} [c] at 10 0
\put {$ \scriptstyle \bullet$} [c] at 10 12
\put {$ \scriptstyle \bullet$} [c] at 12 12
\put {$ \scriptstyle \bullet$} [c] at 14 0
\put {$ \scriptstyle \bullet$} [c] at 14 12
\put {$ \scriptstyle \bullet$} [c] at 16 0
\put {$ \scriptstyle \bullet$} [c] at 16 12
\setlinear \plot  10 12 10 0 12 12 14 0   14 12  /
\setlinear \plot 16 0 16 12   /
\put{$2{,}520$} [c] at 13 -2
\endpicture
\end{minipage}
\begin{minipage}{4cm}
\beginpicture
\setcoordinatesystem units <1.5mm,2mm>
\setplotarea x from 0 to  16, y from -2 to 15
\put{1.924)} [l] at 2 12
\put {$ \scriptstyle \bullet$} [c] at 10 0
\put {$ \scriptstyle \bullet$} [c] at 10 12
\put {$ \scriptstyle \bullet$} [c] at 12 0
\put {$ \scriptstyle \bullet$} [c] at 14 0
\put {$ \scriptstyle \bullet$} [c] at 14 12
\put {$ \scriptstyle \bullet$} [c] at 16 0
\put {$ \scriptstyle \bullet$} [c] at 16 12
\setlinear \plot  10 0 10 12 12 0 14 12   14 0  /
\setlinear \plot 16 0 16 12   /
\put{$2{,}520$} [c] at 13 -2
\endpicture
\end{minipage}
\begin{minipage}{4cm}
\beginpicture
\setcoordinatesystem units <1.5mm,2mm>
\setplotarea x from 0 to 16, y from -2 to 15
\put {${\bf  40}$} [l] at 2 15

\put {1.925)} [l] at 2 12
\put {$ \scriptstyle \bullet$} [c] at  10 0
\put {$ \scriptstyle \bullet$} [c] at  10 12
\put {$ \scriptstyle \bullet$} [c] at  12 12
\put {$ \scriptstyle \bullet$} [c] at  13   0
\put {$ \scriptstyle \bullet$} [c] at  14 12
\put {$ \scriptstyle \bullet$} [c] at  16  0
\put {$ \scriptstyle \bullet$} [c] at  16 12
\setlinear \plot   10 0 10  12  13 0 16 12 16 0   /
\setlinear \plot   12 12 13 0 14 12    /
\put{$1{,}260$} [c] at 13 -2
\endpicture
\end{minipage}
\begin{minipage}{4cm}
\beginpicture
\setcoordinatesystem units <1.5mm,2mm>
\setplotarea x from 0 to 16, y from -2 to 15
\put {1.926)} [l] at 2 12
\put {$ \scriptstyle \bullet$} [c] at  10 0
\put {$ \scriptstyle \bullet$} [c] at  10 12
\put {$ \scriptstyle \bullet$} [c] at  12 0
\put {$ \scriptstyle \bullet$} [c] at  13   12
\put {$ \scriptstyle \bullet$} [c] at  14 0
\put {$ \scriptstyle \bullet$} [c] at  16  0
\put {$ \scriptstyle \bullet$} [c] at  16 12
\setlinear \plot   10 12 10  0  13 12 16 0 16 12   /
\setlinear \plot   12 0 13 12 14 0    /
\put{$1{,}260$} [c] at 13 -2
\endpicture
\end{minipage}
$$
$$
\begin{minipage}{4cm}
\beginpicture
\setcoordinatesystem units <1.5mm,2mm>
\setplotarea x from 0 to 16, y from -2 to 15
\put{1.927)} [l] at 2 12
\put {$ \scriptstyle \bullet$} [c] at 10 0
\put {$ \scriptstyle \bullet$} [c] at  10 12
\put {$ \scriptstyle \bullet$} [c] at  13 0
\put {$ \scriptstyle \bullet$} [c] at  13 6
\put {$ \scriptstyle \bullet$} [c] at  13 12
\put {$ \scriptstyle \bullet$} [c] at  16  0
\put {$ \scriptstyle \bullet$} [c] at  16 12
\setlinear \plot  10 12 10   0 13 12 /
\setlinear \plot  10  0 16 12 16 0 /
\setlinear \plot  13 0 16 12 /
\put{$1{,}260$} [c] at 13 -2
\endpicture
\end{minipage}
\begin{minipage}{4cm}
\beginpicture
\setcoordinatesystem units <1.5mm,2mm>
\setplotarea x from 0 to 16, y from -2 to 15
\put{1.928)} [l] at 2 12
\put {$ \scriptstyle \bullet$} [c] at 10 0
\put {$ \scriptstyle \bullet$} [c] at 10 12
\put {$ \scriptstyle \bullet$} [c] at 14 12
\put {$ \scriptstyle \bullet$} [c] at 14 0
\put {$ \scriptstyle \bullet$} [c] at 14 6
\put {$ \scriptstyle \bullet$} [c] at 12 12
\put{$\scriptstyle \bullet$} [c] at 16  0
\setlinear \plot 10 0 10  12 14 0  14 12 /
\setlinear \plot 12 12 14 0 /
\put{$5{,}040$} [c] at 13 -2
 \endpicture
\end{minipage}
\begin{minipage}{4cm}
\beginpicture
\setcoordinatesystem units <1.5mm,2mm>
\setplotarea x from 0 to 16, y from -2 to 15
\put{1.929)} [l] at 2 12
\put {$ \scriptstyle \bullet$} [c] at 10 0
\put {$ \scriptstyle \bullet$} [c] at 10 12
\put {$ \scriptstyle \bullet$} [c] at 14 12
\put {$ \scriptstyle \bullet$} [c] at 14 0
\put {$ \scriptstyle \bullet$} [c] at 14 6
\put {$ \scriptstyle \bullet$} [c] at 12 0
\put{$\scriptstyle \bullet$} [c] at 16  0
\setlinear \plot 10 12 10  0 14 12  14 0 /
\setlinear \plot 12 0 14 12 /
\put{$5{,}040$} [c] at 13 -2
 \endpicture
\end{minipage}
\begin{minipage}{4cm}
\beginpicture
\setcoordinatesystem units <1.5mm,2mm>
\setplotarea x from 0 to 16, y from -2 to 15
\put{1.930)} [l] at 2 12
\put {$ \scriptstyle \bullet$} [c] at 10 0
\put {$ \scriptstyle \bullet$} [c] at 11 0
\put {$ \scriptstyle \bullet$} [c] at 11 12
\put {$ \scriptstyle \bullet$} [c] at 14 12
\put {$ \scriptstyle \bullet$} [c] at 14 0
\put {$ \scriptstyle \bullet$} [c] at 15 12
\setlinear \plot 10 0 11 12 11 0 14 12  14 0 11 12  /
\setlinear \plot  14 0 15 12  /
\put{$5{,}040$} [c] at 13 -2
\put{$\scriptstyle \bullet$} [c] at 16  0 \endpicture
\end{minipage}
\begin{minipage}{4cm}
\beginpicture
\setcoordinatesystem units <1.5mm,2mm>
\setplotarea x from 0 to 16, y from -2 to 15
\put{1.931)} [l] at 2 12
\put {$ \scriptstyle \bullet$} [c] at 10 0
\put {$ \scriptstyle \bullet$} [c] at 10 12
\put {$ \scriptstyle \bullet$} [c] at 13 0
\put {$ \scriptstyle \bullet$} [c] at 13 12
\put {$ \scriptstyle \bullet$} [c] at 14 12
\put {$ \scriptstyle \bullet$} [c] at 15 0
\put{$\scriptstyle \bullet$} [c] at 16  12
\setlinear \plot 10 12 10  0 13 12 13 0 /
\setlinear \plot 14 12 15 0  16 12 /
\put{$2{,}520$} [c] at 13 -2
 \endpicture
\end{minipage}
\begin{minipage}{4cm}
\beginpicture
\setcoordinatesystem units <1.5mm,2mm>
\setplotarea x from 0 to 16, y from -2 to 15
\put{1.932)} [l] at 2 12
\put {$ \scriptstyle \bullet$} [c] at 10 0
\put {$ \scriptstyle \bullet$} [c] at 10 12
\put {$ \scriptstyle \bullet$} [c] at 13 0
\put {$ \scriptstyle \bullet$} [c] at 13 12
\put {$ \scriptstyle \bullet$} [c] at 14 0
\put {$ \scriptstyle \bullet$} [c] at 15 12
\put{$\scriptstyle \bullet$} [c] at 16  0
\setlinear \plot 10 0 10  12 13 0 13 12 /
\setlinear \plot 14 0 15 12  16 0 /
\put{$2{,}520$} [c] at 13 -2
\endpicture
\end{minipage}
$$

$$
\begin{minipage}{4cm}
\beginpicture
\setcoordinatesystem units <1.5mm,2mm>
\setplotarea x from 0 to 16, y from -2 to 15
\put{1.933)} [l] at 2 12
\put {$ \scriptstyle \bullet$} [c] at 10 0
\put {$ \scriptstyle \bullet$} [c] at 10 12
\put {$ \scriptstyle \bullet$} [c] at 11 12
\put {$ \scriptstyle \bullet$} [c] at 11 6
\put {$ \scriptstyle \bullet$} [c] at 12 0
\put {$ \scriptstyle \bullet$} [c] at 14 0
\setlinear \plot 10 12 10 0 11 6 11 12  /
\setlinear \plot 11 6 12 0  /
\put{$2{,}520$} [c] at 13 -2
\put{$\scriptstyle \bullet$} [c] at 16  0 \endpicture
\end{minipage}
\begin{minipage}{4cm}
\beginpicture
\setcoordinatesystem units <1.5mm,2mm>
\setplotarea x from 0 to 16, y from -2 to 15
\put{1.934)} [l] at 2 12
\put {$ \scriptstyle \bullet$} [c] at 10 0
\put {$ \scriptstyle \bullet$} [c] at 10 12
\put {$ \scriptstyle \bullet$} [c] at 11 0
\put {$ \scriptstyle \bullet$} [c] at 11 6
\put {$ \scriptstyle \bullet$} [c] at 12 12
\put {$ \scriptstyle \bullet$} [c] at 14 0
\setlinear \plot 10 0 10 12 11 6 11 0  /
\setlinear \plot 11 6 12 12  /
\put{$2{,}520$} [c] at 13 -2
\put{$\scriptstyle \bullet$} [c] at 16  0 \endpicture
\end{minipage}
\begin{minipage}{4cm}
\beginpicture
\setcoordinatesystem units <1.5mm,2mm>
\setplotarea x from 0 to 16, y from -2 to 15
\put{1.935)} [l] at 2 12
\put {$ \scriptstyle \bullet$} [c] at 10 12
\put {$ \scriptstyle \bullet$} [c] at 11 0
\put {$ \scriptstyle \bullet$} [c] at 12 12
\put {$ \scriptstyle \bullet$} [c] at 14 0
\put {$ \scriptstyle \bullet$} [c] at 14 6
\put {$ \scriptstyle \bullet$} [c] at 14 12
\setlinear \plot 10 12 11 0 12 12 /
\setlinear \plot 14 12 14 0 /
\put{$2{,}520$} [c] at 13 -2
\put{$\scriptstyle \bullet$} [c] at 16  0 \endpicture
\end{minipage}
\begin{minipage}{4cm}
\beginpicture
\setcoordinatesystem units <1.5mm,2mm>
\setplotarea x from 0 to 16, y from -2 to 15
\put{1.936)} [l] at 2 12
\put {$ \scriptstyle \bullet$} [c] at 10 0
\put {$ \scriptstyle \bullet$} [c] at 11 12
\put {$ \scriptstyle \bullet$} [c] at 12 0
\put {$ \scriptstyle \bullet$} [c] at 14 0
\put {$ \scriptstyle \bullet$} [c] at 14 6
\put {$ \scriptstyle \bullet$} [c] at 14 12
\setlinear \plot 10 0 11 12 12 0 /
\setlinear \plot 14 12 14 0 /
\put{$2{,}520$} [c] at 13 -2
\put{$\scriptstyle \bullet$} [c] at 16  0 \endpicture
\end{minipage}
\begin{minipage}{4cm}
\beginpicture
\setcoordinatesystem units <1.5mm,2mm>
\setplotarea x from 0 to 16, y from -2 to 15
\put{1.937)} [l] at 2 12
\put {$ \scriptstyle \bullet$} [c] at 10 12
\put {$ \scriptstyle \bullet$} [c] at 12 6
\put {$ \scriptstyle \bullet$} [c] at 12.5 0
\put {$ \scriptstyle \bullet$} [c] at 12.5 12
\put {$ \scriptstyle \bullet$} [c] at 13 6
\put {$ \scriptstyle \bullet$} [c] at 15 0
\setlinear \plot 10 12 12.5 0 12 6 12.5 12 13 6 12.5 0 /
\setlinear \plot  12.5 12 15 0 /
\put{$2{,}520$} [c] at 13 -2
\put{$\scriptstyle \bullet$} [c] at 16  0 \endpicture
\end{minipage}
\begin{minipage}{4cm}
\beginpicture
\setcoordinatesystem units <1.5mm,2mm>
\setplotarea x from 0 to 16, y from -2 to 15
\put{1.938)} [l] at 2 12
\put {$ \scriptstyle \bullet$} [c] at 10 0
\put {$ \scriptstyle \bullet$} [c] at 10 6
\put {$ \scriptstyle \bullet$} [c] at 10 12
\put {$ \scriptstyle \bullet$} [c] at 13 0
\put {$ \scriptstyle \bullet$} [c] at 13 12
\put {$ \scriptstyle \bullet$} [c] at 14.5 0
\setlinear \plot 10 12 10 0 13  12 13  0 10 12 /
\put{$2{,}520$} [c] at 13 -2
\put{$\scriptstyle \bullet$} [c] at 16  0 \endpicture
\end{minipage}
$$
$$
\begin{minipage}{4cm}
\beginpicture
\setcoordinatesystem units <1.5mm,2mm>
\setplotarea x from 0 to 16, y from -2 to 15
\put{1.939)} [l] at 2 12
\put {$ \scriptstyle \bullet$} [c] at 10 0
\put {$ \scriptstyle \bullet$} [c] at 10 6
\put {$ \scriptstyle \bullet$} [c] at 11.5 12
\put {$ \scriptstyle \bullet$} [c] at 13 6
\put {$ \scriptstyle \bullet$} [c] at 13 0
\put {$ \scriptstyle \bullet$} [c] at 14.5 0
\setlinear \plot 10 0 10 6 11.5 12 13 6 13 0  /
\put{$1{,}260$} [c] at 13 -2
\put{$\scriptstyle \bullet$} [c] at 16  0 \endpicture
\end{minipage}
\begin{minipage}{4cm}
\beginpicture
\setcoordinatesystem units <1.5mm,2mm>
\setplotarea x from 0 to 16, y from -2 to 15
\put{1.940)} [l] at 2 12
\put {$ \scriptstyle \bullet$} [c] at 10 12
\put {$ \scriptstyle \bullet$} [c] at 10 6
\put {$ \scriptstyle \bullet$} [c] at 11.5 0
\put {$ \scriptstyle \bullet$} [c] at 13 6
\put {$ \scriptstyle \bullet$} [c] at 13 12
\put {$ \scriptstyle \bullet$} [c] at 14.5 0
\setlinear \plot 10 12 10 6 11.5 0 13 6 13 12  /
\put{$1{,}260$} [c] at 13 -2
\put{$\scriptstyle \bullet$} [c] at 16  0 \endpicture
\end{minipage}
\begin{minipage}{4cm}
\beginpicture
\setcoordinatesystem units <1.5mm,2mm>
\setplotarea x from 0 to 16, y from -2 to 15
\put{1.941)} [l] at 2 12
\put {$ \scriptstyle \bullet$} [c] at 10 0
\put {$ \scriptstyle \bullet$} [c] at 10 12
\put {$ \scriptstyle \bullet$} [c] at 11.5  0
\put {$ \scriptstyle \bullet$} [c] at 13.5 0
\put {$ \scriptstyle \bullet$} [c] at 15 12
\put {$ \scriptstyle \bullet$} [c] at 15 0
\setlinear \plot 10 0 10 12 11.5 0 15 12 15 0 10 12 13.5 0 15 12 /
\put{$840$} [c] at 13 -2
\put{$\scriptstyle \bullet$} [c] at 16  0 \endpicture
\end{minipage}
\begin{minipage}{4cm}
\beginpicture
\setcoordinatesystem units <1.5mm,2mm>
\setplotarea x from 0 to 16, y from -2 to 15
\put{1.942)} [l] at 2 12
\put {$ \scriptstyle \bullet$} [c] at 10 0
\put {$ \scriptstyle \bullet$} [c] at 10 12
\put {$ \scriptstyle \bullet$} [c] at 11.5  12
\put {$ \scriptstyle \bullet$} [c] at 13.5 12
\put {$ \scriptstyle \bullet$} [c] at 15 12
\put {$ \scriptstyle \bullet$} [c] at 15 0
\setlinear \plot 10 12 10 0 11.5 12 15 0 15 12 10 0 13.5 12 15 0 /
\put{$840$} [c] at 13 -2
\put{$\scriptstyle \bullet$} [c] at 16  0 \endpicture
\end{minipage}
\begin{minipage}{4cm}
\beginpicture
\setcoordinatesystem units <1.5mm,2mm>
\setplotarea x from 0 to 16, y from -2 to 15
\put{1.943)} [l] at 2 12
\put {$ \scriptstyle \bullet$} [c] at 10 0
\put {$ \scriptstyle \bullet$} [c] at 10 4
\put {$ \scriptstyle \bullet$} [c] at 10 8
\put {$ \scriptstyle \bullet$} [c] at 10 12
\put {$ \scriptstyle \bullet$} [c] at 12 0
\put {$ \scriptstyle \bullet$} [c] at 14 0
\put{$\scriptstyle \bullet$} [c] at 16  0
\setlinear \plot 10 0 10 12 /
\put{$840$} [c] at 13 -2
 \endpicture
\end{minipage}
\begin{minipage}{4cm}
\beginpicture
\setcoordinatesystem units <1.5mm,2mm>
\setplotarea x from 0 to 16, y from -2 to 15
\put{1.944)} [l] at 2 12
\put {$ \scriptstyle \bullet$} [c] at 10 0
\put {$ \scriptstyle \bullet$} [c] at 12 6
\put {$ \scriptstyle \bullet$} [c] at 12 0
\put {$ \scriptstyle \bullet$} [c] at 12 12
\put {$ \scriptstyle \bullet$} [c] at 14 0
\put {$ \scriptstyle \bullet$} [c] at 15 0
\put {$ \scriptstyle \bullet$} [c] at 16 0
\setlinear \plot 10 0  12 6 14 0  /
\setlinear \plot 12 0  12  12   /
\put{$420$} [c] at 13 -2
\endpicture
\end{minipage}
$$

$$
\begin{minipage}{4cm}
\beginpicture
\setcoordinatesystem units <1.5mm,2mm>
\setplotarea x from 0 to 16, y from -2 to 15
\put{1.945)} [l] at 2 12
\put {$ \scriptstyle \bullet$} [c] at 10 12
\put {$ \scriptstyle \bullet$} [c] at 12 6
\put {$ \scriptstyle \bullet$} [c] at 12 0
\put {$ \scriptstyle \bullet$} [c] at 12 12
\put {$ \scriptstyle \bullet$} [c] at 14 12
\put {$ \scriptstyle \bullet$} [c] at 15 0
\put {$ \scriptstyle \bullet$} [c] at 16 0
\setlinear \plot 10 12  12 6 14 12  /
\setlinear \plot 12 0  12  12   /
\put{$420$} [c] at 13 -2
\endpicture
\end{minipage}
\begin{minipage}{4cm}
\beginpicture
\setcoordinatesystem units <1.5mm,2mm>
\setplotarea x from 0 to 16, y from -2 to 15
\put{1.946)} [l] at 2 12
\put {$ \scriptstyle \bullet$} [c] at 10 6
\put {$ \scriptstyle \bullet$} [c] at 12 0
\put {$ \scriptstyle \bullet$} [c] at 12 6
\put {$ \scriptstyle \bullet$} [c] at 12 12
\put {$ \scriptstyle \bullet$} [c] at 14 6
\put {$ \scriptstyle \bullet$} [c] at 15 0
\setlinear \plot 12 0 10 6 12 12  14 6 12 0  /
\setlinear \plot 12 0  12  12   /
\put{$420$} [c] at 13 -2
\put{$\scriptstyle \bullet$} [c] at 16  0 \endpicture
\end{minipage}
\begin{minipage}{4cm}
\beginpicture
\setcoordinatesystem units <1.5mm,2mm>
\setplotarea x from 0 to 16, y from -2 to 15
\put {${\bf  41}$} [l] at 2 15

\put {1.947)} [l] at 2 12
\put {$ \scriptstyle \bullet$} [c] at  10 12
\put {$ \scriptstyle \bullet$} [c] at  11.5  12
\put {$ \scriptstyle \bullet$} [c] at  13 12
\put {$ \scriptstyle \bullet$} [c] at  14.5 12
\put {$ \scriptstyle \bullet$} [c] at  16  12
\put {$ \scriptstyle \bullet$} [c] at  11.5 0
\put {$ \scriptstyle \bullet$} [c] at  14.5 0
\setlinear \plot   10 12 11.5 0  13 12 14.5 0 16 12  /
\setlinear \plot  11.5 0 11.5 12 /
\setlinear  \plot  14.5  0 14.5 12 /
\put{$630$} [c] at 13 -2
\endpicture
\end{minipage}
\begin{minipage}{4cm}
\beginpicture
\setcoordinatesystem units <1.5mm,2mm>
\setplotarea x from 0 to 16, y from -2 to 15
\put {1.948)} [l] at 2 12
\put {$ \scriptstyle \bullet$} [c] at  10 0
\put {$ \scriptstyle \bullet$} [c] at  11.5  0
\put {$ \scriptstyle \bullet$} [c] at  13 0
\put {$ \scriptstyle \bullet$} [c] at  14.5 0
\put {$ \scriptstyle \bullet$} [c] at  16  0
\put {$ \scriptstyle \bullet$} [c] at  11.5 12
\put {$ \scriptstyle \bullet$} [c] at  14.5 12
\setlinear \plot   10 0 11.5 12  13 0 14.5 12 16 0  /
\setlinear \plot  11.5 0 11.5 12 /
\setlinear  \plot  14.5  0 14.5 12 /
\put{$630$} [c] at 13 -2
\endpicture
\end{minipage}
\begin{minipage}{4cm}
\beginpicture
\setcoordinatesystem units <1.5mm,2mm>
\setplotarea x from 0 to 16, y from -2 to 15
\put{1.949)} [l] at 2 12
\put {$ \scriptstyle \bullet$} [c] at 10 12
\put {$ \scriptstyle \bullet$} [c] at 12 12
\put {$ \scriptstyle \bullet$} [c] at 13 12
\put {$ \scriptstyle \bullet$} [c] at 14.5 12
\put {$ \scriptstyle \bullet$} [c] at 16 12
\put {$ \scriptstyle \bullet$} [c] at 14 0
\put {$ \scriptstyle \bullet$} [c] at 11  6
\setlinear \plot 10 12 11 6 14 0 16 12    /
\setlinear \plot 12 12 11 6     /
\setlinear \plot 13 12 14 0 14.5 12    /
\put{$420$} [c] at 13 -2
\endpicture
\end{minipage}
$$
$$
\begin{minipage}{4cm}
\beginpicture
\setcoordinatesystem units <1.5mm,2mm>
\setplotarea x from 0 to 16, y from -2 to 15
\put{1.950)} [l] at 2 12
\put {$ \scriptstyle \bullet$} [c] at 10 0
\put {$ \scriptstyle \bullet$} [c] at 12 0
\put {$ \scriptstyle \bullet$} [c] at 13 0
\put {$ \scriptstyle \bullet$} [c] at 14.5 0
\put {$ \scriptstyle \bullet$} [c] at 16 0
\put {$ \scriptstyle \bullet$} [c] at 14  12
\put {$ \scriptstyle \bullet$} [c] at 11  6
\setlinear \plot 10 0 11 6 14 12 16 0    /
\setlinear \plot 12 0 11 6     /
\setlinear \plot 13 0 14 12 14.5 0    /
\put{$420$} [c] at 13 -2
\endpicture
\end{minipage}
\begin{minipage}{4cm}
\beginpicture
\setcoordinatesystem units <1.5mm,2mm>
\setplotarea x from 0 to 16, y from -2 to 15
\put{1.951)} [l] at 2 12
\put {$ \scriptstyle \bullet$} [c] at 10 6
\put {$ \scriptstyle \bullet$} [c] at 10.5 12
\put {$ \scriptstyle \bullet$} [c] at 10.5 0
\put {$ \scriptstyle \bullet$} [c] at 11 6
\put {$ \scriptstyle \bullet$} [c] at 14 12
\put {$ \scriptstyle \bullet$} [c] at 15  12
\put {$ \scriptstyle \bullet$} [c] at 16  12
\setlinear \plot 10.5 0 10 6 10.5 12 11 6 10.5 0 16 12     /
\setlinear \plot 14 12 10.5 0 15 12     /
\put{$420$} [c] at 13 -2
\endpicture
\end{minipage}
\begin{minipage}{4cm}
\beginpicture
\setcoordinatesystem units <1.5mm,2mm>
\setplotarea x from 0 to 16, y from -2 to 15
\put{1.952)} [l] at 2 12
\put {$ \scriptstyle \bullet$} [c] at 10 6
\put {$ \scriptstyle \bullet$} [c] at 10.5 12
\put {$ \scriptstyle \bullet$} [c] at 10.5 0
\put {$ \scriptstyle \bullet$} [c] at 11 6
\put {$ \scriptstyle \bullet$} [c] at 14 0
\put {$ \scriptstyle \bullet$} [c] at 15  0
\put {$ \scriptstyle \bullet$} [c] at 16  0
\setlinear \plot 10.5 12 10 6 10.5 0 11 6 10.5 12 16 0     /
\setlinear \plot 14 0 10.5 12 15 0     /
\put{$420$} [c] at 13 -2
\endpicture
\end{minipage}
\begin{minipage}{4cm}
\beginpicture
\setcoordinatesystem units <1.5mm,2mm>
\setplotarea x from 0 to 16, y from -2 to 15
\put {${\bf  42}$} [l] at 2 15

\put {1.953)} [l] at 2 12
\put {$ \scriptstyle \bullet$} [c] at  10 12
\put {$ \scriptstyle \bullet$} [c] at  12 12
\put {$ \scriptstyle \bullet$} [c] at  14 12
\put {$ \scriptstyle \bullet$} [c] at  16 12
\put {$ \scriptstyle \bullet$} [c] at  13 0
\put {$ \scriptstyle \bullet$} [c] at  16 0
\put {$ \scriptstyle \bullet$} [c] at  14.8 7
\setlinear \plot  10 12  13 0 16 12 16 0 /
\setlinear \plot  12 12 13 0 14 12  /
\put{$840$} [c] at 13 -2
\endpicture
\end{minipage}
\begin{minipage}{4cm}
\beginpicture
\setcoordinatesystem units <1.5mm,2mm>
\setplotarea x from 0 to 16, y from -2 to 15
\put {1.954)} [l] at 2 12
\put {$ \scriptstyle \bullet$} [c] at  10 0
\put {$ \scriptstyle \bullet$} [c] at  12 0
\put {$ \scriptstyle \bullet$} [c] at  14 0
\put {$ \scriptstyle \bullet$} [c] at  16 0
\put {$ \scriptstyle \bullet$} [c] at  13 12
\put {$ \scriptstyle \bullet$} [c] at  16 12
\put {$ \scriptstyle \bullet$} [c] at  14.3 7
\setlinear \plot  10 0  13 12 16 0 16 12 /
\setlinear \plot  12 0 13 12 14 0  /
\put{$840$} [c] at 13 -2
\endpicture
\end{minipage}
\begin{minipage}{4cm}
\beginpicture
\setcoordinatesystem units <1.5mm,2mm>
\setplotarea x from 0 to 16, y from -2 to 15
\put {1.955)} [l] at 2 12
\put {$ \scriptstyle \bullet$} [c] at  10 12
\put {$ \scriptstyle \bullet$} [c] at  11.5 12
\put {$ \scriptstyle \bullet$} [c] at  13 12
\put {$ \scriptstyle \bullet$} [c] at  14 12
\put {$ \scriptstyle \bullet$} [c] at  16 12
\put {$ \scriptstyle \bullet$} [c] at  14 0
\put {$ \scriptstyle \bullet$} [c] at  16 0
\setlinear \plot  10 12 14 0 14 12 16 0 16 12 14 0 13 12    /
\setlinear  \plot  11.5 12 14 0 /
\put{$420$} [c] at 13 -2
\endpicture
\end{minipage}
$$
$$
\begin{minipage}{4cm}
\beginpicture
\setcoordinatesystem units <1.5mm,2mm>
\setplotarea x from 0 to 16, y from -2 to 15
\put {1.956)} [l] at 2 12
\put {$ \scriptstyle \bullet$} [c] at  10 0
\put {$ \scriptstyle \bullet$} [c] at  11.5 0
\put {$ \scriptstyle \bullet$} [c] at  13 0
\put {$ \scriptstyle \bullet$} [c] at  14 0
\put {$ \scriptstyle \bullet$} [c] at  16 0
\put {$ \scriptstyle \bullet$} [c] at  14 12
\put {$ \scriptstyle \bullet$} [c] at  16 12
\setlinear \plot  10 0 14 12 14 0 16 12 16 0 14 12 13 0    /
\setlinear  \plot  11.5 0 14 12 /
\put{$420$} [c] at 13 -2
\endpicture
\end{minipage}
\begin{minipage}{4cm}
\beginpicture
\setcoordinatesystem units <1.5mm,2mm>
\setplotarea x from 0 to 16, y from -2 to 15
\put{1.957)} [l] at 2 12
\put {$ \scriptstyle \bullet$} [c] at 10 0
\put {$ \scriptstyle \bullet$} [c] at 14 0
\put {$ \scriptstyle \bullet$} [c] at 12 12
\put {$ \scriptstyle \bullet$} [c] at 11 6
\put {$ \scriptstyle \bullet$} [c] at 15 0
\put {$ \scriptstyle \bullet$} [c] at 15 12
\setlinear \plot 10 0 12 12 14 0   /
\setlinear \plot  15 12 15 0   /
\put{$5{,}040$} [c] at 13 -2
\put{$\scriptstyle \bullet$} [c] at 16  0 \endpicture
\end{minipage}
\begin{minipage}{4cm}
\beginpicture
\setcoordinatesystem units <1.5mm,2mm>
\setplotarea x from 0 to 16, y from -2 to 15
\put{1.958)} [l] at 2 12
\put {$ \scriptstyle \bullet$} [c] at 10 12
\put {$ \scriptstyle \bullet$} [c] at 14 12
\put {$ \scriptstyle \bullet$} [c] at 12 0
\put {$ \scriptstyle \bullet$} [c] at 11 6
\put {$ \scriptstyle \bullet$} [c] at 15 0
\put {$ \scriptstyle \bullet$} [c] at 15 12
\setlinear \plot 10 12 12 0 14 12   /
\setlinear \plot  15 12 15 0   /
\put{$5{,}040$} [c] at 13 -2
\put{$\scriptstyle \bullet$} [c] at 16  0
\endpicture
\end{minipage}
\begin{minipage}{4cm}
\beginpicture
\setcoordinatesystem units <1.5mm,2mm>
\setplotarea x from 0 to 16, y from -2 to 15
\put{1.959)} [l] at 2 12
\put {$ \scriptstyle \bullet$} [c] at 10 12
\put {$ \scriptstyle \bullet$} [c] at 10 0
\put {$ \scriptstyle \bullet$} [c] at 11.5 0
\put {$ \scriptstyle \bullet$} [c] at 13 12
\put {$ \scriptstyle \bullet$} [c] at 14 0
\put {$ \scriptstyle \bullet$} [c] at 14 12
\setlinear \plot 10 0 10 12 11.5 0 13  12 14 0 14 12 /
\put{$5{,}040$} [c] at 13 -2
\put{$\scriptstyle \bullet$} [c] at 16  0 \endpicture
\end{minipage}
\begin{minipage}{4cm}
\beginpicture
\setcoordinatesystem units <1.5mm,2mm>
\setplotarea x from 0 to  16, y from -2 to 15
\put{1.960)} [l] at 2 12
\put {$ \scriptstyle \bullet$} [c] at 10 12
\put {$ \scriptstyle \bullet$} [c] at 12 12
\put {$ \scriptstyle \bullet$} [c] at 12 0
\put {$ \scriptstyle \bullet$} [c] at 14 0
\put {$ \scriptstyle \bullet$} [c] at 14 12
\setlinear \plot  10 12   12 0 14 12 14 0    /
\setlinear \plot  12 12   12 0      /
\put {$ \scriptstyle \bullet$} [c] at 16 0
\put {$ \scriptstyle \bullet$} [c] at 16 12
\setlinear \plot 16 0 16 12   /
\put{$2{,}520$} [c] at 13 -2
\endpicture
\end{minipage}
\begin{minipage}{4cm}
\beginpicture
\setcoordinatesystem units <1.5mm,2mm>
\setplotarea x from 0 to  16, y from -2 to 15
\put{1.961)} [l] at 2 12
\put {$ \scriptstyle \bullet$} [c] at 10 0
\put {$ \scriptstyle \bullet$} [c] at 12 0
\put {$ \scriptstyle \bullet$} [c] at 12 12
\put {$ \scriptstyle \bullet$} [c] at 14 0
\put {$ \scriptstyle \bullet$} [c] at 14 12
\setlinear \plot  10 0 12 12   14 0 14 12    /
\setlinear \plot   12 12  12 0    /
\put {$ \scriptstyle \bullet$} [c] at 16 0
\put {$ \scriptstyle \bullet$} [c] at 16 12
\setlinear \plot 16 0 16 12   /
\put{$2{,}520$} [c] at 13 -2
\endpicture
\end{minipage}
$$
$$
\begin{minipage}{4cm}
\beginpicture
\setcoordinatesystem units <1.5mm,2mm>
\setplotarea x from 0 to 16, y from -2 to 15
\put{1.962)} [l] at 2 12
\put {$ \scriptstyle \bullet$} [c] at 10 0
\put {$ \scriptstyle \bullet$} [c] at 11.5 6
\put {$ \scriptstyle \bullet$} [c] at 12 12
\put {$ \scriptstyle \bullet$} [c] at 12 0
\put {$ \scriptstyle \bullet$} [c] at 12.5 6
\put {$ \scriptstyle \bullet$} [c] at 14 0
\setlinear \plot 10 0 12 12 11.5 6 12 0 12.5 6 12 12 14 0 /
\put{$1{,}260$} [c] at 13 -2
\put{$\scriptstyle \bullet$} [c] at 16  0 \endpicture
\end{minipage}
\begin{minipage}{4cm}
\beginpicture
\setcoordinatesystem units <1.5mm,2mm>
\setplotarea x from 0 to 16, y from -2 to 15
\put{1.963)} [l] at 2 12
\put {$ \scriptstyle \bullet$} [c] at 10 12
\put {$ \scriptstyle \bullet$} [c] at 11.5 6
\put {$ \scriptstyle \bullet$} [c] at 12 0
\put {$ \scriptstyle \bullet$} [c] at 12 12
\put {$ \scriptstyle \bullet$} [c] at 12.5 6
\put {$ \scriptstyle \bullet$} [c] at 14 12
\setlinear \plot 10 12 12 0 11.5 6 12 12 12.5 6 12 0 14 12 /
\put{$1{,}260$} [c] at 13 -2
\put{$\scriptstyle \bullet$} [c] at 16  0 \endpicture
\end{minipage}
\begin{minipage}{4cm}
\beginpicture
\setcoordinatesystem units <1.5mm,2mm>
\setplotarea x from 0 to 16, y from -2 to 15
\put{1.964)} [l] at 2 12
\put {$ \scriptstyle \bullet$} [c] at 10 12
\put {$ \scriptstyle \bullet$} [c] at 11 12
\put {$ \scriptstyle \bullet$} [c] at 11 0
\put {$ \scriptstyle \bullet$} [c] at 13 0
\put {$ \scriptstyle \bullet$} [c] at 13 12
\put {$ \scriptstyle \bullet$} [c] at 14 12
\setlinear \plot 10 12 11 0 13 12 13 0 11 12 11 0 /
\setlinear \plot 13 0 14  12     /
\put{$1{,}260$} [c] at 13 -2
\put{$\scriptstyle \bullet$} [c] at 16  0 \endpicture
\end{minipage}
\begin{minipage}{4cm}
\beginpicture
\setcoordinatesystem units <1.5mm,2mm>
\setplotarea x from 0 to 16, y from -2 to 15
\put{1.965)} [l] at 2 12
\put {$ \scriptstyle \bullet$} [c] at 10 0
\put {$ \scriptstyle \bullet$} [c] at 11 12
\put {$ \scriptstyle \bullet$} [c] at 11 0
\put {$ \scriptstyle \bullet$} [c] at 13 0
\put {$ \scriptstyle \bullet$} [c] at 13 12
\put {$ \scriptstyle \bullet$} [c] at 14 0
\setlinear \plot 10 0 11 12 13 0 13 12 11 0 11 12 /
\setlinear \plot 13 12 14  0     /
\put{$1{,}260$} [c] at 13 -2
\put{$\scriptstyle \bullet$} [c] at 16  0 \endpicture
\end{minipage}
\begin{minipage}{4cm}
\beginpicture
\setcoordinatesystem units <1.5mm,2mm>
\setplotarea x from 0 to 16, y from -2 to 15
\put{1.966)} [l] at 2 12
\put {$ \scriptstyle \bullet$} [c] at 10 12
\put {$ \scriptstyle \bullet$} [c] at 12 0
\put {$ \scriptstyle \bullet$} [c] at 12 12
\put {$ \scriptstyle \bullet$} [c] at 14 6
\put {$ \scriptstyle \bullet$} [c] at 13 12
\put {$ \scriptstyle \bullet$} [c] at 15 12
\setlinear \plot 10 12 12 0 12 12   /
\setlinear \plot 12 0 14  6 13  12   /
\setlinear \plot 14 6 15 12   /
\put{$1{,}260$} [c] at 13 -2
\put{$\scriptstyle \bullet$} [c] at 16  0 \endpicture
\end{minipage}
\begin{minipage}{4cm}
\beginpicture
\setcoordinatesystem units <1.5mm,2mm>
\setplotarea x from 0 to 16, y from -2 to 15
\put{1.967)} [l] at 2 12
\put {$ \scriptstyle \bullet$} [c] at 10 0
\put {$ \scriptstyle \bullet$} [c] at 12 0
\put {$ \scriptstyle \bullet$} [c] at 12 12
\put {$ \scriptstyle \bullet$} [c] at 14 6
\put {$ \scriptstyle \bullet$} [c] at 13 0
\put {$ \scriptstyle \bullet$} [c] at 15 0
\setlinear \plot 10 0 12 12 12 0   /
\setlinear \plot 12 12 14  6 13  0   /
\setlinear \plot 14 6 15 0   /
\put{$1{,}260$} [c] at 13 -2
\put{$\scriptstyle \bullet$} [c] at 16  0 \endpicture
\end{minipage}
$$
$$
\begin{minipage}{4cm}
\beginpicture
\setcoordinatesystem units <1.5mm,2mm>
\setplotarea x from 0 to 16, y from -2 to 15
\put{1.968)} [l] at 2 12
\put {$ \scriptstyle \bullet$} [c] at 10 12
\put {$ \scriptstyle \bullet$} [c] at 10 0
\put {$ \scriptstyle \bullet$} [c] at 13 12
\put {$ \scriptstyle \bullet$} [c] at 13 0
\put {$ \scriptstyle \bullet$} [c] at 14  0
\put {$ \scriptstyle \bullet$} [c] at 14 12
\setlinear \plot 10 12 10 0  13 12 13 0 10 12  /
\setlinear \plot 14 0 14 12    /
\put{$1{,}260$} [c] at 13 -2
\put{$\scriptstyle \bullet$} [c] at 16  0
\endpicture
\end{minipage}
\begin{minipage}{4cm}
\beginpicture
\setcoordinatesystem units <1.5mm,2mm>
\setplotarea x from 0 to 16, y from -2 to 15
\put {${\bf  43}$}[l] at 2 15

\put {1.969)}[l] at 2 12
\put {$ \scriptstyle \bullet$} [c] at  10 12
\put {$ \scriptstyle \bullet$} [c] at  11 12
\put {$ \scriptstyle \bullet$} [c] at  13 12
\put {$ \scriptstyle \bullet$} [c] at  14 12
\put {$ \scriptstyle \bullet$} [c] at  16  12
\put {$ \scriptstyle \bullet$} [c] at  12  0
\put {$ \scriptstyle \bullet$} [c] at  15 0
\setlinear \plot   10 12 12 0 14  12 15 0 16 12  /
\setlinear \plot  11 12 12 0 13 12 /
\put{$840$} [c] at 13 -2
\endpicture
\end{minipage}
\begin{minipage}{4cm}
\beginpicture
\setcoordinatesystem units <1.5mm,2mm>
\setplotarea x from 0 to 16, y from -2 to 15
\put {1.970)}[l] at 2 12
\put {$ \scriptstyle \bullet$} [c] at  10 0
\put {$ \scriptstyle \bullet$} [c] at  11 0
\put {$ \scriptstyle \bullet$} [c] at  13 0
\put {$ \scriptstyle \bullet$} [c] at  14 0
\put {$ \scriptstyle \bullet$} [c] at  16  0
\put {$ \scriptstyle \bullet$} [c] at  12  12
\put {$ \scriptstyle \bullet$} [c] at  15 12
\setlinear \plot   10 0 12 12 14  0 15 12 16 0  /
\setlinear \plot  11 0 12 12 13 0 /
\put{$840$} [c] at 13 -2
\endpicture
\end{minipage}
\begin{minipage}{4cm}
\beginpicture
\setcoordinatesystem units <1.5mm,2mm>
\setplotarea x from 0 to 16, y from -2 to 15
\put{${\bf  44}$} [l] at 2 15

\put{1.971)} [l] at 2 12
\put {$ \scriptstyle \bullet$} [c] at  10 0
\put {$ \scriptstyle \bullet$} [c] at  10 12
\put {$ \scriptstyle \bullet$} [c] at  12 0
\put {$ \scriptstyle \bullet$} [c] at  12 12
\put {$ \scriptstyle \bullet$} [c] at  14 12
\put {$ \scriptstyle \bullet$} [c] at  16  0
\put {$ \scriptstyle \bullet$} [c] at  16 12
\setlinear \plot   10 0 10  12  16 0 16 12    /
\setlinear \plot   12 12 16 0  14 12    /
\setlinear \plot   10 12 12 0    /
\put{$420$} [c] at 13 -2
\endpicture
\end{minipage}
\begin{minipage}{4cm}
\beginpicture
\setcoordinatesystem units <1.5mm,2mm>
\setplotarea x from 0 to 16, y from -2 to 15
\put{1.972)} [l] at 2 12
\put {$ \scriptstyle \bullet$} [c] at  10 0
\put {$ \scriptstyle \bullet$} [c] at  10 12
\put {$ \scriptstyle \bullet$} [c] at  12 0
\put {$ \scriptstyle \bullet$} [c] at  12 12
\put {$ \scriptstyle \bullet$} [c] at  14 0
\put {$ \scriptstyle \bullet$} [c] at  16  0
\put {$ \scriptstyle \bullet$} [c] at  16 12
\setlinear \plot   10 12 10  0  16 12 16 0    /
\setlinear \plot   12 0 16 12  14 0    /
\setlinear \plot   10 0 12 12    /
\put{$420$} [c] at 13 -2
\endpicture
\end{minipage}
\begin{minipage}{4cm}
\beginpicture
\setcoordinatesystem units <1.5mm,2mm>
\setplotarea x from 0 to 16, y from -2 to 15
\put{1.973)} [l] at 2 12
\put {$ \scriptstyle \bullet$} [c] at 10 12
\put {$ \scriptstyle \bullet$} [c] at 12 12
\put {$ \scriptstyle \bullet$} [c] at 12 0
\put {$ \scriptstyle \bullet$} [c] at 14 0
\put {$ \scriptstyle \bullet$} [c] at 14 12
\put {$ \scriptstyle \bullet$} [c] at 13 6
\setlinear \plot 10 12 12 0 14  12 14 0  /
\setlinear \plot 12 0 12 12   /
\put{$2{,}520$} [c] at 13 -2
\put{$\scriptstyle \bullet$} [c] at 16  0
\endpicture
\end{minipage}
$$
$$
\begin{minipage}{4cm}
\beginpicture
\setcoordinatesystem units <1.5mm,2mm>
\setplotarea x from 0 to 16, y from -2  to 15
\put{1.974)} [l] at 2 12
\put {$ \scriptstyle \bullet$} [c] at 10 0
\put {$ \scriptstyle \bullet$} [c] at 12 12
\put {$ \scriptstyle \bullet$} [c] at 12 0
\put {$ \scriptstyle \bullet$} [c] at 14 0
\put {$ \scriptstyle \bullet$} [c] at 14 12
\put {$ \scriptstyle \bullet$} [c] at 13 6
\setlinear \plot 10 0 12 12 14  0 14 12  /
\setlinear \plot 12 0 12 12   /
\put{$2{,}520$} [c] at 13 -2
\put{$\scriptstyle \bullet$} [c] at 16  0 \endpicture
\end{minipage}
\begin{minipage}{4cm}
\beginpicture
\setcoordinatesystem units <1.5mm,2mm>
\setplotarea x from 0 to 16, y from -2  to 15
\put{1.975)} [l] at 2 12
\put {$ \scriptstyle \bullet$} [c] at 10 12
\put {$ \scriptstyle \bullet$} [c] at 10 0
\put {$ \scriptstyle \bullet$} [c] at 12 12
\put {$ \scriptstyle \bullet$} [c] at 12 0
\put {$ \scriptstyle \bullet$} [c] at 14  0
\put {$ \scriptstyle \bullet$} [c] at 14 12
\setlinear \plot 10 0 10 12  12 0 14 12 14 0  /
\setlinear \plot 12 0 12 12    /
\put{$2{,}520$} [c] at 13 -2
\put{$\scriptstyle \bullet$} [c] at 16  0 \endpicture
\end{minipage}
\begin{minipage}{4cm}
\beginpicture
\setcoordinatesystem units <1.5mm,2mm>
\setplotarea x from 0 to 16, y from -2  to 15
\put{1.976)} [l] at 2 12
\put {$ \scriptstyle \bullet$} [c] at 10 12
\put {$ \scriptstyle \bullet$} [c] at 10 0
\put {$ \scriptstyle \bullet$} [c] at 12 12
\put {$ \scriptstyle \bullet$} [c] at 12 0
\put {$ \scriptstyle \bullet$} [c] at 14  0
\put {$ \scriptstyle \bullet$} [c] at 14 12
\setlinear \plot 10 12 10 0  12 12 14 0 14 12  /
\setlinear \plot 12 0 12 12    /
\put{$2{,}520$} [c] at 13 -2
\put{$\scriptstyle \bullet$} [c] at 16  0 \endpicture
\end{minipage}
\begin{minipage}{4cm}
\beginpicture
\setcoordinatesystem units <1.5mm,2mm>
\setplotarea x from 0 to 16, y from -2  to 15
\put{1.977)} [l] at 2 12
\put {$ \scriptstyle \bullet$} [c] at 10 0
\put {$ \scriptstyle \bullet$} [c] at 10 6
\put {$ \scriptstyle \bullet$} [c] at 10 12
\put {$ \scriptstyle \bullet$} [c] at 14 12
\put {$ \scriptstyle \bullet$} [c] at 14 0
\put {$ \scriptstyle \bullet$} [c] at 15 0
\setlinear \plot 10 12 10 0  14 12 14 0  /
\put{$2{,}520$} [c] at 13 -2
\put{$\scriptstyle \bullet$} [c] at 16  0 \endpicture
\end{minipage}
\begin{minipage}{4cm}
\beginpicture
\setcoordinatesystem units <1.5mm,2mm>
\setplotarea x from 0 to 16, y from -2  to 15
\put{1.978)} [l] at 2 12
\put {$ \scriptstyle \bullet$} [c] at 10 0
\put {$ \scriptstyle \bullet$} [c] at 10 6
\put {$ \scriptstyle \bullet$} [c] at 10 12
\put {$ \scriptstyle \bullet$} [c] at 14 12
\put {$ \scriptstyle \bullet$} [c] at 14 0
\put {$ \scriptstyle \bullet$} [c] at 15 0
\setlinear \plot 10 0 10 12  14 0 14 12  /
\put{$2{,}520$} [c] at 13 -2
\put{$\scriptstyle \bullet$} [c] at 16  0 \endpicture
\end{minipage}
\begin{minipage}{4cm}
\beginpicture
\setcoordinatesystem units <1.5mm,2mm>
\setplotarea x from 0 to 16, y from -2  to 15
\put{1.979)} [l] at 2 12
\put {$ \scriptstyle \bullet$} [c] at 10 0
\put {$ \scriptstyle \bullet$} [c] at 10 12
\put {$ \scriptstyle \bullet$} [c] at 12 0
\put {$ \scriptstyle \bullet$} [c] at 12 12
\put {$ \scriptstyle \bullet$} [c] at 14 12
\put {$ \scriptstyle \bullet$} [c] at 15 12
\setlinear \plot 15 12 12 0 10 12 10 0 12 12 12 0 14 12 /
\put{$1{,}260$} [c] at 13 -2
\put{$\scriptstyle \bullet$} [c] at 16  0 \endpicture
\end{minipage}
$$
$$
\begin{minipage}{4cm}
\beginpicture
\setcoordinatesystem units <1.5mm,2mm>
\setplotarea x from 0 to 16, y from -2  to 15
\put{1.980)} [l] at 2 12
\put {$ \scriptstyle \bullet$} [c] at 10 0
\put {$ \scriptstyle \bullet$} [c] at 10 12
\put {$ \scriptstyle \bullet$} [c] at 12 0
\put {$ \scriptstyle \bullet$} [c] at 12 12
\put {$ \scriptstyle \bullet$} [c] at 14 0
\put {$ \scriptstyle \bullet$} [c] at 15 0
\setlinear \plot 15 0 12 12 10 0 10 12 12 0 12 12 14 0 /
\put{$1{,}260$} [c] at 13 -2
\put{$\scriptstyle \bullet$} [c] at 16  0 \endpicture
\end{minipage}
\begin{minipage}{4cm}
\beginpicture
\setcoordinatesystem units <1.5mm,2mm>
\setplotarea x from 0 to 16, y from -2  to 15
\put{1.981)} [l] at 2 12
\put {$ \scriptstyle \bullet$} [c] at 10 0
\put {$ \scriptstyle \bullet$} [c] at 12 0
\put {$ \scriptstyle \bullet$} [c] at 14 0
\put {$ \scriptstyle \bullet$} [c] at 13 6
\put {$ \scriptstyle \bullet$} [c] at 13 12
\put {$ \scriptstyle \bullet$} [c] at 15 0
\setlinear \plot 10 0 13  12 13 6 12 0 /
\setlinear \plot 13 6  14 0   /
\put{$1{,}260$} [c] at 13 -2
\put{$\scriptstyle \bullet$} [c] at 16  0 \endpicture
\end{minipage}
\begin{minipage}{4cm}
\beginpicture
\setcoordinatesystem units <1.5mm,2mm>
\setplotarea x from 0 to 16, y from -2  to 15
\put{1.982)} [l] at 2 12
\put {$ \scriptstyle \bullet$} [c] at 10 12
\put {$ \scriptstyle \bullet$} [c] at 12 12
\put {$ \scriptstyle \bullet$} [c] at 14 12
\put {$ \scriptstyle \bullet$} [c] at 13 6
\put {$ \scriptstyle \bullet$} [c] at 13 0
\put {$ \scriptstyle \bullet$} [c] at 15 0
\setlinear \plot 10 12 13 0 13 6 12 12 /
\setlinear \plot 13 6  14 12   /
\put{$1{,}260$} [c] at 13 -2
\put{$\scriptstyle \bullet$} [c] at 16  0 \endpicture
\end{minipage}
\begin{minipage}{4cm}
\beginpicture
\setcoordinatesystem units <1.5mm,2mm>
\setplotarea x from 0 to 16, y from -2  to 15
\put{1.983)} [l] at 2 12
\put {$ \scriptstyle \bullet$} [c] at 10 0
\put {$ \scriptstyle \bullet$} [c] at 12.5 6
\put {$ \scriptstyle \bullet$} [c] at 13 12
\put {$ \scriptstyle \bullet$} [c] at 13 0
\put {$ \scriptstyle \bullet$} [c] at 13.5 6
\put {$ \scriptstyle \bullet$} [c] at 15 0
\setlinear \plot 10 0 13  12 12.5 6 13 0 13.5 6 13 12 /
\put{$1{,}260$} [c] at 13 -2
\put{$\scriptstyle \bullet$} [c] at 16  0
\endpicture
\end{minipage}
\begin{minipage}{4cm}
\beginpicture
\setcoordinatesystem units <1.5mm,2mm>
\setplotarea x from 0 to 16, y from -2  to 15
\put{1.984)} [l] at 2 12
\put {$ \scriptstyle \bullet$} [c] at 10 12
\put {$ \scriptstyle \bullet$} [c] at 12.5 6
\put {$ \scriptstyle \bullet$} [c] at 13 12
\put {$ \scriptstyle \bullet$} [c] at 13 0
\put {$ \scriptstyle \bullet$} [c] at 13.5 6
\put {$ \scriptstyle \bullet$} [c] at 15 0
\setlinear \plot 10 12 13  0 12.5 6 13 12 13.5 6 13 0 /
\put{$1{,}260$} [c] at 13 -2
\put{$\scriptstyle \bullet$} [c] at 16  0
\endpicture
\end{minipage}
\begin{minipage}{4cm}
\beginpicture
\setcoordinatesystem units <1.5mm,2mm>
\setplotarea x from 0 to 16, y from -2  to 15
\put{1.985)} [l] at 2 12
\put {$ \scriptstyle \bullet$} [c] at 10 12
\put {$ \scriptstyle \bullet$} [c] at 10 0
\put {$ \scriptstyle \bullet$} [c] at 12 12
\put {$ \scriptstyle \bullet$} [c] at 14 0
\put {$ \scriptstyle \bullet$} [c] at 14 12
\put {$ \scriptstyle \bullet$} [c] at 15 0
\setlinear \plot 10 12 10 0 14 12 14 0 10 12 /
\setlinear \plot 10 0 12 12 14 0    /
\put{$210$} [c] at 13 -2
\put{$\scriptstyle \bullet$} [c] at 16  0 \endpicture
\end{minipage}
$$
$$
\begin{minipage}{4cm}
\beginpicture
\setcoordinatesystem units <1.5mm,2mm>
\setplotarea x from 0 to 16, y from -2  to 15
\put{1.986)} [l] at 2 12
\put {$ \scriptstyle \bullet$} [c] at 10 12
\put {$ \scriptstyle \bullet$} [c] at 10 0
\put {$ \scriptstyle \bullet$} [c] at 12 0
\put {$ \scriptstyle \bullet$} [c] at 14 0
\put {$ \scriptstyle \bullet$} [c] at 14 12
\put {$ \scriptstyle \bullet$} [c] at 15 0
\setlinear \plot 10 0 10 12 14 0 14 12 10 0 /
\setlinear \plot 10 12 12 0 14 12    /
\put{$210$} [c] at 13 -2
\put{$\scriptstyle \bullet$} [c] at 16  0 \endpicture
\end{minipage}
\begin{minipage}{4cm}
\beginpicture
\setcoordinatesystem units <1.5mm,2mm>
\setplotarea  x from 0 to 16, y from -2 to 15
\put{${\bf  45}$} [l]  at 2 15

\put{1.987)} [l]  at 2 12
\put {$ \scriptstyle \bullet$} [c] at 10 12
\put {$ \scriptstyle \bullet$} [c] at 11 0
\put {$ \scriptstyle \bullet$} [c] at 12 12
\put {$ \scriptstyle \bullet$} [c] at 14 0
\put {$ \scriptstyle \bullet$} [c] at 14 12
\put {$ \scriptstyle \bullet$} [c] at 16 0
\put {$ \scriptstyle \bullet$} [c] at 16 12
\setlinear \plot  10 12 11 0 12 12      /
\setlinear \plot  14 0 14 12    /
\setlinear \plot 16 0 16  12   /
\put{$1{,}260$} [c] at 13 -2
\endpicture
\end{minipage}
\begin{minipage}{4cm}
\beginpicture
\setcoordinatesystem units <1.5mm,2mm>
\setplotarea  x from 0 to 16, y from -2 to 15
\put{1.988)} [l]  at 2 12
\put {$ \scriptstyle \bullet$} [c] at 10 0
\put {$ \scriptstyle \bullet$} [c] at 11 12
\put {$ \scriptstyle \bullet$} [c] at 12 0
\put {$ \scriptstyle \bullet$} [c] at 14 0
\put {$ \scriptstyle \bullet$} [c] at 14 12
\put {$ \scriptstyle \bullet$} [c] at 16 0
\put {$ \scriptstyle \bullet$} [c] at 16 12
\setlinear \plot  10 0 11 12 12 0      /
\setlinear \plot  14 0 14 12    /
\setlinear \plot 16 0 16  12   /
\put{$1{,}260$} [c] at 13 -2
\endpicture
\end{minipage}
\begin{minipage}{4cm}
\beginpicture
\setcoordinatesystem units <1.5mm,2mm>
\setplotarea  x from 0 to 16, y from -2 to 15
\put{1.989)} [l]  at 2 12
\put {$ \scriptstyle \bullet$} [c] at 10 12
\put {$ \scriptstyle \bullet$} [c] at 11 0
\put {$ \scriptstyle \bullet$} [c] at 11 12
\put {$ \scriptstyle \bullet$} [c] at 12 12
\put {$ \scriptstyle \bullet$} [c] at 14 12
\put {$ \scriptstyle \bullet$} [c] at 15 0
\put {$ \scriptstyle \bullet$} [c] at 16 12
\setlinear \plot  10 12 11 0 12 12      /
\setlinear \plot   11 0 11 12    /
\setlinear \plot 14 12 15 0 16 12   /
\put{$420$} [c] at 13 -2
\endpicture
\end{minipage}
\begin{minipage}{4cm}
\beginpicture
\setcoordinatesystem units <1.5mm,2mm>
\setplotarea  x from 0 to 16, y from -2 to 15
\put{1.990)} [l]  at 2 12
\put {$ \scriptstyle \bullet$} [c] at 10 12
\put {$ \scriptstyle \bullet$} [c] at 11 0
\put {$ \scriptstyle \bullet$} [c] at 11 12
\put {$ \scriptstyle \bullet$} [c] at 12 12
\put {$ \scriptstyle \bullet$} [c] at 14 0
\put {$ \scriptstyle \bullet$} [c] at 15 12
\put {$ \scriptstyle \bullet$} [c] at 16 0
\setlinear \plot  10 12 11 0 12 12      /
\setlinear \plot   11 0 11 12    /
\setlinear \plot 14 0 15 12 16 0   /
\put{$420$} [c] at 13 -2
\endpicture
\end{minipage}
\begin{minipage}{4cm}
\beginpicture
\setcoordinatesystem units <1.5mm,2mm>
\setplotarea  x from 0 to 16, y from -2 to 15
\put{1.991)} [l]  at 2 12
\put {$ \scriptstyle \bullet$} [c] at 10 0
\put {$ \scriptstyle \bullet$} [c] at 11 12
\put {$ \scriptstyle \bullet$} [c] at 11 0
\put {$ \scriptstyle \bullet$} [c] at 12 0
\put {$ \scriptstyle \bullet$} [c] at 14 12
\put {$ \scriptstyle \bullet$} [c] at 15 0
\put {$ \scriptstyle \bullet$} [c] at 16 12
\setlinear \plot  10 0 11 12 12 0      /
\setlinear \plot   11 0 11 12    /
\setlinear \plot 14 12 15 0 16 12   /
\put{$420$} [c] at 13 -2
\endpicture
\end{minipage}
$$
$$
\begin{minipage}{4cm}
\beginpicture
\setcoordinatesystem units <1.5mm,2mm>
\setplotarea  x from 0 to 16, y from -2 to 15
\put{1.992)} [l]  at 2 12
\put {$ \scriptstyle \bullet$} [c] at 10 0
\put {$ \scriptstyle \bullet$} [c] at 11 12
\put {$ \scriptstyle \bullet$} [c] at 11 0
\put {$ \scriptstyle \bullet$} [c] at 12 0
\put {$ \scriptstyle \bullet$} [c] at 14 0
\put {$ \scriptstyle \bullet$} [c] at 15 12
\put {$ \scriptstyle \bullet$} [c] at 16 0
\setlinear \plot  10 0 11 12 12 0      /
\setlinear \plot   11 0 11 12    /
\setlinear \plot 14 0 15 12 16 0   /
\put{$420$} [c] at 13 -2
\endpicture
\end{minipage}
\begin{minipage}{4cm}
\beginpicture
\setcoordinatesystem units <1.5mm,2mm>
\setplotarea x from  0 to 16, y from -2 to 15
\put{${\bf  46}$} [l]  at 2 15

\put{1.993)} [l]  at 2 12
\put {$ \scriptstyle \bullet$} [c] at 10 12
\put {$ \scriptstyle \bullet$} [c] at 11 12
\put {$ \scriptstyle \bullet$} [c] at 11 0
\put {$ \scriptstyle \bullet$} [c] at 12 12
\put {$ \scriptstyle \bullet$} [c] at 14 12
\put {$ \scriptstyle \bullet$} [c] at 13 0
\put{$\scriptstyle \bullet$} [c] at 16  0
\setlinear \plot 10 12 11 0  12 12 13 0 14 12 /
\setlinear \plot 11 0 11 12    /
\put{$2{,}520$} [c] at 13 -2
 \endpicture
\end{minipage}
\begin{minipage}{4cm}
\beginpicture
\setcoordinatesystem units <1.5mm,2mm>
\setplotarea x from  0 to 16, y from -2 to 15
\put{1.994)} [l]  at 2 12
\put {$ \scriptstyle \bullet$} [c] at 10 0
\put {$ \scriptstyle \bullet$} [c] at 11 12
\put {$ \scriptstyle \bullet$} [c] at 11 0
\put {$ \scriptstyle \bullet$} [c] at 12 0
\put {$ \scriptstyle \bullet$} [c] at 14 0
\put {$ \scriptstyle \bullet$} [c] at 13 12
\put{$\scriptstyle \bullet$} [c] at 16  0
\setlinear \plot 10 0 11 12  12 0 13 12 14 0 /
\setlinear \plot 11 0 11 12    /
\put{$2{,}520$} [c] at 13 -2
\endpicture
\end{minipage}
\begin{minipage}{4cm}
\beginpicture
\setcoordinatesystem units <1.5mm,2mm>
\setplotarea x from  0 to 16, y from -2 to 15
\put{${\bf  48}$} [l]  at 2 15

\put{1.995)} [l]  at 2 12
\put {$ \scriptstyle \bullet$} [c] at 10 0
\put {$ \scriptstyle \bullet$} [c] at 10 12
\put {$ \scriptstyle \bullet$} [c] at 12 0
\put {$ \scriptstyle \bullet$} [c] at 12 12
\put {$ \scriptstyle \bullet$} [c] at 14 12
\put {$ \scriptstyle \bullet$} [c] at 14 0
\setlinear \plot 10 12 10 0 12 12 12 0    /
\setlinear \plot 14  0 14 12      /
\put{$5{,}040$} [c] at 13 -2
\put{$\scriptstyle \bullet$} [c] at 16  0 \endpicture
\end{minipage}
\begin{minipage}{4cm}
\beginpicture
\setcoordinatesystem units <1.5mm,2mm>
\setplotarea x from  0 to 16, y from -2 to 15
\put{1.996)} [l]  at 2 12
\put {$ \scriptstyle \bullet$} [c] at 10 0
\put {$ \scriptstyle \bullet$} [c] at 10 12
\put {$ \scriptstyle \bullet$} [c] at 11 6
\put {$ \scriptstyle \bullet$} [c] at 12 0
\put {$ \scriptstyle \bullet$} [c] at 12 12
\put {$ \scriptstyle \bullet$} [c] at 14 0
\setlinear \plot 10 12 10 0 12 12 12 0    /
\put{$2{,}520$} [c] at 13 -2
\put{$\scriptstyle \bullet$} [c] at 16  0 \endpicture
\end{minipage}
\begin{minipage}{4cm}
\beginpicture
\setcoordinatesystem units <1.5mm,2mm>
\setplotarea  x from 0 to 16, y from -2 to 15
\put{1.997)} [l]  at 2 12
\put {$ \scriptstyle \bullet$} [c] at 10 0
\put {$ \scriptstyle \bullet$} [c] at 10 6
\put {$ \scriptstyle \bullet$} [c] at 10 12
\put {$ \scriptstyle \bullet$} [c] at 12 0
\put {$ \scriptstyle \bullet$} [c] at 12 12
\put {$ \scriptstyle \bullet$} [c] at 14 0
\put {$ \scriptstyle \bullet$} [c] at 16 0
\setlinear \plot 10 0 10 12   /
\setlinear \plot  12  0  12 12      /
\put{$2{,}520$} [c] at 13 -2
\endpicture
\end{minipage}
$$
$$
\begin{minipage}{4cm}
\beginpicture
\setcoordinatesystem units <1.5mm,2mm>
\setplotarea x from  0 to 16, y from -2 to 15
\put{1.998)} [l]  at 2 12
\put {$ \scriptstyle \bullet$} [c] at 10 0
\put {$ \scriptstyle \bullet$} [c] at 11 12
\put {$ \scriptstyle \bullet$} [c] at 11 0
\put {$ \scriptstyle \bullet$} [c] at 13 12
\put {$ \scriptstyle \bullet$} [c] at 13 0
\put {$ \scriptstyle \bullet$} [c] at 15 0
\setlinear \plot 10 0 11 12 11 0 13 12 13 0 11 12    /
\put{$1{,}260$} [c] at 13 -2
\put{$\scriptstyle \bullet$} [c] at 16  0 \endpicture
\end{minipage}
\begin{minipage}{4cm}
\beginpicture
\setcoordinatesystem units <1.5mm,2mm>
\setplotarea x from  0 to 16, y from -2 to 15
\put{1.999)} [l]  at 2 12
\put {$ \scriptstyle \bullet$} [c] at 10 12
\put {$ \scriptstyle \bullet$} [c] at 11 12
\put {$ \scriptstyle \bullet$} [c] at 11 0
\put {$ \scriptstyle \bullet$} [c] at 13 12
\put {$ \scriptstyle \bullet$} [c] at 13 0
\put {$ \scriptstyle \bullet$} [c] at 15 0
\setlinear \plot 10 12 11 0 11 12 13 0 13 12 11 0    /
\put{$1{,}260$} [c] at 13 -2
\put{$\scriptstyle \bullet$} [c] at 16  0 \endpicture
\end{minipage}
\begin{minipage}{4cm}
\beginpicture
\setcoordinatesystem units <1.5mm,2mm>
\setplotarea x from  0 to 16, y from -2 to 15
\put{2.000)} [l]  at 2 12
\put {$ \scriptstyle \bullet$} [c] at 10 12
\put {$ \scriptstyle \bullet$} [c] at 10 0
\put {$ \scriptstyle \bullet$} [c] at 12 12
\put {$ \scriptstyle \bullet$} [c] at 12 0
\put {$ \scriptstyle \bullet$} [c] at 14 12
\put {$ \scriptstyle \bullet$} [c] at 14 0
\put{$\scriptstyle \bullet$} [c] at 16  0
\setlinear \plot 10 12 10 0 14 12 14 0  /
\setlinear \plot 10 0 12 12    /
\setlinear \plot 12 0 14 12    /
\put{$1{,}260$} [c] at 13 -2
\endpicture
\end{minipage}
\begin{minipage}{4cm}
\beginpicture
\setcoordinatesystem units <1.5mm,2mm>
\setplotarea x from  0 to 16, y from -2 to 15
\put{2.001)} [l]  at 2 12
\put {$ \scriptstyle \bullet$} [c] at 10 12
\put {$ \scriptstyle \bullet$} [c] at 11 0
\put {$ \scriptstyle \bullet$} [c] at 11 6
\put {$ \scriptstyle \bullet$} [c] at 12 12
\put {$ \scriptstyle \bullet$} [c] at 14 0
\put {$ \scriptstyle \bullet$} [c] at 15 0
\put{$\scriptstyle \bullet$} [c] at 16  0
\setlinear \plot 10 12  11 6 11 0    /
\setlinear \plot 12 12 11 6     /
\put{$420$} [c] at 13 -2
 \endpicture
\end{minipage}
\begin{minipage}{4cm}
\beginpicture
\setcoordinatesystem units <1.5mm,2mm>
\setplotarea x from  0 to 16, y from -2 to 15
\put{2.002)} [l]  at 2 12
\put {$ \scriptstyle \bullet$} [c] at 10 0
\put {$ \scriptstyle \bullet$} [c] at 11 12
\put {$ \scriptstyle \bullet$} [c] at 11 6
\put {$ \scriptstyle \bullet$} [c] at 12 0
\put {$ \scriptstyle \bullet$} [c] at 14 0
\put {$ \scriptstyle \bullet$} [c] at 15 0
\put{$\scriptstyle \bullet$} [c] at 16  0
\setlinear \plot 10 0  11 6 11 12    /
\setlinear \plot 12 0 11 6     /
\put{$420$} [c] at 13 -2
\put{$\scriptstyle \bullet$} [c] at 16  0 \endpicture
\end{minipage}
\begin{minipage}{4cm}
\beginpicture
\setcoordinatesystem units <1.5mm,2mm>
\setplotarea x from  0 to 16, y from -2 to 15
\put{2.003)} [l]  at 2 12
\put {$ \scriptstyle \bullet$} [c] at 10 6
\put {$ \scriptstyle \bullet$} [c] at 11 0
\put {$ \scriptstyle \bullet$} [c] at 11 12
\put {$ \scriptstyle \bullet$} [c] at 12 6
\put {$ \scriptstyle \bullet$} [c] at 14 0
\put {$ \scriptstyle \bullet$} [c] at 15 0
\setlinear \plot 11 0 10 6 11 12 12  6 11 0   /
\put{$420$} [c] at 13 -2
\put{$\scriptstyle \bullet$} [c] at 16  0 \endpicture
\end{minipage}
$$
$$
\begin{minipage}{4cm}
\beginpicture
\setcoordinatesystem units <1.5mm,2mm>
\setplotarea x from  0 to 16, y from -2 to 15
\put{${\bf  49}$} [l]  at 2 15

\put{2.004)} [l]  at 2 12
\put {$ \scriptstyle \bullet$} [c] at 10 12
\put {$ \scriptstyle \bullet$} [c] at 12 12
\put {$ \scriptstyle \bullet$} [c] at 13 12
\put {$ \scriptstyle \bullet$} [c] at 14 12
\put {$ \scriptstyle \bullet$} [c] at 16 12
\put {$ \scriptstyle \bullet$} [c] at 13  0
\put {$ \scriptstyle \bullet$} [c] at 11.5  6
\setlinear \plot 10 12 13 0 16 12     /
\setlinear \plot 12 12 13 0 14 12     /
\setlinear \plot 13 0 13 12     /
\put{$210$} [c] at 13 -2
\endpicture
\end{minipage}
\begin{minipage}{4cm}
\beginpicture
\setcoordinatesystem units <1.5mm,2mm>
\setplotarea x from  0 to 16, y from -2 to 15
\put{2.005)} [l]  at 2 12
\put {$ \scriptstyle \bullet$} [c] at 10 0
\put {$ \scriptstyle \bullet$} [c] at 12 0
\put {$ \scriptstyle \bullet$} [c] at 13 0
\put {$ \scriptstyle \bullet$} [c] at 14 0
\put {$ \scriptstyle \bullet$} [c] at 16 0
\put {$ \scriptstyle \bullet$} [c] at 13  12
\put {$ \scriptstyle \bullet$} [c] at 11.5  6
\setlinear \plot 10 0 13 12 16 0     /
\setlinear \plot 12 0 13 12 14 0     /
\setlinear \plot 13 0 13 12     /
\put{$210$} [c] at 13 -2
\endpicture
\end{minipage}
\begin{minipage}{4cm}
\beginpicture
\setcoordinatesystem units <1.5mm,2mm>
\setplotarea x from  0 to 16, y from -2 to 15
\put {${\bf  50}$} [l] at 2 15

\put {2.006)} [l] at 2 12
\put {$ \scriptstyle \bullet$} [c] at  10 0
\put {$ \scriptstyle \bullet$} [c] at  12 0
\put {$ \scriptstyle \bullet$} [c] at  13 0
\put {$ \scriptstyle \bullet$} [c] at  14 0
\put {$ \scriptstyle \bullet$} [c] at  16  0
\put {$ \scriptstyle \bullet$} [c] at  13 12
\put {$ \scriptstyle \bullet$} [c] at  16  12
\setlinear \plot   10 0 13 12 16 0 16 12   /
\setlinear \plot  12 0  13 12  14 0 /
\setlinear \plot  13 0 13 12 /
\put{$210$} [c] at 13 -2
\endpicture
\end{minipage}
\begin{minipage}{4cm}
\beginpicture
\setcoordinatesystem units <1.5mm,2mm>
\setplotarea x from  0 to 16, y from -2 to 15
\put {2.007)} [l] at 2 12
\put {$ \scriptstyle \bullet$} [c] at  10 12
\put {$ \scriptstyle \bullet$} [c] at  12 12
\put {$ \scriptstyle \bullet$} [c] at  13 12
\put {$ \scriptstyle \bullet$} [c] at  14 12
\put {$ \scriptstyle \bullet$} [c] at  16  12
\put {$ \scriptstyle \bullet$} [c] at  13 0
\put {$ \scriptstyle \bullet$} [c] at  16  0
\setlinear \plot   10 12 13 0 16 12 16 0   /
\setlinear \plot  12 12  13 0  14 12 /
\setlinear \plot  13 0 13 12 /
\put{$210$} [c] at 13 -2
\endpicture
\end{minipage}
\begin{minipage}{4cm}
\beginpicture
\setcoordinatesystem units <1.5mm,2mm>
\setplotarea x from  0 to 16, y from -2 to 15
\put{2.008)} [l]  at 2 12
\put {$ \scriptstyle \bullet$} [c] at 10 0
\put {$ \scriptstyle \bullet$} [c] at 10.5 12
\put {$ \scriptstyle \bullet$} [c] at 11 0
\put {$ \scriptstyle \bullet$} [c] at 12 12
\put {$ \scriptstyle \bullet$} [c] at 12.5 0
\put {$ \scriptstyle \bullet$} [c] at 13 12
\setlinear \plot 10 0 10.5 12 11 0    /
\setlinear \plot 12 12 12.5 0 13 12     /
\put{$1{,}260$} [c] at 13 -2
\put{$\scriptstyle \bullet$} [c] at 16  0
\endpicture
\end{minipage}
\begin{minipage}{4cm}
\beginpicture
\setcoordinatesystem units <1.5mm,2mm>
\setplotarea x from  0 to 16, y from -2 to 15
\put{2.009)} [l]  at 2 12
\put {$ \scriptstyle \bullet$} [c] at 10 12
\put {$ \scriptstyle \bullet$} [c] at 12 12
\put {$ \scriptstyle \bullet$} [c] at 13 12
\put {$ \scriptstyle \bullet$} [c] at 15 12
\put {$ \scriptstyle \bullet$} [c] at 12.5 0
\put {$ \scriptstyle \bullet$} [c] at 13.8 6
\setlinear \plot 10 12 12.5 0 15  12   /
\setlinear \plot 12 12 12.5 0 13 12   /
\put{$840$} [c] at 13 -2
\put{$\scriptstyle \bullet$} [c] at 16  0
\endpicture
\end{minipage}
$$

$$
\begin{minipage}{4cm}
\beginpicture
\setcoordinatesystem units <1.5mm,2mm>
\setplotarea x from  0 to 16, y from -2 to 15
\put{2.010)} [l]  at 2 12
\put {$ \scriptstyle \bullet$} [c] at 10 0
\put {$ \scriptstyle \bullet$} [c] at 12 0
\put {$ \scriptstyle \bullet$} [c] at 13 0
\put {$ \scriptstyle \bullet$} [c] at 15 0
\put {$ \scriptstyle \bullet$} [c] at 12.5 12
\put {$ \scriptstyle \bullet$} [c] at 13.8 6
\setlinear \plot 10 0 12.5 12 15  0   /
\setlinear \plot 12 0 12.5 12 13 0   /
\put{$840$} [c] at 13 -2
\put{$\scriptstyle \bullet$} [c] at 16  0 \endpicture
\end{minipage}
\begin{minipage}{4cm}
\beginpicture
\setcoordinatesystem units <1.5mm,2mm>
\setplotarea x from  0 to 16, y from -2 to 15
\put{2.011)} [l]  at 2 12
\put {$ \scriptstyle \bullet$} [c] at 10 12
\put {$ \scriptstyle \bullet$} [c] at 11 0
\put {$ \scriptstyle \bullet$} [c] at 12 12
\put {$ \scriptstyle \bullet$} [c] at 13 12
\put {$ \scriptstyle \bullet$} [c] at 14 0
\put {$ \scriptstyle \bullet$} [c] at 15 12
\setlinear \plot 10 12 11 0 12 12    /
\setlinear \plot 13 12 14 0 15 12     /
\put{$630$} [c] at 13 -2
\put{$\scriptstyle \bullet$} [c] at 16  0 \endpicture
\end{minipage}
\begin{minipage}{4cm}
\beginpicture
\setcoordinatesystem units <1.5mm,2mm>
\setplotarea x from  0 to 16, y from -2 to 15
\put{2.012)} [l]  at 2 12
\put {$ \scriptstyle \bullet$} [c] at 10 0
\put {$ \scriptstyle \bullet$} [c] at 11 12
\put {$ \scriptstyle \bullet$} [c] at 12 0
\put {$ \scriptstyle \bullet$} [c] at 13 0
\put {$ \scriptstyle \bullet$} [c] at 14 12
\put {$ \scriptstyle \bullet$} [c] at 15 0
\setlinear \plot 10 0 11 12 12 0    /
\setlinear \plot 13 0 14 12 15 0     /
\put{$630$} [c] at 13 -2
\put{$\scriptstyle \bullet$} [c] at 16  0 \endpicture
\end{minipage}
\begin{minipage}{4cm}
\beginpicture
\setcoordinatesystem units <1.5mm,2mm>
\setplotarea  x from 0 to 16, y from -2 to 15
\put{${\bf  51}$} [l]  at 2 15

\put{2.013)} [l]  at 2 12
\put {$ \scriptstyle \bullet$} [c] at 10 12
\put {$ \scriptstyle \bullet$} [c] at 11 12
\put {$ \scriptstyle \bullet$} [c] at 14 12
\put {$ \scriptstyle \bullet$} [c] at 15 12
\put {$ \scriptstyle \bullet$} [c] at 12.5 0
\put {$ \scriptstyle \bullet$} [c] at 16 0
\put {$ \scriptstyle \bullet$} [c] at 16 12
\setlinear \plot  10 12   12.5  0 15 12     /
\setlinear \plot  11 12  12.5 0 14 12     /
\setlinear \plot 16 0 16 12   /
\put{$210$} [c] at 13 -2
\endpicture
\end{minipage}
\begin{minipage}{4cm}
\beginpicture
\setcoordinatesystem units <1.5mm,2mm>
\setplotarea  x from 0 to 16, y from -2 to 15
\put{2.014)} [l]  at 2 12
\put {$ \scriptstyle \bullet$} [c] at 10 0
\put {$ \scriptstyle \bullet$} [c] at 11 0
\put {$ \scriptstyle \bullet$} [c] at 14 0
\put {$ \scriptstyle \bullet$} [c] at 15 0
\put {$ \scriptstyle \bullet$} [c] at 12.5 12
\put {$ \scriptstyle \bullet$} [c] at 16 0
\put {$ \scriptstyle \bullet$} [c] at 16 12
\setlinear \plot  10 0   12.5  12 15 0     /
\setlinear \plot  11 0  12.5 12 14 0     /
\setlinear \plot 16 0 16 12   /
\put{$210$} [c] at 8 -2
\endpicture
\end{minipage}
\begin{minipage}{4cm}
\beginpicture
\setcoordinatesystem units <1.5mm,2mm>
\setplotarea x from  0 to 16, y from -2 to 15
\put{${\bf  52}$} [l]  at 2 15

\put{2.015)} [l]  at 2 12
\put {$ \scriptstyle \bullet$} [c] at 10 12
\put {$ \scriptstyle \bullet$} [c] at 12 12
\put {$ \scriptstyle \bullet$} [c] at 12 6
\put {$ \scriptstyle \bullet$} [c] at 12 0
\put {$ \scriptstyle \bullet$} [c] at 14 12
\put {$ \scriptstyle \bullet$} [c] at 15 0
\setlinear \plot 10 12 12 0 14 12     /
\setlinear \plot 12  0 12 12     /
\put{$1{,}260$} [c] at 13 -2
\put{$\scriptstyle \bullet$} [c] at 16  0 \endpicture
\end{minipage}
$$
$$
\begin{minipage}{4cm}
\beginpicture
\setcoordinatesystem units <1.5mm,2mm>
\setplotarea x from  0 to 16, y from -2 to 15
\put{2.016)} [l]  at 2 12
\put {$ \scriptstyle \bullet$} [c] at 10 0
\put {$ \scriptstyle \bullet$} [c] at 12 12
\put {$ \scriptstyle \bullet$} [c] at 12 6
\put {$ \scriptstyle \bullet$} [c] at 12 0
\put {$ \scriptstyle \bullet$} [c] at 14 0
\put {$ \scriptstyle \bullet$} [c] at 15 0
\setlinear \plot 10 0 12 12 14 0     /
\setlinear \plot 12  0 12 12     /
\put{$1{,}260$} [c] at 13 -2
\put{$\scriptstyle \bullet$} [c] at 16  0 \endpicture
\end{minipage}
\begin{minipage}{4cm}
\beginpicture
\setcoordinatesystem units <1.5mm,2mm>
\setplotarea x from  0 to 16, y from -2 to 15
\put{2.017)} [l]  at 2 12
\put {$ \scriptstyle \bullet$} [c] at 10 12
\put {$ \scriptstyle \bullet$} [c] at 10 0
\put {$ \scriptstyle \bullet$} [c] at 12 12
\put {$ \scriptstyle \bullet$} [c] at 14 0
\put {$ \scriptstyle \bullet$} [c] at 14 12
\put {$ \scriptstyle \bullet$} [c] at 15 0
\setlinear \plot 10 12 10 0 12 12 14 0 14 12    /
\put{$1{,}260$} [c] at 13 -2
\put{$\scriptstyle \bullet$} [c] at 16  0 \endpicture
\end{minipage}
\begin{minipage}{4cm}
\beginpicture
\setcoordinatesystem units <1.5mm,2mm>
\setplotarea x from  0 to 16, y from -2 to 15
\put{2.018)} [l]  at 2 12
\put {$ \scriptstyle \bullet$} [c] at 10 12
\put {$ \scriptstyle \bullet$} [c] at 10 0
\put {$ \scriptstyle \bullet$} [c] at 12 0
\put {$ \scriptstyle \bullet$} [c] at 14 0
\put {$ \scriptstyle \bullet$} [c] at 14 12
\put {$ \scriptstyle \bullet$} [c] at 15 0
\setlinear \plot 10 0 10 12 12 0 14 12 14 0    /
\put{$1{,}260$} [c] at 13 -2
\put{$\scriptstyle \bullet$} [c] at 16  0 \endpicture
\end{minipage}
\begin{minipage}{4cm}
\beginpicture
\setcoordinatesystem units <1.5mm,2mm>
\setplotarea x from  0 to 16, y from -2 to 15
\put{2.019)} [l]  at 2 12
\put {$ \scriptstyle \bullet$} [c] at 10 12
\put {$ \scriptstyle \bullet$} [c] at 11 12
\put {$ \scriptstyle \bullet$} [c] at 13 12
\put {$ \scriptstyle \bullet$} [c] at 14 12
\put {$ \scriptstyle \bullet$} [c] at 12 0
\put {$ \scriptstyle \bullet$} [c] at 14 0
\setlinear \plot 10 12 12 0 14  12 14 0  /
\setlinear \plot 11 12 12 0 13 12   /
\put{$840$} [c] at 13 -2
\put{$\scriptstyle \bullet$} [c] at 16  0 \endpicture
\end{minipage}
\begin{minipage}{4cm}
\beginpicture
\setcoordinatesystem units <1.5mm,2mm>
\setplotarea x from  0 to 16, y from -2 to 15
\put{2.020)} [l]  at 2 12
\put {$ \scriptstyle \bullet$} [c] at 10 0
\put {$ \scriptstyle \bullet$} [c] at 11 0
\put {$ \scriptstyle \bullet$} [c] at 13 0
\put {$ \scriptstyle \bullet$} [c] at 14 0
\put {$ \scriptstyle \bullet$} [c] at 12 12
\put {$ \scriptstyle \bullet$} [c] at 14 12
\setlinear \plot 10 0 12 12 14  0 14 12  /
\setlinear \plot 11 0 12 12 13 0   /
\put{$840$} [c] at 13 -2
\put{$\scriptstyle \bullet$} [c] at 16  0 \endpicture
\end{minipage}
\begin{minipage}{4cm}
\beginpicture
\setcoordinatesystem units <1.5mm,2mm>
\setplotarea x from  0 to 16, y from -2 to 15
\put{${\bf  54}$} [l]  at 2 15

\put{2.021)} [l]  at 2 12
\put {$ \scriptstyle \bullet$} [c] at 10 12
\put {$ \scriptstyle \bullet$} [c] at 12 12
\put {$ \scriptstyle \bullet$} [c] at 14 12
\put {$ \scriptstyle \bullet$} [c] at 12 0
\put {$ \scriptstyle \bullet$} [c] at 15 12
\put {$ \scriptstyle \bullet$} [c] at 15 0
\put{$\scriptstyle \bullet$} [c] at 16  0
\setlinear \plot 10 12 12 0 14 12    /
\setlinear \plot 12  0 12 12     /
\setlinear \plot 15  0 15 12     /
\put{$840$} [c] at 13 -2
\endpicture
\end{minipage}
$$
$$
\begin{minipage}{4cm}
\beginpicture
\setcoordinatesystem units <1.5mm,2mm>
\setplotarea x from  0 to 16, y from -2 to 15
\put{2.022)} [l]  at 2 12
\put {$ \scriptstyle \bullet$} [c] at 10 0
\put {$ \scriptstyle \bullet$} [c] at 12 0
\put {$ \scriptstyle \bullet$} [c] at 14 0
\put {$ \scriptstyle \bullet$} [c] at 12 12
\put {$ \scriptstyle \bullet$} [c] at 15 12
\put {$ \scriptstyle \bullet$} [c] at 15 0
\put{$\scriptstyle \bullet$} [c] at 16  0
\setlinear \plot 10 0 12 12 14 0    /
\setlinear \plot 12  0 12 12     /
\setlinear \plot 15  0 15 12     /
\put{$840$} [c] at 13 -2
\endpicture
\end{minipage}
\begin{minipage}{4cm}
\beginpicture
\setcoordinatesystem units <1.5mm,2mm>
\setplotarea x from  0 to 16, y from -2 to 15
\put{2.023)} [l]  at 2 12
\put {$ \scriptstyle \bullet$} [c] at 10 0
\put {$ \scriptstyle \bullet$} [c] at 10 12
\put {$ \scriptstyle \bullet$} [c] at 12 0
\put {$ \scriptstyle \bullet$} [c] at 12 12
\put {$ \scriptstyle \bullet$} [c] at 14 12
\put {$ \scriptstyle \bullet$} [c] at 14 0
\put {$ \scriptstyle \bullet$} [c] at 16 0
\setlinear \plot 10 0 10 12   /
\setlinear \plot 12 0 12 12   /
\setlinear \plot 14 0 14 12   /
\put{$840$} [c] at 13 -2
\endpicture
\end{minipage}
\begin{minipage}{4cm}
\beginpicture
\setcoordinatesystem units <1.5mm,2mm>
\setplotarea x from  0 to 16, y from -2 to 15
\put{${\bf  56}$} [l]  at 2 15

\put{2.024)} [l]  at 2 12
\put {$ \scriptstyle \bullet$} [c] at 10 12
\put {$ \scriptstyle \bullet$} [c] at 12 12
\put {$ \scriptstyle \bullet$} [c] at 14 12
\put {$ \scriptstyle \bullet$} [c] at 14 0
\put {$ \scriptstyle \bullet$} [c] at 12 0
\put {$ \scriptstyle \bullet$} [c] at 15 0
\setlinear \plot 10 12 12 0 14 12 14 0   /
\setlinear \plot 12  0 12 12     /
\put{$1{,}260$} [c] at 13 -2
\put{$\scriptstyle \bullet$} [c] at 16  0 \endpicture
\end{minipage}
\begin{minipage}{4cm}
\beginpicture
\setcoordinatesystem units <1.5mm,2mm>
\setplotarea x from  0 to 16, y from -2 to 15
\put{2.025)} [l]  at 2 12
\put {$ \scriptstyle \bullet$} [c] at 10 0
\put {$ \scriptstyle \bullet$} [c] at 12 0
\put {$ \scriptstyle \bullet$} [c] at 14 0
\put {$ \scriptstyle \bullet$} [c] at 14 12
\put {$ \scriptstyle \bullet$} [c] at 12 12
\put {$ \scriptstyle \bullet$} [c] at 15 0
\setlinear \plot 10 0 12 12 14 0 14 12   /
\setlinear \plot 12  0 12 12     /
\put{$1{,}260$} [c] at 13 -2
\put{$\scriptstyle \bullet$} [c] at 16  0 \endpicture
\end{minipage}
\begin{minipage}{4cm}
\beginpicture
\setcoordinatesystem units <1.5mm,2mm>
\setplotarea x from  0 to 16, y from -2 to 15
\put{2.026)} [l]  at 2 12
\put {$ \scriptstyle \bullet$} [c] at 10 12
\put {$ \scriptstyle \bullet$} [c] at 12 12
\put {$ \scriptstyle \bullet$} [c] at 10.5 6
\put {$ \scriptstyle \bullet$} [c] at 11 0
\put {$ \scriptstyle \bullet$} [c] at 14 0
\put {$ \scriptstyle \bullet$} [c] at 15 0
\put{$\scriptstyle \bullet$} [c] at 16  0
\setlinear \plot 10 12 11 0 12 12    /
\put{$840$} [c] at 13 -2
 \endpicture
\end{minipage}
\begin{minipage}{4cm}
\beginpicture
\setcoordinatesystem units <1.5mm,2mm>
\setplotarea x from  0 to 16, y from -2 to 15
\put{2.027)} [l]  at 2 12
\put {$ \scriptstyle \bullet$} [c] at 10 0
\put {$ \scriptstyle \bullet$} [c] at 12 0
\put {$ \scriptstyle \bullet$} [c] at 10.5 6
\put {$ \scriptstyle \bullet$} [c] at 11 12
\put {$ \scriptstyle \bullet$} [c] at 14 0
\put {$ \scriptstyle \bullet$} [c] at 15 0
\put{$\scriptstyle \bullet$} [c] at 16  0
\setlinear \plot 10 0 11 12 12 0    /
\put{$840$} [c] at 13 -2
\endpicture
\end{minipage}
$$

$$
\begin{minipage}{4cm}
\beginpicture
\setcoordinatesystem units <1.5mm,2mm>
\setplotarea x from  0 to 16, y from -2 to 15
\put{2.028)} [l]  at 2 12
\put {$ \scriptstyle \bullet$} [c] at 10 0
\put {$ \scriptstyle \bullet$} [c] at 10 12
\put {$ \scriptstyle \bullet$} [c] at 12 0
\put {$ \scriptstyle \bullet$} [c] at 12 12
\put {$ \scriptstyle \bullet$} [c] at 14 0
\put {$ \scriptstyle \bullet$} [c] at 15 0
\setlinear \plot 10 0 10 12 12 0 12 12 10 0   /
\put{$210$} [c] at 13 -2
\put{$\scriptstyle \bullet$} [c] at 16  0 \endpicture
\end{minipage}
\begin{minipage}{4cm}
\beginpicture
\setcoordinatesystem units <1.5mm,2mm>
\setplotarea x from 0 to 16, y from -2 to 15
\put{${\bf  60}$} [l] at 2 15

\put{2.029)} [l] at 2 12
\put {$ \scriptstyle \bullet$} [c] at 10 0
\put {$ \scriptstyle \bullet$} [c] at 11 12
\put {$ \scriptstyle \bullet$} [c] at 12 0
\put {$ \scriptstyle \bullet$} [c] at 14 12
\put {$ \scriptstyle \bullet$} [c] at 14 0
\put {$ \scriptstyle \bullet$} [c] at 15 0
\put{$\scriptstyle \bullet$} [c] at 16  0
\setlinear \plot 10 0 11 12 12 0    /
\setlinear \plot 14  0 14 12     /
\put{$1{,}260$} [c] at 13 -2
\endpicture
\end{minipage}
\begin{minipage}{4cm}
\beginpicture
\setcoordinatesystem units <1.5mm,2mm>
\setplotarea x from 0 to 16, y from -2 to 15
\put{2.030)} [l] at 2 12
\put {$ \scriptstyle \bullet$} [c] at 10 12
\put {$ \scriptstyle \bullet$} [c] at 11 0
\put {$ \scriptstyle \bullet$} [c] at 12 12
\put {$ \scriptstyle \bullet$} [c] at 14 12
\put {$ \scriptstyle \bullet$} [c] at 14 0
\put {$ \scriptstyle \bullet$} [c] at 15 0
\put{$\scriptstyle \bullet$} [c] at 16  0
\setlinear \plot 10 12 11 0 12 12    /
\setlinear \plot 14  0 14 12     /
\put{$1{,}260$} [c] at 13 -2
\endpicture
\end{minipage}
\begin{minipage}{4cm}
\beginpicture
\setcoordinatesystem units <1.5mm,2mm>
\setplotarea x from 0 to 16, y from -2 to 15
\put{${\bf  64}$} [l] at 2 15

\put{2.031)} [l] at 2 12
\put {$ \scriptstyle \bullet$} [c] at 10 0
\put {$ \scriptstyle \bullet$} [c] at 10 12
\put {$ \scriptstyle \bullet$} [c] at 12 0
\put {$ \scriptstyle \bullet$} [c] at 12 12
\put {$ \scriptstyle \bullet$} [c] at 14 0
\put {$ \scriptstyle \bullet$} [c] at 15 0
\put{$\scriptstyle \bullet$} [c] at 16  0
\setlinear \plot 10 12 10 0 12 12 12 0    /
\put{$840$} [c] at 13 -2
\endpicture
\end{minipage}
\begin{minipage}{4cm}
\beginpicture
\setcoordinatesystem units <1.5mm,2mm>
\setplotarea x from 0 to 16, y from -2 to 15
\put{2.032)} [l] at 2 12
\put {$ \scriptstyle \bullet$} [c] at 10 0
\put {$ \scriptstyle \bullet$} [c] at 10 6
\put {$ \scriptstyle \bullet$} [c] at 10 12
\put {$ \scriptstyle \bullet$} [c] at 13 0
\put {$ \scriptstyle \bullet$} [c] at 14 0
\put {$ \scriptstyle \bullet$} [c] at 15 0
\put{$\scriptstyle \bullet$} [c] at 16  0
\setlinear \plot 10 0 10 12      /
\put{$210$} [c] at 13 -2
\endpicture
\end{minipage}
\begin{minipage}{4cm}
\beginpicture
\setcoordinatesystem units <1.5mm,2mm>
\setplotarea x from 0 to 16, y from -2 to 15
\put{ ${\bf  65}$} [l] at 2 15

\put{2.033)} [l] at 2 12
\put {$ \scriptstyle \bullet$} [c] at 10 12
\put {$ \scriptstyle \bullet$} [c] at 11 12
\put {$ \scriptstyle \bullet$} [c] at 12.5 12
\put {$ \scriptstyle \bullet$} [c] at 13.5 12
\put {$ \scriptstyle \bullet$} [c] at 15 12
\put {$ \scriptstyle \bullet$} [c] at 16 12
\put{$\scriptstyle \bullet$} [c] at 13  0
\setlinear \plot 10 12 13 0 16 12     /
\setlinear \plot 11 12 13 0 15 12     /
\setlinear \plot 12.5 12  13 0  13.5 12   /
\put{$7$} [c] at 13 -2
 \endpicture
\end{minipage}
$$
$$
\begin{minipage}{4cm}
\beginpicture
\setcoordinatesystem units <1.5mm,2mm>
\setplotarea x from 0 to 16, y from -2 to 15
\put{2.034)} [l] at 2 12
\put {$ \scriptstyle \bullet$} [c] at 10 0
\put {$ \scriptstyle \bullet$} [c] at 11 0
\put {$ \scriptstyle \bullet$} [c] at 12.5 0
\put {$ \scriptstyle \bullet$} [c] at 13.5 0
\put {$ \scriptstyle \bullet$} [c] at 15 0
\put {$ \scriptstyle \bullet$} [c] at 16 0
\put{$\scriptstyle \bullet$} [c] at 13  12
\setlinear \plot 10 0 13 12 16 0     /
\setlinear \plot 11 0 13 12 15 0     /
\setlinear \plot 12.5 0  13 12  13.5 0   /
\put{$7$} [c] at 13 -2
 \endpicture
\end{minipage}
\begin{minipage}{4cm}
\beginpicture
\setcoordinatesystem units <1.5mm,2mm>
\setplotarea x from 0 to 16, y from -2 to 15
\put{${\bf  66}$} [l] at 2 15

\put{2.035)} [l] at 2 12
\put {$ \scriptstyle \bullet$} [c] at 10 12
\put {$ \scriptstyle \bullet$} [c] at 11 12
\put {$ \scriptstyle \bullet$} [c] at 12 12
\put {$ \scriptstyle \bullet$} [c] at 13 12
\put {$ \scriptstyle \bullet$} [c] at 14 12
\put {$ \scriptstyle \bullet$} [c] at 12 0
\put{$\scriptstyle \bullet$} [c] at 16  0
\setlinear \plot 10 12 12 0 14 12     /
\setlinear \plot 11 12 12 0 13 12     /
\setlinear \plot 12 0 12 12     /
\put{$42$} [c] at 13 -2
 \endpicture
\end{minipage}
\begin{minipage}{4cm}
\beginpicture
\setcoordinatesystem units <1.5mm,2mm>
\setplotarea x from 0 to 16, y from -2 to 15
\put{2.036)} [l] at 2 12
\put {$ \scriptstyle \bullet$} [c] at 10 0
\put {$ \scriptstyle \bullet$} [c] at 11 0
\put {$ \scriptstyle \bullet$} [c] at 12 0
\put {$ \scriptstyle \bullet$} [c] at 13 0
\put {$ \scriptstyle \bullet$} [c] at 14 0
\put {$ \scriptstyle \bullet$} [c] at 12 12
\put{$\scriptstyle \bullet$} [c] at 16  0
\setlinear \plot 10 0 12 12 14 0     /
\setlinear \plot 11 0 12 12 13 0     /
\setlinear \plot 12 0 12 12     /
\put{$42$} [c] at 13 -2
 \endpicture
\end{minipage}
\begin{minipage}{4cm}
\beginpicture
\setcoordinatesystem units <1.5mm,2mm>
\setplotarea x from 0 to 16, y from -2 to 15
\put{${\bf  68}$} [l] at 2 15

\put{2.037)} [l] at 2 12
\put {$ \scriptstyle \bullet$} [c] at 10 12
\put {$ \scriptstyle \bullet$} [c] at 11 12
\put {$ \scriptstyle \bullet$} [c] at 13 12
\put {$ \scriptstyle \bullet$} [c] at 14 12
\put {$ \scriptstyle \bullet$} [c] at 12 0
\put {$ \scriptstyle \bullet$} [c] at 15  0
\put{$\scriptstyle \bullet$} [c] at 16  0
\setlinear \plot 10 12 12 0 14 12     /
\setlinear \plot 11 12 12 0 13 12     /
\put{$105$} [c] at 13 -2
 \endpicture
\end{minipage}
\begin{minipage}{4cm}
\beginpicture
\setcoordinatesystem units <1.5mm,2mm>
\setplotarea x from 0 to 16, y from -2 to 15
\put{2.038)} [l] at 2 12
\put {$ \scriptstyle \bullet$} [c] at 10 0
\put {$ \scriptstyle \bullet$} [c] at 11 0
\put {$ \scriptstyle \bullet$} [c] at 13 0
\put {$ \scriptstyle \bullet$} [c] at 14 0
\put {$ \scriptstyle \bullet$} [c] at 12 12
\put {$ \scriptstyle \bullet$} [c] at 15  0
\put{$\scriptstyle \bullet$} [c] at 16  0
\setlinear \plot 10 0 12 12 14 0     /
\setlinear \plot 11 0 12 12 13 0     /
\put{$105$} [c] at 13 -2
 \endpicture
\end{minipage}
\begin{minipage}{4cm}
\beginpicture
\setcoordinatesystem units <1.5mm,2mm>
\setplotarea x from 0 to 16, y from -2 to 15
\put{${\bf  72}$} [l] at 2 15

\put{2.039)} [l] at 2 12
\put {$ \scriptstyle \bullet$} [c] at 10 12
\put {$ \scriptstyle \bullet$} [c] at 10 0
\put {$ \scriptstyle \bullet$} [c] at 12 0
\put {$ \scriptstyle \bullet$} [c] at 12 12
\put {$ \scriptstyle \bullet$} [c] at 14 0
\put {$ \scriptstyle \bullet$} [c] at 15  0
\put{$\scriptstyle \bullet$} [c] at 16  0
\setlinear \plot 12 0 12 12     /
\setlinear \plot 10 0 10 12     /
\put{$420$} [c] at 13 -2
 \endpicture
\end{minipage}
$$
$$
\begin{minipage}{4cm}
\beginpicture
\setcoordinatesystem units <1.5mm,2mm>
\setplotarea x from 0 to 16, y from -2 to 15
\put{2.040)} [l] at 2 12
\put {$ \scriptstyle \bullet$} [c] at 10 12
\put {$ \scriptstyle \bullet$} [c] at 11 12
\put {$ \scriptstyle \bullet$} [c] at 12 12
\put {$ \scriptstyle \bullet$} [c] at 11 0
\put {$ \scriptstyle \bullet$} [c] at 14 0
\put {$ \scriptstyle \bullet$} [c] at 15  0
\put{$\scriptstyle \bullet$} [c] at 16  0
\setlinear \plot 10 12 11 0 12  12     /
\setlinear \plot  11 0 11  12     /
\put{$140$} [c] at 13 -2
 \endpicture
\end{minipage}
\begin{minipage}{4cm}
\beginpicture
\setcoordinatesystem units <1.5mm,2mm>
\setplotarea x from 0 to 16, y from -2 to 15
\put{2.041)} [l] at 2 12
\put {$ \scriptstyle \bullet$} [c] at 10 0
\put {$ \scriptstyle \bullet$} [c] at 11 12
\put {$ \scriptstyle \bullet$} [c] at 12 0
\put {$ \scriptstyle \bullet$} [c] at 11 0
\put {$ \scriptstyle \bullet$} [c] at 14 0
\put {$ \scriptstyle \bullet$} [c] at 15  0
\put{$\scriptstyle \bullet$} [c] at 16  0
\setlinear \plot 10 0 11 12 12  0     /
\setlinear \plot  11 0 11  12     /
\put{$140$} [c] at 13 -2
\endpicture
\end{minipage}
\begin{minipage}{4cm}
\beginpicture
\setcoordinatesystem units <1.5mm,2mm>
\setplotarea x from 0 to 16, y from -2 to 15
\put{${\bf  80}$} [l] at 2 15
\put{2.042)} [l] at 2 12
\put {$ \scriptstyle\bullet$} [c] at 10 12
\put {$\scriptstyle \bullet$} [c] at 11 0
\put {$\scriptstyle \bullet$} [c] at 12 12
\put {$\scriptstyle \bullet$} [c] at 13 0
\put {$\scriptstyle \bullet$} [c] at 14 0
\put {$\scriptstyle \bullet$} [c] at 15  0
\put{$\scriptstyle\bullet$} [c] at 16  0
\setlinear \plot 10 12 11 0 12 12     /
\put{$105$} [c] at 13 -2
 \endpicture
\end{minipage}
\begin{minipage}{4cm}
\beginpicture
\setcoordinatesystem units <1.5mm,2mm>
\setplotarea x from 0 to 16, y from -2 to 15
\put{2.043)} [l] at 2 12
\put {$\scriptstyle \bullet$} [c] at 10 0
\put {$\scriptstyle \bullet$} [c] at 11 12
\put {$\scriptstyle \bullet$} [c] at 12 0
\put {$\scriptstyle \bullet$} [c] at 13 0
\put {$\scriptstyle \bullet$} [c] at 14 0
\put {$\scriptstyle \bullet$} [c] at 15  0
\put{$\scriptstyle\bullet$} [c] at 16  0
\setlinear \plot 10 0 11 12 12 0     /
\put{$105$} [c] at 13 -2
 \endpicture
\end{minipage}
\begin{minipage}{4cm}
\beginpicture
\setcoordinatesystem units <1.5mm,2mm>
\setplotarea x from 0 to 16, y from -2 to 15
\put{${\bf  96}$} [l] at 2 15

\put{2.044)} [l] at 2 12
\put {$ \scriptstyle\bullet$} [c] at 10 0
\put {$ \scriptstyle\bullet$} [c] at 10 12
\put {$ \scriptstyle\bullet$} [c] at 12 0
\put {$ \scriptstyle\bullet$} [c] at 13 0
\put {$ \scriptstyle\bullet$} [c] at 14 0
\put {$ \scriptstyle\bullet$} [c] at 15  0
\put  {$\scriptstyle\bullet$} [c] at 16  0
\setlinear \plot 10 0 10 12     /
\put{$42$} [c] at 13 -2
\endpicture
\end{minipage}
\begin{minipage}{4cm}
\beginpicture
\setcoordinatesystem units <1.5mm,2mm>
\setplotarea x  from 0 to 16, y from 10 to 15
\put{${\bf 128}$} [l] at 2 15

\put{2.045)} [l] at 2 12
\put {$ \scriptstyle\bullet$} [c] at 10 0
\put {$ \scriptstyle\bullet$} [c] at 11 0
\put {$ \scriptstyle\bullet$} [c] at 12 0
\put {$ \scriptstyle\bullet$} [c] at 13 0
\put {$ \scriptstyle\bullet$} [c] at 14 0
\put{$\scriptstyle\bullet$} [c] at 15  0
\put {$ \scriptstyle\bullet$} [c] at 16 0
\put{$1$} [c] at 13 -2
\endpicture
\end{minipage}
$$

\section{Determination of $|P(n)|$ for $1\le n\le 10$}

Summarizing the information accompanying the diagrams in the previous section, we can build  the following table:

\begin{center}
\begin{tabular}{r||r|r|r|r|r}
	                                  {$j$} &                    {$P^{(j)}(1)$} &                    {$P^{(j)}(2)$} &                    {$P^{(j)}(3)$} &                    {$P^{(j)}(4)$} &                     {$P^{(j)}(5)$} \\ \hline\hline
	                                      1 &                               $0$ &                               $0$ &                               $0$ &                               $0$ &                                $0$ \\ \hline
	                                      2 &                               $1$ &                               $0$ &                               $0$ &                               $0$ &                                $0$ \\ \hline
	                                      3 &                                   &                               $2$ &                               $0$ &                               $0$ &                                $0$ \\ \hline
	                                      4 &                                   &                               $1$ &                               $6$ &                               $0$ &                                $0$ \\ \hline
	                                      5 &                                   &                                   &                               $6$ &                              $24$ &                                $0$ \\ \hline
	                                      6 &                                   &                                   &                               $6$ &                              $36$ &                              $120$ \\ \hline
	                                      7 &                                   &                                   &                               $0$ &                              $54$ &                              $240$ \\ \hline
	                                      8 &                                   &                                   &                               $1$ &                              $48$ &                              $450$ \\ \hline
	                                      9 &                                   &                                   &                                   &                              $20$ &                              $600$ \\ \hline
	                                     10 &                                   &                                   &                                   &                              $24$ &                              $660$ \\ \hline
	                                     11 &                                   &                                   &                                   &                               $0$ &                              $500$ \\ \hline
	                                     12 &                                   &                                   &                                   &                              $12$ &                              $540$ \\ \hline
	                                     13 &                                   &                                   &                                   &                               $0$ &                              $240$ \\ \hline
	                                     14 &                                   &                                   &                                   &                               $0$ &                              $390$ \\ \hline
	                                     15 &                                   &                                   &                                   &                               $0$ &                              $120$ \\ \hline
	                                     16 &                                   &                                   &                                   &                               $1$ &                              $180$ \\ \hline
	                                     17 &                                   &                                   &                                   &                                   &                               $10$ \\ \hline
	                                     18 &                                   &                                   &                                   &                                   &                              $100$ \\ \hline
	                                     19 &                                   &                                   &                                   &                                   &                                $0$ \\ \hline
	                                     20 &                                   &                                   &                                   &                                   &                               $60$ \\ \hline
	                               21 to 23 &                                   &                                   &                                   &                                   &                                $0$ \\ \hline
	                                     24 &                                   &                                   &                                   &                                   &                               $20$ \\ \hline
	                               25 to 31 &                                   &                                   &                                   &                                   &                                $0$ \\ \hline
	                                     32 &                                   &                                   &                                   &                                   &                                  1 \\ 
\end{tabular}
\end{center}

 Then, using equation (\ref{PnSum}), we get $|P(2)|= 3$, $|P(3)|= 19$, $|P(4)|= 219$, and $|P(5)|=4{,}231$. We write the values for $n=6$ and $7$ in a separate table:

\begin{center}
\begin{tabular}{r|r||r|r||r|r||r|r}
{$j$} & {$P^{(j)}(6)$} &{$j$} & {$P^{(j)}(6)$} & {$j$} & {$P^{(j)}(7)$} & {$j$} & {$P^{(j)}(7)$} \\ \hline\hline

	                        $1$ &                               $0$ &                      $22$ &                         $4{,}800$ &                         1 &                                 0 &                        37 &                         $10{,}080$ \\ \hline
	                        $2$ &                               $0$ &                      $23$ &                             $720$ &                         2 &                                 0 &                        38 &                         $42{,}120$ \\ \hline
	                        $3$ &                               $0$ &                      $24$ &                         $3{,}600$ &                         3 &                                 0 &                        39 &                         $12{,}600$ \\ \hline
	                        $4$ &                               $0$ &                      $25$ &                             $600$ &                         4 &                                 0 &                        40 &                         $45{,}360$ \\ \hline
	                        $5$ &                               $0$ &                      $26$ &                         $1{,}680$ &                         5 &                                 0 &                        41 &                          $2{,}940$ \\ \hline
	                        $6$ &                               $0$ &                      $27$ &                             $360$ &                         6 &                                 0 &                        42 &                         $31{,}500$ \\ \hline
	                        $7$ &                             $720$ &                      $28$ &                         $1{,}530$ &                         7 &                                 0 &                        43 &                          $1{,}680$ \\ \hline
	                        $8$ &                         $1{,}800$ &                      $29$ &                               $0$ &                         8 &                         $5{,}040$ &                        44 &                         $23{,}940$ \\ \hline
	                        $9$ &                         $3{,}960$ &                      $30$ &                             $720$ &                         9 &                       $ 15{,}120$ &                        45 &                          $4{,}200$ \\ \hline
	                       $10$ &                         $6{,}570$ &                      $31$ &                               $0$ &                        10 &                       $ 37{,}800$ &                        46 &                          $5{,}040$ \\ \hline
	                       $11$ &                         $9{,}480$ &                      $32$ &                             $480$ &                        11 &                        $73{,}080$ &                        47 &                                $0$ \\ \hline
	                       $12$ &                        $11{,}520$ &                      $33$ &                              $12$ &                        12 &                       $123{,}900$ &                        48 &                         $15{,}120$ \\ \hline
	                       $13$ &                        $11{,}760$ &                      $34$ &                              $60$ &                        13 &                       $183{,}960$ &                        49 &                              $420$ \\ \hline
	                       $14$ &                        $12{,}960$ &                      $35$ &                               $0$ &                        14 &                       $245{,}070$ &                        50 &                          $4{,}620$ \\ \hline
	                       $15$ &                        $10{,}820$ &                      $36$ &                             $300$ &                        15 &                       $288{,}120$ &                        51 &                              $420$ \\ \hline
	                       $16$ &                        $12{,}240$ &              $37$ to $39$ &                               $0$ &                        16 &                       $350{,}280$ &                        52 &                          $6{,}720$ \\ \hline
	                       $17$ &                         $8{,}280$ &                      $40$ &                             $120$ &                        17 &                       $355{,}740$ &                        53 &                                $0$ \\ \hline
	                       $18$ &                        $10{,}290$ &              $41$ to $47$ &                               $0$ &                        18 &                       $414{,}540$ &                        54 &                          $2{,}520$ \\ \hline
	                       $19$ &                         $4{,}110$ &                      $48$ &                              $30$ &                        19 &                       $381{,}360$ &                        55 &                                $0$ \\ \hline
	                       $20$ &                         $7{,}080$ &              $49$ to $63$ &                               $0$ &                        20 &                       $430{,}920$ &                        56 &                          $4{,}410$ \\ \hline
	                       $21$ &                         $3{,}420$ &                      $64$ &                               $1$ &                        21 &                       $328{,}020$ &                  57 to 59 &                                $0$ \\ \hline
	                            &                                   &                           &                                   &                        22 &                       $401{,}940$ &                        60 &                          $2{,}520$ \\ \hline
	                            &                                   &                           &                                   &                        23 &                       $290{,}710$ &                  61 to 63 &                                $0$ \\ \hline
	                            &                                   &                           &                                   &                        24 &                       $346{,}500$ &                        64 &                          $1{,}050$ \\ \hline
	                            &                                   &                           &                                   &                        25 &                       $209{,}160$ &                        65 &                               $14$ \\ \hline
	                            &                                   &                           &                                   &                        26 &                       $292{,}740$ &                        66 &                               $84$ \\ \hline
	                            &                                   &                           &                                   &                        27 &                       $150{,}360$ &                        67 &                                $0$ \\ \hline
	                            &                                   &                           &                                   &                        28 &                       $232{,}260$ &                        68 &                              $210$ \\ \hline
	                            &                                   &                           &                                   &                        29 &                        $95{,}760$ &                  69 to 71 &                                $0$ \\ \hline
	                            &                                   &                           &                                   &                        30 &                       $192{,}500$ &                        72 &                              $700$ \\ \hline
	                            &                                   &                           &                                   &                        31 &                        $67{,}200$ &                  73 to 79 &                                $0$ \\ \hline
	                            &                                   &                           &                                   &                        32 &                       $137{,}970$ &                        80 &                              $210$ \\ \hline
	                            &                                   &                           &                                   &                        33 &                        $50{,}400$ &                  81 to 95 &                                $0$ \\ \hline
	                            &                                   &                           &                                   &                        34 &                       $100{,}086$ &                        96 &                               $42$ \\ \hline
	                            &                                   &                           &                                   &                        35 &                        $20{,}622$ &                 97 to 157 &                                $0$ \\ \hline
	                            &                                   &                           &                                   &                        36 &                        $90{,}090$ &                       128 &                                $1$ \\ 
\end{tabular}
\end{center}
Thus  $|P(6)|= 130{,}023$ and $|P(7)|= 6{,}129{,}859.$

 Ern\'e \cite{EM4} determines the numbers $P^{(j)}(n)$, $ 1\leq n\leq 6$, $1\leq j\leq 2^n$, and his tables coincide with the ones above.
To determine $|P(n)|$, $ 2\leq n\leq 9$, Ern\'e \cite{EM4} defines:
$$\tilde{P}^{(j)}(n)=\sum\limits_{r=0}^n {n \choose r}P^{(j)}(n-r)(-j)^r,\;\;1\leq j\leq 2^n,$$
and proves

\begin{equation}
\label{E1}
|P(n)|= \sum\limits_{m=1}^{n} {n \choose m} E(n,m), \;\;{\rm where} \;\; E(n,m)= \sum\limits_{j=1}^{2^{n-m}} \tilde{P}^{(j)}(n-m)\, j^m .
\end{equation}

  Let $P(n,m)$,  $1\leq m\leq n$, be the number of posets with $n$ elements and $m$ minimal elements. Note that $P(n,m)={n \choose m} E(n,m),\; 1\leq m\leq n$. Ern\'e  proves that:

\begin{equation}
\label{E2}
E(n,1)=E(n,2) =|P(n-1)|, \; E(n,n)=1.
\end{equation}
From the values $P^{(j)}(7)$,  $1\leq j\leq 128,$ we obtain:

\begin{center}
\begin{tabular}{r|r||r|r||r|r||r|r}
	                       {$j$} &                    {$\tilde{P}^{(j)}(7)$} &                     {$j$} &                    {$\tilde{P}^{(j)}(7)$} &                     {$j$} &                    {$\tilde{P}^{(j)}(7)$} &                     {$j$} &                     {$\tilde{P}^{(j)}(7)$} \\ \hline\hline
	                           1 &                                        -1 &                        19 &                                  -165.270 &                        37 &                                    10.080 &                        55 &                                          0 \\ \hline
	                           2 &                                       448 &                        20 &                                   -56.280 &                        38 &                                    42.210 &                        56 &                                      4.410 \\ \hline
	                           3 &                                   -10.206 &                        21 &                                  -174.720 &                        39 &                                    12.600 &                  57 to 59 &                                          0 \\ \hline
	                           4 &                                    32.256 &                        22 &                                  -337.260 &                        40 &                                    11.760 &                        60 &                                      2.520 \\ \hline
	                           5 &                                    26.250 &                        23 &                                   174.790 &                        41 &                                     2.940 &                  61 to 63 &                                          0 \\ \hline
	                           6 &                                    90.720 &                        24 &                                   -16.380 &                        42 &                                    31.500 &                        64 &                                        602 \\ \hline
	                           7 &                                  -436.590 &                        25 &                                   104.160 &                        43 &                                     1.680 &                        65 &                                         14 \\ \hline
	                           8 &                                  -207.760 &                        26 &                                   -13.020 &                        44 &                                    23.940 &                        66 &                                         84 \\ \hline
	                           9 &                                   275.940 &                        27 &                                    82.320 &                        45 &                                     4.200 &                        67 &                                          0 \\ \hline
	                          10 &                                   123.900 &                        28 &                                   -67.620 &                        46 &                                     5.040 &                        68 &                                        210 \\ \hline
	                          11 &                                   613.620 &                        29 &                                    95.760 &                        47 &                                         0 &                  69 to 71 &                                          0 \\ \hline
	                          12 &                                    63.420 &                        30 &                                    41.300 &                        48 &                                     5.040 &                        72 &                                        700 \\ \hline
	                          13 &                                   -34.440 &                        31 &                                    67.200 &                        49 &                                       420 &                  73 to 79 &                                          0 \\ \hline
	                          14 &                                   580.230 &                        32 &                                    51.954 &                        50 &                                     4.620 &                        80 &                                        210 \\ \hline
	                          15 &                                  -280.980 &                        33 &                                    47.628 &                        51 &                                       420 &                  81 to 95 &                                          0 \\ \hline
	                          16 &                                  -196.280 &                        34 &                                    85.806 &                        52 &                                     6.720 &                        96 &                                         42 \\ \hline
	                          17 &                                  -568.890 &                        35 &                                    20.622 &                        53 &                                         0 &                 97 to 157 &                                          0 \\ \hline
	                          18 &                                  -201.600 &                        36 &                                    14.490 &                        54 &                                     2.520 &                       128 &                                          1 \\ 
\end{tabular}
\end{center}

  From equation (\ref{E2}) and Ern\'e's results, $E(10,1)=E(10,2)=|P(9)|=44{,}511{,}042{,}511$, $E(10,10)=1$. From the preceding  tables and  Ern\'e's formula (\ref{E1}),  we obtain $E(10,m)$, $3\leq m\leq 9$.

\begin{center}
\begin{tabular}{lr|lr}
$E(10,1)=$&$44{,}511{,}042{,}511$&$E(10,2)=$&$44{,}511{,}042{,}511$\\
$E(10,3)=$&$21{,}724{,}257{,}583$&$E(10,4)=$&$6{,}110{,}970{,}115$\\
$E(10,5)=$&$1{,}004{,}241{,}331$&$E(10,6)=$&$91{,}323{,}031$\\
$E(10,7)=$&$4{,}109{,}383$&$E(10,8)=$&$77{,}635$\\
$E(10,9)=$&$511$&$E(10,10)=$&$1$\\
\end{tabular}
\end{center}

  Then:

$$|P(10)|=\sum\limits_{j=1}^{10} {10 \choose j} E(10,j)=6{,}611{,}065{,}248{,}783$$

\section{Connected and non-connected posets}

Another strategy to find the total number of posets is to count the connected and non-connected posets. If we let $CP(n)$ and $NCP(n)$ be the number of connected  and non-connected posets, respectively, with $n$ elements. It is clear that $|P(n)|=CP(n)+NCP(n)$. 

We further refine the partition considering $CP(n,m)$  and $NCP(n,m)$, $1\leq m\leq n$, the number of connected and non-connected posets with $n$ elements and $m$ minimal elements. Thus $CP(n)=\sum_{m=1}^{n}CP(n,m)$ and $NCP(n)=\sum_{m=1}^{n}NCP(n,m)$.

To determine the numbers $CP(n,n-k)$ for $k=2,3,4,5$, we have made explicit use of the Hasse diagrams of all the posets with $k$ elements. As an example, we sketch how we obtained the formula for $CP(n,n-3)$.

For each of the five posets with three elements, we consider all the posets $X$ with $n-3$ minimal elements that have this poset as the set of non-minimal elements. We consider in each case the restrictions that will ensure that the poset obtained will be connected. 

More specifically, for the antichain with three elements, which we call 1, 2, and 3, we let $P_i$ be the set of minimal elements that are below the element $i$. Next we partition the set of minimal elements in the following seven parts:
$$A_1= P_1\cap \overline{P_2}\cap \overline{P_3},\;\;
A_2= P_1\cap P_2\cap \overline{P_3},\;\;
A_3= P_1\cap P_2\cap P_3,$$
$$A_4= P_1\cap \overline{P_2}\cap P_3,\;\;
A_5= \overline{P_1}\cap P_2\cap \overline{P_3},\;\;
A_6= \overline{P_1}\cap P_2\cap P_3,$$
$$A_7= \overline{P_1}\cap \overline{P_2}\cap P_3.$$
Let  $A=\{A_1,A_2,A_3,A_4,A_5,A_6,A_7\}$. Thus $\bigcup\limits_{i=1}^7A_i=\bigcup\limits_{i=1}^3 P_i$. Each of the minimal elements of the poset must be in one of the elements of $A$, so we count the functions from the set $m(X)$ of the $n-3$ minimal elements of $X$ to $A$. Let  $G=\{A_1,\ldots, A_7\}^{m(X)}$, so $|G|=7^{n-3}$.

Since we are looking for connected posets, we are looking for the set $F$ of functions satisfying:
$$B_1=A_2\cup A_3\cup A_4\neq \emptyset,\;\;B_2= A_2\cup A_3\cup A_6\neq \emptyset,\;\;
B_3= A_3\cup A_4\cup A_6\neq \emptyset.$$
Here with $A_i\neq \emptyset$ we mean the set of functions $g\in G$ such that for some $x \in m(X)$, $g(x) \in A_i$.

The above conditions also imply that 
$$P_1=A_1\cup A_2\cup A_3 \cup A_4 \neq \emptyset,\;\; P_2=A_2\cup A_3\cup A_5 \cup A_6 \neq \emptyset,$$
$$P_3= A_3\cup A_4\cup A_6 \cup A_7 \neq \emptyset.$$

Thus we have 
$$F=\{g\in G:\; {\rm there\ exist}\; x,y,z \in m(X)\;{\rm such\; that}\; g(x)\in B_1,\; g(y)\in B_2,\; g(z)\in B_3\}$$
In order to determine $|F|$, notice that
$$|F|=|G|-|\overline{F}|=7^{n-3}-|\overline{F}|$$
and also $\overline{F}=\bigcup\limits_{i=1}^3 F_i$ where
$$F_1=\{g\in G:\forall x\in m(X), g(x)\in A_1\cup A_5\cup A_6\cup A_7 \}$$
$$F_2=\{g\in G:\forall x\in m(X), g(x)\in A_1\cup A_4\cup A_5\cup A_7\}$$
$$F_3=\{g\in G:\forall x\in m(X), g(x)\in A_1\cup A_2\cup A_5\cup A_7\}.$$
We also observe that $F_1\cap F_2=\{g\in G:\forall x\in m(X), g(x)\in A_1\cup A_5\cup A_7\}$ 
$=F_1\cap F_3$
$=F_2\cap F_3$ $=F_1\cap F_2\cap F_3$, so using the inclusion exclusion principle we get
$$|\overline{F}|=\sum\limits_{i=1}^3|F_i|-|F_1\cap F_2|-|F_1\cap F_3|-|F_2\cap F_3|+
|F_1\cap F_2\cap F_3|$$
and since $|F_1|=|F_2|=|F_3|=4^{n-3}$ and $|F_1\cap F_2|=|F_1\cap  F_3| =|F_2\cap F_3|=|F_1\cap F_2\cap F_3|= 3^{n-3}$
then $|\overline{F}|=3\; 4^{n-3}-2\;3^{n-3}$ so $|F|=7^{n-3}-(3\; 4^{n-3}-2\;3^{n-3})=7^{n-3}-3\; 4^{n-3}+2\;3^{n-3}$.

In a similar way, for the posets with three elements depicted in section \ref{diagrams} we get the formulas:

\[1. {n \choose n-3}\;6(3^{n-3}-2^{n-3})\]
\[2. {n \choose n-3}\;3(4^{n-3}-2\;2^{n-3}+1)\]
\[3.{n \choose n-3}\;3(4^{n-3}-3^{n-3})\]
\[4. {n \choose n-3}\;6(5^{n-3}-2\;3^{n-3}+2^{n-3})\]

Adding these up, we obtain the following result:
\begin{lemma} For $n \geq 4:$
\[CP(n,n-3)={n \choose n-3}\left( 7^{n-3}+6\; 5^{n-3}+ 3\; 4^{n-3}- 7\;3^{n-3}-6\;2^{n-3}+3\right).\]
\end{lemma}

Following this methodology we also obtained:
\begin{theorem}
It is easy to see that if $n\geq 2$ then $CP(n,n-1)=n$ and $NCP(n,n-1)=n(2^{n-1}-2)$.

For $n\geq 3:$
$CP(n,n-2)={n \choose n-2}\left(3^{n-2}+ 2^{n-2}-2\right)$ and $NCP(n,n-2)=$ ${n \choose n-2}(4^{n-2}+ 3^{n-2} -5 \;2^{n-2} +3)$.

For $n\ge 4$: $NCP(n,n-3)= {n\choose{n-3}}(8^{n-3}-7^{n-3}+6\;6^{n-3}-9\;4^{n-3}-11\;3^{n-3}+18\;2^{n-3}-4)$.

For $n\geq 5:$ $CP(n,n-4)={n \choose n-4}(15^{n-4}+12\; 11^{n-4}+24\;9^{n-4}+16\;8^{n-4}+20\;7^{n-4}+$
$27\;6^{n-4}-108\;5^{n-4}-42\;4^{n-4}+30\;3^{n-4}+24\;2^{n-4}-4)$ and 
$NCP(n,n-4)={n\choose{n-4}} (16^{n-4}-15^{n-4}+12\;12^{n-4}-12\;11^{n-4}+24\;10^{n-4}-$
$4\;9^{n-4}+34\;7^{n-4}-135\;6^{n-4}+12\;5^{n-4}+42\;4^{n-4}+78\;3^{n-4}-56\;2^{n-4}+5)$.

For $n\geq 6$:
$CP(n,n-5)={n \choose n-5}(31^{n-5}+20\;23^{n-5}+60\;19^{n-5}+100\;17^{n-5}+$
$5\;16^{n-5}+105\; 15^{n-5}+120\;14^{n-5}+390\;13^{n-5}+ 180\;12^{n-5}$
$-120\;11^{n-5}+370\;10^{n-5}-520\; 9^{n-5}-340\;8^{n-5}-1{,}040\;7^{n-5}-750\;6^{n-5}+$
$1{,}224\;5^{n-5}+360\;4^{n-5}-90\; 3^{n-5}-80\;2^{n-5}+5)$.

\end{theorem}

\section{Other results}

  From the formulas in the previous section and using that $P(n,1)= |P(n-1)|,\; P(n,2)= {n \choose 2}|P(n-1)|,\; P(n,n)=1$ we can also calculate $|P(n)$ for $1\le n\le 8$, without the aid of any computer program. 
  
  As an example of one of the intermediate steps, here we calculate the number of non connected posets with 8 elements, three of which are minimal. We consider how the three minimal elements can be part of different connected components:
\begin{center}
\begin{tabular}{lr}\\
${8 \choose 7}\;CP(7,2)\;CP(1,1)=$&$19{,}862{,}808$\\[2mm]
${8 \choose 6}\;CP(6,2)\;CP(2,1)=$&$3{,}013{,}080$\\[2mm]
${8 \choose 5}\;CP(5,2)\;CP(3,1)=$&$821{,}520$\\[2mm]
${8 \choose 4}\;CP(4,2)\;CP(4,1)=$&$351{,}120$\\[2mm]
${8 \choose 5}\;CP(5,1)\;CP(3,2)=$&$183{,}960$\\[2mm]
${8 \choose 6}\;CP(6,1)\;{2 \choose 1}\;PC(1,1)\;CP(1,1)/2!=$&$710{,}808$\\[2mm]
${8 \choose 5}\;CP(5,1)\;{3 \choose 2}\;PC(2,1)\;CP(1,1)=$&$367{,}920$\\[2mm]
${8 \choose 4}\;CP(4,1)\;{4 \choose 3}\;PC(3,1)\;CP(1,1)=$&$191{,}520$\\[2mm]
${8 \choose 4}\;CP(4,1)\;{4 \choose 2}\;PC(2,1)\;CP(2,1)/2!=$&$63{,}840$\\[2mm]
${8 \choose 3}\;CP(3,1)\;{5 \choose 3}\;CP(3,1)\;CP(2,1)/2!=$&$45{,}360$\\[3mm] \hline
$NCP(8,3)=$ & $25{,}611{,}936$
\end{tabular}
\end{center}

Some further results about $P(n)$ can be found in \cite{MSV}:
 
 \begin{itemize}
\item  We have determined the numbers  $CP^{(j)}(n)$ ($NCP^{(j)}(n)$)  of conected (non-conected) posets on $n$ elements and $j$ antichains for $1\leq n\leq 7$ and $1\leq j\leq 2^n$, which can be seen in  the Hasse diagrams in section \ref{diagrams}. 
    \item We proved that, if $n\geq 5$, $P^{(2^{n-1})}(n)=NCP^{(2^{n-1})}(n)$, and if $n \ge 6$:

\begin{itemize}

\item If   $j\ge 2^{n-1}+2$ then $CP^{(j)}(n)=0$,
\item  $P^{(7\, 2^{n-4})}(n)=NCP^{(7\, 2^{n-4})}(n)$,
\item   $P^{(15\, 2^{n-5})}(n)=NCP^{(15\, 2^{n-5})}(n)$,
\item $P^{(2^{n-1}+1)}(n)=CP^{(2^{n-1}+1)}(n)$,
\item  $NCP^{(2^{n-1}+2)}(n)= 2(n-1)n$ and $CP^{(2^{n-1}+2)}(n)=0,$
\item  $CP^{(7\;2^{n-4})}(n)=CP^{(15\;2^{n-5})}(n)=CP^{(2^{n-1})}(n)=NCP^{(2^{n-1}+1)}(n)=0.$
\end{itemize}

\item  We also obtained the following formula to determine the number of non-connected posets with $n$ elements and $2$ minimal elements, when $n\geq 3$:
%
%

\[
|NCP(n,2)|=\left\{\begin{array}{lr}
\displaystyle{\sum\limits_{k=1}^{\frac{n}{2}-1} {n \choose k} |CP(k,1)|.|CP(n-k,1)|}+\displaystyle{{n \choose \frac{n}{2}} \frac{|PC(\frac{n}{2},1)|^2}{2}}, & \text{if $n$ is even}\\
&\\
	\displaystyle{\sum\limits_{k=1}^{\frac{n-1}{2}} {n \choose k} |CP(k,1)|.|CP(n-k,1)|}, & \text{if $n$ is odd.}
\end{array}\right.
\]

\end{itemize}

%
%
%
%
%
%

\textbf{Authors' contributions} All authors have contributed equally to the production of this manuscript

\end{document}